\documentclass[10pt]{amsart}
\usepackage{amscd,amssymb,amsthm}
\usepackage{amsxtra}
\usepackage{amssymb, amsmath}
\usepackage{amsthm}
\usepackage{multicol}
\usepackage {hyperref}
\usepackage{latexsym}

\usepackage[condensed,math]{iwona} 

\begin{document}
\title[Crazy Sequential Representation - Inder J. Taneja]{}
\begin{center}
\textbf{\huge{Crazy Sequential Representation: Numbers from 0 to 11111 \\ \vspace{6pt} in terms of Increasing and Decreasing Orders of 1 to 9}}
\end{center}

\bigskip
\begin{center}
\textbf{\Large{Inder J. Taneja}}\footnote{\textit{Formerly, Professor of Mathematics, Universidade Federal de Santa Catarina, 88.040-900 Florian\'{o}polis, SC, Brazil. e-mail: ijtaneja@gmail.com.}}
\end{center}

\begin{abstract}
 Natural numbers from 0 to 11111 are written in terms of 1 to 9 in two different ways. The first one in increasing order of 1 to 9, and the second one in decreasing order. This is done by using the operations of \textit{addition, multiplication, subtraction, potentiation}, and \textit{division}. In both the situations there are no missing numbers, except one, i.e., 10958 in the increasing case.
\end{abstract}

\maketitle

\section{\textbf{Introduction}}

Author \cite{t1, t2, t3} wrote the numbers from  44 to 11111 in terms of 1 to 9 in two different ways, one is in increasing order and another in decreasing order. Some comments on this work can be seen in \cite{ma1, ma2, jn1, jn2}. The operations used are only \textit{addition, multiplication}, and \textit{potentiation}. The idea of brackets is also used, i.e., the following operations were used:

\[
[plus, \,product, \, potentiation, \,brackets].
\]

\bigskip
From the mathematical point of view, the brackets are understood as \textit{composite rule}. The operations such as \textit{subtraction} and  \textit{division} are also very important. In this work, the operations of \textit{subtraction} and \textit{division} are also included. This is done to find missing numbers not available in the previous versions. This work is done by using the following operations:

\[
[plus,\,minus, \,product, \,potentiation, \,division, \,brackets].
\]

\bigskip
 In the previous work \cite{t3}, there were approximately 1250 numbers were missing in both the cases. Here, we have found almost all the missing numbers from 0 to 11111, except one, i.e., 10958 in the increasing case. These missing numbers having either \textit{subtraction} and/or \textit{division}, and are written in \textit{italic} forms to identify. Still, there are more operations that can be applied, such as:
 \[
 [factorial, \,decimal, \,square \, root, \,etc.]
 \]

\bigskip
By applying these operations, may be one can find the number 10958. This shall be dealt elsewhere.

\bigskip
The mathematical idea behind this work is based on simple combinations. If we have two different positive natural numbers in a sequence, for example $a$ and $b$, then we can write,
\[
a+b, \, a\times b, \, a^{b} \, \mbox{and} \,\, ab
\]

\bigskip
We have only four ways of writing two numbers, for example if we have $a=2$ and $b=3$, then one can write  $2+3$, $2\times 3$, $2^{3}$ and $23$ in the increasing order, and $3+2$, $3\times 2$,  $3^{2}$ and $32$ in decreasing order.

\bigskip
Again, let us consider three positive natural numbers, $a$, $b$ and $c$ with either $a < b < c$ or $a > b > c$. Following the same procedure for two numbers, here below are 23 possibilities of writing these three numbers:
\[
a+b+c, \,ab+c, \,a+bc, (a+b)\times c, \, a\times(b+c), \, a\times b+c, \, a+b\times c, \,
ab\times c, \, a\times bc, \, a\times b\times c, \, abc,
\]
\[
a^{bc}, \, {a^b}^c, \, (a^b)^c, \, a^{b}\times c, \, a\times b^{c},  \, (a\times b)^{c},  \, (ab)^{c},  \, a^{b}+c, \,
a+b^{c},  \, a^{b+c}, \, a^{b\times c}, \,\, \mbox{and} \,\, (a+b)^{c}.
\]

\bigskip
The expressions $(a^b)^c$ and $a^{b\times c}$ are the same. The expressions $a^{bc}$ and ${a^b}^c$ give very big values except $a=1$.

\bigskip
Imagine if these letters increases from 3 to 4, 5, ... to 9, one may have millions of possibilities of writing these 9 letters either in increasing or in decreasing orders. The above explanation is only for addition, multiplication and brackets. If we allow more operations, such as subtraction, division, etc., these possibilities increases much more.

\bigskip
From first version to this, there is a gap of approximately one year. During this time, I came across, two historical books, \cite{dh, jm}, where these authors specified only the representation of number 100 in different ways including much more opeartions, such as, factorial, decimal, square root etc.

\section{\textbf{Crazy Sequential Representation}}

This is the fifth version of previous works. Readers can see previous versions at \cite{t1, t2, t3, t4}. Here below are \textit{crazy sequential representation} of natural numbers written in terms of 1 to 9 in increasing as well as decreasing orders. The first column represents the increasing order and the second represent the decreasing. Numbers with \textit{subtraction} and/or \textit{division} are written in \textit{Italic} form.

{\footnotesize
\begin{multicols}{2}
\begin{itemize}
\item[]$\mbox{Increasing order}$
\item [] $\mathit{0=12+34-56-7+8+9.}$
\item [] $1=1^{23456789}.$
\item [] $\mathit{2=123+4-56-78+9.}$
\item [] $\mathit{3=123-45-6-78+9.}$
\item [] $\mathit{4=12-34-56-7+89.}$
\item [] $\mathit{5=12-34+5-67+89.}$
\item [] $\mathit{6=12+34+56-7-89.}$
\item [] $\mathit{7=1+23-4+56-78+9.}$
\item [] $\mathit{8=1-23-45+6+78-9.}$
\item [] $9=1^{2345678}\times 9.$
\item [] $10=1^{2345678}+9.$
\item [] $\mathit{11=1+23+4+5+67-89.}$
\item [] $\mathit{12=123+45-67-89.}$
\item [] $\mathit{13=1-23+4-56+78+9.}$
\item [] $\mathit{14=12-3-45+67-8-9.}$
\item [] $\mathit{15=123-45+6-78+9.}$
\item [] $\mathit{16=1-2+34+5+67-89.}$
\item [] $17=1^{234567}\times 8+9.$
\item [] $18=1^{234567}+8+9.$
\item [] $\mathit{19=12+34-5+67-89.}$
\item [] $\mathit{20=12+3-45+67-8-9.}$
\item [] $\mathit{21=1-23-45+6-7+89.}$
\item [] $\mathit{22=1-23+4-56+7+89.}$
\item [] $\mathit{23=1+2-3+45+67-89.}$
\item [] $24=1^{23456}\times 7+8+9.$
\item [] $25=1^{23456} +7+8+9.$
\item [] $\mathit{26=12-3+4-56+78-9.}$
\item [] $\mathit{27=12-3-45-6+78-9.}$
\item [] $\mathit{28=12+3-4-5-67+89.}$
\item [] $\mathit{29=12+34+5+67-89.}$
\item [] $30=1^{2345}\times 6 +7+8+9.$
\item [] $31=1^{2345}+6+7+8+9.$
\item [] $\mathit{32=12-3+45+67-89.}$
\item [] $\mathit{33=12+34+56-78+9.}$
\item [] $\mathit{34=123+4-5-6+7-89.}$
\item [] $35=1^{234}\times 5+6+7+8+9.$
\item [] $36=1^{234}+5+6+7+8+9.$
\item [] $\mathit{37=1+23-4-5-67+89.}$
\item [] $\mathit{38=12+3+45+67-89.}$
\item [] $39=1^{23}\times 4+5+6+7+8+9.$
\item [] $40=1^{23}+4+5+6+7+8+9.$
\item [] $\mathit{41=12-34-5+67-8+9.}$
\item [] $42=1^{2}\times 3+4+5+6+7+8+9.$
\item [] $43=1^{2}+3+4+5+6+7+8+9.$
\item [] $44=1\times 2+3+4+5+6+7+8+9.$
\item [] $45=1+2+3+4+5+6+7+8+9.$
\item [] $46=1+2\times 3+4+5+6+7+8+9.$
\item [] $47=1\times 2^3+4+5+6+7+8+9.$
\item [] $48=1+2^3+4+5+6+7+8+9.$
\item [] $49=1\times 2+3\times 4+5+6+7+8+9.$
\item [] $50=1+2+3\times 4+5+6+7+8+9.$
\item [] $51=1^{23}+4\times 5+6+7+8+9.$
\item [] $\mathit{52=12-3-45+6-7+89.}$
\item [] $53=1^2\times 3+4\times 5+6+7+8+9.$
\item [] $54=12+3+4+5+6+7+8+9.$
\item [] $55=1\times 2+3+4\times 5+6+7+8+9.$
\item [] $56=1+2+3+4\times 5+6+7+8+9.$
\item [] $57=1+2\times 3+4\times 5+6+7+8+9.$
\item [] $58=1\times 2^3+4\times 5+6+7+8+9.$
\item [] $59=1\times 2\times 3\times 4+5+6+7+8+9.$
\item [] $60=1+2\times 3\times 4+5+6+7+8+9.$
\\
\item[]$\mbox{Decreasing order}$
\item [] $\mathit{0=98-7-6-54-32+1.}$
\item [] $\mathit{1=98-76-54+32+1.}$
\item [] $\mathit{2=9+87-65+4-32-1.}$
\item [] $\mathit{3=98-76-5+4+3-21.}$
\item [] $\mathit{4=98-7-65-43+21.}$
\item [] $\mathit{5=98-76+5-43+21.}$
\item [] $\mathit{6=98-7-65+4-3-21.}$
\item [] $\mathit{7=98-7-6-54-3-21.}$
\item [] $\mathit{8=9-8+76-5-43-21.}$
\item [] $\mathit{9=9+87-65-43+21.}$
\item [] $\mathit{10=98-7+6-54-32-1.}$
\item [] $\mathit{11=9+8-7+65-43-21.}$
\item [] $\mathit{12=987-654-321.}$
\item [] $\mathit{13=98-7-6-54+3-21.}$
\item [] $\mathit{14=98+7-6-54-32+1.}$
\item [] $\mathit{15=98-76-5-4+3-2+1.}$
\item [] $\mathit{16=98-7-6-5-43-21.}$
\item [] $\mathit{17=9+87-65+4+3-21.}$
\item [] $\mathit{18=98+7-65-43+21.}$
\item [] $\mathit{19=98-7+6-54-3-21.}$
\item [] $\mathit{20=98+7-65+4-3-21.}$
\item [] $\mathit{21=9+87-6-5-43-21.}$
\item [] $\mathit{22=9-87+65+4+32-1.}$
\item [] $\mathit{23=9+87-65-4-3-2+1.}$
\item [] $\mathit{24=98+7+6-54-32-1.}$
\item [] $\mathit{25=9+8+7+65-43-21.}$
\item [] $\mathit{26=98-7-6+5-43-21.}$
\item [] $\mathit{27=9-87+65+43-2-1.}$
\item [] $\mathit{28=98-7+6-5-43-21.}$
\item [] $\mathit{29=9-87+65+43-2+1.}$
\item [] $\mathit{30=98+7-6-5-43-21.}$
\item [] $\mathit{31=98-76-5-4-3+21.}$
\item [] $\mathit{32=98-7-65+4+3-2+1.}$
\item [] $\mathit{33=98+7+6-54-3-21.}$
\item [] $\mathit{34=9+8+76+5-43-21.}$
\item [] $\mathit{35=98-7-6-54+3+2-1.}$
\item [] $\mathit{36=98-7-6-5-43-2+1.}$
\item [] $\mathit{37=98-76-5-4+3+21.}$
\item [] $\mathit{38=98-7-6-5-43+2-1.}$
\item [] $\mathit{39=98-76-5+43-21.}$
\item [] $\mathit{40=98-7-65-4-3+21.}$
\item [] $\mathit{41=98-76+5-4-3+21.}$
\item [] $\mathit{42=98+7+6-5-43-21.}$
\item [] $\mathit{43=98-76+54-32-1.}$
\item [] $44=9+8+7+6+5+4+3+2\times 1.$
\item [] $45=9+8+7+6+5+4+3+2+1.$
\item [] $46=9+8+7+6+5+4+3\times 2+1.$
\item [] $\mathit{47=98-76+5-4+3+21.}$
\item [] $48=9+8+7+6+5+4+3^2\times 1.$
\item [] $49=9+8+7+6+5+4\times 3+2\times 1.$
\item [] $50=9+8+7+6+5+4\times 3+2+1.$
\item [] $\mathit{51=9+87-65-4+3+21.}$
\item [] $\mathit{52=98-76+54-3-21.}$
\item [] $\mathit{53=9+87-65+43-21.}$
\item [] $54=9+8+7+6+(5+4+3)\times 2\times 1.$
\item [] $55=9+8+7+6+5\times 4+3+2\times 1.$
\item [] $56=9+8+7+6+5\times 4+3+2+1.$
\item [] $57=9+8+7+6+5\times 4+3\times 2+1.$
\item [] $\mathit{58=98-7-6-5-43+21.}$
\item [] $59=9+8+7+6+5+4\times 3\times 2\times 1.$
\item [] $60=9+8+7+6+5\times 4+3^2+1.$
\newpage
\item[]$\mbox{Increasing order}$
\item [] $61=1^2\times 3+4+5\times 6+7+8+9.$
\item [] $62=1\times 23+4+5+6+7+8+9.$
\item [] $63=1+23+4+5+6+7+8+9.$
\item [] $64=1+2+3+4+5\times 6+7+8+9.$
\item [] $65=12+3+4\times 5+6+7+8+9.$
\item [] $66=1\times 2^3+4+5\times 6+7+8+9.$
\item [] $67=1+2^3+4+5\times 6+7+8+9.$
\item [] $68=1\times 2+3\times 4+5\times 6+7+8+9.$
\item [] $69=1+2+3\times 4+5\times 6+7+8+9.$
\item [] $70=1^2+34+5+6+7+8+9.$
\item [] $71=1\times 2+34+5+6+7+8+9.$
\item [] $72=1+2+34+5+6+7+8+9.$
\item [] $73=12+3+4+5\times 6+7+8+9.$
\item [] $74=1+2+3+4+5+6\times 7+8+9.$
\item [] $75=12\times 3+4+5+6+7+8+9.$
\item [] $76=1\times 2^3+4+5+6\times 7+8+9.$
\item [] $77=1^2+3\times 4+5+6\times 7+8+9.$
\item [] $78=12+3\times 4+5\times 6+7+8+9.$
\item [] $79=1+2+3\times 4+5+6\times 7+8+9.$
\item [] $80=1\times 2+3+45+6+7+8+9.$
\item [] $81=1+2+3+45+6+7+8+9.$
\item [] $82=1+2\times 3+45+6+7+8+9.$
\item [] $83=12+3+4+5+6\times 7+8+9.$
\item [] $84=1\times 2+3+4\times 5+6\times 7+8+9.$
\item [] $85=1+2+3+4\times 5+6\times 7+8+9.$
\item [] $86=1+2+3+4+5+6+7\times 8+9.$
\item [] $87=1+2\times 3+4+5+6+7\times 8+9.$
\item [] $88=12+3\times 4+5+6\times 7+8+9.$
\item [] $89=1\times 2+3+4+56+7+8+9.$
\item [] $90=12+3+45+6+7+8+9.$
\item [] $91=1+2+34+5\times 6+7+8+9.$
\item [] $92=1+23+4+5+6\times 7+8+9.$
\item [] $93=1+2+3\times 4\times 5+6+7+8+9.$
\item [] $94=1\times 2+3\times 4+56+7+8+9.$
\item [] $95=12+3+4+5+6+7\times 8+9.$
\item [] $96=1\times 2+3+4\times 5+6+7\times 8+9.$
\item [] $97=1+2+3+4\times 5+6+7\times 8+9.$
\item [] $98=1\times 23+45+6+7+8+9.$
\item [] $99=1+2+3+4+5+67+8+9.$
\item [] $100=1+2+3+4+5+6+7+8\times 9.$
\item [] $101=1+2+34+5+6\times 7+8+9.$
\item [] $102=12+3\times 4\times 5+6+7+8+9.$
\item [] $103=1\times 2\times 34+5+6+7+8+9.$
\item [] $104=1+23+4+5+6+7\times 8+9.$
\item [] $105=1+2\times 3\times 4+56+7+8+9.$
\item [] $106=12+3+4\times 5+6+7\times 8+9.$
\item [] $107=1\times 23+4+56+7+8+9.$
\item [] $108=1+2+3+4+5+6+78+9.$
\item [] $109=1+2\times 3+4+5+6+78+9.$
\item [] $110=12+34+5+6\times 7+8+9.$
\item [] $111=12\times 3+45+6+7+8+9.$
\item [] $112=1\times 2+3\times 4+5+6+78+9.$
\item [] $113=12+3\times 4+5+67+8+9.$
\item [] $114=1+2\times 3\times 4+5+67+8+9.$
\item [] $115=1+23+4\times 5+6+7\times 8+9.$
\item [] $116=1\times 2+34+56+7+8+9.$
\item [] $117=1+2+34+56+7+8+9.$
\item [] $118=1+23+4+5+6+7+8\times 9.$
\item [] $119=1+2+3+4\times 5+6+78+9.$
\item [] $120=12\times 3+4+56+7+8+9.$
\item [] $121=1\times 2+3\times 4+5+6+7+89.$
\item [] $122=1+2+3\times 4+5+6+7+89.$
\item [] $123=1+2\times 3\times 4+5+6+78+9.$
\item [] $124=1+2+3\times 4+5\times 6+7+8\times 9.$
\item [] $125=1\times 2+34+5+67+8+9.$
\item [] $126=12+34+56+7+8+9.$
\item [] $127=1+2+34+5+6+7+8\times 9.$
\item [] $128=1+2+3+4\times 5+6+7+89.$
\item [] $129=12\times 3+4+5+67+8+9.$
\item [] $130=1\times 2+3+4+56+7\times 8+9.$
\item[]$\mbox{Decreasing order}$
\item [] $61=9+8+7+(6+5+4+3)\times 2+1.$
\item [] $\mathit{62=98+7-65+43-21.}$
\item [] $63=9+8+7+6+5+4+3+21.$
\item [] $64=9+8+7+6\times 5+4+3+2+1.$
\item [] $65=9+8+7+6\times 5+4+3\times 2+1.$
\item [] $66=9+8+7+6+(5+4+3)\times (2+1).$
\item [] $67=9+8+7+6\times 5+4+3^2\times 1.$
\item [] $68=9+8+7+6+5+4\times 3+21.$
\item [] $69=9+8+7+6\times 5+4\times 3+2+1.$
\item [] $70=9+8+7+(6+5+4\times 3)\times 2\times 1.$
\item [] $71=9+8+7+6+5+4+32\times 1.$
\item [] $72=9+8+7+6+5+4+32+1.$
\item [] $73=9+8+7\times 6+5+4+3+2\times 1.$
\item [] $74=9+8+7\times 6+5+4+3+2+1.$
\item [] $75=9+8+7\times 6+5+4+3\times 2+1.$
\item [] $76=9+8+7+6+5\times (4+3+2)+1.$
\item [] $77=9+8+7\times 6+5+4+3^2\times 1.$
\item [] $78=9+8+7\times 6+5+4\times 3+2\times 1.$
\item [] $79=9+8+7\times 6+5+4\times 3+2+1.$
\item [] $80=9+8+7+6+5+43+2\times 1.$
\item [] $81=9+8+7+6+5+43+2+1.$
\item [] $82=9+8+7+6\times 5+4+3+21.$
\item [] $83=9+8+7+6+5\times 4+32+1.$
\item [] $84=9+8+7\times 6+5\times 4+3+2\times 1.$
\item [] $85=9+8+7\times 6+5\times 4+3+2+1.$
\item [] $86=9+8\times 7+6+5+4+3+2+1.$
\item [] $87=9+8\times 7+6+5+4+3\times 2+1.$
\item [] $88=9+8+7\times 6+5+4\times 3\times 2\times 1.$
\item [] $89=9+8+7\times 6+5+4\times 3\times 2+1.$
\item [] $90=9+8+7+6+54+3+2+1.$
\item [] $91=9+8+7+6+54+3\times 2+1.$
\item [] $92=9+8+7\times 6+5+4+3+21.$
\item [] $93=9+8+7+6+5\times 4\times 3+2+1.$
\item [] $94=9+8+7+6+54+3^2+1.$
\item [] $95=9+8+(7+6)\times 5+4+3^2\times 1.$
\item [] $96=9+8\times 7+6+5\times 4+3+2\times 1.$
\item [] $97=9+8\times 7+6+5\times 4+3+2+1.$
\item [] $98=9+8+7+65+4+3+2\times 1.$
\item [] $99=9+8+7+65+4+3+2+1.$
\item [] $100=9\times 8+7+6+5+4+3+2+1.$
\item [] $101=9\times 8+7+6+5+4+3\times 2+1.$
\item [] $102=9+8+7+6+5+4^3+2+1.$
\item [] $103=9+8+7\times 6+5\times 4+3+21.$
\item [] $104=9+8+7+65+4\times 3+2+1.$
\item [] $105=9+8\times 7+6\times 5+4+3+2+1.$
\item [] $106=9+8\times 7+6\times 5+4+3\times 2+1.$
\item [] $107=9+8+76+5+4+3+2\times 1.$
\item [] $108=9+8+76+5+4+3+2+1.$
\item [] $109=9+8+76+5+4+3\times 2+1.$
\item [] $110=9+8\times 7+6\times 5+4\times 3+2+1.$
\item [] $111=9\times 8+7+6+5\times 4+3+2+1.$
\item [] $112=9\times 8+7+6+5\times 4+3\times 2+1.$
\item [] $113=9+8+76+5+4\times 3+2+1.$
\item [] $114=9+8+7+65+4\times 3\times 2+1.$
\item [] $115=9\times 8+7+6+5+4\times 3\times 2+1.$
\item [] $116=9+87+6+5+4+3+2\times 1.$
\item [] $117=9+87+6+5+4+3+2+1.$
\item [] $118=9+87+6+5+4+3\times 2+1.$
\item [] $119=9\times 8+7+6\times 5+4+3+2+1.$
\item [] $120=9\times 8+7+6\times 5+4+3\times 2+1.$
\item [] $121=9+8+7\times 6+5\times 4\times 3+2\times 1.$
\item [] $122=9+87+6+5+4\times 3+2+1.$
\item [] $123=9+8+76+5+4\times 3\times 2+1.$
\item [] $124=9\times 8+7+6\times 5+4\times 3+2+1.$
\item [] $125=98+7+6+5+4+3+2\times 1.$
\item [] $126=98+7+6+5+4+3+2+1.$
\item [] $127=98+7+6+5+4+3\times 2+1.$
\item [] $128=9+87+6+5\times 4+3+2+1.$
\item [] $129=9\times 8+7\times 6+5+4+3+2+1.$
\item [] $130=9\times 8+7\times 6+5+4+3\times 2+1.$
\item[]$\mbox{Increasing order}$
\item [] $131=1+2+3+4+56+7\times 8+9.$
\item [] $132=1+2\times 3\times 4+5+6+7+89.$
\item [] $133=1\times 2\times 3\times 4+5\times 6+7+8\times 9.$
\item [] $134=1\times 2+34+5+6+78+9.$
\item [] $135=12+34+5+67+8+9.$
\item [] $136=12+34+5+6+7+8\times 9.$
\item [] $137=1+23+4\times 5+6+78+9.$
\item [] $138=12\times 3+4+5+6+78+9.$
\item [] $139=1\times 23+45+6+7\times 8+9.$
\item [] $140=12+3+4+56+7\times 8+9.$
\item [] $141=1+2+3\times 4+5\times 6+7+89.$
\item [] $142=1+2\times 3\times 4+5\times 6+78+9.$
\item [] $143=1\times 2+3+45+6+78+9.$
\item [] $144=12+34+5+6+78+9.$
\item [] $145=12+3+45+6+7+8\times 9.$
\item [] $146=1+2+3+4+5+6\times 7+89.$
\item [] $147=1+23+4+5+6\times 7+8\times 9.$
\item [] $148=1\times 2\times 34+56+7+8+9.$
\item [] $149=1+23+4+56+7\times 8+9.$
\item [] $150=1+2+3\times 4+56+7+8\times 9.$
\item [] $151=1+2+3\times 4+5+6\times 7+89.$
\item [] $152=1\times 2+3+45+6+7+89.$
\item [] $153=1+23+45+67+8+9.$
\item [] $154=1+2\times 3+4+56+78+9.$
\item [] $155=12+3+4+5+6\times 7+89.$
\item [] $156=12+3\times 4\times 5+67+8+9.$
\item [] $157=1\times 2+3\times 4+56+78+9.$
\item [] $158=1+2\times 34+5+67+8+9.$
\item [] $159=1+2\times 34+5+6+7+8\times 9.$
\item [] $160=12+3\times 4+5+6\times 7+89.$
\item [] $161=1\times 2+3+4+56+7+89.$
\item [] $162=123+4+5+6+7+8+9.$
\item [] $163=12+34+5\times 6+78+9.$
\item [] $164=1+23+4+5+6\times 7+89.$
\item [] $165=12\times 3+45+67+8+9.$
\item [] $166=1\times 2\times 34+5+6+78+9.$
\item [] $167=1+2+3\times 4+56+7+89.$
\item [] $168=1+2+3\times 45+6+7+8+9.$
\item [] $169=1+23\times 4+5+6+7\times 8+9.$
\item [] $170=1\times 23+4+56+78+9.$
\item [] $171=1+23+45+6+7+89.$
\item [] $172=1+23+4+5+67+8\times 9.$
\item [] $173=123+4\times 5+6+7+8+9.$
\item [] $174=12\times 3+45+6+78+9.$
\item [] $175=1\times 2\times 34+5+6+7+89.$
\item [] $176=12+3\times 4+56+7+89.$
\item [] $177=12+3\times 45+6+7+8+9.$
\item [] $178=1+2\times 34+5\times 6+7+8\times 9.$
\item [] $179=1\times 2+34+56+78+9.$
\item [] $180=1+2+34+56+78+9.$
\item [] $181=123+4+5\times 6+7+8+9.$
\item [] $182=1+2+3+4\times 5+67+89.$
\item [] $183=12\times 3+4+56+78+9.$
\item [] $184=12\times 3+4+5+67+8\times 9.$
\item [] $185=12+3\times 4+5+67+89.$
\item [] $186=1+2\times 3\times 4+5+67+89.$
\item [] $187=1\times 2\times 34+5+6\times 7+8\times 9.$
\item [] $188=1\times 2+34+56+7+89.$
\item [] $189=1+2+34+56+7+89.$
\item [] $190=1+2+3+45+67+8\times 9.$
\item [] $191=1+23\times 4+5+6+78+9.$
\item [] $192=12\times 3+4+56+7+89.$
\item [] $193=1\times 2+3\times 4\times 5+6\times 7+89.$
\item [] $194=1+2+3\times 4\times 5+6\times 7+89.$
\item [] $195=1+2\times 34+5\times 6+7+89.$
\item [] $196=1\times 2+3\times 45+6\times 7+8+9.$
\item [] $197=1\times 2+34+5+67+89.$
\item [] $198=12+34+56+7+89.$
\item [] $199=12+3+45+67+8\times 9.$
\item [] $200=1+23\times 4+5+6+7+89.$
\item[]$\mbox{Decreasing order}$
\item [] $131=98+7+6+5+4\times 3+2+1.$
\item [] $132=9+8\times 7+6+54+3\times 2+1.$
\item [] $133=9\times 8+7\times 6+5+4\times 3+2\times 1.$
\item [] $134=9\times 8+7\times 6+5+4\times 3+2+1.$
\item [] $135=9+8+7+65+43+2+1.$
\item [] $136=9\times 8+7+6+5+43+2+1.$
\item [] $137=98+7+6+5\times 4+3+2+1.$
\item [] $138=98+7+6+5\times 4+3\times 2+1.$
\item [] $139=9+8\times 7+65+4+3+2\times 1.$
\item [] $140=9\times 8+7\times 6+5\times 4+3+2+1.$
\item [] $141=9+87+6\times 5+4\times 3+2+1.$
\item [] $142=9\times 8+7+6\times 5+4\times 3+21.$
\item [] $143=9\times 8+7\times 6+5+4\times 3\times 2\times 1.$
\item [] $144=98+7+6+5+4+3+21.$
\item [] $145=98+7+6\times 5+4+3+2+1.$
\item [] $146=9\times 8+7+6+54+3\times 2+1.$
\item [] $147=9\times 8+7\times 6+5+4+3+21.$
\item [] $148=9\times 8+7+6+5\times 4\times 3+2+1.$
\item [] $149=9+8\times 7+6+54+3+21.$
\item [] $150=9+8+7+6\times 5\times 4+3+2+1.$
\item [] $151=9+8+7+6\times 5\times 4+3\times 2+1.$
\item [] $152=9+8+76+54+3+2\times 1.$
\item [] $153=9+8+76+54+3+2+1.$
\item [] $154=9\times 8+7+65+4+3+2+1.$
\item [] $155=9\times 8+7+6\times 5+43+2+1.$
\item [] $156=98+7\times 6+5+4+3\times 2+1.$
\item [] $157=9+8\times 7+6+54+32\times 1.$
\item [] $158=9+8\times 7+65+4+3+21.$
\item [] $159=9\times 8+7+65+4\times 3+2+1.$
\item [] $160=98+7\times 6+5+4\times 3+2+1.$
\item [] $161=9+87+6+54+3+2\times 1.$
\item [] $162=9+87+6+54+3+2+1.$
\item [] $163=9\times 8+76+5+4+3+2+1.$
\item [] $164=9\times 8+76+5+4+3\times 2+1.$
\item [] $165=9\times 8+7\times 6+5+43+2+1.$
\item [] $166=98+7\times 6+5\times 4+3+2+1.$
\item [] $167=98+7\times 6+5\times 4+3\times 2+1.$
\item [] $168=9\times 8+76+5+4\times 3+2+1.$
\item [] $169=9\times 8+7+65+4\times 3\times 2+1.$
\item [] $170=98+7\times 6+5+4\times 3\times 2+1.$
\item [] $171=9+87+65+4+3+2+1.$
\item [] $172=9+87+65+4+3\times 2+1.$
\item [] $173=98+7\times 6+5+4+3+21.$
\item [] $174=9\times 8+76+5\times 4+3+2+1.$
\item [] $175=9\times 8+76+5\times 4+3\times 2+1.$
\item [] $176=9+8\times 7+65+43+2+1.$
\item [] $177=9\times 8+7\times 6+5\times 4\times 3+2+1.$
\item [] $178=9\times 8+76+5+4\times 3\times 2+1.$
\item [] $179=9+8+76+54+32\times 1.$
\item [] $180=98+7+65+4+3+2+1.$
\item [] $181=98+7+6\times 5+43+2+1.$
\item [] $182=98+7\times 6+5+4+32+1.$
\item [] $183=98+7+6+5+4+3\times 21.$
\item [] $184=98+7\times 6+5\times 4+3+21.$
\item [] $185=98+7+65+4\times 3+2+1.$
\item [] $186=9+87+65+4\times 3\times 2+1.$
\item [] $187=(9+8+7+6)\times 5+4+32+1.$
\item [] $188=98+76+5+4+3+2\times 1.$
\item [] $189=98+76+5+4+3+2+1.$
\item [] $190=98+76+5+4+3\times 2+1.$
\item [] $191=98+7\times 6+5+43+2+1.$
\item [] $192=9\times 8+7\times 6+54+3+21.$
\item [] $193=98+76+5+4\times 3+2\times 1.$
\item [] $194=98+76+5+4\times 3+2+1.$
\item [] $195=9+8+7+6+54\times 3+2+1.$
\item [] $196=9+8\times 7+65+4^3+2\times 1.$
\item [] $197=9+87+65+4+32\times 1.$
\item [] $198=98+7+65+4+3+21.$
\item [] $199=9\times 8+76+5+43+2+1.$
\item [] $200=98+7\times 6+54+3+2+1.$
\item[]$\mbox{Increasing order}$
\item [] $201=12\times 3+4+5+67+89.$
\item [] $202=123+4\times 5+6\times 7+8+9.$
\item [] $203=123+4+5+6+7\times 8+9.$
\item [] $204=1+2\times 34+56+7+8\times 9.$
\item [] $205=1+2\times 3\times 4\times 5+67+8+9.$
\item [] $206=1\times 2+3+45+67+89.$
\item [] $207=1+2+3+45+67+89.$
\item [] $208=1+2\times 3+45+67+89.$
\item [] $209=1\times 23\times 4+5\times 6+78+9.$
\item [] $210=1+23\times 4+5\times 6+78+9.$
\item [] $211=1\times 2\times 34+56+78+9.$
\item [] $212=1+2\times 34+56+78+9.$
\item [] $213=1+2\times 34+5+67+8\times 9.$
\item [] $214=123+4\times 5+6+7\times 8+9.$
\item [] $215=1\times 2^3+4\times 5\times 6+78+9.$
\item [] $216=12+3+45+67+89.$
\item [] $217=123+4+5+6+7+8\times 9.$
\item [] $218=12+3\times 45+6+7\times 8+9.$
\item [] $219=1+2+3\times 4\times 5+67+89.$
\item [] $220=1\times 2\times 34+56+7+89.$
\item [] $221=1+2\times 34+56+7+89.$
\item [] $222=1+2+3\times 45+67+8+9.$
\item [] $223=1+2\times 3\times 4\times 5+6+7+89.$
\item [] $224=1\times 23+45+67+89.$
\item [] $225=1+23+45+67+89.$
\item [] $226=1\times 2+3^4+56+78+9.$
\item [] $227=123+4\times 5+67+8+9.$
\item [] $228=1\times 23\times 4+5+6\times 7+89.$
\item [] $229=1\times 2\times 34+5+67+89.$
\item [] $230=1+2\times 34+5+67+89.$
\item [] $231=12+3\times 45+67+8+9.$
\item [] $232=12+3\times 45+6+7+8\times 9.$
\item [] $233=12\times 3\times 4+5+67+8+9.$
\item [] $234=123+4+5+6+7+89.$
\item [] $235=1\times 23\times 4+56+78+9.$
\item [] $236=1+23\times 4+56+78+9.$
\item [] $237=12\times 3+45+67+89.$
\item [] $238=1+2\times 3+4+5\times 6\times 7+8+9.$
\item [] $239=123+45+6+7\times 8+9.$
\item [] $240=1+2+3\times 45+6+7+89.$
\item [] $241=12+34\times 5+6\times 7+8+9.$
\item [] $242=12\times 3\times 4+5+6+78+9.$
\item [] $243=12\times 3+4\times 5\times 6+78+9.$
\item [] $244=123+4+5\times 6+78+9.$
\item [] $245=123+4\times 5+6+7+89.$
\item [] $246=123+4+5+6\times 7+8\times 9.$
\item [] $247=12+3^4+5\times (6+7)+89.$
\item [] $248=123+4+56+7\times 8+9.$
\item [] $249=12+3\times 45+6+7+89.$
\item [] $250=1^2+3\times 45+6\times 7+8\times 9.$
\item [] $251=12\times 3\times 4+5+6+7+89.$
\item [] $252=123+45+67+8+9.$
\item [] $253=123+4+5\times 6+7+89.$
\item [] $254=1+23\times 4+5+67+89.$
\item [] $255=1+2\times 3+4\times 56+7+8+9.$
\item [] $256=1\times 2+34\times 5+67+8+9.$
\item [] $257=1+2+34\times 5+67+8+9.$
\item [] $258=1+2+34\times 5+6+7+8\times 9.$
\item [] $259=1\times 2\times 3\times 4\times 5+67+8\times 9.$
\item [] $260=1+2\times 3\times 4\times 5+67+8\times 9.$
\item [] $261=123+45+6+78+9.$
\item [] $262=123+4+56+7+8\times 9.$
\item [] $263=12+3+4\times 56+7+8+9.$
\item [] $264=1+2+34+5\times 6\times 7+8+9.$
\item [] $265=1\times 2+34\times 5+6+78+9.$
\item [] $266=12+34\times 5+67+8+9.$
\item [] $267=123+4\times 5\times 6+7+8+9.$
\item [] $268=1\times 2+3\times 45+6\times 7+89.$
\item [] $269=1+2+3\times 45+6\times 7+89.$
\item [] $270=123+4+56+78+9.$
\item[]$\mbox{Decreasing order}$
\item [] $201=98+76+5\times 4+3\times 2+1.$
\item [] $202=98+7\times 6+5\times 4\times 3+2\times 1.$
\item [] $203=98+76+5+4\times 3\times 2\times 1.$
\item [] $204=9\times 8+7+6\times 5\times 4+3+2\times 1.$
\item [] $205=9\times 8+7+6\times 5\times 4+3+2+1.$
\item [] $206=9\times 8+7+6+5\times 4\times 3\times 2+1.$
\item [] $207=9+87+65+43+2+1.$
\item [] $208=9\times 8+76+54+3+2+1.$
\item [] $209=9\times 8+76+54+3\times 2+1.$
\item [] $210=9\times 8+76+5\times 4\times 3+2\times 1.$
\item [] $211=9\times 8+76+5\times 4\times 3+2+1.$
\item [] $212=9\times 8+76+54+3^2+1.$
\item [] $213=9+8+76+5\times 4\times 3\times 2\times 1.$
\item [] $214=9+8+76+5\times 4\times 3\times 2+1.$
\item [] $215=98+7+65+43+2\times 1.$
\item [] $216=98+7+65+43+2+1.$
\item [] $217=9\times 8+76+5+43+21.$
\item [] $218=98+76+5\times 4+3+21.$
\item [] $219=9+87+6+54+3\times 21.$
\item [] $220=9\times 8+76+5+4+3\times 21.$
\item [] $221=9+87+6\times 5\times 4+3+2\times 1.$
\item [] $222=9+87+6\times 5\times 4+3+2+1.$
\item [] $223=9+87+6\times 5\times 4+3\times 2+1.$
\item [] $224=9+8+7\times 6+54\times 3+2+1.$
\item [] $225=98+76+5+43+2+1.$
\item [] $226=9\times 8+76+54+3+21.$
\item [] $227=98+76+5\times 4+32+1.$
\item [] $228=9+87+65+4+3\times 21.$
\item [] $229=9\times 8+76+5\times 4\times 3+21.$
\item [] $230=9\times 8+7+65+43\times 2\times 1.$
\item [] $231=98+7+6\times 5\times 4+3+2+1.$
\item [] $232=98+7+6\times 5\times 4+3\times 2+1.$
\item [] $233=98+76+54+3+2\times 1.$
\item [] $234=98+76+54+3+2+1.$
\item [] $235=98+76+54+3\times 2+1.$
\item [] $236=9+8\times 7+6+54\times 3+2+1.$
\item [] $237=98+76+5\times 4\times 3+2+1.$
\item [] $238=9+8+7\times 6\times 5+4+3\times 2+1.$
\item [] $239=9\times 8+76+5+43\times 2\times 1.$
\item [] $240=9+87+6\times 5\times 4+3+21.$
\item [] $241=9+8+7\times 6\times 5+4\times 3+2\times 1.$
\item [] $242=9+8+7\times 6\times 5+4\times 3+2+1.$
\item [] $243=98+76+5+43+21.$
\item [] $244=98+7+6+5+4\times 32\times 1.$
\item [] $245=98+7+6+5+4\times 32+1.$
\item [] $246=98+76+5+4+3\times 21.$
\item [] $247=9+87+65+43\times 2\times 1.$
\item [] $248=9+8+7+6+5\times 43+2+1.$
\item [] $249=98+7+6\times 5\times 4+3+21.$
\item [] $250=9\times 8+7+6+54\times 3+2+1.$
\item [] $251=9+8+7\times 6\times 5+4\times 3\times 2\times 1.$
\item [] $252=98+76+54+3+21.$
\item [] $253=9\times (8+7+6)+54+3^2+1.$
\item [] $254=9+8\times 7+6+54\times 3+21.$
\item [] $255=9+8+7\times 6\times 5+4+3+21.$
\item [] $256=98+7+65+43\times 2\times 1.$
\item [] $257=9+8+76+54\times 3+2\times 1.$
\item [] $258=9+8+76+54\times 3+2+1.$
\item [] $259=9+8\times 7+65+4\times 32+1.$
\item [] $260=98+7\times 6+5\times 4\times 3\times 2\times 1.$
\item [] $261=98+76+54+32+1.$
\item [] $262=9\times 8+7+6\times 5\times 4+3\times 21.$
\item [] $263=9+8+7\times 6\times 5+4+32\times 1.$
\item [] $264=9+8+7\times 6\times 5+4+32+1.$
\item [] $265=98+76+5+43\times 2\times 1.$
\item [] $266=9+87+6+54\times 3+2\times 1.$
\item [] $267=9+87+6+54\times 3+2+1.$
\item [] $268=9\times 8+76+5\times 4\times 3\times 2\times 1.$
\item [] $269=9\times 8+76+5\times 4\times 3\times 2+1.$
\item [] $270=98+7+6\times (5+4)\times 3+2+1.$
\item[]$\mbox{Increasing order}$
\item [] $271=123+4+5+67+8\times 9.$
\item [] $272=1+23+4\times 56+7+8+9.$
\item [] $273=12+34+5\times 6\times 7+8+9.$
\item [] $274=1\times 2+34\times 5+6+7+89.$
\item [] $275=12+34\times 5+6+78+9.$
\item [] $276=1\times 2+3\times 45+67+8\times 9.$
\item [] $277=1+2+3\times 45+67+8\times 9.$
\item [] $278=12+3\times 45+6\times 7+89.$
\item [] $279=123+4+56+7+89.$
\item [] $280=12\times 3\times 4+5+6\times 7+89.$
\item [] $281=1\times 2\times 3^4+5+6\times 7+8\times 9.$
\item [] $282=123+4\times 5+67+8\times 9.$
\item [] $283=1\times 2\times 3^4+56+7\times 8+9.$
\item [] $284=12+34\times 5+6+7+89.$
\item [] $285=1^2+34\times 5+6\times 7+8\times 9.$
\item [] $286=12+3\times 45+67+8\times 9.$
\item [] $287=12\times 3\times 4+56+78+9.$
\item [] $288=123+4+5+67+89.$
\item [] $289=1+234+5\times 6+7+8+9.$
\item [] $290=1^2+3+4+5\times 6\times 7+8\times 9.$
\item [] $291=1\times 2+3+4+5\times 6\times 7+8\times 9.$
\item [] $292=1+2\times 3+4+56\times 7+89.$
\item [] $293=1\times 2+3\times 45+67+89.$
\item [] $294=1+2+3\times 45+67+89.$
\item [] $295=1+2+3+4\times 56+7\times 8+9.$
\item [] $296=12\times 3\times 4+56+7+89.$
\item [] $297=1+2+3\times 4+5\times 6\times 7+8\times 9.$
\item [] $298=1\times 234+5+6\times 7+8+9.$
\item [] $299=123+4\times 5+67+89.$
\item [] $300=1+2+3+45\times 6+7+8+9.$
\item [] $301=1+2\times 3\times 45+6+7+8+9.$
\item [] $302=1^2+34\times 5+6\times 7+89.$
\item [] $303=12+3\times 45+67+89.$
\item [] $304=1+2+34\times 5+6\times 7+89.$
\item [] $305=12\times 3\times 4+5+67+89.$
\item [] $306=12+3\times 4+5\times 6\times 7+8\times 9.$
\item [] $307=123+45+67+8\times 9.$
\item [] $308=123+4\times 5\times 6+7\times 8+9.$
\item [] $309=12+3+45\times 6+7+8+9.$
\item [] $310=1+23+4+5\times 6\times 7+8\times 9.$
\item [] $311=1+234+5+6+7\times 8+9.$
\item [] $312=12+34\times 5+6+7+8+9.$
\item [] $313=12+34\times 5+6\times 7+89.$
\item [] $314=1\times 234+56+7+8+9.$
\item [] $315=1+234+56+7+8+9.$
\item [] $316=1\times 2+3+4\times 56+78+9.$
\item [] $317=1+2+3+4\times 56+78+9.$
\item [] $318=1+23+45\times 6+7+8+9.$
\item [] $319=1\times 23\times 4+5\times 6\times 7+8+9.$
\item [] $320=1+23\times 4+5\times 6\times 7+8+9.$
\item [] $321=12+34\times 5+67+8\times 9.$
\item [] $322=123+4\times 5\times 6+7+8\times 9.$
\item [] $323=1\times 234+5+67+8+9.$
\item [] $324=123+45+67+89.$
\item [] $325=1+234+5+6+7+8\times 9.$
\item [] $326=12+3+4\times 56+78+9.$
\item [] $327=1+23+4+5\times 6\times 7+89.$
\item [] $328=1\times 2+34\times 5+67+89.$
\item [] $329=1+2+34\times 5+67+89.$
\item [] $330=1+234+5\times 6+7\times 8+9.$
\item [] $331=1\times 2^3\times 4+5\times 6\times 7+89.$
\item [] $332=1\times 234+5+6+78+9.$
\item [] $333=1+234+5+6+78+9.$
\item [] $334=(1\times 2\times 3\times 4+5+6)\times 7+89.$
\item [] $335=12+3+4\times 56+7+89.$
\item [] $336=1+2+34+5\times 6\times 7+89.$
\item [] $337=1+(2+34)\times 5+67+89.$
\item [] $338=12+34\times 5+67+89.$
\item [] $339=123+4\times 5\times 6+7+89.$
\item [] $340=1\times 2+3+45\times 6+7\times 8+9.$
\item[]$\mbox{Decreasing order}$
\item [] $271=(9+8)\times 7+65+43\times 2+1.$
\item [] $272=9\times 8+7+65+4\times 32\times 1.$
\item [] $273=9+8+7\times 6\times 5+43+2+1.$
\item [] $274=98+7\times 6+5+4\times 32+1.$
\item [] $275=98+7+6+54\times 3+2\times 1.$
\item [] $276=98+7+6+54\times 3+2+1.$
\item [] $277=9+8+7\times 6+5\times 43+2+1.$
\item [] $278=9\times 8+7\times 6+54\times 3+2\times 1.$
\item [] $279=9\times 8+7\times 6+54\times 3+2+1.$
\item [] $280=98+7\times 6+5\times (4+3+21).$
\item [] $281=9\times 8+76+5+4\times 32\times 1.$
\item [] $282=9\times 8+76+5+4\times 32+1.$
\item [] $283=9+87+6+5\times 4\times 3^2+1.$
\item [] $284=(9+8)\times (7+6)+5\times 4\times 3+2+1.$
\item [] $285=9+87+6+54\times 3+21.$
\item [] $286=98+7\times (6+5\times 4)+3+2+1.$
\item [] $287=9+8+7+6+5+4\times 3\times 21.$
\item [] $288=98+7+6\times 5\times 4+3\times 21.$
\item [] $289=98+7+65\times 4+3+21.$
\item [] $290=9+8+7+65\times 4+3+2+1.$
\item [] $291=9+8+7+65\times 4+3\times 2+1.$
\item [] $292=9\times 8+7\times 6\times 5+4+3+2+1.$
\item [] $293=9\times 8+7\times 6\times 5+4+3\times 2+1.$
\item [] $294=98+76+5\times 4\times 3\times 2\times 1.$
\item [] $295=98+76+5\times 4\times 3\times 2+1.$
\item [] $296=9\times 8+7\times 6\times 5+4\times 3+2\times 1.$
\item [] $297=9\times 8+7\times 6\times 5+4\times 3+2+1.$
\item [] $298=98+7+65+4\times 32\times 1.$
\item [] $299=98+7+65+4\times 32+1.$
\item [] $300=9+87+(6\times 5+4)\times 3\times 2\times 1.$
\item [] $301=9+8+7+6+54\times (3+2)+1.$
\item [] $302=9\times 8+7+6+5\times 43+2\times 1.$
\item [] $303=9\times 8+7+6+5\times 43+2+1.$
\item [] $304=98+7\times 6+54\times 3+2\times 1.$
\item [] $305=98+7\times 6+54\times 3+2+1.$
\item [] $306=9\times 8+7\times 6\times 5+4\times 3\times 2\times 1.$
\item [] $307=9\times 8+7\times 6\times 5+4\times 3\times 2+1.$
\item [] $308=9+8+7+65\times 4+3+21.$
\item [] $309=9+87+6\times 5\times (4+3)+2+1.$
\item [] $310=9\times 8+7\times 6\times 5+4+3+21.$
\item [] $311=9+8+76+5\times 43+2+1.$
\item [] $312=9\times 8+76+54\times 3+2\times 1.$
\item [] $313=9\times 8+76+54\times 3+2+1.$
\item [] $314=98+76+5\times 4\times (3\times 2+1).$
\item [] $315=9\times 8+7\times 6\times 5+4\times 3+21.$
\item [] $316=9+8+7+65\times 4+32\times 1.$
\item [] $317=98+7\times 6\times 5+4+3+2\times 1.$
\item [] $318=98+7\times 6\times 5+4+3+2+1.$
\item [] $319=98+7\times 6\times 5+4+3\times 2+1.$
\item [] $320=9+87+6+5\times 43+2+1.$
\item [] $321=9\times 8+7+6+5\times 43+21.$
\item [] $322=98+7\times 6\times 5+4\times 3+2\times 1.$
\item [] $323=98+7\times 6\times 5+4\times 3+2+1.$
\item [] $324=9\times (8+7)+6+54\times 3+21.$
\item [] $325=98+7\times 6+5\times (4+32+1).$
\item [] $326=98+76+5+(4+3)\times 21.$
\item [] $327=9\times 8+7\times 6\times 5+43+2\times 1.$
\item [] $328=9\times 8+7\times 6\times 5+43+2+1.$
\item [] $329=98+7+6+5\times 43+2+1.$
\item [] $330=9+8\times 7+65\times 4+3+2\times 1.$
\item [] $331=9+8\times 7+65\times 4+3+2+1.$
\item [] $332=9+8\times 7+65\times 4+3\times 2+1.$
\item [] $333=98+7\times 6\times 5+4\times 3\times 2+1.$
\item [] $334=9+8\times 7+65\times 4+3^2\times 1.$
\item [] $335=9+(8+7+65)\times 4+3+2+1.$
\item [] $336=98+7\times 6\times 5+4+3+21.$
\item [] $337=98+7+6+5\times (43+2)+1.$
\item [] $338=98+76+54\times 3+2\times 1.$
\item [] $339=98+76+54\times 3+2+1.$
\item [] $340=9+8+(7+6)\times 5\times 4+3\times 21.$
\item[]$\mbox{Increasing order}$
\item [] $341=1\times 234+5+6+7+89.$
\item [] $342=1+234+5+6+7+89.$
\item [] $343=1\times 23+4\times 56+7+89.$
\item [] $344=1+23+4\times 56+7+89.$
\item [] $345=12+34+5\times 6\times 7+89.$
\item [] $346=1^{2345}+6\times 7\times 8+9.$
\item [] $347=12\times 3+4\times 56+78+9.$
\item [] $348=1\times 234+5\times (6+7+8)+9.$
\item [] $349=1^{23}\times 45\times 6+7+8\times 9.$
\item [] $350=12+3+45\times 6+7\times 8+9.$
\item [] $351=1\times 234+5\times 6+78+9.$
\item [] $352=1+234+5\times 6+78+9.$
\item [] $353=1\times 234+5+6\times 7+8\times 9.$
\item [] $354=123+4+5\times 6\times 7+8+9.$
\item [] $355=1+2+3+45\times 6+7+8\times 9.$
\item [] $356=1+234+56+7\times 8+9.$
\item [] $357=1^2\times 3+4+5+6\times 7\times 8+9.$
\item [] $358=1\times 23+45\times 6+7\times 8+9.$
\item [] $359=1+23+45\times 6+7\times 8+9.$
\item [] $360=1+2+3+4+5+6\times 7\times 8+9.$
\item [] $361=1+234+5\times 6+7+89.$
\item [] $362=1+2+3+4+5\times 67+8+9.$
\item [] $363=1+2+3+45\times 6+78+9.$
\item [] $364=12+3+45\times 6+7+8\times 9.$
\item [] $365=1+2+3\times 4+5+6\times 7\times 8+9.$
\item [] $366=1\times 2+3\times 4+5\times 67+8+9.$
\item [] $367=1\times 2\times 34+5\times 6\times 7+89.$
\item [] $368=1+2\times 34+5\times 6\times 7+89.$
\item [] $369=1\times 234+56+7+8\times 9.$
\item [] $370=1+234+56+7+8\times 9.$
\item [] $371=1+234+5+6\times 7+89.$
\item [] $372=12+3+45\times 6+78+9.$
\item [] $373=1+2\times 3\times 45+6+7+89.$
\item [] $374=12+3\times 4+5+6\times 7\times 8+9.$
\item [] $375=1+23\times 4+5\times 6\times 7+8\times 9.$
\item [] $376=12+3\times 4+5\times 67+8+9.$
\item [] $377=1\times 234+56+78+9.$
\item [] $378=1+234+56+78+9.$
\item [] $379=1+234+5+67+8\times 9.$
\item [] $380=1+23+4+5\times 67+8+9.$
\item [] $381=1+23+45\times 6+78+9.$
\item [] $382=1\times 2+3^4+5\times 6\times 7+89.$
\item [] $383=1+2^3\times 4+5+6\times 7\times 8+9.$
\item [] $384=1\times 2\times 3\times 45+6\times 7+8\times 9.$
\item [] $385=1+2\times 3\times 45+6\times 7+8\times 9.$
\item [] $386=1\times 234+56+7+89.$
\item [] $387=12+345+6+7+8+9.$
\item [] $388=1\times 2+34+5\times 67+8+9.$
\item [] $389=1\times 23+45\times 6+7+89.$
\item [] $390=1+23+45\times 6+7+89.$
\item [] $391=1\times 23\times 4+5\times 6\times 7+89.$
\item [] $392=1+23\times 4+5\times 6\times 7+89.$
\item [] $393=12\times 3+45\times 6+78+9.$
\item [] $394=1^2+3+45+6\times 7\times 8+9.$
\item [] $395=1\times 234+5+67+89.$
\item [] $396=1+234+5+67+89.$
\item [] $397=1+2\times 3+45+6\times 7\times 8+9.$
\item [] $398=12+34+5\times 67+8+9.$
\item [] $399=1\times 2\times 34\times 5+6\times 7+8+9.$
\item [] $400=1+2\times 34\times 5+6\times 7+8+9.$
\item [] $401=1\times 2\times 3\times 45+6\times 7+89.$
\item [] $402=12\times 3+45\times 6+7+89.$
\item [] $403=1+(2\times 3+45)\times 6+7+89.$
\item [] $404=1^2\times 345+6\times 7+8+9.$
\item [] $405=12+3+45+6\times 7\times 8+9.$
\item [] $406=1\times 2+345+6\times 7+8+9.$
\item [] $407=1+2+345+6\times 7+8+9.$
\item [] $408=1+2+3\times 4\times 5+6\times 7\times 8+9.$
\item [] $409=123+4+5\times 6\times 7+8\times 9.$
\item [] $410=1+2\times 3\times 45+67+8\times 9.$
\item[]$\mbox{Decreasing order}$
\item [] $341=98+7\times 6\times 5+4\times 3+21.$
\item [] $342=9\times 8+7+6+5+4\times 3\times 21.$
\item [] $343=(9+8)\times 7+6+5\times 43+2+1.$
\item [] $344=9\times 8+7+65\times 4+3+2\times 1.$
\item [] $345=9\times 8+7+65\times 4+3+2+1.$
\item [] $346=9\times 8+7+65\times 4+3\times 2+1.$
\item [] $347=98+7+6+5\times 43+21.$
\item [] $348=9\times 8+7+65\times 4+3^2\times 1.$
\item [] $349=9+8\times 7+65\times 4+3+21.$
\item [] $350=9+8+76+5+4\times 3\times 21.$
\item [] $351=9+(87+6+5\times 4)\times 3+2+1.$
\item [] $352=9+8+7+6+5\times 4^3+2\times 1.$
\item [] $353=9+8+7+6\times 54+3+2\times 1.$
\item [] $354=9+8+7+6\times 54+3+2+1.$
\item [] $355=9+8+7+6\times 54+3\times 2+1.$
\item [] $356=9+8+7\times 6\times 5+4\times 32+1.$
\item [] $357=98+76+54\times 3+21.$
\item [] $358=98+7\times 6+5\times 43+2+1.$
\item [] $359=9+8\times 7\times 6+5+4+3+2\times 1.$
\item [] $360=9+8+7+6+5+4+321.$
\item [] $361=9+8\times 7\times 6+5+4+3\times 2+1.$
\item [] $362=9+87+65\times 4+3+2+1.$
\item [] $363=9+87+65\times 4+3\times 2+1.$
\item [] $364=9+8\times 7\times 6+5+4\times 3+2\times 1.$
\item [] $365=9+8\times 7\times 6+5+4\times 3+2+1.$
\item [] $366=9\times 8+76+5\times 43+2+1.$
\item [] $367=9\times 8+7\times 6\times 5+4^3+21.$
\item [] $368=98+7+6+5+4\times 3\times 21.$
\item [] $369=9+8+7+6\times (54+3)+2+1.$
\item [] $370=98+7+65\times 4+3+2\times 1.$
\item [] $371=98+7+65\times 4+3+2+1.$
\item [] $372=98+7+65\times 4+3\times 2+1.$
\item [] $373=9+87+6+54\times (3+2)+1.$
\item [] $374=9+8\times 7\times 6+5+4\times 3\times 2\times 1.$
\item [] $375=9+8\times 7\times 6+5+4\times 3\times 2+1.$
\item [] $376=98+7\times 6+5\times 43+21.$
\item [] $377=9+8\times (7+6+5+4+3+21).$
\item [] $378=9+8\times 7\times 6+5+4+3+21.$
\item [] $379=9+8+7+6\times 5+4+321.$
\item [] $380=9+87+65\times 4+3+21.$
\item [] $381=9+8+7+6\times 54+32+1.$
\item [] $382=9+8\times 7+65+4\times 3\times 21.$
\item [] $383=9+8+7\times 6+54\times 3\times 2\times 1.$
\item [] $384=9+8+7\times 6+54\times 3\times 2+1.$
\item [] $385=(9+8)\times 7+65\times 4+3+2+1.$
\item [] $386=9+8+7+6\times 5\times 4\times 3+2\times 1.$
\item [] $387=9+8+7+6\times 5\times 4\times 3+2+1.$
\item [] $388=9+87+65\times 4+32\times 1.$
\item [] $389=9+8\times 7\times 6+5\times 4+3+21.$
\item [] $390=9\times 8+7\times (6+5+4)\times 3+2+1.$
\item [] $391=98+76+5\times 43+2\times 1.$
\item [] $392=98+76+5\times 43+2+1.$
\item [] $393=98+7\times 6\times 5+4^3+21.$
\item [] $394=9+8\times 7+6\times 54+3+2\times 1.$
\item [] $395=9+8\times 7+6\times 54+3+2+1.$
\item [] $396=9+8\times 7+6\times 54+3\times 2+1.$
\item [] $397=98+7+65\times 4+32\times 1.$
\item [] $398=98+7+65\times 4+32+1.$
\item [] $399=9+8\times 7+6\times 54+3^2+1.$
\item [] $400=9+8+7\times 6+5\times 4+321.$
\item [] $401=9+8\times 7+6+5+4+321.$
\item [] $402=9\times 8+7+65\times 4+3\times 21.$
\item [] $403=(9+8)\times 7+65\times 4+3+21.$
\item [] $404=9+8\times 7\times 6+54+3+2\times 1.$
\item [] $405=9+8+7+6+54+321.$
\item [] $406=9+8\times 7\times 6+54+3\times 2+1.$
\item [] $407=9+8+76\times 5+4+3+2+1.$
\item [] $408=9+8+76\times 5+4+3\times 2+1.$
\item [] $409=9\times 8+7+6\times 54+3+2+1.$
\item [] $410=98+76+5\times 43+21.$
\item[]$\mbox{Increasing order}$
\item [] $411=1\times 2\times 34\times 5+6+7\times 8+9.$
\item [] $412=123+4\times 56+7\times 8+9.$
\item [] $413=1\times 23+45+6\times 7\times 8+9.$
\item [] $414=1+23+45+6\times 7\times 8+9.$
\item [] $415=1^2+3+4+5\times 67+8\times 9.$
\item [] $416=12+345+6\times 7+8+9.$
\item [] $417=123+45\times 6+7+8+9.$
\item [] $418=1\times 2+345+6+7\times 8+9.$
\item [] $419=1+2+345+6+7\times 8+9.$
\item [] $420=1+2\times 3+4+56\times 7+8+9.$
\item [] $421=1+2\times 34+5\times 67+8+9.$
\item [] $422=1+2+3\times 4+5\times 67+8\times 9.$
\item [] $423=1\times 2+3\times 4+56\times 7+8+9.$
\item [] $424=1+2+3\times 4+56\times 7+8+9.$
\item [] $425=1+2\times 34\times 5+67+8+9.$
\item [] $426=1\times 2\times 3\times 45+67+89.$
\item [] $427=1+2\times 3\times 45+67+89.$
\item [] $428=12+345+6+7\times 8+9.$
\item [] $429=1^2\times 345+67+8+9.$
\item [] $430=1^2+345+67+8+9.$
\item [] $431=1\times 2+345+67+8+9.$
\item [] $432=1+2+345+67+8+9.$
\item [] $433=1+2+345+6+7+8\times 9.$
\item [] $434=123+4\times 56+78+9.$
\item [] $435=1+23+4+5\times 67+8\times 9.$
\item [] $436=1\times 23+4+56\times 7+8+9.$
\item [] $437=1+23+4+56\times 7+8+9.$
\item [] $438=1\times 2+3\times 4+5\times 67+89.$
\item [] $439=1+2+3\times 4+5\times 67+89.$
\item [] $440=1\times 2+345+6+78+9.$
\item [] $441=12+345+67+8+9.$
\item [] $442=12+345+6+7+8\times 9.$
\item [] $443=123+4\times 56+7+89.$
\item [] $444=1+2+34+5\times 67+8\times 9.$
\item [] $445=1+23\times 4+5\times 67+8+9.$
\item [] $446=1+2+34+56\times 7+8+9.$
\item [] $447=12\times 3+4+5\times 67+8\times 9.$
\item [] $448=12+3\times 4+5\times 67+89.$
\item [] $449=1\times 2+345+6+7+89.$
\item [] $450=12+345+6+78+9.$
\item [] $451=1\times 23+4+5\times 67+89.$
\item [] $452=1+23+4+5\times 67+89.$
\item [] $453=12+34+5\times 67+8\times 9.$
\item [] $454=1\times 2\times 34\times 5+6\times 7+8\times 9.$
\item [] $455=12+34+56\times 7+8+9.$
\item [] $456=1^2\times 3\times 4\times 5\times 6+7+89.$
\item [] $457=1+2^3\times 4+5\times 67+89.$
\item [] $458=123+45\times 6+7\times 8+9.$
\item [] $459=12+345+6+7+89.$
\item [] $460=1\times 2+34+5\times 67+89.$
\item [] $461=1+2+34+5\times 67+89.$
\item [] $462=1+234+5\times 6\times 7+8+9.$
\item [] $463=1+2^3\times 45+6+7+89.$
\item [] $464=12\times 3+4+5\times 67+89.$
\item [] $465=1\times 2\times 3\times 4\times 5+6\times 7\times 8+9.$
\item [] $466=1+2\times 3\times 4\times 5+6\times 7\times 8+9.$
\item [] $467=1+2+3^4\times 5+6\times 7+8+9.$
\item [] $468=12+3\times 4\times 5\times 6+7+89.$
\item [] $469=1^{23}+4+56\times 7+8\times 9.$
\item [] $470=12+34+5\times 67+89.$
\item [] $471=12+345+6\times 7+8\times 9.$
\item [] $472=123+45\times 6+7+8\times 9.$
\item [] $473=1\times 2+3+4+56\times 7+8\times 9.$
\item [] $474=1+2+3+4+56\times 7+8\times 9.$
\item [] $475=1\times 2\times 34+5\times 67+8\times 9.$
\item [] $476=1+2\times 34+5\times 67+8\times 9.$
\item [] $477=123+4+5+6\times 7\times 8+9.$
\item [] $478=1\times 2+345+6\times 7+89.$
\item [] $479=123+4+5\times 67+8+9.$
\item [] $480=123+45\times 6+78+9.$
\item[]$\mbox{Decreasing order}$
\item [] $411=9+8+76\times 5+4\times 3+2\times 1.$
\item [] $412=9+8+76\times 5+4\times 3+2+1.$
\item [] $413=9+8\times 7+6\times 54+3+21.$
\item [] $414=9+8+7+65+4+321.$
\item [] $415=9\times 8+7+6+5+4+321.$
\item [] $416=9+8+76+5\times 4^3+2+1.$
\item [] $417=9+8+76+54\times 3\times 2\times 1.$
\item [] $418=9+8+76+54\times 3\times 2+1.$
\item [] $419=9+87+65\times 4+3\times 21.$
\item [] $420=9+8\times 7+6\times 5+4+321.$
\item [] $421=9+8+76\times 5+4\times 3\times 2\times 1.$
\item [] $422=9+8+76\times 5+4\times 3\times 2+1.$
\item [] $423=9+8+76+5+4+321.$
\item [] $424=9+87+6+5\times 4^3+2\times 1.$
\item [] $425=9+8+76\times 5+4+3+21.$
\item [] $426=9+87+6\times 54+3+2+1.$
\item [] $427=9+87+6\times 54+3\times 2+1.$
\item [] $428=9+8\times 7+6\times 5\times 4\times 3+2+1.$
\item [] $429=9+87+6\times 54+3^2\times 1.$
\item [] $430=9+8+76\times 5+4\times 3+21.$
\item [] $431=98+76+5+4\times 3\times 21.$
\item [] $432=9+87+6+5+4+321.$
\item [] $433=9+8+76\times 5+4+32\times 1.$
\item [] $434=98+7+6\times 54+3+2\times 1.$
\item [] $435=98+7+6\times 54+3+2+1.$
\item [] $436=98+7+6\times 54+3\times 2+1.$
\item [] $437=98+7\times 6\times 5+4\times 32+1.$
\item [] $438=9\times 8+7\times 6+54\times 3\times 2\times 1.$
\item [] $439=9\times 8+7\times 6+54\times 3\times 2+1.$
\item [] $440=(9\times 8+7+6)\times 5+4\times 3+2+1.$
\item [] $441=98+7+6+5+4+321.$
\item [] $442=9\times 8+7+6\times 5\times 4\times 3+2+1.$
\item [] $443=9+8+76\times 5+43+2+1.$
\item [] $444=9+87+6\times 54+3+21.$
\item [] $445=98+76+54\times (3+2)+1.$
\item [] $446=9+8\times 7+6+54+321.$
\item [] $447=(9\times 8+7+65+4)\times 3+2+1.$
\item [] $448=(9+8)\times 7+6\times 54+3+2\times 1.$
\item [] $449=(9+8)\times 7+6\times 54+3+2+1.$
\item [] $450=98+7+6\times (54+3)+2+1.$
\item [] $451=9+87+6\times 5+4+321.$
\item [] $452=98+7+6+5\times 4+321.$
\item [] $453=98+7+6\times 54+3+21.$
\item [] $454=9+(8+76)\times 5+4\times 3\times 2+1.$
\item [] $455=9+8\times 7+65+4+321.$
\item [] $456=9+87+6\times (54+3+2+1).$
\item [] $457=9+(8+76)\times 5+4+3+21.$
\item [] $458=9+87+6\times 5\times 4\times 3+2\times 1.$
\item [] $459=9+87+6\times 5\times 4\times 3+2+1.$
\item [] $460=9\times 8+7+6+54+321.$
\item [] $461=9\times 8+76\times 5+4+3+2\times 1.$
\item [] $462=9\times 8+76\times 5+4+3+2+1.$
\item [] $463=9\times 8+76\times 5+4+3\times 2+1.$
\item [] $464=98+7\times 6+54\times 3\times 2\times 1.$
\item [] $465=98+7\times 6+54\times 3\times 2+1.$
\item [] $466=9\times 8+76\times 5+4\times 3+2\times 1.$
\item [] $467=9\times 8+76\times 5+4\times 3+2+1.$
\item [] $468=9+8+76+54+321.$
\item [] $469=9\times 8+7+65+4+321.$
\item [] $470=98+7\times 6+5+4+321.$
\item [] $471=9\times 8+76+5\times 4^3+2+1.$
\item [] $472=9\times 8+76+54\times 3\times 2\times 1.$
\item [] $473=9\times 8+76+54\times 3\times 2+1.$
\item [] $474=9\times 8+7\times (54+3)+2+1.$
\item [] $475=(9+8)\times 7+6\times 54+32\times 1.$
\item [] $476=9\times 8+76\times 5+4\times 3\times 2\times 1.$
\item [] $477=9+87+6+54+321.$
\item [] $478=9\times 8+76+5+4+321.$
\item [] $479=9+8\times 7\times 6+5+4\times 32+1.$
\item [] $480=9\times 8+76\times 5+4+3+21.$
\item[]$\mbox{Increasing order}$
\item [] $481=1^2+3\times 45+6\times 7\times 8+9.$
\item [] $482=1\times 2+3\times 45+6\times 7\times 8+9.$
\item [] $483=12+3+4+56\times 7+8\times 9.$
\item [] $484=12\times 34+5+6+7\times 8+9.$
\item [] $485=1\times 2+3+456+7+8+9.$
\item [] $486=1+2+3+456+7+8+9.$
\item [] $487=1+2+345+67+8\times 9.$
\item [] $488=12\times 34+56+7+8+9.$
\item [] $489=123+45\times 6+7+89.$
\item [] $490=1\times 2+3+4+56\times 7+89.$
\item [] $491=1+2+3+4+56\times 7+89.$
\item [] $492=1+23+4+56\times 7+8\times 9.$
\item [] $493=1+2\times 34+5\times 67+89.$
\item [] $494=12\times 3\times 4+5+6\times 7\times 8+9.$
\item [] $495=12+3+456+7+8+9.$
\item [] $496=12+345+67+8\times 9.$
\item [] $497=12\times 34+5+67+8+9.$
\item [] $498=12\times 34+5+6+7+8\times 9.$
\item [] $499=1\times 23\times 4+5\times 67+8\times 9.$
\item [] $500=12+3+4+56\times 7+89.$
\item [] $501=12+3+4+5+6\times 78+9.$
\item [] $502=1+23\times 4+56\times 7+8+9.$
\item [] $503=1\times 2+345+67+89.$
\item [] $504=1+2+345+67+89.$
\item [] $505=12+3\times 4+56\times 7+89.$
\item [] $506=12\times 34+5+6+78+9.$
\item [] $507=1+2\times 3\times 4+5+6\times 78+9.$
\item [] $508=1\times 23+4+56\times 7+89.$
\item [] $509=1+23+4+56\times 7+89.$
\item [] $510=12+34+56\times 7+8\times 9.$
\item [] $511=1^{2345}+6+7\times 8\times 9.$
\item [] $512=12+3+4\times 5+6\times 78+9.$
\item [] $513=12+345+67+89.$
\item [] $514=1+2^3\times 4+56\times 7+89.$
\item [] $515=12\times 34+5+6+7+89.$
\item [] $516=12\times 3+456+7+8+9.$
\item [] $517=1\times 2+34+56\times 7+89.$
\item [] $518=1+2+34+56\times 7+89.$
\item [] $519=1+2+34+5+6\times 78+9.$
\item [] $520=1+23\times 4\times 5+6\times 7+8+9.$
\item [] $521=12\times 3+4+56\times 7+89.$
\item [] $522=12\times 3+4+5+6\times 78+9.$
\item [] $523=1^2+3+4+5+6+7\times 8\times 9.$
\item [] $524=1\times 2+3+4+5+6+7\times 8\times 9.$
\item [] $525=1+2+3+4+5+6+7\times 8\times 9.$
\item [] $526=1\times 2+3+456+7\times 8+9.$
\item [] $527=123\times 4+5+6+7+8+9.$
\item [] $528=12+34+5+6\times 78+9.$
\item [] $529=12\times 34+56+7\times 8+9.$
\item [] $530=1+2+3\times 4+5+6+7\times 8\times 9.$
\item [] $531=1\times 23\times 4\times 5+6+7\times 8+9.$
\item [] $532=1\times 2\times 34+56\times 7+8\times 9.$
\item [] $533=12\times 3+4\times 5+6\times 78+9.$
\item [] $534=1+234+5\times 6\times 7+89.$
\item [] $535=1\times 2+3+4\times 5+6+7\times 8\times 9.$
\item [] $536=12+3+456+7\times 8+9.$
\item [] $537=12+3+45+6\times 78+9.$
\item [] $538=1^2+3\times 4\times 5+6\times 78+9.$
\item [] $539=12+3\times 4+5+6+7\times 8\times 9.$
\item [] $540=1\times 2+3+456+7+8\times 9.$
\item [] $541=1+2+3+456+7+8\times 9.$
\item [] $542=1+2\times 3+456+7+8\times 9.$
\item [] $543=12\times 34+56+7+8\times 9.$
\item [] $544=12\times 34+5+6\times 7+89.$
\item [] $545=1+23+456+7\times 8+9.$
\item [] $546=1+23+45+6\times 78+9.$
\item [] $547=1^2+3+456+78+9.$
\item [] $548=1\times 2+3+456+78+9.$
\item [] $549=1+2+3+456+78+9.$
\item [] $550=1+2\times 3+456+78+9.$
\item[]$\mbox{Decreasing order}$
\item [] $481=9+8+7\times 65+4+3+2\times 1.$
\item [] $482=9+8+7\times 65+4+3+2+1.$
\item [] $483=9+8+7\times 65+4+3\times 2+1.$
\item [] $484=9+8+76\times 5+43\times 2+1.$
\item [] $485=9\times 8+76\times 5+4\times 3+21.$
\item [] $486=9+87+65+4+321.$
\item [] $487=9+8+7\times 65+4\times 3+2+1.$
\item [] $488=98+76\times 5+4+3+2+1.$
\item [] $489=98+76\times 5+4+3\times 2+1.$
\item [] $490=(9+8+76)\times 5+4\times 3\times 2+1.$
\item [] $491=98+76\times 5+4+3^2\times 1.$
\item [] $492=98+76\times 5+4\times 3+2\times 1.$
\item [] $493=98+76\times 5+4\times 3+2+1.$
\item [] $494=98+76+5\times (43+21).$
\item [] $495=98+7+65+4+321.$
\item [] $496=9+8+7\times 65+4\times 3\times 2\times 1.$
\item [] $497=9+8+7\times 65+4\times 3\times 2+1.$
\item [] $498=9\times 8+76\times 5+43+2+1.$
\item [] $499=98+76+54\times 3\times 2+1.$
\item [] $500=9+8+7\times 65+4+3+21.$
\item [] $501=9+8\times 7+6+5\times 43\times 2\times 1.$
\item [] $502=98+76\times 5+4\times 3\times 2\times 1.$
\item [] $503=98+76\times 5+4\times 3\times 2+1.$
\item [] $504=98+76+5+4+321.$
\item [] $505=9+8+7\times 65+4\times 3+21.$
\item [] $506=98+76\times 5+4+3+21.$
\item [] $507=(9\times 8+7)\times 6+5+4+3+21.$
\item [] $508=9+8\times 7+6+5+432\times 1.$
\item [] $509=9+8\times 7+6+5+432+1.$
\item [] $510=9+8\times 7\times 6+54\times 3+2+1.$
\item [] $511=98+76\times 5+4\times 3+21.$
\item [] $512=98+7\times (54+3+2)+1.$
\item [] $513=9\times 8+7\times 6\times (5+4)+3\times 21.$
\item [] $514=98+76\times 5+4+32\times 1.$
\item [] $515=98+7\times 6+54+321.$
\item [] $516=9\times 8+76\times 5+43+21.$
\item [] $517=9+8+7\times 65+43+2\times 1.$
\item [] $518=9+8+7\times 65+43+2+1.$
\item [] $519=9\times 8+76\times 5+4+3\times 21.$
\item [] $520=9\times 8+7+6\times 5\times 4+321.$
\item [] $521=9+8+7+65+432\times 1.$
\item [] $522=9+8+7+65+432+1.$
\item [] $523=9\times 8+76+54+321.$
\item [] $524=98+76\times 5+43+2+1.$
\item [] $525=9\times 8\times 7+6+5+4+3+2+1.$
\item [] $526=9\times 8\times 7+6+5+4+3\times 2+1.$
\item [] $527=9+8\times 7+6\times 5+432\times 1.$
\item [] $528=9+8\times 7+6\times 5+432+1.$
\item [] $529=9\times 8\times 7+6+5+4\times 3+2\times 1.$
\item [] $530=9\times 8\times 7+6+5+4\times 3+2+1.$
\item [] $531=9+8+76+5+432+1.$
\item [] $532=9+87+6+5\times 43\times 2\times 1.$
\item [] $533=(9\times 8+7)\times 6+54+3+2\times 1.$
\item [] $534=9\times 8+7\times 6\times 5+4\times 3\times 21.$
\item [] $535=9\times 8\times 7+6+5\times 4+3+2\times 1.$
\item [] $536=9\times 8\times 7+6+5\times 4+3+2+1.$
\item [] $537=9\times 8+7\times 65+4+3+2+1.$
\item [] $538=9\times 8+7\times 65+4+3\times 2+1.$
\item [] $539=9\times 8\times 7+6+5+4\times 3\times 2\times 1.$
\item [] $540=9+87+6+5+432+1.$
\item [] $541=9\times 8+7\times 65+4\times 3+2\times 1.$
\item [] $542=9\times 8+7\times 65+4\times 3+2+1.$
\item [] $543=9\times 8\times 7+6+5+4+3+21.$
\item [] $544=9\times 8\times 7+6\times 5+4+3+2+1.$
\item [] $545=9+87\times 6+5+4+3+2\times 1.$
\item [] $546=9+87\times 6+5+4+3+2+1.$
\item [] $547=9+87\times 6+5+4+3\times 2+1.$
\item [] $548=9\times 8\times 7+6\times 5+4\times 3+2\times 1.$
\item [] $549=98+76+54+321.$
\item [] $550=9+87\times 6+5+4\times 3+2\times 1.$
\item[]$\mbox{Increasing order}$
\item [] $551=12\times 34+56+78+9.$
\item [] $552=12\times 34+5+67+8\times 9.$
\item [] $553=12\times 3\times 4+56\times 7+8+9.$
\item [] $554=1+23\times 4\times 5+6+78+9.$
\item [] $555=12\times 3+4+5+6+7\times 8\times 9.$
\item [] $556=123\times 4+5+6\times 7+8+9.$
\item [] $557=1\times 2+3+456+7+89.$
\item [] $558=12+3+456+78+9.$
\item [] $559=1+23+456+7+8\times 9.$
\item [] $560=12\times 34+56+7+89.$
\item [] $561=1+2+3+45+6+7\times 8\times 9.$
\item [] $562=1+23+4+5\times 6+7\times 8\times 9.$
\item [] $563=1+23\times 4\times 5+6+7+89.$
\item [] $564=1\times 2+3\times 4+5+67\times 8+9.$
\item [] $565=1+2+3\times 4+5+67\times 8+9.$
\item [] $566=1\times 23+456+78+9.$
\item [] $567=1+23+456+78+9.$
\item [] $568=12\times 3\times 4+5\times 67+89.$
\item [] $569=12\times 34+5+67+89.$
\item [] $570=12+3+45+6+7\times 8\times 9.$
\item [] $571=12\times 3+456+7+8\times 9.$
\item [] $572=123\times 4+56+7+8+9.$
\item [] $573=1\times 23\times 4+56\times 7+89.$
\item [] $574=12+3\times 4+5+67\times 8+9.$
\item [] $575=1\times 23+456+7+89.$
\item [] $576=1+23+456+7+89.$
\item [] $577=1\times 23+4+5+67\times 8+9.$
\item [] $578=1+23+4+5+67\times 8+9.$
\item [] $579=12\times 3+456+78+9.$
\item [] $580=12+34+5\times 6+7\times 8\times 9.$
\item [] $581=123\times 4+5+67+8+9.$
\item [] $582=123\times 4+5+6+7+8\times 9.$
\item [] $583=1\times 2\times 34+5+6+7\times 8\times 9.$
\item [] $584=12+3\times 4+56+7\times 8\times 9.$
\item [] $585=1+234+5+6\times 7\times 8+9.$
\item [] $586=1\times 234+5\times 67+8+9.$
\item [] $587=1+234+5\times 67+8+9.$
\item [] $588=12\times 3+456+7+89.$
\item [] $589=1+23+4\times 5+67\times 8+9.$
\item [] $590=123\times 4+5+6+78+9.$
\item [] $591=123+4+56\times 7+8\times 9.$
\item [] $592=1+23\times 4\times 5+6\times 7+89.$
\item [] $593=1\times 2+3+4+567+8+9.$
\item [] $594=1+2+3+4+567+8+9.$
\item [] $595=1+2\times 3+4+567+8+9.$
\item [] $596=12+34+5+67\times 8+9.$
\item [] $597=1+2+34+56+7\times 8\times 9.$
\item [] $598=1\times 2+3\times 4+567+8+9.$
\item [] $599=123\times 4+5+6+7+89.$
\item [] $600=12\times 3+4+56+7\times 8\times 9.$
\item [] $601=12\times 3+4\times 5+67\times 8+9.$
\item [] $602=1\times 2\times 34+5\times 6+7\times 8\times 9.$
\item [] $603=123+456+7+8+9.$
\item [] $604=(1\times 2+3)\times 4+567+8+9.$
\item [] $605=12+3+45+67\times 8+9.$
\item [] $606=12+34+56+7\times 8\times 9.$
\item [] $607=1\times 23\times 4+5+6+7\times 8\times 9.$
\item [] $608=123+4+56\times 7+89.$
\item [] $609=123+4+5+6\times 78+9.$
\item [] $610=1+2\times 3+45+(6+7\times 8)\times 9.$
\item [] $611=1\times 23+4+567+8+9.$
\item [] $612=1+23+4+567+8+9.$
\item [] $613=123\times 4+56+7\times 8+9.$
\item [] $614=1+23+45+67\times 8+9.$
\item [] $615=1+2+3\times 45+6\times 78+9.$
\item [] $616=1\times 23\times 4\times 5+67+89.$
\item [] $617=1+23\times 4\times 5+67+89.$
\item [] $618=123\times 4+5\times 6+7+89.$
\item [] $619=1+2\times 34+5+67\times 8+9.$
\item [] $620=1\times 2+34+567+8+9.$
\item[]$\mbox{Decreasing order}$
\item [] $551=9+87\times 6+5+4\times 3+2+1.$
\item [] $552=9\times 8+7\times 65+4\times 3\times 2+1.$
\item [] $553=98+7\times (6+54+3+2\times 1).$
\item [] $554=9\times 8\times 7+6+5\times 4+3+21.$
\item [] $555=9\times 8+7\times 65+4+3+21.$
\item [] $556=9+87\times 6+5\times 4+3+2\times 1.$
\item [] $557=9+87\times 6+5\times 4+3+2+1.$
\item [] $558=9+87\times 6+5\times 4+3\times 2+1.$
\item [] $559=9\times 8\times 7+6\times 5+4\times 3\times 2+1.$
\item [] $560=9\times 8+7\times 65+4\times 3+21.$
\item [] $561=9+87\times 6+5+4\times 3\times 2+1.$
\item [] $562=98+7\times 65+4+3+2\times 1.$
\item [] $563=98+7\times 65+4+3+2+1.$
\item [] $564=98+7\times 65+4+3\times 2+1.$
\item [] $565=98+76\times 5+43\times 2+1.$
\item [] $566=98+7\times 65+4+3^2\times 1.$
\item [] $567=98+7\times 65+4\times 3+2\times 1.$
\item [] $568=98+7\times 65+4\times 3+2+1.$
\item [] $569=9\times 8\times 7+6+54+3+2\times 1.$
\item [] $570=9\times 8\times 7+6+54+3+2+1.$
\item [] $571=9\times 8\times 7+6+54+3\times 2+1.$
\item [] $572=9\times 8\times 7+6+5\times 4\times 3+2\times 1.$
\item [] $573=9\times 8+7\times 65+43+2+1.$
\item [] $574=9\times 8\times 7+6+54+3^2+1.$
\item [] $575=9+87\times 6+5\times 4+3+21.$
\item [] $576=9+8+7+6+543+2+1.$
\item [] $577=98+7\times 65+4\times 3\times 2\times 1.$
\item [] $578=98+7\times 65+4\times 3\times 2+1.$
\item [] $579=9\times 8\times 7+65+4+3+2+1.$
\item [] $580=9\times 8\times 7+65+4+3\times 2+1.$
\item [] $581=98+7\times 65+4+3+21.$
\item [] $582=9+87\times 6+5+43+2+1.$
\item [] $583=9\times 8\times 7+65+4\times 3+2\times 1.$
\item [] $584=9\times 8\times 7+65+4\times 3+2+1.$
\item [] $585=9\times 8+76+5+432\times 1.$
\item [] $586=9\times 8+76+5+432+1.$
\item [] $587=98+7\times (65+4)+3+2+1.$
\item [] $588=9\times 8\times 7+6+54+3+21.$
\item [] $589=98+7\times 65+4+32\times 1.$
\item [] $590=9+87\times 6+54+3+2\times 1.$
\item [] $591=9+87\times 6+54+3+2+1.$
\item [] $592=9+87\times 6+54+3\times 2+1.$
\item [] $593=9\times 8\times 7+65+4\times 3\times 2\times 1.$
\item [] $594=9+87+65+432+1.$
\item [] $595=9+87\times 6+54+3^2+1.$
\item [] $596=9\times 8\times 7+6+54+32\times 1.$
\item [] $597=9\times 8\times 7+65+4+3+21.$
\item [] $598=98+7\times 65+43+2\times 1.$
\item [] $599=98+7\times 65+43+2+1.$
\item [] $600=9+87\times 6+5+43+21.$
\item [] $601=9+8+7\times 65+4\times 32+1.$
\item [] $602=98+7+65+432\times 1.$
\item [] $603=98+7+65+432+1.$
\item [] $604=98+76+5\times 43\times 2\times 1.$
\item [] $605=9+8+7\times 6+543+2+1.$
\item [] $606=98+76\times 5+4\times 32\times 1.$
\item [] $607=98+76\times 5+4\times 32+1.$
\item [] $608=98+7\times (6+5)+432+1.$
\item [] $609=9+87\times 6+54+3+21.$
\item [] $610=98+76\times 5+4\times (32+1).$
\item [] $611=98+76+5+432\times 1.$
\item [] $612=98+76+5+432+1.$
\item [] $613=9\times 8+7\times 65+43\times 2\times 1.$
\item [] $614=9\times 8\times 7+65+43+2\times 1.$
\item [] $615=9\times 8\times 7+65+43+2+1.$
\item [] $616=9+8\times 7+6+543+2\times 1.$
\item [] $617=9+8\times 7+6+543+2+1.$
\item [] $618=9+87\times 6+54+32+1.$
\item [] $619=9\times 8\times 7+6\times 5+4^3+21.$
\item [] $620=98+7\times 65+4+3\times 21.$
\item[]$\mbox{Increasing order}$
\item [] $621=1+2+34+567+8+9.$
\item [] $622=123\times 4+5+6+7\times (8+9).$
\item [] $623=12+3+4\times (56+7+89).$
\item [] $624=12\times 3+4+567+8+9.$
\item [] $625=12\times 3\times 4+56\times 7+89.$
\item [] $626=12\times 3+45+67\times 8+9.$
\item [] $627=123\times 4+56+7+8\times 9.$
\item [] $628=123\times 4+5+6\times 7+89.$
\item [] $629=1+2\times 34+56+7\times 8\times 9.$
\item [] $630=12+34+567+8+9.$
\item [] $631=1+2\times 3\times 4\times 5+6+7\times 8\times 9.$
\item [] $632=1\times 2^3+4\times 5\times 6+7\times 8\times 9.$
\item [] $633=1+2^3+4\times 5\times 6+7\times 8\times 9.$
\item [] $634=1+2+3^4+5+67\times 8+9.$
\item [] $635=123\times 4+56+78+9.$
\item [] $636=123\times 4+5+67+8\times 9.$
\item [] $637=1\times 2^3\times 4\times 5+6\times 78+9.$
\item [] $638=1+2^3\times 4\times 5+6\times 78+9.$
\item [] $639=12+3+4\times 5\times 6+7\times 8\times 9.$
\item [] $640=1^{234}+567+8\times 9.$
\item [] $641=1\times 234+5\times 67+8\times 9.$
\item [] $642=123+4+5+6+7\times 8\times 9.$
\item [] $643=1\times 234+56\times 7+8+9.$
\item [] $644=123+456+7\times 8+9.$
\item [] $645=123+45+6\times 78+9.$
\item [] $646=1^2\times 3+4+567+8\times 9.$
\item [] $647=1\times 2+3\times 45+6+7\times 8\times 9.$
\item [] $648=1+2+3\times 45+6+7\times 8\times 9.$
\item [] $649=1+2+3+4+567+8\times 9.$
\item [] $650=1+2\times 3+4+567+8\times 9.$
\item [] $651=1\times 2^3+4+567+8\times 9.$
\item [] $652=1\times 2\times 34+567+8+9.$
\item [] $653=123\times 4+5+67+89.$
\item [] $654=1+2+3\times 4+567+8\times 9.$
\item [] $655=1+2+3+4\times 5+6+7\times 89.$
\item [] $656=1+2\times 3+4\times 5+6+7\times 89.$
\item [] $657=12+3\times 45+6+7\times 8\times 9.$
\item [] $658=123+456+7+8\times 9.$
\item [] $659=1+234+5\times 67+89.$
\item [] $660=12\times 3+4\times 5\times 6+7\times 8\times 9.$
\item [] $661=123+4+5\times 6+7\times 8\times 9.$
\item [] $662=1+23+4+5+6+7\times 89.$
\item [] $663=12+3\times 4+567+8\times 9.$
\item [] $664=12+3+4\times 5+6+7\times 89.$
\item [] $665=1\times 2+3+4+567+89.$
\item [] $666=123+456+78+9.$
\item [] $667=1+23+4+567+8\times 9.$
\item [] $668=1+2+3\times 4+5\times 6+7\times 89.$
\item [] $669=1+2^3+4+567+89.$
\item [] $670=1\times 2+3\times 4+567+89.$
\item [] $671=1+2+3\times 4+567+89.$
\item [] $672=1\times 23+4\times 5+6+7\times 89.$
\item [] $673=1+23+4\times 5+6+7\times 89.$
\item [] $674=12\times 3+4+5+6+7\times 89.$
\item [] $675=123+456+7+89.$
\item [] $676=1+2+34+567+8\times 9.$
\item [] $677=1+23\times 4+567+8+9.$
\item [] $678=123+45+6+7\times 8\times 9.$
\item [] $679=12\times 3+4+567+8\times 9.$
\item [] $680=12+3\times 4+567+89.$
\item [] $681=1+23+4+5\times 6+7\times 89.$
\item [] $682=1\times 2+3\times 45+67\times 8+9.$
\item [] $683=1\times 23+4+567+89.$
\item [] $684=1+23+4+567+89.$
\item [] $685=12+34+567+8\times 9.$
\item [] $686=1+2\times 34\times 5+6\times 7\times 8+9.$
\item [] $687=123+4+56+7\times 8\times 9.$
\item [] $688=123+4\times 5+67\times 8+9.$
\item [] $689=12+3+45+6+7\times 89.$
\item [] $690=12\times 34+5\times 6\times 7+8\times 9.$
\item[]$\mbox{Decreasing order}$
\item [] $621=9+87\times 6+5+4^3+21.$
\item [] $622=9+87\times 6+5+43\times 2\times 1.$
\item [] $623=9+8+7\times 6+543+21.$
\item [] $624=(9+8+7+6)\times 5\times 4+3+21.$
\item [] $625=98+76\times 5+(4+3)\times 21.$
\item [] $626=98+7\times 6+54\times 3^2\times 1.$
\item [] $627=9\times 8\times 7+6+54+3\times 21.$
\item [] $628=(9+87)\times 6+5\times 4+32\times 1.$
\item [] $629=9\times 8\times 7+6\times 5\times 4+3+2\times 1.$
\item [] $630=9\times 8\times 7+6+5\times 4\times 3\times 2\times 1.$
\item [] $631=9\times 8+7+6+543+2+1.$
\item [] $632=9+8\times 76+5+4+3+2+1.$
\item [] $633=9+8\times 76+5+4+3\times 2+1.$
\item [] $634=9\times 8+76+54\times 3^2\times 1.$
\item [] $635=9+8\times 7+6+543+21.$
\item [] $636=9+8\times 76+5+4\times 3+2\times 1.$
\item [] $637=9+8\times 76+5+4\times 3+2+1.$
\item [] $638=9+8+76+543+2\times 1.$
\item [] $639=9+8+76+543+2+1.$
\item [] $640=98+7\times (65+4\times 3)+2+1.$
\item [] $641=(9+87)\times 6+5\times (4+3\times (2+1)).$
\item [] $642=9+8\times 76+5\times 4+3+2\times 1.$
\item [] $643=9+8\times 76+5\times 4+3+2+1.$
\item [] $644=9+8\times 76+5\times 4+3\times 2+1.$
\item [] $645=(98+7)\times 6+5+4+3+2+1.$
\item [] $646=9+8\times 76+5+4\times 3\times 2\times 1.$
\item [] $647=9+87+6+543+2\times 1.$
\item [] $648=9+87+6+543+2+1.$
\item [] $649=9\times 8+7+6+543+21.$
\item [] $650=9+8\times 76+5+4+3+21.$
\item [] $651=9+87\times 6+5\times 4\times 3\times 2\times 1.$
\item [] $652=9+87\times 6+5\times 4\times 3\times 2+1.$
\item [] $653=9+8\times 7\times (6+5)+4+3+21.$
\item [] $654=9\times 8\times 7+65+4^3+21.$
\item [] $655=9+8\times 76+5+4\times 3+21.$
\item [] $656=98+7+6+543+2\times 1.$
\item [] $657=98+7+6+543+2+1.$
\item [] $658=9+8\times 76+5+4+32\times 1.$
\item [] $659=9+8\times 76+5+4+32+1.$
\item [] $660=9\times 8+7\times 6+543+2+1.$
\item [] $661=9+8\times 76+5\times 4+3+21.$
\item [] $662=9\times 8\times 7+6\times 5+4\times 32\times 1.$
\item [] $663=9\times 8\times 7+6\times 5+4\times 32+1.$
\item [] $664=9+87\times 6+5+4\times 32\times 1.$
\item [] $665=9+87\times 6+5+4\times 32+1.$
\item [] $666=9+87+6+543+21.$
\item [] $667=9+8\times 76+5+43+2\times 1.$
\item [] $668=9+8\times 76+5+43+2+1.$
\item [] $669=9+8\times 7\times 6+54\times 3\times 2\times 1.$
\item [] $670=9+8\times 7\times 6+54\times 3\times 2+1.$
\item [] $671=9+8+7+6+5\times 4\times 32+1.$
\item [] $672=9+8\times (76+5)+4\times 3+2+1.$
\item [] $673=9+8\times (7+65+4+3\times 2+1).$
\item [] $674=9\times 8\times 7+6+54\times 3+2\times 1.$
\item [] $675=98+7+6+543+21.$
\item [] $676=9+8\times 76+54+3+2\times 1.$
\item [] $677=9+8\times 76+54+3+2+1.$
\item [] $678=9+8\times 76+54+3\times 2+1.$
\item [] $679=9+8\times 76+5\times 4\times 3+2\times 1.$
\item [] $680=9+8\times 76+5\times 4\times 3+2+1.$
\item [] $681=98+7\times 65+4\times 32\times 1.$
\item [] $682=98+7\times 65+4\times 32+1.$
\item [] $683=9+8+7+654+3+2\times 1.$
\item [] $684=9+8+7+654+3+2+1.$
\item [] $685=9+8+7+654+3\times 2+1.$
\item [] $686=98+7\times 6+543+2+1.$
\item [] $687=9\times 8\times 7+6\times 5\times 4+3\times 21.$
\item [] $688=9+8+7+654+3^2+1.$
\item [] $689=9+8\times 76+5+4+3\times 21.$
\item [] $690=9+8+7\times 6+5^4+3\times 2\times 1.$
\item[]$\mbox{Increasing order}$
\item [] $691=1\times 2+3\times 4\times 5+6+7\times 89.$
\item [] $692=1\times 2+34+567+89.$
\item [] $693=1+2+34+567+89.$
\item [] $694=1+2+3\times 4+56+7\times 89.$
\item [] $695=123+4\times (56+78+9).$
\item [] $696=12\times 3+4+567+89.$
\item [] $697=1\times 23+45+6+7\times 89.$
\item [] $698=12+3+4+56+7\times 89.$
\item [] $699=12+34+5\times 6+7\times 89.$
\item [] $700=1^2+3+4+5+678+9.$
\item [] $701=1\times 2+3+4+5+678+9.$
\item [] $702=12+34+567+89.$
\item [] $703=12+3\times 4+56+7\times 89.$
\item [] $704=12\times 3\times 4+56+7\times 8\times 9.$
\item [] $705=1+2^3+4+5+678+9.$
\item [] $706=1\times 23+4+56+7\times 89.$
\item [] $707=1+23+4+56+7\times 89.$
\item [] $708=1+2\times 34+567+8\times 9.$
\item [] $709=1^{2345}+6+78\times 9.$
\item [] $710=12\times 3+45+6+7\times 89.$
\item [] $711=123+4+567+8+9.$
\item [] $712=1\times 2+3+4\times 5+678+9.$
\item [] $713=123+45+67\times 8+9.$
\item [] $714=1+2\times 3+4\times 5+678+9.$
\item [] $715=1\times 234+56\times 7+89.$
\item [] $716=1+234+56\times 7+89.$
\item [] $717=1+234+5+6\times 78+9.$
\item [] $718=1+2+34\times 5+67\times 8+9.$
\item [] $719=123\times 4+5\times 6\times 7+8+9.$
\item [] $720=1+23+4+5+678+9.$
\item [] $721=1+2\times 345+6+7+8+9.$
\item [] $722=12+3+4\times 5+678+9.$
\item [] $723=1+2+3+4+5+6+78\times 9.$
\item [] $724=1\times 2\times 34+567+89.$
\item [] $725=12+34+56+7\times 89.$
\item [] $726=1\times 23\times 4+5+6+7\times 89.$
\item [] $727=12+34\times 5+67\times 8+9.$
\item [] $728=1\times 2+34+5+678+9.$
\item [] $729=1+2+34+5+678+9.$
\item [] $730=1\times 23+4\times 5+678+9.$
\item [] $731=1+23+4\times 5+678+9.$
\item [] $732=12\times 3+4+5+678+9.$
\item [] $733=1\times 2+3+4\times 5+6+78\times 9.$
\item [] $734=1+2+3+4\times 5+6+78\times 9.$
\item [] $735=1+2\times 3+4\times 5+6+78\times 9.$
\item [] $736=1\times 2+3^4+5\times 6+7\times 89.$
\item [] $737=1\times 2+3+45+678+9.$
\item [] $738=12+34+5+678+9.$
\item [] $739=1+2\times 3+45+678+9.$
\item [] $740=1\times 23+4+5+6+78\times 9.$
\item [] $741=1+23+4+5+6+78\times 9.$
\item [] $742=1+2+3+4+5\times 6+78\times 9.$
\item [] $743=12\times 3+4\times 5+678+9.$
\item [] $744=1\times 2\times 3\times 4\times 5\times 6+7+8+9.$
\item [] $745=1\times 23\times 4+5\times 6+7\times 89.$
\item [] $746=1+23\times 4+5\times 6+7\times 89.$
\item [] $747=12+3+45+678+9.$
\item [] $748=1\times 23\times 4+567+89.$
\item [] $749=1+23\times 4+567+89.$
\item [] $750=1+2+34+5+6+78\times 9.$
\item [] $751=12+3+4+5\times 6+78\times 9.$
\item [] $752=1+23+4\times 5+6+78\times 9.$
\item [] $753=12\times 3+4+5+6+78\times 9.$
\item [] $754=1+2+3\times 4\times 56+7+8\times 9.$
\item [] $755=1\times 23+45+678+9.$
\item [] $756=1+23+45+678+9.$
\item [] $757=1+2\times 3\times 4+5\times 6+78\times 9.$
\item [] $758=1\times 2+3+45+6+78\times 9.$
\item [] $759=12+34+5+6+78\times 9.$
\item [] $760=12\times 34+5\times 67+8+9.$
\item[]$\mbox{Decreasing order}$
\item [] $691=9\times 8\times 7+6+5\times 4\times 3^2+1.$
\item [] $692=(98+7)\times 6+5\times 4\times 3+2\times 1.$
\item [] $693=9\times 8+76+543+2\times 1.$
\item [] $694=9\times 8+76+543+2+1.$
\item [] $695=9+87\times 6+54\times 3+2\times 1.$
\item [] $696=9+87\times 6+54\times 3+2+1.$
\item [] $697=9\times 8\times 7+65+4\times 32\times 1.$
\item [] $698=9\times 8\times 7+65+4\times 32+1.$
\item [] $699=9+8+7\times 6+5\times 4\times 32\times 1.$
\item [] $700=9+8\times 76+5\times 4+3\times 21.$
\item [] $701=9+8\times 7+6+5^4+3+2\times 1.$
\item [] $702=9+8+7+654+3+21.$
\item [] $703=9+8\times 7+6+5^4+3\times 2+1.$
\item [] $704=98+7\times 6+543+21.$
\item [] $705=9+8\times 7+6+5^4+3^2\times 1.$
\item [] $706=98\times 7+6+5+4+3+2\times 1.$
\item [] $707=98\times 7+6+5+4+3+2+1.$
\item [] $708=98\times 7+6+5+4+3\times 2+1.$
\item [] $709=(98+7\times 6)\times 5+4+3+2\times 1.$
\item [] $710=9+8\times 7+6\times 54+321.$
\item [] $711=9+8+7+654+32+1.$
\item [] $712=9\times 8+76+543+21.$
\item [] $713=(98+7+65)\times 4+32+1.$
\item [] $714=9+87\times 6+54\times 3+21.$
\item [] $715=9\times 8+7\times 6\times 5+432+1.$
\item [] $716=9+8+7\times 6+5^4+32\times 1.$
\item [] $717=98\times 7+6+5\times 4+3+2\times 1.$
\item [] $718=98\times 7+6+5\times 4+3+2+1.$
\item [] $719=98\times 7+6+5\times 4+3\times 2+1.$
\item [] $720=98+76+543+2+1.$
\item [] $721=98\times 7+6+5+4\times 3\times 2\times 1.$
\item [] $722=98\times 7+6+5+4\times 3\times 2+1.$
\item [] $723=9+8+76+5^4+3+2\times 1.$
\item [] $724=9+8\times 7+654+3+2\times 1.$
\item [] $725=9+8\times 7+654+3+2+1.$
\item [] $726=98\times 7+6\times 5+4+3+2+1.$
\item [] $727=98\times 7+6\times 5+4+3\times 2+1.$
\item [] $728=9\times 8\times 7+6+5\times 43+2+1.$
\item [] $729=9+8\times 7+654+3^2+1.$
\item [] $730=98\times 7+6\times 5+4\times 3+2\times 1.$
\item [] $731=98\times 7+6\times 5+4\times 3+2+1.$
\item [] $732=9+87+6+5^4+3+2\times 1.$
\item [] $733=98\times 7+6+5+4+32\times 1.$
\item [] $734=98\times 7+6+5+4+32+1.$
\item [] $735=9\times (8+7+6)+543+2+1.$
\item [] $736=98\times 7+6+5\times 4+3+21.$
\item [] $737=9+8\times 76+5\times 4\times 3\times 2\times 1.$
\item [] $738=98+76+543+21.$
\item [] $739=9\times 8+7+654+3+2+1.$
\item [] $740=98\times 7+6\times 5+4\times 3\times 2\times 1.$
\item [] $741=98\times 7+6\times 5+4\times 3\times 2+1.$
\item [] $742=98\times 7+6+5+43+2\times 1.$
\item [] $743=98\times 7+6+5+43+2+1.$
\item [] $744=9+8+7+6\times 5\times 4\times 3\times 2\times 1.$
\item [] $745=9+8+7+6\times 5\times 4\times 3\times 2+1.$
\item [] $746=9\times 8\times 7+6+5\times 43+21.$
\item [] $747=9+8+7\times 6+5^4+3\times 21.$
\item [] $748=9+87\times 6+5\times 43+2\times 1.$
\item [] $749=9+87\times 6+5\times 43+2+1.$
\item [] $750=98+7+6\times 54+321.$
\item [] $751=98\times 7+6+54+3+2\times 1.$
\item [] $752=98\times 7+6+54+3+2+1.$
\item [] $753=98\times 7+6+54+3\times 2+1.$
\item [] $754=98\times 7+6+5\times 4\times 3+2\times 1.$
\item [] $755=9+87+654+3+2\times 1.$
\item [] $756=9+87+654+3+2+1.$
\item [] $757=9+87+654+3\times 2+1.$
\item [] $758=(9+8)\times 7\times 6+5\times 4+3+21.$
\item [] $759=9+87+654+3\times (2+1).$
\item [] $760=98\times 7+65+4+3+2\times 1.$
\item[]$\mbox{Increasing order}$
\item [] $761=1+2\times 34+5+678+9.$
\item [] $762=1+2\times 345+6+7\times 8+9.$
\item [] $763=12+3\times 4\times 56+7+8\times 9.$
\item [] $764=12\times 3+4\times 5+6+78\times 9.$
\item [] $765=1+2^3\times 4+5\times 6+78\times 9.$
\item [] $766=123+4+567+8\times 9.$
\item [] $767=1\times 2+3+4+56+78\times 9.$
\item [] $768=12\times 3+45+678+9.$
\item [] $769=1+2+34+5\times 6+78\times 9.$
\item [] $770=1\times 2+3\times 4\times 56+7+89.$
\item [] $771=12+3\times 4\times 56+78+9.$
\item [] $772=123+4\times 5+6+7\times 89.$
\item [] $773=1+2+3\times 4+56+78\times 9.$
\item [] $774=1\times 2\times 345+67+8+9.$
\item [] $775=1+2\times 345+67+8+9.$
\item [] $776=12+3\times 45+6+7\times 89.$
\item [] $777=12+3+4+56+78\times 9.$
\item [] $778=12+34+5\times 6+78\times 9.$
\item [] $779=12\times 3+4\times 5\times 6+7\times 89.$
\item [] $780=123+4+5\times 6+7\times 89.$
\item [] $781=1\times 2\times 34+5+6+78\times 9.$
\item [] $782=12+3\times 4+56+78\times 9.$
\item [] $783=123+4+567+89.$
\item [] $784=1+2\times 345+6+78+9.$
\item [] $785=1+23\times 4+5+678+9.$
\item [] $786=1+23+4+56+78\times 9.$
\item [] $787=(1\times 2+3)\times 4\times 5+678+9.$
\item [] $788=(1+23)\times 4+5+678+9.$
\item [] $789=12\times 3+45+6+78\times 9.$
\item [] $790=1\times 2^3\times 4+56+78\times 9.$
\item [] $791=123\times 4+5\times 6\times 7+89.$
\item [] $792=1\times 2\times 345+6+7+89.$
\item [] $793=1+2\times 345+6+7+89.$
\item [] $794=1\times 2+34+56+78\times 9.$
\item [] $795=1+2+34+56+78\times 9.$
\item [] $796=1\times 2+3^4+5+6+78\times 9.$
\item [] $797=123+45+6+7\times 89.$
\item [] $798=12\times 3+4+56+78\times 9.$
\item [] $799=1\times 2\times 3\times 4\times 5\times 6+7+8\times 9.$
\item [] $800=12\times 3\times 4+567+89.$
\item [] $801=1+2\times 34+5\times 6+78\times 9.$
\item [] $802=1+2+34\times 5+6+7\times 89.$
\item [] $803=12\times 34+5+6\times (7\times 8+9).$
\item [] $804=12+34+56+78\times 9.$
\item [] $805=1+2\times 345+6\times 7+8\times 9.$
\item [] $806=123+4+56+7\times 89.$
\item [] $807=1\times 2\times 3\times 4\times 5+678+9.$
\item [] $808=1+2\times 3\times 4\times 5+678+9.$
\item [] $809=1\times 2+3+4+5+6+789.$
\item [] $810=1+2+3+4+5+6+789.$
\item [] $811=12+34\times 5+6+7\times 89.$
\item [] $812=1\times 2^3+4+5+6+789.$
\item [] $813=12\times 3\times 4\times 5+6+78+9.$
\item [] $814=1\times 2+3\times 4+5+6+789.$
\item [] $815=12\times 34+5\times 67+8\times 9.$
\item [] $816=1+2\times 3\times 45+67\times 8+9.$
\item [] $817=12\times 34+56\times 7+8+9.$
\item [] $818=1+2\times 34\times 5+6\times 78+9.$
\item [] $819=123+4+5+678+9.$
\item [] $820=1\times 2+3+4\times 5+6+789.$
\item [] $821=1+2+3+4\times 5+6+789.$
\item [] $822=1+2\times 345+6\times 7+89.$
\item [] $823=12\times 3\times 4+56+7\times 89.$
\item [] $824=12+3\times 4+5+6+789.$
\item [] $825=1+2+3\times 45+678+9.$
\item [] $826=1\times 2\times 34+56+78\times 9.$
\item [] $827=1+2\times 34+56+78\times 9.$
\item [] $828=1+23+4+5+6+789.$
\item [] $829=1+2+3+4+5\times 6+789.$
\item [] $830=123+4\times 5+678+9.$
\item[]$\mbox{Decreasing order}$
\item [] $761=98\times 7+65+4+3+2+1.$
\item [] $762=98\times 7+65+4+3\times 2+1.$
\item [] $763=98\times 7+6+5+4^3+2\times 1.$
\item [] $764=98+7+654+3+2\times 1.$
\item [] $765=98+7+654+3+2+1.$
\item [] $766=98\times 7+65+4\times 3+2+1.$
\item [] $767=9+87\times 6+5\times 43+21.$
\item [] $768=98+7+6+5^4+32\times 1.$
\item [] $769=9\times 8\times 7+65\times 4+3+2\times 1.$
\item [] $770=9\times 8\times 7+65\times 4+3+2+1.$
\item [] $771=9\times 8\times 7+65\times 4+3\times 2+1.$
\item [] $772=9\times 8+7\times 6+5^4+32+1.$
\item [] $773=98\times 7+6+5\times 4\times 3+21.$
\item [] $774=9+87+654+3+21.$
\item [] $775=98\times 7+65+4\times 3\times 2\times 1.$
\item [] $776=98\times 7+65+4\times 3\times 2+1.$
\item [] $777=9\times 8+76\times 5+4+321.$
\item [] $778=98\times 7+6+54+32\times 1.$
\item [] $779=98\times 7+6+54+32+1.$
\item [] $780=98\times 7+6\times 5+43+21.$
\item [] $781=9+8\times 76+54\times 3+2\times 1.$
\item [] $782=9+8\times 76+54\times 3+2+1.$
\item [] $783=98+7+654+3+21.$
\item [] $784=98\times 7+65+4\times 3+21.$
\item [] $785=9+8\times 7+6\times 5\times 4\times 3\times 2\times 1.$
\item [] $786=9+8\times 7+6\times 5\times 4\times 3\times 2+1.$
\item [] $787=98\times 7+65+4+32\times 1.$
\item [] $788=9\times 8\times 7+65\times 4+3+21.$
\item [] $789=9\times 8+76+5\times 4\times 32+1.$
\item [] $790=9+87+6+5^4+3\times 21.$
\item [] $791=9+8+765+4+3+2\times 1.$
\item [] $792=9+8+765+4+3+2+1.$
\item [] $793=9+8+765+4+3\times 2+1.$
\item [] $794=9+8+(7+6+5)\times 43+2+1.$
\item [] $795=9+8+765+4+3^2\times 1.$
\item [] $796=9+8+765+4\times 3+2\times 1.$
\item [] $797=9+8+765+4\times 3+2+1.$
\item [] $798=98+7\times 6+5^4+32+1.$
\item [] $799=9\times 8+7+6\times 5\times 4\times 3\times 2\times 1.$
\item [] $800=9+8\times 76+54\times 3+21.$
\item [] $801=98\times 7+6\times 5+4^3+21.$
\item [] $802=98\times 7+6\times 5+43\times 2\times 1.$
\item [] $803=9\times 87+6+5+4+3+2\times 1.$
\item [] $804=9\times 87+6+5+4+3+2+1.$
\item [] $805=9\times 87+6+5+4+3\times 2+1.$
\item [] $806=9+8+7+65\times 4\times 3+2\times 1.$
\item [] $807=9+8+7+65\times 4\times 3+2+1.$
\item [] $808=9\times 87+6+5+4\times 3+2\times 1.$
\item [] $809=98\times 7+6+54+3\times 21.$
\item [] $810=9+8+765+4+3+21.$
\item [] $811=98\times 7+6\times 5\times 4+3+2\times 1.$
\item [] $812=98\times 7+6\times 5\times 4+3+2+1.$
\item [] $813=9+87+654+3\times 21.$
\item [] $814=9\times 87+6+5\times 4+3+2\times 1.$
\item [] $815=98\times 7+65+43+21.$
\item [] $816=9+87+6\times 5\times 4\times 3\times 2\times 1.$
\item [] $817=9+87+6\times 5\times 4\times 3\times 2+1.$
\item [] $818=9+8+765+4+32\times 1.$
\item [] $819=9+8+765+4+32+1.$
\item [] $820=98+7+65\times (4+3\times 2+1).$
\item [] $821=9\times 8\times 7+65+4\times 3\times 21.$
\item [] $822=98+7+654+3\times 21.$
\item [] $823=9\times 87+6\times 5+4+3+2+1.$
\item [] $824=9\times 87+6\times 5+4+3\times 2+1.$
\item [] $825=98+7+6\times 5\times 4\times 3\times 2\times 1.$
\item [] $826=98+7+6\times 5\times 4\times 3\times 2+1.$
\item [] $827=9+8+765+43+2\times 1.$
\item [] $828=9+8+765+43+2+1.$
\item [] $829=9+8+76\times 5+432\times 1.$
\item [] $830=98\times 7+6\times 5\times 4+3+21.$
\item[]$\mbox{Increasing order}$
\item [] $831=1+2^3+4\times 5\times 6+78\times 9.$
\item [] $832=12\times 34+5\times 67+89.$
\item [] $833=1\times 2+3\times 4+5\times 6+789.$
\item [] $834=12+345+6\times 78+9.$
\item [] $835=1^2+34+5+6+789.$
\item [] $836=12\times 3\times 4+5+678+9.$
\item [] $837=1+2+34+5+6+789.$
\item [] $838=12+3+4+5\times 6+789.$
\item [] $839=1+23+4\times 5+6+789.$
\item [] $840=12\times 3+4+5+6+789.$
\item [] $841=1\times 2+3^4+56+78\times 9.$
\item [] $842=123\times 4+5+6\times 7\times 8+9.$
\item [] $843=12+3\times 4+5\times 6+789.$
\item [] $844=123\times 4+5\times 67+8+9.$
\item [] $845=1\times 2+3\times 45+6+78\times 9.$
\item [] $846=12+34+5+6+789.$
\item [] $847=1+2\times 345+67+89.$
\item [] $848=1\times 2^3+45+6+789.$
\item [] $849=1+2^3+45+6+789.$
\item [] $850=1\times 23\times 4+56+78\times 9.$
\item [] $851=123+4\times 56+7\times 8\times 9.$
\item [] $852=1\times 2+3+4\times 56+7\times 89.$
\item [] $853=1+2+3+4\times 56+7\times 89.$
\item [] $854=1\times 2+3+4+56+789.$
\item [] $855=123+45+678+9.$
\item [] $856=1+2+34+5\times 6+789.$
\item [] $857=12\times 3\times 4+5+6+78\times 9.$
\item [] $858=1+2+345+6+7\times 8\times 9.$
\item [] $859=12\times 3+4+5\times 6+789.$
\item [] $860=1+2+3\times 4+56+789.$
\item [] $861=1^2+3+4\times 5\times 6\times 7+8+9.$
\item [] $862=12+3+4\times 56+7\times 89.$
\item [] $863=1+2+3+4\times 5\times 6\times 7+8+9.$
\item [] $864=12+3+4+56+789.$
\item [] $865=12+34+5\times 6+789.$
\item [] $866=123+4\times 5\times 6+7\times 89.$
\item [] $867=12+345+6+7\times 8\times 9.$
\item [] $868=1\times 234+5+6+7\times 89.$
\item [] $869=12+3\times 4+56+789.$
\item [] $870=1+2\times 3\times 4+56+789.$
\item [] $871=1+23+4\times 56+7\times 89.$
\item [] $872=12\times 34+56\times 7+8\times 9.$
\item [] $873=1+23+4+56+789.$
\item [] $874=1+234+567+8\times 9.$
\item [] $875=1\times 2\times 3^4+5+6+78\times 9.$
\item [] $876=12\times 3+45+6+789.$
\item [] $877=1\times 2+3\times 45\times 6+7\times 8+9.$
\item [] $878=1+2+3\times 45\times 6+7\times 8+9.$
\item [] $879=1^2+34\times 5+6+78\times 9.$
\item [] $880=1\times 2+34\times 5+6+78\times 9.$
\item [] $881=1\times 2+34+56+789.$
\item [] $882=1+2+34+56+789.$
\item [] $883=12\times 3+4\times 56+7\times 89.$
\item [] $884=1+2+3^4+5+6+789.$
\item [] $885=123+4+56+78\times 9.$
\item [] $886=1+2\times 34\times 5+67\times 8+9.$
\item [] $887=12+3\times 45\times 6+7\times 8+9.$
\item [] $888=1+234+5\times 6+7\times 89.$
\item [] $889=12\times 34+56\times 7+89.$
\item [] $890=12+34\times 5+6+78\times 9.$
\item [] $891=12+34+56+789.$
\item [] $892=1\times 23\times 4+5+6+789.$
\item [] $893=1+23\times 4+5+6+789.$
\item [] $894=12+3^4\times 5+6\times 78+9.$
\item [] $895=1+2\times 3^4+5\times 6+78\times 9.$
\item [] $896=1+2\times 3^4\times 5+6+7+8\times 9.$
\item [] $897=123+45\times 6+7\times 8\times 9.$
\item [] $898=1\times 2+3+45\times 6+7\times 89.$
\item [] $899=123\times 4+5\times 67+8\times 9.$
\item [] $900=1+2\times 3+45\times 6+7\times 89.$
\item[]$\mbox{Decreasing order}$
\item [] $831=9\times 87+6+5+4+32+1.$
\item [] $832=9+8\times 7\times 6+54\times 3^2+1.$
\item [] $833=9\times 8\times 7+6\times 54+3+2\times 1.$
\item [] $834=9\times 8\times 7+6\times 54+3+2+1.$
\item [] $835=9+8\times 76+5\times 43+2+1.$
\item [] $836=98\times 7+65+4^3+21.$
\item [] $837=9\times 87+6\times 5+4\times 3\times 2\times 1.$
\item [] $838=9\times 87+6\times 5+4\times 3\times 2+1.$
\item [] $839=9\times 87+6+5+43+2\times 1.$
\item [] $840=9\times 87+6+5+43+2+1.$
\item [] $841=9\times 87+6\times 5+4+3+21.$
\item [] $842=9\times 87+6+5\times 4+32+1.$
\item [] $843=9\times (8+7)\times 6+5+4+3+21.$
\item [] $844=98\times 7+6\times 5+4\times 32\times 1.$
\item [] $845=98\times 7+6\times 5+4\times 32+1.$
\item [] $846=9+8+765+43+21.$
\item [] $847=9\times 8+765+4+3+2+1.$
\item [] $848=9\times 8+765+4+3\times 2+1.$
\item [] $849=9\times 87+6+54+3+2+1.$
\item [] $850=9\times 87+6+54+3\times 2+1.$
\item [] $851=9\times 8+765+4\times 3+2\times 1.$
\item [] $852=9\times 8+765+4\times 3+2+1.$
\item [] $853=9+8\times 76+5\times 43+21.$
\item [] $854=9+87\times 6+5\times 4^3+2+1.$
\item [] $855=9+87\times 6+54\times 3\times 2\times 1.$
\item [] $856=9+87\times 6+54\times 3\times 2+1.$
\item [] $857=98\times 7+6+54\times 3+2+1.$
\item [] $858=9\times 87+65+4+3+2+1.$
\item [] $859=9\times 87+6\times 5+43+2+1.$
\item [] $860=9\times 8\times 7+6\times 54+32\times 1.$
\item [] $861=9\times 8+765+4\times 3\times 2\times 1.$
\item [] $862=9\times 8+7+65\times 4\times 3+2+1.$
\item [] $863=9\times 87+65+4\times 3+2+1.$
\item [] $864=9+8+7\times 6\times 5\times 4+3\times 2+1.$
\item [] $865=9\times 8+765+4+3+21.$
\item [] $866=9\times 8\times 7+6\times 5\times 4\times 3+2\times 1.$
\item [] $867=9\times 87+6+54+3+21.$
\item [] $868=9+8+765+43\times 2\times 1.$
\item [] $869=98\times 7+6\times 5\times 4+3\times 21.$
\item [] $870=9\times 8+765+4\times 3+21.$
\item [] $871=9\times (8+7)\times 6+54+3\times 2+1.$
\item [] $872=98+765+4+3+2\times 1.$
\item [] $873=98+765+4+3+2+1.$
\item [] $874=98+765+4+3\times 2+1.$
\item [] $875=9\times 87+6+54+32\times 1.$
\item [] $876=9\times 87+65+4+3+21.$
\item [] $877=98+765+4\times 3+2\times 1.$
\item [] $878=98+765+4\times 3+2+1.$
\item [] $879=9+87+65\times 4\times 3+2+1.$
\item [] $880=98\times 7+65+4\times 32+1.$
\item [] $881=9\times 87+65+4\times 3+21.$
\item [] $882=9\times 8+765+43+2\times 1.$
\item [] $883=9\times 8+765+43+2+1.$
\item [] $884=9\times 87+65+4+32\times 1.$
\item [] $885=9\times 87+65+4+32+1.$
\item [] $886=9\times 87+6+(5+43)\times 2+1.$
\item [] $887=98+765+4\times 3\times 2\times 1.$
\item [] $888=98+765+4\times 3\times 2+1.$
\item [] $889=9+8+7\times 6\times 5\times 4+32\times 1.$
\item [] $890=9+8+7\times 6\times 5\times 4+32+1.$
\item [] $891=98+765+4+3+21.$
\item [] $892=9\times (87+6+5)+4+3+2+1.$
\item [] $893=9\times 87+65+43+2\times 1.$
\item [] $894=9\times 87+65+43+2+1.$
\item [] $895=(98+76)\times 5+4\times 3\times 2+1.$
\item [] $896=98+765+4\times 3+21.$
\item [] $897=9+87+65\times 4\times 3+21.$
\item [] $898=9\times 87+6\times 5+4^3+21.$
\item [] $899=9+876+5+4+3+2\times 1.$
\item [] $900=98+765+4+32+1.$
\item[]$\mbox{Increasing order}$
\item [] $901=123\times 4+56\times 7+8+9.$
\item [] $902=12+345+67\times 8+9.$
\item [] $903=1+2+3^4+5\times 6+789.$
\item [] $904=1+2\times 3^4\times 5+6+78+9.$
\item [] $905=1\times 2^3\times 45+67\times 8+9.$
\item [] $906=1+2^3\times 45+67\times 8+9.$
\item [] $907=1^2+3\times 45\times 6+7+89.$
\item [] $908=12+3+45\times 6+7\times 89.$
\item [] $909=12+3\times 45\times 6+78+9.$
\item [] $910=1^{23}+4\times 5\times 6+789.$
\item [] $911=1\times 23\times 4+5\times 6+789.$
\item [] $912=1+23\times 4+5\times 6+789.$
\item [] $913=1\times 2\times 34+56+789.$
\item [] $914=1+234+56+7\times 89.$
\item [] $915=1\times 2\times 3\times 4\times 5+6+789.$
\item [] $916=123\times 4+5\times 67+89.$
\item [] $917=1+23+45\times 6+7\times 89.$
\item [] $918=12+3\times 45\times 6+7+89.$
\item [] $919=1+2\times 3+4\times 5\times 6\times 7+8\times 9.$
\item [] $920=1\times 2\times 3^4+56+78\times 9.$
\item [] $921=1+2^3+4\times 5\times 6\times 7+8\times 9.$
\item [] $922=1^2+345+6\times (7+89).$
\item [] $923=12\times 34+5+6+7\times 8\times 9.$
\item [] $924=12+3+4\times 5\times 6+789.$
\item [] $925=1+2\times 3^4\times 5+6\times 7+8\times 9.$
\item [] $926=1\times 234+5+678+9.$
\item [] $927=123+4+5+6+789.$
\item [] $928=1\times 2+3^4+56+789.$
\item [] $929=12\times 3+45\times 6+7\times 89.$
\item [] $930=1^2+3+4\times 56+78\times 9.$
\item [] $931=1\times 2+3+4\times 56+78\times 9.$
\item [] $932=1\times 23+4\times 5\times 6+789.$
\item [] $933=1+23+4\times 5\times 6+789.$
\item [] $934=1\times 2+3+4\times 5\times 6\times 7+89.$
\item [] $935=1+2+3+4\times 5\times 6\times 7+89.$
\item [] $936=1+2\times 3+4\times 5\times 6\times 7+89.$
\item [] $937=1\times 23\times 4+56+789.$
\item [] $938=123+4\times 5+6+789.$
\item [] $939=1+2+3\times 4\times (5+67)+8\times 9.$
\item [] $940=(1\times 2\times 34+56)\times 7+8\times 9.$
\item [] $941=12+3+4\times 56+78\times 9.$
\item [] $942=12+3\times 45+6+789.$
\item [] $943=(1\times 23\times 4+5\times 6)\times 7+89.$
\item [] $944=12\times 3\times 4+5+6+789.$
\item [] $945=123+4\times 5\times 6+78\times 9.$
\item [] $946=123+4+5\times 6+789.$
\item [] $947=1\times 234+5+6+78\times 9.$
\item [] $948=1+234+5+6+78\times 9.$
\item [] $949=1\times 23+4\times 56+78\times 9.$
\item [] $950=1+23+4\times 56+78\times 9.$
\item [] $951=1^2+3^4\times 5+67\times 8+9.$
\item [] $952=1\times 23+4\times 5\times 6\times 7+89.$
\item [] $953=1+23+4\times 5\times 6\times 7+89.$
\item [] $954=12+3\times (45+6)+789.$
\item [] $955=1\times 2^3\times 4\times 5+6+789.$
\item [] $956=123\times 4+56\times 7+8\times 9.$
\item [] $957=1\times 2\times 3\times 45+678+9.$
\item [] $958=1+2\times 3\times 45+678+9.$
\item [] $959=1\times 23\times 4+(5+6)\times 78+9.$
\item [] $960=1^{23}\times 456+7\times 8\times 9.$
\item [] $961=1^{23}+456+7\times 8\times 9.$
\item [] $962=12\times 3+4\times 56+78\times 9.$
\item [] $963=123+45+6+789.$
\item [] $964=1^2+3+456+7\times 8\times 9.$
\item [] $965=12\times 3+4\times 5\times 6\times 7+89.$
\item [] $966=1+2+3+456+7\times 8\times 9.$
\item [] $967=1+234+5\times 6+78\times 9.$
\item [] $968=12\times 34+56+7\times 8\times 9.$
\item [] $969=1\times 2\times 34\times 5+6+7\times 89.$
\item [] $970=123+4\times 56+7\times 89.$
\item[]$\mbox{Decreasing order}$
\item [] $901=9+876+5+4+3\times 2+1.$
\item [] $902=(9+8)\times 7+65\times 4\times 3+2+1.$
\item [] $903=9+876+5+4+3\times (2+1).$
\item [] $904=9+876+5+4\times 3+2\times 1.$
\item [] $905=9+876+5+4\times 3+2+1.$
\item [] $906=9+87\times 6+54+321.$
\item [] $907=9\times (87+6+5)+4\times 3\times 2+1.$
\item [] $908=98+765+43+2\times 1.$
\item [] $909=98+765+43+2+1.$
\item [] $910=98\times 7+6+5\times 43+2+1.$
\item [] $911=9+876+5\times 4+3+2+1.$
\item [] $912=9+876+5\times 4+3\times 2+1.$
\item [] $913=9\times 87+6\times 5\times 4+3^2+1.$
\item [] $914=9+876+5+4\times 3\times 2\times 1.$
\item [] $915=9+876+5+4\times 3\times 2+1.$
\item [] $916=(98+76)\times 5+43+2+1.$
\item [] $917=9\times 8+7\times 6\times 5\times 4+3+2\times 1.$
\item [] $918=9+876+5+4+3+21.$
\item [] $919=9\times 8+7\times 6\times 5\times 4+3\times 2+1.$
\item [] $920=9+8+7\times 6\times 5\times 4+3\times 21.$
\item [] $921=9\times (8+76)+54\times 3+2+1.$
\item [] $922=9\times 87+6+5+4\times 32\times 1.$
\item [] $923=9+876+5+4\times 3+21.$
\item [] $924=9\times 87+6\times (5\times 4+3)+2+1.$
\item [] $925=9\times 8\times 7+6\times 5\times (4+3)\times 2+1.$
\item [] $926=9+876+5+4\times 3^2\times 1.$
\item [] $927=98+765+43+21.$
\item [] $928=98\times 7+6+5\times 43+21.$
\item [] $929=9+876+5\times 4+3+21.$
\item [] $930=98+765+4+3\times 21.$
\item [] $931=9+876+5\times (4+3+2)+1.$
\item [] $932=(9+8)\times 7\times 6+5\times 43+2+1.$
\item [] $933=9\times 87+65+4^3+21.$
\item [] $934=9\times 87+65+43\times 2\times 1.$
\item [] $935=9+876+5+43+2\times 1.$
\item [] $936=9+876+5+43+2+1.$
\item [] $937=9+876+5\times 4+32\times 1.$
\item [] $938=9+876+5\times 4+32+1.$
\item [] $939=9+8\times 76+5\times 4^3+2\times 1.$
\item [] $940=9\times 8\times 7+6+5\times 43\times 2\times 1.$
\item [] $941=9+8\times 76+54\times 3\times 2\times 1.$
\item [] $942=9+8\times 76+54\times 3\times 2+1.$
\item [] $943=98+7\times 6\times 5\times 4+3+2\times 1.$
\item [] $944=9+876+54+3+2\times 1.$
\item [] $945=9+876+54+3+2+1.$
\item [] $946=9+876+54+3\times 2+1.$
\item [] $947=9+876+5\times 4\times 3+2\times 1.$
\item [] $948=9+876+5\times 4\times 3+2+1.$
\item [] $949=98+765+43\times 2\times 1.$
\item [] $950=98+765+43\times 2+1.$
\item [] $951=98\times 7+65\times 4+3+2\times 1.$
\item [] $952=98\times 7+65\times 4+3+2+1.$
\item [] $953=9\times 87+6+54\times 3+2\times 1.$
\item [] $954=9+876+5+43+21.$
\item [] $955=98\times 7+65\times 4+3^2\times 1.$
\item [] $956=9+876+5+4^3+2\times 1.$
\item [] $957=9+876+5+4+3\times 21.$
\item [] $958=9+8\times 76+5\times 4+321.$
\item [] $959=9\times 8+7\times 65+432\times 1.$
\item [] $960=9\times 8+7\times 65+432+1.$
\item [] $961=9+87\times 6+5\times 43\times 2\times 1.$
\item [] $962=98+7\times 6\times 5\times 4+3+21.$
\item [] $963=9+876+54+3+21.$
\item [] $964=98+(7+65)\times 4\times 3+2\times 1.$
\item [] $965=9\times 8+765+4\times 32\times 1.$
\item [] $966=9\times 8+765+4\times 32+1.$
\item [] $967=9+8+7\times 6+5+43\times 21.$
\item [] $968=9+87\times 6+5+432\times 1.$
\item [] $969=9+87\times 6+5+432+1.$
\item [] $970=98\times 7+65\times 4+3+21.$
\item[]$\mbox{Increasing order}$
\item [] $971=1+23\times 4\times 5+6+7\times 8\times 9.$
\item [] $972=123+4+56+789.$
\item [] $973=123\times 4+56\times 7+89.$
\item [] $974=123\times 4+5+6\times 78+9.$
\item [] $975=12+3+456+7\times 8\times 9.$
\item [] $976=1\times 2+345+6+7\times 89.$
\item [] $977=12+34\times 5+6+789.$
\item [] $978=1\times 2\times 3\times 45+6+78\times 9.$
\item [] $979=1+2\times 3+45\times 6+78\times 9.$
\item [] $980=123+4\times 5\times 6\times 7+8+9.$
\item [] $981=1+2^3+45\times 6+78\times 9.$
\item [] $982=1+2\times 3^4+5\times 6+789.$
\item [] $983=1\times 23+456+7\times 8\times 9.$
\item [] $984=1+23+456+7\times 8\times 9.$
\item [] $985=1\times 2+3\times 4\times 5\times 6+7\times 89.$
\item [] $986=12+345+6+7\times 89.$
\item [] $987=12+3+45\times 6+78\times 9.$
\item [] $988=1^2+3^4\times (5+6)+7+89.$
\item [] $989=12\times 3\times 4+56+789.$
\item [] $990=1+2^3\times 45+6+7\times 89.$
\item [] $991=1+2\times 3\times (4+5+67+89).$
\item [] $992=12\times 34+567+8+9.$
\item [] $993=1+234+56+78\times 9.$
\item [] $994=123+4+(5+6)\times 78+9.$
\item [] $995=12+3\times 4\times 5\times 6+7\times 89.$
\item [] $996=12\times 3+456+7\times 8\times 9.$
\item [] $997=1+2\times (3\times 4\times 5+6)\times 7+8\times 9.$
\item [] $998=1\times 2+3\times 4\times (5+6)\times 7+8\times 9.$
\item [] $999=12\times 3\times (4+5)+(67+8)\times 9.$
\item [] $1000=1+2+34\times (5+6)+7\times 89.$
\item [] $1001=1\times 2+(3+4)\times 5\times 6+789.$
\item [] $1002=1\times 2\times (345+67+89).$
\item [] $1003=1+23\times (4+5)+6+789.$
\item [] $1004=123+4\times (5\times 6\times 7+8)+9.$
\item [] $1005=1\times 23\times 4\times 5+67\times 8+9.$
\item [] $1006=1+23\times 4\times 5+67\times 8+9.$
\item [] $1007=12\times 3^4+5+6+7+8+9.$
\item [] $1008=12\times 3+45\times 6+78\times 9.$
\item [] $1009=12+34\times (5+6)+7\times 89.$
\item [] $1010=1\times (2+3+4+5)\times 67+8\times 9.$
\item [] $1011=12\times 3\times (4+5)+678+9.$
\item [] $1012=1+2\times (3^4+5\times 6)+789.$
\item [] $1013=1^{23}\times 4\times 56+789.$
\item [] $1014=1^{23}+4\times 56+789.$
\item [] $1015=1\times 2\times 34\times 5+(67+8)\times 9.$
\item [] $1016=123+45\times 6+7\times 89.$
\item [] $1017=1^2+3+4\times 56+789.$
\item [] $1018=1\times 2+3+4\times 56+789.$
\item [] $1019=1+2+3+4\times 56+789.$
\item [] $1020=1+2\times 3+4\times 56+789.$
\item [] $1021=1\times 2^3+4\times 56+789.$
\item [] $1022=1+2^3+4\times 56+789.$
\item [] $1023=1+2+345+(67+8)\times 9.$
\item [] $1024=12+3+4\times 5\times (6\times 7+8)+9.$
\item [] $1025=1\times 2^3\times 4\times 5\times 6+7\times 8+9.$
\item [] $1026=123\times 4+5\times 6+7\times 8\times 9.$
\item [] $1027=12+(3+4)\times 56+7\times 89.$
\item [] $1028=12+3+4\times 56+789.$
\item [] $1029=12+(3+4+5+6)\times 7\times 8+9.$
\item [] $1030=1\times 2\times (34+56\times 7+89).$
\item [] $1031=1\times (23+45)\times 6+7\times 89.$
\item [] $1032=123+4\times 5\times 6+789.$
\item [] $1033=1^2+345+678+9.$
\item [] $1034=1\times 234+5+6+789.$
\item [] $1035=1+234+5+6+789.$
\item [] $1036=1\times 23+4\times 56+789.$
\item [] $1037=1+23+4\times 56+789.$
\item [] $1038=(123+4+5+6)\times 7+8\times 9.$
\item [] $1039=1\times 2^3\times 4\times 5\times 6+7+8\times 9.$
\item [] $1040=(1+2)^3+4\times 56+789.$
\item[]$\mbox{Decreasing order}$
\item [] $971=98+7\times 6\times 5\times 4+32+1.$
\item [] $972=9+876+54+32+1.$
\item [] $973=9\times (8+7+6)\times 5+4+3+21.$
\item [] $974=9\times (8+7)\times 6+54\times 3+2\times 1.$
\item [] $975=9\times 8+7\times 6\times 5\times 4+3\times 21.$
\item [] $976=9\times 87+65+4\times 32\times 1.$
\item [] $977=9\times 87+65+4\times 32+1.$
\item [] $978=98\times 7+65\times 4+32\times 1.$
\item [] $979=9+8\times 7+6+5+43\times 21.$
\item [] $980=9+8\times 7\times 6+5^4+3^2+1.$
\item [] $981=9+87\times (6+5)+4\times 3+2+1.$
\item [] $982=9+876+(5+43)\times 2+1.$
\item [] $983=98+(7+65)\times 4\times 3+21.$
\item [] $984=9\times 8+765+(4+3)\times 21.$
\item [] $985=98+7\times 65+432\times 1.$
\item [] $986=98+7\times 65+432+1.$
\item [] $987=(98+7+6)\times 5+432\times 1.$
\item [] $988=9\times 8\times (7+6)+5\times 4+32\times 1.$
\item [] $989=(98+76+5\times 4^3)\times 2+1.$
\item [] $990=9\times 8+7+65\times (4+3)\times 2+1.$
\item [] $991=98+765+4\times 32\times 1.$
\item [] $992=9+8\times 76+54+321.$
\item [] $993=9\times 8+7+6+5+43\times 21.$
\item [] $994=9+8\times 7\times 6+5^4+3+21.$
\item [] $995=9+(8+7)\times 65+4+3\times 2+1.$
\item [] $996=9\times 8\times 7+6+54\times 3^2\times 1.$
\item [] $997=9\times 8\times 7+6+54\times 3^2+1.$
\item [] $998=9+8+7+6\times 54\times 3+2\times 1.$
\item [] $999=9+8+7+654+321.$
\item [] $1000=9+876+(54+3)\times 2+1.$
\item [] $1001=9\times 8\times 7+65+432\times 1.$
\item [] $1002=9\times 8\times 7+65+432+1.$
\item [] $1003=98\times 7+65+4\times 3\times 21.$
\item [] $1004=9+8\times (76+5+43)+2+1.$
\item [] $1005=9+876+5\times 4\times 3\times 2\times 1.$
\item [] $1006=9\times 87+6+5\times 43+2\times 1.$
\item [] $1007=9\times 87+6+5\times 43+2+1.$
\item [] $1008=987+6+5+4+3+2+1.$
\item [] $1009=987+6+5+4+3\times 2+1.$
\item [] $1010=9+87+6+5+43\times 21.$
\item [] $1011=987+6+5+4+3^2\times 1.$
\item [] $1012=9\times 8+7+6\times 5+43\times 21.$
\item [] $1013=987+6+5+4\times 3+2+1.$
\item [] $1014=98\times 7+6+5\times 4^3+2\times 1.$
\item [] $1015=98\times 7+6\times 54+3+2\times 1.$
\item [] $1016=98\times 7+6+54\times 3\times 2\times 1.$
\item [] $1017=9+8+7+6\times 54\times 3+21.$
\item [] $1018=987+6+5\times 4+3+2\times 1.$
\item [] $1019=987+6+5\times 4+3+2+1.$
\item [] $1020=987+6+5\times 4+3\times 2+1.$
\item [] $1021=9+(8+7)\times 65+4+32+1.$
\item [] $1022=98\times 7+6+5+4+321.$
\item [] $1023=987+6+5+4\times 3\times 2+1.$
\item [] $1024=98+7\times (6+5)\times 4\times 3+2\times 1.$
\item [] $1025=9\times 87+6+5\times 43+21.$
\item [] $1026=987+6+5+4+3+21.$
\item [] $1027=987+6\times 5+4+3+2+1.$
\item [] $1028=987+6\times 5+4+3\times 2+1.$
\item [] $1029=9+87+6\times 5+43\times 21.$
\item [] $1030=987+6\times 5+4+3^2\times 1.$
\item [] $1031=9\times 8+7+6+5^4+321.$
\item [] $1032=987+6\times 5+4\times 3+2+1.$
\item [] $1033=98\times 7+6+5\times 4+321.$
\item [] $1034=98\times 7+6\times 54+3+21.$
\item [] $1035=987+6+5+4+32+1.$
\item [] $1036=987+(6+5)\times 4+3+2\times 1.$
\item [] $1037=987+6+5\times 4+3+21.$
\item [] $1038=98+7+6\times 5+43\times 21.$
\item [] $1039=9+8+76+5^4+321.$
\item [] $1040=9+8\times 7+654+321.$
\item[]$\mbox{Increasing order}$
\item [] $1041=(1+2+34+5)\times 6+789.$
\item [] $1042=123\times 4+5+67\times 8+9.$
\item [] $1043=(1+2)\times 3^4+5+6+789.$
\item [] $1044=12+345+678+9.$
\item [] $1045=1^2+34\times 5\times 6+7+8+9.$
\item [] $1046=1\times 2+34\times 5\times 6+7+8+9.$
\item [] $1047=12\times 34+567+8\times 9.$
\item [] $1048=12\times 3^4+5+6+7\times 8+9.$
\item [] $1049=123+4\times 56+78\times 9.$
\item [] $1050=12+(3+4\times 5)\times 6\times 7+8\times 9.$
\item [] $1051=1+2\times (3+45+6\times 78+9).$
\item [] $1052=123+4\times 5\times 6\times 7+89.$
\item [] $1053=1\times 234+5\times 6+789.$
\item [] $1054=1+234+5\times 6+789.$
\item [] $1055=1^{23}+4^5+6+7+8+9.$
\item [] $1056=12+34\times 5\times 6+7+8+9.$
\item [] $1057=1^2\times 3+4^5+6+7+8+9.$
\item [] $1058=1^2+3+4^5+6+7+8+9.$
\item [] $1059=1^{23}\times 45\times 6+789.$
\item [] $1060=1^{23}+45\times 6+789.$
\item [] $1061=12\times 34+5\times 6+7\times 89.$
\item [] $1062=12\times 3^4+5+6+7+8\times 9.$
\item [] $1063=1^2+3+45\times 6+789.$
\item [] $1064=12\times 34+567+89.$
\item [] $1065=12+345+6+78\times 9.$
\item [] $1066=1+23\times 45+6+7+8+9.$
\item [] $1067=12\times 3^4+5\times 6+7\times 8+9.$
\item [] $1068=1+2^3+45\times 6+789.$
\item [] $1069=12+3+4^5+6+7+8+9.$
\item [] $1070=12\times 3^4+5+6+78+9.$
\item [] $1071=12+345+6\times 7\times (8+9).$
\item [] $1072=1+2\times (3\times 45+6)+789.$
\item [] $1073=(123+45)\times 6+7\times 8+9.$
\item [] $1074=12+3+45\times 6+789.$
\item [] $1075=1+(2+3\times 4\times 5)\times 6+78\times 9.$
\item [] $1076=123\times 4+567+8+9.$
\item [] $1077=1\times 23+4^5+6+7+8+9.$
\item [] $1078=1+23+4^5+6+7+8+9.$
\item [] $1079=1\times 234+56+789.$
\item [] $1080=1+234+56+789.$
\item [] $1081=12\times 3^4+5\times 6+7+8\times 9.$
\item [] $1082=1\times 23+45\times 6+789.$
\item [] $1083=123+456+7\times 8\times 9.$
\item [] $1084=1^{23}+4^5+6\times 7+8+9.$
\item [] $1085=1\times 2\times 3+456+7\times 89.$
\item [] $1086=1+2\times 3+456+7\times 89.$
\item [] $1087=12\times 34+56+7\times 89.$
\item [] $1088=1+2^3+456+7\times 89.$
\item [] $1089=12\times 3^4+5\times 6+78+9.$
\item [] $1090=12\times 3+4^5+6+7+8+9.$
\item [] $1091=12\times 3^4+5+6\times 7+8\times 9.$
\item [] $1092=1^2\times 3^4\times 5+678+9.$
\item [] $1093=12\times 3^4+56+7\times 8+9.$
\item [] $1094=12+3+456+7\times 89.$
\item [] $1095=123+45\times 6+78\times 9.$
\item [] $1096=1^{23}+4^5+6+7\times 8+9.$
\item [] $1097=12+34\times 5\times 6+7\times 8+9.$
\item [] $1098=12+3+4^5+6\times 7+8+9.$
\item [] $1099=1^2\times 34\times 5\times 6+7+8\times 9.$
\item [] $1100=12\times 34+5+678+9.$
\item [] $1101=1\times 2+34\times 5\times 6+7+8\times 9.$
\item [] $1102=1\times 23+456+7\times 89.$
\item [] $1103=1+23+456+7\times 89.$
\item [] $1104=12+3^4\times 5+678+9.$
\item [] $1105=1+(2+34)\times 5\times 6+7+8+9.$
\item [] $1106=1\times 23\times 45+6+7\times 8+9.$
\item [] $1107=1+23\times 45+6+7\times 8+9.$
\item [] $1108=12\times 3^4+5+6\times 7+89.$
\item [] $1109=1\times 2+34\times 5\times 6+78+9.$
\item [] $1110=1+2+34\times 5\times 6+78+9.$
\item[]$\mbox{Decreasing order}$
\item [] $1041=98\times 7+6\times 5+4+321.$
\item [] $1042=987+6\times 5+4\times 3\times 2+1.$
\item [] $1043=987+6+5+43+2\times 1.$
\item [] $1044=987+6+5+43+2+1.$
\item [] $1045=987+6\times 5+4+3+21.$
\item [] $1046=987+6+5\times 4+32+1.$
\item [] $1047=9\times 8+7+65+43\times 21.$
\item [] $1048=98+7\times 6+5+43\times 21.$
\item [] $1049=9\times 87+65\times 4+3+2+1.$
\item [] $1050=9+876+54\times 3+2+1.$
\item [] $1051=98+7+(6+5)\times 43\times 2\times 1.$
\item [] $1052=987+6+54+3+2\times 1.$
\item [] $1053=987+6+54+3+2+1.$
\item [] $1054=9\times 8+7+654+321.$
\item [] $1055=987+6+5\times 4\times 3+2\times 1.$
\item [] $1056=987+6+5\times 4\times 3+2+1.$
\item [] $1057=987+6+54+3^2+1.$
\item [] $1058=9+8\times 7+6\times 54\times 3+21.$
\item [] $1059=9\times 87+6+54\times (3+2)\times 1.$
\item [] $1060=9\times 8+7\times 6+5^4+321.$
\item [] $1061=987+65+4+3+2\times 1.$
\item [] $1062=987+65+4+3+2+1.$
\item [] $1063=987+65+4+3\times 2+1.$
\item [] $1064=9+87+65+43\times 21.$
\item [] $1065=987+6+5+4+3\times 21.$
\item [] $1066=987+65+4+3^2+1.$
\item [] $1067=98\times 7+6+54+321.$
\item [] $1068=9+876+54\times 3+21.$
\item [] $1069=9+(8+7)\times 65+4^3+21.$
\item [] $1070=9+87+6\times 54\times 3+2\times 1.$
\item [] $1071=9+87+654+321.$
\item [] $1072=9\times 8+7+6\times 54\times 3+21.$
\item [] $1073=98+7+65+43\times 21.$
\item [] $1074=9\times 8\times 7+6+543+21.$
\item [] $1075=9\times 87+65\times 4+32\times 1.$
\item [] $1076=987+65+4\times 3\times 2\times 1.$
\item [] $1077=9+87\times 6+543+2+1.$
\item [] $1078=9\times 87+6+(5+4)\times 32+1.$
\item [] $1079=987+6+54+32\times 1.$
\item [] $1080=98+7+654+321.$
\item [] $1081=987+6\times 5+43+21.$
\item [] $1082=98+76+5+43\times 21.$
\item [] $1083=987+6+5+4^3+21.$
\item [] $1084=987+6+5+43\times 2\times 1.$
\item [] $1085=987+65+4\times 3+21.$
\item [] $1086=98+7\times 6+5^4+321.$
\item [] $1087=(9+8)\times 7+65+43\times 21.$
\item [] $1088=987+65+4+32\times 1.$
\item [] $1089=987+65+4+32+1.$
\item [] $1090=9+8\times (76+5)+432+1.$
\item [] $1091=98+(76+5)\times 4\times 3+21.$
\item [] $1092=9\times 87+6\times (5+43)+21.$
\item [] $1093=987+6+5\times 4\times (3+2)\times 1.$
\item [] $1094=9\times 8+76+5^4+321.$
\item [] $1095=9+87\times 6+543+21.$
\item [] $1096=9+8+7\times (65+4\times 3)\times 2+1.$
\item [] $1097=987+65+43+2\times 1.$
\item [] $1098=987+65+43+2+1.$
\item [] $1099=(98+76+5+4)\times 3\times 2+1.$
\item [] $1100=98\times 7+(65+4)\times 3\times 2\times 1.$
\item [] $1101=9+8+76+(5+43)\times 21.$
\item [] $1102=987+6\times 5+4^3+21.$
\item [] $1103=9+876+5\times 43+2+1.$
\item [] $1104=987+6\times 5+43\times 2+1.$
\item [] $1105=9+8+7+6\times 5\times 4\times 3^2+1.$
\item [] $1106=9\times 87+65\times 4+3\times 21.$
\item [] $1107=9+8+765+4+321.$
\item [] $1108=987+6+(54+3)\times 2+1.$
\item [] $1109=9\times 87+6+5\times (43+21).$
\item [] $1110=987+6+54+3\times 21.$
\item[]$\mbox{Increasing order}$
\item [] $1111=12+34\times 5\times 6+7+8\times 9.$
\item [] $1112=1^2+3+4^5+67+8+9.$
\item [] $1113=1^2+3+4^5+6+7+8\times 9.$
\item [] $1114=1\times 2+3+4^5+6+7+8\times 9.$
\item [] $1115=12\times 3+456+7\times 89.$
\item [] $1116=12\times 3^4+5+67+8\times 9.$
\item [] $1117=1\times 2^3+4^5+6+7+8\times 9.$
\item [] $1118=1\times 23+4^5+6+7\times 8+9.$
\item [] $1119=12+34\times 5\times 6+78+9.$
\item [] $1120=1+23\times 45+67+8+9.$
\item [] $1121=1^2+3+4^5+6+78+9.$
\item [] $1122=1\times 2+3+4^5+6+78+9.$
\item [] $1123=12+3+4^5+67+8+9.$
\item [] $1124=12\times 3^4+56+7+89.$
\item [] $1125=12+3^4\times 5+6+78\times 9.$
\item [] $1126=123\times 4+5+6+7\times 89.$
\item [] $1127=1^{23}+4^5+6+7+89.$
\item [] $1128=1\times 23\times 45+6+78+9.$
\item [] $1129=1+23\times 45+6+78+9.$
\item [] $1130=1^2+3+4^5+6+7+89.$
\item [] $1131=123\times 4+567+8\times 9.$
\item [] $1132=12+3+4^5+6+78+9.$
\item [] $1133=12\times 3^4+5+67+89.$
\item [] $1134=1\times 2^3+4^5+6+7+89.$
\item [] $1135=1\times 2\times 34\times 5+6+789.$
\item [] $1136=123+4\times 56+789.$
\item [] $1137=1\times 23\times 45+6+7+89.$
\item [] $1138=1+23\times 45+6+7+89.$
\item [] $1139=1^{23}+4^5+6\times 7+8\times 9.$
\item [] $1140=12\times 34+5\times 6+78\times 9.$
\item [] $1141=12+3+4^5+6+7+89.$
\item [] $1142=1\times 2+345+6+789.$
\item [] $1143=1+2+345+6+789.$
\item [] $1144=12\times 3+4^5+67+8+9.$
\item [] $1145=123\times 4+5\times 6+7\times 89.$
\item [] $1146=1\times 2^3+4^5+6\times 7+8\times 9.$
\item [] $1147=1\times 23\times 4\times 5+678+9.$
\item [] $1148=123\times 4+567+89.$
\item [] $1149=1^2\times 3\times 4\times 5\times 6+789.$
\item [] $1150=1+23\times 45+6\times 7+8\times 9.$
\item [] $1151=1\times 2+3\times 4\times 5\times 6+789.$
\item [] $1152=12+345+6+789.$
\item [] $1153=12\times 3+4^5+6+78+9.$
\item [] $1154=1+2+3+4+5+67\times (8+9).$
\item [] $1155=1\times 2^3\times 45+6+789.$
\item [] $1156=1+2^3\times 45+6+789.$
\item [] $1157=1\times 2+345+6\times (7+8)\times 9.$
\item [] $1158=1^{23}\times 456+78\times 9.$
\item [] $1159=1^{23}+456+78\times 9.$
\item [] $1160=1\times 2+3+4^5+6\times 7+89.$
\item [] $1161=12+3\times 4\times 5\times 6+789.$
\item [] $1162=12\times 3+4^5+6+7+89.$
\item [] $1163=1\times 2+3+456+78\times 9.$
\item [] $1164=1\times 2\times 3+456+78\times 9.$
\item [] $1165=1+2\times 3+456+78\times 9.$
\item [] $1166=12\times 34+56+78\times 9.$
\item [] $1167=1\times 2\times 345+6\times 78+9.$
\item [] $1168=1\times 23\times 4\times 5+6+78\times 9.$
\item [] $1169=1+23\times 4\times 5+6+78\times 9.$
\item [] $1170=1+2\times 3+4^5+67+8\times 9.$
\item [] $1171=123\times 4+56+7\times 89.$
\item [] $1172=1+2^3+4^5+67+8\times 9.$
\item [] $1173=12+3+456+78\times 9.$
\item [] $1174=1\times 23\times 45+67+8\times 9.$
\item [] $1175=1+23\times 45+67+8\times 9.$
\item [] $1176=1^2\times 3\times 4\times 56+7\times 8\times 9.$
\item [] $1177=123+4^5+6+7+8+9.$
\item [] $1178=12+3+4^5+67+8\times 9.$
\item [] $1179=1+2+3\times 4\times 56+7\times 8\times 9.$
\item [] $1180=1^{23}\times 4^5+67+89.$
\item[]$\mbox{Decreasing order}$
\item [] $1111=9\times 87+6+5\times 4^3+2\times 1.$
\item [] $1112=9\times 87+6\times 54+3+2\times 1.$
\item [] $1113=987+6+5\times 4\times 3\times 2\times 1.$
\item [] $1114=9\times 87+6+54\times 3\times 2+1.$
\item [] $1115=98+765+4\times 3\times 21.$
\item [] $1116=987+65+43+21.$
\item [] $1117=9+8+7+6+543\times 2+1.$
\item [] $1118=987+65+4^3+2\times 1.$
\item [] $1119=9\times 87+6+5+4+321.$
\item [] $1120=98+76+5^4+321.$
\item [] $1121=9+876+5\times 43+21.$
\item [] $1122=98\times 7+6+5\times 43\times 2\times 1.$
\item [] $1123=98\times 7+6+5\times 43\times 2+1.$
\item [] $1124=9\times 87+6+5\times (4+3\times 21).$
\item [] $1125=9\times 8\times 7+(65+4)\times 3^2\times 1.$
\item [] $1126=987+6+5+4\times 32\times 1.$
\item [] $1127=987+6+5+4\times 32+1.$
\item [] $1128=9\times 87+6\times (54+3)+2+1.$
\item [] $1129=98\times 7+6+5+432\times 1.$
\item [] $1130=98\times 7+6+5+432+1.$
\item [] $1131=987+6\times 5\times 4+3+21.$
\item [] $1132=987+(65+4+3)\times 2+1.$
\item [] $1133=987+6+5\times 4\times (3\times 2+1).$
\item [] $1134=9\times (8+7)\times 6+54\times 3\times 2\times 1.$
\item [] $1135=9\times (8+7)\times 6+54\times 3\times 2+1.$
\item [] $1136=9\times 8+7\times (65+43\times 2+1).$
\item [] $1137=987+65+4^3+21.$
\item [] $1138=987+65+43\times 2\times 1.$
\item [] $1139=987+65+43\times 2+1.$
\item [] $1140=9\times 87+6\times 54+32+1.$
\item [] $1141=9\times 8\times 7+6+5^4+3+2+1.$
\item [] $1142=9+876+5+4\times 3\times 21.$
\item [] $1143=987+(6+5\times 4)\times 3\times 2\times 1.$
\item [] $1144=9\times 8\times 7+6+5^4+3^2\times 1.$
\item [] $1145=987+6\times 5+4\times 32\times 1.$
\item [] $1146=987+6\times 5+4\times 32+1.$
\item [] $1147=(9+8+7+6+543)\times 2+1.$
\item [] $1148=98\times 7+6\times 5+432\times 1.$
\item [] $1149=9\times 8\times 7+6\times 54+321.$
\item [] $1150=9\times 8\times 7+6+5\times 4\times 32\times 1.$
\item [] $1151=9\times 8\times 7+6+5\times 4\times 32+1.$
\item [] $1152=(9+8)\times 7\times 6+5+432+1.$
\item [] $1153=9+8+7\times 6\times (5+4)\times 3+2\times 1.$
\item [] $1154=9\times 8\times (7+6)+5\times 43+2+1.$
\item [] $1155=9+876+54\times (3+2)\times 1.$
\item [] $1156=9+876+54\times (3+2)+1.$
\item [] $1157=987+6+54\times 3+2\times 1.$
\item [] $1158=987+6+54\times 3+2+1.$
\item [] $1159=9\times 8\times 7+6+5^4+3+21.$
\item [] $1160=9\times 8+7+6\times 5\times 4\times 3^2+1.$
\item [] $1161=9+87\times 6+5^4+3+2\times 1.$
\item [] $1162=9\times 8+765+4+321.$
\item [] $1163=9\times 8\times 7+654+3+2\times 1.$
\item [] $1164=9\times 87+6+54+321.$
\item [] $1165=9\times 8\times 7+654+3\times 2+1.$
\item [] $1166=9+87\times 6+5^4+3^2+1.$
\item [] $1167=9\times 8\times 7+6+5^4+32\times 1.$
\item [] $1168=9\times 8\times 7+654+3^2+1.$
\item [] $1169=987+(6+54)\times 3+2\times 1.$
\item [] $1170=987+6\times 5\times 4+3\times 21.$
\item [] $1171=9+87\times 6+5\times 4\times 32\times 1.$
\item [] $1172=9\times 8+7+6+543\times 2+1.$
\item [] $1173=9\times 87+65+4+321.$
\item [] $1174=987+6+5\times 4\times 3^2+1.$
\item [] $1175=(98+7)\times 6+543+2\times 1.$
\item [] $1176=987+6+54\times 3+21.$
\item [] $1177=9+87+6\times 5\times 4\times 3^2+1.$
\item [] $1178=9+8+7\times 6\times 5\times 4+321.$
\item [] $1179=9+8+76+543\times 2\times 1.$
\item [] $1180=987+65+4\times 32\times 1.$
\item[]$\mbox{Increasing order}$
\item [] $1181=1\times 23+456+78\times 9.$
\item [] $1182=123+45\times 6+789.$
\item [] $1183=1^2\times 3+4^5+67+89.$
\item [] $1184=123\times 4+5+678+9.$
\item [] $1185=1\times 2+3+4^5+67+89.$
\item [] $1186=1+2+3+4^5+67+89.$
\item [] $1187=1+2\times 3+4^5+67+89.$
\item [] $1188=1\times 2^3+4^5+67+89.$
\item [] $1189=1+2^3+4^5+67+89.$
\item [] $1190=(1+2)^3+4^5+67+8\times 9.$
\item [] $1191=1\times 23\times 45+67+89.$
\item [] $1192=1+23\times 45+67+89.$
\item [] $1193=12+(3+4)\times 56+789.$
\item [] $1194=12\times 3+456+78\times 9.$
\item [] $1195=12+3+4^5+67+89.$
\item [] $1196=1^2\times 34+5+(6+7)\times 89.$
\item [] $1197=12\times 3\times 4\times 5+6\times 78+9.$
\item [] $1198=1\times 2+34+5+(6+7)\times 89.$
\item [] $1199=12\times 3+4^5+67+8\times 9.$
\item [] $1200=1\times 234\times 5+6+7+8+9.$
\item [] $1201=1+234\times 5+6+7+8+9.$
\item [] $1202=123+456+7\times 89.$
\item [] $1203=1+2+3^4\times 5+6+789.$
\item [] $1204=1+23+4^5+67+89.$
\item [] $1205=123\times 4+5+6+78\times 9.$
\item [] $1206=123+4^5+6\times 7+8+9.$
\item [] $1207=1+2\times 3\times (45+67+89).$
\item [] $1208=12\times 34+5+6+789.$
\item [] $1209=1\times 2\times (3+4)\times 5\times 6+789.$
\item [] $1210=1+2\times (3+4)\times 5\times 6+789.$
\item [] $1211=123\times 4+5+6\times 7\times (8+9).$
\item [] $1212=12+3^4\times 5+6+789.$
\item [] $1213=1+2\times 34+5+67\times (8+9).$
\item [] $1214=1+2+(3\times 4+5)\times 67+8\times 9.$
\item [] $1215=(12+3+45+67+8)\times 9.$
\item [] $1216=12\times 3+4^5+67+89.$
\item [] $1217=12+3+45+(6+7)\times 89.$
\item [] $1218=123+4^5+6+7\times 8+9.$
\item [] $1219=1\times 2\times (34+567)+8+9.$
\item [] $1220=(12+3)\times 45+67\times 8+9.$
\item [] $1221=(12+3\times 4\times 5)\times 6+789.$
\item [] $1222=1+2^3\times (4+5)\times 6+789.$
\item [] $1223=12\times 34+5+6\times (7+8)\times 9.$
\item [] $1224=1\times 2\times 3\times 4\times 5\times 6+7\times 8\times 9.$
\item [] $1225=1+2\times 3\times 4\times 5\times 6+7\times 8\times 9.$
\item [] $1226=1+23+45+(6+7)\times 89.$
\item [] $1227=12\times 34+5\times 6+789.$
\item [] $1228=1+(2\times 34+5)\times 6+789.$
\item [] $1229=1\times 234\times 5+6\times 7+8+9.$
\item [] $1230=1+234\times 5+6\times 7+8+9.$
\item [] $1231=123+4^5+67+8+9.$
\item [] $1232=123+4^5+6+7+8\times 9.$
\item [] $1233=12\times (3+4+5\times 6)+789.$
\item [] $1234=(1+234)\times 5+6\times 7+8+9.$
\item [] $1235=1\times 2\times 345+67\times 8+9.$
\item [] $1236=1+2\times 345+67\times 8+9.$
\item [] $1237=12+3^4+5+67\times (8+9).$
\item [] $1238=123\times (4+5)+6\times 7+89.$
\item [] $1239=1+23\times (45+6)+7\times 8+9.$
\item [] $1240=123+4^5+6+78+9.$
\item [] $1241=1\times 234\times 5+6+7\times 8+9.$
\item [] $1242=1+234\times 5+6+7\times 8+9.$
\item [] $1243=1\times 2\times (3+4\times 56)+789.$
\item [] $1244=1+2\times (3+4\times 56)+789.$
\item [] $1245=1^{23}\times 456+789.$
\item [] $1246=1^{23}+456+789.$
\item [] $1247=1\times 2\times (3\times 4+567)+89.$
\item [] $1248=1^2\times 3+456+789.$
\item [] $1249=123+4^5+6+7+89.$
\item [] $1250=123\times 4+56+78\times 9.$
\item[]$\mbox{Decreasing order}$
\item [] $1181=987+65+4\times 32+1.$
\item [] $1182=9\times 8\times 7+654+3+21.$
\item [] $1183=98\times 7+65+432\times 1.$
\item [] $1184=98\times 7+65+432+1.$
\item [] $1185=9\times 8+7\times 6\times 5+43\times 21.$
\item [] $1186=98+7+6\times 5\times 4\times 3^2+1.$
\item [] $1187=(9+8+7\times 6)\times 5\times 4+3\times 2+1.$
\item [] $1188=98+765+4+321.$
\item [] $1189=9+87+6+543\times 2+1.$
\item [] $1190=9\times 8\times 7+654+32\times 1.$
\item [] $1191=9\times 8\times 7+654+32+1.$
\item [] $1192=987+(6\times 5+4)\times 3\times 2+1.$
\item [] $1193=9\times 8\times (7+6)+5+4\times 3\times 21.$
\item [] $1194=(98+7)\times 6+543+21.$
\item [] $1195=(98+76)\times 5+4+321.$
\item [] $1196=987+(65+4)\times 3+2\times 1.$
\item [] $1197=98+7+6+543\times 2\times 1.$
\item [] $1198=9\times 8\times 7+6+5^4+3\times 21.$
\item [] $1199=9+8\times (7+6)+543\times 2\times 1.$
\item [] $1200=9\times 8+7\times 6+543\times 2\times 1.$
\item [] $1201=9\times 8+7\times 6+543\times 2+1.$
\item [] $1202=9+(8+7\times 6\times 5)\times 4+321.$
\item [] $1203=987+(65+43)\times 2\times 1.$
\item [] $1204=987+(65+43)\times 2+1.$
\item [] $1205=9+876+5\times (43+21).$
\item [] $1206=9+(8+7+6)\times 54+3\times 21.$
\item [] $1207=9+876+5\times 4^3+2\times 1.$
\item [] $1208=9+876+5\times 4^3+2+1.$
\item [] $1209=9+876+54\times 3\times 2\times 1.$
\item [] $1210=9+876+54\times 3\times 2+1.$
\item [] $1211=987+6+5\times 43+2+1.$
\item [] $1212=(9+8)\times 7+6+543\times 2+1.$
\item [] $1213=(9\times 8+76+54)\times 3\times 2+1.$
\item [] $1214=9+8+765+432\times 1.$
\item [] $1215=9+8+765+432+1.$
\item [] $1216=(9+87)\times 6+5\times 4\times 32\times 1.$
\item [] $1217=(9+87)\times 6+5\times 4\times 32+1.$
\item [] $1218=987+6\times 5\times (4+3)+21.$
\item [] $1219=9\times 87+6+5\times 43\times 2\times 1.$
\item [] $1220=9\times 87+6+5\times 43\times 2+1.$
\item [] $1221=9\times 8\times 7+654+3\times 21.$
\item [] $1222=9+(8\times 7+6)\times 5+43\times 21.$
\item [] $1223=987+6+5\times (43+2+1).$
\item [] $1224=9\times 87+6\times 5\times 4+321.$
\item [] $1225=9\times 8\times 7+6\times 5\times 4\times 3\times 2+1.$
\item [] $1226=9+876+5\times 4+321.$
\item [] $1227=9\times 87+6+5+432+1.$
\item [] $1228=9\times (8+7)+6+543\times 2+1.$
\item [] $1229=987+6+5\times 43+21.$
\item [] $1230=9+87+(6+5+43)\times 21.$
\item [] $1231=(9+8+7)\times 6+543\times 2+1.$
\item [] $1232=98+7\times 6\times (5\times 4+3\times 2+1).$
\item [] $1233=9\times 8+7\times 6\times 5\times 4+321.$
\item [] $1234=9\times 8+76+543\times 2\times 1.$
\item [] $1235=9\times 8+76+543\times 2+1.$
\item [] $1236=9+(8+7)\times 65+4\times 3\times 21.$
\item [] $1237=98\times 7+6+543+2\times 1.$
\item [] $1238=98\times 7+6+543+2+1.$
\item [] $1239=987+6\times (5+4+32+1).$
\item [] $1240=9\times 87+65\times (4+3)+2\times 1.$
\item [] $1241=9\times 87+65\times (4+3)+2+1.$
\item [] $1242=9+876+(5+4\times 3)\times 21.$
\item [] $1243=(9+8+7\times 6)\times 5\times 4+3\times 21.$
\item [] $1244=9+(87+6\times 54)\times 3+2\times 1.$
\item [] $1245=9\times 87+6\times 5+432\times 1.$
\item [] $1246=9\times 87+6\times 5+432+1.$
\item [] $1247=9+8\times 76+5^4+3+2\times 1.$
\item [] $1248=9+8\times 76+5^4+3+2+1.$
\item [] $1249=9+8\times 76+5^4+3\times 2+1.$
\item [] $1250=987+6+5+4\times 3\times 21.$
\item[]$\mbox{Increasing order}$
\item [] $1251=1+2+3+456+789.$
\item [] $1252=1+2\times 3+456+789.$
\item [] $1253=12\times 34+56+789.$
\item [] $1254=1\times 234\times 5+67+8+9.$
\item [] $1255=1\times 234\times 5+6+7+8\times 9.$
\item [] $1256=1+23\times 4\times 5+6+789.$
\item [] $1257=12\times 34+56\times (7+8)+9.$
\item [] $1258=1+2\times (34+5)\times 6+789.$
\item [] $1259=1\times 2\times 3\times 4\times 5+67\times (8+9).$
\item [] $1260=12+3+456+789.$
\item [] $1261=123+4^5+6\times 7+8\times 9.$
\item [] $1262=(1+2+3+4)\times 56+78\times 9.$
\item [] $1263=1\times 234\times 5+6+78+9.$
\item [] $1264=1+234\times 5+6+78+9.$
\item [] $1265=12\times 3\times 4\times 5+67\times 8+9.$
\item [] $1266=1\times 23\times (4+5)\times 6+7+8+9.$
\item [] $1267=1+2\times 345+6\times (7+89).$
\item [] $1268=1\times 23+456+789.$
\item [] $1269=1234+5+6+7+8+9.$
\item [] $1270=1\times 23\times 4\times 5+6\times (7+8)\times 9.$
\item [] $1271=12\times 3^4+5\times 6\times 7+89.$
\item [] $1272=1\times 234\times 5+6+7+89.$
\item [] $1273=1+234\times 5+6+7+89.$
\item [] $1274=1\times 2\times (34+567)+8\times 9.$
\item [] $1275=1+2\times (34+567)+8\times 9.$
\item [] $1276=1\times 2+3\times 45+67\times (8+9).$
\item [] $1277=(1+234)\times 5+6+7+89.$
\item [] $1278=123+4^5+6\times 7+89.$
\item [] $1279=1+2\times (3+45)\times 6+78\times 9.$
\item [] $1280=1+2+3\times (4+56)\times 7+8+9.$
\item [] $1281=123+456+78\times 9.$
\item [] $1282=1+(2+3)^4+567+89.$
\item [] $1283=1\times 234\times 5+(6+7)\times 8+9.$
\item [] $1284=1\times 234\times 5+6\times 7+8\times 9.$
\item [] $1285=1+234\times 5+6\times 7+8\times 9.$
\item [] $1286=123+4^5+67+8\times 9.$
\item [] $1287=1\times 2\times 3^4\times 5+6\times 78+9.$
\item [] $1288=1234+5\times 6+7+8+9.$
\item [] $1289=12+3\times (4+56\times 7)+89.$
\item [] $1290=(123+45+6)\times 7+8\times 9.$
\item [] $1291=1\times 2\times (34+567)+89.$
\item [] $1292=123\times 4+5+6+789.$
\item [] $1293=12\times (3\times 4+5\times 6)+789.$
\item [] $1294=1+(2+3+4)\times 56+789.$
\item [] $1295=1^2\times 3\times 4\times 56+7\times 89.$
\item [] $1296=1^2+3\times 4\times 56+7\times 89.$
\item [] $1297=1\times 2+3\times 4\times 56+7\times 89.$
\item [] $1298=1234+5+6\times 7+8+9.$
\item [] $1299=(1+2)^3\times 45+67+8+9.$
\item [] $1300=123+4\times 5+(6+7)\times 89.$
\item [] $1301=1\times 234\times 5+6\times 7+89.$
\item [] $1302=1+234\times 5+6\times 7+89.$
\item [] $1303=123+4^5+67+89.$
\item [] $1304=(12+3)\times 45+6+7\times 89.$
\item [] $1305=(1+2)\times 34\times 5+6+789.$
\item [] $1306=12\times 3\times 4+5+(6+7)\times 89.$
\item [] $1307=12+3\times 4\times 56+7\times 89.$
\item [] $1308=1+23\times (4+5)\times 6+7\times 8+9.$
\item [] $1309=1\times 234\times 5+67+8\times 9.$
\item [] $1310=1234+5+6+7\times 8+9.$
\item [] $1311=123\times 4+5\times 6+789.$
\item [] $1312=1+2+34\times 5+67\times (8+9).$
\item [] $1313=12\times (3\times 4\times 5+6\times 7)+89.$
\item [] $1314=1234+56+7+8+9.$
\item [] $1315=1^2+3\times 45\times 6+7\times 8\times 9.$
\item [] $1316=1\times 2+3\times 45\times 6+7\times 8\times 9.$
\item [] $1317=1+2+3\times 45\times 6+7\times 8\times 9.$
\item [] $1318=1+(2+3)^4+5+678+9.$
\item [] $1319=1\times 2\times 345+6+7\times 89.$
\item [] $1320=1+2\times 345+6+7\times 89.$
\item[]$\mbox{Decreasing order}$
\item [] $1251=9+8\times 76+5^4+3^2\times 1.$
\item [] $1252=987+65\times 4+3+2\times 1.$
\item [] $1253=987+65\times 4+3+2+1.$
\item [] $1254=987+65\times 4+3\times 2+1.$
\item [] $1255=9+(8\times 7+6)\times 5\times 4+3+2+1.$
\item [] $1256=98\times 7+6+543+21.$
\item [] $1257=9+8\times 76+5\times 4\times 32\times 1.$
\item [] $1258=9+8\times 76+5\times 4\times 32+1.$
\item [] $1259=98+7\times 6\times 5\times 4+321.$
\item [] $1260=9+876+54+321.$
\item [] $1261=98+76+543\times 2+1.$
\item [] $1262=9+8\times 76+5\times 43\times (2+1).$
\item [] $1263=987+6+54\times (3+2)\times 1.$
\item [] $1264=987+6+54\times (3+2)+1.$
\item [] $1265=(98+7\times 6\times 5)\times 4+32+1.$
\item [] $1266=9+8\times 76+5^4+3+21.$
\item [] $1267=98\times 7+65\times 4+321.$
\item [] $1268=9+8\times 7+6+(54+3)\times 21.$
\item [] $1269=9\times 8+765+432\times 1.$
\item [] $1270=9\times 8+765+432+1.$
\item [] $1271=987+65\times 4+3+21.$
\item [] $1272=(9+8+76+543)\times 2\times 1.$
\item [] $1273=(9+8+76+543)\times 2+1.$
\item [] $1274=9+8\times 76+5^4+32\times 1.$
\item [] $1275=9+8\times 76+5^4+32+1.$
\item [] $1276=9\times 87+6+54\times 3^2+1.$
\item [] $1277=987+6\times (5+43)+2\times 1.$
\item [] $1278=987+6\times (5+43)+2+1.$
\item [] $1279=987+65\times 4+32\times 1.$
\item [] $1280=987+65\times 4+32+1.$
\item [] $1281=9\times 87+65+432+1.$
\item [] $1282=987+6+(5+4)\times 32+1.$
\item [] $1283=(98+7)\times (6+5)+4\times 32\times 1.$
\item [] $1284=(98+7)\times (6+5)+4\times 32+1.$
\item [] $1285=((9+8\times 7+6)\times (5+4)+3)\times 2+1.$
\item [] $1286=9\times 8\times 7+65\times 4\times 3+2\times 1.$
\item [] $1287=9\times 8\times 7+65\times 4\times 3+2+1.$
\item [] $1288=98\times (7+6)+5+4+3+2\times 1.$
\item [] $1289=98\times (7+6)+5+4+3+2+1.$
\item [] $1290=9+8+7+6+5\times 4\times 3\times 21.$
\item [] $1291=9+8\times 7\times 6+5^4+321.$
\item [] $1292=98\times (7+6)+5+4+3^2\times 1.$
\item [] $1293=98\times (7+6)+5+4\times 3+2\times 1.$
\item [] $1294=98\times (7+6)+5+4\times 3+2+1.$
\item [] $1295=98+765+432\times 1.$
\item [] $1296=98+765+432+1.$
\item [] $1297=9\times (8+7)\times 6+54\times 3^2+1.$
\item [] $1298=9+8+7\times (6+54)\times 3+21.$
\item [] $1299=9+87+6+(54+3)\times 21.$
\item [] $1300=9+8+76\times 5+43\times 21.$
\item [] $1301=98\times (7+6)+5\times 4+3\times 2+1.$
\item [] $1302=(98+76)\times 5+432\times 1.$
\item [] $1303=(98+76)\times 5+432+1.$
\item [] $1304=987+65+4\times 3\times 21.$
\item [] $1305=9\times 8\times 7+65\times 4\times 3+21.$
\item [] $1306=(9+8)\times 76+5+4+3+2\times 1.$
\item [] $1307=98\times (7+6)+5+4+3+21.$
\item [] $1308=98+7+6+(54+3)\times 21.$
\item [] $1309=9+(8+7)\times 65+4+321.$
\item [] $1310=987+65\times 4+3\times 21.$
\item [] $1311=987+6\times (5+4)\times 3\times 2\times 1.$
\item [] $1312=987+6\times (5+4)\times 3\times 2+1.$
\item [] $1313=9+8\times (76+54+32+1).$
\item [] $1314=98+76\times (5+4+3\times 2+1).$
\item [] $1315=9+876+5\times 43\times 2\times 1.$
\item [] $1316=9+876+5\times 43\times 2+1.$
\item [] $1317=987+6+54\times 3\times 2\times 1.$
\item [] $1318=987+6+54\times 3\times 2+1.$
\item [] $1319=9+8+7+6+5+4\times 321.$
\item [] $1320=987+6\times 54+3^2\times 1.$
\item[]$\mbox{Increasing order}$
\item [] $1321=1+2\times 3^4\times 5+6+7\times 8\times 9.$
\item [] $1322=12\times 3^4+5+6\times 7\times 8+9.$
\item [] $1323=1234+5+67+8+9.$
\item [] $1324=1234+5+6+7+8\times 9.$
\item [] $1325=123+45+(6+7)\times 89.$
\item [] $1326=12+3\times 45\times 6+7\times 8\times 9.$
\item [] $1327=1+234\times 5+67+89.$
\item [] $1328=1234+(5+6)\times 7+8+9.$
\item [] $1329=1234+5\times 6+7\times 8+9.$
\item [] $1330=1+23\times (4+5)\times 6+78+9.$
\item [] $1331=(1+234)\times 5+67+89.$
\item [] $1332=1234+5+6+78+9.$
\item [] $1333=12+(34\times 5+6)\times 7+89.$
\item [] $1334=1\times 2+3\times (4+56)\times 7+8\times 9.$
\item [] $1335=(1+2)^3\times 4\times 5+6+789.$
\item [] $1336=1\times 2\times (3\times 4+567+89).$
\item [] $1337=123\times 4+56+789.$
\item [] $1338=1234+5+6\times (7+8)+9.$
\item [] $1339=12+34\times 5+(6+7)\times 89.$
\item [] $1340=1+234+5\times (6+7)\times (8+9).$
\item [] $1341=1234+5+6+7+89.$
\item [] $1342=12+34\times (5\times 6+7)+8\times 9.$
\item [] $1343=1234+5\times 6+7+8\times 9.$
\item [] $1344=1+2\times 3\times 4\times 5\times 6+7\times 89.$
\item [] $1345=1^{234}+56\times (7+8+9).$
\item [] $1346=1+(2+3\times 4+5)\times 67+8\times 9.$
\item [] $1347=1+2\times (34+567+8\times 9).$
\item [] $1348=1234+5\times (6+7+8)+9.$
\item [] $1349=12\times 3\times 4\times 5+6+7\times 89.$
\item [] $1350=1+2\times (34+56)\times 7+89.$
\item [] $1351=1234+5\times 6+78+9.$
\item [] $1352=1234+5+(6+7)\times 8+9.$
\item [] $1353=1234+5+6\times 7+8\times 9.$
\item [] $1354=1\times 2\times 3+4+56\times (7+8+9).$
\item [] $1355=1234+56+7\times 8+9.$
\item [] $1356=1+2\times 3^4\times 5+67\times 8+9.$
\item [] $1357=1^{23}\times 4\times 5\times 67+8+9.$
\item [] $1358=1^{23}+4\times 5\times 67+8+9.$
\item [] $1359=(12+3+4)\times 5\times 6+789.$
\item [] $1360=1234+5\times 6+7+89.$
\item [] $1361=1^2+3+4\times 5\times 67+8+9.$
\item [] $1362=1\times 2+3+4\times 5\times 67+8+9.$
\item [] $1363=1+2+3+4\times 5\times 67+8+9.$
\item [] $1364=1+2\times 3+4\times 5\times 67+8+9.$
\item [] $1365=1\times 2^3+4\times 5\times 67+8+9.$
\item [] $1366=1+2^3+4\times 5\times 67+8+9.$
\item [] $1367=12\times 3^4+5+6\times (7\times 8+9).$
\item [] $1368=123+456+789.$
\item [] $1369=1234+56+7+8\times 9.$
\item [] $1370=1234+5+6\times 7+89.$
\item [] $1371=12\times (3+45)+6+789.$
\item [] $1372=12+3+4\times 5\times 67+8+9.$
\item [] $1373=1^2+3+4^5+6\times 7\times 8+9.$
\item [] $1374=1^2\times 3\times 4\times 56+78\times 9.$
\item [] $1375=1^2+3\times 4\times 56+78\times 9.$
\item [] $1376=1\times 2+3\times 4\times 56+78\times 9.$
\item [] $1377=1234+56+78+9.$
\item [] $1378=1234+5+67+8\times 9.$
\item [] $1379=12\times 3^4+5\times 67+8\times 9.$
\item [] $1380=1\times 23+4\times 5\times 67+8+9.$
\item [] $1381=12\times 3^4+56\times 7+8+9.$
\item [] $1382=1\times 2\times (3\times 4+56+7\times 89).$
\item [] $1383=1234+5+6\times (7+8+9).$
\item [] $1384=12+3+4^5+6\times 7\times 8+9.$
\item [] $1385=(1+2+3)^4+5+67+8+9.$
\item [] $1386=1234+56+7+89.$
\item [] $1387=1+2\times 3^4\times 5+6\times (7+89).$
\item [] $1388=1234+5\times (6+7)+89.$
\item [] $1389=1\times (2+3)\times 4\times 5\times 6+789.$
\item [] $1390=12+34+56\times (7+8+9).$
\item[]$\mbox{Decreasing order}$
\item [] $1321=987+6\times 54+3^2+1.$
\item [] $1322=9+876+5+432\times 1.$
\item [] $1323=9+876+5+432+1.$
\item [] $1324=98\times 7+6+5^4+3\times 2+1.$
\item [] $1325=98\times (7+6)+5+43+2+1.$
\item [] $1326=98\times 7+6+5^4+3^2\times 1.$
\item [] $1327=98\times 7+6+5^4+3^2+1.$
\item [] $1328=(9\times 8+7+6)\times 5+43\times 21.$
\item [] $1329=(9\times 8+7+6\times 5)\times 4\times 3+21.$
\item [] $1330=9+8+76\times (5+4\times 3)+21.$
\item [] $1331=98\times 7+6\times 54+321.$
\item [] $1332=987+6\times (54+3)+2+1.$
\item [] $1333=98\times 7+6+5\times 4\times 32+1.$
\item [] $1334=987+6+5\times 4+321.$
\item [] $1335=987+6\times 54+3+21.$
\item [] $1336=(9+8)\times 76+5\times 4+3+21.$
\item [] $1337=98\times 7+6+5\times 43\times (2+1).$
\item [] $1338=9+8+7+6\times 5+4\times 321.$
\item [] $1339=9+8+7+(654+3)\times 2+1.$
\item [] $1340=9+8+7\times (6+54\times 3+21).$
\item [] $1341=98\times 7+6+5^4+3+21.$
\item [] $1342=987+6\times 5+4+321.$
\item [] $1343=987+6\times 54+32\times 1.$
\item [] $1344=987+6\times 54+32+1.$
\item [] $1345=98\times 7+654+3+2\times 1.$
\item [] $1346=98\times 7+654+3+2+1.$
\item [] $1347=98\times 7+654+3\times 2+1.$
\item [] $1348=9+8+7\times 6+5+4\times 321.$
\item [] $1349=987+6\times 5\times 4\times 3+2\times 1.$
\item [] $1350=987+6\times 5\times 4\times 3+2+1.$
\item [] $1351=(9+8)\times 76+54+3+2\times 1.$
\item [] $1352=987+(6+5)\times 4+321.$
\item [] $1353=9\times 87+6+543+21.$
\item [] $1354=(9+8)\times 7\times 6+5\times 4\times 32\times 1.$
\item [] $1355=9\times 8+76\times 5+43\times 21.$
\item [] $1356=9\times 87+6+(5+4)\times 3\times 21.$
\item [] $1357=9+8\times 7+6\times 5\times 43+2\times 1.$
\item [] $1358=9+8\times 7+6\times 5\times 43+2+1.$
\item [] $1359=9\times 87+6\times (5+43)\times 2\times 1.$
\item [] $1360=98\times (7+6)+54+32\times 1.$
\item [] $1361=(9+8)\times 76+5+43+21.$
\item [] $1362=9+87+6+5\times 4\times 3\times 21.$
\item [] $1363=9\times 8\times (7+6+5)+4+3\times 21.$
\item [] $1364=98\times 7+654+3+21.$
\item [] $1365=9+876+5\times 4\times (3+21).$
\item [] $1366=98\times (7+6)+5+43\times 2+1.$
\item [] $1367=(98+7\times 6+543)\times 2+1.$
\item [] $1368=987+6+54+321.$
\item [] $1369=(98+76+54)\times 3\times 2+1.$
\item [] $1370=(9+8)\times 76+54+3+21.$
\item [] $1371=9\times 8+7+6\times 5\times 43+2\times 1.$
\item [] $1372=98\times 7+654+32\times 1.$
\item [] $1373=98\times 7+654+32+1.$
\item [] $1374=987+6\times 54+3\times 21.$
\item [] $1375=9+8+7\times 65+43\times 21.$
\item [] $1376=9+8\times 7+6\times 5\times 43+21.$
\item [] $1377=987+65+4+321.$
\item [] $1378=(9+8)\times 76+54+32\times 1.$
\item [] $1379=9+8\times 7+6\times 5+4\times 321.$
\item [] $1380=98\times 7+6+5^4+3\times 21.$
\item [] $1381=98+76\times 5+43\times 21.$
\item [] $1382=9+8+76+5+4\times 321.$
\item [] $1383=(9+8)\times 76+5+43\times 2\times 1.$
\item [] $1384=(9+8)\times 76+5+43\times 2+1.$
\item [] $1385=9\times 8+76\times (5+4\times 3)+21.$
\item [] $1386=(9+87\times 6+54\times 3)\times 2\times 1.$
\item [] $1387=(9+87\times 6+54\times 3)\times 2+1.$
\item [] $1388=9+87+6\times 5\times 43+2\times 1.$
\item [] $1389=9+87+6\times 5\times 43+2+1.$
\item [] $1390=9\times 8+7+6\times 5\times 43+21.$
\item[]$\mbox{Increasing order}$
\item [] $1391=12\times (34+5\times 6)+7\times 89.$
\item [] $1392=1^2\times 3\times 456+7+8+9.$
\item [] $1393=12\times 3+4\times 5\times 67+8+9.$
\item [] $1394=1\times 2+3\times 456+7+8+9.$
\item [] $1395=1234+5+67+89.$
\item [] $1396=12\times 3^4+5\times 67+89.$
\item [] $1397=1+234+5+(6+7)\times 89.$
\item [] $1398=1\times 2\times 345+6+78\times 9.$
\item [] $1399=1+2\times 345+6+78\times 9.$
\item [] $1400=1234+(5+6)\times 7+89.$
\item [] $1401=(123+4+5+6\times 7)\times 8+9.$
\item [] $1402=1^2+3\times (45+6+7)\times 8+9.$
\item [] $1403=(1+23+4+5)\times 6\times 7+8+9.$
\item [] $1404=12+3\times 456+7+8+9.$
\item [] $1405=12\times 3+4^5+6\times 7\times 8+9.$
\item [] $1406=1+23\times (4+5\times 6)+7\times 89.$
\item [] $1407=12\times 3\times 4\times 5+678+9.$
\item [] $1408=1+2\times (3+4)\times 56+7\times 89.$
\item [] $1409=1\times 2\times 3\times 4\times 56+7\times 8+9.$
\item [] $1410=1+2\times 3\times 4\times 56+7\times 8+9.$
\item [] $1411=(1^2+3)^4\times 5+6\times 7+89.$
\item [] $1412=1^{23}\times 4\times 5\times 67+8\times 9.$
\item [] $1413=1^{23}+4\times 5\times 67+8\times 9.$
\item [] $1414=1^2+(3+4+5+6)\times 78+9.$
\item [] $1415=1^2\times 3+4\times 5\times 67+8\times 9.$
\item [] $1416=1^2+3+4\times 5\times 67+8\times 9.$
\item [] $1417=1\times 2+3+4\times 5\times 67+8\times 9.$
\item [] $1418=1+2+3+4\times 5\times 67+8\times 9.$
\item [] $1419=1+2\times 3+4\times 5\times 67+8\times 9.$
\item [] $1420=1\times 2^3+4\times 5\times 67+8\times 9.$
\item [] $1421=1+2^3+4\times 5\times 67+8\times 9.$
\item [] $1422=1\times 2\times 3\times 4\times 5\times 6+78\times 9.$
\item [] $1423=1\times 2\times 3\times 4\times 56+7+8\times 9.$
\item [] $1424=1+2\times 3\times 4\times 56+7+8\times 9.$
\item [] $1425=1234+56+(7+8)\times 9.$
\item [] $1426=1+(2+3)^4+5+6+789.$
\item [] $1427=12+3+4\times 5\times 67+8\times 9.$
\item [] $1428=12\times 3\times 4\times 5+6+78\times 9.$
\item [] $1429=1^{23}\times 4\times 5\times 67+89.$
\item [] $1430=1^{23}+4\times 5\times 67+89.$
\item [] $1431=1\times 2\times 3\times 4\times 56+78+9.$
\item [] $1432=1+2\times 3\times 4\times 56+78+9.$
\item [] $1433=1^2\times 3\times 45\times 6+7\times 89.$
\item [] $1434=1\times 2+3+4\times 5\times 67+89.$
\item [] $1435=1\times 23+4\times 5\times 67+8\times 9.$
\item [] $1436=1+2+3\times 45\times 6+7\times 89.$
\item [] $1437=1\times 2^3+4\times 5\times 67+89.$
\item [] $1438=1+2^3+4\times 5\times 67+89.$
\item [] $1439=1\times 2\times 3^4\times 5+6+7\times 89.$
\item [] $1440=1\times 2\times 3\times 4\times 56+7+89.$
\item [] $1441=1+2\times 3\times 4\times 56+7+89.$
\item [] $1442=1\times 2+3\times (456+7+8+9).$
\item [] $1443=1+2+3\times (456+7+8+9).$
\item [] $1444=12+3+4\times 5\times 67+89.$
\item [] $1445=12+3\times 456+7\times 8+9.$
\item [] $1446=1\times 2\times 345+(6+78)\times 9.$
\item [] $1447=1^2\times 3\times 456+7+8\times 9.$
\item [] $1448=12\times 3+4\times 5\times 67+8\times 9.$
\item [] $1449=1\times 2+3\times 456+7+8\times 9.$
\item [] $1450=1+2+3\times 456+7+8\times 9.$
\item [] $1451=1+2\times 34\times (5+6)+78\times 9.$
\item [] $1452=1\times 23+4\times 5\times 67+89.$
\item [] $1453=12\times 3^4+56\times 7+89.$
\item [] $1454=12\times 3^4+5+6\times 78+9.$
\item [] $1455=1^2\times 3\times 456+78+9.$
\item [] $1456=1^2+3\times 456+78+9.$
\item [] $1457=1\times 2+3\times 456+78+9.$
\item [] $1458=1+2+3\times 456+78+9.$
\item [] $1459=12+3\times 456+7+8\times 9.$
\item [] $1460=1234+5+(6+7)\times (8+9).$
\item[]$\mbox{Decreasing order}$
\item [] $1391=9+87+6+5+4\times 321.$
\item [] $1392=98+76\times (5+4\times 3)+2\times 1.$
\item [] $1393=9\times 8+7+6\times 5+4\times 321.$
\item [] $1394=9\times 8+7+(654+3)\times 2+1.$
\item [] $1395=98\times (7+6)+5\times 4\times 3\times 2+1.$
\item [] $1396=9+87+65\times 4\times (3+2)\times 1.$
\item [] $1397=98+7+6\times 5\times 43+2\times 1.$
\item [] $1398=98+7+6\times 5\times 43+2+1.$
\item [] $1399=9+87\times (6+5)+432+1.$
\item [] $1400=98+7+6+5+4\times 321.$
\item [] $1401=987+(65+4)\times 3\times 2\times 1.$
\item [] $1402=987+(65+4)\times 3\times 2+1.$
\item [] $1403=98\times 7+654+3\times 21.$
\item [] $1404=(9+8+7)\times 6+5\times 4\times 3\times 21.$
\item [] $1405=98+7+65\times 4\times (3+2)\times 1.$
\item [] $1406=98\times 7+6\times 5\times 4\times 3\times 2\times 1.$
\item [] $1407=9+87+6\times 5\times 43+21.$
\item [] $1408=9\times 8+76+5\times 4\times 3\times 21.$
\item [] $1409=(9+8)\times 76+54+3\times 21.$
\item [] $1410=9+87+6\times 5+4\times 321.$
\item [] $1411=98+76\times (5+4\times 3)+21.$
\item [] $1412=(9+8)\times 7+6\times 5\times 43+2+1.$
\item [] $1413=(9+87+6)\times 5+43\times 21.$
\item [] $1414=9+8\times 7+65+4\times 321.$
\item [] $1415=9\times 8\times 7+65\times (4+3)\times 2+1.$
\item [] $1416=98+7+6\times 5\times 43+21.$
\item [] $1417=9+(8+7)\times 65+432+1.$
\item [] $1418=9\times 8\times 7+6+5+43\times 21.$
\item [] $1419=98+7+6\times 5+4\times 321.$
\item [] $1420=9\times 87+6+5^4+3+2+1.$
\item [] $1421=9\times 87+6+5^4+3\times 2+1.$
\item [] $1422=(9+8)\times (76+5)+43+2\times 1.$
\item [] $1423=987+6+5\times 43\times 2\times 1.$
\item [] $1424=987+6+5\times 43\times 2+1.$
\item [] $1425=(9+8)\times 76+5+4\times 32\times 1.$
\item [] $1426=(9+8)\times 76+5+4\times 32+1.$
\item [] $1427=98\times 7+6+5\times (4+3)\times 21.$
\item [] $1428=987+6\times 5\times 4+321.$
\item [] $1429=98+7\times 6+5+4\times 321.$
\item [] $1430=9\times 8+7\times 65+43\times 21.$
\item [] $1431=987+6+5+432+1.$
\item [] $1432=9+8\times 7\times 6+543\times 2+1.$
\item [] $1433=9+8+7+65+4^3\times 21.$
\item [] $1434=98+76+5\times 4\times 3\times 21.$
\item [] $1435=(98+76+543)\times 2+1.$
\item [] $1436=98+7\times 6+54\times (3+21).$
\item [] $1437=9\times 8+76+5+4\times 321.$
\item [] $1438=9\times 87+6+5^4+3+21.$
\item [] $1439=9+87\times 6+5+43\times 21.$
\item [] $1440=9+8+7+6\times (5\times 43+21).$
\item [] $1441=9\times (8+7)\times 6+5^4+3\times 2\times 1.$
\item [] $1442=9\times 87+654+3+2\times 1.$
\item [] $1443=9\times 87+654+3+2+1.$
\item [] $1444=9\times 87+654+3\times 2+1.$
\item [] $1445=9+87+65+4\times 321.$
\item [] $1446=9\times 87+654+3^2\times 1.$
\item [] $1447=9\times 87+654+3^2+1.$
\item [] $1448=(98+76\times 5+4)\times 3+2\times 1.$
\item [] $1449=9+876+543+21.$
\item [] $1450=987+6\times 5+432+1.$
\item [] $1451=9+87+6+5+4^3\times 21.$
\item [] $1452=9+876+(5+4)\times 3\times 21.$
\item [] $1453=9\times 8+7+6\times 5+4^3\times 21.$
\item [] $1454=98+7+65+4\times 321.$
\item [] $1455=987+6\times (54+3+21).$
\item [] $1456=98+7\times 65+43\times 21.$
\item [] $1457=98\times (7+6)+54\times 3+21.$
\item [] $1458=9\times 8+(7+6+5+4)\times 3\times 21.$
\item [] $1459=98+7\times (6+5)+4\times 321.$
\item [] $1460=98+7+6+5+4^3\times 21.$
\item[]$\mbox{Increasing order}$
\item [] $1461=1234+5\times 6\times 7+8+9.$
\item [] $1462=1^2+3\times 4\times 56+789.$
\item [] $1463=1\times 2+3\times 4\times 56+789.$
\item [] $1464=1+2+3\times 4\times 56+789.$
\item [] $1465=12\times 3+4\times 5\times 67+89.$
\item [] $1466=1\times 2+3\times 456+7+89.$
\item [] $1467=12+3\times 456+78+9.$
\item [] $1468=1+23\times (4+56)+78+9.$
\item [] $1469=12\times (3+45+67)+89.$
\item [] $1470=(12+3)\times 45+6+789.$
\item [] $1471=123+4+56\times (7+8+9).$
\item [] $1472=(12\times 3\times 4+56)\times 7+8\times 9.$
\item [] $1473=12+3\times 4\times 56+789.$
\item [] $1474=1\times 2\times (3^4+567+89).$
\item [] $1475=1+2\times (3^4+567+89).$
\item [] $1476=12+3\times 456+7+89.$
\item [] $1477=1+23\times (4+56)+7+89.$
\item [] $1478=1+2+3^4\times (5+6+7)+8+9.$
\item [] $1479=1\times 2\times 34\times 5+67\times (8+9).$
\item [] $1480=123+4\times 5\times 67+8+9.$
\item [] $1481=(123+4\times 5)\times 6+7\times 89.$
\item [] $1482=12\times 3^4+5\times (6+7+89).$
\item [] $1483=(1+2\times 3+4)\times (56+78)+9.$
\item [] $1484=1^2\times 345+67\times (8+9).$
\item [] $1485=1\times 2\times 345+6+789.$
\item [] $1486=1+2\times 345+6+789.$
\item [] $1487=12\times 3^4+5+6+7\times 8\times 9.$
\item [] $1488=12\times 3\times 4+56\times (7+8+9).$
\item [] $1489=1\times 2+3\times 456+7\times (8+9).$
\item [] $1490=12+3\times (456+7)+89.$
\item [] $1491=1\times 2\times (345+6)+789.$
\item [] $1492=123+4^5+6\times 7\times 8+9.$
\item [] $1493=1234+5\times (6\times 7+8)+9.$
\item [] $1494=(123+4+5)\times 6+78\times 9.$
\item [] $1495=1+2\times (3\times 4\times 5+678+9).$
\item [] $1496=12+345+67\times (8+9).$
\item [] $1497=1\times 2\times 3^4\times 5+678+9.$
\item [] $1498=1+2\times 3^4\times 5+678+9.$
\item [] $1499=12+3\times 456+7\times (8+9).$
\item [] $1500=1\times 2\times 345+6\times (7+8)\times 9.$
\item [] $1501=1^{23}\times 4^5+6\times 78+9.$
\item [] $1502=1^{23}+4^5+6\times 78+9.$
\item [] $1503=(1+2)\times (345+67+89).$
\item [] $1504=1^2\times 3+4^5+6\times 78+9.$
\item [] $1505=1^2+3+4^5+6\times 78+9.$
\item [] $1506=1\times 2+3+4^5+6\times 78+9.$
\item [] $1507=1+2+3+4^5+6\times 78+9.$
\item [] $1508=1+2\times 3+4^5+6\times 78+9.$
\item [] $1509=1\times 2\times 3\times 4\times 5\times 6+789.$
\item [] $1510=1+2\times 3\times 4\times 5\times 6+789.$
\item [] $1511=1\times (2+3^4)\times (5+6+7)+8+9.$
\item [] $1512=1\times 23\times 45+6\times 78+9.$
\item [] $1513=1+23\times 45+6\times 78+9.$
\item [] $1514=1\times 2+3\times 45\times 6+78\times 9.$
\item [] $1515=12\times 3\times 4\times 5+6+789.$
\item [] $1516=1234+5\times 6\times 7+8\times 9.$
\item [] $1517=1\times 2^3\times 45+(6+7)\times 89.$
\item [] $1518=1\times 2\times 3^4\times 5+6+78\times 9.$
\item [] $1519=1+2\times 3^4\times 5+6+78\times 9.$
\item [] $1520=(1+234)\times 5+6\times 7\times 8+9.$
\item [] $1521=1\times (234+5)\times 6+78+9.$
\item [] $1522=12\times 3^4+5+67\times 8+9.$
\item [] $1523=123+4\times (5+6\times 7\times 8+9).$
\item [] $1524=12+3\times 45\times 6+78\times 9.$
\item [] $1525=1+23+4^5+6\times 78+9.$
\item [] $1526=1\times 2+34\times 5\times 6+7\times 8\times 9.$
\item [] $1527=1+2+34\times 5\times 6+7\times 8\times 9.$
\item [] $1528=12+3+4\times (5+6\times 7)\times 8+9.$
\item [] $1529=1\times 2\times 3\times 4\times (56+7)+8+9.$
\item [] $1530=1\times (234+5)\times 6+7+89.$
\item[]$\mbox{Decreasing order}$
\item [] $1461=9\times 87+654+3+21.$
\item [] $1462=987+(6+5)\times 43+2\times 1.$
\item [] $1463=98+76+5+4\times 321.$
\item [] $1464=9+(8+7+6)\times 54+321.$
\item [] $1465=9+8\times 7\times (6+5+4\times 3+2+1).$
\item [] $1466=(9\times 8+7+6)\times (5+4\times 3)+21.$
\item [] $1467=987+(6+5+4)\times 32\times 1.$
\item [] $1468=98\times 7+65\times 4\times 3+2\times 1.$
\item [] $1469=9\times 87+654+32\times 1.$
\item [] $1470=9\times 87+654+32+1.$
\item [] $1471=98\times (7+6)+5+4^3\times (2+1).$
\item [] $1472=9\times 8\times 7+65+43\times 21.$
\item [] $1473=9+8\times 7+(6+5)\times 4\times 32\times 1.$
\item [] $1474=9+8\times 7+65+4^3\times 21.$
\item [] $1475=(9+8)\times 76+54\times 3+21.$
\item [] $1476=987+6+(5\times 4+3)\times 21.$
\item [] $1477=9+87\times 6+5^4+321.$
\item [] $1478=9\times 8\times 7+6\times 54\times 3+2\times 1.$
\item [] $1479=9\times 8\times 7+654+321.$
\item [] $1480=987+6+54\times 3^2+1.$
\item [] $1481=9\times 8\times (7+6)+543+2\times 1.$
\item [] $1482=9\times 8\times (7+6)+543+2+1.$
\item [] $1483=(9+8\times 7+6)\times 5\times 4+3\times 21.$
\item [] $1484=987+65+432\times 1.$
\item [] $1485=987+65+432+1.$
\item [] $1486=(98+7+6+54)\times 3^2+1.$
\item [] $1487=98\times 7+65\times 4\times 3+21.$
\item [] $1488=9\times 8+7+65+4^3\times 21.$
\item [] $1489=98+7\times 6+5+4^3\times 21.$
\item [] $1490=9\times 87+(6+5)\times 4^3+2+1.$
\item [] $1491=98\times (7+6)+5\times 43+2\times 1.$
\item [] $1492=98\times (7+6)+5\times 43+2+1.$
\item [] $1493=(9+8)\times 7+6\times 5+4^3\times 21.$
\item [] $1494=9\times 8\times 7+6\times (54\times 3+2+1).$
\item [] $1495=9\times (8+7)\times (6+5)+4+3\times 2\times 1.$
\item [] $1496=98+(7\times 65+4)\times 3+21.$
\item [] $1497=9+8\times (76+5+4\times 3)\times 2\times 1.$
\item [] $1498=9\times 87+65\times (4+3\times 2+1).$
\item [] $1499=(98\times 7+6+54+3)\times 2+1.$
\item [] $1500=9\times 87+654+3\times 21.$
\item [] $1501=9\times 8+7\times 6\times (5+4\times 3)\times 2+1.$
\item [] $1502=98\times 7+(6\times 5+4)\times (3+21).$
\item [] $1503=9\times 87+6\times 5\times 4\times 3\times 2\times 1.$
\item [] $1504=9\times 87+6\times 5\times 4\times 3\times 2+1.$
\item [] $1505=9+87+65+4^3\times 21.$
\item [] $1506=(9+87+654+3)\times 2\times 1.$
\item [] $1507=(9+87+654+3)\times 2+1.$
\item [] $1508=(9+8+7+6+5)\times 43+2+1.$
\item [] $1509=987+6\times (54+32+1).$
\item [] $1510=98\times (7+6)+5\times 43+21.$
\item [] $1511=9+8+7\times 6\times 5+4\times 321.$
\item [] $1512=9+87+6\times (5\times 43+21).$
\item [] $1513=98+7+(6+5)\times 4\times 32\times 1.$
\item [] $1514=98+7+65+4^3\times 21.$
\item [] $1515=9+876+5^4+3+2\times 1.$
\item [] $1516=9+876+5^4+3+2+1.$
\item [] $1517=9+876+5^4+3\times 2+1.$
\item [] $1518=98\times 7+(6+5\times 4)\times 32\times 1.$
\item [] $1519=9+876+5^4+3^2\times 1.$
\item [] $1520=9+876+5^4+3^2+1.$
\item [] $1521=9\times 87+6+(5+4)^3+2+1.$
\item [] $1522=9+8\times (7\times 6+5)\times 4+3^2\times 1.$
\item [] $1523=98+76+5+4^3\times 21.$
\item [] $1524=9\times 87+6+5\times (4+3)\times 21.$
\item [] $1525=9+876+5\times 4\times 32\times 1.$
\item [] $1526=9+876+5\times 4\times 32+1.$
\item [] $1527=(98\times 7+65+4\times 3)\times 2+1.$
\item [] $1528=(9+8)\times 76+5\times 43+21.$
\item [] $1529=9+8+7\times (65+43)\times 2\times 1.$
\item [] $1530=9+8+7\times (65+43)\times 2+1.$
\item[]$\mbox{Increasing order}$
\item [] $1531=1+(234+5)\times 6+7+89.$
\item [] $1532=12\times 3^4+56+7\times 8\times 9.$
\item [] $1533=1234+5\times 6\times 7+89.$
\item [] $1534=1^{23}\times 4^5+6+7\times 8\times 9.$
\item [] $1535=123+4\times 5\times 67+8\times 9.$
\item [] $1536=12+34\times 5\times 6+7\times 8\times 9.$
\item [] $1537=12\times 3+4^5+6\times 78+9.$
\item [] $1538=1^2+3+4^5+6+7\times 8\times 9.$
\item [] $1539=1\times 2+3+4^5+6+7\times 8\times 9.$
\item [] $1540=1+2+3+4^5+6+7\times 8\times 9.$
\item [] $1541=1+2\times 3+4^5+6+7\times 8\times 9.$
\item [] $1542=1\times 2^3+4^5+6+7\times 8\times 9.$
\item [] $1543=1+2^3+4^5+6+7\times 8\times 9.$
\item [] $1544=1\times 2+(3+4)\times 5\times 6\times 7+8\times 9.$
\item [] $1545=1\times 23\times 45+6+7\times 8\times 9.$
\item [] $1546=1+23\times 45+6+7\times 8\times 9.$
\item [] $1547=1^{23}+4^5+6\times (78+9).$
\item [] $1548=1+2\times 34\times (5+6)+789.$
\item [] $1549=12+3+4^5+6+7\times 8\times 9.$
\item [] $1550=1^2+3+4^5+6\times (78+9).$
\item [] $1551=12+34\times (5\times 6+7+8)+9.$
\item [] $1552=123+4\times 5\times 67+89.$
\item [] $1553=1+2\times 3+4^5+6\times (78+9).$
\item [] $1554=(12+3)\times (45+6)+789.$
\item [] $1555=12\times (3+4\times 5\times 6)+7+8\times 9.$
\item [] $1556=12\times 3^4+567+8+9.$
\item [] $1557=1\times 23+4^5+6+7\times 8\times 9.$
\item [] $1558=1+23+4^5+6+7\times 8\times 9.$
\item [] $1559=(12+3+4\times 5)\times 6\times 7+89.$
\item [] $1560=(1+2)^3\times 45+6\times 7\times 8+9.$
\item [] $1561=1\times 2+(3+4)\times 5\times 6\times 7+89.$
\item [] $1562=1+2+(3+4)\times 5\times 6\times 7+89.$
\item [] $1563=12\times (3+4\times 5\times 6)+78+9.$
\item [] $1564=1\times 2+3^4\times 5+(6+7)\times 89.$
\item [] $1565=1\times 2+3\times (456+7\times 8+9).$
\item [] $1566=(12\times 3+4+56+78)\times 9.$
\item [] $1567=1+2\times 3^4\times 5+(6+78)\times 9.$
\item [] $1568=(12+3+4)\times 56+7\times 8\times 9.$
\item [] $1569=1^{23}\times 4^5+67\times 8+9.$
\item [] $1570=12\times 3+4^5+6+7\times 8\times 9.$
\item [] $1571=1\times 23\times (4+5\times 6)+789.$
\item [] $1572=1^2\times 3+4^5+67\times 8+9.$
\item [] $1573=1^2+3+4^5+67\times 8+9.$
\item [] $1574=1\times 2+3+4^5+67\times 8+9.$
\item [] $1575=1\times 2\times 3+4^5+67\times 8+9.$
\item [] $1576=1+2\times 3+4^5+67\times 8+9.$
\item [] $1577=1\times 2^3+4^5+67\times 8+9.$
\item [] $1578=1+2^3+4^5+67\times 8+9.$
\item [] $1579=1+234+56\times (7+8+9).$
\item [] $1580=1\times 23\times 45+67\times 8+9.$
\item [] $1581=1+23\times 45+67\times 8+9.$
\item [] $1582=12\times 3+4^5+6\times (78+9).$
\item [] $1583=1\times 2^3\times 4\times 5\times 6+7\times 89.$
\item [] $1584=1234+5+6\times 7\times 8+9.$
\item [] $1585=1^{23}\times 4\times 56\times 7+8+9.$
\item [] $1586=1234+5\times 67+8+9.$
\item [] $1587=1\times 2+3+4^5+(6+7\times 8)\times 9.$
\item [] $1588=1^2\times 3+4\times 56\times 7+8+9.$
\item [] $1589=1^2+3+4\times 56\times 7+8+9.$
\item [] $1590=1\times 2+3+4\times 56\times 7+8+9.$
\item [] $1591=1\times 2\times 3+4\times 56\times 7+8+9.$
\item [] $1592=1+2\times 3+4\times 56\times 7+8+9.$
\item [] $1593=1+23+4^5+67\times 8+9.$
\item [] $1594=1+2^3+4\times 56\times 7+8+9.$
\item [] $1595=(1+2+3)^4+5\times 6\times 7+89.$
\item [] $1596=(1+2)^3+4^5+67\times 8+9.$
\item [] $1597=123\times 4+5\times (6+7)\times (8+9).$
\item [] $1598=(1^2+3)\times 4\times 56+78\times 9.$
\item [] $1599=1^2\times 3\times 45\times 6+789.$
\item [] $1600=12+3+4\times 56\times 7+8+9.$
\item[]$\mbox{Decreasing order}$
\item [] $1531=98\times (7+6)+5+4\times 3\times 21.$
\item [] $1532=(9+8)\times (7+6\times 5)+43\times 21.$
\item [] $1533=(9+87+6\times 5)\times 4\times 3+21.$
\item [] $1534=9+876+5^4+3+21.$
\item [] $1535=(9\times 8\times 7+65\times 4+3)\times 2+1.$
\item [] $1536=9+8+7+6\times (5+4+3)\times 21.$
\item [] $1537=9+8\times (7\times 6+5)\times 4+3+21.$
\item [] $1538=987+6+543+2\times 1.$
\item [] $1539=987+6+543+2+1.$
\item [] $1540=9\times 8\times 7+6+5+4^(3+2)+1.$
\item [] $1541=9+(8+7)\times (6\times 5+4)\times 3+2\times 1.$
\item [] $1542=9+8+76\times 5\times 4+3+2\times 1.$
\item [] $1543=9+876+5^4+32+1.$
\item [] $1544=9+8+76\times 5\times 4+3\times 2+1.$
\item [] $1545=9\times 8+7\times 6\times 5\times (4+3)+2+1.$
\item [] $1546=9+8+76\times 5\times 4+3^2\times 1.$
\item [] $1547=9+8+76\times 5\times 4+3^2+1.$
\item [] $1548=9\times 8\times (7+6+5)+4\times 3\times 21.$
\item [] $1549=(9+8)\times 76+5+4\times 3\times 21.$
\item [] $1550=98+7\times (65+4)\times 3+2+1.$
\item [] $1551=9+87+6+(5+4^3)\times 21.$
\item [] $1552=(9\times 8+76+5^4+3)\times 2\times 1.$
\item [] $1553=(98+7\times 6\times 5)\times 4+321.$
\item [] $1554=(98+7+6)\times (5+4+3+2\times 1).$
\item [] $1555=(98+7+6)\times (5+4+3+2)+1.$
\item [] $1556=9+(8+7\times 6\times 5)\times (4+3)+21.$
\item [] $1557=987+6+543+21.$
\item [] $1558=(9+8\times 76+54\times 3)\times 2\times 1.$
\item [] $1559=(9+8\times 76+54\times 3)\times 2+1.$
\item [] $1560=987+6+(5+4)\times 3\times 21.$
\item [] $1561=9+8+76\times 5\times 4+3+21.$
\item [] $1562=(9+8\times 7+65)\times 4\times 3+2\times 1.$
\item [] $1563=9+8\times 76+5^4+321.$
\item [] $1564=987+6\times (5+43)\times 2+1.$
\item [] $1565=9\times 87+65\times 4\times 3+2\times 1.$
\item [] $1566=9\times 8+7\times 6\times 5+4\times 321.$
\item [] $1567=9+(8+76\times 5)\times 4+3+2+1.$
\item [] $1568=987+65\times 4+321.$
\item [] $1569=9+8+76\times 5\times 4+32\times 1.$
\item [] $1570=9+8+76\times 5\times 4+32+1.$
\item [] $1571=9+8+7\times 6\times 5+4^3\times 21.$
\item [] $1572=987+65\times (4+3+2\times 1).$
\item [] $1573=9+876+5^4+3\times 21.$
\item [] $1574=(9+8+7)\times 65+4\times 3+2\times 1.$
\item [] $1575=(9+8+7)\times 65+4\times 3+2+1.$
\item [] $1576=(98+7)\times 6+5^4+321.$
\item [] $1577=9+8+(7+6)\times 5\times 4\times 3\times 2\times 1.$
\item [] $1578=9+87\times (6+5+4+3)+2+1.$
\item [] $1579=(9+8+7\times 6)\times 5+4\times 321.$
\item [] $1580=(98+7)\times (6+5+4)+3+2\times 1.$
\item [] $1581=(9+8\times 7+65)\times 4\times 3+21.$
\item [] $1582=9\times 87+6\times (5+4\times 32)+1.$
\item [] $1583=(9+8)\times (76+5+4\times 3)+2\times 1.$
\item [] $1584=9\times 87+65\times 4\times 3+21.$
\item [] $1585=9+8+7+65\times 4\times 3\times 2+1.$
\item [] $1586=98+(7\times 6+5\times 4)\times (3+21).$
\item [] $1587=987+6\times 5\times 4\times (3+2)\times 1.$
\item [] $1588=987+6\times 5\times 4\times (3+2)+1.$
\item [] $1589=98+7\times 6\times 5\times (4+3)+21.$
\item [] $1590=9\times 8\times (7+6+5+4)+3\times 2\times 1.$
\item [] $1591=9\times 8\times (7+6+5+4)+3\times 2+1.$
\item [] $1592=98+7\times 6\times 5+4\times 321.$
\item [] $1593=9+(8+76\times 5)\times 4+32\times 1.$
\item [] $1594=(9+8+76\times 5)\times 4+3\times 2\times 1.$
\item [] $1595=(9+8+76\times 5)\times 4+3\times 2+1.$
\item [] $1596=9\times 8\times 7+6+543\times 2\times 1.$
\item [] $1597=9\times 8+76\times 5\times 4+3+2\times 1.$
\item [] $1598=9\times 8+76\times 5\times 4+3+2+1.$
\item [] $1599=9\times 8+76\times 5\times 4+3\times 2+1.$
\item [] $1600=98\times 7+6+5+43\times 21.$
\item[]$\mbox{Increasing order}$
\item [] $1601=1\times 2+3\times 45\times 6+789.$
\item [] $1602=1+2+3\times 45\times 6+789.$
\item [] $1603=1^2\times 3+4^5+6\times (7+89).$
\item [] $1604=1^2+3+4^5+6\times (7+89).$
\item [] $1605=12\times 3+4^5+67\times 8+9.$
\item [] $1606=12\times 3^4+5+6+7\times 89.$
\item [] $1607=1+2\times 3+4^5+6\times (7+89).$
\item [] $1608=1\times 23+4\times 56\times 7+8+9.$
\item [] $1609=1+23+4\times 56\times 7+8+9.$
\item [] $1610=(12+34)\times (5+6+7+8+9).$
\item [] $1611=12+3\times 45\times 6+789.$
\item [] $1612=1+23\times 45+6\times (7+89).$
\item [] $1613=12\times (3\times 4\times 5+67)+89.$
\item [] $1614=1\times 2\times (3+4+5+6+789).$
\item [] $1615=12+3+4^5+6\times (7+89).$
\item [] $1616=1+(2\times 3^4+56)\times 7+89.$
\item [] $1617=(1+2+3\times 45)\times 6+789.$
\item [] $1618=1234+5\times (67+8)+9.$
\item [] $1619=1\times 2+3\times 4\times (56+78)+9.$
\item [] $1620=12\times (3^4+5\times 6+7+8+9).$
\item [] $1621=12\times 3+4\times 56\times 7+8+9.$
\item [] $1622=1\times 2+3\times 4\times (56+7+8\times 9).$
\item [] $1623=1\times 23+4^5+6\times (7+89).$
\item [] $1624=123+4^5+6\times 78+9.$
\item [] $1625=12\times 3^4+5\times 6+7\times 89.$
\item [] $1626=1\times 2\times 345+(6+7)\times 8\times 9.$
\item [] $1627=1+2\times 345+(6+7)\times 8\times 9.$
\item [] $1628=12\times 3^4+567+89.$
\item [] $1629=(1+23+4)\times 5\times 6+789.$
\item [] $1630=1^2+3\times (456+78+9).$
\item [] $1631=1+23\times (4+56+7)+89.$
\item [] $1632=1+2+3\times (456+78+9).$
\item [] $1633=(1+23+4)\times 56+7\times 8+9.$
\item [] $1634=(1+2\times 3+4\times 56)\times 7+8+9.$
\item [] $1635=1\times 23\times 4\times (5+6)+7\times 89.$
\item [] $1636=123\times 4+5+67\times (8+9).$
\item [] $1637=1+2\times (3+4\times 5+6+789).$
\item [] $1638=(1+2)\times (3+456+78+9).$
\item [] $1639=1+2\times (3+4)\times (5\times 6+78+9).$
\item [] $1640=1^{23}\times 4\times 56\times 7+8\times 9.$
\item [] $1641=1234+5\times 67+8\times 9.$
\item [] $1642=12+(3+4\times 5)\times 67+89.$
\item [] $1643=1234+56\times 7+8+9.$
\item [] $1644=1^2+34\times 5\times 6+7\times 89.$
\item [] $1645=1+2\times 3\times 45\times 6+7+8+9.$
\item [] $1646=1\times 2\times 3+4\times 56\times 7+8\times 9.$
\item [] $1647=1\times 234\times 5+6\times 78+9.$
\item [] $1648=1+234\times 5+6\times 78+9.$
\item [] $1649=1+2^3+4\times 56\times 7+8\times 9.$
\item [] $1650=12\times 3\times 45+6+7+8+9.$
\item [] $1651=12\times 3^4+56+7\times 89.$
\item [] $1652=123\times (4+5)+67\times 8+9.$
\item [] $1653=1^{23}\times 4^5+6+7\times 89.$
\item [] $1654=1^{23}+4^5+6+7\times 89.$
\item [] $1655=12+34\times 5\times 6+7\times 89.$
\item [] $1656=1^2\times 3+4^5+6+7\times 89.$
\item [] $1657=123+4^5+6+7\times 8\times 9.$
\item [] $1658=1234+5\times 67+89.$
\item [] $1659=1\times 2\times 3+4^5+6+7\times 89.$
\item [] $1660=1^2\times 3+4\times 56\times 7+89.$
\item [] $1661=1^2+3+4\times 56\times 7+89.$
\item [] $1662=1\times 2+3+4\times 56\times 7+89.$
\item [] $1663=1+2+3+4\times 56\times 7+89.$
\item [] $1664=12\times 3^4+5+678+9.$
\item [] $1665=1+23\times 45+6+7\times 89.$
\item [] $1666=1+2^3+4\times 56\times 7+89.$
\item [] $1667=12\times 3^4+5\times (67+8\times 9).$
\item [] $1668=12+3+4^5+6+7\times 89.$
\item [] $1669=1+2\times (34+5+6+789).$
\item [] $1670=1^2\times 34\times (5+6\times 7)+8\times 9.$
\item[]$\mbox{Decreasing order}$
\item [] $1601=9\times 8+76\times 5\times 4+3^2\times 1.$
\item [] $1602=9\times 8+76\times 5\times 4+3^2+1.$
\item [] $1603=(98+7\times 6)\times 5+43\times 21.$
\item [] $1604=98\times (7+6)+5+4+321.$
\item [] $1605=9+8\times 7\times 6+5\times 4\times 3\times 21.$
\item [] $1606=(9+8+7)\times 65+43+2+1.$
\item [] $1607=(98+7)\times (6+5+4)+32\times 1.$
\item [] $1608=987+(65+4)\times 3^2\times 1.$
\item [] $1609=987+(65+4)\times 3^2+1.$
\item [] $1610=98+7\times (65+43)\times 2\times 1.$
\item [] $1611=(9\times 8\times 7+6+5\times 4)\times 3+21.$
\item [] $1612=(9+8+76\times 5)\times 4+3+21.$
\item [] $1613=(9\times 87+6+5+4\times 3)\times 2+1.$
\item [] $1614=9\times (8+7\times 6+5+4)\times 3+21.$
\item [] $1615=98\times (7+6)+5\times 4+321.$
\item [] $1616=9\times 8+76\times 5\times 4+3+21.$
\item [] $1617=9+87\times 6+543\times 2\times 1.$
\item [] $1618=9+87\times 6+543\times 2+1.$
\item [] $1619=98\times 7+6\times 5+43\times 21.$
\item [] $1620=9+876+5\times (4+3)\times 21.$
\item [] $1621=9\times 8\times (7+6+5)+4+321.$
\item [] $1622=(9+8)\times 76+5+4+321.$
\item [] $1623=987+6+5^4+3+2\times 1.$
\item [] $1624=987+6+5^4+3+2+1.$
\item [] $1625=9\times 8+76\times 5\times 4+32+1.$
\item [] $1626=9+8\times 7+65\times 4\times 3\times 2+1.$
\item [] $1627=987+6+5^4+3^2\times 1.$
\item [] $1628=987+6+5^4+3^2+1.$
\item [] $1629=9+8+(7+6\times 5)\times 43+21.$
\item [] $1630=(9+8+7+6)\times 54+3^2+1.$
\item [] $1631=98\times 7+(6+5+4)\times 3\times 21.$
\item [] $1632=987+6\times 54+321.$
\item [] $1633=987+6+5\times 4\times 32\times 1.$
\item [] $1634=987+6+5\times 4\times 32+1.$
\item [] $1635=9+87\times (6+5+4)+321.$
\item [] $1636=9+8\times 7\times (6+5\times 4+3)+2+1.$
\item [] $1637=987+65\times (4+3+2+1).$
\item [] $1638=98\times 7+6+5^4+321.$
\item [] $1639=9\times 8+7+65\times 4\times 3\times 2\times 1.$
\item [] $1640=9\times 8+7+65\times 4\times 3\times 2+1.$
\item [] $1641=(9+87\times 6+5+4)\times 3+21.$
\item [] $1642=98+76\times 5\times 4+3+21.$
\item [] $1643=(9\times 8+76)\times 5+43\times 21.$
\item [] $1644=(9+8+7+6)\times 54+3+21.$
\item [] $1645=(9+8+7)\times 65+4^3+21.$
\item [] $1646=987+654+3+2\times 1.$
\item [] $1647=987+654+3+2+1.$
\item [] $1648=987+654+3\times 2+1.$
\item [] $1649=98\times (7+6)+54+321.$
\item [] $1650=98+76\times 5\times 4+32\times 1.$
\item [] $1651=987+654+3^2+1.$
\item [] $1652=98+7\times 6\times 5+4^3\times 21.$
\item [] $1653=(9+8+7+6)\times 54+32+1.$
\item [] $1654=98\times 7+65+43\times 21.$
\item [] $1655=9\times 8+76\times 5\times 4+3\times 21.$
\item [] $1656=9+87+65\times 4\times 3\times 2\times 1.$
\item [] $1657=9+87+65\times 4\times 3\times 2+1.$
\item [] $1658=98+(7+6)\times 5\times 4\times 3\times 2\times 1.$
\item [] $1659=98+(7+6)\times 5\times 4\times 3\times 2+1.$
\item [] $1660=98\times 7+6\times 54\times 3+2\times 1.$
\item [] $1661=98\times 7+654+321.$
\item [] $1662=(9+87)\times 6+543\times 2\times 1.$
\item [] $1663=(9+87)\times 6+543\times 2+1.$
\item [] $1664=(9+8\times 76+5\times 43)\times 2\times 1.$
\item [] $1665=987+654+3+21.$
\item [] $1666=98+7+65\times 4\times 3\times 2+1.$
\item [] $1667=(9+8)\times 76+54+321.$
\item [] $1668=9\times 8+7\times 6\times (5+4\times 3+21).$
\item [] $1669=9+8\times (7\times 6+5)+4\times 321.$
\item [] $1670=98\times 7+6\times (54\times 3+2)\times 1.$
\item[]$\mbox{Increasing order}$
\item [] $1671=(12+3\times 45)\times 6+789.$
\item [] $1672=12+3+4\times 56\times 7+89.$
\item [] $1673=1\times 23\times 4\times (5+6+7)+8+9.$
\item [] $1674=1+23\times 4\times (5+6+7)+8+9.$
\item [] $1675=(1\times 2+3+4\times 56)\times 7+8\times 9.$
\item [] $1676=1\times 23+4^5+6+7\times 89.$
\item [] $1677=1+23+4^5+6+7\times 89.$
\item [] $1678=1+23\times (45+6)+7\times 8\times 9.$
\item [] $1679=12\times 3\times 45+6\times 7+8+9.$
\item [] $1680=1\times 23+4\times 56\times 7+89.$
\item [] $1681=1+23+4\times 56\times 7+89.$
\item [] $1682=12+34\times (5+6\times 7)+8\times 9.$
\item [] $1683=(12\times 3\times 4+5)\times 6+789.$
\item [] $1684=(1+2)^3+4\times 56\times 7+89.$
\item [] $1685=12\times 3^4+5+6+78\times 9.$
\item [] $1686=1+2\times 3\times 45\times 6+7\times 8+9.$
\item [] $1687=1+2\times (3+45+6+789).$
\item [] $1688=1\times 2^3+4\times 5\times (67+8+9).$
\item [] $1689=12\times 3+4^5+6+7\times 89.$
\item [] $1690=1^{234}+5\times 6\times 7\times 8+9.$
\item [] $1691=12\times 3\times 45+6+7\times 8+9.$
\item [] $1692=123+4^5+67\times 8+9.$
\item [] $1693=12\times 3+4\times 56\times 7+89.$
\item [] $1694=1^{23}+4+5\times 6\times 7\times 8+9.$
\item [] $1695=(12+3)\times (4\times 5+6+78+9).$
\item [] $1696=1^2\times 3+4+5\times 6\times 7\times 8+9.$
\item [] $1697=1^2+3+4+5\times 6\times 7\times 8+9.$
\item [] $1698=1234+56\times 7+8\times 9.$
\item [] $1699=1\times 2\times 3+4+5\times 6\times 7\times 8+9.$
\item [] $1700=1+2\times 3+4+5\times 6\times 7\times 8+9.$
\item [] $1701=1\times 2^3+4+5\times 6\times 7\times 8+9.$
\item [] $1702=1+2^3+4+5\times 6\times 7\times 8+9.$
\item [] $1703=1\times 2+3\times 4+5\times 6\times 7\times 8+9.$
\item [] $1704=12\times 3^4+5\times 6+78\times 9.$
\item [] $1705=12\times 3\times 45+6+7+8\times 9.$
\item [] $1706=1\times 2\times (34+5\times 6+789).$
\item [] $1707=1\times 2\times 3\times 45\times 6+78+9.$
\item [] $1708=123+4\times 56\times 7+8+9.$
\item [] $1709=(1+23)\times 45+6+7\times 89.$
\item [] $1710=(123+45)\times 6+78\times 9.$
\item [] $1711=1^{23}\times 4^5+678+9.$
\item [] $1712=1^{23}+4^5+678+9.$
\item [] $1713=12\times 3\times 45+6+78+9.$
\item [] $1714=1+(2^3+4\times 56)\times 7+89.$
\item [] $1715=1234+56\times 7+89.$
\item [] $1716=1234+5+6\times 78+9.$
\item [] $1717=1+2+3+4^5+678+9.$
\item [] $1718=1+2\times 3+4^5+678+9.$
\item [] $1719=1\times 2^3+4^5+678+9.$
\item [] $1720=1+2^3+4^5+678+9.$
\item [] $1721=1\times 2^3\times 4+5\times 6\times 7\times 8+9.$
\item [] $1722=12\times 3\times 45+6+7+89.$
\item [] $1723=1+23\times 45+678+9.$
\item [] $1724=1^2+34+5\times 6\times 7\times 8+9.$
\item [] $1725=1\times 2+34+5\times 6\times 7\times 8+9.$
\item [] $1726=12+3+4^5+678+9.$
\item [] $1727=12\times (3^4+5+6)+7\times 89.$
\item [] $1728=12+3\times 4\times (56+78+9).$
\item [] $1729=12\times 3+4+5\times 6\times 7\times 8+9.$
\item [] $1730=1+2+(34+5)\times 6\times 7+89.$
\item [] $1731=12\times (3^4+56)+78+9.$
\item [] $1732=1^{23}\times 4^5+6+78\times 9.$
\item [] $1733=1^{23}+4^5+6+78\times 9.$
\item [] $1734=1\times 23+4^5+678+9.$
\item [] $1735=1+23+4^5+678+9.$
\item [] $1736=1^2+3+4^5+6+78\times 9.$
\item [] $1737=1\times 2+3+4^5+6+78\times 9.$
\item [] $1738=1\times 2\times 3+4^5+6+78\times 9.$
\item [] $1739=1+2\times 3+4^5+6+78\times 9.$
\item [] $1740=1\times 2^3+4^5+6+78\times 9.$
\item[]$\mbox{Decreasing order}$
\item [] $1671=987+6\times (54+3)\times 2\times 1.$
\item [] $1672=987+6\times (54+3)\times 2+1.$
\item [] $1673=987+654+32\times 1.$
\item [] $1674=987+654+32+1.$
\item [] $1675=(9\times 87+6+5+43)\times 2+1.$
\item [] $1676=98\times 7+6\times (54\times 3+2+1).$
\item [] $1677=987+(65+4)\times (3^2+1).$
\item [] $1678=(9\times (8+7)+(6+5)\times 4^3)\times 2\times 1.$
\item [] $1679=98\times 7+6\times 54\times 3+21.$
\item [] $1680=(9+8)\times 7+65\times 4\times 3\times 2+1.$
\item [] $1681=9+8+76\times 5+4\times 321.$
\item [] $1682=9+8\times 7+(65+4\times 3)\times 21.$
\item [] $1683=(9+8+7+6)\times 54+3\times 21.$
\item [] $1684=9\times 8+(7+6\times 5)\times 43+21.$
\item [] $1685=9+8+765+43\times 21.$
\item [] $1686=(9+8)\times (76+5\times 4+3)+2+1.$
\item [] $1687=(9\times 8+7\times 6+(5+4)^3)\times 2+1.$
\item [] $1688=(9+8+7)\times 65+4\times 32\times 1.$
\item [] $1689=(9+8+7)\times 65+4\times 32+1.$
\item [] $1690=9+8\times 7+65\times (4\times 3\times 2+1).$
\item [] $1691=98+(7+6\times 5)\times 43+2\times 1.$
\item [] $1692=(98+7\times 65+4)\times 3+21.$
\item [] $1693=9\times 87+65\times (4+3)\times 2\times 1.$
\item [] $1694=9+8\times 7\times 6+5+4^3\times 21.$
\item [] $1695=9\times (8+7)+65\times 4\times 3\times 2\times 1.$
\item [] $1696=9+(8+76)\times 5\times 4+3\times 2+1.$
\item [] $1697=9\times 87+6+5+43\times 21.$
\item [] $1698=9+8\times 7\times 6\times 5+4+3+2\times 1.$
\item [] $1699=9+8\times 7\times 6\times 5+4+3+2+1.$
\item [] $1700=9+8\times 7\times 6\times 5+4+3\times 2+1.$
\item [] $1701=98+(7+65\times 4)\times 3\times 2+1.$
\item [] $1702=9+8\times 7\times 6\times 5+4+3^2\times 1.$
\item [] $1703=9+8\times 76+543\times 2\times 1.$
\item [] $1704=987+654+3\times 21.$
\item [] $1705=98\times (7+6)+5\times 43\times 2+1.$
\item [] $1706=(9\times 8+7+6)\times 5\times 4+3+2+1.$
\item [] $1707=987+6\times 5\times 4\times 3\times 2\times 1.$
\item [] $1708=987+6\times 5\times 4\times 3\times 2+1.$
\item [] $1709=(9\times 8+7+6)\times 5+4\times 321.$
\item [] $1710=(9\times 8+7+6)\times 5\times 4+3^2+1.$
\item [] $1711=(9\times 87+65+4+3)\times 2+1.$
\item [] $1712=98\times (7+6)+5+432+1.$
\item [] $1713=9+(8+76)\times 5+4\times 321.$
\item [] $1714=9+8\times 7\times 6\times 5+4\times 3\times 2+1.$
\item [] $1715=98+7\times (6+5)\times (4+3)\times (2+1).$
\item [] $1716=9\times 87+6\times 5+43\times 21.$
\item [] $1717=(98+7)\times 6+543\times 2+1.$
\item [] $1718=9\times (8+7)\times 6+5+43\times 21.$
\item [] $1719=(9+8+7)\times (65+4)+3\times 21.$
\item [] $1720=(9\times 87+65+4\times 3)\times 2\times 1.$
\item [] $1721=9+(8+76)\times 5\times 4+32\times 1.$
\item [] $1722=9+8\times 7\times 6\times 5+4\times 3+21.$
\item [] $1723=(9+8)\times 76+5\times 43\times 2+1.$
\item [] $1724=(9\times 8+7+6)\times 5\times 4+3+21.$
\item [] $1725=9+8\times 7\times 6\times 5+4+32\times 1.$
\item [] $1726=9+8\times 7\times 6\times 5+4+32+1.$
\item [] $1727=9+8\times (7\times 6\times 5+4)+3+2+1.$
\item [] $1728=987+6+5\times (4+3)\times 21.$
\item [] $1729=9\times 87+(6+5)\times 43\times 2\times 1.$
\item [] $1730=(9+8)\times 76+5+432+1.$
\item [] $1731=98+(7\times 6+5+4)\times 32+1.$
\item [] $1732=(9\times 8+7+6)\times 5\times 4+32\times 1.$
\item [] $1733=(9\times 8+7+6)\times 5\times 4+32+1.$
\item [] $1734=9+8\times 7\times 6\times 5+43+2\times 1.$
\item [] $1735=9+8\times 7\times 6\times 5+43+2+1.$
\item [] $1736=9\times 8+76\times 5+4\times 321.$
\item [] $1737=9\times (8+7)\times (6+5)+4\times 3\times 21.$
\item [] $1738=(9+8\times 7\times 6)\times 5+4+3^2\times 1.$
\item [] $1739=9+(87\times 6+54)\times 3+2\times 1.$
\item [] $1740=9\times 8+765+43\times 21.$
\item[]$\mbox{Increasing order}$
\item [] $1741=(1+2)\times 3+4^5+6+78\times 9.$
\item [] $1742=12+3\times (4+567)+8+9.$
\item [] $1743=1\times 23\times 45+6+78\times 9.$
\item [] $1744=1+23\times 45+6+78\times 9.$
\item [] $1745=12\times (3^4+56+7)+8+9.$
\item [] $1746=(123+4+5+6+7\times 8)\times 9.$
\item [] $1747=12+3+4^5+6+78\times 9.$
\item [] $1748=1+(2+3+4\times 5)\times 67+8\times 9.$
\item [] $1749=1234+5+6+7\times 8\times 9.$
\item [] $1750=1+2^3\times 4\times 5\times 6+789.$
\item [] $1751=12\times 3\times 45+6\times 7+89.$
\item [] $1752=12\times (3+4\times 5)\times 6+7+89.$
\item [] $1753=12+3+4^5+6\times 7\times (8+9).$
\item [] $1754=(1+2)\times (3\times 4+567)+8+9.$
\item [] $1755=1\times 23+4^5+6+78\times 9.$
\item [] $1756=1+23+4^5+6+78\times 9.$
\item [] $1757=1\times 2\times 34+5\times 6\times 7\times 8+9.$
\item [] $1758=1+2\times 34+5\times 6\times 7\times 8+9.$
\item [] $1759=12\times 3\times 45+67+8\times 9.$
\item [] $1760=(1+2)^3\times 45+67\times 8+9.$
\item [] $1761=(1+2+3)^4+5\times (6+78+9).$
\item [] $1762=(12+3+4\times 56)\times 7+89.$
\item [] $1763=123+4\times 56\times 7+8\times 9.$
\item [] $1764=12\times 3\times 45+6\times (7+8+9).$
\item [] $1765=1+(2+3+4\times 5)\times 67+89.$
\item [] $1766=(12+3+4)\times 56+78\times 9.$
\item [] $1767=12^3+4+5+6+7+8+9.$
\item [] $1768=1234+5\times 6+7\times 8\times 9.$
\item [] $1769=(1^2+3+4)\times 5\times 6\times 7+89.$
\item [] $1770=1^2\times 3^4+5\times 6\times 7\times 8+9.$
\item [] $1771=1^2+3^4+5\times 6\times 7\times 8+9.$
\item [] $1772=12\times 3^4+5+6+789.$
\item [] $1773=1+2+3^4+5\times 6\times 7\times 8+9.$
\item [] $1774=1\times 2+3\times 45\times (6+7)+8+9.$
\item [] $1775=1+2+3\times 45\times (6+7)+8+9.$
\item [] $1776=12\times 3\times 45+67+89.$
\item [] $1777=(1+2+3)^4+56\times 7+89.$
\item [] $1778=12^3+4\times 5+6+7+8+9.$
\item [] $1779=12\times (3\times 45+6)+78+9.$
\item [] $1780=123+4\times 56\times 7+89.$
\item [] $1781=1\times 23\times 4+5\times 6\times 7\times 8+9.$
\item [] $1782=1+23\times 4+5\times 6\times 7\times 8+9.$
\item [] $1783=1+(2+34)\times 5\times 6+78\times 9.$
\item [] $1784=1234+5+67\times 8+9.$
\item [] $1785=(1+23)\times 4+5\times 6\times 7\times 8+9.$
\item [] $1786=12^3+4+5\times 6+7+8+9.$
\item [] $1787=12\times 3^4+5+6\times (7+8)\times 9.$
\item [] $1788=12\times (3\times 45+6)+7+89.$
\item [] $1789=(1+2)\times 3+4^5+(6+78)\times 9.$
\item [] $1790=(1^2+3)^4\times 5+6+7\times 8\times 9.$
\item [] $1791=12\times 3^4+5\times 6+789.$
\item [] $1792=1+(2\times 3^4+5)\times 6+789.$
\item [] $1793=12\times (3\times 45+6+7)+8+9.$
\item [] $1794=1234+56+7\times 8\times 9.$
\item [] $1795=1+2\times (3\times 45\times 6+78+9).$
\item [] $1796=12^3+4+5+6\times 7+8+9.$
\item [] $1797=1234+5+(6+7\times 8)\times 9.$
\item [] $1798=(1\times 2+3\times 45)\times (6+7)+8+9.$
\item [] $1799=1\times 234\times 5+6+7\times 89.$
\item [] $1800=1+234\times 5+6+7\times 89.$
\item [] $1801=1\times 23\times 4\times (5+6)+789.$
\item [] $1802=1^2\times 3\times (4+567)+89.$
\item [] $1803=12^3+45+6+7+8+9.$
\item [] $1804=(1+234)\times 5+6+7\times 89.$
\item [] $1805=1+2+3\times (4+567)+89.$
\item [] $1806=1+2+(3\times 4+5+6)\times 78+9.$
\item [] $1807=12^3+4\times 5+6\times 7+8+9.$
\item [] $1808=12^3+4+5+6+7\times 8+9.$
\item [] $1809=1^2\times 34\times 5\times 6+789.$
\item [] $1810=1^2+34\times 5\times 6+789.$
\item[]$\mbox{Decreasing order}$
\item [] $1741=9+8+76\times 5+4^3\times 21.$
\item [] $1742=(98+7\times 6+5)\times 4\times 3+2\times 1.$
\item [] $1743=9+(8+7)\times 6\times 5+4\times 321.$
\item [] $1744=9+(87+65\times 4)\times (3+2)\times 1.$
\item [] $1745=9+8+(7+65)\times 4\times 3\times 2\times 1.$
\item [] $1746=9+8+(7+65)\times 4\times 3\times 2+1.$
\item [] $1747=(98+7\times 6+54)\times 3^2+1.$
\item [] $1748=9+8\times (7\times 6\times 5+4+3)+2+1.$
\item [] $1749=(9+8+76)\times 5+4\times 321.$
\item [] $1750=(98+765+4\times 3)\times 2\times 1.$
\item [] $1751=9\times 87+65+43\times 21.$
\item [] $1752=9+(8+76)\times 5\times 4+3\times 21.$
\item [] $1753=9+8\times 7\times 6\times 5+43+21.$
\item [] $1754=9+8\times (7\times 6\times 5+4)+32+1.$
\item [] $1755=9+8\times 7\times 6\times 5+4^3+2\times 1.$
\item [] $1756=9+8+7\times 65+4\times 321.$
\item [] $1757=9\times 87+6\times 54\times 3+2\times 1.$
\item [] $1758=9\times 87+654+321.$
\item [] $1759=9+8+7+6+54\times 32+1.$
\item [] $1760=(9\times 8+765+43)\times 2\times 1.$
\item [] $1761=98\times (7+6)+54\times 3^2+1.$
\item [] $1762=98+76\times 5+4\times 321.$
\item [] $1763=(9\times 8+7)\times 6+5+4\times 321.$
\item [] $1764=987+(6\times 5+4+3)\times 21.$
\item [] $1765=(9+8+7+6+5+4+3)^2+1.$
\item [] $1766=98+765+43\times 21.$
\item [] $1767=98\times 7+6\times 5\times 4\times 3^2+1.$
\item [] $1768=9\times 87+6\times (54\times 3+2)+1.$
\item [] $1769=987+65\times 4\times 3+2\times 1.$
\item [] $1770=987+65\times 4\times 3+2+1.$
\item [] $1771=(9+8\times 7\times 6)\times 5+43+2+1.$
\item [] $1772=(9\times 8+7)\times (6+5)+43\times 21.$
\item [] $1773=(98+76)\times 5+43\times 21.$
\item [] $1774=9+8\times 7\times 6\times 5+4^3+21.$
\item [] $1775=9+8\times 7\times 6\times 5+43\times 2\times 1.$
\item [] $1776=9\times 87+6\times 54\times 3+21.$
\item [] $1777=9\times 8+76+543\times (2+1).$
\item [] $1778=98\times 7+6+543\times 2\times 1.$
\item [] $1779=98\times 7+6+543\times 2+1.$
\item [] $1780=(9+8+7+65)\times 4\times (3+2)\times 1.$
\item [] $1781=9+8+7\times 6\times (5+4+32+1).$
\item [] $1782=(9\times 87+65+43)\times 2\times 1.$
\item [] $1783=(9\times 87+65+43)\times 2+1.$
\item [] $1784=9+8\times (7\times 6\times 5+4)+3\times 21.$
\item [] $1785=98\times 7+(6+543)\times 2+1.$
\item [] $1786=9+8+76\times (5\times 4+3)+21.$
\item [] $1787=9+8+7\times 6+54\times 32\times 1.$
\item [] $1788=987+65\times 4\times 3+21.$
\item [] $1789=98\times (7+6+5)+4\times 3\times 2+1.$
\item [] $1790=98+(7\times 6+5)\times 4\times 3^2\times 1.$
\item [] $1791=9+87\times 6+5\times 4\times 3\times 21.$
\item [] $1792=98+7\times (6+5\times 43+21).$
\item [] $1793=9+876+5+43\times 21.$
\item [] $1794=9+8\times 7+6\times (5+4)\times 32+1.$
\item [] $1795=(9+876+5+4+3)\times 2+1.$
\item [] $1796=9\times 8\times 7+6\times 5\times 43+2\times 1.$
\item [] $1797=9\times 8\times 7+6\times 5\times 43+2+1.$
\item [] $1798=9+8\times (7+6)\times (5+4\times 3)+21.$
\item [] $1799=9+8\times 7+6+54\times 32\times 1.$
\item [] $1800=9+8\times 7+6+54\times 32+1.$
\item [] $1801=(9+8)\times 7\times 6+543\times 2+1.$
\item [] $1802=(9+8+7+6)\times 5\times 4\times 3+2\times 1.$
\item [] $1803=9\times 8\times (7+6+5+4+3)+2+1.$
\item [] $1804=9\times 8\times 7+65\times 4\times (3+2)\times 1.$
\item [] $1805=(9+876+5+4\times 3)\times 2+1.$
\item [] $1806=9+8\times (7\times 6\times 5+4\times 3)+21.$
\item [] $1807=9\times 8\times 7+6\times (5\times 43+2)+1.$
\item [] $1808=9\times 8+7+6\times (5+4)\times 32+1.$
\item [] $1809=98\times (7+6+5)+43+2\times 1.$
\item [] $1810=98\times (7+6+5)+43+2+1.$
\item[]$\mbox{Increasing order}$
\item [] $1811=1\times 2+34\times 5\times 6+789.$
\item [] $1812=12^3+4+56+7+8+9.$
\item [] $1813=1+2\times (3\times 45\times 6+7+89).$
\item [] $1814=12+3\times (4+567)+89.$
\item [] $1815=1234+5+6\times (7+89).$
\item [] $1816=123+4+5\times 6\times 7\times 8+9.$
\item [] $1817=12\times 3^4+56+789.$
\item [] $1818=1234+567+8+9.$
\item [] $1819=12^3+4\times 5+6+7\times 8+9.$
\item [] $1820=1^{23}+4^5+6+789.$
\item [] $1821=12+34\times 5\times 6+789.$
\item [] $1822=12^3+4+5+6+7+8\times 9.$
\item [] $1823=1^2+3+4^5+6+789.$
\item [] $1824=1\times 2+3+4^5+6+789.$
\item [] $1825=1\times 2\times 3+4^5+6+789.$
\item [] $1826=1+2\times 3+4^5+6+789.$
\item [] $1827=12^3+4+5\times 6+7\times 8+9.$
\item [] $1828=1+2^3+4^5+6+789.$
\item [] $1829=1\times 2\times 345+67\times (8+9).$
\item [] $1830=12^3+4+5+6+78+9.$
\item [] $1831=1+23\times 45+6+789.$
\item [] $1832=12^3+45+6\times 7+8+9.$
\item [] $1833=12^3+4\times 5+6+7+8\times 9.$
\item [] $1834=123+4^5+678+9.$
\item [] $1835=1^{23}+4^5+6\times (7+8)\times 9.$
\item [] $1836=123\times 4+56\times (7+8+9).$
\item [] $1837=12^3+4\times (5+6)+7\times 8+9.$
\item [] $1838=12+3\times (4+5)\times 67+8+9.$
\item [] $1839=12^3+4+5+6+7+89.$
\item [] $1840=1\times 23\times (4+5+6+7\times 8+9).$
\item [] $1841=12^3+4+5\times 6+7+8\times 9.$
\item [] $1842=1\times 23+4^5+6+789.$
\item [] $1843=1+23+4^5+6+789.$
\item [] $1844=12^3+45+6+7\times 8+9.$
\item [] $1845=12+3\times (4+5+67)\times 8+9.$
\item [] $1846=1+234\times 5+(67+8)\times 9.$
\item [] $1847=1+2+3\times 45\times (6+7)+89.$
\item [] $1848=1\times 2\times 3\times 4\times 56+7\times 8\times 9.$
\item [] $1849=12^3+4+5\times 6+78+9.$
\item [] $1850=12^3+4\times 5+6+7+89.$
\item [] $1851=12^3+4+5+6\times 7+8\times 9.$
\item [] $1852=1+2\times 3^4+5\times 6\times 7\times 8+9.$
\item [] $1853=12^3+4+56+7\times 8+9.$
\item [] $1854=1\times (2+3)\times 45\times 6+7\times 8\times 9.$
\item [] $1855=12\times 3+4^5+6+789.$
\item [] $1856=12+3\times 45\times (6+7)+89.$
\item [] $1857=12^3+45+67+8+9.$
\item [] $1858=12^3+45+6+7+8\times 9.$
\item [] $1859=1234+(5+6)\times 7\times 8+9.$
\item [] $1860=12\times (3\times 4+56+78+9).$
\item [] $1861=1+2\times (3\times 45+6+789).$
\item [] $1862=(1+234)\times 5+678+9.$
\item [] $1863=1+2\times 3+4\times (56\times 7+8\times 9).$
\item [] $1864=1\times 2^3+4\times (56\times 7+8\times 9).$
\item [] $1865=1\times 23\times (4+5)\times 6+7\times 89.$
\item [] $1866=12^3+45+6+78+9.$
\item [] $1867=12^3+4+56+7+8\times 9.$
\item [] $1868=1234+5+6+7\times 89.$
\item [] $1869=12\times (34+56)+789.$
\item [] $1870=1+(2+34)\times 5\times 6+789.$
\item [] $1871=1\times 2^3\times 4\times 56+7+8\times 9.$
\item [] $1872=1^2\times 3\times 456+7\times 8\times 9.$
\item [] $1873=1^2+3\times 456+7\times 8\times 9.$
\item [] $1874=1\times 2+3\times 456+7\times 8\times 9.$
\item [] $1875=12^3+45+6+7+89.$
\item [] $1876=12^3+4+5+67+8\times 9.$
\item [] $1877=12\times 3\times 4\times 5+(6+7)\times 89.$
\item [] $1878=1\times 234\times 5+6+78\times 9.$
\item [] $1879=1+234\times 5+6+78\times 9.$
\item [] $1880=1+2^3\times 4\times 56+78+9.$
\item[]$\mbox{Decreasing order}$
\item [] $1811=9\times 8+7\times 65+4\times 321.$
\item [] $1812=(98+765+43)\times 2\times 1.$
\item [] $1813=9\times 8+7+6+54\times 32\times 1.$
\item [] $1814=9\times 8+7+6+54\times 32+1.$
\item [] $1815=9\times 8\times 7+6\times 5\times 43+21.$
\item [] $1816=9+8+7\times 65+4^3\times 21.$
\item [] $1817=9+8\times 7\times 6\times 5+4\times 32\times 1.$
\item [] $1818=9\times 8\times 7+6\times 5+4\times 321.$
\item [] $1819=98\times (7+6)+543+2\times 1.$
\item [] $1820=9+87\times 6+5+4\times 321.$
\item [] $1821=9+8+76+54\times 32\times 1.$
\item [] $1822=98+76\times 5+4^3\times 21.$
\item [] $1823=9\times 8+76\times (5\times 4+3)+2+1.$
\item [] $1824=9+87+6\times (5+4)\times 32\times 1.$
\item [] $1825=9+87+6\times (5+4)\times 32+1.$
\item [] $1826=98+(7+65)\times 4\times 3\times 2\times 1.$
\item [] $1827=9+(8+7)\times 6+54\times 32\times 1.$
\item [] $1828=9+(8+7)\times 6+54\times 32+1.$
\item [] $1829=9\times 8+7\times 6+5\times (4+3)^{(2+1)}.$
\item [] $1830=9+87+6+54\times 32\times 1.$
\item [] $1831=9+876+5^4+321.$
\item [] $1832=(9\times 8+76\times 5)\times 4+3+21.$
\item [] $1833=98+7+6\times (5+4)\times 32\times 1.$
\item [] $1834=98+7+6\times (5+4)\times 32+1.$
\item [] $1835=9+(87+65)\times 4\times 3+2\times 1.$
\item [] $1836=(9+87)\times 6+5\times 4\times 3\times 21.$
\item [] $1837=98+7\times 65+4\times 321.$
\item [] $1838=98\times (7+6)+543+21.$
\item [] $1839=98+7+6+54\times 32\times 1.$
\item [] $1840=98+7+6+54\times 32+1.$
\item [] $1841=9+8\times (7+6)+54\times 32\times 1.$
\item [] $1842=9\times 8+7\times 6+54\times 32\times 1.$
\item [] $1843=9\times 8+7\times 6+54\times 32+1.$
\item [] $1844=9+8+7\times 65\times 4+3\times 2+1.$
\item [] $1845=9+876+5\times 4^3\times (2+1).$
\item [] $1846=9+8+7\times 65\times 4+3^2\times 1.$
\item [] $1847=9+8+7\times 65\times 4+3^2+1.$
\item [] $1848=9\times (8+7+6)\times 5+43\times 21.$
\item [] $1849=9+8\times (7+65+43)\times 2\times 1.$
\item [] $1850=98\times (7+6+5)+43\times 2\times 1.$
\item [] $1851=98\times (7+6+5)+43\times 2+1.$
\item [] $1852=987+6\times (5+4+3)^2+1.$
\item [] $1853=9\times 8\times 7+65+4\times 321.$
\item [] $1854=9+(87+65)\times 4\times 3+21.$
\item [] $1855=(98+765+4^3)\times 2+1.$
\item [] $1856=(9+8)\times 76+543+21.$
\item [] $1857=(9\times 8+76+5)\times 4\times 3+21.$
\item [] $1858=9+8+76\times 5\times 4+321.$
\item [] $1859=9\times 8\times 7+6+5+4^3\times 21.$
\item [] $1860=9+8+7\times (65\times 4+3)+2\times 1.$
\item [] $1861=9+8+7\times 65\times 4+3+21.$
\item [] $1862=9+(8\times 76+5+4)\times 3+2\times 1.$
\item [] $1863=9\times 87+6\times 5\times 4\times 3^2\times 1.$
\item [] $1864=9\times 87+6\times 5\times 4\times 3^2+1.$
\item [] $1865=(9+87)\times 6+5+4\times 321.$
\item [] $1866=(9+876+5+43)\times 2\times 1.$
\item [] $1867=98+76\times (5\times 4+3)+21.$
\item [] $1868=98+7\times 6+54\times 32\times 1.$
\item [] $1869=98+7\times 6+54\times 32+1.$
\item [] $1870=9+8+7\times 65\times 4+32+1.$
\item [] $1871=9\times 8+7\times 65+4^3\times 21.$
\item [] $1872=(9+8+7)\times 6+54\times 32\times 1.$
\item [] $1873=(9+8+7)\times 6+54\times 32+1.$
\item [] $1874=9+(87+6)\times 5\times 4+3+2\times 1.$
\item [] $1875=9\times 87+6+543\times 2\times 1.$
\item [] $1876=9\times 8+76+54\times 32\times 1.$
\item [] $1877=9+8\times 76+5\times 4\times 3\times 21.$
\item [] $1878=9\times 8\times 7+6\times 5+4^3\times 21.$
\item [] $1879=9+8+7\times (65\times 4+3)+21.$
\item [] $1880=9+87\times 6+5+4^3\times 21.$
\item[]$\mbox{Increasing order}$
\item [] $1881=(12+34\times 5)\times 6+789.$
\item [] $1882=12^3+4+5\times (6+7+8+9).$
\item [] $1883=1\times 23+4\times 5\times (6+78+9).$
\item [] $1884=12+3\times 456+7\times 8\times 9.$
\item [] $1885=1+234\times 5+6\times 7\times (8+9).$
\item [] $1886=12^3+4+5\times (6+7)+89.$
\item [] $1887=1234+5\times 6+7\times 89.$
\item [] $1888=1\times 2^3\times 4\times 56+7+89.$
\item [] $1889=1+2^3\times 4\times 56+7+89.$
\item [] $1890=1234+567+89.$
\item [] $1891=12+34\times (5+6\times 7+8)+9.$
\item [] $1892=12\times 3+4\times (56\times 7+8\times 9).$
\item [] $1893=12^3+4+5+67+89.$
\item [] $1894=1+(2^3+4\times 5)\times 67+8+9.$
\item [] $1895=1\times (23+4\times 5)\times 6\times 7+89.$
\item [] $1896=(12+3)\times 4\times 5\times 6+7+89.$
\item [] $1897=12^3+4+5\times 6+(7+8)\times 9.$
\item [] $1898=12^3+45+6+7\times (8+9).$
\item [] $1899=(1+2+34)\times 5\times 6+789.$
\item [] $1900=1\times 2+3\times (4+5)\times 67+89.$
\item [] $1901=1+2+3\times (4+5)\times 67+89.$
\item [] $1902=123\times (4+5)+6+789.$
\item [] $1903=123+4^5+(6+78)\times 9.$
\item [] $1904=12^3+45+6\times 7+89.$
\item [] $1905=1^2\times (34\times 5+67)\times 8+9.$
\item [] $1906=1\times 2+3+4\times (5+6\times 78)+9.$
\item [] $1907=1^{23}\times 45\times 6\times 7+8+9.$
\item [] $1908=1^{23}+45\times 6\times 7+8+9.$
\item [] $1909=1\times (2+3)\times 4\times 56+789.$
\item [] $1910=1+(2+3)\times 4\times 56+789.$
\item [] $1911=1^2+3+45\times 6\times 7+8+9.$
\item [] $1912=12^3+45+67+8\times 9.$
\item [] $1913=1234+56+7\times 89.$
\item [] $1914=1+2\times 3+45\times 6\times 7+8+9.$
\item [] $1915=1\times 2^3+45\times 6\times 7+8+9.$
\item [] $1916=1+2^3+45\times 6\times 7+8+9.$
\item [] $1917=12^3+45+6\times (7+8+9).$
\item [] $1918=1\times 2+3+(4+5\times 6)\times 7\times 8+9.$
\item [] $1919=1\times 2\times 3+(4+5\times 6)\times 7\times 8+9.$
\item [] $1920=1\times 2\times 3\times 4\times (5+6)\times 7+8\times 9.$
\item [] $1921=1+2\times 3\times 4\times (56+7+8+9).$
\item [] $1922=12+3+45\times 6\times 7+8+9.$
\item [] $1923=1\times 234+5\times 6\times 7\times 8+9.$
\item [] $1924=1+234+5\times 6\times 7\times 8+9.$
\item [] $1925=12\times (3^4+5+67)+89.$
\item [] $1926=1234+5+678+9.$
\item [] $1927=12^3+4\times 5\times 6+7+8\times 9.$
\item [] $1928=1^2\times 34\times 56+7+8+9.$
\item [] $1929=12^3+45+67+89.$
\item [] $1930=1\times 23+45\times 6\times 7+8+9.$
\item [] $1931=1+23+45\times 6\times 7+8+9.$
\item [] $1932=1\times 23\times (4+56+7+8+9).$
\item [] $1933=1+23\times (4+56+7+8+9).$
\item [] $1934=(1+2)^3+45\times 6\times 7+8+9.$
\item [] $1935=12^3+4\times 5\times 6+78+9.$
\item [] $1936=1\times 2^3+4\times (5+6\times 78+9).$
\item [] $1937=(1\times 234+5\times 6)\times 7+89.$
\item [] $1938=1+(234+5\times 6)\times 7+89.$
\item [] $1939=1+(234+5)\times 6+7\times 8\times 9.$
\item [] $1940=12+34\times 56+7+8+9.$
\item [] $1941=12+3\times (4+567+8\times 9).$
\item [] $1942=123+4^5+6+789.$
\item [] $1943=12\times 3+45\times 6\times 7+8+9.$
\item [] $1944=12^3+4\times 5\times 6+7+89.$
\item [] $1945=1+(2+34)\times (5\times 6+7+8+9).$
\item [] $1946=(1\times 2+3+4+5)\times (67+8\times 9).$
\item [] $1947=1234+5+6+78\times 9.$
\item [] $1948=1+23+4\times (56\times 7+89).$
\item [] $1949=1+2\times (345+6+7\times 89).$
\item [] $1950=1+(2\times 3\times 45+6)\times 7+8+9.$
\item[]$\mbox{Decreasing order}$
\item [] $1881=9\times 87+(6+543)\times 2\times 1.$
\item [] $1882=9\times 87+(6+543)\times 2+1.$
\item [] $1883=(98+7\times 6\times 5\times 4+3)\times 2+1.$
\item [] $1884=9+8+7+6+5+43^2\times 1.$
\item [] $1885=9+8+7+6+5+43^2+1.$
\item [] $1886=9+8+7\times (65\times 4+3\times 2+1).$
\item [] $1887=9+(8+7)\times 65+43\times 21.$
\item [] $1888=98+765+4^(3+2)+1.$
\item [] $1889=98\times 7+6+(54+3)\times 21.$
\item [] $1890=(98+7)\times 6+5\times 4\times 3\times 21.$
\item [] $1891=(9+876+5\times 4\times 3)\times 2+1.$
\item [] $1892=98\times (7+6+5)+4\times 32\times 1.$
\item [] $1893=(9+8+76)\times 5\times 4+32+1.$
\item [] $1894=9+(8+7\times 65)\times 4+32+1.$
\item [] $1895=(9+8+7\times 65)\times 4+3\times 2+1.$
\item [] $1896=9\times (8+7)\times 6+543\times 2\times 1.$
\item [] $1897=98+7\times 65+4^3\times 21.$
\item [] $1898=9\times 8+7\times 65\times 4+3+2+1.$
\item [] $1899=9\times 8+7\times 65\times 4+3\times 2+1.$
\item [] $1900=9+8+7\times 65\times 4+3\times 21.$
\item [] $1901=987+6+5+43\times 21.$
\item [] $1902=98+76+54\times 32\times 1.$
\item [] $1903=98+76+54\times 32+1.$
\item [] $1904=9+8+7+6\times 5+43^2+1.$
\item [] $1905=98\times (7+6)+5^4+3\times 2\times 1.$
\item [] $1906=9+8\times 76+5+4\times 321.$
\item [] $1907=9+8+7+6+5^4\times 3+2\times 1.$
\item [] $1908=9+8+7+6+5^4\times 3+2+1.$
\item [] $1909=(9+8+76\times 5)\times 4+321.$
\item [] $1910=9+(8\times 7+6\times 54)\times (3+2)+1.$
\item [] $1911=9\times (8+7\times 6+5\times 4)\times 3+21.$
\item [] $1912=9\times 8\times 7+(6+5)\times 4\times 32\times 1.$
\item [] $1913=9\times 8+76\times 5\times 4+321.$
\item [] $1914=9+8+7\times 6+5+43^2+1.$
\item [] $1915=98\times (7+6)+5\times 4\times 32+1.$
\item [] $1916=9\times 8+7\times 65\times 4+3+21.$
\item [] $1917=9\times 87+(6+5+43)\times 21.$
\item [] $1918=(98+76\times 5)\times 4+3+2+1.$
\item [] $1919=(98+7)\times 6+5+4\times 321.$
\item [] $1920=987+6\times 5+43\times 21.$
\item [] $1921=(9+8+7\times 65)\times 4+32+1.$
\item [] $1922=98+7\times 6+54\times (32+1).$
\item [] $1923=98+7\times 65\times 4+3+2\times 1.$
\item [] $1924=98+7\times 65\times 4+3+2+1.$
\item [] $1925=9\times 8+7\times 65\times 4+32+1.$
\item [] $1926=9+8+7+6+5^4\times 3+21.$
\item [] $1927=98+7\times 65\times 4+3^2\times 1.$
\item [] $1928=98+7\times 65\times 4+3^2+1.$
\item [] $1929=9+(8+7+65)\times 4\times 3\times 2\times 1.$
\item [] $1930=9\times 8+76+54\times (32+1).$
\item [] $1931=98\times (7+6)+5^4+32\times 1.$
\item [] $1932=987+(6+5+4)\times 3\times 21.$
\item [] $1933=(9+8)\times 76+5\times 4\times 32+1.$
\item [] $1934=987+(6+5)\times 43\times 2+1.$
\item [] $1935=(9+87+65)\times 4\times 3+2+1.$
\item [] $1936=9+8+7\times 6+5^4\times 3+2\times 1.$
\item [] $1937=9+8+7\times 6+5^4\times 3+2+1.$
\item [] $1938=9+8+7+65+43^2\times 1.$
\item [] $1939=98+76\times 5\times 4+321.$
\item [] $1940=9\times 8+7+6+5+43^2+1.$
\item [] $1941=9+8\times 7\times 6\times 5+4\times 3\times 21.$
\item [] $1942=98+7\times 65\times 4+3+21.$
\item [] $1943=9+8+7\times (6+5)+43^2\times 1.$
\item [] $1944=9+8\times 7+6\times 5+43^2\times 1.$
\item [] $1945=9+8\times 7+6\times 5+43^2+1.$
\item [] $1946=9+8+7+6\times 5\times 4^3+2\times 1.$
\item [] $1947=9+8+76+5+43^2\times 1.$
\item [] $1948=9+8+76+5+43^2+1.$
\item [] $1949=9+8\times 7+6+5^4\times 3+2+1.$
\item [] $1950=98+7\times 65\times 4+32\times 1.$
\item[]$\mbox{Increasing order}$
\item [] $1951=1\times 23+4\times (5+6\times 78+9).$
\item [] $1952=(1+2+3)^4+567+89.$
\item [] $1953=1234+5+6\times 7\times (8+9).$
\item [] $1954=1234+5\times 6\times (7+8+9).$
\item [] $1955=(1+2\times 3)\times 45\times 6+7\times 8+9.$
\item [] $1956=(1+2\times 3+45\times 6)\times 7+8+9.$
\item [] $1957=123+4^5+6\times (7+8)\times 9.$
\item [] $1958=12^3+4+5+(6+7)\times (8+9).$
\item [] $1959=12^3+4+5\times 6\times 7+8+9.$
\item [] $1960=12\times 3+4\times (56\times 7+89).$
\item [] $1961=1+2\times 3+4+5\times 6\times (7\times 8+9).$
\item [] $1962=1^{23}\times 45\times 6\times 7+8\times 9.$
\item [] $1963=1^{23}+45\times 6\times 7+8\times 9.$
\item [] $1964=12\times 3+4\times (5+6\times 78+9).$
\item [] $1965=1\times 234\times 5+6+789.$
\item [] $1966=1234+5\times 6+78\times 9.$
\item [] $1967=1\times 2+3+45\times 6\times 7+8\times 9.$
\item [] $1968=1+2\times 3\times 4\times 56+7\times 89.$
\item [] $1969=1+2\times 3+45\times 6\times 7+8\times 9.$
\item [] $1970=1\times 2^3+45\times 6\times 7+8\times 9.$
\item [] $1971=1\times 2+34\times 56+7\times 8+9.$
\item [] $1972=1+2+34\times 56+7\times 8+9.$
\item [] $1973=1\times (2+3)\times 45\times 6+7\times 89.$
\item [] $1974=1+(2+3)\times 45\times 6+7\times 89.$
\item [] $1975=(1+2+3)^4+56+7\times 89.$
\item [] $1976=12^3+4\times 56+7+8+9.$
\item [] $1977=12+3+45\times 6\times 7+8\times 9.$
\item [] $1978=1+23+4+5\times 6\times (7\times 8+9).$
\item [] $1979=1^{23}\times 45\times 6\times 7+89.$
\item [] $1980=1^{23}+45\times 6\times 7+89.$
\item [] $1981=12+34\times 56+7\times 8+9.$
\item [] $1982=1^2\times 3+45\times 6\times 7+89.$
\item [] $1983=1^2+3+45\times 6\times 7+89.$
\item [] $1984=1\times 2+3+45\times 6\times 7+89.$
\item [] $1985=1\times 23+45\times 6\times 7+8\times 9.$
\item [] $1986=1+23+45\times 6\times 7+8\times 9.$
\item [] $1987=1\times 2^3+45\times 6\times 7+89.$
\item [] $1988=1+2^3+45\times 6\times 7+89.$
\item [] $1989=(12\times 3+4)\times 5\times 6+789.$
\item [] $1990=12\times 3+4+5\times 6\times (7\times 8+9).$
\item [] $1991=1^2\times 34\times 56+78+9.$
\item [] $1992=1234+56+78\times 9.$
\item [] $1993=1\times 2+34\times 56+78+9.$
\item [] $1994=12+3+45\times 6\times 7+89.$
\item [] $1995=12+34\times 56+7+8\times 9.$
\item [] $1996=12+34+5\times 6\times (7\times 8+9).$
\item [] $1997=12^3+4\times (56+7)+8+9.$
\item [] $1998=12\times 3+45\times 6\times 7+8\times 9.$
\item [] $1999=12^3+(4\times 5+6)\times 7+89.$
\item [] $2000=1^2\times 34\times 56+7+89.$
\item [] $2001=1^2+34\times 56+7+89.$
\item [] $2002=1\times 23+45\times 6\times 7+89.$
\item [] $2003=12+34\times 56+78+9.$
\item [] $2004=1+23\times (4+56)+7\times 89.$
\item [] $2005=12^3+4\times 5\times (6+7)+8+9.$
\item [] $2006=(1+2)^3+45\times 6\times 7+89.$
\item [] $2007=1\times 2\times (3^4+56)\times 7+89.$
\item [] $2008=1+2\times (3^4+56)\times 7+89.$
\item [] $2009=12\times 3\times (4+5)\times 6+7\times 8+9.$
\item [] $2010=12\times 3\times 45+6\times (7\times 8+9).$
\item [] $2011=(1+2\times 3+45\times 6)\times 7+8\times 9.$
\item [] $2012=12+34\times 56+7+89.$
\item [] $2013=12\times (3\times 4+5)\times 6+789.$
\item [] $2014=1234+5\times (67+89).$
\item [] $2015=12\times 3+45\times 6\times 7+89.$
\item [] $2016=12\times (3+4+5+67+89).$
\item [] $2017=12^3+4\times 56+7\times 8+9.$
\item [] $2018=1\times 2\times 34+5\times 6\times (7\times 8+9).$
\item [] $2019=1+2\times 34+5\times 6\times (7\times 8+9).$
\item [] $2020=(123+4)\times (5+6)+7\times 89.$
\item[]$\mbox{Decreasing order}$
\item [] $1951=98+7\times 65\times 4+32+1.$
\item [] $1952=98\times 7+6+5\times 4\times 3\times 21.$
\item [] $1953=(9+87+65)\times 4\times 3+21.$
\item [] $1954=9+(8+7)\times 6+5+43^2+1.$
\item [] $1955=987+65+43\times 21.$
\item [] $1956=9+87+6+5+43^2\times 1.$
\item [] $1957=9+87+6+5+43^2+1.$
\item [] $1958=9\times 8+7+6\times 5+43^2\times 1.$
\item [] $1959=9\times 8+7+6\times 5+43^2+1.$
\item [] $1960=98+7\times (65\times 4+3)+21.$
\item [] $1961=987+6\times 54\times 3+2\times 1.$
\item [] $1962=987+654+321.$
\item [] $1963=9\times 8+7+6+5^4\times 3+2+1.$
\item [] $1964=9+(8+7+6)\times 5+43^2+1.$
\item [] $1965=98+7+6+5+43^2\times 1.$
\item [] $1966=98+7+6+5+43^2+1.$
\item [] $1967=9+8\times 7+6+5^4\times 3+21.$
\item [] $1968=9+8+7+6\times 54\times 3\times 2\times 1.$
\item [] $1969=9+8+7+6\times 54\times 3\times 2+1.$
\item [] $1970=9+8+76+5^4\times 3+2\times 1.$
\item [] $1971=9+8+76+5^4\times 3+2+1.$
\item [] $1972=9+876+543\times 2+1.$
\item [] $1973=9\times 8+76\times (5\times 4+3+2)+1.$
\item [] $1974=(9+8)\times 7\times 6+5\times 4\times 3\times 21.$
\item [] $1975=9+87+6\times 5+43^2\times 1.$
\item [] $1976=9+87+6\times 5+43^2+1.$
\item [] $1977=987+6\times 5\times (4\times 3+21).$
\item [] $1978=98\times 7+6\times 5\times 43+2\times 1.$
\item [] $1979=98\times 7+6\times 5\times 43+2+1.$
\item [] $1980=987+6\times 54\times 3+21.$
\item [] $1981=98+7\times 65\times 4+3\times 21.$
\item [] $1982=98\times 7+6\times (5+4)\times (3+21).$
\item [] $1983=9\times 87+6\times 5\times 4\times (3^2+1).$
\item [] $1984=98+7+6\times 5+43^2\times 1.$
\item [] $1985=98+7+6\times 5+43^2+1.$
\item [] $1986=9+8\times 7+(6+54)\times 32+1.$
\item [] $1987=9+8\times 7+6\times 5\times 4^3+2\times 1.$
\item [] $1988=9+8\times 7+6\times 5\times 4^3+2+1.$
\item [] $1989=98+7+6+5^4\times 3+2+1.$
\item [] $1990=9+8\times (7+6)+5^4\times 3+2\times 1.$
\item [] $1991=9\times 8+7\times 6+5^4\times 3+2\times 1.$
\item [] $1992=9\times 8+7\times 6+5^4\times 3+2+1.$
\item [] $1993=9\times 8+7+65+43^2\times 1.$
\item [] $1994=9\times 8+7+65+43^2+1.$
\item [] $1995=98+7\times 6+5+43^2+1.$
\item [] $1996=(9+8)\times (7+6)\times (5+4)+3\times 2+1.$
\item [] $1997=98\times 7+6\times 5\times 43+21.$
\item [] $1998=9+87+6+5^4\times 3+21.$
\item [] $1999=9\times 8+7+(6+54)\times 32\times 1.$
\item [] $2000=98\times 7+6\times 5+4\times 321.$
\item [] $2001=9\times 8+7+6\times 5\times 4^3+2\times 1.$
\item [] $2002=9\times 8+76+5+43^2\times 1.$
\item [] $2003=9\times 8+76+5+43^2+1.$
\item [] $2004=(9+8+7)\times 6\times 5+4\times 321.$
\item [] $2005=9\times 8+7\times 6\times (5\times 4+3)\times 2+1.$
\item [] $2006=9+8\times 7+6\times 5\times 4^3+21.$
\item [] $2007=9+8+7+654\times 3+21.$
\item [] $2008=(9\times 87+6+5\times 43)\times 2\times 1.$
\item [] $2009=9+8\times 7+6\times 54\times 3\times 2\times 1.$
\item [] $2010=9+8\times 7+6\times 54\times 3\times 2+1.$
\item [] $2011=9+87+65+43^2+1.$
\item [] $2012=9\times (8+7\times 6+5)\times 4+32\times 1.$
\item [] $2013=98+(7+6)\times 5+43^2+1.$
\item [] $2014=9+8\times 7\times 6\times 5+4+321.$
\item [] $2015=(9+8\times 76+54)\times 3+2\times 1.$
\item [] $2016=9+87+(6+54)\times 32\times 1.$
\item [] $2017=98+7\times 6+5^4\times 3+2\times 1.$
\item [] $2018=98+7\times 6+5^4\times 3+2+1.$
\item [] $2019=98+7+65+43^2\times 1.$
\item [] $2020=98+7+65+43^2+1.$
\item[]$\mbox{Increasing order}$
\item [] $2021=12\times (3\times 4+5+6)\times 7+89.$
\item [] $2022=12^3+45\times 6+7+8+9.$
\item [] $2023=12\times 3\times (4+5)\times 6+7+8\times 9.$
\item [] $2024=123+4\times (5+6\times 78)+9.$
\item [] $2025=12\times 3+(4+5)\times (6+7)\times (8+9).$
\item [] $2026=1+2+34\times 56+7\times (8+9).$
\item [] $2027=12^3+4+5\times (6\times 7+8+9).$
\item [] $2028=12\times (34+56+7+8\times 9).$
\item [] $2029=1^2+3+4\times (56+7)\times 8+9.$
\item [] $2030=123+45\times 6\times 7+8+9.$
\item [] $2031=12^3+4+5\times 6\times 7+89.$
\item [] $2032=1+23\times (4+5)\times 6+789.$
\item [] $2033=1+(2\times 3\times 4+5)\times 67+89.$
\item [] $2034=1234+5+6+789.$
\item [] $2035=12+34\times 56+7\times (8+9).$
\item [] $2036=123+(4+5\times 6)\times 7\times 8+9.$
\item [] $2037=(1+2+3\times 4+5+6)\times 78+9.$
\item [] $2038=1^{23}+(4\times 5+6)\times 78+9.$
\item [] $2039=12^3+4\times 56+78+9.$
\item [] $2040=1\times (2+34)\times 56+7+8+9.$
\item [] $2041=1+(2+34)\times 56+7+8+9.$
\item [] $2042=1+2+34\times 56+(7+8)\times 9.$
\item [] $2043=(123+4\times 5+6+78)\times 9.$
\item [] $2044=1+2\times 3+(4\times 5+6)\times 78+9.$
\item [] $2045=12+(3+45)\times 6\times 7+8+9.$
\item [] $2046=1\times 2\times 3\times 4\times 56+78\times 9.$
\item [] $2047=1+2\times 3\times 4\times 56+78\times 9.$
\item [] $2048=12^3+4\times 56+7+89.$
\item [] $2049=1234+5+6\times (7+8)\times 9.$
\item [] $2050=1^2+3\times (4+56+7\times 89).$
\item [] $2051=123+4\times (5+6\times 78+9).$
\item [] $2052=1\times (2+3)\times 45\times 6+78\times 9.$
\item [] $2053=1234+5\times 6+789.$
\item [] $2054=(1+234+56)\times 7+8+9.$
\item [] $2055=12+3+4\times 5\times (6+7+89).$
\item [] $2056=(1^2+34)\times 56+7+89.$
\item [] $2057=(1\times 234+5)\times 6+7\times 89.$
\item [] $2058=1+(234+5)\times 6+7\times 89.$
\item [] $2059=(12+3\times 4+5)\times (6+7\times 8+9).$
\item [] $2060=12^3+4\times 5\times (6+7)+8\times 9.$
\item [] $2061=12+3\times (4+56+7\times 89).$
\item [] $2062=1\times 2^3+(4\times 5+6)\times (7+8\times 9).$
\item [] $2063=12^3+45\times 6+7\times 8+9.$
\item [] $2064=1\times 2\times 34\times 5\times 6+7+8+9.$
\item [] $2065=1+2\times 34\times 5\times 6+7+8+9.$
\item [] $2066=1+(2+3)\times (4+56\times 7+8+9).$
\item [] $2067=(12+3+45\times 6)\times 7+8\times 9.$
\item [] $2068=1\times 2^3+4\times (5+6+7\times 8\times 9).$
\item [] $2069=12^3+4\times (56+7)+89.$
\item [] $2070=12\times 34\times 5+6+7+8+9.$
\item [] $2071=1^2+3\times 456+78\times 9.$
\item [] $2072=1\times 2+3\times 456+78\times 9.$
\item [] $2073=1+2+3\times 456+78\times 9.$
\item [] $2074=1+23\times (45+6\times 7)+8\times 9.$
\item [] $2075=(1+23+45\times 6)\times 7+8+9.$
\item [] $2076=12\times (3\times 4+5+67+89).$
\item [] $2077=12^3+45\times 6+7+8\times 9.$
\item [] $2078=1+23+(4\times 5+6)\times (7+8\times 9).$
\item [] $2079=1234+56+789.$
\item [] $2080=(1+2+3+4\times 5+6)\times (7\times 8+9).$
\item [] $2081=1\times (2+34)\times 56+7\times 8+9.$
\item [] $2082=12+3\times 456+78\times 9.$
\item [] $2083=1234+56\times (7+8)+9.$
\item [] $2084=12^3+4+5\times 67+8+9.$
\item [] $2085=123+45\times 6\times 7+8\times 9.$
\item [] $2086=1+(2\times 3+4+5)\times (67+8\times 9).$
\item [] $2087=12^3+4\times 56+(7+8)\times 9.$
\item [] $2088=12\times 3\times (4+5\times 6+7+8+9).$
\item [] $2089=1^2+3\times (4+5+678+9).$
\item [] $2090=1\times 23\times (45+6\times 7)+89.$
\item[]$\mbox{Decreasing order}$
\item [] $2021=(9\times 8+7+6)\times 5\times 4+321.$
\item [] $2022=9\times 8\times (7+6)+543\times 2\times 1.$
\item [] $2023=9\times 8+7+6\times 54\times 3\times 2\times 1.$
\item [] $2024=9\times 8+7+6\times 54\times 3\times 2+1.$
\item [] $2025=9\times 8+76+5^4\times 3+2\times 1.$
\item [] $2026=9\times 8+76+5^4\times 3+2+1.$
\item [] $2027=98+7+6\times 5\times 4^3+2\times 1.$
\item [] $2028=98+76+5+43^2\times 1.$
\item [] $2029=98+76+5+43^2+1.$
\item [] $2030=9+8\times 7+654\times 3+2+1.$
\item [] $2031=9+87\times (6+5+4\times 3)+21.$
\item [] $2032=(987+6+5\times 4+3)\times 2\times 1.$
\item [] $2033=(987+6+5\times 4+3)\times 2+1.$
\item [] $2034=(9+8\times 76+54)\times 3+21.$
\item [] $2035=98\times 7+65+4\times 321.$
\item [] $2036=98+7\times 6+5^4\times 3+21.$
\item [] $2037=9+87+6\times 5\times 4^3+21.$
\item [] $2038=98+7+6\times (5\times 4^3+2)+1.$
\item [] $2039=9+8+(7\times 6+5^4)\times 3+21.$
\item [] $2040=9+87+6\times 54\times 3\times 2\times 1.$
\item [] $2041=9+87+6\times 54\times 3\times 2+1.$
\item [] $2042=98+(76+5)\times 4\times 3\times 2\times 1.$
\item [] $2043=9\times 8+7+654\times 3+2\times 1.$
\item [] $2044=9\times 8+76+5^4\times 3+21.$
\item [] $2045=(9+87+6)\times 5\times 4+3+2\times 1.$
\item [] $2046=98+7+6\times 5\times 4^3+21.$
\item [] $2047=(9+87+6)\times 5\times 4+3\times 2+1.$
\item [] $2048=9+8\times 7+654\times 3+21.$
\item [] $2049=9\times 87+6+5\times 4\times 3\times 21.$
\item [] $2050=98+7+6\times 54\times 3\times 2+1.$
\item [] $2051=98+76+5^4\times 3+2\times 1.$
\item [] $2052=98+76+5^4\times 3+2+1.$
\item [] $2053=(98+76+54)\times 3^2+1.$
\item [] $2054=9+8+7\times 6\times (5+43)+21.$
\item [] $2055=(9+8+7+654)\times 3+21.$
\item [] $2056=9\times 8+(7\times 6+5\times 4)\times 32\times 1.$
\item [] $2057=9\times 8+(7+654)\times 3+2\times 1.$
\item [] $2058=9\times 8+(7+654)\times 3+2+1.$
\item [] $2059=9+8+(7\times 6+5)\times 43+21.$
\item [] $2060=9+87+654\times 3+2\times 1.$
\item [] $2061=9+87+654\times 3+2+1.$
\item [] $2062=9\times 8+7+654\times 3+21.$
\item [] $2063=(9+8)\times 7+6\times 54\times 3\times 2\times 1.$
\item [] $2064=9\times 8\times 7+65\times 4\times 3\times 2\times 1.$
\item [] $2065=9\times 8\times 7+65\times 4\times 3\times 2+1.$
\item [] $2066=9+8+765+4\times 321.$
\item [] $2067=987+6\times 5\times 4\times 3^2\times 1.$
\item [] $2068=987+6\times 5\times 4\times 3^2+1.$
\item [] $2069=98+7+654\times 3+2\times 1.$
\item [] $2070=98+7+654\times 3+2+1.$
\item [] $2071=9+8+76\times (5+4)\times 3+2\times 1.$
\item [] $2072=(9+87+6)\times 5\times 4+32\times 1.$
\item [] $2073=9+8\times 7\times 6+54\times 32\times 1.$
\item [] $2074=9+8\times 7\times 6+54\times 32+1.$
\item [] $2075=9\times 87+6\times 5\times 43+2\times 1.$
\item [] $2076=9\times 87+6\times 5\times 43+2+1.$
\item [] $2077=9+8+7\times 6\times 5+43^2+1.$
\item [] $2078=9\times 87+6+5+4\times 321.$
\item [] $2079=9+87+654\times 3+21.$
\item [] $2080=987+6+543\times 2+1.$
\item [] $2081=9+8\times (76+54\times 3+21).$
\item [] $2082=9+876+(54+3)\times 21.$
\item [] $2083=98+(7+654)\times 3+2\times 1.$
\item [] $2084=(9+8)\times 7+654\times 3+2+1.$
\item [] $2085=98+(7+6\times 54)\times 3\times 2+1.$
\item [] $2086=987+(6+543)\times 2+1.$
\item [] $2087=(9\times 8+7)\times (6+5\times 4)+32+1.$
\item [] $2088=98+7+654\times 3+21.$
\item [] $2089=98\times (7+6+5)+4+321.$
\item [] $2090=9+8+76\times (5+4)\times 3+21.$
\item[]$\mbox{Increasing order}$
\item [] $2091=1+23\times (45+6\times 7)+89.$
\item [] $2092=1^2\times 3+4\times 5\times (6+7)\times 8+9.$
\item [] $2093=12^3+4\times 5+6\times 7\times 8+9.$
\item [] $2094=12^3+45\times 6+7+89.$
\item [] $2095=1^2+345\times 6+7+8+9.$
\item [] $2096=1\times 2+345\times 6+7+8+9.$
\item [] $2097=12\times 34+5\times 6\times 7\times 8+9.$
\item [] $2098=1^2+3\times (45+6\times 7)\times 8+9.$
\item [] $2099=12\times 34\times 5+6\times 7+8+9.$
\item [] $2100=12+3\times (4+5+678+9).$
\item [] $2101=1234+(5+6)\times 78+9.$
\item [] $2102=123+45\times 6\times 7+89.$
\item [] $2103=1+(234+56)\times 7+8\times 9.$
\item [] $2104=1+(2+34)\times 56+78+9.$
\item [] $2105=1\times 2\times 34\times 5\times 6+7\times 8+9.$
\item [] $2106=12+345\times 6+7+8+9.$
\item [] $2107=1+234\times 5+(6+7)\times 8\times 9.$
\item [] $2108=12^3+4\times (5+6)\times 7+8\times 9.$
\item [] $2109=(1+234+56)\times 7+8\times 9.$
\item [] $2110=1+(2+3+4\times 5)\times (6+78)+9.$
\item [] $2111=12\times 34\times 5+6+7\times 8+9.$
\item [] $2112=1\times 23+4\times 5\times (6+7)\times 8+9.$
\item [] $2113=1+(2+34)\times 56+7+89.$
\item [] $2114=1+2\times (3+4^5)+6\times 7+8+9.$
\item [] $2115=(1+2+3)^4+5\times 6+789.$
\item [] $2116=12^3+4+5\times (67+8)+9.$
\item [] $2117=12+(3+45)\times 6\times 7+89.$
\item [] $2118=12^3+45+6\times 7\times 8+9.$
\item [] $2119=1\times 2\times 34\times 5\times 6+7+8\times 9.$
\item [] $2120=1+2\times 34\times 5\times 6+7+8\times 9.$
\item [] $2121=12\times (3^4+5\times 6)+789.$
\item [] $2122=1^2+3\times (4\times 5+678+9).$
\item [] $2123=1\times 2\times 3\times (4+5\times 67)+89.$
\item [] $2124=12\times 34\times 5+67+8+9.$
\item [] $2125=12\times 34\times 5+6+7+8\times 9.$
\item [] $2126=(1+234+56)\times 7+89.$
\item [] $2127=1\times 2\times 34\times 5\times 6+78+9.$
\item [] $2128=1+2\times 34\times 5\times 6+78+9.$
\item [] $2129=1+2\times (3+4)\times (56+7+89).$
\item [] $2130=12\times 3\times 45+6+7\times 8\times 9.$
\item [] $2131=1\times 2\times (3+4^5+6)+7\times 8+9.$
\item [] $2132=1+2\times (3+4^5+6)+7\times 8+9.$
\item [] $2133=1\times 2\times 3\times 4\times 56+789.$
\item [] $2134=1+2\times 3\times 4\times 56+789.$
\item [] $2135=1^2\times 345\times 6+7\times 8+9.$
\item [] $2136=1^2+345\times 6+7\times 8+9.$
\item [] $2137=1\times 2+345\times 6+7\times 8+9.$
\item [] $2138=1+2+345\times 6+7\times 8+9.$
\item [] $2139=12^3+4+5\times 67+8\times 9.$
\item [] $2140=1+(2+3)\times 45\times 6+789.$
\item [] $2141=12^3+4+56\times 7+8+9.$
\item [] $2142=12\times 34\times 5+6+7+89.$
\item [] $2143=1+2\times 3+4\times (5\times 6+7\times 8\times 9).$
\item [] $2144=1+2\times (34\times 5\times 6+7)+89.$
\item [] $2145=(1+2+3\times 4)\times (56+78+9).$
\item [] $2146=1+2\times 3\times 4\times (5+6+78)+9.$
\item [] $2147=12+345\times 6+7\times 8+9.$
\item [] $2148=123+4\times (56+7)\times 8+9.$
\item [] $2149=1^2\times 345\times 6+7+8\times 9.$
\item [] $2150=1^2+345\times 6+7+8\times 9.$
\item [] $2151=1\times 2+345\times 6+7+8\times 9.$
\item [] $2152=1+2+345\times 6+7+8\times 9.$
\item [] $2153=12\times 34\times 5+(6+7)\times 8+9.$
\item [] $2154=12\times 34\times 5+6\times 7+8\times 9.$
\item [] $2155=(1^2+345)\times 6+7+8\times 9.$
\item [] $2156=12^3+4+5\times 67+89.$
\item [] $2157=1^2\times 345\times 6+78+9.$
\item [] $2158=1^2+345\times 6+78+9.$
\item [] $2159=1\times 2+345\times 6+78+9.$
\item [] $2160=1+2+3\times 456+789.$
\item[]$\mbox{Decreasing order}$
\item [] $2091=9+(87+65\times 4)\times 3\times 2\times 1.$
\item [] $2092=9+(87+65\times 4)\times 3\times 2+1.$
\item [] $2093=(9\times 87+65\times 4+3)\times 2+1.$
\item [] $2094=9\times 87+6\times 5\times 43+21.$
\item [] $2095=98\times 7+65+4^3\times 21.$
\item [] $2096=9+8\times (7+6)\times 5\times 4+3\times 2+1.$
\item [] $2097=9\times 87+6\times 5+4\times 321.$
\item [] $2098=9\times 87+(654+3)\times 2+1.$
\item [] $2099=9\times (8+7)\times 6+5+4\times 321.$
\item [] $2100=9\times (8+7)+654\times 3+2+1.$
\item [] $2101=(987+6+54+3)\times 2+1.$
\item [] $2102=(9+8)\times 7+654\times 3+21.$
\item [] $2103=(9+87+6)\times 5\times 4+3\times 21.$
\item [] $2104=(9\times (87+6)+5\times 43)\times 2\times 1.$
\item [] $2105=(98\times 7+6+5+4)\times 3+2\times 1.$
\item [] $2106=9\times 87+(6+54+3)\times 21.$
\item [] $2107=(987+6+5\times 4\times 3)\times 2+1.$
\item [] $2108=9+(8+7\times 6)\times 5+43^2\times 1.$
\item [] $2109=9\times 8+7\times 6\times (5+43)+21.$
\item [] $2110=(9+8)\times (76+5+43)+2\times 1.$
\item [] $2111=9+8\times 7+6\times (5\times 4+321).$
\item [] $2112=(9+8)\times (7+6\times 5+4)\times 3+21.$
\item [] $2113=9+8\times (7+6)\times 5\times 4+3+21.$
\item [] $2114=(9\times 8+7\times 65)\times 4+3+2+1.$
\item [] $2115=(9\times 8+7\times 65)\times 4+3\times 2+1.$
\item [] $2116=98+7\times 6\times (5+43)+2\times 1.$
\item [] $2117=98+7\times 6\times (5+43)+2+1.$
\item [] $2118=9\times (8+7)+654\times 3+21.$
\item [] $2119=(987+65+4+3)\times 2+1.$
\item [] $2120=(987+6\times 5+43)\times 2\times 1.$
\item [] $2121=9\times 8+765+4\times 321.$
\item [] $2122=9+8\times 7\times 6\times 5+432+1.$
\item [] $2123=(9\times 8\times 7+6+5)\times 4+3\times 21.$
\item [] $2124=(98\times 7+6+5+4)\times 3+21.$
\item [] $2125=(9+8\times 7\times 6+5+4)\times 3\times 2+1.$
\item [] $2126=9+8+765+4^3\times 21.$
\item [] $2127=9+(87\times 6+5)\times 4+3^2+1.$
\item [] $2128=(987+65+4\times 3)\times 2\times 1.$
\item [] $2129=(9\times 8+76\times 5)\times 4+321.$
\item [] $2130=9+(8+7)\times 6\times 5\times 4+321.$
\item [] $2131=9\times 8+7\times 6\times 5+43^2\times 1.$
\item [] $2132=9\times 87+65+4\times 321.$
\item [] $2133=9\times 87+6\times 5\times (43+2)\times 1.$
\item [] $2134=9+8\times ((7+6)\times 5\times 4+3)+21.$
\item [] $2135=98+7\times 6\times (5+43)+21.$
\item [] $2136=(9+8+7+65)\times 4\times 3\times 2\times 1.$
\item [] $2137=(9+8+7+65)\times 4\times 3\times 2+1.$
\item [] $2138=9\times 87+6+5+4^3\times 21.$
\item [] $2139=(98\times 7+6+5\times 4)\times 3+2+1.$
\item [] $2140=98+(7\times 6+5)\times 43+21.$
\item [] $2141=9+(87\times 6+5)\times 4+3+21.$
\item [] $2142=(9+8\times 7+6\times 5+4+3)\times 21.$
\item [] $2143=(9\times 8+76+5)\times (4+3)\times 2+1.$
\item [] $2144=9+8+(7+6)\times 54\times 3+21.$
\item [] $2145=9+876+5\times 4\times 3\times 21.$
\item [] $2146=9+(8+76+5)\times 4\times 3\times 2+1.$
\item [] $2147=98+765+4\times 321.$
\item [] $2148=(9+87\times 6+543)\times 2\times 1.$
\item [] $2149=(9+87\times 6+543)\times 2+1.$
\item [] $2150=(9+87\times 6+5)\times 4+3+2+1.$
\item [] $2151=9+8\times (7+65\times 4)+3+2+1.$
\item [] $2152=98+76\times (5+4)\times 3+2\times 1.$
\item [] $2153=98+76\times (5+4)\times 3+2+1.$
\item [] $2154=(98+76)\times 5+4\times 321.$
\item [] $2155=9+8\times (7+6\times 5)+43^2+1.$
\item [] $2156=9+87\times 6+5\times (4+321).$
\item [] $2157=9\times 87+6\times 5+4^3\times 21.$
\item [] $2158=9+8+7\times 65\times 4+321.$
\item [] $2159=9\times (8+7)\times 6+5+4^3\times 21.$
\item [] $2160=(9+8\times 7+654)\times 3+2+1.$
\item[]$\mbox{Increasing order}$
\item [] $2161=12+345\times 6+7+8\times 9.$
\item [] $2162=1+(23+4+5)\times 67+8+9.$
\item [] $2163=123+4\times 5\times (6+7+89).$
\item [] $2164=1\times 2\times (3+456+7\times 89).$
\item [] $2165=12\times 3\times 45+67\times 8+9.$
\item [] $2166=1^2\times 345\times 6+7+89.$
\item [] $2167=1^2+345\times 6+7+89.$
\item [] $2168=1\times 2+345\times 6+7+89.$
\item [] $2169=12+3\times 456+789.$
\item [] $2170=1+23\times (4+56)+789.$
\item [] $2171=12\times 34\times 5+6\times 7+89.$
\item [] $2172=1+2^3+4^5+67\times (8+9).$
\item [] $2173=1+(2+3+45)\times 6\times 7+8\times 9.$
\item [] $2174=1\times 23\times 45+67\times (8+9).$
\item [] $2175=1234+5+(6+7)\times 8\times 9.$
\item [] $2176=1+2\times 34\times 5\times 6+(7+8)\times 9.$
\item [] $2177=12\times 3\times (45+6+7)+89.$
\item [] $2178=12+345\times 6+7+89.$
\item [] $2179=12\times 34\times 5+67+8\times 9.$
\item [] $2180=1+2\times 3+4\times (5+67\times 8)+9.$
\item [] $2181=123\times 4+5\times 6\times 7\times 8+9.$
\item [] $2182=1+2^3+4\times (5+67\times 8)+9.$
\item [] $2183=123+4\times (5+6+7\times 8\times 9).$
\item [] $2184=1\times 2^3\times 45\times 6+7+8+9.$
\item [] $2185=1+2^3\times 45\times 6+7+8+9.$
\item [] $2186=1\times 23+4^5+67\times (8+9).$
\item [] $2187=1+23+4^5+67\times (8+9).$
\item [] $2188=12+34\times (5+6\times 7+8+9).$
\item [] $2189=(1+2+3+4)\times 5\times 6\times 7+89.$
\item [] $2190=1^2+345\times 6+7\times (8+9).$
\item [] $2191=(1+23+4)\times 56+7\times 89.$
\item [] $2192=1\times 23\times 45+(6+7)\times 89.$
\item [] $2193=1+23\times 45+(6+7)\times 89.$
\item [] $2194=1+2\times (3+4^5)+67+8\times 9.$
\item [] $2195=12^3+(4+5)\times 6\times 7+89.$
\item [] $2196=12\times 34\times 5+67+89.$
\item [] $2197=12^3+4+5\times (6+78+9).$
\item [] $2198=12^3+4\times 5+(6\times 7+8)\times 9.$
\item [] $2199=12\times (34\times 5+6)+78+9.$
\item [] $2200=(1+2)^3+4\times (5+67\times 8)+9.$
\item [] $2201=12+345\times 6+7\times (8+9).$
\item [] $2202=1+(2+345)\times 6+7\times (8+9).$
\item [] $2203=1^2\times 3+4\times (5+67\times 8+9).$
\item [] $2204=1\times 23+4^5+(6+7)\times 89.$
\item [] $2205=1+23+4^5+(6+7)\times 89.$
\item [] $2206=1^2+345\times 6+(7+8)\times 9.$
\item [] $2207=12\times 3\times 4\times (5+6)+7\times 89.$
\item [] $2208=12^3+456+7+8+9.$
\item [] $2209=12\times 3+4\times (5+67\times 8)+9.$
\item [] $2210=1\times 2\times (3+4^5)+67+89.$
\item [] $2211=1+2\times (3+4^5)+67+89.$
\item [] $2212=123+4\times 5\times (6+7)\times 8+9.$
\item [] $2213=12^3+4+56\times 7+89.$
\item [] $2214=12^3+4+5+6\times 78+9.$
\item [] $2215=12+3+4\times (5+67\times 8+9).$
\item [] $2216=1\times (23+4+5)\times 67+8\times 9.$
\item [] $2217=12+345\times 6+(7+8)\times 9.$
\item [] $2218=1+(2+345)\times 6+(7+8)\times 9.$
\item [] $2219=(1+23)\times 45+67\times (8+9).$
\item [] $2220=12^3+(4+56)\times 7+8\times 9.$
\item [] $2221=(12+345)\times 6+7+8\times 9.$
\item [] $2222=1\times 2\times (3+4^5+67+8+9).$
\item [] $2223=(1\times 234+5)\times 6+789.$
\item [] $2224=1+(234+5)\times 6+789.$
\item [] $2225=12^3+4\times 5+6\times 78+9.$
\item [] $2226=1+2^3\times 45\times 6+7\times 8+9.$
\item [] $2227=1^2\times 34\times 5\times (6+7)+8+9.$
\item [] $2228=(1+23+4+5)\times 67+8+9.$
\item [] $2229=(1+234+5)\times 6+789.$
\item [] $2230=1+2+34\times 5\times (6+7)+8+9.$
\item[]$\mbox{Decreasing order}$
\item [] $2161=9+8+(7+6+54)\times 32\times 1.$
\item [] $2162=9+8+(7+6+54)\times 32+1.$
\item [] $2163=(98\times 7+6\times 5+4)\times 3+2+1.$
\item [] $2164=(9+8)\times (7+6\times 5\times 4)+3+2\times 1.$
\item [] $2165=98\times (7+6+5+4)+3^2\times 1.$
\item [] $2166=9+8+7+6\times (5+4\times 3)\times 21.$
\item [] $2167=(9\times 8+7\times 6)\times (5+4\times 3+2)+1.$
\item [] $2168=(9+87\times 6+5)\times 4+3+21.$
\item [] $2169=(98+76+5)\times 4\times 3+21.$
\item [] $2170=9+(8+7)\times (65+4+3)\times 2+1.$
\item [] $2171=(9\times 8+7\times 65)\times 4+3\times 21.$
\item [] $2172=9+8\times (7+65\times 4+3)+2+1.$
\item [] $2173=9\times 8+7\times 6\times (5+43+2)+1.$
\item [] $2174=9+876+5+4\times 321.$
\item [] $2175=(9+8+76+5^4)\times 3+21.$
\item [] $2176=(9+87\times 6+5)\times 4+32\times 1.$
\item [] $2177=(9+87\times 6+5)\times 4+32+1.$
\item [] $2178=(9+8\times 7+654)\times 3+21.$
\item [] $2179=(98+7)\times (6+5)+4^(3+2)\times 1.$
\item [] $2180=98\times (7+6+5+4)+3+21.$
\item [] $2181=9\times 8+765+4^3\times 21.$
\item [] $2182=98\times (7+6)+5+43\times 21.$
\item [] $2183=(9+87+6+5^4)\times 3+2\times 1.$
\item [] $2184=(9+87+6+5^4)\times 3+2+1.$
\item [] $2185=9+8\times (76+5\times 4\times 3)\times 2\times 1.$
\item [] $2186=9+8\times (76+5\times 4\times 3)\times 2+1.$
\item [] $2187=987+6\times 5\times 4\times (3^2+1).$
\item [] $2188=98\times (7+6+5+4)+32\times 1.$
\item [] $2189=98\times (7+6+5+4)+32+1.$
\item [] $2190=987+6+(54+3)\times 21.$
\item [] $2191=9+8+7+6+5\times 432+1.$
\item [] $2192=9\times 87+65+4^3\times 21.$
\item [] $2193=9+8\times 7\times (6+5+4+3+21).$
\item [] $2194=9+8\times 7\times (6\times 5+4+3+2)+1.$
\item [] $2195=9+8+7+6+5\times (432+1).$
\item [] $2196=98\times (7+6+5)+432\times 1.$
\item [] $2197=98\times (7+6+5)+432+1.$
\item [] $2198=98+7\times 6\times 5\times (4+3+2+1).$
\item [] $2199=9+8\times 7\times 6+5+43^2\times 1.$
\item [] $2200=9+8\times 7\times 6+5+43^2+1.$
\item [] $2201=(9+876+5\times 43)\times 2+1.$
\item [] $2202=(9\times 8+7+654)\times 3+2+1.$
\item [] $2203=(9\times 8+7)\times 6+54\times 32+1.$
\item [] $2204=(9+8\times 7+6)\times 5+43^2\times 1.$
\item [] $2205=(9+8+76+5+4+3)\times 21.$
\item [] $2206=98+(7+6)\times 54\times 3+2\times 1.$
\item [] $2207=98+765+4^3\times 21.$
\item [] $2208=9+8\times (7+65\times 4)+3\times 21.$
\item [] $2209=(9+8+7\times 65)\times 4+321.$
\item [] $2210=(98+7+6+5^4)\times 3+2\times 1.$
\item [] $2211=(98+7+6+5^4)\times 3+2+1.$
\item [] $2212=(9+8)\times 7\times (6+5)+43\times 21.$
\item [] $2213=9\times 8+7\times 65\times 4+321.$
\item [] $2214=(98+76)\times 5+4^3\times 21.$
\item [] $2215=9\times 8+7\times (6\times 5+4)\times 3^2+1.$
\item [] $2216=9\times 8+(7+6+54)\times 32\times 1.$
\item [] $2217=9\times 8+(7+6+54)\times 32+1.$
\item [] $2218=(98+7\times 65)\times 4+3+2+1.$
\item [] $2219=9+8+7\times 6+5\times 432\times 1.$
\item [] $2220=9+8+7\times 6+5\times 432+1.$
\item [] $2221=9\times 8+7+(6\times 5+4)\times 3\times 21.$
\item [] $2222=9+8\times 7\times 6+5^4\times 3+2\times 1.$
\item [] $2223=9+8\times 7\times 6+5^4\times 3+2+1.$
\item [] $2224=9+8+7\times 6+5\times (432+1).$
\item [] $2225=98+(7+6)\times 54\times 3+21.$
\item [] $2226=(98+7+6)\times 5\times 4+3\times 2\times 1.$
\item [] $2227=(98+7+6)\times 5\times 4+3\times 2+1.$
\item [] $2228=98\times 7+6\times (5+4\times 3\times 21).$
\item [] $2229=9\times (8+7+6)\times 5+4\times 321.$
\item [] $2230=(98+7+6)\times 5\times 4+3^2+1.$
\item[]$\mbox{Increasing order}$
\item [] $2231=(12\times 3+45\times 6)\times 7+89.$
\item [] $2232=(1+2+34\times 5+67+8)\times 9.$
\item [] $2233=1\times (23+4+5)\times 67+89.$
\item [] $2234=1+2+34\times (56+7)+89.$
\item [] $2235=(1+2+3\times 4)\times (5\times 6+7\times (8+9)).$
\item [] $2236=12\times 3+4\times (5+67\times 8+9).$
\item [] $2237=(1+23)\times 45+(6+7)\times 89.$
\item [] $2238=(12+345)\times 6+7+89.$
\item [] $2239=12\times 3\times (4+56)+7+8\times 9.$
\item [] $2240=1+2^3\times 45\times 6+7+8\times 9.$
\item [] $2241=12+3\times (4\times 5\times 6+7\times 89).$
\item [] $2242=12^3+4+5\times (6+7+89).$
\item [] $2243=1\times 2\times 3\times 45\times 6+7\times 89.$
\item [] $2244=1+2\times 3\times 45\times 6+7\times 89.$
\item [] $2245=1\times 2+3+4\times (56+7\times 8\times 9).$
\item [] $2246=123\times (4+5)+67\times (8+9).$
\item [] $2247=12^3+4+5+6+7\times 8\times 9.$
\item [] $2248=1+2^3\times 45\times 6+78+9.$
\item [] $2249=12^3+456+7\times 8+9.$
\item [] $2250=12^3+45+6\times 78+9.$
\item [] $2251=(12+34)\times (5+6\times 7)+89.$
\item [] $2252=1^{23}\times 4\times (5+(6+7\times 8)\times 9).$
\item [] $2253=(1+2)\times (3\times 4\times 56+7+8\times 9).$
\item [] $2254=1+(2+34\times 5)\times (6+7)+8+9.$
\item [] $2255=12+3+4\times (56+7\times 8\times 9).$
\item [] $2256=1\times 2^3\times 45\times 6+7+89.$
\item [] $2257=1+2^3\times 45\times 6+7+89.$
\item [] $2258=12^3+4\times 5+6+7\times 8\times 9.$
\item [] $2259=12^3+4+5+6\times (78+9).$
\item [] $2260=1\times 2\times (3+4^5+67)+8\times 9.$
\item [] $2261=12^3+4+5\times (6+7)\times 8+9.$
\item [] $2262=(1+2+3\times 4+5+6)\times (78+9).$
\item [] $2263=12^3+456+7+8\times 9.$
\item [] $2264=1+23+4\times (56+7\times 8\times 9).$
\item [] $2265=12\times (3+4\times 5\times 6)+789.$
\item [] $2266=12^3+4+5\times 6+7\times 8\times 9.$
\item [] $2267=12\times (3^4+56)+7\times 89.$
\item [] $2268=(123+45+6+78)\times 9.$
\item [] $2269=12+3+4+5\times (6\times 7+8)\times 9.$
\item [] $2270=(1+23+4)\times 56+78\times 9.$
\item [] $2271=12^3+456+78+9.$
\item [] $2272=12^3+4\times (5+6\times 7+89).$
\item [] $2273=(123+45)\times (6+7)+89.$
\item [] $2274=12+3+45\times (6\times 7+8)+9.$
\item [] $2275=1+2\times 3\times 4+5\times (6\times 7+8)\times 9.$
\item [] $2276=12\times 3+4\times (56+7\times 8\times 9).$
\item [] $2277=(1+2)\times (3\times 4\times 56+78+9).$
\item [] $2278=1+23+4+5\times (6\times 7+8)\times 9.$
\item [] $2279=12\times (3+4\times 5)\times 6+7\times 89.$
\item [] $2280=12^3+456+7+89.$
\item [] $2281=1+2\times (345+6+789).$
\item [] $2282=12^3+4+5+67\times 8+9.$
\item [] $2283=12^3+45+6+7\times 8\times 9.$
\item [] $2284=1\times 2+34\times 5\times (6+7)+8\times 9.$
\item [] $2285=1^{23}\times 4\times 567+8+9.$
\item [] $2286=1^{23}+4\times 567+8+9.$
\item [] $2287=1+2\times 3^4\times (5+6)+7\times 8\times 9.$
\item [] $2288=1^2\times 3+4\times 567+8+9.$
\item [] $2289=1^2+3+4\times 567+8+9.$
\item [] $2290=1\times 2+3+4\times 567+8+9.$
\item [] $2291=1+2+3+4\times 567+8+9.$
\item [] $2292=1+2\times 3+4\times 567+8+9.$
\item [] $2293=12^3+4\times 5+67\times 8+9.$
\item [] $2294=1+2^3+4\times 567+8+9.$
\item [] $2295=(1\times 23+45+67)\times (8+9).$
\item [] $2296=1\times 2^3\times 4\times 56+7\times 8\times 9.$
\item [] $2297=1+2^3\times 4\times 56+7\times 8\times 9.$
\item [] $2298=1\times 2\times (3\times 4\times 5\times 6+789).$
\item [] $2299=1+2\times (3\times 4\times 5\times 6+789).$
\item [] $2300=12+3+4\times 567+8+9.$
\item[]$\mbox{Decreasing order}$
\item [] $2231=9+8\times 7+6+5\times 432\times 1.$
\item [] $2232=9+8\times 7+6+5\times 432+1.$
\item [] $2233=9\times 8\times 7+6\times (5+4)\times 32+1.$
\item [] $2234=9+876+5+4^3\times 21.$
\item [] $2235=9+(87+654)\times 3+2+1.$
\item [] $2236=(98+7\times 65)\times 4+3+21.$
\item [] $2237=((9\times 8+7)\times 6+5)\times 4+321.$
\item [] $2238=9\times 8\times 7+6+54\times 32\times 1.$
\item [] $2239=98+7\times 65\times 4+321.$
\item [] $2240=(98\times 7+6+54)\times 3+2\times 1.$
\item [] $2241=(98\times 7+6+54)\times 3+2+1.$
\item [] $2242=9+8\times 76+5\times (4+321).$
\item [] $2243=98+(7+6+54)\times 32+1.$
\item [] $2244=(98+7\times 65)\times 4+32\times 1.$
\item [] $2245=9\times 8+7+6+5\times 432\times 1.$
\item [] $2246=98\times 7+65\times 4\times 3\times 2\times 1.$
\item [] $2247=98\times 7+65\times 4\times 3\times 2+1.$
\item [] $2248=(9+8)\times (7+65)+4^(3+2)\times 1.$
\item [] $2249=9+(8+7\times 6+5\times 4)\times 32\times 1.$
\item [] $2250=9+87\times (6+5)+4\times 321.$
\item [] $2251=(9\times 8+7\times 6\times (5+4))\times (3+2)+1.$
\item [] $2252=(9+87+654)\times 3+2\times 1.$
\item [] $2253=987+6+5\times 4\times 3\times 21.$
\item [] $2254=9+8+76+5\times 432+1.$
\item [] $2255=987+(6+5^4+3)\times 2\times 1.$
\item [] $2256=987+(6+5^4+3)\times 2+1.$
\item [] $2257=9\times 8+7\times 6\times (5\times 4+32)+1.$
\item [] $2258=9+8+76+5\times (432+1).$
\item [] $2259=9+87\times 6+54\times 32\times 1.$
\item [] $2260=9+87\times 6+54\times 32+1.$
\item [] $2261=(9+8+7+6+5)\times 4^3+21.$
\item [] $2262=9+87+6+5\times 432\times 1.$
\item [] $2263=9+87+6+5\times 432+1.$
\item [] $2264=9\times (8+7\times 6)\times 5+4\times 3+2\times 1.$
\item [] $2265=9\times (8+7\times 6)\times 5+4\times 3+2+1.$
\item [] $2266=(9+8\times (7+6))\times 5\times 4+3+2+1.$
\item [] $2267=(98\times 7+65+4)\times 3+2\times 1.$
\item [] $2268=9+(8+7)\times 65+4\times 321.$
\item [] $2269=(9+8\times (7+6))\times 5\times 4+3\times (2+1).$
\item [] $2270=9\times (8+7+65+4)\times 3+2\times 1.$
\item [] $2271=98+7+6+5\times 432\times 1.$
\item [] $2272=98+7+6+5\times 432+1.$
\item [] $2273=9+8\times (7+6)+5\times 432\times 1.$
\item [] $2274=9\times 8+7\times 6+5\times 432\times 1.$
\item [] $2275=9\times 8+7\times 6+5\times 432+1.$
\item [] $2276=98+7+6+5\times (432+1).$
\item [] $2277=9+87\times (6+5\times 4)+3+2+1.$
\item [] $2278=9+87\times (6+5\times 4)+3\times 2+1.$
\item [] $2279=987+6\times 5\times 43+2\times 1.$
\item [] $2280=987+6\times 5\times 43+2+1.$
\item [] $2281=(9\times 8\times 7+65)\times 4+3+2\times 1.$
\item [] $2282=987+6+5+4\times 321.$
\item [] $2283=9\times (8+7\times 6)\times 5+4\times 3+21.$
\item [] $2284=(9+8\times (7+6))\times 5\times 4+3+21.$
\item [] $2285=9+8+7\times 6\times (5+4)\times 3\times 2\times 1.$
\item [] $2286=(9+8)\times 7+6+5\times 432+1.$
\item [] $2287=987+65\times 4\times (3+2)\times 1.$
\item [] $2288=9\times 8+7+(65+4)\times 32+1.$
\item [] $2289=987+6+54\times (3+21).$
\item [] $2290=9+8+7\times 6\times 54+3+2\times 1.$
\item [] $2291=9+8+7\times 6\times 54+3+2+1.$
\item [] $2292=9+8+7\times 6\times 54+3\times 2+1.$
\item [] $2293=(9\times 8+7)\times (6+5\times 4+3)+2\times 1.$
\item [] $2294=9+8+7\times 6\times 54+3^2\times 1.$
\item [] $2295=9+8+7\times 6\times 54+3^2+1.$
\item [] $2296=9\times (8+7\times 6)\times 5+43+2+1.$
\item [] $2297=(98+7\times 6+5^4)\times 3+2\times 1.$
\item [] $2298=987+6\times 5\times 43+21.$
\item [] $2299=9\times (8+7\times 6)\times 5+(4+3)^2\times 1.$
\item [] $2300=98+7\times 6+5\times 432\times 1.$
\item[]$\mbox{Increasing order}$
\item [] $2301=1\times 2+34\times 5\times (6+7)+89.$
\item [] $2302=1+(23+4)\times 56+789.$
\item [] $2303=12^3+456+7\times (8+9).$
\item [] $2304=123+4^5+(6+7)\times 89.$
\item [] $2305=1+(2+3)\times 456+7+8+9.$
\item [] $2306=12^3+4\times 5+(6+7\times 8)\times 9.$
\item [] $2307=12\times 3\times 45+678+9.$
\item [] $2308=1\times 23+4\times 567+8+9.$
\item [] $2309=1+23+4\times 567+8+9.$
\item [] $2310=1+234\times 5+67\times (8+9).$
\item [] $2311=12+34\times 5\times (6+7)+89.$
\item [] $2312=(1+23+45+67)\times (8+9).$
\item [] $2313=(1+2+34\times 5+6+78)\times 9.$
\item [] $2314=(1^2+3+4+5+6+7)\times 89.$
\item [] $2315=12\times (3\times 45+6)+7\times 89.$
\item [] $2316=12^3+4+567+8+9.$
\item [] $2317=1\times 2^3+4\times (567+8)+9.$
\item [] $2318=12^3+45+67\times 8+9.$
\item [] $2319=12^3+456+(7+8)\times 9.$
\item [] $2320=1+2\times 3+4\times (5+67)\times 8+9.$
\item [] $2321=12\times 3+4\times 567+8+9.$
\item [] $2322=1\times 2\times 3\times 45\times 6+78\times 9.$
\item [] $2323=1+2\times 3\times 45\times 6+78\times 9.$
\item [] $2324=12+3+4\times (567+8)+9.$
\item [] $2325=1+2+3\times (45\times 6+7\times 8\times 9).$
\item [] $2326=1+(2+34\times 5)\times (6+7)+89.$
\item [] $2327=(12\times 3+4)\times 56+78+9.$
\item [] $2328=12\times 3\times 45+6+78\times 9.$
\item [] $2329=(1\times 2+3)\times (45+6+7)\times 8+9.$
\item [] $2330=12^3+45\times (6+7)+8+9.$
\item [] $2331=12^3+45+(6+7\times 8)\times 9.$
\item [] $2332=1\times 23+4\times (567+8)+9.$
\item [] $2333=123+(4+5\times 6)\times (7\times 8+9).$
\item [] $2334=12\times 3\times 45+6\times 7\times (8+9).$
\item [] $2335=1+2\times 3\times (45\times 6+7\times (8+9)).$
\item [] $2336=(12\times 3+4)\times 56+7+89.$
\item [] $2337=1+23+4\times (5+67)\times 8+9.$
\item [] $2338=(1+2+34\times 5)\times (6+7)+89.$
\item [] $2339=1234+5\times (6+7)\times (8+9).$
\item [] $2340=1^{23}\times 4\times 567+8\times 9.$
\item [] $2341=1^{23}+4\times 567+8\times 9.$
\item [] $2342=1+2+3+4\times (567+8+9).$
\item [] $2343=1^2\times 3+4\times 567+8\times 9.$
\item [] $2344=1^2+3+4\times 567+8\times 9.$
\item [] $2345=1\times 2+3+4\times 567+8\times 9.$
\item [] $2346=1+2+3+4\times 567+8\times 9.$
\item [] $2347=1+2\times 3+4\times 567+8\times 9.$
\item [] $2348=1\times 2^3+4\times 567+8\times 9.$
\item [] $2349=1+2^3+4\times 567+8\times 9.$
\item [] $2350=1^{234}+5\times 6\times 78+9.$
\item [] $2351=12+3+4\times (567+8+9).$
\item [] $2352=12^3+4\times 5\times 6+7\times 8\times 9.$
\item [] $2353=1^{23}\times 4+5\times 6\times 78+9.$
\item [] $2354=1^{23}+4+5\times 6\times 78+9.$
\item [] $2355=12+3+4\times 567+8\times 9.$
\item [] $2356=1^2\times 3+4+5\times 6\times 78+9.$
\item [] $2357=1^{23}\times 4\times 567+89.$
\item [] $2358=1\times 2+3+4+5\times 6\times 78+9.$
\item [] $2359=1\times 2\times 3+4+5\times 6\times 78+9.$
\item [] $2360=1+2\times 3+4+5\times 6\times 78+9.$
\item [] $2361=1\times 2^3+4+5\times 6\times 78+9.$
\item [] $2362=1+2^3+4+5\times 6\times 78+9.$
\item [] $2363=1\times 2\times 3+4\times 567+89.$
\item [] $2364=1+23+4\times 567+8\times 9.$
\item [] $2365=1\times 2^3+4\times 567+89.$
\item [] $2366=12^3+4+5+6+7\times 89.$
\item [] $2367=1\times (2+3)\times 456+78+9.$
\item [] $2368=12+3+4+5\times 6\times 78+9.$
\item [] $2369=1\times (2+3)\times 4+5\times 6\times 78+9.$
\item [] $2370=1+(2+3)\times 4+5\times 6\times 78+9.$
\item[]$\mbox{Decreasing order}$
\item [] $2301=98+7\times 6+5\times 432+1.$
\item [] $2302=9\times 8+76\times 5+43^2+1.$
\item [] $2303=98\times 7+(65+4\times 3)\times 21.$
\item [] $2304=9+87+(65+4)\times 32\times 1.$
\item [] $2305=(9+87)\times 6+54\times 32+1.$
\item [] $2306=98\times 7+6\times 54\times (3+2)\times 1.$
\item [] $2307=((9+8)\times 7\times 6+54)\times 3+2+1.$
\item [] $2308=9\times 8+76+5\times 432\times 1.$
\item [] $2309=9\times 8+76+5\times 432+1.$
\item [] $2310=987+(6+54+3)\times 21.$
\item [] $2311=(987+6+54\times 3)\times 2+1.$
\item [] $2312=(9\times 8+7\times 6)\times 5\times 4+32\times 1.$
\item [] $2313=98+7+(65+4)\times 32\times 1.$
\item [] $2314=98+7+(65+4)\times 32+1.$
\item [] $2315=9+8\times (7\times 6+54)\times 3+2\times 1.$
\item [] $2316=9\times 87+(6\times 5+43)\times 21.$
\item [] $2317=9+8+7\times 6\times 54+32\times 1.$
\item [] $2318=9+8+7\times 6\times 54+32+1.$
\item [] $2319=9+8\times (7+65)\times 4+3+2+1.$
\item [] $2320=(9+8\times 76+543)\times 2\times 1.$
\item [] $2321=9+8+7\times 65+43^2\times 1.$
\item [] $2322=9+8+7\times 65+43^2+1.$
\item [] $2323=(9\times 8\times 7+654+3)\times 2+1.$
\item [] $2324=(9+8+7+6\times 5)\times 43+2\times 1.$
\item [] $2325=(9+8+7+6\times 5)\times 43+2+1.$
\item [] $2326=9+8+(765+4)\times 3+2\times 1.$
\item [] $2327=98+76\times 5+43^2\times 1.$
\item [] $2328=98+76\times 5+43^2+1.$
\item [] $2329=(9\times 8+7)\times 6+5+43^2+1.$
\item [] $2330=((9+87)\times 6+5)\times 4+3+2+1.$
\item [] $2331=((9+87)\times 6+5)\times 4+3\times 2+1.$
\item [] $2332=9+8\times 76+5\times (4+3)^{(2+1)}.$
\item [] $2333=9+8\times 7+(65+43)\times 21.$
\item [] $2334=98+76+5\times 432\times 1.$
\item [] $2335=98+76+5\times 432+1.$
\item [] $2336=987+65+4\times 321.$
\item [] $2337=9+8\times (7+65)\times 4+3+21.$
\item [] $2338=987+6\times 5\times (43+2)+1.$
\item [] $2339=98+76+5\times (432+1).$
\item [] $2340=9+8\times (76+5\times 43)+2+1.$
\item [] $2341=9\times 8+7\times 6\times (5+4)\times 3\times 2+1.$
\item [] $2342=987+6+5+4^3\times 21.$
\item [] $2343=9\times 87+65\times 4\times 3\times 2\times 1.$
\item [] $2344=9\times 87+65\times 4\times 3\times 2+1.$
\item [] $2345=9+8\times 76+54\times 32\times 1.$
\item [] $2346=9+8\times 76+54\times 32+1.$
\item [] $2347=9\times 8+7\times 6\times 54+3\times 2+1.$
\item [] $2348=9+8+7\times 6\times 54+3\times 21.$
\item [] $2349=9\times 8+7\times 6\times 54+3^2\times 1.$
\item [] $2350=9\times 8+7\times 6\times 54+3^2+1.$
\item [] $2351=(9\times 8+7)\times 6+5^4\times 3+2\times 1.$
\item [] $2352=(9+8+76+5)\times 4\times 3\times 2\times 1.$
\item [] $2353=(9+8+76+5)\times 4\times 3\times 2+1.$
\item [] $2354=9+((87+6)\times 5+4)\times (3+2)\times 1.$
\item [] $2355=98+7+6\times (54+321).$
\item [] $2356=((9+87)\times 6+5)\times 4+32\times 1.$
\item [] $2357=9+8+(7+6)\times 5\times 4\times 3^2\times 1.$
\item [] $2358=9+8\times (76+5\times 43)+21.$
\item [] $2359=(98+7)\times 6+54\times 32+1.$
\item [] $2360=98\times (7+6)+543\times 2\times 1.$
\item [] $2361=987+6\times 5+4^3\times 21.$
\item [] $2362=9+(87+6+5)\times 4\times 3\times 2+1.$
\item [] $2363=9\times 8+7\times (6\times 54+3)+2\times 1.$
\item [] $2364=9+8+7+65\times 4\times 3^2\times 1.$
\item [] $2365=9+8+7+65\times 4\times 3^2+1.$
\item [] $2366=98+7\times 6\times (5+4)\times 3\times 2\times 1.$
\item [] $2367=987+6\times 5\times (43+2+1).$
\item [] $2368=9\times 8+7\times (6+5\times 4^3+2)\times 1.$
\item [] $2369=9+8+7\times (6+5+4+321).$
\item [] $2370=9+8+7\times (6+54\times 3)\times 2+1.$
\item[]$\mbox{Increasing order}$
\item [] $2371=12^3+4+567+8\times 9.$
\item [] $2372=12+3+4\times 567+89.$
\item [] $2373=12+3\times 4+5\times 6\times 78+9.$
\item [] $2374=1+2\times 3\times 4+5\times 6\times 78+9.$
\item [] $2375=1\times 2345+6+7+8+9.$
\item [] $2376=12\times 3+4\times 567+8\times 9.$
\item [] $2377=12^3+4\times 5+6+7\times 89.$
\item [] $2378=1234+5+67\times (8+9).$
\item [] $2379=1\times 2+3+4+5\times 6\times (7+8\times 9).$
\item [] $2380=1\times 23+4\times 567+89.$
\item [] $2381=1+23+4\times 567+89.$
\item [] $2382=1+2^3\times 4+5\times 6\times 78+9.$
\item [] $2383=1^2\times 34+5\times 6\times 78+9.$
\item [] $2384=1^2+34+5\times 6\times 78+9.$
\item [] $2385=1\times 2+34+5\times 6\times 78+9.$
\item [] $2386=1+2+34+5\times 6\times 78+9.$
\item [] $2387=12^3+4+5\times (6\times 7+89).$
\item [] $2388=12^3+4+567+89.$
\item [] $2389=12\times 3+4+5\times 6\times 78+9.$
\item [] $2390=1^{23}+4+5\times (6\times 78+9).$
\item [] $2391=1+(2+3)\times 4+5\times 6\times (7+8\times 9).$
\item [] $2392=1\times 23\times (45+6\times 7+8+9).$
\item [] $2393=12\times 3+4\times 567+89.$
\item [] $2394=1+2^3\times 4\times (5+67)+89.$
\item [] $2395=12+34+5\times 6\times 78+9.$
\item [] $2396=1234+5+(6+7)\times 89.$
\item [] $2397=1\times 23+4+5\times 6\times (7+8\times 9).$
\item [] $2398=1^2+3\times 4+5\times (6\times 78+9).$
\item [] $2399=1\times 2+3\times 4+5\times (6\times 78+9).$
\item [] $2400=12\times (34\times 5+6+7+8+9).$
\item [] $2401=1\times 23\times 4\times (5+6+7+8)+9.$
\item [] $2402=12^3+45+6+7\times 89.$
\item [] $2403=12^3+(4+5)\times 67+8\times 9.$
\item [] $2404=1\times 2345+6\times 7+8+9.$
\item [] $2405=1+2345+6\times 7+8+9.$
\item [] $2406=1\times 2+34+5\times 6\times (7+8\times 9).$
\item [] $2407=1+2+34+5\times 6\times (7+8\times 9).$
\item [] $2408=123+4\times 567+8+9.$
\item [] $2409=1\times 2\times 3\times 45\times 6+789.$
\item [] $2410=1+2\times 3\times 45\times 6+789.$
\item [] $2411=12^3+4+56+7\times 89.$
\item [] $2412=1^2\times 3\times (4+5+6+789).$
\item [] $2413=1^2+3\times (4+5+6+789).$
\item [] $2414=1\times 2+3\times (4+5+6+789).$
\item [] $2415=1\times 2^3\times 4\times 56+7\times 89.$
\item [] $2416=1\times 2345+6+7\times 8+9.$
\item [] $2417=1\times 2\times 34+5\times 6\times 78+9.$
\item [] $2418=1+2\times 34+5\times 6\times 78+9.$
\item [] $2419=1\times 2+(3+4)\times 5\times 67+8\times 9.$
\item [] $2420=12+34\times 56+7\times 8\times 9.$
\item [] $2421=(1+2)\times (3+4+5+6+789).$
\item [] $2422=1+2+34+5\times (6\times 78+9).$
\item [] $2423=12^3+4\times 5+(67+8)\times 9.$
\item [] $2424=12^3+4+5+678+9.$
\item [] $2425=12\times 3+4+5\times (6\times 78+9).$
\item [] $2426=1+2+34\times (56+7+8)+9.$
\item [] $2427=12^3+4+5\times (67+8\times 9).$
\item [] $2428=1+(2+3+4\times 5+6)\times 78+9.$
\item [] $2429=1\times 2345+67+8+9.$
\item [] $2430=1+2345+67+8+9.$
\item [] $2431=1+2345+6+7+8\times 9.$
\item [] $2432=1\times 2+3^4+5\times 6\times 78+9.$
\item [] $2433=1+2+3^4+5\times 6\times 78+9.$
\item [] $2434=(12+3+4\times 5)\times 67+89.$
\item [] $2435=12^3+4\times 5+678+9.$
\item [] $2436=123+4\times (5+67)\times 8+9.$
\item [] $2437=1+2+(3+4)\times 5\times 67+89.$
\item [] $2438=1\times 2345+6+78+9.$
\item [] $2439=1+2345+6+78+9.$
\item [] $2440=1+(2\times 34\times 5+6)\times 7+8+9.$
\item[]$\mbox{Decreasing order}$
\item [] $2371=98+7\times 6\times 54+3+2\times 1.$
\item [] $2372=9\times 8+7\times 6\times 54+32\times 1.$
\item [] $2373=9\times 8+7\times 6\times 54+32+1.$
\item [] $\mathit{2374=98+7\times 6\times 54+3^2-1.}$
\item [] $2375=98+7\times 6\times 54+3^2\times 1.$
\item [] $2376=98+7\times 6\times 54+3^2+1.$
\item [] $2377=9\times 8+7\times 65+43^2+1.$
\item [] $2378=(9+8)\times 76+543\times 2\times 1.$
\item [] $2379=(9+8)\times 76+543\times 2+1.$
\item [] $2380=(9\times 8+7)\times 6\times 5+4+3+2+1.$
\item [] $2381=(9\times 8\times 7+6+5)\times 4+321.$
\item [] $2382=9\times 8+(765+4)\times 3+2+1.$
\item [] $2383=9\times 8\times 7+6\times 5+43^2\times 1.$
\item [] $2384=9\times 8\times 7+6\times 5+43^2+1.$
\item [] $2385=9+87\times 6+5+43^2\times 1.$
\item [] $2386=9+87\times 6+5+43^2+1.$
\item [] $2387=9\times 8\times 7+6+5^4\times 3+2\times 1.$
\item [] $2388=9\times 8\times 7+6+5^4\times 3+2+1.$
\item [] $2389=98+7\times (6\times 54+3)+2\times 1.$
\item [] $2390=98+7\times 6\times 54+3+21.$
\item [] $2391=(9\times 8+7)\times 6\times 5+4\times (3+2)+1.$
\item [] $2392=9\times (87+65)+4^(3+2)\times 1.$
\item [] $2393=98+7\times 6\times 54+3^{(2+1)}.$
\item [] $2394=(9+87+6+5+4+3)\times 21.$
\item [] $2395=987+(6+5)\times 4\times 32\times 1.$
\item [] $2396=987+65+4^3\times 21.$
\item [] $2397=(9\times 87+6+5+4)\times 3+2+1.$
\item [] $2398=98+7\times 6\times 54+32\times 1.$
\item [] $2399=98+7\times 6\times 54+32+1.$
\item [] $2400=9\times 8+(765+4)\times 3+21.$
\item [] $2401=98+7\times (6\times 54+3+2)\times 1.$
\item [] $2402=98+7\times 65+43^2\times 1.$
\item [] $2403=98+7\times 65+43^2+1.$
\item [] $2404=(98+7+6)\times 5+43^2\times 1.$
\item [] $2405=9+8\times 7+65\times 4\times 3^2\times 1.$
\item [] $2406=9\times 8\times 7+6+5^4\times 3+21.$
\item [] $2407=98+(765+4)\times 3+2\times 1.$
\item [] $2408=9+87\times 6+5^4\times 3+2\times 1.$
\item [] $2409=9+87\times 6+5^4\times 3+2+1.$
\item [] $2410=9+8\times (7+6)\times 5\times 4+321.$
\item [] $2411=(98+76\times 5+4)\times (3+2)+1.$
\item [] $2412=9\times 8+(7+6)\times 5\times 4\times 3^2\times 1.$
\item [] $2413=9\times 8+(7+6)\times 5\times 4\times 3^2+1.$
\item [] $2414=9+8+7\times 6\times (54+3)+2+1.$
\item [] $2415=(9\times 87+6+5+4)\times 3+21.$
\item [] $2416=(987+6+5\times 43)\times 2\times 1.$
\item [] $2417=(987+6+5\times 43)\times 2+1.$
\item [] $2418=(9\times 8)\times 7+65+43^2\times 1.$
\item [] $2419=9\times 8\times 7+65+43^2+1.$
\item [] $2420=98\times 7+6+54\times 32\times 1.$
\item [] $2421=98\times 7+6+54\times 32+1.$
\item [] $2422=9+(8\times 7+6+5)\times 4\times 3^2+1.$
\item [] $\mathit{2423=-98+7\times 6\times 5\times 4\times 3+2-1.}$
\item [] $2424=9\times 8\times 7+(6+54)\times 32\times 1.$
\item [] $2425=9\times 8\times 7+(6+54)\times 32+1.$
\item [] $2426=9\times 8\times 7+6\times 5\times 4^3+2\times 1.$
\item [] $2427=9+87\times 6+5^4\times 3+21.$
\item [] $2428=(98+7\times 6)\times 5+(4\times 3)^{(2+1)}.$
\item [] $2429=98+7\times 6\times 54+3\times 21.$
\item [] $2430=(9+87)\times 6+5+43^2\times 1.$
\item [] $2431=(9+87)\times 6+5+43^2+1.$
\item [] $2432=9+8+7\times 6\times (54+3)+21.$
\item [] $2433=9\times (8\times 7+6\times 5+4)\times 3+2+1.$
\item [] $2434=(9\times 8+7)\times 6\times 5+43+21.$
\item [] $2435=9+8\times (7+65)+43^2+1.$
\item [] $2436=9+87+65\times 4\times 3^2\times 1.$
\item [] $2437=9+87+65\times 4\times 3^2+1.$
\item [] $2438=98+(7+6)\times 5\times 4\times 3^2\times 1.$
\item [] $2439=(98\times 7+6\times 5\times 4)\times 3+21.$
\item [] $2440=9\times 8+(7+6\times 5)\times (43+21).$
\item[]$\mbox{Increasing order}$
\item [] $2441=1\times 23\times 4+5\times 6\times 78+9.$
\item [] $2442=1+23\times 4+5\times 6\times 78+9.$
\item [] $2443=(9+8)\times 7\times 6+54\times 32+1.$
\item [] $2444=1\times 2345+6\times (7+8)+9.$
\item [] $2445=12^3+4+5+6+78\times 9.$
\item [] $2446=12+(3+4)\times 5\times 67+89.$
\item [] $2447=1\times 2345+6+7+89.$
\item [] $2448=1+2345+6+7+89.$
\item [] $2449=1+2\times 3\times (4+56\times 7)+8\times 9.$
\item [] $2450=1+(2+3+4\times 5+6)\times (7+8\times 9).$
\item [] $2451=(1+2)\times 34+5\times 6\times 78+9.$
\item [] $2452=12^3+4+5\times 6\times (7+8+9).$
\item [] $2453=1\times 2+3^4+5\times 6\times (7+8\times 9).$
\item [] $2454=1^2\times 3^4\times 5\times 6+7+8+9.$
\item [] $2455=1^2+3^4\times 5\times 6+7+8+9.$
\item [] $2456=12^3+4\times 56+7\times 8\times 9.$
\item [] $2457=1+2+3^4\times 5\times 6+7+8+9.$
\item [] $2458=1+2\times 3\times (4+5+6\times 7)\times 8+9.$
\item [] $2459=1\times 2345+6\times 7+8\times 9.$
\item [] $2460=12^3+45+678+9.$
\item [] $2461=1+2\times 3\times (4\times 5+6\times (7\times 8+9)).$
\item [] $2462=12^3+4\times 5+6\times 7\times (8+9).$
\item [] $2463=123+4\times 567+8\times 9.$
\item [] $2464=12^3+4+5\times 6+78\times 9.$
\item [] $2465=(1+2\times 3)^4+5+6\times 7+8+9.$
\item [] $2466=12+3^4\times 5\times 6+7+8+9.$
\item [] $2467=1^2+3^4+5\times (6\times 78+9).$
\item [] $2468=123\times (4+5+6)+7\times 89.$
\item [] $2469=1+2+3^4+5\times (6\times 78+9).$
\item [] $2470=1\times 2345+6+7\times (8+9).$
\item [] $2471=12^3+4\times 5\times 6+7\times 89.$
\item [] $2472=12\times (3\times 45+6+7\times 8+9).$
\item [] $2473=1+(23\times 4+5+6)\times (7+8+9).$
\item [] $2474=12^3+4\times (5+6)+78\times 9.$
\item [] $2475=(12+3)\times (4+5+67+89).$
\item [] $2476=123+4+5\times 6\times 78+9.$
\item [] $2477=1+2345+6\times 7+89.$
\item [] $2478=12+3^4+5\times (6\times 78+9).$
\item [] $2479=1\times 2\times 3+4\times (5+6)\times 7\times 8+9.$
\item [] $2480=123+4\times 567+89.$
\item [] $2481=12^3+45+6+78\times 9.$
\item [] $2482=1+2^3+4\times (5+6)\times 7\times 8+9.$
\item [] $2483=(1+2\times 3)^4+5\times (6+7)+8+9.$
\item [] $2484=1\times 2345+67+8\times 9.$
\item [] $2485=1+2345+67+8\times 9.$
\item [] $2486=1\times 2345+6+(7+8)\times 9.$
\item [] $2487=12^3+45+6\times 7\times (8+9).$
\item [] $2488=1\times 2^3\times (4\times 56+78+9).$
\item [] $2489=1\times 2345+6\times (7+8+9).$
\item [] $2490=123\times 4\times 5+6+7+8+9.$
\item [] $2491=(1+2\times 3)^4+5+6+7+8\times 9.$
\item [] $2492=(1+23+4)\times (5+67+8+9).$
\item [] $2493=12\times 3\times 4+5\times 6\times 78+9.$
\item [] $2494=1\times 2^3\times 4\times 56+78\times 9.$
\item [] $2495=1+2^3\times 4\times 56+78\times 9.$
\item [] $2496=1^2+3^4\times 5\times 6+7\times 8+9.$
\item [] $2497=1\times 2+3^4\times 5\times 6+7\times 8+9.$
\item [] $2498=1+2+3^4\times 5\times 6+7\times 8+9.$
\item [] $2499=(1+2\times 3)^4+5+6+78+9.$
\item [] $2500=(1+2)^3+4\times (5+6)\times 7\times 8+9.$
\item [] $2501=1\times 2345+67+89.$
\item [] $2502=1+2345+67+89.$
\item [] $2503=12+34\times (5\times (6+7)+8)+9.$
\item [] $2504=1\times 23\times (4+5+6)\times 7+89.$
\item [] $2505=(1+2+3+4\times 5+6)\times 78+9.$
\item [] $2506=(1^2+3)^4+5\times (6\times 7+8)\times 9.$
\item [] $2507=12+3^4\times 5\times 6+7\times 8+9.$
\item [] $2508=12\times (34+56+7\times (8+9)).$
\item [] $2509=1^2\times 3^4\times 5\times 6+7+8\times 9.$
\item [] $2510=1^2+3^4\times 5\times 6+7+8\times 9.$
\item[]$\mbox{Decreasing order}$
\item [] $2441=9+(8+7)\times 6\times (5+4)\times 3+2\times 1.$
\item [] $2442=(9+8)\times 7\times 6+54\times 32\times 1.$
\item [] $2443=(9+8)\times 7\times 6+54\times 32+1.$
\item [] $2444=9\times (8\times 7+6+5)\times 4+32\times 1.$
\item [] $2445=9\times 8\times 7+6\times 5\times 4^3+21.$
\item [] $2446=98+7+65\times 4\times 3^2+1.$
\item [] $2447=(9+8)\times 7\times 6+5+(4\times 3)^{(2+1)}.$
\item [] $2448=9\times 8\times 7+6\times 54\times 3\times 2\times 1.$
\item [] $2449=9\times 8\times 7+6\times 54\times 3\times 2+1.$
\item [] $2450=98+7\times (6+54\times 3)\times 2\times 1.$
\item [] $2451=98+7\times (6+54\times 3)\times 2+1.$
\item [] $2452=98\times (7+6+5+4+3)+2\times 1.$
\item [] $2453=98\times (7+6+5+4+3)+2+1.$
\item [] $2454=(9+87)\times 6+5^4\times 3+2+1.$
\item [] $2455=(9\times 8+7)\times 6\times 5+4^3+21.$
\item [] $2456=(9\times 8+7)\times 6\times 5+43\times 2\times 1.$
\item [] $2457=(9+8+7+6+5+4)\times 3\times 21.$
\item [] $2458=987+6\times 5\times (4+3)^2+1.$
\item [] $2459=(9+8)\times 7+65\times 4\times 3^2\times 1.$
\item [] $2460=(9+8)\times 7+65\times 4\times 3^2+1.$
\item [] $2461=9\times 8+(7+6\times 5)\times 4^3+21.$
\item [] $2462=9\times (87+6)+5\times (4+321).$
\item [] $2463=(9+8+7)\times 65+43\times 21.$
\item [] $2464=9+(8+7\times (65+4))\times (3+2)\times 1.$
\item [] $2465=(9+87\times 6+5)\times 4+321.$
\item [] $2466=9+8\times (7+65\times 4)+321.$
\item [] $2467=9+(8\times 76+5)\times 4+3\times 2\times 1.$
\item [] $2468=9\times 8\times 7+654\times 3+2\times 1.$
\item [] $2469=9\times 8\times 7+654\times 3+2+1.$
\item [] $2470=(98\times 7+6+543)\times 2\times 1.$
\item [] $2471=9+8\times 76+5+43^2\times 1.$
\item [] $2472=9+8\times 76+5+43^2+1.$
\item [] $2473=9+8\times 7\times (6+5+4\times 3+21).$
\item [] $2474=98\times 7+6+54\times (32+1).$
\item [] $2475=9+8\times 7\times (6+5)+43^2+1.$
\item [] $2476=9+(87+6\times 54)\times 3\times 2+1.$
\item [] $2477=98\times (7+6+5+4)+321.$
\item [] $2478=9+8\times 7\times (6+5)\times 4+3+2\times 1.$
\item [] $2479=9+8\times 7\times (6+5)\times 4+3\times 2\times 1.$
\item [] $2480=9+8\times 7\times (6+5)\times 4+3\times 2+1.$
\item [] $2481=9+8\times (76+(5+4)\times 3)\times (2+1).$
\item [] $2482=9+8\times 7\times (6+5)\times 4+3^2\times 1.$
\item [] $2483=9+8\times 7\times (6+5)\times 4+3^2+1.$
\item [] $2484=(98+7)\times 6+5+43^2\times 1.$
\item [] $2485=9+(8\times 76+5)\times 4+3+21.$
\item [] $2486=((9\times 8)\times 7+6\times 54)\times 3+2\times 1.$
\item [] $2487=9\times 8\times 7+654\times 3+21.$
\item [] $2488=9\times 8+7\times (65+4)\times (3+2)+1.$
\item [] $2489=(9+8)\times 76+(54+3)\times 21.$
\item [] $2490=(9+8\times 76+5^4+3)\times 2\times 1.$
\item [] $2491=(9+8\times 76+5^4+3)\times 2+1.$
\item [] $2492=(9+8+7+65)\times (4+3+21).$
\item [] $2493=(9+8\times 76+5)\times 4+3+2\times 1.$
\item [] $2494=9+8\times 76+5^4\times 3+2\times 1.$
\item [] $2495=9+8\times 76+5^4\times 3+2+1.$
\item [] $2496=9+8+7\times 6\times (54+3+2)+1.$
\item [] $2497=9+8\times 7\times (6+5)\times 4+3+21.$
\item [] $2498=(9\times 8+7)\times 6\times 5+4\times 32\times 1.$
\item [] $2499=987+6\times (5+4+3)\times 21.$
\item [] $2500=(987+65\times 4+3)\times 2\times 1.$
\item [] $2501=(987+65\times 4+3)\times 2+1.$
\item [] $2502=9\times (8+7\times 6)\times 5+4\times 3\times 21.$
\item [] $2503=(9+87)\times (6+5\times 4)+3\times 2+1.$
\item [] $2504=9\times 8+(7+65+4)\times 32\times 1.$
\item [] $2505=9+8\times 7\times 6+5\times 432\times 1.$
\item [] $2506=9+8\times 7\times 6+5\times 432+1.$
\item [] $2507=9+8\times (76+5)+43^2+1.$
\item [] $2508=(98+7)\times 6+5^4\times 3+2+1.$
\item [] $2509=9+8+7\times (6\times 54+32)\times 1.$
\item [] $2510=9+876+5\times (4+321).$
\item[]$\mbox{Increasing order}$
\item [] $2511=1\times 2+3^4\times 5\times 6+7+8\times 9.$
\item [] $2512=1+2\times 3^4+5\times 6\times 78+9.$
\item [] $2513=12\times (3\times 45+67)+89.$
\item [] $2514=12\times 34\times 5+6\times (7+8\times 9).$
\item [] $2515=1+(2+3^4)\times 5\times 6+7+8+9.$
\item [] $2516=(12\times 3+45+67)\times (8+9).$
\item [] $2517=12\times 34\times 5+6\times 78+9.$
\item [] $2518=1^2+3^4\times 5\times 6+78+9.$
\item [] $2519=123\times 4\times 5+6\times 7+8+9.$
\item [] $2520=1+2+3^4\times 5\times 6+78+9.$
\item [] $2521=12+3^4\times 5\times 6+7+8\times 9.$
\item [] $2522=1\times 2+3\times (45+6+789).$
\item [] $2523=1+2+3\times (45+6+789).$
\item [] $2524=1\times 2^3+4\times (5\times 6+7)\times (8+9).$
\item [] $2525=12\times (3+4\times 5+6)\times 7+89.$
\item [] $2526=1^2\times 3^4\times 5\times 6+7+89.$
\item [] $2527=1^2\times 34\times 56+7\times 89.$
\item [] $2528=1^2+34\times 56+7\times 89.$
\item [] $2529=1\times 2+34\times 56+7\times 89.$
\item [] $2530=1+2+34\times 56+7\times 89.$
\item [] $2531=123\times 4\times 5+6+7\times 8+9.$
\item [] $2532=12^3+4+5+6+789.$
\item [] $2533=1\times 2\times 34\times (5\times 6+7)+8+9.$
\item [] $2534=1\times 2345+(6+7+8)\times 9.$
\item [] $2535=1+2345+(6+7+8)\times 9.$
\item [] $2536=(1+2\times 3)^4+56+7+8\times 9.$
\item [] $2537=1^2\times 3\times 4\times 5\times 6\times 7+8+9.$
\item [] $2538=12+3^4\times 5\times 6+7+89.$
\item [] $2539=12+34\times 56+7\times 89.$
\item [] $2540=1+2+3\times 4\times 5\times 6\times 7+8+9.$
\item [] $2541=12+(345+6)\times 7+8\times 9.$
\item [] $2542=1\times 2\times 3+4\times (5+6+7\times 89).$
\item [] $2543=12^3+4\times 5+6+789.$
\item [] $2544=123\times 4\times 5+67+8+9.$
\item [] $2545=123\times 4\times 5+6+7+8\times 9.$
\item [] $2546=1^2\times (345+6)\times 7+89.$
\item [] $2547=1^2\times 3\times (4+56+789).$
\item [] $2548=1\times 2+(345+6)\times 7+89.$
\item [] $2549=12+3\times 4\times 5\times 6\times 7+8+9.$
\item [] $2550=12^3+4\times 5\times 6+78\times 9.$
\item [] $2551=12^3+4+5\times 6+789.$
\item [] $2552=1+(2^3\times 4+5)\times 67+8\times 9.$
\item [] $2553=123\times 4\times 5+6+78+9.$
\item [] $2554=1+23\times (4+5+6+7+89).$
\item [] $2555=(12+3+4)\times (56+78)+9.$
\item [] $2556=(1+2)\times (3+4+56+789).$
\item [] $2557=(12\times 34+5)\times 6+7+8\times 9.$
\item [] $2558=12+(345+6)\times 7+89.$
\item [] $2559=12+3\times (4+56+789).$
\item [] $2560=1+23+4\times (5+6+7\times 89).$
\item [] $2561=12^3+4\times (5+6)+789.$
\item [] $2562=123\times 4\times 5+6+7+89.$
\item [] $2563=1+2\times 3+4\times (567+8\times 9).$
\item [] $2564=1\times 2^3+4\times (567+8\times 9).$
\item [] $2565=(12\times 34+5)\times 6+78+9.$
\item [] $2566=1\times 2345+(6+7)\times (8+9).$
\item [] $2567=1+2345+(6+7)\times (8+9).$
\item [] $2568=12^3+45+6+789.$
\item [] $2569=1+2\times 3\times 4\times (5+6+7+89).$
\item [] $2570=1+2^3\times 4\times (5+67+8)+9.$
\item [] $2571=12^3+(4+5)\times 6+789.$
\item [] $2572=1+2\times 3^4\times (5+6)+789.$
\item [] $2573=123\times 4\times 5+(6+7)\times 8+9.$
\item [] $2574=123\times 4\times 5+6\times 7+8\times 9.$
\item [] $2575=12^3+4\times 56+7\times 89.$
\item [] $2576=1\times 2+345\times 6+7\times 8\times 9.$
\item [] $2577=1+2+345\times 6+7\times 8\times 9.$
\item [] $2578=1234+56\times (7+8+9).$
\item [] $2579=1\times 23+4\times (567+8\times 9).$
\item [] $2580=1+23+4\times (567+8\times 9).$
\item[]$\mbox{Decreasing order}$
\item [] $2511=9\times 87+6\times (5+4)\times 32\times 1.$
\item [] $2512=9\times 87+6\times (5+4)\times 32+1.$
\item [] $2513=9+8\times 76+5^4\times 3+21.$
\item [] $2514=9+876+543\times (2+1).$
\item [] $2515=((9+8)\times 7\times 6+543)\times 2+1.$
\item [] $2516=(9+8)\times (76+5+4+3\times 21).$
\item [] $2517=9\times 87+6+54\times 32\times 1.$
\item [] $2518=9\times 87+6+54\times 32+1.$
\item [] $2519=(9\times 8+76)\times (5+4\times 3)+2+1.$
\item [] $2520=(9+8\times 76+5)\times 4+32\times 1.$
\item [] $2521=(9+8\times 76+5)\times 4+32+1.$
\item [] $2522=(9+8)\times (7+6\times 5)\times 4+3+2+1.$
\item [] $2523=(9+8)\times (7+6\times 5)\times 4+3\times 2+1.$
\item [] $2524=9+(8\times 76+5)\times 4+3\times 21.$
\item [] $2525=(9\times 8+765+4)\times 3+2\times 1.$
\item [] $2526=(9\times 8+765+4)\times 3+2+1.$
\item [] $2527=98+7\times (6+5\times 4+321).$
\item [] $2528=(9+8+7\times 6+5\times 4)\times 32\times 1.$
\item [] $2529=(9+8+7\times 6+5\times 4)\times 32+1.$
\item [] $2530=98+(7+65+4)\times 32\times 1.$
\item [] $2531=(9+8)\times (76+54)+321.$
\item [] $2532=(9\times 87+6+54)\times 3+2+1.$
\item [] $2533=(98+7\times 65)\times 4+321.$
\item [] $2534=98\times (7+6)+5\times 4\times 3\times 21.$
\item [] $2535=9+87\times (6+5\times 4+3)+2+1.$
\item [] $2536=9+8\times 7\times (6+5)\times 4+3\times 21.$
\item [] $2537=9+8+7\times (6+54)\times 3\times 2\times 1.$
\item [] $2538=9\times (8+7)\times 6+54\times 32\times 1.$
\item [] $2539=9+8+7\times 6\times 5\times 4\times 3+2\times 1.$
\item [] $2540=9+8+7\times 6\times 5\times 4\times 3+2+1.$
\item [] $2541=(98+7+6)\times 5\times 4+321.$
\item [] $2542=(9+8+(7+6)\times 5)\times (4+3^{(2+1)}).$
\item [] $2543=(9+8)\times (7+6\times 5)\times 4+3^{(2+1)}.$
\item [] $2544=(9\times 8+765+4)\times 3+21.$
\item [] $2545=((98+7)\times 6+5)\times 4+3+2\times 1.$
\item [] $2546=98\times 7+6+5+43^2\times 1.$
\item [] $2547=987+65\times 4\times 3\times 2\times 1.$
\item [] $2548=987+65\times 4\times 3\times 2+1.$
\item [] $2549=(98+7\times 6)\times 5+43^2\times 1.$
\item [] $2550=(9\times 87+6+54)\times 3+21.$
\item [] $2551=(9+8\times 76+5)\times 4+3\times 21.$
\item [] $2552=(9+8)\times 76+5\times 4\times 3\times 21.$
\item [] $2553=9+87\times (6+5\times 4+3)+21.$
\item [] $2554=(9\times (8+7)+6)\times 5+43^2\times 1.$
\item [] $2555=(9+8\times 7)\times 6+5\times (432+1).$
\item [] $2556=9+(8+7\times 6\times 5\times 4)\times 3+2+1.$
\item [] $2557=(9+8+7+65\times 4)\times 3^2+1.$
\item [] $2558=9+8+7\times 6\times 5\times 4\times 3+21.$
\item [] $2559=(9\times 87+65+4)\times 3+2+1.$
\item [] $2560=(9+8\times 7+6+5+4)\times 32\times 1.$
\item [] $2561=(9+8\times 7+6+5+4)\times 32+1.$
\item [] $2562=987+(6+5+4^3)\times 21.$
\item [] $2563=98\times (7+6)+5+4\times 321.$
\item [] $2564=(9\times 8\times 7+6)\times 5+4\times 3+2\times 1.$
\item [] $2565=98\times 7+6\times 5+43^2\times 1.$
\item [] $2566=98\times 7+6\times 5+43^2+1.$
\item [] $2567=(9+8)\times (7+6\times 5\times 4+3+21).$
\item [] $2568=(9+87+6+5)\times 4\times 3\times 2\times 1.$
\item [] $2569=98\times 7+6+5^4\times 3+2\times 1.$
\item [] $2570=98\times 7+6+5^4\times 3+2+1.$
\item [] $2571=9\times 87+6+54\times (32+1).$
\item [] $2572=((98+7)\times 6+5)\times 4+32\times 1.$
\item [] $2573=((98+7)\times 6+5)\times 4+32+1.$
\item [] $2574=(98+7+6\times 54)\times 3\times 2\times 1.$
\item [] $2575=9\times (8+7\times 6)\times 5+4+321.$
\item [] $2576=98+7\times 6\times (54+3+2)\times 1.$
\item [] $2577=(9\times 87+65+4)\times 3+21.$
\item [] $2578=(9\times 8\times 7+6)\times 5+4+3+21.$
\item [] $2579=9+8+7\times 6\times (54+3\times 2+1).$
\item [] $2580=9\times 8\times (7+6+5)+4\times 321.$
\item[]$\mbox{Increasing order}$
\item [] $2581=1\times 2^3\times 4\times 56+789.$
\item [] $2582=1+2^3\times 4\times 56+789.$
\item [] $2583=1\times 234+5\times 6\times 78+9.$
\item [] $2584=1+234+5\times 6\times 78+9.$
\item [] $2585=12\times 34\times 5+67\times 8+9.$
\item [] $2586=12+345\times 6+7\times 8\times 9.$
\item [] $2587=1+(2+345)\times 6+7\times 8\times 9.$
\item [] $2588=1\times 2\times 34\times (5\times 6+7)+8\times 9.$
\item [] $2589=(12+3)\times 4\times 5\times 6+789.$
\item [] $2590=12\times 3+4+5\times (6+7\times 8\times 9).$
\item [] $2591=123\times 4\times 5+6\times 7+89.$
\item [] $2592=1^2\times 3\times 4\times 5\times 6\times 7+8\times 9.$
\item [] $2593=1^2+3\times 4\times 5\times 6\times 7+8\times 9.$
\item [] $2594=1\times 2+3\times 4\times 5\times 6\times 7+8\times 9.$
\item [] $2595=1+2+3\times 4\times 5\times 6\times 7+8\times 9.$
\item [] $2596=12+34\times (5+6+7\times 8+9).$
\item [] $2597=(12\times 34+5)\times 6+7\times (8+9).$
\item [] $2598=1\times 234\times (5+6)+7+8+9.$
\item [] $2599=123\times 4\times 5+67+8\times 9.$
\item [] $2600=(12+34)\times 56+7+8+9.$
\item [] $2601=123\times 4\times 5+6+(7+8)\times 9.$
\item [] $2602=1+2\times 3\times 4\times (5\times 6+78)+9.$
\item [] $2603=1\times 2+3\times (4+5+6\times 7)\times (8+9).$
\item [] $2604=12+3\times 4\times 5\times 6\times 7+8\times 9.$
\item [] $2605=1\times 2\times 34\times (5\times 6+7)+89.$
\item [] $2606=1^2\times 34\times 56+78\times 9.$
\item [] $2607=1^2+34\times 56+78\times 9.$
\item [] $2608=1\times 2+34\times 56+78\times 9.$
\item [] $2609=1+2+34\times 56+78\times 9.$
\item [] $2610=1^2+3\times 4\times 5\times 6\times 7+89.$
\item [] $2611=1\times 2+3\times 4\times 5\times 6\times 7+89.$
\item [] $2612=1+2+3\times 4\times 5\times 6\times 7+89.$
\item [] $2613=(1^2+3)\times 456+789.$
\item [] $2614=1^{23}\times 4+5\times 6\times (78+9).$
\item [] $2615=1^2\times 3+4\times (5\times 6+7\times 89).$
\item [] $2616=123\times 4\times 5+67+89.$
\item [] $2617=1^2\times 3+4+5\times 6\times (78+9).$
\item [] $2618=12+34\times 56+78\times 9.$
\item [] $2619=1\times 234+5\times (6\times 78+9).$
\item [] $2620=1+234+5\times (6\times 78+9).$
\item [] $2621=12^3+45\times 6+7\times 89.$
\item [] $2622=1^2\times 3\times 4+5\times 6\times (78+9).$
\item [] $2623=1+2^3+4+5\times 6\times (78+9).$
\item [] $2624=1\times 2+3\times 4+5\times 6\times (78+9).$
\item [] $2625=12\times 3\times (45+6)+789.$
\item [] $2626=1\times 2\times 3+4\times 5\times (6\times 7+89).$
\item [] $2627=1+2\times 3+4\times 5\times (6\times 7+89).$
\item [] $2628=1^2+3+4\times (567+89).$
\item [] $2629=12+3+4+5\times 6\times (78+9).$
\item [] $2630=(12+345+6)\times 7+89.$
\item [] $2631=1+2\times 3+4\times (567+89).$
\item [] $2632=1\times 2+(34+5)\times 67+8+9.$
\item [] $2633=12\times (34\times 5+6\times 7)+89.$
\item [] $2634=123\times (4+5+6)+789.$
\item [] $2635=12+3+4\times 5\times (6\times 7+89).$
\item [] $2636=1+23+4\times (5\times 6+7\times 89).$
\item [] $2637=12^3+4\times 5\times 6+789.$
\item [] $2638=1+23+4+5\times 6\times (78+9).$
\item [] $2639=1\times 234\times (5+6)+7\times 8+9.$
\item [] $2640=12^3+4\times 5\times 6\times 7+8\times 9.$
\item [] $2641=(12+34)\times 56+7\times 8+9.$
\item [] $2642=12+(34+5)\times 67+8+9.$
\item [] $2643=1+2^3\times 4+5\times 6\times (78+9).$
\item [] $2644=1^2\times 34+5\times 6\times (78+9).$
\item [] $2645=1^2+34+5\times 6\times (78+9).$
\item [] $2646=(12+345)\times 6+7\times 8\times 9.$
\item [] $2647=1\times 23+4\times (567+89).$
\item [] $2648=1+23+4\times (567+89).$
\item [] $2649=1+23\times (45+67)+8\times 9.$
\item [] $2650=12\times 3+4+5\times 6\times (78+9).$
\item[]$\mbox{Decreasing order}$
\item [] $2581=(9+8)\times 76+5+4\times 321.$
\item [] $2582=98\times 7+(6+5^4)\times 3+2+1.$
\item [] $2583=(98+7+6+5+4+3)\times 21.$
\item [] $2584=(9+8+7)\times 65+4^(3+2)\times 1.$
\item [] $2585=9+8\times 7\times (6+5+4\times 3)\times 2\times 1.$
\item [] $2586=(9\times 8\times 7+6)\times 5+4+32\times 1.$
\item [] $2587=(9\times 8\times 7+6)\times 5+4\times 3^2+1.$
\item [] $2588=98\times 7+6+5^4\times 3+21.$
\item [] $2589=(9\times 8+76)\times 5+43^2\times 1.$
\item [] $2590=(9\times 8+76)\times 5+43^2+1.$
\item [] $2591=(9+8)\times 7\times 6+5^4\times 3+2\times 1.$
\item [] $2592=9+8\times 7\times 6\times 5+43\times 21.$
\item [] $2593=9\times 8+7\times (6+54)\times 3\times 2+1.$
\item [] $2594=9\times 8+7\times 6\times 5\times 4\times 3+2\times 1.$
\item [] $2595=(9\times 8\times 7+6)\times 5+43+2\times 1.$
\item [] $2596=(9\times 8\times 7+6)\times 5+43+2+1.$
\item [] $2597=(9\times 8\times 7+65)\times 4+321.$
\item [] $2598=(9+8)\times 7\times 6+(5^4+3)\times (2+1).$
\item [] $2599=(9\times 8\times 7+6)\times 5+(4+3)^2\times 1.$
\item [] $2600=98\times 7+65+43^2\times 1.$
\item [] $2601=98\times 7+65+43^2+1.$
\item [] $2602=9+8+76\times (5+4\times 3)\times 2+1.$
\item [] $2603=(98+765+4)\times 3+2\times 1.$
\item [] $2604=9+8+7+6\times 5\times 43\times 2\times 1.$
\item [] $2605=9+8+7+6\times 5\times 43\times 2+1.$
\item [] $2606=9+8+7\times 6\times 54+321.$
\item [] $2607=987+6\times 54\times (3+2)\times 1.$
\item [] $2608=98\times 7+6\times 5\times 4^3+2\times 1.$
\item [] $2609=98\times 7+6\times 5\times 4^3+2+1.$
\item [] $2610=(9+8)\times 7\times 6+5^4\times 3+21.$
\item [] $2611=9\times (8+7\times 6)+5\times 432+1.$
\item [] $2612=987+65\times (4\times 3\times 2+1).$
\item [] $2613=9+876+54\times 32\times 1.$
\item [] $2614=9+876+54\times 32+1.$
\item [] $2615=9+8\times 7+6\times 5\times (4^3+21).$
\item [] $2616=(9\times 8+7+6\times 5)\times 4\times 3\times 2\times 1.$
\item [] $2617=(9\times 8\times 7+6)\times 5+4+3\times 21.$
\item [] $2618=987+6+5\times (4+321).$
\item [] $2619=98+7\times (6+54)\times 3\times 2+1.$
\item [] $2620=98+7\times 6\times 5\times 4\times 3+2\times 1.$
\item [] $2621=98+7\times 6\times 5\times 4\times 3+2+1.$
\item [] $2622=(98+765+4)\times 3+21.$
\item [] $2623=98\times (7+6)+5+4^3\times 21.$
\item [] $2624=9\times (8\times 7+6\times 5)+43^2+1.$
\item [] $2625=987+(6+5\times 4)\times 3\times 21.$
\item [] $2626=9+(8\times (7+6)+5)\times 4\times 3\times 2+1.$
\item [] $2627=98\times 7+6\times 5\times 4^3+21.$
\item [] $2628=9+87\times 6\times 5+4+3+2\times 1.$
\item [] $2629=9+87\times 6\times 5+4+3+2+1.$
\item [] $2630=98\times 7+6\times 54\times 3\times 2\times 1.$
\item [] $2631=98\times 7+6\times 54\times 3\times 2+1.$
\item [] $2632=9+8+765+43^2+1.$
\item [] $2633=9+87\times 6\times 5+4\times 3+2\times 1.$
\item [] $2634=9+87\times 6\times 5+4\times 3+2+1.$
\item [] $2635=(9\times 8+7)\times 6+5\times 432+1.$
\item [] $2636=(9\times 8\times 7+6)\times 5+43\times 2\times 1.$
\item [] $2637=(9\times 8\times 7+6)\times 5+43\times 2+1.$
\item [] $2638=9\times (8+(7+6)\times 5)\times 4+3^2+1.$
\item [] $2639=98+7\times 6\times 5\times 4\times 3+21.$
\item [] $2640=9+87\times 6\times 5+4\times (3+2)+1.$
\item [] $2641=(9+8)\times 76+5+4^3\times 21.$
\item [] $2642=9+8+7\times (6\times (5+4)+321).$
\item [] $2643=9+87\times 6\times 5+4\times 3\times 2\times 1.$
\item [] $2644=9+87\times 6\times 5+4\times 3\times 2+1.$
\item [] $2645=9+8\times 7+6\times 5\times 43\times 2\times 1.$
\item [] $2646=9+8\times 7+6\times 5\times 43\times 2+1.$
\item [] $2647=9+87\times 6\times 5+4+3+21.$
\item [] $2648=(9+8)\times (7\times 6+5)+43^2\times 1.$
\item [] $2649=9\times 8\times 7+65\times (4\times 3+21).$
\item [] $2650=98\times 7+654\times 3+2\times 1.$
\item[]$\mbox{Increasing order}$
\item [] $2651=(1+2)^3+4\times (567+89).$
\item [] $2652=(1+2+3+4\times 5)\times (6+7+89).$
\item [] $2653=1\times 234\times (5+6)+7+8\times 9.$
\item [] $2654=12^3+4\times 56+78\times 9.$
\item [] $2655=(12+34)\times 56+7+8\times 9.$
\item [] $2656=12+34+5\times 6\times (78+9).$
\item [] $2657=12^3+4\times 5\times 6\times 7+89.$
\item [] $2658=(1+2^3\times 45+6)\times 7+89.$
\item [] $2659=123+4\times (5+6+7\times 89).$
\item [] $2660=12\times 3+4\times (567+89).$
\item [] $2661=12\times 3^4+5\times 6\times 7\times 8+9.$
\item [] $2662=1+234\times (5+6)+78+9.$
\item [] $2663=1\times 2\times 34\times 5\times 6+7\times 89.$
\item [] $2664=1+2\times 34\times 5\times 6+7\times 89.$
\item [] $2665=1+2^3\times 45\times 6+7\times 8\times 9.$
\item [] $2666=1+23\times (45+67)+89.$
\item [] $2667=1+2+3+(4+5\times 6)\times 78+9.$
\item [] $2668=1\times 23\times (45+6+7\times 8+9).$
\item [] $2669=12\times 34\times 5+6+7\times 89.$
\item [] $2670=1\times 234\times (5+6)+7+89.$
\item [] $2671=1+234\times (5+6)+7+89.$
\item [] $2672=(12+34)\times 56+7+89.$
\item [] $2673=12^3+4+5+(6+7)\times 8\times 9.$
\item [] $2674=1+(2+3+4)\times (5\times 6+7)\times 8+9.$
\item [] $2675=(1\times 2+3)\times (456+7+8\times 9).$
\item [] $2676=12+3+(4+5\times 6)\times 78+9.$
\item [] $2677=123+4+5\times (6+7\times 8\times 9).$
\item [] $2678=1\times 23+45\times (6\times 7+8+9).$
\item [] $2679=123+4\times (567+8\times 9).$
\item [] $2680=1^2+3\times (45\times 6+7\times 89).$
\item [] $2681=1\times 2+3\times (45\times 6+7\times 89).$
\item [] $2682=1+2+3\times (45\times 6+7\times 89).$
\item [] $2683=1\times 2\times (3+4^5)+6+7\times 89.$
\item [] $2684=1\times 23+(4+5\times 6)\times 78+9.$
\item [] $2685=1+23+(4+5\times 6)\times 78+9.$
\item [] $2686=1^2+(34+5)\times 67+8\times 9.$
\item [] $2687=1\times 2+(34+5)\times 67+8\times 9.$
\item [] $2688=12^3+456+7\times 8\times 9.$
\item [] $2689=1\times (2+3)\times 4\times (56+78)+9.$
\item [] $2690=1\times 2345+6\times 7\times 8+9.$
\item [] $2691=1+2345+6\times 7\times 8+9.$
\item [] $2692=1^2+3^4+5\times 6\times (78+9).$
\item [] $2693=1^2\times 34\times 56+789.$
\item [] $2694=1^2+34\times 56+789.$
\item [] $2695=1\times 2+34\times 56+789.$
\item [] $2696=1^2\times 3+4+5\times 67\times 8+9.$
\item [] $2697=1^2+3+4+5\times 67\times 8+9.$
\item [] $2698=1\times 2+3+4+5\times 67\times 8+9.$
\item [] $2699=1+2+3+4+5\times 67\times 8+9.$
\item [] $2700=12^3+45\times 6+78\times 9.$
\item [] $2701=1\times 2^3+4+5\times 67\times 8+9.$
\item [] $2702=1+2^3+4+5\times 67\times 8+9.$
\item [] $2703=1\times 2+3\times 4+5\times 67\times 8+9.$
\item [] $2704=1+2+3\times 4+5\times 67\times 8+9.$
\item [] $2705=12+345\times 6+7\times 89.$
\item [] $2706=1+(2+345)\times 6+7\times 89.$
\item [] $2707=1^2\times 34\times (5+6)\times 7+89.$
\item [] $2708=12+3+4+5\times 67\times 8+9.$
\item [] $2709=1\times 2+34\times (5+6)\times 7+89.$
\item [] $2710=1+(2+3)\times 4+5\times 67\times 8+9.$
\item [] $2711=(1+2+345)\times 6+7\times 89.$
\item [] $2712=(1+23)\times (4\times 5+6+78+9).$
\item [] $2713=12+3\times 4+5\times 67\times 8+9.$
\item [] $2714=1+2\times 3\times 4+5\times 67\times 8+9.$
\item [] $2715=12\times 34\times 5+(67+8)\times 9.$
\item [] $2716=1\times 23+4+5\times 67\times 8+9.$
\item [] $2717=1+23+4+5\times 67\times 8+9.$
\item [] $2718=1\times (2+34)\times 56+78\times 9.$
\item [] $2719=12+34\times (5+6)\times 7+89.$
\item [] $2720=(1+2)^3+4+5\times 67\times 8+9.$
\item[]$\mbox{Decreasing order}$
\item [] $2651=98\times 7+654\times 3+2+1.$
\item [] $2652=9+87\times 6\times 5+4\times 3+21.$
\item [] $2653=(9+8)\times 7\times (6+5)+4^3\times 21.$
\item [] $2654=(9+8)\times (7+6)\times (5+4+3)+2\times 1.$
\item [] $2655=9+87\times 6\times 5+4+32\times 1.$
\item [] $2656=9+87\times 6\times 5+4+32+1.$
\item [] $2657=9\times 8+76\times (5+4\times 3)\times 2+1.$
\item [] $2658=((9+8)\times 7+6\times 54)\times 3\times 2\times 1.$
\item [] $2659=9\times 8+7+6\times 5\times 43\times 2\times 1.$
\item [] $2660=9\times 8+7+6\times 5\times 43\times 2+1.$
\item [] $2661=9\times 8+7\times 6\times 54+321.$
\item [] $2662=9\times 87+6\times 5+43^2\times 1.$
\item [] $2663=9\times 87+6\times 5+43^2+1.$
\item [] $2664=9+87\times 6\times 5+43+2\times 1.$
\item [] $2665=9+87\times 6\times 5+43+2+1.$
\item [] $2666=9+(876+5+4)\times 3+2\times 1.$
\item [] $2667=9+876+54\times (32+1).$
\item [] $2668=(9+87\times 6)\times 5+4+3^2\times 1.$
\item [] $2669=98\times 7+654\times 3+21.$
\item [] $2670=9\times 8\times 7+6+5\times 432\times 1.$
\item [] $2671=9\times 8\times 7+6+5\times 432+1.$
\item [] $2672=9\times 8+(7+6)\times 5\times 4\times (3^2+1).$
\item [] $2673=9\times 8\times (7+6\times 5)+4+3+2\times 1.$
\item [] $2674=9\times 8\times (7+6\times 5)+4+3+2+1.$
\item [] $2675=9\times 8\times 7+6+5\times (432+1).$
\item [] $2676=9+87+6\times 5\times 43\times 2\times 1.$
\item [] $2677=9+87+6\times 5\times 43\times 2+1.$
\item [] $2678=(9\times 8\times 7+6)\times 5+4\times 32\times 1.$
\item [] $2679=(9\times 8\times 7+6)\times 5+4\times 32+1.$
\item [] $2680=9+8+76\times 5\times (4+3)+2+1.$
\item [] $2681=98+(7+6\times 5+4)\times 3\times 21.$
\item [] $2682=98+76\times (5+4\times 3)\times 2\times 1.$
\item [] $2683=9+87\times 6\times 5+43+21.$
\item [] $2684=(9+876+5+4)\times 3+2\times 1.$
\item [] $2685=98+7+6\times 5\times 43\times 2\times 1.$
\item [] $2686=9+87\times 6\times 5+4+3\times 21.$
\item [] $2687=98+7\times 6\times 54+321.$
\item [] $2688=(9+8+7+6+54)\times 32\times 1.$
\item [] $2689=(9+8+7+6+54)\times 32+1.$
\item [] $2690=98+(7+65)\times 4\times 3^2\times 1.$
\item [] $2691=9+87\times 6+5\times 432\times 1.$
\item [] $2692=9+87\times 6+5\times 432+1.$
\item [] $2693=9+8\times 7+6\times (5+432+1).$
\item [] $2694=9\times (8+7)\times 6+(5^4+3)\times (2+1).$
\item [] $2695=(9\times 8+7)\times 6\times 5+4+321.$
\item [] $2696=98\times 7+6\times 5\times (4+3\times 21).$
\item [] $2697=9\times 87+65+43^2\times 1.$
\item [] $2698=9\times 87+65+43^2+1.$
\item [] $2699=(9+8)\times 7+6\times 5\times 43\times 2\times 1.$
\item [] $2700=(9+8)\times 7+6\times 5\times 43\times 2+1.$
\item [] $2701=9\times 8\times (7+6\times 5)+4+32+1.$
\item [] $2702=9\times 8+7+6\times (5+432)+1.$
\item [] $2703=9\times 87+(6+54)\times 32\times 1.$
\item [] $2704=9+87\times 6\times 5+4^3+21.$
\item [] $2705=9+87\times 6\times 5+43\times 2\times 1.$
\item [] $2706=9+87\times 6\times 5+43\times 2+1.$
\item [] $2707=9\times 8+7+6\times (5+432+1).$
\item [] $2708=987+6+5\times (4+3)^{(2+1)}.$
\item [] $2709=9\times 8\times (7+6\times 5)+43+2\times 1.$
\item [] $2710=9\times 8\times (7+6\times 5)+43+2+1.$
\item [] $2711=9+(8+7)\times (6+54)\times 3+2\times 1.$
\item [] $2712=98+765+43^2\times 1.$
\item [] $2713=98+765+43^2+1.$
\item [] $2714=9+(8+7\times 6)\times 54+3+2\times 1.$
\item [] $2715=987+6\times (5+4)\times 32\times 1.$
\item [] $2716=987+6\times (5+4)\times 32+1.$
\item [] $2717=(9+876+5\times 4)\times 3+2\times 1.$
\item [] $2718=9+(876+5\times 4)\times 3+21.$
\item [] $2719=(9+87\times 6)\times 5+43+21.$
\item [] $2720=9+8\times (7+6\times 54)+3\times 21.$
\item[]$\mbox{Increasing order}$
\item [] $2721=1\times 2^3\times 4+5\times 67\times 8+9.$
\item [] $2722=1+2^3\times 4+5\times 67\times 8+9.$
\item [] $2723=1^2\times 34+5\times 67\times 8+9.$
\item [] $2724=1^2+34+5\times 67\times 8+9.$
\item [] $2725=1\times 2+34+5\times 67\times 8+9.$
\item [] $2726=1+2+34+5\times 67\times 8+9.$
\item [] $2727=12\times 34\times 5+678+9.$
\item [] $2728=1^2+3\times (4\times 5\times 6+789).$
\item [] $2729=12\times 3+4+5\times 67\times 8+9.$
\item [] $2730=1^2+34\times (5+67+8)+9.$
\item [] $2731=12+3+4\times (56+7\times 89).$
\item [] $2732=1+2+34\times (5+67+8)+9.$
\item [] $2733=12\times 3\times (4+5)\times 6+789.$
\item [] $2734=1\times 2+3+4+5\times (67\times 8+9).$
\item [] $2735=12+34+5\times 67\times 8+9.$
\item [] $2736=12\times (3\times 45+6+78+9).$
\item [] $2737=123+4+5\times 6\times (78+9).$
\item [] $2738=1^2+3\times 4+5\times (67\times 8+9).$
\item [] $2739=1\times 23+4\times (56+7\times 89).$
\item [] $2740=1+23+4\times (56+7\times 89).$
\item [] $2741=12^3+4\times 56+789.$
\item [] $2742=1\times 2\times 34\times 5\times 6+78\times 9.$
\item [] $2743=1+2\times 34\times 5\times 6+78\times 9.$
\item [] $2744=(12\times 3+4)\times 56+7\times 8\times 9.$
\item [] $2745=1^2+3+4\times (5+678)+9.$
\item [] $2746=1\times 2+3+4\times (5+678)+9.$
\item [] $2747=123+4\times (567+89).$
\item [] $2748=12\times 34\times 5+6+78\times 9.$
\item [] $2749=(1^2+34)\times 56+789.$
\item [] $2750=1+2\times 3\times 4+5\times (67\times 8+9).$
\item [] $2751=(1+2\times 3)^4+5+6\times 7\times 8+9.$
\item [] $2752=12\times 3+4\times (56+7\times 89).$
\item [] $2753=12\times 3\times (4+5+67)+8+9.$
\item [] $2754=12\times 34\times 5+6\times 7\times (8+9).$
\item [] $2755=1+2\times (34+5+6\times 7)\times (8+9).$
\item [] $2756=12+3+4\times (5+678)+9.$
\item [] $2757=12\times 34+5\times 6\times 78+9.$
\item [] $2758=1+2\times 34+5\times 67\times 8+9.$
\item [] $2759=12\times 3\times 45+67\times (8+9).$
\item [] $2760=1\times 2\times 3\times 456+7+8+9.$
\item [] $2761=1+2\times 3\times 456+7+8+9.$
\item [] $2762=1+2+34+5\times (67\times 8+9).$
\item [] $2763=1\times 2+(3+4)\times 56\times 7+8+9.$
\item [] $2764=1+2+(3+4)\times 56\times 7+8+9.$
\item [] $2765=(12+345)\times 6+7\times 89.$
\item [] $2766=12+3\times (4+5)\times (6+7+89).$
\item [] $2767=1\times 2\times (3\times 456+7)+8+9.$
\item [] $2768=(123+45\times 6)\times 7+8+9.$
\item [] $2769=(1^2+34+5)\times 67+89.$
\item [] $2770=1^2\times 3^4+5\times 67\times 8+9.$
\item [] $2771=1^2+3^4+5\times 67\times 8+9.$
\item [] $2772=1\times 2+3^4+5\times 67\times 8+9.$
\item [] $2773=1^2+345\times 6+78\times 9.$
\item [] $2774=1\times 2+345\times 6+78\times 9.$
\item [] $2775=1+2+345\times 6+78\times 9.$
\item [] $2776=1\times 2^3+4\times (5+678+9).$
\item [] $2777=12\times 3\times 45+(6+7)\times 89.$
\item [] $2778=123+45\times (6\times 7+8+9).$
\item [] $2779=1^2+3\times (4\times 56+78\times 9).$
\item [] $2780=1^{23}\times 4\times 5\times (67+8\times 9).$
\item [] $2781=1\times 23\times 4+5\times 67\times 8+9.$
\item [] $2782=12^3+4^5+6+7+8+9.$
\item [] $2783=1\times 2^3\times 45\times 6+7\times 89.$
\item [] $2784=12+345\times 6+78\times 9.$
\item [] $2785=1+23\times 4\times 5\times 6+7+8+9.$
\item [] $2786=1+2+3+4\times 5\times (67+8\times 9).$
\item [] $2787=12^3+45\times 6+789.$
\item [] $2788=1\times 2^3+4\times 5\times (67+8\times 9).$
\item [] $2789=12\times (3+4\times 56)+7\times 8+9.$
\item [] $2790=(1+2+345)\times 6+78\times 9.$
\item[]$\mbox{Decreasing order}$
\item [] $2721=987+6+54\times 32\times 1.$
\item [] $2722=987+6+54\times 32+1.$
\item [] $2723=(9+8\times 7+6\times 54)\times (3\times 2+1).$
\item [] $2724=9\times 87+6\times 5\times 4^3+21.$
\item [] $2725=(9\times 8+7+6\times 5)\times (4\times 3\times 2+1).$
\item [] $2726=9+8+7\times (6\times 54+3\times 21).$
\item [] $2727=9\times 87+6\times 54\times 3\times 2\times 1.$
\item [] $2728=9\times 87+6\times 54\times 3\times 2+1.$
\item [] $2729=9+8\times (7+6\times 54+3^2\times 1).$
\item [] $2730=9+(8+7)\times (6+54)\times 3+21.$
\item [] $2731=9\times 8\times (7+6\times 5)+4+3\times 21.$
\item [] $2732=98\times 7+6\times (5\times 4+321).$
\item [] $2733=9\times (8\times 7+6+5)\times 4+321.$
\item [] $2734=9\times 8+76\times 5\times (4+3)+2\times 1.$
\item [] $2735=9\times 8+76\times 5\times (4+3)+2+1.$
\item [] $2736=(9+87)\times 6+5\times 432\times 1.$
\item [] $2737=(9+87)\times 6+5\times 432+1.$
\item [] $2738=9+8+(76+5+4)\times 32+1.$
\item [] $2739=9+876+5+43^2\times 1.$
\item [] $2740=9+876+5+43^2+1.$
\item [] $2741=9+(8+7\times 6)\times 54+32\times 1.$
\item [] $2742=9+(8+7\times 6)\times 54+32+1.$
\item [] $2743=9+((8\times 7\times 6+5)\times 4+3)\times 2\times 1.$
\item [] $2744=98+7\times (6+5+4+3)\times 21.$
\item [] $2745=(9+8+7\times 6\times 5)\times 4\times 3+21.$
\item [] $2746=9+(8+76\times 5)\times (4+3)+21.$
\item [] $2747=9\times 87+654\times 3+2\times 1.$
\item [] $2748=9\times 87+654\times 3+2+1.$
\item [] $2749=9\times 8\times (7+6\times 5)+4^3+21.$
\item [] $2750=9\times 8\times (7+6\times 5)+43\times 2\times 1.$
\item [] $2751=9\times 8\times (7+6\times 5)+43\times 2+1.$
\item [] $2752=(9+8+(7+6)\times 5+4)\times 32\times 1.$
\item [] $2753=9\times 8+76\times 5\times (4+3)+21.$
\item [] $2754=9\times 8\times 7+6\times (54+321).$
\item [] $2755=(9+8+7\times 6+5)\times 43+2+1.$
\item [] $2756=9\times (87+6+5+4)\times 3+2\times 1.$
\item [] $2757=9\times (87+6+5+4)\times 3+2+1.$
\item [] $2758=(9+8)\times 7\times (6+5+4\times 3)+21.$
\item [] $2759=9+(8+7+65\times 4)\times (3^2+1).$
\item [] $2760=98+76\times 5\times (4+3)+2\times 1.$
\item [] $2761=98+76\times 5\times (4+3)+2+1.$
\item [] $2762=9+876+5^4\times 3+2\times 1.$
\item [] $2763=9+876+5^4\times 3+2+1.$
\item [] $2764=9+87\times 6\times 5+(4\times 3)^2+1.$
\item [] $2765=98+7\times (6+54+321).$
\item [] $2766=9\times 87+654\times 3+21.$
\item [] $2767=((9+8\times 7)\times 6+5)\times (4+3)+2\times 1.$
\item [] $2768=98+(7+65\times 4)\times (3^2+1).$
\item [] $2769=9\times (8+7)\times (6+5)+4\times 321.$
\item [] $2770=(9+8+(7+6)\times 5\times 4)\times (3^2+1).$
\item [] $2771=(9+8)\times (76+54+32+1).$
\item [] $2772=(9+8+7+65+43)\times 21.$
\item [] $2773=(9+8+7\times 6+5)\times 43+21.$
\item [] $2774=9+8\times 7+(65+4^3)\times 21.$
\item [] $2775=987+6+54\times (32+1).$
\item [] $2776=(98+7+6)\times (5\times 4+3+2)+1.$
\item [] $2777=9+8\times 76+5\times 432\times 1.$
\item [] $2778=9+8\times 76+5\times 432+1.$
\item [] $2779=98+76\times 5\times (4+3)+21.$
\item [] $2780=987+65+(4\times 3)^{(2+1)}.$
\item [] $2781=9+876+5^4\times 3+21.$
\item [] $2782=9+(8\times 76+5)\times 4+321.$
\item [] $2783=(9+87\times 6)\times 5+4\times 32\times 1.$
\item [] $2784=(9+87\times 6)\times 5+4\times 32+1.$
\item [] $2785=(98+7+6+5)\times 4\times 3\times 2+1.$
\item [] $2786=98+7\times 6\times (54+3^2+1).$
\item [] $2787=(9+87)\times (6+5\times 4+3)+2+1.$
\item [] $2788=9\times 8+7+(65+4^3)\times 21.$
\item [] $2789=9+8+7\times (6+5)\times 4\times 3^2\times 1.$
\item [] $2790=(98+7)\times 6+5\times 432\times 1.$
\item[]$\mbox{Increasing order}$
\item [] $2791=1\times 23+4\times (5+678+9).$
\item [] $2792=1+23+4\times (5+678+9).$
\item [] $2793=12\times 34+5\times (6\times 78+9).$
\item [] $2794=1+2\times 34+5\times (67\times 8+9).$
\item [] $2795=123\times (4+5+6+7)+89.$
\item [] $2796=1+2345+(6\times 7+8)\times 9.$
\item [] $2797=(1+2)^3\times 4+5\times 67\times 8+9.$
\item [] $2798=1^2+3+4+5\times (6+7\times 8)\times 9.$
\item [] $2799=(1+2)\times 3\times (4\times 56+78+9).$
\item [] $2800=1\times 2\times 3+4+5\times (6+7\times 8)\times 9.$
\item [] $2801=1\times 2\times 3\times 456+7\times 8+9.$
\item [] $2802=1+2\times 3\times 456+7\times 8+9.$
\item [] $2803=12\times (3+4\times 56)+7+8\times 9.$
\item [] $2804=12\times 3+4\times (5+678+9).$
\item [] $2805=123\times 4\times 5+6\times 7\times 8+9.$
\item [] $2806=1+(2+34)\times 56+789.$
\item [] $2807=12^3+456+7\times 89.$
\item [] $2808=12\times 3\times (4+5+67)+8\times 9.$
\item [] $2809=1+2\times 3\times 4\times (5\times 6+78+9).$
\item [] $2810=(1+2\times 3)^4+56\times 7+8+9.$
\item [] $2811=12^3+4^5+6\times 7+8+9.$
\item [] $2812=(1+2+34)\times (5+6+7\times 8+9).$
\item [] $2813=1+2\times (3\times (456+7)+8+9).$
\item [] $2814=12+3+45\times (6+7\times 8)+9.$
\item [] $2815=1\times 2\times 3\times 456+7+8\times 9.$
\item [] $2816=123+4+5\times 67\times 8+9.$
\item [] $2817=1\times 23+4+5\times (6+7\times 8)\times 9.$
\item [] $2818=1+23\times 4+5\times (67\times 8+9).$
\item [] $2819=1\times 2345+6\times (7+8\times 9).$
\item [] $2820=1+2345+6\times (7+8\times 9).$
\item [] $2821=(1+23)\times 4+5\times (67\times 8+9).$
\item [] $2822=1\times 2345+6\times 78+9.$
\item [] $2823=1\times 2\times 3\times 456+78+9.$
\item [] $2824=1+2\times 3\times 456+78+9.$
\item [] $2825=1\times 23\times 4\times 5\times 6+7\times 8+9.$
\item [] $2826=1+23\times 4\times 5\times 6+7\times 8+9.$
\item [] $2827=1+2+34+5\times (6+7\times 8)\times 9.$
\item [] $2828=12+(3+4)\times 56\times 7+8\times 9.$
\item [] $2829=1\times 2\times 34\times 5\times 6+789.$
\item [] $2830=1+2\times 34\times 5\times 6+789.$
\item [] $2831=(1+2+34+5)\times 67+8+9.$
\item [] $2832=1\times 2\times 3\times 456+7+89.$
\item [] $2833=12\times 3\times 4+5\times 67\times 8+9.$
\item [] $2834=1^2+(3+4)\times 56\times 7+89.$
\item [] $2835=12\times 34\times 5+6+789.$
\item [] $2836=12^3+4^5+67+8+9.$
\item [] $2837=12^3+4^5+6+7+8\times 9.$
\item [] $2838=1+2+3^4\times (5+6+7+8+9).$
\item [] $2839=1\times 23\times 4\times 5\times 6+7+8\times 9.$
\item [] $2840=(123+45\times 6)\times 7+89.$
\item [] $2841=123\times 4+5\times 6\times 78+9.$
\item [] $2842=1+2\times (3\times 45+6\times 7)\times 8+9.$
\item [] $2843=12\times (3+4\times 56)+7\times (8+9).$
\item [] $2844=1\times 234+5\times 6\times (78+9).$
\item [] $2845=12^3+4^5+6+78+9.$
\item [] $2846=1\times 2+3\times 45\times (6+7+8)+9.$
\item [] $2847=1\times 23\times 4\times 5\times 6+78+9.$
\item [] $2848=1+23\times 4\times 5\times 6+78+9.$
\item [] $2849=1+2^3\times 4\times (5+67+8+9).$
\item [] $2850=12\times 34\times 5+6\times (7+8)\times 9.$
\item [] $2851=1\times 2\times 3^4+5\times 67\times 8+9.$
\item [] $2852=1+2\times 3^4+5\times 67\times 8+9.$
\item [] $2853=(12+3+4\times 56+78)\times 9.$
\item [] $2854=12^3+4^5+6+7+89.$
\item [] $2855=1\times 2345+6+7\times 8\times 9.$
\item [] $2856=1+2345+6+7\times 8\times 9.$
\item [] $2857=1+23\times 4\times 5\times 6+7+89.$
\item [] $2858=1+2+3+4\times (5+6+78\times 9).$
\item [] $2859=1^2\times 345\times 6+789.$
\item [] $2860=1^2+345\times 6+789.$
\item[]$\mbox{Decreasing order}$
\item [] $2791=(98+7)\times 6+5\times 432+1.$
\item [] $2792=9\times 8\times (7+6\times 5)+4\times 32\times 1.$
\item [] $2793=(9\times 8+7+6+5+43)\times 21.$
\item [] $2794=9+8\times 7\times (6+5)\times 4+321.$
\item [] $2795=(98\times 7+6+5)\times 4+3\times 2+1.$
\item [] $2796=9+8+7+(6+5)\times 4\times 3\times 21.$
\item [] $2797=(98\times 7+6+5)\times 4+3^2\times 1.$
\item [] $2798=(98\times 7+6+5)\times 4+3^2+1.$
\item [] $2799=9\times (87+6+5\times 43+2+1).$
\item [] $2800=9\times (8\times 7+6)\times 5+4+3\times 2\times 1.$
\item [] $2801=9+(876+54)\times 3+2\times 1.$
\item [] $2802=9+(876+54)\times 3+2+1.$
\item [] $2803=(9\times 8+7)\times 6\times 5+432+1.$
\item [] $2804=9\times (8\times 7+6)\times 5+4\times 3+2\times 1.$
\item [] $2805=9+87+(65+4^3)\times 21.$
\item [] $2806=(98+7\times 6)\times 5\times 4+3\times 2\times 1.$
\item [] $2807=987+65\times (4+3+21).$
\item [] $2808=(98+7\times 6\times 5+4)\times 3^2\times 1.$
\item [] $2809=(9+8\times 76+5)\times 4+321.$
\item [] $2810=9\times (8+7\times 6+54)\times 3+2\times 1.$
\item [] $2811=9\times (8+76+5\times 4)\times 3+2+1.$
\item [] $2812=9+(8+7\times 65+4)\times 3\times 2+1.$
\item [] $2813=9\times 8\times (7+6)+5^4\times 3+2\times 1.$
\item [] $2814=9+8+(7+6)\times 5\times 43+2\times 1.$
\item [] $2815=9+8+(7+6)\times 5\times 43+2+1.$
\item [] $2816=9+87\times (6+5)+43^2+1.$
\item [] $2817=9+(87+6\times 5)\times 4\times 3\times 2\times 1.$
\item [] $2818=98+(76+5+4)\times 32\times 1.$
\item [] $2819=(9+876+54)\times 3+2\times 1.$
\item [] $2820=9+(876+54)\times 3+21.$
\item [] $2821=9+8+7+65\times 43+2\times 1.$
\item [] $2822=9+8+7+65\times 43+2+1.$
\item [] $2823=9\times (8\times 7+6)\times 5+4\times 3+21.$
\item [] $2824=(98+7\times 6)\times 5\times 4+3+21.$
\item [] $2825=(9\times (8+7)+6)\times 5\times 4+3+2\times 1.$
\item [] $2826=(9+876+54+3)\times (2+1).$
\item [] $2827=9\times 8+(7\times 65+4)\times 3\times 2+1.$
\item [] $2828=98+7\times (65+4+321).$
\item [] $2829=9\times (8+7\times 6+54)\times 3+21.$
\item [] $2830=9\times (8\times 7+6)\times 5+4\times (3^2+1).$
\item [] $2831=98\times 7+65\times (4\times 3+21).$
\item [] $2832=(98+7\times 6)\times 5\times 4+32\times 1.$
\item [] $2833=9+8+(7+6)\times 5\times 43+21.$
\item [] $2834=9+(8+7)\times 65+43^2+1.$
\item [] $2835=9\times (8\times 7+6)\times 5+43+2\times 1.$
\item [] $2836=9\times (8\times 7+6)\times 5+43+2+1.$
\item [] $2837=9+8\times 7+(6+5)\times 4\times 3\times 21.$
\item [] $2838=(9+876+54)\times 3+21.$
\item [] $2839=(9+8)\times (76+5+43\times 2\times 1).$
\item [] $2840=9+8+7+65\times 43+21.$
\item [] $2841=(98\times 7+65\times 4)\times 3+2+1.$
\item [] $2842=98\times (7+6+5+4+3\times 2+1).$
\item [] $2843=9+(8+7\times 6\times 5)\times (4+3^2\times 1).$
\item [] $2844=9\times 8\times 7+65\times 4\times 3^2\times 1.$
\item [] $2845=9\times 8\times 7+65\times 4\times 3^2+1.$
\item [] $2846=((9+8\times 7+6)\times 5\times 4+3)\times 2\times 1.$
\item [] $2847=987+6+5+43^2\times 1.$
\item [] $2848=987+6+5+43^2+1.$
\item [] $2849=9+8\times (7+6\times 54+3+21).$
\item [] $2850=(9\times (8+7+6)\times 5+4)\times 3+2+1.$
\item [] $2851=(98\times 7+6+5)\times 4+3\times 21.$
\item [] $2852=98\times 7+6+5\times 432\times 1.$
\item [] $2853=98\times 7+6+5\times 432+1.$
\item [] $2854=9\times (8\times 7+6)\times 5+43+21.$
\item [] $2855=9\times ((8\times 7+6)\times 5+4+3)+2\times 1.$
\item [] $2856=(9+876+543)\times 2\times 1.$
\item [] $2857=(9+876+543)\times 2+1.$
\item [] $2858=9+(8+7\times (6+5)+4)\times 32+1.$
\item [] $2859=(98\times 7+65\times 4)\times 3+21.$
\item [] $2860=(9+8\times 7)\times (6+5+4\times 3+21).$
\item[]$\mbox{Increasing order}$
\item [] $2861=1\times 2+345\times 6+789.$
\item [] $2862=1\times 2^3\times 45\times 6+78\times 9.$
\item [] $2863=1+2^3\times 45\times 6+78\times 9.$
\item [] $2864=123+4\times (5+678)+9.$
\item [] $2865=(12\times 3+4+5+6)\times 7\times 8+9.$
\item [] $2866=12^3+4^5+6\times 7+8\times 9.$
\item [] $2867=1\times 2345+6\times (78+9).$
\item [] $2868=1+2345+6\times (78+9).$
\item [] $2869=12\times 3\times 4+5\times (67\times 8+9).$
\item [] $2870=1\times 2+3+(45+6)\times 7\times 8+9.$
\item [] $2871=12+345\times 6+789.$
\item [] $2872=1+(2+345)\times 6+789.$
\item [] $2873=(1\times 23+45)\times 6\times 7+8+9.$
\item [] $2874=1+2^3+(45+6)\times 7\times 8+9.$
\item [] $2875=1\times 23\times (4+56+7\times 8+9).$
\item [] $2876=12^3+4+5+67\times (8+9).$
\item [] $2877=(1+2+345)\times 6+789.$
\item [] $2878=1+2+(3+4\times 5)\times (6+7\times (8+9)).$
\item [] $2879=1\times 23\times 4\times 5\times 6+7\times (8+9).$
\item [] $2880=(12+3+4+5+6)\times (7+89).$
\item [] $2881=1^2+3\times (456+7\times 8\times 9).$
\item [] $2882=(1+2\times 3)^4+56\times 7+89.$
\item [] $2883=12^3+4^5+6\times 7+89.$
\item [] $2884=1^2+3+4\times 5\times 6\times (7+8+9).$
\item [] $2885=1\times 2+3+4\times 5\times 6\times (7+8+9).$
\item [] $2886=12^3+456+78\times 9.$
\item [] $2887=12\times (34+5)\times 6+7+8\times 9.$
\item [] $2888=(12+3+4)\times (56+7+89).$
\item [] $2889=(1\times 2+345+6+7)\times 8+9.$
\item [] $2890=1\times 2345+67\times 8+9.$
\item [] $2891=1+2345+67\times 8+9.$
\item [] $2892=12+3\times (456+7\times 8\times 9).$
\item [] $2893=1^2+3\times 4+5\times 6\times (7+89).$
\item [] $2894=12^3+4+5+(6+7)\times 89.$
\item [] $2895=12\times (34+5)\times 6+78+9.$
\item [] $2896=1+23\times 4\times 5\times 6+(7+8)\times 9.$
\item [] $2897=12\times (3+4\times 56+7)+89.$
\item [] $2898=(1\times 23+4\times 5)\times 67+8+9.$
\item [] $2899=12+3+4+5\times 6\times (7+89).$
\item [] $2900=(1\times 2+3)\times 4+5\times 6\times (7+89).$
\item [] $2901=12\times (34\times 5+6)+789.$
\item [] $2902=1+2+34\times ((5+6)\times 7+8)+9.$
\item [] $2903=(1+2+34+5)\times 67+89.$
\item [] $2904=1+(2+3)\times 456+7\times 89.$
\item [] $2905=12^3+4\times 5+(6+7)\times 89.$
\item [] $2906=12+(3^4\times 5+6)\times 7+8+9.$
\item [] $2907=12+(3+4+5\times 6)\times 78+9.$
\item [] $2908=12^3+4^5+67+89.$
\item [] $2909=(12+3^4)\times 5\times 6+7\times (8+9).$
\item [] $2910=1\times 2\times (3\times 456+78+9).$
\item [] $2911=1+2\times (3\times 456+78+9).$
\item [] $2912=12^3+45+67\times (8+9).$
\item [] $2913=1+2^3\times 4+5\times 6\times (7+89).$
\item [] $2914=1^2\times 34+5\times 6\times (7+89).$
\item [] $2915=1^2+34+5\times 6\times (7+89).$
\item [] $2916=12\times (34\times 5+67)+8\times 9.$
\item [] $2917=1+2+34+5\times 6\times (7+89).$
\item [] $2918=1\times 2^{(3+4)}+5\times (6+7\times 8)\times 9.$
\item [] $2919=1+2^{(3+4)}+5\times (6+7\times 8)\times 9.$
\item [] $2920=12\times 3+4+5\times 6\times (7+89).$
\item [] $2921=1\times 2345+6\times (7+89).$
\item [] $2922=1+2345+6\times (7+89).$
\item [] $2923=1234+5\times 6\times 7\times 8+9.$
\item [] $2924=1+234+5\times 67\times 8+9.$
\item [] $2925=(1\times 23+4\times 56+78)\times 9.$
\item [] $2926=12+34+5\times 6\times (7+89).$
\item [] $2927=12\times (34+5)\times 6+7\times (8+9).$
\item [] $2928=(1\times 23+45)\times 6\times 7+8\times 9.$
\item [] $2929=1+(23+45)\times 6\times 7+8\times 9.$
\item [] $2930=12^3+45+(6+7)\times 89.$
\item[]$\mbox{Decreasing order}$
\item [] $2861=((98+7)\times 6+5)\times 4+321.$
\item [] $2862=9+8\times 7+65\times 43+2\times 1.$
\item [] $2863=9+8\times 7+65\times 43+2+1.$
\item [] $2864=(9\times (87+65)+4^3)\times 2\times 1.$
\item [] $2865=(9\times 8\times 7+65+4)\times (3+2)\times 1.$
\item [] $2866=987+6\times 5+43^2\times 1.$
\item [] $2867=987+6\times 5+43^2+1.$
\item [] $2868=9+87+(6+5)\times 4\times 3\times 21.$
\item [] $2869=9\times 8+(7+6)\times 5\times 43+2\times 1.$
\item [] $2870=987+6+5^4\times 3+2\times 1.$
\item [] $2871=9+87\times 6\times 5+4\times 3\times 21.$
\item [] $2872=(98\times (7+6)+54\times 3)\times 2\times 1.$
\item [] $2873=(98\times 7+6\times 5)\times 4+3\times (2+1).$
\item [] $2874=(9+8)\times 7\times 6+5\times 432\times 1.$
\item [] $2875=(9+8)\times 7\times 6+5\times 432+1.$
\item [] $2876=9\times 8+7+65\times 43+2\times 1.$
\item [] $2877=9\times 8+7+65\times 43+2+1.$
\item [] $2878=((9+8)\times 7+6)\times (5\times 4+3)+2+1.$
\item [] $2879=(9+8)\times 7\times 6+5\times (432+1).$
\item [] $2880=(9\times 87+654+3)\times 2\times 1.$
\item [] $2881=9+8\times 7+65\times 43+21.$
\item [] $2882=987+(6+5^4)\times 3+2\times 1.$
\item [] $2883=987+(6+5^4)\times 3+2+1.$
\item [] $2884=((9+8\times 7)\times (6+5)\times 4+3+21).$
\item [] $2885=9\times (8+7+65)\times 4+3+2\times 1.$
\item [] $2886=9\times (8+7+65)\times 4+3\times 2\times 1.$
\item [] $2887=(9+8+7)\times 6\times 5\times 4+3\times 2+1.$
\item [] $2888=9\times 8+(7+6)\times 5\times 43+21.$
\item [] $2889=987+6+5^4\times 3+21.$
\item [] $2890=(9+87)\times 6\times 5+4+3\times 2\times 1.$
\item [] $2891=(9+87)\times 6\times 5+4+3\times 2+1.$
\item [] $2892=(98+76\times 5+4)\times 3\times 2\times 1.$
\item [] $2893=9+87+65\times 43+2\times 1.$
\item [] $2894=9+87+65\times 43+2+1.$
\item [] $2895=9\times 8+7+65\times 43+21.$
\item [] $2896=98+(7+6)\times 5\times 43+2+1.$
\item [] $2897=(98\times 7+6\times 5)\times 4+32+1.$
\item [] $2898=(9+8\times 7+6\times 5+43)\times 21.$
\item [] $2899=98\times (7+6)+5\times (4+321).$
\item [] $2900=(98+7+6+5)\times (4\times 3\times 2+1).$
\item [] $2901=987+65+43^2\times 1.$
\item [] $2902=98+7+65\times 43+2\times 1.$
\item [] $2903=98+7+65\times 43+2+1.$
\item [] $2904=9\times (8+7+65)\times 4+3+21.$
\item [] $2905=(9+87)\times 6\times 5+4\times 3\times 2+1.$
\item [] $2906=(9+87)\times (6+5)+43^2+1.$
\item [] $2907=987+(6+54)\times 32\times 1.$
\item [] $2908=987+(6+54)\times 32+1.$
\item [] $2909=987+6\times 5\times 4^3+2\times 1.$
\item [] $2910=987+6\times 5\times 4^3+2+1.$
\item [] $2911=9+(8\times 7+6+5)\times 43+21.$
\item [] $2912=9+87+65\times 43+21.$
\item [] $2913=(9+87)\times 6\times 5+4\times 3+21.$
\item [] $2914=98+(7+6)\times 5\times 43+21.$
\item [] $2915=9+8+7\times (65+4)\times 3\times 2\times 1.$
\item [] $2916=(9+87)\times 6\times 5+4+32\times 1.$
\item [] $2917=(9+87)\times 6\times 5+4+32+1.$
\item [] $2918=9\times (8\times 7+6)\times 5+4\times 32\times 1.$
\item [] $2919=9+8+7+6+(5+4)\times 321.$
\item [] $2920=(9+87)\times 6\times 5+4\times (3^2+1).$
\item [] $2921=98+7+65\times 43+21.$
\item [] $2922=9+8\times (76+5)\times 4+321.$
\item [] $2923=((9+8\times 7)\times (6+5)\times 4+3\times 21).$
\item [] $2924=(9+8)\times ((7+6\times 5)\times 4+3+21).$
\item [] $2925=9\times 87+(6\times 5+4)\times 3\times 21.$
\item [] $2926=(9+87)\times 6\times 5+43+2+1.$
\item [] $2927=(98\times 7+6\times 5)\times 4+3\times 21.$
\item [] $2928=987+6\times 5\times 4^3+21.$
\item [] $2929=9+8\times (7\times 6+5\times 4^3+2+1).$
\item [] $2930=(9+87)\times 6\times 5+(4+3)^2+1.$
\item[]$\mbox{Increasing order}$
\item [] $2931=(12+345)\times 6+789.$
\item [] $2932=(1+2)\times 3^4+5\times 67\times 8+9.$
\item [] $2933=12\times (34\times 5+67)+89.$
\item [] $2934=1^2\times 3^4\times 5\times 6+7\times 8\times 9.$
\item [] $2935=1^2+3^4\times 5\times 6+7\times 8\times 9.$
\item [] $2936=1\times 2+3^4\times 5\times 6+7\times 8\times 9.$
\item [] $2937=123\times 4\times 5+6\times 78+9.$
\item [] $2938=(1+2+34)\times (5+6)\times 7+89.$
\item [] $2939=1+23\times (4\times 5\times 6+7)+8+9.$
\item [] $2940=12\times (3+4)\times (5+6+7+8+9).$
\item [] $2941=1+(2+3)\times (4+567+8+9).$
\item [] $2942=(12\times 3+4)\times 56+78\times 9.$
\item [] $2943=12+3+4\times (5\times 6+78\times 9).$
\item [] $2944=(12+34)\times (5+6\times 7+8+9).$
\item [] $2945=(1\times 23+45)\times 6\times 7+89.$
\item [] $2946=1+(23+45)\times 6\times 7+89.$
\item [] $2947=1+(2+3^4\times 5)\times 6+7\times 8\times 9.$
\item [] $2948=1\times 2\times 34+5\times 6\times (7+89).$
\item [] $2949=1\times 2^3\times 45\times 6+789.$
\item [] $2950=1+2^3\times 45\times 6+789.$
\item [] $2951=(1+2\times 3)^4+5+67\times 8+9.$
\item [] $2952=1+23+4\times (5\times 6+78\times 9).$
\item [] $2953=(1\times 23+4\times 5)\times 67+8\times 9.$
\item [] $2954=1+(23+4\times 5)\times 67+8\times 9.$
\item [] $2955=(1+2^3\times 45)\times 6+789.$
\item [] $2956=(1^2+3^4\times 5+6)\times 7+8\times 9.$
\item [] $2957=1\times 2\times (3+4)\times 5\times 6\times 7+8+9.$
\item [] $2958=(1\times 2^3+4\times 5+6)\times (78+9).$
\item [] $2959=1\times 234+5\times (67\times 8+9).$
\item [] $2960=1+234+5\times (67\times 8+9).$
\item [] $2961=(1+2\times 3)^4+56+7\times 8\times 9.$
\item [] $2962=1^2+3^4+5\times 6\times (7+89).$
\item [] $2963=1\times 2+3^4+5\times 6\times (7+89).$
\item [] $2964=12\times 3+4\times (5\times 6+78\times 9).$
\item [] $2965=(1+23+4\times 5)\times 67+8+9.$
\item [] $2966=(1^2\times 3^4\times 5+6)\times 7+89.$
\item [] $2967=1\times 23\times (45+67+8+9).$
\item [] $2968=1+23\times (45+67+8+9).$
\item [] $2969=(12\times 3+4)\times (5+67)+89.$
\item [] $2970=123\times 4\times 5+6+7\times 8\times 9.$
\item [] $2971=1+(23+4\times 5)\times 67+89.$
\item [] $2972=1\times 23\times 4+5\times 6\times (7+89).$
\item [] $2973=12^3+456+789.$
\item [] $2974=1\times 2345+6+7\times 89.$
\item [] $2975=1+2345+6+7\times 89.$
\item [] $2976=12\times 34\times 5+(6+7)\times 8\times 9.$
\item [] $2977=12^3+4\times 5\times (6+7\times 8)+9.$
\item [] $2978=12+(3^4\times 5+6)\times 7+89.$
\item [] $2979=12^3+(4+5)\times (67+8\times 9).$
\item [] $2980=12+(3+4)\times (5\times 67+89).$
\item [] $2981=1\times 23+(4+5\times 6)\times (78+9).$
\item [] $2982=1\times (2+3)\times 456+78\times 9.$
\item [] $2983=1+(2+3)\times 456+78\times 9.$
\item [] $2984=((1^2+3^4)\times 5+6)\times 7+8\times 9.$
\item [] $2985=(1+2\times 3)^4+567+8+9.$
\item [] $2986=1+(2\times 3)^4+5\times 6\times 7\times 8+9.$
\item [] $2987=(12+34)\times (56+7)+89.$
\item [] $2988=123+(45+6)\times 7\times 8+9.$
\item [] $2989=1+2\times 3^4\times (5+6+7)+8\times 9.$
\item [] $2990=1\times 23\times (45+6+7+8\times 9).$
\item [] $2991=1+23\times (45+6+7+8\times 9).$
\item [] $2992=1\times 2^3\times (4+5+6+7)\times (8+9).$
\item [] $2993=1\times 23\times (4\times 5\times 6+7)+8\times 9.$
\item [] $2994=1+23\times (4\times 5\times 6+7)+8\times 9.$
\item [] $2995=1\times (2+3)^4+5\times 6\times (7+8\times 9).$
\item [] $2996=(1+23+4)\times (5+6+7+89).$
\item [] $2997=12+(3+45)\times (6+7\times 8)+9.$
\item [] $2998=1+(2+3)\times 45\times (6+7)+8\times 9.$
\item [] $2999=(1+23)\times 4\times 5\times 6+7\times (8+9).$
\item [] $3000=12\times (34+5\times 6\times 7)+8\times 9.$
\item[]$\mbox{Decreasing order}$
\item [] $2931=987+6\times 54\times 3\times 2\times 1.$
\item [] $2932=987+6\times 54\times 3\times 2+1.$
\item [] $2933=9\times (8+7)+65\times 43+2+1.$
\item [] $2934=9+8+(76+5)\times 4\times 3^2+1.$
\item [] $2935=(9+8)\times 7+65\times 43+21.$
\item [] $2936=98\times 7+6\times (54+321).$
\item [] $2937=987+6\times (54\times 3\times 2+1).$
\item [] $2938=((9+8)\times (7+6)+5)\times (4+3^2\times 1).$
\item [] $2939=((9+8)\times 7\times 6+5)\times 4+3\times 21.$
\item [] $2940=9+8+7+6\times 54\times 3^2\times 1.$
\item [] $2941=9+8+7+6\times 54\times 3^2+1.$
\item [] $2942=98\times (7+6+5+4\times 3)+2\times 1.$
\item [] $2943=(9+8+7)\times 6\times 5\times 4+3\times 21.$
\item [] $2944=9+87\times 6\times 5+4+321.$
\item [] $2945=9+8\times 7+6\times 5\times 4\times (3+21).$
\item [] $2946=(9+87)\times 6\times 5+4^3+2\times 1.$
\item [] $2947=(9+87)\times 6\times 5+4+3\times 21.$
\item [] $2948=9+8+7\times 6+(5+4)\times 321.$
\item [] $2949=9\times 87+6+5\times 432\times 1.$
\item [] $2950=9\times 87+6+5\times 432+1.$
\item [] $2951=987+654\times 3+2\times 1.$
\item [] $2952=987+654\times 3+2+1.$
\item [] $2953=(9\times 8+7\times (6+54))\times 3\times 2+1.$
\item [] $2954=9\times 87+6+5\times (432+1).$
\item [] $2955=(9+8)\times (7+6)\times 5+43^2+1.$
\item [] $2956=9+(8+7\times (65+4))\times 3\times 2+1.$
\item [] $2957=9+8+7\times 6\times 5\times (4+3)\times 2\times 1.$
\item [] $2958=9+8+7\times 6\times 5\times (4+3)\times 2+1.$
\item [] $2959=9\times 8+7+6\times 5\times 4\times (3+21).$
\item [] $2960=9\times 8+76\times (5+4\times 3+21).$
\item [] $2961=(9+8+76+5+43)\times 21.$
\item [] $2962=9\times 8\times (7+6\times 5+4)+3^2+1.$
\item [] $2963=9+8\times 7+6\times (5\times 4+3)\times 21.$
\item [] $2964=(9\times 8+7\times 6)\times (5\times 4+3+2+1).$
\item [] $2965=(9+87)\times 6\times 5+4^3+21.$
\item [] $2966=(9\times 8+76)\times 5\times 4+3\times 2\times 1.$
\item [] $2967=(9+87)\times 6\times 5+43\times 2+1.$
\item [] $2968=987+6\times 5\times (4^3+2)+1.$
\item [] $2969=(9\times 8+76)\times 5\times 4+3^2\times 1.$
\item [] $2970=987+654\times 3+21.$
\item [] $2971=9\times (8+7)\times 6+5\times 432+1.$
\item [] $2972=9+87\times (6\times 5+4)+3+2\times 1.$
\item [] $2973=9+8\times 7\times 6\times 5+4\times 321.$
\item [] $2974=9\times 8+7+6+(5+4)\times 321.$
\item [] $2975=(9+8)\times (7\times 6+5+4\times 32\times 1).$
\item [] $2976=9\times 8\times (7+6\times 5+4)+3+21.$
\item [] $2977=(9+8+7+65+4)\times 32+1.$
\item [] $2978=9+8\times (7+6\times 54)+321.$
\item [] $2979=(98+76)\times (5+4\times 3)+21.$
\item [] $2980=(9+87\times 6)\times 5+4+321.$
\item [] $2981=9+8\times 7+6\times 54\times 3^2\times 1.$
\item [] $2982=9+8\times 7+6\times 54\times 3^2+1.$
\item [] $2983=(9\times 8\times 7+6)\times 5+432+1.$
\item [] $2984=(9\times 8+76)\times 5\times 4+3+21.$
\item [] $2985=9+(8+76+5+4)\times 32\times 1.$
\item [] $2986=98+76\times (5+4\times 3+21).$
\item [] $2987=9+(8\times 7+6)\times (5+43)+2\times 1.$
\item [] $2988=9+(8+7)\times 6+(5+4)\times 321.$
\item [] $2989=9\times 8\times (7+6\times 5)+4+321.$
\item [] $2990=9+8\times 7+65\times (43+2)\times 1.$
\item [] $2991=9+87+6+(5+4)\times 321.$
\item [] $2992=9\times 87+(65+4)\times 32+1.$
\item [] $2993=(9\times 8+76)\times 5\times 4+32+1.$
\item [] $2994=9+87+6\times (5\times 4+3)\times 21.$
\item [] $2995=9\times 8+7+6\times 54\times 3^2\times 1.$
\item [] $2996=9\times 8+7+6\times 54\times 3^2+1.$
\item [] $2997=9\times (87+6)+5\times 432\times 1.$
\item [] $2998=9\times (87+6)+5\times 432+1.$
\item [] $2999=9+87\times (6\times 5+4)+32\times 1.$
\item [] $3000=98+7+6+(5+4)\times 321.$
\item[]$\mbox{Increasing order}$
\item [] $3001=(1+2^3\times 45+6+7)\times 8+9.$
\item [] $3002=1\times 2\times (3+4^5+6\times (7+8\times 9)).$
\item [] $3003=123+4\times (5+67+8)\times 9.$
\item [] $3004=12^3+4\times (5\times (6+7\times 8)+9).$
\item [] $3005=123\times 4\times 5+67\times 8+9.$
\item [] $3006=(123\times 4+5)\times 6+7+8+9.$
\item [] $3007=123+4+5\times 6\times (7+89).$
\item [] $3008=1\times 2\times (3+4^5+6\times 78+9).$
\item [] $3009=1^2\times 3\times 4\times 5\times (6\times 7+8)+9.$
\item [] $3010=1\times 23\times (4\times 5\times 6+7)+89.$
\item [] $3011=1+23\times (4\times 5\times 6+7)+89.$
\item [] $3012=12\times (3+4\times 56+7+8+9).$
\item [] $3013=1+2\times (3+4)\times 5\times 6\times 7+8\times 9.$
\item [] $3014=1^2\times 3^4\times (5\times 6+7)+8+9.$
\item [] $3015=1+(2+3)\times 45\times (6+7)+89.$
\item [] $3016=1\times 2+3^4\times (5\times 6+7)+8+9.$
\item [] $3017=12\times (34+5\times 6\times 7)+89.$
\item [] $3018=12\times 34+5\times 6\times (78+9).$
\item [] $3019=(1^2+34)\times (5\times 6+7\times 8)+9.$
\item [] $3020=1\times 2345+(67+8)\times 9.$
\item [] $3021=1+2345+(67+8)\times 9.$
\item [] $3022=(12\times 34+5+6)\times 7+89.$
\item [] $3023=(1+2\times 3^4)\times (5+6+7)+89.$
\item [] $3024=12\times 3\times (4+56+7+8+9).$
\item [] $3025=1+234+5\times (6+7\times 8)\times 9.$
\item [] $3026=1^2\times 34\times (5+67+8+9).$
\item [] $3027=1^2+34\times (5+67+8+9).$
\item [] $3028=1\times 2+34\times (5+67+8+9).$
\item [] $3029=1\times 2\times (3+4)\times 5\times 6\times 7+89.$
\item [] $3030=1+2\times (3+4)\times 5\times 6\times 7+89.$
\item [] $3031=(1+2\times 3)\times (4+5\times (6+78)+9).$
\item [] $3032=1\times 2345+678+9.$
\item [] $3033=1+2345+678+9.$
\item [] $3034=1^{23}+4+5+6\times 7\times 8\times 9.$
\item [] $3035=1^2\times 3+45\times 67+8+9.$
\item [] $3036=1^2+3+45\times 67+8+9.$
\item [] $3037=1\times 2+3+45\times 67+8+9.$
\item [] $3038=1\times 2\times 3+45\times 67+8+9.$
\item [] $3039=1+2+3+4+5+6\times 7\times 8\times 9.$
\item [] $3040=1+2\times 3+4+5+6\times 7\times 8\times 9.$
\item [] $3041=1+2^3+45\times 67+8+9.$
\item [] $3042=1+2^3+4+5+6\times 7\times 8\times 9.$
\item [] $3043=1\times 2+3\times 4+5+6\times 7\times 8\times 9.$
\item [] $3044=1+2+3\times 4+5+6\times 7\times 8\times 9.$
\item [] $3045=1^{23}+4\times 5+6\times 7\times 8\times 9.$
\item [] $3046=1+(2\times 3\times 4+5+6)\times (78+9).$
\item [] $3047=12+3+45\times 67+8+9.$
\item [] $3048=12+3+4+5+6\times 7\times 8\times 9.$
\item [] $3049=1\times 2+3+4\times 5+6\times 7\times 8\times 9.$
\item [] $3050=1+2+3+4\times 5+6\times 7\times 8\times 9.$
\item [] $3051=1+2\times 3+4\times 5+6\times 7\times 8\times 9.$
\item [] $3052=1\times 2^3+4\times 5+6\times 7\times 8\times 9.$
\item [] $3053=1\times 2345+6+78\times 9.$
\item [] $3054=1+2\times 3\times 4+5+6\times 7\times 8\times 9.$
\item [] $3055=1\times 23+45\times 67+8+9.$
\item [] $3056=1+23+45\times 67+8+9.$
\item [] $3057=1+23+4+5+6\times 7\times 8\times 9.$
\item [] $3058=1\times 2\times (3\times 4+5)+6\times 7\times 8\times 9.$
\item [] $3059=12+3+4\times 5+6\times 7\times 8\times 9.$
\item [] $3060=1+2345+6\times 7\times (8+9).$
\item [] $3061=1\times 2^3\times 4+5+6\times 7\times 8\times 9.$
\item [] $3062=1+2^3\times 4+5+6\times 7\times 8\times 9.$
\item [] $3063=1^2\times 34+5+6\times 7\times 8\times 9.$
\item [] $3064=1^2+34+5+6\times 7\times 8\times 9.$
\item [] $3065=1\times 2+34+5+6\times 7\times 8\times 9.$
\item [] $3066=1+2+34+5+6\times 7\times 8\times 9.$
\item [] $3067=1\times 23+4\times 5+6\times 7\times 8\times 9.$
\item [] $3068=12\times 3+45\times 67+8+9.$
\item [] $3069=12\times 3+4+5+6\times 7\times 8\times 9.$
\item [] $3070=1+(2+3)\times 456+789.$
\item[]$\mbox{Decreasing order}$
\item [] $3001=(9+8+7+6)\times 5\times 4\times (3+2)+1.$
\item [] $3002=98\times (7+6)+54\times 32\times 1.$
\item [] $3003=98\times (7+6)+54\times 32+1.$
\item [] $3004=9\times 8+7+65\times (43+2)\times 1.$
\item [] $3005=9\times 8+7+65\times (43+2)+1.$
\item [] $3006=9+(8\times 7+6)\times (5+43)+21.$
\item [] $3007=(9\times 8+7)\times (6\times 5+4)+321.$
\item [] $3008=(9+87)\times 6\times 5+4\times 32\times 1.$
\item [] $3009=(987+6+5+4)\times 3+2+1.$
\item [] $3010=(98\times 7+65)\times 4+3+2+1.$
\item [] $3011=(9+8+7+6+5)\times 43\times 2+1.$
\item [] $3012=9+87+6\times 54\times 3^2\times 1.$
\item [] $3013=9+87+6\times 54\times 3^2+1.$
\item [] $3014=9+8+7+65\times (43+2+1).$
\item [] $3015=98+(76+5)\times 4\times 3^2+1.$
\item [] $3016=9\times ((8+7)\times 6+5\times (4+3)^2)+1.$
\item [] $3017=98+7\times 6\times (5+4^3)+21.$
\item [] $3018=9+87+6\times (54\times 3^2+1).$
\item [] $3019=(9+8)\times 7\times 6+(5+43)^2+1.$
\item [] $3020=(9+8)\times 76+54\times 32\times 1.$
\item [] $3021=98+7+6\times 54\times 3^2\times 1.$
\item [] $3022=98+7+6\times 54\times 3^2+1.$
\item [] $3023=(9\times 8+76)\times 5\times 4+3\times 21.$
\item [] $3024=(9+87+6\times 5)\times 4\times 3\times 2\times 1.$
\item [] $3025=(9+87+6\times 5)\times 4\times 3\times 2+1.$
\item [] $3026=98\times 7+65\times 4\times 3^2\times 1.$
\item [] $3027=98\times 7+65\times 4\times 3^2+1.$
\item [] $3028=(98\times 7+65)\times 4+3+21.$
\item [] $3029=98+7\times 6+(5+4)\times 321.$
\item [] $3030=9+(8+7\times 6)\times 54+321.$
\item [] $3031=98+7+65\times (43+2)+1.$
\item [] $3032=(98\times 7+6\times 54)\times 3+2\times 1.$
\item [] $3033=9+8\times 7\times 6\times 5+4^3\times 21.$
\item [] $3034=9+(8+7+6)\times (5+4+3)^2+1.$
\item [] $3035=(9+8)\times 7+6\times 54\times 3^2\times 1.$
\item [] $3036=(98\times 7+65)\times 4+32\times 1.$
\item [] $3037=9\times 8+76+(5+4)\times 321.$
\item [] $3038=9\times 8\times 7\times 6+5+4+3+2\times 1.$
\item [] $3039=9\times 8\times 7\times 6+5+4+3\times 2\times 1.$
\item [] $3040=9\times 8\times 7\times 6+5+4+3\times 2+1.$
\item [] $3041=9+8+7\times 6\times (5+4+3\times 21).$
\item [] $3042=9\times 8\times 7\times 6+5+4+3^2\times 1.$
\item [] $3043=9\times 8\times 7\times 6+5+4\times 3+2\times 1.$
\item [] $3044=9\times 8\times 7\times 6+5+4\times 3+2+1.$
\item [] $3045=9+876+5\times 432\times 1.$
\item [] $3046=9+876+5\times 432+1.$
\item [] $3047=(98+7)\times (6+5\times 4+3)+2\times 1.$
\item [] $3048=98\times (7+6+5)+4\times 321.$
\item [] $3049=9\times 8\times 7\times 6+5\times 4+3+2\times 1.$
\item [] $3050=9\times 8\times 7\times 6+5\times 4+3\times 2\times 1.$
\item [] $3051=9+87\times 6\times 5+432\times 1.$
\item [] $3052=9+87\times 6\times 5+432+1.$
\item [] $3053=9\times 8\times 7\times 6+5+4\times 3\times 2\times 1.$
\item [] $3054=9\times 8\times 7\times 6+5+4\times 3\times 2+1.$
\item [] $3055=9+8\times 7+65\times (43+2+1).$
\item [] $3056=98\times (7+6)+54\times (32+1).$
\item [] $3057=9\times 8\times 7\times 6+5+4+3+21.$
\item [] $3058=9+8\times 76\times 5+4+3+2\times 1.$
\item [] $3059=9+8\times 76\times 5+4+3\times 2\times 1.$
\item [] $3060=9+8\times 76\times 5+4+3\times 2+1.$
\item [] $3061=9+(8+7+6\times 54)\times 3^2+1.$
\item [] $3062=9\times 8\times 7\times 6+5+4\times 3+21.$
\item [] $3063=9+8\times 76\times 5+4\times 3+2\times 1.$
\item [] $3064=9+8\times 76\times 5+4\times 3+2+1.$
\item [] $3065=9\times 8\times 7\times 6+5+4+32\times 1.$
\item [] $3066=9\times 8\times 7\times 6+5+4+32+1.$
\item [] $3067=(98\times 7+65)\times 4+3\times 21.$
\item [] $3068=9\times 8\times 7\times 6+5\times 4+3+21.$
\item [] $3069=9+8\times 76\times 5+4\times (3+2)\times 1.$
\item [] $3070=9+8\times 76\times 5+4\times (3+2)+1.$
\item[]$\mbox{Increasing order}$
\item [] $3071=12+(3+4)\times 5+6\times 7\times 8\times 9.$
\item [] $3072=1^2\times 3+45+6\times 7\times 8\times 9.$
\item [] $3073=1^2+3+45+6\times 7\times 8\times 9.$
\item [] $3074=1\times 2+3+45+6\times 7\times 8\times 9.$
\item [] $3075=1+2+3+45+6\times 7\times 8\times 9.$
\item [] $3076=1+2\times 3+45+6\times 7\times 8\times 9.$
\item [] $3077=1\times 2^3+45+6\times 7\times 8\times 9.$
\item [] $3078=1+2^3+45+6\times 7\times 8\times 9.$
\item [] $3079=(12+34)\times 5\times (6+7)+89.$
\item [] $3080=12\times 3+4\times 5+6\times 7\times 8\times 9.$
\item [] $3081=123+(4+5\times 6)\times (78+9).$
\item [] $3082=1\times 23\times (4\times 5+6\times 7+8\times 9).$
\item [] $3083=1\times (2\times 34+5)\times 6\times 7+8+9.$
\item [] $3084=12+3+45+6\times 7\times 8\times 9.$
\item [] $3085=12^3+4\times 5\times 67+8+9.$
\item [] $3086=1\times 2+3\times 4\times 5+6\times 7\times 8\times 9.$
\item [] $3087=1+2+3\times 4\times 5+6\times 7\times 8\times 9.$
\item [] $3088=1^{23}+45\times 67+8\times 9.$
\item [] $3089=123\times 4\times 5+6+7\times 89.$
\item [] $3090=1^2\times 3+45\times 67+8\times 9.$
\item [] $3091=1^2+3+45\times 67+8\times 9.$
\item [] $3092=1\times 23+45+6\times 7\times 8\times 9.$
\item [] $3093=1+23+45+6\times 7\times 8\times 9.$
\item [] $3094=1+2\times 3+45\times 67+8\times 9.$
\item [] $3095=1\times 2^3+45\times 67+8\times 9.$
\item [] $3096=12+3\times 4\times 5+6\times 7\times 8\times 9.$
\item [] $3097=12\times 34+5\times 67\times 8+9.$
\item [] $3098=1+2\times 34+5+6\times 7\times 8\times 9.$
\item [] $3099=(1+2+3\times 4)\times 5+6\times 7\times 8\times 9.$
\item [] $3100=1+2\times (3+4\times 5)\times 67+8+9.$
\item [] $3101=1\times 2345+(6+78)\times 9.$
\item [] $3102=12+3+45\times 67+8\times 9.$
\item [] $3103=1^2+3\times 4^5+6+7+8+9.$
\item [] $3104=1^{23}\times 45\times 67+89.$
\item [] $3105=12\times 3+45+6\times 7\times 8\times 9.$
\item [] $3106=1+(2^3+4+5+6)\times (7+8)\times 9.$
\item [] $3107=1^2\times 3+45\times 67+89.$
\item [] $3108=1^2+3+45\times 67+89.$
\item [] $3109=1\times 2+3+45\times 67+89.$
\item [] $3110=1+2+3+45\times 67+89.$
\item [] $3111=1+23+45\times 67+8\times 9.$
\item [] $3112=1\times 2^3+45\times 67+89.$
\item [] $3113=1+2^3+45\times 67+89.$
\item [] $3114=12+3\times 4^5+6+7+8+9.$
\item [] $3115=1+234+5\times 6\times (7+89).$
\item [] $3116=(1+2\times 3)^4+(5+6)\times (7\times 8+9).$
\item [] $3117=(123\times 4+5)\times 6+(7+8)\times 9.$
\item [] $3118=1+23+(4\times 5+6)\times 7\times (8+9).$
\item [] $3119=12+3+45\times 67+89.$
\item [] $3120=1\times 2\times (3+45)+6\times 7\times 8\times 9.$
\item [] $3121=1\times 23\times 4+5+6\times 7\times 8\times 9.$
\item [] $3122=12+3^4+5+6\times 7\times 8\times 9.$
\item [] $3123=12\times 3+45\times 67+8\times 9.$
\item [] $3124=1\times (2+3)\times 4\times 5+6\times 7\times 8\times 9.$
\item [] $3125=(1+2)\times 34\times 5\times 6+7\times 8+9.$
\item [] $3126=1+2+3+4\times 5\times (67+89).$
\item [] $3127=1\times 23+45\times 67+89.$
\item [] $3128=1+23+45\times 67+89.$
\item [] $3129=(1+2)\times (3+4)\times 5+6\times 7\times 8\times 9.$
\item [] $3130=1+(2\times 3+4+5\times 6)\times 78+9.$
\item [] $3131=1^2\times 3\times 4^5+6\times 7+8+9.$
\item [] $3132=1^2+3\times 4^5+6\times 7+8+9.$
\item [] $3133=1\times 2+3\times 4^5+6\times 7+8+9.$
\item [] $3134=1\times 2+3^4\times 5\times 6+78\times 9.$
\item [] $3135=1+2+3^4\times 5\times 6+78\times 9.$
\item [] $3136=12^3+4\times (5\times 67+8+9).$
\item [] $3137=(1+2)^3\times 4+5+6\times 7\times 8\times 9.$
\item [] $3138=(12\times 34+5\times 6)\times 7+8\times 9.$
\item [] $3139=(1+2)\times 34\times 5\times 6+7+8\times 9.$
\item [] $3140=12\times 3+45\times 67+89.$
\item[]$\mbox{Decreasing order}$
\item [] $3071=9\times 8\times 7\times 6+(5\times 4+3)\times 2+1.$
\item [] $3072=9\times 8\times 7\times 6+(5+4)\times 3+21.$
\item [] $3073=9+8\times 76\times 5+4\times 3\times 2\times 1.$
\item [] $3074=9\times 8\times 7\times 6+5+43+2\times 1.$
\item [] $3075=9\times 8\times 7\times 6+5+43+2+1.$
\item [] $3076=9\times 8\times 7\times 6+5\times 4+32\times 1.$
\item [] $3077=9+8\times 76\times 5+4+3+21.$
\item [] $3078=9+8\times (7+6\times 5\times 4)\times 3+21.$
\item [] $3079=9\times 8\times 7\times 6+5+(4+3)^2+1.$
\item [] $3080=(9+8+76\times 5\times 4+3)\times 2\times 1.$
\item [] $3081=(9+8+76\times 5\times 4+3)\times 2+1.$
\item [] $3082=9+8\times 76\times 5+4\times 3+21.$
\item [] $3083=9\times 8\times 7\times 6+54+3+2\times 1.$
\item [] $3084=9\times 8\times 7\times 6+54+3\times 2\times 1.$
\item [] $3085=9+8\times 76\times 5+4+32\times 1.$
\item [] $3086=9+8\times 76\times 5+4+32+1.$
\item [] $3087=9\times 8\times 7\times 6+5\times 4\times 3+2+1.$
\item [] $3088=9\times 8\times 7\times 6+54+3^2+1.$
\item [] $3089=9+8+(7\times 6+54)\times 32\times 1.$
\item [] $3090=9+8+(7\times 6+54)\times 32+1.$
\item [] $3091=98\times 7+65\times (4+32+1).$
\item [] $3092=9\times (8+7\times (6+5))\times 4+32\times 1.$
\item [] $3093=9\times 8\times 7\times 6+5+43+21.$
\item [] $3094=9+8\times 76\times 5+43+2\times 1.$
\item [] $3095=9+8\times 76\times 5+43+2+1.$
\item [] $3096=9\times 8\times 7\times 6+5+4+3\times 21.$
\item [] $3097=9\times 8\times (7+6\times 5)+432+1.$
\item [] $3098=9+8\times 76\times 5+(4+3)^2\times 1.$
\item [] $3099=9+8\times 76\times 5+(4+3)^2+1.$
\item [] $3100=(9+8\times 76)\times 5+4\times 3+2+1.$
\item [] $3101=9+8+765\times 4+3+21.$
\item [] $3102=9\times 8\times 7\times 6+54+3+21.$
\item [] $3103=9\times (8\times 7+6\times 5)\times 4+3\times 2+1.$
\item [] $3104=9+8+765\times 4+3^{(2+1)}.$
\item [] $3105=9\times 8\times 7\times 6+5\times 4\times 3+21.$
\item [] $3106=9+(8+765)\times 4+3+2\times 1.$
\item [] $3107=9\times 8\times 7\times 6+5\times 4+3\times 21.$
\item [] $3108=9+(8+765)\times 4+3\times 2+1.$
\item [] $3109=9+8+765\times 4+32\times 1.$
\item [] $3110=9\times 8\times 7\times 6+54+32\times 1.$
\item [] $3111=9\times 8\times 7\times 6+54+32+1.$
\item [] $3112=9+8+7\times (6+5\times 43)\times 2+1.$
\item [] $3113=9+8\times 76\times 5+43+21.$
\item [] $3114=9\times 8\times 7\times 6+5+4^3+21.$
\item [] $3115=9\times 8\times 7\times 6+5+43\times 2\times 1.$
\item [] $3116=9+8\times 76\times 5+4+3\times 21.$
\item [] $3117=9+(8+7)\times (65+4)\times 3+2+1.$
\item [] $3118=9+8+7\times (6+5+432\times 1).$
\item [] $3119=9+8+(7\times 6+5)\times (4^3+2)\times 1.$
\item [] $3120=(9+8\times 7+65)\times 4\times 3\times 2\times 1.$
\item [] $3121=(98+7\times 6)\times 5\times 4+321.$
\item [] $3122=98+7\times 6\times (5+4+3\times 21).$
\item [] $3123=9\times 87+65\times 4\times 3^2\times 1.$
\item [] $3124=9\times 87+65\times 4\times 3^2+1.$
\item [] $3125=9+(8+765)\times 4+3+21.$
\item [] $3126=9+8\times (76\times 5+4+3)+21.$
\item [] $3127=98+(7\times 6+5)\times 4^3+21.$
\item [] $3128=98\times (7+6)+5+43^2\times 1.$
\item [] $3129=987+6\times (5+4\times 3)\times 21.$
\item [] $3130=(9+8\times 76)\times 5+43+2\times 1.$
\item [] $3131=9+8\times (76+54)\times 3+2\times 1.$
\item [] $3132=987+65\times (4\times 3+21).$
\item [] $3133=9+(8+765)\times 4+32\times 1.$
\item [] $3134=9+8\times 76\times 5+4^3+21.$
\item [] $3135=9+8\times 76\times 5+43\times 2\times 1.$
\item [] $3136=9+8\times 76\times 5+43\times 2+1.$
\item [] $3137=9\times 8+765\times 4+3+2\times 1.$
\item [] $3138=9\times 8+765\times 4+3+2+1.$
\item [] $3139=9\times 8+765\times 4+3\times 2+1.$
\item [] $3140=9+8+765\times 4+3\times 21.$
\item[]$\mbox{Increasing order}$
\item [] $3141=1+2345+6+789.$
\item [] $3142=12^3+4^5+6\times (7\times 8+9).$
\item [] $3143=12+3\times 4^5+6\times 7+8+9.$
\item [] $3144=1\times 2\times 3\times 4\times 5+6\times 7\times 8\times 9.$
\item [] $3145=1+2\times 3\times 4\times 5+6\times 7\times 8\times 9.$
\item [] $3146=1+2+3\times 4^5+6+7\times 8+9.$
\item [] $3147=123\times 4\times 5+678+9.$
\item [] $3148=12^3+4\times 5\times (6+7\times 8+9).$
\item [] $3149=1^{23}\times 4+56\times 7\times 8+9.$
\item [] $3150=1^{23}+4+56\times 7\times 8+9.$
\item [] $3151=1+(2+3+4)\times (5+6\times 7\times 8+9).$
\item [] $3152=1^2\times 3+4+56\times 7\times 8+9.$
\item [] $3153=1^2+3+4+56\times 7\times 8+9.$
\item [] $3154=1\times 2+3+4+56\times 7\times 8+9.$
\item [] $3155=123+45\times 67+8+9.$
\item [] $3156=123+4+5+6\times 7\times 8\times 9.$
\item [] $3157=12^3+4\times 5\times 67+89.$
\item [] $3158=1+2^3+4+56\times 7\times 8+9.$
\item [] $3159=1+2+3\times 4^5+67+8+9.$
\item [] $3160=1+2+3\times 4+56\times 7\times 8+9.$
\item [] $3161=1\times 2+3\times 45+6\times 7\times 8\times 9.$
\item [] $3162=1+2+3\times 45+6\times 7\times 8\times 9.$
\item [] $3163=1+2\times (3\times 4+5)\times (6+78+9).$
\item [] $3164=12+3+4+56\times 7\times 8+9.$
\item [] $3165=1^2\times 3\times 4^5+6+78+9.$
\item [] $3166=1^2+3\times 4^5+6+78+9.$
\item [] $3167=123+4\times 5+6\times 7\times 8\times 9.$
\item [] $3168=123\times 4\times 5+6+78\times 9.$
\item [] $3169=12+3\times 4+56\times 7\times 8+9.$
\item [] $3170=1+2\times 3\times 4+56\times 7\times 8+9.$
\item [] $3171=12+3\times 45+6\times 7\times 8\times 9.$
\item [] $3172=1\times 23+4+56\times 7\times 8+9.$
\item [] $3173=12\times 3\times 4+5+6\times 7\times 8\times 9.$
\item [] $3174=1^2\times 3\times 4^5+6+7+89.$
\item [] $3175=1^2+3\times 4^5+6+7+89.$
\item [] $3176=1\times 2+3\times 4^5+6+7+89.$
\item [] $3177=12+3\times 4^5+6+78+9.$
\item [] $3178=1+2^3\times 4+56\times 7\times 8+9.$
\item [] $3179=1^2\times 34+56\times 7\times 8+9.$
\item [] $3180=1^2+34+56\times 7\times 8+9.$
\item [] $3181=1\times 2+34+56\times 7\times 8+9.$
\item [] $3182=1+2+34+56\times 7\times 8+9.$
\item [] $3183=12+3\times 4^5+6\times (7+8)+9.$
\item [] $3184=1\times 2^3\times 4\times 5+6\times 7\times 8\times 9.$
\item [] $3185=1+2^3\times 4\times 5+6\times 7\times 8\times 9.$
\item [] $3186=12+3\times 4^5+6+7+89.$
\item [] $3187=1^2+3\times 4^5+6\times 7+8\times 9.$
\item [] $3188=1\times 2+3\times 4^5+6\times 7+8\times 9.$
\item [] $3189=1+2+3\times 4^5+6\times 7+8\times 9.$
\item [] $3190=(1\times 2+3)\times (4+5+6+7\times 89).$
\item [] $3191=12+34+56\times 7\times 8+9.$
\item [] $3192=123+45+6\times 7\times 8\times 9.$
\item [] $3193=1+(2+3^4)\times 5\times 6+78\times 9.$
\item [] $3194=1^2\times 34\times 5+6\times 7\times 8\times 9.$
\item [] $3195=1^2+34\times 5+6\times 7\times 8\times 9.$
\item [] $3196=1\times 2+34\times 5+6\times 7\times 8\times 9.$
\item [] $3197=1+2+34\times 5+6\times 7\times 8\times 9.$
\item [] $3198=12+3\times 4^5+6\times 7+8\times 9.$
\item [] $3199=(12+34)\times 56+7\times 89.$
\item [] $3200=1\times 2+3+45\times (6+7\times 8+9).$
\item [] $3201=(12+34+5+6)\times 7\times 8+9.$
\item [] $3202=1+(23+4+5\times 6)\times 7\times 8+9.$
\item [] $3203=1^2\times 3\times 4^5+6\times 7+89.$
\item [] $3204=1^2+3\times 4^5+6\times 7+89.$
\item [] $3205=1\times 2+3\times 4^5+6\times 7+89.$
\item [] $3206=12+34\times 5+6\times 7\times 8\times 9.$
\item [] $3207=(12+3+4\times 5+6)\times 78+9.$
\item [] $3208=(1^2+3+4)\times 56\times 7+8\times 9.$
\item [] $3209=1^{23}\times 456\times 7+8+9.$
\item [] $3210=123+45\times 67+8\times 9.$
\item[]$\mbox{Decreasing order}$
\item [] $3141=9\times 8+765\times 4+3^2\times 1.$
\item [] $3142=9\times 8+765\times 4+3^2+1.$
\item [] $3143=(987+6+54)\times 3+2\times 1.$
\item [] $3144=9\times 8\times 7\times 6+5\times 4\times 3\times 2\times 1.$
\item [] $3145=9\times 8\times 7\times 6+5\times 4\times 3\times 2+1.$
\item [] $3146=(9+8)\times 76+5+43^2\times 1.$
\item [] $3147=(9+8)\times 76+5+43^2+1.$
\item [] $3148=9+8\times (76\times 5+4\times 3)+2+1.$
\item [] $3149=9\times 8\times 7\times 6+5\times (4\times 3\times 2+1).$
\item [] $3150=(9+8+76+54+3)\times 21.$
\item [] $3151=98\times (7+6)+5^4\times 3+2\times 1.$
\item [] $3152=(9+8+765)\times 4+3+21.$
\item [] $3153=987+6+5\times 432\times 1.$
\item [] $3154=987+6+5\times 432+1.$
\item [] $3155=9+8+7+6+5^4\times (3+2\times 1).$
\item [] $3156=9\times 8+765\times 4+3+21.$
\item [] $3157=9\times 8\times 7\times 6+5+4^3\times 2\times 1.$
\item [] $3158=9\times 8\times 7\times 6+5+4^3\times 2+1.$
\item [] $3159=(98+7)\times 6\times 5+4+3+2\times 1.$
\item [] $3160=(98+7)\times 6\times 5+4+3+2+1.$
\item [] $3161=(98+7)\times 6\times 5+4+3\times 2+1.$
\item [] $3162=(987+6+54)\times 3+21.$
\item [] $3163=(98+7)\times 6\times 5+4+3^2\times 1.$
\item [] $3164=9\times 8+765\times 4+32\times 1.$
\item [] $3165=98+765\times 4+3\times 2+1.$
\item [] $3166=9+8\times (76\times 5+4\times 3)+21.$
\item [] $3167=98+765\times 4+3^2\times 1.$
\item [] $3168=98+765\times 4+3^2+1.$
\item [] $3169=(9+8)\times 76+5^4\times 3+2\times 1.$
\item [] $3170=(9+8)\times 76+5^4\times 3+2+1.$
\item [] $3171=(987+65+4)\times 3+2+1.$
\item [] $3172=(9+8\times 76)\times 5+43\times 2+1.$
\item [] $3173=9\times 8+7\times (6+5+432\times 1).$
\item [] $3174=(98+7)\times 6\times 5+4\times 3\times 2\times 1.$
\item [] $3175=(98+7)\times 6\times 5+4\times 3\times 2+1.$
\item [] $3176=9\times 8\times 7\times 6+5+(4+3)\times 21.$
\item [] $3177=9+8\times 76\times 5+4\times 32\times 1.$
\item [] $3178=9+8\times 76\times 5+4^3\times 2+1.$
\item [] $3179=(9+8+7\times 6\times 5)\times (4+3)\times 2+1.$
\item [] $3180=9\times 8+7\times (6+5+432+1).$
\item [] $3181=9+8\times 76\times 5+4\times (32+1).$
\item [] $3182=98+765\times 4+3+21.$
\item [] $3183=(98+7)\times 6\times 5+4\times 3+21.$
\item [] $3184=9+8+7\times 6+5^4\times (3+2)\times 1.$
\item [] $3185=(98\times 7+6\times 5)\times 4+321.$
\item [] $3186=(98+7)\times 6\times 5+4+32\times 1.$
\item [] $3187=(98+7)\times 6\times 5+4+32+1.$
\item [] $3188=9\times 8\times 7\times 6+54\times 3+2\times 1.$
\item [] $3189=9\times 8\times 7\times 6+54\times 3+2+1.$
\item [] $3190=98+765\times 4+32\times 1.$
\item [] $3191=98+765\times 4+32+1.$
\item [] $3192=98+7\times (6+5\times 43)\times 2\times 1.$
\item [] $3193=(9+8+7\times 6)\times 54+3\times 2+1.$
\item [] $3194=9+8\times 76\times 5+(4\times 3)^2+1.$
\item [] $3195=9\times 8+765\times 4+3\times 21.$
\item [] $3196=987+(65+4)\times 32+1.$
\item [] $3197=98+(7+65)\times 43+2+1.$
\item [] $3198=9+8+7+6\times (5\times 4+3)^2\times 1.$
\item [] $3199=98+7\times (6+5+432)\times 1.$
\item [] $3200=(9\times 87+6+5)\times 4+3+21.$
\item [] $3201=9\times (8+7+65)\times 4+321.$
\item [] $3202=(9+8)\times (7\times 6+5)\times 4+3+2+1.$
\item [] $3203=(9+8)\times (7\times 6+5)\times 4+3\times 2+1.$
\item [] $3204=9\times 8\times 7\times 6+5\times 4\times 3^2\times 1.$
\item [] $3205=(9+87)\times 6\times 5+4+321.$
\item [] $3206=98+7\times (6+5+432+1).$
\item [] $3207=9\times 8\times 7\times 6+54\times 3+21.$
\item [] $3208=(9\times 87+6+5)\times 4+32\times 1.$
\item [] $3209=9+8+76\times (5+4+32+1).$
\item [] $3210=9+8+7+65\times (4+3)^2+1.$
\item[]$\mbox{Increasing order}$
\item [] $3211=1^2\times 3\times 4^5+67+8\times 9.$
\item [] $3212=1^2+3\times 4^5+67+8\times 9.$
\item [] $3213=1\times 2\times 34+56\times 7\times 8+9.$
\item [] $3214=1+2\times 34+56\times 7\times 8+9.$
\item [] $3215=12+3\times 4^5+6\times 7+89.$
\item [] $3216=1+2\times 3+456\times 7+8+9.$
\item [] $3217=1\times 2^3+456\times 7+8+9.$
\item [] $3218=1+2^3+456\times 7+8+9.$
\item [] $3219=1^2\times 3^4\times 5\times 6+789.$
\item [] $3220=1^2+3^4\times 5\times 6+789.$
\item [] $3221=1\times 2+3^4\times 5\times 6+789.$
\item [] $3222=1+2+3^4\times 5\times 6+789.$
\item [] $3223=12+3\times 4^5+67+8\times 9.$
\item [] $3224=12+3+456\times 7+8+9.$
\item [] $3225=1\times 2\times 3\times 4\times (56+78)+9.$
\item [] $3226=1^2\times 3^4+56\times 7\times 8+9.$
\item [] $3227=123+45\times 67+89.$
\item [] $3228=1^2\times 3\times 4^5+67+89.$
\item [] $3229=12^3+4^5+6\times 78+9.$
\item [] $3230=1\times 2+3\times 4^5+67+89.$
\item [] $3231=12+3^4\times 5\times 6+789.$
\item [] $3232=1\times 23+456\times 7+8+9.$
\item [] $3233=1+23+456\times 7+8+9.$
\item [] $3234=1^2+(3+45)\times 67+8+9.$
\item [] $3235=1\times 2+(3+45)\times 67+8+9.$
\item [] $3236=1+2+(3+45)\times 67+8+9.$
\item [] $3237=1\times 23\times 4+56\times 7\times 8+9.$
\item [] $3238=12+3^4+56\times 7\times 8+9.$
\item [] $3239=1\times 2+3\times (456+7\times 89).$
\item [] $3240=12+3\times 4^5+67+89.$
\item [] $3241=1+2\times 3\times 456+7\times 8\times 9.$
\item [] $3242=12+(3+456)\times 7+8+9.$
\item [] $3243=123+4\times 5\times (67+89).$
\item [] $3244=(1\times 2+3+456)\times 7+8+9.$
\item [] $3245=12\times 3+456\times 7+8+9.$
\item [] $3246=1+(2+3)\times (4\times 5+6+7\times 89).$
\item [] $3247=(1+2)\times 34+56\times 7\times 8+9.$
\item [] $3248=(1+2)^3\times (4+5)\times (6+7)+89.$
\item [] $3249=12+3\times (456+7\times 89).$
\item [] $3250=1+(2+3)\times 45+6\times 7\times 8\times 9.$
\item [] $3251=(1+2+3+456)\times 7+8+9.$
\item [] $3252=1\times 2\times 3\times (4\times 5+6\times (78+9)).$
\item [] $3253=(1+2)^3\times 4+56\times 7\times 8+9.$
\item [] $3254=(12+34)\times 5+6\times 7\times 8\times 9.$
\item [] $3255=123\times 4\times 5+6+789.$
\item [] $3256=1\times (2+3)^4\times 5+6\times 7+89.$
\item [] $3257=(1+2\times 3)\times 456+7\times 8+9.$
\item [] $3258=1\times 2\times 3\times (456+78+9).$
\item [] $3259=(1+2^3+4)\times 5\times (6\times 7+8)+9.$
\item [] $3260=(1+2\times 34)\times (5+6\times 7)+8+9.$
\item [] $3261=1^2\times 3\times 4^5+(6+7+8)\times 9.$
\item [] $3262=12^3+4^5+6+7\times 8\times 9.$
\item [] $3263=1\times 234+5+6\times 7\times 8\times 9.$
\item [] $3264=1\times 23\times 4\times 5\times 6+7\times 8\times 9.$
\item [] $3265=1+23\times 4\times 5\times 6+7\times 8\times 9.$
\item [] $3266=1+(2^3+456)\times 7+8+9.$
\item [] $3267=1^2\times 3+456\times 7+8\times 9.$
\item [] $3268=1^2+3+456\times 7+8\times 9.$
\item [] $3269=1\times 2+3+456\times 7+8\times 9.$
\item [] $3270=1\times 2\times 3+456\times 7+8\times 9.$
\item [] $3271=1+2\times 3+456\times 7+8\times 9.$
\item [] $3272=123+4+56\times 7\times 8+9.$
\item [] $3273=1+2^3+456\times 7+8\times 9.$
\item [] $3274=12^3+4^5+6\times (78+9).$
\item [] $3275=(12+34)\times (56+7+8)+9.$
\item [] $3276=1\times 234\times (5+6)+78\times 9.$
\item [] $3277=123\times (4\times 5+6)+7+8\times 9.$
\item [] $3278=(12+34)\times 56+78\times 9.$
\item [] $3279=12+3+456\times 7+8\times 9.$
\item [] $3280=1+(23\times 4\times 5+6)\times 7+8+9.$
\item[]$\mbox{Decreasing order}$
\item [] $3211=9\times (8+76+5)\times 4+3\times 2+1.$
\item [] $\mathit{3212=9-8\times 7+6\times 543+2-1.}$
\item [] $3213=(9+8\times 76)\times 5+4\times 32\times 1.$
\item [] $3214=(98+7)\times 6\times 5+43+21.$
\item [] $3215=98+(7+65)\times 43+21.$
\item [] $3216=(98+7)\times 6\times 5+4^3+2\times 1.$
\item [] $3217=(98+7)\times 6\times 5+4+3\times 21.$
\item [] $3218=(9+8+7\times 6)\times 54+32\times 1.$
\item [] $3219=(9+8+7\times 6)\times 54+32+1.$
\item [] $3220=(9+87+65)\times 4\times (3+2)\times 1.$
\item [] $3221=98+765\times 4+3\times 21.$
\item [] $3222=9\times (8\times 7+6)\times 5+432\times 1.$
\item [] $3223=9+8+7\times 65\times (4+3)+21.$
\item [] $3224=9\times 8\times 7\times 6+5\times 4\times (3^2+1).$
\item [] $3225=9+8\times (7+6+54)\times 3\times 2\times 1.$
\item [] $3226=(9+8)\times (76+5)+43^2\times 1.$
\item [] $3227=9+87+6+5^4\times (3+2)\times 1.$
\item [] $3228=9\times (8+76+5)\times 4+3+21.$
\item [] $3229=(9\times 8\times 7+6\times 5+4)\times 3\times 2+1.$
\item [] $3230=9+87\times (6\times 5+4+3)+2\times 1.$
\item [] $3231=9+87\times (6\times 5+4+3)+2+1.$
\item [] $3232=(9+8\times 76)\times 5+(4+3)\times 21.$
\item [] $3233=9+(8\times 7+6)\times (5\times 4+32\times 1).$
\item [] $3234=((9+8\times 7)\times 6)\times 5+4\times 321.$
\item [] $3235=(98+7)\times 6\times 5+4^3+21.$
\item [] $3236=(98+7)\times 6\times 5+43\times 2\times 1.$
\item [] $3237=987+6\times (54+321).$
\item [] $3238=9+8\times (7+6)+5^4\times (3+2)\times 1.$
\item [] $3239=98\times (7+6+5\times 4)+3+2\times 1.$
\item [] $3240=(98+7+6\times 5)\times 4\times 3\times 2\times 1.$
\item [] $3241=(98+7+6\times 5)\times 4\times 3\times 2+1.$
\item [] $3242=9\times 8\times 7\times 6+5\times 43+2+1.$
\item [] $3243=(98+76\times 5\times 4+3)\times 2+1.$
\item [] $3244=98\times (7+6+5\times 4)+3^2+1.$
\item [] $3245=(9+8+7\times 6)\times 5\times (4+3\times 2+1).$
\item [] $3246=(98+7)\times 6\times 5+4\times (3+21).$
\item [] $3247=98+(7\times 6+5)\times (4+3\times 21).$
\item [] $3248=98+7\times 6\times 5\times (4\times 3+2+1).$
\item [] $3249=(9+8+7\times 6)\times 54+3\times 21.$
\item [] $3250=9\times 8\times 7\times 6+5\times (43+2)+1.$
\item [] $3251=9+8+7\times (6\times 5+432\times 1).$
\item [] $3252=9+(8+7\times 65)\times (4+3)+2\times 1.$
\item [] $3253=9+(8+7\times 65)\times (4+3)+2+1.$
\item [] $3254=9\times 8+(7+6\times 5)\times 43\times 2\times 1.$
\item [] $3255=987+(65+43)\times 21.$
\item [] $3256=(98+7)\times (6+5\times 4+3+2)+1.$
\item [] $3257=9+8\times (76+5+4+321).$
\item [] $3258=(9\times 87+6\times 5)\times 4+3\times 2\times 1.$
\item [] $3259=9\times 8+7\times 65\times (4+3)+2\times 1.$
\item [] $3260=9\times 8\times 7\times 6+5\times 43+21.$
\item [] $3261=9+(87\times 6+5\times 4)\times 3\times 2\times 1.$
\item [] $3262=9+(87\times 6+5\times 4)\times 3\times 2+1.$
\item [] $3263=9\times 8+(7+6+5^4)\times (3+2)+1.$
\item [] $3264=(9+8+76+5+4)\times 32\times 1.$
\item [] $3265=(9+8+76+5+4)\times 32+1.$
\item [] $3266=98\times 7+6\times 5\times 43\times 2\times 1.$
\item [] $3267=98\times 7+6\times 5\times 43\times 2+1.$
\item [] $3268=(9+8)\times ((7\times 6+5)\times 4+3)+21.$
\item [] $3269=9\times (87+6\times 5+4)\times 3+2\times 1.$
\item [] $3270=9+87+6\times (5\times 4+3)^2\times 1.$
\item [] $3271=(9+8\times 7+6+5)\times 43+2+1.$
\item [] $3272=98\times 7+6\times (5\times 43\times 2+1).$
\item [] $3273=9+8\times (7\times 6+54\times 3)\times 2\times 1.$
\item [] $3274=9+(87+6+5+4)\times 32+1.$
\item [] $3275=(9+8+7+6+5^4)\times (3+2)\times 1.$
\item [] $3276=9\times 8\times 7+(6+5)\times 4\times 3\times 21.$
\item [] $3277=9\times (8+7+6+5)\times (4+3)\times 2+1.$
\item [] $3278=9\times 8+7\times 65\times (4+3)+21.$
\item [] $3279=(98+7)\times 6\times 5+4\times 32+1.$
\item [] $3280=9+8+7+6+(54+3)^2+1.$
\item[]$\mbox{Increasing order}$
\item [] $3281=1^{23}\times 456\times 7+89.$
\item [] $3282=1^{23}+456\times 7+89.$
\item [] $3283=1+2\times 3+4\times (5\times 6+789).$
\item [] $3284=1^2\times 3+456\times 7+89.$
\item [] $3285=1^2+3+456\times 7+89.$
\item [] $3286=1\times 2+3+456\times 7+89.$
\item [] $3287=1\times 2\times 3+456\times 7+89.$
\item [] $3288=1+23+456\times 7+8\times 9.$
\item [] $3289=12\times 3\times 4+56\times 7\times 8+9.$
\item [] $3290=1+2^3+456\times 7+89.$
\item [] $3291=(1+2)^3+456\times 7+8\times 9.$
\item [] $3292=(1^2+3+456)\times 7+8\times 9.$
\item [] $3293=(12+3+4+5+6+7)\times 89.$
\item [] $3294=1\times 2\times 3\times 45+6\times 7\times 8\times 9.$
\item [] $3295=1+2\times 3\times 45+6\times 7\times 8\times 9.$
\item [] $3296=12+3+456\times 7+89.$
\item [] $3297=12^3+4^5+67\times 8+9.$
\item [] $3298=1\times 2\times (34+56+7)\times (8+9).$
\item [] $3299=1\times 23+4\times (5\times 6+789).$
\item [] $3300=12\times 3+456\times 7+8\times 9.$
\item [] $3301=1+(2+3)\times (4+567+89).$
\item [] $3302=1^2\times (3+456)\times 7+89.$
\item [] $3303=(12+3\times 4+5\times 67+8)\times 9.$
\item [] $3304=1\times 23+456\times 7+89.$
\item [] $3305=1+23+456\times 7+89.$
\item [] $3306=(1+2+3+456)\times 7+8\times 9.$
\item [] $3307=1\times 2\times 3^4+56\times 7\times 8+9.$
\item [] $3308=1+2\times 3^4+56\times 7\times 8+9.$
\item [] $3309=12\times (3+4)\times 5\times 6+789.$
\item [] $3310=12^3+4^5+(6+7\times 8)\times 9.$
\item [] $3311=(123+4)\times (5+6+7+8)+9.$
\item [] $3312=(1\times 234+56+78)\times 9.$
\item [] $3313=12^3+4\times 56\times 7+8+9.$
\item [] $3314=(12+3+456)\times 7+8+9.$
\item [] $3315=1+(2+3)^4+5\times 67\times 8+9.$
\item [] $3316=(1\times 2+3+456)\times 7+89.$
\item [] $3317=12\times 3+456\times 7+89.$
\item [] $3318=123+45\times (6+7\times 8+9).$
\item [] $3319=123+4\times (5+6\times 7)\times (8+9).$
\item [] $3320=(12\times 34+56)\times 7+8\times 9.$
\item [] $3321=12\times 3^4+5\times 6\times 78+9.$
\item [] $3322=1\times 2+(3+4)\times (5+6\times 78)+9.$
\item [] $3323=(1+2+3+456)\times 7+89.$
\item [] $3324=(12+3)\times 4\times 5+6\times 7\times 8\times 9.$
\item [] $3325=1+2\times 3\times (4+5+67\times 8+9).$
\item [] $3326=1\times 2+3\times (4^5+67+8+9).$
\item [] $3327=(1+2)^3\times 4\times 5\times 6+78+9.$
\item [] $3328=12^3+4^5+6\times (7+89).$
\item [] $3329=(1+2+345+67)\times 8+9.$
\item [] $3330=1\times 2\times 3\times (45+6+7\times 8\times 9).$
\item [] $3331=1+2\times 3\times (45+6+7\times 8\times 9).$
\item [] $3332=123+456\times 7+8+9.$
\item [] $3333=1^2+34\times (5+6+78+9).$
\item [] $3334=(1\times 23\times 4\times 5+6)\times 7+8\times 9.$
\item [] $3335=1+2+34\times (5+6+78+9).$
\item [] $3336=(1+2)^3\times 4\times 5\times 6+7+89.$
\item [] $3337=(12\times 34+56)\times 7+89.$
\item [] $3338=1+(2^3+456)\times 7+89.$
\item [] $3339=(1\times 23+45\times 6+78)\times 9.$
\item [] $3340=1+(23+45\times 6+78)\times 9.$
\item [] $3341=12\times (3+45\times 6)+7\times 8+9.$
\item [] $3342=12\times 3^4+5\times 6\times (7+8\times 9).$
\item [] $3343=(1^2+3)^4+(5\times 67+8)\times 9.$
\item [] $3344=12+34\times (5+6+78+9).$
\item [] $3345=(1+2\times 3\times 4+56\times 7)\times 8+9.$
\item [] $3346=1\times (2+3)^4\times 5+(6+7)\times (8+9).$
\item [] $3347=12\times (3+4\times 56)+7\times 89.$
\item [] $3348=12\times 3\times (4+5+67+8+9).$
\item [] $3349=1+2\times (34+5)\times 6\times 7+8\times 9.$
\item [] $3350=1\times (2+3)^4+5\times (67\times 8+9).$
\item[]$\mbox{Decreasing order}$
\item [] $3281=(9\times 8+76)\times 5\times 4+321.$
\item [] $3282=9+87+65\times (4+3)^2+1.$
\item [] $3283=98+(7+6)\times 5\times (4+3)^2\times 1.$
\item [] $3284=9+8+7+6\times 543+2\times 1.$
\item [] $3285=9+8+7+6\times 543+2+1.$
\item [] $3286=98+7\times 65\times (4+3)+2+1.$
\item [] $3287=(9\times (8+7)\times 6+5)\times 4+3^{(2+1)}.$
\item [] $3288=(987+654+3)\times 2\times 1.$
\item [] $3289=(987+654+3)\times 2+1.$
\item [] $3290=98+76\times (5+4+32+1).$
\item [] $3291=98+7+65\times (4+3)^2+1.$
\item [] $3292=9+(8\times 7+6+5)\times (4+3)^2\times 1.$
\item [] $3293=9+8+7\times 6\times (54+3+21).$
\item [] $3294=9\times 8\times 7\times 6+54\times (3+2)\times 1.$
\item [] $3295=9+8+7+654\times (3+2)+1.$
\item [] $3296=98\times 7+6\times 5\times (43\times 2+1).$
\item [] $3297=98\times (7+6+5\times 4)+3\times 21.$
\item [] $3298=987+6+(5+43)^2+1.$
\item [] $3299=98+76+5^4\times (3+2)\times 1.$
\item [] $3300=98+76+5^4\times (3+2)+1.$
\item [] $3301=9+8\times 76\times 5+4\times 3\times 21.$
\item [] $3302=9\times 8\times 7+65\times 43+2+1.$
\item [] $3303=9+8+7+6\times 543+21.$
\item [] $3304=98+7\times 65\times (4+3)+21.$
\item [] $3305=9+8\times 7+6\times 54\times (3^2+1).$
\item [] $3306=(9+87\times 6+5\times 4)\times 3\times 2\times 1.$
\item [] $3307=(9+87\times 6+5\times 4)\times 3\times 2+1.$
\item [] $3308=98\times 7+6\times (5+432\times 1).$
\item [] $3309=98\times 7+6\times (5+432)+1.$
\item [] $3310=9+(8+7)\times (6+5)\times 4\times (3+2)+1.$
\item [] $3311=98+(7\times 6+5+4)\times 3\times 21.$
\item [] $3312=(9+87)\times 6\times 5+432\times 1.$
\item [] $3313=(9+87)\times 6\times 5+432+1.$
\item [] $3314=98\times 7+6\times (5+432+1).$
\item [] $3315=(9\times 87+6\times 5)\times 4+3\times 21.$
\item [] $3316=(9+8)\times (7+6)\times (5+4+3\times 2)+1.$
\item [] $3317=98+(7+6\times 5)\times (43\times 2+1).$
\item [] $3318=(9\times (8+7)+6+5+4\times 3)\times 21.$
\item [] $3319=9+(8\times 76+54)\times (3+2\times 1).$
\item [] $3320=9\times 8\times 7+65\times 43+21.$
\item [] $3321=9\times (8+7+6+5\times 4)\times 3^2\times 1.$
\item [] $3322=9\times (8+7+6+5\times 4)\times 3^2+1.$
\item [] $3323=(9\times 87+6\times 54)\times 3+2\times 1.$
\item [] $3324=(9\times 87+6\times 54)\times 3+2+1.$
\item [] $3325=9+8\times 7+6\times 543+2\times 1.$
\item [] $3326=9+8\times 7+6\times 543+2+1.$
\item [] $3327=987+65\times 4\times 3^2\times 1.$
\item [] $3328=987+65\times 4\times 3^2+1.$
\item [] $3329=9+8+(7+65)\times (43+2+1).$
\item [] $3330=9+8+7\times (6+5)\times 43+2\times 1.$
\item [] $3331=9+8+7\times (6+5)\times 43+2+1.$
\item [] $3332=98+7\times (6\times 5+432\times 1).$
\item [] $3333=9\times 87+6\times 5\times (4^3+21).$
\item [] $3334=9\times 8+7+6+(54+3)^2\times 1.$
\item [] $3335=9+8\times 7+654\times (3+2)\times 1.$
\item [] $3336=9+8\times 7+654\times (3+2)+1.$
\item [] $3337=(9+8\times 76)\times 5+4\times 3\times 21.$
\item [] $3338=9+(8+7\times 6+54)\times 32+1.$
\item [] $3339=9\times 8+7+6\times 543+2\times 1.$
\item [] $3340=9\times 8+7+6\times 543+2+1.$
\item [] $3341=9+8\times 7+6\times (543+2+1).$
\item [] $3342=(98+7\times 65+4)\times 3\times 2\times 1.$
\item [] $3343=(98+7\times 65+4)\times 3\times 2+1.$
\item [] $3344=9+8\times 7+6\times 543+21.$
\item [] $3345=98+7+6\times 54\times (3^2+1).$
\item [] $3346=9\times 8\times 7\times 6+5\times 4^3+2\times 1.$
\item [] $3347=9\times 8\times 7\times 6+5\times 4^3+2+1.$
\item [] $3348=9\times 8\times 7\times 6+54\times 3\times 2\times 1.$
\item [] $3349=9\times 8\times 7\times 6+54\times 3\times 2+1.$
\item [] $3350=9\times 8+7+654\times (3+2)+1.$
\item[]$\mbox{Increasing order}$
\item [] $3351=(1\times 23\times 4\times 5+6)\times 7+89.$
\item [] $3352=1+(23\times 4\times 5+6)\times 7+89.$
\item [] $3353=(1^2\times 3^4\times 5+6+7)\times 8+9.$
\item [] $3354=(1+2)\times (3+4^5+67)+8\times 9.$
\item [] $3355=12\times (3+45\times 6)+7+8\times 9.$
\item [] $3356=12^3+4\times (5\times 67+8\times 9).$
\item [] $3357=12\times 3^4+5\times (6\times 78+9).$
\item [] $3358=(1+23\times 4\times 5+6)\times 7+89.$
\item [] $3359=1\times 2\times 3\times 456+7\times 89.$
\item [] $3360=1+2\times 3\times 456+7\times 89.$
\item [] $3361=(1^2+3^4\times 5+6+7)\times 8+9.$
\item [] $3362=1+2\times 34+(5\times 6+7)\times 89.$
\item [] $3363=1\times 234\times (5+6)+789.$
\item [] $3364=1\times 2\times 34\times 5+6\times 7\times 8\times 9.$
\item [] $3365=(12+34)\times 56+789.$
\item [] $3366=1+2\times (34+5)\times 6\times 7+89.$
\item [] $3367=(1+2+3+4)\times 5\times 67+8+9.$
\item [] $3368=12^3+4\times 56\times 7+8\times 9.$
\item [] $3369=1^2\times 345+6\times 7\times 8\times 9.$
\item [] $3370=1^2+345+6\times 7\times 8\times 9.$
\item [] $3371=1\times 2+345+6\times 7\times 8\times 9.$
\item [] $3372=1+2+345+6\times 7\times 8\times 9.$
\item [] $3373=1^2+3+(4+56)\times 7\times 8+9.$
\item [] $3374=1\times 2+3+(4+56)\times 7\times 8+9.$
\item [] $3375=1\times 2\times 3+(4+56)\times 7\times 8+9.$
\item [] $3376=1+2\times 3+4\times 56\times (7+8)+9.$
\item [] $3377=(1+23+456)\times 7+8+9.$
\item [] $3378=1+2^3\times (4+56)\times 7+8+9.$
\item [] $3379=1\times 234+56\times 7\times 8+9.$
\item [] $3380=1+234+56\times 7\times 8+9.$
\item [] $3381=12+345+6\times 7\times 8\times 9.$
\item [] $3382=1+23\times (45+6+7+89).$
\item [] $3383=1\times 23\times 4\times 5\times 6+7\times 89.$
\item [] $3384=1+23\times 4\times 5\times 6+7\times 89.$
\item [] $3385=12^3+4\times 56\times 7+89.$
\item [] $3386=(12+3+456)\times 7+89.$
\item [] $3387=123+456\times 7+8\times 9.$
\item [] $3388=(1+2)\times 3^4+56\times 7\times 8+9.$
\item [] $3389=1\times 2+3+45\times (67+8)+9.$
\item [] $3390=1\times 2\times 3+45\times (67+8)+9.$
\item [] $3391=1+2\times 3+45\times (67+8)+9.$
\item [] $3392=1\times 23+(4+56)\times 7\times 8+9.$
\item [] $3393=1+23+(4+56)\times 7\times 8+9.$
\item [] $3394=12+3+4+5\times (67+8)\times 9.$
\item [] $3395=12+3+4\times (56+789).$
\item [] $3396=123\times 4\times 5+(6+7)\times 8\times 9.$
\item [] $3397=(1+2^3+4+5\times 6)\times (7+8\times 9).$
\item [] $3398=((1+2)^3+456)\times 7+8+9.$
\item [] $3399=1^{234}\times 5\times 678+9.$
\item [] $3400=1^{234}+5\times 678+9.$
\item [] $3401=(12+34)\times (5+67)+89.$
\item [] $3402=1\times 23+4+5\times (67+8)\times 9.$
\item [] $3403=1^{23}\times 4+5\times 678+9.$
\item [] $3404=123+456\times 7+89.$
\item [] $3405=12\times 3+4\times 56\times (7+8)+9.$
\item [] $3406=1^2\times 3+4+5\times 678+9.$
\item [] $3407=1^2+3+4+5\times 678+9.$
\item [] $3408=1\times 2+3+4+5\times 678+9.$
\item [] $3409=1\times 2\times 3+4+5\times 678+9.$
\item [] $3410=1+2\times 3+4+5\times 678+9.$
\item [] $3411=1^2\times 3\times 4+5\times 678+9.$
\item [] $3412=1^2+3\times 4+5\times 678+9.$
\item [] $3413=1\times 2+3\times 4+5\times 678+9.$
\item [] $3414=1+2+3\times 4+5\times 678+9.$
\item [] $3415=12\times 3+4+5\times (67+8)\times 9.$
\item [] $3416=12\times 3+4\times (56+789).$
\item [] $3417=1^2\times 3\times 4^5+6\times 7\times 8+9.$
\item [] $3418=12+3+4+5\times 678+9.$
\item [] $3419=1\times 2+3\times 4^5+6\times 7\times 8+9.$
\item [] $3420=1+2+3\times 4^5+6\times 7\times 8+9.$
\item[]$\mbox{Decreasing order}$
\item [] $3351=9\times (8\times (7+6)+5\times 4)\times 3+2+1.$
\item [] $3352=9+87+6+(54+3)^2+1.$
\item [] $3353=(9\times 8+765)\times 4+3+2\times 1.$
\item [] $3354=9\times 8\times 7\times 6+5+4+321.$
\item [] $3355=(9\times 8+765)\times 4+3\times 2+1.$
\item [] $3356=9+87+6\times 543+2\times 1.$
\item [] $3357=9+87+6\times 543+2+1.$
\item [] $3358=9\times 8+7+6\times 543+21.$
\item [] $3359=9\times 8\times 7\times 6+5\times (4+3\times 21).$
\item [] $3360=98+7+6+(54+3)^2\times 1.$
\item [] $3361=98+7+6+(54+3)^2+1.$
\item [] $3362=9+(8\times 7+6)\times 54+3+2\times 1.$
\item [] $3363=9\times 87+6\times 5\times 43\times 2\times 1.$
\item [] $3364=9\times 87+6\times 5\times 43\times 2+1.$
\item [] $3365=98+7+6\times 543+2\times 1.$
\item [] $3366=98+7+6\times 543+2+1.$
\item [] $3367=9+87+654\times (3+2)+1.$
\item [] $3368=(9+(8+7\times 6)\times 5)\times (4+3^2)+1.$
\item [] $3369=9\times 87+6\times (5\times 43\times 2+1).$
\item [] $3370=9+8\times 7\times (6\times (5+4)+3\times 2)+1.$
\item [] $\mathit{3371=9+8\times 7\times (6+54)+3-2+1.}$
\item [] $3372=(9\times 8+765)\times 4+3+21.$
\item [] $3373=(9\times (87+6)+5)\times 4+3+2\times 1.$
\item [] $3374=9+8\times 76\times 5+4+321.$
\item [] $3375=9+87+6\times 543+21.$
\item [] $3376=98+7+654\times (3+2)+1.$
\item [] $3377=9+8+7\times (6+5+4)\times 32\times 1.$
\item [] $3378=9+8+7\times (6+5+4)\times 32+1.$
\item [] $3379=(9+8)\times 7+6\times 543+2\times 1.$
\item [] $3380=(9+8)\times 7+6\times 543+2+1.$
\item [] $3381=98+7+6\times (543+2+1).$
\item [] $3382=9\times 8+(7+6\times 54)\times (3^2+1).$
\item [] $3383=9\times 87+65\times 4\times (3^2+1).$
\item [] $3384=98+7+6\times 543+21.$
\item [] $3385=9\times 8+7\times (6+5)\times 43+2\times 1.$
\item [] $3386=9\times 8+7\times (6+5)\times 43+2+1.$
\item [] $3387=(9\times 8+7\times 6\times 5)\times 4\times 3+2+1.$
\item [] $3388=(98+7\times 6\times 5)\times (4+3\times 2+1).$
\item [] $3389=9+(8\times 7+6)\times 54+32\times 1.$
\item [] $3390=9+(8\times 7+6)\times 54+32+1.$
\item [] $3391=(9+8+7+654)\times (3+2)+1.$
\item [] $3392=987+65\times (4+32+1).$
\item [] $3393=9\times 87+6\times 5\times (43\times 2+1).$
\item [] $3394=9\times (8\times 7\times 6+5)+4+321.$
\item [] $3395=9+8\times (76\times 5+43)+2\times 1.$
\item [] $3396=9\times (8\times 7+65+4)\times 3+21.$
\item [] $3397=(9\times 87+65)\times 4+3+2\times 1.$
\item [] $3398=9+8+765\times 4+321.$
\item [] $3399=9\times 8\times 7\times 6+54+321.$
\item [] $3400=(9+8)\times (7+65+4\times 32\times 1).$
\item [] $3401=(9\times 87+65)\times 4+3^2\times 1.$
\item [] $3402=9+8\times (76\times 5+4)+321.$
\item [] $3403=98+(7+654)\times (3+2)\times 1.$
\item [] $3404=9\times 8+7\times (6+5)\times 43+21.$
\item [] $3405=(9\times 8+7\times 6\times 5)\times 4\times 3+21.$
\item [] $3406=(9\times 8+76)\times (5\times 4+3)+2\times 1.$
\item [] $3407=(9\times 8\times 7+6+5^4)\times 3+2\times 1.$
\item [] $3408=9+8+7+6\times (543+21).$
\item [] $3409=(9+8+7)\times 65+43^2\times 1.$
\item [] $3410=(9+8\times 76)\times 5+4+321.$
\item [] $3411=98+7\times (6+5)\times 43+2\times 1.$
\item [] $3412=98+7\times (6+5)\times 43+2+1.$
\item [] $3413=9+(8+7+6)\times 54\times 3+2\times 1.$
\item [] $3414=9\times (8+7)+6\times 543+21.$
\item [] $3415=(9+8)\times 7\times (6+5\times 4)+321.$
\item [] $3416=(9\times 87+65)\times 4+3+21.$
\item [] $3417=9\times (8\times 7+6\times 5)\times 4+321.$
\item [] $3418=9+8\times (7\times (6+54)+3\times 2)+1.$
\item [] $3419=9+(8\times 7+6)\times 5\times (4+3\times 2+1).$
\item [] $3420=9\times 8\times 7+6\times 54\times 3^2\times 1.$
\item[]$\mbox{Increasing order}$
\item [] $3421=12^3+4+5\times 6\times 7\times 8+9.$
\item [] $3422=(1+2+3+4)\times 5\times 67+8\times 9.$
\item [] $3423=12+3\times 4+5\times 678+9.$
\item [] $3424=1+2\times 3\times 4+5\times 678+9.$
\item [] $3425=(12+3)\times 4\times 56+7\times 8+9.$
\item [] $3426=1\times 23+4+5\times 678+9.$
\item [] $3427=1+23+4+5\times 678+9.$
\item [] $3428=12^3+4\times 5\times (6+7+8\times 9).$
\item [] $3429=12+3\times 4^5+6\times 7\times 8+9.$
\item [] $3430=1^2+3^4\times 5+6\times 7\times 8\times 9.$
\item [] $3431=1\times 2^3\times 4+5\times 678+9.$
\item [] $3432=1+2^3\times 4+5\times 678+9.$
\item [] $3433=1^2\times 34+5\times 678+9.$
\item [] $3434=1^2+34+5\times 678+9.$
\item [] $3435=1\times 2+34+5\times 678+9.$
\item [] $3436=1+2+34+5\times 678+9.$
\item [] $3437=12\times 34+5+6\times 7\times 8\times 9.$
\item [] $3438=1\times 2\times 3\times 456+78\times 9.$
\item [] $3439=12\times 3+4+5\times 678+9.$
\item [] $3440=1+(2+3+45)\times 67+89.$
\item [] $3441=12+3^4\times 5+6\times 7\times 8\times 9.$
\item [] $3442=1+(2+3\times 4+5\times 6)\times 78+9.$
\item [] $3443=1+(23+456)\times 7+89.$
\item [] $3444=1+2\times 34+5\times (67+8)\times 9.$
\item [] $3445=12+34+5\times 678+9.$
\item [] $3446=1+2\times 3+4+5\times (678+9).$
\item [] $3447=(12+3)\times 4\times 56+78+9.$
\item [] $3448=1+2\times 3+4\times (5+6)\times 78+9.$
\item [] $3449=(1+23+456)\times 7+89.$
\item [] $3450=1+2^3\times (4+56)\times 7+89.$
\item [] $3451=(1^2+3)\times 4+5\times (678+9).$
\item [] $3452=1\times 2\times 34+(5+6\times 7)\times 8\times 9.$
\item [] $3453=1+2\times 34+(5+6\times 7)\times 8\times 9.$
\item [] $3454=12+3+4+5\times (678+9).$
\item [] $3455=1\times (2+3)\times 4+5\times (678+9).$
\item [] $3456=12+3+4\times (5+6)\times 78+9.$
\item [] $3457=1+2\times 3\times 4\times (5+67+8\times 9).$
\item [] $3458=1\times 2+3^4+5\times (67+8)\times 9.$
\item [] $3459=(12+3)\times 4+5\times 678+9.$
\item [] $3460=12^3+4^5+6+78\times 9.$
\item [] $3461=(12\times 3+456)\times 7+8+9.$
\item [] $3462=1\times 23\times 4\times 5\times 6+78\times 9.$
\item [] $3463=1+23\times 4\times 5\times 6+78\times 9.$
\item [] $3464=1\times 23+4\times (5+6)\times 78+9.$
\item [] $3465=1+23+4\times (5+6)\times 78+9.$
\item [] $3466=12^3+4^5+6\times 7\times (8+9).$
\item [] $3467=1\times 2\times 34+5\times 678+9.$
\item [] $3468=1+2\times 34+5\times 678+9.$
\item [] $3469=1^2\times 34+5\times (678+9).$
\item [] $3470=1^2+34+5\times (678+9).$
\item [] $3471=1\times 2+34+5\times (678+9).$
\item [] $3472=1+2+34+5\times (678+9).$
\item [] $3473=1\times 23\times (4\times 5+6\times 7+89).$
\item [] $3474=(12\times 3+45)\times 6\times 7+8\times 9.$
\item [] $3475=12\times 3+4+5\times (678+9).$
\item [] $3476=1\times 23\times 4\times (5\times 6+7)+8\times 9.$
\item [] $3477=12\times 3+4\times (5+6)\times 78+9.$
\item [] $3478=1^2+3\times (4\times 5+67\times (8+9)).$
\item [] $3479=(12+3)\times 4\times 56+7\times (8+9).$
\item [] $3480=1^2\times 3456+7+8+9.$
\item [] $3481=1^2+3456+7+8+9.$
\item [] $3482=1\times 2+3456+7+8+9.$
\item [] $3483=1+2+3456+7+8+9.$
\item [] $3484=1\times 23\times 4\times 5+6\times 7\times 8\times 9.$
\item [] $3485=1+23\times 4\times 5+6\times 7\times 8\times 9.$
\item [] $3486=(123\times 4+5)\times 6+7\times 8\times 9.$
\item [] $\mathit{3487=-1\times 2-3^4+5\times 6\times 7\times (8+9).}$
\item [] $3488=1\times 2^3\times 4\times (5\times 6+7+8\times 9).$
\item [] $3489=(12+34+5)\times 67+8\times 9.$
\item [] $3490=1+(2\times 3+45)\times 67+8\times 9.$
\item[]$\mbox{Decreasing order}$
\item [] $3421=9\times 8\times 7+6\times 54\times 3^2+1.$
\item [] $3422=9+(8+765)\times 4+321.$
\item [] $3423=98+76+(54+3)^2\times 1.$
\item [] $3424=98+76+(54+3)^2+1.$
\item [] $3425=(9\times 87+65)\times 4+32+1.$
\item [] $3426=9+8+7+6\times (5+4)\times 3\times 21.$
\item [] $3427=9\times (8\times 7+6\times 54)+3\times 2+1.$
\item [] $3428=(9+8\times 76)\times 5+(4+3)^{(2+1)}.$
\item [] $3429=(9+8\times 7+6)\times (5+43)+21.$
\item [] $3430=98+7\times (6+5)\times 43+21.$
\item [] $3431=(9\times (87+6)+5)\times 4+3\times 21.$
\item [] $3432=9+8+7+6+54\times 3\times 21.$
\item [] $3433=9\times 8+7\times (6+5+4)\times 32+1.$
\item [] $3434=98\times (7+6)+5\times 432\times 1.$
\item [] $3435=98\times (7+6)+5\times 432+1.$
\item [] $\mathit{3436=9+8+76\times 5\times (4+3+2)-1.}$
\item [] $3437=9+8+76\times 5\times (4+3+2)\times 1.$
\item [] $3438=(9\times 8\times 7+65+4)\times 3\times 2\times 1.$
\item [] $3439=(9\times 8\times 7+65+4)\times 3\times 2+1.$
\item [] $3440=9+8+7\times (6+(5\times 4+3)\times 21).$
\item [] $3441=9+8\times (7+6)\times (5+4+3+21).$
\item [] $3442=98+76\times (5\times 4+3+21).$
\item [] $3443=9+8\times 7+(6+5+4)^3+2+1.$
\item [] $3444=9\times (8\times 7+6\times 54)+3+21.$
\item [] $3445=(987+6+(5+4)^3)\times 2+1.$
\item [] $3446=9+(8+7\times (65+4))\times (3\times 2+1).$
\item [] $3447=9\times (8\times 7+6)+(5+4)\times 321.$
\item [] $3448=9\times 8\times (7\times 6+5)+43+21.$
\item [] $3449=(9+8+765)\times 4+321.$
\item [] $3450=9\times 8\times (7\times 6+5)+4^3+2\times 1.$
\item [] $3451=9+(8+7+65)\times 43+2\times 1.$
\item [] $3452=(9+8)\times 76+5\times 432\times 1.$
\item [] $3453=9\times 8+765\times 4+321.$
\item [] $3454=9\times 8\times 7\times 6+5\times 43\times 2\times 1.$
\item [] $3455=9\times 8\times 7\times 6+5\times 43\times 2+1.$
\item [] $3456=(9+8+7+6\times 5)\times (43+21).$
\item [] $3457=(98+765)\times 4+3+2\times 1.$
\item [] $3458=(98+765)\times 4+3+2+1.$
\item [] $3459=(98+765)\times 4+3\times 2+1.$
\item [] $\mathit{3460=98\times 7+65\times 43-21.}$
\item [] $3461=9+8+7\times 6+54\times 3\times 21.$
\item [] $3462=9\times 8\times 7\times 6+5+432+1.$
\item [] $3463=9\times 8+7+6\times (543+21).$
\item [] $3464=(9+8\times 7\times 6\times 5+43)\times 2\times 1.$
\item [] $3465=9+8\times (7+6+5)\times 4\times 3\times 2\times 1.$
\item [] $3466=9+(87\times 6+54)\times 3\times 2+1.$
\item [] $3467=9+8\times 7+6\times (5+4)\times 3\times 21.$
\item [] $3468=9+8+7\times (6+54\times 3^2+1).$
\item [] $3469=(98\times 7+6)\times 5+4+3+2\times 1.$
\item [] $3470=9+(8+7+65)\times 43+21.$
\item [] $3471=9\times 8\times (7\times 6+5)+43\times 2+1.$
\item [] $3472=(9+8)\times 7\times (6+5\times 4+3)+21.$
\item [] $3473=9+8\times 7+6+54\times 3\times 21.$
\item [] $3474=(98\times 7+6)\times 5+4\times 3+2\times 1.$
\item [] $3475=(98+7)\times 6\times 5+4+321.$
\item [] $3476=(98+765)\times 4+3+21.$
\item [] $3477=(9\times 8\times 7+654)\times 3+2+1.$
\item [] $3478=(9+8+7\times (6+5))\times (4+32+1).$
\item [] $3479=98+765\times 4+321.$
\item [] $3480=9+87+6\times (543+21).$
\item [] $3481=9+8\times 76\times 5+432\times 1.$
\item [] $3482=9+8\times 76\times 5+432+1.$
\item [] $3483=98\times 7+65\times 43+2\times 1.$
\item [] $3484=98\times 7+65\times 43+2+1.$
\item [] $3485=(98+765)\times 4+32+1.$
\item [] $3486=(98+76)\times 5\times 4+3+2+1.$
\item [] $3487=9\times 8+7+6+54\times 3\times 21.$
\item [] $3488=(98\times 7+6)\times 5+4+3+21.$
\item [] $3489=(98+76)\times 5\times 4+3^2\times 1.$
\item [] $3490=9+(8\times (7+65)+4)\times 3\times 2+1.$
\item[]$\mbox{Increasing order}$
\item [] $3491=1\times 23\times 4+5\times 678+9.$
\item [] $3492=12+3456+7+8+9.$
\item [] $3493=1\times 23\times 4\times (5\times 6+7)+89.$
\item [] $3494=1+23\times 4\times (5\times 6+7)+89.$
\item [] $3495=(1+23)\times 4+5\times 678+9.$
\item [] $3496=(12+34)\times (5+6+7\times 8+9).$
\item [] $3497=1\times 2^3\times 4\times (5+(6+7)\times 8)+9.$
\item [] $3498=1\times 2\times 3\times (4+567)+8\times 9.$
\item [] $3499=1+2\times 3\times (4+567)+8\times 9.$
\item [] $3500=(1+2+3+4)\times (5+6\times 7\times 8+9).$
\item [] $3501=(1+2\times 3+45)\times 67+8+9.$
\item [] $3502=123+4+5\times (67+8)\times 9.$
\item [] $3503=123+4\times (56+789).$
\item [] $3504=1+2\times 34+5\times (678+9).$
\item [] $3505=1\times (2+3)^4+5\times 6\times (7+89).$
\item [] $3506=(12+34+5)\times 67+89.$
\item [] $3507=123+45\times (67+8)+9.$
\item [] $3508=12^3+4^5+(6+78)\times 9.$
\item [] $3509=12\times (3^4+5\times 6\times 7)+8+9.$
\item [] $3510=12\times (34+5)\times 6+78\times 9.$
\item [] $3511=1^2+(3+45+6)\times (7\times 8+9).$
\item [] $3512=12\times 3+4\times (5+6)\times (7+8\times 9).$
\item [] $3513=12\times (3+4\times 56)+789.$
\item [] $3514=1^2+3+(4+5)\times 6\times (7\times 8+9).$
\item [] $3515=1\times 2\times 3\times (4+567)+89.$
\item [] $3516=1+2\times 3\times (4+567)+89.$
\item [] $3517=1^2+3^4+5\times (678+9).$
\item [] $3518=1\times (2+3^4)+5\times (678+9).$
\item [] $3519=(1^2\times 34+5+6)\times 78+9.$
\item [] $3520=1^2+(34+5+6)\times 78+9.$
\item [] $3521=123\times 4+5+6\times 7\times 8\times 9.$
\item [] $3522=1^2+3456+7\times 8+9.$
\item [] $3523=1\times 2+3456+7\times 8+9.$
\item [] $3524=1+2+3456+7\times 8+9.$
\item [] $3525=1\times 2\times 3\times 456+789.$
\item [] $3526=123+4+5\times 678+9.$
\item [] $3527=1\times 23\times 4+5\times (678+9).$
\item [] $3528=1+23\times 4+5\times (678+9).$
\item [] $3529=1+2\times 3\times (4+567+8+9).$
\item [] $3530=(12+3^4)\times (5\times 6+7)+89.$
\item [] $3531=12+(34+5+6)\times 78+9.$
\item [] $3532=1+2\times (3^4+5\times 6\times 7\times 8)+9.$
\item [] $3533=12+3456+7\times 8+9.$
\item [] $3534=(1+2)\times 34\times 5+6\times 7\times 8\times 9.$
\item [] $3535=1^2\times 3456+7+8\times 9.$
\item [] $3536=1^2+3456+7+8\times 9.$
\item [] $3537=1\times 2+3456+7+8\times 9.$
\item [] $3538=1+2+3456+7+8\times 9.$
\item [] $3539=1\times 2+(3+4+56)\times 7\times 8+9.$
\item [] $3540=1+2+(3+4+56)\times 7\times 8+9.$
\item [] $3541=1+2\times 3\times (45+67\times 8+9).$
\item [] $3542=1\times 2+3\times 4\times 5\times (6\times 7+8+9).$
\item [] $3543=12\times 3\times 4+5\times 678+9.$
\item [] $3544=1^2+3456+78+9.$
\item [] $3545=1\times 2+3456+78+9.$
\item [] $3546=1+2+3456+78+9.$
\item [] $3547=12+3456+7+8\times 9.$
\item [] $3548=1\times 2+3\times 4^5+6\times (7+8\times 9).$
\item [] $3549=1\times 23\times 4\times 5\times 6+789.$
\item [] $3550=1+23\times 4\times 5\times 6+789.$
\item [] $3551=1\times 2+3\times 4^5+6\times 78+9.$
\item [] $3552=1+2+3\times 4^5+6\times 78+9.$
\item [] $3553=12\times 34+56\times 7\times 8+9.$
\item [] $3554=1\times 2+3456+7+89.$
\item [] $3555=12+3456+78+9.$
\item [] $3556=(1+2\times 3+45)\times 67+8\times 9.$
\item [] $3557=1\times 2+3\times (4+5+6)\times (7+8\times 9).$
\item [] $3558=12+3\times 4^5+6\times (7+8\times 9).$
\item [] $3559=12+3+4+5\times (6+78\times 9).$
\item [] $3560=(12\times 3+4)\times (5+67+8+9).$
\item[]$\mbox{Decreasing order}$
\item [] $3491=(98\times 7+6)\times 5+4+3^{(2+1)}.$
\item [] $3492=9\times 87+(65+4^3)\times 21.$
\item [] $3493=(98\times 7+6)\times 5+4\times 3+21.$
\item [] $3494=9\times 8\times 7+65\times (43+2+1).$
\item [] $3495=9+8+76+54\times 3\times 21.$
\item [] $3496=(98\times 7+6)\times 5+4+32\times 1.$
\item [] $3497=(9\times 87+6+5)\times 4+321.$
\item [] $3498=9+87+6\times (5+4)\times 3\times 21.$
\item [] $3499=9+8+(7\times 6+5+4\times 3)^2+1.$
\item [] $3500=98+7\times 6\times (5+4)\times 3^2\times 1.$
\item [] $3501=9+(8+7)\times 6+54\times 3\times 21.$
\item [] $3502=98\times 7+65\times 43+21.$
\item [] $3503=9+8+(76+5)\times 43+2+1.$
\item [] $3504=9+87+6+54\times 3\times 21.$
\item [] $3505=(98\times 7+6)\times 5+43+2\times 1.$
\item [] $3506=9+8\times 76+(5+4)\times 321.$
\item [] $3507=98+7+6\times (5+4)\times 3\times 21.$
\item [] $3508=(9\times 8+7)\times (6+5)\times 4+32\times 1.$
\item [] $3509=(9\times 8+7)\times (6+5)\times 4+32+1.$
\item [] $3510=9\times 8\times 7\times 6+54\times 3^2\times 1.$
\item [] $3511=9\times 8\times 7\times 6+54\times 3^2+1.$
\item [] $3512=(98+76)\times 5\times 4+32\times 1.$
\item [] $3513=98+7+6+54\times 3\times 21.$
\item [] $3514=9+8+76\times (5\times 4+3)\times 2+1.$
\item [] $3515=(98+765)\times 4+3\times 21.$
\item [] $3516=9\times 8+7\times 6+54\times 3\times 21.$
\item [] $3517=(9+8\times 76)\times 5+432\times 1.$
\item [] $3518=98+76\times 5\times (4+3+2)\times 1.$
\item [] $3519=(98+7)\times 6+(5+4)\times 321.$
\item [] $3520=9\times (8+76\times 5)+4+3+21.$
\item [] $3521=(9+8)\times 7+6\times (5+4)\times 3\times 21.$
\item [] $3522=9+87\times 6\times 5+43\times 21.$
\item [] $3523=9\times 8+7\times (6+54\times 3^2+1).$
\item [] $3524=(98\times 7+6)\times 5+43+21.$
\item [] $3525=9\times (8+76+5)\times 4+321.$
\item [] $3526=(98\times 7+6)\times 5+4^3+2\times 1.$
\item [] $3527=(98\times 7+6)\times 5+4+3\times 21.$
\item [] $3528=(98+7+6+54+3)\times 21.$
\item [] $3529=(9+8+76+5)\times 4\times 3^2+1.$
\item [] $3530=98\times (7+6+5\times 4+3)+2\times 1.$
\item [] $3531=98\times (7+6+5\times 4+3)+2+1.$
\item [] $3532=9+(8+7\times 6+5)\times 4^3+2+1.$
\item [] $3533=9\times (87+6+5)\times 4+3+2\times 1.$
\item [] $3534=9\times (8+7\times 6)\times 5+4\times 321.$
\item [] $3535=9\times (87+6+5)\times 4+3\times 2+1.$
\item [] $3536=(9+8)\times (7+6)\times (5+4+3\times 2+1).$
\item [] $3537=9\times (8+76\times 5)+43+2\times 1.$
\item [] $3538=9+8\times 7\times 6\times 5+43^2\times 1.$
\item [] $3539=9+8\times 7\times 6\times 5+43^2+1.$
\item [] $3540=9+(876+5)\times 4+3\times 2+1.$
\item [] $3541=(9+8+7\times 6)\times (54+3\times 2)+1.$
\item [] $3542=98+7\times 6+54\times 3\times 21.$
\item [] $3543=9+(876+5)\times 4+3^2+1.$
\item [] $3544=9\times 8+7\times ((6+5)\times (43+2)+1).$
\item [] $3545=(98\times 7+6)\times 5+4^3+21.$
\item [] $3546=(98\times 7+6)\times 5+43\times 2\times 1.$
\item [] $3547=(98\times 7+6)\times 5+43\times 2+1.$
\item [] $\mathit{3548=(9+8)\times 7\times 6\times 5-43+21.}$
\item [] $3549=(9+87+6\times 5+43)\times 21.$
\item [] $3550=9\times 8+76+54\times 3\times 21.$
\item [] $3551=9\times 8+7\times (65+432)\times 1.$
\item [] $3552=(9+87+6+5+4)\times 32\times 1.$
\item [] $3553=(9+87+6+5+4)\times 32+1.$
\item [] $3554=98+(76\times 5+4)\times 3^2\times 1.$
\item [] $3555=9\times 87+(6+5)\times 4\times 3\times 21.$
\item [] $3556=9\times (8+76\times 5)+43+21.$
\item [] $3557=9+(876+5)\times 4+3+21.$
\item [] $3558=9+8\times 7\times (6+54+3)+21.$
\item [] $3559=9\times (8+76\times 5)+4+3\times 21.$
\item [] $3560=9\times (87+6+5)\times 4+32\times 1.$
\item[]$\mbox{Increasing order}$
\item [] $3561=12+3\times 4^5+6\times 78+9.$
\item [] $3562=1+2\times 3^4+5\times 678+9.$
\item [] $3563=1\times 2^{(3+4)}+5\times (678+9).$
\item [] $3564=12+3456+7+89.$
\item [] $3565=1+2\times 3\times 4+5\times (6+78\times 9).$
\item [] $3566=1\times 2\times (3+4^5+(6+78)\times 9).$
\item [] $3567=1\times 23+4+5\times (6+78\times 9).$
\item [] $3568=1+23+4+5\times (6+78\times 9).$
\item [] $3569=1+(2^3+45)\times 67+8+9.$
\item [] $3570=(123+45+6\times 7)\times (8+9).$
\item [] $3571=(1+2)^3+4+5\times (6+78\times 9).$
\item [] $3572=1\times 2^3\times 4+5\times (6+78\times 9).$
\item [] $3573=12+3\times 4\times (5\times 6+7)\times 8+9.$
\item [] $3574=1^2\times 34+5\times (6+78\times 9).$
\item [] $3575=1^2\times 3456+7\times (8+9).$
\item [] $3576=1\times 2+34+5\times (6+78\times 9).$
\item [] $3577=1\times 2+3456+7\times (8+9).$
\item [] $3578=1+2+3456+7\times (8+9).$
\item [] $3579=12\times 3\times 4+5\times (678+9).$
\item [] $3580=1+2+3+4+5\times 6\times 7\times (8+9).$
\item [] $3581=1+2\times 3+4+5\times 6\times 7\times (8+9).$
\item [] $3582=1^2\times 3\times 4^5+6+7\times 8\times 9.$
\item [] $3583=1^2+3\times 4^5+6+7\times 8\times 9.$
\item [] $3584=1\times 2+3\times 4^5+6+7\times 8\times 9.$
\item [] $3585=1+2+3\times 4^5+6+7\times 8\times 9.$
\item [] $3586=(1^2+3)\times 4+5\times 6\times 7\times (8+9).$
\item [] $3587=12+3456+7\times (8+9).$
\item [] $3588=1\times 23\times (4+56+7+89).$
\item [] $3589=12+3+4+5\times 6\times 7\times (8+9).$
\item [] $3590=1\times 2+3\times 4\times (5\times 6\times 7+89).$
\item [] $3591=1^2\times 3456+(7+8)\times 9.$
\item [] $3592=1^2+3456+(7+8)\times 9.$
\item [] $3593=1\times 2+3456+(7+8)\times 9.$
\item [] $3594=1+2+3456+(7+8)\times 9.$
\item [] $3595=1+2\times 3\times 4+5\times 6\times 7\times (8+9).$
\item [] $3596=1+2+(34+5\times 6)\times 7\times 8+9.$
\item [] $3597=12\times (34+5)\times 6+789.$
\item [] $3598=1+23+4+5\times 6\times 7\times (8+9).$
\item [] $3599=1\times 2\times 3\times 45\times (6+7)+89.$
\item [] $3600=12\times (3+45)+6\times 7\times 8\times 9.$
\item [] $3601=1+2\times 3\times 4\times 5\times (6+7+8+9).$
\item [] $3602=(1\times 2+3)^4\times 5+6\times 78+9.$
\item [] $3603=12+3456+(7+8)\times 9.$
\item [] $3604=1234+5\times 6\times (7+8\times 9).$
\item [] $3605=(123\times 4+5)\times 6+7\times 89.$
\item [] $3606=1\times 2+34+5\times 6\times 7\times (8+9).$
\item [] $3607=1+2+34+5\times 6\times 7\times (8+9).$
\item [] $3608=1\times 2\times 34+5\times (6+78\times 9).$
\item [] $3609=1\times 234+5\times (67+8)\times 9.$
\item [] $3610=1+234+5\times (67+8)\times 9.$
\item [] $3611=1\times 2+(3+45)\times (67+8)+9.$
\item [] $3612=12\times (34\times 5+6\times 7+89).$
\item [] $3613=1+(23+4\times 5)\times (67+8+9).$
\item [] $3614=(12+34)\times (5+6)\times 7+8\times 9.$
\item [] $3615=1+(2\times 3+4\times 5)\times (67+8\times 9).$
\item [] $3616=12+34+5\times 6\times 7\times (8+9).$
\item [] $3617=1^2\times 3\times 4^5+67\times 8+9.$
\item [] $3618=1^2+3\times 4^5+67\times 8+9.$
\item [] $3619=1\times 2+3\times 4^5+67\times 8+9.$
\item [] $3620=1+2+3\times 4^5+67\times 8+9.$
\item [] $3621=12+(3+45)\times (67+8)+9.$
\item [] $3622=1^2+3^4+5\times (6+78\times 9).$
\item [] $3623=(1\times 2^3+45)\times 67+8\times 9.$
\item [] $3624=1+(2^3+45)\times 67+8\times 9.$
\item [] $3625=((1^2+3^4)\times 5+6\times 7)\times 8+9.$
\item [] $3626=(1+2+34)\times (5+6+78+9).$
\item [] $3627=(123+4\times 56+7\times 8)\times 9.$
\item [] $3628=1+(23+4)\times (56+78)+9.$
\item [] $3629=12+3\times 4^5+67\times 8+9.$
\item [] $3630=1+2+(34+5)\times (6+78+9).$
\item[]$\mbox{Decreasing order}$
\item [] $3561=(9+8+7\times 6)\times 5\times 4\times 3+21.$
\item [] $3562=9+8\times (7\times 6\times 5+4\times 3)\times 2+1.$
\item [] $3563=9+8\times (7+6\times 5)\times 4\times 3+2\times 1.$
\item [] $3564=9+8\times (7+6\times 5)\times 4\times 3+2+1.$
\item [] $3565=9+(876+5)\times 4+32\times 1.$
\item [] $3566=9+(876+5)\times 4+32+1.$
\item [] $3567=987+6\times 5\times 43\times 2\times 1.$
\item [] $3568=987+6\times 5\times 43\times 2+1.$
\item [] $3569=9\times 8\times 7\times 6+543+2\times 1.$
\item [] $3570=9\times 8\times 7\times 6+543+2+1.$
\item [] $3571=(9+8)\times 7\times (6+(5+4+3)\times 2)+1.$
\item [] $3572=(9+87+6)\times 5\times (4+3)+2\times 1.$
\item [] $3573=(9\times 87+6\times 5)\times 4+321.$
\item [] $3574=(9+8\times 7\times 6)\times 5+43^2\times 1.$
\item [] $3575=(9+8\times 7\times 6)\times 5+43^2+1.$
\item [] $3576=98+76+54\times 3\times 21.$
\item [] $3577=9\times (8+76\times 5)+4^3+21.$
\item [] $3578=9\times (8+76\times 5)+43\times 2\times 1.$
\item [] $3579=9\times (8+76\times 5)+43\times 2+1.$
\item [] $3580=9\times 87+65\times 43+2\times 1.$
\item [] $3581=9\times 87+65\times 43+2+1.$
\item [] $3582=(98+7)\times 6\times 5+432\times 1.$
\item [] $3583=(98+7)\times 6\times 5+432+1.$
\item [] $3584=(9+876+5)\times 4+3+21.$
\item [] $3585=(9+8)\times 7\times 6\times 5+4\times 3+2+1.$
\item [] $3586=9+(8+(7+6)\times 5)\times (4+3)^2\times 1.$
\item [] $3587=987+65\times 4\times (3^2+1).$
\item [] $3588=9\times 8\times 7\times 6+543+21.$
\item [] $3589=(98\times 7+6)\times 5+4\times 32+1.$
\item [] $3590=(9+8)\times 7\times 6\times 5+4\times (3+2)\times 1.$
\item [] $3591=(9\times 8+76+5\times 4+3)\times 21.$
\item [] $3592=(9+876+5)\times 4+32\times 1.$
\item [] $3593=(9+876+5)\times 4+32+1.$
\item [] $3594=98+76\times (5\times 4+3)\times 2\times 1.$
\item [] $3595=98+76\times (5\times 4+3)\times 2+1.$
\item [] $3596=9+(876+5)\times 4+3\times 21.$
\item [] $3597=987+6\times 5\times (43\times 2+1).$
\item [] $3598=(9+8)\times 7\times 6\times 5+4+3+21.$
\item [] $3599=9\times 87+65\times 43+21.$
\item [] $3600=(9+8+7+6)\times 5\times 4\times 3\times 2\times 1.$
\item [] $3601=9\times 8\times (7+6+5+4+3)\times 2+1.$
\item [] $3602=98\times 7+6\times 54\times 3^2\times 1.$
\item [] $3603=98\times 7+6\times 54\times 3^2+1.$
\item [] $3604=(98\times 7+6)\times 5+(4\times 3)^2\times 1.$
\item [] $3605=(98\times 7+6)\times 5+(4\times 3)^2+1.$
\item [] $3606=(9+8)\times 7\times 6\times 5+4+32\times 1.$
\item [] $3607=(9+8)\times 7\times 6\times 5+4+32+1.$
\item [] $3608=98+(7+6)\times 54\times (3+2)\times 1.$
\item [] $3609=987+6\times (5+432)\times 1.$
\item [] $3610=987+6\times (5+432)+1.$
\item [] $3611=98\times 7+65\times (43+2)\times 1.$
\item [] $3612=98\times 7+65\times (43+2)+1.$
\item [] $3613=98\times (7+6+5)+43^2\times 1.$
\item [] $3614=98\times (7+6+5)+43^2+1.$
\item [] $3615=987+6\times (5+432+1).$
\item [] $3616=(9+8+76+5\times 4)\times 32\times 1.$
\item [] $3617=9+8\times (76+54+321).$
\item [] $3618=9\times (8+7+6\times 54+3\times 21).$
\item [] $3619=(9+8)\times 7\times 6\times 5+(4+3)^2\times 1.$
\item [] $3620=9\times (8+76\times 5)+4\times 32\times 1.$
\item [] $3621=9\times (8+76\times 5)+4\times 32+1.$
\item [] $3622=((9\times 8+76\times 5)\times 4+3)\times 2\times 1.$
\item [] $3623=(9+876+5)\times 4+3\times 21.$
\item [] $3624=9+87+(6+54\times 3)\times 21.$
\item [] $3625=9+(87+6+5\times 4)\times 32\times 1.$
\item [] $3626=98+7\times (6+54\times 3)\times (2+1).$
\item [] $3627=9\times (8\times 7+6+5\times 4+321).$
\item [] $3628=(9\times 8+7+6\times 54)\times 3^2+1.$
\item [] $3629=9+8+7\times 6\times (54+32\times 1).$
\item [] $3630=9+8+7\times 6\times (54+32)+1.$
\item[]$\mbox{Increasing order}$
\item [] $3631=1\times 2+(3^4+5)\times 6\times 7+8+9.$
\item [] $3632=1\times 23\times 4+5\times (6+78\times 9).$
\item [] $3633=1\times 234+5\times 678+9.$
\item [] $3634=1+234+5\times 678+9.$
\item [] $3635=12^3+45\times 6\times 7+8+9.$
\item [] $3636=1+2\times 3\times (4+5)\times 67+8+9.$
\item [] $3637=123\times 4+56\times 7\times 8+9.$
\item [] $3638=1\times 2\times 34+5\times 6\times 7\times (8+9).$
\item [] $3639=1+2\times 34+5\times 6\times 7\times (8+9).$
\item [] $3640=1\times 2+34\times (5+6+7+89).$
\item [] $3641=1+(2^3+45)\times 67+89.$
\item [] $3642=(1+2)\times 3^4+5\times 678+9.$
\item [] $3643=(1\times 2+(3^4+5)\times 6)\times 7+8+9.$
\item [] $3644=1\times 2\times (3+4^5+6+789).$
\item [] $3645=(1+2+3)^4+5\times 6\times 78+9.$
\item [] $3646=1+(23+4)\times (56+7+8\times 9).$
\item [] $3647=1^2\times 3+4+56\times (7\times 8+9).$
\item [] $3648=1\times 2\times 3\times 4\times (56+7+89).$
\item [] $3649=1\times 2+3+4+56\times (7\times 8+9).$
\item [] $3650=12+34\times (5+6+7+89).$
\item [] $3651=1+2\times 3+4+56\times (7\times 8+9).$
\item [] $3652=12^3+4\times (56\times 7+89).$
\item [] $3653=1+2^3+4+56\times (7\times 8+9).$
\item [] $3654=1\times 2+3\times 4+56\times (7\times 8+9).$
\item [] $3655=1+(2+3)^4+5+6\times 7\times 8\times 9.$
\item [] $3656=12^3+4\times (5+6\times 78+9).$
\item [] $3657=1\times 23\times (45+6\times 7+8\times 9).$
\item [] $3658=1+23\times (4\times 5+67+8\times 9).$
\item [] $3659=(123+4)\times 5+6\times 7\times 8\times 9.$
\item [] $3660=12+3\times 4^5+6\times (7+89).$
\item [] $3661=12\times 3^4+5\times 67\times 8+9.$
\item [] $3662=1\times 23\times 4+5\times 6\times 7\times (8+9).$
\item [] $3663=12+3^4+5\times 6\times 7\times (8+9).$
\item [] $3664=12+3\times 4+56\times (7\times 8+9).$
\item [] $3665=1+2\times 3\times 4+56\times (7\times 8+9).$
\item [] $3666=12+3^4\times (5\times 6+7+8)+9.$
\item [] $3667=1\times 23+4+56\times (7\times 8+9).$
\item [] $3668=1+23+4+56\times (7\times 8+9).$
\item [] $3669=1\times 234+5\times (678+9).$
\item [] $3670=1+234+5\times (678+9).$
\item [] $3671=(123\times 4+5\times 6)\times 7+8+9.$
\item [] $3672=1\times 2^3\times 456+7+8+9.$
\item [] $3673=1+2^3\times 456+7+8+9.$
\item [] $3674=1^2\times 34+56\times (7\times 8+9).$
\item [] $3675=(1\times 2+34+5+6)\times 78+9.$
\item [] $3676=1+(2+34+5+6)\times 78+9.$
\item [] $3677=1+2+34+56\times (7\times 8+9).$
\item [] $3678=(1+2)\times 3^4+5\times (678+9).$
\item [] $3679=1+2\times 3+(4+5+6\times 7)\times 8\times 9.$
\item [] $3680=12\times 3+4+56\times (7\times 8+9).$
\item [] $3681=(12+3^4\times 5+6\times 7)\times 8+9.$
\item [] $3682=12^3+4+5\times 6\times (7\times 8+9).$
\item [] $3683=(1+2)\times 34\times 5\times 6+7\times 89.$
\item [] $3684=12\times 3\times 4+5\times (6+78\times 9).$
\item [] $3685=1^2+(3^4+5)\times 6\times 7+8\times 9.$
\item [] $3686=12+34+56\times (7\times 8+9).$
\item [] $3687=1+2+(3^4+5)\times 6\times 7+8\times 9.$
\item [] $3688=12^3+4^5+(6+7)\times 8\times 9.$
\item [] $3689=(12+34)\times (5+67+8)+9.$
\item [] $3690=12^3+45\times 6\times 7+8\times 9.$
\item [] $3691=1+2\times 3\times (4+5)\times 67+8\times 9.$
\item [] $3692=(1+2\times 3+45)\times (6+7\times 8+9).$
\item [] $3693=12+34\times (5\times 6+78)+9.$
\item [] $3694=12^3+4+(5\times 6\times 7+8)\times 9.$
\item [] $3695=1\times 23+(4+5+6\times 7)\times 8\times 9.$
\item [] $3696=12+(3^4+5)\times 6\times 7+8\times 9.$
\item [] $3697=123+4+5\times 6\times 7\times (8+9).$
\item [] $3698=(1\times 2+(3^4+5)\times 6)\times 7+8\times 9.$
\item [] $3699=(12+3)\times 45+6\times 7\times 8\times 9.$
\item [] $3700=(12+3)\times 4+56\times (7\times 8+9).$
\item[]$\mbox{Decreasing order}$
\item [] $3631=9+8+7+6+(5\times 4\times 3)^2+1.$
\item [] $\mathit{3632=-9+8\times 7\times 65-4+3\times 2-1.}$
\item [] $3633=(9+8\times 7+65+43)\times 21.$
\item [] $3634=(9+8)\times 7\times 6\times 5+43+21.$
\item [] $3635=98\times (7+6\times 5)+4+3+2\times 1.$
\item [] $3636=98\times (7+6\times 5)+4+3+2+1.$
\item [] $3637=(9+8)\times 7\times 6\times 5+4+3\times 21.$
\item [] $3638=(9+8)\times (7\times (6+5\times 4)+32\times 1).$
\item [] $3639=98\times (7+6\times 5)+4+3^2\times 1.$
\item [] $3640=98\times (7+6\times 5)+4+3^2+1.$
\item [] $3641=98\times (7+6\times 5)+4\times 3+2+1.$
\item [] $3642=(9\times (8\times 7+6+5)+4)\times 3\times 2\times 1.$
\item [] $3643=(9+8)\times (7\times 6\times 5+4)+3+2\times 1.$
\item [] $3644=(9+8)\times (7\times 6\times 5+4)+3\times 2\times 1.$
\item [] $3645=9\times (8+7+65+4+321).$
\item [] $3646=98\times (7+6\times 5)+4\times (3+2)\times 1.$
\item [] $3647=98\times (7+6\times 5)+4\times (3+2)+1.$
\item [] $3648=(9+8+7)\times (65+43\times 2+1).$
\item [] $3649=9+8\times 7\times (6+54+3+2\times 1).$
\item [] $3650=98\times (7+6\times 5)+4\times 3\times 2\times 1.$
\item [] $3651=98\times (7+6\times 5)+4\times 3\times 2+1.$
\item [] $3652=(98\times 7+6)\times 5+4^3\times (2+1).$
\item [] $\mathit{3653=9+8\times 7\times 65+4-3+2+1.}$
\item [] $3654=9\times 8\times 7\times 6+5^4+3+2\times 1.$
\item [] $3655=9\times 8\times 7\times 6+5^4+3\times 2\times 1.$
\item [] $3656=9\times 8\times 7\times 6+5^4+3\times 2+1.$
\item [] $3657=(9+8)\times 7\times 6\times 5+43\times 2+1.$
\item [] $3658=9+8\times 7\times 65+4+3+2\times 1.$
\item [] $3659=9+8\times 7\times 65+4+3\times 2\times 1.$
\item [] $3660=9+8\times 7\times 65+4+3\times 2+1.$
\item [] $\mathit{3661=9+8\times 7\times 65+4+3^2-1.}$
\item [] $3662=9+8\times 7\times 65+4+3^2\times 1.$
\item [] $3663=9+8\times 7\times 65+4+3^2+1.$
\item [] $3664=9\times 8\times 7\times 6+5\times 4\times 32\times 1.$
\item [] $3665=9\times 8\times 7\times 6+5\times 4\times 32+1.$
\item [] $3666=9\times (8+7+6\times 5\times 4)\times 3+21.$
\item [] $3667=9+8+76\times (5+43)+2\times 1.$
\item [] $3668=9+8+76\times (5+43)+2+1.$
\item [] $3669=(9\times 8+765)\times 4+321.$
\item [] $3670=9+8\times 7\times 65+4\times (3+2)+1.$
\item [] $3671=98\times (7+6\times 5)+43+2\times 1.$
\item [] $3672=(9\times 8+76+5)\times 4\times 3\times 2\times 1.$
\item [] $3673=9\times 8\times 7\times 6+5^4+3+21.$
\item [] $3674=9+8\times 7\times 65+4\times 3\times 2+1.$
\item [] $3675=(98+7)\times (6+5+4\times 3\times 2\times 1).$
\item [] $3676=98\times 7+65\times (43+2+1).$
\item [] $3677=9+8\times 7\times 65+4+3+21.$
\item [] $3678=9+(8\times 7+6)\times 54+321.$
\item [] $3679=9\times 8\times 7+6\times (5\times 4+3)^2+1.$
\item [] $3680=(9+8+7\times 65\times 4+3)\times 2\times 1.$
\item [] $3681=(9+8+7\times 65\times 4+3)\times 2+1.$
\item [] $3682=9+8\times 7\times 65+4\times 3+21.$
\item [] $3683=9\times (8\times 7+6)+5^4\times (3+2)\times 1.$
\item [] $3684=9\times 8+7\times 6\times (54+32)\times 1.$
\item [] $3685=9+8\times 7\times 65+4+32\times 1.$
\item [] $3686=9+8\times 7\times 65+4+32+1.$
\item [] $3687=9+8\times (7\times 65+4)+3+2+1.$
\item [] $3688=9+8\times (7\times 65+4)+3\times 2+1.$
\item [] $3689=9\times 8\times 7+65\times (4+3)^2\times 1.$
\item [] $3690=98\times (7+6\times 5)+43+21.$
\item [] $3691=9+8\times (7\times 65+4)+3^2+1.$
\item [] $3692=98\times (7+6\times 5)+4^3+2\times 1.$
\item [] $3693=98\times (7+6\times 5)+4+3\times 21.$
\item [] $3694=9+8\times 7\times 65+43+2\times 1.$
\item [] $3695=9+8\times 7\times 65+43+2+1.$
\item [] $3696=987+(65+4^3)\times 21.$
\item [] $3697=(9+8+7)\times (6+5)\times (4+3)\times 2+1.$
\item [] $3698=(98+7\times 6\times 5)\times 4\times 3+2\times 1.$
\item [] $3699=9\times 87+6\times 54\times 3^2\times 1.$
\item [] $3700=9\times 87+6\times 54\times 3^2+1.$
\item[]$\mbox{Increasing order}$
\item [] $3701=1^2\times 3\times 4^5+6+7\times 89.$
\item [] $3702=1^2+3\times 4^5+6+7\times 89.$
\item [] $3703=1\times 2+3\times 4^5+6+7\times 89.$
\item [] $3704=1+2+3\times 4^5+6+7\times 89.$
\item [] $3705=(12+3+45+6)\times 7\times 8+9.$
\item [] $3706=1^2\times 34\times (5\times 6+7+8\times 9).$
\item [] $3707=12^3+45\times 6\times 7+89.$
\item [] $3708=1\times 2\times 34+56\times (7\times 8+9).$
\item [] $3709=1+2\times 34+56\times (7\times 8+9).$
\item [] $3710=(1+2)\times (3+4^5)+6+7\times 89.$
\item [] $3711=12+(3+(4\times (5+6)+7)\times 8)\times 9.$
\item [] $3712=(1^2+3+4)\times (56\times 7+8\times 9).$
\item [] $3713=12+3\times 4^5+6+7\times 89.$
\item [] $3714=1+2^3\times 456+7\times 8+9.$
\item [] $3715=1+2\times 345+6\times 7\times 8\times 9.$
\item [] $3716=1+(2+3)\times (4\times 5\times 6+7\times 89).$
\item [] $3717=(1+2)\times 3\times (4+56\times 7+8+9).$
\item [] $3718=12+34\times (5\times 6+7+8\times 9).$
\item [] $3719=1+2+3+(456+7)\times 8+9.$
\item [] $3720=(12+3)\times (4\times 56+7+8+9).$
\item [] $3721=1\times 2^3\times (456+7)+8+9.$
\item [] $3722=1+2^3\times (456+7)+8+9.$
\item [] $3723=1\times 2+3^4+56\times (7\times 8+9).$
\item [] $3724=1+2+3^4+56\times (7\times 8+9).$
\item [] $3725=(1+2\times 3\times 4)\times (5\times 6+7\times (8+9)).$
\item [] $3726=(123\times 4+5\times 6)\times 7+8\times 9.$
\item [] $3727=1\times 2^3\times 456+7+8\times 9.$
\item [] $3728=1+2^3\times 456+7+8\times 9.$
\item [] $3729=1^2\times 3\times 4\times 5\times (6+7\times 8)+9.$
\item [] $3730=1^2+3\times 4\times 5\times (6+7\times 8)+9.$
\item [] $3731=1\times 2+3\times 4\times 5\times (6+7\times 8)+9.$
\item [] $3732=1\times 23\times 4+56\times (7\times 8+9).$
\item [] $3733=1+23\times 4+56\times (7\times 8+9).$
\item [] $\mathit{3734=1-2\times 3-4+5+6\times 7\times 89.}$
\item [] $3735=1\times 2^3\times 456+78+9.$
\item [] $3736=1+2^3\times 456+78+9.$
\item [] $3737=1\times 2+3\times (456+789).$
\item [] $3738=1+2+3\times (456+789).$
\item [] $3739=1^2+(3+4)\times (5\times 6+7\times 8\times 9).$
\item [] $3740=(1+23+4\times 5)\times (6+7+8\times 9).$
\item [] $3741=12+3\times 4\times 5\times (6+7\times 8)+9.$
\item [] $3742=(1+2)\times 34+56\times (7\times 8+9).$
\item [] $3743=(123\times 4+5\times 6)\times 7+89.$
\item [] $3744=1\times 2^3\times 456+7+89.$
\item [] $3745=1+2^3\times 456+7+89.$
\item [] $3746=1\times 2+(34+5+(6+7))\times 8\times 9.$
\item [] $3747=12+3\times (456+789).$
\item [] $3748=1^{23}+4+5+6\times 7\times 89.$
\item [] $3749=1\times 2+3\times 4^5+(67+8)\times 9.$
\item [] $3750=1^2\times 3+4+5+6\times 7\times 89.$
\item [] $3751=1^2+3+4+5+6\times 7\times 89.$
\item [] $3752=1\times 2+3+4+5+6\times 7\times 89.$
\item [] $3753=1+2+3+4+5+6\times 7\times 89.$
\item [] $3754=1+2\times 3+4+5+6\times 7\times 89.$
\item [] $3755=1\times 2^3+4+5+6\times 7\times 89.$
\item [] $3756=1^2+3\times 4+5+6\times 7\times 89.$
\item [] $3757=1\times 2+3\times 4+5+6\times 7\times 89.$
\item [] $3758=1+2+3\times 4+5+6\times 7\times 89.$
\item [] $3759=1^2\times 3\times 4^5+678+9.$
\item [] $3760=1^2+3\times 4^5+678+9.$
\item [] $3761=1\times 2+3\times 4^5+678+9.$
\item [] $3762=12+3+4+5+6\times 7\times 89.$
\item [] $3763=1\times 2+3+4\times 5+6\times 7\times 89.$
\item [] $3764=1\times 2\times 3+4\times 5+6\times 7\times 89.$
\item [] $3765=1+2\times 3+4\times 5+6\times 7\times 89.$
\item [] $3766=1\times 2^3+4\times 5+6\times 7\times 89.$
\item [] $3767=1\times 2\times 3\times 4+5+6\times 7\times 89.$
\item [] $3768=1+2\times 3\times 4+5+6\times 7\times 89.$
\item [] $3769=(1+2\times 3+456+7)\times 8+9.$
\item [] $3770=1\times 23+4+5+6\times 7\times 89.$
\item[]$\mbox{Decreasing order}$
\item [] $3701=(9+8)\times (7\times 6\times 5+4)+3\times 21.$
\item [] $3702=9+87+6+(5\times 4\times 3)^2\times 1.$
\item [] $3703=9+87+6+(5\times 4\times 3)^2+1.$
\item [] $3704=9\times 8\times (7\times 6+5+4)+32\times 1.$
\item [] $3705=9+8\times 7\times (6+54+3+2+1).$
\item [] $3706=9\times (87+6\times 54)+3\times 2+1.$
\item [] $3707=9+8\times (7\times 65+4+3)+2\times 1.$
\item [] $3708=9\times 87+65\times (43+2)\times 1.$
\item [] $3709=9\times 8\times (7\times 6+5)+4+321.$
\item [] $3710=98+7\times 6\times (54+32\times 1).$
\item [] $3711=98+7\times 6\times (54+32)+1.$
\item [] $3712=9\times 8\times 7\times 6+5^4+3\times 21.$
\item [] $3713=9+8\times 7\times 65+43+21.$
\item [] $3714=9+8\times (7\times 65+4)+32+1.$
\item [] $3715=9+8\times 7\times 65+4^3+2\times 1.$
\item [] $3716=9+8\times 7\times 65+4+3\times 21.$
\item [] $3717=(98+7+65+4+3)\times 21.$
\item [] $3718=(9+8+7\times 6)\times (54+3^2)+1.$
\item [] $3719=9+(87\times 6+5)\times (4+3)+21.$
\item [] $3720=9+8+7\times (6+5+4\times 3)^2\times 1.$
\item [] $3721=9+(8\times 7+6+54)\times 32\times 1.$
\item [] $3722=9\times 8+76\times (5+43)+2\times 1.$
\item [] $3723=9\times (87+6\times 54)+3+21.$
\item [] $3724=98\times (7+6+5\times 4+3+2\times 1).$
\item [] $3725=(9+8\times 7+6+5)\times (4+3)^2+1.$
\item [] $3726=9+8\times (7\times 65+4+3)+21.$
\item [] $3727=((9+8+76)\times 5\times 4+3)\times 2+1.$
\item [] $3728=9\times (8+76+54)\times 3+2\times 1.$
\item [] $3729=9\times (8+76+54)\times 3+2+1.$
\item [] $3730=(98\times 7+6+54)\times (3+2)\times 1.$
\item [] $3731=9\times (87+6\times 54)+32\times 1.$
\item [] $3732=9\times (87+6\times 54)+32+1.$
\item [] $3733=(9\times (87+6\times 5\times 4)+3)\times 2+1.$
\item [] $3734=9+8\times 7\times 65+4^3+21.$
\item [] $3735=9+8\times 7\times 65+43\times 2\times 1.$
\item [] $3736=9+8\times 7\times 65+43\times 2+1.$
\item [] $3737=9+8\times (76\times 5+43\times 2\times 1).$
\item [] $3738=(98+7+6\times 5+43)\times 21.$
\item [] $3739=(9+(8\times 7\times (6+5)+4)\times 3)\times 2+1.$
\item [] $3740=98+7\times 6+(5\times 4\times 3)^2\times 1.$
\item [] $3741=9\times 8+76\times (5+43)+21.$
\item [] $3742=9+8\times 76+5^4\times (3+2)\times 1.$
\item [] $3743=(987+65\times 4)\times 3+2\times 1.$
\item [] $3744=(9\times 8\times 7+6\times 5\times 4)\times 3\times 2\times 1.$
\item [] $3745=(9\times 8\times 7+6\times 5\times 4)\times 3\times 2+1.$
\item [] $3746=98+76\times ((5+4)\times 3+21).$
\item [] $3747=9+8\times 7\times 6+54\times 3\times 21.$
\item [] $3748=98+76\times (5+43)+2\times 1.$
\item [] $3749=98+76\times (5+43)+2+1.$
\item [] $3750=9+(8\times 7+6)\times 5\times 4\times 3+21.$
\item [] $3751=(9+87+654)\times (3+2)+1.$
\item [] $3752=98+7\times 6\times (54+32+1).$
\item [] $3753=9+8\times (7\times 65+4+3^2\times 1).$
\item [] $3754=98\times (7+6\times 5)+4\times 32\times 1.$
\item [] $3755=98\times (7+6\times 5)+4\times 32+1.$
\item [] $3756=(9+8\times 76+5+4)\times 3\times 2\times 1.$
\item [] $3757=(9+8\times 76+5+4)\times 3\times 2+1.$
\item [] $3758=98\times (7+6\times 5)+4\times (32+1).$
\item [] $3759=9\times 8\times 7\times 6+5\times (4+3)\times 21.$
\item [] $3760=9\times 8\times 7+6+(54+3)^2+1.$
\item [] $3761=9+8\times 7\times (6+54+3\times 2+1).$
\item [] $3762=(987+65\times 4)\times 3+21.$
\item [] $3763=(9+8\times 7+6)\times (5\times 4+32+1).$
\item [] $3764=9\times 8\times 7+6\times 543+2\times 1.$
\item [] $3765=9\times 8\times 7+6\times 543+2+1.$
\item [] $3766=9+8\times (7\times 65+4\times 3)+21.$
\item [] $3767=98+76\times (5+43)+21.$
\item [] $3768=(9\times 8+76+5+4)\times (3+21).$
\item [] $3769=9+8\times (7+6\times 5+432+1).$
\item [] $3770=(9+8\times 7)\times (6+5\times 4+32\times 1).$
\item[]$\mbox{Increasing order}$
\item [] $3771=12+3\times 4^5+678+9.$
\item [] $3772=1\times 2\times (3\times 4+5)+6\times 7\times 89.$
\item [] $3773=12+3+4\times 5+6\times 7\times 89.$
\item [] $3774=(1+2\times 3+4)\times 5\times 67+89.$
\item [] $3775=1\times 2^3\times 4+5+6\times 7\times 89.$
\item [] $3776=1+2^3\times 4+5+6\times 7\times 89.$
\item [] $3777=1^2\times 34+5+6\times 7\times 89.$
\item [] $3778=1^2+34+5+6\times 7\times 89.$
\item [] $3779=1\times 2+34+5+6\times 7\times 89.$
\item [] $3780=1^2\times 3\times 4^5+6+78\times 9.$
\item [] $3781=1\times 23+4\times 5+6\times 7\times 89.$
\item [] $3782=1+23+4\times 5+6\times 7\times 89.$
\item [] $3783=12\times 3+4+5+6\times 7\times 89.$
\item [] $3784=1^{23}+45+6\times 7\times 89.$
\item [] $3785=12+(3+4)\times 5+6\times 7\times 89.$
\item [] $3786=1^2\times 3+45+6\times 7\times 89.$
\item [] $3787=1^2+3+45+6\times 7\times 89.$
\item [] $3788=1\times 2+3+45+6\times 7\times 89.$
\item [] $3789=12+34+5+6\times 7\times 89.$
\item [] $3790=1+2\times 3+45+6\times 7\times 89.$
\item [] $3791=1\times 2^3+45+6\times 7\times 89.$
\item [] $3792=1+2^3+45+6\times 7\times 89.$
\item [] $3793=1\times 2^3\times (456+7)+89.$
\item [] $3794=12\times 3+4\times 5+6\times 7\times 89.$
\item [] $3795=1\times 23\times (4+5+67+89).$
\item [] $3796=1+23\times (4+5+67+89).$
\item [] $3797=12\times 3\times (4+5+6)\times 7+8+9.$
\item [] $3798=12+3+45+6\times 7\times 89.$
\item [] $3799=12+3+4+5\times (6+78)\times 9.$
\item [] $3800=1\times 2+3\times 4\times 5+6\times 7\times 89.$
\item [] $3801=1+2+3\times 4\times 5+6\times 7\times 89.$
\item [] $3802=1\times 2\times 3^4+56\times (7\times 8+9).$
\item [] $3803=1\times 23+45\times (67+8+9).$
\item [] $3804=12+3+45\times (6+78)+9.$
\item [] $3805=1+234+5\times 6\times 7\times (8+9).$
\item [] $3806=1\times 23+45+6\times 7\times 89.$
\item [] $3807=12\times 34+5\times 678+9.$
\item [] $3808=1\times 2\times (3+4)\times 5+6\times 7\times 89.$
\item [] $3809=1+2\times (3+4)\times 5+6\times 7\times 89.$
\item [] $3810=12+3\times 4\times 5+6\times 7\times 89.$
\item [] $3811=1\times 2\times 34+5+6\times 7\times 89.$
\item [] $3812=1+2\times 34+5+6\times 7\times 89.$
\item [] $3813=1+23+45\times (6+78)+9.$
\item [] $3814=1+(2+34+5)\times (6+78+9).$
\item [] $3815=(1+2\times 3)\times 456+7\times 89.$
\item [] $3816=1\times 2\times (34+5)+6\times 7\times 89.$
\item [] $3817=1+2+34+5\times (6+78)\times 9.$
\item [] $3818=1^2+34\times (56+7\times 8)+9.$
\item [] $3819=12\times 3+45+6\times 7\times 89.$
\item [] $3820=12\times 3+4+5\times (6+78)\times 9.$
\item [] $3821=123\times (4\times 5+6)+7\times 89.$
\item [] $3822=(1+2)\times (3\times 45+67\times (8+9)).$
\item [] $3823=1\times 2\times (3+(4+5\times 6)\times 7\times 8)+9.$
\item [] $3824=1^2\times 3^4+5+6\times 7\times 89.$
\item [] $3825=(123+4\times 56+78)\times 9.$
\item [] $3826=1\times 2+3^4+5+6\times 7\times 89.$
\item [] $3827=1+2+3^4+5+6\times 7\times 89.$
\item [] $3828=1\times 2\times 3\times (4+5+6+7\times 89).$
\item [] $3829=12+34\times (56+7\times 8)+9.$
\item [] $3830=1\times 2+3+45\times (6+7+8\times 9).$
\item [] $3831=1\times 2\times 3+45\times (6+7+8\times 9).$
\item [] $3832=1\times 2\times 34\times 56+7+8+9.$
\item [] $3833=1+2\times 34\times 56+7+8+9.$
\item [] $3834=1\times 2\times 3^4\times 5+6\times 7\times 8\times 9.$
\item [] $3835=1\times 23\times 4+5+6\times 7\times 89.$
\item [] $3836=12+3^4+5+6\times 7\times 89.$
\item [] $3837=1+2^3+4\times (5+6)\times (78+9).$
\item [] $3838=1\times (2+3)\times 4\times 5+6\times 7\times 89.$
\item [] $3839=(1+23)\times 4+5+6\times 7\times 89.$
\item [] $3840=12+3\times 4^5+(6+78)\times 9.$
\item[]$\mbox{Decreasing order}$
\item [] $3771=9\times (8\times 7\times 6+5\times 4+3\times 21).$
\item [] $3772=(9+8+(7+6)\times 5)\times (43+2+1).$
\item [] $3773=(98+765)\times 4+321.$
\item [] $3774=9+876+(5+4)\times 321.$
\item [] $3775=9\times 8\times 7+654\times (3+2)+1.$
\item [] $3776=9\times 8+7\times (6+5+4\times 3)^2+1.$
\item [] $3777=9+8\times 7\times 65+4^3\times 2\times 1.$
\item [] $3778=9+8\times 7\times 65+4\times 32+1.$
\item [] $3779=9+8+7+6\times 5^4+3+2\times 1.$
\item [] $3780=9+8+7+6\times 5^4+3+2+1.$
\item [] $3781=9+8+7+6\times 5^4+3\times 2+1.$
\item [] $3782=(98+7\times 6)\times (5+4)\times 3+2\times 1.$
\item [] $3783=9\times 8\times 7+6\times 543+21.$
\item [] $3784=987+65\times 43+2\times 1.$
\item [] $3785=987+65\times 43+2+1.$
\item [] $3786=9\times (8+7+6)\times 5\times 4+3\times 2\times 1.$
\item [] $3787=9\times (8+7+6)\times 5\times 4+3\times 2+1.$
\item [] $3788=(9\times 8\times (7+6)+5)\times 4+3+21.$
\item [] $3789=9\times (8+76\times 5+4\times 3+21).$
\item [] $3790=(9\times 8+7\times 65\times 4+3)\times 2\times 1.$
\item [] $3791=(9\times 8+7\times 65\times 4+3)\times 2+1.$
\item [] $3792=(9+8\times 7)\times 6+54\times 3\times 21.$
\item [] $3793=9+8\times 7\times 65+(4\times 3)^2\times 1.$
\item [] $3794=9+8\times 7\times 65+(4\times 3)^2+1.$
\item [] $3795=9\times (8+76)\times 5+4\times 3+2+1.$
\item [] $3796=9+8\times 7\times 65+(4+3)\times 21.$
\item [] $3797=(9\times 8+7)\times (6+5)\times 4+321.$
\item [] $3798=9+8+7+6\times 5^4+3+21.$
\item [] $3799=(9+8\times 7)\times 6\times 5+43^2\times 1.$
\item [] $3800=(9+8\times 7)\times 6\times 5+43^2+1.$
\item [] $3801=(98+76)\times 5\times 4+321.$
\item [] $3802=98+7\times (6+5+4\times 3)^2+1.$
\item [] $3803=987+65\times 43+21.$
\item [] $3804=9\times (8+76)\times 5+4\times 3\times 2\times 1.$
\item [] $3805=9\times (8+76)\times 5+4\times 3\times 2+1.$
\item [] $3806=9+8+7+6\times 5^4+32\times 1.$
\item [] $3807=9+8+7+6\times 5^4+32+1.$
\item [] $3808=9\times (8+76)\times 5+4+3+21.$
\item [] $3809=9+8+7\times 6+5^4\times 3\times 2\times 1.$
\item [] $3810=9+8+7\times 6+5^4\times 3\times 2+1.$
\item [] $3811=9+8+7+(6+5^4)\times 3\times 2+1.$
\item [] $3812=9\times (8+7+6)\times 5\times 4+32\times 1.$
\item [] $3813=9\times (8+76)\times 5+4\times 3+21.$
\item [] $3814=9\times (8\times 7+6)\times 5+4^(3+2)\times 1.$
\item [] $3815=9\times (8\times 7+6)\times 5+4^(3+2)+1.$
\item [] $3816=9\times (8+76)\times 5+4+32\times 1.$
\item [] $3817=9\times 8\times (7\times 6+5)+432+1.$
\item [] $3818=9+8\times (76+54\times 3)\times 2+1.$
\item [] $3819=98+(7+6+5+43)^2\times 1.$
\item [] $3820=9+8\times 7+6\times 5^4+3+2\times 1.$
\item [] $3821=9+8\times 7+6+5^4\times 3\times 2\times 1.$
\item [] $3822=9+8\times 7+6+5^4\times 3\times 2+1.$
\item [] $3823=(9+8\times 76+5\times 4)\times 3\times 2+1.$
\item [] $3824=9+8\times 7+6\times 5^4+3^2\times 1.$
\item [] $3825=9+8\times 7+6\times 5^4+3^2+1.$
\item [] $3826=9\times (8+76)\times 5+43+2+1.$
\item [] $3827=(9\times 8\times (7+6)+5)\times 4+3\times 21.$
\item [] $3828=9+8+7+6\times (5^4+3^2\times 1).$
\item [] $3829=(9+8+7+65)\times 43+2\times 1.$
\item [] $3830=(9+8+7+65)\times 43+2+1.$
\item [] $3831=9\times (8\times 7\times 6+54)+321.$
\item [] $3832=(9+8\times 7\times (6\times 5+4)+3)\times 2\times 1.$
\item [] $3833=(9+8\times 7\times (6\times 5+4)+3)\times 2+1.$
\item [] $3834=9\times 8+7+6\times 5^4+3+2\times 1.$
\item [] $3835=9\times 8+7+6+5^4\times 3\times 2\times 1.$
\item [] $3836=9\times 8+7+6+5^4\times 3\times 2+1.$
\item [] $3837=9+8+7+6\times 5^4+3\times 21.$
\item [] $3838=9\times 8+7+6\times 5^4+3^2\times 1.$
\item [] $3839=9+8\times 7+6\times 5^4+3+21.$
\item [] $3840=(9+8\times 7+6)\times 54+3+2+1.$
\item[]$\mbox{Increasing order}$
\item [] $3841=1+2^3\times (456+7+8+9).$
\item [] $3842=(1+2\times 3^4+56+7)\times (8+9).$
\item [] $3843=12\times 34+5\times (678+9).$
\item [] $3844=1234+5\times 6\times (78+9).$
\item [] $3845=(1+2)\times 34+5+6\times 7\times 89.$
\item [] $3846=(1+2)\times 34\times (5\times 6+7)+8\times 9.$
\item [] $3847=1+(2\times 3)^4+5\times (6+7\times 8\times 9).$
\item [] $3848=1\times 23+45\times (6+7+8\times 9).$
\item [] $3849=1+23+45\times (6+7+8\times 9).$
\item [] $3850=(1+2\times 3+4)\times (5+6\times 7\times 8+9).$
\item [] $3851=(1+2)^3\times 4+5+6\times 7\times 89.$
\item [] $3852=12\times 3\times (4+5+6)\times 7+8\times 9.$
\item [] $3853=(123\times 4+56)\times 7+8+9.$
\item [] $3854=1+(2+3\times 4\times 5)\times (6+7\times 8)+9.$
\item [] $3855=(1+2)\times (34+5)+6\times 7\times 89.$
\item [] $3856=1\times 2\times (34\times 56+7+8+9).$
\item [] $3857=(1\times 2^3\times 4+5)\times (6+7)\times 8+9.$
\item [] $3858=1\times 2\times 3\times 4\times 5+6\times 7\times 89.$
\item [] $3859=1+2\times 3\times 4\times 5+6\times 7\times 89.$
\item [] $3860=1+((2+3)^4+5)\times 6+7+8\times 9.$
\item [] $3861=12\times 3+45\times (6+7+8\times 9).$
\item [] $3862=1+(23+4)\times (56+78+9).$
\item [] $3863=1\times 2+3^4+5\times (6+78)\times 9.$
\item [] $3864=1\times 2\times 3\times 4\times (5+67+89).$
\item [] $3865=1+2\times 3\times 4\times (5+67+89).$
\item [] $3866=1+2\times ((34+5)\times 6+7)\times 8+9.$
\item [] $3867=1^2\times 3\times 4^5+6+789.$
\item [] $3868=1^2+3\times 4^5+6+789.$
\item [] $3869=1\times 2+3\times 4^5+6+789.$
\item [] $3870=123+4+5+6\times 7\times 89.$
\item [] $3871=(12+3+4+5\times 6)\times (7+8\times 9).$
\item [] $3872=1\times 23\times 4+5\times (6+78)\times 9.$
\item [] $3873=1\times 2\times 34\times 56+7\times 8+9.$
\item [] $3874=1+2\times 34\times 56+7\times 8+9.$
\item [] $3875=1\times 2+3\times 45+6\times 7\times 89.$
\item [] $3876=1+2+3\times 45+6\times 7\times 89.$
\item [] $3877=1+2\times (3\times 4+5)\times (6\times 7+8\times 9).$
\item [] $3878=(1+23+4)\times 5+6\times 7\times 89.$
\item [] $3879=12+3\times 4^5+6+789.$
\item [] $3880=1\times 2^3\times (4+56\times 7+89).$
\item [] $3881=123+4\times 5+6\times 7\times 89.$
\item [] $3882=12\times (3+4+5)+6\times 7\times 89.$
\item [] $3883=(1+2)\times 3^4+56\times (7\times 8+9).$
\item [] $3884=1\times 2+3\times 4^5+6\times (7+8)\times 9.$
\item [] $3885=12+3\times 45+6\times 7\times 89.$
\item [] $3886=(1+2\times 3)^4+(5+6)\times (7+8)\times 9.$
\item [] $3887=12\times 3\times 4+5+6\times 7\times 89.$
\item [] $3888=1+2\times 34\times 56+7+8\times 9.$
\item [] $3889=(123+4)\times 5\times 6+7+8\times 9.$
\item [] $3890=1\times 2+3\times (4+5)\times 6\times (7+8+9).$
\item [] $3891=123\times 4+5\times 678+9.$
\item [] $3892=1\times (2^3+4\times 5)\times (67+8\times 9).$
\item [] $3893=1+(2^3+4\times 5)\times (67+8\times 9).$
\item [] $3894=12+3\times 4^5+6\times (7+8)\times 9.$
\item [] $3895=1\times 2\times 34\times 56+78+9.$
\item [] $3896=1+2\times 34\times 56+78+9.$
\item [] $3897=(123+4)\times 5\times 6+78+9.$
\item [] $3898=1\times 2^3\times 4\times 5+6\times 7\times 89.$
\item [] $3899=12\times (3+45\times 6)+7\times 89.$
\item [] $3900=(1\times 2+3+4\times 5)\times (67+89).$
\item [] $3901=12^3+4\times (5+67\times 8)+9.$
\item [] $3902=1\times 2+3\times (4+(5+6+7)\times 8\times 9).$
\item [] $3903=123+45\times (67+8+9).$
\item [] $3904=1\times 2\times 34\times 56+7+89.$
\item [] $3905=1+2\times 34\times 56+7+89.$
\item [] $3906=123+45+6\times 7\times 89.$
\item [] $3907=123+4+5\times (6+78)\times 9.$
\item [] $3908=1^2\times 34\times 5+6\times 7\times 89.$
\item [] $3909=1^2+34\times 5+6\times 7\times 89.$
\item [] $3910=1\times 2+34\times 5+6\times 7\times 89.$
\item[]$\mbox{Decreasing order}$
\item [] $3841=(98+7+6+5+4)\times 32+1.$
\item [] $3842=9+87\times (6+5)\times 4+3+2\times 1.$
\item [] $3843=9+8+76+5^4\times 3\times 2\times 1.$
\item [] $3844=9+8+76+5^4\times 3\times 2+1.$
\item [] $3845=9+8\times 7+(6+54)\times 3\times 21.$
\item [] $3846=9\times (8+76)\times 5+4^3+2\times 1.$
\item [] $3847=9+8\times 7+6\times 5^4+32\times 1.$
\item [] $3848=9+8\times 7+6\times 5^4+32+1.$
\item [] $3849=9+8\times (7\times 6+5+432+1).$
\item [] $3850=9+(8+7)\times 6+5^4\times 3\times 2+1.$
\item [] $3851=9+87+6\times 5^4+3+2\times 1.$
\item [] $3852=9+87+6\times 5^4+3\times 2\times 1.$
\item [] $3853=9\times 8+7+6\times 5^4+3+21.$
\item [] $3854=9+(876+5)\times 4+321.$
\item [] $3855=9+87+6\times 5^4+3^2\times 1.$
\item [] $3856=9+87+6\times 5^4+3^2+1.$
\item [] $3857=9+8\times 7+6\times (5^4+3\times 2+1).$
\item [] $3858=(9+8\times 7+6)\times 54+3+21.$
\item [] $3859=9\times 8+7+(6+54)\times 3\times 21.$
\item [] $3860=98+7+6\times 5^4+3+2\times 1.$
\item [] $3861=9\times 8+7+6\times 5^4+32\times 1.$
\item [] $3862=98+7+6+5^4\times 3\times 2+1.$
\item [] $3863=9+8+7\times (6+543)+2+1.$
\item [] $3864=9+8+7+6\times 5\times 4\times 32\times 1.$
\item [] $3865=9+8+7+6\times 5\times 4\times 32+1.$
\item [] $3866=9\times (8+76)\times 5+43\times 2\times 1.$
\item [] $3867=987+6\times 5\times 4\times (3+21).$
\item [] $3868=9\times 8+7+6\times (5^4+3)+21.$
\item [] $3869=9+87\times (6+5)\times 4+32\times 1.$
\item [] $3870=9+87+6\times 5^4+3+21.$
\item [] $3871=98\times 7+65\times (4+3)^2\times 1.$
\item [] $3872=98\times 7+65\times (4+3)^2+1.$
\item [] $3873=(9\times 8+7+6+5)\times 43+2+1.$
\item [] $3874=(9+8)\times 7+6\times 5^4+3+2\times 1.$
\item [] $3875=(9+8)\times 7+6+5^4\times 3\times 2\times 1.$
\item [] $3876=9+87+(6+54)\times 3\times 21.$
\item [] $3877=9+87+6\times (5^4+3+2)+1.$
\item [] $3878=9+8\times 7+6\times 5^4+3\times 21.$
\item [] $3879=9+87+6\times 5^4+32+1.$
\item [] $3880=9+(8+7+6\times 5)\times 43\times 2+1.$
\item [] $3881=9+(87+6\times 5+4)\times 32\times 1.$
\item [] $3882=987+6+(5+4)\times 321.$
\item [] $3883=9+87+(6+5^4)\times 3\times 2+1.$
\item [] $3884=(98\times (7+6)+5\times 4)\times 3+2\times 1.$
\item [] $3885=987+6\times (5\times 4+3)\times 21.$
\item [] $3886=98+7+6\times (5^4+3+2)+1.$
\item [] $3887=98+7+6\times 5^4+32\times 1.$
\item [] $3888=98+7+6\times 5^4+32+1.$
\item [] $3889=9+(8+76\times 5)\times (4+3+2+1).$
\item [] $3890=98+7\times 6+5^4\times 3\times 2\times 1.$
\item [] $3891=98+7\times 6+5^4\times 3\times 2+1.$
\item [] $3892=9\times 8+7+6\times 5^4+3\times 21.$
\item [] $3893=(98\times 7+6)\times 5+432+1.$
\item [] $3894=98+7+6\times (5^4+3)+21.$
\item [] $3895=(9+8)\times 7\times 6\times 5+4+321.$
\item [] $3896=(9+8)\times 7+6\times 5^4+3^{(2+1)}.$
\item [] $3897=(9+8\times 7+6)\times 54+3\times 21.$
\item [] $3898=9\times 8+76+5^4\times 3\times 2\times 1.$
\item [] $3899=9\times 8+76+5^4\times 3\times 2+1.$
\item [] $3900=9+87\times (6+5)\times 4+3\times 21.$
\item [] $3901=9+8\times 7\times 65+4\times 3\times 21.$
\item [] $3902=(9+8)\times 7+6\times 5^4+32+1.$
\item [] $3903=9+87\times 6\times 5+4\times 321.$
\item [] $3904=987+6\times 54\times 3^2+1.$
\item [] $3905=9+8\times 7+6\times 5\times 4\times 32\times 1.$
\item [] $3906=9+8\times 7+6\times 5\times 4\times 32+1.$
\item [] $3907=(9+8\times (7\times (6+5)+4)\times 3)\times 2+1.$
\item [] $3908=9\times (8+76)\times 5+4\times 32\times 1.$
\item [] $3909=9+87+6\times 5^4+3\times 21.$
\item [] $3910=(9+8)\times (7+65+43)\times 2\times 1.$
\item[]$\mbox{Increasing order}$
\item [] $3911=1+2+34\times 5+6\times 7\times 89.$
\item [] $3912=123+45\times (6+78)+9.$
\item [] $3913=1\times 2\times (34+5\times 6\times 7)\times 8+9.$
\item [] $3914=1+2\times (34+5\times 6\times 7)\times 8+9.$
\item [] $3915=12+3+(4+56)\times (7\times 8+9).$
\item [] $3916=1^2\times (3+4+5\times 6+7)\times 89.$
\item [] $3917=1\times 2+3\times (4+5+6)\times (78+9).$
\item [] $3918=1\times (2+34)\times 5+6\times 7\times 89.$
\item [] $3919=1+(2+34)\times 5+6\times 7\times 89.$
\item [] $3920=12+34\times 5+6\times 7\times 89.$
\item [] $3921=1+(2+3)^4\times 5+6+789.$
\item [] $3922=1\times 2+(3+4)\times (56+7\times 8\times 9).$
\item [] $3923=1234+5\times 67\times 8+9.$
\item [] $3924=12\times 3\times 4+5\times (6+78)\times 9.$
\item [] $3925=(123\times 4+56)\times 7+89.$
\item [] $3926=(1^2+3^4)\times (5+6\times 7)+8\times 9.$
\item [] $3927=123\times 4+5\times (678+9).$
\item [] $3928=1+2\times 34\times 56+7\times (8+9).$
\item [] $3929=(1+2\times 34)\times 56+7\times 8+9.$
\item [] $3930=1\times 2\times (3+45\times 6\times 7+8\times 9).$
\item [] $3931=1+2\times (3+45\times 6\times 7+8\times 9).$
\item [] $3932=12+(3+4)\times (56+7\times 8\times 9).$
\item [] $3933=(1+2+34+56\times 7+8)\times 9.$
\item [] $3934=1+(2+34)\times (5+(6+7)\times 8)+9.$
\item [] $3935=1\times (2+3)^4\times 5+6\times (7+8)\times 9.$
\item [] $3936=12\times 3+(4+56)\times (7\times 8+9).$
\item [] $3937=(1+2+3+4)\times 56\times 7+8+9.$
\item [] $3938=(12\times 3+4)\times 5+6\times 7\times 89.$
\item [] $3939=1+2\times (34\times 56+7\times 8+9).$
\item [] $3940=(1+(2+3)^4)\times 5+6\times (7+8)\times 9.$
\item [] $3941=(1^2\times 3+4)\times (5+(6+7\times 8)\times 9).$
\item [] $3942=12\times (3\times 4+5)+6\times 7\times 89.$
\item [] $3943=1\times 2\times 34\times 56+(7+8)\times 9.$
\item [] $3944=1+2\times 34\times 56+(7+8)\times 9.$
\item [] $3945=1\times 23\times (4+5)+6\times 7\times 89.$
\item [] $3946=1+23\times (4+5)+6\times 7\times 89.$
\item [] $3947=(1\times 23+4\times (5+(6+7)\times 8)\times 9).$
\item [] $3948=123+45\times (6+7+8\times 9).$
\item [] $3949=1+2\times (3+4)\times (5\times 6\times 7+8\times 9).$
\item [] $3950=(1+23+4\times 5+6)\times (7+8\times 9).$
\item [] $3951=(1+2\times 34)\times 56+78+9.$
\item [] $3952=1\times 2^3\times (4\times 5+6\times (7+8\times 9)).$
\item [] $3953=((1+23)\times 4\times 5+6+7)\times 8+9.$
\item [] $3954=(1+23)\times (4+5)+6\times 7\times 89.$
\item [] $3955=(12+3+4\times 5)\times ((6+7)\times 8+9).$
\item [] $3956=1^2+(3+4)\times 5\times ((6+7)\times 8+9).$
\item [] $3957=((1+2)^3+4\times 5)\times (6+78)+9.$
\item [] $3958=1+2+(3+4)\times 5\times ((6+7)\times 8+9).$
\item [] $3959=1234+5\times (67\times 8+9).$
\item [] $3960=1^2\times 3456+7\times 8\times 9.$
\item [] $3961=1^2+3456+7\times 8\times 9.$
\item [] $3962=1\times 2+3456+7\times 8\times 9.$
\item [] $3963=1+2+3456+7\times 8\times 9.$
\item [] $3964=1+(2+3)\times 45+6\times 7\times 89.$
\item [] $3965=1+2\times (3+45\times 6\times 7+89).$
\item [] $3966=1\times 2\times (34\times 56+7+8\times 9).$
\item [] $3967=1+2\times (34\times 56+7+8\times 9).$
\item [] $3968=(12+34)\times 5+6\times 7\times 89.$
\item [] $3969=(12+3)\times (4\times 5+6+7)\times 8+9.$
\item [] $3970=1+2\times 3\times 4\times (5+6)\times (7+8)+9.$
\item [] $\mathit{3971=-1+2+(3+4)\times 567-8+9.}$
\item [] $3972=12+3456+7\times 8\times 9.$
\item [] $3973=(1\times 2+3^4)\times (5+6\times 7)+8\times 9.$
\item [] $3974=1+(2+3^4)\times (5+6\times 7)+8\times 9.$
\item [] $3975=1^{234}\times 5\times (6+789).$
\item [] $3976=1^{234}+5\times (6+789).$
\item [] $3977=1\times 234+5+6\times 7\times 89.$
\item [] $3978=1+234+5+6\times 7\times 89.$
\item [] $3979=1^{23}\times 4+5\times (6+789).$
\item [] $3980=1\times 2+34\times (5\times 6+78+9).$
\item[]$\mbox{Decreasing order}$
\item [] $3911=9+8\times 7+6\times (5\times 4\times 32+1).$
\item [] $3912=9\times 8\times 7+6+54\times 3\times 21.$
\item [] $3913=987+65\times (43+2)+1.$
\item [] $3914=9+(8+7)\times 65\times 4+3+2\times 1.$
\item [] $3915=9+(8+7)\times 65\times 4+3\times 2\times 1.$
\item [] $3916=9+(8+7)\times 65\times 4+3\times 2+1.$
\item [] $3917=9\times 8+7\times (6+543)+2\times 1.$
\item [] $3918=98+7+6\times 5^4+3\times 21.$
\item [] $3919=9\times 8+7+6\times 5\times 4^3\times 2\times 1.$
\item [] $3920=9\times 8+7+6\times 5\times 4\times 32+1.$
\item [] $3921=98+7\times 6\times (5+43\times 2)+1.$
\item [] $3922=9\times (8+7)+(6+5^4)\times 3\times 2+1.$
\item [] $3923=9+8+(7\times 6+5\times 4)\times 3\times 21.$
\item [] $3924=98+76+5^4\times 3\times 2\times 1.$
\item [] $3925=98+76+5^4\times 3\times 2+1.$
\item [] $3926=9+87\times (6+5+4)\times 3+2\times 1.$
\item [] $3927=(9\times 8+7+65+43)\times 21.$
\item [] $3928=9+(8+76+5^4\times 3)\times 2+1.$
\item [] $3929=9\times 8+7\times (6+543+2)\times 1.$
\item [] $3930=(9+8+7+6+5^4)\times 3\times 2\times 1.$
\item [] $3931=(9+8+7+6+5^4)\times 3\times 2+1.$
\item [] $3932=9\times 8\times 7\times 6+5+43\times 21.$
\item [] $3933=9+87\times 6+54\times 3\times 21.$
\item [] $3934=9\times (8\times (7+6)+5)\times 4+3^2+1.$
\item [] $3935=(987+6\times 54)\times 3+2\times 1.$
\item [] $3936=9+87+6\times 5\times 4\times 32\times 1.$
\item [] $3937=9+87+6\times 5\times 4\times 32+1.$
\item [] $3938=9\times 87+(6+5^4)\times (3+2)\times 1.$
\item [] $3939=(9+87\times 6)\times 5+4\times 321.$
\item [] $3940=9\times (8+7+6)+5^4\times 3\times 2+1.$
\item [] $3941=9+(8+7)\times 65\times 4+32\times 1.$
\item [] $3942=9+(8+7)\times 65\times 4+32+1.$
\item [] $3943=98+7\times (6+543)+2\times 1.$
\item [] $3944=98+7\times (6+543)+2+1.$
\item [] $3945=98+7+6\times 5\times 4\times 32\times 1.$
\item [] $3946=98\times 7+6\times 543+2\times 1.$
\item [] $3947=98\times 7+6\times 543+2+1.$
\item [] $3948=9+8+7+654\times 3\times 2\times 1.$
\item [] $3949=9+8+7+654\times 3\times 2+1.$
\item [] $3950=(9\times 8+7)\times (6+5\times 4+3+21).$
\item [] $3951=98\times (7+6\times 5)+4+321.$
\item [] $3952=9+8\times 76\times 5+43\times 21.$
\item [] $3953=(98\times 7+6+5^4)\times 3+2\times 1.$
\item [] $3954=(987+6\times 54)\times 3+21.$
\item [] $3955=98+7\times (6+543+2\times 1).$
\item [] $3956=98\times 7+654\times (3+2)\times 1.$
\item [] $3957=98\times 7+654\times (3+2)+1.$
\item [] $3958=(9\times 87+6)\times 5+4+3^2\times 1.$
\item [] $3959=(9+8)\times 7+6\times 5\times 4\times 32\times 1.$
\item [] $3960=(9+8)\times 7+6\times 5\times 4^3\times 2+1.$
\item [] $3961=9+8\times (7+6)\times (5+4\times 3+21).$
\item [] $3962=98+7\times (6+543)+21.$
\item [] $3963=9+87\times 6\times 5+4^3\times 21.$
\item [] $3964=(9+8)\times 7\times 6+(54+3)^2+1.$
\item [] $3965=98\times 7+6\times 543+21.$
\item [] $3966=9+87+6\times 5\times 43\times (2+1).$
\item [] $3967=9+8+7+6\times (5^4+32)+1.$
\item [] $3968=9\times 87+65\times (4+3)^2\times 1.$
\item [] $3969=(9\times 87+6)\times 5+4\times 3\times 2\times 1.$
\item [] $3970=9\times 8\times 7\times 6+5^4+321.$
\item [] $3971=9\times (87+6+54)\times 3+2\times 1.$
\item [] $3972=9\times (87+6+54)\times 3+2+1.$
\item [] $3973=(9\times 87+6)\times 5+4+3+21.$
\item [] $3974=9+8\times 7\times 65+4+321.$
\item [] $3975=98+7+6\times 5\times 43\times (2+1).$
\item [] $3976=9\times (8+7)+6\times 5\times 4\times 32+1.$
\item [] $3977=987+65\times (43+2+1).$
\item [] $3978=(9+87)\times 6+54\times 3\times 21.$
\item [] $3979=(9\times 8+7\times 6+5^4\times 3)\times 2+1.$
\item [] $\mathit{3980=-98+7+6\times 5^4+321.}$
\item[]$\mbox{Increasing order}$
\item [] $3981=1+2+34\times (5\times 6+78+9).$
\item [] $3982=1\times 2\times (34\times 56+78+9).$
\item [] $3983=1+2\times (34\times 56+78+9).$
\item [] $3984=1\times 2+3+4+5\times (6+789).$
\item [] $3985=1\times 2\times 3+4+5\times (6+789).$
\item [] $3986=1+2\times 3+4+5\times (6+789).$
\item [] $3987=123\times (4\times 5+6)+789.$
\item [] $3988=1\times 2+(3+4)\times 567+8+9.$
\item [] $3989=1+2+(3+4)\times 567+8+9.$
\item [] $3990=12+34\times (5\times 6+78+9).$
\item [] $3991=1^2+3+(45+6)\times 78+9.$
\item [] $3992=1\times 2+3+(45+6)\times 78+9.$
\item [] $3993=1+2+3+(45+6)\times 78+9.$
\item [] $3994=12+3+4+5\times (6+789).$
\item [] $3995=(123+45+67)\times (8+9).$
\item [] $3996=12\times (34+5\times 6\times 7+89).$
\item [] $3997=1+(2^3\times 45+6+78)\times 9.$
\item [] $3998=12+(3+4)\times 567+8+9.$
\item [] $3999=12+3\times 4+5\times (6+789).$
\item [] $4000=1+2\times 3\times 4+5\times (6+789).$
\item [] $4001=12\times 3^4+5+6\times 7\times 8\times 9.$
\item [] $4002=12+3+(45+6)\times 78+9.$
\item [] $4003=1+23+4+5\times (6+789).$
\item [] $4004=(1+23+4)\times (56+78+9).$
\item [] $4005=(1\times 23+4+5+6+7)\times 89.$
\item [] $4006=(1+2)^3+4+5\times (6+789).$
\item [] $4007=1\times 2^3\times 4+5\times (6+789).$
\item [] $4008=1+2^3\times 4+5\times (6+789).$
\item [] $4009=1+2\times 3\times 45+6\times 7\times 89.$
\item [] $4010=1\times 23+(45+6)\times 78+9.$
\item [] $4011=1\times 2+34+5\times (6+789).$
\item [] $4012=1+2+34+5\times (6+789).$
\item [] $4013=12^3+4\times 567+8+9.$
\item [] $4014=12\times (3+4\times 5)+6\times 7\times 89.$
\item [] $4015=12\times 3+4+5\times (6+789).$
\item [] $4016=(1+(2\times 3+4)\times 56)\times 7+89.$
\item [] $4017=(1+2)\times (3+4^5)+(6+7)\times 8\times 9.$
\item [] $4018=(1+2\times 3)\times (4+5\times (6\times 7+8\times 9)).$
\item [] $\mathit{4019=-1\times 2+3\times 4\times 5\times 67-8+9.}$
\item [] $4020=12+3\times 4^5+(6+7)\times 8\times 9.$
\item [] $4021=12+34+5\times (6+789).$
\item [] $4022=1+(2\times 3)^4+5\times (67\times 8+9).$
\item [] $4023=12\times 3+(45+6)\times 78+9.$
\item [] $4024=1234+5\times (6+7\times 8)\times 9.$
\item [] $4025=1+(2+3)^4+5\times 678+9.$
\item [] $4026=1+(23\times 4\times 5+6\times 7)\times 8+9.$
\item [] $\mathit{4027=1-2\times 3\times 4+5\times 6\times (7+8)\times 9.}$
\item [] $\mathit{4028=12-34+5\times 6\times (7+8)\times 9.}$
\item [] $4029=(1+2)^3\times 4\times 5\times 6+789.$
\item [] $4030=(12\times 3+4\times 5+6)\times (7\times 8+9).$
\item [] $4031=(12+3\times 4+5)\times (67+8\times 9).$
\item [] $4032=123\times 4+5\times (6+78\times 9).$
\item [] $4033=1+2^3\times 4\times (5\times 6+7+89).$
\item [] $4034=1\times 2+(3\times 4+5\times 6)\times (7+89).$
\item [] $4035=1+2+3+(45+6)\times (7+8\times 9).$
\item [] $4036=1+2\times 3+(45+6)\times (7+8\times 9).$
\item [] $4037=1^2\times 3\times 4\times 5\times 67+8+9.$
\item [] $4038=(12+3)\times 4\times 5+6\times 7\times 89.$
\item [] $4039=1\times 2+3\times 4\times 5\times 67+8+9.$
\item [] $4040=1+2+3\times 4\times 5\times 67+8+9.$
\item [] $4041=1^2\times (3+4)\times 567+8\times 9.$
\item [] $4042=1^2+(3+4)\times 567+8\times 9.$
\item [] $4043=1\times 2+(3+4)\times 567+8\times 9.$
\item [] $4044=1+2\times 34+5\times (6+789).$
\item [] $4045=1+2\times 3\times (45+6+7\times 89).$
\item [] $4046=1^2\times 34\times (5+6\times 7+8\times 9).$
\item [] $4047=1^2+34\times (5+6\times 7+8\times 9).$
\item [] $4048=12\times 34+56\times (7\times 8+9).$
\item [] $4049=12+3\times 4\times 5\times 67+8+9.$
\item [] $4050=(12+3\times 4\times 5\times 6+78)\times 9.$
\item[]$\mbox{Decreasing order}$
\item [] $3981=(9\times 87+6)\times 5+4+32\times 1.$
\item [] $3982=(9\times 87+6)\times 5+4+32+1.$
\item [] $3983=9+8+(7+654)\times 3\times 2\times 1.$
\item [] $3984=9+8+7+6\times 5\times 4\times (32+1).$
\item [] $3985=(9\times 8\times 7+65)\times (4+3)+2\times 1.$
\item [] $3986=(9\times 8\times 7+65)\times (4+3)+2+1.$
\item [] $3987=9\times (87+6\times 54+32\times 1).$
\item [] $3988=(9+8\times 76)\times 5+43\times 21.$
\item [] $3989=9+8\times 7+654\times 3\times 2\times 1.$
\item [] $3990=9+8\times 7+654\times 3\times 2+1.$
\item [] $3991=(9\times 87+6)\times 5+43+2+1.$
\item [] $3992=(9\times 8+7\times 6)\times 5\times (4+3)+2\times 1.$
\item [] $3993=9\times 8\times (7\times 6+5+4)+321.$
\item [] $3994=(9\times 87+6)\times 5+(4+3)^2\times 1.$
\item [] $3995=9+8\times (7\times 65+43)+2\times 1.$
\item [] $3996=9\times (87+6\times 54+32+1).$
\item [] $3997=(987+6+5)\times 4+3+2\times 1.$
\item [] $3998=(987+6+5)\times 4+3+2+1.$
\item [] $3999=(987+6+5)\times 4+3\times 2+1.$
\item [] $4000=(9+8\times 7+6+54)\times 32\times 1.$
\item [] $4001=(987+6+5)\times 4+3^2\times 1.$
\item [] $4002=(9+8)\times 7\times 6\times 5+432\times 1.$
\item [] $4003=9\times 8+7+654\times 3\times 2\times 1.$
\item [] $4004=9\times 8+7+654\times 3\times 2+1.$
\item [] $4005=(98+76)\times (5\times 4+3)+2+1.$
\item [] $4006=(9+8+7+65)\times (43+2)+1.$
\item [] $4007=9+8+7\times (6+543+21).$
\item [] $4008=9\times 8\times (7+6\times 5)+4^3\times 21.$
\item [] $4009=(9\times 87+6)\times 5+43+21.$
\item [] $4010=9+(8\times 7+65+4)\times 32+1.$
\item [] $4011=(98+76+5+4\times 3)\times 21.$
\item [] $4012=(9\times 87+6)\times 5+4+3\times 21.$
\item [] $4013=9+8\times 7+6\times (5^4+32+1).$
\item [] $4014=9+8\times (7\times 65+43)+21.$
\item [] $4015=((9+8)\times (7\times 6+5)+4)\times (3+2)\times 1.$
\item [] $4016=(987+6+5)\times 4+3+21.$
\item [] $4017=9+8\times (7\times 65+43+2+1).$
\item [] $4018=(9+8+(7+6)\times 5)\times (4+3)^2\times 1.$
\item [] $4019=9+8\times 76+54\times 3\times 21.$
\item [] $4020=9+87+654\times 3\times 2\times 1.$
\item [] $4021=9+87+654\times 3\times 2+1.$
\item [] $4022=(98\times 7+654)\times 3+2\times 1.$
\item [] $4023=(98+76)\times (5\times 4+3)+21.$
\item [] $4024=(987+6+5)\times 4+32\times 1.$
\item [] $4025=(987+6+5)\times 4+32+1.$
\item [] $4026=(9+8\times 76+54)\times 3\times 2\times 1.$
\item [] $4027=(9+8\times 76+54)\times 3\times 2+1.$
\item [] $4028=98\times (7+6\times 5+4)+3^2+1.$
\item [] $4029=98+7+654\times 3\times 2\times 1.$
\item [] $4030=98+7+654\times 3\times 2+1.$
\item [] $4031=(9\times 87+6)\times 5+43\times 2\times 1.$
\item [] $4032=(9\times 87+6)\times 5+43\times 2+1.$
\item [] $4033=(9\times 8+7+6\times 5)\times (4+32+1).$
\item [] $4034=(9+87)\times (6\times 5+4\times 3)+2\times 1.$
\item [] $4035=(9+87)\times (6\times 5+4\times 3)+2+1.$
\item [] $4036=(9+8+(7\times 6+5^4)\times 3)\times 2\times 1.$
\item [] $4037=(9+8+(7\times 6+5^4)\times 3)\times 2+1.$
\item [] $4038=9\times 8+(7+654)\times 3\times 2\times 1.$
\item [] $4039=9\times 8+(7+654)\times 3\times 2+1.$
\item [] $4040=((9+8)\times 76+54)\times 3+2\times 1.$
\item [] $4041=((9+8)\times 76+54)\times 3+2+1.$
\item [] $4042=98\times (7+6\times 5+4)+3+21.$
\item [] $4043=9\times 87+6\times 543+2\times 1.$
\item [] $4044=9\times 87+6\times 543+2+1.$
\item [] $4045=9+8+76\times (5\times 4+32+1).$
\item [] $4046=(9\times 8+76+5^4\times 3)\times 2\times 1.$
\item [] $4047=98+7+6\times (5^4+32)\times 1.$
\item [] $4048=98+7+6\times (5^4+32)+1.$
\item [] $4049=9+8\times (7+65+432+1).$
\item [] $4050=98+76\times (5\times 4+32)\times 1.$
\item[]$\mbox{Increasing order}$
\item [] $4051=1^2\times 3+4^5+6\times 7\times 8\times 9.$
\item [] $4052=1^2+3+4^5+6\times 7\times 8\times 9.$
\item [] $4053=12+(3+4)\times 567+8\times 9.$
\item [] $4054=1+2+3+4^5+6\times 7\times 8\times 9.$
\item [] $4055=1+2\times 3+4^5+6\times 7\times 8\times 9.$
\item [] $4056=1\times 2^3+4^5+6\times 7\times 8\times 9.$
\item [] $4057=1+2^3+4^5+6\times 7\times 8\times 9.$
\item [] $4058=1^2\times (3+4)\times 567+89.$
\item [] $4059=1\times 23\times 45+6\times 7\times 8\times 9.$
\item [] $4060=1\times 2+(3+4)\times 567+89.$
\item [] $4061=1+2+(3+4)\times 567+89.$
\item [] $4062=(12+3)\times 4\times 56+78\times 9.$
\item [] $4063=12+3+4^5+6\times 7\times 8\times 9.$
\item [] $4064=12^3+4\times (567+8+9).$
\item [] $4065=12\times (3+45\times 6)+789.$
\item [] $4066=1+2\times 3+45\times 6\times (7+8)+9.$
\item [] $4067=1\times 23\times 4+5\times (6+789).$
\item [] $4068=12^3+4\times 567+8\times 9.$
\item [] $4069=12+3+4+5\times 6\times (7+8)\times 9.$
\item [] $4070=(123+456)\times 7+8+9.$
\item [] $4071=1\times 23+4^5+6\times 7\times 8\times 9.$
\item [] $4072=1+23+4^5+6\times 7\times 8\times 9.$
\item [] $4073=1\times 23+(4+5)\times (6\times 7+8)\times 9.$
\item [] $4074=12+3\times 4+5\times 6\times (7+8)\times 9.$
\item [] $4075=(1+2)^3+4^5+6\times 7\times 8\times 9.$
\item [] $4076=123\times (4\times 5+6+7)+8+9.$
\item [] $4077=1\times 23+4+5\times 6\times (7+8)\times 9.$
\item [] $4078=1\times 2\times 34\times 5+6\times 7\times 89.$
\item [] $4079=1+2\times 34\times 5+6\times 7\times 89.$
\item [] $4080=1^2+3456+7\times 89.$
\item [] $4081=12^3+4+5\times 6\times 78+9.$
\item [] $4082=1+2+3456+7\times 89.$
\item [] $4083=1+23+45\times 6\times (7+8)+9.$
\item [] $4084=12\times 3+4^5+6\times 7\times 8\times 9.$
\item [] $4085=1\times 2+345+6\times 7\times 89.$
\item [] $4086=1+2+345+6\times 7\times 89.$
\item [] $4087=1+2+34+5\times 6\times (7+8)\times 9.$
\item [] $4088=(1^2\times 3+4)\times (567+8+9).$
\item [] $4089=(123\times 4+5+6+7)\times 8+9.$
\item [] $4090=(1+2\times 3)^4+5\times 6\times 7\times 8+9.$
\item [] $4091=12+3456+7\times 89.$
\item [] $4092=1^2\times 3\times 4\times 5\times 67+8\times 9.$
\item [] $4093=1^2+3\times 4\times 5\times 67+8\times 9.$
\item [] $4094=1\times 2+3\times 4\times 5\times 67+8\times 9.$
\item [] $4095=12+345+6\times 7\times 89.$
\item [] $4096=12+34+5\times 6\times (7+8)\times 9.$
\item [] $4097=1\times 2+(3+4+56)\times (7\times 8+9).$
\item [] $4098=1\times 2^3\times 45+6\times 7\times 89.$
\item [] $4099=1+2^3\times 45+6\times 7\times 89.$
\item [] $4100=12+(3+4)\times (567+8+9).$
\item [] $4101=(1+2)\times (3+4\times 5\times 67)+8\times 9.$
\item [] $4102=123+4+5\times (6+789).$
\item [] $4103=(12+34)\times (5+6+78)+9.$
\item [] $4104=12+3\times 4\times 5\times 67+8\times 9.$
\item [] $4105=1+2\times (3+45)\times 6\times 7+8\times 9.$
\item [] $4106=1+2^3\times (456+7\times 8)+9.$
\item [] $4107=12+(3+4+56)\times (7\times 8+9).$
\item [] $4108=(1^{23}+45+6)\times (7+8\times 9).$
\item [] $4109=1^2\times 3\times 4\times 5\times 67+89.$
\item [] $4110=1^2+3\times 4\times 5\times 67+89.$
\item [] $4111=1\times 2+3\times 4\times 5\times 67+89.$
\item [] $4112=1+2+3\times 4\times 5\times 67+89.$
\item [] $4113=12+3\times 4\times (5+6\times 7\times 8)+9.$
\item [] $4114=1234+5\times 6\times (7+89).$
\item [] $4115=(12+3)\times 45\times 6+7\times 8+9.$
\item [] $4116=1\times 2+34\times (56+7\times 8+9).$
\item [] $4117=12\times 3^4+56\times 7\times 8+9.$
\item [] $4118=1\times 2\times 34+5\times 6\times (7+8)\times 9.$
\item [] $4119=12\times 3\times 4+5\times (6+789).$
\item [] $4120=(1^2+3+4)\times (5+6+7\times 8\times 9).$
\item[]$\mbox{Decreasing order}$
\item [] $4051=98+76\times (5\times 4+32)+1.$
\item [] $4052=(9+8)\times 7\times (6\times 5+4)+3\times 2\times 1.$
\item [] $4053=(98+7)\times 6\times 5+43\times 21.$
\item [] $4054=9\times 87+654\times (3+2)+1.$
\item [] $4055=(987+6+5)\times 4+3\times 21.$
\item [] $4056=9+87+6\times 5\times 4\times (32+1).$
\item [] $4057=9+8\times (7\times 6\times 5+43)\times 2\times 1.$
\item [] $4058=98\times (7+6\times 5)+432\times 1.$
\item [] $4059=98\times (7+6\times 5)+432+1.$
\item [] $4060=9\times (8+7)\times 6\times 5+4+3\times 2\times 1.$
\item [] $4061=9\times (8+7)\times 6\times 5+4+3\times 2+1.$
\item [] $4062=9\times 87+6\times 543+21.$
\item [] $4063=9\times (8+7)\times 6\times 5+4+3^2\times 1.$
\item [] $4064=98+(7+654)\times 3\times 2\times 1.$
\item [] $4065=98+(7+654)\times 3\times 2+1.$
\item [] $4066=(9+8+7\times 6\times (5+43))\times 2\times 1.$
\item [] $4067=98+(76+5)\times (4+3)^2\times 1.$
\item [] $4068=(9+8+7+654)\times 3\times 2\times 1.$
\item [] $4069=(9+8+7+654)\times 3\times 2+1.$
\item [] $4070=98\times 7+6\times (543+21).$
\item [] $4071=(9+8+7\times 6)\times (5+43+21).$
\item [] $\mathit{4072=9+876\times 5+4-321.}$
\item [] $4073=(987+6\times 5)\times 4+3+2\times 1.$
\item [] $4074=9\times (8+7)\times 6\times 5+4\times 3\times 2\times 1.$
\item [] $4075=9\times (8+7)\times 6\times 5+4\times 3\times 2+1.$
\item [] $4076=(9+8+(7\times 6+5)\times 43)\times 2\times 1.$
\item [] $4077=(987+6\times 5)\times 4+3^2\times 1.$
\item [] $4078=9\times (8+7)\times 6\times 5+4+3+21.$
\item [] $4079=9\times 8\times 7+(6+5)\times (4+321).$
\item [] $4080=(98+7+65)\times 4\times 3\times 2\times 1.$
\item [] $4081=9+8\times 7\times 65+432\times 1.$
\item [] $4082=9+8\times 7\times 65+432+1.$
\item [] $4083=9\times (8+7)\times 6\times 5+4\times 3+21.$
\item [] $4084=9+8+7\times (65\times 4+321).$
\item [] $4085=(9\times 8\times (7+6)+5)\times 4+321.$
\item [] $4086=9\times (87+6)+(54+3)^2\times 1.$
\item [] $4087=9\times (8+7)\times 6\times 5+4+32+1.$
\item [] $4088=98+7\times (6+543+21).$
\item [] $4089=9+8\times (76+5+4)\times 3\times 2\times 1.$
\item [] $4090=9+8\times (76+5+4)\times 3\times 2+1.$
\item [] $4091=9+8+7\times 6\times ((5+43)\times 2+1).$
\item [] $4092=(987+6\times 5)\times 4+3+21.$
\item [] $4093=9\times (8\times 7\times 6+5)+4^(3+2)\times 1.$
\item [] $4094=98\times 7+6+54\times 3\times 21.$
\item [] $4095=9+8+7+6\times 5^4+321.$
\item [] $4096=9+8\times 7\times 6+5^4\times 3\times 2+1.$
\item [] $4097=(9+8)\times (76+54\times 3+2+1).$
\item [] $4098=(98+76+5^4\times 3)\times 2\times 1.$
\item [] $4099=9\times (8+7\times 6)\times 5+43^2\times 1.$
\item [] $4100=(987+6\times 5)\times 4+32\times 1.$
\item [] $4101=(987+6\times 5)\times 4+32+1.$
\item [] $\mathit{4102=98\times 7\times 6-5\times 4+3\times 2\times 1.}$
\item [] $4103=9+(8\times 7\times 6+5)\times 4\times 3+2\times 1.$
\item [] $4104=9\times 8+7\times 6\times (5+43)\times 2\times 1.$
\item [] $4105=9\times 8+7\times 6\times (5+43)\times 2+1.$
\item [] $4106=(9+8\times 7+6\times 5)\times 43+21.$
\item [] $4107=(98+7+6)\times (5\times (4+3)+2\times 1).$
\item [] $4108=9+(87+654\times 3)\times 2+1.$
\item [] $4109=(9+8)\times 7\times (6\times 5+4)+3\times 21.$
\item [] $4110=9\times 8\times 7\times 6+543\times 2\times 1.$
\item [] $4111=9\times 8\times 7\times 6+543\times 2+1.$
\item [] $4112=9+8+7\times 65\times (4+3+2)\times 1.$
\item [] $4113=9+8\times (76+5+432\times 1).$
\item [] $4114=9\times (8+7)\times 6\times 5+43+21.$
\item [] $4115=9\times 8+(7\times 6+5)\times 43\times 2+1.$
\item [] $4116=9\times (8+7)\times 6\times 5+4^3+2\times 1.$
\item [] $4117=9\times (8+7)\times 6\times 5+4+3\times 21.$
\item [] $4118=987+6+5^4\times (3+2)\times 1.$
\item [] $4119=987+6+5^4\times (3+2)+1.$
\item [] $\mathit{4120=98\times 7\times 6-5\times 4+3+21.}$
\item[]$\mbox{Increasing order}$
\item [] $4121=12+3\times 4\times 5\times 67+89.$
\item [] $4122=1+2\times (3+45)\times 6\times 7+89.$
\item [] $4123=(1+2\times 3)\times (4+(5+67)\times 8+9).$
\item [] $\mathit{4124=12\times 345-6+7-8-9.}$
\item [] $4125=(123+456)\times 7+8\times 9.$
\item [] $4126=12+34\times (56+7\times 8+9).$
\item [] $4127=1\times 2+3\times 4\times (5\times 67+8)+9.$
\item [] $4128=(1+2)\times 3\times 456+7+8+9.$
\item [] $4129=(12+3)\times 45\times 6+7+8\times 9.$
\item [] $4130=1\times 2\times (3+4)\times 5\times (6\times 7+8+9).$
\item [] $4131=123\times (4+5)+6\times 7\times 8\times 9.$
\item [] $4132=123\times 4+56\times (7\times 8+9).$
\item [] $4133=(12+3^4+5)\times 6\times 7+8+9.$
\item [] $4134=12\times (3+4)+5\times 6\times (7+8)\times 9.$
\item [] $4135=1+(2+3+45\times 6)\times (7+8)+9.$
\item [] $4136=1\times 2^3\times (4+((56+7)\times 8+9)).$
\item [] $4137=(12+3)\times 45\times 6+78+9.$
\item [] $4138=1+2\times 3^4+5\times (6+789).$
\item [] $\mathit{4139=-1\times 2\times 34-5+6\times 78\times 9.}$
\item [] $4140=(12+34)\times (5+6+7+8\times 9).$
\item [] $4141=1^2+(3+4+5)\times (6\times 7\times 8+9).$
\item [] $4142=(123+456)\times 7+89.$
\item [] $4143=1^2\times 3^4\times 5+6\times 7\times 89.$
\item [] $4144=1^2+3^4\times 5+6\times 7\times 89.$
\item [] $4145=1\times 2+3^4\times 5+6\times 7\times 89.$
\item [] $4146=1+2+3^4\times 5+6\times 7\times 89.$
\item [] $4147=1+(23\times 4+5)\times 6\times 7+8\times 9.$
\item [] $4148=123\times (4\times 5+6+7)+89.$
\item [] $4149=(12+3)\times 4\times 56+789.$
\item [] $4150=1+23\times 4\times (5\times 6+7+8)+9.$
\item [] $4151=12\times 34+5+6\times 7\times 89.$
\item [] $4152=1\times 2^3\times 456+7\times 8\times 9.$
\item [] $4153=1+2^3\times 456+7\times 8\times 9.$
\item [] $4154=1+(23+45+6)\times 7\times 8+9.$
\item [] $4155=12+3^4\times 5+6\times 7\times 89.$
\item [] $4156=1+2^{(3+4+5)}+6\times 7+8+9.$
\item [] $\mathit{4157=1-2+3456+78\times 9.}$
\item [] $4158=1^2\times 3456+78\times 9.$
\item [] $4159=1^2+3456+78\times 9.$
\item [] $4160=1\times 2+3456+78\times 9.$
\item [] $4161=1+2+3456+78\times 9.$
\item [] $4162=1\times 2+(34+5\times 6)\times (7\times 8+9).$
\item [] $4163=1\times (23\times 4+5)\times 6\times 7+89.$
\item [] $4164=1\times 2\times 345\times 6+7+8+9.$
\item [] $4165=1+2\times 345\times 6+7+8+9.$
\item [] $4166=1+(2+3)^4+5\times (6+78\times 9).$
\item [] $4167=(12\times 34+5+6\times 7+8)\times 9.$
\item [] $4168=1+2^{(3+4+5)}+6+7\times 8+9.$
\item [] $4169=(1+2)\times 3\times 456+7\times 8+9.$
\item [] $4170=12+3456+78\times 9.$
\item [] $4171=123+4^5+6\times 7\times 8\times 9.$
\item [] $4172=1+2\times (345\times 6+7)+8+9.$
\item [] $4173=123+(4+5)\times (6\times 7+8)\times 9.$
\item [] $\mathit{4174=1+2\times 3-45+6\times 78\times 9.}$
\item [] $4175=1\times 2\times (3+4\times 5\times (6+7)\times 8)+9.$
\item [] $4176=12\times (3+4+5\times 67)+8\times 9.$
\item [] $4177=123+4+5\times 6\times (7+8)\times 9.$
\item [] $4178=1\times 2^{(3+4)}+5\times 6\times (7+8)\times 9.$
\item [] $4179=1+2^{(3\times 4)}+5\times (6+7)+8+9.$
\item [] $4180=1\times 2^{(3+4+5)}+67+8+9.$
\item [] $4181=12\times (3\times 4+5\times 67)+8+9.$
\item [] $4182=(12+3+4\times 56+7)\times (8+9).$
\item [] $4183=(1+2)\times 3\times 456+7+8\times 9.$
\item [] $4184=(1+2)\times 3\times (456+7)+8+9.$
\item [] $4185=1\times 2\times 3\times (45+6\times 7)\times 8+9.$
\item [] $4186=1+2\times 3\times (45+6\times 7)\times 8+9.$
\item [] $4187=1+2\times (3+4)\times (5\times 6\times 7+89).$
\item [] $4188=12\times 34+5\times (6+78)\times 9.$
\item [] $4189=1+2\times (345\times 6+7+8+9).$
\item [] $4190=1\times 2+3+45\times (6+78+9).$
\item[]$\mbox{Decreasing order}$
\item [] $4121=9+8\times (76+5+432+1).$
\item [] $4122=9+8+76\times (5+4)\times 3\times 2+1.$
\item [] $4123=(9+8\times (7+65)+4)\times (3\times 2+1).$
\item [] $\mathit{4124=98\times 7\times 6-5+4\times 3+2-1.}$
\item [] $4125=(9+8\times 7+6+54)\times (32+1).$
\item [] $4126=9+8+76\times 54+3+2\times 1.$
\item [] $4127=9+8+76\times 54+3+2+1.$
\item [] $4128=9+8+76\times 54+3\times 2+1.$
\item [] $4129=(9\times (8+7+6)+5^4\times 3)\times 2+1.$
\item [] $4130=98\times 7\times 6+5+4+3+2\times 1.$
\item [] $4131=9+8+76\times 54+3^2+1.$
\item [] $4132=98\times 7\times 6+5+4+3\times 2+1.$
\item [] $4133=(9+8)\times (7\times (6+5)+4)\times 3+2\times 1.$
\item [] $4134=98\times 7\times 6+5+4+3^2\times 1.$
\item [] $4135=98\times 7\times 6+5+4+3^2+1.$
\item [] $4136=9+8\times 7+6\times 5^4+321.$
\item [] $4137=(98+76+5\times 4+3)\times 21.$
\item [] $4138=(9+8+76\times (5+4)\times 3)\times 2\times 1.$
\item [] $4139=9\times 8+7\times (65\times 4+321).$
\item [] $4140=98+(7\times 6+5)\times 43\times 2\times 1.$
\item [] $4141=98\times 7\times 6+5\times 4+3+2\times 1.$
\item [] $4142=98\times 7\times 6+5\times 4+3+2+1.$
\item [] $4143=98\times 7\times 6+5\times 4+3\times 2+1.$
\item [] $\mathit{4144=98\times 7\times 6-5+4\times 3+21.}$
\item [] $4145=9+8+76\times 54+3+21.$
\item [] $4146=98\times 7\times 6+5+4\times 3\times 2+1.$
\item [] $4147=(9+(8\times 7+6)\times 5)\times (4+3^2)\times 1.$
\item [] $4148=9+8+76\times 54+3^{(2+1)}.$
\item [] $4149=98\times 7\times 6+5+4+3+21.$
\item [] $4150=98\times 7\times 6+(5+4\times 3)\times 2\times 1.$
\item [] $4151=98\times 7\times 6+(5+4\times 3)\times 2+1.$
\item [] $4152=9+8+7+6\times (5^4+3\times 21).$
\item [] $4153=9+8+76\times 54+32\times 1.$
\item [] $4154=98\times 7\times 6+5+4\times 3+21.$
\item [] $4155=(9+8\times 7+6)\times 54+321.$
\item [] $\mathit{4156=9\times 87\times 6-543+2-1.}$
\item [] $4157=98\times 7\times 6+5+4+32\times 1.$
\item [] $4158=98\times 7\times 6+5+4+32+1.$
\item [] $4159=9+(8+7\times 6)\times (5\times 4+3\times 21).$
\item [] $4160=98\times 7\times 6+5\times 4+3+21.$
\item [] $4161=9+8\times (7\times 65+43+21).$
\item [] $4162=98+(7+6\times 5\times 4)\times 32\times 1.$
\item [] $4163=98+(7+6\times 5\times 4)\times 32+1.$
\item [] $4164=(9+87)\times 6\times 5+4\times 321.$
\item [] $4165=98+7\times (65\times 4+321).$
\item [] $4166=98\times 7\times 6+5+43+2\times 1.$
\item [] $4167=98\times 7\times 6+5+43+2+1.$
\item [] $4168=98\times 7\times 6+5\times 4+32\times 1.$
\item [] $4169=98\times 7\times 6+5\times 4+32+1.$
\item [] $4170=98\times 7\times 6+(5+4)\times 3\times 2\times 1.$
\item [] $4171=98\times 7\times 6+5+(4+3)^2+1.$
\item [] $4172=987+65\times (4+3)^2\times 1.$
\item [] $4173=987+65\times (4+3)^2+1.$
\item [] $4174=(98+7)\times 6\times 5+4^(3+2)\times 1.$
\item [] $4175=98\times 7\times 6+54+3+2\times 1.$
\item [] $4176=98+7+6\times 5^4+321.$
\item [] $4177=98\times 7\times 6+54+3\times 2+1.$
\item [] $4178=98\times 7\times 6+5\times 4\times 3+2\times 1.$
\item [] $4179=98\times 7\times 6+54+3^2\times 1.$
\item [] $4180=98\times 7\times 6+54+3^2+1.$
\item [] $4181=9\times 8+76\times 54+3+2\times 1.$
\item [] $4182=9\times 8+76\times 54+3+2+1.$
\item [] $4183=9\times 8+76\times 54+3\times 2+1.$
\item [] $4184=9+8+76\times 54+3\times 21.$
\item [] $4185=98\times 7\times 6+5+43+21.$
\item [] $4186=9+8+7+65\times 4^3+2\times 1.$
\item [] $4187=98\times 7\times 6+5+4^3+2\times 1.$
\item [] $4188=98\times 7\times 6+5+4+3\times 21.$
\item [] $4189=(9+8\times 7+6)\times (54+3+2\times 1).$
\item [] $4190=(9+8)\times 7+6\times 5^4+321.$
\item[]$\mbox{Increasing order}$
\item [] $4191=(1+2)\times 3\times 456+78+9.$
\item [] $4192=1+2\times 3+45\times (6+78+9).$
\item [] $4193=12\times (3+4+5\times 67)+89.$
\item [] $4194=1\times 234\times 5+6\times 7\times 8\times 9.$
\item [] $4195=1+234\times 5+6\times 7\times 8\times 9.$
\item [] $4196=1+(2+3)^4+5\times 6\times 7\times (8+9).$
\item [] $4197=(1+2)\times (3\times 456+7)+8\times 9.$
\item [] $4198=1\times 23\times 4\times 5+6\times 7\times 89.$
\item [] $4199=12\times 345+6\times 7+8+9.$
\item [] $4200=(1+2^3)\times 456+7+89.$
\item [] $4201=1+(2+3)\times (45+6+789).$
\item [] $4202=1\times 2+3\times 4\times (5+6\times 7\times 8+9).$
\item [] $4203=(12+3^4)\times 5+6\times 7\times 89.$
\item [] $4204=1+2^{(3\times 4)}+5+6+7+89.$
\item [] $4205=1\times 2\times 345\times 6+7\times 8+9.$
\item [] $4206=123\times (4+5\times 6)+7+8+9.$
\item [] $4207=(1^2+3)\times (4^5+6)+78+9.$
\item [] $4208=1\times 23+45\times (6+78+9).$
\item [] $4209=1+23+45\times (6+78+9).$
\item [] $4210=1+234+5\times (6+789).$
\item [] $4211=12\times 345+6+7\times 8+9.$
\item [] $4212=12+3\times 4\times (5+6\times 7\times 8+9).$
\item [] $4213=1\times 2+3\times 4^5+67\times (8+9).$
\item [] $4214=1+2+3\times 4^5+67\times (8+9).$
\item [] $4215=1+2+3\times (4+5)\times (67+89).$
\item [] $4216=1+2^{(3\times 4)}+5+6\times 7+8\times 9.$
\item [] $4217=1^{234}\times 5+6\times 78\times 9.$
\item [] $4218=(1+23)\times 4\times 5+6\times 7\times 89.$
\item [] $4219=1\times 2\times 345\times 6+7+8\times 9.$
\item [] $4220=1+2\times 345\times 6+7+8\times 9.$
\item [] $4221=1^{23}\times 4+5+6\times 78\times 9.$
\item [] $4222=1^{23}+4+5+6\times 78\times 9.$
\item [] $4223=12+3\times 4^5+67\times (8+9).$
\item [] $4224=12\times 345+67+8+9.$
\item [] $4225=12\times 345+6+7+8\times 9.$
\item [] $4226=1\times 2+3+4+5+6\times 78\times 9.$
\item [] $4227=1+2+3+4+5+6\times 78\times 9.$
\item [] $4228=1+2\times 345\times 6+78+9.$
\item [] $4229=1^2\times 3\times 4+5+6\times 78\times 9.$
\item [] $4230=1+2^3+4+5+6\times 78\times 9.$
\item [] $4231=1\times 2+3\times 4+5+6\times 78\times 9.$
\item [] $4232=1+2+3\times 4+5+6\times 78\times 9.$
\item [] $4233=12\times 345+6+78+9.$
\item [] $4234=1+2^{(3\times 4)}+5\times (6+7)+8\times 9.$
\item [] $4235=123\times 4+5+6\times 7\times 89.$
\item [] $4236=12+3+4+5+6\times 78\times 9.$
\item [] $4237=1+2\times 345\times 6+7+89.$
\item [] $4238=1+2+3+4\times 5+6\times 78\times 9.$
\item [] $4239=1+2\times 3+4\times 5+6\times 78\times 9.$
\item [] $4240=1\times 2^3+4\times 5+6\times 78\times 9.$
\item [] $4241=12+3\times 4+5+6\times 78\times 9.$
\item [] $4242=12\times 345+6+7+89.$
\item [] $4243=1\times 2\times (345\times 6+7)+89.$
\item [] $4244=1\times 23+4+5+6\times 78\times 9.$
\item [] $4245=1^2\times 3456+789.$
\item [] $4246=1^2+3456+789.$
\item [] $4247=1\times 2+3456+789.$
\item [] $4248=1+2+3456+789.$
\item [] $4249=1\times 2^3\times 4+5+6\times 78\times 9.$
\item [] $4250=1+2^3\times 4+5+6\times 78\times 9.$
\item [] $4251=1^2\times 34+5+6\times 78\times 9.$
\item [] $4252=1^2+34+5+6\times 78\times 9.$
\item [] $4253=1\times 2+34+5+6\times 78\times 9.$
\item [] $4254=12\times 345+6\times 7+8\times 9.$
\item [] $4255=1\times 23+4\times 5+6\times 78\times 9.$
\item [] $4256=1+23+4\times 5+6\times 78\times 9.$
\item [] $4257=12+3456+789.$
\item [] $4258=1^{23}+45+6\times 78\times 9.$
\item [] $4259=12+(3+4)\times 5+6\times 78\times 9.$
\item [] $4260=1^2\times 3+45+6\times 78\times 9.$
\item[]$\mbox{Decreasing order}$
\item [] $4191=9\times 87+6+54\times 3\times 21.$
\item [] $4192=(98+7+6+5\times 4)\times 32\times 1.$
\item [] $4193=(98+7+6+5\times 4)\times 32+1.$
\item [] $4194=98\times 7\times 6+54+3+21.$
\item [] $4195=9\times (87+6)\times 5+4+3+2+1.$
\item [] $4196=9\times (87+6)\times 5+4+3\times 2+1.$
\item [] $4197=98\times 7\times 6+5\times 4\times 3+21.$
\item [] $4198=9+8+(7+6)\times 5\times 4^3+21.$
\item [] $4199=98\times 7\times 6+5\times 4+3\times 21.$
\item [] $4200=9\times 8+76\times 54+3+21.$
\item [] $4201=9\times (8+7\times 6)+5^4\times 3\times 2+1.$
\item [] $4202=98\times 7\times 6+54+32\times 1.$
\item [] $4203=98\times 7\times 6+54+32+1.$
\item [] $4204=9\times (8+7\times 65)+4+32+1.$
\item [] $4205=9+8+7+65\times 4^3+21.$
\item [] $4206=98\times 7\times 6+5+4^3+21.$
\item [] $4207=98\times 7\times 6+5+43\times 2\times 1.$
\item [] $4208=98+76\times 54+3+2+1.$
\item [] $4209=98+76\times 54+3\times 2+1.$
\item [] $4210=9\times (87+6)\times 5+4\times 3\times 2+1.$
\item [] $4211=98+76\times 54+3^2\times 1.$
\item [] $4212=98+76\times 54+3^2+1.$
\item [] $4213=9\times (8+76)\times 5+432+1.$
\item [] $4214=(987+65)\times 4+3+2+1.$
\item [] $4215=(987+65)\times 4+3\times 2+1.$
\item [] $4216=(9+8+76+5)\times 43+2\times 1.$
\item [] $4217=(987+65)\times 4+3^2\times 1.$
\item [] $4218=(987+65)\times 4+3^2+1.$
\item [] $4219=(9\times 8+7)\times 6\times 5+43^2\times 1.$
\item [] $4220=(9\times 8+7)\times 6\times 5+43^2+1.$
\item [] $4221=9\times (87+6)\times 5+4+32\times 1.$
\item [] $4222=9\times (87+6)\times 5+4+32+1.$
\item [] $4223=9+87+6\times 5+4^{(3\times 2)}+1.$
\item [] $4224=9+87+6\times (5^4+3\times 21).$
\item [] $4225=9+8\times 7+65\times (43+21).$
\item [] $4226=98+76\times 54+3+21.$
\item [] $4227=9+8\times 7+65\times 4^3+2\times 1.$
\item [] $4228=9+8\times 7+65\times 4^3+2+1.$
\item [] $4229=9+8+(7+6)\times 54\times 3\times 2\times 1.$
\item [] $4230=9+(8+7)\times 65\times 4+321.$
\item [] $4231=9\times (8+7\times 65)+43+21.$
\item [] $4232=(987+65)\times 4+3+21.$
\item [] $4233=98\times 7\times 6+54+3\times 21.$
\item [] $4234=98+76\times 54+32\times 1.$
\item [] $4235=98+76\times 54+32+1.$
\item [] $4236=98\times 7\times 6+5\times 4\times 3\times 2\times 1.$
\item [] $4237=98\times 7\times 6+5\times 4\times 3\times 2+1.$
\item [] $4238=9+8+(7+6+54)\times 3\times 21.$
\item [] $4239=9\times 8+76\times 54+3\times 21.$
\item [] $4240=(987+65)\times 4+32\times 1.$
\item [] $4241=(987+65)\times 4+32+1.$
\item [] $4242=9\times 8+7+65\times 4^3+2+1.$
\item [] $4243=987+6+(54+3)^2+1.$
\item [] $4244=9\times (87+6\times 5)\times 4+32\times 1.$
\item [] $4245=9\times (87+6\times 5)\times 4+32+1.$
\item [] $4246=9+8\times 7+65\times 4^3+21.$
\item [] $4247=987+6\times 543+2\times 1.$
\item [] $4248=987+6\times 543+2+1.$
\item [] $4249=98\times 7\times 6+5+4\times 32\times 1.$
\item [] $4250=98\times 7\times 6+5+4^3\times 2+1.$
\item [] $4251=(9\times 87+6+5^4+3)\times (2+1).$
\item [] $4252=9\times (87+6)\times 5+4+3\times 21.$
\item [] $4253=(9+8\times 7)\times 65+4+3+21.$
\item [] $4254=9\times (8+7\times 65)+43\times 2+1.$
\item [] $4255=98\times 7\times 6+(5+4^3)\times 2+1.$
\item [] $4256=9+87+65\times (43+21).$
\item [] $4257=987+6\times (543+2)\times 1.$
\item [] $4258=9+87+65\times 4^3+2\times 1.$
\item [] $4259=9+87+65\times 4^3+2+1.$
\item [] $4260=9\times 8+7+65\times 4^3+21.$
\item[]$\mbox{Increasing order}$
\item [] $4261=1^2+3+45+6\times 78\times 9.$
\item [] $4262=1\times 2+3+45+6\times 78\times 9.$
\item [] $4263=1234+5+6\times 7\times 8\times 9.$
\item [] $4264=1+2\times 3+45+6\times 78\times 9.$
\item [] $4265=1\times 2^3+45+6\times 78\times 9.$
\item [] $4266=1+2^3+45+6\times 78\times 9.$
\item [] $4267=1+2\times 3\times (4+5)+6\times 78\times 9.$
\item [] $4268=12\times 3+4\times 5+6\times 78\times 9.$
\item [] $4269=123\times (4+5\times 6)+78+9.$
\item [] $4270=1^2+3+(4+5)\times 6\times (7+8\times 9).$
\item [] $4271=12\times 345+6\times 7+89.$
\item [] $4272=12+3+45+6\times 78\times 9.$
\item [] $4273=1^2+3\times 4\times 5+6\times 78\times 9.$
\item [] $4274=1\times 2+3\times 4\times 5+6\times 78\times 9.$
\item [] $4275=1+2+3\times 4\times 5+6\times 78\times 9.$
\item [] $4276=1+2\times 345\times 6+(7+8)\times 9.$
\item [] $4277=12\times (345+6)+7\times 8+9.$
\item [] $4278=123\times (4+5\times 6)+7+89.$
\item [] $4279=12\times 345+67+8\times 9.$
\item [] $4280=1\times 23+45+6\times 78\times 9.$
\item [] $4281=1+23+45+6\times 78\times 9.$
\item [] $4282=1\times 2\times (3+4)\times 5+6\times 78\times 9.$
\item [] $4283=1+2\times (3+4)\times 5+6\times 78\times 9.$
\item [] $4284=12+3\times 4\times 5+6\times 78\times 9.$
\item [] $4285=1\times 2\times 34+5+6\times 78\times 9.$
\item [] $4286=1+2\times 34+5+6\times 78\times 9.$
\item [] $4287=1+2+34\times (5\times 6+7+89).$
\item [] $4288=1+2^{(3\times 4)}+56+(7+8)\times 9.$
\item [] $4289=(1\times 2+3)\times 4\times 5\times 6\times 7+89.$
\item [] $4290=1+(2+3)\times 4\times 5\times 6\times 7+89.$
\item [] $4291=1+2\times (34+5)+6\times 78\times 9.$
\item [] $4292=(1^2+3)\times 4\times 5+6\times 78\times 9.$
\item [] $4293=12\times 3+45+6\times 78\times 9.$
\item [] $4294=1+(2\times 3+45)\times (6+78)+9.$
\item [] $4295=1\times (2+3)\times (4+(5+6\times (7+8))\times 9).$
\item [] $4296=12\times 345+67+89.$
\item [] $4297=1\times 2^3\times 4\times (56+78)+9.$
\item [] $4298=1^2\times 3^4+5+6\times 78\times 9.$
\item [] $4299=1^2+3^4+5+6\times 78\times 9.$
\item [] $4300=1\times 2+3^4+5+6\times 78\times 9.$
\item [] $4301=1+2+3^4+5+6\times 78\times 9.$
\item [] $4302=1+234\times (5+6+7)+89.$
\item [] $4303=1234+(5+6\times 7\times 8)\times 9.$
\item [] $4304=(1^2+3)^4\times 5+6\times 7\times 8\times 9.$
\item [] $4305=1+2^3\times (4+5\times 6+7\times 8\times 9).$
\item [] $4306=(1^2+3)^4+5\times 6\times (7+8)\times 9.$
\item [] $4307=(1+2^3\times (4+5))\times (6\times 7+8+9).$
\item [] $4308=123+45\times (6+78+9).$
\item [] $4309=1\times (23\times 4+5)+6\times 78\times 9.$
\item [] $4310=12+3^4+5+6\times 78\times 9.$
\item [] $4311=(1+2)^3+4\times (56+7)\times (8+9).$
\item [] $4312=1\times 2\times 34\times 56+7\times 8\times 9.$
\item [] $4313=1+(2+3)\times 4\times 5+6\times 78\times 9.$
\item [] $4314=(123+4)\times 5\times 6+7\times 8\times 9.$
\item [] $4315=1+2\times (3\times 456+789).$
\item [] $4316=(1^2+3)\times (456+7\times 89).$
\item [] $4317=1+23+(4+5)\times (6\times 78+9).$
\item [] $4318=(1\times 2\times (3\times 4\times 5+67))\times (8+9).$
\item [] $4319=(1+2)\times 34+5+6\times 78\times 9.$
\item [] $4320=12\times (3+4+5)\times (6+7+8+9).$
\item [] $4321=1+2\times (34+56)\times (7+8+9).$
\item [] $4322=1\times 2+3\times (4+5+6)\times (7+89).$
\item [] $4323=1+2^{(3\times 4)}+5+(6+7)\times (8+9).$
\item [] $4324=1+2^{(3\times 4)}+5\times 6\times 7+8+9.$
\item [] $4325=(1+2)^3\times 4+5+6\times 78\times 9.$
\item [] $4326=1+(2\times 3)^4+5+6\times 7\times 8\times 9.$
\item [] $4327=12\times 3\times 4+(5+6\times 7)\times 89.$
\item [] $\mathit{4328=1+2-34+56\times 78-9.}$
\item [] $4329=1^2\times (3+45)\times 6\times (7+8)+9.$
\item [] $4330=1+(2+3+4)\times (56\times 7+89).$
\item[]$\mbox{Decreasing order}$
\item [] $4261=9\times 8\times 7+6+5^4\times 3\times 2+1.$
\item [] $4262=(9+8\times 7)\times 65+4+32+1.$
\item [] $4263=9\times 8\times 7+6\times 5^4+3^2\times 1.$
\item [] $4264=9\times 8\times 7+6\times 5^4+3^2+1.$
\item [] $4265=98+76\times 54+3\times 21.$
\item [] $4266=987+6\times 543+21.$
\item [] $4267=98+7+65\times 4^3+2\times 1.$
\item [] $4268=98+7+65\times 4^3+2+1.$
\item [] $4269=9\times (8\times (7+6)+54)\times 3+2+1.$
\item [] $4270=(9\times 87+6)\times 5+4+321.$
\item [] $4271=9\times (87+6)\times 5+43\times 2\times 1.$
\item [] $4272=9\times 8+7\times 6\times 5\times 4\times (3+2)\times 1.$
\item [] $4273=9\times 8+7\times 6\times 5\times 4\times (3+2)+1.$
\item [] $4274=9\times 8\times 7+6\times (5^4+3)+2\times 1.$
\item [] $4275=9\times (8+76\times 5+43\times 2+1).$
\item [] $4276=98+76+5+4^{(3\times 2)}+1.$
\item [] $4277=9+87+65\times 4^3+21.$
\item [] $4278=9\times 8\times 7+6\times 5^4+3+21.$
\item [] $4279=98+(7+6)\times 5\times 4^3+21.$
\item [] $4280=98\times 7\times 6+54\times 3+2\times 1.$
\item [] $4281=9+87\times 6+5^4\times 3\times 2\times 1.$
\item [] $4282=9+87\times 6+5^4\times 3\times 2+1.$
\item [] $\mathit{4283=-9+8-7\times 6+5+4321.}$
\item [] $4284=(9+87+65+43)\times 21.$
\item [] $4285=9\times 8+(7+6)\times 54\times 3\times 2+1.$
\item [] $4286=98+7+65\times 4^3+21.$
\item [] $4287=9+876+54\times 3\times 21.$
\item [] $4288=(9+8\times 7+65+4)\times 32\times 1.$
\item [] $4289=(9+8\times 7)\times 65+43+21.$
\item [] $4290=(9+87+6\times 5+4)\times (32+1).$
\item [] $4291=((9+8\times 7)\times 65+4^3+2)\times 1.$
\item [] $4292=(9+8\times 7)\times 65+4+3\times 21.$
\item [] $4293=9\times (87+65+4+321).$
\item [] $4294=((9+87\times 6+5)\times 4+3)\times 2\times 1.$
\item [] $4295=9\times (8+7\times 65)+4\times 32\times 1.$
\item [] $4296=98\times 7\times 6+5\times (4+32)\times 1.$
\item [] $4297=98\times 7\times 6+5\times 4\times 3^2+1.$
\item [] $4298=98+7\times 6\times 5\times 4\times (3+2)\times 1.$
\item [] $4299=98\times 7\times 6+54\times 3+21.$
\item [] $4300=(9+8)\times 7+65\times 4^3+21.$
\item [] $4301=98\times 7\times 6+5\times (4+32+1).$
\item [] $4302=9\times (8+7)\times 6\times 5+4\times 3\times 21.$
\item [] $4303=9+(8\times 76+5)\times (4+3)+2+1.$
\item [] $4304=(9+8)\times (7\times 6\times 5+43)+2+1.$
\item [] $4305=(9\times 8+76+54+3)\times 21.$
\item [] $\mathit{4306=-9+876\times 5-4^3-2+1.}$
\item [] $4307=9+8+(76+54)\times (32+1).$
\item [] $4308=9\times (8\times 7+6)+5^4\times 3\times 2\times 1.$
\item [] $4309=(9+8+76+5^4)\times 3\times 2+1.$
\item [] $4310=98+(7+6)\times 54\times 3\times 2\times 1.$
\item [] $4311=98+(7+6)\times 54\times 3\times 2+1.$
\item [] $4312=(9+8\times 7)\times 65+43\times 2+1.$
\item [] $4313=9\times 8\times 7\times 6+5+4\times 321.$
\item [] $4314=(9+8\times 7+654)\times 3\times 2\times 1.$
\item [] $4315=98\times (7+6\times 5+4+3)+2+1.$
\item [] $4316=9\times (8+7)+65\times 4^3+21.$
\item [] $4317=9\times 8\times 7+6\times 5^4+3\times 21.$
\item [] $4318=9+(8\times 7+6+5)\times 4^3+21.$
\item [] $4319=98+(7+6+54)\times 3\times 21.$
\item [] $4320=9\times 8\times 7\times 6+54\times (3+21).$
\item [] $4321=(98\times 7+6\times 5+4)\times 3\times 2+1.$
\item [] $4322=(9+8)\times (7\times 6\times 5+43)+21.$
\item [] $4323=(9+8+7)\times (6+54)\times 3+2+1.$
\item [] $4324=9+8\times 7\times (65+4\times 3)+2+1.$
\item [] $4325=((9+8)\times (7+6\times 5\times 4)+3)\times 2+1.$
\item [] $4326=(9+87)\times 6+5^4\times 3\times 2\times 1.$
\item [] $4327=(9+87)\times 6+5^4\times 3\times 2+1.$
\item [] $4328=9+(8\times 76+5+4)\times (3\times 2+1).$
\item [] $4329=9+(8+7+6\times 5\times 4)\times 32\times 1.$
\item [] $4330=9+8\times (7+65\times 4+3)\times 2+1.$
\item[]$\mbox{Increasing order}$
\item [] $4331=(1+(2+3)\times 4\times 5)\times 6\times 7+89.$
\item [] $4332=1\times 2\times (3+4^5+67\times (8+9)).$
\item [] $4333=1+2\times 3\times 4\times 5+6\times 78\times 9.$
\item [] $\mathit{4334=123+4-5+6\times 78\times 9.}$
\item [] $4335=(12+34+5)\times ((6+7)+8\times 9).$
\item [] $4336=1+(2\times 3+45)\times ((6+7)+8\times 9).$
\item [] $4337=(1+2\times 3\times 4)\times 5+6\times 78\times 9.$
\item [] $4338=(1\times 2+3^4\times 5+67+8)\times 9.$
\item [] $4339=1+(2+3+4)\times (5+6\times 78+9).$
\item [] $4340=12^3+4\times (5\times 6+7\times 89).$
\item [] $4341=(1+2)\times (3\times 456+7+8\times 9).$
\item [] $4342=1+23\times 4\times (5+6\times 7)+8+9.$
\item [] $4343=1\times 2\times (3+4\times (5+67\times 8))+9.$
\item [] $4344=123+(4+5)\times 6\times 78+9.$
\item [] $4345=1+2^{(3+4)}\times 5\times 6+7\times 8\times 9.$
\item [] $4346=1+(2^{(3+4)}+5)+6\times 78\times 9.$
\item [] $4347=1^2\times 3\times 45+6\times 78\times 9.$
\item [] $4348=1^2+3\times 45+6\times 78\times 9.$
\item [] $4349=1\times 2+3\times 45+6\times 78\times 9.$
\item [] $4350=(1+23+4\times 5+6)\times (78+9).$
\item [] $4351=(1+2\times 3)^4+5\times 6\times (7\times 8+9).$
\item [] $4352=1\times 2^3\times 4\times (5+6\times 7+89).$
\item [] $4353=1+2\times (3+4\times (5+67\times 8)+9).$
\item [] $\mathit{4354=1\times 2-3-4+56\times 78-9.}$
\item [] $4355=123+4\times 5+6\times 78\times 9.$
\item [] $4356=(1+2\times 3+4+5+6\times 78)\times 9.$
\item [] $4357=1+(2+34)\times (56+7\times 8+9).$
\item [] $4358=1\times 2+(34+5\times 6\times (7+8))\times 9.$
\item [] $4359=12+3\times 45+6\times 78\times 9.$
\item [] $4360=(12\times 3+4)\times (5\times 6+7+8\times 9).$
\item [] $4361=12\times 3\times 4+5+6\times 78\times 9.$
\item [] $4362=123+(456+7+8)\times 9.$
\item [] $4363=(1+2\times 3)^4+(5\times 6\times 7+8)\times 9.$
\item [] $4364=(1+2\times 34)\times (56+7)+8+9.$
\item [] $4365=(1+2)\times (3\times 456+78+9).$
\item [] $4366=(1+(2\times 34+5))\times (6\times 7+8+9).$
\item [] $4367=((1+2)^3+4)\times 5+6\times 78\times 9.$
\item [] $4368=1\times 2\times 3\times (4\times 56+7\times 8\times 9).$
\item [] $4369=1+(2+3)^4+5+6\times 7\times 89.$
\item [] $4370=(12+34)\times (5\times 6+7\times 8+9).$
\item [] $4371=12\times 3^4+5\times 678+9.$
\item [] $4372=1\times 2^3\times 4\times 5+6\times 78\times 9.$
\item [] $4373=1+2^3\times 4\times 5+6\times 78\times 9.$
\item [] $4374=(1+2^{(3+4)})\times 5\times 6+7\times 8\times 9.$
\item [] $4375=1^2\times (3+4)\times 5\times (6+7\times (8+9)).$
\item [] $4376=1\times 2+3^4\times (5\times 6+7+8+9).$
\item [] $4377=1^{234}\times 56\times 78+9.$
\item [] $4378=1^{234}+56\times 78+9.$
\item [] $4379=1234+56\times 7\times 8+9.$
\item [] $4380=123+45+6\times 78\times 9.$
\item [] $4381=1^{23}\times 4+56\times 78+9.$
\item [] $4382=1^{23}+4+56\times 78+9.$
\item [] $4383=1^2+34\times 5+6\times 78\times 9.$
\item [] $4384=1^2\times 3+4+56\times 78+9.$
\item [] $4385=12\times 3\times 4\times 5\times 6+7\times 8+9.$
\item [] $4386=1\times (23+4\times 5)\times ((6+7)+89).$
\item [] $4387=1+2+3+4+56\times 78+9.$
\item [] $4388=(12+3^4)\times (5+6\times 7)+8+9.$
\item [] $4389=1^2\times (3\times 4+56\times 78+9).$
\item [] $4390=1^2+3\times 4+56\times 78+9.$
\item [] $4391=1\times 2+3\times 4+56\times 78+9.$
\item [] $4392=1+2^{(3\times 4)}+5\times (6\times 7+8+9).$
\item [] $4393=1+(2+34)\times 5+6\times 78\times 9.$
\item [] $4394=12+34\times 5+6\times 78\times 9.$
\item [] $4395=1\times 2^{(3\times 4)}+5\times 6\times 7+89.$
\item [] $4396=12+3+4+56\times 78+9.$
\item [] $4397=(1+2+34)\times 5+6\times 78\times 9.$
\item [] $4398=1+(2+3)\times 4+56\times 78+9.$
\item [] $4399=12\times 3\times 4\times 5\times 6+7+8\times 9.$
\item [] $4400=(1^2+3+4)\times (5+67\times 8+9).$
\item[]$\mbox{Decreasing order}$
\item [] $4331=9+(8+7)\times 6\times (5+43)+2\times 1.$
\item [] $4332=(9\times 87+654)\times 3+21.$
\item [] $4333=98\times 7\times 6+5\times 43+2\times 1.$
\item [] $4334=98\times 7\times 6+5\times 43+2+1.$
\item [] $4335=(98+765+4)\times (3+2)\times 1.$
\item [] $4336=(98+765+4)\times (3+2)+1.$
\item [] $4337=9+8\times (7\times 65+43\times 2\times 1).$
\item [] $4338=9\times 8\times (7+6)+54\times 3\times 21.$
\item [] $4339=98\times (7+6\times 5+4)+321.$
\item [] $4340=98+7\times 6\times (5+4\times (3+21)).$
\item [] $4341=98\times 7\times 6+5\times (43+2)\times 1.$
\item [] $4342=98\times 7\times 6+5\times (43+2)+1.$
\item [] $\mathit{4343=9\times 8\times 7+6\times 5\times 4\times 32-1.}$
\item [] $4344=9\times 8\times 7+6\times 5\times 4\times 32\times 1.$
\item [] $4345=9\times 8\times 7+6\times 5\times 4\times 32+1.$
\item [] $4346=98\times 7\times 6+5\times (43+2+1).$
\item [] $4347=(9\times (8+7)+65+4+3)\times 21.$
\item [] $4348=((9\times 8+7)\times 6+5+4)\times 3^2+1.$
\item [] $4349=9\times (8+7+6)\times (5\times 4+3)+2\times 1.$
\item [] $4350=9\times 8\times 7+6\times (5\times 4\times 32+1).$
\item [] $4351=9+8+76\times (54+3)+2\times 1.$
\item [] $4352=98\times 7\times 6+5\times 43+21.$
\item [] $4353=9\times (8\times 7\times 6+5)+4\times 321.$
\item [] $4354=98+76\times (5\times (4+3)+21).$
\item [] $4355=9+8\times 7+65\times (4^3+2)\times 1.$
\item [] $4356=9+8+7+6+5+4321.$
\item [] $4357=(9+8+7+6\times 5+4\times 3)^2+1.$
\item [] $4358=9\times 87+(6+5)\times (4+321).$
\item [] $4359=9+87\times (6+5\times 4+3+21).$
\item [] $4360=(9\times 8+7+6\times 5)\times 4\times (3^2+1).$
\item [] $4361=98\times 7\times 6+5\times (4+3)^2\times 1.$
\item [] $4362=98\times 7\times 6+5\times (4+3)^2+1.$
\item [] $4363=(98\times 7+6)\times 5+43\times 21.$
\item [] $4364=987+(6+5+4)^3+2\times 1.$
\item [] $4365=987+(6+5+4)^3+2+1.$
\item [] $4366=9+(8\times 7+65)\times 4\times 3^2+1.$
\item [] $4367=9+8\times 76+5^4\times 3\times 2\times 1.$
\item [] $4368=9+8\times 76+5^4\times 3\times 2+1.$
\item [] $4369=(9+8\times 76)\times 5+4\times 321.$
\item [] $4370=9+8+76\times (54+3)+21.$
\item [] $4371=987+6\times (543+21).$
\item [] $4372=(9+8\times 7)\times 65+(4+3)\times 21.$
\item [] $4373=98\times 7\times 6+5+4\times 3\times 21.$
\item [] $4374=9\times 8\times 7+6\times 5\times 43\times (2+1).$
\item [] $4375=9+8+7+6\times 5+4321.$
\item [] $\mathit{4376=9\times 87\times 6-5+4-321.}$
\item [] $4377=(9\times 87+6)\times 5+432\times 1.$
\item [] $4378=(9\times 87+6)\times 5+432+1.$
\item [] $4379=9+8+7+65\times (4+3\times 21).$
\item [] $4380=9\times (8\times 7+65)\times 4+3+21.$
\item [] $4381=(98+7)\times 6+5^4\times 3\times 2+1.$
\item [] $4382=(9+8\times 76+5+4)\times (3\times 2+1).$
\item [] $4383=9\times (8+7\times 6+5+432\times 1).$
\item [] $4384=9+(8\times (7+6)+5^4)\times 3\times 2+1.$
\item [] $4385=9+8+7\times 6+5+4321.$
\item [] $4386=98\times 7\times 6+54\times (3+2\times 1).$
\item [] $4387=98\times 7\times 6+54\times (3+2)+1.$
\item [] $4388=98+(76+54)\times (32+1).$
\item [] $4389=(987+6\times 5)\times 4+321.$
\item [] $4390=(9+8)\times (7\times 6+5\times 43)+21.$
\item [] $4391=(9+8+7\times 6)\times 5+4^{(3+2+1)}.$
\item [] $4392=(98+76+5+4)\times (3+21).$
\item [] $4393=9+8\times 76\times 5+4^3\times 21.$
\item [] $4394=9\times 8\times (7+6+5+43)+2\times 1.$
\item [] $4395=987+6+54\times 3\times 21.$
\item [] $4396=98+7+65\times (4^3+2)+1.$
\item [] $4397=9+8\times 7+6+5+4321.$
\item [] $4398=9+876\times 5+4+3+2\times 1.$
\item [] $4399=9+876\times 5+4+3+2+1.$
\item [] $4400=9+876\times 5+4+3\times 2+1.$
\item[]$\mbox{Increasing order}$
\item [] $4401=1\times 2\times 3\times 4+56\times 78+9.$
\item [] $4402=1\times (2+3\times 4\times 5)\times (6+7\times 8+9).$
\item [] $4403=(1^2\times 3+4+5\times 6)\times 7\times (8+9).$
\item [] $4404=1\times (23+4)+56\times 78+9.$
\item [] $4405=1\times 2+(3+4+5\times 6)\times 7\times (8+9).$
\item [] $4406=1+(2+3)^4+5\times (6+78)\times 9.$
\item [] $4407=12\times 3\times 4\times 5\times 6+78+9.$
\item [] $4408=(1+2)^3+4+56\times 78+9.$
\item [] $4409=1\times 2^3\times 4+56\times 78+9.$
\item [] $4410=1+2^3\times 4+56\times 78+9.$
\item [] $4411=1^2\times 34+56\times 78+9.$
\item [] $4412=1^2+34+56\times 78+9.$
\item [] $4413=1\times 2+34+56\times 78+9.$
\item [] $4414=1+2+34+56\times 78+9.$
\item [] $4415=1\times 2^{(3\times 4)}+5\times (6+7\times 8)+9.$
\item [] $4416=12\times 3\times 4\times 5\times 6+7+89.$
\item [] $4417=12\times 3+4+56\times 78+9.$
\item [] $4418=1+234+(5+6\times 7)\times 89.$
\item [] $4419=(1+(2+3^4)\times 5+67+8)\times 9.$
\item [] $4420=1^{23}\times 4\times 5\times (6+7)\times (8+9).$
\item [] $4421=1+(2+3\times 4\times 5+6)\times (7\times 8+9).$
\item [] $4422=1\times 2+(3\times 4+56)\times (7\times 8+9).$
\item [] $4423=1+2+34\times (5+6+7\times (8+9)).$
\item [] $4424=((1+2+3+4)\times 5+6)\times (7+8\times 9).$
\item [] $4425=1^{234}+56\times (7+8\times 9).$
\item [] $4426=1\times 2\times 3+4\times 5\times (6+7)\times (8+9).$
\item [] $4427=(1+2)\times (3+4)\times 5\times 6\times 7+8+9.$
\item [] $4428=(12+3+4+5+6\times 78)\times 9.$
\item [] $4429=1^{23}+4+56\times (7+8\times 9).$
\item [] $\mathit{4430=-1+2\times 34\times 56+7\times 89.}$
\item [] $4431=(1\times 2\times 34)\times 56+7\times 89.$
\item [] $4432=1^2+3+4+56\times (7+8\times 9).$
\item [] $4433=1\times (2+3+4)+56\times (7+8\times 9).$
\item [] $4434=1+(2\times 34+5+6)\times 7\times 8+9.$
\item [] $4435=12+3+4\times (5\times (6+7)\times (8+9)).$
\item [] $4436=1^2\times (3\times 4+56\times (7+8\times 9)).$
\item [] $4437=1\times 2\times 34\times 5\times (6+7)+8+9.$
\item [] $4438=1+2\times 34\times 5\times (6+7)+8+9.$
\item [] $4439=1+2+3\times 4+56\times (7+8\times 9).$
\item [] $4440=(1^2+3)\times 4+56\times (7+8\times 9).$
\item [] $4441=(1^2+3)\times 4^5+6\times 7\times 8+9.$
\item [] $4442=(1^2+3\times 4)\times (5+6\times 7\times 8)+9.$
\item [] $4443=1+2+3+(45+6)\times (78+9).$
\item [] $4444=1+23+4\times 5\times (6+7)\times (8+9).$
\item [] $4445=1\times 2\times 34+56\times 78+9.$
\item [] $4446=1+2\times 34+56\times 78+9.$
\item [] $4447=1^2+(34+5)\times (6\times 7+8\times 9).$
\item [] $4448=12+3\times 4+56\times (7+8\times 9).$
\item [] $4449=1+2\times 3\times 4+56\times (7+8\times 9).$
\item [] $4450=12+(3+4)\times (5+6+7\times 89).$
\item [] $4451=1\times 234+5+6\times 78\times 9.$
\item [] $4452=1+234+5+6\times 78\times 9.$
\item [] $4453=1+(2^3+45)\times (67+8+9).$
\item [] $4454=1\times 2\times (3\times 4+5)\times (6\times 7+89).$
\item [] $4455=(1+23+456+7+8)\times 9.$
\item [] $4456=1+(2\times 3+45+6)\times 78+9.$
\item [] $4457=12^3+4+5\times (67\times 8+9).$
\item [] $4458=12\times 3\times 4\times 5+6\times 7\times 89.$
\item [] $4459=1^2+34+56\times (7+8\times 9).$
\item [] $4460=1\times 2+3^4+56\times 78+9.$
\item [] $4461=1+2+3^4+56\times 78+9.$
\item [] $4462=1+(2^3+45)\times (6+78)+9.$
\item [] $4463=1\times 2^{(3+4)}\times 5\times 6+7\times 89.$
\item [] $4464=12\times 3+4+56\times (7+8\times 9).$
\item [] $4465=1^2+3^4\times (5+6\times 7+8)+9.$
\item [] $4466=1\times 2+(3+45)\times (6+78+9).$
\item [] $4467=123\times 4+5\times (6+789).$
\item [] $4468=(1^2+3)\times (4^5+6+78+9).$
\item [] $4469=1\times 23\times 4+56\times 78+9.$
\item [] $4470=12+3^4+56\times 78+9.$
\item[]$\mbox{Decreasing order}$
\item [] $4401=9+8+7+6\times (5+4)^3+2+1.$
\item [] $4402=9+876\times 5+4+3^2\times 1.$
\item [] $4403=9+876\times 5+4+3^2+1.$
\item [] $4404=9+876\times 5+4\times 3+2+1.$
\item [] $4405=98\times 7\times 6+(5+4)\times 32+1.$
\item [] $4406=9\times 8+76\times (54+3)+2\times 1.$
\item [] $4407=9\times 8+76\times (54+3)+2+1.$
\item [] $4408=9\times 8\times (7\times 6+5)+4^(3+2)\times 1.$
\item [] $4409=9+876\times 5+4\times (3+2)\times 1.$
\item [] $4410=9+8+7+65+4321.$
\item [] $4411=9\times 8+7+6+5+4321.$
\item [] $4412=9+8\times 7+(65+4)\times 3\times 21.$
\item [] $4413=9+876\times 5+4\times 3\times 2^1.$
\item [] $4414=9+876\times 5+4\times 3\times 2+1.$
\item [] $4415=9+8+7\times (6+5)+4321.$
\item [] $4416=9+8\times 7+6\times 5+4321.$
\item [] $4417=9+876\times 5+4+3+21.$
\item [] $4418=(9+8)\times 76+5^4\times (3+2)+1.$
\item [] $4419=9+8+76+5+4321.$
\item [] $4420=9\times 8+7\times (65+4)\times 3^2+1.$
\item [] $4421=(9\times 8+7+6)\times (5\times 4+32)+1.$
\item [] $4422=9+876\times 5+4\times 3+21.$
\item [] $\mathit{4423=98-7+6+5+4321.}$
\item [] $4424=(9\times 8+7)\times (6+5+43+2\times 1).$
\item [] $4425=9+876\times 5+4+32\times 1.$
\item [] $4426=9+876\times 5+4+32+1.$
\item [] $4427=9+8+7\times 6\times 5\times (4+3)\times (2+1).$
\item [] $4428=9\times 8\times 7+654\times 3\times 2\times 1.$
\item [] $4429=9\times 8\times 7+654\times 3\times 2+1.$
\item [] $4430=9\times 8+7+6\times 5+4321.$
\item [] $4431=(98+7)\times (6\times 5+4\times 3)+21.$
\item [] $4432=98+76\times (54+3)+2\times 1.$
\item [] $4433=98+76\times (54+3)+2+1.$
\item [] $4434=9+876\times 5+43+2^1.$
\item [] $4435=9+876\times 5+43+2+1.$
\item [] $4436=98\times 7\times 6+5\times (43+21).$
\item [] $4437=98+7+6+5+4321.$
\item [] $4438=98\times 7\times 6+5\times 4^3+2\times 1.$
\item [] $4439=98\times 7\times 6+5\times 4^3+2+1.$
\item [] $4440=98\times 7\times 6+54\times 3\times 2\times 1.$
\item [] $4441=98\times 7\times 6+54\times 3\times 2+1.$
\item [] $4442=9+8+76\times 54+321.$
\item [] $4443=98\times 7+6+5^4\times 3\times 2+1.$
\item [] $4444=(9+8)\times (7+6)\times 5\times 4+3+21.$
\item [] $4445=98\times 7+6\times 5^4+3^2\times 1.$
\item [] $4446=98\times 7\times 6+5+4+321.$
\item [] $4447=9+87+6\times 5+4321.$
\item [] $4448=(98+7+6\times 5+4)\times 32\times 1.$
\item [] $4449=(9\times 8+7+6+54)\times 32+1.$
\item [] $4450=(9+876)\times 5+4\times 3\times 2+1.$
\item [] $4451=98+76\times (54+3)+21.$
\item [] $4452=98+7+(65+4)\times 3\times 21.$
\item [] $4453=9+876\times 5+43+21.$
\item [] $4454=(9+8)\times (7\times (6\times 5+4)+3+21).$
\item [] $4455=9+876\times 5+4^3+2\times 1.$
\item [] $4456=98+7+6\times 5+4321.$
\item [] $4457=98\times 7\times 6+5\times 4+321.$
\item [] $4458=9\times 8+(7+6)\times 5+4321.$
\item [] $4459=9+8\times (7+6+543)+2\times 1.$
\item [] $4460=98\times 7+6\times 5^4+3+21.$
\item [] $4461=(9+876)\times 5+4+32\times 1.$
\item [] $4462=(9+876)\times 5+4+32+1.$
\item [] $4463=98\times 7+6\times 5^4+3^{(2+1)}.$
\item [] $4464=(9+8)\times 7\times 6+5^4\times 3\times 2\times 1.$
\item [] $4465=9\times 8+7+65+4321.$
\item [] $4466=98+7\times 6+5+4321.$
\item [] $4467=9\times (8+7)+6+5+4321.$
\item [] $4468=98\times 7+6\times 5^4+32\times 1.$
\item [] $4469=98\times 7+6\times 5^4+32+1.$
\item [] $4470=(9+8)\times 7+6\times 5+4321.$
\item[]$\mbox{Increasing order}$
\item [] $4471=(1\times 2+3\times (4\times 5+67))\times (8+9).$
\item [] $4472=1+(2+3\times (4\times 5+67))\times (8+9).$
\item [] $4473=(1+23)\times 4+56\times 78+9.$
\item [] $4474=1^2+(3+4)\times (567+8\times 9).$
\item [] $4475=1\times 2+(3+4)\times (567+8\times 9).$
\item [] $4476=12+(3+45)\times (6+78+9).$
\item [] $4477=(1+2+34)\times (56+7\times 8+9).$
\item [] $4478=1\times 23+45\times (6\times (7+8)+9).$
\item [] $4479=(1+2)\times 34+56\times 78+9.$
\item [] $4480=(1^2+3+4)\times (56+7\times 8\times 9).$
\item [] $4481=(1\times 23+4\times 5)\times (6+7)\times 8+9.$
\item [] $4482=1\times 2\times 3\times 45+6\times 78\times 9.$
\item [] $4483=1+2\times 3\times 45+6\times 78\times 9.$
\item [] $4484=(12+3)\times 4+56\times (7+8\times 9).$
\item [] $4485=12\times 345+6\times 7\times 8+9.$
\item [] $4486=1+23\times (4+56+(7+8)\times 9).$
\item [] $4487=1+2^{(3+4+5)}+6\times (7\times 8+9).$
\item [] $4488=12\times (3+4\times 5)+6\times 78\times 9.$
\item [] $4489=(12+3\times 4+56)\times 7\times 8+9.$
\item [] $4490=1+(2\times 3\times 4+56)\times 7\times 8+9.$
\item [] $4491=12\times 3+45\times (6\times (7+8)+9).$
\item [] $4492=1\times 2\times 34+56\times (7+8\times 9).$
\item [] $4493=12\times (3\times 4\times 5\times 6+7)+89.$
\item [] $4494=1^2\times 3\times (4^5+6\times (7+8\times 9)).$
\item [] $4495=1^2+3\times (4^5+6\times (7+8\times 9)).$
\item [] $4496=12^3+4\times (5+678+9).$
\item [] $4497=1+2+3\times (4^5+6\times (7+8\times 9)).$
\item [] $4498=(1+2+3^4+7\times (5\times 6))\times (8+9).$
\item [] $4499=(1+2)\times (3+4)\times 5\times 6\times 7+89.$
\item [] $4500=12\times (345+6+7+8+9).$
\item [] $4501=1+(23+4+5+6\times 78)\times 9.$
\item [] $4502=(1+2\times 34)\times 5\times (6+7)+8+9.$
\item [] $4503=1^2+345\times (6+7)+8+9.$
\item [] $4504=123+4+56\times 78+9.$
\item [] $4505=1+2+345\times (6+7)+8+9.$
\item [] $4506=1+2^{(3+4)}+56\times 78+9.$
\item [] $4507=1\times 2+3^4+56\times (7+8\times 9).$
\item [] $4508=12^3+4\times 5\times (67+8\times 9).$
\item [] $4509=1\times 2\times 34\times 5\times (6+7)+89.$
\item [] $4510=1\times 2\times 34\times 56+78\times 9.$
\item [] $4511=1+2\times 34\times 56+78\times 9.$
\item [] $4512=12\times 34\times (5+6)+7+8+9.$
\item [] $4513=1+2^3\times (4+56+7\times 8\times 9).$
\item [] $4514=12+345\times (6+7)+8+9.$
\item [] $4515=(1^2+345)\times (6+7)+8+9.$
\item [] $4516=1\times 23\times 4+56\times (7+8\times 9).$
\item [] $4517=12\times (34+5\times 67)+89.$
\item [] $4518=((1+2)\times 34+56\times 7+8)\times 9.$
\item [] $4519=(1^2+3^4)\times (5+6\times 7+8)+9.$
\item [] $4520=(1\times 2^3)^4+5\times 67+89.$
\item [] $4521=12\times 3\times 4+56\times 78+9.$
\item [] $4522=12^3+4+5\times (6+7\times 8)\times 9.$
\item [] $4523=1\times 2+3\times 4\times (5+6\times 7)\times 8+9.$
\item [] $4524=(1^{23}+45+6)\times (78+9).$
\item [] $4525=1+(2\times 3\times 4+5)\times (67+89).$
\item [] $4526=(1+2)\times 34+56\times (7+8\times 9).$
\item [] $4527=12^3+45\times (6+7\times 8)+9.$
\item [] $4528=1\times (2+345)\times (6+7)+8+9.$
\item [] $4529=1+(2+345)\times (6+7)+8+9.$
\item [] $4530=12\times 345+6\times (7\times 8+9).$
\item [] $4531=(1+(2+3)^4+5+6)\times 7+8\times 9.$
\item [] $4532=(1+2)^3\times 4+56\times (7+8\times 9).$
\item [] $4533=(1+23+4+5\times 6)\times 78+9.$
\item [] $4534=1+2\times (3+4\times 5+6)\times 78+9.$
\item [] $4535=1\times (2\times 34+5)\times (6+7\times 8)+9.$
\item [] $4536=12\times 3\times (4+5)+6\times 78\times 9.$
\item [] $4537=1+(2+34)\times (5\times 6+7+89).$
\item [] $4538=1+(2+3\times 4\times (5+6\times 7))\times 8+9.$
\item [] $4539=1\times 2\times 3^4+56\times 78+9.$
\item [] $4540=1+2\times 3^4+56\times 78+9.$
\item[]$\mbox{Decreasing order}$
\item [] $4471=(9+8+7+6)\times 5+4321.$
\item [] $4472=98\times 7+(6+5^4)\times 3\times 2\times 1.$
\item [] $4473=(98+7+65+43)\times 21.$
\item [] $4474=9\times 8+76+5+4321.$
\item [] $4475=9+876\times 5+43\times 2\times 1.$
\item [] $4476=9+876\times 5+43\times 2+1.$
\item [] $4477=9+87+6+5^4\times (3\times 2+1).$
\item [] $4478=9+8\times (7+6+543)+21.$
\item [] $4479=98+7+6\times ((5+4)\times 3)^2\times 1.$
\item [] $4480=(9+8+7+6+5)\times 4\times 32\times 1.$
\item [] $4481=(9+8+7+6+5)\times 4\times 32+1.$
\item [] $4482=9+87+65+4321.$
\item [] $4483=9\times (8+7)\times 6\times 5+432+1.$
\item [] $4484=98+(7+6)\times 5+4321.$
\item [] $4485=9+876\times 5+4\times (3+21).$
\item [] $4486=9\times (8+7)+6\times 5+4321.$
\item [] $\mathit{4487=-98+7\times 654+3\times 2+1.}$
\item [] $4488=(9\times 8)\times (7\times 6+5\times 4)+3+21.$
\item [] $4489=(9+876)\times 5+43+21.$
\item [] $4490=9\times (8+7)+65\times (4+3\times 21).$
\item [] $4491=98+7+65+4321.$
\item [] $4492=(9+876)\times 5+4+3\times 21.$
\item [] $4493=9+8+76\times 5+4^{(3+2+1)}.$
\item [] $4494=98\times 7\times 6+54\times (3\times 2+1).$
\item [] $4495=9+(8+7)\times (6+5)+4321.$
\item [] $4496=98+7\times (6+5)+4321.$
\item [] $4497=9\times 8+76\times 54+321.$
\item [] $4498=9\times 8+7\times (6+5^4)+3^2\times 1.$
\item [] $4499=98\times 7+6\times 5^4+3\times 21.$
\item [] $4500=98+76+5+4321.$
\item [] $4501=(9+87+654)\times 3\times 2+1.$
\item [] $4502=9+8+76\times (54+3+2)+1.$
\item [] $4503=(9\times 8+7)\times (6+5+43+2+1).$
\item [] $4504=9+8+7\times (6+5^4+3^2+1).$
\item [] $4505=(9+8)\times 7+65+4321.$
\item [] $4506=9+87+6\times 5\times (4+3)\times 21.$
\item [] $4507=9\times (8+76)+5^4\times 3\times 2+1.$
\item [] $4508=98\times (7+6+5+4+3+21).$
\item [] $4509=9+(8+7\times 6)\times (5+4^3+21).$
\item [] $4510=9\times (87+6)\times 5+4+321.$
\item [] $4511=(9+876)\times 5+43\times 2\times 1.$
\item [] $4512=(9+876)\times 5+43\times 2+1.$
\item [] $4513=9+8\times (76+54\times 3^2+1).$
\item [] $4514=9+(876+5^4)\times 3+2\times 1.$
\item [] $4515=9\times (8+7+6)+5+4321.$
\item [] $4516=(9\times 87+6\times 5\times 4)\times (3+2)+1.$
\item [] $4517=9+876\times 5+4^3\times 2\times 1.$
\item [] $4518=9+876\times 5+4\times 32+1.$
\item [] $4519=(9+8+7)\times 6+5^4\times (3\times 2+1).$
\item [] $4520=98+7\times (6+5^4)+3+2\times 1.$
\item [] $4521=9\times (8+7)+65+4321.$
\item [] $4522=9\times 8+7\times (6+5^4)+32+1.$
\item [] $4523=98+76\times 54+321.$
\item [] $4524=98+7\times (6+5^4)+3^2\times 1.$
\item [] $4525=98+7\times (6+5^4)+3^2+1.$
\item [] $4526=98\times 7+6\times 5\times 4\times 32\times 1.$
\item [] $4527=98\times 7+6\times 5\times 4^3\times 2+1.$
\item [] $4528=9\times (87+6)\times 5+(4+3)^{(2+1)}.$
\item [] $4529=(987+65)\times 4+321.$
\item [] $4530=(98\times 7+65+4)\times 3\times 2\times 1.$
\item [] $4531=(98\times 7+65+4)\times 3\times 2+1.$
\item [] $4532=9+8+7\times (6\times 54+321).$
\item [] $4533=9\times (87+6\times 5)\times 4+321.$
\item [] $4534=9+876\times 5+(4\times 3)^2+1.$
\item [] $\mathit{4535=9\times 87\times 6-54\times 3-2+1.}$
\item [] $4536=9+876\times 5+(4+3)\times 21.$
\item [] $4537=(9+87+6\times 5)\times 4\times 3^2+1.$
\item [] $4538=9\times 87+6\times 5^4+3+2\times 1.$
\item [] $4539=9\times 87+6\times 5^4+3+2+1.$
\item [] $4540=9\times 87+6+5^4\times 3\times 2+1.$
\item[]$\mbox{Increasing order}$
\item [] $4541=1+(2+3)\times 4\times (5\times 6\times 7+8+9).$
\item [] $4542=123\times 4+5\times 6\times (7+8)\times 9.$
\item [] $4543=123+4\times 5\times (6+7)\times (8+9).$
\item [] $4544=1\times 2+3+(4+5+6\times 7)\times 89.$
\item [] $4545=1^{234}\times 567\times 8+9.$
\item [] $4546=1^{234}+567\times 8+9.$
\item [] $4547=1+2^{(3+4+5)}+(6\times 7+8)\times 9.$
\item [] $4548=1\times 2\times 3^4\times 5+6\times 7\times 89.$
\item [] $4549=1^{23}\times 4+567\times 8+9.$
\item [] $4550=1^{23}+4+567\times 8+9.$
\item [] $4551=123+4+56\times (7+8\times 9).$
\item [] $4552=1^2\times 3+4+567\times 8+9.$
\item [] $4553=1+2\times 34\times 5+6\times 78\times 9.$
\item [] $4554=1\times 2+3+4+567\times 8+9.$
\item [] $4555=1+2+3+4+567\times 8+9.$
\item [] $4556=1+2\times 3+4+567\times 8+9.$
\item [] $4557=1^2\times 345+6\times 78\times 9.$
\item [] $4558=1+2^3+4+567\times 8+9.$
\item [] $4559=1\times 2+345+6\times 78\times 9.$
\item [] $4560=1+2+345+6\times 78\times 9.$
\item [] $4561=1^2+3^4\times 56+7+8+9.$
\item [] $4562=1\times 2+3^4\times 56+7+8+9.$
\item [] $4563=1+2+3^4\times 56+7+8+9.$
\item [] $4564=12+3+4+567\times 8+9.$
\item [] $4565=1\times (2+3)\times 4+567\times 8+9.$
\item [] $4566=1+(2+3)\times 4+567\times 8+9.$
\item [] $4567=12\times 34\times (5+6)+7+8\times 9.$
\item [] $4568=12\times 3\times 4+56\times (7+8\times 9).$
\item [] $4569=12+345+6\times 78\times 9.$
\item [] $4570=1+2\times 3\times 4+567\times 8+9.$
\item [] $4571=1^2\times (3+4)\times (5\times 6+7\times 89).$
\item [] $4572=1\times 23+4+567\times 8+9.$
\item [] $4573=1+23+4+567\times 8+9.$
\item [] $4574=1+(23+45)\times 67+8+9.$
\item [] $4575=12\times 34\times (5+6)+78+9.$
\item [] $4576=1\times 2+345\times (6+7)+89.$
\item [] $4577=1\times 2^3\times 4+567\times 8+9.$
\item [] $4578=1+2^3\times 4+567\times 8+9.$
\item [] $4579=1^2\times 34+567\times 8+9.$
\item [] $4580=1^2+34+567\times 8+9.$
\item [] $4581=1\times 2+34+567\times 8+9.$
\item [] $4582=1+2+34+567\times 8+9.$
\item [] $4583=12+(3+4)\times (5\times 6+7\times 89).$
\item [] $4584=1^{23}\times 4567+8+9.$
\item [] $4585=1^{23}+4567+8+9.$
\item [] $4586=12+345\times (6+7)+89.$
\item [] $4587=1^2\times 3+4567+8+9.$
\item [] $4588=1^2+3+4567+8+9.$
\item [] $4589=1\times 2+3+4567+8+9.$
\item [] $4590=1+2+3+4567+8+9.$
\item [] $4591=1+2\times 3+4567+8+9.$
\item [] $4592=1\times 2^3+4567+8+9.$
\item [] $4593=1+2^3+4567+8+9.$
\item [] $4594=1^2+3+45\times (6+7+89).$
\item [] $4595=1\times 2+3+45\times (6+7+89).$
\item [] $4596=1+2+3+45\times (6+7+89).$
\item [] $4597=1+2\times 3+45\times (6+7+89).$
\item [] $4598=1+2\times 34\times 56+789.$
\item [] $4599=12+3+4567+8+9.$
\item [] $4600=(1\times 2+345)\times (6+7)+89.$
\item [] $4601=1^2\times 3^4\times 56+7\times 8+9.$
\item [] $4602=12+34\times (56+7+8\times 9).$
\item [] $4603=1\times 2+3^4\times 56+7\times 8+9.$
\item [] $4604=1+2+3^4\times 56+7\times 8+9.$
\item [] $4605=12+3+45\times (6+7+89).$
\item [] $4606=1\times 2^{(3+4+5)}+6+7\times 8\times 9.$
\item [] $4607=1\times 23+4567+8+9.$
\item [] $4608=1+23+4567+8+9.$
\item [] $4609=1234+5\times (67+8)\times 9.$
\item [] $4610=1\times 2+3\times 4\times (5\times (67+8)+9).$
\item[]$\mbox{Decreasing order}$
\item [] $4541=98\times 7\times 6+5\times (4^3+21).$
\item [] $4542=9\times 87+6\times 5^4+3^2\times 1.$
\item [] $4543=9\times 87+6\times 5^4+3^2+1.$
\item [] $4544=(9+8\times 7+6)\times (54+3^2+1).$
\item [] $4545=9+8+7+6+5\times 43\times 21.$
\item [] $4546=98\times 7\times 6+5\times 43\times 2\times 1.$
\item [] $4547=98\times 7\times 6+5\times 43\times 2+1.$
\item [] $4548=9+8+7\times 6\times 5+4321.$
\item [] $4549=98+76+5^4\times (3\times 2+1).$
\item [] $4550=(9+8\times 7)\times 65+4+321.$
\item [] $4551=(9+876+5^4)\times 3+21.$
\item [] $4552=9+8\times 7\times 65+43\times 21.$
\item [] $4553=98\times 7\times 6+5+432\times 1.$
\item [] $4554=98\times 7\times 6+5+432+1.$
\item [] $4555=9\times (8+7+6+5)+4321.$
\item [] $4556=98\times 7+6\times 5\times 43\times (2+1).$
\item [] $4557=9\times 87+6\times 5^4+3+21.$
\item [] $4558=((9\times 8\times 7+65)\times 4+3)\times 2\times 1.$
\item [] $4559=9\times 8+7\times (6+5^4+3^2+1).$
\item [] $4560=9\times (8+7)\times 6+5^4\times 3\times 2\times 1.$
\item [] $4561=9\times (8+7)\times 6+5^4\times 3\times 2+1.$
\item [] $4562=987+(6+5)\times (4+321).$
\item [] $4563=9\times 87+(6+54)\times 3\times 21.$
\item [] $4564=9\times 87+6\times (5^4+3+2)+1.$
\item [] $4565=9\times 87+6\times 5^4+32\times 1.$
\item [] $4566=9\times 87+6\times 5^4+32+1.$
\item [] $4567=9+8+7\times 65\times (4+3+2+1).$
\item [] $4568=9+8+7\times 65\times (4+3\times 2)+1.$
\item [] $4569=9+(8+76)\times 54+3+21.$
\item [] $4570=9\times 87+(6+5^4)\times 3\times 2+1.$
\item [] $4571=(9+8)\times (7+65\times 4)+32\times 1.$
\item [] $4572=9\times 87+6\times (5^4+3)+21.$
\item [] $4573=9+(8\times (76+5)+4)\times (3\times 2+1).$
\item [] $4574=9+8+7\times 6+5\times 43\times 21.$
\item [] $4575=9\times 8\times 7+6\times 5^4+321.$
\item [] $4576=(9+8+7\times 6\times 54+3)\times 2\times 1.$
\item [] $4577=9+(8+76)\times 54+32\times 1.$
\item [] $4578=9+(8+76)\times 54+32+1.$
\item [] $4579=9+8+76\times 5\times 4\times 3+2\times 1.$
\item [] $4580=9+8+76\times 5\times 4\times 3+2+1.$
\item [] $4581=9+876\times 5+4^3\times (2+1).$
\item [] $4582=98+76\times (54+3+2)\times 1.$
\item [] $4583=98+76\times (54+3+2)+1.$
\item [] $4584=(9\times 8\times 7+65\times 4)\times 3\times 2\times 1.$
\item [] $4585=(9\times 8\times 7+65\times 4)\times 3\times 2+1.$
\item [] $4586=9+8\times 7+6+5\times 43\times 21.$
\item [] $4587=9\times 8+7\times (6\times 54+321).$
\item [] $4588=9\times (87+6)+5^4\times 3\times 2+1.$
\item [] $4589=9+(8+7\times (6+5\times 4^3))\times 2\times 1.$
\item [] $4590=(98+7\times 6+5^4)\times 3\times 2\times 1.$
\item [] $4591=(98+7\times 6+5^4)\times 3\times 2+1.$
\item [] $4592=(9+8)\times (7+65\times 4+3)+2\times 1.$
\item [] $4593=9\times 87+6\times (5^4+3^2+1).$
\item [] $4594=9\times 8+7\times (6+5\times 4^3\times 2\times 1).$
\item [] $4595=9+8+7\times (6\times 54+3)\times 2\times 1.$
\item [] $4596=9\times 87+6\times 5^4+3\times 21.$
\item [] $4597=(9\times 8+7+6)\times 54+3\times 2+1.$
\item [] $4598=9+8+76\times 5\times 4\times 3+21.$
\item [] $4599=98\times 7\times 6+(5\times 4+3)\times 21.$
\item [] $4600=9\times 8+7+6+5\times 43\times 21.$
\item [] $4601=9+8+7\times 654+3+2+1.$
\item [] $4602=9+8+7\times 654+3\times 2+1.$
\item [] $4603=9\times 8+7\times 6\times 5+4321.$
\item [] $4604=9+8+7\times 654+3^2\times 1.$
\item [] $4605=9+8+7\times 654+3^2+1.$
\item [] $4606=98\times (7+6\times 5+4+3+2+1).$
\item [] $4607=9+8+(7\times 65+4)\times (3^2+1).$
\item [] $4608=9+8+76+5\times 43\times 21.$
\item [] $4609=(9+8+7+6\times 5\times 4)\times 32+1.$
\item [] $4610=98\times 7+654\times 3\times 2\times 1.$
\item[]$\mbox{Increasing order}$
\item [] $4611=1\times 234+56\times 78+9.$
\item [] $4612=1+234+56\times 78+9.$
\item [] $4613=1\times 2\times 34+567\times 8+9.$
\item [] $4614=1+2\times 34+567\times 8+9.$
\item [] $4615=1^2\times 3^4\times 56+7+8\times 9.$
\item [] $4616=1^2+3^4\times 56+7+8\times 9.$
\item [] $4617=12\times 345+6\times 78+9.$
\item [] $4618=1+2+3^4\times 56+7+8\times 9.$
\item [] $4619=1\times 2+3^4\times 5+6\times 78\times 9.$
\item [] $4620=12\times 3+4567+8+9.$
\item [] $4621=1^2+3+(45+6\times 78)\times 9.$
\item [] $4622=(1^2+3^4)\times 5+6\times 78\times 9.$
\item [] $4623=1^2\times 3^4\times 56+78+9.$
\item [] $4624=1^2+3^4\times 56+78+9.$
\item [] $4625=12\times 34+5+6\times 78\times 9.$
\item [] $4626=1+2+3^4\times 56+78+9.$
\item [] $4627=12+3^4\times 56+7+8\times 9.$
\item [] $4628=1\times 2+3^4+567\times 8+9.$
\item [] $4629=1+2+3^4+567\times 8+9.$
\item [] $4630=1+2^{(3+4)}\times 5\times 6+789.$
\item [] $4631=1\times 2\times (3+4\times 567)+89.$
\item [] $4632=1^2\times 3^4\times 56+7+89.$
\item [] $4633=1234+5\times 678+9.$
\item [] $4634=1\times 2+3^4\times 56+7+89.$
\item [] $4635=12+3^4\times 56+78+9.$
\item [] $4636=12+34\times (5+6\times 7+89).$
\item [] $4637=1\times 23\times 4+567\times 8+9.$
\item [] $4638=12+3^4+567\times 8+9.$
\item [] $4639=1^{23}\times 4567+8\times 9.$
\item [] $4640=1^{23}+4567+8\times 9.$
\item [] $4641=(1+23)\times 4+567\times 8+9.$
\item [] $4642=1^2\times 3+4567+8\times 9.$
\item [] $4643=1^2+3+4567+8\times 9.$
\item [] $4644=1\times 2+3+4567+8\times 9.$
\item [] $4645=1\times 2\times 3+4567+8\times 9.$
\item [] $4646=1+2\times 3+4567+8\times 9.$
\item [] $4647=1\times 2^3+4567+8\times 9.$
\item [] $4648=1+2^3+4567+8\times 9.$
\item [] $4649=1^{23}+4\times (5+(6+7)\times 89).$
\item [] $4650=12\times 345+6+7\times 8\times 9.$
\item [] $4651=1+(2+3+45)\times (6+78+9).$
\item [] $4652=(1^2+3)\times (4^5+67+8\times 9).$
\item [] $4653=(1+2\times 34)\times 56+789.$
\item [] $4654=12+3+4567+8\times 9.$
\item [] $4655=1^2\times 3^4\times 56+7\times (8+9).$
\item [] $4656=1^{23}\times 4567+89.$
\item [] $4657=1^{23}+4567+89.$
\item [] $4658=1\times 234+56\times (7+8\times 9).$
\item [] $4659=1^2\times 3+4567+89.$
\item [] $4660=1^2+3+4567+89.$
\item [] $4661=1\times 2+3+4567+89.$
\item [] $4662=1\times 23+4567+8\times 9.$
\item [] $4663=1+23+4567+8\times 9.$
\item [] $4664=1\times 2^3+4567+89.$
\item [] $4665=1+2^3+4567+89.$
\item [] $4666=(1+2)^3+4567+8\times 9.$
\item [] $4667=12+3^4\times 56+7\times (8+9).$
\item [] $4668=12\times (34+5\times (6+7\times 8+9)).$
\item [] $4669=1234+5\times (678+9).$
\item [] $4670=12+34\times (5\times (6+7)+8\times 9).$
\item [] $4671=12+3+4567+89.$
\item [] $4672=123+4+567\times 8+9.$
\item [] $4673=1+23\times 4\times 5+6\times 78\times 9.$
\item [] $4674=1+2+3^4\times 56+(7+8)\times 9.$
\item [] $4675=12\times 3+4567+8\times 9.$
\item [] $4676=1+(2+3)^4+5\times 6\times (7+8)\times 9.$
\item [] $4677=(12+3^4)\times 5+6\times 78\times 9.$
\item [] $4678=1+2^{(3\times 4)}+5+6\times (7+89).$
\item [] $4679=1\times 23+4567+89.$
\item [] $4680=1+23+4567+89.$
\item[]$\mbox{Decreasing order}$
\item [] $4611=98\times 7+654\times 3\times 2+1.$
\item [] $4612=9\times 8\times 7+6+5+4^{(3\times 2)}+1.$
\item [] $4613=98+7\times (6\times 54+321).$
\item [] $4614=9+(8+7)\times 6+5\times 43\times 21.$
\item [] $4615=98\times (7\times 6+5)+4+3+2\times 1.$
\item [] $4616=98\times (7\times 6+5)+4+3+2+1.$
\item [] $4617=9+87+6+5\times 43\times 21.$
\item [] $4618=9\times (87+6)\times 5+432+1.$
\item [] $4619=9+8+7\times 654+3+21.$
\item [] $4620=98+7\times (6+5\times 4\times 32\times 1).$
\item [] $4621=98\times (7\times 6+5)+4\times 3+2+1.$
\item [] $4622=(9+87+6+5)\times 43+21.$
\item [] $4623=9\times 87+6\times 5\times 4^3\times 2\times 1.$
\item [] $4624=9\times 87+6\times 5\times 4^3\times 2+1.$
\item [] $4625=9+8\times (7+6+543+21).$
\item [] $4626=98+7+6+5\times 43\times 21.$
\item [] $4627=9+8+7\times 654+32\times 1.$
\item [] $4628=9+8+7\times 654+32+1.$
\item [] $4629=98+7\times 6\times 5+4321.$
\item [] $4630=98\times (7\times 6+5)+4\times 3\times 2\times 1.$
\item [] $4631=98\times (7\times 6+5)+4\times 3\times 2+1.$
\item [] $4632=9+8+(765+4)\times 3\times 2+1.$
\item [] $4633=9\times 8+76\times (54+3\times 2)+1.$
\item [] $4634=98+(7+6+5)\times 4\times 3\times 21.$
\item [] $4635=9+876+5^4\times 3\times 2\times 1.$
\item [] $4636=9+876+5^4\times 3\times 2+1.$
\item [] $4637=9+8+7\times (654+3+2+1).$
\item [] $4638=(9\times 8+76+5^4)\times 3\times 2\times 1.$
\item [] $4639=98\times (7\times 6+5)+4\times 3+21.$
\item [] $4640=9+(8\times 7+6)\times 5+4321.$
\item [] $4641=9+876\times 5+4\times 3\times 21.$
\item [] $4642=98\times (7\times 6+5)+4+32\times 1.$
\item [] $4643=98\times (7\times 6+5)+4+32+1.$
\item [] $4644=9+8+7\times (654+3\times 2+1).$
\item [] $4645=(9+8+7+6\times 5)\times 43\times 2+1.$
\item [] $4646=98\times 7+6\times 5\times 4\times (32+1).$
\item [] $4647=(98+7)\times (6+5)\times 4+3^{(2+1)}.$
\item [] $4648=98+7\times 65\times (4+3+2+1).$
\item [] $4649=98+7\times 65\times (4+3\times 2)+1.$
\item [] $4650=9\times 8+7\times (6\times 54+3)\times 2\times 1.$
\item [] $4651=98\times (7\times 6+5)+43+2\times 1.$
\item [] $4652=98\times (7\times 6+5)+43+2+1.$
\item [] $4653=9\times 8+76\times 5\times 4\times 3+21.$
\item [] $4654=(((9+87)\times 6+5)\times 4+3)\times 2\times 1.$
\item [] $4655=98+7\times 6+5\times 43\times 21.$
\item [] $4656=9\times 8+7\times 654+3+2+1.$
\item [] $4657=9\times 8+7\times 654+3\times 2+1.$
\item [] $4658=9+8+7\times 654+3\times 21.$
\item [] $4659=9\times 8+7\times 654+3^2\times 1.$
\item [] $4660=98+76\times 5\times 4\times 3+2\times 1.$
\item [] $4661=98+76\times 5\times 4\times 3+2+1.$
\item [] $4662=98\times 7\times 6+543+2+1.$
\item [] $4663=9\times 8+76+5\times 43\times 21.$
\item [] $4664=9\times 8\times 7+65\times (43+21).$
\item [] $4665=9+8+7+(6+5\times 43)\times 21.$
\item [] $4666=9\times 8\times 7+65\times 4^3+2\times 1.$
\item [] $4667=9\times 8\times 7+65\times 4^3+2+1.$
\item [] $4668=9\times 8\times (7\times 6+5)+4\times 321.$
\item [] $4669=(9\times (8\times 7+6\times 5)+4)\times 3\times 2+1.$
\item [] $4670=98\times (7\times 6+5)+43+21.$
\item [] $4671=9+8\times 7\times 6+5+4321.$
\item [] $4672=9+(8+765+4)\times 3\times 2+1.$
\item [] $4673=98\times (7\times 6+5)+4+3\times 21.$
\item [] $4674=9\times 8+7\times 654+3+21.$
\item [] $4675=(9\times 8+7\times 6)\times (5+4+32)+1.$
\item [] $4676=(9+8\times 7+6)\times 5+4321.$
\item [] $4677=(9+876)\times 5+4\times 3\times 21.$
\item [] $4678=(9+87)\times 6+5+4^{(3\times 2)}+1.$
\item [] $4679=98+76\times 5\times 4\times 3+21.$
\item [] $4680=98\times 7\times 6+543+21.$
\item[]$\mbox{Increasing order}$
\item [] $4681=1+2^{(3\times 4)}+567+8+9.$
\item [] $4682=1\times 2+(3+4+5)\times 6\times (7\times 8+9).$
\item [] $4683=(1+2)^3+4567+89.$
\item [] $4684=1^{23}\times 4+5\times (6+7)\times 8\times 9.$
\item [] $4685=12\times 345+67\times 8+9.$
\item [] $4686=123\times (4+5\times 6)+7\times 8\times 9.$
\item [] $4687=1+2\times (3+4\times 567+8\times 9).$
\item [] $4688=(1^2+3^4)\times 56+7+89.$
\item [] $4689=12\times 3\times 4+567\times 8+9.$
\item [] $4690=1^{23}+(4+56)\times 78+9.$
\item [] $4691=1+2\times 3+4+5\times (6+7)\times 8\times 9.$
\item [] $4692=12\times 3+4567+89.$
\item [] $4693=1^2+3+45\times (6+7)\times 8+9.$
\item [] $4694=1\times 2+3+45\times (6+7)\times 8+9.$
\item [] $4695=(1+23+45)\times 67+8\times 9.$
\item [] $4696=1+2\times 3+45\times (6+7)\times 8+9.$
\item [] $4697=(123+456+7)\times 8+9.$
\item [] $4698=12\times 345+(6+7\times 8)\times 9.$
\item [] $4699=12+3+4+5\times (6+7)\times 8\times 9.$
\item [] $4700=1\times 2+(3+45+6)\times (78+9).$
\item [] $4701=1+2+(3+45+6)\times (78+9).$
\item [] $4702=1^2+3+(4+5)\times 6\times (78+9).$
\item [] $4703=1\times 2+3+(4+5)\times 6\times (78+9).$
\item [] $4704=12+3+45\times (6+7)\times 8+9.$
\item [] $4705=1+2^3\times (4+567+8+9).$
\item [] $4706=1\times 2^3+(4+5)\times 6\times (78+9).$
\item [] $4707=123+4567+8+9.$
\item [] $4708=1+2\times 3^4+567\times 8+9.$
\item [] $4709=123\times 4+5+6\times 78\times 9.$
\item [] $4710=12+(3+45+6)\times (78+9).$
\item [] $4711=(1^2+3^4)\times 56+7\times (8+9).$
\item [] $4712=(1+23+45)\times 67+89.$
\item [] $4713=123+45\times (6+7+89).$
\item [] $4714=1+(2+3^4)\times 56+7\times 8+9.$
\item [] $4715=12\times 3^4+5+6\times 7\times 89.$
\item [] $4716=12\times 345+6\times (7+89).$
\item [] $4717=(1^2+34+5+6+7)\times 89.$
\item [] $4718=(1+2\times 3)\times (45+6+7\times 89).$
\item [] $4719=(1+2^3\times 4)\times (56+78+9).$
\item [] $4720=1\times 2\times (3+4\times 567+89).$
\item [] $4721=1+2\times (3+4\times 567+89).$
\item [] $4722=1^2+3\times 4\times 56\times 7+8+9.$
\item [] $4723=1\times 2+3\times 4\times 56\times 7+8+9.$
\item [] $4724=1+2+3\times 4\times 56\times 7+8+9.$
\item [] $4725=12\times 3+45\times (6+7)\times 8+9.$
\item [] $4726=12+34+5\times (6+7)\times 8\times 9.$
\item [] $4727=(1+2)\times 3\times 456+7\times 89.$
\item [] $4728=12\times (3+4)\times 56+7+8+9.$
\item [] $4729=1+2^3\times (456+(7+8)\times 9).$
\item [] $4730=(12+345)\times (6+7)+89.$
\item [] $4731=(1+2+3)^4+5\times (678+9).$
\item [] $4732=1+(2\times 3)^4+5\times (678+9).$
\item [] $4733=12+3\times 4\times 56\times 7+8+9.$
\item [] $4734=1^{2345}\times 6\times 789.$
\item [] $4735=1^{2345}+6\times 789.$
\item [] $4736=1+(2+3^4)\times 56+78+9.$
\item [] $4737=(1\times 2\times 3\times 4+567)\times 8+9.$
\item [] $4738=1+(2\times 3\times 4+567)\times 8+9.$
\item [] $4739=1^{234}\times 5+6\times 789.$
\item [] $4740=1^{234}+5+6\times 789.$
\item [] $4741=1^{23}+(4+56)\times (7+8\times 9).$
\item [] $4742=(1+2+3\times 4\times 56)\times 7+8+9.$
\item [] $4743=1^{23}\times 4+5+6\times 789.$
\item [] $4744=1^{23}+4+5+6\times 789.$
\item [] $4745=1+(2+3^4)\times 56+7+89.$
\item [] $4746=1^2\times 3+4+5+6\times 789.$
\item [] $4747=1^2+3+4+5+6\times 789.$
\item [] $4748=1\times 2+3+4+5+6\times 789.$
\item [] $4749=1+2+3+4+5+6\times 789.$
\item [] $4750=1+2\times 3+4+5+6\times 789.$
\item[]$\mbox{Decreasing order}$
\item [] $4681=98+7\times 654+3+2\times 1.$
\item [] $4682=98+7\times 654+3+2+1.$
\item [] $4683=9\times 8+7\times 654+32+1.$
\item [] $4684=9+(8+(7+6)\times 5)\times 4^3+2+1.$
\item [] $4685=98+7\times 654+3^2\times 1.$
\item [] $4686=98+7\times 654+3^2+1.$
\item [] $4687=9\times 8+(765+4)\times 3\times 2+1.$
\item [] $4688=98+(7\times 65+4)\times (3^2+1).$
\item [] $4689=98+76+5\times 43\times 21.$
\item [] $4690=9\times 8\times (7+6)\times 5+4+3+2+1.$
\item [] $4691=98\times (7\times 6+5)+4^3+21.$
\item [] $4692=98\times (7\times 6+5)+43\times 2\times 1.$
\item [] $4693=98\times (7\times 6+5)+43\times 2+1.$
\item [] $4694=(9+8)\times 76+54\times 3\times 21.$
\item [] $4695=(9+876+54)\times (3+2)\times 1.$
\item [] $4696=(9+876+54)\times (3+2)+1.$
\item [] $4697=(9+8\times 76+54)\times (3\times 2+1).$
\item [] $4698=9\times (8+76+5+432+1).$
\item [] $4699=98+7\times (654+3)+2\times 1.$
\item [] $4700=98+7\times 654+3+21.$
\item [] $4701=9\times 8\times (7\times 6+5\times 4+3)+21.$
\item [] $4702=98\times (7\times 6+5)+4\times (3+21).$
\item [] $4703=98+7\times 654+3^{(2+1)}.$
\item [] $4704=(9+8+76+54)\times 32\times 1.$
\item [] $4705=(9+8+76+54)\times 32+1.$
\item [] $4706=9+8\times (7\times 6+5)+4321.$
\item [] $4707=9\times 87+654\times 3\times 2\times 1.$
\item [] $4708=98+7\times 654+32\times 1.$
\item [] $4709=98+7\times 654+32+1.$
\item [] $4710=9+8\times (7\times 6+543)+21.$
\item [] $4711=98+7\times (654+3+2\times 1).$
\item [] $4712=9\times 87\times 6+5+4+3+2\times 1.$
\item [] $4713=9\times 8+7\times 654+3\times 21.$
\item [] $4714=9+876\times 5+4+321.$
\item [] $\mathit{4715=9\times 87\times 6-5+4-3+21.}$
\item [] $4716=9\times 87\times 6+5+4+3^2\times 1.$
\item [] $4717=9\times 87\times 6+5+4+3^2+1.$
\item [] $4718=9+8+76\times 5+4321.$
\item [] $4719=(98+76)\times (5+4)\times 3+21.$
\item [] $4720=9\times 8+7+(6+5\times 43)\times 21.$
\item [] $4721=(9+8+7\times 65)\times (4+3\times 2)+1.$
\item [] $4722=9\times 87\times 6+(5+4+3)\times 2\times 1.$
\item [] $4723=9\times 87\times 6+5\times 4+3+2\times 1.$
\item [] $4724=9\times 87\times 6+5\times 4+3+2+1.$
\item [] $4725=9\times 87\times 6+5\times 4+3\times 2+1.$
\item [] $4726=9\times (8+7+6\times 5)+4321.$
\item [] $4727=98+(7+65)\times 4^3+21.$
\item [] $4728=9\times 87\times 6+5+4\times 3\times 2+1.$
\item [] $4729=9+8+76\times (5\times 4\times 3+2)\times 1.$
\item [] $4730=9+8+76\times (5\times 4\times 3+2)+1.$
\item [] $4731=9\times 87\times 6+5+4+3+21.$
\item [] $4732=9+876\times 5+(4+3)^{(2+1)}.$
\item [] $4733=9\times 87\times 6+(5+4\times 3)\times 2+1.$
\item [] $4734=98+76\times (54+3\times 2+1).$
\item [] $4735=98\times (7\times 6+5)+4\times 32+1.$
\item [] $4736=9\times 87\times 6+5+4\times 3+21.$
\item [] $4737=9+87+(6+5\times 43)\times 21.$
\item [] $4738=(98+7\times 6\times 54+3)\times 2\times 1.$
\item [] $4739=98+7\times 654+3\times 21.$
\item [] $4740=9\times 87\times 6+5+4+32+1.$
\item [] $4741=(9+8)\times (7+6)\times 5\times 4+321.$
\item [] $4742=9\times 87\times 6+5\times 4+3+21.$
\item [] $4743=987+6+5^4\times 3\times 2\times 1.$
\item [] $4744=987+6+5^4\times 3\times 2+1.$
\item [] $4745=9\times 87\times 6+(5\times 4+3)\times 2+1.$
\item [] $4746=987+6\times 5^4+3^2\times 1.$
\item [] $4747=987+6\times 5^4+3^2+1.$
\item [] $4748=98\times 7\times 6+5^4+3\times 2+1.$
\item [] $4749=9\times 87\times 6+5+43+2+1.$
\item [] $4750=987+(6+5^4\times 3)\times 2+1.$
\item[]$\mbox{Increasing order}$
\item [] $4751=1^2\times 3\times 4+5+6\times 789.$
\item [] $4752=1+2^3+4+5+6\times 789.$
\item [] $4753=1\times 2+3\times 4+5+6\times 789.$
\item [] $4754=1+2+3\times 4+5+6\times 789.$
\item [] $4755=1^{23}+4\times 5+6\times 789.$
\item [] $4756=1+2+(3+4)\times (56+7\times 89).$
\item [] $4757=1^2\times 3+4\times 5+6\times 789.$
\item [] $4758=12+3+4+5+6\times 789.$
\item [] $4759=1\times 2+3+4\times 5+6\times 789.$
\item [] $4760=12^3+45\times 67+8+9.$
\item [] $4761=12^3+4+5+6\times 7\times 8\times 9.$
\item [] $4762=123+4567+8\times 9.$
\item [] $4763=12+3\times 4+5+6\times 789.$
\item [] $4764=1+2\times 345\times 6+7\times 89.$
\item [] $4765=1^2\times 3+4^5+6\times 7\times 89.$
\item [] $4766=1\times 23+4+5+6\times 789.$
\item [] $4767=1+23+4+5+6\times 789.$
\item [] $4768=1+2+3+4^5+6\times 7\times 89.$
\item [] $4769=12+3+4\times 5+6\times 789.$
\item [] $4770=(1+2)^3+4+5+6\times 789.$
\item [] $4771=1+2^3+4^5+6\times 7\times 89.$
\item [] $4772=12^3+4\times 5+6\times 7\times 8\times 9.$
\item [] $4773=1^2\times 34+5+6\times 789.$
\item [] $4774=1^2+34+5+6\times 789.$
\item [] $4775=1\times 2+34+5+6\times 789.$
\item [] $4776=1+2+34+5+6\times 789.$
\item [] $4777=12+3+4^5+6\times 7\times 89.$
\item [] $4778=1+23+4\times 5+6\times 789.$
\item [] $4779=123+4567+89.$
\item [] $4780=1+234+567\times 8+9.$
\item [] $4781=12+(3+4)\times 5+6\times 789.$
\item [] $4782=1^2\times 3+45+6\times 789.$
\item [] $4783=1^2+3+45+6\times 789.$
\item [] $4784=1\times 2+3+45+6\times 789.$
\item [] $4785=12+34+5+6\times 789.$
\item [] $4786=1+2\times 3+45+6\times 789.$
\item [] $4787=1\times 2^3+45+6\times 789.$
\item [] $4788=12+3\times 4\times 56\times 7+8\times 9.$
\item [] $4789=(1+2)^3+4^5+6\times 7\times 89.$
\item [] $4790=12\times 3+4\times 5+6\times 789.$
\item [] $4791=12\times (3+4)\times 56+78+9.$
\item [] $4792=1\times 2+(3+4)\times (5+678)+9.$
\item [] $4793=1^2\times 3\times 4\times 56\times 7+89.$
\item [] $4794=12+3+45+6\times 789.$
\item [] $4795=1^2+3\times 4\times 5+6\times 789.$
\item [] $4796=1+2+3\times 4\times 56\times 7+89.$
\item [] $4797=12^3+45+6\times 7\times 8\times 9.$
\item [] $4798=12\times 3+4^5+6\times 7\times 89.$
\item [] $4799=(12+3)\times 4+5+6\times 789.$
\item [] $4800=12\times (3+4)\times 56+7+89.$
\item [] $4801=1+2^3\times 4\times 5\times (6+7+8+9).$
\item [] $4802=1\times 23+45+6\times 789.$
\item [] $4803=1+23+45+6\times 789.$
\item [] $4804=1234+5\times 6\times 7\times (8+9).$
\item [] $4805=123\times (4+5\times 6)+7\times 89.$
\item [] $4806=12+3\times 4\times 5+6\times 789.$
\item [] $4807=1\times 2\times 34+5+6\times 789.$
\item [] $4808=1+2\times 34+5+6\times 789.$
\item [] $4809=(1+2+3\times 4)\times 5+6\times 789.$
\item [] $4810=(1+2+3+4)\times (56\times 7+89).$
\item [] $4811=1\times 2+3\times 4\times (56\times 7+8)+9.$
\item [] $4812=123+45\times (6+7)\times 8+9.$
\item [] $4813=1+2\times (34+5)+6\times 789.$
\item [] $4814=(1+2+3\times 4\times 56)\times 7+89.$
\item [] $4815=12^3+45\times 67+8\times 9.$
\item [] $4816=1\times 2\times (34\times 56+7\times 8\times 9).$
\item [] $4817=1+2\times (34\times 56+7\times 8\times 9).$
\item [] $4818=(1+23)\times 45+6\times 7\times 89.$
\item [] $4819=(1^{23}+4+56)\times (7+8\times 9).$
\item [] $4820=1^2\times 3^4+5+6\times 789.$
\item[]$\mbox{Decreasing order}$
\item [] $4751=98\times 7\times 6+5^4+3^2+1.$
\item [] $4752=9\times 8\times 7\times 6+54\times 32\times 1.$
\item [] $4753=9\times 8\times 7\times 6+54\times 32+1.$
\item [] $4754=9\times 87\times 6+5\times (4+3)+21.$
\item [] $4755=9+8+7\times (6+5^4)+321.$
\item [] $4756=98\times 7\times 6+5\times 4^3\times 2\times 1.$
\item [] $4757=98\times 7+6\times 5^4+321.$
\item [] $4758=9\times 87\times 6+54+3+2+1.$
\item [] $4759=9\times 87\times 6+54+3\times 2+1.$
\item [] $4760=9\times 87\times 6+5\times 4\times 3+2\times 1.$
\item [] $4761=987+6\times 5^4+3+21.$
\item [] $4762=9\times 87\times 6+54+3^2+1.$
\item [] $4763=9+8+7\times (654+3+21).$
\item [] $4764=(9\times 8+7)\times (6+54)+3+21.$
\item [] $4765=98\times 7\times 6+5^4+3+21.$
\item [] $4766=9\times 8\times (7+6)\times 5+43\times 2\times 1.$
\item [] $4767=9\times 87\times 6+5+43+21.$
\item [] $4768=987+6\times (5^4+3+2)+1.$
\item [] $4769=9\times 87\times 6+5+4^3+2\times 1.$
\item [] $4770=987+6\times 5^4+32+1.$
\item [] $4771=9\times (87\times 6+5)+4+3+21.$
\item [] $4772=(9\times 8+7)\times (6+54)+32\times 1.$
\item [] $4773=9\times 8+76\times 5+4321.$
\item [] $4774=98\times 7\times 6+5^4+32+1.$
\item [] $4775=(98+7\times (6\times 54+3))\times 2+1.$
\item [] $4776=9\times 87\times 6+54+3+21.$
\item [] $4777=9+8+7+(6+5)\times 432+1.$
\item [] $4778=(9\times 8+76\times 5\times 4)\times 3+2\times 1.$
\item [] $4779=9\times 87\times 6+5\times 4\times 3+21.$
\item [] $4780=9\times 87\times 6+(5+4)\times 3^2+1.$
\item [] $4781=9\times 87\times 6+5\times 4+3\times 21.$
\item [] $4782=(98\times (7+6)+5\times 4^3)\times (2+1).$
\item [] $4783=9\times (87\times 6+5)+4\times (3^2+1).$
\item [] $4784=9\times 87\times 6+54+32\times 1.$
\item [] $4785=9\times 87\times 6+54+32+1.$
\item [] $4786=(9+8+76)\times 5+4321.$
\item [] $4787=9+8+7+(6+5)\times (432+1).$
\item [] $4788=9\times 87\times 6+5+4^3+21.$
\item [] $4789=9\times 87\times 6+5+43\times 2\times 1.$
\item [] $4790=9\times 87\times 6+5+43\times 2+1.$
\item [] $4791=9+8+7+6+(5+4^3)^2\times 1.$
\item [] $4792=9+8+7+6+(5+4^3)^2+1.$
\item [] $4793=9+8+7\times 65+4321.$
\item [] $4794=(98+76+5^4)\times 3\times 2\times 1.$
\item [] $4795=9+(87+6)\times 5+4321.$
\item [] $4796=(98+7\times 6)\times 5+4^{(3+2+1)}.$
\item [] $4797=987+6\times (5^4+3^2+1).$
\item [] $4798=9\times 87\times 6+5\times 4\times (3+2)\times 1.$
\item [] $4799=98+76\times 5+4321.$
\item [] $4800=987+6\times 5^4+3\times 21.$
\item [] $4801=(98+(7+6)\times 54)\times 3\times 2+1.$
\item [] $4802=(9+8+76+5)\times (4+3)^2\times 1.$
\item [] $4803=(9\times 8+7)\times (6+54)+3\times 21.$
\item [] $4804=98\times 7\times 6+5^4+3\times 21.$
\item [] $4805=98\times (7+6\times 5+4\times 3)+2+1.$
\item [] $4806=9+8+76\times (54+3^2)+1.$
\item [] $4807=9\times (87\times 6+5)+43+21.$
\item [] $4808=9\times 8\times (7+6)\times 5+4\times 32\times 1.$
\item [] $4809=9\times 8\times (7+6)\times 5+4\times 32+1.$
\item [] $4810=9\times 8+7\times (6+5^4)+321.$
\item [] $4811=98+76\times (5\times 4\times 3+2)+1.$
\item [] $4812=9\times 87\times 6+(54+3)\times 2\times 1.$
\item [] $4813=9\times 87\times 6+(54+3)\times 2+1.$
\item [] $4814=98\times 7+6\times (5^4+3\times 21).$
\item [] $4815=9\times 87\times 6+54+3\times 21.$
\item [] $4816=(9+87+6+5)\times (43+2)+1.$
\item [] $4817=9+8\times 7+(6+5)\times 432\times 1.$
\item [] $4818=9\times 87\times 6+5\times 4\times 3\times 2\times 1.$
\item [] $4819=9\times 87\times 6+5\times 4\times 3\times 2+1.$
\item [] $4820=(98+76\times 5+4)\times (3^2+1).$
\item[]$\mbox{Increasing order}$
\item [] $4821=1^2+3^4+5+6\times 789.$
\item [] $4822=1\times 2+3^4+5+6\times 789.$
\item [] $4823=1+2+3^4+5+6\times 789.$
\item [] $4824=1^{2345}\times 67\times 8\times 9.$
\item [] $4825=1^{2345}+67\times 8\times 9.$
\item [] $4826=1+2^{(3+4)}\times (5\times 6+7)+89.$
\item [] $4827=12\times 345+678+9.$
\item [] $4828=1\times (23+45)\times (6+7\times 8+9).$
\item [] $4829=1^{234}\times 5+67\times 8\times 9.$
\item [] $4830=1^{234}+5+67\times 8\times 9.$
\item [] $4831=1\times 23\times 4+5+6\times 789.$
\item [] $4832=12+3^4+5+6\times 789.$
\item [] $4833=1^{23}\times 4+5+67\times 8\times 9.$
\item [] $4834=1^{23}+4+5+67\times 8\times 9.$
\item [] $4835=12\times (345+6)+7\times 89.$
\item [] $4836=1^2\times 3+4+5+67\times 8\times 9.$
\item [] $4837=1^2+3+4+5+67\times 8\times 9.$
\item [] $4838=1\times 2+3+4+5+67\times 8\times 9.$
\item [] $4839=(12+3)\times 45\times 6+789.$
\item [] $4840=1+2\times 3+4+5+67\times 8\times 9.$
\item [] $4841=1\times 2^3+4+5+67\times 8\times 9.$
\item [] $4842=1\times 2\times 345\times 6+78\times 9.$
\item [] $4843=1+2\times 345\times 6+78\times 9.$
\item [] $4844=1+2+3\times 4+5+67\times 8\times 9.$
\item [] $4845=1^{23}+4\times 5+67\times 8\times 9.$
\item [] $4846=1\times 2+(3+4)\times (5+678+9).$
\item [] $4847=1^2\times 3+4\times 5+67\times 8\times 9.$
\item [] $4848=12\times 345+6+78\times 9.$
\item [] $4849=1\times 2+3+4\times 5+67\times 8\times 9.$
\item [] $4850=1+2+3+4\times 5+67\times 8\times 9.$
\item [] $4851=1+2\times 3+4\times 5+67\times 8\times 9.$
\item [] $4852=1\times 2^3+4\times 5+67\times 8\times 9.$
\item [] $4853=12+3\times 4+5+67\times 8\times 9.$
\item [] $4854=1\times 2\times 3\times 4\times 5+6\times 789.$
\item [] $4855=1+2\times 3\times 4\times 5+6\times 789.$
\item [] $4856=1\times 23+4+5+67\times 8\times 9.$
\item [] $4857=1+23+4+5+67\times 8\times 9.$
\item [] $4858=1\times 2\times (3\times 4+5)+67\times 8\times 9.$
\item [] $4859=12+3+4\times 5+67\times 8\times 9.$
\item [] $4860=12\times (3+4+56\times 7)+8\times 9.$
\item [] $4861=1\times 2^3\times 4+5+67\times 8\times 9.$
\item [] $4862=1+2^3\times 4+5+67\times 8\times 9.$
\item [] $4863=1^2\times 34+5+67\times 8\times 9.$
\item [] $4864=1^2+34+5+67\times 8\times 9.$
\item [] $4865=1\times 2+34+5+67\times 8\times 9.$
\item [] $4866=123+4+5+6\times 789.$
\item [] $4867=1\times 23+4\times 5+67\times 8\times 9.$
\item [] $4868=1+23+4\times 5+67\times 8\times 9.$
\item [] $4869=123\times 4+56\times 78+9.$
\item [] $4870=1^{23}+45+67\times 8\times 9.$
\item [] $4871=1\times 2+3\times 45+6\times 789.$
\item [] $4872=1+2+3\times 45+6\times 789.$
\item [] $4873=1^2+3+45+67\times 8\times 9.$
\item [] $4874=1\times 2+3+45+67\times 8\times 9.$
\item [] $4875=12+34+5+67\times 8\times 9.$
\item [] $4876=1+2\times 3+45+67\times 8\times 9.$
\item [] $4877=123+4\times 5+6\times 789.$
\item [] $4878=1+2^3+45+67\times 8\times 9.$
\item [] $4879=1+2\times 3\times (4+5)+67\times 8\times 9.$
\item [] $4880=12\times 3+4\times 5+67\times 8\times 9.$
\item [] $4881=12+3\times 45+6\times 789.$
\item [] $4882=1+2+3+4+56\times (78+9).$
\item [] $4883=12\times 3\times 4+5+6\times 789.$
\item [] $4884=12+3+45+67\times 8\times 9.$
\item [] $4885=123+4^5+6\times 7\times 89.$
\item [] $4886=1\times 2+3\times 4\times 5+67\times 8\times 9.$
\item [] $4887=1+2+3\times 4\times 5+67\times 8\times 9.$
\item [] $4888=(1^2+3)\times 4+56\times (78+9).$
\item [] $4889=(12+3)\times 4+5+67\times 8\times 9.$
\item [] $4890=12\times 3^4\times 5+6+7+8+9.$
\item[]$\mbox{Decreasing order}$
\item [] $4821=9+876\times 5+432\times 1.$
\item [] $4822=9+876\times 5+432+1.$
\item [] $4823=98\times (7+6\times 5+4\times 3)+21.$
\item [] $4824=9\times (87\times 6+5+4+3+2\times 1).$
\item [] $4825=9\times 8\times 7+6\times 5\times (4\times 3)^2+1.$
\item [] $4826=9+8+7\times (654+32+1).$
\item [] $4827=987+6\times 5\times 4^3\times 2\times 1.$
\item [] $4828=987+6\times 5\times 4^3\times 2+1.$
\item [] $4829=9\times 8\times (7+6+54)+3+2\times 1.$
\item [] $4830=9\times (87\times 6+5)+43\times 2+1.$
\item [] $4831=9\times 87\times 6+5+4^3\times 2\times 1.$
\item [] $4832=9\times 87\times 6+5+4^3\times 2+1.$
\item [] $4833=9\times 8\times (7+6+54)+3^2\times 1.$
\item [] $4834=98+(7+6\times 5)\times 4^3\times 2\times 1.$
\item [] $4835=98+(7+6\times 5)\times 4\times 32+1.$
\item [] $4836=9\times 8\times 7+6+5+4321.$
\item [] $4837=98\times 7\times 6+5\times (4\times 3)^2+1.$
\item [] $4838=9\times 87\times 6+5\times 4\times (3\times 2+1).$
\item [] $4839=9+8+7+(6+5+4)\times 321.$
\item [] $\mathit{4840=9\times 8\times 76-5^4-3\times 2-1.}$
\item [] $4841=9+8+(7+65)\times (4+3\times 21).$
\item [] $4842=(9\times 8\times 7+6\times 5+4)\times 3^2\times 1.$
\item [] $4843=9\times 87\times 6+(5+4+3)^2+1.$
\item [] $4844=98+7\times (654+3+21).$
\item [] $4845=9+(8\times 7+6)\times (54+3+21).$
\item [] $4846=98\times 7+65\times (43+21).$
\item [] $4847=(9\times 8+7+6)\times (54+3)+2\times 1.$
\item [] $4848=9\times 8+7\times 65+4321.$
\item [] $4849=98\times 7+65\times 4^3+2+1.$
\item [] $4850=9+8\times (7+6)\times 5+4321.$
\item [] $4851=(98+76+54+3)\times 21.$
\item [] $\mathit{4852=9\times 8\times 76-5^4+3\times 2-1.}$
\item [] $4853=9\times 87\times 6+5\times (4+3^{(2+1)}).$
\item [] $4854=9\times 87+6\times 5^4+321.$
\item [] $4855=9\times 8\times 7+6\times 5+4321.$
\item [] $4856=9\times 8\times (7+6+54)+32\times 1.$
\item [] $4857=9+87\times 6+5+4321.$
\item [] $4858=98+7+(6+5)\times 432+1.$
\item [] $4859=9\times 8\times 7+65\times (4+3\times 21).$
\item [] $4860=9+8\times 7\times 6+5\times 43\times 21.$
\item [] $4861=9\times 8+7\times 6\times (54+3)\times 2+1.$
\item [] $4862=9\times 87\times 6+54\times 3+2\times 1.$
\item [] $4863=9\times 87\times 6+54\times 3+2+1.$
\item [] $4864=9+87+6+(5+4^3)^2+1.$
\item [] $4865=98+7\times (654+3^{(2+1)}).$
\item [] $4866=9+(8+76)\times 54+321.$
\item [] $4867=98\times 7+65\times 4^3+21.$
\item [] $4868=9+8+7\times (6+(5+4)\times 3)\times 21.$
\item [] $4869=9\times 87\times 6+(54+3)\times (2+1).$
\item [] $4870=9+(8+7)\times 6\times (5+4)\times 3\times 2+1.$
\item [] $4871=9\times (87\times 6+5)+4\times 32\times 1.$
\item [] $4872=9\times (87\times 6+5)+4\times 32+1.$
\item [] $4873=98+7+6+(5+4^3)^2+1.$
\item [] $4874=98+7\times 65+4321.$
\item [] $4875=9+(8+7)\times 6\times 54+3\times 2\times 1.$
\item [] $4876=(98+7+6)\times 5+4321.$
\item [] $4877=(9+8\times 7)\times (6+5+4^3)+2\times 1.$
\item [] $4878=9\times 8\times 7\times 6+5+43^2\times 1.$
\item [] $4879=9\times 87\times 6+5\times 4\times 3^2+1.$
\item [] $4880=9+8\times 7+(6+5+4)\times 321.$
\item [] $4881=9\times 87\times 6+54\times 3+21.$
\item [] $4882=(98\times 7+6+5)\times (4+3)+2+1.$
\item [] $4883=9\times 87\times 6+5\times (4+32+1).$
\item [] $4884=9\times (8\times 7+6)+5+4321.$
\item [] $4885=(9+8\times 7+6+5)\times 4^3+21.$
\item [] $4886=98+76\times (54+3^2\times 1).$
\item [] $4887=98+76\times (54+3^2)+1.$
\item [] $4888=9+(87\times 6+5\times 4)\times 3^2+1.$
\item [] $4889=9+8\times (7\times (6+(5+4)\times 3^2)+1).$
\item [] $4890=9\times 87+6+5+4^{(3+2+1)}.$
\item[]$\mbox{Increasing order}$
\item [] $4891=12+3+4+56\times (78+9).$
\item [] $4892=1\times 23+45+67\times 8\times 9.$
\item [] $4893=1+23+45+67\times 8\times 9.$
\item [] $4894=1\times 2^3\times 4\times 5+6\times 789.$
\item [] $4895=1+2^3\times 4\times 5+6\times 789.$
\item [] $4896=12+3\times 4+56\times (78+9).$
\item [] $4897=1\times 2\times 34+5+67\times 8\times 9.$
\item [] $4898=1+2\times 34+5+67\times 8\times 9.$
\item [] $4899=1\times 23+4+56\times (78+9).$
\item [] $4900=1+23+4+56\times (78+9).$
\item [] $4901=1\times 2\times 3^4+5+6\times 789.$
\item [] $4902=123+45+6\times 789.$
\item [] $4903=1+2\times 345+6\times 78\times 9.$
\item [] $4904=1^2\times 34\times 5+6\times 789.$
\item [] $4905=12\times 3+45+67\times 8\times 9.$
\item [] $4906=1\times 2+34\times 5+6\times 789.$
\item [] $4907=1+2+34\times 5+6\times 789.$
\item [] $4908=1\times 234\times 5+6\times 7\times 89.$
\item [] $4909=1+234\times 5+6\times 7\times 89.$
\item [] $4910=1^2\times 3^4+5+67\times 8\times 9.$
\item [] $4911=1^2+3^4+5+67\times 8\times 9.$
\item [] $4912=1\times 2+3^4+5+67\times 8\times 9.$
\item [] $4913=1+2+3^4+5+67\times 8\times 9.$
\item [] $4914=1\times (2+34)\times 5+6\times 789.$
\item [] $4915=1+234+5\times (6+7)\times 8\times 9.$
\item [] $4916=12+34\times 5+6\times 789.$
\item [] $4917=12+34\times (5+6+7)\times 8+9.$
\item [] $4918=12+34+56\times (78+9).$
\item [] $4919=12\times 3^4\times 5+6\times 7+8+9.$
\item [] $4920=1\times 2\times (3+45)+67\times 8\times 9.$
\item [] $4921=1\times 23\times 4+5+67\times 8\times 9.$
\item [] $4922=1+23\times 4+5+67\times 8\times 9.$
\item [] $4923=(1^2\times 3+4+56)\times 78+9.$
\item [] $4924=1^2+(3+4+56)\times 78+9.$
\item [] $4925=1\times 2\times 3^4\times 5\times 6+7\times 8+9.$
\item [] $4926=1+2\times 3^4\times 5\times 6+7\times 8+9.$
\item [] $\mathit{4927=123-4\times 5+67\times 8\times 9.}$
\item [] $4928=12^3+4\times (5+6+789).$
\item [] $4929=1\times 2\times 345\times 6+789.$
\item [] $4930=1+2\times 345\times 6+789.$
\item [] $4931=12\times 3^4\times 5+6+7\times 8+9.$
\item [] $4932=12\times 3\times 4\times 5+6\times 78\times 9.$
\item [] $4933=1+(2+34)\times (5\times (6+7)+8\times 9).$
\item [] $4934=(12\times 3+4)\times 5+6\times 789.$
\item [] $4935=12\times 345+6+789.$
\item [] $4936=(1+2+3)^4+56\times (7\times 8+9).$
\item [] $4937=12^3+456\times 7+8+9.$
\item [] $4938=1\times 2\times 3\times (4+5\times 6+789).$
\item [] $4939=1\times 2\times 3^4\times 5\times 6+7+8\times 9.$
\item [] $4940=1+2\times 3^4\times 5\times 6+7+8\times 9.$
\item [] $4941=1\times 23\times (4+5)+6\times 789.$
\item [] $4942=1+23\times (4+5)+6\times 789.$
\item [] $4943=12\times 3\times 4\times 5\times 6+7\times 89.$
\item [] $4944=12\times 3^4\times 5+67+8+9.$
\item [] $4945=12\times 3^4\times 5+6+7+8\times 9.$
\item [] $4946=1+2^{(3\times 4)}+56\times (7+8)+9.$
\item [] $4947=1\times 2\times 3^4\times 5\times 6+78+9.$
\item [] $4948=1+2\times 3^4\times 5\times 6+78+9.$
\item [] $4949=(1+2\times 3\times 4)\times 5+67\times 8\times 9.$
\item [] $4950=12\times 345+6\times (7+8)\times 9.$
\item [] $4951=(1+2\times 345+6)\times 7+8\times 9.$
\item [] $4952=1\times 2\times (3+4\times (5+6)\times 7\times 8+9).$
\item [] $4953=12\times 34+567\times 8+9.$
\item [] $4954=1^2+3^4+56\times (78+9).$
\item [] $4955=1\times 2+3^4+56\times (78+9).$
\item [] $4956=123+4+5+67\times 8\times 9.$
\item [] $4957=1+2\times 3^4\times 5\times 6+7+89.$
\item [] $4958=1+2^{(3+4)}+5+67\times 8\times 9.$
\item [] $4959=1^2\times 3\times 45+67\times 8\times 9.$
\item [] $4960=1^2+3\times 45+67\times 8\times 9.$
\item[]$\mbox{Decreasing order}$
\item [] $4891=(9\times 8+7\times 6)\times 5+4321.$
\item [] $4892=9+8+(7+6)\times (54+321).$
\item [] $4893=9+(8+7)\times 6\times 54+3+21.$
\item [] $4894=9\times 8+7+(6+5+4)\times 321.$
\item [] $4895=9\times 87\times 6+5+4^3\times (2+1).$
\item [] $4896=9\times 8\times (7\times 6+5\times 4+3\times 2\times 1).$
\item [] $4897=((9+8+7)\times 6+5+4)\times 32+1.$
\item [] $4898=9+8\times 76\times 5+43^2\times 1.$
\item [] $4899=9+8\times 76\times 5+43^2+1.$
\item [] $4900=98+7\times (654+32)\times 1.$
\item [] $4901=9\times 8\times 7\times 6+5^4\times 3+2\times 1.$
\item [] $4902=9\times 8\times 7\times 6+5^4\times 3+2+1.$
\item [] $4903=(9\times 87+6\times 5+4)\times 3\times 2+1.$
\item [] $4904=(9+(8+7+6)\times 5)\times 43+2\times 1.$
\item [] $4905=(9+8\times 7)\times 6+5\times 43\times 21.$
\item [] $4906=9+8\times (7+65)+4321.$
\item [] $4907=98+7\times (654+32+1).$
\item [] $4908=9\times 8\times 7\times 6+(5^4+3)\times (2+1).$
\item [] $4909=9\times 8+76+(5+4^3)^2\times 1.$
\item [] $4910=98\times (7+6\times 5)+4\times 321.$
\item [] $4911=987+654\times 3\times 2\times 1.$
\item [] $4912=987+654\times 3\times 2+1.$
\item [] $4913=(9+8)\times (7+6\times 5+4\times 3\times 21).$
\item [] $4914=(98+76+5\times 4\times 3)\times 21.$
\item [] $4915=9\times 87\times 6+5\times 43+2\times 1.$
\item [] $4916=9+8+7\times 654+321.$
\item [] $4917=(9+87+6)\times (5+43)+21.$
\item [] $4918=9\times (8\times 7\times 6+5)+43^2\times 1.$
\item [] $4919=9\times (8\times 7\times 6+5)+43^2+1.$
\item [] $4920=98+7+(6+5+4)\times 321.$
\item [] $4921=(98+7\times 6)\times 5\times (4+3)+21.$
\item [] $4922=98+(7+65)\times (4+3\times 21).$
\item [] $4923=9\times 87\times 6+5\times (43+2\times 1).$
\item [] $4924=9\times 87\times 6+5\times (43+2)+1.$
\item [] $4925=(987+654)\times 3+2\times 1.$
\item [] $4926=(987+654)\times 3+2+1.$
\item [] $\mathit{4927=9\times 8\times 76-543-2\times 1.}$
\item [] $4928=9\times 87\times 6+5\times (43+2+1).$
\item [] $4929=987+6\times (5^4+32\times 1).$
\item [] $4930=987+6\times (5^4+32)+1.$
\item [] $4931=98\times (7\times 6+5)+4+321.$
\item [] $4932=9+(8+7)\times 6\times 54+3\times 21.$
\item [] $4933=9+8\times 7\times 65+4\times 321.$
\item [] $4934=9\times 87\times 6+5\times 43+21.$
\item [] $4935=(9+8\times 76)\times 5+43^2+1.$
\item [] $4936=98+76+(5+4^3)^2+1.$
\item [] $4937=9+8\times (76+54\times (3^2+1)).$
\item [] $4938=9+(87+6)\times (5\times 4+32+1).$
\item [] $4939=9\times (87+6)+5+4^{(3\times 2)}+1.$
\item [] $4940=(9+(8+7)\times 65+4)\times (3+2)\times 1.$
\item [] $4941=9\times (87+6\times 5+432\times 1).$
\item [] $4942=9\times (87+(6\times 5+432))+1.$
\item [] $4943=9+8\times 76+5+4321.$
\item [] $4944=(9+8+7)\times (6+5\times 4\times (3^2+1)).$
\item [] $4945=9\times 87+65\times 4^3+2\times 1.$
\item [] $4946=9+8\times 7\times (6+5)+4321.$
\item [] $4947=987+6\times 5\times 4\times (32+1).$
\item [] $4948=9+8+7\times (6+5)\times 4^3+2+1.$
\item [] $4949=(9+8)\times (76+5\times 43)+2\times 1.$
\item [] $4950=(9+8)\times ((7+65)\times 4+3)+2+1.$
\item [] $4951=(9+87)\times 6+5^4\times (3\times 2+1).$
\item [] $4952=(9+8\times 7\times (6+5)\times 4+3)\times 2\times 1.$
\item [] $4953=9\times (8+7)\times 6\times 5+43\times 21.$
\item [] $\mathit{4954=-9+8+7\times 6+(5+4\times 3)^{(2+1)}.}$
\item [] $4955=9\times 87\times 6+5+4\times 3\times 21.$
\item [] $4956=(98+7)\times 6+5+4321.$
\item [] $4957=9+8+76\times 5\times (4+3^2)\times 1.$
\item [] $4958=9+8+76\times 5\times (4+3^2)+1.$
\item [] $4959=98+765+4^{(3+2+1)}.$
\item [] $4960=(9+87\times 6+5\times 4)\times 3^2+1.$
\item[]$\mbox{Increasing order}$
\item [] $4961=1\times 2+3\times 45+67\times 8\times 9.$
\item [] $4962=12\times 3^4\times 5+6+7+89.$
\item [] $4963=1\times (2\times 34+5)\times 67+8\times 9.$
\item [] $4964=(12+34)\times 5+6\times 789.$
\item [] $4965=12+3^4+56\times (78+9).$
\item [] $4966=1\times 2+34\times (5+6+(7+8)\times 9).$
\item [] $4967=123+4\times 5+67\times 8\times 9.$
\item [] $4968=12\times 3\times (45+6+78+9).$
\item [] $4969=1+2\times 34\times (5+67)+8\times 9.$
\item [] $4970=1\times 2\times (3+4)\times 5\times (6+7\times 8+9).$
\item [] $4971=12+3\times 45+67\times 8\times 9.$
\item [] $4972=1^2+3\times (4\times 56\times 7+89).$
\item [] $4973=1\times 234+5+6\times 789.$
\item [] $4974=1+234+5+6\times 789.$
\item [] $4975=(1+2\times 34+5)\times 67+8+9.$
\item [] $4976=12+34\times (5+6+(7+8)\times 9).$
\item [] $4977=1234+5+6\times 7\times 89.$
\item [] $4978=1^2+(3+4+56)\times (7+8\times 9).$
\item [] $4979=1\times 2\times 3^4\times 5\times 6+7\times (8+9).$
\item [] $4980=1\times (2\times 34+5)\times 67+89.$
\item [] $4981=1+(2\times 34+5)\times 67+89.$
\item [] $4982=(1+2)\times 3^4+5+6\times 789.$
\item [] $4983=12+3\times (4\times 56\times 7+89).$
\item [] $4984=1\times 2^3\times 4\times 5+67\times 8\times 9.$
\item [] $4985=1+2^3\times 4\times 5+67\times 8\times 9.$
\item [] $4986=1\times 2\times (345+6)\times 7+8\times 9.$
\item [] $4987=1+2\times (345+6)\times 7+8\times 9.$
\item [] $4988=12^3+4\times (5+6\times (7+8)\times 9).$
\item [] $4989=12+(3+4+56)\times (7+8\times 9).$
\item [] $4990=1\times 2\times (3^4\times 5\times 6+7\times 8+9).$
\item [] $4991=12\times 3^4\times 5+6\times 7+89.$
\item [] $4992=123+45+67\times 8\times 9.$
\item [] $4993=1+(23+4+5)\times (67+89).$
\item [] $4994=1^2\times 34\times 5+67\times 8\times 9.$
\item [] $4995=1^2+34\times 5+67\times 8\times 9.$
\item [] $4996=1\times 2+34\times 5+67\times 8\times 9.$
\item [] $4997=1+2+34\times 5+67\times 8\times 9.$
\item [] $4998=(1^2+3+45)\times (6+7+89).$
\item [] $4999=12\times 3^4\times 5+67+8\times 9.$
\item [] $5000=1\times 2+(3\times 4+5\times 6)\times 7\times (8+9).$
\item [] $5001=12\times (345+6)+789.$
\item [] $5002=1^2+(34+5\times 6)\times 78+9.$
\item [] $5003=1\times 2\times (345+6)\times 7+89.$
\item [] $5004=1\times 2\times 3\times 45+6\times 789.$
\item [] $5005=1+2\times 3\times 45+6\times 789.$
\item [] $5006=12+34\times 5+67\times 8\times 9.$
\item [] $5007=1^2\times 3+(4\times 5+67\times 8)\times 9.$
\item [] $5008=1^2+3+(4\times 5+67\times 8)\times 9.$
\item [] $5009=12^3+456\times 7+89.$
\item [] $5010=12\times (3+4\times 5)+6\times 789.$
\item [] $5011=1\times 23\times (4+5\times 6\times 7)+89.$
\item [] $5012=1+23\times (4+5\times 6\times 7)+89.$
\item [] $5013=12+(34+5\times 6)\times 78+9.$
\item [] $5014=1234+5\times (6+78)\times 9.$
\item [] $5015=1^2\times ((3+45)\times 6+7)\times (8+9).$
\item [] $5016=12\times 3^4\times 5+67+89.$
\item [] $5017=(1\times 234+56\times 7)\times 8+9.$
\item [] $5018=(1^2+3)^4\times 5+6\times 7\times 89.$
\item [] $5019=123+(45+6)\times (7+89).$
\item [] $\mathit{5020=-1+(2+3)\times 4^5-6\times (7+8)-9.}$
\item [] $\mathit{5021=-1+2\times 3^4\times 5+6\times 78\times 9.}$
\item [] $5022=12\times 3\times 4\times 5\times 6+78\times 9.$
\item [] $5023=1+2\times 3^4\times 5+6\times 78\times 9.$
\item [] $5024=(12\times 3+4)\times 5+67\times 8\times 9.$
\item [] $5025=1^2\times 3+(4+5)\times (6+7\times 8)\times 9.$
\item [] $5026=1^2+3+(4+5)\times (6+7\times 8)\times 9.$
\item [] $5027=(1+2\times 3^4)\times 5+6\times 78\times 9.$
\item [] $5028=12\times (3\times 4+5)+67\times 8\times 9.$
\item [] $5029=1+(234+5)\times (6+7+8)+9.$
\item [] $5030=(1+2\times 34+5)\times 67+8\times 9.$
\item[]$\mbox{Decreasing order}$
\item [] $4961=9\times (8\times 7\times 6+5\times 43)+2\times 1.$
\item [] $4962=98+76\times (54+3^2+1).$
\item [] $4963=(9\times 87+(6+5)\times 4)\times 3\times 2+1.$
\item [] $4964=9\times 87+65\times 4^3+21.$
\item [] $4965=9\times (8+7\times 6)+5\times 43\times 21.$
\item [] $4966=9+8+7\times (6+5)\times 4^3+21.$
\item [] $4967=(98+76)\times 5+4^{(3\times 2)}+1.$
\item [] $4968=9+87\times (6+5+43+2+1).$
\item [] $4969=9\times 87\times 6+54\times (3+2)+1.$
\item [] $4970=9\times 8\times (7\times 6+(5+4)\times 3)+2\times 1.$
\item [] $4971=9\times 8+7\times 654+321.$
\item [] $4972=(9\times 8+76\times 5)\times (4+3\times 2+1).$
\item [] $4973=98+(7+6)\times (54+321).$
\item [] $4974=(987+6)\times 5+4+3+2\times 1.$
\item [] $4975=(987+6)\times 5+4+3+2+1.$
\item [] $4976=(987+6)\times 5+4+3\times 2+1.$
\item [] $4977=9+8\times (76+543+2\times 1).$
\item [] $4978=9+8\times (76+5)+4321.$
\item [] $4979=(987+6)\times 5+4+3^2+1.$
\item [] $4980=(987+6)\times 5+4\times 3+2+1.$
\item [] $4981=(9+8)\times ((7+65)\times 4+3+2\times 1).$
\item [] $4982=((9+8\times 76+5)\times 4+3)\times 2\times 1.$
\item [] $4983=((9+8\times 76+5)\times 4+3)\times 2+1.$
\item [] $4984=(9+8\times 76+5^4\times 3)\times 2\times 1.$
\item [] $4985=(9+8\times (7\times (6+5)\times 4+3)\times 2\times 1).$
\item [] $4986=9\times (8\times 7\times 6+5\times 43+2+1).$
\item [] $4987=9\times 87\times 6+(5+4)\times 32+1.$
\item [] $4988=9\times 87\times 6+(5+4\times 3)^2+1.$
\item [] $4989=(987+6)\times 5+4\times 3\times 2\times 1.$
\item [] $4990=(987+6)\times 5+4\times 3\times 2+1.$
\item [] $4991=(98+7+6+5)\times 43+2+1.$
\item [] $4992=(9+87+6+54)\times 32\times 1.$
\item [] $4993=9+8\times 7\times 65+4^3\times 21.$
\item [] $4994=9+8\times 7\times (65+4\times 3\times 2)+1.$
\item [] $4995=(98+7+6)\times 5\times (4+3+2)\times 1.$
\item [] $4996=(98+7+6)\times 5\times (4+3+2)+1.$
\item [] $4997=98+7\times 654+321.$
\item [] $4998=(987+6)\times 5+4\times 3+21.$
\item [] $4999=(98+7)\times 6\times 5+43^2\times 1.$
\item [] $5000=(98+7)\times 6\times 5+43^2+1.$
\item [] $5001=9+(87+65+4)\times 32\times 1.$
\item [] $5002=(987+6)\times 5+4+32+1.$
\item [] $5003=9\times 8+7\times (6+5)\times 4^3+2+1.$
\item [] $5004=98\times (7\times 6+5+4)+3+2+1.$
\item [] $5005=9\times 8\times (7+6)\times 5+4+321.$
\item [] $5006=98\times 7+6\times 5\times (4\times 3)^2\times 1.$
\item [] $5007=98\times 7+6\times 5\times (4\times 3)^2+1.$
\item [] $5008=98\times (7\times 6+5+4)+3^2+1.$
\item [] $5009=(98+7+6+5)\times 43+21.$
\item [] $5010=(987+6)\times 5+43+2\times 1.$
\item [] $5011=(987+6)\times 5+43+2+1.$
\item [] $5012=(98+76+5)\times (4+3+21).$
\item [] $5013=(98+7\times 65+4)\times 3^2\times 1.$
\item [] $5014=(98+7\times 65+4)\times 3^2+1.$
\item [] $5015=(9\times 8+7+6)\times (54+3+2\times 1).$
\item [] $5016=9\times (8+7)\times (6\times 5+4+3)+21.$
\item [] $5017=(9+8\times 7+6+5)\times (4^3+2)+1.$
\item [] $5018=98\times 7+6+5+4321.$
\item [] $5019=(9+8+7\times 6\times 5+4\times 3)\times 21.$
\item [] $5020=9\times 87\times 6+5\times 4^3+2\times 1.$
\item [] $5021=9\times 87\times 6+5\times 4^3+2+1.$
\item [] $5022=9\times 87\times 6+54\times 3\times 2\times 1.$
\item [] $5023=9\times 87\times 6+54\times 3\times 2+1.$
\item [] $5024=98\times 7\times 6+5+43\times 21.$
\item [] $5025=9\times 8\times 7+6+5\times 43\times 21.$
\item [] $5026=98+7\times (6+5)\times (43+21).$
\item [] $5027=9\times (8\times 7+6)\times (5+4)+3+2\times 1.$
\item [] $5028=9+87\times 6\times (5+4)+321.$
\item [] $5029=(987+6)\times 5+43+21.$
\item [] $5030=98\times (7\times 6+5+4)+32\times 1.$
\item[]$\mbox{Increasing order}$
\item [] $5031=1\times 23\times (4+5)+67\times 8\times 9.$
\item [] $5032=1+23\times (4+5)+67\times 8\times 9.$
\item [] $5033=12\times (345+67)+89.$
\item [] $5034=(12+3)\times 4\times 5+6\times 789.$
\item [] $5035=1+2\times 3^4+56\times (78+9).$
\item [] $5036=1+2+(3+4)\times (5+6\times 7\times (8+9)).$
\item [] $5037=123\times 4+567\times 8+9.$
\item [] $5038=1+2^{(3\times 4)}+5+(6+7)\times 8\times 9.$
\item [] $5039=1\times (2\times 3)^4+5+6\times 7\times 89.$
\item [] $5040=1^2\times 3^4\times 56+7\times 8\times 9.$
\item [] $5041=1^2+3^4\times 56+7\times 8\times 9.$
\item [] $5042=1\times 2+3^4\times 56+7\times 8\times 9.$
\item [] $5043=1+2+3^4\times 56+7\times 8\times 9.$
\item [] $\mathit{5044=-1+(2+3)\times 4^5-6-78+9.}$
\item [] $5045=1\times 23+(4+5)\times (6+7\times 8)\times 9.$
\item [] $5046=(1+2\times 3+45+6)\times (78+9).$
\item [] $5047=(1+2\times 34+5)\times 67+89.$
\item [] $\mathit{5048=1-2+3\times 4\times 5\times (6+78)+9.}$
\item [] $5049=(1+2\times 3)\times 45+6\times 789.$
\item [] $5050=1+(2+3)\times 45+67\times 8\times 9.$
\item [] $5051=1\times 2+(34+56)\times 7\times 8+9.$
\item [] $5052=12+3^4\times 56+7\times 8\times 9.$
\item [] $5053=1+2\times (3^4\times 5\times 6+7+89).$
\item [] $5054=(12+34)\times 5+67\times 8\times 9.$
\item [] $5055=1+2\times (34\times 56+7\times 89).$
\item [] $5056=1^2\times (34+5\times 6)\times (7+8\times 9).$
\item [] $5057=1\times 2\times 3\times 4\times 5\times 6\times 7+8+9.$
\item [] $5058=1+2\times 3\times 4\times 5\times 6\times 7+8+9.$
\item [] $5059=1+2+(34+5\times 6)\times (7+8\times 9).$
\item [] $5060=123\times (4+5\times 6+7)+8+9.$
\item [] $5061=12+(34+56)\times 7\times 8+9.$
\item [] $5062=1^2+3\times (4+5\times 6\times 7\times 8)+9.$
\item [] $5063=1\times 234+5+67\times 8\times 9.$
\item [] $5064=12\times 34\times 5+6\times 7\times 8\times 9.$
\item [] $5065=(1\times 2+3)\times (4\times 56+789).$
\item [] $5066=(1+234+56+7)\times (8+9).$
\item [] $5067=(123\times 4+56+7+8)\times 9.$
\item [] $5068=12+(34+5\times 6)\times (7+8\times 9).$
\item [] $5069=1\times 23\times 4\times (5+6\times 7+8)+9.$
\item [] $5070=(12+3\times 4\times 5+6)\times (7\times 8+9).$
\item [] $5071=1^2+3\times (4\times 5+6)\times (7\times 8+9).$
\item [] $5072=(1+2)\times 3^4+5+67\times 8\times 9.$
\item [] $5073=(1+2+3+4+5+6\times 7)\times 89.$
\item [] $5074=1\times 2\times 34\times 5+6\times 789.$
\item [] $5075=1+2\times 34\times 5+6\times 789.$
\item [] $5076=12\times 345+(6+7)\times 8\times 9.$
\item [] $5077=1+2+(3^4+5)\times (6\times 7+8+9).$
\item [] $5078=12+34\times (5+6\times (7+8+9)).$
\item [] $5079=1^2\times 345+6\times 789.$
\item [] $5080=1^2+345+6\times 789.$
\item [] $5081=1\times 2+345+6\times 789.$
\item [] $5082=1+2+345+6\times 789.$
\item [] $5083=1+2\times 3\times (4\times 56+7\times 89).$
\item [] $5084=1^2+(3+4\times 5)\times (6+7)\times (8+9).$
\item [] $5085=(1\times 23\times 4+5+6\times 78)\times 9.$
\item [] $5086=12+(3^4+5)\times (6\times 7+8+9).$
\item [] $\mathit{5087=1\times 23\times 4\times 56-7\times 8-9.}$
\item [] $5088=12\times 34+5\times (6+7)\times 8\times 9.$
\item [] $5089=1+(2^3+4)\times (5\times 67+89).$
\item [] $5090=1\times 2+3\times 4\times (5\times 67+89).$
\item [] $5091=12+345+6\times 789.$
\item [] $5092=1+2\times 3+45\times ((6+7)\times 8+9).$
\item [] $5093=(1+2)^3\times 4\times (5+6\times 7)+8+9.$
\item [] $5094=1\times 2^3\times 45+6\times 789.$
\item [] $5095=1+2^3\times 45+6\times 789.$
\item [] $5096=(1^2+3^4)\times 56+7\times 8\times 9.$
\item [] $5097=12^3+4\times 56\times (7+8)+9.$
\item [] $5098=1\times 2\times (3^4\times 5\times 6+7\times (8+9)).$
\item [] $5099=(12\times 3\times 4\times 5+6)\times 7+8+9.$
\item [] $5100=(12+3+45)\times (6+7+8\times 9).$
\item[]$\mbox{Decreasing order}$
\item [] $5031=(987+6)\times 5+4^3+2\times 1.$
\item [] $5032=(987+6)\times 5+4+3\times 21.$
\item [] $5033=98\times 7+(65+4)\times 3\times 21.$
\item [] $5034=9\times (87\times 6+5\times (4+3))+21.$
\item [] $5035=9+(87+6+5^4)\times (3\times 2+1).$
\item [] $5036=9+8+7\times (654+3\times 21).$
\item [] $5037=98\times 7+6\times 5+4321.$
\item [] $5038=98+76\times 5\times (4+3^2\times 1).$
\item [] $5039=9\times 87\times 6+5\times 4+321.$
\item [] $5040=(98+7)\times (6+5+4+32+1).$
\item [] $5041=9\times (8+7+65)+4321.$
\item [] $5042=9+(87+6\times 5)\times 43+2\times 1.$
\item [] $5043=9+(87+6\times 5)\times 43+2+1.$
\item [] $\mathit{5044=9\times 8\times 76+5-432-1.}$
\item [] $5045=(9+(8+7\times 6)\times 5\times 4)\times (3+2)\times 1.$
\item [] $5046=9+87\times 6+5\times 43\times 21.$
\item [] $5047=98+7\times (6+5)\times 4^3+21.$
\item [] $\mathit{5048=9+(8+76)\times 5\times 4\times 3-2+1.}$
\item [] $5049=9+8\times 7\times 6\times (5+4+3+2+1).$
\item [] $5050=(987+6)\times 5+4^3+21.$
\item [] $5051=(987+6)\times 5+43\times 2\times 1.$
\item [] $5052=(987+6)\times 5+43\times 2+1.$
\item [] $5053=98+7\times 6+(5+4\times 3)^(2+1).$
\item [] $5054=(9+8+76)\times 54+32\times 1.$
\item [] $5055=9+(87+6)\times 54+3+21.$
\item [] $5056=9+87\times (6+5\times 4+32)+1.$
\item [] $5057=9+8+7\times 6\times 5\times 4\times 3\times 2\times 1.$
\item [] $5058=987+6\times 5^4+321.$
\item [] $5059=(9\times 87+6+54)\times 3\times 2+1.$
\item [] $5060=98\times 7+6\times ((5+4)\times 3)^2\times 1.$
\item [] $5061=(9\times 8+76)\times 5+4321.$
\item [] $5062=98\times 7\times 6+5^4+321.$
\item [] $5063=98\times 7+6\times (5+4)^3+2+1.$
\item [] $5064=9\times (8+76)\times 5+4\times 321.$
\item [] $5065=9+8\times (76\times 5+4\times 3\times 21).$
\item [] $5066=(987+6+5\times 4)\times (3+2)+1.$
\item [] $5067=98\times 7+6+5^4\times (3\times 2+1).$
\item [] $5068=9\times (87\times 6+5)+4+321.$
\item [] $5069=(9+(8+76)\times 5\times 4)\times 3+2\times 1.$
\item [] $5070=9+(8+76)\times 5\times 4\times 3+21.$
\item [] $5071=(9+8\times 7)\times (65+4+3^2)+1.$
\item [] $5072=98\times 7+65+4321.$
\item [] $5073=9\times 87\times 6+54+321.$
\item [] $5074=9\times 87+65\times (4^3+2)+1.$
\item [] $5075=(9+8+7\times 6\times 5\times 4\times 3)\times 2+1.$
\item [] $5076=9\times (87\times 6+5+4+32+1).$
\item [] $5077=9\times 8+7\times 65\times (4+3\times 2+1).$
\item [] $5078=98\times 7+6\times ((5+4)^3+2+1).$
\item [] $5079=((9+8\times 7)\times (6+5\times 4)+3)\times (2+1).$
\item [] $5080=9+(8+7)\times 65+4^{(3+2+1)}.$
\item [] $5081=98\times 7+6\times (5+4)^3+21.$
\item [] $5082=9\times (8+76)+5+4321.$
\item [] $5083=(9+8)\times (7+65\times 4+32\times 1).$
\item [] $5084=(9+8)\times (7+(65\times 4+32))+1.$
\item [] $5085=(9+8+76)\times 54+3\times 21.$
\item [] $5086=9\times (87+6+5\times 4)\times (3+2)+1.$
\item [] $5087=98+76+(5+4\times 3)^(2+1).$
\item [] $5088=9\times (87+6)\times 5+43\times 21.$
\item [] $5089=9+8\times (7+6\times 5\times 4)\times (3+2\times 1).$
\item [] $5090=9+8\times (7+6\times 5\times 4)\times (3+2)+1.$
\item [] $5091=(9+87)\times 6+5\times 43\times 21.$
\item [] $5092=(9+8\times 7+6+5)\times (4+3\times 21).$
\item [] $5093=(987+6)\times 5+4\times 32\times 1.$
\item [] $5094=(987+6)\times 5+4\times 32+1.$
\item [] $5095=9\times (8\times 7+6\times 5)+4321.$
\item [] $5096=98\times 7+6\times 5\times (4+3)\times 21.$
\item [] $5097=(9+(8+7)\times 6+5)\times (4+3)^2+1.$
\item [] $5098=98\times (7+(6+5+4)\times 3)+2\times 1.$
\item [] $5099=9+8+7\times 6\times (5\times 4\times 3\times 2+1).$
\item [] $5100=9\times (87\times 6+5+4)+321.$
\item[]$\mbox{Increasing order}$
\item [] $5101=1^2+34\times 5\times (6+7+8+9).$
\item [] $5102=1\times 2+34\times 5\times (6+7+8+9).$
\item [] $5103=1+2+34\times 5\times (6+7+8+9).$
\item [] $5104=1\times 2^3\times (4+5+6+7\times 89).$
\item [] $5105=(1^2+34+56)\times 7\times 8+9.$
\item [] $5106=1\times 234+56\times (78+9).$
\item [] $5107=1+234+56\times (78+9).$
\item [] $5108=12^3+4\times (56+789).$
\item [] $5109=12\times 3\times 4\times 5\times 6+789.$
\item [] $5110=1+(23+45)\times (67+8)+9.$
\item [] $5111=12\times 34\times (5+6)+7\times 89.$
\item [] $5112=1\times 2\times 3\times 4\times 5\times 6\times 7+8\times 9.$
\item [] $5113=1+2\times 3\times 4\times 5\times 6\times 7+8\times 9.$
\item [] $5114=(1+2+3\times 4)\times 5\times 67+89.$
\item [] $5115=123\times (4+5\times 6+7)+8\times 9.$
\item [] $5116=12^3+4+(5+6\times 7)\times 8\times 9.$
\item [] $5117=(1+23+45\times 6+7)\times (8+9).$
\item [] $5118=1+(2+(3+4)\times 5+6)\times 7\times (8+9).$
\item [] $5119=1+2\times 3\times (4+56\times (7+8)+9).$
\item [] $5120=1\times (2+3+4)\times 567+8+9.$
\item [] $5121=1+(2+3+4)\times 567+8+9.$
\item [] $5122=1\times 2+3^4\times (56+7)+8+9.$
\item [] $5123=1+2+3^4\times (56+7)+8+9.$
\item [] $5124=(12+3)\times 4\times 5+67\times 8\times 9.$
\item [] $5125=1^2+(3+4)\times (5\times 6+78\times 9).$
\item [] $5126=(1+2\times 3)^4+5\times (67\times 8+9).$
\item [] $5127=123+(4\times 5+67\times 8)\times 9.$
\item [] $5128=(1^2+3)^4+56\times (78+9).$
\item [] $5129=1\times 2\times 3\times 4\times 5\times 6\times 7+89.$
\item [] $5130=1+2\times 3\times 4\times 5\times 6\times 7+89.$
\item [] $5131=12^3+4+5\times 678+9.$
\item [] $5132=123\times (4+5\times 6+7)+89.$
\item [] $5133=1^2\times 3\times (4^5+678+9).$
\item [] $5134=1^2+3\times (4^5+678+9).$
\item [] $5135=1\times 2+3\times (4^5+678+9).$
\item [] $5136=1+2+3\times (4^5+678+9).$
\item [] $5137=1+2\times 3+45\times (6\times 7+8\times 9).$
\item [] $5138=1\times 2^3+45\times (6\times 7+8\times 9).$
\item [] $5139=1^2\times 3^4\times 5+6\times 789.$
\item [] $5140=1^2+3^4\times 5+6\times 789.$
\item [] $5141=1\times 2+3^4\times 5+6\times 789.$
\item [] $5142=1+2+3^4\times 5+6\times 789.$
\item [] $5143=1+2\times 3\times (4\times 5\times 6\times 7+8+9).$
\item [] $5144=1\times 2^3\times (4+567+8\times 9).$
\item [] $5145=12+3+45\times (6\times 7+8\times 9).$
\item [] $5146=1^2+3\times (4+5\times 6\times 7)\times 8+9.$
\item [] $5147=12\times 34+5+6\times 789.$
\item [] $5148=12\times (345+67+8+9).$
\item [] $5149=1+(2+34)\times (56+78+9).$
\item [] $5150=1+(2+3^4)\times 5+6\times 789.$
\item [] $5151=12+3^4\times 5+6\times 789.$
\item [] $5152=1\times (2+3^4)\times 56+7\times 8\times 9.$
\item [] $5153=1+(2+3^4)\times 56+7\times 8\times 9.$
\item [] $5154=12\times (3+4)\times 5+6\times 789.$
\item [] $5155=12\times 3^4+(5+6\times 7)\times 89.$
\item [] $5156=(1+2)\times 3\times (4+567)+8+9.$
\item [] $5157=12\times (345+6+78)+9.$
\item [] $5158=1+(2+34+5\times 6)\times 78+9.$
\item [] $5159=1^2\times 3^4\times 56+7\times 89.$
\item [] $5160=1^2+3^4\times 56+7\times 89.$
\item [] $5161=1\times 2+3^4\times 56+7\times 89.$
\item [] $5162=1+2+3^4\times 56+7\times 89.$
\item [] $5163=1^{23}+(45+6+7)\times 89.$
\item [] $5164=1\times 2\times 34\times 5+67\times 8\times 9.$
\item [] $5165=1+2\times 34\times 5+67\times 8\times 9.$
\item [] $5166=1^2+3+(45+6+7)\times 89.$
\item [] $5167=12^3+4+5\times (678+9).$
\item [] $5168=1^2\times 34\times (56+7+89).$
\item [] $5169=1^2\times 345+67\times 8\times 9.$
\item [] $5170=1^2+345+67\times 8\times 9.$
\item[]$\mbox{Decreasing order}$
\item [] $5101=(9+87+6)\times 5\times (4+3\times 2)+1.$
\item [] $5102=(9\times 8+7+6)\times 5\times 4\times 3+2\times 1.$
\item [] $5103=9+8+765+4321.$
\item [] $5104=9\times 87+6\times 5\times (4\times 3)^2+1.$
\item [] $5105=(987+6\times 5+4)\times (3+2)\times 1.$
\item [] $5106=(987+6\times 5+4)\times (3+2)+1.$
\item [] $5107=(98+7+6)\times (5\times 4+3)\times 2+1.$
\item [] $\mathit{5108=(9\times 8+7)\times 65+4-32+1.}$
\item [] $5109=(987+6)\times 5+(4\times 3)^2\times 1.$
\item [] $5110=(987+6)\times 5+(4\times 3)^2+1.$
\item [] $\mathit{5111=9\times 8+7\times 6\times 5\times 4\times 3\times 2-1.}$
\item [] $5112=9\times 8+7\times 6\times 5\times 4\times 3\times 2\times 1.$
\item [] $5113=9\times 8+7\times 6\times 5\times 4\times 3\times 2+1.$
\item [] $5114=987+6\times 5+4^{(3\times 2)}+1.$
\item [] $5115=9\times 87+6+5+4321.$
\item [] $\mathit{5116=-9\times 8\times 7-6+5^4\times 3^2+1.}$
\item [] $5117=98+7\times (654+3\times 21).$
\item [] $5118=9+8\times (7+6+5^4)+3+2\times 1.$
\item [] $5119=(9\times 8+7\times 6+5)\times 43+2\times 1.$
\item [] $5120=(98+7\times 6+5\times 4)\times 32\times 1.$
\item [] $5121=(9\times 8+7+6)\times 5\times 4\times 3+21.$
\item [] $5122=9\times (8+76+5)+4321.$
\item [] $5123=9\times 87\times 6+5\times (4^3+21).$
\item [] $5124=9\times (8+76)\times 5+4^3\times 21.$
\item [] $5125=(98\times (7+6)+5)\times 4+3^2\times 1.$
\item [] $5126=(98\times (7+6)+5)\times 4+3^2+1.$
\item [] $\mathit{5127=9\times 87\times 6+5\times 43\times 2-1.}$
\item [] $5128=9\times 87\times 6+5\times 43\times 2\times 1.$
\item [] $5129=9\times 87\times 6+5\times 43\times 2+1.$
\item [] $5130=9\times 87+(65+4)\times 3\times 21.$
\item [] $5131=9+(8+7+65)\times 4^3+2\times 1.$
\item [] $5132=9+8\times 76+5\times 43\times 21.$
\item [] $5133=9\times 87\times 6+5\times (43\times 2+1).$
\item [] $5134=9\times 87+6\times 5+4321.$
\item [] $5135=9\times 87\times 6+5+432\times 1.$
\item [] $5136=9\times 87\times 6+5+432+1.$
\item [] $5137=9+8\times (7+6+5^4)+3+21.$
\item [] $5138=98+7\times 6\times 5\times 4\times 3\times 2\times 1.$
\item [] $5139=98+7\times 6\times 5\times 4\times 3\times 2+1.$
\item [] $5140=(98\times (7+6)+5)\times 4+3+21.$
\item [] $5141=9+(8+7)\times 6\times (54+3)+2\times 1.$
\item [] $5142=9+87\times (6+5\times 4+32+1).$
\item [] $5143=(9+8+7\times 6\times 5\times 4)\times 3\times 2+1.$
\item [] $5144=(9\times 8+7)\times 65+4+3+2\times 1.$
\item [] $5145=9\times 8\times (7+6+54)+321.$
\item [] $5146=(9\times 8+7)\times 65+4+3\times 2+1.$
\item [] $5147=987+65\times (43+21).$
\item [] $5148=987+65+4^{(3+2+1)}.$
\item [] $5149=987+65\times 4^3+2\times 1.$
\item [] $5150=987+65\times 4^3+2+1.$
\item [] $5151=9+8+7+6+5\times 4^(3+2)+1.$
\item [] $5152=(98+7+6+5^4)\times (3\times 2+1).$
\item [] $5153=9+8\times (7\times 6\times 5+432+1).$
\item [] $5154=9+(8+7+6)\times 5\times (4+3)^2\times 1.$
\item [] $5155=9+(8\times 7+6)\times (5\times 4+3\times 21).$
\item [] $5156=(9\times 8+7)\times 65+4\times (3+2)+1.$
\item [] $5157=9\times (87\times 6+5+43+2+1).$
\item [] $5158=9\times 8+765+4321.$
\item [] $5159=(9\times 8+7)\times 65+4\times 3\times 2\times 1.$
\item [] $5160=9\times 87+6\times (5+4)^3+2+1.$
\item [] $5161=9+8\times 7\times (6+54+32)\times 1.$
\item [] $5162=9+8\times 7\times (6+54+32)+1.$
\item [] $5163=9\times (87+6)+5+4321.$
\item [] $5164=9\times 87+6+5^4\times (3\times 2+1).$
\item [] $5165=9+8+7\times (6+(5+4)^3)+2+1.$
\item [] $5166=9+8\times (7\times 6\times 5+4)\times 3+21.$
\item [] $5167=((9\times 8+7\times 6)\times 5+4)\times 3^2+1.$
\item [] $5168=987+65\times 4^3+21.$
\item [] $5169=9\times 87+65+4321.$
\item [] $5170=9\times 87+6\times ((5+4)^3+2)+1.$
\item[]$\mbox{Increasing order}$
\item [] $5171=12+3^4\times 56+7\times 89.$
\item [] $5172=12\times (3+4+5\times 67+89).$
\item [] $5173=12+(3^4+5+6)\times 7\times 8+9.$
\item [] $5174=(12+3)\times (4+5\times 67)+89.$
\item [] $5175=(1\times 2+3+4)\times 567+8\times 9.$
\item [] $5176=1\times 23\times 4\times 56+7+8+9.$
\item [] $5177=1+23\times 4\times 56+7+8+9.$
\item [] $5178=1+2+3^4\times (56+7)+8\times 9.$
\item [] $5179=1^{23}\times 4+(567+8)\times 9.$
\item [] $5180=12+34\times (56+7+89).$
\item [] $5181=12+345+67\times 8\times 9.$
\item [] $5182=1^2\times 3+4+(567+8)\times 9.$
\item [] $5183=1^2+3+4+(567+8)\times 9.$
\item [] $5184=1\times 2^3\times 45+67\times 8\times 9.$
\item [] $5185=1+2^3\times 45+67\times 8\times 9.$
\item [] $5186=1+23+(45+6+7)\times 89.$
\item [] $5187=12+3^4\times (56+7)+8\times 9.$
\item [] $5188=1^2+3\times 4+(567+8)\times 9.$
\item [] $5189=12\times 3^4+5+6\times 78\times 9.$
\item [] $5190=12\times 34\times (5+6)+78\times 9.$
\item [] $5191=1^2\times 3+4+(5+67)\times 8\times 9.$
\item [] $5192=1^2\times 3^4\times (56+7)+89.$
\item [] $5193=1+(2+3+4)\times 567+89.$
\item [] $5194=1\times 23\times 4\times 5+6\times 789.$
\item [] $5195=1+23\times 4\times 5+6\times 789.$
\item [] $5196=12\times (3^4+5\times 67+8+9).$
\item [] $5197=1+2^3+4+(5+67)\times 8\times 9.$
\item [] $5198=12\times 3+(45+6+7)\times 89.$
\item [] $5199=(12+3^4)\times 5+6\times 789.$
\item [] $5200=1+2\times 3\times 4+(567+8)\times 9.$
\item [] $5201=12\times (34+56\times 7)+89.$
\item [] $5202=12\times (34+5)+6\times 789.$
\item [] $5203=1+23+4+(567+8)\times 9.$
\item [] $5204=12+3^4\times (56+7)+89.$
\item [] $5205=12\times 3^4\times 5+6\times 7\times 8+9.$
\item [] $5206=(1+2)^3+4+(567+8)\times 9.$
\item [] $5207=1\times 23+4\times (5+6+7)\times 8\times 9.$
\item [] $5208=1+2^3\times 4+(567+8)\times 9.$
\item [] $5209=1234+5\times (6+789).$
\item [] $5210=1^2+34+(567+8)\times 9.$
\item [] $5211=1\times 23+4+(5+67)\times 8\times 9.$
\item [] $5212=1+2+34+(567+8)\times 9.$
\item [] $5213=1+2\times (34\times 56+78\times 9).$
\item [] $5214=(1+23)\times 4\times 5+6\times 789.$
\item [] $5215=1^2+3+(4+567+8)\times 9 .$
\item [] $5216=1\times 2+3+(4+567+8)\times 9 .$
\item [] $5217=1\times 23\times 4\times 56+7\times 8+9.$
\item [] $5218=1+23\times 4\times 56+7\times 8+9.$
\item [] $5219=1\times 2^3+(4+567+8)\times 9 .$
\item [] $5220=12\times 3+4\times (5+6+7)\times 8\times 9.$
\item [] $5221=12+34+(567+8)\times 9.$
\item [] $5222=1\times (2+3)\times 4^5+6+7+89.$
\item [] $5223=1+(2+3)\times 4^5+6+7+89.$
\item [] $5224=12\times 3+4+(5+67)\times 8\times 9.$
\item [] $5225=1\times 2+3+(4+56)\times (78+9).$
\item [] $5226=12+3+(4+567+8)\times 9 .$
\item [] $5227=1+2\times 3+(4+56)\times (78+9).$
\item [] $5228=(1+2)\times 3\times (4+567)+89.$
\item [] $5229=1^2\times 3^4\times 5+67\times 8\times 9.$
\item [] $5230=1^2+3^4\times 5+67\times 8\times 9.$
\item [] $5231=123\times 4+5+6\times 789.$
\item [] $5232=1+23\times 4\times 56+7+8\times 9.$
\item [] $5233=1^2+3+(45+67\times 8)\times 9.$
\item [] $5234=1\times 23+(4+567+8)\times 9 .$
\item [] $5235=1+23+(4+567+8)\times 9 .$
\item [] $5236=1^{23}\times 4^5+6\times 78\times 9.$
\item [] $5237=1^{23}+4^5+6\times 78\times 9.$
\item [] $5238=1^2\times 3^4\times 56+78\times 9.$
\item [] $5239=1\times 23\times 4\times 56+78+9.$
\item [] $5240=1+23\times 4\times 56+78+9.$
\item[]$\mbox{Decreasing order}$
\item [] $5171=(9\times 8+7)\times 65+4+32\times 1.$
\item [] $5172=(9\times 8+7)\times 65+4+32+1.$
\item [] $5173=(9\times 8+7\times 6+5^4)\times (3\times 2+1).$
\item [] $5174=98+(7\times 6+5)\times 4\times 3^{(2+1)}.$
\item [] $5175=9+(8+76+54\times 3)\times 21.$
\item [] $5176=9\times (87\times 6+5)+432+1.$
\item [] $5177=9+8+7\times (6+(5+4)^3+2)+1.$
\item [] $5178=9\times 87+6\times (5+4)^3+21.$
\item [] $5179=9+8+7\times 6+5\times 4^(3+2)\times 1.$
\item [] $5180=(9\times 8+7)\times 65+43+2\times 1.$
\item [] $5181=(9\times 8+7)\times 65+43+2+1.$
\item [] $5182=(9\times 8+76)\times 5\times (4+3)+2\times 1.$
\item [] $5183=(9\times 8+76)\times 5\times (4+3)+2+1.$
\item [] $5184=98+765+4321.$
\item [] $5185=9\times 87\times 6+54\times 3^2+1.$
\item [] $5186=(9+87)\times (6+5+43)+2\times 1.$
\item [] $5187=(9\times 8+7+6+54\times 3)\times 21.$
\item [] $5188=(98\times 7+6)\times 5+(4\times 3)^{(2+1)}.$
\item [] $5189=9\times (87\times 6+54)+3+2\times 1.$
\item [] $5190=9+(8+7)\times 6\times 54+321.$
\item [] $5191=(98+76)\times 5+4321.$
\item [] $5192=9+8\times 7+6+5\times 4^(3+2)+1.$
\item [] $5193=9\times (87\times 6+54)+3^2\times 1.$
\item [] $5194=9\times (87\times 6+54)+3^2+1.$
\item [] $5195=9\times 8\times (7+65)+4+3\times 2+1.$
\item [] $5196=(9+8\times 7\times 6\times 5+43)\times (2+1).$
\item [] $5197=9\times 8\times (7+65)+4+3^2\times 1.$
\item [] $5198=9\times 8\times (7+65)+4+3^2+1.$
\item [] $5199=(9\times 8+7)\times 65+43+21.$
\item [] $5200=(9+8\times 7)\times (65+4\times 3+2+1).$
\item [] $5201=(9\times 8+76)\times 5\times (4+3)+21.$
\item [] $5202=98\times 7\times 6+543\times 2\times 1.$
\item [] $5203=98\times 7\times 6+543\times 2+1.$
\item [] $5204=9+8+(76+5)\times 4^3+2+1.$
\item [] $5205=9\times 8\times (7+65)+4\times (3+2)+1.$
\item [] $5206=9\times 8+7+6+5\times 4^(3+2)+1.$
\item [] $5207=98\times 7+6+5\times 43\times 21.$
\item [] $5208=9\times (87\times 6+54)+3+21.$
\item [] $5209=9+8+7+(65+4+3)^2+1.$
\item [] $5210=9+8\times (76+54)\times (3+2)+1.$
\item [] $5211=9+876+5+4321.$
\item [] $5212=9\times 8\times (7+65)+4+3+21.$
\item [] $5213=9+8+76+5\times 4^(3+2)\times 1.$
\item [] $5214=9+8+76+5\times 4^(3+2)+1.$
\item [] $5215=9+(8\times 7+65)\times 43+2+1.$
\item [] $5216=9\times (87\times 6+54)+32\times 1.$
\item [] $5217=(987+6)\times 5+4\times 3\times 21.$
\item [] $5218=9+8+76+5\times (4^(3+2)+1).$
\item [] $5219=9+(8+7)\times 6+5\times 4^(3+2)\times 1.$
\item [] $5220=9\times (87+6+54\times 3^2+1).$
\item [] $5221=9\times 8\times (7+65)+4+32+1.$
\item [] $5222=9+8+(76+5)\times 4^3+21.$
\item [] $5223=9+87+6+5\times 4^(3+2)+1.$
\item [] $5224=9\times 8\times (7+65)+4\times (3^2+1).$
\item [] $5225=9+8\times (7+6\times 54+321).$
\item [] $5226=9\times (8\times (7+65)+4)+3\times 2\times 1.$
\item [] $5227=9\times 87\times 6+(5\times 4+3)^2\times 1.$
\item [] $5228=9\times 87\times 6+(5\times 4+3)^2+1.$
\item [] $5229=(9\times 87+6)\times 5+4\times 321.$
\item [] $5230=9\times (87\times 6+54+3+2)+1.$
\item [] $5231=98+7+6+5\times 4^(3+2)\times 1.$
\item [] $5232=98+7+6+5\times 4^(3+2)+1.$
\item [] $5233=9+(8\times 7+65)\times 43+21.$
\item [] $5234=9+87\times (6+54)+3+2\times 1.$
\item [] $5235=9+87\times (6+54)+3+2+1.$
\item [] $5236=9+87\times (6+54)+3\times 2+1.$
\item [] $5237=(98+7\times 6\times 5\times 4\times 3)\times 2+1.$
\item [] $5238=9\times (87\times 6+54+3+2+1).$
\item [] $5239=9+87\times (6+54)+3^2+1.$
\item [] $5240=9\times (87\times 6+5\times 4\times 3)+2\times 1.$
\item[]$\mbox{Increasing order}$
\item [] $5241=12+3^4\times 5+67\times 8\times 9.$
\item [] $5242=1+2+3+4^5+6\times 78\times 9.$
\item [] $5243=1+2\times 3+4^5+6\times 78\times 9.$
\item [] $5244=12\times (3+4)\times 5+67\times 8\times 9.$
\item [] $5245=1+2^3+4^5+6\times 78\times 9.$
\item [] $5246=1+(2+3)\times 4^5+6+7\times (8+9).$
\item [] $5247=1\times 23\times 45+6\times 78\times 9.$
\item [] $5248=1+23\times 45+6\times 78\times 9.$
\item [] $5249=1+23\times 4\times 56+7+89.$
\item [] $5250=12+3^4\times 56+78\times 9.$
\item [] $5251=12+3+4^5+6\times 78\times 9.$
\item [] $5252=1\times 23+(45+67\times 8)\times 9.$
\item [] $5253=123+45\times (6\times 7+8\times 9).$
\item [] $5254=1+2\times 34\times (5+6)\times 7+8+9.$
\item [] $5255=(1^2+3^4)\times (56+7)+89.$
\item [] $5256=12\times (345+6+78+9).$
\item [] $5257=1+2^3\times (4+5\times 6+7\times 89).$
\item [] $5258=1\times 2+3^4+(567+8)\times 9.$
\item [] $5259=1\times 23+4^5+6\times 78\times 9.$
\item [] $5260=1+23+4^5+6\times 78\times 9.$
\item [] $5261=1\times (2+3)\times 4^5+6+(7+8)\times 9.$
\item [] $5262=1+(2+3)\times 4^5+6+(7+8)\times 9.$
\item [] $5263=(1+2)^3+4^5+6\times 78\times 9.$
\item [] $5264=(1\times 2+3)\times 4^5+6\times (7+8+9).$
\item [] $5265=(1+2+3+4+567+8)\times 9 .$
\item [] $5266=1^2+3^4+(5+67)\times 8\times 9.$
\item [] $5267=1\times 23\times 4+(567+8)\times 9.$
\item [] $5268=12+3^4+(567+8)\times 9.$
\item [] $5269=1\times (2+3)\times (4^5+6)+7\times (8+9).$
\item [] $5270=(1\times 2+3\times 4\times 5)\times (6+7+8\times 9).$
\item [] $5271=1\times 23\times 4\times 56+7\times (8+9).$
\item [] $5272=12\times 3+4^5+6\times 78\times 9.$
\item [] $5273=(12+3^4)\times 56+7\times 8+9.$
\item [] $5274=(1+2)^3\times 4\times 5+6\times 789.$
\item [] $5275=1^2+(345+6)\times (7+8)+9.$
\item [] $5276=1\times (2+3)\times 4^5+67+89.$
\item [] $5277=12\times 34\times (5+6)+789.$
\item [] $\mathit{5278=12+3^4\times 5\times (6+7)-8+9.}$
\item [] $5279=12\times 345+67\times (8+9).$
\item [] $5280=12\times 34+56\times (78+9).$
\item [] $5281=(1\times 23\times 4+567)\times 8+9.$
\item [] $5282=1^2\times 3^4\times 5\times (6+7)+8+9.$
\item [] $5283=1^2\times (3\times 4+567+8)\times 9.$
\item [] $5284=1\times 23\times 4\times 5+67\times 8\times 9.$
\item [] $5285=1+23\times 4\times 5+67\times 8\times 9.$
\item [] $5286=12+(345+6)\times (7+8)+9.$
\item [] $5287=1\times 23\times 4\times 56+(7+8)\times 9.$
\item [] $5288=1+23\times 4\times 56+(7+8)\times 9.$
\item [] $5289=(1+23\times 4)\times 5+67\times 8\times 9.$
\item [] $5290=1+2\times (34\times 5+6)\times (7+8)+9.$
\item [] $5291=(1+2+34)\times (56+78+9).$
\item [] $5292=12\times 3\times (45+6+7+89).$
\item [] $5293=(1+2+34+5\times 6)\times (7+8\times 9).$
\item [] $5294=(1^2+3^4)\times 56+78\times 9.$
\item [] $5295=(12+3^4)\times 56+78+9.$
\item [] $5296=1+(2+3)\times (45\times 6+789).$
\item [] $5297=12\times 345+(6+7)\times 89.$
\item [] $5298=1\times 2\times (34+5)\times 67+8\times 9.$
\item [] $5299=1+2\times (34+5)\times 67+8\times 9.$
\item [] $5300=(1\times 2+3)\times (4+(5+6)\times (7+89)).$
\item [] $5301=(1\times 2+3\times 4+567+8)\times 9 .$
\item [] $5302=123+4+(567+8)\times 9.$
\item [] $5303=1+2^{(3\times 4)}+(56+78)\times 9.$
\item [] $5304=(12+3^4)\times 56+7+89.$
\item [] $5305=1+2\times (3\times 4+5)\times (67+89).$
\item [] $5306=1^2\times (3+4)\times (56+78\times 9).$
\item [] $5307=123+4\times (5+6+7)\times 8\times 9.$
\item [] $5308=1\times 2\times 34\times (5+6)\times 7+8\times 9.$
\item [] $5309=12\times (3+4+56)\times 7+8+9.$
\item [] $5310=12\times (3+45)+6\times 789.$
\item[]$\mbox{Decreasing order}$
\item [] $5241=9+8\times (7+6+5\times 4\times 32+1).$
\item [] $5242=9+(8+7\times 6\times 5)\times 4\times 3\times 2+1.$
\item [] $5243=9\times 87\times 6+543+2\times 1.$
\item [] $5244=9\times 87\times 6+543+2+1.$
\item [] $5245=(9+8)\times 7\times (6+5)\times 4+3^2\times 1.$
\item [] $5246=(9+8)\times 7\times (6+5)\times 4+3^2+1.$
\item [] $5247=9\times (87\times 6+54)+3\times 21.$
\item [] $5248=9\times 8\times (7+65)+43+21.$
\item [] $5249=9+8\times 7+(65+4+3)^2\times 1.$
\item [] $5250=(9+87+654)\times (3\times 2+1).$
\item [] $5251=(9\times 8)\times (7+65)+4+3\times 21.$
\item [] $5252=9\times (8\times (7+65)+4)+32\times 1.$
\item [] $5253=9+87\times (6+54)+3+21.$
\item [] $5254=9+(87\times 6\times 5+4\times 3)\times 2+1.$
\item [] $\mathit{5255=9\times 8\times 76-5\times 43-2\times 1.}$
\item [] $5256=9\times (87+65+432\times 1).$
\item [] $5257=9\times 8\times 7+(6+5)\times 432+1.$
\item [] $5258=9\times 8+(76+5)\times 4^3+2\times 1.$
\item [] $5259=9\times 8+(76+5)\times 4^3+2+1.$
\item [] $5260=98+7\times 6+5\times 4^(3+2)\times 1.$
\item [] $5261=9+87\times (6+54)+32\times 1.$
\item [] $5262=9\times 87\times 6+543+21.$
\item [] $5263=(9\times 8+7)\times 65+4^3\times 2\times 1.$
\item [] $5264=98\times 7+654\times (3\times 2+1).$
\item [] $5265=9\times (87+65+432+1).$
\item [] $5266=(9+87\times 6+54)\times 3^2+1.$
\item [] $5267=9+8+7\times 6\times 5\times (4\times 3\times 2+1).$
\item [] $5268=9\times 8+76+5\times 4^(3+2)\times 1.$
\item [] $5269=9\times 8\times (7+65)+4^3+21.$
\item [] $5270=9\times 8\times (7+65)+43\times 2\times 1.$
\item [] $5271=9\times (8+76)+5\times 43\times 21.$
\item [] $5272=9\times 8\times 7+6+(5+4^3)^2+1.$
\item [] $5273=9+8\times 7\times (6\times 5+43+21).$
\item [] $5274=9\times (87\times 6+54+3^2+1).$
\item [] $5275=9+8\times 7\times (6\times 5+4^3)+2\times 1.$
\item [] $5276=9+8\times 7\times (6\times 5+4^3)+2+1.$
\item [] $5277=9\times 8+(76+5)\times 4^3+21.$
\item [] $5278=987+65\times (4^3+2)+1.$
\item [] $5279=9+(87\times 6+5)\times (4+3+2+1).$
\item [] $5280=(9+87+65+4)\times 32\times 1.$
\item [] $5281=(9+87+65+4)\times 32+1.$
\item [] $5282=9+8+76\times (5+4^3)+21.$
\item [] $5283=9\times (8\times (7+65)+4)+3\times 21.$
\item [] $5284=98+(76+5)\times 4^3+2\times 1.$
\item [] $5285=98+(76+5)\times 4^3+2+1.$
\item [] $5286=(9+(8+7\times 6\times 5)\times 4)\times 3\times 2\times 1.$
\item [] $5287=9+87\times (6+5)+4321.$
\item [] $\mathit{5288=98\times 76-5\times 432\times 1.}$
\item [] $5289=(9\times 87+6)\times 5+4^3\times 21.$
\item [] $5290=(987+6)\times 5+4+321.$
\item [] $5291=(9+8\times (7+6)\times 5)\times (4+3\times 2)+1.$
\item [] $5292=9+876\times 5+43\times 21.$
\item [] $5293=9+87\times 6+(5+4^3)^2+1.$
\item [] $5294=9+8\times 7\times (6\times 5+4^3)+21.$
\item [] $5295=98+76+5\times 4^(3+2)+1.$
\item [] $5296=(98\times 7+654\times 3)\times 2\times 1.$
\item [] $5297=(98\times 7+654\times 3)\times 2+1.$
\item [] $5298=(9+8)\times 7\times 6\times 5+(4\times 3)^{(2+1)}.$
\item [] $5299=(9+8)\times 7\times (6+5)\times 4+3\times 21.$
\item [] $5300=(9\times 8\times 7+6+5\times 4)\times (3^2+1).$
\item [] $5301=9+(8+7+65+4)\times 3\times 21.$
\item [] $5302=9+8\times (7+654)+3+2\times 1.$
\item [] $5303=98+(76+5)\times 4^3+21.$
\item [] $5304=9\times 87+6+5\times 43\times 21.$
\item [] $5305=9+(8+7)\times 65+4321.$
\item [] $5306=9+8\times (7+654)+3^2\times 1.$
\item [] $5307=987+6\times 5\times (4\times 3)^2\times 1.$
\item [] $5308=987+6\times 5\times (4\times 3)^2+1.$
\item [] $5309=(98\times 7+6)\times 5+43^2\times 1.$
\item [] $5310=(98\times 7+6)\times 5+43^2+1.$
\item[]$\mbox{Increasing order}$
\item [] $5311=123+4+(5+67)\times 8\times 9.$
\item [] $5312=1\times 2^{(3+4)}+(5+67)\times 8\times 9.$
\item [] $5313=1^2\times (3\times 4+56)\times 78+9.$
\item [] $5314=1+23\times (4+5\times 6\times 7+8+9).$
\item [] $5315=1\times 2\times (34+5)\times 67+89.$
\item [] $5316=12\times (34+56\times 7+8+9).$
\item [] $5317=(1+2^3+4)\times (56\times 7+8+9).$
\item [] $5318=12+(3+4)\times (56+78\times 9).$
\item [] $5319=123\times (4+5)+6\times 78\times 9.$
\item [] $5320=(1^2+34)\times (56+7+89).$
\item [] $5321=123\times 4+5+67\times 8\times 9.$
\item [] $5322=(1+2\times 3\times 4)\times 5\times 6\times 7+8\times 9.$
\item [] $\mathit{5323=-1\times 2+3^4\times 56+789.}$
\item [] $\mathit{5324=1-2+3^4\times 56+789.}$
\item [] $5325=1^2\times 3^4\times 56+789.$
\item [] $5326=1^2+3^4\times 56+789.$
\item [] $5327=1\times 2+3^4\times 56+789.$
\item [] $5328=1+2+3^4\times 56+789.$
\item [] $5329=1+(2^3\times 4+5)\times 6\times (7+8+9).$
\item [] $5330=1\times 23\times (4\times 56+7)+8+9.$
\item [] $5331=1+23\times (4\times 56+7)+8+9.$
\item [] $5332=12^3+4+(56\times 7+8)\times 9.$
\item [] $\mathit{5333=1-23+4^5\times 6-789.}$
\item [] $5334=123+(4+567+8)\times 9 .$
\item [] $5335=1+2\times (3+(4\times (5+67)+8)\times 9).$
\item [] $5336=1\times 2^3\times (4\times (5+6)+7\times 89).$
\item [] $5337=12+3^4\times 56+789.$
\item [] $5338=1+2\times 3^4+(567+8)\times 9.$
\item [] $5339=1\times 2+3^4\times 5\times (6+7)+8\times 9.$
\item [] $5340=(1+2^3+4+5+6\times 7)\times 89.$
\item [] $5341=1\times (2+3)\times 4^5+(6+7)\times (8+9).$
\item [] $5342=1\times 2+3+(45\times (6+7)+8)\times 9.$
\item [] $5343=(12+3^4)\times 56+(7+8)\times 9.$
\item [] $5344=1+2+(3^4+5)\times (6+7\times 8)+9.$
\item [] $5345=12\times 3\times 4\times (5\times 6+7)+8+9.$
\item [] $5346=(123+456+7+8)\times 9.$
\item [] $5347=1+2\times 3^4+(5+67)\times 8\times 9.$
\item [] $5348=1\times 2+3\times (4\times 5\times 6+78)\times 9.$
\item [] $5349=12\times 3^4+56\times 78+9.$
\item [] $5350=(1\times 2+3^4)\times 56+78\times 9.$
\item [] $5351=1+(2+3^4)\times 56+78\times 9.$
\item [] $5352=123+(45+67\times 8)\times 9.$
\item [] $5353=12+(3^4+5)\times (6+7\times 8)+9.$
\item [] $5354=1^2\times 3^4\times 5\times (6+7)+89.$
\item [] $5355=1^2+3^4\times 5\times (6+7)+89.$
\item [] $5356=1^2+(34+5+6)\times 7\times (8+9).$
\item [] $5357=1+2+3^4\times 5\times (6+7)+89.$
\item [] $5358=12\times 3\times 45+6\times 7\times 89.$
\item [] $5359=123+4^5+6\times 78\times 9.$
\item [] $5360=1\times 23+(45\times (6+7)+8)\times 9.$
\item [] $5361=123\times 4+(5+67\times 8)\times 9.$
\item [] $\mathit{5362=1+2\times 3\times 4^5+6-789.}$
\item [] $5363=(1\times 2+3^4\times 5)\times (6+7)+8\times 9.$
\item [] $5364=1\times 2\times 3^4\times 5\times 6+7\times 8\times 9.$
\item [] $5365=1+2\times 3^4\times 5\times 6+7\times 8\times 9.$
\item [] $5366=12+3^4\times 5\times (6+7)+89.$
\item [] $5367=12+(34+5+6)\times 7\times (8+9).$
\item [] $\mathit{5368=-1+2345+6\times 7\times 8\times 9.}$
\item [] $5369=1\times 2345+6\times 7\times 8\times 9.$
\item [] $5370=1+2345+6\times 7\times 8\times 9.$
\item [] $5371=(1\times 2+34+5)\times (6\times 7+89).$
\item [] $5372=12^3+4+56\times (7\times 8+9).$
\item [] $5373=12\times 3+(45\times (6+7)+8)\times 9.$
\item [] $5374=1\times 2^{(3+4)}\times 5+6\times 789.$
\item [] $5375=1+2^{(3+4)}\times 5+6\times 789.$
\item [] $5376=12\times (3+4)\times (5+6\times 7+8+9).$
\item [] $5377=1+(2+3+45+6)\times (7+89).$
\item [] $5378=1\times 2\times (3+(4+5\times 6)\times (7+8\times 9)).$
\item [] $5379=1^2\times 3+4\times 56\times (7+8+9).$
\item [] $5380=1^{23}\times 4+56\times (7+89).$
\item[]$\mbox{Decreasing order}$
\item [] $5311=9\times 8\times (7+6)+5^4\times (3\times 2+1).$
\item [] $5312=9\times 8\times (7+65)+4\times 32\times 1.$
\item [] $5313=98\times 7\times 6+(54+3)\times 21.$
\item [] $5314=9\times (8+7+6)+5\times (4^(3+2)+1).$
\item [] $5315=9+(87\times 6\times 5+43)\times 2\times 1.$
\item [] $5316=9\times 8+76\times (5+43+21).$
\item [] $5317=(9\times (87+6+5)+4)\times 3\times 2+1.$
\item [] $5318=9\times 8+76\times (5+4^3)+2\times 1.$
\item [] $5319=987+6+5+4321.$
\item [] $5320=9+(876+5+4)\times 3\times 2+1.$
\item [] $5321=9+8\times (7+654)+3+21.$
\item [] $5322=9\times 8+7\times 6\times 5\times (4\times 3\times 2+1).$
\item [] $5323=9+8\times (7+654+3)+2\times 1.$
\item [] $5324=(9+87\times 6\times 5+43)\times 2\times 1.$
\item [] $5325=9\times (8+7)\times 6+5\times 43\times 21.$
\item [] $\mathit{5326=9+87\times (65-4)+3^2+1.}$
\item [] $5327=98\times 7+(6+5\times 43)\times 21.$
\item [] $5328=9\times 87\times 6+5^4+3+2\times 1.$
\item [] $5329=9\times 87\times 6+5^4+3+2+1.$
\item [] $5330=9\times 87\times 6+5^4+3\times 2+1.$
\item [] $5331=9+(87+6)\times (54+3)+21.$
\item [] $5332=9\times 87\times 6+5^4+3^2\times 1.$
\item [] $5333=9\times 87\times 6+5^4+3^2+1.$
\item [] $5334=987+(65+4)\times 3\times 21.$
\item [] $5335=((9+87\times 6)\times 5+4\times 3)\times 2+1.$
\item [] $5336=(98+7+6+5)\times (43+2+1).$
\item [] $5337=9+8+76\times 5\times (4+3)\times 2\times 1.$
\item [] $5338=987+6\times 5+4321.$
\item [] $5339=9\times 87\times 6+5\times 4^3\times 2+1.$
\item [] $5340=(9\times 8\times 7+6\times 5)\times (4+3\times 2)\times 1.$
\item [] $5341=9\times (8+76\times 5)+43^2\times 1.$
\item [] $5342=987+65\times (4+3\times 21).$
\item [] $5343=(9+8+76)\times 54+321.$
\item [] $5344=987+(6+5\times 4\times 3)^2+1.$
\item [] $5345=98+76\times (5+4^3)+2+1.$
\item [] $5346=9\times (87\times 6+5+4+3\times 21).$
\item [] $5347=9\times 87\times 6+5^4+3+21.$
\item [] $5348=98+7\times 6\times 5\times (4\times 3\times 2+1).$
\item [] $5349=(98+7+6)\times (5+43)+21.$
\item [] $5350=9+8\times (7\times 6+5^4)+3+2\times 1.$
\item [] $5351=(9+(8+7)\times 6)\times 54+3+2\times 1.$
\item [] $5352=9+(87+6)\times 54+321.$
\item [] $5353=9+8\times (7+654+3\times 2+1).$
\item [] $5354=9+8\times (7+6\times 54+3)\times 2+1.$
\item [] $5355=9\times 87\times 6+5^4+32\times 1.$
\item [] $5356=9\times 87\times 6+5^4+32+1.$
\item [] $5357=9+(87\times 6\times 5+4^3)\times 2\times 1.$
\item [] $5358=9\times 87\times 6+5\times 4\times (32+1).$
\item [] $\mathit{5359=9-87+6+5432-1.}$
\item [] $5360=9+8\times (7+654)+3\times 21.$
\item [] $5361=9\times 87+654\times (3\times 2+1).$
\item [] $5362=987+6\times ((5+4)\times 3)^2+1.$
\item [] $5363=987+6\times (5+4)^3+2\times 1.$
\item [] $5364=987+6\times (5+4)^3+2+1.$
\item [] $5365=(9+876+5+4)\times 3\times 2+1.$
\item [] $5366=(9+87\times 6\times 5+4^3)\times 2^1.$
\item [] $5367=(9+87\times 6\times 5+4^3)\times 2+1.$
\item [] $5368=987+6+5^4\times (3\times 2+1).$
\item [] $5369=9+8\times (7+654+3^2\times 1).$
\item [] $5370=(9+(8+7)\times 6)\times 54+3+21.$
\item [] $5371=9+8\times (7\times 6+5^4+3)+2\times 1.$
\item [] $5372=98\times 7\times 6+(5^4+3)\times 2\times 1.$
\item [] $5373=987+65+4321.$
\item [] $5374=9+8+765\times (4+3)+2^1.$
\item [] $5375=9+8+765\times (4+3)+2+1.$
\item [] $5376=98\times 7\times 6+5\times 4\times 3\times 21.$
\item [] $5377=(9+87)\times (6+5)+4321.$
\item [] $5378=9+8\times 76+(5+4^3)^2\times 1.$
\item [] $5379=9+8\times 76+(5+4^3)^2+1.$
\item [] $5380=(9\times 8\times 7+6\times 5+4)\times (3^2+1).$
\item[]$\mbox{Increasing order}$
\item [] $5381=(1^2+3^4)\times 56+789.$
\item [] $5382=1\times 234\times 5+6\times 78\times 9.$
\item [] $5383=1+234\times 5+6\times 78\times 9.$
\item [] $5384=1^2+3+4+56\times (7+89).$
\item [] $5385=(12\times 3+4+56)\times 7\times 8+9.$
\item [] $5386=1\times 2\times (34\times 56+789).$
\item [] $5387=1+2\times (34\times 56+789).$
\item [] $5388=1\times 2^3+4+56\times (7+89).$
\item [] $5389=1+2^3+4+56\times (7+89).$
\item [] $5390=1\times 2+3\times 4+56\times (7+89).$
\item [] $5391=12+3+4\times 56\times (7+8+9).$
\item [] $5392=1\times 2\times (3+4+5\times 67\times 8+9).$
\item [] $5393=1+2\times (3+4+5\times 67\times 8+9).$
\item [] $5394=(12\times 3+4\times 5+6)\times (78+9).$
\item [] $5395=12+3+4+56\times (7+89).$
\item [] $5396=12\times 3^4+56\times (7+8\times 9).$
\item [] $5397=1+(2+3)\times 4+56\times (7+89).$
\item [] $\mathit{5398=-123+4^5\times 6-7\times 89.}$
\item [] $5399=1\times 23+4\times 56\times (7+8+9).$
\item [] $5400=12\times (3+45)+67\times 8\times 9.$
\item [] $5401=1+2\times 3\times 4+56\times (7+89).$
\item [] $5402=1\times 23\times (4\times 56+7)+89.$
\item [] $5403=1\times 23+4+56\times (7+89).$
\item [] $5404=1+23+4+56\times (7+89).$
\item [] $5405=12\times 3^4\times 5+67\times 8+9.$
\item [] $5406=12\times (3+4)\times 56+78\times 9.$
\item [] $5407=1+(2^3+45)\times (6+7+89).$
\item [] $5408=1\times 2^3\times 4+56\times (7+89).$
\item [] $5409=1\times 234+(567+8)\times 9.$
\item [] $5410=1+234+(567+8)\times 9.$
\item [] $5411=1^2+34+56\times (7+89).$
\item [] $5412=1\times 2+34+56\times (7+89).$
\item [] $5413=1+2+34+56\times (7+89).$
\item [] $5414=1\times 2\times (34\times (5+6)\times 7+89).$
\item [] $5415=(1\times 2+3)\times (4^5+6\times 7+8+9).$
\item [] $5416=12\times 3+4+56\times (7+89).$
\item [] $5417=1234+(5+6\times 7)\times 89.$
\item [] $5418=(1\times 23+4+567+8)\times 9.$
\item [] $5419=1+234+(5+67)\times 8\times 9.$
\item [] $5420=1+(2+3)\times (4^5+6\times 7)+89.$
\item [] $5421=12+3\times 4\times 5\times 6\times (7+8)+9.$
\item [] $5422=12+34+56\times (7+89).$
\item [] $5423=1\times 2+(34+5)\times (67+8\times 9).$
\item [] $5424=1\times 2\times 345+6\times 789.$
\item [] $5425=1+2\times 345+6\times 789.$
\item [] $5426=1\times 2+(3+45)\times ((6+7)\times 8+9).$
\item [] $5427=(1+23+4+567+8)\times 9.$
\item [] $5428=(12+34)\times (5+(6+7)\times 8+9).$
\item [] $5429=(1^2\times 3+45+6+7)\times 89.$
\item [] $5430=(1+2\times 3)^4+5+6\times 7\times 8\times 9.$
\item [] $5431=(1^2+3)^4+(567+8)\times 9.$
\item [] $5432=(1^2+3+4)\times (56+7\times 89).$
\item [] $5433=12+(34+5)\times (67+8\times 9).$
\item [] $5434=1\times 2+3+((4+5)\times 6+7)\times 89.$
\item [] $5435=1+2+3+((4+5)\times 6+7)\times 89.$
\item [] $5436=123\times 4\times (5+6)+7+8+9.$
\item [] $5437=1\times (2+3^4)\times 56+789.$
\item [] $5438=1+(2+3^4)\times 56+789.$
\item [] $5439=1+2+3^4\times (5+6+7\times 8)+9.$
\item [] $5440=1\times 2\times 34\times (56+7+8+9).$
\item [] $5441=1+2\times 34\times (56+7+8+9).$
\item [] $5442=1\times 2\times (3+(4\times 56+78)\times 9).$
\item [] $5443=1+2\times (3+(4\times 56+78)\times 9).$
\item [] $5444=1\times 2\times 34+56\times (7+89).$
\item [] $5445=(1\times 2+3\times 45+6\times 78)\times 9.$
\item [] $5446=1\times 2\times (34+5\times 67\times 8+9).$
\item [] $5447=1+2\times (34+5\times 67\times 8+9).$
\item [] $5448=12\times (34+5\times 6)\times 7+8\times 9.$
\item [] $5449=(1^2+3\times 4\times 56+7)\times 8+9.$
\item [] $5450=1+2\times 34\times (5+67+8)+9.$
\item[]$\mbox{Decreasing order}$
\item [] $\mathit{5381=987\times 6-543+2\times 1.}$
\item [] $5382=987+6\times (5+4)^3+21.$
\item [] $5383=9\times (8+76+5\times 43)\times 2+1.$
\item [] $\mathit{5384=(9+8)\times 76-5+4^(3\times 2)+1.}$
\item [] $5385=9+(876+5\times 4)\times 3\times 2\times 1.$
\item [] $5386=9+(876+5\times 4)\times 3\times 2+1.$
\item [] $5387=(9\times 8+7)\times 65+4\times 3\times 21.$
\item [] $\mathit{5388=-9+8-7\times 6+5432-1.}$
\item [] $5389=(9+8)\times (7\times (6+5+4)\times 3+2\times 1).$
\item [] $5390=98\times (7+6+5+4+32+1).$
\item [] $5391=9+(87+6\times 5)\times (43+2+1).$
\item [] $5392=9\times 8+76\times 5\times (4+3)\times 2\times 1.$
\item [] $5393=9+8+765\times (4+3)+21.$
\item [] $5394=9+8\times 7+(6\times 5+43)^2\times 1.$
\item [] $5395=9+8\times 7+(6\times 5+43)^2+1.$
\item [] $5396=((9+87\times 6)\times 5+43)\times 2\times 1.$
\item [] $5397=(987+6)\times 5+432\times 1.$
\item [] $5398=(987+6)\times 5+432+1.$
\item [] $\mathit{5399=9\times 8\times (76+5)-432-1.}$
\item [] $5400=9+876+5\times 43\times 21.$
\item [] $5401=(9+8+7+6)\times 5\times 4\times 3^2+1.$
\item [] $5402=9\times (8+7\times 6)\times (5+4+3)+2\times 1.$
\item [] $5403=9+(8\times 7+6)\times (54+32+1).$
\item [] $5404=987+(6+5^4)\times (3\times 2+1).$
\item [] $5405=98\times 7\times 6+5+4\times 321.$
\item [] $5406=9+(8+7+6)\times (5+4\times 3\times 21).$
\item [] $5407=(9\times 87+6\times 5\times 4^3)\times 2+1.$
\item [] $5408=9+8\times (7\times 6+5^4)+3\times 21.$
\item [] $5409=9+(8+7)\times 6\times (54+3+2+1).$
\item [] $5410=9\times (8\times 7+65)+4321.$
\item [] $5411=9+(8+7)\times 6\times 5\times 4\times 3+2\times 1.$
\item [] $5412=98\times 7\times 6+54\times (3+21).$
\item [] $5413=9+876\times 5+4^(3+2)\times 1.$
\item [] $5414=9+876\times 5+4^(3+2)+1.$
\item [] $5415=9+((8+7\times 6)\times 54+3)\times 2\times 1.$
\item [] $5416=9+((8+7\times 6)\times 54+3)\times 2+1.$
\item [] $5417=9+8\times (7+6)\times (5\times 4+32\times 1).$
\item [] $5418=98+76\times 5\times (4+3^2+1).$
\item [] $5419=98+76\times 5\times (4+3)\times 2+1.$
\item [] $5420=(9+87+6\times 5)\times 43+2\times 1.$
\item [] $5421=(9+87+6\times 5)\times 43+2+1.$
\item [] $5422=9+(8+765)\times (4+3)+2\times 1.$
\item [] $5423=9+(8+765)\times (4+3)+2+1.$
\item [] $5424=9\times 87+(6+5\times 43)\times 21.$
\item [] $5425=9+8\times (7\times 6+5^4+3^2+1).$
\item [] $5426=(9+8)\times (7+6)\times 5+4321.$
\item [] $5427=(9\times 8+76\times 5)\times 4\times 3+2+1.$
\item [] $5428=(9+8+7\times 6)\times (5+43\times 2+1).$
\item [] $5429=9\times 8+765\times (4+3)+2\times 1.$
\item [] $5430=9+(8+7)\times 6\times 5\times 4\times 3+21.$
\item [] $5431=(9+876+5\times 4)\times 3\times 2+1.$
\item [] $5432=(9\times (8+7+6)+5)\times (4+3+21).$
\item [] $5433=9\times 87\times 6+5\times (4+3)\times 21.$
\item [] $5434=9+8\times (7+6+5^4)+321.$
\item [] $5435=98+7+(6\times 5+43)^2+1.$
\item [] $5436=9\times 8\times (7+65)+4\times 3\times 21.$
\item [] $5437=(98\times (7+6)+5)\times 4+321.$
\item [] $5438=98\times 7+(6+5)\times 432\times 1.$
\item [] $5439=98\times 7+(6+5)\times 432+1.$
\item [] $5440=(9+8+7\times 6+5)\times (4^3+21).$
\item [] $5441=9+(8+765)\times (4+3)+21.$
\item [] $5442=9+8\times 7\times (6+5+43\times 2)+1.$
\item [] $\mathit{5443=-9+8+7+6+5432-1.}$
\item [] $5444=98+(76+5)\times (4^3+2)\times 1.$
\item [] $5445=(9\times 8+76\times 5)\times 4\times 3+21.$
\item [] $5446=9\times (87+6\times 5+4)\times (3+2)+1.$
\item [] $5447=9\times (8\times 7+6+543)+2\times 1.$
\item [] $5448=9\times 8+765\times (4+3)+21.$
\item [] $5449=(9+8+7\times 6\times 5)\times 4\times 3\times 2+1.$
\item [] $5450=(9+876)\times 5+4^(3+2)+1.$
\item[]$\mbox{Increasing order}$
\item [] $5451=1234+5+6\times 78\times 9.$
\item [] $5452=1+23\times (4\times (5\times 6+7)+89).$
\item [] $5453=12+(3\times 4\times 56+7)\times 8+9.$
\item [] $5454=12\times 3\times 4\times 5+6\times 789.$
\item [] $5455=(1+23)\times 4\times 56+7+8\times 9.$
\item [] $5456=(1+2)^3+((4+5)\times 6+7)\times 89.$
\item [] $5457=1^2\times 3^4+56\times (7+89).$
\item [] $5458=1^2+3^4+56\times (7+89).$
\item [] $5459=1\times 2+3^4+56\times (7+89).$
\item [] $5460=12\times (3+4)+56\times (7+89).$
\item [] $5461=1+2\times (3\times 4+5\times 6)\times (7\times 8+9).$
\item [] $5462=1\times 2+(3+4)\times 5\times (67+89).$
\item [] $5463=(1+23)\times 4\times 56+78+9.$
\item [] $5464=1\times 2^3\times (4+56+7\times 89).$
\item [] $5465=12\times (34+5\times 6)\times 7+89.$
\item [] $5466=1+(2+3)\times 4^5+6\times 7\times 8+9.$
\item [] $5467=1\times (2+3^4)\times 5\times (6+7)+8\times 9.$
\item [] $5468=1\times 23\times 4+56\times (7+89).$
\item [] $5469=12+3^4+56\times (7+89).$
\item [] $5470=1+(2+3\times 4+56)\times 78+9.$
\item [] $\mathit{5471=1\times 2\times (3+4)\times 56\times 7-8-9.}$
\item [] $5472=(1^2+3\times 4+56+7)\times 8\times 9.$
\item [] $5473=1+(2+34)\times (56+7+89).$
\item [] $5474=1^2\times 34\times (5+67+89).$
\item [] $5475=12^3+4+5+6\times 7\times 89.$
\item [] $5476=1^2+3+(4+5+67)\times 8\times 9.$
\item [] $5477=123\times 4\times (5+6)+7\times 8+9.$
\item [] $5478=(1+2)\times 34+56\times (7+89).$
\item [] $5479=1\times 2\times (3+4\times (5+678))+9.$
\item [] $5480=1\times 2^3+(4+5+67)\times 8\times 9.$
\item [] $5481=1^2\times (34+567+8)\times 9.$
\item [] $5482=1+2\times (3+4+5\times 67)\times 8+9.$
\item [] $5483=1\times 2\times 3^4\times 5\times 6+7\times 89.$
\item [] $5484=1+2\times 3^4\times 5\times 6+7\times 89.$
\item [] $5485=1+(2+3^4)\times 5\times (6+7)+89.$
\item [] $5486=12+34\times (5+67+89).$
\item [] $5487=12+3+(4+5+67)\times 8\times 9.$
\item [] $5488=1\times 2\times (3+4\times (5+678)+9).$
\item [] $5489=12\times 3^4\times 5+6+7\times 89.$
\item [] $5490=1\times 2\times (3+4\times 56+78)\times 9.$
\item [] $5491=1234+(5+6\times 78)\times 9.$
\item [] $5492=1+23\times (4+5\times 6)\times 7+8+9.$
\item [] $5493=12\times (3+4)\times 56+789.$
\item [] $5494=1\times (2+3)^4+(5+67\times 8)\times 9.$
\item [] $5495=1\times 23+(4+5+67)\times 8\times 9.$
\item [] $5496=12\times (34+5\times 67+89).$
\item [] $5497=1\times (2+3)^4+56\times (78+9).$
\item [] $5498=1+(2+3)^4+56\times (78+9).$
\item [] $5499=(12+3)\times 45+67\times 8\times 9.$
\item [] $5500=1\times (2\times 3+4)\times (5+67\times 8+9).$
\item [] $5501=1+(2\times 3+4)\times (5+67\times 8+9).$
\item [] $5502=(1+2+34+5)\times (6\times 7+89).$
\item [] $5503=123+4+56\times (7+89).$
\item [] $5504=1\times 2^{(3+4)}+56\times (7+89).$
\item [] $5505=1\times 2^3\times (4+5+678)+9.$
\item [] $5506=1+2\times (3+4)\times 56\times 7+8+9.$
\item [] $5507=1+23\times (4\times 56+7+8)+9.$
\item [] $5508=(1+2+34+567+8)\times 9 .$
\item [] $5509=1+2\times 3\times (4+5)\times (6+7+89).$
\item [] $5510=(12+3^4\times 5)\times (6+7)+89.$
\item [] $5511=12^3+45+6\times 7\times 89.$
\item [] $5512=12^3+4+5\times (6+78)\times 9.$
\item [] $5513=(1+2+3)^4+5+6\times 78\times 9.$
\item [] $5514=1\times 2\times 345+67\times 8\times 9.$
\item [] $5515=1+2\times 345+67\times 8\times 9.$
\item [] $5516=(12\times 3+45)\times 67+89.$
\item [] $5517=12\times 3^4+567\times 8+9.$
\item [] $5518=(1+2\times 3\times 4+5\times 6+7)\times 89.$
\item [] $5519=1^{23}+(4\times 5+6\times 7)\times 89.$
\item [] $5520=12\times 3\times 4+56\times (7+89).$
\item[]$\mbox{Decreasing order}$
\item [] $5451=9\times (8+7\times 65)+4\times 321.$
\item [] $5452=9+(8+7\times (6+5))\times 4^3+2+1.$
\item [] $5453=98+(76+5+4)\times 3\times 21.$
\item [] $5454=98\times 7+6+(5+4^3)^2+1.$
\item [] $5455=98+765\times (4+3)+2\times 1.$
\item [] $5456=98+765\times (4+3)+2+1.$
\item [] $5457=(9\times 8+7)\times (65+4)+3+2+1.$
\item [] $5458=(9\times 8+7)\times (65+4)+3\times 2+1.$
\item [] $5459=9+(8+7\times 6\times 5)\times (4\times 3\times 2+1).$
\item [] $5460=(9\times 8+7)\times 65+4+321.$
\item [] $5461=(9\times 8+7)\times (65+4)+3^2+1.$
\item [] $5462=9+8+7+6+5432^1.$
\item [] $5463=9+8+7+6+5432+1.$
\item [] $5464=9\times (8+7)+(6\times 5+43)^2\times 1.$
\item [] $5465=98\times 7\times 6+5+4^3\times 21.$
\item [] $5466=9+8\times 7\times 6+5\times 4^(3+2)+1.$
\item [] $5467=((98+7)\times (6+5\times 4)+3)\times 2+1.$
\item [] $5468=9\times 8+76\times (5+4^3+2)\times 1.$
\item [] $5469=9\times (87+6)\times 5+4\times 321.$
\item [] $5470=9+8\times 7\times 6+5\times (4^(3+2)+1).$
\item [] $5471=9\times (8\times 7+6)+(5+4\times 3)^(2+1).$
\item [] $5472=9\times (8\times 7+6+543+2+1).$
\item [] $5473=9+8\times (7\times 6+5\times 4\times 32+1).$
\item [] $5474=98+765\times (4+3)+21.$
\item [] $5475=98\times (7+6\times 5)+43^2\times 1.$
\item [] $5476=98\times (7+6\times 5)+43^2+1.$
\item [] $5477=9\times (87+65)\times 4+3+2\times 1.$
\item [] $5478=9\times (87+65)\times 4+3\times 2\times 1.$
\item [] $5479=9+8+7\times 65\times 4\times 3+2\times 1.$
\item [] $5480=9+8+7\times 65\times 4\times 3+2+1.$
\item [] $5481=9\times (87\times 6+54+32+1).$
\item [] $5482=9+(87+65)\times 4\times 3^2+1.$
\item [] $5483=(9\times 8+7)\times (65+4)+32\times 1.$
\item [] $5484=(9\times 8+7)\times (65+4)+32+1.$
\item [] $5485=9+(8\times 7+6+5+4+3)^2\times 1.$
\item [] $5486=9\times 8\times 76+5+4+3+2\times 1.$
\item [] $5487=9\times 8\times 76+5+4+3\times 2\times 1.$
\item [] $5488=9\times 8\times 76+5+4+3\times 2+1.$
\item [] $5489=9+8\times (7+654+3+21).$
\item [] $5490=9\times 8\times 76+5+4+3^2\times 1.$
\item [] $5491=9+8+7\times 6+5432\times 1.$
\item [] $5492=9+8+7\times 6+5432+1.$
\item [] $5493=9+87\times (6+54+3)+2+1.$
\item [] $5494=98+76\times (5+4^3+2\times 1).$
\item [] $5495=98+76\times (5+4^3+2)+1.$
\item [] $5496=9\times 8\times (7+65+4)+3+21.$
\item [] $5497=9\times 8\times 76+5\times 4+3+2\times 1.$
\item [] $5498=9\times 8\times 76+5\times 4+3+2+1.$
\item [] $5499=9\times 8\times 76+5\times 4+3\times 2+1.$
\item [] $5500=((9+8)\times 7+6)\times (5\times 4+3+21).$
\item [] $5501=9\times 8\times 76+5\times 4+3^2\times 1.$
\item [] $5502=9\times 8\times 76+5+4\times 3\times 2+1.$
\item [] $5503=9+8\times 7+6+5432\times 1.$
\item [] $5504=9+8\times 7+6+5432+1.$
\item [] $5505=9\times 8\times 76+5+4+3+21.$
\item [] $5506=9\times 8\times 76+(5+4\times 3)\times 2\times 1.$
\item [] $5507=9\times 8\times 76+(5+4\times 3)\times 2+1.$
\item [] $5508=987+6+5\times 43\times 21.$
\item [] $5509=9\times 8\times (7+65)+4+321.$
\item [] $5510=9\times 8\times 76+5+4\times 3+21.$
\item [] $5511=9\times (8+7\times 65)+4^3\times 21.$
\item [] $5512=9+(876+5^4\times 3)\times 2+1.$
\item [] $5513=(9+87+6)\times 54+3+2\times 1.$
\item [] $5514=9\times 8\times 76+5+4+32+1.$
\item [] $5515=(9+87+6)\times 54+3\times 2+1.$
\item [] $5516=9\times 8\times 76+5\times 4+3+21.$
\item [] $5517=9\times 8+7+6+5432\times 1.$
\item [] $5518=9\times 8+7+6+5432+1.$
\item [] $5519=9\times 8\times 76+5\times 4+3^{(2+1)}.$
\item [] $5520=9\times 8\times 76+(5+4)\times 3+21.$
\item[]$\mbox{Increasing order}$
\item [] $5521=1^2\times 3+(4\times 5+6\times 7)\times 89.$
\item [] $5522=1^2+3+(4\times 5+6\times 7)\times 89.$
\item [] $5523=(1^{2345}+6)\times 789.$
\item [] $5524=1+2+3+(4\times 5+6\times 7)\times 89.$
\item [] $5525=1+2\times 3+(4\times 5+6\times 7)\times 89.$
\item [] $5526=1\times 2^3+(4\times 5+6\times 7)\times 89.$
\item [] $5527=(1\times 2+3)\times (4^5+67)+8\times 9.$
\item [] $5528=1+(2+3)\times (4^5+67)+8\times 9.$
\item [] $5529=(1+2\times 34)\times (5+67+8)+9.$
\item [] $5530=1\times (2+3\times 4+56)\times (7+8\times 9).$
\item [] $5531=123\times 4\times (5+6)+7\times (8+9).$
\item [] $5532=12\times (3+456)+7+8+9.$
\item [] $5533=12+3+(4\times 5+6\times 7)\times 89.$
\item [] $\mathit{5534=-1-2-3-4+(5+6)\times 7\times 8\times 9.}$
\item [] $5535=(12+3\times 45+6\times 78)\times 9.$
\item [] $5536=(1^2+3+4)\times (5+678+9).$
\item [] $5537=(12+3\times 4\times 56+7)\times 8+9.$
\item [] $5538=1\times 2\times 3^4+56\times (7+89).$
\item [] $5539=1+2\times 3^4+56\times (7+89).$
\item [] $5540=1\times 2\times (3^4+5\times 67\times 8+9).$
\item [] $5541=1\times 23+(4\times 5+6\times 7)\times 89.$
\item [] $5542=1+23+(4\times 5+6\times 7)\times 89.$
\item [] $5543=1+2\times (3+4\times (5+678+9)).$
\item [] $5544=12\times 3\times 4\times 5+67\times 8\times 9.$
\item [] $5545=1+2\times 3^4\times 5+6\times 789.$
\item [] $5546=1\times 23\times (4+5\times 6)\times 7+8\times 9.$
\item [] $5547=12\times 3^4\times 5+678+9.$
\item [] $5548=1^{23}\times 4+(5+6)\times 7\times 8\times 9.$
\item [] $5549=1^{23}+4+(5+6)\times 7\times 8\times 9.$
\item [] $5550=(1+2+34)\times 5\times (6+7+8+9).$
\item [] $5551=1^2\times 3+4+(5+6)\times 7\times 8\times 9.$
\item [] $5552=1^2+3+4+(5+6)\times 7\times 8\times 9.$
\item [] $5553=(12\times 3+45+67\times 8)\times 9.$
\item [] $5554=12\times 3+(4\times 5+6\times 7)\times 89.$
\item [] $5555=1\times 2+3\times (4\times 56+7)\times 8+9.$
\item [] $5556=1+2+3\times (4\times 56+7)\times 8+9.$
\item [] $5557=1+2\times 3\times (4\times 56+78\times 9).$
\item [] $5558=1\times 2+3\times 4+(5+6)\times 7\times 8\times 9.$
\item [] $5559=1+2+3\times 4+(5+6)\times 7\times 8\times 9.$
\item [] $5560=1\times (2+3\times 4)\times 56\times 7+8\times 9.$
\item [] $5561=(123+4+567)\times 8+9.$
\item [] $5562=1\times 2\times 3^4\times 5\times 6+78\times 9.$
\item [] $5563=1+2\times 3^4\times 5\times 6+78\times 9.$
\item [] $5564=1+23\times (4+5\times 6)\times 7+89.$
\item [] $5565=123\times 45+6+7+8+9.$
\item [] $5566=(12+34)\times (56+7\times 8+9).$
\item [] $5567=(1+(2+3\times 4)\times 56)\times 7+8\times 9.$
\item [] $5568=(12\times 3^4)\times 5+6+78\times 9.$
\item [] $5569=1+2^3\times (4+5+678+9).$
\item [] $5570=1\times 2+3\times 4\times (56\times 7+8\times 9).$
\item [] $5571=1\times 23+4+(5+6)\times 7\times 8\times 9.$
\item [] $5572=1+23+4+(5+6)\times 7\times 8\times 9.$
\item [] $5573=12\times (3+456)+7\times 8+9.$
\item [] $5574=12\times 3^4\times 5+6\times 7\times (8+9).$
\item [] $5575=1+2\times 3\times (4\times 5\times 6\times 7+89).$
\item [] $5576=1\times 2^3\times 4+(5+6)\times 7\times 8\times 9.$
\item [] $5577=1\times 2\times (3+4)\times 56\times 7+89.$
\item [] $5578=1+2\times (3+4)\times 56\times 7+89.$
\item [] $5579=1^2+34+(5+6)\times 7\times 8\times 9.$
\item [] $5580=(12+3+45)\times (6+78+9).$
\item [] $5581=1^2+3\times 4\times 5\times (6+78+9).$
\item [] $5582=1\times 2+3\times 4\times 5\times (6+78+9).$
\item [] $5583=1+2+3\times 4\times 5\times (6+78+9).$
\item [] $5584=12\times 3+4+(5+6)\times 7\times 8\times 9.$
\item [] $5585=(1\times 2+3)\times (4^5+6+78+9).$
\item [] $5586=1+(2+3)\times (4^5+6+78+9).$
\item [] $5587=12\times (3+456)+7+8\times 9.$
\item [] $5588=1\times 2^3+(4+(5+6)\times 7\times 8)\times 9.$
\item [] $5589=(12+34+567+8)\times 9 .$
\item [] $5590=12+34+(5+6)\times 7\times 8\times 9.$
\item[]$\mbox{Decreasing order}$
\item [] $5521=(9+876+5^4\times 3)\times 2+1.$
\item [] $5522=9\times 8\times 76+5+43+2\times 1.$
\item [] $5523=9\times 8\times 76+5+43+2+1.$
\item [] $5524=9\times 8\times 76+5\times 4+32\times 1.$
\item [] $5525=9+8+76+5432\times 1.$
\item [] $5526=9+8+76+5432+1.$
\item [] $5527=9\times 8\times 76+5+(4+3)^2+1.$
\item [] $5528=9\times 8\times 76+5\times (4+3)+21.$
\item [] $5529=9\times (87+6)\times 5+4^3\times 21.$
\item [] $5530=(9\times 8+7)\times (6+54+3^2+1).$
\item [] $5531=9\times 8\times 76+54+3+2\times 1.$
\item [] $5532=9\times 8\times 76+54+3+2+1.$
\item [] $5533=9\times 8\times 76+54+3\times 2+1.$
\item [] $5534=9+87+6+5432\times 1.$
\item [] $5535=9+87+6+5432+1.$
\item [] $5536=9\times 87+(6+5)\times 432+1.$
\item [] $5537=98+7\times (6\times 5+4+3)\times 21.$
\item [] $5538=9\times 8\times 76+5\times (4+3^2)+1.$
\item [] $\mathit{5539=9\times 8\times 76+5+4^3-2\times 1.}$
\item [] $5540=(9+87+6)\times 54+32\times 1.$
\item [] $5541=9\times 8\times 76+5+43+21.$
\item [] $5542=9\times 8\times 76+5\times (4+3)\times 2\times 1.$
\item [] $5543=98+7+6+5432\times 1.$
\item [] $5544=98+7+6+5432+1.$
\item [] $5545=9+8\times (7+6)+5432\times 1.$
\item [] $5546=9+8\times (7+6)+5432+1.$
\item [] $5547=9\times 8+7\times 6+5432+1.$
\item [] $5548=9\times (8\times 76+5)+4+3^{(2+1)}.$
\item [] $5549=9+8+7+65\times (4^3+21).$
\item [] $5550=9\times 8\times 76+54+3+21.$
\item [] $5551=9\times 87+6+(5+4^3)^2+1.$
\item [] $\mathit{5552=9\times 8\times 76-5+43\times 2-1.}$
\item [] $5553=9\times 8\times 76+5\times 4\times 3+21.$
\item [] $5554=9\times 8\times 76+(5+4)\times 3^2+1.$
\item [] $5555=9\times 8\times 76+5\times 4+3\times 21.$
\item [] $5556=9\times (8+7)\times (6+5\times (4+3))+21.$
\item [] $5557=(9+8)\times 7+6+5432\times 1.$
\item [] $5558=9\times 8\times 76+54+32\times 1.$
\item [] $5559=9\times 8\times 76+54+32+1.$
\item [] $5560=98+7\times 65\times 4\times 3+2\times 1.$
\item [] $5561=98+7\times 65\times 4\times 3+2+1.$
\item [] $5562=9\times 8\times 76+5+4^3+21.$
\item [] $5563=9\times 8\times 76+5+43\times 2\times 1.$
\item [] $5564=9\times 8\times 76+5+43\times 2+1.$
\item [] $5565=987+654\times (3\times 2+1).$
\item [] $5566=9\times (8\times 76+5)+(4+3)^2\times 1.$
\item [] $5567=9+(8+7\times 65)\times 4\times 3+2\times 1.$
\item [] $5568=(9\times 8+7)\times 65+432+1.$
\item [] $5569=(98+7+65+4)\times 32+1.$
\item [] $5570=98+76\times (5+4+3\times 21).$
\item [] $5571=(9+87+6)\times 54+3\times 21.$
\item [] $5572=98+7\times 6+5432\times 1.$
\item [] $5573=98+7\times 6+5432+1.$
\item [] $5574=98\times 7\times 6+54\times 3^{(2+1)}.$
\item [] $5575=98\times 7\times 6+(5+4)^3\times 2+1.$
\item [] $5576=(9+8+7)\times 6+5432\times 1.$
\item [] $5577=(9+8+7)\times 6+5432+1.$
\item [] $5578=(9+(8+76)\times 5)\times (4+3^2)+1.$
\item [] $5579=98+7\times 65\times 4\times 3+21.$
\item [] $5580=9\times 8+76+5432\times 1.$
\item [] $5581=9\times 8+76+5432+1.$
\item [] $5582=(9+8+76)\times 5\times 4\times 3+2\times 1.$
\item [] $5583=9\times 8\times (7+6)\times 5+43\times 21.$
\item [] $5584=9\times (8\times 76+5)+4+3\times 21.$
\item [] $5585=9\times (8\times 76+5+4)+32\times 1.$
\item [] $5586=98\times (7+6+5\times 4+3+21).$
\item [] $5587=9\times 8\times 76+(54+3)\times 2+1.$
\item [] $5588=9+8\times 7+(65\times 4+3)\times 21.$
\item [] $5589=9\times 8\times 76+54+3\times 21.$
\item [] $5590=9+(876+54)\times 3\times 2+1.$
\item[]$\mbox{Increasing order}$
\item [] $\mathit{5591=1+2-3^4+5678-9.}$
\item [] $5592=12+3\times 4\times 5\times (6+78+9).$
\item [] $5593=1\times 2^3\times (4\times 5+678)+9.$
\item [] $5594=123\times 45+6\times 7+8+9.$
\item [] $5595=12\times (3+456)+78+9.$
\item [] $5596=1+2\times (3+45\times (6+7\times 8))+9.$
\item [] $5597=1\times (2+3)\times 4^5+6\times 78+9.$
\item [] $5598=1+(2+3)\times 4^5+6\times 78+9.$
\item [] $5599=1+(2+34\times (5+6+7)+8)\times 9.$
\item [] $5600=1^2\times (3+4)\times (5+6+789).$
\item [] $5601=1^2+(3+4)\times (5+6+789).$
\item [] $5602=1\times 2+(3+4)\times (5+6+789).$
\item [] $5603=1+2+(3+4)\times (5+6+789).$
\item [] $5604=12\times (3+456)+7+89.$
\item [] $5605=(12+3+4)\times 5\times (6\times 7+8+9).$
\item [] $5606=123\times 45+6+7\times 8+9.$
\item [] $5607=1^{234}\times (56+7)\times 89.$
\item [] $5608=1^2+(3+4)\times (5+6+78)\times 9.$
\item [] $5609=12\times (3+456+7)+8+9.$
\item [] $5610=1\times 234+56\times (7+89).$
\item [] $5611=1234+56\times 78+9.$
\item [] $5612=12+(3+4)\times (5+6+789).$
\item [] $5613=1+2\times 34+(5+6)\times 7\times 8\times 9.$
\item [] $5614=1^2\times 3+4+(56+7)\times 89.$
\item [] $5615=1\times (2+3)\times (4^5+6\times (7+8)+9).$
\item [] $5616=(123+4+5)\times 6\times 7+8\times 9.$
\item [] $5617=1+2\times 3\times (4\times 5+6+78)\times 9.$
\item [] $5618=1\times 2+(34+5)\times 6\times (7+8+9).$
\item [] $5619=123\times 45+67+8+9.$
\item [] $5620=123\times 45+6+7+8\times 9.$
\item [] $5621=1\times 2+34\times (5+6)\times (7+8)+9.$
\item [] $5622=12+34\times (5\times 6+(7+8)\times 9).$
\item [] $5623=(1^2+3)\times 4+(56+7)\times 89.$
\item [] $5624=(1+2+34)\times (56+7+89).$
\item [] $5625=(123\times 4+5\times 6\times 7)\times 8+9.$
\item [] $5626=1^2+(3+4+5)\times 6\times 78+9.$
\item [] $5627=1\times 2+(3+4+5)\times 6\times 78+9.$
\item [] $5628=123\times 45+6+78+9.$
\item [] $5629=1^2+3+45\times (6+7\times (8+9)).$
\item [] $5630=(1\times 2+3)\times 4^5+6+7\times 8\times 9.$
\item [] $5631=1+(2+3)\times 4^5+6+7\times 8\times 9.$
\item [] $5632=(1^2+3)^4+56\times (7+89).$
\item [] $5633=(123+4+5)\times 6\times 7+89.$
\item [] $5634=1\times 2\times 3^4\times 5+67\times 8\times 9.$
\item [] $5635=1+2\times 3^4\times 5+67\times 8\times 9.$
\item [] $5636=1\times 23\times 4+(5+6)\times 7\times 8\times 9.$
\item [] $5637=123\times 45+6+7+89.$
\item [] $5638=(1+2)^3+4+(56+7)\times 89.$
\item [] $5639=(1+2\times 3^4)\times 5+67\times 8\times 9.$
\item [] $5640=12\times (345+6+7\times (8+9)).$
\item [] $5641=123+(4\times 5+6\times 7)\times 89.$
\item [] $5642=(1\times 2+3)\times 4^5+6\times (78+9).$
\item [] $5643=12\times (3+456)+(7+8)\times 9.$
\item [] $5644=(1+2\times 3\times (4+5)\times 6+7)\times (8+9).$
\item [] $5645=(1+2^{(3+4)}+5)\times 6\times 7+8+9.$
\item [] $5646=(1+2)\times 34+(5+6)\times 7\times 8\times 9.$
\item [] $5647=1+2\times (3+4\times 5\times (6+(7+8)\times 9)).$
\item [] $5648=123\times 45+(6+7)\times 8+9.$
\item [] $5649=123\times 45+6\times 7+8\times 9.$
\item [] $5650=1+2\times 3^4\times 5\times 6+789.$
\item [] $5651=1+(2+3+45)\times ((6+7)\times 8+9).$
\item [] $5652=12\times 3^4+5\times (6+7)\times 8\times 9.$
\item [] $5653=12+34+(56+7)\times 89.$
\item [] $5654=(1\times 2+3)\times (4^5+6)+7\times 8\times 9.$
\item [] $5655=12\times 3^4\times 5+6+789.$
\item [] $5656=1\times 23\times 4\times 56+7\times 8\times 9.$
\item [] $5657=1+23\times 4\times 56+7\times 8\times 9.$
\item [] $5658=1234+56\times (7+8\times 9).$
\item [] $5659=1+(2^{(3+4)}+5)\times 6\times 7+8\times 9.$
\item [] $5660=123\times 45+6+7\times (8+9).$
\item[]$\mbox{Decreasing order}$
\item [] $5591=9+(87+6)\times 5\times 4\times 3+2\times 1.$
\item [] $5592=9\times 8\times 76+5\times 4\times 3\times 2\times 1.$
\item [] $5593=9\times 8\times 76+5\times 4\times 3\times 2+1.$
\item [] $5594=(9+8)\times 7\times (6+5+4+32)+1.$
\item [] $5595=9\times 8+7\times (65\times 4+3)\times (2+1).$
\item [] $5596=(9+8)\times (7\times 6+5)\times (4+3)+2+1.$
\item [] $5597=9\times 8\times 76+5\times (4\times 3\times 2+1).$
\item [] $5598=9\times (8\times 76+5+4+3+2\times 1).$
\item [] $5599=(9+8\times 76+5)\times (4+3+2)+1.$
\item [] $5600=98\times (7+6)+5+4321.$
\item [] $5601=(9+8+76)\times 5\times 4\times 3+21.$
\item [] $5602=9\times 8+7+(65\times 4+3)\times 21.$
\item [] $5603=9\times (8\times 76+5)+43\times 2\times 1.$
\item [] $5604=9\times (8\times 76+5)+43\times 2+1.$
\item [] $5605=9\times 8\times 76+5+4\times 32\times 1.$
\item [] $5606=9\times 87\times 6+5+43\times 21.$
\item [] $5607=98+76+5432+1.$
\item [] $5608=9\times 8\times 7\times (6+5)+43+21.$
\item [] $5609=9\times 8\times 76+5+4\times (32+1).$
\item [] $5610=9\times (87+6\times 5\times 4)\times 3+21.$
\item [] $5611=9\times 8\times 7\times (6+5)+4+3\times 21.$
\item [] $5612=9\times 8\times 76+5\times (4+3+21).$
\item [] $5613=9\times (8\times 76+5)+4\times (3+21).$
\item [] $5614=9+8+7+65\times 43\times 2\times 1.$
\item [] $5615=9+8+7+65\times 43\times 2+1.$
\item [] $5616=9\times 8\times (7+65)+432\times 1.$
\item [] $5617=9\times 8\times (7+65)+432+1.$
\item [] $5618=9+8\times (7+654)+321.$
\item [] $5619=9+87+(65\times 4+3)\times 21.$
\item [] $5620=9+8+(7+6)\times (5\times 43\times 2+1).$
\item [] $5621=9\times (8+7+6)+5432\times 1.$
\item [] $5622=9\times (8+7+6)+5432+1.$
\item [] $5623=9+8\times (76+5^4)+3\times 2\times 1.$
\item [] $5624=9\times 8\times 76+5+(4+3)\times 21.$
\item [] $5625=9+8\times (7+6)\times (5+4)\times 3\times 2\times 1.$
\item [] $5626=9+8\times (76+5^4)+3^2\times 1.$
\item [] $5627=9+8\times (76+5^4)+3^2+1.$
\item [] $5628=987+(6+5\times 43)\times 21.$
\item [] $5629=9+(8+7+65\times 43)\times 2\times 1.$
\item [] $5630=9+8\times (7+6)\times 54+3+2\times 1.$
\item [] $5631=9+8\times (7+6)\times 54+3\times 2\times 1.$
\item [] $5632=9+8\times (7+6)\times 54+3\times 2+1.$
\item [] $5633=9+(87+65)\times (4+32+1).$
\item [] $5634=(9+876+54)\times 3\times 2\times 1.$
\item [] $5635=(9+876+54)\times 3\times 2+1.$
\item [] $5636=9\times 8\times 76+54\times 3+2\times 1.$
\item [] $5637=9\times 8\times 76+54\times 3+2+1.$
\item [] $5638=(9+8+7+65\times 43)\times 2\times 1.$
\item [] $5639=(9+8+7+65\times 43)\times 2+1.$
\item [] $5640=9\times 8\times 7\times (6+5)+4\times (3+21).$
\item [] $5641=9+8\times (76+5^4)+3+21.$
\item [] $5642=(9+8)\times 7+(65\times 4+3)\times 21.$
\item [] $5643=9+8\times (76+5^4+3)+2\times 1.$
\item [] $5644=9\times 87\times 6+5^4+321.$
\item [] $5645=9\times (8\times 76+5)+4^3\times 2\times 1.$
\item [] $5646=9\times (87\times 6+5)+43\times 21.$
\item [] $5647=9+876+(5+4^3)^2+1.$
\item [] $\mathit{5648=-9+87\times 65-4+3\times 2\times 1.}$
\item [] $5649=9+8\times (76\times 5+4+321).$
\item [] $5650=9+8\times (76+5^4)+32+1.$
\item [] $5651=(9+8)\times (7+6\times 54)+3+21.$
\item [] $5652=9\times 8\times 76+5\times 4\times 3^2\times 1.$
\item [] $5653=9\times 8\times 76+5\times 4\times 3^2+1.$
\item [] $5654=(9+8)\times (7+6)+5432+1.$
\item [] $5655=9\times 8\times 76+54\times 3+21.$
\item [] $5656=9+8\times 7+65\times 43\times 2+1.$
\item [] $5657=9+8\times (7+6)\times 54+32\times 1.$
\item [] $5658=9+8\times (7+6)\times 54+32+1.$
\item [] $5659=9\times 8+7\times 6\times (5+4^3\times 2)+1.$
\item [] $5660=(9+8)\times (7+6\times 54)+32+1.$
\item[]$\mbox{Increasing order}$
\item [] $5661=12\times (3\times 45+6\times 7\times 8)+9.$
\item [] $5662=1+(2+3+4)\times (5\times 6+7)\times (8+9).$
\item [] $\mathit{5663=-1+2-3-4+5678-9.}$
\item [] $5664=12\times (3+456+7)+8\times 9.$
\item [] $5665=1\times (2+3)\times 4^5+67\times 8+9.$
\item [] $5666=123\times 45+6\times 7+89.$
\item [] $5667=123\times 4+(567+8)\times 9.$
\item [] $5668=1+2\times 3+(4+5)\times (6+7\times 89).$
\item [] $5669=1\times 2^3+(4+5)\times (6+7\times 89).$
\item [] $5670=12\times 3^4\times 5+6\times (7+8)\times 9.$
\item [] $5671=123+4+(5+6)\times 7\times 8\times 9.$
\item [] $5672=1\times 2+(3\times 4+5\times 6)\times (7+8)\times 9.$
\item [] $5673=1\times (2\times 3)^4+56\times 78+9.$
\item [] $5674=123\times 45+67+8\times 9.$
\item [] $5675=(1+2\times 3\times 4)\times (5\times 6\times 7+8+9).$
\item [] $5676=123\times 45+6+(7+8)\times 9.$
\item [] $\mathit{5677=-1+2+3+4+5678-9.}$
\item [] $5678=1\times (2+3)\times 4^5+(6+7\times 8)\times 9.$
\item [] $5679=(1+2\times 3^4)\times 5\times 6+789.$
\item [] $5680=1+2\times 3\times 45\times (6+7+8)+9.$
\item [] $5681=12\times (3+456+7)+89.$
\item [] $5682=12+(3\times 4+5\times 6)\times (7+8)\times 9.$
\item [] $\mathit{5683=1+2-3-4+5678+9.}$
\item [] $5684=1\times 23+(4+5)\times (6+7\times 89).$
\item [] $5685=1+23+(4+5)\times (6+7\times 89).$
\item [] $5686=1+(2^3+4)\times (5+6\times 78)+9.$
\item [] $5687=1^{234}\times 5678+9.$
\item [] $5688=1^{234}+5678+9.$
\item [] $5689=1\times 2+3\times 45\times 6\times 7+8+9.$
\item [] $5690=1+2+3\times 45\times 6\times 7+8+9.$
\item [] $5691=123\times 45+67+89.$
\item [] $5692=1^{23}+4+5678+9.$
\item [] $5693=12\times (3^4+56\times 7)+8+9.$
\item [] $5694=1^2\times 3+4+5678+9.$
\item [] $5695=1^2+3+4+5678+9.$
\item [] $5696=1\times 2+3+4+5678+9.$
\item [] $5697=1+2+3+4+5678+9.$
\item [] $5698=1+2\times 3+4+5678+9.$
\item [] $5699=1^2\times 3\times 4+5678+9.$
\item [] $5700=1+2^3+4+5678+9.$
\item [] $5701=1\times 2+3\times 4+5678+9.$
\item [] $5702=1+2+3\times 4+5678+9.$
\item [] $5703=(1^2+3)\times 4+5678+9.$
\item [] $5704=1\times 23\times (4\times 56+7+8+9).$
\item [] $5705=1+23\times (4\times 56+7+8+9).$
\item [] $5706=12+3+4+5678+9.$
\item [] $5707=12^3+4+5\times (6+789).$
\item [] $5708=1+(2+3)\times 4+5678+9.$
\item [] $5709=1\times 2+3\times 4+5\times 67\times (8+9).$
\item [] $5710=1+2+3\times 4+5\times 67\times (8+9).$
\item [] $5711=12+3\times 4+5678+9.$
\item [] $5712=1+2\times 3\times 4+5678+9.$
\item [] $5713=1+(23+45)\times (67+8+9).$
\item [] $5714=1\times 23+4+5678+9.$
\item [] $5715=1+23+4+5678+9.$
\item [] $5716=1+(2+3)\times 4+5\times 67\times (8+9).$
\item [] $5717=(1+2^{(3+4)}+5)\times 6\times 7+89.$
\item [] $5718=(1+2)^3+4+5678+9.$
\item [] $5719=1\times 2^3\times 4+5678+9.$
\item [] $5720=1+2^3\times 4+5678+9.$
\item [] $5721=1^2\times 34+5678+9.$
\item [] $5722=1^2+34+5678+9.$
\item [] $5723=1\times 2+34+5678+9.$
\item [] $5724=1+2+34+5678+9.$
\item [] $5725=1^2+(3+4+5)\times (6\times 78+9).$
\item [] $5726=1\times 2+(3+4+5)\times (6\times 78+9).$
\item [] $5727=12\times 3+4+5678+9.$
\item [] $5728=1+2^3\times 4+5\times 67\times (8+9).$
\item [] $5729=12\times (3\times 4+56)\times 7+8+9.$
\item [] $5730=1^2+34+5\times 67\times (8+9).$
\item[]$\mbox{Decreasing order}$
\item [] $5661=9+(8\times 76+5\times 4)\times 3^2\times 1.$
\item [] $5662=9\times 8+(7+6)\times 5\times 43\times 2\times 1.$
\item [] $5663=9\times 8+(7+6)\times 5\times 43\times 2+1.$
\item [] $5664=9+87\times (6+54+3+2\times 1).$
\item [] $5665=9+87\times (6+54+3+2)+1.$
\item [] $5666=9+8\times (7\times 6+5^4)+321.$
\item [] $5667=(9+8+7\times 65)\times 4\times 3+2+1.$
\item [] $\mathit{5668=9\times 8+7+65\times 43\times 2-1.}$
\item [] $5669=9\times 8+7+65\times 43\times 2\times 1.$
\item [] $5670=9\times 8+7+65\times 43\times 2+1.$
\item [] $5671=(9+87+6\times 5)\times (43+2)+1.$
\item [] $5672=9\times 8\times 7\times (6+5)+4\times 32\times 1.$
\item [] $5673=9+876\times 5+4\times 321.$
\item [] $5674=9+87\times 65+4+3+2+1.$
\item [] $5675=9+87\times 65+4+3\times 2+1.$
\item [] $5676=(98\times 7+65\times 4)\times 3\times 2\times 1.$
\item [] $5677=9+87\times 65+4+3^2\times 1.$
\item [] $5678=9+87\times 65+4+3^2+1.$
\item [] $5679=9+87\times 65+4\times 3+2+1.$
\item [] $5680=9+8\times (76+5^4)+3\times 21.$
\item [] $5681=9+8\times (7+6\times (54+3\times 21)).$
\item [] $5682=9+(87+6)\times (54+3\times 2+1).$
\item [] $5683=9\times (8\times 7+6)+5\times (4^(3+2)+1).$
\item [] $5684=9+8+7\times 6+5^4\times 3^2\times 1.$
\item [] $5685=9+8+7\times 6+5^4\times 3^2+1.$
\item [] $5686=9+87+65\times 43\times 2\times 1.$
\item [] $5687=9+87+65\times 43\times 2+1.$
\item [] $5688=9+87\times 65+4\times 3\times 2\times 1.$
\item [] $5689=9\times 8\times 76+5\times 43+2\times 1.$
\item [] $5690=9\times 8\times 76+5\times 43+2+1.$
\item [] $5691=(98+7)\times (6+5+43)+21.$
\item [] $5692=9+87\times 65+4+3+21.$
\item [] $5693=(9\times 87+6\times 5)\times (4+3)+2\times 1.$
\item [] $5694=(98+7)\times 6\times (5+4)+3+21.$
\item [] $5695=98+7+65\times 43\times 2\times 1.$
\item [] $5696=98+7+65\times 43\times 2+1.$
\item [] $5697=9+87\times 65+4\times 3+21.$
\item [] $5698=9\times 8\times 76+5\times (43+2)+1.$
\item [] $5699=(9+8+7+65)\times 4^3+2+1.$
\item [] $5700=9+87\times 65+4+32\times 1.$
\item [] $5701=9+87\times 65+4+32+1.$
\item [] $5702=9\times 8\times 76+5\times (43+2+1).$
\item [] $5703=(98+7)\times 6\times (5+4)+32+1.$
\item [] $5704=9+87\times 65+4\times (3^2+1).$
\item [] $5705=98+7\times (65\times 4\times 3+21).$
\item [] $5706=9\times 87\times 6+(5+43)\times 21.$
\item [] $5707=9+(8+76+5)\times 4^3+2\times 1.$
\item [] $5708=9\times 8\times 76+5\times 43+21.$
\item [] $5709=9+87\times 65+43+2\times 1.$
\item [] $5710=9+87\times 65+43+2+1.$
\item [] $5711=9\times 8+7+6+5^4\times 3^2+1.$
\item [] $5712=9+(8\times 7+65\times 43)\times 2+1.$
\item [] $5713=9+87\times 65+(4+3)^2\times 1.$
\item [] $5714=9+87\times 65+(4+3)^2+1.$
\item [] $5715=9\times (8\times 76+5\times 4)+3\times 21.$
\item [] $5716=((98+7)\times 6+5)\times (4+3+2)+1.$
\item [] $5717=9+8+76\times 5\times (4\times 3+2+1).$
\item [] $5718=9+8+76+5^4\times 3^2\times 1.$
\item [] $5719=9+8+76+5^4\times 3^2+1.$
\item [] $5720=9+8\times 7+65\times (43\times 2+1).$
\item [] $5721=9+8\times 7\times (65+4+32+1).$
\item [] $5722=9+87+(6+5+4^3)^2+1.$
\item [] $5723=9+8\times 7\times (6\times 5+4)\times 3+2\times 1.$
\item [] $5724=9+8\times 7\times 6\times (5+4\times 3)+2+1.$
\item [] $5725=9\times (8+7)+65\times 43\times 2\times 1.$
\item [] $5726=98+7\times 6\times (5+4\times 32+1).$
\item [] $5727=9+87+6+5^4\times 3^2\times 1.$
\item [] $5728=9+87\times 65+43+21.$
\item [] $5729=9\times 8\times 76+5+4\times 3\times 21.$
\item [] $5730=9+87\times 65+4^3+2\times 1.$
\item[]$\mbox{Increasing order}$
\item [] $5731=1\times 2+34+5\times 67\times (8+9).$
\item [] $5732=1+2+34+5\times 67\times (8+9).$
\item [] $5733=12+34+5678+9.$
\item [] $5734=1^2+(3+4)\times (5\times 6+789).$
\item [] $5735=12\times 3+4+5\times 67\times (8+9).$
\item [] $5736=1+2+(3+4)\times (5\times 6+789).$
\item [] $5737=1\times 2\times (345+6+7)\times 8+9.$
\item [] $5738=1+2\times (345+6+7)\times 8+9.$
\item [] $5739=1^2\times 3\times ((4+5\times 6)\times 7\times 8+9).$
\item [] $5740=1\times (2+3)\times (4+5+67\times (8+9)).$
\item [] $5741=12+34+5\times 67\times (8+9).$
\item [] $5742=1^2\times 3\times 45\times 6\times 7+8\times 9.$
\item [] $5743=1^2+3\times 45\times 6\times 7+8\times 9.$
\item [] $5744=1\times 2+3\times 45\times 6\times 7+8\times 9.$
\item [] $5745=1+2+3\times 45\times 6\times 7+8\times 9.$
\item [] $5746=(1+2\times 3\times 45+67)\times (8+9).$
\item [] $5747=(12+3)\times 4+5678+9.$
\item [] $5748=12\times (3^4+56\times 7)+8\times 9.$
\item [] $5749=1\times (2+3)\times 4^5+6+7\times 89.$
\item [] $5750=(1+2)\times (34\times 56+7)+8+9.$
\item [] $5751=(12\times 3+45)\times (6+7\times 8+9).$
\item [] $5752=1+(2+3+4)\times (567+8\times 9).$
\item [] $5753=1+2^3\times (4+(5+6)\times (7\times 8+9)).$
\item [] $5754=12+3\times 45\times 6\times 7+8\times 9.$
\item [] $5755=1\times 2\times 34+5678+9.$
\item [] $5756=1+2\times 34+5678+9.$
\item [] $5757=1+(2+3\times 45\times 6)\times 7+8\times 9.$
\item [] $5758=1^{23}\times 4^5+6\times 789.$
\item [] $5759=1^{23}+4^5+6\times 789.$
\item [] $5760=1^2+3\times 45\times 6\times 7+89.$
\item [] $5761=1^2\times 3+4^5+6\times 789.$
\item [] $5762=1+2+3\times 45\times 6\times 7+89.$
\item [] $5763=1\times 2+3+4^5+6\times 789.$
\item [] $5764=1+2+3+4^5+6\times 789.$
\item [] $5765=1+2\times 3+4^5+6\times 789.$
\item [] $5766=1\times 2^3+4^5+6\times 789.$
\item [] $5767=1+2^3+4^5+6\times 789.$
\item [] $5768=1^2\times 3^4+5678+9.$
\item [] $5769=1^2+3^4+5678+9.$
\item [] $5770=1+23\times 45+6\times 789.$
\item [] $5771=12+3\times 45\times 6\times 7+89.$
\item [] $5772=12+3^4\times (56+7+8)+9.$
\item [] $5773=12+3+4^5+6\times 789.$
\item [] $5774=1+(2+3\times 45\times 6)\times 7+89.$
\item [] $5775=1\times 23\times 4\times 56+7\times 89.$
\item [] $5776=1+23\times 4\times 56+7\times 89.$
\item [] $5777=1^2+3^4+5\times 67\times (8+9).$
\item [] $5778=12\times 34\times 5+6\times 7\times 89.$
\item [] $5779=1234+567\times 8+9.$
\item [] $5780=12+3^4+5678+9.$
\item [] $5781=1\times 23+4^5+6\times 789.$
\item [] $5782=1+23+4^5+6\times 789.$
\item [] $5783=(1+23)\times 4+5678+9.$
\item [] $5784=12\times 34+56\times (7+89).$
\item [] $5785=1^2\times (3+4\times 5+6\times 7)\times 89.$
\item [] $5786=1^{234}+5\times (6+7)\times 89.$
\item [] $5787=12^3+45\times 6\times (7+8)+9.$
\item [] $5788=12+3^4+5\times 67\times (8+9).$
\item [] $5789=(1+2)\times 34+5678+9.$
\item [] $5790=1^{23}+4+5\times (6+7)\times 89.$
\item [] $5791=12+(3^4+5)\times 67+8+9.$
\item [] $5792=1^2\times 3+4+5\times (6+7)\times 89.$
\item [] $5793=1^2+3+4+5\times (6+7)\times 89.$
\item [] $5794=12\times 3+4^5+6\times 789.$
\item [] $5795=(1+2)^3\times 4+5678+9.$
\item [] $5796=12\times (3+456+7+8+9).$
\item [] $5797=1^2\times 3\times 4+5\times (6+7)\times 89.$
\item [] $5798=1+2^3+4+5\times (6+7)\times 89.$
\item [] $5799=(1+2)\times 34\times 56+78+9.$
\item [] $5800=1+2+3\times 4+5\times (6+7)\times 89.$
\item[]$\mbox{Decreasing order}$
\item [] $5731=9+87\times 65+4+3\times 21.$
\item [] $\mathit{5732=-9\times 87+6\times 543\times 2-1.}$
\item [] $5733=9+876\times 5+4^3\times 21.$
\item [] $5734=9\times 8+7+65\times (43\times 2+1).$
\item [] $5735=((98\times 7+6\times 5)\times 4+3)\times 2+1.$
\item [] $5736=98+7+6+5^4\times 3^2\times 1.$
\item [] $5737=98+7+6+5^4\times 3^2+1.$
\item [] $5738=9+8\times 76+5\times 4^(3+2)+1.$
\item [] $5739=9\times 8+7\times 6+5^4\times 3^2\times 1.$
\item [] $5740=987+(6+5)\times 432+1.$
\item [] $5741=98\times 7\times 6+5\times (4+321).$
\item [] $5742=9+8\times 7\times 6\times (5+4\times 3)+21.$
\item [] $5743=9\times 8\times 76+54\times (3+2)+1.$
\item [] $5744=9+8\times 7+(6+5^4)\times 3^2\times 1.$
\item [] $5745=98\times 7\times 6+543\times (2+1).$
\item [] $5746=(9+8)\times (7+6+54\times 3\times 2+1).$
\item [] $5747=(9+8)\times (7+6\times (5+4)\times 3)\times 2+1.$
\item [] $5748=(9\times 8+7+65\times 43)\times 2\times 1.$
\item [] $5749=9+87\times 65+4^3+21.$
\item [] $5750=9+87\times 65+43\times 2\times 1.$
\item [] $5751=9+87\times 65+43\times 2+1.$
\item [] $5752=9+87\times (6+54+3\times 2)+1.$
\item [] $5753=9+87\times (6+5\times 4\times 3)+2\times 1.$
\item [] $5754=(9+8+7\times 6+5\times 43)\times 21.$
\item [] $5755=9\times 8+7+6\times (5^4+321).$
\item [] $5756=(98+7\times 65\times 4)\times 3+2\times 1.$
\item [] $5757=(98+76\times 5)\times 4\times 3+21.$
\item [] $5758=9\times 8+7+(6+5^4)\times 3^2\times 1.$
\item [] $5759=9\times 8+7+(6+5^4)\times 3^2+1.$
\item [] $5760=98+7+65\times (43\times 2+1).$
\item [] $5761=9\times 8\times 76+(5+4)\times 32+1.$
\item [] $5762=(9\times 8+7+6+5)\times 4^3+2\times 1.$
\item [] $5763=9\times 8\times (7+6\times 5+43)+2+1.$
\item [] $\mathit{5764=98+7\times 6+5^4\times 3^2-1.}$
\item [] $5765=98+7\times 6+5^4\times 3^2\times 1.$
\item [] $5766=98+7\times 6+5^4\times 3^2+1.$
\item [] $5767=9\times (8+7)+6+5^4\times 3^2+1.$
\item [] $5768=9\times (87\times 6+5)+4^(3+2)+1.$
\item [] $5769=(9+876)\times 5+4^3\times 21.$
\item [] $5770=9+(8+7+6\times 5)\times 4\times 32+1.$
\item [] $5771=9+(8\times 7+6+5)\times 43\times 2\times 1.$
\item [] $5772=9+87\times (6+5\times 4\times 3)+21.$
\item [] $5773=9\times 8+76+5^4\times 3^2\times 1.$
\item [] $5774=9\times 8+76+5^4\times 3^2+1.$
\item [] $5775=(98+7\times 65\times 4)\times 3+21.$
\item [] $5776=9+87+(6+5^4)\times 3^2+1.$
\item [] $5777=9+8\times 7\times 6+5432\times 1.$
\item [] $5778=9+8\times 7\times 6+5432+1.$
\item [] $5779=98+(7+6)\times (5+432\times 1).$
\item [] $5780=(9+8)\times (7+6\times 54+3^2\times 1).$
\item [] $5781=(9\times 8+7+6+5)\times 4^3+21.$
\item [] $5782=(9+87+65\times 43)\times 2\times 1.$
\item [] $5783=(9+87+65\times 43)\times 2+1.$
\item [] $5784=9\times 87\times 6+543\times 2\times 1.$
\item [] $5785=9\times 87\times 6+543\times 2+1.$
\item [] $5786=9+8\times 7\times 6\times 5+4^{(3\times 2)}+1.$
\item [] $5787=(98+(7+6)\times 5\times 43)\times 2+1.$
\item [] $5788=(9\times 8+7)\times (6\times 5+43)+21.$
\item [] $5789=98\times (7+6)+5\times 43\times 21.$
\item [] $5790=9\times (8+7)+65\times (43\times 2+1).$
\item [] $5791=(9\times 8\times 7+654)\times (3+2)+1.$
\item [] $5792=9+87\times 65+4\times 32\times 1.$
\item [] $5793=9+87\times 65+4\times 32+1.$
\item [] $5794=9\times 8\times 76+5\times 4^3+2\times 1.$
\item [] $5795=9\times 8\times 76+5\times 4^3+2+1.$
\item [] $5796=9\times 8\times 76+54\times 3\times 2\times 1.$
\item [] $5797=9\times 8\times 76+54\times 3\times 2+1.$
\item [] $5798=98+76\times 5\times (4\times 3+2+1).$
\item [] $5799=98+76+5^4\times 3^2\times 1.$
\item [] $5800=98+76+5^4\times 3^2+1.$
\item[]$\mbox{Increasing order}$
\item [] $5801=(1^2+3\times 45)\times 6\times 7+89.$
\item [] $5802=1+(2\times 3^4\times 5+6)\times 7+89.$
\item [] $5803=(1+2)^3\times 4+5\times 67\times (8+9).$
\item [] $5804=12+3+4+5\times (6+7)\times 89.$
\item [] $5805=(1+23+45)\times (6+78)+9.$
\item [] $5806=1+(2+3)\times 4+5\times (6+7)\times 89.$
\item [] $5807=1\times (2+3)\times 4^5+678+9.$
\item [] $5808=(1+2)\times 34\times 56+7+89.$
\item [] $5809=12+3\times 4+5\times (6+7)\times 89.$
\item [] $5810=1+2\times 3\times 4+5\times (6+7)\times 89.$
\item [] $5811=1+(2\times 3+4)\times (5+6\times (7+89)).$
\item [] $5812=1\times 23+4+5\times (6+7)\times 89.$
\item [] $5813=1+23+4+5\times (6+7)\times 89.$
\item [] $5814=123+4+5678+9.$
\item [] $5815=1\times 2^{(3+4)}+5678+9.$
\item [] $5816=1+2^{(3+4)}+5678+9.$
\item [] $5817=1\times 2^3\times 4+5\times (6+7)\times 89.$
\item [] $5818=1+2^3\times 4+5\times (6+7)\times 89.$
\item [] $5819=1^2\times 34+5\times (6+7)\times 89.$
\item [] $5820=1^2+34+5\times (6+7)\times 89.$
\item [] $5821=1\times 2+34+5\times (6+7)\times 89.$
\item [] $5822=123+4+5\times 67\times (8+9).$
\item [] $5823=1\times 2^{(3+4)}+5\times 67\times (8+9).$
\item [] $5824=1\times 2^3\times (4\times 56+7\times 8\times 9).$
\item [] $5825=12\times 3+4+5\times (6+7)\times 89.$
\item [] $5826=(12+3\times 45\times 6)\times 7+8\times 9.$
\item [] $5827=1+(2+3\times 45)\times 6\times 7+8\times 9.$
\item [] $5828=1\times (2+3)\times 4^5+6+78\times 9.$
\item [] $5829=1+(2+3)\times 4^5+6+78\times 9.$
\item [] $5830=1+(23+4\times (5+6))\times (78+9).$
\item [] $5831=12\times 3\times 4+5678+9.$
\item [] $5832=12\times 3\times 45+6\times 78\times 9.$
\item [] $5833=1+2\times 3\times (45\times 6+78\times 9).$
\item [] $5834=1^2\times (3^4+5)\times 67+8\times 9.$
\item [] $5835=1^2+(3^4+5)\times 67+8\times 9.$
\item [] $5836=1\times 2+(3^4+5)\times 67+8\times 9.$
\item [] $5837=1+2+(3^4+5)\times 67+8\times 9.$
\item [] $5838=(1+2+34+5)\times (67+8\times 9).$
\item [] $5839=12\times 3\times 4+5\times 67\times (8+9).$
\item [] $5840=(1+2+3+4)\times (567+8+9).$
\item [] $5841=123\times (4+5)+6\times 789.$
\item [] $5842=1+(2\times 3)^4+567\times 8+9.$
\item [] $5843=(1\times 2+3\times 45)\times 6\times 7+89.$
\item [] $5844=12\times 3^4+56\times (78+9).$
\item [] $5845=(12+3)\times 4+5\times (6+7)\times 89.$
\item [] $5846=12+(3^4+5)\times 67+8\times 9.$
\item [] $5847=1+(23+45+6)\times (7+8\times 9).$
\item [] $5848=1^{23}\times 4^5+67\times 8\times 9.$
\item [] $5849=1^{23}+4^5+67\times 8\times 9.$
\item [] $5850=1+2\times 3^4+5678+9.$
\item [] $5851=1^2\times 3+4^5+67\times 8\times 9.$
\item [] $5852=1^2+3+4^5+67\times 8\times 9.$
\item [] $5853=1\times 2\times 34+5\times (6+7)\times 89.$
\item [] $5854=1\times 23\times 4\times 56+78\times 9.$
\item [] $5855=1+23\times 4\times 56+78\times 9.$
\item [] $5856=1\times 2^3+4^5+67\times 8\times 9.$
\item [] $5857=1+2^3+4^5+67\times 8\times 9.$
\item [] $5858=1+2\times 3^4+5\times 67\times (8+9).$
\item [] $5859=1\times 23\times 45+67\times 8\times 9.$
\item [] $5860=1+23\times 45+67\times 8\times 9.$
\item [] $5861=12+3^4\times (5+67)+8+9.$
\item [] $5862=12+(34+56)\times (7\times 8+9).$
\item [] $5863=12+3+4^5+67\times 8\times 9.$
\item [] $5864=1\times 2+3\times (4+5\times 6\times (7\times 8+9)).$
\item [] $5865=12\times (3+45+6+7)\times 8+9.$
\item [] $5866=1^2\times 3^4+5\times (6+7)\times 89.$
\item [] $5867=1^2+3^4+5\times (6+7)\times 89.$
\item [] $5868=123\times 4+56\times (7+89).$
\item [] $5869=1+23\times 4\times (56+7)+8\times 9.$
\item [] $5870=1\times 2+((3^4+5+6)\times 7+8)\times 9.$
\item[]$\mbox{Decreasing order}$
\item [] $5801=(98+7+65\times 43)\times 2+1.$
\item [] $5802=9\times 8\times 76+5+4+321.$
\item [] $5803=98\times (7\times 6+5+4\times 3)+21.$
\item [] $\mathit{5804=9+8\times 765-4-321.}$
\item [] $5805=9+(8+76)\times (5+43+21).$
\item [] $5806=9\times (8+7)\times (6+5)+4321.$
\item [] $5807=(9+8)\times 76+5\times 43\times 21.$
\item [] $5808=9+87\times 65+(4\times 3)^2\times 1.$
\item [] $5809=9+87\times 65+(4\times 3)^2+1.$
\item [] $5810=98+7\times (6\times 5+4)\times (3+21).$
\item [] $5811=9+87\times 65+(4+3)\times 21.$
\item [] $5812=98\times 7+6+5\times 4^(3+2)\times 1.$
\item [] $5813=9\times 8\times 76+5\times 4+321.$
\item [] $5814=9\times (8\times 76+5+4\times 3+21).$
\item [] $5815=9\times 8+(7+6+5^4)\times 3^2+1.$
\item [] $5816=(9+87+6)\times (54+3)+2\times 1.$
\item [] $5817=(9+87+6)\times (54+3)+2+1.$
\item [] $\mathit{5818=9\times 8\times (76+5)+4+3-21.}$
\item [] $5819=9\times (8+7+6+5^4)+3+2\times 1.$
\item [] $5820=(9+8\times 7\times 6+5^4)\times 3\times 2\times 1.$
\item [] $5821=9\times (8+7+6+5^4)+3\times 2+1.$
\item [] $5822=(9+8\times 7\times 6)\times 5+4^{(3\times 2)}+1.$
\item [] $5823=(9+8\times 7)\times 6+5432+1.$
\item [] $5824=9+(8+7+6+5^4)\times 3^2+1.$
\item [] $5825=9+8\times (7+6\times 5\times 4\times 3\times 2\times 1).$
\item [] $5826=(98+7+6\times 5)\times 43+21.$
\item [] $\mathit{5827=987\times 6-(5+43)\times 2+1.}$
\item [] $5828=((9+8)\times 7+65\times 43)\times 2\times 1.$
\item [] $5829=(9+87+6)\times 54+321.$
\item [] $5830=(9+8+7\times 6\times (5+4^3))\times 2\times 1.$
\item [] $5831=(9+8)\times (7+6+5+4+321).$
\item [] $5832=9\times (8+76+543+21).$
\item [] $5833=9+8\times (7+6\times 5\times 4\times 3\times 2+1).$
\item [] $5834=9\times 8\times 7+(6\times 5+43)^2+1.$
\item [] $5835=(9+87+6)\times (54+3)+21.$
\item [] $\mathit{5836=-98\times 7+6543-21.}$
\item [] $5837=9\times 8\times (7\times (6+5)+4)+3+2\times 1.$
\item [] $5838=9+87\times (6+54+3\times 2+1).$
\item [] $5839=9\times 8+7+6\times 5\times 4^3\times (2+1).$
\item [] $5840=98+(7+6+5^4)\times 3^2\times 1.$
\item [] $5841=98+(7+6+5^4)\times 3^2+1.$
\item [] $5842=9\times (8\times 76+5)+4+321.$
\item [] $5843=9\times 8\times (76+5)+4+3\times 2+1.$
\item [] $5844=98\times 7\times 6+54\times 32\times 1.$
\item [] $5845=98\times 7\times 6+54\times 32+1.$
\item [] $5846=9\times 8\times (76+5)+4+3^2+1.$
\item [] $5847=9\times 8\times 76+54+321.$
\item [] $5848=(9+8)\times (7+(6+54\times 3)\times 2+1).$
\item [] $5849=9+8\times (7\times 6+5^4+3\times 21).$
\item [] $5850=9\times 8\times 76+54\times (3\times 2+1).$
\item [] $5851=(9+8\times 7+65)\times (43+2)+1.$
\item [] $5852=9\times 8\times (76+5)+4\times (3+2)\times 1.$
\item [] $5853=9\times 8\times (76+5)+4\times (3+2)+1.$
\item [] $\mathit{5854=-98\times 7+6543-2-1.}$
\item [] $5855=9+8+7\times 6\times ((5+4^3)\times 2+1).$
\item [] $5856=(98+76+5+4)\times 32\times 1.$
\item [] $5857=(98+76+5+4)\times 32+1.$
\item [] $\mathit{5858=987\times 6-5+4-3\times 21.}$
\item [] $5859=(9+8+7+65+4)\times 3\times 21.$
\item [] $5860=9\times 8\times (76+5)+4+3+21.$
\item [] $5861=(9\times 8+765)\times (4+3)+2\times 1.$
\item [] $5862=(9\times 8+765)\times (4+3)+2+1.$
\item [] $5863=9\times 8\times (76+5)+4+3^{(2+1)}.$
\item [] $5864=9\times 8\times (7\times (6+5)+4)+32\times 1.$
\item [] $5865=(9\times 8)\times (76+5)+4\times 3+21.$
\item [] $5866=(9+8)\times ((7+6)\times 5+4)\times (3+2)+1.$
\item [] $5867=9+8+(7+6+5)\times (4+321).$
\item [] $5868=(987+6)\times 5+43\times 21.$
\item [] $5869=9\times 8\times 7\times (6+5)+4+321.$
\item [] $5870=98\times 7+(65+4+3)^2\times 1.$
\item[]$\mbox{Increasing order}$
\item [] $5871=1\times 23+4^5+67\times 8\times 9.$
\item [] $5872=1+23+4^5+67\times 8\times 9.$
\item [] $5873=(1+2\times 345+6\times 7)\times 8+9.$
\item [] $5874=(12+3+4+5+6\times 7)\times 89.$
\item [] $5875=1^2+3\times (4+5+6+7)\times 89.$
\item [] $5876=(1\times 2+3)\times 4^5+(6+78)\times 9.$
\item [] $5877=1\times 23\times 4+5\times (6+7)\times 89.$
\item [] $5878=1+23\times 4+5\times (6+7)\times 89.$
\item [] $\mathit{5879=1\times 2-3^4-5+67\times 89.}$
\item [] $5880=123\times 45+6\times 7\times 8+9.$
\item [] $5881=123+4^5+6\times 789.$
\item [] $5882=1+2\times (3\times 4\times 5\times 6+7)\times 8+9.$
\item [] $5883=123+(4+56)\times (7+89).$
\item [] $5884=12\times 3+4^5+67\times 8\times 9.$
\item [] $5885=1\times 23\times 4\times (56+7)+89.$
\item [] $5886=1+23\times 4\times (56+7)+89.$
\item [] $5887=(1+2)\times 34+5\times (6+7)\times 89.$
\item [] $5888=1\times 23\times 4\times (5+6\times 7+8+9).$
\item [] $5889=(1+2)\times (34\times 56+7\times 8)+9.$
\item [] $5890=1^2+(3+4)\times 56\times (7+8)+9.$
\item [] $5891=1\times 2+(3+4)\times 56\times (7+8)+9.$
\item [] $5892=1+2+(3+4)\times 56\times (7+8)+9.$
\item [] $5893=(1+2)^3\times 4+5\times (6+7)\times 89.$
\item [] $\mathit{5894=-1-23-45+67\times 89.}$
\item [] $5895=1^{23}\times 45\times (6\times 7+89).$
\item [] $5896=1^{23}+45\times (6\times 7+89).$
\item [] $5897=(1+23+4)\times 5\times 6\times 7+8+9.$
\item [] $5898=12+3\times (45\times 6\times 7+8\times 9).$
\item [] $5899=1^2+3+45\times (6\times 7+89).$
\item [] $5900=1\times 2+3+45\times (6\times 7+89).$
\item [] $5901=1+2+3+45\times (6\times 7+89).$
\item [] $5902=1+2\times 3+45\times (6\times 7+89).$
\item [] $5903=1\times 2^3+45\times (6\times 7+89).$
\item [] $5904=1\times 234\times 5+6\times 789.$
\item [] $5905=1+234\times 5+6\times 789.$
\item [] $5906=1\times 2+3^4\times (5+67)+8\times 9.$
\item [] $5907=1+2+3^4\times (5+67)+8\times 9.$
\item [] $\mathit{5908=1\times 23\times 4\times 5\times (6+7)-8\times 9.}$
\item [] $5909=(1+234)\times 5+6\times 789.$
\item [] $5910=(12+3^4)\times 56+78\times 9.$
\item [] $5911=1+2\times 3\times (4+(5+(6+7)\times 8)\times 9).$
\item [] $5912=123+4+5\times (6+7)\times 89.$
\item [] $5913=12^3+45\times (6+78+9).$
\item [] $5914=1234+5\times (6+7)\times 8\times 9.$
\item [] $5915=(1^2\times 3+4)\times (56+789).$
\item [] $5916=123\times 4\times (5+6)+7\times 8\times 9.$
\item [] $5917=1\times 2+(3+4)\times (56+789).$
\item [] $5918=1\times 23+45\times (6\times 7+89).$
\item [] $5919=1+23+45\times (6\times 7+89).$
\item [] $5920=(1\times 2+3)\times (45+67\times (8+9)).$
\item [] $5921=1\times 234+5678+9.$
\item [] $5922=1+234+5678+9.$
\item [] $5923=1\times 2+3^4\times (5+67)+89.$
\item [] $5924=1+2+3^4\times (5+67)+89.$
\item [] $5925=123\times 45+6\times (7\times 8+9).$
\item [] $5926=1+(2+3)\times (4+5+6)\times (7+8\times 9).$
\item [] $5927=12+(3+4)\times (56+789).$
\item [] $5928=12+(3\times 4+56)\times (78+9).$
\item [] $5929=12\times 3\times 4+5\times (6+7)\times 89.$
\item [] $5930=1+234+5\times 67\times (8+9).$
\item [] $5931=123\times (4+5)+67\times 8\times 9.$
\item [] $5932=1\times 2\times ((3^4\times 5+6)\times 7+89).$
\item [] $5933=12+3^4\times (5+67)+89.$
\item [] $5934=12+3^4\times (5\times (6+7)+8)+9.$
\item [] $\mathit{5935=1+2-34\times 5+678\times 9.}$
\item [] $5936=1\times 2\times (3+4)\times (5\times 67+89).$
\item [] $5937=(12+34+5\times 6)\times 78+9.$
\item [] $5938=1^2+3\times (45\times 6\times 7+89).$
\item [] $5939=(1\times 2+3)\times (4^5+6)+789.$
\item [] $5940=12\times 3\times (4+5+67+89).$
\item[]$\mbox{Decreasing order}$
\item [] $5871=98\times 7+(65+4+3)^2+1.$
\item [] $5872=9\times 8\times (76+5)+4\times (3^2+1).$
\item [] $5873=9+8\times (76+5^4+32\times 1).$
\item [] $5874=9\times (8\times 76+5+4)+321.$
\item [] $5875=(9+8+7+65)\times (4^3+2)+1.$
\item [] $5876=(9\times 8+76\times 5)\times (4+3^2)\times 1.$
\item [] $5877=9\times 8\times (76+5)+43+2\times 1.$
\item [] $5878=9\times 8\times (76+5)+43+2+1.$
\item [] $5879=(9+8)\times 7+6\times 5\times 4^3\times (2+1).$
\item [] $5880=98\times (7+6+5+4\times 3)\times 2\times 1.$
\item [] $5881=9+8\times (76+5^4+32+1).$
\item [] $5882=9\times (8+7\times 6)+5432\times 1.$
\item [] $5883=9\times (8+7\times 6)+5432+1.$
\item [] $5884=9+(8+76+5)\times (4^3+2)+1.$
\item [] $\mathit{5885=987\times 6+5-43+2-1.}$
\item [] $5886=(9+8\times 7\times 6)\times (5+4\times 3)+21.$
\item [] $5887=9\times ((8+7+6)\times 5+4)\times 3\times 2+1.$
\item [] $5888=((9+87)\times 6\times 5+4^3)\times 2\times 1.$
\item [] $5889=9+8\times 7\times (6+5+4)\times (3\times 2+1).$
\item [] $5890=98\times (7\times 6+5)+4\times 321.$
\item [] $5891=(9+8\times (7+6\times 5\times 4\times 3))\times 2+1.$
\item [] $5892=9\times (8\times (76+5)+4)+3+21.$
\item [] $5893=(9+8\times 7+6)\times (5\times 4+3\times 21).$
\item [] $5894=98+7\times 6\times (5+4^3)\times 2\times 1.$
\item [] $5895=9\times (87\times 6+5+4\times 32\times 1).$
\item [] $5896=9\times 8\times (76+5)+43+21.$
\item [] $5897=9\times 8+7\times (6+5\times 4)\times 32+1.$
\item [] $5898=(987+654\times 3)\times 2\times 1.$
\item [] $5899=9\times 8\times (76+5)+4+3\times 21.$
\item [] $5900=9\times (8+7)\times 6\times 5+43^2+1.$
\item [] $5901=(9\times 87+6+54)\times (3\times 2+1).$
\item [] $5902=9\times 8\times 76+5\times 43\times 2\times 1.$
\item [] $5903=9\times 8\times 76+5\times 43\times 2+1.$
\item [] $5904=9\times 8\times (7+65+4+3\times 2\times 1).$
\item [] $5905=9+8+7\times (6+5\times 4+3)^2+1.$
\item [] $5906=(9\times 8+7)\times 6+5432\times 1.$
\item [] $5907=(9\times 8+7)\times 6+5432+1.$
\item [] $\mathit{5908=9\times 8\times 76+5+432-1.}$
\item [] $5909=9\times 8\times 76+5+432\times 1.$
\item [] $5910=9\times 8\times 76+5+432+1.$
\item [] $5911=9+8+7+654\times 3^2+1.$
\item [] $5912=(9\times 8+(7+6)\times 5)\times 43+21.$
\item [] $5913=9\times 8\times 7\times 6+(5+4)\times 321.$
\item [] $5914=9\times 87+6+5\times (4^(3+2)+1).$
\item [] $5915=(9\times (87+6)+5)\times (4+3)+21.$
\item [] $5916=9+87\times 65+4\times 3\times 21.$
\item [] $5917=9\times 8\times (76+5)+4^3+21.$
\item [] $5918=9\times 8\times (76+5)+43\times 2\times 1.$
\item [] $5919=9\times 8\times (76+5)+43\times 2+1.$
\item [] $5920=(9+8\times 7+6\times 5\times 4)\times 32\times 1.$
\item [] $5921=(9+8\times 7+6\times 5\times 4)\times 32+1.$
\item [] $5922=(9\times 8+7+6+5+4)\times 3\times 21.$
\item [] $5923=98+7\times (6+5\times 4)\times 32+1.$
\item [] $5924=9\times 87\times 6+(5\times (4+3))^2+1.$
\item [] $5925=9+87\times (6+5\times 4\times 3+2\times 1).$
\item [] $5926=9+87\times (6+5\times 4\times 3+2)+1.$
\item [] $5927=(9\times 8+7)\times (6+5+4^3)+2\times 1.$
\item [] $5928=9\times 8\times (76+5)+4\times (3+21).$
\item [] $5929=9+8\times (7+6\times 5)\times 4\times (3+2\times 1).$
\item [] $5930=9+8\times (7+6\times 5)\times 4\times (3+2)+1.$
\item [] $5931=9\times (8\times 76+5+43+2+1).$
\item [] $5932=9+8+7\times 65\times (4+3^2\times 1).$
\item [] $5933=9+8+7\times 65\times (4+3^2)+1.$
\item [] $5934=9\times (8\times 7\times (6+5)+43)+2+1.$
\item [] $5935=9\times (8+7)\times 6+5\times (4^(3+2)+1).$
\item [] $5936=987\times 6+5+4+3+2\times 1.$
\item [] $5937=987\times 6+5+4+3+2+1.$
\item [] $5938=987\times 6+5+4+3\times 2+1.$
\item [] $5939=9+8\times (7+6)\times (54+3)+2\times 1.$
\item [] $5940=987\times 6+5+4+3^2\times 1.$
\item[]$\mbox{Increasing order}$
\item [] $5941=1\times 23\times 4\times 56+789.$
\item [] $5942=1+23\times 4\times 56+789.$
\item [] $5943=(1^2+3)^4+5678+9.$
\item [] $5944=1^2+3+4\times (5+6)\times (7+8)\times 9.$
\item [] $5945=1+2^3\times (4\times 5\times 6+7\times 89).$
\item [] $5946=(1+2)\times (3+45\times 6\times 7+89).$
\item [] $5947=1\times 2\times 3^4+5\times (6+7)\times 89.$
\item [] $5948=(1+23\times 4)\times (56+7)+89.$
\item [] $5949=(1+2)^3\times 45+6\times 789.$
\item [] $5950=(1+23+4\times 5+6)\times 7\times (8+9).$
\item [] $5951=1+2\times (3+4)\times 5\times (6+7+8\times 9).$
\item [] $5952=12\times 34+(5+6)\times 7\times 8\times 9.$
\item [] $5953=(12\times 34+5\times 67)\times 8+9.$
\item [] $\mathit{5954=12\times 3-45+67\times 89.}$
\item [] $5955=(1+2)\times (34\times 56+78)+9.$
\item [] $5956=(1^2+3)\times (4+(5+6)\times (7+8)\times 9).$
\item [] $5957=(1+2+34)\times (5+67+89).$
\item [] $5958=(12\times 3^4+5)\times 6+7+89.$
\item [] $5959=((1+23)\times 4+5)\times (6\times 7+8+9).$
\item [] $5960=12^3+4\times 5+6\times 78\times 9.$
\item [] $5961=1\times 2\times (3+45)\times (6+7\times 8)+9.$
\item [] $5962=12+34\times (56+7\times (8+9)).$
\item [] $5963=1^{2345}\times 67\times 89.$
\item [] $5964=1^{2345}+67\times 89.$
\item [] $5965=1\times (2+3)\times (4\times (5\times 6+7)\times 8+9).$
\item [] $5966=1^2\times 3+(4+56+7)\times 89.$
\item [] $5967=1^2+3+(4+56+7)\times 89.$
\item [] $5968=1^{234}\times 5+67\times 89.$
\item [] $5969=1^{234}+5+67\times 89.$
\item [] $5970=1+2\times 3+(4+56+7)\times 89.$
\item [] $5971=123+4^5+67\times 8\times 9.$
\item [] $5972=1^{23}\times 4+5+67\times 89.$
\item [] $5973=1234+5+6\times 789.$
\item [] $\mathit{5974=1-2+3+4+5+67\times 89.}$
\item [] $5975=1^2\times 3+4+5+67\times 89.$
\item [] $5976=1^2+3+4+5+67\times 89.$
\item [] $5977=1\times 2+3+4+5+67\times 89.$
\item [] $5978=1+2+3+4+5+67\times 89.$
\item [] $5979=1+2\times 3+4+5+67\times 89.$
\item [] $5980=1^2\times 3\times 4+5+67\times 89.$
\item [] $5981=1+2^3+4+5+67\times 89.$
\item [] $5982=1\times 2+3\times 4+5+67\times 89.$
\item [] $5983=1+2+3\times 4+5+67\times 89.$
\item [] $5984=1^{23}+4\times 5+67\times 89.$
\item [] $5985=12^3+45+6\times 78\times 9.$
\item [] $5986=1^2\times 3+4\times 5+67\times 89.$
\item [] $5987=12+3+4+5+67\times 89.$
\item [] $5988=1+2\times (3+4+5)+67\times 89.$
\item [] $5989=1+2+3+4\times 5+67\times 89.$
\item [] $5990=1+2\times 3+4\times 5+67\times 89.$
\item [] $5991=1^2+3\times (4+5)+67\times 89.$
\item [] $5992=12+3\times 4+5+67\times 89.$
\item [] $5993=1+2+3\times (4+5)+67\times 89.$
\item [] $5994=1\times 234\times 5+67\times 8\times 9.$
\item [] $5995=1+234\times 5+67\times 8\times 9.$
\item [] $5996=1+23+4+5+67\times 89.$
\item [] $5997=(12+3^4)\times 56+789.$
\item [] $5998=12+3+4\times 5+67\times 89.$
\item [] $5999=(1+2)^3+4+5+67\times 89.$
\item [] $6000=1\times 2^3\times 4+5+67\times 89.$
\item [] $6001=1+2+(3+4)\times 5+67\times 89.$
\item [] $6002=1^2\times 34+5+67\times 89.$
\item [] $6003=1^2+34+5+67\times 89.$
\item [] $6004=1\times 2+34+5+67\times 89.$
\item [] $6005=1+2+34+5+67\times 89.$
\item [] $6006=1\times 23+4\times 5+67\times 89.$
\item [] $6007=1+23+4\times 5+67\times 89.$
\item [] $6008=12\times 3+4+5+67\times 89.$
\item [] $6009=1^{23}+45+67\times 89.$
\item [] $6010=12+(3+4)\times 5+67\times 89.$
\item[]$\mbox{Decreasing order}$
\item [] $5941=987\times 6+5+4+3^2+1.$
\item [] $5942=9\times 8\times 7+6+5432\times 1.$
\item [] $5943=9\times 8\times 7+6+5432+1.$
\item [] $\mathit{5944=987\times 6+54-32\times 1.}$
\item [] $5945=9+8+76\times (54+3+21).$
\item [] $5946=9+8\times (7+6)\times 54+321.$
\item [] $5947=987\times 6+5\times 4+3+2\times 1.$
\item [] $5948=987\times 6+5\times 4+3+2+1.$
\item [] $5949=987\times 6+5\times 4+3\times 2+1.$
\item [] $5950=9\times (8\times 76+5)+432+1.$
\item [] $5951=987\times 6+5+4\times 3\times 2\times 1.$
\item [] $5952=987\times 6+5+4\times 3\times 2+1.$
\item [] $5953=9+8+7+(65+4\times 3)^2\times 1.$
\item [] $5954=9+8+7+(65+4\times 3)^2+1.$
\item [] $5955=987\times 6+5+4+3+21.$
\item [] $5956=987\times 6+(5+4\times 3)\times 2\times 1.$
\item [] $5957=987\times 6+(5+4\times 3)\times 2+1.$
\item [] $5958=9\times 87\times 6+5\times 4\times 3\times 21.$
\item [] $5959=9\times 8\times 76+54\times 3^2+1.$
\item [] $5960=987\times 6+5+4\times 3+21.$
\item [] $5961=9\times 8\times (76+5)+4\times 32+1.$
\item [] $5962=9+(8\times 7+6)\times (5+43)\times 2+1.$
\item [] $5963=9+87\times 6+5432\times 1.$
\item [] $5964=9+87\times 6+5432+1.$
\item [] $5965=9\times 8+7+654\times 3^2\times 1.$
\item [] $5966=987\times 6+5\times 4+3+21.$
\item [] $5967=9+(8\times 76+54)\times 3^2\times 1.$
\item [] $5968=987\times 6+5\times (4+3+2)+1.$
\item [] $5969=987\times 6+(5\times 4+3)\times 2+1.$
\item [] $5970=987\times 6+(5+4)\times 3+21.$
\item [] $5971=98\times 7\times 6+5+43^2+1.$
\item [] $5972=987\times 6+5+43+2\times 1.$
\item [] $5973=987\times 6+5+43+2+1.$
\item [] $5974=987\times 6+5\times 4+32\times 1.$
\item [] $5975=987\times 6+5\times 4+32+1.$
\item [] $5976=9\times 8\times 7\times (6+5)+432\times 1.$
\item [] $5977=987\times 6+5+(4+3)^2+1.$
\item [] $5978=987\times 6+5\times (4+3)+21.$
\item [] $5979=9\times 8\times (76+5)+(4+3)\times 21.$
\item [] $5980=98\times (7+6+5+43)+2\times 1.$
\item [] $5981=987\times 6+54+3+2\times 1.$
\item [] $5982=987\times 6+54+3+2+1.$
\item [] $5983=987\times 6+54+3\times 2+1.$
\item [] $5984=987\times 6+5\times 4\times 3+2\times 1.$
\item [] $5985=987\times 6+54+3^2\times 1.$
\item [] $5986=987\times 6+54+3^2+1.$
\item [] $5987=9\times 87\times 6+5+4\times 321.$
\item [] $5988=987\times 6+5\times (4+3^2)+1.$
\item [] $5989=9+87\times 65+4+321.$
\item [] $5990=9\times (8\times 76+54)+32\times 1.$
\item [] $5991=987\times 6+5+43+21.$
\item [] $5992=987\times 6+5\times (4+3)\times 2\times 1.$
\item [] $5993=98\times 7\times 6+5^4\times 3+2\times 1.$
\item [] $5994=987\times 6+5+4+3\times 21.$
\item [] $5995=(98+7+6)\times (5+4)\times 3\times 2+1.$
\item [] $5996=9\times (8+7\times 6\times 5+4)\times 3+2\times 1.$
\item [] $5997=987\times 6+5\times (4\times 3+2+1).$
\item [] $\mathit{5998=-9+87\times (65+4)+3+2-1.}$
\item [] $5999=98\times (7+6+5+43)+21.$
\item [] $6000=987\times 6+54+3+21.$
\item [] $6001=(98+7+6)\times 54+3\times 2+1.$
\item [] $6002=9\times 8\times 76+(5\times 4+3)^2+1.$
\item [] $6003=987\times 6+5\times 4\times 3+21.$
\item [] $6004=987\times 6+(5+4)\times 3^2+1.$
\item [] $6005=987\times 6+5\times 4+3\times 21.$
\item [] $6006=9+876+5\times 4^(3+2)+1.$
\item [] $6007=9+87\times 65+(4+3)^{(2+1)}.$
\item [] $6008=987\times 6+54+32\times 1.$
\item [] $6009=987\times 6+54+32+1.$
\item [] $6010=9+8\times 7\times 6\times 5+4321.$
\item[]$\mbox{Increasing order}$
\item [] $6011=1^2\times 3+45+67\times 89.$
\item [] $6012=123\times 45+6\times 78+9.$
\item [] $6013=1\times 2+3+45+67\times 89.$
\item [] $6014=12+34+5+67\times 89.$
\item [] $6015=1+2\times 3+45+67\times 89.$
\item [] $6016=1\times 2^3+45+67\times 89.$
\item [] $6017=1+2^3+45+67\times 89.$
\item [] $6018=123+45\times (6\times 7+89).$
\item [] $6019=12\times 3+4\times 5+67\times 89.$
\item [] $6020=1+234+5\times (6+7)\times 89.$
\item [] $6021=12^3+(4+5)\times (6\times 78+9).$
\item [] $6022=1+(2^{(3+4)}+5+67\times 8)\times 9.$
\item [] $6023=12+3+45+67\times 89.$
\item [] $6024=1^2+3\times 4\times 5+67\times 89.$
\item [] $6025=1\times 2+3\times 4\times 5+67\times 89.$
\item [] $6026=1+2+3\times 4\times 5+67\times 89.$
\item [] $6027=12+(3+4)\times (5+6)\times 78+9.$
\item [] $6028=(12+3)\times 4+5+67\times 89.$
\item [] $6029=(12+3)\times (4+56\times 7)+89.$
\item [] $6030=(1+2+3+4)\times (5+6+7\times 8)\times 9.$
\item [] $6031=1\times 23+45+67\times 89.$
\item [] $6032=1+23+45+67\times 89.$
\item [] $6033=1\times 2\times (3+4)\times 5+67\times 89.$
\item [] $6034=1+2\times (3+4)\times 5+67\times 89.$
\item [] $6035=12+3\times 4\times 5+67\times 89.$
\item [] $6036=1\times 2\times 34+5+67\times 89.$
\item [] $6037=1+2\times 34+5+67\times 89.$
\item [] $6038=(1+2+3\times 4)\times 5+67\times 89.$
\item [] $6039=(1+2)^3\times 45+67\times 8\times 9.$
\item [] $6040=1^2+(3\times 45+67\times 8)\times 9.$
\item [] $6041=1\times 2\times (34+5)+67\times 89.$
\item [] $6042=1+2\times (34+5)+67\times 89.$
\item [] $6043=(1^2+3)\times 4\times 5+67\times 89.$
\item [] $6044=12\times 3+45+67\times 89.$
\item [] $6045=123\times 45+6+7\times 8\times 9.$
\item [] $6046=1\times 2^{(3\times 4)}+5\times 6\times (7\times 8+9).$
\item [] $6047=1+2^{(3\times 4)}+5\times 6\times (7\times 8+9).$
\item [] $6048=1\times 2\times 3\times (4+5\times 6+78)\times 9.$
\item [] $6049=1^2\times 3^4+5+67\times 89.$
\item [] $6050=1^2+3^4+5+67\times 89.$
\item [] $6051=1\times 2+3^4+5+67\times 89.$
\item [] $6052=1+2+3^4+5+67\times 89.$
\item [] $6053=1+23\times 4\times 5\times (6+7)+8\times 9.$
\item [] $6054=1+2\times (3+45\times 67)+8+9.$
\item [] $6055=1+2\times 3\times (4\times 5\times (6\times 7+8)+9).$
\item [] $6056=1\times (2+3)\times 4^5+(6+7)\times 8\times 9.$
\item [] $6057=123\times 45+6\times (78+9).$
\item [] $6058=(12+3+4)\times 5+67\times 89.$
\item [] $6059=1\times 2\times (3+45)+67\times 89.$
\item [] $6060=1\times 23\times 4+5+67\times 89.$
\item [] $6061=12+3^4+5+67\times 89.$
\item [] $6062=1+2\times 34\times (5+6+78)+9.$
\item [] $6063=1234+5+67\times 8\times 9.$
\item [] $6064=(1+23)\times 4+5+67\times 89.$
\item [] $6065=12\times (3+4+5)\times 6\times 7+8+9.$
\item [] $6066=(1+2+3\times 45+67\times 8)\times 9.$
\item [] $6067=1+2\times 3^4\times (5\times 6+7)+8\times 9.$
\item [] $6068=(1+2)\times (3+4)\times 5+67\times 89.$
\item [] $6069=1\times 23\times 4\times 5\times (6+7)+89.$
\item [] $6070=(1+2)\times 34+5+67\times 89.$
\item [] $6071=1+2\times (3+45\times 67+8+9).$
\item [] $6072=(1+2)^3\times 4\times 56+7+8+9.$
\item [] $6073=(1+2\times 345+67)\times 8+9.$
\item [] $6074=1\times 2+3+(45+6)\times 7\times (8+9).$
\item [] $6075=1+2+3+(45+6)\times 7\times (8+9).$
\item [] $6076=(1+2)^3\times 4+5+67\times 89.$
\item [] $6077=1\times 2^3+(45+6)\times 7\times (8+9).$
\item [] $6078=1\times 2\times 3\times (4\times 56+789).$
\item [] $6079=1+2\times 3\times (4\times 56+789).$
\item [] $6080=123\times 45+67\times 8+9.$
\item[]$\mbox{Decreasing order}$
\item [] $6011=98\times (7+6\times (5+4))+32+1.$
\item [] $6012=987\times 6+5+4^3+21.$
\item [] $6013=987\times 6+5+43\times 2\times 1.$
\item [] $6014=987\times 6+5+43\times 2+1.$
\item [] $6015=(9+8\times 7)\times 6+5^4\times 3^2\times 1.$
\item [] $6016=98\times 7+(6\times 5+43)^2+1.$
\item [] $6017=9\times 8\times 76+543+2\times 1.$
\item [] $6018=9\times 8\times 76+543+2+1.$
\item [] $6019=987\times 6+(5+43)\times 2+1.$
\item [] $6020=9+8+(7\times 6+5^4)\times 3^2\times 1.$
\item [] $6021=9+87\times (65+4)+3^2\times 1.$
\item [] $6022=9+87\times (65+4)+3^2+1.$
\item [] $6023=987\times 6+5\times 4\times (3+2)+1.$
\item [] $6024=9\times 8\times (7+6)\times 5+4^3\times 21.$
\item [] $6025=9+87+(65+4\times 3)^2\times 1.$
\item [] $6026=(98+7+6)\times 54+32\times 1.$
\item [] $6027=9\times (8+7+654)+3+2+1.$
\item [] $6028=9\times (8+7+654)+3\times 2+1.$
\item [] $6029=9\times 8\times 7+65\times (4^3+21).$
\item [] $6030=9+(8+7+654)\times 3^2\times 1.$
\item [] $6031=9+(8+7+654)\times 3^2+1.$
\item [] $\mathit{6032=-9+8\times 76+5432+1.}$
\item [] $6033=9+8+(7\times 6+5)\times 4^3\times 2\times 1.$
\item [] $6034=9\times (87+6)\times 5+43^2\times 1.$
\item [] $6035=9\times (87+6)\times 5+43^2+1.$
\item [] $6036=9\times 8\times 76+543+21.$
\item [] $6037=987\times 6+(54+3)\times 2+1.$
\item [] $6038=(9\times 8+7)\times 65+43\times 21.$
\item [] $6039=987\times 6+54+3\times 21.$
\item [] $6040=(9+8\times 76+54)\times 3^2+1.$
\item [] $6041=98\times (7+6\times (5+4))+3\times 21.$
\item [] $6042=987\times 6+5\times 4\times 3\times 2\times 1.$
\item [] $6043=987\times 6+5\times 4\times 3\times 2+1.$
\item [] $6044=9+87\times (65+4)+32\times 1.$
\item [] $6045=9+87\times (65+4)+32+1.$
\item [] $6046=(9+8\times 7\times 6)\times 5+4321.$
\item [] $6047=9\times 87\times 6+5+4^3\times 21.$
\item [] $6048=98+(7+654)\times 3^2+1.$
\item [] $6049=9+8\times 76+5432\times 1.$
\item [] $6050=9+8\times 76+5432+1.$
\item [] $6051=9\times 8\times (7+65+4\times 3)+2+1.$
\item [] $6052=(9+8)\times ((76+5)\times 4+32\times 1).$
\item [] $6053=9\times (8+7+654)+32\times 1.$
\item [] $6054=9\times (8+7+654)+32+1.$
\item [] $6055=987\times 6+5+4^3\times 2\times 1.$
\item [] $6056=987\times 6+5+4^3\times 2+1.$
\item [] $6057=(98+7+6)\times 54+3\times 21.$
\item [] $6058=9+8\times 7\times (6+5+43)\times 2+1.$
\item [] $6059=987\times 6+5+4\times (32+1).$
\item [] $6060=987\times 6+(5+4^3)\times 2\times 1.$
\item [] $6061=(98\times 7+6\times 54)\times 3\times 2+1.$
\item [] $6062=(98+7)\times 6+5432\times 1.$
\item [] $6063=(98+7)\times 6+5432+1.$
\item [] $6064=9\times (8+7)+(65+4\times 3)^2\times 1.$
\item [] $6065=9+8+(76+5\times 4)\times 3\times 21.$
\item [] $6066=987\times 6+(5+4+3)^2\times 1.$
\item [] $6067=987\times 6+(5+4+3)^2+1.$
\item [] $6068=(9\times 8+76)\times (5+4+32\times 1).$
\item [] $6069=9\times 8\times (7+65+4\times 3)+21.$
\item [] $6070=(9+8)\times 7\times (6+5\times (4+3+2))+1.$
\item [] $6071=987\times 6+5+(4\times 3)^2\times 1.$
\item [] $6072=987\times 6+5+(4\times 3)^2+1.$
\item [] $6073=(9\times 8\times 7\times 6+5+4+3)\times 2+1.$
\item [] $6074=987\times 6+5+(4+3)\times 21.$
\item [] $6075=9+87\times (65+4)+3\times 21.$
\item [] $6076=98\times (7+6\times 5+4\times 3\times 2+1).$
\item [] $6077=987\times 6+5\times (4+3^{(2+1)}).$
\item [] $6078=(987+6+5\times 4)\times 3\times 2\times 1.$
\item [] $6079=(987+6+5\times 4)\times 3\times 2+1.$
\item [] $6080=(9+8\times 7+6\times 5)\times (43+21).$
\item[]$\mbox{Increasing order}$
\item [] $6081=1\times 23\times (4\times 5+6+7)\times 8+9.$
\item [] $6082=1+23\times (4\times 5+6+7)\times 8+9.$
\item [] $6083=1\times 2345+6\times 7\times 89.$
\item [] $6084=1+2345+6\times 7\times 89.$
\item [] $6085=1^2+(34+5)\times (67+89).$
\item [] $6086=123+(4+56+7)\times 89.$
\item [] $6087=12+(34+5+6)\times (7+8)\times 9.$
\item [] $6088=(1+2\times 3\times 4)\times 5+67\times 89.$
\item [] $6089=1+2^{(3+4)}\times (5+6\times 7)+8\times 9.$
\item [] $6090=(1\times 2+3\times 4+56)\times (78+9).$
\item [] $6091=1+(2+3\times 4+56)\times (78+9).$
\item [] $6092=1\times 23+(45+6)\times 7\times (8+9).$
\item [] $6093=123\times 45+(6+7\times 8)\times 9.$
\item [] $6094=1+234\times (5+6+7+8)+9.$
\item [] $6095=12\times 34+5678+9.$
\item [] $6096=1^2\times 3\times 4^5+6\times 7\times 8\times 9.$
\item [] $6097=1+2^3\times (4+56+78\times 9).$
\item [] $6098=1\times 2+3\times 4^5+6\times 7\times 8\times 9.$
\item [] $6099=1^2+3\times 45+67\times 89.$
\item [] $6100=1\times 2+3\times 45+67\times 89.$
\item [] $6101=1+2+3\times 45+67\times 89.$
\item [] $6102=1^{2345}\times 678\times 9.$
\item [] $6103=1^{2345}+678\times 9.$
\item [] $6104=(1^2+3)^4\times 5+67\times 8\times 9.$
\item [] $6105=12+3\times (4\times 5+6)\times 78+9.$
\item [] $6106=1234+56\times (78+9).$
\item [] $6107=1^{234}\times 5+678\times 9.$
\item [] $6108=12+3\times 4^5+6\times 7\times 8\times 9.$
\item [] $6109=12^3+4+56\times 78+9.$
\item [] $6110=12+3\times 45+67\times 89.$
\item [] $6111=1^{23}\times 4+5+678\times 9.$
\item [] $6112=12\times 3\times 4+5+67\times 89.$
\item [] $6113=(1+2)^3\times 4\times 56+7\times 8+9.$
\item [] $6114=1^2\times 3+4+5+678\times 9.$
\item [] $6115=1^2+3+4+5+678\times 9.$
\item [] $6116=1\times 2+3+4+5+678\times 9.$
\item [] $6117=1+2+3+4+5+678\times 9.$
\item [] $6118=1+2\times 3+4+5+678\times 9.$
\item [] $6119=1^2\times 3\times 4+5+678\times 9.$
\item [] $6120=1+2^3+4+5+678\times 9.$
\item [] $6121=1\times 2+3\times 4+5+678\times 9.$
\item [] $6122=1+2+3\times 4+5+678\times 9.$
\item [] $6123=1^{23}+4\times 5+678\times 9.$
\item [] $6124=1+2^3\times 4\times 5+67\times 89.$
\item [] $6125=1^2\times 3+4\times 5+678\times 9.$
\item [] $6126=12+3+4+5+678\times 9.$
\item [] $6127=1\times 2+3+4\times 5+678\times 9.$
\item [] $6128=1+2+3+4\times 5+678\times 9.$
\item [] $6129=1+2\times 3+4\times 5+678\times 9.$
\item [] $6130=1\times 2^3+4\times 5+678\times 9.$
\item [] $6131=123+45+67\times 89.$
\item [] $6132=1+2\times 3\times 4+5+678\times 9.$
\item [] $6133=1^2\times 34\times 5+67\times 89.$
\item [] $6134=1^2+34\times 5+67\times 89.$
\item [] $6135=1\times 2+34\times 5+67\times 89.$
\item [] $6136=1+2+34\times 5+67\times 89.$
\item [] $6137=12+3+4\times 5+678\times 9.$
\item [] $6138=(1+2)^3+4+5+678\times 9.$
\item [] $6139=1\times 2^3\times 4+5+678\times 9.$
\item [] $6140=1+2^3\times 4+5+678\times 9.$
\item [] $6141=1^2\times 34+5+678\times 9.$
\item [] $6142=1^2+34+5+678\times 9.$
\item [] $6143=1\times 2+34+5+678\times 9.$
\item [] $6144=1+2+34+5+678\times 9.$
\item [] $6145=1\times 23+4\times 5+678\times 9.$
\item [] $6146=1+23+4\times 5+678\times 9.$
\item [] $6147=12\times 3+4+5+678\times 9.$
\item [] $6148=1^{23}+45+678\times 9.$
\item [] $6149=12+(3+4)\times 5+678\times 9.$
\item [] $6150=1^2\times 3+45+678\times 9.$
\item[]$\mbox{Decreasing order}$
\item [] $6081=98\times (7\times 6+5\times 4)+3+2\times 1.$
\item [] $6082=98\times (7\times 6+5\times 4)+3+2+1.$
\item [] $6083=(9\times 8\times 7\times 6+5+4\times 3)\times 2+1.$
\item [] $6084=9\times 8\times (76+5)+4\times 3\times 21.$
\item [] $6085=98\times (7+6+5)+4321.$
\item [] $6086=987\times 6+54\times 3+2\times 1.$
\item [] $6087=987\times 6+54\times 3+2+1.$
\item [] $6088=9\times 8+(7\times 6+5)\times 4\times 32\times 1.$
\item [] $6089=9\times 8+(7\times 6+5)\times 4\times 32+1.$
\item [] $6090=(98+7)\times (6+5\times 4+32\times 1).$
\item [] $6091=(98+7)\times (6+5\times 4+32)+1.$
\item [] $6092=(98+76)\times 5\times (4+3)+2\times 1.$
\item [] $6093=9\times (8\times 76+5+43+21).$
\item [] $6094=9\times 8\times 7+65\times 43\times 2\times 1.$
\item [] $6095=9\times 8\times 7+65\times 43\times 2+1.$
\item [] $6096=9+87\times 65+432\times 1.$
\item [] $6097=9+87\times 65+432+1.$
\item [] $\mathit{6098=-9+8\times 765-4\times 3-2+1.}$
\item [] $6099=(9\times 8+7)\times 6+5^4\times 3^2\times 1.$
\item [] $6100=(9\times 8+7)\times 6+5^4\times 3^2+1.$
\item [] $6101=98+(7\times 6+5^4)\times 3^2\times 1.$
\item [] $6102=987\times 6+5\times 4\times 3^2\times 1.$
\item [] $6103=987\times 6+5\times 4\times 3^2+1.$
\item [] $6104=9\times 8\times 76+5^4+3\times 2+1.$
\item [] $6105=987\times 6+54\times 3+21.$
\item [] $6106=9\times 8\times 76+5^4+3^2\times 1.$
\item [] $6107=9\times 8\times 76+5^4+3^2+1.$
\item [] $6108=98\times (7\times 6+5\times 4)+32\times 1.$
\item [] $6109=98\times (7\times 6+5\times 4)+32+1.$
\item [] $6110=(9+8\times 7)\times (6\times 5+43+21).$
\item [] $6111=9\times (8+7+654+3^2+1).$
\item [] $6112=9\times 8\times 76+5\times 4^3\times 2\times 1.$
\item [] $6113=9\times 8\times 76+5\times 4^3\times 2+1.$
\item [] $6114=98+(7\times 6+5)\times 4^3\times 2\times 1.$
\item [] $6115=98+(7\times 6+5)\times 4^3\times 2+1.$
\item [] $\mathit{6116=9+8\times 765-4\times 3-2+1.}$
\item [] $6117=9\times 8\times 76+5\times 43\times (2+1).$
\item [] $6118=((98+76)\times 5+4)\times (3\times 2+1).$
\item [] $6119=987\times 6+5+4^3\times (2+1).$
\item [] $6120=9\times 8+7\times 6\times (5+43)\times (2+1).$
\item [] $6121=9\times 8\times 76+5^4+3+21.$
\item [] $6122=(9+87+6)\times 5\times 4\times 3+2\times 1.$
\item [] $6123=(9+87+6)\times 5\times 4\times 3+2+1.$
\item [] $6124=98\times 7+6+5432\times 1.$
\item [] $6125=98\times 7+6+5432+1.$
\item [] $6126=(987+6\times 5+4)\times 3\times 2\times 1.$
\item [] $6127=9\times 8+7+6\times (5+43)\times 21.$
\item [] $6128=9+8+(76+5\times 43)\times 21.$
\item [] $6129=9\times 8\times 76+5^4+32\times 1.$
\item [] $6130=9\times 8\times 76+5^4+32+1.$
\item [] $6131=(9+8\times 7)\times (6\times 5+4^3)+21.$
\item [] $6132=9\times 8\times 76+5\times 4\times (32+1).$
\item [] $\mathit{6133=9+8\times 765+4+3-2-1.}$
\item [] $6134=98\times 7+6\times (5+43\times 21).$
\item [] $6135=9\times 8\times 7+6+5^4\times 3^2\times 1.$
\item [] $6136=9\times 8\times 7+6+5^4\times 3^2+1.$
\item [] $6137=9+8+(7+65)\times (4^3+21).$
\item [] $6138=9+8\times 765+4+3+2\times 1.$
\item [] $6139=987\times 6+5\times 43+2\times 1.$
\item [] $6140=987\times 6+5\times 43+2+1.$
\item [] $6141=(9+87+6)\times 5\times 4\times 3+21.$
\item [] $6142=9+8\times 765+4+3^2\times 1.$
\item [] $6143=9+8\times 765+4+3^2+1.$
\item [] $6144=9+8\times 765+4\times 3+2+1.$
\item [] $6145=(9\times 8\times 7\times 6+5+43)\times 2+1.$
\item [] $6146=(9+8)\times 7\times 6+5432\times 1.$
\item [] $6147=(9+8)\times 7\times 6+5432+1.$
\item [] $6148=987\times 6+5\times (43+2)+1.$
\item [] $6149=9+8\times 765+4\times (3+2\times 1).$
\item [] $6150=9+8\times 765+4\times (3+2)+1.$
\item[]$\mbox{Increasing order}$
\item [] $6151=1^2+3+45+678\times 9.$
\item [] $6152=1\times 2+3+45+678\times 9.$
\item [] $6153=12+34+5+678\times 9.$
\item [] $6154=1+2\times 3+45+678\times 9.$
\item [] $6155=1\times 2^3+45+678\times 9.$
\item [] $6156=1+2^3+45+678\times 9.$
\item [] $6157=(1+2\times 3+4)\times 5+678\times 9.$
\item [] $6158=12\times 3+4\times 5+678\times 9.$
\item [] $6159=1^2\times 3\times 4+(5+678)\times 9.$
\item [] $6160=1+2^3+4+(5+678)\times 9.$
\item [] $6161=1\times 2+3\times 4+(5+678)\times 9.$
\item [] $6162=12+3+45+678\times 9.$
\item [] $6163=1^2+3\times 4\times 5+678\times 9.$
\item [] $6164=123\times 45+6+7\times 89.$
\item [] $6165=1+2+3\times 4\times 5+678\times 9.$
\item [] $6166=12+3+4+(5+678)\times 9.$
\item [] $6167=12\times (3\times 4+5)+67\times 89.$
\item [] $6168=1^{23}\times 4^5\times 6+7+8+9.$
\item [] $6169=1^{23}+4^5\times 6+7+8+9.$
\item [] $6170=1\times 23+45+678\times 9.$
\item [] $6171=1+23+45+678\times 9.$
\item [] $6172=1^2+3+4^5\times 6+7+8+9.$
\item [] $6173=1\times 2+3+4^5\times 6+7+8+9.$
\item [] $6174=12+3\times 4\times 5+678\times 9.$
\item [] $6175=1\times 2\times 34+5+678\times 9.$
\item [] $6176=1+2\times 34+5+678\times 9.$
\item [] $6177=1+2^3+4^5\times 6+7+8+9.$
\item [] $6178=(1+2)^3+4+(5+678)\times 9.$
\item [] $6179=123\times 4+5678+9.$
\item [] $6180=1\times 2\times (34+5)+678\times 9.$
\item [] $6181=1^2\times 34+(5+678)\times 9.$
\item [] $6182=1^2+34+(5+678)\times 9.$
\item [] $6183=12\times 3+45+678\times 9.$
\item [] $6184=1+2+34+(5+678)\times 9.$
\item [] $6185=((12+3)\times (45+6)+7)\times 8+9.$
\item [] $6186=1^2\times 3+(4+5+678)\times 9.$
\item [] $6187=123\times 4+5\times 67\times (8+9).$
\item [] $6188=1^2\times 3^4+5+678\times 9.$
\item [] $6189=1^2+3^4+5+678\times 9.$
\item [] $6190=1\times 2+3^4+5+678\times 9.$
\item [] $6191=1+2+3^4+5+678\times 9.$
\item [] $6192=1+23+4^5\times 6+7+8+9.$
\item [] $6193=(12+34)\times 5+67\times 89.$
\item [] $6194=1\times 2\times (3+(4\times 5+6)\times 7\times (8+9)).$
\item [] $6195=(1+2)^3+4^5\times 6+7+8+9.$
\item [] $6196=1+(2\times (3\times 4^5+6+7+8)+9).$
\item [] $6197=(12+3+4)\times 5+678\times 9.$
\item [] $6198=123\times 4\times 5+6\times 7\times 89.$
\item [] $6199=1\times 23\times 4+5+678\times 9.$
\item [] $6200=1234\times 5+6+7+8+9.$
\item [] $6201=123\times 4\times (5+6)+789.$
\item [] $6202=1\times 234+5+67\times 89.$
\item [] $6203=1+234+5+67\times 89.$
\item [] $6204=12\times 3+4^5\times 6+7+8+9.$
\item [] $6205=1+(2\times 3+4^5)\times 6+7+8+9.$
\item [] $6206=1\times 23+(4+5+678)\times 9.$
\item [] $6207=1+23+(4+5+678)\times 9.$
\item [] $6208=1\times 2^3\times (4\times 5+(6+78)\times 9).$
\item [] $6209=1^{23}\times 4^5\times 6+7\times 8+9.$
\item [] $6210=1^{23}+4^5\times 6+7\times 8+9.$
\item [] $6211=(1+2)\times 3^4+5+67\times 89.$
\item [] $6212=1^2\times 3+4^5\times 6+7\times 8+9.$
\item [] $6213=1^2+3+4^5\times 6+7\times 8+9.$
\item [] $6214=1\times 2+3+4^5\times 6+7\times 8+9.$
\item [] $6215=1+2+3+4^5\times 6+7\times 8+9.$
\item [] $6216=1+2\times 3+4^5\times 6+7\times 8+9.$
\item [] $6217=1\times 2^3+4^5\times 6+7\times 8+9.$
\item [] $6218=1+2^3+4^5\times 6+7\times 8+9.$
\item [] $6219=(1^2+3+4+5+678)\times 9.$
\item [] $6220=1\times 2\times (3^4+5+6\times 7\times 8\times 9).$
\item[]$\mbox{Decreasing order}$
\item [] $6151=9+(8+765\times 4+3)\times 2\times 1.$
\item [] $6152=987\times 6+5\times (43+2+1).$
\item [] $6153=9+8\times 765+4\times 3\times 2\times 1.$
\item [] $6154=9+8\times 765+4\times 3\times 2+1.$
\item [] $\mathit{6155=9\times 876-54\times 32-1.}$
\item [] $6156=9+87\times 6+5^4\times 3^2\times 1.$
\item [] $6157=9+8\times 765+4+3+21.$
\item [] $6158=987\times 6+5\times 43+21.$
\item [] $6159=9\times 8\times 7+65\times (43\times 2+1).$
\item [] $6160=9\times 8\times 76+5^4+3\times 21.$
\item [] $6161=(9+8+765\times 4+3)\times 2+1.$
\item [] $6162=9+8\times 765+4\times 3+21.$
\item [] $6163=(9\times 8+7\times 6)\times 54+3\times 2+1.$
\item [] $6164=(9\times 8+7)\times (6+5\times 4)\times 3+2\times 1.$
\item [] $6165=9+8\times 765+4+32\times 1.$
\item [] $6166=9+8\times 765+4+32+1.$
\item [] $6167=987\times 6+5\times (4+3)^2\times 1.$
\item [] $6168=987\times 6+5\times (4+3)^2+1.$
\item [] $6169=9+8\times 765+4\times (3^2+1).$
\item [] $6170=9+8\times 7\times (65+43+2)+1.$
\item [] $6171=987+(65+4+3)^2\times 1.$
\item [] $6172=987+(65+4+3)^2+1.$
\item [] $6173=9+8+76\times (5+4)\times 3^2\times 1.$
\item [] $6174=9+8\times 765+43+2\times 1.$
\item [] $6175=9+8\times 765+43+2+1.$
\item [] $6176=9+(8\times 76\times 5+43)\times 2+1.$
\item [] $6177=9+8\times (7\times 6+5\times 43)\times (2+1).$
\item [] $6178=9+8\times 765+(4+3)^2\times 1.$
\item [] $6179=987\times 6+5+4\times 3\times 21.$
\item [] $6180=9\times 8\times 7+6\times (5^4+321).$
\item [] $6181=9+87+((6+5\times 4)\times 3)^2+1.$
\item [] $6182=(9\times 8\times 7+6+5)\times 4\times 3+2\times 1.$
\item [] $6183=9\times 8\times (76+5+4)+3\times 21.$
\item [] $6184=9\times 8\times 7+(6+5^4)\times 3^2+1.$
\item [] $6185=9+8\times (765+4)+3+21.$
\item [] $6186=(9\times 8\times 7\times 6+5+4^3)\times 2\times 1.$
\item [] $6187=9+8\times (765+4+3)+2\times 1.$
\item [] $6188=9\times (8+76)+5432\times 1.$
\item [] $6189=9\times (8+76)+5432+1.$
\item [] $6190=9\times (8\times 7+6+5^4)+3\times 2+1.$
\item [] $6191=9+8+7\times (6\times 5+4\times 3)\times 21.$
\item [] $6192=987\times 6+54\times (3+2)\times 1.$
\item [] $6193=9+8\times 765+43+21.$
\item [] $6194=9+8\times (765+4)+32+1.$
\item [] $6195=9+8\times 765+4^3+2\times 1.$
\item [] $6196=9+8\times 765+4+3\times 21.$
\item [] $6197=9+(876+5)\times (4+3)+21.$
\item [] $6198=9+8\times (7\times 6+(5+4)^3)+21.$
\item [] $6199=(9\times (8\times 7+6\times 5)\times 4+3)\times 2+1.$
\item [] $6200=(9\times (8+7)+65)\times (4+3^{(2+1)}).$
\item [] $6201=(9\times 8\times 7+6+5)\times 4\times 3+21.$
\item [] $6202=9+8\times (76\times 5+4+3)\times 2+1.$
\item [] $6203=9+8\times (7+6+5)\times 43+2\times 1.$
\item [] $6204=9+8\times (7+6+5)\times 43+2+1.$
\item [] $6205=(9+8)\times (7\times 6+5\times 4^3+2+1).$
\item [] $6206=9+8\times (765+4+3)+21.$
\item [] $6207=9\times 8\times 76+5\times (4+3)\times 21.$
\item [] $6208=(98+76+5\times 4)\times 32\times 1.$
\item [] $6209=(98+76+5\times 4)\times 32+1.$
\item [] $6210=987\times 6+(5+4)\times 32\times 1.$
\item [] $6211=987\times 6+(5+4)\times 32+1.$
\item [] $6212=987\times 6+(5+4\times 3)^2+1.$
\item [] $6213=(9\times 8+7+65)\times 43+21.$
\item [] $6214=9+8\times 765+4^3+21.$
\item [] $6215=9+8\times 765+43\times 2\times 1.$
\item [] $6216=9+8\times 765+43\times 2+1.$
\item [] $6217=9+8\times (765+4+3\times 2+1).$
\item [] $6218=9+(8\times 76\times 5+4^3)\times 2+1.$
\item [] $6219=(9\times 8+7\times 6)\times 54+3\times 21.$
\item [] $6220=(9+8\times 76+5)\times (4+3+2+1).$
\item[]$\mbox{Increasing order}$
\item [] $6221=(1+2)\times (3^4+5)+67\times 89.$
\item [] $6222=123\times 45+678+9.$
\item [] $6223=1+2\times 3\times 4\times 5+678\times 9.$
\item [] $6224=12+3+4^5\times 6+7\times 8+9.$
\item [] $6225=(1\times 2+3)\times (456+789).$
\item [] $6226=1^2\times 3+4^5\times 6+7+8\times 9.$
\item [] $6227=1^2+3+4^5\times 6+7+8\times 9.$
\item [] $6228=1\times 2+3+4^5\times 6+7+8\times 9.$
\item [] $6229=1\times 2\times 3+4^5\times 6+7+8\times 9.$
\item [] $6230=1+2\times 3+4^5\times 6+7+8\times 9.$
\item [] $6231=1^{23}\times 4^5\times 6+78+9.$
\item [] $6232=1^{23}+4^5\times 6+78+9.$
\item [] $6233=1\times 2\times 3\times 45+67\times 89.$
\item [] $6234=123+4+5+678\times 9.$
\item [] $6235=1+23\times 45\times 6+7+8+9.$
\item [] $6236=1\times 2+3+4^5\times 6+78+9.$
\item [] $6237=1^2\times 3\times 45+678\times 9.$
\item [] $6238=12+3+4^5\times 6+7+8\times 9.$
\item [] $6239=1\times 2+3\times 45+678\times 9.$
\item [] $6240=1+2+3\times 45+678\times 9.$
\item [] $6241=1234\times 5+6+7\times 8+9.$
\item [] $6242=(1+23+4)\times 5+678\times 9.$
\item [] $6243=123\times 45+6+78\times 9.$
\item [] $6244=1^2+3+4^5\times 6+7+89.$
\item [] $6245=123+4\times 5+678\times 9.$
\item [] $6246=12+3+4^5\times 6+78+9.$
\item [] $6247=1+23+4^5\times 6+7+8\times 9.$
\item [] $6248=1\times 2^3+4^5\times 6+7+89.$
\item [] $6249=12+3\times 45+678\times 9.$
\item [] $6250=(1+2)^3+4^5\times 6+7+8\times 9.$
\item [] $6251=12\times 3\times 4+5+678\times 9.$
\item [] $6252=12\times 34\times 5+6\times 78\times 9.$
\item [] $6253=12+(3+4^5)\times 6+7+8\times 9.$
\item [] $6254=1234\times 5+67+8+9.$
\item [] $6255=1234\times 5+6+7+8\times 9.$
\item [] $6256=1^{23}+45\times (67+8\times 9).$
\item [] $6257=12\times (3^4+5)\times 6+7\times 8+9.$
\item [] $6258=1\times 2\times 3\times 4^5+6\times 7+8\times 9.$
\item [] $6259=12\times 3+4^5\times 6+7+8\times 9.$
\item [] $6260=1\times 2+3+45\times (67+8\times 9).$
\item [] $6261=12+3\times 4\times 5\times (6+7)\times 8+9.$
\item [] $6262=1+2\times 3+45\times (67+8\times 9).$
\item [] $6263=1234\times 5+6+78+9.$
\item [] $6264=1+23+4^5\times 6+7+89.$
\item [] $6265=1^{23}+(45+6\times 7)\times 8\times 9.$
\item [] $6266=12+3^4\times (5+6)\times 7+8+9.$
\item [] $6267=12\times 3+4^5\times 6+78+9.$
\item [] $6268=1+(2\times 3+4^5)\times 6+78+9.$
\item [] $6269=1\times 2\times 3^4+5+678\times 9.$
\item [] $6270=123+45+678\times 9.$
\item [] $6271=1+(2+3+4^5)\times 6+7+89.$
\item [] $6272=1234\times 5+6+7+89.$
\item [] $6273=1^2+34\times 5+678\times 9.$
\item [] $6274=1\times 2+34\times 5+678\times 9.$
\item [] $6275=1+2+34\times 5+678\times 9.$
\item [] $6276=12\times 3+4^5\times 6+7+89.$
\item [] $6277=12^3+4+567\times 8+9.$
\item [] $6278=1\times 23+45\times (67+8\times 9).$
\item [] $6279=1+23+45\times (67+8\times 9).$
\item [] $6280=1+(2^3+4^5)\times 6+78+9.$
\item [] $6281=(1+23\times 45)\times 6+7\times 8+9.$
\item [] $6282=(1^2+3)\times 45+678\times 9.$
\item [] $6283=1\times 2\times 3\times 4^5+67+8\times 9.$
\item [] $6284=12+34\times 5+678\times 9.$
\item [] $6285=1^2\times 3+(4\times 5+678)\times 9.$
\item [] $6286=1^2+3+(4\times 5+678)\times 9.$
\item [] $6287=12\times 3\times (4+5)+67\times 89.$
\item [] $6288=12\times (3+456+7\times 8+9).$
\item [] $6289=1\times 23\times 45\times 6+7+8\times 9.$
\item [] $6290=1+23\times 45\times 6+7+8\times 9.$
\item[]$\mbox{Decreasing order}$
\item [] $6221=9\times 87+6+5432\times 1.$
\item [] $6222=9\times 87+6+5432+1.$
\item [] $6223=(9+8)\times (7+6\times (5+4))\times 3\times 2+1.$
\item [] $6224=9+8\times (765+4)+3\times 21.$
\item [] $6225=9+8\times 765+4\times (3+21).$
\item [] $6226=(9+8\times 76\times 5+4^3)\times 2\times 1.$
\item [] $6227=(9+8\times 76\times 5+4^3)\times 2+1.$
\item [] $6228=9+8\times (765+4\times 3)+2+1.$
\item [] $6229=9\times 8+76\times (5+4)\times 3^2+1.$
\item [] $\mathit{6230=(9+87)\times 65-4\times 3+2\times 1.}$
\item [] $6231=9\times 87+6\times (5+43\times 21).$
\item [] $6232=(9+876+5)\times (4+3)+2\times 1.$
\item [] $6233=9+8\times (765+4+3^2\times 1).$
\item [] $6234=(98\times 7+6)\times (5+4)+3+2+1.$
\item [] $6235=(987+65\times 4)\times (3+2\times 1).$
\item [] $6236=(987+65\times 4)\times (3+2)+1.$
\item [] $6237=9\times (8+7+654+3+21).$
\item [] $6238=9+876\times 5+43^2\times 1.$
\item [] $6239=9+876\times 5+43^2+1.$
\item [] $6240=(9+87)\times (6+54+3+2\times 1).$
\item [] $6241=9+8\times (765+4+3^2+1).$
\item [] $6242=9\times (8+7)\times 6+5432\times 1.$
\item [] $6243=9+8\times 76+5^4\times 3^2+1.$
\item [] $6244=987\times 6+5\times 4^3+2^1.$
\item [] $6245=987\times 6+5\times 4^3+2+1.$
\item [] $6246=987\times 6+54\times 3\times 2\times 1.$
\item [] $6247=987\times 6+54\times 3\times 2+1.$
\item [] $\mathit{6248=9+8\times (7+6)\times 5\times 4\times 3-2+1.}$
\item [] $6249=(9+87)\times 65+4+3+2\times 1.$
\item [] $6250=(9+87)\times 65+4+3+2+1.$
\item [] $6251=(9+87)\times 65+4+3\times 2+1.$
\item [] $6252=987\times 6+5+4+321.$
\item [] $6253=(9+87)\times 65+4+3^2\times 1.$
\item [] $6254=(9+87)\times 65+4+3^2+1.$
\item [] $6255=(9+87)\times 65+4\times 3+2+1.$
\item [] $6256=(98+7)\times 6+5^4\times 3^2+1.$
\item [] $6257=9+8\times 765+4^3\times 2\times 1.$
\item [] $6258=9+8\times 765+4\times 32+1.$
\item [] $6259=(9\times 87+65\times 4)\times 3\times 2+1.$
\item [] $6260=(98\times 7+6)\times (5+4)+32\times 1.$
\item [] $6261=9+8\times 765+4\times (32+1).$
\item [] $6262=((9+8+765)\times 4+3)\times 2\times 1.$
\item [] $6263=987\times 6+5\times 4+321.$
\item [] $6264=(9+8+7)\times 65\times 4+3+21.$
\item [] $6265=(9+87)\times 65+4\times 3\times 2+1.$
\item [] $6266=9+8\times ((7+6)\times 5\times 4\times 3+2)+1.$
\item [] $6267=(9+8+7)\times 65\times 4+3^{(2+1)}.$
\item [] $6268=(9+87)\times 65+4+3+21.$
\item [] $6269=9\times (87+6)+5432\times 1.$
\item [] $6270=9\times (87+6)+5432+1.$
\item [] $6271=(9+8\times 7)\times 6\times 5+4321.$
\item [] $6272=(9+8+7)\times 65\times 4+32\times 1.$
\item [] $6273=(9+87)\times 65+4\times 3+21.$
\item [] $6274=(9+876)\times 5+43^2\times 1.$
\item [] $6275=(9+876)\times 5+43^2+1.$
\item [] $6276=98\times 7+65\times 43\times 2\times 1.$
\item [] $6277=98\times 7+65\times 43\times 2+1.$
\item [] $\mathit{6278=(9\times (8+7)+6+5)\times 43\times (2-1).}$
\item [] $6279=987\times 6+(5+4\times 3)\times 21.$
\item [] $6280=(9+87)\times 65+4\times (3^2+1).$
\item [] $6281=98\times 7\times 6+5\times (432+1).$
\item [] $6282=(987+6+54)\times 3\times 2\times 1.$
\item [] $6283=(987+6+54)\times 3\times 2+1.$
\item [] $6284=9+(87+6+5)\times 4^3+2+1.$
\item [] $6285=(9+87)\times 65+43+2\times 1.$
\item [] $6286=(9+87)\times 65+43+2+1.$
\item [] $6287=9+(8+(7+6)\times 5)\times 43\times 2\times 1.$
\item [] $6288=9+(8+76+5\times 43)\times 21.$
\item [] $6289=(9+87)\times 65+(4+3)^2\times 1.$
\item [] $6290=98+(7+65)\times 43\times 2\times 1.$
\item[]$\mbox{Increasing order}$
\item [] $6291=123+4^5\times 6+7+8+9.$
\item [] $6292=1+2\times (3\times 4^5+6)+(7+8)\times 9.$
\item [] $6293=12+(3+4^5)\times 6+7\times (8+9).$
\item [] $6294=12+3+4^5\times 6+(7+8)\times 9.$
\item [] $6295=1234\times 5+6+7\times (8+9).$
\item [] $6296=1+2\times (3+4+56\times 7\times 8)+9.$
\item [] $6297=1\times 23\times 45\times 6+78+9.$
\item [] $6298=1+23\times 45\times 6+78+9.$
\item [] $6299=12\times 3+4^5\times 6+7\times (8+9).$
\item [] $6300=1\times 2\times 3\times 4^5+67+89.$
\item [] $6301=1234\times 5+6\times 7+89.$
\item [] $6302=(12\times 3+4)\times 5+678\times 9.$
\item [] $6303=1\times 2\times 34\times 5+67\times 89.$
\item [] $6304=1+2\times 34\times 5+67\times 89.$
\item [] $6305=1\times 23+(4\times 5+678)\times 9.$
\item [] $6306=1\times 23\times 45\times 6+7+89.$
\item [] $6307=1+23\times 45\times 6+7+89.$
\item [] $6308=1^2\times 345+67\times 89.$
\item [] $6309=1234\times 5+67+8\times 9.$
\item [] $6310=1\times 2+345+67\times 89.$
\item [] $6311=1+2+345+67\times 89.$
\item [] $6312=12^3+4567+8+9.$
\item [] $6313=1+(2+3)^4+5678+9.$
\item [] $6314=1234\times 5+6\times (7+8+9).$
\item [] $6315=12\times 3+4^5\times 6+(7+8)\times 9.$
\item [] $6316=1+(2\times 3+4^5)\times 6+(7+8)\times 9.$
\item [] $6317=(12+3)\times (4+56)\times 7+8+9.$
\item [] $6318=12\times 3+(4\times 5+678)\times 9.$
\item [] $6319=1^2\times (34+5\times 6+7)\times 89.$
\item [] $6320=(1+2)\times (345\times 6+7)+89.$
\item [] $6321=(12+3+4^5)\times 6+78+9.$
\item [] $6322=1+(23\times (4+5\times 6)+7)\times 8+9.$
\item [] $6323=1\times 2^3\times 45+67\times 89.$
\item [] $6324=1+2^3\times 45+67\times 89.$
\item [] $6325=1+(23+45)\times (6+78+9).$
\item [] $6326=1234\times 5+67+89.$
\item [] $6327=1\times (2+3)\times 45+678\times 9.$
\item [] $6328=1+(2+3)\times 45+678\times 9.$
\item [] $6329=1\times 23\times 45\times 6+7\times (8+9).$
\item [] $6330=123\times 45+6+789.$
\item [] $6331=(1^2+3^4)\times (5+6)\times 7+8+9.$
\item [] $6332=123+4^5\times 6+7\times 8+9.$
\item [] $6333=1\times 2\times 3\times 4^5+(6+7+8)\times 9.$
\item [] $6334=1+2\times 3\times 4^5+(6+7+8)\times 9.$
\item [] $6335=(1+2)\times 34\times 56+7\times 89.$
\item [] $6336=12\times 3\times (4\times 5+67+89).$
\item [] $6337=1+2\times (3+4+5\times 6+7)\times 8\times 9.$
\item [] $6338=123\times (45+6)+7\times 8+9.$
\item [] $6339=1+2+(3\times 4\times 5+6)\times (7+89).$
\item [] $6340=(1+2+3+4)\times (5+6+7\times 89).$
\item [] $6341=1\times 234+5+678\times 9.$
\item [] $6342=1+234+5+678\times 9.$
\item [] $6343=1+2\times (3\times 4^5+6\times (7+8)+9).$
\item [] $6344=(1^2+3)\times 4\times 56\times 7+8\times 9.$
\item [] $6345=123\times 45+6\times (7+8)\times 9.$
\item [] $6346=123+4^5\times 6+7+8\times 9.$
\item [] $6347=(12\times 3\times 4+5)\times 6\times 7+89.$
\item [] $6348=12\times 3^4+56\times (7+89).$
\item [] $6349=1\times 2\times (34+56\times 7\times 8)+9.$
\item [] $6350=1+2\times (34+56\times 7\times 8)+9.$
\item [] $6351=(1+23\times 45)\times 6+(7+8)\times 9.$
\item [] $6352=123\times (45+6)+7+8\times 9.$
\item [] $6353=(1+23+4^5)\times 6+7\times 8+9.$
\item [] $6354=12\times 3\times 45+6\times 789.$
\item [] $6355=1+2\times 3\times (45\times 6+789).$
\item [] $6356=(1+23+4)\times (5\times 6\times 7+8+9).$
\item [] $6357=12+(3+4\times (5+6))\times (7+8)\times 9.$
\item [] $6358=1\times 2\times (34+56\times 7\times 8+9).$
\item [] $6359=1234\times 5+(6+7+8)\times 9.$
\item [] $6360=123\times (45+6)+78+9.$
\item[]$\mbox{Decreasing order}$
\item [] $6291=98+(7+65)\times 43\times 2+1.$
\item [] $6292=(9\times (8+7+6\times (54+3)\times 2)+1).$
\item [] $6293=(9+8+76+5)\times 4^3+21.$
\item [] $6294=9+87\times (65+4+3)+21.$
\item [] $6295=(9+8\times (76+54))\times 3\times 2+1.$
\item [] $6296=9+8+(7+6)\times (5\times 4+3)\times 21.$
\item [] $6297=987\times 6+54+321.$
\item [] $6298=(98\times (7+6)+5^4\times 3)\times 2\times 1.$
\item [] $6299=(9\times (8+7)+6+5)\times 43+21.$
\item [] $6300=(9\times 8+7+6+5\times 43)\times 21.$
\item [] $6301=9\times (87+65\times 4+3)\times 2+1.$
\item [] $6302=9+(87+6+5)\times 4^3+21.$
\item [] $6303=(9+8+7)\times 65\times 4+3\times 21.$
\item [] $6304=(9+87)\times 65+43+21.$
\item [] $6305=9+8\times 7+65\times 4\times (3+21).$
\item [] $6306=(9+87)\times 65+4^3+2\times 1.$
\item [] $6307=(9+87)\times 65+4+3\times 21.$
\item [] $6308=9\times 87+65\times (4^3+21).$
\item [] $6309=(987+6)\times 5+4^3\times 21.$
\item [] $6310=9+(8+7)\times 6\times 5\times (4+3)\times 2+1.$
\item [] $6311=98\times 7+(6+5+4^3)^2\times 1.$
\item [] $6312=98\times 7+(6+5+4^3)^2+1.$
\item [] $6313=9\times 8+(7+65+4+3)^2\times 1.$
\item [] $6314=(9+8+7)\times (65\times 4+3)+2\times 1.$
\item [] $6315=(98+7+6)\times 54+321.$
\item [] $6316=987+(6\times 5+43)^2\times 1.$
\item [] $6317=9+876+5432\times 1.$
\item [] $6318=9+876+5432+1.$
\item [] $6319=9\times 8+7+65\times 4\times (3+21).$
\item [] $6320=(9\times 8+7)\times (65+4\times 3+2+1).$
\item [] $6321=9+8\times (765+4\times 3\times 2\times 1).$
\item [] $6322=(98+765\times 4+3)\times 2\times 1.$
\item [] $6323=(98+765\times 4+3)\times 2+1.$
\item [] $6324=(98+7)\times (6+54)+3+21.$
\item [] $6325=(9+87)\times 65+4^3+21.$
\item [] $6326=(9+87)\times 65+43\times 2\times 1.$
\item [] $6327=(9+87)\times 65+43\times 2+1.$
\item [] $6328=(9\times 8+76\times 5)\times (4\times 3+2)\times 1.$
\item [] $6329=9+8\times (765+4\times 3\times 2+1).$
\item [] $6330=98+76\times ((5+4)\times 3^2+1).$
\item [] $6331=9\times 8+7\times 6\times (5+(4\times 3)^2)+1.$
\item [] $6332=(98+7)\times (6+54)+32\times 1.$
\item [] $6333=9+87\times (65+4)+321.$
\item [] $6334=9+8+7+(6+5^4)\times (3^2+1).$
\item [] $6335=9+(87\times 6+5)\times 4\times 3+2\times 1.$
\item [] $6336=(987+65+4)\times 3\times 2\times 1.$
\item [] $6337=(987+65+4)\times 3\times 2+1.$
\item [] $6338=9\times 8\times (76+5+4+3)+2\times 1.$
\item [] $6339=9\times 8\times (76+5+4+3)+2+1.$
\item [] $6340=(9+8)\times 7\times 6+5^4\times 3^2+1.$
\item [] $6341=98\times 7+65\times (43\times 2+1).$
\item [] $6342=9\times (8+7+654)+321.$
\item [] $6343=9+8+76+5^4\times (3^2+1).$
\item [] $\mathit{6344=(-9+8)\times 76+5\times 4\times 321.}$
\item [] $6345=(9\times 8+7\times 65)\times 4\times 3+21.$
\item [] $6346=9+8\times (76+5\times 4^3)\times 2+1.$
\item [] $6347=987\times 6+5\times (4^3+21).$
\item [] $6348=(98+7+6)\times (54+3)+21.$
\item [] $6349=9+(8+7)\times 6+5^4\times (3^2+1).$
\item [] $6350=((98+7)\times 6+5)\times (4+3+2+1).$
\item [] $6351=(9+8\times (7+6)\times 5)\times 4\times 3+2+1.$
\item [] $6352=987\times 6+5\times 43\times 2\times 1.$
\item [] $6353=987\times 6+5\times 43\times 2+1.$
\item [] $6354=9+(87\times 6+5)\times 4\times 3+21.$
\item [] $6355=9\times (8\times 7\times 6+5+4\times 3)\times 2+1.$
\item [] $6356=98+7\times 6\times (5+(4\times 3)^2)\times 1.$
\item [] $6357=987\times 6+5\times (43\times 2+1).$
\item [] $6358=((9\times 87+6+5)\times 4+3)\times 2\times 1.$
\item [] $6359=987\times 6+5+432\times 1.$
\item [] $6360=987\times 6+5+432+1.$
\item[]$\mbox{Increasing order}$
\item [] $6361=(1^2+3)\times 4\times 56\times 7+89.$
\item [] $6362=123\times (4+5+6\times 7)+89.$
\item [] $6363=123+4^5\times 6+7+89.$
\item [] $6364=1+(234+5+6\times 78)\times 9.$
\item [] $6365=1\times 2\times 3\times 4^5+(6+7)\times (8+9).$
\item [] $6366=(12\times 3^4+5)\times 6+7\times 8\times 9.$
\item [] $6367=12^3+4567+8\times 9.$
\item [] $6368=1^2\times 3^4\times 5+67\times 89.$
\item [] $6369=1^2+3^4\times 5+67\times 89.$
\item [] $6370=1\times 2+3^4\times 5+67\times 89.$
\item [] $6371=1+2+3^4\times 5+67\times 89.$
\item [] $6372=1\times 2\times 3\times 45+678\times 9.$
\item [] $6373=1+2\times 3\times 45+678\times 9.$
\item [] $\mathit{6374=-1^2\times 34+(5+67)\times 89.}$
\item [] $6375=(1+23+4^5)\times 6+78+9.$
\item [] $6376=12\times 34+5+67\times 89.$
\item [] $6377=1\times 2+3+(4+5)\times (6+78\times 9).$
\item [] $6378=123+45\times (67+8\times 9).$
\item [] $6379=1+(23+4^5)\times 6+7+89.$
\item [] $6380=12+3^4\times 5+67\times 89.$
\item [] $6381=1\times 234+(5+678)\times 9.$
\item [] $6382=1+234+(5+678)\times 9.$
\item [] $6383=12\times (3+4)\times 5+67\times 89.$
\item [] $6384=12^3+4567+89.$
\item [] $6385=((1+2)^3+4^5)\times 6+7+8\times 9.$
\item [] $6386=123+4^5\times 6+7\times (8+9).$
\item [] $6387=123+(45+6\times 7)\times 8\times 9.$
\item [] $6388=1\times 2\times (34\times 5+6\times 7\times 8\times 9).$
\item [] $6389=1+23\times (45\times 6+7)+8+9.$
\item [] $6390=(1\times 23+4+5+678)\times 9.$
\item [] $6391=1234\times 5+(6+7)\times (8+9).$
\item [] $6392=123\times (45+6)+7\times (8+9).$
\item [] $6393=1\times 2\times (3+4+56\times 7)\times 8+9.$
\item [] $6394=1+2\times (3+4+56\times 7)\times 8+9.$
\item [] $6395=1\times 23+(4+5)\times (6+78\times 9).$
\item [] $6396=12\times (3+4\times 5+6+7\times 8\times 9).$
\item [] $6397=1+(2+34+5)\times (67+89).$
\item [] $6398=1\times 2\times (3+4\times (5+6\times 7)\times (8+9)).$
\item [] $6399=(1+23+4+5+678)\times 9.$
\item [] $6400=(1^2+3+4)\times (5+6+789).$
\item [] $6401=(1\times 23+4^5)\times 6+7\times (8+9).$
\item [] $6402=(12+3)\times 4\times 5+678\times 9.$
\item [] $6403=(1^2+3^4)\times (5+6)\times 7+89.$
\item [] $6404=(1+2)^(3+4)+5+6\times 78\times 9.$
\item [] $6405=123+(4\times 5+678)\times 9.$
\item [] $6406=1\times 2\times (3\times 4^5+6\times 7+89).$
\item [] $6407=1\times 2\times (3+456\times 7)+8+9.$
\item [] $6408=(1\times 2+3+4+56+7)\times 89.$
\item [] $6409=1234+(567+8)\times 9.$
\item [] $6410=1\times 2+3\times 4\times (5\times 6+7\times 8\times 9).$
\item [] $6411=1^2\times 3+4\times (5+6+7)\times 89.$
\item [] $6412=12^3+4+5\times (6+7)\times 8\times 9.$
\item [] $6413=1\times 2+3+4\times (5+6+7)\times 89.$
\item [] $6414=(1+2)\times 34\times 56+78\times 9.$
\item [] $6415=1\times 2\times (3+456\times 7+8)+9.$
\item [] $6416=1\times 2^3+4\times (5+6+7)\times 89.$
\item [] $6417=12^3+(4+56)\times 78+9.$
\item [] $6418=1234+(5+67)\times 8\times 9.$
\item [] $6419=1+2\times 3+4+(5+67)\times 89.$
\item [] $6420=12+3\times 4\times (5\times 6+7\times 8\times 9).$
\item [] $6421=1+2^3+4+(5+67)\times 89.$
\item [] $6422=1\times 2+3\times 4+(5+67)\times 89.$
\item [] $6423=1\times 23\times 4\times 5+67\times 89.$
\item [] $6424=1+23\times 4\times 5+67\times 89.$
\item [] $6425=(1+234+567)\times 8+9.$
\item [] $6426=12\times 3\times (4+5)+678\times 9.$
\item [] $6427=12+3+4+(5+67)\times 89.$
\item [] $6428=(1+23\times 4)\times 5+67\times 89.$
\item [] $6429=1+(2+3)\times 4+(5+67)\times 89.$
\item [] $6430=1^2+3+(4+5)\times 6\times 7\times (8+9).$
\item[]$\mbox{Decreasing order}$
\item [] $6361=98+7+6+5^4\times (3^2+1).$
\item [] $6362=98\times 7+6\times (5^4+321).$
\item [] $6363=9+87\times (6\times 5+43)+2+1.$
\item [] $6364=9\times 8+7\times 6+5^4\times (3^2+1).$
\item [] $6365=98\times 7+(6+5^4)\times 3^2\times 1.$
\item [] $6366=98\times 7+(6+5^4)\times 3^2+1.$
\item [] $\mathit{6367=9-8\times 7-6+5\times 4\times 321.}$
\item [] $6368=(9+87)\times 65+4\times 32\times 1.$
\item [] $6369=9\times 8\times (7+6)+5432+1.$
\item [] $6370=98\times (7+6+5\times 4+32\times 1).$
\item [] $6371=(9+8\times 7+65)\times (4+3)^2+1.$
\item [] $6372=98\times (7\times 6+5\times 4+3)+2\times 1.$
\item [] $6373=9\times 87+65\times 43\times 2\times 1.$
\item [] $6374=9\times 87+65\times 43\times 2+1.$
\item [] $6375=9\times (8+7)+65\times 4\times (3+21).$
\item [] $\mathit{6376=98\times (7+6)\times 5+4+3-2+1.}$
\item [] $6377=9+8\times (76+5\times (4\times 3)^2\times 1).$
\item [] $6378=((9+8+7\times 6)\times 54+3)\times 2\times 1.$
\item [] $6379=98\times (7+6)\times 5+4+3+2\times 1.$
\item [] $6380=9\times 8\times 76+5+43\times 21.$
\item [] $6381=9+8\times 765+4\times 3\times 21.$
\item [] $6382=9\times (8+76)+5^4\times 3^2+1.$
\item [] $6383=98\times (7+6)\times 5+4+3^2\times 1.$
\item [] $6384=98\times (7+6)\times 5+4+3^2+1.$
\item [] $6385=98\times (7+6)\times 5+4\times 3+2+1.$
\item [] $6386=((98+7)\times 6\times 5+43)\times 2\times 1.$
\item [] $6387=(9+87)\times 65+(4+3)\times 21.$
\item [] $6388=9+8+7\times 65\times (4+3)\times 2+1.$
\item [] $6389=(98+(7+65)\times 43)\times 2+1.$
\item [] $6390=9\times 8\times 7+654\times 3^2\times 1.$
\item [] $6391=9\times 8\times 7+654\times 3^2+1.$
\item [] $6392=(9+8)\times (7+(6\times 5\times 4+3)\times (2+1)).$
\item [] $6393=9+8\times 7\times 6\times (5+4+3^2+1).$
\item [] $6394=98\times (7+6)\times 5+4\times 3\times 2\times 1.$
\item [] $6395=9\times (8\times 7+654)+3+2\times 1.$
\item [] $6396=9\times (8+7\times 6+5^4)+321.$
\item [] $6397=98\times (7\times 6+5\times 4)+321.$
\item [] $6398=98\times (7+6)\times 5+4+3+21.$
\item [] $6399=9\times (8\times 7\times 6+54+321).$
\item [] $6400=9+(8+(7+6)\times 54)\times 3^2+1.$
\item [] $6401=98\times (7+6)\times 5+4+3^{(2+1)}.$
\item [] $6402=987\times 6+5\times 4\times (3+21).$
\item [] $6403=98\times (7+6)\times 5+4\times 3+21.$
\item [] $6404=(9+8+7\times 65\times (4+3))\times 2\times 1.$
\item [] $6405=987\times 6+(5\times 4+3)\times 21.$
\item [] $6406=98\times (7+6)\times 5+4+32\times 1.$
\item [] $6407=98\times (7+6)\times 5+4+32+1.$
\item [] $6408=987\times 6+54\times 3^2\times 1.$
\item [] $6409=987\times 6+54\times 3^2+1.$
\item [] $6410=(9\times 8\times 7+6\times 5)\times 4\times 3+2\times 1.$
\item [] $6411=(9\times 8\times 7+6\times 5)\times 4\times 3+2+1.$
\item [] $6412=(9+8)\times 76+5\times 4^(3+2)\times 1.$
\item [] $6413=9\times (8+76+5^4)+32\times 1.$
\item [] $6414=9\times 87+6+5^4\times 3^2\times 1.$
\item [] $6415=9\times 87+6+5^4\times 3^2+1.$
\item [] $6416=98\times (7+6)\times 5+43+2+1.$
\item [] $6417=98+(7+6)\times 54\times 3^2+1.$
\item [] $6418=9\times 8\times 76+5^4+321.$
\item [] $6419=9+8\times (7+65\times 4)\times 3+2\times 1.$
\item [] $6420=9\times (8\times 76+5)+43\times 21.$
\item [] $6421=98\times 7\times 6+(5+43)^2+1.$
\item [] $6422=9\times (8\times 7+654)+32\times 1.$
\item [] $6423=9\times (8\times 7+654)+32+1.$
\item [] $6424=98+76+5^4\times (3^2+1).$
\item [] $6425=987+6+5432\times 1.$
\item [] $6426=987+6+5432+1.$
\item [] $6427=9\times 87\times 6+54\times 32+1.$
\item [] $6428=9+(87+6)\times (5+4^3)+2\times 1.$
\item [] $6429=(9\times 8\times 7+6\times 5)\times 4\times 3+21.$
\item [] $\mathit{6430=9\times 8\times 7+6^5-43^2-1.}$
\item[]$\mbox{Increasing order}$
\item [] $6431=12\times (34+5)+67\times 89.$
\item [] $6432=12+3\times 4+(5+67)\times 89.$
\item [] $6433=1+2^3\times (4+5+6+789).$
\item [] $6434=1\times 2^3+(4+5)\times 6\times 7\times (8+9).$
\item [] $6435=1\times 23+4+(5+67)\times 89.$
\item [] $6436=1+23+4+(5+67)\times 89.$
\item [] $6437=(12+3+4)\times 5\times 67+8\times 9.$
\item [] $6438=1\times 2\times (3^4\times 5\times 6+789).$
\item [] $6439=1+2\times (3^4\times 5\times 6+789).$
\item [] $6440=1\times 2^3\times 4+(5+67)\times 89.$
\item [] $6441=1+2^3\times 4+(5+67)\times 89.$
\item [] $6442=1\times 2\times 34\times 5+678\times 9.$
\item [] $6443=1+2\times 34\times 5+678\times 9.$
\item [] $6444=12\times 3\times 45+67\times 8\times 9.$
\item [] $6445=1+2+34+(5+67)\times 89.$
\item [] $6446=1+2^{(3\times 4)}+5\times 6\times 78+9.$
\item [] $6447=1^2\times 345+678\times 9.$
\item [] $6448=1^2+345+678\times 9.$
\item [] $6449=1\times 2+345+678\times 9.$
\item [] $6450=1+2+345+678\times 9.$
\item [] $6451=(1+2\times 3)^4+5\times 6\times (7+8)\times 9.$
\item [] $6452=1\times 2\times (3^4+56\times 7\times 8+9).$
\item [] $6453=12\times (3+456+78)+9.$
\item [] $6454=12+34+(5+67)\times 89.$
\item [] $6455=1\times 2+3\times (4\times 56+7+8)\times 9.$
\item [] $6456=12\times (3+456+7+8\times 9).$
\item [] $6457=1+2\times (3\times 4^5+67+89).$
\item [] $6458=1+(2+3\times 4\times 5)\times (6+7)\times 8+9.$
\item [] $6459=12+345+678\times 9.$
\item [] $6460=123\times 4+5+67\times 89.$
\item [] $6461=1+23\times (45\times 6+7)+89.$
\item [] $6462=1\times 2^3\times 45+678\times 9.$
\item [] $6463=1+2^3\times 45+678\times 9.$
\item [] $6464=1+(2+3^4)\times (5+6)\times 7+8\times 9.$
\item [] $6465=12+3\times (4\times 56+7+8)\times 9.$
\item [] $6466=1\times 2^{(3\times 4)}+5\times 6\times (7+8\times 9).$
\item [] $6467=1+2^{(3\times 4)}+5\times 6\times (7+8\times 9).$
\item [] $6468=12^3+(4+56)\times (7+8\times 9).$
\item [] $6469=1\times 2\times 34\times (5+6\times (7+8))+9.$
\item [] $6470=1+2\times 34\times (5+6\times (7+8))+9.$
\item [] $6471=12^3+4+5+6\times 789.$
\item [] $6472=1+(2\times 3)^4+(567+8)\times 9.$
\item [] $6473=(1+2)\times 34\times 5+67\times 89.$
\item [] $6474=1\times 2\times 3\times (456+7\times 89).$
\item [] $6475=1+2\times 3\times (456+7\times 89).$
\item [] $6476=1\times 2\times 34+(5+67)\times 89.$
\item [] $6477=1+2\times 34+(5+67)\times 89.$
\item [] $6478=1^2+3\times (4\times 5\times 6+7)\times (8+9).$
\item [] $6479=1\times 2\times (3+456\times 7)+89.$
\item [] $6480=(1+2+34+5+678)\times 9.$
\item [] $6481=1^{23}+45\times 6\times (7+8+9).$
\item [] $6482=12^3+4\times 5+6\times 789.$
\item [] $6483=1\times (23+4+56)\times 78+9.$
\item [] $6484=1+(23+4+56)\times 78+9.$
\item [] $6485=1\times 2+3+45\times 6\times (7+8+9).$
\item [] $6486=1+2+3+45\times 6\times (7+8+9).$
\item [] $6487=1+2\times 3+45\times 6\times (7+8+9).$
\item [] $6488=1\times 2^3+45\times 6\times (7+8+9).$
\item [] $6489=1\times 2\times 3\times 4^5+6\times 7\times 8+9.$
\item [] $6490=12^3+4^5+6\times 7\times 89.$
\item [] $6491=1\times 2+3^4+(5+67)\times 89.$
\item [] $6492=1+2+3^4+(5+67)\times 89.$
\item [] $\mathit{6493=1+(2\times 3)^4\times 5+6+7+8-9.}$
\item [] $6494=(1\times 2\times 34\times 5+6\times 7)\times (8+9).$
\item [] $6495=12+3+45\times 6\times (7+8+9).$
\item [] $6496=1\times 2+34\times (56+(7+8)\times 9).$
\item [] $6497=(12+3+45+6+7)\times 89.$
\item [] $6498=(1+23+4\times 5+678)\times 9.$
\item [] $6499=1+2\times (3+456)\times 7+8\times 9.$
\item [] $6500=1\times 23\times 4+(5+67)\times 89.$
\item[]$\mbox{Decreasing order}$
\item [] $6431=(9+8)\times 7\times 6\times (5+4)+3+2\times 1.$
\item [] $6432=9\times 8\times 76+5\times 4^3\times (2+1).$
\item [] $6433=9\times 8\times 7+(65+4\times 3)^2\times 1.$
\item [] $6434=98\times (7+6)\times 5+43+21.$
\item [] $6435=987+6\times (5+43\times 21).$
\item [] $6436=98\times (7+6)\times 5+4^3+2\times 1.$
\item [] $6437=98\times (7+6)\times 5+4+3\times 21.$
\item [] $6438=9+8\times (7+65\times 4)\times 3+21.$
\item [] $6439=9\times (8+7+6)+5^4\times (3^2+1).$
\item [] $6440=(9+87+65)\times 4\times (3^2+1).$
\item [] $6441=9\times 8\times (76+5+4)+321.$
\item [] $6442=9\times 8+7\times 65\times (4+3^2+1).$
\item [] $6443=9\times 8+7\times 65\times (4+3)\times 2+1.$
\item [] $6444=(98+76+5)\times 4\times 3^2\times 1.$
\item [] $6445=(98+76+5)\times 4\times 3^2+1.$
\item [] $6446=98\times 7+6\times 5\times 4^3\times (2+1).$
\item [] $6447=9\times 8\times 7\times (6+5)+43\times 21.$
\item [] $6448=9+87\times (6\times 5+4+3)\times 2+1.$
\item [] $6449=9+8\times 7\times (6\times 5+4^3+21).$
\item [] $6450=9+8+7+6+5\times 4\times 321.$
\item [] $6451=987\times 6+(5\times 4+3)^2\times 1.$
\item [] $6452=987\times 6+(5\times 4+3)^2+1.$
\item [] $6453=9\times (8\times 7+654)+3\times 21.$
\item [] $6454=9+8\times 765+4+321.$
\item [] $6455=98\times (7+6)\times 5+4^3+21.$
\item [] $6456=98\times (7+6)\times 5+43\times 2\times 1.$
\item [] $6457=98\times (7+6)\times 5+43\times 2+1.$
\item [] $6458=(9+8)\times 7\times 6\times (5+4)+32\times 1.$
\item [] $6459=9\times 87+6\times (5^4+321).$
\item [] $6460=(9+8\times 7+6+5)\times (4^3+21).$
\item [] $6461=(9+8\times 7+6)\times (5+43\times 2\times 1).$
\item [] $6462=(9+8+76+5^4)\times 3^2\times 1.$
\item [] $6463=9\times 87+(6+5^4)\times 3^2+1.$
\item [] $6464=(9\times 8+76+54)\times 32\times 1.$
\item [] $6465=(9\times 8+76+54)\times 32+1.$
\item [] $6466=98\times (7+6)\times 5+4\times (3+21).$
\item [] $6467=987\times 6+543+2\times 1.$
\item [] $6468=(9+8+76+5\times 43)\times 21.$
\item [] $6469=(9+8)\times 76\times 5+4+3+2\times 1.$
\item [] $6470=(9+8)\times 76\times 5+4+3+2+1.$
\item [] $6471=(9+8)\times 76\times 5+4+3\times 2+1.$
\item [] $6472=(9+8\times 7+654)\times 3^2+1.$
\item [] $6473=(9+8)\times 76\times 5+4+3^2\times 1.$
\item [] $6474=(9+8)\times 76\times 5+4+3^2+1.$
\item [] $6475=(9+8)\times 76\times 5+4\times 3+2+1.$
\item [] $6476=9+8\times (765+43)+2+1.$
\item [] $6477=9+(87+6+5\times 43)\times 21.$
\item [] $6478=(9\times 8\times 7\times 6+5\times 43)\times 2\times 1.$
\item [] $6479=(9\times 8+7)\times 65+4^3\times 21.$
\item [] $6480=9\times 87\times 6+54\times (32+1).$
\item [] $6481=(98\times 7+6\times 5+4)\times 3^2+1.$
\item [] $6482=9+8\times (765+4)+321.$
\item [] $6483=9\times 8\times (7\times 6+5+43)+2+1.$
\item [] $6484=(9+8)\times 76\times 5+4\times 3\times 2\times 1.$
\item [] $6485=(9+8)\times 76\times 5+4\times 3\times 2+1.$
\item [] $6486=987\times 6+543+21.$
\item [] $6487=(9\times (8+7)+6)\times (5\times 4+3)\times 2+1.$
\item [] $6488=(9+8)\times 76\times 5+4+3+21.$
\item [] $6489=987\times 6+(5+4)\times 3\times 21.$
\item [] $6490=9+(8+7+6\times 5)\times (4\times 3)^2+1.$
\item [] $6491=9+8\times 7+6+5\times 4\times 321.$
\item [] $6492=(9+87)\times 65+4\times 3\times 21.$
\item [] $6493=(9+8)\times 76\times 5+4\times 3+21.$
\item [] $6494=9+8\times (765+43)+21.$
\item [] $6495=9\times (87+6+5^4)+32+1.$
\item [] $6496=98+7\times (6+5+43\times 21).$
\item [] $6497=9+8\times (765+43+2+1).$
\item [] $6498=98\times (7+6)\times 5+4^3\times 2\times 1.$
\item [] $6499=98\times (7+6)\times 5+4\times 32+1.$
\item [] $6500=(9+8)\times 76\times 5+4\times (3^2+1).$
\item[]$\mbox{Increasing order}$
\item [] $6501=(1+2)\times 34\times 56+789.$
\item [] $6502=1+(2^3+4)\times (5+67\times 8)+9.$
\item [] $6503=(1+2)^3\times 4\times 5+67\times 89.$
\item [] $6504=1+23+45\times 6\times (7+8+9).$
\item [] $6505=1+2\times (3+45)\times 67+8\times 9.$
\item [] $6506=12+34\times (56+(7+8)\times 9).$
\item [] $6507=12^3+45+6\times 789.$
\item [] $6508=1^2+3^4\times 5+678\times 9.$
\item [] $6509=1\times 2+3^4\times 5+678\times 9.$
\item [] $6510=1+2+3^4\times 5+678\times 9.$
\item [] $6511=1^2+3+(45+678)\times 9.$
\item [] $6512=1\times 2+3+(45+678)\times 9.$
\item [] $6513=12+3\times 4\times (5+67\times 8)+9.$
\item [] $6514=1+2\times 3+(45+678)\times 9.$
\item [] $6515=1234\times 5+6\times 7\times 8+9.$
\item [] $6516=1+2\times (3+456)\times 7+89.$
\item [] $6517=1+(23\times 4+5)\times 67+8+9.$
\item [] $6518=1+(2+3^4)\times 5+678\times 9.$
\item [] $6519=12+3^4\times 5+678\times 9.$
\item [] $6520=1\times 2^3\times (4\times 5+6+789).$
\item [] $6521=1\times 2\times (3+45)\times 67+89.$
\item [] $6522=12+3+(45+678)\times 9.$
\item [] $6523=1+2\times (3\times 4^5+(6+7+8)\times 9).$
\item [] $\mathit{6524=1+(2\times 3)^4\times 5+6\times 7-8+9.}$
\item [] $6525=(12+3+4+56)\times (78+9).$
\item [] $6526=(12+3+4)\times (5\times 67+8)+9.$
\item [] $6527=((1+2)^3+4)\times 5\times 6\times 7+8+9.$
\item [] $6528=(1+2)\times 34\times (5+6\times 7+8+9).$
\item [] $6529=1+(2^3+4+56)\times (7+89).$
\item [] $6530=1\times 23+(45+678)\times 9.$
\item [] $6531=1+23+(45+678)\times 9.$
\item [] $\mathit{6532=12^3-4\times 5+67\times 8\times 9.}$
\item [] $6533=(123\times 4+5)\times (6+7)+8\times 9.$
\item [] $6534=(1+2)^3+(45+678)\times 9.$
\item [] $6535=123+4+(5+67)\times 89.$
\item [] $6536=1\times 2^{(3+4)}+(5+67)\times 89.$
\item [] $6537=1+2^{(3+4)}+(5+67)\times 89.$
\item [] $\mathit{6538=(1+2^3)^4+56-7-8\times 9.}$
\item [] $6539=12\times (3+45)+67\times 89.$
\item [] $6540=12\times (3^4+56\times 7+8\times 9).$
\item [] $6541=1\times 23\times 4\times (56+7+8)+9.$
\item [] $6542=1+23\times 4\times (56+7+8)+9.$
\item [] $6543=12\times 3+(45+678)\times 9.$
\item [] $6544=1^2+3\times (4^5+(6+7)\times 89).$
\item [] $6545=(1+23)\times 45\times 6+7\times 8+9.$
\item [] $6546=1+2+3\times (4^5+(6+7)\times 89).$
\item [] $6547=(1\times 2+(3\times 45\times 6+7)\times 8)+9.$
\item [] $6548=(1+2+(3\times 45\times 6+7)\times 8)+9.$
\item [] $6549=123+(4+5)\times 6\times 7\times (8+9).$
\item [] $6550=(123\times 4+5)\times (6+7)+89.$
\item [] $6551=1\times (2\times 3)^4\times 5+6+7\times 8+9.$
\item [] $6552=(1\times 2+3+45+678)\times 9.$
\item [] $6553=1+2^3\times (4+5+6\times (7+8)\times 9).$
\item [] $6554=1\times 2+(3+45\times 6)\times (7+8+9).$
\item [] $6555=12\times 34+(5+678)\times 9.$
\item [] $6556=(1+(2\times 3)^4)\times 5+6+7\times 8+9.$
\item [] $6557=1\times 2345+6\times 78\times 9.$
\item [] $6558=1+2345+6\times 78\times 9.$
\item [] $6559=(1+23)\times 45\times 6+7+8\times 9.$
\item [] $6560=1234\times 5+6\times (7\times 8+9).$
\item [] $6561=12^3+4+5+67\times 8\times 9.$
\item [] $6562=1\times 23\times 4\times 5+678\times 9.$
\item [] $6563=1+23\times 4\times 5+678\times 9.$
\item [] $6564=1\times 23\times (45\times 6+7+8)+9.$
\item [] $6565=1+23\times (45\times 6+7+8)+9.$
\item [] $6566=1+(2\times 3)^4\times 5+6+7+8\times 9.$
\item [] $6567=(12+3^4)\times 5+678\times 9.$
\item [] $6568=1\times 2\times (3+456\times 7+89).$
\item [] $6569=1+2\times (3+456\times 7+89).$
\item [] $6570=12\times (34+5)+678\times 9.$
\item[]$\mbox{Decreasing order}$
\item [] $6501=9\times 8\times (7\times 6+5+43)+21.$
\item [] $6502=98\times (7+6)\times 5+4\times (32+1).$
\item [] $\mathit{6503=(9+8)\times 76\times 5+4^3-21.}$
\item [] $6504=9\times (8\times 7+6+5^4)+321.$
\item [] $6505=9\times 8+7+6+5\times 4\times 321.$
\item [] $6506=(9+8)\times 76\times 5+43+2+1.$
\item [] $6507=9\times (87+6+5^4+3+2\times 1).$
\item [] $6508=((9+8)\times 7\times 6+5+4)\times 3^2+1.$
\item [] $6509=(9+8)\times 76\times 5+(4+3)^2\times 1.$
\item [] $6510=987+(65\times 4+3)\times 21.$
\item [] $6511=9+876+5^4\times 3^2+1.$
\item [] $6512=987+65\times (4^3+21).$
\item [] $6513=9+8+76+5\times 4\times 321.$
\item [] $6514=98\times (7+6)\times 5+(4\times 3)^2\times 1.$
\item [] $6515=98\times (7+6)\times 5+(4\times 3)^2+1.$
\item [] $6516=9\times (87+6+5^4+3+2+1).$
\item [] $6517=98\times (7+6)\times 5+(4+3)\times 21.$
\item [] $6518=98\times 7+(6+5+4+3)^(2+1).$
\item [] $6519=9+(8+7)\times 6+5\times 4\times 321.$
\item [] $6520=(98+(7+6)\times 5)\times 4\times (3^2+1).$
\item [] $6521=9+8\times (765+(4+3)^2\times 1).$
\item [] $6522=9+87+6+5\times 4\times 321.$
\item [] $\mathit{6523=-98+7\times (6+5)\times 43\times 2-1.}$
\item [] $6524=(9+8)\times 76\times 5+43+21.$
\item [] $6525=9\times (8\times 76+54+3\times 21).$
\item [] $6526=(9+8)\times 76\times 5+4^3+2\times 1.$
\item [] $6527=(9+8)\times 76\times 5+4+3\times 21.$
\item [] $6528=9\times 8\times (7+65)+4^3\times 21.$
\item [] $6529=9\times 8\times (7+6)\times 5+43^2\times 1.$
\item [] $6530=9\times 8\times (7+6)\times 5+43^2+1.$
\item [] $6531=98+7+6+5\times 4\times 321.$
\item [] $6532=(9+8\times 7+6)\times (5+43\times 2+1).$
\item [] $6533=9+8\times (7+6)+5\times 4\times 321.$
\item [] $6534=9\times 8+7\times 6+5\times 4\times 321.$
\item [] $6535=9+87\times (65+4+3\times 2)+1.$
\item [] $6536=9+87\times (6+5+4^3)+2\times 1.$
\item [] $6537=9+87\times (6+5+4^3)+2+1.$
\item [] $6538=(9+8)\times (76\times 5+4)+3^2+1.$
\item [] $\mathit{6539=9+8+7+6\times 543\times 2-1.}$
\item [] $6540=9+8+7+6\times 543\times 2\times 1.$
\item [] $6541=9+8+7+6\times 543\times 2+1.$
\item [] $6542=9\times (8\times 76+5)+4^(3+2)+1.$
\item [] $6543=9\times 87+6\times 5\times 4^3\times (2+1).$
\item [] $6544=(9+87+6+5^4)\times 3^2+1.$
\item [] $6545=(9+8)\times 76\times 5+4^3+21.$
\item [] $6546=9+8+7+6\times (543\times 2+1).$
\item [] $6547=9+(87+65)\times 43+2\times 1.$
\item [] $6548=9+8+7\times (6\times 5+43\times 21).$
\item [] $6549=(98+7+6)\times (54+3+2\times 1).$
\item [] $6550=(98+7+6)\times (54+3+2)+1.$
\item [] $\mathit{6551=9+8-7+6543-2\times 1.}$
\item [] $6552=987\times 6+5^4+3+2\times 1.$
\item [] $6553=987\times 6+5^4+3+2+1.$
\item [] $6554=987\times 6+5^4+3\times 2+1.$
\item [] $6555=9+87\times (6+5+4^3)+21.$
\item [] $6556=987\times 6+5^4+3^2\times 1.$
\item [] $6557=987\times 6+5^4+3^2+1.$
\item [] $6558=9\times 8\times 76+543\times 2\times 1.$
\item [] $6559=9\times 8\times 76+543\times 2+1.$
\item [] $6560=98+7\times 6+5\times 4\times 321.$
\item [] $6561=9+8\times 765+432\times 1.$
\item [] $6562=9+8\times 765+432+1.$
\item [] $6563=987\times 6+5\times 4^3\times 2+1.$
\item [] $6564=(9+8+7)\times 6+5\times 4\times 321.$
\item [] $6565=(9+87)\times 65+4+321.$
\item [] $6566=9+(87+65)\times 43+21.$
\item [] $6567=9+87\times 65+43\times 21.$
\item [] $6568=9\times 8+76+5\times 4\times 321.$
\item [] $6569=9+8+7+6543+2\times 1.$
\item [] $6570=9+8+7+6543+2+1.$
\item[]$\mbox{Increasing order}$
\item [] $6571=1\times (23\times 4+5)\times 67+8\times 9.$
\item [] $6572=12^3+4\times 5+67\times 8\times 9.$
\item [] $6573=(1+2+3)^4\times 5+6+78+9.$
\item [] $6574=1+(2\times 3)^4\times 5+6+78+9.$
\item [] $6575=1+2\times (3\times (4^5+6\times 7)+89).$
\item [] $6576=(1+23)\times 45\times 6+7+89.$
\item [] $6577=1+2^3\times (4\times 5\times 6+78\times 9).$
\item [] $6578=(12+34)\times (56+78+9).$
\item [] $6579=1+23\times (4+5\times 6\times 7+8\times 9).$
\item [] $6580=1+(2\times 34+5)\times 6\times (7+8)+9.$
\item [] $6581=1\times 2+(345+6\times 7)\times (8+9).$
\item [] $6582=(1+23)\times 4\times 5+678\times 9.$
\item [] $6583=1+(2\times 3)^4\times 5+6+7+89.$
\item [] $6584=(1\times 2)^3\times (4+5\times 6+789).$
\item [] $6585=1+2^3\times (4+5\times 6+789).$
\item [] $6586=(1+2+34+5\times 6+7)\times 89.$
\item [] $6587=(1+(2\times 3)^4)\times 5+6+7+89.$
\item [] $6588=12\times (34+5+6+7\times 8\times 9).$
\item [] $6589=1+(23\times 4+5)\times 67+89.$
\item [] $6590=1+2+(3+4)\times (5+(6+7)\times 8\times 9).$
\item [] $6591=12+(345+6\times 7)\times (8+9).$
\item [] $\mathit{6592=12^3\times 4-56\times 7+8\times 9.}$
\item [] $6593=1\times (2+3)^4+5+67\times 89.$
\item [] $6594=1+(2+3)^4+5+67\times 89.$
\item [] $6595=1+(2\times 3)^4\times 5+6\times 7+8\times 9.$
\item [] $6596=(1+2^3)^4+5+6+7+8+9.$
\item [] $6597=12^3+45+67\times 8\times 9.$
\item [] $6598=(123+4)\times 5+67\times 89.$
\item [] $6599=123\times 4+5+678\times 9.$
\item [] $6600=1^2\times 3\times 4\times (5+67\times 8+9).$
\item [] $6601=12^3+4+(5+67\times 8)\times 9.$
\item [] $6602=1+2\times (345+67)\times 8+9.$
\item [] $6603=123+45\times 6\times (7+8+9).$
\item [] $6604=12^3+4+56\times (78+9).$
\item [] $6605=1\times (2\times 3)^4\times 5+6+7\times (8+9).$
\item [] $6606=(12\times 3+4\times 5+678)\times 9.$
\item [] $6607=1+((2+3)^4+5+(6+7)\times 8)\times 9.$
\item [] $6608=12+34\times (5+(6+7+8)\times 9).$
\item [] $6609=(123+4+5)\times (6\times 7+8)+9.$
\item [] $6610=1234+56\times (7+89).$
\item [] $6611=(1+2+3)^4\times 5+6\times 7+89.$
\item [] $6612=(1+2)\times 34\times 5+678\times 9.$
\item [] $6613=(1\times 2+345+6\times 7)\times (8+9).$
\item [] $6614=1+(2+345+6\times 7)\times (8+9).$
\item [] $6615=(12+3+4+5\times 6)\times (7+8)\times 9.$
\item [] $6616=1+(23+4\times 5+6)\times (7+8)\times 9.$
\item [] $6617=1\times 2+(3+4)\times 5\times (6+7+8)\times 9.$
\item [] $6618=1\times 2\times 3\times 4^5+6\times (7+8\times 9).$
\item [] $6619=(1+2+3)^4\times 5+67+8\times 9.$
\item [] $6620=1234\times 5+(6\times 7+8)\times 9.$
\item [] $6621=1\times 2\times 3\times 4^5+6\times 78+9.$
\item [] $6622=1+2\times 3\times 4^5+6\times 78+9.$
\item [] $\mathit{6623=1-2\times 3^4-5+6789.}$
\item [] $6624=12\times 3\times (45+67+8\times 9).$
\item [] $6625=(123\times 4+5\times 67)\times 8+9.$
\item [] $6626=(1+(2\times 3)^4)\times 5+6+(7+8)\times 9.$
\item [] $6627=12+(3+4)\times 5\times (6+7+8)\times 9.$
\item [] $6628=(1^2+3)\times (4\times 56\times 7+89).$
\item [] $6629=(1+(2\times 3)^4)\times 5+6\times (7+8+9).$
\item [] $6630=123+(45+678)\times 9.$
\item [] $6631=1+2\times (34+5)\times (6+7+8\times 9).$
\item [] $6632=1\times 2+(3\times 4+5)\times 6\times (7\times 8+9).$
\item [] $6633=(1\times 2\times 3^4+567+8)\times 9.$
\item [] $6634=1+((2+3)^4+56+7\times 8)\times 9.$
\item [] $6635=(1+2+3)\times (4^5+67)+89.$
\item [] $6636=(1\times 2\times 3)^4\times 5+67+89.$
\item [] $6637=1+(2\times 3)^4\times 5+67+89.$
\item [] $6638=(12+3)\times 45+67\times 89.$
\item [] $6639=123\times 4+(5+678)\times 9.$
\item [] $6640=1\times (2+3^4)\times (56+7+8+9).$
\item[]$\mbox{Decreasing order}$
\item [] $6571=987\times 6+5^4+3+21.$
\item [] $6572=98\times 7+654\times 3^2\times 1.$
\item [] $6573=98\times 7+654\times 3^2+1.$
\item [] $6574=987\times 6+5^4+3^{(2+1)}.$
\item [] $6575=9\times 87\times 6+5^4\times 3+2\times 1.$
\item [] $6576=9\times 87\times 6+5^4\times 3+2+1.$
\item [] $6577=987+65\times 43\times 2\times 1.$
\item [] $6578=987+65\times 43\times 2+1.$
\item [] $6579=987\times 6+5^4+32\times 1.$
\item [] $6580=987\times 6+5^4+32+1.$
\item [] $6581=9+8\times 7+6\times 543\times 2\times 1.$
\item [] $6582=9+8\times 7+6\times 543\times 2+1.$
\item [] $6583=(9+87)\times 65+(4+3)^{(2+1)}.$
\item [] $\mathit{6584=98\times (7+6+54)-3+21.}$
\item [] $6585=(98\times 7+6+5^4)\times (3+2)\times 1.$
\item [] $6586=(98\times 7+6+5^4)\times (3+2)+1.$
\item [] $6587=9+8\times 7+6\times (543\times 2+1).$
\item [] $6588=9+8+7+6543+21.$
\item [] $6589=(9+8)\times 76\times 5+4^3\times 2+1.$
\item [] $6590=98\times (7+6+54)+3+21.$
\item [] $6591=(9+8)\times (76\times 5+4)+3\times 21.$
\item [] $6592=9\times (87\times 6+5)+43^2\times 1.$
\item [] $6593=9\times (87\times 6+5)+43^2+1.$
\item [] $6594=98+76+5\times 4\times 321.$
\item [] $6595=9\times 8+7+6\times 543\times 2\times 1.$
\item [] $6596=9\times 8+7+6\times 543\times 2+1.$
\item [] $6597=9\times (8+76+5^4+3+21).$
\item [] $6598=98\times (7+6+54)+32\times 1.$
\item [] $6599=98\times (7+6+54)+32+1.$
\item [] $6600=(9\times 8+76+5)\times 43+21.$
\item [] $6601=9\times 8+7+6\times (543\times 2+1).$
\item [] $6602=9\times 8+(7+6\times 543)\times 2\times 1.$
\item [] $6603=9\times 8+7\times (6\times 5+43\times 21).$
\item [] $6604=(9+8)\times 76\times 5+(4\times 3)^2\times 1.$
\item [] $6605=9+8\times 7+654\times (3^2+1).$
\item [] $6606=9\times (8\times 7+654+3+21).$
\item [] $6607=(9+8)\times 76\times 5+(4+3)\times 21.$
\item [] $6608=9\times 8+76\times (54+32\times 1).$
\item [] $6609=9\times 8+76\times (54+32)+1.$
\item [] $6610=987\times 6+5^4+3\times 21.$
\item [] $6611=9+8\times 7+6543+2+1.$
\item [] $6612=9+87+6\times 543\times 2\times 1.$
\item [] $6613=9+87+6\times 543\times 2+1.$
\item [] $6614=(9+8)\times (76\times 5+4+3+2)+1.$
\item [] $6615=98\times 7+(65+4\times 3)^2\times 1.$
\item [] $6616=98\times 7+(65+4\times 3)^2+1.$
\item [] $6617=9\times (8+7+6)\times 5\times (4+3)+2\times 1.$
\item [] $6618=987+6+5^4\times 3^2\times 1.$
\item [] $6619=987+6+5^4\times 3^2+1.$
\item [] $\mathit{6620=-9+87+6543-2+1.}$
\item [] $6621=98+7+6\times 543\times 2\times 1.$
\item [] $6622=98+7+6\times 543\times 2+1.$
\item [] $6623=(98+76+5)\times (4+32+1).$
\item [] $6624=9\times 8+7+6543+2\times 1.$
\item [] $6625=9\times 8+7+6543+2+1.$
\item [] $\mathit{6626=9876-(54+3)^2-1.}$
\item [] $6627=98+7+6\times (543\times 2+1).$
\item [] $6628=98+(7+6\times 543)\times 2\times 1.$
\item [] $6629=9+8\times 7+6543+21.$
\item [] $6630=(9+87)\times (65+4)+3+2+1.$
\item [] $6631=(9+87)\times (65+4)+3\times 2+1.$
\item [] $6632=9+8+7\times (6+5+4)\times 3\times 21.$
\item [] $6633=9+(87+6\times 5\times 4)\times 32\times 1.$
\item [] $6634=98+76\times (54+32\times 1).$
\item [] $6635=98+76\times (54+32)+1.$
\item [] $6636=(9\times 8+7)\times (6+54+3+21).$
\item [] $6637=9+(8\times 7+6\times 543)\times 2\times 1.$
\item [] $6638=(98+7\times 65)\times 4\times 3+2\times 1.$
\item [] $6639=(98+7\times 65)\times 4\times 3+2+1.$
\item [] $6640=9+8+7\times (6+5)\times 43\times 2+1.$
\item[]$\mbox{Increasing order}$
\item [] $6641=12\times 3\times (4\times 5+6)\times 7+89.$
\item [] $6642=1\times 234+(5+67)\times 89.$
\item [] $6643=1+234+(5+67)\times 89.$
\item [] $6644=1234\times 5+6\times (7+8\times 9).$
\item [] $6645=12\times (3\times 4+5+67\times 8)+9.$
\item [] $6646=1\times 2^{(3\times 4)}+5\times (6+7\times 8\times 9).$
\item [] $6647=1234\times 5+6\times 78+9.$
\item [] $6648=1^{23}\times 4^5\times 6+7\times 8\times 9.$
\item [] $6649=1^{23}+4^5\times 6+7\times 8\times 9.$
\item [] $6650=(1+2^3)^4+5+67+8+9.$
\item [] $6651=(12\times 3^4+5)\times 6+789.$
\item [] $6652=1^2+3+4^5\times 6+7\times 8\times 9.$
\item [] $6653=1\times 2\times 345+67\times 89.$
\item [] $6654=1+2\times 345+67\times 89.$
\item [] $6655=1+2\times 3+4^5\times 6+7\times 8\times 9.$
\item [] $6656=1\times 2^3+4^5\times 6+7\times 8\times 9.$
\item [] $6657=1+2^3+4^5\times 6+7\times 8\times 9.$
\item [] $6658=1^2+3\times (45\times 6+7)\times 8+9.$
\item [] $6659=12\times 3^4+5678+9.$
\item [] $6660=12\times (3+456+7+89).$
\item [] $6661=1+2\times (3\times 4^5+6)+7\times 8\times 9.$
\item [] $\mathit{6662=1-2+3-4+56\times 7\times (8+9).}$
\item [] $6663=12+3+4^5\times 6+7\times 8\times 9.$
\item [] $6664=1^{234}\times 56\times 7\times (8+9).$
\item [] $6665=1+2\times 34\times (5+6+78+9).$
\item [] $6666=1\times 2\times 3\times 4^5+6\times (78+9).$
\item [] $6667=12\times 3^4+5\times 67\times (8+9).$
\item [] $6668=1^{23}\times 4+56\times 7\times (8+9).$
\item [] $6669=12+3\times (45\times 6+7)\times 8+9.$
\item [] $6670=1+(2\times 3)^4\times 5+(6+7+8)\times 9.$
\item [] $6671=1\times 23+4^5\times 6+7\times 8\times 9.$
\item [] $6672=123\times 4\times 5+6\times 78\times 9.$
\item [] $6673=1+(2\times 3)^4+56\times (7+89).$
\item [] $6674=123\times 45+67\times (8+9).$
\item [] $6675=1+2\times 3+4+56\times 7\times (8+9).$
\item [] $6676=1^2\times 3\times 4+56\times 7\times (8+9).$
\item [] $6677=1+2^3+4+56\times 7\times (8+9).$
\item [] $6678=12\times (3+45)+678\times 9.$
\item [] $6679=1+2+3\times 4+56\times 7\times (8+9).$
\item [] $6680=1234\times 5+6+7\times 8\times 9.$
\item [] $6681=(12+34+5)\times (6\times 7+89).$
\item [] $6682=(1+2^3)^4+56+7\times 8+9.$
\item [] $6683=12\times 3\times 4\times 5+67\times 89.$
\item [] $6684=12\times 3+4^5\times 6+7\times 8\times 9.$
\item [] $6685=1+(2\times 3+4^5)\times 6+7\times 8\times 9.$
\item [] $\mathit{6686=-123+4\times 5+6789.}$
\item [] $6687=(123+4+(5+6)\times 7\times 8)\times 9.$
\item [] $6688=12+3\times 4+56\times 7\times (8+9).$
\item [] $6689=1\times 2\times 3\times 4^5+67\times 8+9.$
\item [] $6690=1+2\times 3\times 4^5+67\times 8+9.$
\item [] $6691=1\times 23+4+56\times 7\times (8+9).$
\item [] $6692=1234\times 5+6\times (78+9).$
\item [] $6693=1\times 2\times 3\times (4^5+6\times (7+8))+9.$
\item [] $6694=1+2\times 3\times (4^5+6\times (7+8))+9.$
\item [] $6695=(1+2)^3+4+56\times 7\times (8+9).$
\item [] $6696=1\times 23\times 4\times (5+67)+8\times 9.$
\item [] $6697=1+23\times 4\times (5+67)+8\times 9.$
\item [] $6698=1^2\times 34+56\times 7\times (8+9).$
\item [] $6699=1^2+34+56\times 7\times (8+9).$
\item [] $6700=1\times 2+34+56\times 7\times (8+9).$
\item [] $6701=1+2+34+56\times 7\times (8+9).$
\item [] $6702=(1+2^3+4^5)\times 6+7\times 8\times 9.$
\item [] $6703=1+2\times 3\times 4^5+(6+7\times 8)\times 9.$
\item [] $6704=(1+2^3)^4+56+78+9.$
\item [] $6705=(1+2^3)^4+5+67+8\times 9.$
\item [] $6706=(1\times 2^3)^4+5\times 6\times (78+9).$
\item [] $6707=123\times (4+5)\times 6+7\times 8+9.$
\item [] $6708=12+(3+4+5)\times (6+7\times 8)\times 9.$
\item [] $6709=1+(23+4\times 5)\times (67+89).$
\item [] $6710=12+34+56\times 7\times (8+9).$
\item[]$\mbox{Decreasing order}$
\item [] $6641=9+87+6543+2\times 1.$
\item [] $6642=9+87+6543+2+1.$
\item [] $6643=9\times 8+7+6543+21.$
\item [] $6644=9+8\times (765+4^3)+2+1.$
\item [] $6645=98+7+654\times (3^2+1).$
\item [] $6646=(9+8\times 7+6\times 543)\times 2\times 1.$
\item [] $6647=(9+8\times 7+6\times 543)\times 2+1.$
\item [] $6648=(9+87)\times (65+4)+3+21.$
\item [] $6649=(9+(8+7)\times (6+5\times 43))\times 2+1.$
\item [] $6650=98+7+6543+2\times 1.$
\item [] $6651=98+7+6543+2+1.$
\item [] $6652=9\times (8+7)+6\times 543\times 2+1.$
\item [] $6653=987\times 6+(5+4)^3+2\times 1.$
\item [] $6654=987\times 6+(5+4)^3+2+1.$
\item [] $6655=((9+8)\times (7+6)\times 5+4)\times 3\times 2+1.$
\item [] $6656=(9+87)\times (65+4)+32\times 1.$
\item [] $6657=(98+7\times 65)\times 4\times 3+21.$
\item [] $6658=(9+(8+7)\times 6+5)\times 4^3+2\times 1.$
\item [] $6659=(9+8)\times 7+654\times (3^2+1).$
\item [] $6660=9+87+6543+21.$
\item [] $6661=9\times (8+7\times 6+5\times 4^3)\times 2+1.$
\item [] $6662=9+8\times (765+4^3)+21.$
\item [] $6663=(98+7+6)\times 5\times 4\times 3+2+1.$
\item [] $6664=(9+8)\times 7+6543+2\times 1.$
\item [] $6665=(9+8)\times 7+6543+2+1.$
\item [] $6666=987+(6+5^4)\times 3^2\times 1.$
\item [] $6667=987+(6+5^4)\times 3^2+1.$
\item [] $6668=(9+8\times (7+6))\times (54+3+2)+1.$
\item [] $6669=98+7+6543+21.$
\item [] $6670=9\times 87+654\times 3^2+1.$
\item [] $\mathit{6671=(9+87)\times 65+432-1.}$
\item [] $6672=(9+87)\times 65+432\times 1.$
\item [] $6673=(9+87)\times 65+432+1.$
\item [] $6674=9\times (87+654)+3+2\times 1.$
\item [] $6675=9\times (87+654)+3+2+1.$
\item [] $6676=9\times (87+654)+3\times 2+1.$
\item [] $6677=(9+(8+7)\times 6+5)\times 4^3+21.$
\item [] $6678=9\times (87+654)+3^2\times 1.$
\item [] $6679=9+(87+654)\times 3^2+1.$
\item [] $6680=9\times (8+7)+6543+2\times 1.$
\item [] $6681=9\times (8+7)+6543+2+1.$
\item [] $6682=9\times 8+(7+654)\times (3^2+1).$
\item [] $6683=(9+8)\times 7+6543+21.$
\item [] $6684=9\times 8+76\times (54+32+1).$
\item [] $6685=(9+8)\times (76\times 5+4\times 3)+21.$
\item [] $\mathit{6686=-9-8\times 7+(6+5+4)^3\times 2+1.}$
\item [] $6687=9\times 8+7\times (6+5+4)\times 3\times 21.$
\item [] $6688=9+87\times 65+4^(3+2\times 1).$
\item [] $6689=9+87\times 65+4^(3+2)+1.$
\item [] $6690=(9+(8+7)\times (6+5)\times 4)\times (3^2+1).$
\item [] $6691=(9\times 8+7)\times 6\times 5+4321.$
\item [] $6692=(98+76\times 5)\times (4+3)\times 2\times 1.$
\item [] $6693=9\times (87+654)+3+21.$
\item [] $6694=9\times 8+7\times (6+5)\times 43\times 2\times 1.$
\item [] $6695=98\times (7+6)\times 5+4+321.$
\item [] $6696=(9+8+76)\times (5+4+3\times 21).$
\item [] $6697=((9\times 8+76)\times 5+4)\times 3^2+1.$
\item [] $6698=9\times 8\times 76+(5\times (4+3))^2+1.$
\item [] $6699=9\times (8+7)+6543+21.$
\item [] $6700=9+(87+6\times 543)\times 2+1.$
\item [] $6701=9\times (87+654)+32\times 1.$
\item [] $6702=9\times (87+654)+32+1.$
\item [] $6703=((9\times 8+765)\times 4+3)\times 2+1.$
\item [] $\mathit{6704=1-(2-34)\times 5\times 6\times 7-(8+9).}$
\item [] $6705=9+(87+6)\times (5+4+3\times 21).$
\item [] $6706=98\times (7+6)+5432\times 1.$
\item [] $6707=98\times (7+6)+5432+1.$
\item [] $6708=(9+87+6\times 543)\times 2\times 1.$
\item [] $6709=(9+87+6\times 543)\times 2+1.$
\item [] $6710=98+76\times (54+32+1).$
\item[]$\mbox{Increasing order}$
\item [] $6711=12+(3+4)\times (5+6)\times (78+9).$
\item [] $\mathit{6712=-1-2+3\times 4\times 567-89.}$
\item [] $6713=1\times 23\times 4\times (5+67)+89.$
\item [] $6714=1\times 23\times 45\times 6+7\times 8\times 9.$
\item [] $6715=1234\times 5+67\times 8+9.$
\item [] $6716=(12+34)\times (5+6+(7+8)\times 9).$
\item [] $6717=1\times (2+3)\times 4\times 5\times 67+8+9.$
\item [] $6718=1+(2+3)\times 4\times 5\times 67+8+9.$
\item [] $6719=(1+2)\times 34\times 5\times (6+7)+89.$
\item [] $6720=1\times 2^3\times (45+6+789).$
\item [] $6721=1+2^3\times (45+6+789).$
\item [] $6722=1\times 2+3\times 4\times (56+7\times 8\times 9).$
\item [] $6723=(1+23+45+678)\times 9.$
\item [] $6724=(12+3)\times 4+56\times 7\times (8+9).$
\item [] $6725=1\times 23\times 4\times (5\times (6+7)+8)+9.$
\item [] $6726=1+23\times 4\times (5\times (6+7)+8)+9.$
\item [] $6727=(1+2^3)^4+(5+6)\times 7+89.$
\item [] $6728=1234\times 5+(6+7\times 8)\times 9.$
\item [] $6729=1^{23}\times 4\times 5\times 6\times 7\times 8+9.$
\item [] $6730=1^{23}+4\times 5\times 6\times 7\times 8+9.$
\item [] $\mathit{6731=1-2+3\times 4\times 567-8\times 9.}$
\item [] $6732=1\times 2\times 34+56\times 7\times (8+9).$
\item [] $6733=1+2\times 34+56\times 7\times (8+9).$
\item [] $6734=1\times 2+3+4\times 5\times 6\times 7\times 8+9.$
\item [] $6735=1+2+3+4\times 5\times 6\times 7\times 8+9.$
\item [] $6736=1+2\times 3+4\times 5\times 6\times 7\times 8+9.$
\item [] $6737=1\times 2^3+4\times 5\times 6\times 7\times 8+9.$
\item [] $6738=1+2^3+4\times 5\times 6\times 7\times 8+9.$
\item [] $6739=1+2\times (345+6\times 7\times 8\times 9).$
\item [] $6740=1+23\times (4\times (56+7+8)+9).$
\item [] $6741=(12\times 34+5+6\times 7\times 8)\times 9.$
\item [] $6742=1\times 2^{(3+4)}\times 5+678\times 9.$
\item [] $6743=1+2^{(3+4)}\times 5+678\times 9.$
\item [] $6744=12+3+4\times 5\times 6\times 7\times 8+9.$
\item [] $6745=1^2\times 3^4+56\times 7\times (8+9).$
\item [] $6746=1234\times 5+6\times (7+89).$
\item [] $6747=1\times 2+3^4+56\times 7\times (8+9).$
\item [] $6748=12\times (3+4)+56\times 7\times (8+9).$
\item [] $6749=(1^{23}+4+56\times 7)\times (8+9).$
\item [] $6750=1^2\times (3\times 4\times 56+78)\times 9.$
\item [] $6751=1+2\times (34+5+6\times 7\times 8)\times 9.$
\item [] $6752=1\times 23+4\times 5\times 6\times 7\times 8+9.$
\item [] $6753=1+23+4\times 5\times 6\times 7\times 8+9.$
\item [] $6754=1^2+(3+4\times 5\times 6\times 7)\times 8+9.$
\item [] $6755=1\times 23+(4+56\times 7)\times (8+9).$
\item [] $6756=1\times 23\times 4+56\times 7\times (8+9).$
\item [] $6757=12+3^4+56\times 7\times (8+9).$
\item [] $6758=1\times 2+3\times 4\times (5+(6+7\times 8)\times 9).$
\item [] $6759=(12\times 34+5\times 67+8)\times 9.$
\item [] $6760=1^2+3\times 45\times (6\times 7+8)+9.$
\item [] $6761=123\times (4+5)\times 6+7\times (8+9).$
\item [] $6762=1\times 23\times (45\times 6+7+8+9).$
\item [] $6763=1+23\times (45\times 6+7+8+9).$
\item [] $6764=1^{23}\times (4+5+67)\times 89.$
\item [] $6765=12\times 3+4\times 5\times 6\times 7\times 8+9.$
\item [] $6766=(1+2)\times 34+56\times 7\times (8+9).$
\item [] $6767=1^{23}\times 4^5\times 6+7\times 89.$
\item [] $6768=1^{23}+4^5\times 6+7\times 89.$
\item [] $6769=1\times 2+3+(4+5+67)\times 89.$
\item [] $6770=1^2\times 3+4^5\times 6+7\times 89.$
\item [] $6771=123+4^5\times 6+7\times 8\times 9.$
\item [] $6772=1\times 2+3+4^5\times 6+7\times 89.$
\item [] $6773=1+2+3+4^5\times 6+7\times 89.$
\item [] $6774=12\times 34\times 5+6\times 789.$
\item [] $6775=1\times 2^3+4^5\times 6+7\times 89.$
\item [] $6776=1+2^3+4^5\times 6+7\times 89.$
\item [] $6777=(12+3)\times 45+678\times 9.$
\item [] $6778=1234+(5+6)\times 7\times 8\times 9.$
\item [] $6779=12+3+(4+5+67)\times 89.$
\item [] $6780=1+23+4\times (5\times 6\times 7\times 8+9).$
\item[]$\mbox{Decreasing order}$
\item [] $6711=9+87\times (65+4\times 3)+2+1.$
\item [] $6712=(9+8)\times 76\times 5+4\times 3\times 21.$
\item [] $6713=9\times 87+(65+4\times 3)^2+1.$
\item [] $6714=(98\times 7+6+54)\times 3^2\times 1.$
\item [] $6715=(98\times 7+6+54)\times 3^2+1.$
\item [] $6716=9+87\times 6\times 5+4^{(3\times 2)}+1.$
\item [] $6717=9\times 8\times (76+5+4\times 3)+21.$
\item [] $\mathit{6718=-9+8\times 7\times 6\times 5\times 4+3\times 2+1.}$
\item [] $\mathit{6719=9+8\times 7\times 6\times 5\times 4-3^2-1.}$
\item [] $6720=98+7\times (6+5)\times 43\times 2\times 1.$
\item [] $6721=98+7\times (6+5)\times 43\times 2+1.$
\item [] $6722=(98+7\times 6)\times (5+43)+2\times 1.$
\item [] $6723=9\times (87+654+3+2+1).$
\item [] $6724=(9+8)\times 76+5432\times 1.$
\item [] $6725=(9+8)\times 76+5432+1.$
\item [] $6726=(98+7+6\times 543)\times 2\times 1.$
\item [] $6727=(98+7+6\times 543)\times 2+1.$
\item [] $6728=9\times 8\times 76+(5^4+3)\times 2\times 1.$
\item [] $6729=9+87\times (65+4\times 3)+21.$
\item [] $6730=9+8\times 7\times (6+(54+3)\times 2)+1.$
\item [] $6731=9+(8+7+6)\times 5\times 4^3+2\times 1.$
\item [] $6732=9\times 8\times 76+5\times 4\times 3\times 21.$
\item [] $6733=9\times ((8\times 7+6)\times 5+4^3)\times 2+1.$
\item [] $6734=9+8\times 7\times 6\times 5\times 4+3+2\times 1.$
\item [] $6735=9+8\times 7\times 6\times 5\times 4+3\times 2\times 1.$
\item [] $6736=9+8\times 7\times 6\times 5\times 4+3\times 2+1.$
\item [] $6737=9+8\times (7\times (6+(54+3)\times 2)+1).$
\item [] $6738=9+8\times 7\times 6\times 5\times 4+3^2\times 1.$
\item [] $6739=9+8\times 7\times 6\times 5\times 4+3^2+1.$
\item [] $\mathit{6740=(98+7)\times 65-43\times 2+1.}$
\item [] $6741=(98+7\times 6)\times (5+43)+21.$
\item [] $6742=9\times 8+(7\times 6+5^4)\times (3^2+1).$
\item [] $6743=((9\times (87+6)+5)\times 4+3)\times 2+1.$
\item [] $6744=9\times 8\times 7+65\times 4\times (3+21).$
\item [] $6745=(9+8\times 7\times (6+54)+3)\times 2+1.$
\item [] $6746=(9+8\times (7+6\times (5+4^3))\times 2+1).$
\item [] $6747=987+6\times 5\times 4^3\times (2+1).$
\item [] $\mathit{6748=-987+6^5-43+2\times 1.}$
\item [] $6749=(9+8)\times (7+65+4+321).$
\item [] $6750=9+8+7+6+5\times 4^3\times 21.$
\item [] $6751=(9+87+654)\times 3^2+1.$
\item [] $6752=9+(8\times 76+5)\times (4+3\times 2+1).$
\item [] $6753=9+8\times 7\times 6\times 5\times 4+3+21.$
\item [] $6754=((9+8)\times 7+6\times 543)\times 2\times 1.$
\item [] $6755=((9+8)\times 7+6\times 543)\times 2+1.$
\item [] $6756=9+8\times (7\times 6\times 5\times 4+3)+2+1.$
\item [] $6757=((9+8)\times 7+6)\times 54+3\times 2+1.$
\item [] $6758=(9\times (8\times 7+6)+5)\times 4\times 3+2\times 1.$
\item [] $6759=9\times (87+654+3^2+1).$
\item [] $6760=(98\times 7+65)\times (4+3+2)+1.$
\item [] $6761=9\times 8\times 76+5+4\times 321.$
\item [] $6762=9+8\times 7\times 6\times 5\times 4+32+1.$
\item [] $6763=(9+8\times 7\times 6\times 5)\times 4+3\times 2+1.$
\item [] $6764=98\times (7\times 6+(5+4)\times 3)+2\times 1.$
\item [] $6765=9+8\times 7\times 6+5\times 4\times 321.$
\item [] $6766=(9+8\times 7\times 6\times 5)\times 4+3^2+1.$
\item [] $6767=(9\times 8+7\times (6+5)\times 43)\times 2+1.$
\item [] $6768=9\times 8\times 76+54\times (3+21).$
\item [] $6769=(9\times 8+7\times 6\times 5)\times 4\times 3\times 2+1.$
\item [] $6770=(9\times (8+7)+6)\times (5+43)+2\times 1.$
\item [] $6771=(98+7+6)\times (54+3\times 2+1).$
\item [] $6772=9+(8+7+6)\times (5\times 4^3+2)+1.$
\item [] $\mathit{6773=9\times 8\times 7\times 6+5^4\times 3\times 2-1.}$
\item [] $6774=9\times 8\times 7\times 6+5^4\times 3\times 2\times 1.$
\item [] $6775=9\times 8\times 7\times 6+5^4\times 3\times 2+1.$
\item [] $\mathit{6776=98-7\times 6+5\times 4^3\times 21.}$
\item [] $6777=9+8\times (76\times 5+43)\times 2\times 1.$
\item [] $6778=9+8\times (76\times 5+43)\times 2+1.$
\item [] $6779=9+8+7\times 6+5\times 4^3\times 21.$
\item [] $6780=(9+8\times 7\times 6\times 5)\times 4+3+21.$
\item[]$\mbox{Increasing order}$
\item [] $6781=1+(2^3+4)\times 5\times ((6+7)\times 8+9).$
\item [] $6782=12+3+4^5\times 6+7\times 89.$
\item [] $6783=(12+345+6\times 7)\times (8+9).$
\item [] $6784=1+(2\times 3+45+6)\times 7\times (8+9).$
\item [] $6785=1^2\times (3+4^5)\times 6+7\times 89.$
\item [] $6786=(1\times 23+4^5)\times 6+7\times 8\times 9.$
\item [] $6787=1\times 23+(4+5+67)\times 89.$
\item [] $6788=1\times 2\times (3+4^5)+6\times 789.$
\item [] $6789=1^{2345}\times 6789.$
\item [] $6790=1^{2345}+6789.$
\item [] $6791=1+23+4^5\times 6+7\times 89.$
\item [] $6792=1\times 2\times 345+678\times 9.$
\item [] $6793=1+2\times 345+678\times 9.$
\item [] $6794=1^{234}\times 5+6789.$
\item [] $6795=1^{234}+5+6789.$
\item [] $6796=(1^2+3)\times (4^5+(67+8)\times 9).$
\item [] $6797=12+(3+4^5)\times 6+7\times 89.$
\item [] $6798=1^{23}\times 4+5+6789.$
\item [] $6799=1234\times 5+6+7\times 89.$
\item [] $6800=(123+45\times 6+7)\times (8+9).$
\item [] $6801=1^2\times 3+4+5+6789.$
\item [] $6802=1^2+3+4+5+6789.$
\item [] $6803=1\times 2+3+4+5+6789.$
\item [] $6804=1+2+3+4+5+6789.$
\item [] $6805=1+2\times 3+4+5+6789.$
\item [] $6806=1\times 2^3+4+5+6789.$
\item [] $6807=1+2^3+4+5+6789.$
\item [] $6808=1\times 2+3\times 4+5+6789.$
\item [] $6809=1+2+3\times 4+5+6789.$
\item [] $6810=1^{23}+4\times 5+6789.$
\item [] $6811=1^2+3\times 4^5+6\times 7\times 89.$
\item [] $6812=1^2\times 3+4\times 5+6789.$
\item [] $6813=12+3+4+5+6789.$
\item [] $6814=1\times 2+3+4\times 5+6789.$
\item [] $6815=1+2+3+4\times 5+6789.$
\item [] $6816=1+2\times 3+4\times 5+6789.$
\item [] $6817=1\times 2^3+4\times 5+6789.$
\item [] $6818=12+3\times 4+5+6789.$
\item [] $6819=1+2\times 3\times 4+5+6789.$
\item [] $6820=1+2\times 3\times 4^5+(67+8)\times 9.$
\item [] $6821=1\times 23+4+5+6789.$
\item [] $6822=1+23+4+5+6789.$
\item [] $6823=1\times 2\times (3\times 4+5)+6789.$
\item [] $6824=12+3+4\times 5+6789.$
\item [] $6825=(1+2)^3+4+5+6789.$
\item [] $6826=1\times 2^3\times 4+5+6789.$
\item [] $6827=1+2^3\times 4+5+6789.$
\item [] $6828=1^2\times 34+5+6789.$
\item [] $6829=1^2+34+5+6789.$
\item [] $6830=1\times 2+34+5+6789.$
\item [] $6831=1+2+34+5+6789.$
\item [] $6832=1\times 23+4\times 5+6789.$
\item [] $6833=1+23+4\times 5+6789.$
\item [] $6834=12\times 3+4+5+6789.$
\item [] $6835=1^{23}+45+6789.$
\item [] $6836=12+(3+4)\times 5+6789.$
\item [] $6837=1^2\times 3+45+6789.$
\item [] $6838=1^2+3+45+6789.$
\item [] $6839=1\times 2+3+45+6789.$
\item [] $6840=12+34+5+6789.$
\item [] $6841=1+2\times 3+45+6789.$
\item [] $6842=1\times 2^3+45+6789.$
\item [] $6843=1+2^3+45+6789.$
\item [] $6844=1+2\times 3\times (4+5)+6789.$
\item [] $6845=12\times 3+4\times 5+6789.$
\item [] $6846=1^{23}\times 4^5\times 6+78\times 9.$
\item [] $6847=1^{23}+4^5\times 6+78\times 9.$
\item [] $\mathit{6848=1-2+3\times 4\times 5+6789.}$
\item [] $6849=12+3+45+6789.$
\item [] $6850=1+(2^3+4)\times 5+6789.$
\item[]$\mbox{Decreasing order}$
\item [] $6781=9+87\times 6+5^4\times (3^2+1).$
\item [] $6782=((9+8)\times 7+6)\times 54+32\times 1.$
\item [] $6783=9\times (87+6+5^4)+321.$
\item [] $\mathit{6784=9-8+7\times 6\times 54\times 3-21.}$
\item [] $6785=(9+8)\times 76\times 5+4+321.$
\item [] $6786=98\times ((7+6)\times 5+4)+3+21.$
\item [] $6787=(9\times (8+7)+6\times 543)\times 2+1.$
\item [] $6788=(9+8\times 7\times 6\times 5)\times 4+32\times 1.$
\item [] $6789=(9+8\times 7\times 6\times 5)\times 4+32+1.$
\item [] $6790=((9\times 87+65)\times 4+3)\times 2\times 1.$
\item [] $6791=9+8\times 7+6+5\times 4^3\times 21.$
\item [] $6792=9+8\times 7\times 6\times 5\times 4+3\times 21.$
\item [] $6793=9+8+7\times (65+43\times 21).$
\item [] $6794=98\times ((7+6)\times 5+4)+32\times 1.$
\item [] $6795=(98\times 7+65+4)\times 3^2\times 1.$
\item [] $6796=9+87\times (65+4+3^2)+1.$
\item [] $6797=98+7\times (6+5)\times (43\times 2+1).$
\item [] $6798=9+87\times (6+5\times 4)\times 3+2+1.$
\item [] $6799=(9+8+7+(6+5+4)^3)\times 2+1.$
\item [] $6800=(98+7+65)\times 4\times (3^2+1).$
\item [] $6801=9\times (8\times 76+5)+4\times 321.$
\item [] $6802=98\times (7+6)\times 5+432\times 1.$
\item [] $6803=98\times (7+6)\times 5+432+1.$
\item [] $6804=9\times (8+7+65+4)\times 3^2\times 1.$
\item [] $6805=9\times 8+7+6+5\times 4^3\times 21.$
\item [] $\mathit{6806=9-8+7\times 6\times 54\times 3+2-1.}$
\item [] $\mathit{6807=98\times 76-5\times 4\times 32-1.}$
\item [] $6808=((9+8)\times 7+65)\times (4+32+1).$
\item [] $6809=9+8\times (765+4^3+21).$
\item [] $6810=(9+8\times 7)\times 6+5\times 4\times 321.$
\item [] $6811=9\times (8+76)\times (5+4)+3\times 2+1.$
\item [] $\mathit{6812=-987+6^5+4\times 3\times 2-1.}$
\item [] $6813=9+8+76+5\times 4^3\times 21.$
\item [] $6814=(987+6)\times 5+43^2\times 1.$
\item [] $6815=(987+6)\times 5+43^2+1.$
\item [] $6816=9+87\times (6+5\times 4)\times 3+21.$
\item [] $6817=9+8\times (765+43\times 2\times 1).$
\item [] $6818=(98+7\times (6+5)\times 43)\times 2\times 1.$
\item [] $6819=9+(8+7)\times 6+5\times 4^3\times 21.$
\item [] $\mathit{6820=987\times 6-5+43\times 21.}$
\item [] $6821=9\times 8\times 76+5+4^3\times 21.$
\item [] $6822=9+87+6+5\times 4^3\times 21.$
\item [] $6823=9+8+7\times 6\times 54\times 3+2\times 1.$
\item [] $6824=9+8+7\times 6\times 54\times 3+2+1.$
\item [] $6825=(9\times 8+7\times 6\times 5+43)\times 21.$
\item [] $6826=(98+7)\times (6+54+3+2)+1.$
\item [] $\mathit{6827=98\times 76-5^4+3+2-1.}$
\item [] $6828=9\times 8\times 7\times (6+5)+4\times 321.$
\item [] $6829=9\times 8+7+(6+5+4)^3\times 2\times 1.$
\item [] $6830=987\times 6+5+43\times 21.$
\item [] $6831=98+7+6+5\times 4^3\times 21.$
\item [] $6832=(98+7+654)\times 3^2+1.$
\item [] $6833=9+8\times (7+6)+5\times 4^3\times 21.$
\item [] $6834=9\times 8+7\times 6+5\times 4^3\times 21.$
\item [] $6835=(98+7)\times 65+4+3+2+1.$
\item [] $6836=(98+7)\times 65+4+3\times 2+1.$
\item [] $6837=9\times (8+76)\times (5+4)+32+1.$
\item [] $6838=(98+7)\times 65+4+3^2\times 1.$
\item [] $6839=(98+7)\times 65+4+3^2+1.$
\item [] $6840=(98+7)\times 65+4\times 3+2+1.$
\item [] $6841=9\times 8+(7\times 6+5)\times (4\times 3)^2+1.$
\item [] $6842=9+8+7\times 65\times (4\times 3+2+1).$
\item [] $6843=(9\times 8+7\times 6)\times 5\times 4\times 3+2+1.$
\item [] $\mathit{6844=98-7+6^5-4^(3+2)+1.}$
\item [] $6845=(9+8)\times 7+6+5\times 4^3\times 21.$
\item [] $6846=(98+7+6+5\times 43)\times 21.$
\item [] $6847=9+87+(6+5+4)^3\times 2+1.$
\item [] $6848=9\times 8+7\times (65+43\times 21).$
\item [] $6849=(98+7)\times 65+4\times 3\times 2\times 1.$
\item [] $6850=(98+7)\times 65+4\times 3\times 2+1.$
\item[]$\mbox{Increasing order}$
\item [] $6851=1\times 2+3\times 4\times 5+6789.$
\item [] $6852=1+2+3\times 4\times 5+6789.$
\item [] $6853=1+2\times 3+4^5\times 6+78\times 9.$
\item [] $6854=(12+3)\times 4+5+6789.$
\item [] $6855=1+2^3+4^5\times 6+78\times 9.$
\item [] $6856=1+2\times 3+456\times (7+8)+9.$
\item [] $6857=1\times 23+45+6789.$
\item [] $6858=1+23+45+6789.$
\item [] $6859=1\times 2\times (3+4)\times 5+6789.$
\item [] $6860=1+2\times (3+4)\times 5+6789.$
\item [] $6861=12+3\times 4\times 5+6789.$
\item [] $6862=1\times 2\times 34+5+6789.$
\item [] $6863=1+2\times 34+5+6789.$
\item [] $6864=(1+2+3\times 4)\times 5+6789.$
\item [] $6865=1^2+(3+4^5)\times 6+78\times 9.$
\item [] $6866=1\times 2\times (34+5\times 678+9).$
\item [] $6867=1\times 2\times (34+5)+6789.$
\item [] $6868=1+2\times (34+5)+6789.$
\item [] $6869=1\times 23+4^5\times 6+78\times 9.$
\item [] $6870=12\times 3+45+6789.$
\item [] $6871=1+(2\times 3)^4\times 5+6\times (7\times 8+9).$
\item [] $6872=1\times 23+456\times (7+8)+9.$
\item [] $6873=(1+2)^3+4^5\times 6+78\times 9.$
\item [] $6874=1+(2^3\times 4+56)\times 78+9.$
\item [] $6875=1^2\times 3^4+5+6789.$
\item [] $6876=1^2+3^4+5+6789.$
\item [] $6877=1\times 2+3^4+5+6789.$
\item [] $6878=1+2+3^4+5+6789.$
\item [] $6879=1+2+3\times 4\times 567+8\times 9.$
\item [] $6880=1\times 23+4+(5+6)\times 7\times 89.$
\item [] $6881=1+23+4+(5+6)\times 7\times 89.$
\item [] $6882=12\times 3+4^5\times 6+78\times 9.$
\item [] $6883=1+2\times 3\times (4^5+6)+78\times 9.$
\item [] $6884=1234\times 5+6\times 7\times (8+9).$
\item [] $6885=1\times 2\times (3+45)+6789.$
\item [] $6886=1\times 23\times 4+5+6789.$
\item [] $6887=12+3^4+5+6789.$
\item [] $6888=12+3\times 4\times 567+8\times 9.$
\item [] $6889=1\times (2+3)\times 4\times 5+6789.$
\item [] $6890=(1+23)\times 4+5+6789.$
\item [] $6891=1+2+3+(45+6)\times (7+8)\times 9.$
\item [] $6892=1+2\times 3+(45+6)\times (7+8)\times 9.$
\item [] $6893=1^2\times 3\times 4\times 567+89.$
\item [] $6894=(1+2)\times (3+4)\times 5+6789.$
\item [] $6895=1\times 2+3\times 4\times 567+89.$
\item [] $6896=(1+2)\times 34+5+6789.$
\item [] $6897=(12+3^4+5\times 6)\times 7\times 8+9.$
\item [] $6898=1\times 234+56\times 7\times (8+9).$
\item [] $6899=12+34+(5+6)\times 7\times 89.$
\item [] $6900=123\times 4+(5+67)\times 89.$
\item [] $6901=1+2\times 3\times 4^5+(6+78)\times 9.$
\item [] $6902=(1+2)^3\times 4+5+6789.$
\item [] $6903=1^2\times 3+4\times 5\times (6\times 7\times 8+9).$
\item [] $6904=1^2+3+4\times 5\times (6\times 7\times 8+9).$
\item [] $6905=12+3\times 4\times 567+89.$
\item [] $6906=1+(23+4^5)\times 6+7\times 89.$
\item [] $6907=12^3+4+(567+8)\times 9.$
\item [] $6908=1\times 23+(45+6)\times (7+8)\times 9.$
\item [] $6909=1\times 2\times 3\times 4\times 5+6789.$
\item [] $6910=1+2\times 3\times 4\times 5+6789.$
\item [] $6911=1\times 2+3\times 4\times (567+8)+9.$
\item [] $6912=1\times 23\times 45\times 6+78\times 9.$
\item [] $6913=1+2\times 3^4\times 5+678\times 9.$
\item [] $6914=(1+2\times 3\times 4)\times 5+6789.$
\item [] $6915=12+3+4\times 5\times (6\times 7\times 8+9).$
\item [] $6916=12^3+4+(5+67)\times 8\times 9.$
\item [] $6917=(1+2\times 3^4)\times 5+678\times 9.$
\item [] $6918=(1+23\times 45)\times 6+78\times 9.$
\item [] $6919=(12+3)\times 456+7+8\times 9.$
\item [] $6920=1\times (2\times 3+4)\times (5+678+9).$
\item[]$\mbox{Decreasing order}$
\item [] $6851=(9+87+6+5)\times 4^3+2+1.$
\item [] $6852=9\times ((8+76)\times (5+4)+3)+21.$
\item [] $6853=(98+7)\times 65+4+3+21.$
\item [] $6854=(9\times 8+7\times (6+5))\times (43+2+1).$
\item [] $6855=98+7+(6+5+4)^3\times 2\times 1.$
\item [] $6856=98+7+(6+5+4)^3\times 2+1.$
\item [] $6857=9+8+76\times (5+4^3+21).$
\item [] $6858=9\times 87\times 6+5\times 432\times 1.$
\item [] $6859=9\times 87\times 6+5\times 432+1.$
\item [] $6860=98+7\times 6+5\times 4^3\times 21.$
\item [] $6861=(98+7)\times 65+4+32\times 1.$
\item [] $6862=(98+7)\times 65+4+32+1.$
\item [] $6863=98\times (7+6+54+3)+2+1.$
\item [] $6864=(9+8+7)\times 6+5\times 4^3\times 21.$
\item [] $6865=9+8+(7\times 6\times 5+4)\times 32\times 1.$
\item [] $6866=9+8+(7\times 6\times 5+4)\times 32+1.$
\item [] $6867=9+8\times 76+5^4\times (3^2+1).$
\item [] $6868=987\times 6+5^4+321.$
\item [] $6869=(9+87+6+5)\times 4^3+21.$
\item [] $6870=9\times (8+7\times 6)+5\times 4\times 321.$
\item [] $6871=(98+7)\times 65+43+2+1.$
\item [] $6872=9+(8\times 7+(6+5+4)^3)\times 2+1.$
\item [] $6873=987+654\times 3^2\times 1.$
\item [] $6874=987+654\times 3^2+1.$
\item [] $6875=(98+7)\times 65+(4+3)^2+1.$
\item [] $6876=(9+8+7\times 6\times 54)\times 3+21.$
\item [] $6877=(9\times 8\times 7+65\times 4)\times 3^2+1.$
\item [] $6878=9\times 8+7\times 6\times 54\times 3+2\times 1.$
\item [] $6879=9\times 8+7\times 6\times 54\times 3+2+1.$
\item [] $6880=(98+7)\times 6+5^4\times (3^2+1).$
\item [] $6881=98\times (7+6+54+3)+21.$
\item [] $6882=9+87\times (65+4+3^2+1).$
\item [] $6883=9+87\times (65+4\times 3+2)+1.$
\item [] $6884=9+8+7\times (6\times 54+3)\times (2+1).$
\item [] $6885=(98+7\times 6+5^4)\times 3^2\times 1.$
\item [] $6886=(98+7\times 6+5^4)\times 3^2+1.$
\item [] $6887=9\times (8\times 7+6)\times 5+4^{(3\times 2)}+1.$
\item [] $6888=(9+8\times 7+65\times 4+3)\times 21.$
\item [] $6889=(98+7)\times 65+43+21.$
\item [] $6890=9+(8+7+65)\times 43\times 2+1.$
\item [] $6891=9+8+7+(6\times 54+3)\times 21.$
\item [] $6892=(98+7)\times 65+4+3\times 21.$
\item [] $6893=(9+8)\times 76\times 5+432+1.$
\item [] $6894=(9\times 8+7)\times 6+5\times 4\times 321.$
\item [] $6895=(9\times (8+7)+6+5^4)\times 3^2+1.$
\item [] $6896=(9\times 8\times (7\times 6+5)+4^3)\times 2\times 1.$
\item [] $6897=9\times 8+7\times 6\times 54\times 3+21.$
\item [] $6898=(9+(8+7+65)\times 43)\times 2\times 1.$
\item [] $6899=98\times (7+6)+5^4\times 3^2\times 1.$
\item [] $6900=98\times (7+6)+5^4\times 3^2+1.$
\item [] $6901=(9+(87\times 6+5^4)\times 3)\times 2+1.$
\item [] $6902=9+8+765\times (4+3+2\times 1).$
\item [] $6903=9\times (8\times 7\times 6+5\times 43\times 2+1).$
\item [] $6904=98+7\times 6\times 54\times 3+2\times 1.$
\item [] $6905=98+7\times 6\times 54\times 3+2+1.$
\item [] $6906=(9+8\times 7\times 6)\times 5\times 4+3+2+1.$
\item [] $6907=(9+8\times 7\times 6)\times 5\times 4+3\times 2+1.$
\item [] $6908=(9\times 8+7+(6+5+4)^3)\times 2\times 1.$
\item [] $6909=9\times (8+7+6)+5\times 4^3\times 21.$
\item [] $6910=(98+7)\times 65+4^3+21.$
\item [] $6911=(98+7)\times 65+43\times 2\times 1.$
\item [] $6912=(9+8+7+6\times 5)\times 4\times 32\times 1.$
\item [] $6913=(9+87+6\times 5\times 4)\times 32+1.$
\item [] $6914=(9+87)\times (65+4+3)+2\times 1.$
\item [] $6915=(9+87)\times (65+4+3)+2+1.$
\item [] $6916=987+(65+4\times 3)^2\times 1.$
\item [] $6917=(9+8)\times 76+5^4\times 3^2\times 1.$
\item [] $6918=9\times 8\times (76+5\times 4)+3+2+1.$
\item [] $6919=9\times 8\times (76+5\times 4)+3\times 2+1.$
\item [] $6920=9\times 8+(7\times 6\times 5+4)\times 32\times 1.$
\item[]$\mbox{Increasing order}$
\item [] $6921=1234+5678+9.$
\item [] $6922=1\times 2^{(3+4)}+5+6789.$
\item [] $6923=1+2^{(3+4)}+5+6789.$
\item [] $6924=1^2\times 3\times 45+6789.$
\item [] $6925=1^2+3\times 45+6789.$
\item [] $6926=1\times 2+3\times 45+6789.$
\item [] $6927=1+2+3\times 45+6789.$
\item [] $6928=(1^2+3)\times (4^5+6+78\times 9).$
\item [] $6929=(1+2\times 3)\times 4\times 5+6789.$
\item [] $6930=1\times 2\times 3\times (4^5+6\times 7+89).$
\item [] $6931=1+2\times 3\times (4^5+6\times 7+89).$
\item [] $6932=123+4\times 5+6789.$
\item [] $6933=12\times (3+4+5)+6789.$
\item [] $6934=1^{23}+4^5\times 6+789.$
\item [] $6935=1^2+3^4+(5+6)\times 7\times 89.$
\item [] $6936=12+3\times 45+6789.$
\item [] $6937=1+2\times 3456+7+8+9.$
\item [] $6938=12\times 3\times 4+5+6789.$
\item [] $6939=1+2+3+4^5\times 6+789.$
\item [] $6940=12\times 3^4+5+67\times 89.$
\item [] $6941=1\times 2^3+4^5\times 6+789.$
\item [] $6942=1+2^3+4^5\times 6+789.$
\item [] $6943=1\times 2\times (3456+7)+8+9.$
\item [] $6944=1+2\times (3456+7)+8+9.$
\item [] $6945=1\times 23\times 4+(5+6)\times 7\times 89.$
\item [] $6946=12+3^4+(5+6)\times 7\times 89.$
\item [] $6947=12^3\times 4+5+6+7+8+9.$
\item [] $6948=12+3+4^5\times 6+789.$
\item [] $6949=1\times 2^3\times 4\times 5+6789.$
\item [] $6950=1+2^3\times 4\times 5+6789.$
\item [] $6951=1^2\times (3+4^5)\times 6+789.$
\item [] $6952=1+2\times (3456+7+8)+9.$
\item [] $6953=1\times 2+(3+4^5)\times 6+789.$
\item [] $6954=(1+2^3\times 4)\times 5+6789.$
\item [] $6955=1\times 23\times (4\times 56+78)+9.$
\item [] $6956=1\times 2\times 3^4+5+6789.$
\item [] $6957=123+45+6789.$
\item [] $6958=1+(2\times 3)^4\times 5+6\times 78+9.$
\item [] $6959=1^2\times 34\times 5+6789.$
\item [] $6960=1^2+34\times 5+6789.$
\item [] $6961=1\times 2+34\times 5+6789.$
\item [] $6962=1+2+34\times 5+6789.$
\item [] $6963=12+(3+4^5)\times 6+789.$
\item [] $6964=(1^2+34)\times 5+6789.$
\item [] $6965=1234\times 5+6+789.$
\item [] $6966=12^3\times 4+5\times 6+7+8+9.$
\item [] $6967=1+2\times 3\times (45+6+78)\times 9.$
\item [] $6968=(1+2^3)^4+5\times 67+8\times 9.$
\item [] $6969=(1^2+3)\times 45+6789.$
\item [] $6970=1+(2+34)\times 5+6789.$
\item [] $6971=12+34\times 5+6789.$
\item [] $6972=123+456\times (7+8)+9.$
\item [] $6973=1^2+3\times 4\times (5+6\times (7+89)).$
\item [] $6974=(1+2+34)\times 5+6789.$
\item [] $6975=(1+2\times 345+6+78)\times 9.$
\item [] $6976=12^3\times 4+5+6\times 7+8+9.$
\item [] $6977=1\times 2\times 3456+7\times 8+9.$
\item [] $6978=1+2\times 3456+7\times 8+9.$
\item [] $6979=1+2\times 3\times (4^5+67+8\times 9).$
\item [] $6980=1234\times 5+6\times (7+8)\times 9.$
\item [] $6981=1\times (2^3+4^5)\times 6+789.$
\item [] $6982=1+(2^3+4^5)\times 6+789.$
\item [] $6983=1\times (2\times 3)^4+5678+9.$
\item [] $6984=1+(2\times 3)^4+5678+9.$
\item [] $6985=1+(23+4^5)\times 6+78\times 9.$
\item [] $6986=1+2^3\times 4\times (5\times 6\times 7+8)+9.$
\item [] $6987=1^{23}\times 4^5+67\times 89.$
\item [] $6988=1^{23}+4^5+67\times 89.$
\item [] $6989=(12\times 3+4)\times 5+6789.$
\item [] $6990=1^2\times 3+4^5+67\times 89.$
\item[]$\mbox{Decreasing order}$
\item [] $6921=9\times 8\times 76+(5+4^3)\times 21.$
\item [] $6922=9\times 8\times (7\times 6+54)+3^2+1.$
\item [] $6923=98+7\times 6\times 54\times 3+21.$
\item [] $6924=9+8\times (7+65)\times 4\times 3+2+1.$
\item [] $6925=(9+87+65)\times 43+2\times 1.$
\item [] $6926=98\times 7+65\times 4\times (3+21).$
\item [] $6927=(9+8\times 7\times 6)\times 5\times 4+3^{(2+1)}.$
\item [] $\mathit{6928=-9\times 87+6^5-4^3-2+1.}$
\item [] $6929=9+8+(7+65)\times 4\times (3+21).$
\item [] $6930=9\times 8\times 7+6+5\times 4\times 321.$
\item [] $6931=9\times 8\times 76+(5+4)^3\times 2+1.$
\item [] $6932=(9+8\times 7\times 6)\times 5\times 4+32\times 1.$
\item [] $6933=9+8+76\times (5+43\times 2\times 1).$
\item [] $6934=9+8+76\times (5+43\times 2)+1.$
\item [] $6935=((98\times 7+6)\times 5+4+3)\times 2+1.$
\item [] $6936=9\times 8\times (76+5\times 4)+3+21.$
\item [] $6937=(9+87\times 6+5^4)\times 3\times 2+1.$
\item [] $6938=98+76\times (5+4^3+21).$
\item [] $6939=9+(8+7+6)\times (5+4+321).$
\item [] $6940=9+87\times 6\times 5+4321.$
\item [] $6941=(9+8)\times (7+6)+5\times 4^3\times 21.$
\item [] $6942=9+87+(6+5\times 4^3)\times 21.$
\item [] $6943=(9+87+(6+5+4)^3)\times 2+1.$
\item [] $6944=9\times 8\times (76+5\times 4)+32\times 1.$
\item [] $6945=9\times 8\times (76+5\times 4)+32+1.$
\item [] $6946=98+(7\times 6\times 5+4)\times 32\times 1.$
\item [] $6947=98+(7\times 6\times 5+4)\times 32+1.$
\item [] $6948=9+87\times 65+4\times 321.$
\item [] $6949=(9\times 8\times 7+654)\times 3\times 2+1.$
\item [] $6950=(9\times 8+(7+6)\times (5\times 4+3)^2+1).$
\item [] $6951=9+87\times 6+5\times 4\times 321.$
\item [] $6952=987\times 6+5+4^(3+2)+1.$
\item [] $6953=(98+7)\times 65+4^3\times 2\times 1.$
\item [] $6954=(98+7)\times 65+4\times 32+1.$
\item [] $6955=(9+8\times 7)\times (6+5+4\times (3+21)).$
\item [] $6956=(9\times 8+76)\times (5\times 4+3^{(2+1)}).$
\item [] $6957=9\times 8+765\times (4+3+2\times 1).$
\item [] $6958=9\times 8+765\times (4+3+2)+1.$
\item [] $6959=((9\times 8+7)\times (6+5)\times 4+3)\times 2+1.$
\item [] $6960=9\times 8+7\times 6\times (54\times 3+2\times 1).$
\item [] $6961=9\times 8+7\times 6\times (54\times 3+2)+1.$
\item [] $6962=(98\times 7+65\times 43)\times 2\times 1.$
\item [] $6963=(98\times 7+65\times 43)\times 2+1.$
\item [] $6964=(9+8)\times 7\times 6+5^4\times (3^2+1).$
\item [] $6965=98+7\times (6\times 54+3)\times (2+1).$
\item [] $6966=9\times (87+654+32+1).$
\item [] $6967=9\times (8+765)+4+3+2+1.$
\item [] $6968=9+8+7\times (6\times 54\times 3+21).$
\item [] $6969=9+87\times (65+4\times 3+2+1).$
\item [] $6970=(98+7)\times 65+(4\times 3)^2+1.$
\item [] $6971=9\times (8+765)+4+3^2+1.$
\item [] $6972=98+7+(6\times 54+3)\times 21.$
\item [] $\mathit{6973=-9\times 87+6^5+4-3-21.}$
\item [] $6974=((9+87)\times 6+5)\times 4\times 3+2\times 1.$
\item [] $6975=9\times 8\times (76+5\times 4)+3\times 21.$
\item [] $6976=(9+87\times 6)\times 5+4321.$
\item [] $6977=9\times (8+765)+4\times (3+2)\times 1.$
\item [] $6978=9\times (8\times 7+6)+5\times 4\times 321.$
\item [] $6979=(9\times 8+7+6\times 5)\times 4^3+2+1.$
\item [] $\mathit{6980=-9\times 87+6^5-4\times 3-2+1.}$
\item [] $6981=9\times (8+765)+4\times 3\times 2\times 1.$
\item [] $6982=9\times (8+765)+4\times 3\times 2+1.$
\item [] $6983=9+8+(76+5)\times 43\times 2\times 1.$
\item [] $6984=98+765\times (4+3+2)+1.$
\item [] $6985=(9\times 8+7)\times 65+43^2+1.$
\item [] $6986=98+7\times 6\times (54\times 3+2\times 1).$
\item [] $6987=98+7\times 6\times (54\times 3+2)+1.$
\item [] $6988=9\times 8+76\times (5+43\times 2\times 1).$
\item [] $6989=9\times 8+76\times (5+43\times 2)+1.$
\item [] $6990=9\times (87+654)+321.$
\item[]$\mbox{Increasing order}$
\item [] $6991=1\times 2\times 3456+7+8\times 9.$
\item [] $6992=1\times 2+3+4^5+67\times 89.$
\item [] $6993=1+2+3+4^5+67\times 89.$
\item [] $6994=1+2\times 3+4^5+67\times 89.$
\item [] $6995=1\times 2^3+4^5+67\times 89.$
\item [] $6996=1+2^3+4^5+67\times 89.$
\item [] $6997=1+23\times (4+5)+6789.$
\item [] $6998=1\times 23\times 45+67\times 89.$
\item [] $6999=1\times 23\times 45\times 6+789.$
\item [] $7000=1+23\times 45\times 6+789.$
\item [] $7001=12^3\times 4+5+67+8+9.$
\item [] $7002=12+3+4^5+67\times 89.$
\item [] $7003=1+(2\times 3)^4\times 5+6\times (78+9).$
\item [] $7004=1^2\times (345+67)\times (8+9).$
\item [] $7005=(1+23\times 45)\times 6+789.$
\item [] $7006=12^3\times 4+(5+6)\times 7+8+9.$
\item [] $7007=12^3\times 4+5\times 6+7\times 8+9.$
\item [] $7008=1\times 2\times 3456+7+89.$
\item [] $7009=1+2\times 3456+7+89.$
\item [] $7010=12^3\times 4+5+6+78+9.$
\item [] $7011=1+23+4^5+67\times 89.$
\item [] $7012=1+(23+(4+5)\times (6+78))\times 9.$
\item [] $\mathit{7013=1+234\times 5\times 6-7+8-9.}$
\item [] $7014=1\times (2+3)\times 45+6789.$
\item [] $7015=1+(2+3)\times 45+6789.$
\item [] $7016=1+2\times (3456+7)+89.$
\item [] $7017=12\times (3+45+67\times 8)+9.$
\item [] $\mathit{7018=1\times 234-5+6789.}$
\item [] $7019=12^3\times 4+5+6+7+89.$
\item [] $7020=12+3\times 4\times (567+8+9).$
\item [] $7021=12^3\times 4+5\times 6+7+8\times 9.$
\item [] $7022=1+(2^3+45+6)\times 7\times (8+9).$
\item [] $7023=12\times 3+4^5+67\times 89.$
\item [] $7024=1^2+3+45\times (67+89).$
\item [] $7025=1\times 2+3+45\times (67+89).$
\item [] $7026=1+2+3+45\times (67+89).$
\item [] $7027=1+2\times 3+45\times (67+89).$
\item [] $7028=1\times 234+5+6789.$
\item [] $7029=1+234+5+6789.$
\item [] $7030=12^3\times 4+5+(6+7)\times 8+9.$
\item [] $7031=12^3\times 4+5+6\times 7+8\times 9.$
\item [] $7032=1+2+(34+56)\times 78+9.$
\item [] $7033=12^3\times 4+56+7\times 8+9.$
\item [] $7034=1+2\times (3456+7\times 8)+9.$
\item [] $7035=12+3+45\times (67+89).$
\item [] $\mathit{7036=12^3\times 4-5-6+(7+8)\times 9.}$
\item [] $7037=(1+2)\times 3^4+5+6789.$
\item [] $7038=12^3\times 4+5\times 6+7+89.$
\item [] $7039=1+(2\times 3)^4\times 5+(6+7\times 8)\times 9.$
\item [] $7040=1\times (2+3)\times 4\times (5\times 67+8+9).$
\item [] $7041=12+(34+56)\times 78+9.$
\item [] $7042=1\times 2\times (3456+7\times 8+9).$
\item [] $7043=(1+23)\times 45+67\times 89.$
\item [] $7044=1\times 234\times 5\times 6+7+8+9.$
\item [] $7045=1+234\times 5\times 6+7+8+9.$
\item [] $\mathit{7046=12^3\times 4-5+67+8\times 9.}$
\item [] $7047=12^3\times 4+56+7+8\times 9.$
\item [] $7048=12^3\times 4+5+6\times 7+89.$
\item [] $7049=12^3\times 4+5\times (6+7)+8\times 9.$
\item [] $7050=(1^2+3)^4+5+6789.$
\item [] $7051=1+(2+3+45)\times (6+(7+8)\times 9).$
\item [] $7052=(1+2)\times (3+4)\times 5\times 67+8+9.$
\item [] $7053=12\times (3\times 4+567+8)+9.$
\item [] $\mathit{7054=1+2+34\times 5\times 6\times 7-89.}$
\item [] $7055=12^3\times 4+56+78+9.$
\item [] $7056=123+4^5\times 6+789.$
\item [] $7057=1+2\times (3\times 4+5\times 6+7)\times 8\times 9.$
\item [] $7058=12^3\times 4+5+6+(7+8)\times 9.$
\item [] $7059=1\times 2\times 3\times 45+6789.$
\item [] $7060=1+2\times 3\times 45+6789.$
\item[]$\mbox{Decreasing order}$
\item [] $6991=9+(8+(76+5)\times 43)\times 2\times 1.$
\item [] $6992=9+(8+(76+5)\times 43)\times 2+1.$
\item [] $6993=(9+87+6+5+4)\times 3\times 21.$
\item [] $6994=9\times 8+(765+4)\times 3^2+1.$
\item [] $\mathit{6995=-9\times 87+6^5+4-3+2-1.}$
\item [] $6996=(9+87)\times 6+5\times 4\times 321.$
\item [] $6997=(9\times 8+7+6\times 5)\times 4^3+21.$
\item [] $6998=9\times (8+765+4)+3+2\times 1.$
\item [] $6999=9\times (8+765+4)+3+2+1.$
\item [] $7000=9\times (8+765+4)+3\times 2+1.$
\item [] $7001=9+(87+65)\times (43+2+1).$
\item [] $7002=9+(8+765+4)\times 3^2\times 1.$
\item [] $7003=9+(8+765+4)\times 3^2+1.$
\item [] $7004=(9+8)\times (7+(6+5+4)\times 3^{(2+1)}).$
\item [] $7005=98\times 7\times 6+(5+4)\times 321.$
\item [] $7006=9\times (8+76)+5^4\times (3^2+1).$
\item [] $7007=(987+6+5)\times (4+3)+21.$
\item [] $7008=987\times 6+543\times 2\times 1.$
\item [] $7009=987\times 6+543\times 2+1.$
\item [] $7010=98+(7+65)\times 4\times (3+21).$
\item [] $7011=9\times (87\times 6+5+4\times 3\times 21).$
\item [] $7012=(98+(7+6)\times 5)\times 43+2+1.$
\item [] $7013=9\times (8+7\times 6+(5+4)^3)+2\times 1.$
\item [] $7014=98+76\times (5+43\times 2\times 1).$
\item [] $7015=98+76\times (5+43\times 2)+1.$
\item [] $\mathit{7016=-9\times 87+6^5+4\times 3\times 2-1.}$
\item [] $7017=9\times (8+765+4)+3+21.$
\item [] $\mathit{7018=98\times 76-5\times 43\times 2\times 1.}$
\item [] $7019=98+(765+4)\times 3^2\times 1.$
\item [] $7020=9\times 8\times 7+6\times 543\times 2\times 1.$
\item [] $7021=9\times 8\times 7+6\times 543\times 2+1.$
\item [] $7022=(9\times 8+7\times 6\times 54)\times 3+2\times 1.$
\item [] $7023=9\times 87+65\times 4\times (3+21).$
\item [] $7024=9\times (8+765)+4+3\times 21.$
\item [] $7025=9\times (8+765+4)+32\times 1.$
\item [] $7026=9\times (8+765+4)+32+1.$
\item [] $7027=(9\times (8+7)\times (6+5\times 4)+3)\times 2+1.$
\item [] $7028=98+7\times 6\times (54\times 3+2+1).$
\item [] $7029=9\times (87+6+5^4+3\times 21).$
\item [] $7030=(98+(7+6)\times 5)\times 43+21.$
\item [] $7031=(9\times 8+7)\times (65+4\times 3\times 2\times 1).$
\item [] $7032=9+8\times 765+43\times 21.$
\item [] $7033=9\times 8\times (7+65)+43^2\times 1.$
\item [] $7034=9\times 8\times (7+65)+43^2+1.$
\item [] $7035=987+6\times (5+43)\times 21.$
\item [] $\mathit{7036=-9\times 87+6^5+4^3-21.}$
\item [] $7037=9+8\times 76+5\times 4\times 321.$
\item [] $7038=9\times 8+(76+5)\times 43\times 2\times 1.$
\item [] $7039=9\times 8+(76+5)\times 43\times 2+1.$
\item [] $7040=(9+87+6)\times (5+4^3)+2\times 1.$
\item [] $7041=(9\times 8+7\times 6\times 54)\times 3+21.$
\item [] $7042=9\times (8+765)+4^3+21.$
\item [] $7043=9\times (8+765)+43\times 2\times 1.$
\item [] $7044=9\times (8+765)+43\times 2+1.$
\item [] $\mathit{7045=9\times 8\times 7+6543-2\times 1.}$
\item [] $\mathit{7046=9\times 8\times 7+6543-2+1.}$
\item [] $7047=9\times (8+765+4+3+2+1).$
\item [] $7048=((98\times 7+6)\times 5+4^3)\times 2\times 1.$
\item [] $7049=9\times 8\times 7+6543+2\times 1.$
\item [] $7050=9\times 8\times 7+6543+2+1.$
\item [] $7051=(9\times (8+7)+6)\times (5+43+2)+1.$
\item [] $7052=(9+8+(7+6)\times 5)\times 43\times 2\times 1.$
\item [] $7053=9\times (8+765)+4\times (3+21).$
\item [] $\mathit{7054=-98\times 7+6^5-4-32\times 1.}$
\item [] $7055=(9\times 8+7+6)\times (5\times 4+3\times 21).$
\item [] $7056=98\times (7+6+54+3+2\times 1).$
\item [] $7057=98\times (7+6+5\times 4+3)\times 2+1.$
\item [] $7058=(9+(8+7\times 6+5)\times 4^3)\times 2\times 1.$
\item [] $7059=(9+87+6)\times (5+4^3)+21.$
\item [] $7060=9\times (8+7)\times 6+5^4\times (3^2+1).$
\item[]$\mbox{Increasing order}$
\item [] $7061=12^3\times 4+5+6\times (7+8+9).$
\item [] $7062=123\times (45+6)+789.$
\item [] $7063=1+2\times 3\times (4\times 5+(6+7)\times 89).$
\item [] $7064=12^3\times 4+56+7+89.$
\item [] $7065=12\times (3+4\times 5)+6789.$
\item [] $7066=12^3\times 4+5\times (6+7)+89.$
\item [] $7067=(1+2+34)\times (56+(7+8)\times 9).$
\item [] $7068=12+(3+4)\times (56+7\times 8)\times 9.$
\item [] $7069=1+(2+3\times 4\times 5)\times (6\times 7+8\times 9).$
\item [] $7070=123\times (4+5)+67\times 89.$
\item [] $7071=1+2\times (3456+7+8\times 9).$
\item [] $7072=12\times 34+56\times 7\times (8+9).$
\item [] $7073=12^3\times 4+5+67+89.$
\item [] $7074=(1+234)\times 5\times 6+7+8+9.$
\item [] $7075=1+2\times 3\times (4+5)\times (6\times 7+89).$
\item [] $7076=(1+2^3)^4+5+6+7\times 8\times 9.$
\item [] $7077=1\times 2\times (3456+78)+9.$
\item [] $7078=1+2\times (3456+78)+9.$
\item [] $7079=1\times 2345+6\times 789.$
\item [] $7080=1+2345+6\times 789.$
\item [] $7081=1+2\times 3\times 4^5+(6+7)\times 8\times 9.$
\item [] $7082=1+(23+45)\times (6+7)\times 8+9.$
\item [] $7083=12+3\times (4\times 567+89).$
\item [] $7084=(12+34)\times (5\times (6+7)+89).$
\item [] $7085=1\times 234\times 5\times 6+7\times 8+9.$
\item [] $7086=1+234\times 5\times 6+7\times 8+9.$
\item [] $7087=12^3\times 4+56+7\times (8+9).$
\item [] $7088=1+234+(5+6)\times 7\times 89.$
\item [] $7089=(12+3)\times 4\times 5+6789.$
\item [] $7090=1+(2\times 3+45)\times (67+8\times 9).$
\item [] $7091=(1+234\times 5)\times 6+7\times 8+9.$
\item [] $7092=12\times (3+4+567+8+9).$
\item [] $7093=(1+(2\times 3^4+5)\times 6)\times 7+8\times 9.$
\item [] $7094=1\times 2\times (3+4+5\times (6+78\times 9)).$
\item [] $7095=(1+2^3)^4+5\times 6+7\times 8\times 9.$
\item [] $7096=(1+2)\times 3^4+(5+6)\times 7\times 89.$
\item [] $\mathit{7097=(1+2\times 3)\times 4^5-6+7-8\times 9.}$
\item [] $7098=1\times 2\times (3\times 4^5+6\times 78+9).$
\item [] $7099=1\times 234\times 5\times 6+7+8\times 9.$
\item [] $7100=1+234\times 5\times 6+7+8\times 9.$
\item [] $7101=1\times 23\times 4\times (5+6)\times 7+8+9.$
\item [] $7102=1+23\times 4\times (5+6)\times 7+8+9.$
\item [] $7103=(1+23)\times 45\times 6+7\times 89.$
\item [] $7104=1+(2\times 3^4+5)\times 6\times 7+89.$
\item [] $7105=1+2\times (3456+7+89).$
\item [] $7106=1234\times 5+(6+7)\times 8\times 9.$
\item [] $7107=1\times 234\times 5\times 6+78+9.$
\item [] $7108=1+234\times 5\times 6+78+9.$
\item [] $7109=1\times (2\times 3)^4\times 5+6+7\times 89.$
\item [] $7110=123+4^5+67\times 89.$
\item [] $7111=1+2\times (345+6\times 7+8)\times 9.$
\item [] $7112=1\times 2+(34+56)\times (7+8\times 9).$
\item [] $7113=12\times 3\times (4+5)+6789.$
\item [] $7114=1+2\times 3\times 4\times (5\times 6+7)\times 8+9.$
\item [] $7115=(1+234)\times 5\times 6+7\times 8+9.$
\item [] $7116=1\times 234\times 5\times 6+7+89.$
\item [] $7117=1+234\times 5\times 6+7+89.$
\item [] $7118=12+(3^4\times 5+6+7)\times (8+9).$
\item [] $7119=12\times 3^4+(5+678)\times 9.$
\item [] $7120=(1^2+3+4+5+67)\times 89.$
\item [] $7121=(1+2^3)^4+56+7\times 8\times 9.$
\item [] $7122=(1+234\times 5)\times 6+7+89.$
\item [] $7123=(1\times 23+4+56\times 7)\times (8+9).$
\item [] $7124=(1+2)\times (3+4)\times 5\times 67+89.$
\item [] $7125=1\times 2^{(3\times 4)}+5+6\times 7\times 8\times 9.$
\item [] $7126=1^{23}\times 4^5+678\times 9.$
\item [] $7127=1^{23}+4^5+678\times 9.$
\item [] $7128=(123+45)\times 6\times 7+8\times 9.$
\item [] $7129=1\times 2\times 34\times 5+6789.$
\item [] $7130=1+2\times 34\times 5+6789.$
\item[]$\mbox{Decreasing order}$
\item [] $7061=98\times (7+6+5)\times 4+3+2\times 1.$
\item [] $7062=98\times (7+6+5)\times 4+3+2+1.$
\item [] $7063=98\times (7+6+5)\times 4+3\times 2+1.$
\item [] $7064=98+(76+5)\times 43\times 2\times 1.$
\item [] $7065=9+8\times 7\times 6+5\times 4^3\times 21.$
\item [] $7066=98\times (7+65)+4+3+2+1.$
\item [] $7067=98\times (7+65)+4+3\times 2+1.$
\item [] $7068=9\times 8\times 7+6543+21.$
\item [] $7069=98\times (7+65)+4+3^2\times 1.$
\item [] $7070=98\times (7+65)+4+3^2+1.$
\item [] $7071=98\times (7+65)+4\times 3+2+1.$
\item [] $7072=(9+8)\times (76\times 5+4+32\times 1).$
\item [] $7073=(9+(876+5)\times 4+3)\times 2+1.$
\item [] $7074=(9+8+765+4)\times 3^2\times 1.$
\item [] $7075=9\times (8\times 7\times 6+54+3)\times 2+1.$
\item [] $7076=98\times (7+65)+4\times (3+2)\times 1.$
\item [] $7077=(98+7)\times 65+4\times 3\times 21.$
\item [] $7078=9+8+(7+6)\times 543+2\times 1.$
\item [] $7079=9+8+(7+6+5+4)\times 321.$
\item [] $7080=98\times (7+65)+4\times 3\times 2\times 1.$
\item [] $7081=98\times (7+65)+4\times 3\times 2+1.$
\item [] $7082=9\times 8+(76+5^4)\times (3^2+1).$
\item [] $7083=9\times (8+765+4\times 3+2\times 1).$
\item [] $7084=98\times (7+65)+4+3+21.$
\item [] $7085=9\times (8+765)+4^3\times 2\times 1.$
\item [] $7086=9\times (8+765)+4\times 32+1.$
\item [] $7087=9\times (87+6)+5^4\times (3^2+1).$
\item [] $7088=98\times (7+6+5)\times 4+32\times 1.$
\item [] $7089=98\times (7+65)+4\times 3+21.$
\item [] $7090=98+76\times (5+43\times 2+1).$
\item [] $7091=(987+6+5\times 4)\times (3\times 2+1).$
\item [] $7092=98\times (7+65)+4+32\times 1.$
\item [] $7093=98\times (7+65)+4+32+1.$
\item [] $7094=9\times (8+(7+6)\times 5\times 4\times 3)+2\times 1.$
\item [] $7095=(98+(7+6)\times (5+4))\times (32+1).$
\item [] $7096=98\times (7+65)+4\times (3^2+1).$
\item [] $7097=9+8+(7+6)\times 543+21.$
\item [] $7098=(9\times (8+7)+6\times 5)\times 43+2+1.$
\item [] $7099=9\times 8+7+65\times 4\times 3^{(2+1)}.$
\item [] $7100=(98+7\times 6\times 54)\times 3+2\times 1.$
\item [] $7101=98\times (7+65)+43+2\times 1.$
\item [] $7102=98\times (7+65)+43+2+1.$
\item [] $7103=9+8+(7+6)\times (543+2)+1.$
\item [] $7104=9\times (8+765)+(4+3)\times 21.$
\item [] $7105=(9+87)\times (65+4+3+2)+1.$
\item [] $7106=(9\times 87+6)\times (5+4)+3+2\times 1.$
\item [] $7107=(98+7\times 6\times 54+3)\times (2+1).$
\item [] $7108=(9\times 87+6)\times (5+4)+3\times 2+1.$
\item [] $7109=9+(8+(7+6)\times 54)\times (3^2+1).$
\item [] $7110=9\times (8+76\times 5+4+3)\times 2\times 1.$
\item [] $7111=9\times (8\times 7+6)\times 5+4321.$
\item [] $7112=98\times 7+6+5\times 4\times 321.$
\item [] $7113=9+8\times (7\times 65+432+1).$
\item [] $7114=9+(8+7\times 6\times 5+4)\times 32+1.$
\item [] $7115=9+8+7\times (6+(5+43)\times 21).$
\item [] $7116=9\times 8\times (76+5)+4\times 321.$
\item [] $\mathit{7117=-98\times 7+6^5-4+32-1.}$
\item [] $7118=98+(7+6)\times 54\times (3^2+1).$
\item [] $7119=987\times 6+(54+3)\times 21.$
\item [] $7120=98\times (7+65)+43+21.$
\item [] $7121=9+8\times 7\times (6+5\times 4\times 3\times 2+1).$
\item [] $7122=98\times (7+65)+4^3+2\times 1.$
\item [] $7123=98\times (7+65)+4+3\times 21.$
\item [] $7124=((9\times 8+7)\times 6\times 5+4)\times 3+2\times 1.$
\item [] $7125=9\times (8+76)\times (5+4)+321.$
\item [] $7126=9\times 8\times 7\times 6+5+4^{(3\times 2)}+1.$
\item [] $7127=((9+876+5)\times 4+3)\times 2+1.$
\item [] $7128=9+(87+6+5\times 4)\times 3\times 21.$
\item [] $7129=9\times (8\times 7\times 6+5\times 4\times 3)\times 2+1.$
\item [] $7130=9\times 8\times (76+5\times 4+3)+2\times 1.$
\item[]$\mbox{Increasing order}$
\item [] $7131=1\times 2+3+4^5+678\times 9.$
\item [] $7132=1+2+3+4^5+678\times 9.$
\item [] $7133=1+2\times 3+4^5+678\times 9.$
\item [] $7134=1+234\times 5+67\times 89.$
\item [] $7135=1^2+345+6789.$
\item [] $7136=1\times 2+345+6789.$
\item [] $7137=1+2+345+6789.$
\item [] $7138=1+23\times 45+678\times 9.$
\item [] $7139=12^3\times 4+5\times 6\times 7+8+9.$
\item [] $7140=1+234\times 5\times 6+7\times (8+9).$
\item [] $7141=12+3+4^5+678\times 9.$
\item [] $7142=(1+2^3)^4+5+6\times (7+89).$
\item [] $7143=123+45\times (67+89).$
\item [] $7144=1\times 2^3\times (45\times 6+7\times 89).$
\item [] $7145=(123+45)\times 6\times 7+89.$
\item [] $7146=12+345+6789.$
\item [] $7147=1+2\times 3+(4+56)\times 7\times (8+9).$
\item [] $7148=1\times 2^3+(4+56)\times 7\times (8+9).$
\item [] $7149=1\times 23+4^5+678\times 9.$
\item [] $7150=1+23+4^5+678\times 9.$
\item [] $7151=123\times (45+6+7)+8+9.$
\item [] $7152=12\times (3\times 4+567+8+9).$
\item [] $7153=(1+2)^3+4^5+678\times 9.$
\item [] $7154=1+23\times (4\times 56+78+9).$
\item [] $7155=1\times 234\times 5\times 6+(7+8)\times 9.$
\item [] $7156=123\times 4+56\times 7\times (8+9).$
\item [] $7157=1^2\times 34\times 5\times 6\times 7+8+9.$
\item [] $7158=1^2+34\times 5\times 6\times 7+8+9.$
\item [] $7159=1\times 2+34\times 5\times 6\times 7+8+9.$
\item [] $7160=1+2+34\times 5\times 6\times 7+8+9.$
\item [] $7161=1+2+3+(4+5)\times (6+789).$
\item [] $7162=12\times 3+4^5+678\times 9.$
\item [] $7163=1\times 23+(4+56)\times 7\times (8+9).$
\item [] $7164=12\times 345+6\times 7\times 8\times 9.$
\item [] $7165=1+2\times (3\times 4+5\times 6\times 7\times (8+9)).$
\item [] $7166=1+2^{(3\times 4)}+(5+6\times 7\times 8)\times 9.$
\item [] $7167=1\times (2\times 3)^4\times 5+678+9.$
\item [] $7168=1+(2\times 3)^4\times 5+678+9.$
\item [] $7169=1\times 2345+67\times 8\times 9.$
\item [] $7170=1+2345+67\times 8\times 9.$
\item [] $7171=1\times (2+34\times 5\times 6)\times 7+8+9.$
\item [] $7172=1+(2+34\times 5\times 6)\times 7+8+9.$
\item [] $7173=1\times 23\times 4\times (5+6)\times 7+89.$
\item [] $7174=1+23\times 4\times (5+6)\times 7+89.$
\item [] $7175=(1+2\times 3)\times ((4\times 5\times 6+7)\times 8+9).$
\item [] $7176=1\times 2\times 3\times 4\times (5\times 6\times 7+89).$
\item [] $7177=(123+4)\times 56+7\times 8+9.$
\item [] $7178=(1+2)^3\times 45+67\times 89.$
\item [] $7179=1+23+(4+5)\times (6+789).$
\item [] $7180=(1+23\times 4\times (5+6))\times 7+89.$
\item [] $\mathit{7181=(1+2\times 3)^4\times 5-67\times 8\times 9.}$
\item [] $7182=(1+23)\times 45+678\times 9.$
\item [] $7183=1+2\times (3456+(7+8)\times 9).$
\item [] $7184=1+2^{(3\times 4)}+(5\times 67+8)\times 9.$
\item [] $7185=(1\times 2+34+56)\times 78+9.$
\item [] $7186=1+(2+34+56)\times 78+9.$
\item [] $7187=1\times 2+(3^4+5+6)\times 78+9.$
\item [] $7188=1\times (2\times 3)^4\times 5+6+78\times 9.$
\item [] $7189=1+(2\times 3)^4\times 5+6+78\times 9.$
\item [] $7190=1\times (2\times 3+4)\times (5+6\times 7\times (8+9)).$
\item [] $7191=(123+4)\times 56+7+8\times 9.$
\item [] $7192=1\times 2^{(3+4)}\times 56+7+8+9.$
\item [] $7193=1+2^{(3+4)}\times 56+7+8+9.$
\item [] $7194=123\times 4\times 5+6\times 789.$
\item [] $7195=1^2+3^4\times 5+6789.$
\item [] $7196=1\times 2+3^4\times 5+6789.$
\item [] $7197=1+2+3^4\times 5+6789.$
\item [] $7198=(1+2\times 3)\times 4^5+6+7+8+9.$
\item [] $7199=(123+4)\times 56+78+9.$
\item [] $7200=(12+3+4+56)\times (7+89).$
\item[]$\mbox{Decreasing order}$
\item [] $7131=9\times 8\times (76+5\times 4+3)+2+1.$
\item [] $\mathit{7132=-98\times 7+6^5+43-2+1.}$
\item [] $7133=9\times 8+(7+6)\times 543+2\times 1.$
\item [] $7134=(9+8)\times 7\times 6+5\times 4\times 321.$
\item [] $7135=9+876+5^4\times (3^2+1).$
\item [] $\mathit{7136=98\times 76+5+4-321.}$
\item [] $7137=9+8\times (76+5)\times (4+3\times 2+1).$
\item [] $7138=9\times (8+765+4\times (3+2))+1.$
\item [] $7139=(9+8)\times 7+65\times 4\times 3^{(2+1)}.$
\item [] $7140=(987+6\times 5)\times (4+3)+21.$
\item [] $7141=98\times (7+65)+4^3+21.$
\item [] $7142=98\times (7+65)+43\times 2\times 1.$
\item [] $7143=(9+87)\times 65+43\times 21.$
\item [] $7144=9+87\times (6+5\times (4+3))\times 2+1.$
\item [] $7145=(9+8)\times 7\times (6+54)+3+2\times 1.$
\item [] $7146=(9+8)\times 7\times (6+54)+3\times 2\times 1.$
\item [] $7147=(9+8)\times 7\times (6+54)+3\times 2+1.$
\item [] $7148=987\times 6+(5\times (4+3))^2+1.$
\item [] $7149=9\times 8\times (76+5\times 4+3)+21.$
\item [] $7150=(98+7)\times 65+4+321.$
\item [] $7151=(9+8\times 7)\times (65+43+2)+1.$
\item [] $7152=9\times 8+(7+6)\times 543+21.$
\item [] $7153=9+8\times (765+4\times 32\times 1).$
\item [] $7154=98\times (7+6+54+3\times 2\times 1).$
\item [] $7155=(9\times 8\times 7+6+5^4\times 3)\times (2+1).$
\item [] $7156=(9\times 87+65\times 43)\times 2\times 1.$
\item [] $7157=(9\times 87+65\times 43)\times 2+1.$
\item [] $7158=(9+8\times (7+6\times 5)\times 4)\times 3\times 2\times 1.$
\item [] $7159=98+(7+6)\times 543+2\times 1.$
\item [] $7160=98+(7+6+5+4)\times 321.$
\item [] $7161=9+8\times (765+4\times 32+1).$
\item [] $7162=(98+(76+5)\times 43)\times 2\times 1.$
\item [] $7163=(98+(76+5)\times 43)\times 2+1.$
\item [] $7164=(9+8)\times 7\times (6+54)+3+21.$
\item [] $7165=9\times (8\times 7+6\times (54+3))\times 2+1.$
\item [] $7166=9\times (8\times 7\times 6+5)+4^{(3\times 2)}+1.$
\item [] $7167=987+6\times (5+4^(3+2)+1).$
\item [] $7168=(98+7)\times 65+(4+3)^{(2+1)}.$
\item [] $7169=(9+87+6+5)\times (4+3\times 21).$
\item [] $7170=9+(8+7+6)\times (5\times 4+321).$
\item [] $7171=(9+8\times 7+6)\times (5+4\times (3+21)).$
\item [] $7172=(9+8)\times 7\times (6+54)+32\times 1.$
\item [] $7173=9\times (8+765+4\times 3\times 2\times 1).$
\item [] $7174=(9+8)\times (7+(65+4)\times 3\times 2+1).$
\item [] $7175=98\times (7+6+5\times 4\times 3)+21.$
\item [] $7176=9\times 8\times (76+5)+4^3\times 21.$
\item [] $\mathit{7177=-98\times 7+6^5+43\times 2+1.}$
\item [] $7178=98+(7+6)\times 543+21.$
\item [] $7179=987\times 6+(5^4+3)\times 2+1.$
\item [] $7180=(9+8+76+5^4)\times (3^2+1).$
\item [] $7181=(9+8\times 76)\times 5+4^{(3+2+1)}.$
\item [] $7182=987\times 6+5\times 4\times 3\times 21.$
\item [] $7183=(9\times 87+6+5+4)\times 3^2+1.$
\item [] $7184=98\times (7+65)+4^3\times 2\times 1.$
\item [] $7185=98\times (7+65)+4\times 32+1.$
\item [] $7186=(9\times 87+654)\times (3+2)+1.$
\item [] $7187=9+8\times (7+6)\times (5+4^3)+2\times 1.$
\item [] $7188=98\times (7+65)+4\times (32+1).$
\item [] $7189=(98+7\times 65)\times (4+3^2\times 1).$
\item [] $7190=(9+8\times 7+654)\times (3^2+1).$
\item [] $7191=9+(8+7+6\times 54+3)\times 21.$
\item [] $7192=(98+76+5^4)\times 3^2+1.$
\item [] $7193=(9+8)\times (76\times 5+43)+2\times 1.$
\item [] $7194=(9\times 8+7)\times 6+5\times 4^3\times 21.$
\item [] $7195=(9\times 8+7+6\times 5)\times (4^3+2)+1.$
\item [] $7196=98+7\times (6+(5+43)\times 21).$
\item [] $\mathit{7197=9\times 87-6+5\times 4\times 321.}$
\item [] $\mathit{7198=-98+76\times (5+43)\times 2\times 1.}$
\item [] $7199=9+8+7\times 6\times (54+3)\times (2+1).$
\item [] $7200=9\times 8\times 76+54\times 32\times 1.$
\item[]$\mbox{Increasing order}$
\item [] $7201=1+(2+3+45)\times 6\times (7+8+9).$
\item [] $7202=1234+5+67\times 89.$
\item [] $7203=(1+2^3\times 4)\times (5\times 6\times 7+8)+9.$
\item [] $7204=1\times (2+3^4)\times 5+6789.$
\item [] $7205=1+(2+3^4)\times 5+6789.$
\item [] $7206=12+3^4\times 5+6789.$
\item [] $7207=12^3\times 4+5\times (6\times 7+8+9).$
\item [] $7208=(123+4)\times 56+7+89.$
\item [] $7209=123\times (4+5)+678\times 9.$
\item [] $7210=1\times 2\times (3+4)\times (5+6+7\times 8\times 9).$
\item [] $7211=12^3\times 4+5\times 6\times 7+89.$
\item [] $7212=1^2\times 34\times 5\times 6\times 7+8\times 9.$
\item [] $7213=1^2+34\times 5\times 6\times 7+8\times 9.$
\item [] $7214=1\times 2+34\times 5\times 6\times 7+8\times 9.$
\item [] $7215=1+2+34\times 5\times 6\times 7+8\times 9.$
\item [] $7216=(1+2\times 3+4)\times (567+89).$
\item [] $7217=(1+2^3)^4+567+89.$
\item [] $7218=1^2\times 3^4\times (5+6+78)+9.$
\item [] $7219=1^2+3^4\times (5+6+78)+9.$
\item [] $7220=1\times 2+3^4\times (5+6+78)+9.$
\item [] $7221=12+3^4\times (5+67+8+9).$
\item [] $7222=1+(23+4+56)\times (78+9).$
\item [] $7223=123\times (45+6+7)+89.$
\item [] $7224=12+34\times 5\times 6\times 7+8\times 9.$
\item [] $7225=(1+2^3\times 4+56\times 7)\times (8+9).$
\item [] $7226=1\times (2+34\times 5\times 6)\times 7+8\times 9.$
\item [] $7227=1+(2+34\times 5\times 6)\times 7+8\times 9.$
\item [] $7228=1^{23}+(4^5+6)\times 7+8+9.$
\item [] $7229=12\times (34+567)+8+9.$
\item [] $7230=1^2+34\times 5\times 6\times 7+89.$
\item [] $7231=1\times 2+34\times 5\times 6\times 7+89.$
\item [] $7232=1+2+34\times 5\times 6\times 7+89.$
\item [] $7233=1\times 2^{(3+4)}\times 56+7\times 8+9.$
\item [] $7234=1+2^{(3+4)}\times 56+7\times 8+9.$
\item [] $7235=1+2\times (3\times 4^5+67\times 8+9).$
\item [] $7236=12\times 3\times (45+67+89).$
\item [] $7237=1+(2\times 3)^4\times 5+(6+78)\times 9.$
\item [] $7238=1\times 2+3\times 4\times (5+6+7\times 8)\times 9.$
\item [] $7239=1+2+3\times 4\times (5+6+7\times 8)\times 9.$
\item [] $7240=(1+2^3)^4+56+7\times 89.$
\item [] $7241=12+34\times 5\times 6\times 7+89.$
\item [] $7242=12+3+(4^5+6)\times 7+8+9.$
\item [] $7243=1\times (2+34\times 5\times 6)\times 7+89.$
\item [] $7244=1+(2+34\times 5\times 6)\times 7+89.$
\item [] $7245=12\times (3\times 45+6\times 78)+9.$
\item [] $7246=12^3+(4\times 5+6\times 7)\times 89.$
\item [] $7247=1\times 2^{(3+4)}\times 56+7+8\times 9.$
\item [] $7248=1+2^{(3+4)}\times 56+7+8\times 9.$
\item [] $7249=123+4^5+678\times 9.$
\item [] $7250=1+23\times 4\times 5+6789.$
\item [] $7251=1+23+(4^5+6)\times 7+8+9.$
\item [] $7252=(1+2\times 3)\times 4^5+67+8+9.$
\item [] $7253=(1+2^3)^4+5+678+9.$
\item [] $7254=(12+3^4)\times 5+6789.$
\item [] $7255=1^2+345\times (6+7+8)+9.$
\item [] $7256=(1+2^3)^4+5\times (67+8\times 9).$
\item [] $7257=12\times (34+5)+6789.$
\item [] $7258=((1+2)\times 34+5)\times 67+89.$
\item [] $7259=(1^2+34+56\times 7)\times (8+9).$
\item [] $7260=12+(3+4^5+6)\times 7+8+9.$
\item [] $7261=12\times 34+(5+6)\times 7\times 89.$
\item [] $7262=12^3\times 4+5+6\times 7\times 8+9.$
\item [] $7263=(1+2+34+56)\times 78+9.$
\item [] $7264=12^3\times 4+5\times 67+8+9.$
\item [] $7265=123\times (4+5)\times 6+7\times 89.$
\item [] $7266=123\times (4+5+6\times 7+8)+9.$
\item [] $7267=12^3\times 4+5\times (6+7\times 8+9).$
\item [] $7268=(12\times 3^4+56)\times 7+8\times 9.$
\item [] $7269=(1+23)\times 45\times 6+789.$
\item [] $7270=(1+2\times 3)\times 4^5+6+7+89.$
\item[]$\mbox{Decreasing order}$
\item [] $7201=9\times 8\times 76+54\times 32+1.$
\item [] $7202=98\times 7+6\times 543\times 2\times 1.$
\item [] $7203=98\times 7+6\times 543\times 2+1.$
\item [] $7204=9+(8\times (7+6)+5)\times (4^3+2)+1.$
\item [] $7205=9\times 8\times 76+5+(4\times 3)^{(2+1)}.$
\item [] $7206=9+8+7+6\times (54+3)\times 21.$
\item [] $7207=((9+87)\times 6+5^4)\times 3\times 2+1.$
\item [] $7208=98\times 7+6\times (543\times 2+1).$
\item [] $7209=9\times 87+6+5\times 4\times 321.$
\item [] $7210=9+(8+7)\times (6+5+4)\times 32+1.$
\item [] $7211=987\times 6+5+4\times 321.$
\item [] $7212=(9+8)\times (76\times 5+43)+21.$
\item [] $\mathit{7213=-9\times 8\times 7+6^5+4-3\times 21.}$
\item [] $\mathit{7214=-9\times 87-6+(5\times 4)^3+2+1.}$
\item [] $7215=(98+7+6)\times 5\times (4+3^2\times 1).$
\item [] $7216=(98+7+6)\times 5\times (4+3^2)+1.$
\item [] $7217=(987+(6+5)\times 4)\times (3\times 2+1).$
\item [] $7218=987\times 6+54\times (3+21).$
\item [] $\mathit{7219=-98\times 7+6^5+4\times 32+1.}$
\item [] $\mathit{7220=-987+6^5+432-1.}$
\item [] $7221=(9+8\times 7\times 6)\times 5\times 4+321.$
\item [] $\mathit{7222=-987+6^5+432+1.}$
\item [] $\mathit{7223=-9\times 8+76\times (5+43)\times 2-1.}$
\item [] $7224=(9+8+(76+5)\times 4+3)\times 21.$
\item [] $7225=(9+8)\times (76\times 5+43+2\times 1).$
\item [] $7226=98\times 7+654\times (3^2+1).$
\item [] $7227=987+65\times 4\times (3+21).$
\item [] $7228=((9+8)\times (7\times 6+5)+4)\times 3^2+1.$
\item [] $\mathit{7229=-9\times 8\times 7+6^5-4^3+21.}$
\item [] $7230=9\times 8\times 7+6+5\times 4^3\times 21.$
\item [] $7231=98\times 7+6543+2\times 1.$
\item [] $7232=98\times 7+6543+2+1.$
\item [] $7233=9+8\times 7\times (65+43+21).$
\item [] $7234=9+(8\times 7+6+5\times 4+3)^2\times 1.$
\item [] $7235=9+8\times 7\times (65+4^3)+2\times 1.$
\item [] $7236=9\times (87+654+3\times 21).$
\item [] $7237=9\times (8+76+5\times (4\times 3)^2)+1.$
\item [] $7238=9\times (8\times 7+6+5)\times 4\times 3+2\times 1.$
\item [] $7239=9\times (8\times 7+6+5)\times 4\times 3+2+1.$
\item [] $7240=(9+(8+7)\times 6+5^4)\times (3^2+1).$
\item [] $7241=98\times 7\times 6+5^4\times (3+2\times 1).$
\item [] $7242=98\times 7\times 6+5^4\times (3+2)+1.$
\item [] $7243=987+6+5^4\times (3^2+1).$
\item [] $\mathit{7244=-9\times 8\times 7+6^5+4-32\times 1.}$
\item [] $7245=(9+87+65)\times (43+2\times 1).$
\item [] $7246=(9+87+65)\times (43+2)+1.$
\item [] $7247=9+8\times 7+6\times (54+3)\times 21.$
\item [] $7248=98\times (7+65)+4^3\times (2+1).$
\item [] $7249=9+8\times (7\times (65+4^3)+2\times 1).$
\item [] $7250=98\times 7+6543+21.$
\item [] $7251=9+87\times 6+5\times 4^3\times 21.$
\item [] $7252=98\times (7+6+54+3\times 2+1).$
\item [] $7253=98\times ((7+6)\times 5+4+3+2)+1.$
\item [] $7254=9\times 8\times 76+54\times (32+1).$
\item [] $7255=9\times 8\times 7+(6+5+4)^3\times 2+1.$
\item [] $\mathit{7256=(98+7)\times 65+432-1.}$
\item [] $7257=(98+7)\times 65+432\times 1.$
\item [] $7258=(98+7)\times 65+432+1.$
\item [] $7259=(9+8)\times 7\times (6+5\times (4+3\times 2+1)).$
\item [] $7260=(9\times 8+(7+6\times 5)\times 4)\times (32+1).$
\item [] $7261=9\times 8+7+6\times (54+3)\times 21.$
\item [] $7262=9+8+7\times (6+5+4^(3+2)\times 1).$
\item [] $7263=9+(87+6)\times (54+3+21).$
\item [] $7264=(9+87)\times 65+4^(3+2\times 1).$
\item [] $7265=(9+87)\times 65+4^(3+2)+1.$
\item [] $7266=(9\times 8+7\times 6\times 5+4^3)\times 21.$
\item [] $7267=(98\times (7+6\times 5)+4+3)\times 2+1.$
\item [] $7268=(9\times 8+7)\times (6+54+32\times 1).$
\item [] $7269=(98+7)\times (65+4)+3+21.$
\item [] $7270=(9+87+6+5^4)\times (3^2+1).$
\item[]$\mbox{Increasing order}$
\item [] $7271=(1^2+34\times 5)\times 6\times 7+89.$
\item [] $7272=1\times 234\times 5+678\times 9.$
\item [] $7273=1+234\times 5+678\times 9.$
\item [] $7274=(1+2^3)^4+5+6+78\times 9.$
\item [] $7275=1\times (2\times 3)^4\times 5+6+789.$
\item [] $7276=1+(2\times 3)^4\times 5+6+789.$
\item [] $7277=(1+234)\times 5+678\times 9.$
\item [] $7278=123+(4+5)\times (6+789).$
\item [] $7279=1+2\times 3\times (4^5+(6+7+8)\times 9).$
\item [] $7280=(1+(2\times 3)^4)\times 5+6+789.$
\item [] $7281=(1\times 234+567+8)\times 9.$
\item [] $7282=1^{23}\times (4^5+6)\times 7+8\times 9.$
\item [] $7283=1\times 2\times 3\times 4^5+67\times (8+9).$
\item [] $7284=123\times 4\times 5+67\times 8\times 9.$
\item [] $7285=1^2+3\times 4^5+6\times 78\times 9.$
\item [] $7286=123\times 4+5+6789.$
\item [] $7287=1+2+3\times 4^5+6\times 78\times 9.$
\item [] $7288=1+2+3+(4^5+6)\times 7+8\times 9.$
\item [] $7289=1\times (2+3)^4+56\times 7\times (8+9).$
\item [] $7290=(123+4+5+678)\times 9.$
\item [] $7291=(1+2)\times 3+(4^5+6)\times 7+8\times 9.$
\item [] $7292=1\times 2+3^4\times (5+6+7+8\times 9).$
\item [] $7293=(1+2^3)^4+5\times 6+78\times 9.$
\item [] $7294=1^2+3+(4+5)\times 6\times (7+8)\times 9.$
\item [] $7295=1\times 2+3+(4+5)\times 6\times (7+8)\times 9.$
\item [] $7296=12+3\times 4^5+6\times 78\times 9.$
\item [] $7297=12+3+(4^5+6)\times 7+8\times 9.$
\item [] $7298=(1\times 2\times 3+4+5+67)\times 89.$
\item [] $7299=(1+2)\times 34\times 5+6789.$
\item [] $7300=1^{23}+(4^5+6)\times 7+89.$
\item [] $7301=12\times (34+567)+89.$
\item [] $7302=1^2\times 3+(4^5+6)\times 7+89.$
\item [] $7303=1^2+3+(4^5+6)\times 7+89.$
\item [] $7304=(1\times 23\times 45+6)\times 7+8+9.$
\item [] $7305=1+2+3+(4^5+6)\times 7+89.$
\item [] $7306=1+23+(4^5+6)\times 7+8\times 9.$
\item [] $7307=12^3\times 4+5+6\times (7\times 8+9).$
\item [] $7308=(1+23+4+56)\times (78+9).$
\item [] $7309=1234\times 5+67\times (8+9).$
\item [] $7310=(1^2+3+4^5+6)\times 7+8\times 9.$
\item [] $7311=(1+23\times 45+6)\times 7+8+9.$
\item [] $7312=1\times 2+(3^4+5)\times (6+7+8\times 9).$
\item [] $7313=1\times (2+34\times 5)\times 6\times 7+89.$
\item [] $7314=12+3+(4^5+6)\times 7+89.$
\item [] $7315=12+(3+4^5+6)\times 7+8\times 9.$
\item [] $\mathit{7316=123+(4^5+6)\times 7-8-9.}$
\item [] $7317=(1+2)^3\times 45+678\times 9.$
\item [] $7318=12\times 3+(4^5+6)\times 7+8\times 9.$
\item [] $7319=12^3\times 4+5\times 67+8\times 9.$
\item [] $7320=(1+2+3+4)\times (5\times 6+78\times 9).$
\item [] $7321=12^3\times 4+56\times 7+8+9.$
\item [] $7322=1\times 23+(4^5+6)\times 7+89.$
\item [] $7323=1+23+(4^5+6)\times 7+89.$
\item [] $7324=(1+2+3+4^5+6)\times 7+8\times 9.$
\item [] $7325=1+(2\times 3+4^5+6)\times 7+8\times 9.$
\item [] $7326=(1^2+345+6\times 78)\times 9.$
\item [] $7327=1234\times 5+(6+7)\times 89.$
\item [] $7328=(1+2)\times (3^4\times 5\times 6+7)+8+9.$
\item [] $7329=(1+2)^3\times 4\times 5+6789.$
\item [] $\mathit{7330=(1+2+3+4)\times (-56+789).}$
\item [] $7331=(1+2\times 3+4^5+6)\times 7+8\times 9.$
\item [] $7332=(12+3+4^5+6)\times 7+8+9.$
\item [] $\mathit{7333=123\times (4+56)-7\times 8+9.}$
\item [] $7334=1\times (2+3+4^5+6)\times 7+89.$
\item [] $7335=(1\times 2+345+6\times 78)\times 9.$
\item [] $7336=12^3\times 4+5\times 67+89.$
\item [] $7337=(1^{23}+4)^5+6\times 78\times 9.$
\item [] $7338=1+(2+3)^4\times 5+6\times 78\times 9.$
\item [] $7339=1+(2^3+4^5+6)\times 7+8\times 9.$
\item [] $\mathit{7340=-1+(2\times 3+4^5+6)\times 7+89.}$
\item[]$\mbox{Decreasing order}$
\item [] $7271=987\times 6+5+4^3\times 21.$
\item [] $7272=9\times 8\times (76+5\times 4+3+2\times 1).$
\item [] $7273=98\times (7+6)\times 5+43\times 21.$
\item [] $7274=9\times 8\times (7+6\times 5+4^3)+2\times 1.$
\item [] $7275=9\times 8\times (7+6\times 5+4^3)+2+1.$
\item [] $7276=(9+8)\times (7+6\times 5\times (4+3)\times 2+1).$
\item [] $7277=(98+7)\times (65+4)+32\times 1.$
\item [] $7278=9+87+6\times (54+3)\times 21.$
\item [] $\mathit{7279=-9-8\times 7+6^5-432\times 1.}$
\item [] $7280=98+7\times 6\times (54+3)\times (2+1).$
\item [] $7281=9\times (8+765+4+32\times 1).$
\item [] $7282=9\times (8+765)+4+321.$
\item [] $7283=9\times (8+(7+65\times 4)\times 3)+2\times 1.$
\item [] $7284=9\times (8+(7+65\times 4)\times 3)+2+1.$
\item [] $\mathit{7285=98\times 76-54\times 3-2+1.}$
\item [] $\mathit{7286=-9\times 8\times 7+6^5-4-3+21.}$
\item [] $7287=98+7+6\times (54+3)\times 21.$
\item [] $\mathit{7288=98\times 76-54\times 3+2\times 1.}$
\item [] $7289=9+8\times (76+54)\times (3\times 2+1).$
\item [] $7290=9\times (8+765+4+32+1).$
\item [] $7291=9\times (8+7+6\times 5\times 4)\times 3\times 2+1.$
\item [] $7292=9\times (8+7)\times (6+5+43)+2\times 1.$
\item [] $7293=9\times 8\times (7+6\times 5+4^3)+21.$
\item [] $7294=(9+8)\times (76\times 5+(4+3)^2)+1.$
\item [] $7295=9\times (8+7)\times 6\times (5+4)+3+2\times 1.$
\item [] $7296=(98+76+54)\times 32\times 1.$
\item [] $7297=(98+76+54)\times 32+1.$
\item [] $7298=(9+(8+7+6)\times 5)\times 4^3+2\times 1.$
\item [] $7299=9\times 87+6\times 543\times 2\times 1.$
\item [] $7300=9\times 87+6\times 543\times 2+1.$
\item [] $7301=(9+8)\times 7+6\times (54+3)\times 21.$
\item [] $7302=(9\times 8+7\times (6+5))\times (4+3)^2+1.$
\item [] $7303=9+(8\times 7\times 65+4+3)\times 2\times 1.$
\item [] $7304=(9\times 8\times 7\times 6+5^4+3)\times 2\times 1.$
\item [] $7305=9+876+5\times 4\times 321.$
\item [] $\mathit{7306=98\times 76+5-(4+3)\times 21.}$
\item [] $7307=9+8\times 76\times (5+4+3)+2\times 1.$
\item [] $7308=98\times (7+65)+4\times 3\times 21.$
\item [] $7309=9\times ((8+7)\times (6+5+43)+2)+1.$
\item [] $7310=(9\times 8+7+6)\times (54+32\times 1).$
\item [] $7311=9\times (8+7)\times (6+5+43)+21.$
\item [] $7312=(98+7+65)\times 43+2\times 1.$
\item [] $7313=(98+7+65)\times 43+2+1.$
\item [] $7314=9+8+76\times (5+43)\times 2+1.$
\item [] $\mathit{7315=98\times 76-5-4\times 32\times 1.}$
\item [] $\mathit{7316=-9-8+7654-321.}$
\item [] $7317=9+87\times (6+54+3+21).$
\item [] $7318=9+87\times 6\times (5+4+3+2)+1.$
\item [] $7319=9+(8\times 7+6\times 5)\times (4^3+21).$
\item [] $7320=(9+8\times (7+6\times 5))\times 4\times 3\times 2\times 1.$
\item [] $7321=9+8\times (76\times (5+4+3)+2\times 1).$
\item [] $7322=(9+8\times 7\times 65+4\times 3)\times 2\times 1.$
\item [] $7323=9\times 87+654\times (3^2+1).$
\item [] $7324=9\times 8+7\times (6+5+4^(3+2)+1).$
\item [] $7325=(9+(8\times 7)\times (6+5\times 4))\times (3+2)\times 1.$
\item [] $7326=9\times 8\times 76+5+43^2\times 1.$
\item [] $7327=9\times 8\times 76+5+43^2+1.$
\item [] $7328=9\times 87+6543+2\times 1.$
\item [] $7329=9\times 87+6543+2+1.$
\item [] $7330=(9\times 8+7+654)\times (3^2+1).$
\item [] $7331=(98+7+65)\times 43+21.$
\item [] $7332=(9\times (8+7)+6)\times (5\times 4+32\times 1).$
\item [] $7333=(9\times (8+7)+6)\times (5\times 4+32)+1.$
\item [] $\mathit{7334=9-8+7654-321.}$
\item [] $7335=(98+(7+6)\times 5)\times (43+2\times 1).$
\item [] $7336=(9\times (8+7)\times 6+5)\times (4+3+2)+1.$
\item [] $7337=9+8\times 76+5\times 4^3\times 21.$
\item [] $7338=(98\times (7+6\times 5)+43)\times 2\times 1.$
\item [] $7339=(98\times (7+6\times 5)+43)\times 2+1.$
\item [] $7340=((9+8)\times 7\times 6+5\times 4)\times (3^2+1).$
\item[]$\mbox{Increasing order}$
\item [] $7341=1234+5+678\times 9.$
\item [] $7342=1+(2+3^4+5+6)\times 78+9.$
\item [] $7343=(1+2^{(3+4)})\times 56+7\times (8+9).$
\item [] $7344=123\times (4+5)\times 6+78\times 9.$
\item [] $7345=123\times 4+(5+6)\times 7\times 89.$
\item [] $7346=1+(2+3^4+5\times 6)\times (7\times 8+9).$
\item [] $7347=(1+2+34+56)\times (7+8\times 9).$
\item [] $7348=(1+2\times 3+4^5+6)\times 7+89.$
\item [] $7349=1+2\times (34+56\times (7\times 8+9)).$
\item [] $7350=123+(4^5+6)\times 7+8+9.$
\item [] $7351=1+2\times ((3+4\times (5+6))\times 78+9).$
\item [] $\mathit{7352=(1234+5)\times 6+7-89.}$
\item [] $7353=1\times 2\times 34\times (5\times 6+78)+9.$
\item [] $7354=1+2\times 34\times (5\times 6+78)+9.$
\item [] $7355=(1+2)^3\times 45\times 6+7\times 8+9.$
\item [] $7356=1+(2^3+4^5+6)\times 7+89.$
\item [] $7357=(1+2\times 3)\times 4^5+(6+7+8)\times 9.$
\item [] $\mathit{7358=-1+(23\times 45+6)\times 7+8\times 9.}$
\item [] $7359=1\times (23\times 45+6)\times 7+8\times 9.$
\item [] $7360=1\times 23\times (4\times 56+7+89).$
\item [] $7361=12\times 34\times (5+6+7)+8+9.$
\item [] $7362=(1+2^3+4^5+6)\times 7+89.$
\item [] $7363=1+2\times (34+5\times (67+8))\times 9.$
\item [] $\mathit{7364=(-1+23+4^5+6)\times 7\times (-8+9).}$
\item [] $7365=12\times (3+45)+6789.$
\item [] $7366=(1+23\times 45+6)\times 7+8\times 9.$
\item [] $7367=(1^2+34)\times 5\times 6\times 7+8+9.$
\item [] $7368=1+(2^3+45)\times (67+8\times 9).$
\item [] $7369=1\times 23\times 4\times (5+67+8)+9.$
\item [] $7370=1+23\times 4\times (5+67+8)+9.$
\item [] $7371=12^3\times 4+5\times 6\times (7+8)+9.$
\item [] $7372=1+(2+3+4)\times (5\times 6+789).$
\item [] $7373=1\times 2+(34+5)\times (6+7+8)\times 9.$
\item [] $7374=1+2+(34+5)\times (6+7+8)\times 9.$
\item [] $7375=(1+2\times 3\times 4)\times 5\times (6\times 7+8+9).$
\item [] $7376=12^3\times 4+56\times 7+8\times 9.$
\item [] $7377=1+(23\times 45+6)\times 7+89.$
\item [] $7378=1\times 2\times (3+4+5\times 6\times 7)\times (8+9).$
\item [] $7379=1+2\times (3+4+5\times 6\times 7)\times (8+9).$
\item [] $7380=12\times 3^4+(5+67)\times 89.$
\item [] $7381=1234+(5+678)\times 9.$
\item [] $7382=(1^2+3)^4\times 5+678\times 9.$
\item [] $7383=(1+23\times 45+6)\times 7+89.$
\item [] $\mathit{7384=12^3\times 4-5+6\times 78+9.}$
\item [] $\mathit{7385=1\times 2+3\times 4\times (5+6)\times 7\times 8-9.}$
\item [] $7386=(1+2)^3\times 45\times 6+7+89.$
\item [] $7387=(12+34+5\times 6+7)\times 89.$
\item [] $7388=(1+2^3+4)\times 567+8+9.$
\item [] $7389=(123+4\times 5+678)\times 9.$
\item [] $7390=1+2^{(3\times 4)}+(5\times 6+7)\times 89.$
\item [] $7391=12^3\times 4+5+6\times (7+8\times 9).$
\item [] $7392=1^2\times (3+4)\times (5+6)\times (7+89).$
\item [] $7393=12^3\times 4+56\times 7+89.$
\item [] $7394=12^3\times 4+5+6\times 78+9.$
\item [] $7395=(1+23+4^5+6)\times 7+8+9.$
\item [] $7396=(1\times 2+3^4)\times (5+6+78)+9.$
\item [] $7397=1+(2+3^4)\times (5+6+78)+9.$
\item [] $7398=(1+2)\times 3\times (4\times 5\times 6+78\times 9).$
\item [] $7399=1+2\times (3+(4+5+6\times 7)\times 8)\times 9.$
\item [] $7400=(1+2)\times (3^4\times 5\times 6+7)+89.$
\item [] $7401=1^2\times 3\times 4\times (5+6)\times 7\times 8+9.$
\item [] $7402=1\times 2\times (3\times 4^5+6+7\times 89).$
\item [] $7403=1\times (2\times 3)^4+5+678\times 9.$
\item [] $7404=123\times (4+56)+7+8+9.$
\item [] $7405=123+(4^5+6)\times 7+8\times 9.$
\item [] $7406=(12+34)\times (5+67+89).$
\item [] $7407=12^3\times 4+(5+6\times 7+8)\times 9.$
\item [] $7408=1\times 2^3\times (4\times 56+78\times 9).$
\item [] $7409=(1\times 2\times (3+456)+7)\times 8+9.$
\item [] $7410=1\times 2\times (3+(4+5)\times 6)\times (7\times 8+9).$
\item[]$\mbox{Decreasing order}$
\item [] $\mathit{7341=-98\times 7+6+(5\times 4)^3+21.}$
\item [] $\mathit{7342=9\times 876-543+2-1.}$
\item [] $7343=98+7\times (6+5+4^(3+2\times 1)).$
\item [] $7344=9\times 8\times (76+5\times 4+3+2+1).$
\item [] $7345=9+8\times 7\times (65+4^3+2\times 1).$
\item [] $7346=9\times 8\times (7\times 6+5\times 4\times 3)+2\times 1.$
\item [] $7347=9\times 87+6543+21.$
\item [] $\mathit{7348=9\times 8\times 76+5^4\times 3+2-1.}$
\item [] $7349=9\times 8\times 76+5^4\times 3+2\times 1.$
\item [] $7350=9\times 8\times 7\times 6+5+4321.$
\item [] $7351=(9+8+7+6)\times 5\times (4+3)^2+1.$
\item [] $7352=98\times ((7+6+5)\times 4+3)+2\times 1.$
\item [] $7353=(9\times 87+6\times 5+4)\times 3^2\times 1.$
\item [] $7354=(9\times 87+6\times 5+4)\times 3^2+1.$
\item [] $\mathit{7355=-9-87+6^5-4-321.}$
\item [] $7356=9\times 8\times (7+6)+5\times 4\times 321.$
\item [] $\mathit{7357=98\times 76-5-43\times 2\times 1.}$
\item [] $\mathit{7358=98\times 76-5-43\times 2+1.}$
\item [] $7359=9+(87+65\times 4+3)\times 21.$
\item [] $7360=(98+7+6+5^4)\times (3^2+1).$
\item [] $7361=((9+8)\times (7+6)+5+4)\times 32+1.$
\item [] $7362=9\times (8+765+43+2\times 1).$
\item [] $7363=(9+8)\times 76\times 5+43\times 21.$
\item [] $\mathit{7364=98\times 7\times 6+(54+3)^2-1.}$
\item [] $7365=98\times 7\times 6+(54+3)^2\times 1.$
\item [] $7366=98\times 7\times 6+(54+3)^2+1.$
\item [] $7367=9+(8\times 76+5)\times 4\times 3+2\times 1.$
\item [] $7368=9\times 8\times 76+5^4\times 3+21.$
\item [] $7369=9\times 8+76\times (5+43)\times 2+1.$
\item [] $7370=9+8\times 76\times 5+4321.$
\item [] $7371=(98+7\times 6\times 5+43)\times 21.$
\item [] $7372=(9\times (8+7)\times 6+5+4)\times 3^2+1.$
\item [] $7373=9\times (87+6\times 5)\times (4+3)+2\times 1.$
\item [] $7374=9+8+7+6\times (5\times (4+3))^2\times 1.$
\item [] $7375=9+(8\times 7\times 65+43)\times 2\times 1.$
\item [] $7376=9+(8\times 7\times 65+43)\times 2+1.$
\item [] $7377=98\times (7+6+5)\times 4+321.$
\item [] $7378=(9\times 8+7\times 65)\times (4+3^2+1).$
\item [] $7379=(9\times 8+7\times 65)\times (4+3)\times 2+1.$
\item [] $7380=987\times 6+(5+4)^3\times 2\times 1.$
\item [] $7381=98\times (7+65)+4+321.$
\item [] $\mathit{7382=98\times 76-5-4^3+2+1.}$
\item [] $\mathit{7383=98\times 76-5\times (4+3^2)\times 1.}$
\item [] $7384=(9+8\times 7\times 65+43)\times 2\times 1.$
\item [] $7385=(9+8\times 7\times 65+43)\times 2+1.$
\item [] $7386=9+(8\times 76+5)\times 4\times 3+21.$
\item [] $7387=9+(87\times 6+5)\times (4+3)\times 2\times 1.$
\item [] $7388=9+8+(7+6)\times (5+4)\times 3\times 21.$
\item [] $7389=9\times (8+765)+432\times 1.$
\item [] $7390=9\times (8+765)+432+1.$
\item [] $\mathit{7391=98\times 76-5\times 4\times 3+2+1.}$
\item [] $7392=(9+8+7\times 6\times 5+4)\times 32\times 1.$
\item [] $7393=9\times 8\times 7\times (6+5)+43^2\times 1.$
\item [] $7394=98+76\times (5+43)\times 2\times 1.$
\item [] $7395=98+76\times (5+43)\times 2+1.$
\item [] $7396=(9+8+7\times 6+(5+4)\times 3)^2\times 1.$
\item [] $7397=(9\times 8\times 7+65)\times (4+3^2\times 1).$
\item [] $7398=(9\times 8\times 7+654\times 3)\times (2+1).$
\item [] $7399=98\times (7+65)+(4+3)^{(2+1)}.$
\item [] $7400=(9\times 8+76)\times 5\times (4+3+2+1).$
\item [] $7401=9+8\times 7\times (65+4+3\times 21).$
\item [] $7402=9+8\times (7\times 65+4+3)\times 2+1.$
\item [] $7403=9+8\times 7\times (6+5)\times 4\times 3+2\times 1.$
\item [] $7404=9+8\times 7\times (6+5)\times 4\times 3+2+1.$
\item [] $7405=9+(8\times 7+6\times 5)\times 43\times 2\times 1.$
\item [] $7406=(9+8\times 76)\times 5+4321.$
\item [] $7407=9+(8\times 76+5^4)\times 3\times 2\times 1.$
\item [] $7408=9+(8\times 76+5^4)\times 3\times 2+1.$
\item [] $7409=9+8+7\times (6+5)\times 4\times (3+21).$
\item [] $7410=(9+8\times 7)\times 6\times (5+4\times 3+2\times 1).$
\item[]$\mbox{Increasing order}$
\item [] $7411=1+(2+3\times 4+5)\times 6\times (7\times 8+9).$
\item [] $7412=1\times 2\times 34\times (5\times 6+7+8\times 9).$
\item [] $7413=123+(4+5)\times 6\times (7+8)\times 9.$
\item [] $\mathit{7414=(1+2\times 3+4)\times (5+678-9).}$
\item [] $\mathit{7415=-1+2\times 3456+7\times 8\times 9.}$
\item [] $7416=1\times 2\times 3456+7\times 8\times 9.$
\item [] $7417=1+2\times 3456+7\times 8\times 9.$
\item [] $7418=1+(2+3\times 4\times (5+6)\times 7)\times 8+9.$
\item [] $7419=12^3+4+5678+9.$
\item [] $7420=1+(2+3)^4+5+6789.$
\item [] $7421=1^2\times 34\times (5\times 6\times 7+8)+9.$
\item [] $7422=123+(4^5+6)\times 7+89.$
\item [] $7423=1\times 2+34\times (5\times 6\times 7+8)+9.$
\item [] $7424=(123+4)\times 5+6789.$
\item [] $7425=(12+345+6\times 78)\times 9.$
\item [] $7426=(1\times 2+3^4+5+6)\times (7+8\times 9).$
\item [] $7427=12^3+4+5\times 67\times (8+9).$
\item [] $7428=(1+2^3)^4+(5+6)\times 78+9.$
\item [] $7429=1\times 2^{(3+4)}\times 5+6789.$
\item [] $7430=1+2^{(3+4)}\times 5+6789.$
\item [] $7431=123\times (4+5)\times 6+789.$
\item [] $7432=1\times 2^3\times (4\times 5\times 6\times 7+89).$
\item [] $7433=12+34\times (5\times 6\times 7+8)+9.$
\item [] $7434=(1+2^{(3+4)})\times 5+6789.$
\item [] $7435=1+(2\times 34\times (5+6)+78)\times 9.$
\item [] $\mathit{7436=1\times 2^3\times 4^5-(6+78)\times 9.}$
\item [] $7437=12+((3\times 45\times 6+7)+8)\times 9.$
\item [] $\mathit{7438=12^3\times 4-5+(67-8)\times 9.}$
\item [] $7439=12^3\times 4+5+6\times (78+9).$
\item [] $7440=(12+3^4)\times (56+7+8+9).$
\item [] $7441=12^3\times 4+5\times (6+7)\times 8+9.$
\item [] $7442=1\times 2\times (3^4+56\times (7\times 8+9)).$
\item [] $7443=(12\times 3\times 4+5+678)\times 9.$
\item [] $7444=1+(23+4^5+6)\times 7+8\times 9.$
\item [] $7445=123\times (4+56)+7\times 8+9.$
\item [] $7446=12^3\times 4+5\times 6+7\times 8\times 9.$
\item [] $7447=1+(2\times 34+5)\times (6+7+89).$
\item [] $\mathit{7448=-1+2\times 3\times 4\times 5\times (6+7\times 8)+9.}$
\item [] $7449=1\times 2\times 3\times 4\times 5\times (6+7\times 8)+9.$
\item [] $7450=1+2\times 3\times 4\times 5\times (6+7\times 8)+9.$
\item [] $\mathit{7451=1\times 2^3\times (4^5+6)-789.}$
\item [] $7452=12\times 3\times (4\times 5\times 6+78+9).$
\item [] $7453=1+2\times 3\times (4+56+78)\times 9.$
\item [] $\mathit{7454=1\times 2+(3\times 4\times ((5+6)\times 7-8)\times 9).}$
\item [] $7455=(12+3)\times (4\times 5+6\times 78+9).$
\item [] $\mathit{7456=(1234-5)\times 6-7+89.}$
\item [] $7457=(12\times (3+4)\times (5+6)+7)\times 8+9.$
\item [] $7458=(1234+5)\times 6+7+8+9.$
\item [] $7459=123\times (4+56)+7+8\times 9.$
\item [] $7460=(1+2^3+4)\times 567+89.$
\item [] $7461=1+(23+4^5+6)\times 7+89.$
\item [] $7462=12^3\times 4+5+67\times 8+9.$
\item [] $7463=(12+3)\times 456+7\times 89.$
\item [] $7464=(12+3)\times 45+6789.$
\item [] $7465=(12+3+4)\times 56\times 7+8+9.$
\item [] $7466=1+2\times (3+456+7)\times 8+9.$
\item [] $7467=123\times (4+56)+78+9.$
\item [] $\mathit{7468=-1\times 23\times 4+56\times (7+8)\times 9.}$
\item [] $7469=(1+2\times 3+4)\times (56+7\times 89).$
\item [] $7470=1\times 2\times 3\times (456+789).$
\item [] $7471=1+2\times 3\times (456+789).$
\item [] $7472=12^3\times 4+56+7\times 8\times 9.$
\item [] $7473=12\times (3^4+5+67\times 8)+9.$
\item [] $7474=1+(2^3+45)\times (6+(7+8)\times 9).$
\item [] $7475=12^3\times 4+5+(6+7\times 8)\times 9.$
\item [] $7476=123\times (4+56)+7+89.$
\item [] $7477=1+(23+(4+5)\times 6+7)\times 89.$
\item [] $7478=(1\times 2+3)^4+(5+6)\times 7\times 89.$
\item [] $7479=1\times 2\times 345+6789.$
\item [] $7480=1+2\times 345+6789.$
\item[]$\mbox{Decreasing order}$
\item [] $7411=(9+8\times 7)\times 6\times (5+4\times 3+2)+1.$
\item [] $7412=98\times 7+6+5\times 4^3\times 21.$
\item [] $7413=987+6+5\times 4\times 321.$
\item [] $7414=(9+(8\times 7+6\times 5)\times 43)\times 2\times 1.$
\item [] $7415=(9+(8\times 7+6\times 5)\times 43)\times 2+1.$
\item [] $7416=9\times 8\times (7\times 6+54+3\times 2+1).$
\item [] $7417=9+(8\times 7\times 65+4^3)\times 2\times 1.$
\item [] $7418=9+(8\times 7\times 65+4^3)\times 2+1.$
\item [] $7419=9\times 8\times (76+(5+4)\times 3)+2+1.$
\item [] $7420=(98+7\times 6)\times (5\times 4+32+1).$
\item [] $7421=(9+8\times 7\times (6+5)\times 4)\times 3+2\times 1.$
\item [] $7422=9+8\times 7\times (6+5)\times 4\times 3+21.$
\item [] $\mathit{7423=98\times 76-5+4-3-21.}$
\item [] $7424=(98+7+6+5)\times (43+21).$
\item [] $7425=9\times (8\times 76+5\times 43+2\times 1).$
\item [] $7426=(98+7+6+5)\times 4^3+2\times 1.$
\item [] $7427=9\times (8+7+65\times 4)\times 3+2\times 1.$
\item [] $7428=9\times (8+7+65\times 4)\times 3+2+1.$
\item [] $7429=(9\times 8+7)\times (6\times 5+4^3)+2+1.$
\item [] $7430=9\times 8+7+6\times (5\times (4+3))^2+1.$
\item [] $\mathit{7431=98+7654-321.}$
\item [] $\mathit{7432=98\times 76-5-4-3\times 2-1.}$
\item [] $\mathit{7433=98\times 76+5+4-3-21.}$
\item [] $7434=(9+8)\times 7\times 6+5\times 4^3\times 21.$
\item [] $7435=9\times (87+6+5\times 4^3)\times 2+1.$
\item [] $7436=98\times 7+(6+5+4)^3\times 2\times 1.$
\item [] $7437=98\times 7+(6+5+4)^3\times 2+1.$
\item [] $\mathit{7438=98\times 76+5-4\times 3-2-1.}$
\item [] $\mathit{7439=98\times 76+54-3\times 21.}$
\item [] $7440=(9+8\times 7\times (6+5)\times 4)\times 3+21.$
\item [] $7441=(9\times 8+76\times (5+43))\times 2+1.$
\item [] $\mathit{7442=98\times 76-5+4-3-2\times 1.}$
\item [] $7443=(9\times 87+(6+5)\times 4)\times 3^2\times 1.$
\item [] $7444=(9\times 87+(6+5)\times 4)\times 3^2+1.$
\item [] $7445=(98+7+6+5)\times 4^3+21.$
\item [] $7446=9\times (8+7+65\times 4)\times 3+21.$
\item [] $7447=(9\times 8+7)\times (6\times 5+4^3)+21.$
\item [] $7448=98\times (7+6+54+3^2\times 1).$
\item [] $7449=9+(8\times 7+6)\times 5\times 4\times 3\times 2\times 1.$
\item [] $7450=9+(8\times 7+6)\times 5\times 4\times 3\times 2+1.$
\item [] $7451=9\times 8+7\times (6\times 5+4^(3+2))+1.$
\item [] $7452=9\times (87+6\times 54+3)\times 2\times 1.$
\item [] $7453=98\times (7+65+4)+3+2\times 1.$
\item [] $7454=98\times (7+65+4)+3+2+1.$
\item [] $7455=(98+7\times 6+5\times 43)\times 21.$
\item [] $7456=98+7+6\times (5\times (4+3))^2+1.$
\item [] $7457=98\times (7+65+4)+3^2\times 1.$
\item [] $7458=98\times (7+65+4)+3^2+1.$
\item [] $7459=(9+(8\times 7+6)\times 5\times 4\times 3)\times 2+1.$
\item [] $7460=(98\times 7+6+54)\times (3^2+1).$
\item [] $7461=(9+8)\times 7\times (6+54)+321.$
\item [] $7462=98\times 76+5+4+3+2\times 1.$
\item [] $7463=98\times 76+5+4+3+2+1.$
\item [] $7464=98\times 76+5+4+3\times 2+1.$
\item [] $7465=9+(8+7)\times (65+432)+1.$
\item [] $7466=98\times 76+5+4+3^2\times 1.$
\item [] $7467=98\times 76+5+4+3^2+1.$
\item [] $7468=98\times 76+5+4\times 3+2+1.$
\item [] $7469=98+(7+6)\times (5+4)\times 3\times 21.$
\item [] $7470=98+76\times ((5+43)\times 2+1).$
\item [] $7471=(98+7)\times 6\times 5+4321.$
\item [] $7472=98\times (7+65+4)+3+21.$
\item [] $7473=98\times 76+5\times 4+3+2\times 1.$
\item [] $7474=98\times 76+5\times 4+3+2+1.$
\item [] $7475=98\times 76+5\times 4+3\times 2+1.$
\item [] $7476=9\times (8+76)+5\times 4^3\times 21.$
\item [] $7477=98\times 76+5+4\times 3\times 2\times 1.$
\item [] $7478=98\times 76+5+4\times 3\times 2+1.$
\item [] $7479=9+(8+7)\times (65+432+1).$
\item [] $7480=98\times (7+65+4)+32\times 1.$
\item[]$\mbox{Increasing order}$
\item [] $7481=1+2^{(3\times 4)}+(5+6\times 7)\times 8\times 9.$
\item [] $7482=12^3\times 4+5\times (6\times 7+8\times 9).$
\item [] $\mathit{7483=-1+2^3\times 4^5-6-78\times 9.}$
\item [] $7484=(1+2^3+4)\times (567+8)+9.$
\item [] $7485=12\times (3+4)\times (5+6+78)+9.$
\item [] $7486=12^3+4^5+6\times 789.$
\item [] $\mathit{7487=-1-234+(5+6)\times 78\times 9.}$
\item [] $7488=12^3+(4+56)\times (7+89).$
\item [] $7489=1+2\times (34+5+6+7)\times 8\times 9.$
\item [] $7490=1\times 2+(3+45)\times (67+89).$
\item [] $7491=1+2+3\times (4\times 5+6)\times (7+89).$
\item [] $\mathit{7492=-1\times 2\times 34+56\times (7+8)\times 9.}$
\item [] $7493=12^3\times 4+5+6\times (7+89).$
\item [] $7494=1\times 2\times 3\times (4\times 5\times (6+7\times 8)+9).$
\item [] $7495=1\times 2^{(3\times 4)}+5\times 678+9.$
\item [] $7496=12^3\times 4+567+8+9.$
\item [] $7497=(12\times 3+4+56)\times 78+9.$
\item [] $7498=1+2^3\times (4+5)\times (6+7)\times 8+9.$
\item [] $7499=(1234+5)\times 6+7\times 8+9.$
\item [] $7500=12+(3+45)\times (67+89).$
\item [] $7501=1+2\times (3+4+5+6\times 7\times 89).$
\item [] $7502=(1+2^3)^4+5+(6+7)\times 8\times 9.$
\item [] $7503=1+2+3\times 4\times 5\times (6+7\times (8+9)).$
\item [] $\mathit{7504=-1+2^3\times 4^5-678-9.}$
\item [] $7505=(123+4^5)\times 6+7\times 89.$
\item [] $7506=(12+345)\times (6+7+8)+9.$
\item [] $7507=1\times 23\times (45\times 6+7\times 8)+9.$
\item [] $7508=1+23\times (45\times 6+7\times 8)+9.$
\item [] $7509=12\times 3\times 4\times 5+6789.$
\item [] $7510=1\times 2\times (3\times 4+5+6\times 7\times 89).$
\item [] $7511=1+2\times (3\times 4+5+6\times 7\times 89).$
\item [] $7512=12+3\times 4\times 5\times (6+7\times (8+9)).$
\item [] $7513=(1234+5)\times 6+7+8\times 9.$
\item [] $7514=1+(2+3+4+5)\times 67\times 8+9.$
\item [] $7515=(123\times 4+5\times 67+8)\times 9.$
\item [] $\mathit{7516=-1+2^3\times 4^5-(67+8)\times 9.}$
\item [] $7517=12^3+4+5\times (6+7)\times 89.$
\item [] $7518=1\times 2\times (3\times 4^5+678+9).$
\item [] $7519=1+2\times (3\times 4^5+678+9).$
\item [] $7520=(12+3+4)\times 56\times 7+8\times 9.$
\item [] $7521=(1234+5)\times 6+78+9.$
\item [] $7522=1\times 2\times (3+4\times 5+6\times 7\times 89).$
\item [] $7523=1+2\times (3+4\times 5+6\times 7\times 89).$
\item [] $7524=12\times (3+4\times 5\times 6+7\times 8\times 9).$
\item [] $7525=1+234\times 5\times 6+7\times 8\times 9.$
\item [] $\mathit{7526=1+2\times 3^4\times (5+6\times 7)-89.}$
\item [] $7527=(1+2)\times (3^4\times 5\times 6+7+8\times 9).$
\item [] $7528=(1^2+3+4)\times (5+(6+7)\times 8\times 9).$
\item [] $7529=1\times (2+3)\times 4\times (5+6\times 7)\times 8+9.$
\item [] $7530=1\times 2\times (3\times (4+5)+6\times 7\times 89).$
\item [] $7531=1\times 2^{(3\times 4)}+5\times (678+9).$
\item [] $7532=1+2^{(3\times 4)}+5\times (678+9).$
\item [] $7533=(12+3\times 45\times 6+7+8)\times 9.$
\item [] $7534=(12\times 3+4^5+6)\times 7+8\times 9.$
\item [] $7535=1\times 2\times 3456+7\times 89.$
\item [] $7536=1+2\times 3456+7\times 89.$
\item [] $7537=12^3\times 4+(5+6)\times 7\times 8+9.$
\item [] $\mathit{7538=12-34+56\times (7+8)\times 9.}$
\item [] $\mathit{7539=-1-(2+3)\times 4+56\times (7+8)\times 9.}$
\item [] $7540=((12+3)\times 4+56)\times (7\times 8+9).$
\item [] $7541=12^3\times 4+(5\times 6+7)\times (8+9).$
\item [] $7542=(12+3)\times 456+78\times 9.$
\item [] $7543=1+((2^3\times 4)\times 5+678)\times 9.$
\item [] $7544=1\times 23\times 4\times (5\times (6+7)+8+9).$
\item [] $7545=12\times (3\times 4+(5+6)\times 7\times 8)+9.$
\item [] $7546=12^3\times 4+5+6+7\times 89.$
\item [] $7547=1+2\times ((3+4)\times 5+6\times 7\times 89).$
\item [] $7548=(1\times 234+5\times 6\times 7)\times (8+9).$
\item [] $7549=1+2\times (3\times 4+5\times 6\times 7)\times (8+9).$
\item [] $7550=1\times 2+3\times 4\times (5\times 6+7)\times (8+9).$
\item[]$\mbox{Decreasing order}$
\item [] $7481=98\times 76+5+4+3+21.$
\item [] $7482=98\times 76+(5+4\times 3)\times 2\times 1.$
\item [] $7483=(987+65\times 4)\times 3\times 2+1.$
\item [] $7484=98\times 76+5+4+3^{(2+1)}.$
\item [] $7485=(9+8\times 76+5)\times 4\times 3+21.$
\item [] $7486=98\times 76+5+4\times 3+21.$
\item [] $\mathit{7487=98\times 76+5\times 4\times 3-21.}$
\item [] $7488=98\times (7+65)+432\times 1.$
\item [] $7489=98\times 76+5+4+32\times 1.$
\item [] $7490=98\times 76+5+4+32+1.$
\item [] $7491=(9+87)\times (6+5\times 4)\times 3+2+1.$
\item [] $7492=98\times 76+5\times 4+3+21.$
\item [] $7493=98\times 76+5+4\times (3^2+1).$
\item [] $7494=98\times 76+(5\times 4+3)\times 2\times 1.$
\item [] $7495=98\times 76+(5\times 4+3)\times 2+1.$
\item [] $7496=98\times 76+(5+4)\times 3+21.$
\item [] $7497=(9+8\times 7\times 6+5+4+3)\times 21.$
\item [] $7498=98\times 76+5+43+2\times 1.$
\item [] $7499=98\times 76+5+43+2+1.$
\item [] $7500=98\times 76+5\times 4+32\times 1.$
\item [] $7501=98\times 76+5\times 4+32+1.$
\item [] $7502=98\times 76+5+(4+3)^2\times 1.$
\item [] $7503=987+6\times 543\times 2\times 1.$
\item [] $7504=987+6\times 543\times 2+1.$
\item [] $7505=(9+87\times 6+5)\times (4+3)\times 2+1.$
\item [] $7506=9+(8+7+6)\times (5+4\times 3)\times 21.$
\item [] $7507=98\times 76+54+3+2^1.$
\item [] $7508=98\times 76+54+3+2+1.$
\item [] $7509=9\times 87+6+5\times 4^3\times 21.$
\item [] $7510=98\times 76+5\times 4\times 3+2\times 1.$
\item [] $7511=98\times 76+54+3^2\times 1.$
\item [] $7512=98\times 76+54+3^2+1.$
\item [] $7513=9+87\times 65+43^2\times 1.$
\item [] $7514=9+87\times 65+43^2+1.$
\item [] $7515=9+(876+5^4)\times (3+2)+1.$
\item [] $7516=9\times ((8\times (7+6)+(5+4)^3)+2)+1.$
\item [] $7517=98\times 76+5+43+21.$
\item [] $7518=98\times 7\times 6+54\times 3\times 21.$
\item [] $7519=98\times 76+5+4^3+2\times 1.$
\item [] $7520=98\times 76+5+4+3\times 21.$
\item [] $7521=(9+8\times 7\times (6+5))\times 4\times 3+21.$
\item [] $\mathit{7522=-98\times 7+6^5+432\times 1.}$
\item [] $7523=98\times 76+5\times (4\times 3+2+1).$
\item [] $7524=(9+87)\times 65+4\times 321.$
\item [] $7525=(9\times 8\times 7+6\times 543)\times 2+1.$
\item [] $7526=98\times 76+54+3+21.$
\item [] $7527=987+654\times (3^2+1).$
\item [] $7528=(98+7\times (6+5))\times 43+2+1.$
\item [] $7529=98\times 76+5\times 4\times 3+21.$
\item [] $7530=98\times 76+(5+4)\times 3^2+1.$
\item [] $7531=98\times 76+5\times 4+3\times 21.$
\item [] $7532=987+6543+2\times 1.$
\item [] $7533=987+6543+2+1.$
\item [] $7534=98\times 76+54+32\times 1.$
\item [] $7535=98\times 76+54+32+1.$
\item [] $7536=(9\times (8+(765+4^3)))+(2+1).$
\item [] $7537=9+8+(7+6\times 5^4+3)\times 2\times 1.$
\item [] $7538=98\times 76+5+4^3+21.$
\item [] $7539=98\times 76+5+43\times 2\times 1.$
\item [] $7540=98\times 76+5+43\times 2+1.$
\item [] $7541=(9+8\times 7)\times (6\times 5+43\times 2)+1.$
\item [] $7542=9+(87+6)\times (5\times 4\times 3+21).$
\item [] $7543=9+(87+6)\times (5+4)\times 3^2+1.$
\item [] $7544=98\times 76+(5+43)\times 2\times 1.$
\item [] $7545=98\times 76+(5+43)\times 2+1.$
\item [] $7546=98\times (7+6+54+3^2+1).$
\item [] $7547=987\times 6+5\times (4+321).$
\item [] $7548=98\times 76+5\times 4\times (3+2\times 1).$
\item [] $7549=98\times 76+5+4\times (3+21).$
\item [] $7550=(98\times 7+65+4)\times (3^2+1).$
\item[]$\mbox{Increasing order}$
\item [] $7551=12^3\times 4+567+8\times 9.$
\item [] $7552=1^2+(3\times (45\times 6+7)+8)\times 9.$
\item [] $7553=1\times 23\times (4+5\times 6+7)\times 8+9.$
\item [] $7554=1\times 2\times (34+5+6\times 7\times 89).$
\item [] $7555=1+2\times (34+5+6\times 7\times 89).$
\item [] $\mathit{7556=1-2-3+(4+5+6)\times 7\times 8\times 9.}$
\item [] $7557=(1+2)\times (3+4\times (5\times 6+7)\times (8+9)).$
\item [] $7558=(1+2\times 3)\times 4^5+6\times (7\times 8+9).$
\item [] $\mathit{7559=1-2\times 3+4+56\times (7+8)\times 9.}$
\item [] $7560=1^{234}\times 56\times (7+8)\times 9.$
\item [] $7561=1+2\times (345+67+8)\times 9.$
\item [] $7562=1\times 2+3\times (4\times 56+7\times 8)\times 9.$
\item [] $7563=1+2+3\times (4\times 56+7\times 8)\times 9.$
\item [] $7564=1^{23}\times 4+56\times (7+8)\times 9.$
\item [] $7565=12^3\times 4+5\times 6+7\times 89.$
\item [] $7566=1+2+3+(4+5+6)\times 7\times 8\times 9.$
\item [] $7567=12^3\times 4+5\times (6\times 7+89).$
\item [] $7568=12^3\times 4+567+89.$
\item [] $7569=1\times 2+3+4+56\times (7+8)\times 9.$
\item [] $7570=1+2+3+4+56\times (7+8)\times 9.$
\item [] $7571=1+2\times 3+4+56\times (7+8)\times 9.$
\item [] $7572=12+3\times (4\times 56+7\times 8)\times 9.$
\item [] $7573=1^2+3\times 4+56\times (7+8)\times 9.$
\item [] $7574=1\times 2+3\times 4+56\times (7+8)\times 9.$
\item [] $7575=1+2+3\times 4+56\times (7+8)\times 9.$
\item [] $7576=12^3+4^5+67\times 8\times 9.$
\item [] $7577=12\times (34+56)\times 7+8+9.$
\item [] $7578=1+(2+34)\times 5\times 6\times 7+8+9.$
\item [] $7579=12+3+4+56\times (7+8)\times 9.$
\item [] $7580=1\times (2+3)\times 4+56\times (7+8)\times 9.$
\item [] $7581=(1+2)\times (34\times 56+7\times 89).$
\item [] $7582=(1+(2\times 34+5)\times 6+7)\times (8+9).$
\item [] $7583=12\times 3\times 45+67\times 89.$
\item [] $7584=12+3\times 4+56\times (7+8)\times 9.$
\item [] $7585=1+2\times 3\times 4+56\times (7+8)\times 9.$
\item [] $7586=1\times 2\times ((3^4+56\times 7)\times 8+9).$
\item [] $7587=1\times 23+4+56\times (7+8)\times 9.$
\item [] $7588=1+23+4+56\times (7+8)\times 9.$
\item [] $7589=1\times 2+(3+(4+5+6)\times 7\times 8)\times 9.$
\item [] $7590=(12+34)\times (5\times 6+(7+8)\times 9).$
\item [] $7591=12^3\times 4+56+7\times 89.$
\item [] $7592=12^3\times 4+5+(67+8)\times 9.$
\item [] $7593=1+2^3\times 4+56\times (7+8)\times 9.$
\item [] $7594=1^2\times 34+56\times (7+8)\times 9.$
\item [] $7595=1^2+34+56\times (7+8)\times 9.$
\item [] $7596=1\times 2+34+56\times (7+8)\times 9.$
\item [] $7597=1+2+34+56\times (7+8)\times 9.$
\item [] $\mathit{7598=-1+2\times 3^4\times 5+6789.}$
\item [] $7599=1\times 2\times 3^4\times 5+6789.$
\item [] $7600=1+2\times 3^4\times 5+6789.$
\item [] $7601=1\times (2\times 34+5)\times (6+7)\times 8+9.$
\item [] $7602=1+(2\times 34+5)\times (6+7)\times 8+9.$
\item [] $7603=1+2\times 3+(4+56\times (7+8))\times 9.$
\item [] $7604=12^3\times 4+5+678+9.$
\item [] $7605=1\times (2+3+4)\times (56+789).$
\item [] $7606=12+34+56\times (7+8)\times 9.$
\item [] $7607=12^3\times 4+5\times (67+8\times 9).$
\item [] $7608=1^2\times 3\times 4\times (5+6+7\times 89).$
\item [] $7609=1^2+3\times 4\times (5+6+7\times 89).$
\item [] $7610=1\times 2+3\times 4\times (5+6+7\times 89).$
\item [] $7611=1+2+3\times 4\times (5+6+7\times 89).$
\item [] $7612=(1+2\times 3+4)\times (5+678+9).$
\item [] $\mathit{7613=-1+2\times 3456+78\times 9.}$
\item [] $7614=1\times 2\times 3456+78\times 9.$
\item [] $7615=1+2\times 3456+78\times 9.$
\item [] $7616=(123+4)\times 56+7\times 8\times 9.$
\item [] $7617=1^2+(34+5\times 6)\times 7\times (8+9).$
\item [] $7618=1\times 2+(34+5\times 6)\times 7\times (8+9).$
\item [] $7619=1+2+(34+5\times 6)\times 7\times (8+9).$
\item [] $7620=12+3\times 4\times (5+6+7\times 89).$
\item[]$\mbox{Decreasing order}$
\item [] $7551=987+6543+21.$
\item [] $7552=(9+8\times (7+6)+5)\times (43+21).$
\item [] $7553=98\times 7+(6\times 54+3)\times 21.$
\item [] $7554=(9+8+7+6\times 5^4+3)\times 2\times 1.$
\item [] $7555=98\times 7\times (6+5)+4+3+2\times 1.$
\item [] $7556=98\times 7\times (6+5)+4+3+2+1.$
\item [] $7557=98\times 7\times (6+5)+4+3\times 2+1.$
\item [] $\mathit{7558=9\times 8\times 7\times (6+5+4)-3+2-1.}$
\item [] $7559=98\times 7\times (6+5)+4+3^2\times 1.$
\item [] $7560=98\times 7\times (6+5)+4+3^2+1.$
\item [] $7561=98\times 7\times (6+5)+4\times 3+2+1.$
\item [] $7562=98\times 76+(54+3)\times 2\times 1.$
\item [] $7563=98\times 76+(54+3)\times 2+1.$
\item [] $7564=(987+65\times 43)\times 2\times 1.$
\item [] $7565=98\times 76+54+3\times 21.$
\item [] $7566=(98+7\times 6)\times 54+3+2+1.$
\item [] $7567=(98+7\times 6)\times 54+3\times 2+1.$
\item [] $7568=98\times 76+5\times 4\times 3\times 2\times 1.$
\item [] $7569=98\times 76+5\times 4\times 3\times 2+1.$
\item [] $7570=98\times 7\times (6+5)+4\times 3\times 2\times 1.$
\item [] $7571=98\times 7\times (6+5)+4\times 3\times 2+1.$
\item [] $7572=9+8\times 7+(6\times 5^4+3)\times 2+1.$
\item [] $7573=98\times 76+5\times (4\times 3\times 2+1).$
\item [] $7574=98\times 7\times (6+5)+4+3+21.$
\item [] $7575=(9\times (8+76)\times 5+4+3)\times 2+1.$
\item [] $7576=9+(8+7+6\times (5^4+3))\times 2+1.$
\item [] $7577=9+8+7\times 6\times 5\times 4\times 3^2\times 1.$
\item [] $7578=9+8+7\times 6\times 5\times 4\times 3^2+1.$
\item [] $7579=98\times 7\times (6+5)+4\times 3+21.$
\item [] $7580=9\times (87\times 6+5\times 4^3)+2\times 1.$
\item [] $7581=98\times 76+5+4^3\times 2\times 1.$
\item [] $7582=98\times 76+5+4\times 32+1.$
\item [] $7583=98\times 7\times (6+5)+4+32+1.$
\item [] $7584=9+8+7+6\times 5\times 4\times 3\times 21.$
\item [] $7585=98\times 76+5+4\times (32+1).$
\item [] $7586=98\times 76+(5+4^3)\times 2\times 1.$
\item [] $7587=9\times 87\times 6+(5+4)\times 321.$
\item [] $7588=98\times 76+5\times (4+3+21).$
\item [] $7589=9\times (8+76\times 5)+4^{(3\times 2)}+1.$
\item [] $7590=(98+7+654)\times (3^2+1).$
\item [] $7591=98\times 7\times (6+5)+43+2\times 1.$
\item [] $7592=98\times 7\times (6+5)+43+2+1.$
\item [] $7593=(98+7\times 6)\times 54+32+1.$
\item [] $\mathit{7594=-9-87+6^5-43\times 2\times 1.}$
\item [] $7595=98\times 7\times (6+5)+(4+3)^2\times 1.$
\item [] $7596=9\times (8\times 76+5\times 43+21).$
\item [] $7597=98\times 76+5+(4\times 3)^2\times 1.$
\item [] $7598=98\times 76+5+(4\times 3)^2+1.$
\item [] $7599=(9+8)\times (76\times 5+4+3\times 21).$
\item [] $7600=98\times 76+5+(4+3)\times 21.$
\item [] $7601=9+8\times 7+6\times (5^4+3)\times 2\times 1.$
\item [] $7602=(9+8\times 7\times 6+5+4\times 3)\times 21.$
\item [] $7603=9+87+(6\times 5^4+3)\times 2+1.$
\item [] $\mathit{7604=-9-8+7654-32-1.}$
\item [] $7605=9+876+5\times 4^3\times 21.$
\item [] $7606=9\times (8\times 7\times (6+5+4)+3+2)+1.$
\item [] $7607=9+8\times 7+6\times ((5^4+3)\times 2+1).$
\item [] $7608=(9+8+7)\times (65+4\times 3\times 21).$
\item [] $7609=9+8\times (7\times 6+5+43\times 21).$
\item [] $7610=(9\times (8+76)+5)\times (4+3+2+1).$
\item [] $7611=9\times (8+7)\times 6\times (5+4)+321.$
\item [] $7612=98\times 76+54\times 3+2\times 1.$
\item [] $7613=98\times 76+54\times 3+2+1.$
\item [] $7614=9\times (87\times 6+54\times 3\times 2\times 1).$
\item [] $7615=9\times 8+7+6\times (5^4+3)\times 2\times 1.$
\item [] $7616=9\times 8+7+6\times (5^4+3)\times 2+1.$
\item [] $7617=9+8+76\times 5\times 4\times (3+2\times 1).$
\item [] $7618=9+8+76\times 5\times 4\times (3+2)+1.$
\item [] $7619=(9\times 8+7\times 6+5)\times 4^3+2+1.$
\item [] $7620=(9\times (8+7)+6)\times 54+3+2+1.$
\item[]$\mbox{Increasing order}$
\item [] $\mathit{7621=12^3\times 4-5+6\times 7\times (8+9).}$
\item [] $\mathit{7622=-1^2+(3+4+56\times  (7+8))\times 9.}$
\item [] $7623=12^3+45\times (6\times 7+89).$
\item [] $7624=(1+(2\times 3)^4)\times 5+67\times (8+9).$
\item [] $7625=12^3\times 4+5+6+78\times 9.$
\item [] $7626=1+2\times 34\times (56+7\times 8)+9.$
\item [] $7627=12^3\times 4+(5+6)\times (7\times 8+9).$
\item [] $7628=1\times 2\times 34+56\times (7+8)\times 9.$
\item [] $7629=(12+3)\times 456+789.$
\item [] $7630=(1\times 2+3+4+5)\times (67\times 8+9).$
\item [] $7631=12^3\times 4+5+6\times 7\times (8+9).$
\item [] $7632=1^2\times (34\times 5+678)\times 9.$
\item [] $7633=1+(2+34)\times 5\times 6\times 7+8\times 9.$
\item [] $7634=1+(2+3^4\times 5+6\times 7)\times (8+9).$
\item [] $7635=1^2\times 3+(4\times 5\times 6\times 7+8)\times 9.$
\item [] $7636=12\times 3^4+56\times 7\times (8+9).$
\item [] $7637=1+2^{(3\times 4)}+5\times (6+78\times 9).$
\item [] $7638=1+(2\times 3)^4\times 5+(6+7)\times 89.$
\item [] $7639=(1+(2+34)\times 5\times 6)\times 7+8\times 9.$
\item [] $7640=1\times 2^3+(4\times 5\times 6\times 7+8)\times 9.$
\item [] $7641=(1+2)\times 3\times (4+56+789).$
\item [] $7642=1234+(5+67)\times 89.$
\item [] $7643=1\times 234\times 5\times 6+7\times 89.$
\item [] $7644=1+234\times 5\times 6+7\times 89.$
\item [] $7645=(1+2\times 3)\times 4^5+6\times 78+9.$
\item [] $\mathit{7646=12^3-45+67\times 89.}$
\item [] $7647=12+3+(4\times 5\times 6\times 7+8)\times 9.$
\item [] $7648=1\times 2\times (3^4+5+6\times 7\times 89).$
\item [] $7649=12\times (34+56)\times 7+89.$
\item [] $7650=1+(2+34)\times 5\times 6\times 7+89.$
\item [] $7651=1+(2\times 3+4+5)\times (6+7\times 8\times 9).$
\item [] $7652=1\times 23\times 4+56\times (7+8)\times 9.$
\item [] $7653=12+3^4+56\times (7+8)\times 9.$
\item [] $7654=(1\times 2\times 34+5+6+7)\times 89.$
\item [] $7655=1\times 23+(4\times 5\times 6\times 7+8)\times 9.$
\item [] $7656=1+23+(4\times 5\times 6\times 7+8)\times 9.$
\item [] $7657=1+(2^3\times 4+56)\times (78+9).$
\item [] $\mathit{7658=-1+23\times (-456+789).}$
\item [] $7659=(1+2+34\times 5+678)\times 9.$
\item [] $7660=((1+2)^3\times 4+5)\times 67+89.$
\item [] $7661=(12+34\times 5)\times 6\times 7+8+9.$
\item [] $7662=12+(3\times 4+5)\times (6\times 7+8)\times 9.$
\item [] $7663=1+2\times (3+4\times (5+6)\times (78+9)).$
\item [] $7664=(1\times 2+3)^4\times (5+6)+789.$
\item [] $7665=12\times (34\times 5+6\times 78)+9.$
\item [] $7666=1\times 2^{(3\times 4)}+5\times 6\times 7\times (8+9).$
\item [] $7667=1+2^{(3\times 4)}+5\times 6\times 7\times (8+9).$
\item [] $7668=1^2\times 3\times 4\times (567+8\times 9).$
\item [] $7669=1^2+3\times 4\times (567+8\times 9).$
\item [] $7670=12^3\times 4+56+78\times 9.$
\item [] $7671=(123+4^5)\times 6+789.$
\item [] $7672=1\times 2^{(3+4)}\times 56+7\times 8\times 9.$
\item [] $7673=(1+234)\times 5\times 6+7\times 89.$
\item [] $7674=(1+(2+34)\times 5)\times 6\times 7+8\times 9.$
\item [] $7675=(1+(2+3)^4)\times (5+6)+789.$
\item [] $\mathit{7676=-12-34+(5+6)\times 78\times 9.}$
\item [] $7677=12^3\times 4+((5+6)\times 7+8)\times 9.$
\item [] $7678=(1+2\times 3)\times 4^5+6+7\times 8\times 9.$
\item [] $\mathit{7679=123-4+56\times (7+8)\times 9.}$
\item [] $7680=12+3\times 4\times (567+8\times 9).$
\item [] $7681=1+2^3\times (456+7\times 8\times 9).$
\item [] $7682=1^2+(3^4+56)\times 7\times 8+9.$
\item [] $7683=123+(4+5+6)\times 7\times 8\times 9.$
\item [] $7684=1+2+(3^4+56)\times 7\times 8+9.$
\item [] $7685=1^2+34\times (5+(6+7)\times (8+9)).$
\item [] $7686=1\times 2\times 3^4\times (5+6\times 7)+8\times 9.$
\item [] $7687=123+4+56\times (7+8)\times 9.$
\item [] $7688=1\times 2^{(3+4)}+56\times (7+8)\times 9.$
\item [] $7689=(12+3)\times (456+7\times 8)+9.$
\item [] $7690=(1+2\times 3)\times 4^5+6\times (78+9).$
\item[]$\mbox{Decreasing order}$
\item [] $7621=(9\times (8+7)+6)\times 54+3\times 2+1.$
\item [] $7622=((98+7)\times 6+5)\times 4\times 3+2\times 1.$
\item [] $7623=(9\times 8+76+5\times 43)\times 21.$
\item [] $7624=(9\times (8+7)+6)\times 54+3^2+1.$
\item [] $7625=9+8\times 7+6\times 5\times 4\times 3\times 21.$
\item [] $7626=9\times (8\times 7+65)\times (4+3)+2+1.$
\item [] $7627=9+(8\times 7+6\times 5^4+3)\times 2\times 1.$
\item [] $7628=98\times 76+5\times 4\times 3^2\times 1.$
\item [] $7629=98\times 76+5\times 4\times 3^2+1.$
\item [] $\mathit{7630=-9-8+7654-3\times 2-1.}$
\item [] $7631=98\times 76+54\times 3+21.$
\item [] $7632=9\times 8\times 76+5\times 432\times 1.$
\item [] $7633=9\times 8\times 76+5\times 432+1.$
\item [] $7634=9+8+7\times (6\times 5+4)\times 32+1.$
\item [] $7635=9+(87+6)\times ((5+4)\times 3^2+1).$
\item [] $7636=(9+8\times 7+6\times 5^4+3)\times 2\times 1.$
\item [] $7637=(9\times 8+7\times 6+5)\times 4^3+21.$
\item [] $7638=9\times (8+7\times 6\times 5\times 4)+3+2+1.$
\item [] $7639=9\times 8+7+6\times 5\times 4\times 3\times 21.$
\item [] $7640=9+8+7\times (65+4^(3+2\times 1)).$
\item [] $7641=98+7+6\times (5^4+3)\times 2\times 1.$
\item [] $7642=98+7+6\times (5^4+3)\times 2+1.$
\item [] $\mathit{7643=-9-8+7654+3+2+1.}$
\item [] $7644=9\times (8\times 7+65)\times (4+3)+21.$
\item [] $7645=98\times 76+5+4^3\times (2+1).$
\item [] $7646=(9\times (8+76)\times 5+43)\times 2\times 1.$
\item [] $7647=(9\times (8+7)+6)\times 54+32+1.$
\item [] $7648=98\times 76+5\times 4\times (3^2+1).$
\item [] $7649=9+(8+76\times 5\times 4)\times (3+2)\times 1.$
\item [] $7650=987\times 6+54\times 32\times 1.$
\item [] $7651=987\times 6+54\times 32+1.$
\item [] $7652=9\times (8+7\times 6)\times (5+4\times 3)+2\times 1.$
\item [] $7653=9+(8+76)\times (5+43\times 2\times 1).$
\item [] $7654=98\times (7+6)\times 5+4\times 321.$
\item [] $7655=(9+8+7+65)\times 43\times 2+1.$
\item [] $7656=9+87+6\times 5\times 4\times 3\times 21.$
\item [] $7657=(98+7+6+5)\times (4^3+2)+1.$
\item [] $7658=98+7\times 6\times 5\times 4\times 3^2\times 1.$
\item [] $7659=98+7\times 6\times 5\times 4\times 3^2+1.$
\item [] $7660=(9\times (8+7)+6+5^4)\times (3^2+1).$
\item [] $7661=9+8+7\times (6+543\times 2\times 1).$
\item [] $7662=9+8+7\times (6+543\times 2)+1.$
\item [] $7663=9+(8+76+5)\times 43\times 2\times 1.$
\item [] $7664=9+(8+76+5)\times 43\times 2+1.$
\item [] $7665=98+7+6\times 5\times 4\times 3\times 21.$
\item [] $7666=98\times 76+5\times 43+2+1.$
\item [] $7667=(9+8)\times (76+54+321).$
\item [] $7668=(9\times 87+65+4)\times 3^2\times 1.$
\item [] $7669=(9\times 87+65+4)\times 3^2+1.$
\item [] $\mathit{7670=-9-8+7654+32+1.}$
\item [] $7671=9\times (8+7\times 6)\times (5+4\times 3)+21.$
\item [] $7672=9\times 8+76\times 5\times 4\times (3+2\times 1).$
\item [] $7673=98\times 76+5\times (43+2\times 1).$
\item [] $7674=98\times 76+5\times (43+2)+1.$
\item [] $7675=98\times 7\times (6+5)+43\times (2+1).$
\item [] $7676=9+8+7654+3+2\times 1.$
\item [] $7677=9+8+7654+3+2+1.$
\item [] $7678=9+8+7654+3\times 2+1.$
\item [] $7679=(9+8)\times 7+6\times 5\times 4\times 3\times 21.$
\item [] $7680=9+8+7654+3^2\times 1.$
\item [] $7681=9+8+7654+3^2+1.$
\item [] $7682=9\times 8\times (76+5)+43^2+1.$
\item [] $7683=(9\times 8\times 7+6)\times 5+4\times 3+21.$
\item [] $7684=98\times 76+5\times 43+21.$
\item [] $7685=(9+8+76\times 5\times 4)\times (3+2\times 1).$
\item [] $7686=9\times 8\times 7+6\times (54+3)\times 21.$
\item [] $7687=9\times (8+(76\times 5+43)\times 2)+1.$
\item [] $7688=9\times 8+7\times (6\times 5+4)\times 32\times 1.$
\item [] $7689=9\times 8+7\times (6\times 5+4)\times 32+1.$
\item [] $7690=9+(87+6\times 5^4+3)\times 2+1.$
\item[]$\mbox{Increasing order}$
\item [] $7691=((1+23)\times 45+6)\times 7+89.$
\item [] $7692=12^3\times 4+5\times (67+89).$
\item [] $7693=12+(3^4+56)\times 7\times 8+9.$
\item [] $\mathit{7694=1\times 2^3\times 4^5+6-7\times 8\times 9.}$
\item [] $7695=(1+2\times 3+4\times 5\times 6\times 7+8)\times 9.$
\item [] $7696=1^2+3^4\times (5\times 6+7\times 8+9).$
\item [] $7697=1\times 2+3^4\times (5\times 6+7\times 8+9).$
\item [] $7698=123\times (4\times 5+6\times 7)+8\times 9.$
\item [] $\mathit{7699=1-2\times 3\times 4+(5+6)\times 78\times 9.}$
\item [] $7700=12^3+4+5+67\times 89.$
\item [] $7701=1\times 2\times 3456+789.$
\item [] $7702=1+2\times 3456+789.$
\item [] $7703=1\times 2\times 3^4\times (5+6\times 7)+89.$
\item [] $7704=(1^2+3)^4\times 5\times 6+7+8+9.$
\item [] $7705=1\times 23\times (45\times 6+7\times 8+9).$
\item [] $7706=1+23\times (45\times 6+7\times 8+9).$
\item [] $7707=12+3^4\times (5\times 6+7\times 8+9).$
\item [] $7708=(1^2+3^4)\times ((5+6)\times 7+8+9).$
\item [] $7709=(1+2\times 3)\times (4^5+67)+8\times 9.$
\item [] $7710=(12+3)\times (4+5\times (6+7+89)).$
\item [] $7711=12^3+4\times 5+67\times 89.$
\item [] $7712=12^3\times 4+5+6+789.$
\item [] $7713=12^3\times 4+(5+6+78)\times 9.$
\item [] $7714=(1+2\times 3)\times (4^5+6)+7\times 8\times 9.$
\item [] $7715=123\times (4\times 5+6\times 7)+89.$
\item [] $7716=(12+34\times 5)\times 6\times 7+8\times 9.$
\item [] $\mathit{7717=-12+3+4+(5+6)\times 78\times 9.}$
\item [] $7718=1^2\times 34\times (5\times 6\times 7+8+9).$
\item [] $7719=123+(4+56\times (7+8))\times 9.$
\item [] $7720=1\times 2+34\times (5\times 6\times 7+8+9).$
\item [] $7721=1+2+34\times (5\times 6\times 7+8+9).$
\item [] $7722=12\times 3\times 45+678\times 9.$
\item [] $7723=1+234\times 5\times 6+78\times 9.$
\item [] $\mathit{7724=1+2+3-4+(5+6)\times 78\times 9.}$
\item [] $7725=(1+2+3\times 4)\times (5+6+7\times 8\times 9).$
\item [] $7726=1^{23}\times 4+(5+6)\times 78\times 9.$
\item [] $7727=1^{23}+4+(5+6)\times 78\times 9.$
\item [] $7728=(1+234\times 5)\times 6+78\times 9.$
\item [] $7729=1^2\times 3+4+(5+6)\times 78\times 9.$
\item [] $7730=12+34\times (5\times 6\times 7+8+9).$
\item [] $7731=12^3\times 4+5\times 6+789.$
\item [] $7732=1\times 2\times 3+4+(5+6)\times 78\times 9.$
\item [] $7733=(12+34\times 5)\times 6\times 7+89.$
\item [] $7734=1^2\times 3\times 4+(5+6)\times 78\times 9.$
\item [] $7735=(123+4)\times 56+7\times 89.$
\item [] $7736=12^3+45+67\times 89.$
\item [] $7737=(123+4+5+6)\times 7\times 8+9.$
\item [] $7738=1^2+(3+4\times 5)\times 6\times 7\times 8+9.$
\item [] $7739=1\times 2+(3+4\times 5)\times 6\times 7\times 8+9.$
\item [] $7740=(12+34\times 5+678)\times 9.$
\item [] $7741=12+3+4+(5+6)\times 78\times 9.$
\item [] $7742=1\times (2+3)\times 4+(5+6)\times 78\times 9.$
\item [] $7743=1^{23}\times (45+6\times 7)\times 89.$
\item [] $7744=1^{23}+(4\times 5+67)\times 89.$
\item [] $7745=12\times (3^4+5+6)\times 7+8+9.$
\item [] $7746=12+3\times 4+(5+6)\times 78\times 9.$
\item [] $7747=1+2\times 3\times 4+(5+6)\times 78\times 9.$
\item [] $7748=1\times 2+3+(45+6\times 7)\times 89.$
\item [] $7749=12+(3+4\times 5)\times 6\times 7\times 8+9.$
\item [] $7750=1+23+4+(5+6)\times 78\times 9.$
\item [] $7751=1\times 2^3+(45+6\times 7)\times 89.$
\item [] $7752=12\times (3+4+567+8\times 9).$
\item [] $7753=(1+2^3\times 4\times 5\times 6+7)\times 8+9.$
\item [] $7754=1\times 2^3\times 4+(5+6)\times 78\times 9.$
\item [] $7755=123+(4\times 5\times 6\times 7+8)\times 9.$
\item [] $7756=1^2\times 34+(5+6)\times 78\times 9.$
\item [] $7757=12^3\times 4+56+789.$
\item [] $7758=12+3+(45+6\times 7)\times 89.$
\item [] $7759=1+2+34+(5+6)\times 78\times 9.$
\item [] $7760=(12\times 3+4)\times (5+(6+7+8)\times 9).$
\item[]$\mbox{Decreasing order}$
\item [] $7691=98\times 7\times (6+5)+(4\times 3)^2+1.$
\item [] $\mathit{7692=9-8\times 7+6^5-4-32-1.}$
\item [] $7693=98\times 76+5\times (4+3)^2\times 1.$
\item [] $7694=98\times (7+6)+5\times 4\times 321.$
\item [] $7695=9+8+7654+3+21.$
\item [] $7696=(9\times 8+76)\times (5\times 4+32\times 1).$
\item [] $7697=(9\times 8+76)\times (5\times 4+32)+1.$
\item [] $7698=(9+87+6\times 5^4+3)\times 2\times 1.$
\item [] $7699=(98+76+5)\times 43+2\times 1.$
\item [] $7700=98\times 76+(5+4+3)\times 21.$
\item [] $7701=(98\times 7+6+5^4\times 3)\times (2+1).$
\item [] $7702=9\times 8+7\times (65+4^(3+2)+1).$
\item [] $7703=9+8+7654+32\times 1.$
\item [] $7704=9+8+7654+32+1.$
\item [] $7705=98\times 76+5+4\times 3\times 21.$
\item [] $7706=98\times 7+65\times 4\times 3^{(2+1)}.$
\item [] $7707=9+8+(765+4)\times (3^2+1).$
\item [] $\mathit{7708=9\times 8+7654+3-21.}$
\item [] $\mathit{7709=9-8\times 7+6^5+4-3-21.}$
\item [] $7710=(9+8+7+6)\times (5+4\times 3\times 21).$
\item [] $7711=9+(8+7\times (6+543))\times 2\times 1.$
\item [] $7712=(9+8)\times 76+5\times 4\times 321.$
\item [] $7713=987+6+5\times 4^3\times 21.$
\item [] $7714=98\times (7+6)\times 5+4^3\times 21.$
\item [] $7715=98+7\times (6\times 5+4)\times 32+1.$
\item [] $7716=9\times 8+7\times (6+543\times 2\times 1).$
\item [] $7717=98\times 7\times 6+(5\times 4\times 3)^2+1.$
\item [] $7718=(98+76+5)\times 43+21.$
\item [] $7719=98\times 76+54\times (3+2)+1.$
\item [] $7720=(9+8+7\times (6+543))\times 2\times 1.$
\item [] $7721=98+7\times (6+(5+4)\times 3)^2\times 1.$
\item [] $7722=9\times (8+765+4^3+21).$
\item [] $7723=9\times 8+765\times (4+3\times 2)+1.$
\item [] $7724=9+8+(7+6\times 5\times 4\times 3)\times 21.$
\item [] $7725=(9\times 8\times 7+6+5)\times (4\times 3+2+1).$
\item [] $\mathit{7726=9\times 8+7654-3+2+1.}$
\item [] $\mathit{7727=9\times 8+7654+3-2\times 1.}$
\item [] $7728=(98+7)\times 65+43\times 21.$
\item [] $7729=(9+87+6\times (5^4+3))\times 2+1.$
\item [] $7730=(9\times 8+76+5^4)\times (3^2+1).$
\item [] $7731=9\times 8+7654+3+2\times 1.$
\item [] $7732=9\times 8+7654+3+2+1.$
\item [] $7733=9\times 8+7654+3\times 2+1.$
\item [] $7734=9+8+7654+3\times 21.$
\item [] $7735=9\times 8+7654+3^2\times 1.$
\item [] $7736=9\times 8+7654+3^2+1.$
\item [] $7737=987+(6+5+4)^3\times 2\times 1.$
\item [] $7738=987+(6+5+4)^3\times 2+1.$
\item [] $7739=9+8\times 7\times 6\times (5\times 4+3)+2\times 1.$
\item [] $7740=9\times (8+765+43\times 2+1).$
\item [] $7741=(9\times 8+7+6+5)\times 43\times 2+1.$
\item [] $7742=98+7\times (6+543\times 2\times 1).$
\item [] $7743=98+7\times (6+543\times 2)+1.$
\item [] $7744=(9+8)\times 76\times 5+4\times 321.$
\item [] $7745=98\times (7+65+4+3)+2+1.$
\item [] $7746=(9+8)\times 7\times 65+4+3\times 2+1.$
\item [] $7747=(98+7+6\times (5^4+3))\times 2+1.$
\item [] $7748=98+765\times (4+3+2+1).$
\item [] $7749=(9+8)\times 7\times 65+4+3^2+1.$
\item [] $7750=9\times 8+7654+3+21.$
\item [] $7751=((9+8)\times 7+6)\times (5\times 4\times 3+2)+1.$
\item [] $7752=9+87\times (65+4\times 3\times 2\times 1).$
\item [] $7753=9+87\times (65+4\times 3\times 2)+1.$
\item [] $7754=9+(8+7+6\times 5+43)^2+1.$
\item [] $7755=9+(8\times 7+65)\times 4^3+2\times 1.$
\item [] $7756=9+(8\times 7+65)\times 4^3+2+1.$
\item [] $7757=98+7654+3+2\times 1.$
\item [] $7758=98+7654+3+2+1.$
\item [] $7759=9\times 8+7654+32+1.$
\item [] $7760=(9+8)\times 7\times 65+4\times 3\times 2+1.$
\item[]$\mbox{Increasing order}$
\item [] $7761=12^3\times 4+56\times (7+8)+9.$
\item [] $7762=12\times 3+4+(5+6)\times 78\times 9.$
\item [] $7763=1\times 2+3+(4+(5+6)\times 78)\times 9.$
\item [] $7764=12\times (34\times 5+6\times 78+9).$
\item [] $7765=1+2\times 3+(4+(5+6)\times 78)\times 9.$
\item [] $7766=12\times 3^4+5+6789.$
\item [] $7767=1+23+(45+6\times 7)\times 89.$
\item [] $7768=12+34+(5+6)\times 78\times 9.$
\item [] $7769=(1+2\times 3\times (4+5+67))\times (8+9).$
\item [] $7770=(1+2)^3+(4\times 5+67)\times 89.$
\item [] $7771=(12+3+4)\times (56\times 7+8+9).$
\item [] $\mathit{7772=(1-2+3+4)^5+6+7-8-9.}$
\item [] $7773=12+3+(4+(5+6)\times 78)\times 9.$
\item [] $\mathit{7774=1\times 23\times (4+5\times 67+8-9).}$
\item [] $\mathit{7775=-1+2\times 3^4\times (5+6\times 7-8+9).}$
\item [] $7776=12\times 3\times (4+5\times 6\times 7)+8\times 9.$
\item [] $7777=1+2^3\times (45\times 6+78\times 9).$
\item [] $7778=1\times 2+(3+(4+5+6)\times 7)\times 8\times 9.$
\item [] $7779=12\times 3+(45+6\times 7)\times 89.$
\item [] $7780=1+(2+(3^4+5)\times 6)\times (7+8)+9.$
\item [] $7781=12^3\times 4+(5+6)\times (7+8\times 9).$
\item [] $7782=(12+3)\times 4+(5+6)\times 78\times 9.$
\item [] $7783=1+2\times (3+4\times (5\times 6+78)\times 9).$
\item [] $7784=(12\times 3+4\times 5)\times (67+8\times 9).$
\item [] $7785=12\times (34+5+6\times 7)\times 8+9.$
\item [] $7786=(1+(2+3^4)\times 5+6\times 7)\times (8+9).$
\item [] $7787=(1+2+34)\times 5\times 6\times 7+8+9.$
\item [] $7788=(123+4+5)\times (6\times 7+8+9).$
\item [] $\mathit{7789=-1-23+4^5+6789.}$
\item [] $7790=1\times 2\times 34+(5+6)\times 78\times 9.$
\item [] $7791=1+2\times 34+(5+6)\times 78\times 9.$
\item [] $7792=1+2^{(3+4)}\times 56+7\times 89.$
\item [] $7793=12\times 3\times (4+5\times 6\times 7)+89.$
\item [] $7794=1\times 234+56\times (7+8)\times 9.$
\item [] $7795=1+234+56\times (7+8)\times 9.$
\item [] $7796=1\times 2+(3\times 45\times 6+7\times 8)\times 9.$
\item [] $7797=(1+2\times 3)\times 4^5+6+7\times 89.$
\item [] $\mathit{7798=-12-3+4^5+6789.}$
\item [] $7799=(1^2+3)^4\times 5\times 6+7\times (8+9).$
\item [] $7800=12\times (3^4+5+6)\times 7+8\times 9.$
\item [] $7801=1^{23}+4\times 5\times 6\times (7\times 8+9).$
\item [] $7802=(1\times 2+3^4)\times ((5+6)\times 7+8+9).$
\item [] $7803=1^2\times 3+4\times 5\times 6\times (7\times 8+9).$
\item [] $7804=1^2+3+4\times 5\times 6\times (7\times 8+9).$
\item [] $7805=1\times 2+3+4\times 5\times 6\times (7\times 8+9).$
\item [] $7806=1^2\times 3\times 4^5+6\times 789.$
\item [] $7807=1^2+3\times 4^5+6\times 789.$
\item [] $7808=1\times 2+3\times 4^5+6\times 789.$
\item [] $7809=1\times 234\times 5\times 6+789.$
\item [] $7810=1+234\times 5\times 6+789.$
\item [] $\mathit{7811=1^2-3+4^5+6789.}$
\item [] $7812=12\times (3\times 4+567+8\times 9).$
\item [] $7813=1^{23}\times 4^5+6789.$
\item [] $7814=1^{23}+4^5+6789.$
\item [] $7815=(1+234\times 5)\times 6+789.$
\item [] $7816=1^2\times 3+4^5+6789.$
\item [] $7817=1^2+3+4^5+6789.$
\item [] $7818=12+3\times 4^5+6\times 789.$
\item [] $7819=1+2+3+4^5+6789.$
\item [] $7820=1+2\times 3+4^5+6789.$
\item [] $7821=1\times 2^3+4^5+6789.$
\item [] $7822=1+2^3+4^5+6789.$
\item [] $7823=1\times 23+4\times 5\times 6\times (7\times 8+9).$
\item [] $7824=1\times 23\times 45+6789.$
\item [] $7825=1+23\times 45+6789.$
\item [] $\mathit{7826=1\times 2\times (3^4\times 56-7\times 89).}$
\item [] $7827=(1+2)^3+4\times 5\times 6\times (7\times 8+9).$
\item [] $7828=12+3+4^5+6789.$
\item [] $7829=1\times 23\times 4\times ((5+6)\times 7+8)+9.$
\item [] $7830=(12+3)\times (45+6\times 78+9).$
\item[]$\mbox{Decreasing order}$
\item [] $7761=98+7654+3^2\times 1.$
\item [] $7762=98+7654+3^2+1.$
\item [] $7763=98\times (7+65+4+3)+21.$
\item [] $7764=(98\times (7+6)+5\times 4)\times (3+2+1).$
\item [] $7765=(98\times (7+6)+5\times 4)\times 3\times 2+1.$
\item [] $7766=(9+8)\times 7\times 65+4+3^{(2+1)}.$
\item [] $7767=9\times (87\times 6+5\times 4+321).$
\item [] $7768=98\times 76+5\times (43+21).$
\item [] $7769=98\times (7+65+4)+321.$
\item [] $7770=98\times 76+5\times 4^3+2\times 1.$
\item [] $7771=98\times 76+5\times 4^3+2+1.$
\item [] $7772=98\times 76+54\times 3\times 2\times 1.$
\item [] $7773=98\times 76+54\times 3\times 2+1.$
\item [] $7774=9+(8\times 7+65)\times 4^3+21.$
\item [] $7775=(9+8)\times 7\times 65+4\times (3^2+1).$
\item [] $7776=98+7654+3+21.$
\item [] $7777=987\times 6+5+43^2+1.$
\item [] $7778=98\times 76+5+4+321.$
\item [] $7779=98+7654+3^{(2+1)}.$
\item [] $7780=(9+8)\times 7\times 65+43+2\times 1.$
\item [] $7781=(9+8)\times 7\times 65+43+2+1.$
\item [] $7782=(9+8+7)\times 6\times 54+3+2+1.$
\item [] $7783=98\times 76+5\times (4+3\times 21).$
\item [] $7784=98+7654+32\times 1.$
\item [] $7785=98+7654+32+1.$
\item [] $7786=(9+8+7)\times 6\times 54+3^2+1.$
\item [] $7787=9+8+7\times 6\times 5\times (4+32+1).$
\item [] $7788=9+8\times (76+5)\times 4\times 3+2+1.$
\item [] $7789=98\times 76+5\times 4+321.$
\item [] $\mathit{7790=9-8-7+6^5-4+3+21.}$
\item [] $7791=(9\times 8+(7+6)\times (5\times 4+3))\times 21.$
\item [] $\mathit{7792=9\times 876-5-43\times 2-1.}$
\item [] $7793=9+8+(7+6+5)\times 432\times 1.$
\item [] $7794=9+8+(7+6+5)\times 432+1.$
\item [] $7795=(9+8\times (7\times (65+4)+3))\times 2+1.$
\item [] $7796=98\times 76+5+(4+3)^{(2+1)}.$
\item [] $7797=(9+8\times (7+6))\times (5+43+21).$
\item [] $7798=98\times 7\times (6+5)+4\times 3\times 21.$
\item [] $7799=987\times 6+5^4\times 3+2\times 1.$
\item [] $7800=987\times 6+5^4\times 3+2+1.$
\item [] $7801=(9+8)\times 7\times 65+4^3+2\times 1.$
\item [] $7802=(9+8)\times 7\times 65+4+3\times 21.$
\item [] $7803=(98+765+4)\times 3^2\times 1.$
\item [] $7804=(98+765+4)\times 3^2+1.$
\item [] $7805=98\times 76+(5+4\times 3)\times 21.$
\item [] $7806=9+8\times (76+5)\times 4\times 3+21.$
\item [] $7807=(9+8\times 7)\times 6\times 5\times 4+3\times 2+1.$
\item [] $7808=(9+8+7)\times 6\times 54+32\times 1.$
\item [] $7809=9+8+7+6^5+4+3+2\times 1.$
\item [] $7810=9+8+7+6^5+4+3+2+1.$
\item [] $7811=9+8+7+6^5+4+3\times 2+1.$
\item [] $7812=(98+7\times 6\times 5+4^3)\times 21.$
\item [] $7813=9+8+7+6^5+4+3^2\times 1.$
\item [] $7814=9+8+7+6^5+4+3^2+1.$
\item [] $7815=98+7654+3\times 21.$
\item [] $7816=9+((8+7)\times 65\times 4+3)\times 2+1.$
\item [] $7817=9\times 8+(76+5+4+3)^2+1.$
\item [] $7818=987\times 6+5^4\times 3+21.$
\item [] $\mathit{7819=9+8\times 7+6^5-43+21.}$
\item [] $7820=(9+8)\times 7\times 65+4^3+21.$
\item [] $7821=(9+8)\times 7\times 65+43\times 2\times 1.$
\item [] $7822=(9+8)\times 7\times 65+43\times 2+1.$
\item [] $7823=98\times 76+54+321.$
\item [] $7824=9+8+7+6^5+4\times 3\times 2\times 1.$
\item [] $7825=9+8+7+6^5+4\times 3\times 2+1.$
\item [] $7826=98\times 76+54\times (3\times 2+1).$
\item [] $7827=(9+8)\times (7\times 65+4)+3+21.$
\item [] $7828=9+8+7+6^5+4+3+21.$
\item [] $\mathit{7829=9\times 876-54-3+2\times 1.}$
\item [] $7830=9+87+6\times (5+4\times 321).$
\item[]$\mbox{Increasing order}$
\item [] $7831=1^2+(34+56)\times (78+9).$
\item [] $7832=(1+2\times 3\times 4+56+7)\times 89.$
\item [] $7833=1+2+(34+56)\times (78+9).$
\item [] $7834=(1^2+3)\times 4^5+6\times 7\times 89.$
\item [] $7835=1+2^{(3+4+5)}+6\times 7\times 89.$
\item [] $7836=1\times 23+4^5+6789.$
\item [] $7837=1+23+4^5+6789.$
\item [] $7838=1\times 2+3\times 4\times (5\times 6+7\times 89).$
\item [] $7839=12^3+4+5+678\times 9.$
\item [] $7840=(1+2)^3+4^5+6789.$
\item [] $7841=1+2\times (3+4)\times (56+7\times 8\times 9).$
\item [] $7842=12+(34+56)\times (78+9).$
\item [] $7843=(1+2\times 3+4)\times (5+6+78\times 9).$
\item [] $7844=1+23\times (4\times (56+7)+89).$
\item [] $7845=(1\times 2+3)\times (4^5+67\times 8+9).$
\item [] $7846=1+(2+3)\times (4^5+67\times 8+9).$
\item [] $7847=(1+2^{(3+4)})\times 56+7\times 89.$
\item [] $7848=12\times (3^4+567)+8\times 9.$
\item [] $7849=12\times 3+4^5+6789.$
\item [] $7850=12^3+4\times 5+678\times 9.$
\item [] $7851=1+2^{(3+4)}+(5+6)\times 78\times 9.$
\item [] $7852=1^2+3+4\times (5\times 6\times 7+8)\times 9.$
\item [] $7853=12^3\times 4+5+(6+7)\times 8\times 9.$
\item [] $7854=(123+4+5\times 67)\times (8+9).$
\item [] $7855=(1+2\times 3)\times 4^5+678+9.$
\item [] $7856=1\times 2^3+4\times (5\times 6\times 7+8)\times 9.$
\item [] $7857=1\times (2+3)\times 4\times 56\times 7+8+9.$
\item [] $7858=1+(2+3)\times 4\times 56\times 7+8+9.$
\item [] $7859=(1+2+34)\times 5\times 6\times 7+89.$
\item [] $7860=1+(2+3)^4\times 5+6\times 789.$
\item [] $7861=1^2+3\times 4\times 5\times (6\times 7+89).$
\item [] $7862=1\times 2+3\times 4\times 5\times (6\times 7+89).$
\item [] $7863=12+3+4\times (5\times 6\times 7+8)\times 9.$
\item [] $7864=(1+(2+3)^4)\times 5+6\times 789.$
\item [] $7865=12\times (3^4+567)+89.$
\item [] $7866=123+(45+6\times 7)\times 89.$
\item [] $7867=1^2+(34+56\times (7+8))\times 9.$
\item [] $7868=1\times 2+(34+56\times (7+8))\times 9.$
\item [] $7869=(1+23)\times 45+6789.$
\item [] $7870=1+23\times (4+5\times 67)+8\times 9.$
\item [] $7871=1+2^{(3+4)}\times 56+78\times 9.$
\item [] $7872=12+3\times 4\times 5\times (6\times 7+89).$
\item [] $7873=1^2+3\times 4\times (567+89).$
\item [] $7874=1\times 2+3\times 4\times (567+89).$
\item [] $7875=12^3+45+678\times 9.$
\item [] $7876=1\times 2+3+(456+7)\times (8+9).$
\item [] $7877=1\times 2\times 3+(456+7)\times (8+9).$
\item [] $7878=12\times 345+6\times 7\times 89.$
\item [] $7879=12^3+4+(5+678)\times 9.$
\item [] $7880=(123\times (4+5)+6)\times 7+89.$
\item [] $7881=12\times (3^4+567+8)+9.$
\item [] $7882=(1+2\times 3)\times (4^5+6+7+89).$
\item [] $7883=1+2\times (3+4)\times (5+(6+7\times 8)\times 9).$
\item [] $7884=12\times 3^4\times 5+6\times 7\times 8\times 9.$
\item [] $7885=1+2\times 3^4+(5+6)\times 78\times 9.$
\item [] $7886=1\times 23\times (4+5\times 67)+89.$
\item [] $7887=1+23\times (4+5\times 67)+89.$
\item [] $7888=(1^{23}+456+7)\times (8+9).$
\item [] $7889=1+2^3\times (45+6+7)\times (8+9).$
\item [] $7890=1+23\times (4\times 56+7\times (8+9)).$
\item [] $7891=(1+(2\times 3)^4+5)\times 6+7+8\times 9.$
\item [] $\mathit{7892=1-2+3\times (-4+5\times (67\times 8-9)).}$
\item [] $7893=(1+2\times 34\times 5+67\times 8)\times 9.$
\item [] $7894=1\times 23+(456+7)\times (8+9).$
\item [] $7895=1+23+(456+7)\times (8+9).$
\item [] $7896=123\times (4+5)+6789.$
\item [] $7897=1^2+3\times 4^5+67\times 8\times 9.$
\item [] $7898=1234+56\times 7\times (8+9).$
\item [] $7899=1+2+3\times 4^5+67\times 8\times 9.$
\item [] $7900=1\times (2+3)\times 4\times (5+6\times (7\times 8+9)).$
\item[]$\mbox{Decreasing order}$
\item [] $7831=(98+76)\times 5\times (4+3+2)+1.$
\item [] $7832=(9+8\times 7)\times 6\times 5\times 4+32\times 1.$
\item [] $7833=987+(6+5\times 4^3)\times 21.$
\item [] $\mathit{7834=9\times 876-5-43-2\times 1.}$
\item [] $7835=(9+8)\times (7\times 65+4)+32\times 1.$
\item [] $7836=9+8+7+6^5+4+32\times 1.$
\item [] $7837=9+8+7+6^5+4+32+1.$
\item [] $\mathit{7838=9\times 876-5-43+2\times 1.}$
\item [] $7839=98+7+6\times (5+4\times 321).$
\item [] $7840=98\times (7\times 6+5+4\times 3+21).$
\item [] $7841=9+8\times 7+6\times 54\times (3+21).$
\item [] $7842=98\times (7+6\times 5+43)+2\times 1.$
\item [] $7843=98\times (7+6\times 5+43)+2+1.$
\item [] $7844=9+8\times (7+6\times 54\times 3)+2+1.$
\item [] $7845=9+8+7+6^5+43+2\times 1.$
\item [] $7846=9+8+7+6^5+43+2+1.$
\item [] $7847=(9+8+7\times 6)\times (5+4\times 32\times 1).$
\item [] $7848=9\times 8\times (7+65+4+32+1).$
\item [] $7849=9+8\times 7\times (6+5+4\times 32+1).$
\item [] $7850=9+8\times 7+6^5+4+3+2\times 1.$
\item [] $7851=9+8\times 7+6^5+4+3+2+1.$
\item [] $7852=9+8\times 7+6^5+4+3\times 2+1.$
\item [] $7853=(9+8)\times 7+6\times (5+4\times 321).$
\item [] $7854=987+(6\times 54+3)\times 21.$
\item [] $7855=9+8\times 7+6^5+4+3^2+1.$
\item [] $7856=9+8\times 7+6^5+4\times 3+2+1.$
\item [] $7857=9+8\times (7+6\times 54\times 3+2\times 1).$
\item [] $7858=9+(8+7\times 6\times 5)\times 4\times 3^2+1.$
\item [] $\mathit{7859=9\times 876-5+4-3-21.}$
\item [] $7860=9\times (8+765)+43\times 21.$
\item [] $7861=9+8\times 7+6^5+4\times (3+2)\times 1.$
\item [] $7862=9+8\times (7+6\times 54\times 3)+21.$
\item [] $7863=(9+8\times 7)\times 6\times 5\times 4+3\times 21.$
\item [] $7864=9+8+7+6^5+43+21.$
\item [] $7865=9+8\times 7+6^5+4\times 3\times 2\times 1.$
\item [] $7866=98\times 7\times 6+5^4\times 3\times 2\times 1.$
\item [] $7867=98\times 7\times 6+5^4\times 3\times 2+1.$
\item [] $7868=9\times 8+7+6^5+4+3^2\times 1.$
\item [] $7869=9\times 8+7+6^5+4+3^2+1.$
\item [] $7870=9\times 8+7+6^5+4\times 3+2+1.$
\item [] $7871=98\times 7\times (6+5)+4+321.$
\item [] $7872=(9+87\times 6\times 5+4)\times 3+2+1.$
\item [] $7873=98\times 76+5\times (4^3+21).$
\item [] $7874=9+8\times 7+6^5+4\times 3+21.$
\item [] $7875=98+(7+6+5)\times 432+1.$
\item [] $7876=9\times 8+7+6^5+4\times (3+2)+1.$
\item [] $7877=9+8\times 7+6^5+4+32\times 1.$
\item [] $7878=98\times 76+5\times 43\times 2\times 1.$
\item [] $7879=98\times 76+5\times 43\times 2+1.$
\item [] $7880=9\times 8+7+6^5+4\times 3\times 2+1.$
\item [] $7881=9+87+6^5+4+3+2\times 1.$
\item [] $7882=9+87+6^5+4+3+2+1.$
\item [] $7883=9\times 8+7+6^5+4+3+21.$
\item [] $7884=9+(8\times 7+65+4)\times 3\times 21.$
\item [] $7885=98\times 76+5+432\times 1.$
\item [] $7886=98\times 76+5+432+1.$
\item [] $7887=9+87+6^5+4\times 3+2+1.$
\item [] $7888=9\times 8+7+6^5+4\times 3+21.$
\item [] $7889=(98\times 7+6\times 543)\times 2+1.$
\item [] $7890=98+7+6^5+4+3+2\times 1.$
\item [] $7891=9\times 8+7+6^5+4+32\times 1.$
\item [] $7892=98+7+6^5+4+3\times 2+1.$
\item [] $7893=9+87+6^5+4\times (3+2)+1.$
\item [] $7894=98+7+6^5+4+3^2\times 1.$
\item [] $7895=98+7+6^5+4+3^2+1.$
\item [] $7896=9+87+6^5+4\times 3\times 2\times 1.$
\item [] $7897=9+87+6^5+4\times 3\times 2+1.$
\item [] $7898=9\times 876+5+4+3+2\times 1.$
\item [] $7899=9\times 876+5+4+3+2+1.$
\item [] $7900=9+87+6^5+4+3+21.$
\item[]$\mbox{Increasing order}$
\item [] $7901=(123+4)\times 56+789.$
\item [] $7902=(12+3\times 45\times 6+7\times 8)\times 9.$
\item [] $7903=1+((2\times 3)^4+5)\times 6+7+89.$
\item [] $\mathit{7904=1-2+(3\times 45+6)\times 7\times 8+9.}$
\item [] $7905=1^2\times (3\times 45+6)\times 7\times 8+9.$
\item [] $7906=1^2+(3\times 45+6)\times 7\times 8+9.$
\item [] $7907=12\times 3+(456+7)\times (8+9).$
\item [] $7908=12+3\times 4^5+67\times 8\times 9.$
\item [] $7909=(1+2\times 3+4)\times (5+6\times 7\times (8+9)).$
\item [] $7910=1\times 2\times (3+4)\times 5\times ((6+7)\times 8+9).$
\item [] $7911=(1+2)\times (34+5)\times 67+8\times 9.$
\item [] $7912=1\times (2+3)\times 4\times 56\times 7+8\times 9.$
\item [] $7913=1+(2+3)\times 4\times 56\times 7+8\times 9.$
\item [] $7914=(1+23\times 4)\times ((5+6)\times 7+8)+9.$
\item [] $\mathit{7915=-12^3+4+567\times (8+9).}$
\item [] $7916=(1^2+3)\times (45\times 6\times 7+89).$
\item [] $7917=12+(3\times 45+6)\times 7\times 8+9.$
\item [] $\mathit{7918=-1+2+3\times (4+5\times (67\times 8-9)).}$
\item [] $7919=(1+(2+3)\times 4\times 56)\times 7+8\times 9.$
\item [] $7920=(1\times 23+45+6\times 7)\times 8\times 9.$
\item [] $7921=1\times 23\times 4\times (5\times 6+7\times 8)+9.$
\item [] $7922=1+23\times 4\times (5\times 6+7\times 8)+9.$
\item [] $7923=123+4\times 5\times 6\times (7\times 8+9).$
\item [] $7924=(1+2\times 3)\times 4^5+(6+78)\times 9.$
\item [] $7925=1\times (2+3)\times (4\times 56\times 7+8+9).$
\item [] $7926=(1+2^{(3+4)})\times 56+78\times 9.$
\item [] $\mathit{7927=1-2\times (3\times 4-5\times (6+789)).}$
\item [] $7928=(1+2)\times (34+5)\times 67+89.$
\item [] $7929=(1\times 2+3)\times 4\times 56\times 7+89.$
\item [] $7930=1+(2+3)\times 4\times 56\times 7+89.$
\item [] $7931=1+(23\times 4+5\times 6)\times (7\times 8+9).$
\item [] $7932=(1+2)\times (34+5\times 6\times (78+9)).$
\item [] $7933=12+(3\times 4+(5+6)\times 7)\times 89.$
\item [] $7934=12+(3+456+7)\times (8+9).$
\item [] $7935=(1+2+3\times 4)\times (5\times (6+7)\times 8+9).$
\item [] $7936=123+4^5+6789.$
\item [] $7937=12^3+4^5\times 6+7\times 8+9.$
\item [] $7938=(1234+5)\times 6+7\times 8\times 9.$
\item [] $7939=1^2+3^4\times (5+6+78+9).$
\item [] $7940=1\times 2+3^4\times (5+6+78+9).$
\item [] $7941=1+2+3^4\times (5+6+78+9).$
\item [] $7942=1+((2\times 3)^4+5)\times 6+(7+8)\times 9.$
\item [] $7943=(1+2)\times 34\times (5+6)\times 7+89.$
\item [] $7944=1\times 23\times (4+5+6\times 7\times 8)+9.$
\item [] $7945=1+23\times (4+5+6\times 7\times 8)+9.$
\item [] $\mathit{7946=1-2+3\times (4+5\times 6)\times 78-9.}$
\item [] $7947=(1\times 2+345+67\times 8)\times 9.$
\item [] $7948=1+(2+345+67\times 8)\times 9.$
\item [] $7949=(1^{23}+4)^5+67\times 8\times 9.$
\item [] $7950=12+3^4\times (5+6+78+9).$
\item [] $7951=12^3+4^5\times 6+7+8\times 9.$
\item [] $\mathit{7952=1+2\times 3\times 4\times 5\times 67-89.}$
\item [] $\mathit{7953=-12+3\times (4+5\times 6)\times 78+9.}$
\item [] $7954=(1+(2+3)^4)\times 5+67\times 8\times 9.$
\item [] $7955=1\times 2\times (3+4)\times 567+8+9.$
\item [] $7956=12\times (3+4+567+89).$
\item [] $7957=1+234+(5+6)\times 78\times 9.$
\item [] $7958=1+2^{(3+4)}\times 56+789.$
\item [] $7959=1\times 234\times 5+6789.$
\item [] $7960=1+234\times 5+6789.$
\item [] $7961=1+(2\times 3)^4+56\times 7\times (8+9).$
\item [] $\mathit{7962=-1-2+3\times (4+5\times 6)\times 78+9.}$
\item [] $7963=(1+2\times 3)\times 4^5+6+789.$
\item [] $7964=(1+234)\times 5+6789.$
\item [] $7965=(12+34+56)\times 78+9.$
\item [] $7966=1+(2^3+45+6)\times (7+8)\times 9.$
\item [] $7967=1\times 2+(3\times 4+5)\times 6\times 78+9.$
\item [] $7968=12^3+4^5\times 6+7+89.$
\item [] $7969=1+(23+4+56)\times (7+89).$
\item [] $\mathit{7970=-1+2^3\times 4^5-(6+7)\times (8+9).}$
\item[]$\mbox{Decreasing order}$
\item [] $7901=9\times 8+7+6^5+43+2+1.$
\item [] $7902=9\times 876+5+4+3^2\times 1.$
\item [] $7903=9\times 876+5+4+3^2+1.$
\item [] $7904=9\times 876+5+4\times 3+2+1.$
\item [] $7905=9+87+6^5+4\times 3+21.$
\item [] $7906=98+7+6^5+4\times 3\times 2+1.$
\item [] $7907=9+8\times 7+6^5+4^3+2\times 1.$
\item [] $7908=9+8\times 7+6^5+4+3\times 21.$
\item [] $7909=9+87+6^5+4+32+1.$
\item [] $7910=9\times 876+5\times 4+3+2+1.$
\item [] $7911=9\times 876+5\times 4+3\times 2+1.$
\item [] $7912=9+87+6^5+4\times (3^2+1).$
\item [] $7913=9\times 876+5+4\times 3\times 2\times 1.$
\item [] $7914=9\times 876+5+4\times 3\times 2+1.$
\item [] $7915=((9+8)\times 7+65)\times 43+2+1.$
\item [] $7916=(9+8)\times 7+6^5+4\times (3+2)+1.$
\item [] $7917=98+7+6^5+4+32\times 1.$
\item [] $7918=9+87+6^5+43+2+1.$
\item [] $7919=9\times 8+7+6^5+43+21.$
\item [] $7920=9\times (8+7+65)\times (4+3\times 2+1).$
\item [] $7921=98+7+6^5+4\times (3^2+1).$
\item [] $7922=9\times 876+5+4\times 3+21.$
\item [] $7923=(9+8)\times 7+6^5+4+3+21.$
\item [] $7924=9\times (8+7)+6^5+4+3^2\times 1.$
\item [] $7925=9\times 876+5+4+32\times 1.$
\item [] $7926=9\times 876+5+4+32+1.$
\item [] $7927=98+7+6^5+43+2+1.$
\item [] $7928=9\times 876+5\times 4+3+21.$
\item [] $7929=9+8+7+6^5+4\times 32+1.$
\item [] $7930=98+7+6^5+(4+3)^2\times 1.$
\item [] $7931=98+7+6^5+(4+3)^2+1.$
\item [] $7932=9\times 876+(5+4)\times 3+21.$
\item [] $7933=((9+8)\times 7+65)\times 43+21.$
\item [] $7934=9\times 876+5+43+2\times 1.$
\item [] $7935=98\times 76+54\times 3^2+1.$
\item [] $7936=9+87+6^5+43+21.$
\item [] $7937=9\times 876+5\times 4+32+1.$
\item [] $7938=9+87+6^5+4^3+2\times 1.$
\item [] $7939=9+87+6^5+4+3\times 21.$
\item [] $7940=9\times 8+7+6^5+4^3+21.$
\item [] $7941=9\times 8+7+6^5+43\times 2\times 1.$
\item [] $7942=9\times 8+7+6^5+43\times 2+1.$
\item [] $7943=9\times 876+54+3+2\times 1.$
\item [] $7944=9\times 876+54+3+2+1.$
\item [] $7945=9\times 876+54+3\times 2+1.$
\item [] $7946=9\times 876+5\times 4\times 3+2\times 1.$
\item [] $7947=9\times 876+5\times 4\times 3+2+1.$
\item [] $7948=98+7+6^5+4+3\times 21.$
\item [] $7949=98\times (76+5)+4+3\times 2+1.$
\item [] $7950=9\times 876+5\times (4+3^2)+1.$
\item [] $7951=98\times (76+5)+4+3^2\times 1.$
\item [] $7952=98\times (76+5)+4\times 3+2\times 1.$
\item [] $7953=9\times 876+5+43+21.$
\item [] $7954=9\times 876+5\times (4+3)\times 2\times 1.$
\item [] $7955=9\times 876+5+4^3+2\times 1.$
\item [] $7956=9\times 876+5+4+3\times 21.$
\item [] $7957=9+87+6^5+4^3+21.$
\item [] $7958=9+87+6^5+43\times 2\times 1.$
\item [] $7959=(9+8)\times 7+6^5+43+21.$
\item [] $7960=(9+8)\times (7\times 65+4\times 3)+21.$
\item [] $7961=(9+8)\times 7+6^5+4^3+2\times 1.$
\item [] $7962=9\times 876+54+3+21.$
\item [] $7963=98\times (76+5)+4\times 3\times 2+1.$
\item [] $\mathit{7964=9\times 876-5+43\times 2-1.}$
\item [] $7965=9\times 876+5\times 4\times 3+21.$
\item [] $7966=98+7+6^5+4^3+21.$
\item [] $7967=9\times 876+5\times 4+3\times 21.$
\item [] $7968=98+7+6^5+43\times 2+1.$
\item [] $7969=9+8\times 7+6^5+4^3\times 2\times 1.$
\item [] $7970=9+8\times 7\times 65+4321.$
\item[]$\mbox{Increasing order}$
\item [] $7971=123+4\times (5\times 6\times 7+8)\times 9.$
\item [] $7972=1+2\times (3+(45+6)\times 78)+9.$
\item [] $7973=(1+2+34+5\times 6)\times 7\times (8+9).$
\item [] $7974=1\times 2\times (3\times 4+5\times (6+789)).$
\item [] $7975=1+2\times (3\times 4+5\times (6+789)).$
\item [] $\mathit{7976=1+2^{(3+4)}\times (56+7)-89.}$
\item [] $7977=12+3\times 45\times (6\times 7+8+9).$
\item [] $7978=(1^2+3)^4+(5+6)\times 78\times 9.$
\item [] $7979=((12+3+4)\times 5+6)\times (7+8\times 9).$
\item [] $7980=12\times (3^4+567+8+9).$
\item [] $7981=1+2\times (3+(45+6)\times 78+9).$
\item [] $\mathit{7982=12\times 3\times 4\times 56+7-89.}$
\item [] $7983=12^3+45\times (67+8\times 9).$
\item [] $7984=1^2+3\times ((4+5\times 6)\times 78+9).$
\item [] $7985=((1+2^3\times 4)\times 5\times 6+7)\times 8+9.$
\item [] $7986=1+2+3\times ((4+5\times 6)\times 78+9).$
\item [] $\mathit{7987=(123+4^5-6)\times 7\times (-8+9).}$
\item [] $\mathit{7988=(123+4^5-6)\times 7-8+9.}$
\item [] $7989=(12+3+4)\times 5\times (6+78)+9.$
\item [] $7990=(1+2\times 3+456+7)\times (8+9).$
\item [] $7991=12^3+4^5\times 6+7\times (8+9).$
\item [] $7992=12\times (3+4+5+6+7\times 8)\times 9.$
\item [] $7993=1+2\times 3\times (4+(5+6+7)\times 8)\times 9.$
\item [] $7994=123+(456+7)\times (8+9).$
\item [] $7995=(12+3^4+5\times 6)\times (7\times 8+9).$
\item [] $7996=1^2+(3+4\times 5\times 6)\times (7\times 8+9).$
\item [] $7997=1\times 2+(3+4\times 5\times 6)\times (7\times 8+9).$
\item [] $7998=1^2\times (3^4+5)\times (6+78+9).$
\item [] $7999=1^2+(3^4+5)\times (6+78+9).$
\item [] $8000=(1+2+3+4)\times (5+6+789).$
\item [] $8001=12\times 3\times (4+5\times 6\times 7+8)+9.$
\item [] $8002=1+(23+4)\times (5\times 6+7)\times 8+9.$
\item [] $8003=12\times 34\times 5+67\times 89.$
\item [] $8004=(1+2)^3\times 45+6789.$
\item [] $8005=1+(2+34+56)\times (78+9).$
\item [] $8006=1\times 2+(3^4+5+6)\times (78+9).$
\item [] $8007=12^3+4^5\times 6+(7+8)\times 9.$
\item [] $8008=1^2\times (3+4)\times (5+67\times (8+9)).$
\item [] $8009=1\times 2^3\times 4\times 5\times (6\times 7+8)+9.$
\item [] $8010=1\times 2\times (3+4)\times 567+8\times 9.$
\item [] $8011=1+2\times (3+4)\times 567+8\times 9.$
\item [] $8012=1\times 2+(3+4\times 5+67)\times 89.$
\item [] $8013=1\times 23\times (45\times 6+78)+9.$
\item [] $8014=1+23\times (45\times 6+78)+9.$
\item [] $\mathit{8015=(1-23+4^5-6+7)\times 8-9.}$
\item [] $8016=12\times (3\times 4+567+89).$
\item [] $8017=1\times 2\times (3+4^5)+67\times 89.$
\item [] $8018=1+2\times (3+4^5)+67\times 89.$
\item [] $8019=(1+2+3\times 45\times 6+78)\times 9.$
\item [] $8020=12+(3+4)\times (5+67\times (8+9)).$
\item [] $\mathit{8021=1\times 2+3^4\times (5\times 6+78-9).}$
\item [] $8022=12+(3+4\times 5+67)\times 89.$
\item [] $8023=((1+2)^3\times 4+5)\times (6+7\times 8+9).$
\item [] $8024=(1+2^3+456+7)\times (8+9).$
\item [] $8025=(12+3)\times (456+7+8\times 9).$
\item [] $8026=1\times 2+(3^4\times 5+67)\times (8+9).$
\item [] $8027=1\times 2\times (3+4)\times 567+89.$
\item [] $8028=1234+5+6789.$
\item [] $8029=((1+2)^3+4)\times (5\times (6\times 7+8)+9).$
\item [] $8030=(1\times 2+3)\times (4+(5+6+7)\times 89).$
\item [] $8031=1+2+(34+(5+6)\times 78)\times 9.$
\item [] $\mathit{8032=1+(2\times 3+4+5)\times 67\times 8-9.}$
\item [] $\mathit{8033=((-1+23)\times 45+6+7)\times 8+9.}$
\item [] $8034=((1+2)\times 34+5)\times (67+8)+9.$
\item [] $\mathit{8035=-1+2^3\times 4^5-67-89.}$
\item [] $8036=(1^2+3^4)\times (5+6+78+9).$
\item [] $8037=(12+345+67\times 8)\times 9.$
\item [] $\mathit{8038=-1+2\times 3\times 4\times 5\times 67+8-9.}$
\item [] $\mathit{8039=1\times 2\times 3\times 4\times 5\times 67+8-9.}$
\item [] $8040=(1+23)\times (45\times 6+7\times 8+9).$
\item[]$\mbox{Decreasing order}$
\item [] $7971=9\times 876+54+32+1.$
\item [] $7972=9\times (876+5+4)+3\times 2+1.$
\item [] $7973=9+8\times 7+6^5+4\times (32+1).$
\item [] $7974=9\times 876+5+4^3+21.$
\item [] $7975=9\times 876+5+43\times 2\times 1.$
\item [] $7976=9\times 876+5+43\times 2+1.$
\item [] $7977=9\times (8+7)+6^5+4^3+2\times 1.$
\item [] $7978=9+8\times 765+43^2\times 1.$
\item [] $7979=9+8\times 765+43^2+1.$
\item [] $7980=9\times 876+(5+43)\times 2\times 1.$
\item [] $7981=(9+8)\times 7+6^5+43\times 2\times 1.$
\item [] $7982=(9+8)\times 7+6^5+43\times 2+1.$
\item [] $7983=9\times 8+7+6^5+4\times 32\times 1.$
\item [] $7984=9\times 8+7+6^5+4^3\times 2+1.$
\item [] $7985=9\times 876+5+4\times (3+21).$
\item [] $7986=9+8\times 7+6^5+(4\times 3)^2+1.$
\item [] $7987=(9+8)\times 7\times 65+4\times 3\times 21.$
\item [] $7988=98\times 76+54\times (3^2+1).$
\item [] $7989=9\times (876+5+4)+3+21.$
\item [] $7990=(9+8)\times (7\times 65+4\times 3+2+1).$
\item [] $7991=(9+8)\times 7+6^5+4\times (3+21).$
\item [] $7992=9+8+7654+321.$
\item [] $7993=98\times 76+543+2\times 1.$
\item [] $7994=98\times 76+543+2+1.$
\item [] $7995=(9+8\times 7)\times (6+54+3\times 21).$
\item [] $7996=9\times (876+5)+4+3\times 21.$
\item [] $7997=9\times (8+7)+6^5+43\times 2\times 1.$
\item [] $7998=9\times 876+(54+3)\times 2\times 1.$
\item [] $7999=9\times 876+(54+3)\times 2+1.$
\item [] $8000=9+87+6^5+4^3\times 2\times 1.$
\item [] $8001=9\times 876+54+3\times 21.$
\item [] $8002=98\times (76+5)+43+21.$
\item [] $8003=(9+8)\times 7+6^5+4\times 3^{(2+1)}.$
\item [] $8004=9\times 876+5\times 4\times 3\times 2\times 1.$
\item [] $8005=9\times 876+5\times 4\times 3\times 2+1.$
\item [] $8006=(9\times 8\times (7+6\times 5)+4)\times 3+2\times 1.$
\item [] $8007=9+(87+6)\times (54+32\times 1).$
\item [] $8008=9+(87+6)\times (54+32)+1.$
\item [] $8009=98+7+6^5+4^3\times 2\times 1.$
\item [] $8010=98+7+6^5+4^3\times 2+1.$
\item [] $8011=(9+876+5)\times (4+3+2)+1.$
\item [] $8012=98\times 76+543+21.$
\item [] $8013=9+87\times (6+54+32\times 1).$
\item [] $8014=9+87\times (6+54+32)+1.$
\item [] $8015=98\times 76+(5+4)\times 3\times 21.$
\item [] $8016=9\times (876+5)+43\times 2+1.$
\item [] $8017=9\times 876+5+4^3\times 2\times 1.$
\item [] $8018=9\times 876+5+4^3\times 2+1.$
\item [] $8019=9+87+6^5+(4+3)\times 21.$
\item [] $8020=(9\times (87+6)+54)\times 3^2+1.$
\item [] $8021=9\times 876+5+4\times (32+1).$
\item [] $8022=9\times 876+(5+4^3)\times 2\times 1.$
\item [] $8023=(9+8)\times 7+6^5+4^3\times 2\times 1.$
\item [] $8024=(9\times 8+76)\times 54+32\times 1.$
\item [] $8025=98\times (76+5)+43\times 2+1.$
\item [] $8026=98+7+6^5+(4\times 3)^2+1.$
\item [] $8027=(9+8)\times 7+6^5+4\times (32+1).$
\item [] $8028=98+7+6^5+(4+3)\times 21.$
\item [] $8029=9\times 876+(5+4+3)^2+1.$
\item [] $8030=((9+8)\times (7\times 6+5)+4)\times (3^2+1).$
\item [] $8031=9+(87\times 6\times 5+4^3)\times (2+1).$
\item [] $8032=9+8+7+6+(5\times 4)^3+2\times 1.$
\item [] $8033=9\times 876+5+(4\times 3)^2\times 1.$
\item [] $8034=9\times 876+5+(4\times 3)^2+1.$
\item [] $8035=((9+8)\times 7\times 6+5^4)\times 3\times 2+1.$
\item [] $8036=9\times 876+5+(4+3)\times 21.$
\item [] $8037=98\times 7+6\times (5\times (4+3))^2+1.$
\item [] $8038=98\times (7+6+5+4^3)+2\times 1.$
\item [] $8039=9\times (8+7)+6^5+4\times 32\times 1.$
\item [] $8040=(98\times 7+654)\times 3\times 2\times 1.$
\item[]$\mbox{Increasing order}$
\item [] $8041=(1\times 23\times 4\times 5+6+7)\times (8+9).$
\item [] $8042=1^2+(3^4+56\times 7)\times (8+9).$
\item [] $8043=(1\times 23\times 4+5+6)\times 78+9.$
\item [] $8044=1+(23\times 4+5+6)\times 78+9.$
\item [] $8045=(1\times 2+3)\times (4\times (56\times 7+8)+9).$
\item [] $8046=(1\times 2\times 3^4\times 5+6+78)\times 9.$
\item [] $8047=1+(2\times 3^4\times 5+6+78)\times 9.$
\item [] $\mathit{8048=1\times 2^3\times 4^5-6\times (7+8+9).}$
\item [] $8049=(1+2+3+4+5)\times 67\times 8+9.$
\item [] $8050=1+(2\times 3+4+5)\times 67\times 8+9.$
\item [] $\mathit{8051=1\times 2^3\times 4^5-6-(7+8)\times 9.}$
\item [] $8052=123\times 4+56\times (7+8)\times 9.$
\item [] $8053=12+(((3^4)+(56\times 7))\times (8+9)).$
\item [] $\mathit{8054=12\times 3\times 4\times 56+7-8-9.}$
\item [] $8055=(1+2\times 3^4\times 5+6+78)\times 9.$
\item [] $8056=12^3\times 4+5+67\times (8+9).$
\item [] $8057=1\times 2\times 3\times 4\times 5\times 67+8+9.$
\item [] $8058=1+2\times 3\times 4\times 5\times 67+8+9.$
\item [] $8059=1\times 2\times (3+4)\times (567+8)+9.$
\item [] $8060=1+2\times (3+4)\times (567+8)+9.$
\item [] $8061=12\times (3\times 45+67\times 8)+9.$
\item [] $8062=1^2+3\times (4+5\times 67\times 8)+9.$
\item [] $8063=1\times 2+3\times (4+5\times 67\times 8)+9.$
\item [] $8064=1^{23}\times (45+67)\times 8\times 9.$
\item [] $8065=1^{23}+(45+67)\times 8\times 9.$
\item [] $8066=(1+2\times (3\times 4\times 5\times 67+8))+9.$
\item [] $8067=1^2\times 3+(45+67)\times 8\times 9.$
\item [] $8068=1^2+3+(45+67)\times 8\times 9.$
\item [] $8069=(1^2+3)^4\times 5+6789.$
\item [] $8070=1\times 2\times 3+(45+67)\times 8\times 9.$
\item [] $8071=1+2\times 3+(45+67)\times 8\times 9.$
\item [] $8072=1\times 2^3+(45+67)\times 8\times 9.$
\item [] $8073=12+3\times (4+5\times 67\times 8)+9.$
\item [] $8074=12^3\times 4+5+(6+7)\times 89.$
\item [] $8075=(1\times 2+3^4+56\times 7)\times (8+9).$
\item [] $8076=12\times (34+567+8\times 9).$
\item [] $\mathit{8077=-1+2^3\times 4^5-6\times 7-8\times 9.}$
\item [] $\mathit{8078=1\times 2^3\times 4^5-6\times 7-8\times 9.}$
\item [] $8079=(1+2)^3\times 45\times 6+789.$
\item [] $8080=1^2+3\times (4+5\times 67\times 8+9).$
\item [] $8081=1\times 2+3\times (4+5\times 67\times 8+9).$
\item [] $8082=123\times (4+56)+78\times 9.$
\item [] $8083=1\times (2+3\times 45)\times (6\times 7+8+9).$
\item [] $8084=1+(2+3\times 45)\times (6\times 7+8+9).$
\item [] $8085=(1+2\times 3)\times (4^5+6\times 7+89).$
\item [] $\mathit{8086=-1\times 2+3\times 4\times (5+678-9).}$
\item [] $8087=1234+(5+6)\times 7\times 89.$
\item [] $8088=12\times 3\times 4\times 56+7+8+9.$
\item [] $8089=1+2\times 3\times (4+56\times (7+8+9)).$
\item [] $8090=1\times (2\times 3)^4+5+6789.$
\item [] $8091=1+(2\times 3)^4+5+6789.$
\item [] $8092=1\times 2\times 34\times (5+6\times 7+8\times 9).$
\item [] $8093=1+2\times 34\times (5+6\times 7+8\times 9).$
\item [] $8094=1\times 2+(3\times 4+56)\times 7\times (8+9).$
\item [] $8095=1+2+(3\times 4+56)\times 7\times (8+9).$
\item [] $8096=(1+2\times 3)^4+5\times 67\times (8+9).$
\item [] $8097=1+2^3\times (4+(56+7\times 8)\times 9).$
\item [] $8098=1\times 2\times (3+(4+5\times 6)\times 7\times (8+9)).$
\item [] $8099=(12+3+4+5+67)\times 89.$
\item [] $8100=1\times 2\times 3\times 45\times (6+7+8+9).$
\item [] $8101=1+2\times 3\times 45\times (6+7+8+9).$
\item [] $8102=1\times 2\times (3+4^5+6\times 7\times 8\times 9).$
\item [] $8103=1+2\times (3+4^5+6\times 7\times 8\times 9).$
\item [] $8104=1\times 2^3\times (4\times 56+789).$
\item [] $8105=1+2^3\times (4\times 56+789).$
\item [] $8106=(1+2\times 3)\times (456+78\times 9).$
\item [] $8107=1+2\times 3+(4+56)\times (7+8)\times 9.$
\item [] $8108=1\times 2^3+(4+56)\times (7+8)\times 9.$
\item [] $8109=(1\times 2^3+4+5)\times (6\times 78+9).$
\item [] $8110=1+(2^3+4+5)\times (6\times 78+9).$
\item[]$\mbox{Decreasing order}$
\item [] $8041=(98\times 7+654)\times 3\times 2+1.$
\item [] $8042=98+(7+6\times 54)\times (3+21).$
\item [] $8043=9\times (8+7)+6^5+4\times (32+1).$
\item [] $\mathit{8044=9\times 876+54\times 3-2\times 1.}$
\item [] $\mathit{8045=9\times 876+54\times 3-2+1.}$
\item [] $8046=(9+876+5+4)\times 3^2\times 1.$
\item [] $8047=9\times 8+7654+321.$
\item [] $8048=9\times 876+54\times 3+2\times 1.$
\item [] $8049=9\times 876+54\times 3+2+1.$
\item [] $8050=(9+87+65)\times ((4+3)^2+1).$
\item [] $8051=9+8+7+6+(5\times 4)^3+21.$
\item [] $8052=9+8+7+6^5+4\times 3\times 21.$
\item [] $\mathit{8053=98-7\times 6+(5\times 4)^3-2-1.}$
\item [] $\mathit{8054=9-8\times 7+6^5+4+321.}$
\item [] $8055=(9\times 8+76)\times 54+3\times 21.$
\item [] $8056=(98+76+5)\times (43+2)+1.$
\item [] $8057=98\times (7+6+5+4^3)+21.$
\item [] $8058=9\times (876+5)+4\times 32+1.$
\item [] $8059=(9\times 8+7)\times (6+(5+43)\times 2)+1.$
\item [] $8060=(9+8)\times 7\times 65+4+321.$
\item [] $8061=9\times (876+5)+4\times (32+1).$
\item [] $8062=9+8+7\times 6+(5\times 4)^3+2+1.$
\item [] $\mathit{8063=9\times 876+5\times 4\times 3^2-1.}$
\item [] $8064=9\times 8\times 7+6\times 5\times 4\times 3\times 21.$
\item [] $8065=9\times 876+5\times 4\times 3^2+1.$
\item [] $8066=98\times (76+5)+4\times 32\times 1.$
\item [] $8067=9\times 876+54\times 3+21.$
\item [] $\mathit{8068=98\times 76+5^4-3-2\times 1.}$
\item [] $8069=9\times 876+5\times (4+32+1).$
\item [] $8070=9\times (876+5\times 4)+3\times 2\times 1.$
\item [] $8071=9\times (876+5\times 4)+3\times 2+1.$
\item [] $\mathit{8072=98\times 76+5^4-3+2\times 1.}$
\item [] $8073=98+7654+321.$
\item [] $8074=9+8\times 7+6+(5\times 4)^3+2+1.$
\item [] $8075=(9+8\times 7+6\times 5)\times (4^3+21).$
\item [] $8076=9\times (876+5)+(4+3)\times 21.$
\item [] $8077=((9+8)\times 76+54)\times 3\times 2+1.$
\item [] $8078=98\times 76+5^4+3+2\times 1.$
\item [] $8079=98\times 76+5^4+3+2+1.$
\item [] $8080=98\times 76+5^4+3\times 2+1.$
\item [] $8081=9\times 876+5+4^3\times (2+1).$
\item [] $8082=987\times 6+5\times 432\times 1.$
\item [] $8083=987\times 6+5\times 432+1.$
\item [] $8084=9\times 876+5\times 4\times (3^2+1).$
\item [] $8085=(9+87+6\times 5)\times 4^3+21.$
\item [] $8086=(98+7)\times (6+5+4^3+2)+1.$
\item [] $8087=987\times 6+5\times (432+1).$
\item [] $8088=98\times 76+5\times 4^3\times 2\times 1.$
\item [] $8089=98\times 76+5\times 4\times 32+1.$
\item [] $8090=(9+87)\times 65+43^2+1.$
\item [] $8091=(9+8+76)\times (54+32+1).$
\item [] $8092=9+8\times 7+6+(5\times 4)^3+21.$
\item [] $8093=9+8\times 7+6^5+4\times 3\times 21.$
\item [] $8094=(9+8\times 7+6)\times (54+3)\times 2\times 1.$
\item [] $8095=9+8+76+(5\times 4)^3+2\times 1.$
\item [] $8096=9\times (876+5\times 4)+32\times 1.$
\item [] $8097=98\times 76+5^4+3+21.$
\item [] $8098=98\times (7+6)\times 5+(4\times 3)^{(2+1)}.$
\item [] $8099=98+(7+6\times 5\times 4)\times 3\times 21.$
\item [] $8100=9\times 87\times 6+54\times 3\times 21.$
\item [] $8101=9\times 876+5\times 43+2\times 1.$
\item [] $8102=9\times 876+5\times 43+2+1.$
\item [] $8103=9\times (8+7)+6^5+4^3\times (2+1).$
\item [] $8104=9+87+6+(5\times 4)^3+2\times 1.$
\item [] $8105=98\times 76+5^4+32\times 1.$
\item [] $8106=98\times 76+5^4+32+1.$
\item [] $8107=9\times 8+7+6^5+4\times 3\times 21.$
\item [] $8108=98\times 76+5\times 4\times (32+1).$
\item [] $8109=(98+7)\times 65+4\times 321.$
\item [] $8110=9\times 876+5\times (43+2)+1.$
\item[]$\mbox{Increasing order}$
\item [] $8111=1\times 2+(3\times 4+5)\times (6\times 78+9).$
\item [] $8112=1\times 2\times 3\times 4\times 5\times 67+8\times 9.$
\item [] $8113=1+2\times 3\times 4\times 5\times 67+8\times 9.$
\item [] $8114=1\times 2\times (3+4+5\times 6\times (7+8)\times 9).$
\item [] $8115=12+3+(4+56)\times (7+8)\times 9.$
\item [] $8116=1\times 2\times ((3+4)\times 567+89).$
\item [] $8117=1+2\times ((3+4)\times 567+89).$
\item [] $8118=12^3\times 4+(56+78)\times 9.$
\item [] $8119=1\times 23\times (4\times (5\times 6+7\times 8)+9).$
\item [] $8120=1+23\times (4\times (5\times 6+7\times 8)+9).$
\item [] $8121=(12+3^4+5+6)\times 78+9.$
\item [] $8122=1+2\times (34+5)\times (6+7)\times 8+9.$
\item [] $8123=1\times 23+(4+56)\times (7+8)\times 9.$
\item [] $8124=12\times (3+45+6+7\times 89).$
\item [] $8125=(1+2\times 34+56)\times (7\times 8+9).$
\item [] $8126=(12+3+456+7)\times (8+9).$
\item [] $8127=123\times (45+6+7+8)+9.$
\item [] $8128=1+(23+4\times 5)\times (6+7+8)\times 9.$
\item [] $8129=1\times 2\times 3\times 4\times 5\times 67+89.$
\item [] $8130=1+2\times 3\times 4\times 5\times 67+89.$
\item [] $8131=1+2\times ((34+5)\times (6+7)\times 8+9).$
\item [] $\mathit{8132=-1+2^3\times 4^5-6\times 7-8-9.}$
\item [] $\mathit{8133=12345-6\times 78\times 9.}$
\item [] $8134=1\times (2+3^4)\times (5+6+78+9).$
\item [] $8135=1+(2+3^4)\times (5+6+78+9).$
\item [] $8136=(1234+5)\times 6+78\times 9.$
\item [] $8137=1\times 2^3\times (4\times 5\times 6+7)\times 8+9.$
\item [] $8138=1+2^3\times (4\times 5\times 6+7)\times 8+9.$
\item [] $8139=12+(3+(4+56)\times (7+8))\times 9.$
\item [] $8140=12^3+4+(5+67)\times 89.$
\item [] $8141=1+(2+3)\times 4\times (5\times 67+8\times 9).$
\item [] $8142=12\times 34\times 5+678\times 9.$
\item [] $8143=12\times 3\times 4\times 56+7+8\times 9.$
\item [] $8144=1\times (2+3)\times 4^5+6\times 7\times 8\times 9.$
\item [] $8145=(12+3)\times (456+78+9).$
\item [] $8146=1^2+(3+4+5)\times 678+9.$
\item [] $8147=1\times 2+(3+4+5)\times 678+9.$
\item [] $8148=1+2+(3+4+5)\times 678+9.$
\item [] $8149=1^2+3\times 4\times (56+7\times 89).$
\item [] $8150=1\times 2+3\times 4\times (56+7\times 89).$
\item [] $8151=12\times 3\times 4\times 56+78+9.$
\item [] $\mathit{8152=1+2^3\times 4^5-6\times 7-8+9.}$
\item [] $8153=(1+23)\times (4+5\times 67)+8+9.$
\item [] $8154=12^3+(4+5)\times 6\times 7\times (8+9).$
\item [] $8155=(1+2)^(3+4)+5+67\times 89.$
\item [] $8156=1\times 2+3\times (4\times 56+78)\times 9.$
\item [] $8157=12+(3+4+5)\times 678+9.$
\item [] $8158=1\times 2\times (3456+7\times 89).$
\item [] $8159=1+2\times (3456+7\times 89).$
\item [] $8160=(123+4^5+6)\times 7+89.$
\item [] $8161=1+2\times (3+45)\times (6+7+8\times 9).$
\item [] $8162=1\times 2+3\times 4\times (5+(67+8)\times 9).$
\item [] $8163=1+2+3\times 4\times (5+(67+8)\times 9).$
\item [] $8164=1^2+3\times ((4+5\times 67)\times 8+9).$
\item [] $8165=12\times (3\times 4\times 56+7)+8+9.$
\item [] $8166=12+3\times (4\times 56+78)\times 9.$
\item [] $8167=1+2\times (345+6\times 7\times 89).$
\item [] $8168=1\times 2\times (34+5\times 6\times (7+8)\times 9).$
\item [] $8169=123\times (4+56)+789.$
\item [] $8170=1+2\times 3\times (4^5+6\times 7\times 8)+9.$
\item [] $\mathit{8171=1+2^3\times 4^5+67-89.}$
\item [] $8172=1\times (2+34)\times (5\times 6\times 7+8+9).$
\item [] $8173=1+(2+34)\times (5\times 6\times 7+8+9).$
\item [] $\mathit{8174=1\times 2^3\times 4^5+6-7-8-9.}$
\item [] $8175=(1+2+3+4+5)\times (67\times 8+9).$
\item [] $8176=1\times 2\times (3+4)\times (567+8+9).$
\item [] $8177=1+2\times (3+4)\times (567+8+9).$
\item [] $8178=(1\times 2+3^4+5+6)\times (78+9).$
\item [] $8179=1+2\times (3+4\times (5+6))\times (78+9).$
\item [] $8180=1\times (2+3)\times 4\times (56\times 7+8+9).$
\item[]$\mbox{Decreasing order}$
\item [] $8111=9+(8+7\times 6)\times 54\times 3+2\times 1.$
\item [] $8112=9+(8+7\times 6)\times 54\times 3+2+1.$
\item [] $8113=98+7+6+(5\times 4)^3+2\times 1.$
\item [] $8114=9\times 876+5\times (43+2+1).$
\item [] $8115=(9\times (8+7)\times 6\times 5+4+3)\times 2+1.$
\item [] $8116=9+(8\times 7+65)\times (4+3\times 21).$
\item [] $8117=9\times 8+7\times 6+(5\times 4)^3+2+1.$
\item [] $8118=9\times (876+5\times 4+3+2+1).$
\item [] $8119=9\times (876+5\times 4+3\times 2)+1.$
\item [] $8120=9\times 876+5\times 43+21.$
\item [] $8121=(9+8\times 7)\times 6\times 5\times 4+321.$
\item [] $8122=(98\times 7+(6+5+4)^3)\times 2\times 1.$
\item [] $8123=9+87+6+(5\times 4)^3+21.$
\item [] $8124=9+87+6^5+4\times 3\times 21.$
\item [] $8125=9+8+7+6^5+4+321.$
\item [] $8126=(9+8\times 7)\times (6\times 5\times 4+3+2)+1.$
\item [] $8127=9\times (876+5\times 4)+3\times 21.$
\item [] $8128=(9\times 87+6\times 5\times 4)\times 3^2+1.$
\item [] $8129=9\times 876+5\times (4+3)^2\times 1.$
\item [] $8130=9+(8+7\times 6)\times 54\times 3+21.$
\item [] $8131=(9\times 87+6\times 5)\times (4+3\times 2)+1.$
\item [] $8132=98+7+6+(5\times 4)^3+21.$
\item [] $8133=98+7+6^5+4\times 3\times 21.$
\item [] $8134=98\times (7+65+4+3\times 2+1).$
\item [] $8135=9\times 8+7\times 6+(5\times 4)^3+21.$
\item [] $8136=98\times 76+5^4+3\times 21.$
\item [] $8137=9\times (8+76\times 5+4^3)\times 2+1.$
\item [] $8138=9+((8+7\times 6)\times 5+4)\times 32+1.$
\item [] $8139=9\times 87+6\times ((5\times (4+3))^2+1).$
\item [] $8140=(9\times 8+76)\times 5\times (4+3\times 2+1).$
\item [] $8141=9\times 876+5+4\times 3\times 21.$
\item [] $8142=98+7\times 6+(5\times 4)^3+2\times 1.$
\item [] $8143=98+7\times 6+(5\times 4)^3+2+1.$
\item [] $8144=9\times 8\times 7\times 6+5\times 4^(3+2)\times 1.$
\item [] $8145=9+8\times (765+4\times 3\times 21).$
\item [] $8146=(9+876+5\times 4)\times 3^2+1.$
\item [] $8147=(9+8)\times 7+6^5+4\times 3\times 21.$
\item [] $8148=9+8+(7+6)\times 5^4+3+2+1.$
\item [] $8149=9+8+(7+6)\times 5^4+3\times 2+1.$
\item [] $8150=9\times 8+76+(5\times 4)^3+2\times 1.$
\item [] $8151=9+8+(7+6)\times 5^4+3^2\times 1.$
\item [] $8152=9+8+(7+6)\times 5^4+3^2+1.$
\item [] $\mathit{8153=9\times 8\times 7+6^5-4\times 32+1.}$
\item [] $8154=9\times 876+54\times (3+2\times 1).$
\item [] $8155=9\times 876+54\times (3+2)+1.$
\item [] $\mathit{8156=98\times 76+(5+4)^3-21.}$
\item [] $8157=9+(8+76\times 5)\times (4+3)\times (2+1).$
\item [] $\mathit{8158=9\times 8-7+6^5-4+321.}$
\item [] $8159=9+8\times 7+6\times (5+4^3\times 21).$
\item [] $8160=(9+8+7\times (6\times 5+4))\times 32\times 1.$
\item [] $8161=98+7\times 6+(5\times 4)^3+21.$
\item [] $8162=9\times (8+7)+6+(5\times 4)^3+21.$
\item [] $8163=9\times (8+7)+6^5+4\times 3\times 21.$
\item [] $8164=9\times (8+(7+6)\times (5+4^3)+2)+1.$
\item [] $8165=(9+8+7)\times 6+(5\times 4)^3+21.$
\item [] $8166=9+8\times 7+6^5+4+321.$
\item [] $8167=(9+8)\times 7\times 65+432\times 1.$
\item [] $8168=(9+8)\times 7\times 65+432+1.$
\item [] $8169=(98+76+5\times 43)\times 21.$
\item [] $8170=(9+8\times 7+6\times 5)\times 43\times 2\times 1.$
\item [] $8171=(9+8\times 7+6\times 5)\times 43\times 2+1.$
\item [] $8172=9\times 876+(5+4)\times 32\times 1.$
\item [] $8173=9\times 876+(5+4)\times 32+1.$
\item [] $8174=9\times 876+(5+4\times 3)^2+1.$
\item [] $8175=9+8+(7+6)\times 5^4+32+1.$
\item [] $8176=98+76+(5\times 4)^3+2\times 1.$
\item [] $8177=98+76+(5\times 4)^3+2+1.$
\item [] $8178=9+(8+7\times 6\times (5+4)+3)\times 21.$
\item [] $8179=98\times 76+(5+4)^3+2\times 1.$
\item [] $8180=9\times 8+7+6^5+4+321.$
\item[]$\mbox{Increasing order}$
\item [] $8181=12\times (345+6\times 7\times 8)+9.$
\item [] $\mathit{8182=1\times 2-3-4^5\times (6-7)\times 8-9.}$
\item [] $8183=12\times 3\times 4\times 56+7\times (8+9).$
\item [] $8184=(1^2+3)^4\times 5\times 6+7\times 8\times 9.$
\item [] $8185=1+(2+3^4+5)\times (6+78+9).$
\item [] $8186=1+(2+3)^4+56\times (7+8)\times 9.$
\item [] $8187=123+(45+67)\times 8\times 9.$
\item [] $8188=1\times 23\times 4\times (5+67+8+9).$
\item [] $8189=1+23\times 4\times (5+67+8+9).$
\item [] $8190=1\times 234\times (5+6+7+8+9).$
\item [] $8191=1+234\times (5+6+7+8+9).$
\item [] $8192=1\times 2^{(3+4)}\times (5+6\times 7+8+9).$
\item [] $8193=1\times 2\times 3\times 4\times (5+6\times 7\times 8)+9.$
\item [] $8194=1+2\times 3\times 4\times (5+6\times 7\times 8)+9.$
\item [] $8195=1+2\times ((34+5)\times 6+7)\times (8+9).$
\item [] $8196=(1+2\times 3\times 4\times 5)\times 67+89.$
\item [] $8197=1\times 23\times 4\times (5+6+78)+9.$
\item [] $8198=1+23\times 4\times (5+6+78)+9.$
\item [] $8199=12\times 3\times 4\times 56+(7+8)\times 9.$
\item [] $8200=1+2\times (3+45\times 6)\times (7+8)+9.$
\item [] $8201=1+(2+3)\times (4\times 56\times 7+8\times 9).$
\item [] $8202=1\times 2\times (3\times 4\times (5+6\times 7\times 8)+9).$
\item [] $8203=1+2\times (3\times 4\times (5+6\times 7\times 8)+9).$
\item [] $\mathit{8204=1\times 2^3\times 4^5+6+7+8-9.}$
\item [] $8205=1^2\times 3\times 4\times (5+678)+9.$
\item [] $8206=1+(2^3+4)\times (5+678)+9.$
\item [] $8207=1\times 2+3\times 4\times (5+678)+9.$
\item [] $8208=12^3+45\times 6\times (7+8+9).$
\item [] $8209=1+2^3\times (4+5)\times (6\times 7+8\times 9).$
\item [] $\mathit{8210=1\times 2^3\times 4^5-6+7+8+9.}$
\item [] $8211=1\times 23\times (45\times 6+78+9).$
\item [] $8212=1+23\times (45\times 6+78+9).$
\item [] $\mathit{8213=-1+2^3\times 4^5-67+89.}$
\item [] $8214=123\times 4+(5+6)\times 78\times 9.$
\item [] $8215=1+2\times 3\times (4^5+6\times 7\times 8+9).$
\item [] $8216=(12+3^4+5+6)\times (7+8\times 9).$
\item [] $8217=12+3\times 4\times (5+678)+9.$
\item [] $8218=1\times 2\times (3\times 4\times 5\times 67+89).$
\item [] $8219=1+2\times (3\times 4\times 5\times 67+89).$
\item [] $8220=12\times (3\times 4\times 56+7)+8\times 9.$
\item [] $8221=12^3\times 4+(5+6)\times 7\times (8+9).$
\item [] $8222=1\times 2^3\times 4^5+6+7+8+9.$
\item [] $8223=1+2^3\times 4^5+6+7+8+9.$
\item [] $8224=1^2+3\times (4\times (5+678)+9).$
\item [] $8225=(1+23)\times (4+5\times 67)+89.$
\item [] $8226=1+2+3\times (4\times (5+678)+9).$
\item [] $8227=1\times 2+(34\times 5\times 6+7)\times 8+9.$
\item [] $8228=1\times 23\times (45+6)\times 7+8+9.$
\item [] $8229=1+23\times (45+6)\times 7+8+9.$
\item [] $\mathit{8230=1\times 2^3\times (4^5+6)+7-8-9.}$
\item [] $\mathit{8231=1+2^3\times (4^5+6)+7-8-9.}$
\item [] $8232=12\times (3+4+56+7\times 89).$
\item [] $8233=(1^2+34\times 5\times 6+7)\times 8+9.$
\item [] $\mathit{8234=1+2^3\times 4^5-6+7\times 8-9.}$
\item [] $8235=(1+2)\times (3+4\times 56+78)\times 9.$
\item [] $8236=1+(23+4)\times ((5\times 6+7)\times 8+9).$
\item [] $8237=12\times (3\times 4\times 56+7)+89.$
\item [] $8238=1\times 2\times 3\times (4\times (5+6\times 7\times 8)+9).$
\item [] $8239=1+2\times 3\times (4\times (5+6\times 7\times 8)+9).$
\item [] $8240=(1^2+3)\times 4\times (5+6+7\times 8\times 9).$
\item [] $8241=1\times 2\times 3\times 4\times (5\times 67+8)+9.$
\item [] $8242=1+2\times 3\times 4\times (5\times 67+8)+9.$
\item [] $\mathit{8243=1+2^3\times 4^5+67-8-9.}$
\item [] $8244=12\times (34+5\times 6+7\times 89).$
\item [] $8245=(12+3^4+56\times 7)\times (8+9).$
\item [] $8246=1\times 2+(3+4+5)\times (678+9).$
\item [] $8247=1+2+(3+4+5)\times (678+9).$
\item [] $\mathit{8248=1\times 2^3\times (4^5+6)+7-8+9.}$
\item [] $8249=(1\times 234\times 5+6)\times 7+8+9.$
\item [] $8250=(1+2^3)^4+5\times 6\times 7\times 8+9.$
\item[]$\mbox{Decreasing order}$
\item [] $8181=9\times (876+5)+4\times 3\times 21.$
\item [] $8182=(9+(8+7)\times (6+54))\times 3^2+1.$
\item [] $8183=98\times 76+5\times (4+3)\times 21.$
\item [] $8184=9+8+(7+6)\times (5^4+3)+2+1.$
\item [] $8185=9+8\times (76+5^4+321).$
\item [] $8186=9+8\times 7\times (6\times 5+43)\times 2+1.$
\item [] $8187=9+87\times (6\times 5+43+21).$
\item [] $8188=9+87\times ((6+5)\times 4+3)\times 2+1.$
\item [] $8189=9+87\times (6\times 5+4^3)+2\times 1.$
\item [] $8190=(9+87+6\times 5+4)\times 3\times 21.$
\item [] $8191=(98+7)\times (65+4+3^2)+1.$
\item [] $8192=98\times 7+(6\times 5^4+3)\times 2\times 1.$
\item [] $8193=(9+8+7\times 6+5)\times 4\times 32+1.$
\item [] $8194=9+(8\times 7\times 6+5)\times 4\times 3\times 2+1.$
\item [] $8195=98+76+(5\times 4)^3+21.$
\item [] $8196=(9\times 8+76\times 5\times (4+3))\times (2+1).$
\item [] $8197=9+87+6^5+4+321.$
\item [] $8198=98\times 76+(5+4)^3+21.$
\item [] $8199=98+7+6\times (5+4^3\times 21).$
\item [] $8200=9+(8+7)\times 6\times (5+43\times 2)+1.$
\item [] $8201=(9+8+7\times 6)\times ((5+4^3)\times 2+1).$
\item [] $8202=9\times 8+(7+6)\times 5^4+3+2\times 1.$
\item [] $8203=9\times 8+(7+6)\times 5^4+3+2+1.$
\item [] $8204=9\times 876+5\times (43+21).$
\item [] $8205=(987+654)\times (3+2)\times 1.$
\item [] $8206=98+7+6^5+4+321.$
\item [] $8207=9\times 876+5\times 4^3+2+1.$
\item [] $8208=9\times 876+54\times 3\times 2\times 1.$
\item [] $8209=9\times 876+54\times 3\times 2+1.$
\item [] $8210=9\times (8+7+6)+(5\times 4)^3+21.$
\item [] $8211=(98+7)\times (6+5\times 4)\times 3+21.$
\item [] $8212=(9+8)\times 7\times (6+54+3^2)+1.$
\item [] $8213=(9+8)\times 7+6\times (5+4^3\times 21).$
\item [] $8214=9\times 876+5+4+321.$
\item [] $8215=9\times 876+5\times (4^3+2)+1.$
\item [] $8216=(9+8)\times 7\times (65+4)+3+2\times 1.$
\item [] $8217=(9+8)\times 7\times (65+4)+3\times 2\times 1.$
\item [] $8218=(9+8)\times 7\times (65+4)+3\times 2+1.$
\item [] $8219=98\times (7+6)\times 5+43^2\times 1.$
\item [] $8220=(9+8)\times 7+6^5+4+321.$
\item [] $8221=9\times 8+(7+6)\times 5^4+3+21.$
\item [] $8222=98\times 7+6\times (5^4+3)\times 2\times 1.$
\item [] $8223=98\times 7+6\times (5^4+3)\times 2+1.$
\item [] $8224=98+7+6^5+(4+3)^{(2+1)}.$
\item [] $8225=9\times 876+5\times 4+321.$
\item [] $8226=987\times 6+(5+43)^2\times 1.$
\item [] $8227=987\times 6+(5+43)^2+1.$
\item [] $8228=98+(7+6)\times 5^4+3+2\times 1.$
\item [] $8229=98+(7+6)\times 5^4+3+2+1.$
\item [] $8230=98+(7+6)\times 5^4+3\times 2+1.$
\item [] $8231=9+8+(76\times 54+3)\times 2\times 1.$
\item [] $8232=9+8+7+6^5+432\times 1.$
\item [] $8233=9+8+7+6^5+432+1.$
\item [] $8234=98\times (7+65+4\times 3)+2\times 1.$
\item [] $8235=(9\times 87+654\times 3)\times (2+1).$
\item [] $8236=9\times (8+7)+6^5+4+321.$
\item [] $\mathit{8237=9\times 8\times 7+6^5-4^3+21.}$
\item [] $8238=9\times 8\times 7+6\times (5+4\times 321).$
\item [] $8239=9+(8+76\times 54+3)\times 2\times 1.$
\item [] $8240=9+(8+76\times 54+3)\times 2+1.$
\item [] $8241=9\times 876+(5+4\times 3)\times 21.$
\item [] $8242=(9+8)\times (7+6)+(5\times 4)^3+21.$
\item [] $8243=(9+8)\times 7\times (65+4)+32\times 1.$
\item [] $8244=(9+8)\times 7\times (65+4)+32+1.$
\item [] $8245=(9\times 8+7+6)\times ((5+43)\times 2+1).$
\item [] $8246=98\times 7+6\times 5\times 4\times 3\times 21.$
\item [] $8247=98+(7+6)\times 5^4+3+21.$
\item [] $8248=(9+8+76\times 54+3)\times 2\times 1.$
\item [] $8249=(9+8+76\times 54+3)\times 2+1.$
\item [] $8250=98+(7+6)\times 5^4+3^{(2+1)}.$
\item[]$\mbox{Increasing order}$
\item [] $8251=1\times 2^3\times 4^5+6\times 7+8+9.$
\item [] $8252=1+2^3\times 4^5+6\times 7+8+9.$
\item [] $8253=(1\times 234+5+678)\times 9.$
\item [] $8254=1+(234+5+678)\times 9.$
\item [] $8255=(1+2\times 3+4\times 5\times 6)\times (7\times 8+9).$
\item [] $8256=12+(3+4+5)\times (678+9).$
\item [] $8257=1\times 23\times (4+5\times (6+7\times 8+9)).$
\item [] $8258=123\times (4+56+7)+8+9.$
\item [] $8259=(1+2^3\times 4)\times 5\times (6\times 7+8)+9.$
\item [] $8260=(1+2\times 3)\times (4^5+67+89).$
\item [] $\mathit{8261=1+2^3\times 4^5+67-8+9.}$
\item [] $8262=(1+234+5+678)\times 9.$
\item [] $8263=1\times 2^3\times 4^5+6+7\times 8+9.$
\item [] $8264=1+2^3\times 4^5+6+7\times 8+9.$
\item [] $8265=1+2^3\times (4^5+6)+7+8+9.$
\item [] $8266=1+(2+3)\times (4^5+6+7\times 89).$
\item [] $\mathit{8267=1\times 2^3\times 4^5+6+78-9.}$
\item [] $8268=12\times (3\times 4\times 5+6+7\times 89).$
\item [] $8269=1+(2^3+45)\times (67+89).$
\item [] $\mathit{8270=1\times 2\times (3+4^{(5-6+7)})+8\times 9.}$
\item [] $8271=(1+2\times (3+4)\times 5\times (6+7)+8)\times 9.$
\item [] $8272=(123+4)\times 5\times (6+7)+8+9.$
\item [] $8273=((1^2+34\times 5)\times 6+7)\times 8+9.$
\item [] $\mathit{8274=1+2^3\times 4^5-6+78+9.}$
\item [] $8275=(1+(2+3)^4+5)\times (6+7)+8\times 9.$
\item [] $8276=1\times 2^3\times 4^5+67+8+9.$
\item [] $8277=1+2^3\times 4^5+67+8+9.$
\item [] $8278=1+2^3\times 4^5+6+7+8\times 9.$
\item [] $8279=(1+23+456+7)\times (8+9).$
\item [] $8280=12\times (34+567+89).$
\item [] $8281=1+23\times 4\times (5+6+7+8\times 9).$
\item [] $8282=1\times 2+(3+45+67)\times 8\times 9.$
\item [] $8283=1\times 23\times (45+6)\times 7+8\times 9.$
\item [] $8284=1+23\times (45+6)\times 7+8\times 9.$
\item [] $8285=1\times 2^3\times 4^5+6+78+9.$
\item [] $8286=1+2^3\times 4^5+6+78+9.$
\item [] $8287=1+2\times (3^4\times 5+6\times 7\times 89).$
\item [] $8288=(1+2\times 3)\times (45+67\times (8+9)).$
\item [] $8289=((1+2)^3\times (4+5)+678)\times 9.$
\item [] $8290=(1+23\times (45+6))\times 7+8\times 9.$
\item [] $8291=1\times 2^3\times 4^5+6\times (7+8)+9.$
\item [] $8292=12\times (3\times 4+56+7\times 89).$
\item [] $\mathit{8293=-1+2^3\times 4^5+6+7+89.}$
\item [] $8294=1\times 2^3\times 4^5+6+7+89.$
\item [] $8295=1+2^3\times 4^5+6+7+89.$
\item [] $8296=1\times 2\times (34+5\times 6\times 7)\times (8+9).$
\item [] $8297=(1\times 23\times 4+56)\times 7\times 8+9.$
\item [] $8298=1+(23\times 4+56)\times 7\times 8+9.$
\item [] $8299=1+(2+(3+45+67)\times 8)\times 9.$
\item [] $8300=1\times 23\times (45+6)\times 7+89.$
\item [] $8301=1+23\times (45+6)\times 7+89.$
\item [] $\mathit{8302=-1+2^{(3\times 4)}-5+6\times 78\times 9.}$
\item [] $8303=(1^2+3)^4\times 5\times 6+7\times 89.$
\item [] $8304=1^2\times 3\times 4\times (5+678+9).$
\item [] $8305=1^2+3\times 4\times (5+678+9).$
\item [] $8306=1\times 2^3\times 4^5+6\times 7+8\times 9.$
\item [] $8307=1+2^3\times 4^5+6\times 7+8\times 9.$
\item [] $8308=1\times 2345+67\times 89.$
\item [] $8309=1+2345+67\times 89.$
\item [] $8310=(12+3)\times (4+5+67\times 8+9).$
\item [] $8311=(1+234\times 5+6)\times 7+8\times 9.$
\item [] $8312=1+2\times 3+(4^5+6+7)\times 8+9.$
\item [] $8313=123\times (4+56+7)+8\times 9.$
\item [] $8314=1+2^{(3\times 4)}+5+6\times 78\times 9.$
\item [] $\mathit{8315=(1+2)^3\times 4\times (5+6)\times 7+8-9.}$
\item [] $8316=1\times 2\times (3456+78\times 9).$
\item [] $8317=1+2\times (3456+78\times 9).$
\item [] $8318=1+2^3\times 4^5+6+7\times (8+9).$
\item [] $8319=1\times 2^3\times (4^5+6)+7+8\times 9.$
\item [] $8320=12+3+(4^5+6+7)\times 8+9.$
\item[]$\mbox{Decreasing order}$
\item [] $\mathit{8251=9\times 8\times 7+6^5+4-32-1.}$
\item [] $\mathit{8252=9\times 8\times 7+6^5-4-3-21.}$
\item [] $8253=(98+7+6+5\times 4)\times 3\times 21.$
\item [] $8254=9\times (876+5)+4+321.$
\item [] $8255=98+(7+6)\times 5^4+32\times 1.$
\item [] $8256=98+(7+6)\times 5^4+32+1.$
\item [] $8257=(98\times 7\times 6+5+4+3)\times 2+1.$
\item [] $\mathit{8258=9\times 8\times 7+6^5-43+21.}$
\item [] $8259=9\times 876+54+321.$
\item [] $8260=9\times 8+(7+6)\times 5^4+3\times 21.$
\item [] $8261=(98\times 7+65)\times (4+3\times 2+1).$
\item [] $8262=9\times (876+5+4+32+1).$
\item [] $8263=98\times (76+5)+4+321.$
\item [] $8264=98+(7+6)\times (5^4+3)+2\times 1.$
\item [] $8265=9+(8+7)\times (6+543)+21.$
\item [] $8266=(9\times 87+6)\times 5+4321.$
\item [] $8267=(98\times 7\times 6+5+4\times 3)\times 2+1.$
\item [] $\mathit{8268=9\times 8\times 7+6^5-4-3^2+1.}$
\item [] $\mathit{8269=9\times 8\times 7+6^5-4-3\times 2-1.}$
\item [] $8270=(9\times 87+(6+5)\times 4)\times (3^2+1).$
\item [] $8271=9\times (8+(7+6)\times 5\times (4+3)\times 2+1).$
\item [] $8272=9+8+(7+6)\times (5^4+3^2+1).$
\item [] $8273=9+8\times 7+6^5+432\times 1.$
\item [] $8274=9+8\times 7+6^5+432+1.$
\item [] $8275=9\times 8+(7+6)\times (5^4+3+2+1).$
\item [] $8276=9\times 8+(7+6)\times (5^4+3\times 2)+1.$
\item [] $\mathit{8277=9\times 8\times 7+6^5+4-3\times 2-1.}$
\item [] $8278=(98\times 7\times 6+5\times 4+3)\times 2\times 1.$
\item [] $8279=(98\times 7\times 6+5\times 4+3)\times 2+1.$
\item [] $8280=9\times 8\times 7+6\times 54\times (3+21).$
\item [] $8281=98\times (76+5)+(4+3)^{(2+1)}.$
\item [] $8282=9\times 8\times (7+65+43)+2\times 1.$
\item [] $8283=9\times 8\times (7+65+43)+2+1.$
\item [] $\mathit{8284=9\times 8\times 7+6^5+4+3-2-1.}$
\item [] $\mathit{8285=9\times 8\times 7+6^5+4+3-2\times 1.}$
\item [] $8286=9\times (876+5+4)+321.$
\item [] $8287=9\times 8+7+6^5+432\times 1.$
\item [] $8288=9\times 8+7+6^5+432+1.$
\item [] $8289=9\times 8\times 7+6^5+4+3+2\times 1.$
\item [] $8290=9\times 87+(6\times 5^4+3)\times 2+1.$
\item [] $8291=9\times 8\times 7+6^5+4+3\times 2+1.$
\item [] $\mathit{8292=9\times 8\times 7+6^5+4+3^2-1.}$
\item [] $8293=9\times 8\times 7+6^5+4+3^2\times 1.$
\item [] $8294=9\times 8\times 7+6^5+4\times 3+2\times 1.$
\item [] $8295=9\times 8\times 7+6^5+4\times 3+2+1.$
\item [] $8296=(9+8)\times (7\times 65+4\times 3+21).$
\item [] $8297=9+8\times (76+5\times 4^3\times (2+1)).$
\item [] $8298=9\times 87\times 6+(5\times 4\times 3)^2\times 1.$
\item [] $8299=9\times 87\times 6+(5\times 4\times 3)^2+1.$
\item [] $8300=9\times 8\times 7+6^5+4\times (3+2)\times 1.$
\item [] $8301=9\times 8\times (7+65+43)+21.$
\item [] $8302=(98\times 7\times 6+5\times (4+3))\times 2\times 1.$
\item [] $8303=(98\times 7\times 6+5\times (4+3))\times 2+1.$
\item [] $8304=9+87+6^5+432\times 1.$
\item [] $8305=9+87+6^5+432+1.$
\item [] $8306=(98\times 7+6)\times (5+4+3)+2\times 1.$
\item [] $8307=(98\times 7+6)\times (5+4+3)+2+1.$
\item [] $8308=9\times 8\times 7+6^5+4+3+21.$
\item [] $8309=9\times 876+5\times (4^3+21).$
\item [] $8310=(9+8)\times 76\times 5+43^2+1.$
\item [] $8311=9\times 8\times 7+6^5+4+3^{(2+1)}.$
\item [] $8312=98+(76\times 54+3)\times 2\times 1.$
\item [] $8313=98+7+6^5+432\times 1.$
\item [] $8314=98+7+6^5+432+1.$
\item [] $8315=9\times 876+5\times 43\times 2+1.$
\item [] $8316=9\times 8\times 7+6^5+4+32\times 1.$
\item [] $8317=9\times 8\times 7+6^5+4+32+1.$
\item [] $8318=9\times (876+5+43)+2\times 1.$
\item [] $8319=9\times 876+5\times (43\times 2+1).$
\item [] $8320=9\times 87+6\times (5^4+3)\times 2+1.$
\item[]$\mbox{Increasing order}$
\item [] $8321=1\times (234\times 5+6)\times 7+89.$
\item [] $8322=1+(234\times 5+6)\times 7+89.$
\item [] $8323=1\times 2^3\times 4^5+6\times 7+89.$
\item [] $8324=1+2^3\times 4^5+6\times 7+89.$
\item [] $8325=(12+3)\times (45+6+7\times 8\times 9).$
\item [] $\mathit{8326=-1+2^3\times (4^5+6)+78+9.}$
\item [] $8327=(123+4)\times 5\times (6+7)+8\times 9.$
\item [] $8328=(1+234\times 5+6)\times 7+89.$
\item [] $8329=1+23+(4^5+6+7)\times 8+9.$
\item [] $8330=123\times (4+56+7)+89.$
\item [] $8331=1\times 2^3\times 4^5+67+8\times 9.$
\item [] $8332=1+2^3\times 4^5+67+8\times 9.$
\item [] $8333=12\times 3\times (4\times 56+7)+8+9.$
\item [] $8334=1+2^3\times 4^5+6+(7+8)\times 9.$
\item [] $8335=1+(2\times 34+(5+6)\times 78)\times 9.$
\item [] $8336=1\times 2^3\times 4^5+6\times (7+8+9).$
\item [] $8337=1+2^3\times 4^5+6\times (7+8+9).$
\item [] $8338=(1+2\times 3+4)\times (56+78\times 9).$
\item [] $8339=((1+234)\times 5+6)\times 7+8\times 9.$
\item [] $8340=1^2\times 3\times 4\times 5\times (67+8\times 9).$
\item [] $8341=12\times 3+(4^5+6+7)\times 8+9.$
\item [] $8342=1\times 2+3\times 4\times 5\times (67+8\times 9).$
\item [] $8343=1+2+3\times 4\times 5\times (67+8\times 9).$
\item [] $8344=(123+4)\times 5\times (6+7)+89.$
\item [] $8345=(1\times 2+3+4^5+6+7)\times 8+9.$
\item [] $8346=1+(2+3+4^5+6+7)\times 8+9.$
\item [] $8347=(1\times 2+3)^4+(5+6)\times 78\times 9.$
\item [] $8348=1\times 2^3\times 4^5+67+89.$
\item [] $8349=1+2^3\times 4^5+67+89.$
\item [] $8350=(1+2^{(3+4)}\times 5)\times (6+7)+8+9.$
\item [] $8351=(1+2\times 3)\times (4\times (5\times 6+7)\times 8+9).$
\item [] $8352=12\times 345+6\times 78\times 9.$
\item [] $8353=(123+4\times 5+6)\times 7\times 8+9.$
\item [] $8354=1+2^{(3\times 4)}+(5+6\times 78)\times 9.$
\item [] $8355=(1+(2+3)\times 4\times 5+6)\times 78+9.$
\item [] $8356=((1+234)\times 5+6)\times 7+89.$
\item [] $\mathit{8357=1+2\times ((3-4+5)^6-7+89).}$
\item [] $8358=1\times 2\times (3+(45+6+7)\times 8\times 9).$
\item [] $8359=(1\times 2\times 3^4+5)\times (6\times 7+8)+9.$
\item [] $8360=1+2^3\times (4^5+6)+7\times (8+9).$
\item [] $8361=(12\times 3^4+5+67)\times 8+9.$
\item [] $8362=(1+2+34)\times (5+(6+7)\times (8+9)).$
\item [] $\mathit{8363=1\times 2+3\times 45\times (6+7\times 8)-9.}$
\item [] $8364=123\times (4+5+6\times 7+8+9).$
\item [] $\mathit{8365=1+2\times (-3+45\times (6+78+9)).}$
\item [] $8366=(12\times 3+45+6+7)\times 89.$
\item [] $8367=1^2+(3\times (4+5)+67)\times 89.$
\item [] $8368=1\times 2+(3\times (4+5)+67)\times 89.$
\item [] $8369=(1+2\times 3)^4+5+67\times 89.$
\item [] $8370=(12+3^4)\times (5+6+7+8\times 9).$
\item [] $8371=1+(2^3\times 4+5\times 6)\times (7+8)\times 9.$
\item [] $8372=(1+23+4)\times (5\times 6\times 7+89).$
\item [] $8373=12\times (3^4+(5+6)\times 7\times 8)+9.$
\item [] $8374=1\times ((2+3)\times 4\times 5+6)\times (7+8\times 9).$
\item [] $8375=1+((2+3)\times 4\times 5+6)\times (7+8\times 9).$
\item [] $8376=12^3+4^5\times 6+7\times 8\times 9.$
\item [] $8377=1+2\times (3+45\times (6+78+9)).$
\item [] $8378=1+2\times ((3^4+5)\times 6+7)\times 8+9.$
\item [] $8379=1^2\times 3\times 45\times (6+7\times 8)+9.$
\item [] $8380=1^2+3\times 45\times (6+7\times 8)+9.$
\item [] $8381=1\times 2+3\times 45\times (6+7\times 8)+9.$
\item [] $8382=1+2+3\times 45\times (6+7\times 8)+9.$
\item [] $8383=1^2+3\times (4+5\times (6+7\times 8)\times 9).$
\item [] $8384=1\times 2+3\times (4+5\times (6+7\times 8)\times 9).$
\item [] $8385=1+2+3\times (4+5\times (6+7\times 8)\times 9).$
\item [] $8386=1+2^3\times (4^5+6+7)+89.$
\item [] $8387=1+2\times (((3^4+5)\times 6+7)\times 8+9).$
\item [] $8388=12\times 3\times (4\times 56+7)+8\times 9.$
\item [] $8389=1^2+(3\times 4\times (5+6)\times 7+8)\times 9.$
\item [] $8390=1\times 2+(3\times 4\times (5+6)\times 7+8)\times 9.$
\item[]$\mbox{Decreasing order}$
\item [] $8321=9\times 876+5+432\times 1.$
\item [] $8322=9\times 876+5+432+1.$
\item [] $8323=(9+8\times 7+65)\times 4^3+2+1.$
\item [] $\mathit{8324=-9+8765-432\times 1.}$
\item [] $8325=9\times 8\times 7+6^5+43+2\times 1.$
\item [] $8326=9\times 8\times 7+6^5+43+2+1.$
\item [] $8327=(9+8)\times 7+6^5+432\times 1.$
\item [] $8328=(9+8)\times 7+6^5+432+1.$
\item [] $8329=(98\times 7\times 6+5+43)\times 2+1.$
\item [] $8330=98\times (7+65+4+3^2\times 1).$
\item [] $8331=(98+7+65)\times (4+3)^2+1.$
\item [] $8332=98\times (7+(6+5\times 4)\times 3)+2\times 1.$
\item [] $8333=98\times (7+(6+5\times 4)\times 3)+2+1.$
\item [] $8334=9\times (876+5+43+2\times 1).$
\item [] $8335=98\times (76+5+4)+3+2\times 1.$
\item [] $8336=98\times (76+5+4)+3+2+1.$
\item [] $8337=98\times (76+5+4)+3\times 2+1.$
\item [] $8338=9+8+(7+6)\times 5\times 4^3\times 2+1.$
\item [] $8339=98\times (76+5+4)+3^2\times 1.$
\item [] $8340=98\times (7+65)+4\times 321.$
\item [] $8341=(9+8\times 7+65)\times 4^3+21.$
\item [] $\mathit{8342=9+8765-432\times 1.}$
\item [] $8343=9\times 87+6\times 5\times 4\times 3\times 21.$
\item [] $8344=9\times 8\times 7+6^5+43+21.$
\item [] $8345=9+8+7+65\times 4\times 32+1.$
\item [] $8346=9\times 8\times 7+6^5+4^3+2\times 1.$
\item [] $8347=9\times 8\times 7+6^5+4+3\times 21.$
\item [] $8348=9+8\times 7\times 6+(5\times 4)^3+2+1.$
\item [] $8349=(9\times (8+7\times 65)+4+3)\times 2+1.$
\item [] $8350=9+8+(7+6)\times (5\times 4\times 32+1).$
\item [] $8351=98\times (7+(6+5\times 4)\times 3)+21.$
\item [] $8352=9\times (876+5\times 4+32\times 1).$
\item [] $8353=(98\times 7\times 6+5\times 4\times 3)\times 2+1.$
\item [] $8354=98\times (76+5+4)+3+21.$
\item [] $8355=(98+76)\times (5+43)+2+1.$
\item [] $8356=98\times 76+5+43\times 21.$
\item [] $8357=98\times (76+5+4)+3^{(2+1)}.$
\item [] $8358=(9\times 8+76\times 54+3)\times 2\times 1.$
\item [] $8359=(9\times 8+76\times 54+3)\times 2+1.$
\item [] $8360=9+(8+7+65\times 4^3)\times 2+1.$
\item [] $8361=9\times (876+5)+432\times 1.$
\item [] $8362=9\times (876+5)+432+1.$
\item [] $8363=98\times (76+5+4)+32+1.$
\item [] $8364=(9+8)\times (7\times 65+4+32+1).$
\item [] $8365=9\times 8\times 7+6^5+4^3+21.$
\item [] $8366=9\times 8\times 7+6^5+43\times 2\times 1.$
\item [] $8367=9\times 8\times 7+6^5+43\times 2+1.$
\item [] $8368=(9+8+7+65\times 4^3)\times 2\times 1.$
\item [] $8369=(9+8+7+65\times 4^3)\times 2+1.$
\item [] $8370=9\times 876+54\times 3^2\times 1.$
\item [] $8371=9\times 876+54\times 3^2+1.$
\item [] $8372=98+7\times 6\times (5+4^3\times (2+1)).$
\item [] $8373=(98+76)\times (5+43)+21.$
\item [] $8374=(9\times 8+7)\times (6+5\times 4\times (3+2\times 1)).$
\item [] $8375=9\times (876+54)+3+2\times 1.$
\item [] $8376=9\times (876+54)+3+2+1.$
\item [] $8377=9\times (876+54)+3\times 2+1.$
\item [] $\mathit{8378=9-876+5\times 43^2\times 1.}$
\item [] $8379=9+(876+54)\times 3^2\times 1.$
\item [] $8380=9\times (876+54)+3^2+1.$
\item [] $8381=(9+8)\times (7+6+5\times 4\times (3+21)).$
\item [] $8382=(9\times 8+7\times (6+5\times 4))\times (32+1).$
\item [] $\mathit{8383=-9-8+7\times 6\times 5\times 4\times (3^2+1).}$
\item [] $8384=(9\times (87+6)\times 5+4+3)\times 2\times 1.$
\item [] $8385=9+8\times 7+65\times 4\times 32\times 1.$
\item [] $8386=9+8\times 7+65\times 4\times 32+1.$
\item [] $8387=(9+8\times 7)\times (65+4^3)+2\times 1.$
\item [] $8388=9\times (8\times 76+54\times 3\times 2\times 1).$
\item [] $8389=(9\times (8+7\times (6+5\times 4\times 3)\times 2)+1).$
\item [] $8390=9\times (8+7\times (6+5)\times 4\times 3)+2\times 1.$
\item[]$\mbox{Increasing order}$
\item [] $8391=12+3\times 45\times (6+7\times 8)+9.$
\item [] $8392=(1+2\times 3\times 4)\times 5\times 67+8+9.$
\item [] $8393=(1\times 23\times 45+6+7)\times 8+9.$
\item [] $8394=12+3\times (4+5\times (6+7\times 8)\times 9).$
\item [] $8395=1\times 23\times (4\times 5+6\times 7\times 8+9).$
\item [] $8396=12^3+4+56\times 7\times (8+9).$
\item [] $8397=(1+2)^3\times (4\times 56+78+9).$
\item [] $8398=1^2+3\times (45\times (6+7\times 8)+9).$
\item [] $8399=(1+2+34)\times (5\times 6\times 7+8+9).$
\item [] $8400=1\times 2\times 3\times 4\times (5+6\times 7\times 8+9).$
\item [] $8401=1+2\times 3\times 4\times (5+6\times 7\times 8+9).$
\item [] $8402=(1+2^{(3+4)})\times 5\times (6+7)+8+9.$
\item [] $\mathit{8403=123+4\times 5\times 6\times (78-9).}$
\item [] $8404=1\times 23\times ((45+6)\times 7+8)+9.$
\item [] $8405=12\times 3\times (4\times 56+7)+89.$
\item [] $8406=(1+2)\times (3+45\times (6+7\times 8)+9).$
\item [] $8407=1+(2\times (3+(45+6+7)\times 8))\times 9.$
\item [] $\mathit{8408=-1+2^{(3+4)}\times 5\times (6+7)+89.}$
\item [] $8409=12\times 3\times 45+6789.$
\item [] $8410=1+(2+3+4\times 5)\times 6\times 7\times 8+9.$
\item [] $\mathit{8411=-1+2\times (3+(4+5)\times 6\times 78-9).}$
\item [] $8412=12\times (3\times (4\times 5+6)+7\times 89).$
\item [] $8413=1\times 2^3\times 4^5+(6+7)\times (8+9).$
\item [] $8414=1+2^3\times 4^5+(6+7)\times (8+9).$
\item [] $8415=(12+3)\times (4\times 5+6+7)\times (8+9).$
\item [] $8416=1+((2+3)^4+5\times (6+7\times 8))\times 9.$
\item [] $8417=(12\times 3+4)\times 5\times 6\times 7+8+9.$
\item [] $8418=1\times 23\times (45\times 6+7+89).$
\item [] $8419=1+23\times (45\times 6+7+89).$
\item [] $\mathit{8420=1-2-3+(4+5)\times (6+7)\times 8\times 9.}$
\item [] $8421=12\times (3+4\times 5+678)+9.$
\item [] $8422=1\times 2\times (3\times 4^5+67\times (8+9)).$
\item [] $8423=123\times 4\times 5+67\times 89.$
\item [] $8424=(1\times 2+3+45+67)\times 8\times 9.$
\item [] $8425=1+2\times 3\times (4+5)\times (67+89).$
\item [] $8426=1\times 2+(345+6)\times (7+8+9).$
\item [] $8427=1+2+(345+6)\times (7+8+9).$
\item [] $8428=123+(4^5+6+7)\times 8+9.$
\item [] $8429=1\times 2+3+(4+5)\times (6+7)\times 8\times 9.$
\item [] $8430=1\times 2\times 3+(4+5)\times (6+7)\times 8\times 9.$
\item [] $8431=1+2\times 3+(4+5)\times (6+7)\times 8\times 9.$
\item [] $8432=1\times 2^3+(4+5)\times (6+7)\times 8\times 9.$
\item [] $8433=1\times 2\times (3+45+6)\times 78+9.$
\item [] $8434=1+2\times (3+45+6)\times 78+9.$
\item [] $8435=(1+(2\times 3)^4+5)\times 6+7\times 89.$
\item [] $8436=12+(345+6)\times (7+8+9).$
\item [] $8437=(123+4\times 5)\times (6\times 7+8+9).$
\item [] $\mathit{8438=1\times 2\times (3\times 4-5+6\times 78\times 9).}$
\item [] $8439=12+3+(4+5)\times (6+7)\times 8\times 9.$
\item [] $8440=1+2\times (3+(4+5)\times 6\times 78)+9.$
\item [] $8441=(1^2+3\times 45)\times (6+7\times 8)+9.$
\item [] $8442=1^2\times (3+4)\times (56+78)\times 9.$
\item [] $8443=1^2+(3+4)\times (56+78)\times 9.$
\item [] $8444=1\times 2+(3+4)\times (56+78)\times 9.$
\item [] $8445=1+2+(3+4)\times (56+78)\times 9.$
\item [] $8446=1+(2+3)\times (4\times 5\times (6+78)+9).$
\item [] $8447=1\times 2345+678\times 9.$
\item [] $8448=1+2345+678\times 9.$
\item [] $8449=1+2\times 3\times 4\times (5\times 67+8+9).$
\item [] $8450=(1\times 2\times 3+4)\times (56+789).$
\item [] $8451=1+(2\times 3+4)\times (56+789).$
\item [] $8452=1^2+(3^4+(5+6)\times 78)\times 9.$
\item [] $8453=((1+23)\times 4+5+6)\times (7+8\times 9).$
\item [] $8454=12+(3+4)\times (56+78)\times 9.$
\item [] $8455=(1\times 2^3+45+6\times 7)\times 89.$
\item [] $8456=1+(23+4\times (5+6+7))\times 89.$
\item [] $8457=12^3+4\times 5\times 6\times 7\times 8+9.$
\item [] $8458=1\times 2\times (3\times 4+5+6\times 78\times 9).$
\item [] $8459=1+2\times (3\times 4+5+6\times 78\times 9).$
\item [] $8460=(12+3)\times (4+56+7\times 8\times 9).$
\item[]$\mbox{Decreasing order}$
\item [] $8391=9+((8+7\times 6)\times 5+4)\times (32+1).$
\item [] $8392=9\times 8+(7+6)\times 5\times 4^3\times 2\times 1.$
\item [] $8393=98\times (76+5+4)+3\times 21.$
\item [] $8394=98\times 76+5^4+321.$
\item [] $8395=(9\times (87+6)\times 5+4\times 3)\times 2+1.$
\item [] $8396=9+8+7\times (6\times (5+4)+3)\times 21.$
\item [] $8397=9\times (8\times 76+54\times 3\times 2+1).$
\item [] $8398=(9+8)\times (7+6)\times (5+4\times 3+21).$
\item [] $8399=9\times 8+7+65\times 4\times 32\times 1.$
\item [] $8400=9\times 8+7+65\times 4\times 32+1.$
\item [] $8401=(9+8\times 7\times 6+5)\times 4\times 3\times 2+1.$
\item [] $8402=9\times (876+54)+32\times 1.$
\item [] $8403=9\times (876+54)+32+1.$
\item [] $\mathit{8404=9\times 8\times 7+6^5+4\times (32-1).}$
\item [] $8405=9\times 8+(7+6)\times (5\times 4\times 32+1).$
\item [] $8406=(9+8\times 7)\times (65+4^3)+21.$
\item [] $8407=((9+87)\times 6+5^4)\times (3\times 2+1).$
\item [] $8408=9\times 8\times 7+6^5+4\times 32\times 1.$
\item [] $8409=9\times 8\times 7+6^5+4\times 32+1.$
\item [] $8410=(98+76\times 54+3)\times 2\times 1.$
\item [] $8411=9+8\times 7+(6+5\times 4)\times 321.$
\item [] $8412=9\times 8\times 7+6^5+4\times (32+1).$
\item [] $8413=9\times 876+(5\times 4+3)^2\times 1.$
\item [] $8414=9\times 876+(5\times 4+3)^2+1.$
\item [] $8415=9\times (876+54+3+2\times 1).$
\item [] $8416=9+87+65\times 4\times 32\times 1.$
\item [] $8417=9+87+65\times 4\times 32+1.$
\item [] $8418=98+(7+6)\times 5\times 4^3\times 2\times 1.$
\item [] $8419=98+(7+6)\times 5\times 4^3\times 2+1.$
\item [] $8420=98\times 7+6\times (5+4\times 321).$
\item [] $8421=9\times (8+7)\times (6+54)+321.$
\item [] $8422=((987+65)\times 4+3)\times 2\times 1.$
\item [] $8423=((987+65)\times 4+3)\times 2+1.$
\item [] $8424=9\times (876+54+3+2+1).$
\item [] $8425=98+7+65\times 4\times 32\times 1.$
\item [] $8426=98+7+65\times 4\times 32+1.$
\item [] $8427=9\times 8\times 7+6^5+(4+3)\times 21.$
\item [] $8428=98\times (7+65+4+3^2+1).$
\item [] $8429=9\times 876+543+2\times 1.$
\item [] $8430=9\times 876+543+2+1.$
\item [] $8431=9\times 8\times (7+6)\times (5+4)+3\times 2+1.$
\item [] $8432=9+8+765\times (4+3\times 2+1).$
\item [] $8433=9\times (876+54)+3\times 21.$
\item [] $8434=9+8\times (7+6)\times (5+4)\times 3^2+1.$
\item [] $\mathit{8435=98\times 7+6^5+4-32+1.}$
\item [] $8436=(9+8+(7+6)\times 5\times 43)\times (2+1).$
\item [] $8437=9+(87+6+5)\times 43\times 2\times 1.$
\item [] $8438=9+(87+6+5)\times 43\times 2+1.$
\item [] $8439=(9+8)\times 7+65\times 4\times 32\times 1.$
\item [] $8440=(9+8)\times 7+65\times 4\times 32+1.$
\item [] $8441=9+8\times ((76+5)\times (4+3^2)+1).$
\item [] $8442=98\times 7\times 6+5+4321.$
\item [] $8443=(98+7\times 6\times 5\times 4)\times 3^2+1.$
\item [] $\mathit{8444=(9+(8+76)\times 5\times 4)\times (3+2)-1.}$
\item [] $8445=(9+(8+76)\times 5\times 4)\times (3+2\times 1).$
\item [] $8446=(9+(87+6+5)\times 43)\times 2\times 1.$
\item [] $8447=(9+(87+6+5)\times 43)\times 2+1.$
\item [] $8448=9\times 876+543+21.$
\item [] $8449=9\times 87\times 6+5^4\times 3\times 2+1.$
\item [] $8450=9+8\times 7+65\times 43\times (2+1).$
\item [] $8451=98+7+(6+5\times 4)\times 321.$
\item [] $8452=(9+876+54)\times 3^2+1.$
\item [] $8453=9\times (8+7\times 6)+(5\times 4)^3+2+1.$
\item [] $8454=9+8+(7+6)\times (5^4+3+21).$
\item [] $8455=9\times (8+7)+65\times 4\times 32\times 1.$
\item [] $8456=9\times 8\times 7\times 6+5432\times 1.$
\item [] $8457=9\times 8\times 7\times 6+5432+1.$
\item [] $8458=9+8\times (7+6+5\times 4)\times 32+1.$
\item [] $8459=9+(8+7+6+5)\times (4+321).$
\item [] $8460=9\times (876+54+3^2+1).$
\item[]$\mbox{Increasing order}$
\item [] $8461=1^2+3\times 4\times 5\times (6+(7+8)\times 9).$
\item [] $8462=1\times 2+3\times 4\times 5\times (6+(7+8)\times 9).$
\item [] $8463=1+2+3\times 4\times 5\times (6+(7+8)\times 9).$
\item [] $8464=(1+2\times 3\times 4)\times 5\times 67+89.$
\item [] $8465=(1\times 2+3)\times (4+5\times 6\times 7\times 8+9).$
\item [] $8466=1+(2+3)\times (4+5\times 6\times 7\times 8+9).$
\item [] $\mathit{8467=1+2\times 3+(4^5-6-78)\times 9.}$
\item [] $\mathit{8468=1\times 2^3+(4^5-6-78)\times 9.}$
\item [] $8469=(1^2+3)^4\times 5\times 6+789.$
\item [] $8470=1\times 2\times (3+4\times 5+6\times 78\times 9).$
\item [] $8471=1+2\times (3+4\times 5+6\times 78\times 9).$
\item [] $8472=1\times 2^3\times (45\times 6+789).$
\item [] $8473=1+2^3\times (45\times 6+789).$
\item [] $8474=1+2^{(3\times 4)}+56\times 78+9.$
\item [] $8475=(12+3)\times (4\times 5+67\times 8+9).$
\item [] $\mathit{8476=(12+3)\times (4+567)-89.}$
\item [] $\mathit{8477=1-2+3\times (4+5\times (6+7\times 8))\times 9.}$
\item [] $8478=1\times 2\times (3\times 45+6\times 7\times 8)\times 9.$
\item [] $8479=(1^2+3\times 4\times 5)\times (67+8\times 9).$
\item [] $8480=1\times (2+3)\times 4\times (5\times 67+89).$
\item [] $8481=(123\times 4+567)\times 8+9.$
\item [] $\mathit{8482=-1\times 23+45\times (6+7+8)\times 9.}$
\item [] $8483=(12\times 3+456+7)\times (8+9).$
\item [] $8484=12^3+4\times (5\times 6\times 7\times 8+9).$
\item [] $8485=1+2\times 3\times (4^5+6\times (7\times 8+9)).$
\item [] $\mathit{8486=-1+23\times (4\times 5+6+7+8)\times 9.}$
\item [] $8487=1\times 23\times (4\times 5+6+7+8)\times 9.$
\item [] $8488=1+23\times (4\times 5+6+7+8)\times 9.$
\item [] $8489=(12\times 3+4)\times 5\times 6\times 7+89.$
\item [] $8490=1\times 2\times (3456+789).$
\item [] $8491=1+2\times (3456+789).$
\item [] $8492=12^3+(4+5+67)\times 89.$
\item [] $8493=((1+23)\times 4+5)\times (6+78)+9.$
\item [] $8494=1\times 2\times ((3+4)\times 5+6\times 78\times 9).$
\item [] $8495=12^3+4^5\times 6+7\times 89.$
\item [] $8496=(1\times 2\times 3+45+67)\times 8\times 9.$
\item [] $8497=(1+23+4^5+6+7)\times 8+9.$
\item [] $8498=1\times 2+(3+4+5)\times (6+78\times 9).$
\item [] $8499=1+2+(3+4+5)\times (6+78\times 9).$
\item [] $8500=1\times (2+3)\times 4\times 5\times (6+7+8\times 9).$
\item [] $8501=1+(2+3)\times 4\times 5\times (6+7+8\times 9).$
\item [] $8502=1\times 2\times (34+5+6\times 78\times 9).$
\item [] $8503=1+2\times (34+5+6\times 78\times 9).$
\item [] $8504=1+(2+3\times 45)\times (6+7\times 8)+9.$
\item [] $8505=1^{23}\times 45\times (6+7+8)\times 9.$
\item [] $8506=1^{23}+45\times (6+7+8)\times 9.$
\item [] $8507=1\times 2+(3+4+56)\times (7+8)\times 9.$
\item [] $8508=(1+2\times 3)^4+5+678\times 9.$
\item [] $8509=1^2\times 34\times 5\times (6\times 7+8)+9.$
\item [] $8510=1\times 2+3+45\times (6+7+8)\times 9.$
\item [] $8511=1+2+3+45\times (6+7+8)\times 9.$
\item [] $8512=1+2\times 3+45\times (6+7+8)\times 9.$
\item [] $8513=1\times 2^3+45\times (6+7+8)\times 9.$
\item [] $8514=12^3\times 4+(5+6+7)\times 89.$
\item [] $8515=1^2+3^4\times 5\times (6+7+8)+9.$
\item [] $8516=1\times 2+3^4\times 5\times (6+7+8)+9.$
\item [] $8517=1+2+3^4\times 5\times (6+7+8)+9.$
\item [] $8518=(123+4)\times (5+6+7\times 8)+9.$
\item [] $8519=(1\times 2+3)\times (4^5+678)+9.$
\item [] $8520=12+3+45\times (6+7+8)\times 9.$
\item [] $8521=12+34\times 5\times (6\times 7+8)+9.$
\item [] $8522=(1+2+3\times 4)\times 567+8+9.$
\item [] $8523=1^2\times (3^4\times (5+6)+7\times 8)\times 9.$
\item [] $8524=1+(2+(3+4+56)\times (7+8))\times 9.$
\item [] $8525=1\times 2+(3^4\times (5+6)+7\times 8)\times 9.$
\item [] $8526=12^3+4+5+6789.$
\item [] $8527=1+(2\times 34+5\times 6)\times (78+9).$
\item [] $8528=1\times 23+45\times (6+7+8)\times 9.$
\item [] $8529=1+23+45\times (6+7+8)\times 9.$
\item [] $\mathit{8530=1-2^3+(4^5+6\times 7)\times 8+9.}$
\item[]$\mbox{Decreasing order}$
\item [] $8461=9\times (8+7\times 65+4+3)\times 2+1.$
\item [] $8462=98\times 7+6\times 54\times (3+21).$
\item [] $8463=9+8+(7+6)\times 5^4+321.$
\item [] $8464=9\times 8+7+65\times 43\times (2+1).$
\item [] $8465=(9+8)\times 7+(6+5\times 4)\times 321.$
\item [] $8466=9+8\times 7\times (65+43\times 2)+1.$
\item [] $\mathit{8467=-9\times 87+6+5\times 43^2-1.}$
\item [] $8468=(9+8)\times (7\times 65+43)+2\times 1.$
\item [] $8469=(9+8)\times (7\times 65+43)+2+1.$
\item [] $8470=(9\times 8\times (7+6)+5)\times (4+3+2)+1.$
\item [] $8471=98\times 7+6^5+4+3+2\times 1.$
\item [] $8472=98\times 7+6^5+4+3+2+1.$
\item [] $8473=98\times 7+6^5+4+3\times 2+1.$
\item [] $8474=9+(8+7+65+4\times 3)^2+1.$
\item [] $8475=98\times 7+6^5+4+3^2\times 1.$
\item [] $8476=98\times 7+6^5+4+3^2+1.$
\item [] $8477=98\times 7+6^5+4\times 3+2+1.$
\item [] $8478=(9\times 8+7+65\times 4^3)\times 2\times 1.$
\item [] $8479=(9\times 8+7+65\times 4^3)\times 2+1.$
\item [] $8480=(9\times 87+65)\times (4+3+2+1).$
\item [] $8481=9+87+65\times 43\times (2+1).$
\item [] $8482=98\times 7+6^5+4\times (3+2\times 1).$
\item [] $8483=98\times 7+6^5+4\times (3+2)+1.$
\item [] $8484=(9+8+76\times 5+4+3)\times 21.$
\item [] $8485=(9\times 87+6+5^4)\times 3\times 2+1.$
\item [] $8486=98\times 7+6^5+4\times 3\times 2\times 1.$
\item [] $8487=98\times 7+6^5+4\times 3\times 2+1.$
\item [] $8488=9\times (8+7\times 65)+4321.$
\item [] $8489=9+(8+7\times 6\times 5\times 4)\times (3^2+1).$
\item [] $8490=(987+6\times 543)\times 2\times 1.$
\item [] $8491=98\times 7\times 6+5^4\times (3\times 2+1).$
\item [] $\mathit{8492=9\times 87+6^5-4-3\times 21.}$
\item [] $8493=987+(6\times 5^4+3)\times 2\times 1.$
\item [] $8494=987+(6\times 5^4+3)\times 2+1.$
\item [] $8495=98\times 7+6^5+4\times 3+21.$
\item [] $8496=9\times 8\times (76+5+4+32+1).$
\item [] $8497=9\times 8\times (7\times 6+5+4\times 3)\times 2+1.$
\item [] $8498=98\times 7+6^5+4+32\times 1.$
\item [] $8499=98\times 7+6^5+4+32+1.$
\item [] $8500=(9\times 8+7+6)\times 5\times 4\times (3+2)\times 1.$
\item [] $8501=(9\times 8+7+6)\times 5\times 4\times (3+2)+1.$
\item [] $8502=98\times 7+6^5+4\times (3^2+1).$
\item [] $8503=9+(8\times 7+6)\times (5+4\times (32+1)).$
\item [] $8504=(9+8)\times 7+65\times 43\times (2+1).$
\item [] $8505=9\times (876+5+43+21).$
\item [] $8506=9\times (87+6)\times 5+4321.$
\item [] $8507=98\times 7+6^5+43+2\times 1.$
\item [] $8508=98\times 7+6^5+43+2+1.$
\item [] $8509=9\times 8+(7+6)\times (5^4+3+21).$
\item [] $\mathit{8510=9\times 876+5^4+3-2\times 1.}$
\item [] $8511=98\times 7+6^5+(4+3)^2\times 1.$
\item [] $8512=98\times 7+6^5+(4+3)^2+1.$
\item [] $8513=98+765\times (4+3\times 2+1).$
\item [] $8514=9\times 876+5^4+3+2\times 1.$
\item [] $8515=9\times 876+5^4+3+2+1.$
\item [] $8516=9\times 876+5^4+3\times 2+1.$
\item [] $8517=9\times 87+6\times (5+4\times 321).$
\item [] $8518=9\times 876+5^4+3^2\times 1.$
\item [] $8519=9\times 876+5^4+3^2+1.$
\item [] $8520=9\times (8+7)+65\times 43\times (2+1).$
\item [] $8521=(9+8\times 7+6)\times 5\times 4\times 3\times 2+1.$
\item [] $8522=9+8\times 76\times (5+4+3+2)+1.$
\item [] $8523=987+6\times (5^4+3)\times 2\times 1.$
\item [] $8524=9\times 876+5\times 4\times 32\times 1.$
\item [] $8525=9\times 876+5\times 4\times 32+1.$
\item [] $8526=98\times 7+6^5+43+21.$
\item [] $8527=(9+(8+7)\times (6+5))\times (4+3)^2+1.$
\item [] $8528=98\times 7+6^5+4^3+2\times 1.$
\item [] $8529=98\times 7+6^5+4+3\times 21.$
\item [] $8530=(98+7+65\times 4^3)\times 2\times 1.$
\item[]$\mbox{Increasing order}$
\item [] $\mathit{8531=(12^3+4)\times 5+6-(7+8)\times 9.}$
\item [] $8532=12^3+(4+5)\times (6+78)\times 9.$
\item [] $8533=1+2\times (3+456+7+8)\times 9.$
\item [] $8534=(1\times 23\times 4\times 5+6\times 7)\times (8+9).$
\item [] $8535=(1^2+3^4\times 5)\times (6+7+8)+9.$
\item [] $\mathit{8536=-1+2^3\times 4^5+6\times 7\times 8+9.}$
\item [] $8537=12^3+4\times 5+6789.$
\item [] $8538=1+2^3\times 4^5+6\times 7\times 8+9.$
\item [] $8539=1+2\times (3+(4+5)\times 6\times (7+8\times 9)).$
\item [] $8540=1^2\times 3+(4^5+6\times 7)\times 8+9.$
\item [] $8541=12\times 3+45\times (6+7+8)\times 9.$
\item [] $8542=1\times 2+3+(4^5+6\times 7)\times 8+9.$
\item [] $8543=1\times 2\times 3+(4^5+6\times 7)\times 8+9.$
\item [] $8544=1\times 2\times (3\times 4\times 5+6\times 78\times 9).$
\item [] $8545=1+2\times (3\times 4\times 5+6\times 78\times 9).$
\item [] $8546=1+2^3+(4^5+6\times 7)\times 8+9.$
\item [] $8547=123+(4+5)\times (6+7)\times 8\times 9.$
\item [] $8548=(1+2\times 3)^4+(5+678)\times 9.$
\item [] $\mathit{8549=1234\times (-5+6)\times 7-89.}$
\item [] $8550=(1+2+3\times 4)\times 5\times (6\times 7+8\times 9).$
\item [] $8551=(1+23\times 4\times 5+6\times 7)\times (8+9).$
\item [] $8552=12+3+(4^5+6\times 7)\times 8+9.$
\item [] $8553=(1+23)\times 4\times (5+6+78)+9.$
\item [] $\mathit{8554=-1\times 2+3\times 4\times (5+6+78\times 9).}$
\item [] $8555=1\times (2+3)\times (4^5+678+9).$
\item [] $8556=12\times (34+56+7\times 89).$
\item [] $8557=1^2+3\times 4\times (5+6+78\times 9).$
\item [] $8558=1\times 2+3\times 4\times (5+6+78\times 9).$
\item [] $8559=123\times 45+6\times 7\times 8\times 9.$
\item [] $8560=1\times 23+(4^5+6\times 7)\times 8+9.$
\item [] $8561=1+23+(4^5+6\times 7)\times 8+9.$
\item [] $8562=123\times 4\times 5+678\times 9.$
\item [] $8563=1\times 2+(3+4^5+6\times 7)\times 8+9.$
\item [] $8564=(1+2)^3+(4^5+6\times 7)\times 8+9.$
\item [] $8565=(1+2+3\times 45)\times (6+7\times 8)+9.$
\item [] $\mathit{8566=1234\times (-5+6)\times 7-8\times 9.}$
\item [] $\mathit{8567=-1+2\times 3^4\times 56-7\times 8\times 9.}$
\item [] $8568=12\times 3\times 4\times 56+7\times 8\times 9.$
\item [] $8569=1+2\times 34\times (5\times 6+7+89).$
\item [] $8570=1+(2+3)\times (4+5\times 6\times 7)\times 8+9.$
\item [] $8571=1+2+(3+4+5)\times 6\times 7\times (8+9).$
\item [] $\mathit{8572=(12^3+4)\times 5-6+7-89.}$
\item [] $8573=12\times 3+(4^5+6\times 7)\times 8+9.$
\item [] $8574=12^3+4^5\times 6+78\times 9.$
\item [] $8575=1+2\times 3\times (4\times 5\times 67+89).$
\item [] $8576=123\times (4+5\times (6+7))+89.$
\item [] $8577=12^3+456\times (7+8)+9.$
\item [] $8578=1^2+3\times (45+6)\times 7\times 8+9.$
\item [] $8579=1\times 2+3\times (45+6)\times 7\times 8+9.$
\item [] $8580=12\times (34\times 5+67\times 8+9).$
\item [] $8581=1^2+3\times 4\times (5+6)\times (7\times 8+9).$
\item [] $8582=1\times 2^3\times 4^5+6\times (7\times 8+9).$
\item [] $8583=1+2^3\times 4^5+6\times (7\times 8+9).$
\item [] $8584=(1+2\times 3\times 4)\times (5\times 67+8)+9.$
\item [] $8585=12^3+4+(5+6)\times 7\times 89.$
\item [] $8586=(1+2^3+4+5)\times (6\times 78+9).$
\item [] $8587=1+(2\times (3^4+56\times 7)+8)\times 9.$
\item [] $\mathit{8588=1-23\times 4+(5+6)\times 789.}$
\item [] $8589=12+3\times (45+6)\times 7\times 8+9.$
\item [] $8590=1+(2\times 3+4)\times (5+6)\times 78+9.$
\item [] $8591=(1+2\times 3\times 4\times 5)\times (6+7\times 8+9).$
\item [] $8592=12+3\times 4\times (5+6)\times (7\times 8+9).$
\item [] $8593=(12\times 3+4^5+6+7)\times 8+9.$
\item [] $8594=(1+2+3\times 4)\times 567+89.$
\item [] $8595=((1\times 2\times 3)^4+5)\times 6+789.$
\item [] $8596=1\times 2\times (3^4+5+6\times 78\times 9).$
\item [] $8597=1+2\times (3^4+5+6\times 78\times 9).$
\item [] $8598=12\times 3^4\times 5+6\times 7\times 89.$
\item [] $\mathit{8599=(1+2)\times 3\times 4^5+6-7\times 89.}$
\item [] $8600=1\times 2^3\times (4^5+6\times 7)+8\times 9.$
\item[]$\mbox{Decreasing order}$
\item [] $8531=(98+7+65\times 4^3)\times 2+1.$
\item [] $8532=9\times 8\times 7+6^5+4\times 3\times 21.$
\item [] $8533=9\times 876+5^4+3+21.$
\item [] $8534=98\times 76+543\times 2\times 1.$
\item [] $8535=98\times 76+543\times 2+1.$
\item [] $8536=9\times 876+5^4+3^{(2+1)}.$
\item [] $8537=((9+8\times 7)\times 65+43)\times 2+1.$
\item [] $8538=((9\times 8+7)\times 6\times (5+4)+3)\times 2\times 1.$
\item [] $8539=((9\times 8+7)\times 6\times (5+4)+3)\times 2+1.$
\item [] $8540=(98+7\times 6)\times (54+3\times 2+1).$
\item [] $8541=9\times 876+5^4+32\times 1.$
\item [] $8542=9\times 876+5^4+32+1.$
\item [] $\mathit{8543=-9-8\times 7-6^5+4^(3\times 2+1).}$
\item [] $8544=98+(7+6)\times 5^4+321.$
\item [] $8545=(9+8+7)\times (6\times 54+32)+1.$
\item [] $8546=(9+8\times 7)\times 65+4321.$
\item [] $8547=98\times 7+6^5+4^3+21.$
\item [] $8548=98\times 7+6^5+43\times 2\times 1.$
\item [] $8549=98\times 7+6^5+43\times 2+1.$
\item [] $8550=9\times (8+7\times 65+4\times 3)\times 2\times 1.$
\item [] $8551=9\times (8+7\times 65+4\times 3)\times 2+1.$
\item [] $8552=9+87\times 6+(5\times 4)^3+21.$
\item [] $8553=9\times (8\times 7+65\times 4)\times 3+21.$
\item [] $\mathit{8554=9\times 87+6^5-4-3+2\times 1.}$
\item [] $8555=(9+8+7\times 6)\times ((5+4+3)^2+1).$
\item [] $8556=(98\times 7\times 6+54\times 3)\times 2\times 1.$
\item [] $8557=(98\times 7\times 6+54\times 3)\times 2+1.$
\item [] $8558=98\times 7+6^5+4\times (3+21).$
\item [] $8559=9\times 87+6\times 54\times (3+21).$
\item [] $8560=9\times (8\times 7+6)+(5\times 4)^3+2\times 1.$
\item [] $8561=9+8+(7+65\times 4)\times 32\times 1.$
\item [] $8562=9+8+(7+65\times 4)\times 32+1.$
\item [] $\mathit{8563=9\times 87+6^5+4+3-2-1.}$
\item [] $\mathit{8564=9\times 87+6^5+4+3-2\times 1.}$
\item [] $8565=9+(87+6)\times (5+43\times 2+1).$
\item [] $8566=9+(8\times 7+6)\times (5+4^3)\times 2+1.$
\item [] $\mathit{8567=9\times 87+6^5+4+3+2-1.}$
\item [] $8568=9\times 87+6^5+4+3+2\times 1.$
\item [] $8569=9\times 87+6^5+4+3+2+1.$
\item [] $8570=9\times 87+6^5+4+3\times 2+1.$
\item [] $8571=9+8+(7+6)\times (5^4+32+1).$
\item [] $8572=9\times 876+5^4+3\times 21.$
\item [] $8573=9\times 87+6^5+4+3^2+1.$
\item [] $8574=9\times 87+6^5+4\times 3+2+1.$
\item [] $8575=(98+7\times (6+5))\times (4+3)^2\times 1.$
\item [] $8576=(98+7\times (6+5))\times (4+3)^2+1.$
\item [] $8577=9\times (87\times 6+5\times 43\times 2+1).$
\item [] $8578=(9+87)\times 6+(5\times 4)^3+2\times 1.$
\item [] $8579=9\times 87+6^5+4\times (3+2)\times 1.$
\item [] $8580=9+8\times 7+65\times 43\times (2+1).$
\item [] $8581=(9+8\times 7+65)\times (4^3+2)+1.$
\item [] $8582=(9+8\times 7)\times (6+5)\times 4\times 3+2\times 1.$
\item [] $8583=9\times 87+6^5+4\times 3\times 2\times 1.$
\item [] $8584=9\times 87+6^5+4\times 3\times 2+1.$
\item [] $8585=9+(8\times 7+6+5)\times 4\times 32\times 1.$
\item [] $8586=9\times (876+54+3+21).$
\item [] $8587=9\times 87+6^5+4+3+21.$
\item [] $\mathit{8588=9\times 87+6^5-4+32+1.}$
\item [] $8589=(9+8+76\times 5+4\times 3)\times 21.$
\item [] $8590=98\times 7+6^5+4\times 32\times 1.$
\item [] $8591=98\times 7+6^5+4^3\times 2+1.$
\item [] $8592=9\times 87+6^5+4\times 3+21.$
\item [] $\mathit{8593=9\times 87\times (6+5)+4-3-21.}$
\item [] $8594=98\times 7+6^5+4\times (32+1).$
\item [] $8595=9\times 87+6^5+4+32\times 1.$
\item [] $8596=9\times 87+6^5+4+32+1.$
\item [] $8597=9\times 8\times 76+5^4\times (3+2)\times 1.$
\item [] $8598=9\times 8\times 76+5^4\times (3+2)+1.$
\item [] $8599=9\times 87+6^5+4\times (3^2+1).$
\item [] $8600=(9+(8\times 76+5)\times (4+3))\times 2\times 1.$
\item[]$\mbox{Increasing order}$
\item [] $8601=12^3\times 4+5\times 6\times 7\times 8+9.$
\item [] $8602=1\times 23\times (4+5+6+7)\times (8+9).$
\item [] $8603=1+23\times (4+5+6+7)\times (8+9).$
\item [] $8604=(1^2+3)\times (4\times 56+7+8)\times 9.$
\item [] $8605=(1\times 2+3)\times ((4+5\times 6\times 7)\times 8+9).$
\item [] $8606=1+(2+3)\times ((4+5\times 6\times 7)\times 8+9).$
\item [] $8607=12+3\times ((45+6)\times 7\times 8+9).$
\item [] $\mathit{8608=(12^3-4)\times 5-6-7-8+9.}$
\item [] $8609=(1+2^3+4^5+6\times 7)\times 8+9.$
\item [] $8610=1+(2+34\times 5)\times (6\times 7+8)+9.$
\item [] $8611=(1+2\times (3+45+6))\times (7+8\times 9).$
\item [] $\mathit{8612=1-2\times 34+(5+6)\times 789.}$
\item [] $8613=12\times (34+5+678)+9.$
\item [] $8614=1+(2+3+4)\times (5+6)\times (78+9).$
\item [] $\mathit{8615=1\times 2\times (3+4)\times (5+6)\times 7\times 8-9.}$
\item [] $8616=(1+23)\times (4\times 56+(7+8)\times 9).$
\item [] $8617=1\times 2^3\times (4^5+6\times 7)+89.$
\item [] $8618=1\times (2+3\times 4\times 5)\times (67+8\times 9).$
\item [] $8619=123\times (4\times 5+6\times 7+8)+9.$
\item [] $8620=1^2+(34+5)\times (6+7)\times (8+9).$
\item [] $8621=1\times 2+(34+5)\times (6+7)\times (8+9).$
\item [] $8622=1+2+(34+5)\times (6+7)\times (8+9).$
\item [] $\mathit{8623=1+23\times (4+5)\times 6\times 7-8\times 9.}$
\item [] $\mathit{8624=1+2^{(3\times 4)}+567\times 8-9.}$
\item [] $8625=(1\times 23\times 45+6\times 7)\times 8+9.$
\item [] $8626=1+(2+3+4\times 5)\times (6\times 7\times 8+9).$
\item [] $\mathit{8627=1-2+3\times 4\times (5+6\times 7\times (8+9)).}$
\item [] $8628=123+45\times (6+7+8)\times 9.$
\item [] $8629=1^2+3\times 4\times (5+6\times 7\times (8+9)).$
\item [] $8630=1\times 2+3\times 4\times (5+6\times 7\times (8+9)).$
\item [] $8631=12+(34+5)\times (6+7)\times (8+9).$
\item [] $8632=1\times 2^3\times (456+7\times 89).$
\item [] $8633=1+2^3\times (456+7\times 89).$
\item [] $8634=(1+2+3\times 4)\times (567+8)+9.$
\item [] $8635=1\times 2+(34+56+7)\times 89.$
\item [] $8636=((1+2)^3\times 45+6)\times 7+89.$
\item [] $8637=12^3\times 4+5\times (6\times 7\times 8+9).$
\item [] $\mathit{8638=-1\times 2+3\times 4\times 5\times 6\times (7+8+9).}$
\item [] $\mathit{8639=-12\times 3-4+(5+6)\times 789.}$
\item [] $8640=12\times 3\times (4+5\times 6)\times 7+8\times 9.$
\item [] $8641=1^2+3\times 4\times 5\times 6\times (7+8+9).$
\item [] $8642=1\times 2+3\times 4\times 5\times 6\times (7+8+9).$
\item [] $8643=1+2+3\times 4\times 5\times 6\times (7+8+9).$
\item [] $\mathit{8644=-12^3\times 4+5^6-78+9.}$
\item [] $8645=(12+3+4)\times (5+(6\times 7+8)\times 9).$
\item [] $8646=1+(2+3)\times (4+5\times (6\times 7\times 8+9)).$
\item [] $\mathit{8647=1\times 2-34+(5+6)\times 789.}$
\item [] $8648=(1+2)\times (3^4\times 5+6)\times 7+8+9.$
\item [] $8649=12\times (3+45+6\times 7)\times 8+9.$
\item [] $8650=1+2\times (3+45)\times 6\times (7+8)+9.$
\item [] $\mathit{8651=1\times 2^3\times 4^5+6\times 78-9.}$
\item [] $8652=12+3\times 4\times 5\times 6\times (7+8+9).$
\item [] $8653=1^2+3\times (4+5\times 6\times (7+89)).$
\item [] $8654=(12+3)\times (4+567)+89.$
\item [] $8655=1+2+3\times (4+5\times 6\times (7+89)).$
\item [] $\mathit{8656=1-2\times 3\times 4+(5+6)\times 789.}$
\item [] $8657=12\times 3\times (4+5\times 6)\times 7+89.$
\item [] $8658=1\times 2\times (3^4+56\times 7+8)\times 9.$
\item [] $8659=1+2\times (3^4+56\times 7+8)\times 9.$
\item [] $8660=123+(4^5+6\times 7)\times 8+9.$
\item [] $8661=12^3+4^5\times 6+789.$
\item [] $\mathit{8662=(12^3+4)\times 5-6+7-8+9.}$
\item [] $\mathit{8663=(1-2-3)\times 4+(5+6)\times 789.}$
\item [] $8664=(1^2+3\times 4\times 5\times 6)\times (7+8+9).$
\item [] $8665=(1^2+3)\times 4\times (5+67\times 8)+9.$
\item [] $8666=1\times 2^3\times 4^5+6\times (7+8\times 9).$
\item [] $8667=1+2^3\times 4^5+6\times (7+8\times 9).$
\item [] $8668=1^2+3^4\times (5+6+7+89).$
\item [] $8669=1\times 2^3\times 4^5+6\times 78+9.$
\item [] $8670=1+2^3\times 4^5+6\times 78+9.$
\item[]$\mbox{Decreasing order}$
\item [] $8601=(9+8\times 7)\times (6+5)\times 4\times 3+21.$
\item [] $8602=(9+8+765)\times (4+3\times 2+1).$
\item [] $8603=(9\times (8+7)+65)\times 43+2+1.$
\item [] $8604=9\times 87+6^5+43+2\times 1.$
\item [] $8605=9\times 8\times 7+6^5+4+321.$
\item [] $8606=98\times 7+6^5+(4\times 3)^2\times 1.$
\item [] $8607=98\times 7+6^5+(4\times 3)^2+1.$
\item [] $8608=9\times 87+6^5+(4+3)^2\times 1.$
\item [] $8609=9\times 87+6^5+(4+3)^2+1.$
\item [] $8610=(9\times (8+7)\times 6+5^4)\times 3\times 2\times 1.$
\item [] $8611=(9\times (8+7)\times 6+5^4)\times 3\times 2+1.$
\item [] $\mathit{8612=9\times 87\times (6+5)+4-3-2\times 1.}$
\item [] $8613=(98\times 7+6\times 5)\times 4\times 3+21.$
\item [] $8614=9\times 876+((5+4)\times 3)^2+1.$
\item [] $8615=9\times 876+(5+4)^3+2\times 1.$
\item [] $8616=9\times 876+(5+4)^3+2+1.$
\item [] $8617=9\times 8+(7+65\times 4)\times 32+1.$
\item [] $\mathit{8618=9\times 87+6^5-4+3\times 21.}$
\item [] $8619=9\times 876+5\times (4+3)\times 21.$
\item [] $8620=9+8\times 76+(5\times 4)^3+2+1.$
\item [] $8621=(9\times (8+7)+65)\times 43+21.$
\item [] $8622=(9\times 87+654)\times 3\times 2\times 1.$
\item [] $8623=9\times 87+6^5+43+21.$
\item [] $8624=9\times 87\times (6+5)+4+3\times 2+1.$
\item [] $8625=9\times 87+6^5+4^3+2\times 1.$
\item [] $8626=9\times 87+6^5+4+3\times 21.$
\item [] $8627=9\times 87\times (6+5)+4\times 3+2\times 1.$
\item [] $8628=9\times 87\times (6+5)+4\times 3+2+1.$
\item [] $\mathit{8629=-9+8765-4\times 32+1.}$
\item [] $8630=(98+765)\times (4+3+2+1).$
\item [] $8631=98\times 7\times 6+5\times 43\times 21.$
\item [] $8632=(98+7)\times 6+(5\times 4)^3+2\times 1.$
\item [] $8633=(98+7)\times 6+(5\times 4)^3+2+1.$
\item [] $8634=9\times 876+(5+4)^3+21.$
\item [] $8635=9+8+7\times (6+(5\times (4+3))^2)+1.$
\item [] $8636=(9+8)\times (76\times 5+4\times 32\times 1).$
\item [] $8637=9\times 87\times (6+5)+4\times 3\times 2\times 1.$
\item [] $8638=9+8\times 76+(5\times 4)^3+21.$
\item [] $8639=98+(7+6)\times (5^4+32\times 1).$
\item [] $8640=9\times (87\times 6+5+432+1).$
\item [] $8641=9\times 87\times (6+5)+4+3+21.$
\item [] $8642=98+(7+65\times 4)\times 32\times 1.$
\item [] $8643=98+(7+65\times 4)\times 32+1.$
\item [] $8644=9\times 87+6^5+4^3+21.$
\item [] $8645=9\times 87+6^5+43\times 2\times 1.$
\item [] $8646=9\times 87+6^5+43\times 2+1.$
\item [] $\mathit{8647=9+8765-4\times 32+1.}$
\item [] $\mathit{8648=9\times (8+7)\times 65-4\times 32+1.}$
\item [] $8649=9\times 8\times 7\times 6+5^4\times 3^2\times 1.$
\item [] $8650=9\times 87\times (6+5)+4+32+1.$
\item [] $8651=98\times (76+5+4)+321.$
\item [] $8652=9\times 8+(7+6)\times 5\times 4\times (32+1).$
\item [] $8653=(9+8)\times (76\times 5+4\times 32+1).$
\item [] $8654=98\times 7+6^5+4^3\times (2+1).$
\item [] $8655=9\times 87+6^5+4\times (3+21).$
\item [] $8656=9\times (8\times 7\times 6+5^4)+3\times 2+1.$
\item [] $\mathit{8657=9\times 87\times (6+5)+43+2-1.}$
\item [] $8658=9\times (876+54+32\times 1).$
\item [] $8659=9\times 87\times (6+5)+43+2+1.$
\item [] $8660=9\times 8+76\times (5+4\times 3^{(2+1)}).$
\item [] $8661=9\times (8+7+65)\times 4\times 3+21.$
\item [] $8662=9\times 87\times (6+5)+(4+3)^2\times 1.$
\item [] $8663=(98\times 7\times 6+5\times 43)\times 2+1.$
\item [] $8664=(9+8+7)\times (6\times (54+3\times 2)+1).$
\item [] $\mathit{8665=9-8+76\times (54+3)\times 2\times 1.}$
\item [] $8666=9+8+(76+5+4\times 3)^2\times 1.$
\item [] $8667=9\times (876+54+32+1).$
\item [] $\mathit{8668=-9\times 8-7+6\times (5+4)^3\times 2-1.}$
\item [] $\mathit{8669=-9+8765-43\times 2-1.}$
\item [] $8670=(98+765+4)\times (3^2+1).$
\item[]$\mbox{Increasing order}$
\item [] $8671=1+(2^3+4+5)\times (6+7\times 8\times 9).$
\item [] $8672=1+23\times (4\times (5+67)+89).$
\item [] $8673=(12+3\times 45)\times (6\times 7+8+9).$
\item [] $\mathit{8674=-12+3+4+(5+6)\times 789.}$
\item [] $8675=1\times 234\times (5\times 6+7)+8+9.$
\item [] $8676=1+234\times (5\times 6+7)+8+9.$
\item [] $8677=1+2\times (3^4+(5+6\times 78)\times 9).$
\item [] $\mathit{8678=1-2\times 3+4+(5+6)\times 789.}$
\item [] $8679=12+3^4\times (5+6+7+89).$
\item [] $8680=1^{234}+(5+6)\times 789.$
\item [] $\mathit{8681=1+2+3-4+(5+6)\times 789.}$
\item [] $8682=(12\times 34+5)\times (6+7+8)+9.$
\item [] $8683=1^{23}\times 4+(5+6)\times 789.$
\item [] $8684=1^{23}+4+(5+6)\times 789.$
\item [] $8685=(123\times 4+5+6\times 78)\times 9.$
\item [] $8686=1^2\times 3+4+(5+6)\times 789.$
\item [] $8687=12\times 3\times 4\times 56+7\times 89.$
\item [] $8688=1\times 2+3+4+(5+6)\times 789.$
\item [] $8689=1+2+3+4+(5+6)\times 789.$
\item [] $8690=1+2\times 3+4+(5+6)\times 789.$
\item [] $8691=1^2\times 3\times 4+(5+6)\times 789.$
\item [] $8692=(1+2)\times 3+4+(5+6)\times 789.$
\item [] $8693=1\times 2+3\times 4+(5+6)\times 789.$
\item [] $8694=12\times 3^4+(5+6)\times 78\times 9.$
\item [] $8695=1+2\times (3\times 45+6\times 78\times 9).$
\item [] $\mathit{8696=1+23\times (4+5)\times 6\times 7-8+9.}$
\item [] $8697=(1+(2+3\times (45+6))\times 7)\times 8+9.$
\item [] $8698=12+3+4+(5+6)\times 789.$
\item [] $8699=1\times (2+3)\times 4+(5+6)\times 789.$
\item [] $8700=12\times (3\times 4+5+6+78\times 9).$
\item [] $\mathit{8701=-1+2^3\times 4^5+6+7\times 8\times 9.}$
\item [] $8702=1\times 2^3\times 4^5+6+7\times 8\times 9.$
\item [] $8703=12+3\times 4+(5+6)\times 789.$
\item [] $8704=1+2\times 3\times 4+(5+6)\times 789.$
\item [] $8705=(12\times (34+56)+7)\times 8+9.$
\item [] $8706=1\times 23+4+(5+6)\times 789.$
\item [] $8707=1+23+4+(5+6)\times 789.$
\item [] $\mathit{8708=(12^3-4)\times 5+6-7+89.}$
\item [] $8709=(12+3)\times (4+(5+67)\times 8)+9.$
\item [] $8710=(1+2)^3+4+(5+6)\times 789.$
\item [] $8711=1\times 23\times (4+5)\times 6\times 7+8+9.$
\item [] $8712=12\times (34+5+678+9).$
\item [] $8713=1^2\times 34+(5+6)\times 789.$
\item [] $8714=1^2+34+(5+6)\times 789.$
\item [] $8715=12^3+4^5+67\times 89.$
\item [] $8716=1+2+34+(5+6)\times 789.$
\item [] $\mathit{8717=1\times 23\times (4+56\times 7-8-9).}$
\item [] $8718=(1+23\times (4+5)\times 6)\times 7+8+9.$
\item [] $8719=12\times 3+4+(5+6)\times 789.$
\item [] $8720=(1+2\times 3+4+5)\times (67\times 8+9).$
\item [] $8721=12\times (3+45+678)+9.$
\item [] $8722=(1^2+34+56+7)\times 89.$
\item [] $8723=1+2\times (3\times 4+5\times 6+7)\times 89.$
\item [] $8724=(12+3)\times (45+67\times 8)+9.$
\item [] $8725=12+34+(5+6)\times 789.$
\item [] $8726=1\times 23\times (4+5\times (67+8))+9.$
\item [] $8727=1+23\times (4+5\times (67+8))+9.$
\item [] $\mathit{8728=1+2\times 3\times (4\times 5+6)\times 7\times 8-9.}$
\item [] $8729=(1+23+4^5+6\times 7)\times 8+9.$
\item [] $8730=1\times 234\times (5\times 6+7)+8\times 9.$
\item [] $8731=1+234\times (5\times 6+7)+8\times 9.$
\item [] $8732=(1234+5+6)\times 7+8+9.$
\item [] $8733=12+(3^4\times (5+6)+78)\times 9.$
\item [] $\mathit{8734=(1-2)\times 3+(4^5+67)\times 8+9.}$
\item [] $\mathit{8735=(12^3+4)\times 5+6+78-9.}$
\item [] $8736=(12\times 3+4\times 5)\times (67+89).$
\item [] $8737=1\times 2^3\times 4^5+67\times 8+9.$
\item [] $8738=1+2^3\times 4^5+67\times 8+9.$
\item [] $8739=(123+4\times 5\times 6\times 7+8)\times 9.$
\item [] $8740=1\times 23\times 4\times (5\times 6+7\times 8+9).$
\item[]$\mbox{Decreasing order}$
\item [] $8671=(9+8)\times (76+5+4)\times 3\times 2+1.$
\item [] $8672=(9\times 8\times 7+6)\times (5+4\times 3)+2\times 1.$
\item [] $8673=(9+8+76+5\times 4^3)\times 21.$
\item [] $8674=(98+7)\times 65+43^2\times 1.$
\item [] $8675=(98+7)\times 65+43^2+1.$
\item [] $8676=9+87+65\times 4\times (32+1).$
\item [] $8677=9\times 87\times (6+5)+43+21.$
\item [] $8678=98+(7+6)\times 5\times 4\times (32+1).$
\item [] $8679=9\times 87\times (6+5)+4^3+2\times 1.$
\item [] $8680=9\times 87\times (6+5)+4^3+2+1.$
\item [] $8681=9+8+76\times (54+3)\times 2\times 1.$
\item [] $8682=9+8+76\times (54+3)\times 2+1.$
\item [] $\mathit{8683=-9\times 8+7+6\times (5+4)^3\times 2\times 1.}$
\item [] $\mathit{8684=(9+8)\times 765-4321.}$
\item [] $8685=98+7+65\times 4\times (32+1).$
\item [] $8686=98+76\times (5+4\times 3^{(2+1)}).$
\item [] $8687=9\times 87+6^5+4\times 32\times 1.$
\item [] $8688=9\times 87+6^5+4\times 32+1.$
\item [] $8689=9+(8\times 7+6)\times 5\times (4+3+21).$
\item [] $8690=(9\times 8+7)\times (65+43+2\times 1).$
\item [] $8691=9\times (876+54)+321.$
\item [] $\mathit{8692=-9+876\times 5+4321.}$
\item [] $\mathit{8693=-9+8765-4^3+2-1.}$
\item [] $8694=9\times (876+5+4^3+21).$
\item [] $8695=98\times 7+6+(5\times 4)^3+2+1.$
\item [] $8696=9\times 8+7\times (6+(5\times (4+3))^2+1).$
\item [] $8697=9+8\times (7\times 6+5\times 4^3)\times (2+1).$
\item [] $8698=9\times 87\times (6+5)+4^3+21.$
\item [] $8699=9\times 87\times (6+5)+43\times 2\times 1.$
\item [] $8700=(98+7+65\times 43)\times (2+1).$
\item [] $8701=(98+76)\times (5+43+2)+1.$
\item [] $8702=98\times (7\times 6+5)+4^{(3+2+1)}.$
\item [] $8703=9+8+7+6^5+43\times 21.$
\item [] $8704=98\times 76+(5^4+3)\times 2\times 1.$
\item [] $8705=98\times 76+(5^4+3)\times 2+1.$
\item [] $8706=9\times 87+6^5+(4+3)\times 21.$
\item [] $\mathit{8707=9+8765-4-3\times 21.}$
\item [] $8708=98\times 76+5\times 4\times 3\times 21.$
\item [] $8709=9\times 87\times (6+5)+4\times (3+21).$
\item [] $8710=9+876\times 5+4321.$
\item [] $8711=9+8+7\times 6\times (5+4^3)\times (2+1).$
\item [] $8712=9\times 8\times 7+6^5+432\times 1.$
\item [] $8713=9\times 8\times 7+6^5+432+1.$
\item [] $8714=98\times 7+6^5+4\times 3\times 21.$
\item [] $8715=9\times (8+7)+65\times 4\times (32+1).$
\item [] $8716=(9+8)\times 7\times 6+(5\times 4)^3+2\times 1.$
\item [] $8717=(9+8)\times 7\times 6+(5\times 4)^3+2+1.$
\item [] $8718=9+8+7+6\times (5+4^3)\times 21.$
\item [] $8719=(9\times (8\times 7+65)\times 4+3)\times 2+1.$
\item [] $\mathit{8720=-9+8765-4-32\times 1.}$
\item [] $8721=987+6\times (5+4\times 321).$
\item [] $8722=9\times 8\times 76+(54+3)^2+1.$
\item [] $8723=98\times (7+(6+5\times (4+3))\times 2)+1.$
\item [] $8724=(987+(6+5+4)^3)\times 2\times 1.$
\item [] $8725=(987+(6+5+4)^3)\times 2+1.$
\item [] $\mathit{8726=987+6^5-4-32-1.}$
\item [] $8727=9+(8+7\times 6\times (5+4^3))\times (2+1).$
\item [] $\mathit{8728=9+8765-43-2-1.}$
\item [] $8729=9+8\times (765+4+321).$
\item [] $8730=(9+87\times (6+5)+4)\times 3^2\times 1.$
\item [] $8731=(9+87\times (6+5)+4)\times 3^2+1.$
\item [] $8732=(9\times 8+76)\times (54+3+2\times 1).$
\item [] $8733=(9\times 8+76)\times (54+3+2)+1.$
\item [] $8734=(9\times 87+6+5)\times (4+3\times 2+1).$
\item [] $8735=(9+8)\times 7\times 6+(5\times 4)^3+21.$
\item [] $8736=9\times 8+76\times (54+3)\times 2\times 1.$
\item [] $8737=98\times 76+5+4\times 321.$
\item [] $8738=(9+8)\times (76+5+432+1).$
\item [] $8739=9\times (87\times (6+5)+4\times 3+2\times 1).$
\item [] $8740=((98+76)\times 5+4)\times (3^2+1).$
\item[]$\mbox{Increasing order}$
\item [] $8741=1+23\times 4\times (5\times 6+7\times 8+9).$
\item [] $8742=123\times (4+5+6+7\times 8)+9.$
\item [] $8743=1+2+3+(4^5+67)\times 8+9.$
\item [] $8744=(12^3+4)\times 5+67+8+9.$
\item [] $8745=(12^3+4)\times 5+6+7+8\times 9.$
\item [] $8746=1+2^3+(4^5+67)\times 8+9.$
\item [] $8747=1\times 234\times (5\times 6+7)+89.$
\item [] $8748=12^3+45\times (67+89).$
\item [] $8749=1+2\times 3^4\times (5\times 6+7+8+9).$
\item [] $8750=1\times 2^3\times 4^5+(6+7\times 8)\times 9.$
\item [] $8751=1+2^3\times 4^5+(6+7\times 8)\times 9.$
\item [] $8752=12+3+(4^5+67)\times 8+9.$
\item [] $8753=(12^3+4)\times 5+6+78+9.$
\item [] $8754=(1+(23\times 4+5)\times 6)\times (7+8)+9.$
\item [] $8755=1+2\times (3\times (4\times 5+6)\times 7\times 8+9).$
\item [] $\mathit{8756=1\times 2\times (-3+4+56\times 78+9).}$
\item [] $8757=1^2\times 3^4\times (5\times 6+78)+9.$
\item [] $8758=1^2+3^4\times (5\times 6+78)+9.$
\item [] $8759=1\times 2+3^4\times (5\times 6+78)+9.$
\item [] $8760=1^2\times 3^4+(5+6)\times 789.$
\item [] $8761=1^2+3^4+(5+6)\times 789.$
\item [] $8762=1\times 2+3^4+(5+6)\times 789.$
\item [] $8763=12\times (3+4)+(5+6)\times 789.$
\item [] $8764=1\times 2\times (34\times 5+6\times 78\times 9).$
\item [] $8765=1+2\times (34\times 5+6\times 78\times 9).$
\item [] $8766=12\times 3\times 4\times 56+78\times 9.$
\item [] $8767=1+23\times (4+5)\times 6\times 7+8\times 9.$
\item [] $8768=1\times 2\times (3+4+56\times 78+9).$
\item [] $8769=12+3^4\times (5\times 6+78)+9.$
\item [] $8770=1+2\times (3\times 4+56\times 78)+9.$
\item [] $8771=1\times 23\times 4+(5+6)\times 789.$
\item [] $8772=12+3^4+(5+6)\times 789.$
\item [] $8773=12\times 3+(4^5+67)\times 8+9.$
\item [] $8774=(12^3+4)\times 5+6\times 7+8\times 9.$
\item [] $8775=1+2+(3^4+5)\times (6+7+89).$
\item [] $8776=1+(2+3+4+56)\times (7+8)\times 9.$
\item [] $8777=(1\times 2+3+4^5+67)\times 8+9.$
\item [] $8778=1\times 2\times (3\times 4+56\times 78+9).$
\item [] $8779=1+2\times (3\times 4+56\times 78+9).$
\item [] $\mathit{8780=1\times 2+3\times (45\times 67-89).}$
\item [] $8781=(1+2)\times 34+(5+6)\times 789.$
\item [] $\mathit{8782=-12^3\times 4+5^6+78-9.}$
\item [] $8783=1\times 23\times (4+5)\times 6\times 7+89.$
\item [] $8784=1+23\times (4+5)\times 6\times 7+89.$
\item [] $8785=1^2+3\times 4\times (5\times 6+78\times 9).$
\item [] $8786=1\times 2+3\times 4\times (5\times 6+78\times 9).$
\item [] $8787=(1234+5+6)\times 7+8\times 9.$
\item [] $\mathit{8788=-1+23+(4^5-6\times 7-8)\times 9.}$
\item [] $8789=(1+2\times 3+4)\times (5+6\times 7)\times (8+9).$
\item [] $8790=(1+23\times (4+5)\times 6)\times 7+89.$
\item [] $8791=(12^3+4)\times 5+6\times 7+89.$
\item [] $8792=(12+3)\times 45\times (6+7)+8+9.$
\item [] $8793=(1+2\times 3+4^5+67)\times 8+9.$
\item [] $8794=1234+56\times (7+8)\times 9.$
\item [] $\mathit{8795=-12^3\times 4+5^6-7+89.}$
\item [] $8796=12+3\times 4\times (5\times 6+78\times 9).$
\item [] $\mathit{8797=1^2-3\times 4\times (56-789).}$
\item [] $\mathit{8798=123-4+(5+6)\times 789.}$
\item [] $8799=(12^3+4)\times 5+67+8\times 9.$
\item [] $8800=(1+2\times 3+4)\times (5+6+789).$
\item [] $8801=(123+4+5\times 6)\times 7\times 8+9.$
\item [] $8802=1^2\times 3\times (45\times 6+7\times 8)\times 9.$
\item [] $8803=1^2+3\times (45\times 6+7\times 8)\times 9.$
\item [] $8804=(1234+5+6)\times 7+89.$
\item [] $8805=1+2+3\times (45\times 6+7\times 8)\times 9.$
\item [] $8806=123+4+(5+6)\times 789.$
\item [] $8807=1+(23+45+6)\times 7\times (8+9).$
\item [] $8808=12\times (3^4+5\times 6+7\times 89).$
\item [] $8809=(1+2^3+4^5+67)\times 8+9.$
\item [] $\mathit{8810=1+2^3\times 4^5-6+7\times 89.}$
\item[]$\mbox{Decreasing order}$
\item [] $8741=9\times 87\times (6+5)+4\times 32\times 1.$
\item [] $8742=9\times 87\times (6+5)+4\times 32+1.$
\item [] $8743=98\times (7\times (6+5)+4\times 3)+21.$
\item [] $8744=9+8\times 7+6^5+43\times 21.$
\item [] $8745=9\times 8\times (7+6)\times (5+4)+321.$
\item [] $8746=(9+876)\times 5+4321.$
\item [] $8747=98+(76+5+4\times 3)^2\times 1.$
\item [] $8748=9\times (87\times (6+5)+(4+3)\times 2+1).$
\item [] $8749=9\times (876+(5+43)\times 2)+1.$
\item [] $8750=(9+8\times 7\times 6+5)\times (4\times 3\times 2+1).$
\item [] $8751=9\times 87+6^5+4^3\times (2+1).$
\item [] $\mathit{8752=9+8765-43+21.}$
\item [] $8753=9+8\times (7\times 6\times (5\times 4+3\times 2)+1).$
\item [] $8754=9\times 8+7\times 6+5\times (4\times 3)^{(2+1)}.$
\item [] $8755=(9+8)\times (7\times (65+4)+32\times 1).$
\item [] $8756=(9+8)\times (7\times (65+4)+32)+1.$
\item [] $8757=(98+7+6\times 5+4)\times 3\times 21.$
\item [] $8758=9\times 8+7+6^5+43\times 21.$
\item [] $8759=9+8\times 7+6\times (5+4^3)\times 21.$
\item [] $8760=(98+7+65\times 4)\times (3+21).$
\item [] $\mathit{8761=9+8765-4\times 3-2+1.}$
\item [] $8762=98+76\times (54+3)\times 2\times 1.$
\item [] $8763=987+6\times 54\times (3+21).$
\item [] $\mathit{8764=9+8765-4-3-2-1.}$
\item [] $8765=9+8+(76+5)\times 4\times 3^{(2+1)}.$
\item [] $8766=9\times 8+7\times 6\times (5+4^3)\times (2+1).$
\item [] $8767=9\times ((8\times 7+65)\times 4+3)\times 2+1.$
\item [] $8768=(9\times 8+(7+6)\times 5)\times (43+21).$
\item [] $8769=(9\times 8+7)\times (65+43+2+1).$
\item [] $8770=(9\times 8+(7+6)\times 5)\times 4^3+2\times 1.$
\item [] $8771=(98+76+5)\times (4+3)^2\times 1.$
\item [] $8772=987+6^5+4+3+2\times 1.$
\item [] $8773=987+6^5+4+3+2+1.$
\item [] $8774=987+6^5+4+3\times 2+1.$
\item [] $8775=9+87+6^5+43\times 21.$
\item [] $8776=987+6^5+4+3^2\times 1.$
\item [] $8777=987+6^5+4+3^2+1.$
\item [] $8778=987+6^5+4\times 3+2+1.$
\item [] $8779=9+8+(7+6\times (5+4)^3)\times 2\times 1.$
\item [] $8780=98+7\times 6+5\times (4\times 3)^{(2+1)}.$
\item [] $8781=9\times (8+7)+6+5\times (4\times 3)^{(2+1)}.$
\item [] $\mathit{8782=9+8765+4+3+2-1.}$
\item [] $8783=9+8765+4+3+2\times 1.$
\item [] $8784=98+7+6^5+43\times 21.$
\item [] $8785=9+8765+4+3\times 2+1.$
\item [] $8786=9\times (8+7)\times 65+4+3\times 2+1.$
\item [] $8787=98\times 7+6^5+4+321.$
\item [] $8788=987+6^5+4\times 3\times 2+1.$
\item [] $8789=9+8765+4\times 3+2+1.$
\item [] $8790=9\times (8+7)\times 65+4\times 3+2+1.$
\item [] $8791=987+6^5+4+3+21.$
\item [] $8792=9\times 876+5+43\times 21.$
\item [] $8793=(9+876\times 5+4+3)\times 2+1.$
\item [] $8794=9+8765+4\times (3+2\times 1).$
\item [] $8795=9+8765+4\times (3+2)+1.$
\item [] $8796=987+6^5+4\times 3+21.$
\item [] $8797=98\times 76+5+4^3\times 21.$
\item [] $8798=9+8765+4\times 3\times 2\times 1.$
\item [] $8799=9+8765+4\times 3\times 2+1.$
\item [] $8800=987+6^5+4+32+1.$
\item [] $8801=9+8\times (7+6+543\times 2\times 1).$
\item [] $8802=9+8765+4+3+21.$
\item [] $8803=9\times (8+7)\times 65+4+3+21.$
\item [] $\mathit{8804=987+6^5+43-2\times 1.}$
\item [] $8805=9+8765+4+3^{(2+1)}.$
\item [] $8806=9\times (8+765)+43^2\times 1.$
\item [] $8807=9+8765+4\times 3+21.$
\item [] $8808=987+6^5+43+2\times 1.$
\item [] $8809=987+6^5+43+2+1.$
\item [] $8810=9+8765+4+32\times 1.$
\item[]$\mbox{Increasing order}$
\item [] $8811=(1\times 23+4+5+67)\times 89.$
\item [] $8812=1^2+3\times (4\times 5+6+7)\times 89.$
\item [] $8813=1\times 2\times (34+56\times 78)+9.$
\item [] $8814=12+3\times (45\times 6+7\times 8)\times 9.$
\item [] $8815=1+2\times (34+5)\times ((6+7)\times 8+9).$
\item [] $8816=(12^3+4)\times 5+67+89.$
\item [] $8817=(1+2)\times (3\times 4\times 5\times 6+7)\times 8+9.$
\item [] $8818=12+34\times (5\times (6\times 7+8)+9).$
\item [] $\mathit{8819=-1+((2+3)^4+5)\times (6+7-8+9).}$
\item [] $8820=12\times (3+45+678+9).$
\item [] $8821=1\times 2^3\times 4^5+6+7\times 89.$
\item [] $8822=1+2^3\times 4^5+6+7\times 89.$
\item [] $8823=(1^2+3+4+5)\times 678+9.$
\item [] $8824=1+(2+3^4+5\times 6)\times 78+9.$
\item [] $8825=(12\times 3+4^5+6\times 7)\times 8+9.$
\item [] $8826=1+(23\times 45+67)\times 8+9.$
\item [] $8827=(1+2^3+4)\times (56+7\times 89).$
\item [] $\mathit{8828=1-2-3^4\times (5-6\times 7-8\times 9).}$
\item [] $8829=12\times 34\times 5+6789.$
\item [] $8830=1\times 2^{(3+4+5)}+6\times 789.$
\item [] $8831=1\times 2+3^4\times (5\times 6+7+8\times 9).$
\item [] $8832=1+2+3^4\times (5\times 6+7+8\times 9).$
\item [] $8833=(1+23\times 45+67)\times 8+9.$
\item [] $8834=1\times 2+(3^4+5+6)\times (7+89).$
\item [] $8835=1\times 2^{(3\times 4)}+5+6\times 789.$
\item [] $8836=1+2^{(3\times 4)}+5+6\times 789.$
\item [] $\mathit{8837=12^3+(4^5-6)\times 7-8-9.}$
\item [] $8838=1\times 2\times (3+4\times 5+6\times 78)\times 9.$
\item [] $8839=1^2+3^4\times (5+(6+7)\times 8)+9.$
\item [] $8840=1\times 2\times (3\times 4+56)\times (7\times 8+9).$
\item [] $8841=1\times 2\times 3^4+(5+6)\times 789.$
\item [] $8842=1+2\times 3^4+(5+6)\times 789.$
\item [] $8843=1\times 2\times (3+4^5)+6789.$
\item [] $8844=1+2\times (3+4^5)+6789.$
\item [] $8845=1+2\times 3\times (4^5+(6\times 7+8)\times 9).$
\item [] $8846=1\times 2\times (3+4\times 5\times (6+7)\times (8+9)).$
\item [] $8847=(12\times 34+567+8)\times 9.$
\item [] $8848=(1^2+3^4+5\times 6)\times (7+8\times 9).$
\item [] $8849=(12^3+4)\times 5+(6+7+8)\times 9.$
\item [] $8850=(12+3)\times (45+67\times 8+9).$
\item [] $8851=(1+2)^(3+4)+56\times 7\times (8+9).$
\item [] $\mathit{8852=1+(2\times 3^4+5)\times (6+7\times 8-9).}$
\item [] $8853=12\times 3\times 4\times 56+789.$
\item [] $8854=12^3+4^5+678\times 9.$
\item [] $\mathit{8855=123\times 4\times (5+6+7)+8-9.}$
\item [] $8856=(12\times 3+45+6\times 7)\times 8\times 9.$
\item [] $8857=1+(2\times 3)^4+56\times (7+8)\times 9.$
\item [] $8858=(1\times 2+3)\times 4^5+6\times 7\times 89.$
\item [] $8859=1+(2+3)\times 4^5+6\times 7\times 89.$
\item [] $8860=123+(4^5+67)\times 8+9.$
\item [] $8861=(123+4+5)\times 67+8+9.$
\item [] $8862=12^3\times 4+5\times 6\times (7\times 8+9).$
\item [] $8863=1\times 2^3\times (4^5+6)+7\times 89.$
\item [] $8864=(12+3)\times 45\times (6+7)+89.$
\item [] $8865=(1234+5\times 6)\times 7+8+9.$
\item [] $\mathit{8866=12^3+(4^5+6)\times 7-8\times 9.}$
\item [] $8867=1\times 2^3\times 4^5+(67+8)\times 9.$
\item [] $8868=12+3\times (4+5\times 6+7)\times 8\times 9.$
\item [] $\mathit{8869=(1+(2\times 3)^4-5\times 6)\times 7\times (-8+9).}$
\item [] $\mathit{8870 = (1+(2\times 3)^4-5\times 6)\times 7-8+9.}$
\item [] $\mathit{8871=(1+2^3)\times 4^5-6\times 7\times 8-9.}$
\item [] $8872=1\times 2^3\times (4^5+6+7+8\times 9).$
\item [] $8873=123\times 4\times (5+6+7)+8+9.$
\item [] $8874=12\times 345+6\times 789.$
\item [] $8875=1+2\times 3\times (45+6\times 7)\times (8+9).$
\item [] $8876=1\times 2+3\times (4+5\times 6)\times (78+9).$
\item [] $8877=1+2+3\times (4+5\times 6)\times (78+9).$
\item [] $\mathit{8878=-1+2^3\times 4^5+678+9.}$
\item [] $8879=1\times 2^3\times 4^5+678+9.$
\item [] $8880=1+2^3\times 4^5+678+9.$
\item[]$\mbox{Decreasing order}$
\item [] $8811=9\times 87+6^5+4\times 3\times 21.$
\item [] $8812=9\times (8+7)\times 65+4+32+1.$
\item [] $8813=987+6^5+(4+3)^2+1.$
\item [] $8814=9\times (8+7)+6^5+43\times 21.$
\item [] $8815=9\times (8+7)\times 65+4\times (3^2+1).$
\item [] $8816=9\times ((8+7)\times 65+4)+3+2\times 1.$
\item [] $8817=9\times ((8+7)\times 65+4)+3\times 2\times 1.$
\item [] $8818=9\times ((8+7)\times 65+4)+3\times 2+1.$
\item [] $8819=9+8765+43+2\times 1.$
\item [] $8820=9+8765+43+2+1.$
\item [] $8821=9\times (8+7)\times 65+43+2+1.$
\item [] $8822=98\times (7\times 6+5+43)+2\times 1.$
\item [] $8823=9+8765+(4+3)^2\times 1.$
\item [] $8824=9\times 8\times 7+65\times 4\times 32\times 1.$
\item [] $8825=9\times 8\times 7+65\times 4\times 32+1.$
\item [] $8826=(9+8)\times (7\times 65+4^3)+2+1.$
\item [] $8827=987+6^5+43+21.$
\item [] $8828=9\times 8+7+6\times (5+4)^3\times 2+1.$
\item [] $8829=987+6^5+4^3+2\times 1.$
\item [] $8830=987+6^5+4+3\times 21.$
\item [] $8831=9\times (8+7)\times 6+(5\times 4)^3+21.$
\item [] $8832=9\times (876+5)+43\times 21.$
\item [] $8833=(9+87)\times (6+54+32)+1.$
\item [] $8834=9\times 8+(7+6\times (5+4)^3)\times 2\times 1.$
\item [] $8835=9\times ((8+7)\times 65+4)+3+21.$
\item [] $8836=(9+8\times 7+6+5\times 4+3)^2\times 1.$
\item [] $8837=(9+8\times 7+6+5\times 4+3)^2+1.$
\item [] $8838=9+8765+43+21.$
\item [] $8839=9\times (8+7)\times 65+43+21.$
\item [] $8840=9+8765+4^3+2\times 1.$
\item [] $8841=9+8765+4+3\times 21.$
\item [] $8842=9\times (8+7)\times 65+4+3\times 21.$
\item [] $8843=9\times ((8+7)\times 65+4)+32\times 1.$
\item [] $8844=9+87+6\times (5+4)^3\times 2\times 1.$
\item [] $8845=9+87+6\times (5+4)^3\times 2+1.$
\item [] $8846=98+(76+5)\times 4\times 3^{(2+1)}.$
\item [] $8847=9\times (8\times 76+54+321).$
\item [] $8848=987+6^5+4^3+21.$
\item [] $8849=987+6^5+43\times 2\times 1.$
\item [] $8850=987+6^5+43\times 2+1.$
\item [] $8851=9+8+7\times (6+(5^4+3)\times 2)\times 1.$
\item [] $8852=9+8+7\times (6+(5^4+3)\times 2)+1.$
\item [] $8853=98+7+6\times (5+4)^3\times 2\times 1.$
\item [] $8854=98+7+6\times (5+4)^3\times 2+1.$
\item [] $8855=9+(876\times 5+43)\times 2\times 1.$
\item [] $8856=9+(876\times 5+43)\times 2+1.$
\item [] $8857=(9+8)\times (7\times 65+4^3+2\times 1).$
\item [] $8858=9+8\times 7\times (6\times 5+4\times 32)+1.$
\item [] $8859=9+8765+4^3+21.$
\item [] $8860=9+8765+43\times 2\times 1.$
\item [] $8861=9+8765+43\times 2+1.$
\item [] $8862=9\times (8+7)\times 65+43\times 2+1.$
\item [] $\mathit{8863=-9+87\times (6\times 5+4)\times 3-2\times 1.}$
\item [] $8864=(9+876\times 5+43)\times 2\times 1.$
\item [] $8865=9\times 87\times (6+5)+4\times 3\times 21.$
\item [] $8866=9+8\times 7+6^5+4^(3+2)+1.$
\item [] $8867=(9+8)\times 7+6\times 54\times 3^{(2+1)}.$
\item [] $8868=(9+8)\times 7+6\times (5+4)^3\times 2+1.$
\item [] $8869=9+(8\times 7+6\times (5+4)^3)\times 2\times 1.$
\item [] $8870=9+8765+4\times (3+21).$
\item [] $8871=987+6^5+4\times 3^{(2+1)}.$
\item [] $8872=(98\times 7\times 6+5\times 4^3)\times 2\times 1.$
\item [] $8873=(98\times 7\times 6+5\times 4^3)\times 2+1.$
\item [] $8874=9\times 8\times 76+54\times 3\times 21.$
\item [] $8875=(9+8\times 7+6)\times 5\times (4\times 3\times 2+1).$
\item [] $8876=((9+8)\times 7+6)\times (5+4^3+2)+1.$
\item [] $8877=98\times 7\times 6+(5+4^3)^2\times 1.$
\item [] $8878=98\times 7\times 6+(5+4^3)^2+1.$
\item [] $8879=9+8+7\times (6+5\times 4\times 3\times 21).$
\item [] $8880=(9\times 8+76)\times (54+3+2+1).$
\item[]$\mbox{Increasing order}$
\item [] $8881=1\times (2+3^4)\times (5+6+7+89).$
\item [] $8882=1+(2+3^4)\times (5+6+7+89).$
\item [] $8883=12^3+(4+5)\times (6+789).$
\item [] $\mathit{8884=1+(2+3\times 4+5)\times 6\times 78-9.}$
\item [] $\mathit{8885=-1+2\times 3\times (-4+(5+6)\times (7+8)\times 9).}$
\item [] $8886=12+3\times (4+5\times 6)\times (78+9).$
\item [] $8887=1\times 23\times ((4+5)\times 6\times 7+8)+9.$
\item [] $8888=(1+2^3+4)\times (5+678)+9.$
\item [] $\mathit{8889=1+2^3\times 4^5-6+78\times 9.}$
\item [] $\mathit{8890=(123+4)\times 5\times (6+7-8+9).}$
\item [] $8891=((1+2)\times 34\times 5+6+7)\times (8+9).$
\item [] $8892=1\times 2\times (34+5)\times (6\times 7+8\times 9).$
\item [] $8893=1+2\times (34+5)\times (6\times 7+8\times 9).$
\item [] $8894=1+2+((3^4+5)\times 6+7)\times (8+9).$
\item [] $\mathit{8895=(-1+23+4^5+67)\times 8-9.}$
\item [] $\mathit{8896=1\times 2^3\times (4^5+6-7+89).}$
\item [] $8897=(12\times (3^4+5+6)+7)\times 8+9.$
\item [] $8898=(1+2)\times ((3^4\times 5+6)\times 7+89).$
\item [] $\mathit{8899=-1+2^3\times 4^5+6+78\times 9.}$
\item [] $8900=1\times 2^3\times 4^5+6+78\times 9.$
\item [] $8901=1+2^3\times 4^5+6+78\times 9.$
\item [] $8902=1+2\times (3+(4+5)\times 6)\times 78+9.$
\item [] $8903=12+((3^4+5)\times 6+7)\times (8+9).$
\item [] $8904=(1+2)\times (3+4)\times (5\times 67+89).$
\item [] $8905=1^2\times (3^4+56)\times (7\times 8+9).$
\item [] $8906=1\times 2^3\times 4^5+6\times 7\times (8+9).$
\item [] $8907=1+2^3\times 4^5+6\times 7\times (8+9).$
\item [] $8908=1+2\times (3^4+56\times 78)+9.$
\item [] $8909=1+(23+45)\times (6\times 7+89).$
\item [] $8910=(12+3+45+6)\times (7+8)\times 9.$
\item [] $8911=1+(2+34+5\times 6)\times (7+8)\times 9.$
\item [] $8912=1\times 2+(3\times 4\times 5+6)\times (7+8)\times 9.$
\item [] $8913=1\times 234+(5+6)\times 789.$
\item [] $8914=1+234+(5+6)\times 789.$
\item [] $\mathit{8915=1\times 2^{(3\times 4)}-5+67\times 8\times 9.}$
\item [] $8916=(123+4+5)\times 67+8\times 9.$
\item [] $8917=12+(3^4+56)\times (7\times 8+9).$
\item [] $\mathit{8918=1\times 2^3+(4^5-6\times 7+8)\times 9.}$
\item [] $8919=1\times 2\times 3^4\times (5+6\times 7+8)+9.$
\item [] $8920=(1234+5\times 6)\times 7+8\times 9.$
\item [] $8921=(1\times 23+4^5+67)\times 8+9.$
\item [] $8922=12+(3\times 4\times 5+6)\times (7+8)\times 9.$
\item [] $\mathit{8923=-1+23\times (4+5\times (67+8)+9).}$
\item [] $8924=1\times 23\times (4+5\times (67+8)+9).$
\item [] $8925=1\times 2^{(3\times 4)}+5+67\times 8\times 9.$
\item [] $8926=1+2^{(3\times 4)}+5+67\times 8\times 9.$
\item [] $8927=1\times (2+3^4+5\times 6)\times (7+8\times 9).$
\item [] $8928=123\times 4\times (5+6+7)+8\times 9.$
\item [] $8929=1+2\times (3+45)\times (6+78+9).$
\item [] $\mathit{8930=(-1+(23-4)\times 5)\times ((6+7)\times 8-9).}$
\item [] $8931=(1+2^3)^4+5\times 6\times (7+8\times 9).$
\item [] $8932=(1+23+4)\times (5\times (6+7\times 8)+9).$
\item [] $8933=(123+4+5)\times 67+89.$
\item [] $8934=1\times 2\times 3\times (4+(5+6)\times (7+8)\times 9).$
\item [] $8935=(1^2+3)^4+(5+6)\times 789.$
\item [] $8936=1\times 2^3\times (4^5+6+78+9).$
\item [] $8937=(1234+5\times 6)\times 7+89.$
\item [] $8938=(1^2+3^4)\times (5\times 6+7+8\times 9).$
\item [] $\mathit{8939=12^3+(4^5+6)\times 7-8+9.}$
\item [] $8940=(12+3)\times 4\times (5+6\times (7+8+9)).$
\item [] $\mathit{8941=-1+2^3\times (4^5+6)+78\times 9.}$
\item [] $8942=1\times 2^3\times (4^5+6)+78\times 9.$
\item [] $8943=1+2^3\times (4^5+6)+78\times 9.$
\item [] $\mathit{8944=-12+3\times 45\times 67-89.}$
\item [] $8945=123\times 4\times (5+6+7)+89.$
\item [] $8946=1\times 2\times (3+4)\times (567+8\times 9).$
\item [] $8947=1+2\times (3+4)\times (567+8\times 9).$
\item [] $8948=1\times 2^3\times 4^5+(6+78)\times 9.$
\item [] $8949=1+2^3\times 4^5+(6+78)\times 9.$
\item [] $\mathit{8950=-1-2+(3^4+5)\times (6+7)\times 8+9.}$
\item[]$\mbox{Decreasing order}$
\item [] $8881=(9\times 8+76)\times (54+3\times 2)+1.$
\item [] $8882=(9\times 8+76)\times 5\times 4\times 3+2\times 1.$
\item [] $8883=(9+87+6\times 54+3)\times 21.$
\item [] $8884=9\times 87+6^5+4+321.$
\item [] $8885=9+87\times 6\times (5+4\times 3)+2\times 1.$
\item [] $8886=9+87\times 6\times (5+4\times 3)+2+1.$
\item [] $8887=9+876+(5\times 4)^3+2\times 1.$
\item [] $8888=9+876+(5\times 4)^3+2+1.$
\item [] $8889=9\times 8\times 7+65\times 43\times (2+1).$
\item [] $8890=98\times 7\times (6+5)+4^3\times 21.$
\item [] $8891=987+6^5+4\times 32\times 1.$
\item [] $8892=987+6^5+4^3\times 2+1.$
\item [] $8893=9+8+7\times (6+5^4+3)\times 2\times 1.$
\item [] $8894=98\times 7+6^5+432\times 1.$
\item [] $8895=98\times 7+6^5+432+1.$
\item [] $8896=9+87+6^5+4^(3+2)\times 1.$
\item [] $8897=98\times 76+(5+4^3)\times 21.$
\item [] $8898=9+(876\times 5+4^3)\times 2+1.$
\item [] $\mathit{8899=-9-8\times 7\times 6+5\times 43^2-1.}$
\item [] $8900=(9+876+5)\times (4+3+2+1).$
\item [] $8901=(9\times 8+76)\times 5\times 4\times 3+21.$
\item [] $8902=9+8765+4\times 32\times 1.$
\item [] $8903=9+8765+4\times 32+1.$
\item [] $8904=9\times (8+7)\times 65+4\times 32+1.$
\item [] $8905=98\times (7+65)+43^2\times 1.$
\item [] $8906=98\times 76+54\times 3^{(2+1)}.$
\item [] $8907=(9+876\times 5+4^3)\times 2+1.$
\item [] $8908=987+6^5+(4\times 3)^2+1.$
\item [] $8909=9+8+76\times (54+3\times 21).$
\item [] $8910=987+6^5+(4+3)\times 21.$
\item [] $8911=(9\times 8\times (7+6)+54)\times 3^2+1.$
\item [] $8912=9\times (8+7)\times (6+5\times 4\times 3)+2\times 1.$
\item [] $8913=9+8\times (7\times 6\times 5+43\times 21).$
\item [] $8914=9\times 876+5+4^(3+2)+1.$
\item [] $8915=9\times 8+7+(6\times 5+4^3)^2\times 1.$
\item [] $8916=9\times 8+7+(6\times 5+4^3)^2+1.$
\item [] $\mathit{8917=9-8\times 7\times 6+5\times 43^2-1.}$
\item [] $8918=9+8765+(4\times 3)^2\times 1.$
\item [] $8919=9+8765+(4\times 3)^2+1.$
\item [] $8920=9\times (8+7)\times 65+(4\times 3)^2+1.$
\item [] $8921=9+8765+(4+3)\times 21.$
\item [] $8922=9\times (8+7)\times 65+(4+3)\times 21.$
\item [] $\mathit{8923=-9\times 8+7\times (6-5+4\times 321).}$
\item [] $8924=(9\times (8+7+6)+5)\times (43+2+1).$
\item [] $8925=(9+8)\times 7\times (6+5+43+21).$
\item [] $8926=(9+8)\times 7\times (65+4+3\times 2)+1.$
\item [] $8927=98\times (7\times 6+5)+4321.$
\item [] $8928=(9+87)\times (6+54+32+1).$
\item [] $8929=(9+8+76)\times (5+43)\times 2+1.$
\item [] $8930=9\times 8\times (76+5+43)+2\times 1.$
\item [] $8931=9+(87+6\times (5+4)^3)\times 2\times 1.$
\item [] $8932=9+87+(6\times 5+4^3)^2\times 1.$
\item [] $8933=9+87+(6\times 5+4^3)^2+1.$
\item [] $8934=9\times 8+7\times (6+5\times 4\times 3\times 21).$
\item [] $8935=(9\times (8+7)\times (6+5)+4)\times 3\times 2+1.$
\item [] $8936=((9+876)\times 5+43)\times 2\times 1.$
\item [] $8937=9\times (876+54+3\times 21).$
\item [] $8938=9\times 87\times (6+5)+4+321.$
\item [] $8939=9\times 8\times (7+6)+(5\times 4)^3+2+1.$
\item [] $8940=(9+876+5+4)\times (3^2+1).$
\item [] $8941=98+7+(6\times 5+4^3)^2\times 1.$
\item [] $8942=(987+6)\times (5+4)+3+2\times 1.$
\item [] $8943=(9\times 87+6\times 5)\times (4+3\times 2+1).$
\item [] $8944=(987+6)\times (5+4)+3\times 2+1.$
\item [] $8945=(9+(8+7)\times 6+5)\times 43\times 2+1.$
\item [] $8946=(987+6)\times (5+4)+3^2\times 1.$
\item [] $8947=(987+6)\times (5+4)+3^2+1.$
\item [] $8948=9\times 8+7\times (6+5^4+3)\times 2\times 1.$
\item [] $8949=9\times 8\times (76+5+43)+21.$
\item [] $8950=(98+76+5)\times ((4+3)^2+1).$
\item[]$\mbox{Increasing order}$
\item [] $\mathit{8951=-1\times 2+(3^4+5)\times (6+7)\times 8+9.}$
\item [] $8952=12\times (34\times 5+6\times (7+89)).$
\item [] $8953=((1+2)^3+4^5+67)\times 8+9.$
\item [] $8954=1^2+(3^4+5)\times (6+7)\times 8+9.$
\item [] $8955=12^3+(4^5+6)\times 7+8+9.$
\item [] $8956=1234+(5+6)\times 78\times 9.$
\item [] $8957=1\times 2+((3\times 45+6)\times 7+8)\times 9.$
\item [] $8958=1+2+((3\times 45+6)\times 7+8)\times 9.$
\item [] $8959=(1\times 23\times 4\times 5+67)\times (8+9).$
\item [] $8960=(1^2+3)^4\times (5+6+7+8+9).$
\item [] $8961=1\times (23\times 4+5+6)\times (78+9).$
\item [] $8962=1+(23\times 4+5+6)\times (78+9).$
\item [] $\mathit{8963=1+2\times (3+4567-89).}$
\item [] $8964=12\times 345+67\times 8\times 9.$
\item [] $8965=12+(3^4+5)\times (6+7)\times 8+9.$
\item [] $8966=1+2^{(3\times 4)}+(5+67\times 8)\times 9.$
\item [] $8967=12+((3\times 45+6)\times 7+8)\times 9.$
\item [] $8968=1\times 2^{(3\times 4)}+56\times (78+9).$
\item [] $8969=1+2^{(3\times 4)}+56\times (78+9).$
\item [] $8970=1\times 23\times (45+6\times 7\times 8+9).$
\item [] $8971=1+23\times (45+6\times 7\times 8+9).$
\item [] $8972=1\times 2+(3+4\times 5)\times 6\times (7\times 8+9).$
\item [] $8973=1+2+(3+4\times 5)\times 6\times (7\times 8+9).$
\item [] $8974=1+(2+3^4)\times (5\times 6+78)+9.$
\item [] $\mathit{8975=1\times 2+3\times 45\times 67-8\times 9.}$
\item [] $8976=1^2\times 34\times (5+6)\times (7+8+9).$
\item [] $8977=1+(2+3^4+5)\times (6+7+89).$
\item [] $8978=1\times 2+34\times (5+6)\times (7+8+9).$
\item [] $8979=1+2+34\times (5+6)\times (7+8+9).$
\item [] $8980=1+23\times (4\times 5+6)\times (7+8)+9.$
\item [] $8981=(1+2)^(3+4)+5+6789.$
\item [] $8982=(12\times 3^4+5+6+7+8)\times 9.$
\item [] $8983=(1\times 2^{(3+4)}+5)\times 67+8\times 9.$
\item [] $8984=1+(2^{(3+4)}+5)\times 67+8\times 9.$
\item [] $8985=1+2^3\times (4^5+6\times (7+8)+9).$
\item [] $\mathit{8986=1+2\times 3^4\times 56-78-9.}$
\item [] $8987=1\times 2^3\times 4^5+6+789.$
\item [] $8988=1+2^3\times 4^5+6+789.$
\item [] $8989=(1\times 2+3^4+5+6+7)\times 89.$
\item [] $8990=1+(2+34+5\times (6+7))\times 89.$
\item [] $8991=(1+2\times (3+(4\times 5+6\times 7)\times 8))\times 9.$
\item [] $\mathit{8992=-1+2\times 3^4\times 56-7-8\times 9.}$
\item [] $8993=(123\times 4+5\times 6+7)\times (8+9).$
\item [] $8994=(1+23\times (4\times 5+6))\times (7+8)+9.$
\item [] $8995=(1+2^{(3+4)}+5)\times 67+8+9.$
\item [] $8996=(12+3+4)\times (5+6\times 78)+9.$
\item [] $8997=((1+2)\times 34+5)\times (6+78)+9.$
\item [] $\mathit{8998=-1+2-3+4\times 5\times (6\times 7+8)\times 9.}$
\item [] $\mathit{8999=1\times 2-3+4\times 5\times (6\times 7+8)\times 9.}$
\item [] $9000=1^{23}\times 4\times 5\times (6\times 7+8)\times 9.$
\item [] $9001=1+(2^{(3+4)}+5)\times 67+89.$
\item [] $9002=1\times 2^3\times 4^5+6\times (7+8)\times 9.$
\item [] $9003=1+2^3\times 4^5+6\times (7+8)\times 9.$
\item [] $9004=1^2+3+4\times 5\times (6\times 7+8)\times 9.$
\item [] $9005=1\times 2+3+4\times 5\times (6\times 7+8)\times 9.$
\item [] $9006=1+2+3+4\times 5\times (6\times 7+8)\times 9.$
\item [] $9007=1+2\times 3+4\times 5\times (6\times 7+8)\times 9.$
\item [] $9008=1\times 2^3+4\times 5\times (6\times 7+8)\times 9.$
\item [] $9009=1\times 2\times 3\times 4\times 5\times (67+8)+9.$
\item [] $9010=1+2\times 3\times 4\times 5\times (67+8)+9.$
\item [] $9011=1+2\times (3^4+56\times (7+8\times 9)).$
\item [] $9012=12\times (3\times 4\times 56+7+8\times 9).$
\item [] $\mathit{9013=-1-(2+3)^4+567\times (8+9).}$
\item [] $\mathit{9014=-1+2\times 3\times 4\times (5+6\times 7)\times 8-9.}$
\item [] $9015=12+3+4\times 5\times (6\times 7+8)\times 9.$
\item [] $9016=1\times 23\times 4\times (5+6+78+9).$
\item [] $9017=1+23\times 4\times (5+6+78+9).$
\item [] $9018=12^3+(4+5)\times 6\times (7+8)\times 9.$
\item [] $9019=1+(2\times 3)^4+(5+6)\times 78\times 9.$
\item [] $\mathit{9020=(1-2+3^4+5\times 6)\times (-7+89).}$
\item[]$\mbox{Decreasing order}$
\item [] $\mathit{8951=-9-8+76\times (-5+4^3)\times 2\times 1.}$
\item [] $8952=(9+8+7)\times (6\times (5\times 4\times 3+2)+1).$
\item [] $8953=9+8\times (7+6)\times (54+32)\times 1.$
\item [] $8954=9+8\times (7+6)\times (54+32)+1.$
\item [] $8955=987+6^5+4^3\times (2+1).$
\item [] $8956=(9+8)\times 7+(6\times 5+4^3)^2+1.$
\item [] $8957=9\times 8\times (7+6)+(5\times 4)^3+21.$
\item [] $8958=(98+7+6\times (5+4)^3)\times 2\times 1.$
\item [] $8959=(98+7+6\times (5+4)^3)\times 2+1.$
\item [] $8960=98+7\times (6+5\times 4\times 3\times 21).$
\item [] $8961=(987+6)\times (5+4)+3+21.$
\item [] $8962=98\times (76+5)+4^(3+2)\times 1.$
\item [] $8963=98\times (76+5)+4^(3+2)+1.$
\item [] $8964=9\times 8+76\times (54+3\times 21).$
\item [] $8965=9\times (8\times 7\times 6+54\times 3)\times 2+1.$
\item [] $8966=9+8765+4^3\times (2+1).$
\item [] $8967=9\times (8+7)\times 65+4^3\times (2+1).$
\item [] $8968=(9+8+7\times 6)\times (5+(4+3)\times 21).$
\item [] $8969=(987+6)\times (5+4)+32\times 1.$
\item [] $8970=9\times 876+543\times 2\times 1.$
\item [] $8971=9\times 876+543\times 2+1.$
\item [] $8972=(9+8\times 7)\times 6\times (5\times 4+3)+2\times 1.$
\item [] $8973=9\times (876+5\times 4\times 3\times 2+1).$
\item [] $8974=98+7\times (6+5^4+3)\times 2\times 1.$
\item [] $8975=98+7\times (6+5^4+3)\times 2+1.$
\item [] $8976=(9+8)\times (7\times 6+54\times 3^2\times 1).$
\item [] $8977=((9+8\times 7)\times (65+4)+3)\times 2+1.$
\item [] $8978=((9+876)\times 5+4^3)\times 2\times 1.$
\item [] $8979=((9+876)\times 5+4^3)\times 2+1.$
\item [] $\mathit{8980=987-6+(5\times 4)^3-2+1.}$
\item [] $8981=(98\times (7+6)+5+4)\times (3\times 2+1).$
\item [] $8982=(987+6+5)\times (4+3+2)\times 1.$
\item [] $8983=(987+6+5)\times (4+3+2)+1.$
\item [] $8984=(9\times 8+7\times (6+5^4)+3)\times 2\times 1.$
\item [] $8985=9\times (8+(76+5)\times 4)\times 3+21.$
\item [] $8986=((9+8)\times 7+6\times (5+4)^3)\times 2\times 1.$
\item [] $8987=((9+8)\times 7+6\times (5+4)^3)\times 2+1.$
\item [] $8988=(9+8\times (7\times 6+5)+43)\times 21.$
\item [] $\mathit{8989=9-8+7(6+5)\times 4\times 321.}$
\item [] $8990=98+76\times (54+3\times 21).$
\item [] $8991=9\times 87+6^5+432\times 1.$
\item [] $8992=9\times 87+6^5+432+1.$
\item [] $8993=(9+8)\times (7\times 6+54\times 3^2+1).$
\item [] $\mathit{8994=9\times (8+7\times 6)\times 5\times 4-3\times 2\times 1.}$
\item [] $8995=987+6+(5\times 4)^3+2\times 1.$
\item [] $8996=987+6+(5\times 4)^3+2+1.$
\item [] $\mathit{8997=9\times (8+7\times 6)\times 5\times 4-3\times (2-1).}$
\item [] $\mathit{8998=9\times (8+7\times 6)\times 5\times 4-3+2-1.}$
\item [] $8999=9+(8\times 7+6)\times ((5+4+3)^2+1).$
\item [] $9000=9\times 8\times (7\times 6+5\times 4+3\times 21).$
\item [] $9001=9\times 8\times (7+6)\times 5+4321.$
\item [] $9002=9\times (87\times (6+5)+43)+2\times 1.$
\item [] $9003=9\times (87\times (6+5)+43)+2+1.$
\item [] $\mathit{9004=9\times (8+7\times 6)\times 5\times 4+3+2-1.}$
\item [] $9005=9\times (8+7\times 6)\times 5\times 4+3+2\times 1.$
\item [] $9006=98\times 7+65\times 4^3\times 2\times 1.$
\item [] $9007=98\times 7+65\times 4\times 32+1.$
\item [] $\mathit{9008=987\times (6+5)-43^2\times 1.}$
\item [] $9009=9\times (87+6+5+43\times 21).$
\item [] $9010=9\times (8+7\times 6)\times 5\times 4+3^2+1.$
\item [] $9011=9\times (8+(7+6\times 54)\times 3)+2\times 1.$
\item [] $9012=9\times (8+(7+6\times 54)\times 3)+2+1.$
\item [] $\mathit{9013=9-8\times 7+6^5+4\times 321.}$
\item [] $9014=987+6+(5\times 4)^3+21.$
\item [] $9015=987+6^5+4\times 3\times 21.$
\item [] $9016=98\times (76+5+4+3\times 2+1).$
\item [] $9017=((9+8)\times 7+65)\times (4+3)^2+1.$
\item [] $9018=(987+6+5+4)\times 3^2\times 1.$
\item [] $9019=(9+8)\times 7\times 65+4\times 321.$
\item [] $9020=(9\times 8+7\times (6+5^4+3))\times 2\times 1.$
\item[]$\mbox{Increasing order}$
\item [] $9021=1\times (23\times 4+5)\times (6+78+9).$
\item [] $9022=1+(23\times 4+5)\times (6+78+9).$
\item [] $9023=1\times 23+4\times 5\times (6\times 7+8)\times 9.$
\item [] $9024=1+23+4\times 5\times (6\times 7+8)\times 9.$
\item [] $9025=1+2^3\times 4\times (5\times 6\times 7+8\times 9).$
\item [] $9026=(123+4)\times (56+7+8)+9.$
\item [] $9027=12^3+(4^5+6)\times 7+89.$
\item [] $9028=1^2+(3+4\times 5\times (6\times 7+8))\times 9.$
\item [] $9029=1\times 2^3\times (4^5+6)+789.$
\item [] $9030=1+2^3\times (4^5+6)+789.$
\item [] $\mathit{9031=1+2+3\times 45\times 67-8-9.}$
\item [] $\mathit{9032=-12+3\times 45\times 67+8-9.}$
\item [] $9033=1\times 2\times 3\times 4\times (5+6\times 7)\times 8+9.$
\item [] $9034=1+2\times 3\times 4\times (5+6\times 7)\times 8+9.$
\item [] $9035=1\times (2+3^4+56)\times (7\times 8+9).$
\item [] $9036=1^2+3\times 4^5+67\times 89.$
\item [] $9037=1\times 2+3\times 4^5+67\times 89.$
\item [] $9038=1+2+3\times 4^5+67\times 89.$
\item [] $9039=(1+23+45)\times (6\times 7+89).$
\item [] $9040=(1+2)^(3+4)+(5+6)\times 7\times 89.$
\item [] $9041=((1+2)\times 34\times (5+6)+7)\times 8+9.$
\item [] $9042=1\times 2\times 345\times (6+7)+8\times 9.$
\item [] $9043=1+2\times 345\times (6+7)+8\times 9.$
\item [] $9044=(12+34+5\times 6)\times 7\times (8+9).$
\item [] $9045=(1+2+34+5\times 6)\times (7+8)\times 9.$
\item [] $9046=1+(23+4\times (5+6))\times (7+8)\times 9.$
\item [] $9047=12+3\times 4^5+67\times 89.$
\item [] $9048=(12+3^4+5+6)\times (78+9).$
\item [] $9049=1+(2+3\times (4+5\times 6))\times (78+9).$
\item [] $9050=(12^3+4)\times 5+6\times (7\times 8+9).$
\item [] $\mathit{9051=-12+3\times (4+5)\times 6\times 7\times 8-9.}$
\item [] $\mathit{9052=1+2\times (3+4567)-89.}$
\item [] $\mathit{9053=(1+2)\times (3+45\times 67)+8-9.}$
\item [] $9054=1\times 2\times 3\times (4\times 5\times (67+8)+9).$
\item [] $9055=1+2\times 345\times (6+7)+8\times 9.$
\item [] $9056=1\times (2+3^4)\times (5+(6+7)\times 8)+9.$
\item [] $9057=((12+3)\times 4+56)\times 78+9.$
\item [] $\mathit{9058=12+3\times 45\times 67-8+9.}$
\item [] $9059=1\times 2\times 345\times (6+7)+89.$
\item [] $9060=1+2\times 345\times (6+7)+89.$
\item [] $9061=1\times (2+34+5)\times (6+7)\times (8+9).$
\item [] $9062=1^2\times 3\times 45\times 67+8+9.$
\item [] $9063=1^2+3\times 45\times 67+8+9.$
\item [] $9064=1\times 2+3\times 45\times 67+8+9.$
\item [] $9065=1+2+3\times 45\times 67+8+9.$
\item [] $\mathit{9066=1\times 2\times 3^4\times 56-7-8+9.}$
\item [] $9067=(1+2^{(3+4)}+5)\times 67+89.$
\item [] $\mathit{9068=1+((2\times 3)^4-5-6)\times 7+8\times 9.}$
\item [] $\mathit{9069=-12+3\times (4+5)\times 6\times 7\times 8+9.}$
\item [] $\mathit{9070=-(1+2^3)^4+5^6+7+8-9.}$
\item [] $9071=(1+2)\times (3+45\times 67)+8+9.$
\item [] $9072=12\times 3^4\times 5+6\times 78\times 9.$
\item [] $9073=1+(2+3+4)\times (56+7\times 8)\times 9.$
\item [] $9074=12+3\times 45\times 67+8+9.$
\item [] $9075=1+2+(3+45)\times (6+7+8)\times 9.$
\item [] $\mathit{9076=-1+((2\times 3)^4+5+6)\times 7-8\times 9.}$
\item [] $\mathit{9077=-1+2\times 3^4\times 56+7+8-9.}$
\item [] $9078=(12+3+45+6\times 7)\times 89.$
\item [] $9079=1+2\times 3\times (4\times (5+6\times 7)\times 8+9).$
\item [] $9080=(1+2)\times 3^4\times (5\times 6+7)+89.$
\item [] $9081=12\times (3+4)\times (5\times 6+78)+9.$
\item [] $9082=1^2+3^4\times (56+7\times 8)+9.$
\item [] $9083=1\times 2+3^4\times (56+7\times 8)+9.$
\item [] $9084=12+(3+45)\times (6+7+8)\times 9.$
\item [] $9085=(1+(2+3\times 4+5)\times 6)\times (7+8\times 9).$
\item [] $9086=1\times 2\times (3+4+(56+7)\times 8\times 9).$
\item [] $9087=12\times 34+(5+6)\times 789.$
\item [] $9088=(1^{23}+4)^5+67\times 89.$
\item [] $9089=12\times 3\times 4\times (56+7)+8+9.$
\item [] $9090=(1+2)\times (3+(4+5)\times 6\times 7\times 8)+9.$
\item[]$\mbox{Decreasing order}$
\item [] $9021=((98\times 7+65)\times 4+3)\times (2+1).$
\item [] $\mathit{9022=-98+76\times 5\times 4\times 3\times 2\times 1.}$
\item [] $\mathit{9023=-98+76\times 5\times 4\times 3\times 2+1.}$
\item [] $9024=9\times 87\times 6+5+4321.$
\item [] $9025=(9+(876+5^4)\times 3)\times 2+1.$
\item [] $9026=9+8765+4\times 3\times 21.$
\item [] $9027=9\times (8+7)\times 65+4\times 3\times 21.$
\item [] $9028=(9\times 8+76)\times (54+3\times 2+1).$
\item [] $9029=(9+87\times 6)\times (5+4\times 3)+2\times 1.$
\item [] $9030=(9\times 87+6\times 5\times 4)\times (3^2+1).$
\item [] $9031=(98+7)\times (6\times (5+4)+32)+1.$
\item [] $9032=98\times 7+(6+5\times 4)\times 321.$
\item [] $9033=(98\times 7+65)\times 4\times 3+21.$
\item [] $9034=9+8\times (7\times 6+5)\times 4\times 3\times 2+1.$
\item [] $9035=(9+8\times 7)\times (6+5+4\times 32\times 1).$
\item [] $9036=(9+8\times 7)\times (6+5+4\times 32)+1.$
\item [] $9037=(98+7\times (6+5^4)+3)\times 2+1.$
\item [] $\mathit{9038=9+(8+7+6)\times 5\times 43\times 2-1.}$
\item [] $9039=9+(8+7+6)\times 5\times 43\times 2\times 1.$
\item [] $9040=9+(8+7+6)\times 5\times 43\times 2+1.$
\item [] $9041=9+8\times (7\times 6+543\times 2+1).$
\item [] $9042=98+(7+6)\times (5^4+3\times 21).$
\item [] $9043=(9+8\times (7\times 6+5)\times 4\times 3)\times 2+1.$
\item [] $9044=(9+8)\times 7\times (65+4+3\times 2+1).$
\item [] $9045=9\times 87\times (6+5)+432\times 1.$
\item [] $9046=9\times 87\times (6+5)+432+1.$
\item [] $9047=987\times 6+5^4\times (3+2\times 1).$
\item [] $9048=987\times 6+5^4\times (3+2)+1.$
\item [] $9049=9+8+7\times 6\times 5\times 43+2\times 1.$
\item [] $9050=9+8+7\times 6\times 5\times 43+2+1.$
\item [] $9051=(98+7+6+5\times 4^3)\times 21.$
\item [] $\mathit{9052=9-8+7\times 6\times 5\times 43+21.}$
\item [] $\mathit{9053=(-9+8)\times 7+6^5+4\times 321.}$
\item [] $9054=9\times (87\times (6+5)+(4+3)^2\times 1).$
\item [] $9055=9\times (87\times (6+5)+(4+3)^2)+1.$
\item [] $\mathit{9056=-9\times (8+7+6)+5\times 43^2\times 1.}$
\item [] $9057=9+8\times (7+6)\times (54+32+1).$
\item [] $9058=9+87\times (6\times (5+4\times 3)+2)+1.$
\item [] $9059=9+(8+7\times 6)\times (5\times (4+32)+1).$
\item [] $9060=9+(8+76\times 5+43)\times 21.$
\item [] $9061=(9+876+5^4)\times 3\times 2+1.$
\item [] $9062=9+8+7\times (6\times 5\times 43+2)+1.$
\item [] $9063=9\times (8+7\times 6)\times 5\times 4+3\times 21.$
\item [] $9064=9\times (87\times 6+5)+4321.$
\item [] $\mathit{9065=-9+8\times 7\times 6\times (5+4)\times 3+2\times 1.}$
\item [] $\mathit{9066=-9+8+7+6^5+4\times 321.}$
\item [] $\mathit{9067=(9-8)\times 7+6^5+4\times 321.}$
\item [] $9068=9+8+7\times 6\times 5\times 43+21.$
\item [] $\mathit{9069=9\times 8\times 7\times 6(5+4)\times 3-2-1.}$
\item [] $\mathit{9070=9+8-7+6^5+4\times 321.}$
\item [] $9071=98\times 7+65\times 43\times (2+1).$
\item [] $9072=(98+7+6\times 54+3)\times 21.$
\item [] $9073=98\times 76+5\times (4+321).$
\item [] $9074=9\times 8\times 7\times (6+5+4+3)+2\times 1.$
\item [] $9075=9\times 8\times 7\times (6+5+4+3)+2+1.$
\item [] $\mathit{9076=(9\times 8\times 7\times 6+5-4)\times 3+2-1.}$
\item [] $9077=98\times 76+543\times (2+1).$
\item [] $9078=(9+8\times (7\times 6+5)\times 4)\times 3\times 2\times 1.$
\item [] $9079=(9+8)\times 7\times 65+4^3\times 21.$
\item [] $9080=(9+8+7\times 6\times 5)\times 4\times (3^2+1).$
\item [] $9081=9\times 876+(54+3)\times 21.$
\item [] $9082=9+8+7\times (6+5+4\times 321).$
\item [] $9083=9+8\times 7\times 6\times (5+4)\times 3+2\times 1.$
\item [] $9084=9+8+7+6^5+4\times 321.$
\item [] $9085=(9\times 8+7)\times (6\times 5+4^3+21).$
\item [] $\mathit{9086=(9-8+7+6)\times (5^4+3+21).}$
\item [] $9087=9+((8+76)\times 54+3)\times 2\times 1.$
\item [] $9088=987+6^5+4+321.$
\item [] $9089=(9+8+7+65\times 4)\times 32+1.$
\item [] $9090=9\times (876+5+4\times 32+1).$
\item[]$\mbox{Increasing order}$
\item [] $9091=((1^2+3)^4\times 5+6)\times 7+89.$
\item [] $\mathit{9092=1+23\times (4+56\times 7)-8-9.}$
\item [] $9093=12+3^4\times (56+7\times 8)+9.$
\item [] $\mathit{9094=-1\times 2+3\times 4\times (56+78\times 9).}$
\item [] $9095=1\times 2\times (3+4+567\times 8)+9.$
\item [] $9096=1\times 2\times 3^4\times 56+7+8+9.$
\item [] $9097=1+2\times 3^4\times 56+7+8+9.$
\item [] $9098=1\times 2+3\times 4\times (56+78\times 9).$
\item [] $9099=1+2+3\times 4\times (56+78\times 9).$
\item [] $9100=1^2+3\times (4+5+6\times 7\times 8\times 9).$
\item [] $9101=1\times 2+3\times (4+5+6\times 7\times 8\times 9).$
\item [] $9102=(1234+56)\times 7+8\times 9.$
\item [] $9103=1\times 2\times (3^4\times 56+7)+8+9.$
\item [] $9104=1+2\times (3^4\times 56+7)+8+9.$
\item [] $9105=1+2\times (3+4+567\times 8+9).$
\item [] $9106=1+2\times (3\times 4+567\times 8)+9.$
\item [] $\mathit{9107=(1\times 2+3\times 45)\times 67-8\times 9.}$
\item [] $9108=12+3\times 4\times (56+78\times 9).$
\item [] $9109=1+2\times (34\times 5+6\times 7\times 8)\times 9.$
\item [] $9110=(12^3+4)\times 5+(6\times 7+8)\times 9.$
\item [] $9111=12+3\times (4+5+6\times 7\times 8\times 9).$
\item [] $9112=1+2\times (3^4\times 56+7+8)+9.$
\item [] $9113=1+2^3\times (4+56+7)\times (8+9).$
\item [] $9114=1\times 2\times (345+6\times 78\times 9).$
\item [] $9115=1+2\times (345+6\times 78\times 9).$
\item [] $\mathit{9116=1-2+3\times 45\times 67+8\times 9.}$
\item [] $9117=1^2\times 3\times 45\times 67+8\times 9.$
\item [] $9118=1^2+3\times 45\times 67+8\times 9.$
\item [] $9119=1\times 2+3\times 45\times 67+8\times 9.$
\item [] $9120=1+2+3\times 45\times 67+8\times 9.$
\item [] $9121=1\times 2\times 34\times (56+78)+9.$
\item [] $9122=1+2\times 34\times (56+78)+9.$
\item [] $9123=123+4\times 5\times (6\times 7+8)\times 9.$
\item [] $9124=1+2+(3\times 4+5)\times 67\times 8+9.$
\item [] $9125=1\times 23\times (4+56\times 7)+8+9.$
\item [] $9126=1\times 2\times (34+5+6\times 78)\times 9.$
\item [] $9127=1+2\times (34+5+6\times 78)\times 9.$
\item [] $9128=1\times 2^3\times 4^5+(6+7)\times 8\times 9.$
\item [] $9129=12+3\times 45\times 67+8\times 9.$
\item [] $9130=1\times 2\times (34\times (56+78)+9).$
\item [] $9131=1\times 23\times (4\times (5+6)\times 7+89).$
\item [] $9132=1^2\times 3\times (4\times 5+6\times 7\times 8\times 9).$
\item [] $9133=12+(3\times 4+5)\times 67\times 8+9.$
\item [] $9134=1\times 2345+6789.$
\item [] $9135=1+2345+6789.$
\item [] $9136=1\times 2+3\times 45\times 67+89.$
\item [] $9137=1+2+3\times 45\times 67+89.$
\item [] $9138=1+2\times 3^4\times 56+7\times 8+9.$
\item [] $9139=(12+3+4)\times (56\times 7+89).$
\item [] $9140=1\times 2\times (34+(56+7)\times 8\times 9).$
\item [] $9141=(1+2)\times (3+4\times 5+6\times 7\times 8\times 9).$
\item [] $9142=1\times 2\times (3+4)\times (5\times 6+7\times 89).$
\item [] $9143=1+2\times (3+4)\times (5\times 6+7\times 89).$
\item [] $9144=12\times 3\times 4\times (56+7)+8\times 9.$
\item [] $9145=1^{23}+(4\times 5\times 6+7)\times 8\times 9.$
\item [] $9146=12+3\times 45\times 67+89.$
\item [] $9147=1^2\times 3+(4\times 5\times 6+7)\times 8\times 9.$
\item [] $9148=1^2+3+(4\times 5\times 6+7)\times 8\times 9.$
\item [] $9149=1\times 2\times (34+567\times 8)+9.$
\item [] $9150=1+2\times (34+567\times 8)+9.$
\item [] $9151=1\times 2\times 3^4\times 56+7+8\times 9.$
\item [] $9152=1+2\times 3^4\times 56+7+8\times 9.$
\item [] $9153=(12\times 3^4+5\times 6+7+8)\times 9.$
\item [] $9154=1+((2\times 3^4+5)\times 6+7+8)\times 9.$
\item [] $\mathit{9155=1+2\times 3^4\times 56-7+89.}$
\item [] $9156=12\times (3+4)\times (5\times 6+7+8\times 9).$
\item [] $9157=1\times 2\times (3+4567)+8+9.$
\item [] $9158=1+2\times (3+4567)+8+9.$
\item [] $9159=1\times 2\times 3^4\times 56+78+9.$
\item [] $9160=1+2\times 3^4\times 56+78+9.$
\item[]$\mbox{Decreasing order}$
\item [] $9091=(98\times 7+6\times 54)\times 3^2+1.$
\item [] $\mathit{9092=9\times (8+7)\times 65-4+321.}$
\item [] $9093=9\times 8\times 7\times (6+5+4+3)+21.$
\item [] $\mathit{9094=9-8+7\times (65\times 4\times (3+2)-1).}$
\item [] $9095=(9+87+6+5)\times (4^3+21).$
\item [] $9096=9\times 8+(7\times 6+5)\times 4^3\times (2+1).$
\item [] $9097=9+8\times (7\times 6\times (5+4)\times 3+2\times 1).$
\item [] $9098=9+8\times (7\times 6\times (5+4)\times 3+2)+1.$
\item [] $9099=9+8765+4+321.$
\item [] $9100=9\times (8+7)\times 65+4+321.$
\item [] $9101=(9\times 8\times 7\times 6+5+4)\times 3+2\times 1.$
\item [] $9102=9+8\times 7\times 6\times (5+4)\times 3+21.$
\item [] $9103=9\times 87+65\times 4\times 32\times 1.$
\item [] $9104=9\times 87+65\times 4\times 32+1.$
\item [] $9105=9\times 8+7\times 6\times 5\times 43+2+1.$
\item [] $9106=987+6^5+(4+3)^{(2+1)}.$
\item [] $9107=((9+8)\times 76+5+4)\times (3\times 2+1).$
\item [] $9108=9\times (8+7\times 65+43)\times 2\times 1.$
\item [] $9109=9\times (8+7\times 65+43)\times 2+1.$
\item [] $9110=9\times 876+(5\times (4+3))^2+1.$
\item [] $\mathit{9111=-9\times 8+765\times 4\times 3+2+1.}$
\item [] $9112=(9+8)\times (7+(6+5+4\times 3)^2)\times 1.$
\item [] $9113=(9+8)\times (7+(6+5+4\times 3)^2)+1.$
\item [] $9114=(9\times 8+7\times 6+5\times 4^3)\times 21.$
\item [] $9115=98\times (7+6\times (5+4)+32)+1.$
\item [] $9116=98\times (76+5+4\times 3)+2\times 1.$
\item [] $9117=9+8+7\times 65\times 4\times (3+2)\times 1.$
\item [] $9118=(987+6+5\times 4)\times 3^2+1.$
\item [] $\mathit{9119=-9+8+76\times 5\times 4\times 3\times 2\times 1.}$
\item [] $9120=(9\times 8\times 7\times 6+5+4)\times 3+21.$
\item [] $\mathit{9121=9-8+76\times 5\times 4\times 3\times 2\times 1.}$
\item [] $9122=98+(7\times 6+5)\times 4^3\times (2+1).$
\item [] $9123=9\times 8+7\times 6\times 5\times 43+21.$
\item [] $9124=9+8+7\times (65\times 4\times (3+2)+1).$
\item [] $9125=9+8\times 7+6^5+4\times 321.$
\item [] $9126=9\times ((8+7)\times (6+5)+4)\times 3\times 2\times 1.$
\item [] $9127=9\times ((8+7)\times (6+5)+4)\times 3\times 2+1.$
\item [] $\mathit{9128=9-8+7+6^5+4^3\times 21.}$
\item [] $9129=9\times 87+(6+5\times 4)\times 321.$
\item [] $9130=98+7\times 6\times 5\times 43+2\times 1.$
\item [] $9131=98+7\times 6\times 5\times 43+2+1.$
\item [] $9132=9+8+7\times 6\times (5\times 43+2)+1.$
\item [] $\mathit{9133=-98-7-6+5\times 43^2-1.}$
\item [] $9134=(9\times 8\times 7\times 6+5\times 4)\times 3+2\times 1.$
\item [] $9135=(9\times 8\times 7\times 6+5\times 4)\times 3+2+1.$
\item [] $9136=(98+7)\times (6+(5+4)\times 3^2)+1.$
\item [] $9137=9+8+76\times 5\times 4\times 3\times 2\times 1.$
\item [] $9138=9+8+76\times 5\times 4\times 3\times 2+1.$
\item [] $9139=9\times 8+7+6^5+4\times 321.$
\item [] $9140=9\times 876+(5^4+3)\times 2\times 1.$
\item [] $9141=9\times 876+(5^4+3)\times 2+1.$
\item [] $9142=98+7\times (6\times 5\times 43+2\times 1).$
\item [] $9143=98+7\times (6\times 5\times 43+2)+1.$
\item [] $9144=9\times 876+5\times 4\times 3\times 21.$
\item [] $9145=(9+8+7\times 6)\times 5\times (4+3^{(2+1)}).$
\item [] $9146=9+(8+76\times 5\times 4\times 3)\times 2+1.$
\item [] $\mathit{9147=-98+7-6+5\times 43^2-1.}$
\item [] $\mathit{9148=-98+7-6+5\times 43^2\times 1.}$
\item [] $9149=98+7\times 6\times 5\times 43+21.$
\item [] $9150=9\times 8\times (7+6\times 5\times 4)+3+2+1.$
\item [] $9151=9\times 8\times (7+6\times 5\times 4)+3\times 2+1.$
\item [] $9152=(98+(7\times 6+5)\times 4)\times 32\times 1.$
\item [] $9153=(9\times 8+7\times 6\times 5+4)\times 32+1.$
\item [] $9154=(9+8+76\times 5\times 4\times 3)\times 2\times 1.$
\item [] $9155=9\times (8+7+6\times 54)\times 3+2\times 1.$
\item [] $9156=9+87+6^5+4\times 321.$
\item [] $\mathit{9157=9876-5\times (4\times 3)^2+1.}$
\item [] $\mathit{9158=-9+8+765\times 4\times 3-21.}$
\item [] $9159=9+(8+7\times 6)\times (54\times 3+21).$
\item [] $\mathit{9160=9-8+765\times 4\times 3-21.}$
\item[]$\mbox{Increasing order}$
\item [] $9161=12\times 3\times 4\times (56+7)+89.$
\item [] $9162=12^3\times 4+5\times (6\times 7+8)\times 9.$
\item [] $9163=(123\times 4+5+6\times 7)\times (8+9).$
\item [] $9164=1+(23+(4+5)\times 6)\times 7\times (8+9).$
\item [] $9165=1\times 2\times (3+4567+8)+9.$
\item [] $9166=1+2\times (3+4567+8)+9.$
\item [] $9167=(12\times 3+4+56+7)\times 89.$
\item [] $9168=1\times 2\times 3^4\times 56+7+89.$
\item [] $9169=1+2\times 3^4\times 56+7+89.$
\item [] $9170=(12^3+4)\times 5+6+7\times 8\times 9.$
\item [] $9171=123\times 4+(5+6)\times 789.$
\item [] $9172=1^2+(3+(4\times 5\times 6+7)\times 8)\times 9.$
\item [] $9173=(1+(2\times 3)^4+5+6)\times 7+8+9.$
\item [] $9174=1^2\times 3\times 4^5+678\times 9.$
\item [] $9175=1+2\times (3+4567+8+9).$
\item [] $9176=1+2\times (3^4\times 56+7)+89.$
\item [] $9177=1+2+3\times 4^5+678\times 9.$
\item [] $9178=1+23\times (4+5+6\times (7\times 8+9)).$
\item [] $\mathit{9179=(1\times 2+3\times 45)\times 67\times (-8+9).}$
\item [] $9180=12\times (3+4+56+78\times 9).$
\item [] $9181=1+23\times (4+56\times 7)+8\times 9.$
\item [] $9182=(12^3+4)\times 5+6\times (78+9).$
\item [] $9183=(1+2)\times (3+4^5)+678\times 9.$
\item [] $9184=1\times 2\times (3+4)\times (567+89).$
\item [] $9185=1+2\times (3+4)\times (567+89).$
\item [] $9186=12+3\times 4^5+678\times 9.$
\item [] $9187=1+2\times (3+45\times (6+7+89)).$
\item [] $\mathit{9188=1\times 2^3+4\times 5\times (6\times 78-9).}$
\item [] $9189=(1+2)\times (34+5+6\times 7\times 8\times 9).$
\item [] $9190=1+(2+3+(4\times 5\times 6+7)\times 8)\times 9.$
\item [] $9191=1\times 2\times 3^4\times 56+7\times (8+9).$
\item [] $9192=12\times (34+5\times 6+78\times 9).$
\item [] $9193=1\times 2\times (3+4+567)\times 8+9.$
\item [] $9194=1234\times 5+6\times 7\times 8\times 9.$
\item [] $9195=(1+2\times 3)^4+5+6789.$
\item [] $9196=(1\times 2+3\times 45)\times 67+8+9.$
\item [] $9197=1\times 23\times (4+56\times 7)+89.$
\item [] $9198=1+23\times (4+56\times 7)+89.$
\item [] $9199=1+2\times (3+4)\times (5\times (6+7)+8)\times 9.$
\item [] $\mathit{9200=(1+2)\times 3\times 4^5-6+7-8-9.}$
\item [] $9201=(1^2+3\times 45)\times 67+89.$
\item [] $9202=1\times 2\times (3^4\times 56+7\times 8+9).$
\item [] $9203=1+2\times (3^4\times 56+7\times 8+9).$
\item [] $9204=(1+2+3)\times (4^5+6+7\times 8\times 9).$
\item [] $9205=(12^3+4)\times 5+67\times 8+9.$
\item [] $\mathit{9206=-1+2\times 3^4\times 56+(7+8)\times 9.}$
\item [] $9207=1^2\times (345+678)\times 9.$
\item [] $9208=1+2\times 3^4\times 56+(7+8)\times 9.$
\item [] $9209=1\times 2+3\times (45+6\times 7\times 8\times 9).$
\item [] $9210=1+2+3\times (45+6\times 7\times 8\times 9).$
\item [] $\mathit{9211=(1+2)\times 3\times 4^5+67-8\times 9.}$
\item [] $9212=1\times 2\times (3+4567)+8\times 9.$
\item [] $9213=1+2\times (3+4567)+8\times 9.$
\item [] $9214=(1+2+(3+4)\times (5+6)\times 7)\times (8+9).$
\item [] $9215=(1+2\times 3^4)\times 56+78+9.$
\item [] $9216=12\times (3\times 4\times 56+7+89).$
\item [] $9217=1+2\times 3\times 4\times (5\times (67+8)+9).$
\item [] $9218=(12^3+4)\times 5+(6+7\times 8)\times 9.$
\item [] $9219=12+3\times (45+6\times 7\times 8\times 9).$
\item [] $9220=1^2+3\times (4+(5+6\times 7\times 8)\times 9).$
\item [] $9221=(1\times (2\times 3)^4+5+6)\times 7+8\times 9.$
\item [] $9222=1+2+3\times (4+(5+6\times 7\times 8)\times 9).$
\item [] $9223=1+((2+3)\times 4\times 5+6)\times (78+9).$
\item [] $9224=(1+2\times 3^4)\times 56+7+89.$
\item [] $9225=123\times (45+6+7+8+9).$
\item [] $9226=1+(2+345+678)\times 9.$
\item [] $9227=(1^{23}+4)^5+678\times 9.$
\item [] $9228=1+(2+3)^4\times 5+678\times 9.$
\item [] $9229=1\times 2\times (3+4567)+89.$
\item [] $9230=1+2\times (3+4567)+89.$
\item[]$\mbox{Decreasing order}$
\item [] $9161=(9+8\times 76\times 5+4)\times 3+2\times 1.$
\item [] $9162=9+(8\times 76\times 5+4)\times 3+21.$
\item [] $9163=98+7\times (6+5+4\times 321).$
\item [] $9164=(9\times 8+7)\times (6\times 5+43\times 2\times 1).$
\item [] $9165=98+7+6^5+4\times 321.$
\item [] $9166=(9+8)\times 7\times (65+4\times 3)+2+1.$
\item [] $\mathit{9167=-9+8-76+5\times 43^2-1.}$
\item [] $9168=9\times 87+65\times 43\times (2+1).$
\item [] $\mathit{9169=9-8-76+5\times 43^2-1.}$
\item [] $9170=98+7\times 6\times (5+4)\times (3+21).$
\item [] $9171=987\times 6+(54+3)^2\times 1.$
\item [] $9172=987\times 6+(54+3)^2+1.$
\item [] $9173=9\times 876+5+4\times 321.$
\item [] $9174=9\times (8+7+6\times 54)\times 3+21.$
\item [] $\mathit{9175=98\times 76+54\times 32-1.}$
\item [] $9176=98\times 76+54\times 32\times 1.$
\item [] $9177=98\times 76+54\times 32+1.$
\item [] $9178=9+(8+76\times 5\times 4)\times 3\times 2+1.$
\item [] $9179=(9+8)\times 7+6^5+4\times 321.$
\item [] $9180=9\times 876+54\times (3+21).$
\item [] $9181=(9+8)\times (7+65\times 4+3)\times 2+1.$
\item [] $9182=9\times (8+7\times (6+5))\times 4\times 3+2\times 1.$
\item [] $9183=9\times 8\times 7+6^5+43\times 21.$
\item [] $9184=(9+8)\times 7\times (65+4\times 3)+21.$
\item [] $9185=9+8\times 7+6^5+4^3\times 21.$
\item [] $9186=9\times 8+7\times 6\times (5\times 43+2\times 1).$
\item [] $9187=9+(8+7+6)\times (5+432)+1.$
\item [] $9188=9+(8+7\times 654+3)\times 2+1.$
\item [] $9189=(987+6\times 5+4)\times 3^2\times 1.$
\item [] $9190=(987+6\times 5+4)\times 3^2+1.$
\item [] $9191=((9+8)\times 7\times (6+5)+4)\times (3\times 2+1).$
\item [] $9192=9\times 8+76\times 5\times 4\times 3\times 2\times 1.$
\item [] $9193=9\times 8+76\times 5\times 4\times 3\times 2+1.$
\item [] $9194=9+8+7\times (6\times 5\times 43+21).$
\item [] $9195=987+6^5+432\times 1.$
\item [] $9196=987+6^5+432+1.$
\item [] $9197=(9+8+7\times 654+3)\times 2+1.$
\item [] $9198=98+7\times 65\times 4\times (3+2)\times 1.$
\item [] $9199=9+8+765\times 4\times 3+2\times 1.$
\item [] $9200=9+8+765\times 4\times 3+2+1.$
\item [] $9201=9\times (8+7\times (6+5))\times 4\times 3+21.$
\item [] $9202=(9+87+6+5)\times 43\times 2\times 1.$
\item [] $9203=(9+87+6+5)\times 43\times 2+1.$
\item [] $9204=9+(8\times 76+5)\times (4\times 3+2+1).$
\item [] $9205=98+7\times (65\times 4\times (3+2)+1).$
\item [] $9206=9+8765+432\times 1.$
\item [] $9207=9+8765+432+1.$
\item [] $9208=9\times (8+7)\times 65+432+1.$
\item [] $\mathit{9209=-9+8\times (7+6+5)\times 4^3+2\times 1.}$
\item [] $9210=9\times (87\times (6+5)+4^3)+21.$
\item [] $\mathit{9211=98-7+6^5+4^3\times 21.}$
\item [] $9212=98+7\times 6\times (5\times 43+2\times 1).$
\item [] $9213=9\times 87\times 6+5\times 43\times 21.$
\item [] $9214=(9+8)\times (7\times 65+43\times 2+1).$
\item [] $9215=9+8+7\times (654+3)\times 2\times 1.$
\item [] $9216=9+87+6^5+4^3\times 21.$
\item [] $9217=(9+8+7+65+4+3)^2+1.$
\item [] $9218=98+76\times 5\times 4\times 3\times 2\times 1.$
\item [] $9219=98+76\times 5\times 4\times 3\times 2+1.$
\item [] $\mathit{9220=9876-5^4-32+1.}$
\item [] $9221=(9\times (8\times 7\times 6+5)+4)\times 3+2\times 1.$
\item [] $9222=9\times 8\times 76+5^4\times 3\times 2\times 1.$
\item [] $9223=9\times 8\times 76+5^4\times 3\times 2+1.$
\item [] $9224=9+(8+7\times (654+3))\times 2+1.$
\item [] $9225=98+7+6^5+4^3\times 21.$
\item [] $9226=9+(8+76+5+4+3)^2+1.$
\item [] $9227=9+8\times (76\times 5+4)\times 3+2\times 1.$
\item [] $9228=9+8\times (76\times 5+4)\times 3+2+1.$
\item [] $\mathit{9229=-9-8+7-6+5\times 43^2\times 1.}$
\item [] $9230=98\times 76+54\times (32+1).$
\item[]$\mbox{Increasing order}$
\item [] $9231=(123+(4+56)\times 7)\times (8+9).$
\item [] $9232=(1+(2+3)^4)\times 5+678\times 9.$
\item [] $9233=(123+4)\times (5+67)+89.$
\item [] $9234=(1+2+345+678)\times 9.$
\item [] $9235=1+2\times (3^4\times 5+6\times 78\times 9).$
\item [] $9236=(12^3+4)\times 5+6\times (7+89).$
\item [] $9237=1\times 2\times (3^4\times 56+78)+9.$
\item [] $9238=(1\times (2\times 3)^4+5+6)\times 7+89.$
\item [] $9239=1+((2\times 3)^4+5+6)\times 7+89.$
\item [] $9240=12\times (3\times 4+56+78\times 9).$
\item [] $9241=1+2^3\times (4^5+6\times 7+89).$
\item [] $9242=(1+2)^3\times (4+5\times 67)+89.$
\item [] $9243=1\times 2\times (3^4+567\times 8)+9.$
\item [] $9244=1+2\times (3^4+567\times 8)+9.$
\item [] $9245=(1+(2\times 3)^4+5+6)\times 7+89.$
\item [] $9246=1\times 2\times (3^4\times 56+78+9).$
\item [] $9247=1+2\times (3^4\times 56+78+9).$
\item [] $9248=1\times 2\times 34\times (5+6\times 7+89).$
\item [] $9249=123\times 4\times 5+6789.$
\item [] $9250=(1+2^3)^4+5\times 67\times 8+9.$
\item [] $9251=1\times (2+3\times 45)\times 67+8\times 9.$
\item [] $9252=1+(2+3\times 45)\times 67+8\times 9.$
\item [] $9253=1+2\times (3^4+567\times 8+9).$
\item [] $9254=(1+2\times 3)^4+(5+6)\times 7\times 89.$
\item [] $9255=(1+2\times 34)\times (56+78)+9.$
\item [] $9256=1\times 2\times (34+5+6+7)\times 89.$
\item [] $9257=1+2\times (34+5+6+7)\times 89.$
\item [] $9258=(1+2)\times (3\times 4^5+6)+7+8+9.$
\item [] $9259=(1+2+34)\times 5\times (6\times 7+8)+9.$
\item [] $\mathit{9260=1-2+3\times (45\times 67+8\times 9).}$
\item [] $9261=12^3\times 4+5\times 6\times 78+9.$
\item [] $9262=1^2+3\times (45\times 67+8\times 9).$
\item [] $9263=(1+2+3\times 45)\times 67+8+9.$
\item [] $9264=1\times 2^3\times (456+78\times 9).$
\item [] $9265=1+2^3\times (456+78\times 9).$
\item [] $9266=1^2+(3\times 4+5)\times (67\times 8+9).$
\item [] $9267=123+(4\times 5\times 6+7)\times 8\times 9.$
\item [] $9268=1\times (2+3\times 45)\times 67+89.$
\item [] $9269=1+(2+3\times 45)\times 67+89.$
\item [] $9270=1\times 2\times (3+456+7\times 8)\times 9.$
\item [] $9271=1+2\times (3+456+7\times 8)\times 9.$
\item [] $9272=1+2^{(3\times 4)}+(567+8)\times 9.$
\item [] $9273=123\times 45+6\times 7\times 89.$
\item [] $9274=1+2\times (3\times 4+567)\times 8+9.$
\item [] $9275=(1+2^3)\times 4^5+6\times 7+8+9.$
\item [] $9276=12\times (3^4+5+678+9).$
\item [] $9277=1+2\times 3\times (4^5+6\times (78+9)).$
\item [] $\mathit{9278=(1-2+3)\times (4567+8\times 9).}$
\item [] $9279=(1+2\times (3+456+7\times 8))\times 9.$
\item [] $9280=(1\times 2^3)^4+(5+67)\times 8\times 9.$
\item [] $9281=1+(2+3)\times 4\times (56\times 7+8\times 9).$
\item [] $9282=12^3\times 4+5\times 6\times (7+8\times 9).$
\item [] $9283=1^2+3\times (4\times 5+6)\times 7\times (8+9).$
\item [] $9284=1\times 2\times (3+4567+8\times 9).$
\item [] $9285=1+2\times (3+4567+8\times 9).$
\item [] $9286=(1+2^3)^4+5\times (67\times 8+9).$
\item [] $9287=(1+2^3)\times 4^5+6+7\times 8+9.$
\item [] $9288=12^3+(4+5+6)\times 7\times 8\times 9.$
\item [] $9289=1+2^3\times (45+6+78)\times 9.$
\item [] $9290=1\times 2+3\times 4\times (5\times 6+7\times 8)\times 9.$
\item [] $9291=1+2+3\times 4\times (5\times 6+7\times 8)\times 9.$
\item [] $9292=12^3+4+56\times (7+8)\times 9.$
\item [] $\mathit{9293=1\times 23+(4^5+6)\times (-7+8)\times 9.}$
\item [] $9294=12+3\times (4\times 5+6)\times 7\times (8+9).$
\item [] $9295=(1+2\times 3+4)\times (56+789).$
\item [] $9296=1+(23+4\times 5\times 6)\times (7\times 8+9).$
\item [] $9297=12^3\times 4+5\times (6\times 78+9).$
\item [] $\mathit{9298=1^2+(3+4^5+6)\times (-7+8)\times 9.}$
\item [] $9299=(12\times (34+5+6)+7)\times (8+9).$
\item [] $9300=12+3\times 4\times (5\times 6+7\times 8)\times 9.$
\item[]$\mbox{Decreasing order}$
\item [] $9231=(9\times 8+7)\times 65+4^{(3+2+1)}.$
\item [] $9232=9+87\times (6+5\times 4\times (3+2))+1.$
\item [] $9233=9\times 876+5+4^3\times 21.$
\item [] $9234=9+8+(7+65)\times 4^3\times 2+1.$
\item [] $9235=(9\times 8+7\times 6)\times (5+4)\times 3^2+1.$
\item [] $9236=(9\times 8\times 7\times 6+54)\times 3+2\times 1.$
\item [] $9237=(9\times 8+7+65)\times 4^3+21.$
\item [] $\mathit{9238=((9+8+7)\times 6+5)\times (4^3-2\times 1).}$
\item [] $9239=(9+8)\times 7+6^5+4^3\times 21.$
\item [] $9240=(9+8+76\times 5+43)\times 21.$
\item [] $9241=(9+8\times (7\times 6+5))\times 4\times 3\times 2+1.$
\item [] $9242=9+8\times ((76\times 5+4)\times 3+2)+1.$
\item [] $9243=(9\times 8+7)\times (6\times 5+43\times 2+1).$
\item [] $\mathit{9244=9876-5^4-3\times 2-1.}$
\item [] $9245=(9+8\times (76\times 5+4))\times 3+2\times 1.$
\item [] $9246=9+8\times (76\times 5+4)\times 3+21.$
\item [] $9247=((98+7)\times (6+5)\times 4+3)\times 2+1.$
\item [] $9248=9+(8\times 7+6)\times (5+(4\times 3)^2)+1.$
\item [] $9249=9\times 8+7\times (6\times 5\times 43+21).$
\item [] $9250=(9+8+(7+65)\times 4^3)\times 2\times 1.$
\item [] $9251=(9+8+(7+65)\times 4^3)\times 2+1.$
\item [] $9252=(9+8+765\times 4)\times 3+21.$
\item [] $9253=9\times 8\times 7+6\times (5+4)^3\times 2+1.$
\item [] $9254=9\times 8+765\times 4\times 3+2\times 1.$
\item [] $9255=9\times 8+765\times 4\times 3+2+1.$
\item [] $9256=(9+8\times 76)\times (5+4+3\times 2)+1.$
\item [] $9257=9+8\times 76+5\times (4\times 3)^{(2+1)}.$
\item [] $9258=9+(8\times 76\times 5+43)\times (2+1).$
\item [] $9259=9+(8+7\times 6)\times 5\times (4+32+1).$
\item [] $9260=98+(7\times 654+3)\times 2\times 1.$
\item [] $9261=(9+8+76+54)\times 3\times 21.$
\item [] $9262=(9\times (87+6)+5)\times (4+3\times 2+1).$
\item [] $\mathit{9263=9+8+7-6+5\times 43^2\times 1.}$
\item [] $9264=(9\times 8+76\times 5\times 4\times 3)\times 2\times 1.$
\item [] $9265=(9\times 8+76\times 5\times 4\times 3)\times 2+1.$
\item [] $9266=98\times 7+65\times 4\times (32+1).$
\item [] $\mathit{9267=98-76+5\times 43^2\times 1.}$
\item [] $9268=9\times 8+76\times (5\times 4\times 3\times 2+1).$
\item [] $9269=9+(8+7\times 65)\times 4\times (3+2\times 1).$
\item [] $9270=9+(87+6+54)\times 3\times 21.$
\item [] $9271=9\times 8+7\times (654+3)\times 2+1.$
\item [] $\mathit{9272=98\times (7+6)+(5\times 4)^3-2\times 1.}$
\item [] $9273=9\times 8+765\times 4\times 3+21.$
\item [] $9274=9\times 8\times 7\times 6+5^4\times (3^2+1).$
\item [] $9275=9+8+7+6+5\times 43^2\times 1.$
\item [] $9276=9+8+7+6+5\times 43^2+1.$
\item [] $9277=98\times (7+6)+(5\times 4)^3+2+1.$
\item [] $9278=9+8+7\times (6+54+3)\times 21.$
\item [] $9279=9\times (8\times 7+654+321).$
\item [] $9280=98+765\times 4\times 3+2\times 1.$
\item [] $9281=98+765\times 4\times 3+2+1.$
\item [] $9282=98\times (76+5)+4^3\times 21.$
\item [] $9283=(98+7\times 6+5)\times 4^3+2+1.$
\item [] $9284=(9+8)\times 7\times (6+5\times 4)\times 3+2\times 1.$
\item [] $9285=(9+8)\times 7\times (6+5\times 4)\times 3+2+1.$
\item [] $9286=(987+6)\times 5+4321.$
\item [] $9287=9+(8+765)\times 4\times 3+2\times 1.$
\item [] $9288=9+(8+765)\times 4\times 3+2+1.$
\item [] $9289=9\times 8+(7+65)\times 4^3\times 2+1.$
\item [] $9290=9\times (8\times 7+6\times 5)\times 4\times 3+2\times 1.$
\item [] $9291=9\times (8\times 7+6\times 5)\times 4\times 3+2+1.$
\item [] $\mathit{9292=98\times 76-5+43^2\times 1.}$
\item [] $\mathit{9293=98\times 76-5+43^2+1.}$
\item [] $9294=98+76\times (5\times 4\times 3\times 2+1).$
\item [] $9295=98\times (7+6)+(5\times 4)^3+21.$
\item [] $9296=98+7\times (654+3)\times 2\times 1.$
\item [] $9297=98+7\times (654+3)\times 2+1.$
\item [] $9298=9\times (87+(6+5)\times 43\times 2)+1.$
\item [] $9299=98+765\times 4\times 3+21.$
\item [] $9300=(9+8+76)\times 5\times 4\times (3+2)\times 1.$
\item[]$\mbox{Increasing order}$
\item [] $9301=1\times 23\times (4+56\times 7+8)+9.$
\item [] $9302=1+23\times (4+56\times 7+8)+9.$
\item [] $9303=1+2\times (3+4\times (5+(6+7)\times 89)).$
\item [] $9304=1\times 2^3\times (4^5+67+8\times 9).$
\item [] $9305=1+2^3\times (4^5+67+8\times 9).$
\item [] $9306=1+2^{(3+4)}\times (5+67)+89.$
\item [] $9307=1+(2+3\times 4\times (5\times 6+7\times 8))\times 9.$
\item [] $9308=(1+2)\times (3\times (4^5+6)+7)+8+9.$
\item [] $9309=(1+2^3)\times 4^5+6+78+9.$
\item [] $9310=1\times 2\times (3^4\times 56+7\times (8+9)).$
\item [] $9311=1+2\times (3^4\times 56+7\times (8+9)).$
\item [] $9312=12\times (3\times (45+6)+7\times 89).$
\item [] $9313=1^2+3\times (45\times 67+89).$
\item [] $9314=1\times 2+3\times (45\times 67+89).$
\item [] $9315=(12+345+678)\times 9.$
\item [] $9316=1+23\times (45\times 6+(7+8)\times 9).$
\item [] $9317=1\times 2+3\times (4+5+6\times 7\times 8)\times 9.$
\item [] $9318=1\times 2\times (3+4567+89).$
\item [] $9319=1+2\times (3+4567+89).$
\item [] $9320=1\times 2^3\times (4^5+6+(7+8)\times 9).$
\item [] $9321=12\times (34+56+7)\times 8+9.$
\item [] $\mathit{9322=-1\times 23+(4+5+6)\times 7\times 89.}$
\item [] $\mathit{9323=1-23+(4+5+6)\times 7\times 89.}$
\item [] $9324=12+3\times (45\times 67+89).$
\item [] $9325=1+(2+34)\times (5\times (6\times 7+8)+9).$
\item [] $\mathit{9326=(1+2\times 3\times 4\times 56)\times 7-89.}$
\item [] $9327=12+3\times (4+5+6\times 7\times 8)\times 9.$
\item [] $9328=1\times 2^3\times (4+5+(6+7)\times 89).$
\item [] $9329=1+2^3\times (4+5+(6+7)\times 89).$
\item [] $9330=(1+2^3)\times 4^5+6\times 7+8\times 9.$
\item [] $9331=1\times 2^3\times 4^5+67\times (8+9).$
\item [] $9332=1+2^3\times 4^5+67\times (8+9).$
\item [] $9333=1+(2+3)\times 4^5+6\times 78\times 9.$
\item [] $9334=1+((2+34\times 5\times 6+7)+8)\times 9.$
\item [] $9335=(1+2+3\times 45)\times 67+89.$
\item [] $9336=(1+23)\times (45\times 6+7\times (8+9)).$
\item [] $\mathit{9337=1+2\times 3\times 4\times 56\times 7-8\times 9.}$
\item [] $\mathit{9338=1-2^3+(4+5+6)\times 7\times 89.}$
\item [] $9339=(1+2\times 3+4)\times (56\times (7+8)+9).$
\item [] $\mathit{9340=1-2\times 3+(4+5+6)\times 7\times 89.}$
\item [] $9341=12\times (3^4+5\times 6)\times 7+8+9.$
\item [] $9342=(1\times 2^3\times 45+678)\times 9.$
\item [] $9343=1+2\times (3^4\times 56+(7+8)\times 9).$
\item [] $9344=(1^2+3)\times 4\times (567+8+9).$
\item [] $9345=(1\times 2\times 34+5\times 6+7)\times 89.$
\item [] $9346=1^{23}+(4+5+6)\times 7\times 89.$
\item [] $9347=(12^3+4)\times 5+678+9.$
\item [] $9348=1^2\times 3+(4+5+6)\times 7\times 89.$
\item [] $9349=1^2+3+(4+5+6)\times 7\times 89.$
\item [] $9350=1+2^3\times 4^5+(6+7)\times 89.$
\item [] $9351=(1+2^3\times 45+678)\times 9.$
\item [] $9352=1+2\times 3+(4+5+6)\times 7\times 89.$
\item [] $9353=1\times 2^3+(4+5+6)\times 7\times 89.$
\item [] $9354=1+2^3+(4+5+6)\times 7\times 89.$
\item [] $9355=(1+2^3)\times 4^5+67+8\times 9.$
\item [] $\mathit{9356=1\times 2+3+4\times 5\times 6\times 78-9.}$
\item [] $9357=(1+2^3)\times 4^5+6+(7+8)\times 9.$
\item [] $\mathit{9358=1+2\times 3+4\times 5\times 6\times 78-9.}$
\item [] $9359=(1+2)\times (3\times 4^5+6\times 7)+8+9.$
\item [] $9360=12+3+(4+5+6)\times 7\times 89.$
\item [] $9361=1+2\times (3+4\times 5+6\times 7)\times 8\times 9.$
\item [] $9362=1\times 2+3\times 4\times 5\times (67+89).$
\item [] $9363=1+2+3\times 4\times 5\times (67+89).$
\item [] $\mathit{9364=1-2\times 3+4\times 5\times 6\times 78+9.}$
\item [] $9365=123\times (4+5+67)+8+9.$
\item [] $9366=(1+2^3)\times (4^5+6)+7+89.$
\item [] $9367=(12+(3+4)\times (5+6)\times 7)\times (8+9).$
\item [] $9368=1\times 23+(4+5+6)\times 7\times 89.$
\item [] $9369=1^{23}\times 4\times 5\times 6\times 78+9.$
\item [] $9370=1^{23}+4\times 5\times 6\times 78+9.$
\item[]$\mbox{Decreasing order}$
\item [] $9301=(98+7\times 6+5)\times 4^3+21.$
\item [] $9302=98\times 76+5+43^2\times 1.$
\item [] $9303=98\times 76+5+43^2+1.$
\item [] $9304=9+8+7\times 6+5\times 43^2\times 1.$
\item [] $9305=9+8+7\times 6+5\times 43^2+1.$
\item [] $9306=9+(8+765)\times 4\times 3+21.$
\item [] $9307=987+65\times 4^3\times 2\times 1.$
\item [] $9308=987+65\times 4\times 32+1.$
\item [] $9309=9+8+7\times 6+5\times (43^2+1).$
\item [] $9310=98\times (76+5+4+3^2+1).$
\item [] $9311=98\times (76+5+4\times 3+2)+1.$
\item [] $9312=(9+87)\times (6+5+43\times 2\times 1).$
\item [] $9313=9+8\times (76+543\times 2+1).$
\item [] $9314=98+(7+65)\times 4^3\times 2\times 1.$
\item [] $9315=98+(7+65)\times 4\times 32+1.$
\item [] $9316=9+8\times 7+6+5\times 43^2\times 1.$
\item [] $9317=9+8\times 7+6+5\times 43^2+1.$
\item [] $9318=(98\times 7\times 6+543)\times 2\times 1.$
\item [] $9319=(98\times 7\times 6+543)\times 2+1.$
\item [] $9320=9\times (8+7)\times (65+4)+3+2\times 1.$
\item [] $9321=9\times (8+7\times 6)\times 5\times 4+321.$
\item [] $9322=9+(8+76\times 5)\times 4\times 3\times 2+1.$
\item [] $\mathit{9323=9-8+76+5\times 43^2+1.}$
\item [] $9324=987\times 6+54\times 3\times 21.$
\item [] $9325=98\times 76+5^4\times 3+2\times 1.$
\item [] $9326=98\times 76+5^4\times 3+2+1.$
\item [] $\mathit{9327=-9+(8+7)\times 6+5\times 43^2+1.}$
\item [] $9328=(9\times 87+65)\times (4+3\times 2+1).$
\item [] $9329=(9\times 8\times 7+65\times 4^3)\times 2+1.$
\item [] $9330=9\times 8+7+6+5\times 43^2\times 1.$
\item [] $9331=9\times 8+7+6+5\times 43^2+1.$
\item [] $9332=98\times 76+(5^4+3)\times (2+1).$
\item [] $9333=987+(6+5\times 4)\times 321.$
\item [] $9334=(9+8)\times (7+6\times (5+4))\times 3^2+1.$
\item [] $9335=9\times 8+7+6+5\times (43^2+1).$
\item [] $9336=9\times (8\times 7\times 6+5+4)\times 3+21.$
\item [] $\mathit{9337=9+8+76+5\times 43^2-1.}$
\item [] $9338=9+8+76+5\times 43^2\times 1.$
\item [] $9339=9+8+76+5\times 43^2+1.$
\item [] $9340=9\times 8\times 7+(6\times 5+4^3)^2\times 1.$
\item [] $9341=9+8+(7+6\times 5)\times 4\times 3\times 21.$
\item [] $9342=9\times (8\times 76+5\times 43\times 2\times 1).$
\item [] $9343=9\times 876+(5+4)^3\times 2+1.$
\item [] $9344=98\times 76+5^4\times 3+21.$
\item [] $9345=(98+7)\times (65+4\times 3\times 2)\times 1.$
\item [] $9346=(98+7)\times (65+4\times 3\times 2)+1.$
\item [] $9347=9+87+6+5\times 43^2\times 1.$
\item [] $9348=9+87+6+5\times 43^2+1.$
\item [] $9349=9+(8+7)\times 6+5\times (43^2+1).$
\item [] $9350=(9+8)\times ((7+6\times (5+4))\times 3^2+1).$
\item [] $9351=9\times (87+6+5^4+321).$
\item [] $9352=9+87+6+5\times (43^2+1).$
\item [] $9353=9+(8+(7+6)\times 5)\times 4\times 32\times 1.$
\item [] $9354=9+(8+(7+6)\times 5)\times 4\times 32+1.$
\item [] $\mathit{9355=98+7+6+5\times 43^2-1.}$
\item [] $9356=98+7+6+5\times 43^2\times 1.$
\item [] $9357=98+7+6+5\times 43^2+1.$
\item [] $9358=(98+7\times 654+3)\times 2\times 1.$
\item [] $9359=9\times 8+7\times 6+5\times 43^2\times 1.$
\item [] $9360=9\times 8+7\times 6+5\times 43^2+1.$
\item [] $9361=98+7+6+5\times (43^2+1).$
\item [] $9362=9\times 8\times (7+6\times 5\times 4+3)+2\times 1.$
\item [] $9363=9\times 87+65\times 4\times (32+1).$
\item [] $9364=9\times 8+7\times 6+5\times (43^2+1).$
\item [] $9365=98\times 7+6^5+43\times 21.$
\item [] $9366=9\times 8\times (76+54)+3+2+1.$
\item [] $9367=9\times 8\times (76+54)+3\times 2+1.$
\item [] $\mathit{9368=9-8+7+65\times (4\times 3)^2\times 1.}$
\item [] $9369=9+8\times (76+54)\times 3^2\times 1.$
\item [] $9370=9+8\times (76+54)\times 3^2+1.$
\item[]$\mbox{Increasing order}$
\item [] $9371=1\times 2+3\times (4+5\times 67+8)\times 9.$
\item [] $9372=1^2\times 3+4\times 5\times 6\times 78+9.$
\item [] $9373=1^2+3+4\times 5\times 6\times 78+9.$
\item [] $9374=1\times 2+3+4\times 5\times 6\times 78+9.$
\item [] $9375=1+2+3+4\times 5\times 6\times 78+9.$
\item [] $9376=1+2\times 3+4\times 5\times 6\times 78+9.$
\item [] $9377=1\times 2^3+4\times 5\times 6\times 78+9.$
\item [] $9378=1+2^3+4\times 5\times 6\times 78+9.$
\item [] $9379=1+2\times (3\times (4^5+67\times 8)+9).$
\item [] $9380=(1+2)\times (3\times (4^5+6)+7)+89.$
\item [] $9381=12+3\times (4+5\times 67+8)\times 9.$
\item [] $\mathit{9382=-1\times 23+(4^5+6+7+8)\times 9.}$
\item [] $\mathit{9383=(123\times (4+5)+67)\times 8-9.}$
\item [] $9384=12+3+4\times 5\times 6\times 78+9.$
\item [] $9385=1+2\times (3+(4+56)\times 78+9).$
\item [] $\mathit{9386=1-2-3\times (4-56\times 7\times 8)-9.}$
\item [] $9387=(12\times 3^4+56+7+8)\times 9.$
\item [] $9388=1+(2+3\times (4+5\times 67+8))\times 9.$
\item [] $9389=(1+2^3)\times (4^5+6)+7\times (8+9).$
\item [] $\mathit{9390=-1+2\times 3\times 4\times 56\times 7-8-9.}$
\item [] $\mathit{9391=1\times 2\times 3\times 4\times 56\times 7-8-9.}$
\item [] $9392=1\times 23+4\times 5\times 6\times 78+9.$
\item [] $9393=1+23+4\times 5\times 6\times 78+9.$
\item [] $9394=1+23\times (4+5+6\times 7)\times 8+9.$
\item [] $\mathit{9395=1\times 2-3+4\times (5\times 6\times 78+9).}$
\item [] $9396=(1+2)^3+4\times 5\times 6\times 78+9.$
\item [] $9397=(1+23+4)\times 5\times 67+8+9.$
\item [] $9398=1\times 2+3\times (45\times 6+78)\times 9.$
\item [] $9399=1+2+3\times (45\times 6+78)\times 9.$
\item [] $9400=1^2+3+4\times (5\times 6\times 78+9).$
\item [] $9401=1\times 2+3+4\times (5\times 6\times 78+9).$
\item [] $9402=1+2+3+4\times (5\times 6\times 78+9).$
\item [] $9403=1+2\times 3+4\times (5\times 6\times 78+9).$
\item [] $9404=1\times 2^3+4\times (5\times 6\times 78+9).$
\item [] $9405=12\times 3+4\times 5\times 6\times 78+9.$
\item [] $9406=1^{23}+(4^5+6+7+8)\times 9.$
\item [] $9407=1\times 23\times (4+(5\times 6+7+8)\times 9).$
\item [] $9408=12+3\times (45\times 6+78)\times 9.$
\item [] $9409=1+(2\times 34+5\times 6)\times (7+89).$
\item [] $9410=1\times 2+3+(4^5+6+7+8)\times 9.$
\item [] $9411=12+3+4\times (5\times 6\times 78+9).$
\item [] $9412=1+2\times 3+(4^5+6+7+8)\times 9.$
\item [] $9413=12\times (3^4+5\times 6)\times 7+89.$
\item [] $9414=1\times 2\times 3\times (4^5+67\times 8+9).$
\item [] $9415=1+2\times 3\times (4^5+67\times 8+9).$
\item [] $9416=(12^3+4)\times 5+(6+78)\times 9.$
\item [] $9417=12\times (3+4)\times (56+7\times 8)+9.$
\item [] $9418=1+(2^3+4\times 5)\times 6\times 7\times 8+9.$
\item [] $9419=1\times 23+4\times (5\times 6\times 78+9).$
\item [] $9420=123\times (4+5+67)+8\times 9.$
\item [] $\mathit{9421=1+2-34\times (5+6\times (-7\times 8+9)).}$
\item [] $9422=(1+2)\times 3\times (4^5+6+7)+89.$
\item [] $9423=(1+2)^3+4\times (5\times 6\times 78+9).$
\item [] $9424=1\times 2\times (3+(45\times 6+7)\times (8+9)).$
\item [] $9425=1\times 2\times 3\times 4\times 56\times 7+8+9.$
\item [] $9426=1+2\times 3\times 4\times 56\times 7+8+9.$
\item [] $9427=1+2\times ((3+45\times (6+7))\times 8+9).$
\item [] $9428=1\times 23+(4^5+6+7+8)\times 9.$
\item [] $9429=1+23+(4^5+6+7+8)\times 9.$
\item [] $9430=1^2+3\times (4+56\times 7\times 8)+9.$
\item [] $9431=1+23\times (4\times 5+6\times (7\times 8+9)).$
\item [] $9432=(123+45\times 6)\times (7+8+9).$
\item [] $9433=1+2^3\times (4+5)\times (6\times 7+89).$
\item [] $9434=1^2\times (34+5+67)\times 89.$
\item [] $9435=(123\times 4+56+7)\times (8+9).$
\item [] $9436=1+(2\times 3+4+5)\times (6+7\times 89).$
\item [] $9437=123\times (4+5+67)+89.$
\item [] $9438=(1+2)\times (3+4+56\times 7\times 8)+9.$
\item [] $9439=1+2\times (3+(4+5\times (6+7)\times 8)\times 9).$
\item [] $9440=1\times 2^3\times (4^5+67+89).$
\item[]$\mbox{Decreasing order}$
\item [] $9371=9\times (87+65\times 4)\times 3+2\times 1.$
\item [] $9372=987+65\times 43\times (2+1).$
\item [] $9373=((9+8)\times 7\times 6+5^4)\times (3\times 2+1).$
\item [] $9374=(9\times 8+7+6\times 5)\times 43\times 2\times 1.$
\item [] $9375=(9\times 8+7+6\times 5)\times 43\times 2+1.$
\item [] $\mathit{9376=-98+(7+6)\times (5+4)^3-2-1.}$
\item [] $9377=9+8+(7+6)\times 5\times (4\times 3)^2\times 1.$
\item [] $9378=9\times (8+7)\times (65+4)+3\times 21.$
\item [] $9379=(9+8\times (7+6))\times (5\times 4+3\times 21).$
\item [] $9380=98\times 7+6\times (5+4^3)\times 21.$
\item [] $9381=9\times 8\times (7+6\times 5\times 4+3)+21.$
\item [] $\mathit{9382=-987+6\times 54\times 32+1.}$
\item [] $9383=9+(8\times (7+6)+5)\times 43\times 2\times 1.$
\item [] $9384=9+8+7+65\times (4\times 3)^2\times 1.$
\item [] $9385=98+7\times 6+5\times 43^2\times 1.$
\item [] $9386=98+7\times 6+5\times 43^2+1.$
\item [] $9387=(9\times 87+65\times 4)\times 3^2\times 1.$
\item [] $9388=(9\times 87+65\times 4)\times 3^2+1.$
\item [] $9389=(9+8+7)\times 6+5\times 43^2\times 1.$
\item [] $9390=9\times (87+65\times 4)\times 3+21.$
\item [] $9391=9\times (8+7)+6+5\times (43^2+1).$
\item [] $9392=9\times 8\times (76+54)+32\times 1.$
\item [] $9393=9\times 8+76+5\times 43^2\times 1.$
\item [] $9394=9\times 8+76+5\times 43^2+1.$
\item [] $9395=98\times 7\times (6+5)+43^2\times 1.$
\item [] $9396=98\times 7\times (6+5)+43^2+1.$
\item [] $9397=9\times (8\times 7+6+54)\times 3^2+1.$
\item [] $9398=(9\times 8+765\times 4)\times 3+2\times 1.$
\item [] $9399=(9\times 8+765\times 4)\times 3+2+1.$
\item [] $\mathit{9400=-98+(7+6)\times (5+4)^3+21.}$
\item [] $9401=(98+76)\times 54+3+2\times 1.$
\item [] $9402=(98+76)\times 54+3+2+1.$
\item [] $9403=(98+76)\times 54+3\times 2+1.$
\item [] $9404=9+(8+7\times 6\times 5)\times 43+21.$
\item [] $9405=(9+8+765)\times 4\times 3+21.$
\item [] $9406=(98+76)\times 54+3^2+1.$
\item [] $9407=9+87\times (65+43)+2\times 1.$
\item [] $9408=9+87\times (65+43)+2+1.$
\item [] $9409=(9+8+7+6\times 5+43)^2\times 1.$
\item [] $9410=(9+8+7+6\times 5+43)^2+1.$
\item [] $9411=9+(87\times 6\times (5+4)+3)\times 2\times 1.$
\item [] $9412=9+(87\times 6\times (5+4)+3)\times 2+1.$
\item [] $9413=98\times (76+5\times 4)+3+2\times 1.$
\item [] $9414=98\times (76+5\times 4)+3+2+1.$
\item [] $9415=98\times (76+5\times 4)+3\times 2+1.$
\item [] $\mathit{9416=98\times (76+5\times 4)+3^2-1.}$
\item [] $9417=98\times (76+5\times 4)+3^2\times 1.$
\item [] $9418=98\times (76+5\times 4)+3^2+1.$
\item [] $9419=98+76+5\times 43^2\times 1.$
\item [] $9420=98+76+5\times 43^2+1.$
\item [] $9421=(9\times 87\times 6+5+4+3)\times 2+1.$
\item [] $9422=98+(7+6\times 5)\times 4\times 3\times 21.$
\item [] $9423=9\times 8\times (76+54)+3\times 21.$
\item [] $9424=(987+6+54)\times 3^2+1.$
\item [] $9425=9+8\times 7+65\times (4\times 3)^2\times 1.$
\item [] $9426=9+87\times (65+43)+21.$
\item [] $\mathit{9427=-98\times 76+5^4\times 3^{(2+1)}.}$
\item [] $9428=(98+76)\times 54+32\times 1.$
\item [] $9429=(98+76)\times 54+32+1.$
\item [] $9430=(9\times 87\times 6+5+4\times 3)\times 2\times 1.$
\item [] $9431=(9\times 87\times 6+5+4\times 3)\times 2+1.$
\item [] $9432=98\times (76+5\times 4)+3+21.$
\item [] $9433=9\times 8+(7+6)\times 5\times (4\times 3)^2+1.$
\item [] $9434=98\times 7+6\times (5+4)^3\times 2\times 1.$
\item [] $9435=98\times 7+6\times (5+4)^3\times 2+1.$
\item [] $\mathit{9436=-(9+8)\times 7+65\times (4+3)\times 21.}$
\item [] $\mathit{9437=-9+87+65\times (4\times 3)^2-1.}$
\item [] $9438=9+8\times 7\times (6+54\times 3)+21.$
\item [] $9439=9\times 8+7+65\times (4\times 3)^2\times 1.$
\item [] $9440=98\times (76+5\times 4)+32\times 1.$
\item[]$\mbox{Increasing order}$
\item [] $9441=1+2^3\times (4^5+67+89).$
\item [] $9442=1\times 2\times ((3\times 4\times 56\times 7+8)+9).$
\item [] $9443=(1\times 2^{(3+4)}+5)\times (6+7\times 8+9).$
\item [] $9444=1+(2^{(3+4)}+5)\times (6+7\times 8+9).$
\item [] $9445=(1+2\times 3)\times (4\times 5\times 67+8)+9.$
\item [] $9446=12+(34+5+67)\times 89.$
\item [] $9447=(1^{23}+4\times 5\times 6)\times 78+9.$
\item [] $9448=1^2+3\times (4+56\times 7\times 8+9).$
\item [] $9449=12\times 3\times 4\times 5\times (6+7)+89.$
\item [] $9450=1+2+3\times (4+56\times 7\times 8+9).$
\item [] $9451=1+(2+3+45)\times (6+7+8)\times 9.$
\item [] $9452=(1+23+4)\times 5\times 67+8\times 9.$
\item [] $9453=1\times 23\times (4+5\times 67+8\times 9).$
\item [] $9454=12^3+4+(5+6)\times 78\times 9.$
\item [] $9455=(12^3+4)\times 5+6+789.$
\item [] $9456=(1+2)\times (3+4+56\times 7\times 8+9).$
\item [] $9457=1+(2+(3+4)\times 56)\times (7+8+9).$
\item [] $\mathit{9458=1-23+4\times 5\times 6\times (7+8\times 9).}$
\item [] $9459=(1\times 2\times 3+4^5+6+7+8)\times 9.$
\item [] $9460=1+(2\times 3+4^5+6+7+8)\times 9.$
\item [] $\mathit{9461=1-2\times (3-4+5-6\times 789).}$
\item [] $9462=12^3\times 4+5\times (6+7\times 8\times 9).$
\item [] $9463=1+(2+3^4)\times (5\times (6+7+8)+9).$
\item [] $\mathit{9464=1+2^3\times 4\times (5\times 6+7)\times 8-9.}$
\item [] $9465=(12\times 3^4+5\times 6\times 7)\times 8+9.$
\item [] $9466=1+2\times 3\times (4\times 56\times 7+8)+9.$
\item [] $9467=1\times (2+3)\times 45\times 6\times 7+8+9.$
\item [] $9468=(1^{234}+5+6)\times 789.$
\item [] $9469=(1+23+4)\times 5\times 67+89.$
\item [] $9470=1+(2+3)\times 4\times (5+6\times 78)+9.$
\item [] $9471=(1+2)\times (3\times 4+56\times 7\times 8+9).$
\item [] $9472=1\times 2^3\times (45+67\times (8+9)).$
\item [] $9473=(1\times 234\times 5+6+7)\times 8+9.$
\item [] $9474=(1+(2+3)\times 45\times 6)\times 7+8+9.$
\item [] $9475=1+2\times (3\times (4\times 56\times 7+8)+9).$
\item [] $\mathit{9476=1-2-3+4\times 5\times 6\times (7+8\times 9).}$
\item [] $9477=1^2\times 3^4\times (5\times 6+78+9).$
\item [] $9478=1^2+3^4\times (5\times 6+78+9).$
\item [] $9479=1\times 2+3^4\times (5\times 6+78+9).$
\item [] $9480=1\times 2\times 3\times 4\times 56\times 7+8\times 9.$
\item [] $9481=1+2\times 3\times 4\times 56\times 7+8\times 9.$
\item [] $9482=1+2^3\times 4\times (5\times 6+7)\times 8+9.$
\item [] $9483=1^2\times 3+4\times 5\times 6\times (7+8\times 9).$
\item [] $9484=1^2+3+4\times 5\times 6\times (7+8\times 9).$
\item [] $9485=1\times 2+3+4\times 5\times 6\times (7+8\times 9).$
\item [] $9486=1\times 2\times 3+4\times 5\times 6\times (7+8\times 9).$
\item [] $9487=1+2\times 3+4\times 5\times 6\times (7+8\times 9).$
\item [] $9488=1\times 2^3+4\times 5\times 6\times (7+8\times 9).$
\item [] $9489=12+3^4\times (5\times 6+78+9).$
\item [] $9490=((1+23+4)\times 5+6)\times (7\times 8+9).$
\item [] $\mathit{9491=(-12+3\times 456)\times 7+8-9.}$
\item [] $9492=123+4\times 5\times 6\times 78+9.$
\item [] $9493=1+2\times (3+4+5+6\times 789).$
\item [] $\mathit{9494=(1^2+3\times 456)\times 7-89.}$
\item [] $9495=12+3+4\times 5\times 6\times (7+8\times 9).$
\item [] $9496=(1^2+3)\times (4+5\times 6\times (7+8\times 9)).$
\item [] $9497=1\times 2\times 3\times 4\times 56\times 7+89.$
\item [] $9498=1+2\times 3\times 4\times 56\times 7+89.$
\item [] $9499=1\times 23\times (4+56\times 7+8+9).$
\item [] $9500=1+23\times (4+56\times 7+8+9).$
\item [] $9501=1\times (2+3+4+5)\times 678+9.$
\item [] $9502=1+(2+3+4+5)\times 678+9.$
\item [] $9503=1+2\times (3\times 4+5+6\times 789).$
\item [] $9504=12\times (34+56+78\times 9).$
\item [] $9505=1+(2\times 3)^4\times 5+6\times 7\times 8\times 9.$
\item [] $9506=(1\times 2+3\times 4)\times (56+7\times 89).$
\item [] $9507=1+2\times (3+4)\times (56+7\times 89).$
\item [] $\mathit{9508=(1-2+3)\times (4\times 5+6\times 789).}$
\item [] $9509=(1+(2+3)\times 45)\times 6\times 7+8+9.$
\item [] $9510=1\times 2\times 3\times ((4\times 56\times 7+8)+9).$
\item[]$\mbox{Decreasing order}$
\item [] $9441=98\times (76+5\times 4)+32+1.$
\item [] $9442=(9\times 87\times 6+5\times 4+3)\times 2\times 1.$
\item [] $9443=(9\times 87\times 6+5\times 4+3)\times 2+1.$
\item [] $9444=9\times 8\times 7\times 6+5\times 4\times 321.$
\item [] $9445=9\times 8+(7+6)\times (5\times (4\times 3)^2+1).$
\item [] $9446=(9\times 8\times (7+6)\times 5+43)\times 2\times 1.$
\item [] $9447=(9\times 8\times (7+6)\times 5+43)\times 2+1.$
\item [] $\mathit{9448=9876+5-432-1.}$
\item [] $9449=9+8+7+65\times ((4\times 3)^2+1).$
\item [] $9450=9\times (87+65\times 4+3)\times (2+1).$
\item [] $9451=9\times (8+7+6)\times (5+43+2)+1.$
\item [] $9452=(9+8)\times (76+5\times 4\times (3+21)).$
\item [] $\mathit{9453=(98+765\times 4)\times 3-21.}$
\item [] $9454=(9+8)\times (7+6+543)+2\times 1.$
\item [] $9455=(9+8)\times (7+6+543)+2+1.$
\item [] $9456=(9\times 8+7)\times 65+4321.$
\item [] $9457=9+87+65\times (4\times 3)^2+1.$
\item [] $9458=98+(7+6)\times 5\times (4\times 3)^2\times 1.$
\item [] $9459=(98+76)\times 54+3\times 21.$
\item [] $9460=9+8+7\times (65+4\times 321).$
\item [] $9461=(9\times 8+7\times 65\times 4)\times (3+2)+1.$
\item [] $9462=9\times 87+6^5+43\times 21.$
\item [] $9463=(9+8)\times 7\times 65+(4\times 3)^{(2+1)}.$
\item [] $9464=((98+7)\times 6\times 5+4)\times 3+2\times 1.$
\item [] $9465=98+7+65\times (4\times 3)^2\times 1.$
\item [] $9466=98+7+65\times (4\times 3)^2+1.$
\item [] $9467=9+8+7\times 6\times 5\times (43+2)\times 1.$
\item [] $9468=(987+65)\times (4+3+2\times 1).$
\item [] $9469=(987+65)\times (4+3+2)+1.$
\item [] $9470=(9\times 87+6)\times (5+4+3)+2\times 1.$
\item [] $9471=98\times (76+5\times 4)+3\times 21.$
\item [] $9472=(9\times 8+76)\times (54+3^2+1).$
\item [] $9473=(9+8)\times (7+6+543)+21.$
\item [] $9474=9+8\times (7+6)\times (5+43\times 2)+1.$
\item [] $\mathit{9475=(9\times 8+7)\times 6\times 5\times 4-3-2\times 1.}$
\item [] $9476=(98+765\times 4)\times 3+2\times 1.$
\item [] $9477=(98+765\times 4)\times 3+2+1.$
\item [] $9478=9\times (87+6\times 5)\times (4+3+2)+1.$
\item [] $9479=(9+8)\times 7+65\times (4\times 3)^2\times 1.$
\item [] $9480=(987+6\times 5^4+3)\times 2\times 1.$
\item [] $9481=(987+6\times 5^4+3)\times 2+1.$
\item [] $9482=9+(8+(7+65)\times 4)\times 32+1.$
\item [] $9483=(98+765\times 4+3)\times (2+1).$
\item [] $9484=9+(8+7)\times (6+5^4)+3^2+1.$
\item [] $9485=(9\times 8+7)\times 6\times 5\times 4+3+2\times 1.$
\item [] $9486=(9\times 8+7)\times 6\times 5\times 4+3+2+1.$
\item [] $9487=(9\times 8+7)\times 6\times 5\times 4+3\times 2+1.$
\item [] $9488=(98\times 7\times 6+5^4+3)\times 2\times 1.$
\item [] $9489=(98\times 7\times 6+5^4+3)\times 2+1.$
\item [] $9490=(9\times 8+7)\times 6\times 5\times 4+3^2+1.$
\item [] $9491=(9+8\times 7)\times (6\times 5+43)\times 2+1.$
\item [] $9492=(9\times 87\times 6+5+43)\times 2\times 1.$
\item [] $9493=(98+765)\times (4+3\times 2+1).$
\item [] $9494=9+8+(7+6)\times ((5+4)\times 3)^2\times 1.$
\item [] $9495=(98+765\times 4)\times 3+21.$
\item [] $9496=9+8+(7+6)\times (5+4)^3+2\times 1.$
\item [] $9497=9+8+(7+6)\times (5+4)^3+2+1.$
\item [] $9498=9+(8+7)\times (6+5^4)+3+21.$
\item [] $\mathit{9499=9-8+(7+6)\times (5+4)^3+21.}$
\item [] $9500=(9\times (87\times 6+5)+4+3)\times 2\times 1.$
\item [] $9501=9+(8+76\times 5+4^3)\times 21.$
\item [] $9502=9+8+7\times (6+5+4^3\times 21).$
\item [] $9503=(9+8)\times (7+6+543+2+1).$
\item [] $9504=9\times 8\times (7+6+5+4)\times 3\times 2\times 1.$
\item [] $9505=9\times 8\times (7+65)+4321.$
\item [] $9506=98\times (7+65+4\times 3\times 2+1).$
\item [] $9507=9\times (8\times 7+6)\times (5+4\times 3)+21.$
\item [] $\mathit{9508=9-8\times 7+65\times (4+3)\times 21.}$
\item [] $9509=9\times 876+5\times (4+321).$
\item [] $9510=(9\times 87\times 6+54+3)\times 2\times 1.$
\item[]$\mbox{Increasing order}$
\item [] $9511=1+2\times 3\times (4\times 56\times 7+8+9).$
\item [] $\mathit{9512=-123-4+567\times (8+9).}$
\item [] $9513=12\times 3\times (4\times 5+6+7)\times 8+9.$
\item [] $9514=1^2+3\times (4+56\times 7)\times 8+9.$
\item [] $9515=1+2\times (3+4\times 5+6\times 789).$
\item [] $9516=12\times 3+4\times 5\times 6\times (7+8\times 9).$
\item [] $\mathit{9517=-1\times 23+4\times 5\times (6\times 78+9).}$
\item [] $\mathit{9518=(1\times 2+3\times 456)\times 7-8\times 9.}$
\item [] $9519=123+4\times (5\times 6\times 78+9).$
\item [] $9520=(12+3\times 4+56)\times 7\times (8+9).$
\item [] $9521=1+2\times (3+45\times 6+7)\times (8+9).$
\item [] $9522=12^3\times 4+5\times 6\times (78+9).$
\item [] $9523=(1^2+34+5+67)\times 89.$
\item [] $\mathit{9524=-1+(23\times 4+5\times 6)\times 78+9.}$
\item [] $9525=12+3\times (4+56\times 7)\times 8+9.$
\item [] $9526=1+(23\times 4+5\times 6)\times 78+9.$
\item [] $\mathit{9527=((1-2\times 3)^4+567)\times 8-9.}$
\item [] $9528=123+(4^5+6+7+8)\times 9.$
\item [] $9529=1\times (2+3)\times (4+5\times 6)\times 7\times 8+9.$
\item [] $9530=1\times 2\times (3+4^5+6\times 7\times 89).$
\item [] $9531=(1+2)\times 3\times (45\times 6+789).$
\item [] $9532=1^2+3\times ((4+56\times 7)\times 8+9).$
\item [] $9533=1\times 2+3\times ((4+56\times 7)\times 8+9).$
\item [] $9534=1+2+3\times ((4+56\times 7)\times 8+9).$
\item [] $9535=1\times (2+3)\times (45\times 6\times 7+8+9).$
\item [] $9536=(1\times (2\times 3)^4+56)\times 7+8\times 9.$
\item [] $9537=(1+2)\times (34+56\times 7\times 8+9).$
\item [] $9538=1\times 2\times ((3+4)\times 5+6\times 789).$
\item [] $9539=1\times (2+3)\times 45\times 6\times 7+89.$
\item [] $9540=1+(2+3)\times 45\times 6\times 7+89.$
\item [] $9541=12^3+4^5+6789.$
\item [] $9542=1\times 2+(3+4+5)\times (6+789).$
\item [] $9543=1+2+(3+4+5)\times (6+789).$
\item [] $9544=1^2+3+4\times 5\times (6\times 78+9).$
\item [] $9545=1\times 2+3+4\times 5\times (6\times 78+9).$
\item [] $9546=1+2+3+4\times 5\times (6\times 78+9).$
\item [] $9547=1+2\times (34+5+6\times 789).$
\item [] $9548=1\times 2^3+4\times 5\times (6\times 78+9).$
\item [] $9549=1+2^3+4\times 5\times (6\times 78+9).$
\item [] $9550=1+((2+3)^4+5+6)\times (7+8)+9.$
\item [] $\mathit{9551=(1-23)\times 4+567\times (8+9).}$
\item [] $9552=12+(3+4+5)\times (6+789).$
\item [] $9553=((1+2+3)^4+56)\times 7+89.$
\item [] $9554=1+((2\times 3)^4+56)\times 7+89.$
\item [] $9555=12+3+4\times 5\times (6\times 78+9).$
\item [] $9556=(1^2+3)\times (4+5\times (6\times 78+9)).$
\item [] $9557=(1+23+4)\times (5+6\times 7\times 8)+9.$
\item [] $9558=1^2\times 3^4\times (5+(6+7)\times 8+9).$
\item [] $9559=(1^{23}+4\times 5\times 6)\times (7+8\times 9).$
\item [] $9560=(1+(2\times 3)^4+56)\times 7+89.$
\item [] $9561=(1+2^3)\times 4^5+6\times 7\times 8+9.$
\item [] $\mathit{9562=1+2+3\times 456\times 7-8-9.}$
\item [] $9563=1\times 23+4\times 5\times (6\times 78+9).$
\item [] $9564=1\times 2\times (3+45+6\times 789).$
\item [] $9565=1+2\times (3+45+6\times 789).$
\item [] $9566=1+(2+3)\times ((4+5\times 6)\times 7\times 8+9).$
\item [] $9567=(1+2)^3+4\times 5\times (6\times 78+9).$
\item [] $9568=1\times 2^3\times 4\times (5\times 6\times 7+89).$
\item [] $9569=1+2^3\times 4\times (5\times 6\times 7+89).$
\item [] $9570=(12+3)\times (4+5+6+7\times 89).$
\item [] $9571=1\times 2\times (3+4)\times (5+678)+9.$
\item [] $9572=1+2\times (3+4)\times (5+678)+9.$
\item [] $\mathit{9573=-1\times 2+3\times 456\times 7+8-9.}$
\item [] $\mathit{9574=-1-2+3\times 456\times 7-8+9.}$
\item [] $\mathit{9575=-1\times 2+3\times 456\times 7-8+9.}$
\item [] $9576=1\times 2\times 3^4\times 56+7\times 8\times 9.$
\item [] $9577=1+2\times 3^4\times 56+7\times 8\times 9.$
\item [] $\mathit{9578=1+2+3\times 456\times 7+8-9.}$
\item [] $\mathit{9579=1\times 2+3\times 456\times 7-8+9.}$
\item [] $9580=1\times (2+3)\times 4\times (5+6\times (7+8\times 9)).$
\item[]$\mbox{Decreasing order}$
\item [] $9511=(9\times 87\times 6+54+3)\times 2+1.$
\item [] $9512=(9\times 8+7)\times 6\times 5\times 4+32\times 1.$
\item [] $9513=9\times 876+543\times (2+1).$
\item [] $9514=(9+8\times 7+6)\times (5+4\times 32+1).$
\item [] $9515=9\times 8+7\times (65+4\times 321).$
\item [] $9516=9\times 87\times (6+5)+43\times 21.$
\item [] $9517=(9\times 87\times 6+5\times 4\times 3)\times 2+1.$
\item [] $\mathit{9518=9-8+7-6\times 5\times (4-321).}$
\item [] $\mathit{9519=-98-7+6^5+43^2-1.}$
\item [] $9520=(9+8)\times 7\times (65+4\times 3+2+1).$
\item [] $9521=9+8+(7+65)\times 4\times (32+1).$
\item [] $9522=987\times 6+(5\times 4\times 3)^2\times 1.$
\item [] $9523=987\times 6+(5\times 4\times 3)^2+1.$
\item [] $\mathit{9524=9\times (8\times 76-5)+4^(3\times 2)+1.}$
\item [] $9525=9+876+5\times (4\times 3)^{(2+1)}.$
\item [] $\mathit{9526=(9\times 87+6+5)\times 4\times 3-2\times 1.}$
\item [] $\mathit{9527=(9+8+76\times 5)\times 4\times 3\times 2-1.}$
\item [] $9528=(9+8+76\times 5)\times 4\times 3\times 2\times 1.$
\item [] $9529=9+8\times 7\times (6+54\times 3+2\times 1).$
\item [] $9530=(9\times 87+6+5)\times 4\times 3+2\times 1.$
\item [] $9531=9\times (876+54\times 3+21).$
\item [] $9532=9\times 87+6\times (5+4)^3\times 2+1.$
\item [] $9533=9\times (87+6\times 54\times 3)+2\times 1.$
\item [] $9534=(9\times 87\times 6+5+4^3)\times 2\times 1.$
\item [] $9535=(9\times 87\times 6+5+4^3)\times 2+1.$
\item [] $9536=((9+8)\times 7+6\times 5)\times (43+21).$
\item [] $9537=9+8\times (7\times (6+54\times 3+2)+1).$
\item [] $9538=(9\times 8+7\times (6+5))\times 4^3+2\times 1.$
\item [] $9539=(9\times 8+7\times (6+5))\times 4^3+2+1.$
\item [] $9540=9\times ((8+7)\times 65+4^3+21).$
\item [] $9541=98+7\times (65+4\times 321).$
\item [] $\mathit{9542=(9\times (87-6)+5)\times (4+3^2\times 1).}$
\item [] $9543=(9\times 8+7)\times 6\times 5\times 4+3\times 21.$
\item [] $9544=(9+8)\times 7+65\times ((4\times 3)^2+1).$
\item [] $\mathit{9545=-9\times 8-7+6^5+43^2-1.}$
\item [] $9546=(98+7+6)\times (54+32\times 1).$
\item [] $9547=(98+7+6)\times (54+32)+1.$
\item [] $9548=98\times 7\times 6+5432\times 1.$
\item [] $9549=98\times 7\times 6+5432+1.$
\item [] $9550=9\times (87+6\times 54\times 3+2)+1.$
\item [] $9551=9\times 8+(7+6)\times (5+4)^3+2\times 1.$
\item [] $9552=(9\times 8+76\times 5\times 4)\times 3\times 2\times 1.$
\item [] $9553=(9\times 8+76\times 5\times 4)\times 3\times 2+1.$
\item [] $9554=(9+8)\times (76+54\times 3^2\times 1).$
\item [] $9555=(98+7)\times (6\times (5+4+3\times 2)+1).$
\item [] $\mathit{9556=9876+5-4-321.}$
\item [] $9557=9\times 8+7\times (6+5+4^3\times 21).$
\item [] $9558=9\times (87+654+321).$
\item [] $9559=(9\times 8+7)\times (6+5)\times (4+3\times 2+1).$
\item [] $9560=(9+8+7\times 6)\times 54\times 3+2\times 1.$
\item [] $9561=(9+8+7\times 6)\times 54\times 3+2+1.$
\item [] $9562=9\times 8+(7+6)\times (((5+4)\times 3)^2+1).$
\item [] $\mathit{9563=9-8+7+65\times (4+3)\times 21.}$
\item [] $9564=9\times 8\times 7+6^5+4\times 321.$
\item [] $9565=((9+87\times 6)\times (5+4)+3)\times 2+1.$
\item [] $9566=(9+8\times 7\times (6\times 5+4))\times (3+2)+1.$
\item [] $9567=987+65\times 4\times (32+1).$
\item [] $9568=((9+8+7\times 6)\times 5+4)\times 32\times 1.$
\item [] $9569=((9+8+7\times 6)\times 5+4)\times 32+1.$
\item [] $9570=(98+76)\times 5\times (4+3\times 2+1).$
\item [] $9571=(9+8)\times (76+54\times 3^2+1).$
\item [] $9572=9\times 8+76\times 5\times (4\times 3\times 2+1).$
\item [] $9573=9\times 8\times 76+5+4^{(3\times 2)}\times 1.$
\item [] $9574=9\times 8\times 76+5+4^{(3\times 2)}+1.$
\item [] $9575=98+(7+6)\times ((5+4)\times 3)^2\times 1.$
\item [] $9576=9\times 8\times (7+6+5\times 4\times 3\times 2\times 1).$
\item [] $9577=98+(7+6)\times (5+4)^3+2\times 1.$
\item [] $9578=9\times (87+65)\times (4+3)+2\times 1.$
\item [] $9579=9+8+7+65\times (4+3)\times 21.$
\item [] $9580=9+87\times (65+43+2)+1.$
\item[]$\mbox{Increasing order}$
\item [] $9581=(1+(2+3)\times 45)\times 6\times 7+89.$
\item [] $9582=(1+2)\times (34\times 5+6\times 7\times 8\times 9).$
\item [] $9583=(1+2+34)\times (5\times (6\times 7+8)+9).$
\item [] $\mathit{9584=(-1+2+3\times 456)\times 7-8+9.}$
\item [] $9585=1^2\times 3\times 45\times (6+7\times 8+9).$
\item [] $9586=1+(23+4)\times 5\times (6+7\times 8+9).$
\item [] $9587=1\times 2+3\times 45\times (6+7\times 8+9).$
\item [] $9588=1\times 2\times (3\times 4\times 5+6\times 789).$
\item [] $9589=1+2\times (3\times 4\times 5+6\times 789).$
\item [] $9590=1\times 2+34\times (5\times 6\times 7+8\times 9).$
\item [] $9591=1+2+34\times (5\times 6\times 7+8\times 9).$
\item [] $9592=1+23\times ((4\times (5+6)+7)\times 8+9).$
\item [] $9593=1^2\times 3\times 456\times 7+8+9.$
\item [] $9594=12\times 3^4\times 5+6\times 789.$
\item [] $9595=1\times 2+3\times 456\times 7+8+9.$
\item [] $9596=1+2+3\times 456\times 7+8+9.$
\item [] $9597=12+3\times 45\times (6+7\times 8+9).$
\item [] $9598=(123+4\times 5)\times 67+8+9.$
\item [] $9599=12^3+(456+7)\times (8+9).$
\item [] $9600=1^2\times 3\times 4\times (5+6+789).$
\item [] $9601=12^3\times 4+5\times 67\times 8+9.$
\item [] $9602=1\times 2+3\times 4\times (5+6+789).$
\item [] $9603=(123\times 4+567+8)\times 9.$
\item [] $9604=1^2+(3+4\times 5\times 6)\times 78+9.$
\item [] $9605=12+3\times 456\times 7+8+9.$
\item [] $9606=1+2+(3+4\times 5\times 6)\times 78+9.$
\item [] $9607=1\times (2+3\times 456)\times 7+8+9.$
\item [] $9608=1+(2+3\times 456)\times 7+8+9.$
\item [] $9609=1\times 2\times 3\times 4\times (56\times 7+8)+9.$
\item [] $9610=1+2\times 3\times 4\times (56\times 7+8)+9.$
\item [] $9611=(1+2\times 3)\times (4\times (5+6\times 7\times 8)+9).$
\item [] $9612=12+3\times 4\times (5+6+789).$
\item [] $9613=1+2\times (3+4+5+6\times 7)\times 89.$
\item [] $9614=(1+2+3\times 456)\times 7+8+9.$
\item [] $9615=12+(3+4\times 5\times 6)\times 78+9.$
\item [] $9616=1\times 2^3\times (45+(6+7)\times 89).$
\item [] $9617=1+2^3\times (45+(6+7)\times 89).$
\item [] $9618=1\times (2+3+4+5)\times (678+9).$
\item [] $9619=1+(2+3+4+5)\times (678+9).$
\item [] $9620=1\times (23\times 4+56)\times (7\times 8+9).$
\item [] $9621=1+(23\times 4+56)\times (7\times 8+9).$
\item [] $9622=(1\times 2+3\times 4\times (5+6\times 7))\times (8+9).$
\item [] $9623=1+(2+3\times 4\times (5+6\times 7))\times (8+9).$
\item [] $9624=12+3\times 4\times (5+6+78)\times 9.$
\item [] $9625=1\times 2\times (34+567)\times 8+9.$
\item [] $9626=1+2\times (34+567)\times 8+9.$
\item [] $9627=1^2\times 3\times (456\times 7+8+9).$
\item [] $9628=1^2+3\times (456\times 7+8+9).$
\item [] $9629=1\times 2+3\times (456\times 7+8+9).$
\item [] $9630=1+2+3\times (456\times 7+8+9).$
\item [] $9631=1+(2+3\times 4\times (5+6+78))\times 9.$
\item [] $9632=(1+2\times 3^4)\times 56+7\times 8\times 9.$
\item [] $9633=(1+2)\times (3\times 4^5+67+8\times 9).$
\item [] $9634=1\times 2\times ((34+567)\times 8+9).$
\item [] $9635=1+2\times ((34+567)\times 8+9).$
\item [] $9636=(1+2)\times (3+456\times 7+8+9).$
\item [] $9637=12^3\times 4+5\times (67\times 8+9).$
\item [] $9638=(1\times 23\times 4+5\times 6)\times (7+8\times 9).$
\item [] $9639=12+3\times (456\times 7+8+9).$
\item [] $9640=1^{234}+567\times (8+9).$
\item [] $9641=1+2\times (3^4+5+6\times 789).$
\item [] $9642=1+2+3^4\times (5+6\times 7+8\times 9).$
\item [] $9643=1^{23}\times 4+567\times (8+9).$
\item [] $9644=1^{23}+4+567\times (8+9).$
\item [] $9645=(12+3)\times (4+567+8\times 9).$
\item [] $9646=1^2\times 3+4+567\times (8+9).$
\item [] $9647=1^2+3+4+567\times (8+9).$
\item [] $9648=1^2\times 3\times 456\times 7+8\times 9.$
\item [] $9649=1^2+3\times 456\times 7+8\times 9.$
\item [] $9650=1\times 2+3\times 456\times 7+8\times 9.$
\item[]$\mbox{Decreasing order}$
\item [] $9581=((9\times 8+7)\times 6+5)\times 4\times (3+2)+1.$
\item [] $9582=(98+(7+65)\times 43)\times (2+1).$
\item [] $9583=98+7\times (6+5+4^3\times 21).$
\item [] $9584=(9+8)\times 7\times 65+43^2\times 1.$
\item [] $9585=(9+8)\times 7\times 65+43^2+1.$
\item [] $9586=9+(8+76)\times (54+3)\times 2+1.$
\item [] $9587=9\times (87\times 6+543)+2\times 1.$
\item [] $9588=9\times (87\times 6+543)+2+1.$
\item [] $\mathit{9589=9+8\times 7\times 6+5\times 43^2-1.}$
\item [] $9590=9+8\times 7\times 6+5\times 43^2\times 1.$
\item [] $9591=9+8\times 7\times 6+5\times 43^2+1.$
\item [] $\mathit{9592=(-9+876+5)\times (4+3\times 2+1).}$
\item [] $9593=9+8\times (765+432+1).$
\item [] $9594=9\times (8+7\times (65+43\times 2)+1).$
\item [] $9595=9+8\times 7\times 6+5\times (43^2+1).$
\item [] $9596=98+(7+6)\times (5+4)^3+21.$
\item [] $9597=9\times 8\times (76+54+3)+21.$
\item [] $9598=98+76\times 5\times (4\times 3\times 2+1).$
\item [] $9599=9\times 876+5\times (4+3)^{(2+1)}.$
\item [] $9600=9+87+6^5+(4\times 3)^{(2+1)}.$
\item [] $9601=9+8+7\times (6\times 5+4+3)^2+1.$
\item [] $9602=(9+8+7+6)\times 5\times 4^3+2\times 1.$
\item [] $9603=9\times (87\times 6+543+2\times 1).$
\item [] $9604=98\times (7+6\times 5+4\times 3)\times 2\times 1.$
\item [] $9605=98\times (7+6\times 5+4\times 3)\times 2+1.$
\item [] $9606=(9+8+7\times 65\times (4+3))\times (2+1).$
\item [] $9607=9+8+7\times ((6\times 5+4+3)^2+1).$
\item [] $9608=98\times 76+5\times 432\times 1.$
\item [] $9609=98\times 76+5\times 432+1.$
\item [] $9610=9+(8\times (7+6\times 5)+4)\times 32+1.$
\item [] $9611=9\times 87\times 6+(5+4\times 3)^(2+1).$
\item [] $9612=9\times 876+54\times 32\times 1.$
\item [] $9613=9\times 876+54\times 32+1.$
\item [] $9614=9+(8+7\times 6+5+43)^2+1.$
\item [] $9615=9\times (8+76+5)\times 4\times 3+2+1.$
\item [] $\mathit{9616=-9+8-7+6^5+43^2-1.}$
\item [] $9617=9\times 876+5+(4\times 3)^{(2+1)}.$
\item [] $9618=(9+8+7\times (6+54+3))\times 21.$
\item [] $9619=9\times 87+(6\times 5+4^3)^2\times 1.$
\item [] $9620=9+8\times 7+65\times (4+3)\times 21.$
\item [] $9621=(9+8+7+6)\times 5\times 4^3+21.$
\item [] $9622=(9+8)\times ((7\times 6+5)\times 4\times 3+2\times 1).$
\item [] $9623=(9+8)\times 7+6^5+(4\times 3)^{(2+1)}.$
\item [] $9624=9\times 8\times 7+6^5+4^3\times 21.$
\item [] $9625=(9+8\times (7\times 6+5))\times (4\times 3\times 2+1).$
\item [] $\mathit{9626=9876-5\times ((4+3)^2+1).}$
\item [] $9627=9\times 8+7\times 65\times (4+3)\times (2+1).$
\item [] $\mathit{9628=(9-8\times 7\times 6-5)\times (4-32-1).}$
\item [] $\mathit{9629=9876+5-4\times 3\times 21.}$
\item [] $9630=9\times (8+7+6\times 5+4^(3+2)+1).$
\item [] $9631=9\times ((8+76+5)\times 4\times 3+2)+1.$
\item [] $9632=(9+(8+(7+6)\times 5)\times 4)\times 32\times 1.$
\item [] $9633=9\times (8+76+5)\times 4\times 3+21.$
\item [] $9634=9\times 8+7+65\times (4+3)\times 21.$
\item [] $9635=9+8+7\times (6\times 5+4^3\times 21).$
\item [] $9636=(9+8\times 7)\times 6+5\times 43^2+1.$
\item [] $9637=9\times 87\times (6+5)+4^(3+2)\times 1.$
\item [] $9638=9\times 87\times (6+5)+4^(3+2)+1.$
\item [] $9639=(9\times 8+76\times 5+4+3)\times 21.$
\item [] $9640=(9+8\times 7)\times 6+5\times (43^2+1).$
\item [] $9641=9\times 8\times 7\times (6+5)+4^{(3\times 2)}+1.$
\item [] $9642=(9+8)\times (76+5)\times (4+3)+2+1.$
\item [] $\mathit{9643=-9-8+7\times 6\times 5\times (43+2+1).}$
\item [] $\mathit{9644=((9\times 8+7)\times (65-4)+3)\times 2\times 1.}$
\item [] $9645=9+(8+(7+6)\times 5)\times 4\times (32+1).$
\item [] $\mathit{9646=98-7+65\times (4+3)\times 21.}$
\item [] $\mathit{9647=-(9-87+6+5)\times (4\times 3)^2-1.}$
\item [] $9648=(9\times 8+7+65)\times (4+3\times 21).$
\item [] $9649=9+8+7+6^5+43^2\times 1.$
\item [] $9650=9+8+7+6^5+43^2+1.$
\item[]$\mbox{Increasing order}$
\item [] $9651=12\times 3^4+(5+6)\times 789.$
\item [] $9652=1+2^3+4+567\times (8+9).$
\item [] $9653=1\times 2+3\times 4+567\times (8+9).$
\item [] $9654=1+2+3\times 4+567\times (8+9).$
\item [] $9655=(1^2+3\times 456)\times 7+8\times 9.$
\item [] $9656=(1+2)^3\times (45+6)\times 7+8+9.$
\item [] $9657=12\times 3\times 4\times (5+6+7\times 8)+9.$
\item [] $9658=12+3+4+567\times (8+9).$
\item [] $9659=1\times 2+(3^4+5\times 6)\times (78+9).$
\item [] $9660=12+3\times 456\times 7+8\times 9.$
\item [] $9661=1+(2^3+4\times 5)\times (6\times 7\times 8+9).$
\item [] $9662=1\times (2+3\times 456)\times 7+8\times 9.$
\item [] $9663=12+3\times 4+567\times (8+9).$
\item [] $9664=1+2\times 3\times 4+567\times (8+9).$
\item [] $9665=1^2\times 3\times 456\times 7+89.$
\item [] $9666=1^2+3\times 456\times 7+89.$
\item [] $9667=1\times 2+3\times 456\times 7+89.$
\item [] $9668=1+2+3\times 456\times 7+89.$
\item [] $9669=1\times 23\times 4\times 5\times (6+7+8)+9.$
\item [] $9670=1+23\times 4\times 5\times (6+7+8)+9.$
\item [] $9671=1\times 2^3\times 4+567\times (8+9).$
\item [] $9672=1+2^3\times 4+567\times (8+9).$
\item [] $9673=1^2\times 34+567\times (8+9).$
\item [] $9674=1^2+34+567\times (8+9).$
\item [] $9675=1\times 2+34+567\times (8+9).$
\item [] $9676=1+2+34+567\times (8+9).$
\item [] $9677=12+3\times 456\times 7+89.$
\item [] $9678=1+23\times (4+56)\times 7+8+9.$
\item [] $9679=12\times 3+4+567\times (8+9).$
\item [] $9680=1+(2+3\times 456)\times 7+89.$
\item [] $9681=(1\times 2\times 34+56)\times 78+9.$
\item [] $9682=1+(2\times 34+56)\times 78+9.$
\item [] $9683=1+2\times (3\times 4+5+67\times 8\times 9).$
\item [] $9684=12\times 3^4\times 5+67\times 8\times 9.$
\item [] $9685=12+34+567\times (8+9).$
\item [] $9686=(1+2+3\times 456)\times 7+89.$
\item [] $9687=1+2+(34\times 5\times 6+7\times 8)\times 9.$
\item [] $9688=1\times 2\times (3+4)\times (5+678+9).$
\item [] $9689=1+2\times (3+4)\times (5+678+9).$
\item [] $9690=(1+234+5\times 67)\times (8+9).$
\item [] $9691=1+(2+3\times 4+5)\times (6+7\times 8\times 9).$
\item [] $\mathit{9692=-12-3+(4+567)\times (8+9).}$
\item [] $9693=(1+2)\times 3\times 4^5+6\times 78+9.$
\item [] $9694=1\times 2\times (3+4\times 5+67\times 8\times 9).$
\item [] $9695=1\times 2\times 3^4\times 56+7\times 89.$
\item [] $9696=1+2\times 3^4\times 56+7\times 89.$
\item [] $9697=(1\times 2\times 3^4+5+6)\times 7\times 8+9.$
\item [] $9698=1+(2\times 3^4+5+6)\times 7\times 8+9.$
\item [] $9699=(12+3)\times 4+567\times (8+9).$
\item [] $9700=1+(2^{(3+4)}\times 5+6)\times (7+8)+9.$
\item [] $9701=(12+34+56+7)\times 89.$
\item [] $9702=12^3\times 4+5\times (6+7\times 8)\times 9.$
\item [] $9703=1+2\times (3+4\times (56+78))\times 9.$
\item [] $9704=1\times 2^3\times (4^5+(6+7+8)\times 9).$
\item [] $9705=(1\times 234\times 5+6\times 7)\times 8+9.$
\item [] $9706=(1+2^3)^4+56\times 7\times 8+9.$
\item [] $9707=1\times 2\times 34+567\times (8+9).$
\item [] $9708=1+2\times 34+567\times (8+9).$
\item [] $\mathit{9709=1-2+3+(4+567)\times (8+9).}$
\item [] $9710=1^2\times 3+(4+567)\times (8+9).$
\item [] $9711=1^2+3+(4+567)\times (8+9).$
\item [] $9712=1+(2+3^4)\times (5\times 6+78+9).$
\item [] $9713=1\times 2\times 3+(4+567)\times (8+9).$
\item [] $9714=1+2\times 3+(4+567)\times (8+9).$
\item [] $9715=1\times 2^3+(4+567)\times (8+9).$
\item [] $9716=1+2^3+(4+567)\times (8+9).$
\item [] $9717=123\times (4\times 5+6\times 7+8+9).$
\item [] $9718=1^2+(3+4\times 5\times 6)\times (7+8\times 9).$
\item [] $9719=1\times 2+(3+4\times 5\times 6)\times (7+8\times 9).$
\item [] $9720=12\times 3\times (4+5+6+7+8)\times 9.$
\item[]$\mbox{Decreasing order}$
\item [] $9651=9+87+65\times (4+3)\times 21.$
\item [] $\mathit{9652=9876-5\times (43+2)+1.}$
\item [] $9653=98+7\times 65\times (4+3)\times (2+1).$
\item [] $9654=(9\times 8\times (7+6+54)+3)\times 2\times 1.$
\item [] $9655=9\times 8+7\times (6\times 5+4+3)^2\times 1.$
\item [] $9656=9\times 8+7\times (6\times 5+4+3)^2+1.$
\item [] $9657=(98+7+6)\times (54+32+1).$
\item [] $9658=9+(8\times 7+6+5)\times (4\times 3)^2+1.$
\item [] $\mathit{9659=9876-5\times 43-2\times 1.}$
\item [] $9660=98+7+65\times (4+3)\times 21.$
\item [] $9661=(98+7)\times (6+54+32)+1.$
\item [] $9662=(98+7\times 6)\times (5+4^3)+2\times 1.$
\item [] $9663=(98+7\times 6)\times (5+4^3)+2+1.$
\item [] $\mathit{9664=9876-5\times 43+2+1.}$
\item [] $9665=9+8\times (7+6\times 5\times 4\times (3^2+1)).$
\item [] $9666=987+6^5+43\times 21.$
\item [] $9667=(9\times (8+76)+5^4)\times (3\times 2+1).$
\item [] $\mathit{9668=-9-8-(7-6^5/4)\times (3+2)\times 1.}$
\item [] $9669=9+(8+76)\times ((54+3)\times 2+1).$
\item [] $\mathit{9670=9876-5\times (43-2)-1.}$
\item [] $\mathit{9671=987\times 6+5^4\times 3\times 2-1.}$
\item [] $9672=987\times 6+5^4\times 3\times 2\times 1.$
\item [] $9673=987\times 6+5^4\times 3\times 2+1.$
\item [] $9674=(9+8)\times 7+65\times (4+3)\times 21.$
\item [] $9675=(9\times (8+7\times 6)+5^4)\times 3^2\times 1.$
\item [] $9676=(9+8+7\times 6)\times (54\times 3+2\times 1).$
\item [] $9677=9+8765+43\times 21.$
\item [] $9678=9\times (8+7)\times 65+43\times 21.$
\item [] $\mathit{9679=9876-5-4^3\times (2+1).}$
\item [] $\mathit{9680=(9+8-7)\times (65+43\times 21).}$
\item [] $9681=987+6\times (5+4^3)\times 21.$
\item [] $9682=98+7\times (6\times 5+4+3)^2+1.$
\item [] $\mathit{9683=-98\times 7+6\times 54\times 32+1.}$
\item [] $9684=9+(8+7)\times (6\times 54+321).$
\item [] $\mathit{9685=-9-8\times 7+6\times 5\times (4+321).}$
\item [] $\mathit{9686=9-8-(7-6^5/4)\times (3+2)\times 1.}$
\item [] $\mathit{9687=(9\times 8\times 76-5^4-3)\times 2-1.}$
\item [] $9688=(98\times 7+6)\times (5+4+3+2\times 1).$
\item [] $9689=(98\times 7+6)\times (5+4+3+2)+1.$
\item [] $9690=9+8\times 7+6^5+43^2\times 1.$
\item [] $9691=9+8\times 7+6^5+43^2+1.$
\item [] $9692=(98\times 7+65\times 4^3)\times 2\times 1.$
\item [] $9693=(98\times 7+65\times 4^3)\times 2+1.$
\item [] $\mathit{9694=9\times (8+7\times 6)+5\times 43^2-1.}$
\item [] $9695=9\times (8+7\times 6)+5\times 43^2\times 1.$
\item [] $9696=(9+87)\times (65+4+32\times 1).$
\item [] $9697=(9+87)\times (65+4+32)+1.$
\item [] $\mathit{9698=9\times 87\times 6+5^4\times (3^2-1).}$
\item [] $9699=9+(8+7)\times (6+5\times 4\times 32\times 1).$
\item [] $9700=9+(876+5)\times (4+3\times 2+1).$
\item [] $\mathit{9701=98\times (7+6+5\times 4)\times 3-2+1.}$
\item [] $9702=98\times (7+6+54+32\times 1).$
\item [] $9703=9\times (87\times 6+5+4\times 3)\times 2+1.$
\item [] $9704=9\times 8+7+6^5+43^2\times 1.$
\item [] $9705=9\times 8+7+6^5+43^2+1.$
\item [] $\mathit{9706=-98-76\times (-5-4\times (32-1)).}$
\item [] $9707=(9+8)\times (7\times 6+(5\times 4+3)^2\times 1).$
\item [] $9708=(987+6+5^4)\times 3\times 2\times 1.$
\item [] $9709=(987+6+5^4)\times 3\times 2+1.$
\item [] $\mathit{9710=(9\times 8+7-6)\times (5+4\times 32)+1.}$
\item [] $9711=9\times (8+7\times 6+5+4^(3+2)\times 1).$
\item [] $9712=9\times (8+7\times 6+5+4^(3+2))+1.$
\item [] $\mathit{9713=9876-54\times 3-2+1.}$
\item [] $9714=9+(8+7)\times (6+5\times 4\times 32+1).$
\item [] $9715=(98+7\times 6+5)\times (4+3\times 21).$
\item [] $9716=98+7\times (6\times 5+4^3\times 21).$
\item [] $9717=(98+76)\times 54+321.$
\item [] $9718=(9+8\times (7+6))\times (54+32\times 1).$
\item [] $9719=9+8+7\times (6+5\times 4\times 3)\times 21.$
\item [] $9720=9\times (87+6\times 54\times 3+21).$
\item[]$\mbox{Increasing order}$
\item [] $9721=(123+4^5+67)\times 8+9.$
\item [] $9722=12+3+(4+567)\times (8+9).$
\item [] $9723=1+2+3^4+567\times (8+9).$
\item [] $9724=1\times 2\times 34\times (56+78+9).$
\item [] $9725=1+2\times 34\times (56+78+9).$
\item [] $9726=1\times 2\times (34+5+67\times 8\times 9).$
\item [] $9727=1+2\times (34+5+67\times 8\times 9).$
\item [] $9728=(1+2)^3\times (45+6)\times 7+89.$
\item [] $9729=(12\times 3+4^5+6+7+8)\times 9.$
\item [] $9730=1\times 23+(4+567)\times (8+9).$
\item [] $9731=1\times 23\times 4+567\times (8+9).$
\item [] $9732=1+23\times 4+567\times (8+9).$
\item [] $9733=1+23\times (4+56)\times 7+8\times 9.$
\item [] $9734=(1+2)^3+(4+567)\times (8+9).$
\item [] $9735=(1+23)\times 4+567\times (8+9).$
\item [] $\mathit{9736=1+(2\times 3\times 4+5)\times 6\times 7\times 8-9.}$
\item [] $9737=12\times (3+4+5)\times 67+89.$
\item [] $9738=1\times 2\times (3\times 45+6\times 789).$
\item [] $9739=1+2\times (3\times 45+6\times 789).$
\item [] $\mathit{9740=-1+(2+3+4\times 5\times 6)\times 78-9.}$
\item [] $9741=(1+2)\times 34+567\times (8+9).$
\item [] $\mathit{9742=1+(2+3+4\times 5\times 6)\times 78-9.}$
\item [] $9743=12\times 3+(4+567)\times (8+9).$
\item [] $9744=12\times 3\times 45\times 6+7+8+9.$
\item [] $9745=1+2\times (3+45+67\times 8\times 9).$
\item [] $\mathit{9746=(1-2+3)\times (4+(5+67\times 8)\times 9).}$
\item [] $9747=123\times 45+6\times 78\times 9.$
\item [] $9748=1^2+(3^4\times 5+678)\times 9.$
\item [] $9749=1\times 23\times (4+56)\times 7+89.$
\item [] $9750=1+23\times (4+56)\times 7+89.$
\item [] $9751=(1+2\times 3^4)\times 56+7\times 89.$
\item [] $9752=1\times 23\times (4+5\times (67+8+9)).$
\item [] $9753=(123+45+6)\times 7\times 8+9.$
\item [] $9754=1+(2+3\times (4+5))\times 6\times 7\times 8+9.$
\item [] $\mathit{9755=-1-2+34\times (5+6\times (7\times 8-9)).}$
\item [] $9756=(1^2+3^4\times 5+678)\times 9.$
\item [] $\mathit{9757=1-2+34\times (5+6\times (7\times 8-9)).}$
\item [] $9758=(1^2\times 3+4+567)\times (8+9).$
\item [] $9759=(1+2\times 34+56)\times 78+9.$
\item [] $9760=1+(2+3+4\times 5\times 6)\times 78+9.$
\item [] $9761=(1+2^3)\times 4^5+67\times 8+9.$
\item [] $9762=1+23\times (4+5\times (6+78))+9.$
\item [] $9763=1+2\times (3\times 4+(5+67\times 8)\times 9).$
\item [] $\mathit{9764=-1+2^{(3\times 4)}+5678-9.}$
\item [] $9765=12\times (345+6\times 78)+9.$
\item [] $9766=123+4+567\times (8+9).$
\item [] $9767=1\times 2^{(3+4)}+567\times (8+9).$
\item [] $9768=1\times 2\times (3\times 4\times 5+67\times 8\times 9).$
\item [] $9769=1+2\times 3\times 4\times (5\times 67+8\times 9).$
\item [] $9770=1\times (2+3)\times (4+5\times 6\times (7\times 8+9)).$
\item [] $9771=1+(2+3)\times (4+5\times 6\times (7\times 8+9)).$
\item [] $\mathit{9772=1^2-3+(4^5+6+7\times 8)\times 9.}$
\item [] $\mathit{9773=-1+2\times 3^4\times 56+78\times 9.}$
\item [] $9774=1\times 2\times 3^4\times 56+78\times 9.$
\item [] $9775=1+2\times 3^4\times 56+78\times 9.$
\item [] $\mathit{9776=1-2+3+(4^5+6+7\times 8)\times 9.}$
\item [] $9777=1^2\times 3+(4^5+6+7\times 8)\times 9.$
\item [] $9778=1^2+3+(4^5+6+7\times 8)\times 9.$
\item [] $9779=1\times 2+3+(4^5+6+7\times 8)\times 9.$
\item [] $9780=1\times 2\times 3+(4^5+6+7\times 8)\times 9.$
\item [] $9781=1+2\times 3+(4^5+6+7\times 8)\times 9.$
\item [] $9782=1\times 2^3+(4^5+6+7\times 8)\times 9.$
\item [] $9783=12\times 3\times 4+567\times (8+9).$
\item [] $9784=1+2^{(3\times 4)}+5678+9.$
\item [] $9785=12\times 3\times 45\times 6+7\times 8+9.$
\item [] $\mathit{9786=1\times 23\times 456-78\times 9.}$
\item [] $\mathit{9787=1+23\times 456-78\times 9.}$
\item [] $\mathit{9788=(123+4)\times (5-6+78)+9.}$
\item [] $9789=123+(4^5+6\times 7+8)\times 9.$
\item [] $9790=(1\times 23+45+6\times 7)\times 89.$
\item[]$\mbox{Decreasing order}$
\item [] $9721=9+87+6^5+43^2\times 1.$
\item [] $9722=9+87+6^5+43^2+1.$
\item [] $9723=98\times (76+5\times 4+3)+21.$
\item [] $9724=(9\times 8+7)\times 6+5\times (43^2+1).$
\item [] $\mathit{9725=((987-6)\times 5-43)\times 2+1.}$
\item [] $\mathit{9726=9876-5-(4\times 3)^2-1.}$
\item [] $9727=(9+8\times 7+6)\times (5+4\times (32+1)).$
\item [] $9728=(9+8\times 7+6+5)\times 4\times 32\times 1.$
\item [] $9729=98\times (76+5\times 4)+321.$
\item [] $9730=98+7+6^5+43^2\times 1.$
\item [] $9731=98+7+6^5+43^2+1.$
\item [] $9732=9\times 8+7\times 6\times 5\times (43+2+1).$
\item [] $\mathit{9733=9876-(5+4+3)^2+1.}$
\item [] $\mathit{9734=9876+5-(4+3)\times 21.}$
\item [] $9735=987+6\times (5+4)^3\times 2\times 1.$
\item [] $9736=987+6\times (5+4)^3\times 2+1.$
\item [] $9737=9+(87+65)\times (43+21).$
\item [] $9738=9\times 876+5+43^2\times 1.$
\item [] $9739=9\times 876+5+43^2+1.$
\item [] $9740=9+(87+65)\times 4^3+2+1.$
\item [] $9741=98\times 7\times 6+5^4\times 3^2\times 1.$
\item [] $9742=98\times 7\times 6+5^4\times 3^2+1.$
\item [] $\mathit{9743=9876-5-4\times 32\times 1.}$
\item [] $9744=9\times 8\times 7\times 6+5\times 4^3\times 21.$
\item [] $9745=(9+8)\times 7+6^5+43^2+1.$
\item [] $9746=98\times 7+6^5+4\times 321.$
\item [] $9747=9\times ((8+7)\times 65+4\times 3^{(2+1)}).$
\item [] $\mathit{9748=9\times 8\times 7-6+5\times (43^2+1).}$
\item [] $\mathit{9749=9876+5-4\times (32+1).}$
\item [] $9750=(9+8\times 7)\times (65+4^3+21).$
\item [] $9751=(9+8\times 7)\times 6\times (5\times 4+3+2)+1.$
\item [] $9752=98\times 76+(5+43)^2\times 1.$
\item [] $9753=98\times 76+(5+43)^2+1.$
\item [] $9754=9+8\times 7\times 6\times (5+4\times 3\times 2)+1.$
\item [] $9755=9\times 8\times 7+6+5\times 43^2\times 1.$
\item [] $9756=9\times 8\times 7+6+5\times 43^2+1.$
\item [] $9757=9\times (87+65\times (4+3))\times 2+1.$
\item [] $9758=9+(87+65)\times 4^3+21.$
\item [] $9759=(9\times 87+6\times 5)\times 4\times 3+2+1.$
\item [] $9760=9\times (8+7)+6^5+43^2\times 1.$
\item [] $9761=9\times 876+5^4\times 3+2\times 1.$
\item [] $9762=9\times 876+5^4\times 3+2+1.$
\item [] $9763=(98\times 7+65)\times (4+3^2\times 1).$
\item [] $9764=(98\times 7+65)\times (4+3^2)+1.$
\item [] $9765=(98+7)\times (6+54+32+1).$
\item [] $9766=9\times (87\times (6+5)+4\times 32)+1.$
\item [] $9767=9+8+(7+6)\times ((5+4)^3+21).$
\item [] $9768=9\times 876+(5^4+3)\times (2+1).$
\item [] $\mathit{9769=-9+(8+7)\times 654-32\times 1.}$
\item [] $9770=((9+8\times 7)\times 6\times 5+4)\times (3+2)\times 1.$
\item [] $9771=((9+8\times 7)\times 6\times 5+4)\times (3+2)+1.$
\item [] $\mathit{9772=(9+8+7-6)\times 543-2\times 1.}$
\item [] $\mathit{9773=9876+5-4\times 3^{(2+1)}.}$
\item [] $9774=9\times (87\times 6+543+21).$
\item [] $9775=(9\times 8+7+6)\times ((54+3)\times 2+1).$
\item [] $9776=9+87\times 6+5\times 43^2\times 1.$
\item [] $9777=9+87\times 6+5\times 43^2+1.$
\item [] $9778=9\times (876+5)+43^2\times 1.$
\item [] $9779=9\times (876+5)+43^2+1.$
\item [] $9780=9\times 876+5^4\times 3+21.$
\item [] $9781=9+87\times 6+5\times (43^2+1).$
\item [] $9782=(9+8+7\times 6\times 5)\times 43+21.$
\item [] $9783=(9\times (8+7)\times 6+5)\times 4\times 3+2+1.$
\item [] $\mathit{9784=9876-5-43\times 2-1.}$
\item [] $\mathit{9785=9876-5-43\times 2\times 1.}$
\item [] $9786=(9+(8+7)\times 6\times 5+4+3)\times 21.$
\item [] $9787=98\times (76+5)+43^2\times 1.$
\item [] $9788=98\times (76+5)+43^2+1.$
\item [] $9789=((9\times (8+7)\times 6+5)\times 4+3)\times (2+1).$
\item [] $9790=(9+876+5)\times (4+3\times 2+1).$
\item[]$\mbox{Increasing order}$
\item [] $9791=1+(2+3)\times (4+5+6+7)\times 89.$
\item [] $9792=12^3\times 4+5\times 6\times (7+89).$
\item [] $9793=1^2+3\times (456\times 7+8\times 9).$
\item [] $9794=1\times 2+3\times (456\times 7+8\times 9).$
\item [] $9795=1+2+3\times (456\times 7+8\times 9).$
\item [] $9796=(123+4)\times (5+6)\times 7+8+9.$
\item [] $9797=1+(2\times 34+56)\times (7+8\times 9).$
\item [] $9798=1+23+(4^5+6+7\times 8)\times 9.$
\item [] $9799=12\times 3\times 45\times 6+7+8\times 9.$
\item [] $9800=(1+23+4)\times (5+6\times 7\times 8+9).$
\item [] $9801=1\times 2\times 3^4+567\times (8+9).$
\item [] $9802=1+2\times 3^4+567\times (8+9).$
\item [] $9803=1\times 2+3^4\times (56+7\times 8+9).$
\item [] $9804=12+3\times (456\times 7+8\times 9).$
\item [] $9805=1^2+(3^4+5)\times (6\times 7+8\times 9).$
\item [] $9806=1\times 2+(3^4+5)\times (6\times 7+8\times 9).$
\item [] $9807=12\times 3\times 45\times 6+78+9.$
\item [] $9808=1\times 2\times (34\times 5+6\times 789).$
\item [] $9809=1+2\times (34\times 5+6\times 789).$
\item [] $9810=1\times (2+3)\times (45\times 6\times 7+8\times 9).$
\item [] $9811=1+(2\times 3+4)\times (5+(6+7)\times 8)\times 9.$
\item [] $9812=1\times 2\times (34+56\times (78+9)).$
\item [] $9813=12+3^4\times (56+7\times 8+9).$
\item [] $\mathit{9814=1\times 2\times ((3+4)\times (5-6+78\times 9)).}$
\item [] $9815=(1+(2+3+4\times 5)\times 6)\times (7\times 8+9).$
\item [] $9816=12\times 3\times 45\times 6+7+89.$
\item [] $9817=1+2\times 3\times 4\times (56\times 7+8+9).$
\item [] $\mathit{9818=-1+(2\times 3+4\times 5\times 6)\times 78-9.}$
\item [] $9819=(12\times 34+5+678)\times 9.$
\item [] $9820=1\times 2\times (3^4+5+67\times 8\times 9).$
\item [] $9821=12\times (3\times 45\times 6+7)+8+9.$
\item [] $\mathit{9822=-123+45\times (6+7)\times (8+9).}$
\item [] $\mathit{9823=(-1-2+(34\times 5+6)\times 7)\times 8-9.}$
\item [] $9824=1\times 2^3\times (4+(5+67)\times (8+9)).$
\item [] $9825=(1+2+3\times 4)\times 5\times (6\times 7+89).$
\item [] $9826=(1+2\times 3+4+567)\times (8+9).$
\item [] $\mathit{9827=1-2+3\times 4\times (5\times 6+789).}$
\item [] $9828=12^3+(4+56)\times (7+8)\times 9.$
\item [] $9829=1+(2^3+4)\times (5\times 6+789).$
\item [] $9830=123+(4+567)\times (8+9).$
\item [] $9831=1+2+3\times 4\times (5\times 6+789).$
\item [] $9832=1+(2+3^4+5\times 6)\times (78+9).$
\item [] $9833=((1+2)^3\times 45+6+7)\times 8+9.$
\item [] $9834=(1\times (2+3)^4+5\times 6)\times (7+8)+9.$
\item [] $9835=1+((2+3)^4+5\times 6)\times (7+8)+9.$
\item [] $9836=(1+2)\times 3\times (4^5+67)+8+9.$
\item [] $9837=(1\times 2\times 3+4\times 5\times 6)\times 78+9.$
\item [] $9838=1+2\times (3+4+56)\times 78+9.$
\item [] $9839=12\times 3\times 45\times 6+7\times (8+9).$
\item [] $9840=123\times (4+5+6+7\times 8+9).$
\item [] $9841=1+2\times 3\times (4\times 56\times 7+8\times 9).$
\item [] $\mathit{9842=1-2+3\times (456\times 7+89).}$
\item [] $9843=1^2\times 3\times (456\times 7+89).$
\item [] $9844=1^2+3\times (456\times 7+89).$
\item [] $9845=1+23\times 4\times (5+6+7+89).$
\item [] $9846=1+2+3\times (456\times 7+89).$
\item [] $9847=1+2\times ((3+4+56)\times 78+9).$
\item [] $\mathit{9848=(12+3\times 45)\times 67+8-9.}$
\item [] $9849=(1+(2+3)^4+5\times 6)\times (7+8)+9.$
\item [] $\mathit{9850=(12+3\times 45)\times 67-8+9.}$
\item [] $9851=(123+4)\times (5+6)\times 7+8\times 9.$
\item [] $9852=(1+2)\times (3+456\times 7+89).$
\item [] $9853=(1\times 2\times 3^4+5)\times (6\times 7+8+9).$
\item [] $9854=1\times (2+3)\times 4^5+6\times 789.$
\item [] $9855=12+3\times (456\times 7+89).$
\item [] $9856=(1+23+4)\times (5\times 67+8+9).$
\item [] $9857=(12\times 3\times (4+5\times 6)+7)\times 8+9.$
\item [] $9858=1+2+3\times ((45+6)\times 7+8)\times 9.$
\item [] $\mathit{9859=-1\times 2+3\times 4^5+6789.}$
\item [] $9860=(1+2^3+4+567)\times (8+9).$
\item[]$\mbox{Decreasing order}$
\item [] $\mathit{9791=9876-54-32+1.}$
\item [] $9792=9\times 8\times (76+54+3+2+1).$
\item [] $9793=(9+87+6)\times (5+43)\times 2+1.$
\item [] $9794=(9\times 8+76+5)\times 4^3+2\times 1.$
\item [] $9795=(9\times 8+76+5)\times 4^3+2+1.$
\item [] $\mathit{9796=9876+5-43\times 2+1.}$
\item [] $\mathit{9797=-9+(8+7)\times 654-3-2+1.}$
\item [] $9798=9\times 8\times 76+5+4321.$
\item [] $9799=(9\times 8\times 7\times 6+5^4\times 3)\times 2+1.$
\item [] $9800=98\times (7+65+4+3+21).$
\item [] $9801=(9\times 8+7)\times 6\times 5\times 4+321.$
\item [] $9802=9\times (87+6\times 5+4)\times 3^2+1.$
\item [] $9803=9\times (8\times 7+6)+5\times 43^2\times 1.$
\item [] $9804=9\times (8\times 7+6)+5\times 43^2+1.$
\item [] $9805=(9+(8+7+6)\times 5)\times 43\times 2+1.$
\item [] $9806=98\times 7+6^5+4^3\times 21.$
\item [] $9807=((98+7+6+5)\times 4+3)\times 21.$
\item [] $9808=9\times (8\times 7+6)+5\times (43^2+1).$
\item [] $9809=(9+8)\times (7+6+543+21).$
\item [] $9810=9\times (87\times 6+5\times 4+3)\times 2\times 1.$
\item [] $9811=9\times (87\times 6+5\times 4+3)\times 2+1.$
\item [] $9812=9+8+(7+6\times 543)\times (2+1).$
\item [] $9813=(9\times 8+76+5)\times 4^3+21.$
\item [] $\mathit{9814=9876+5-4-3\times 21.}$
\item [] $9815=9+8\times 7+6\times 5\times (4+321).$
\item [] $9816=(9+8\times 7)\times (65+43\times 2)+1.$
\item [] $\mathit{9817=9+(8+7)\times 654-3+2-1.}$
\item [] $9818=9\times 87\times 6+5\times 4^(3+2)\times 1.$
\item [] $9819=9+(8+7\times 6\times 5)\times (43+2)\times 1.$
\item [] $9820=9+(8+7)\times (6\times 54+3)\times 2+1.$
\item [] $9821=(9+87)\times 6+5\times 43^2\times 1.$
\item [] $9822=(9+87)\times 6+5\times 43^2+1.$
\item [] $9823=987+(6\times 5+4^3)^2\times 1.$
\item [] $9824=9+(8+7)\times 654+3+2\times 1.$
\item [] $9825=9+(8+7)\times 654+3\times 2\times 1.$
\item [] $9826=9+(8+7)\times 654+3\times 2+1.$
\item [] $9827=(9\times 87\times 6+5\times 43)\times 2+1.$
\item [] $9828=(9+87+6+54)\times 3\times 21.$
\item [] $9829=9+(8+7)\times 654+3^2+1.$
\item [] $9830=(9\times (8+7\times 6\times 5)+4)\times (3+2)\times 1.$
\item [] $9831=(9+8\times (7+6))\times (54+32+1).$
\item [] $\mathit{9832=9876-5\times 4-3-21.}$
\item [] $\mathit{9833=9876+5\times 4-3\times 21.}$
\item [] $9834=9+(8+7)\times (6+5^4+3+21).$
\item [] $9835=((9+8+7)\times 6+5)\times (4^3+2)+1.$
\item [] $\mathit{9836=9876+5-43-2\times 1.}$
\item [] $9837=9+(87+65+4)\times 3\times 21.$
\item [] $9838=9\times (8\times 76+5)+4321.$
\item [] $9839=9+8\times 7+6\times 543\times (2+1).$
\item [] $9840=(9+(8+7)\times 65)\times (4+3+2+1).$
\item [] $9841=(9+(8+7)\times 65)\times (4+3\times 2)+1.$
\item [] $\mathit{9842=9876-5+4-32-1.}$
\item [] $9843=9\times 87+6^5+4\times 321.$
\item [] $\mathit{9844=9876+5-4-32-1.}$
\item [] $\mathit{9845=9876+5-4-32\times 1.}$
\item [] $9846=(987+654)\times 3\times 2\times 1.$
\item [] $9847=(987+654)\times 3\times 2+1.$
\item [] $9848=98+(7+6)\times ((5+4)^3+21).$
\item [] $9849=(98+7\times 65\times (4+3))\times (2+1).$
\item [] $9850=(9\times (8\times (7+6)+5)+4)\times (3^2+1).$
\item [] $9851=9+(8+7)\times 654+32\times 1.$
\item [] $9852=9+(8+7)\times 654+32+1.$
\item [] $9853=9\times 8+7+6\times 543\times (2+1).$
\item [] $\mathit{9854=9876+5+4-32+1.}$
\item [] $9855=9\times (8\times 76+54\times 3^2+1).$
\item [] $9856=9\times (8+7+6\times 5\times (4+32))+1.$
\item [] $9857=9\times (8+7)\times (6\times 5+43)+2\times 1.$
\item [] $9858=9+(8\times 7+6+5)\times (4+3)\times 21.$
\item [] $9859=9+(8+7\times 6)\times (5+4^3\times (2+1)).$
\item [] $9860=(98+7+6+5)\times (4^3+21).$
\item[]$\mbox{Increasing order}$
\item [] $9861=1^2\times 3\times 4^5+6789.$
\item [] $9862=1+2\times 3^4\times 56+789.$
\item [] $9863=1\times 2+3\times 4^5+6789.$
\item [] $9864=1+2+3\times 4^5+6789.$
\item [] $9865=1^2\times (34\times 5+6)\times 7\times 8+9.$
\item [] $9866=(12+3\times 45)\times 67+8+9.$
\item [] $9867=1\times 2+(34\times 5+6)\times 7\times 8+9.$
\item [] $9868=1+2+(34\times 5+6)\times 7\times 8+9.$
\item [] $9869=(1+2\times (3+4)\times 5)\times (67+8\times 9).$
\item [] $9870=(1+2)\times (3+4^5)+6789.$
\item [] $9871=1+2\times (3+4)\times 5\times (6+(7+8)\times 9).$
\item [] $9872=(1+2\times 3\times 4)\times 56\times 7+8\times 9.$
\item [] $9873=12+3\times 4^5+6789.$
\item [] $9874=1+234+567\times (8+9).$
\item [] $9875=(1+2\times 34+56)\times (7+8\times 9).$
\item [] $9876=12\times (3\times 45\times 6+7)+8\times 9.$
\item [] $9877=12+(34\times 5+6)\times 7\times 8+9.$
\item [] $9878=1+(23+4+56)\times 7\times (8+9).$
\item [] $9879=(1+23+45+6\times 7)\times 89.$
\item [] $9880=((1+23)\times 4+56)\times (7\times 8+9).$
\item [] $9881=1\times 2^{(3\times 4)}+5\times (6+7)\times 89.$
\item [] $9882=1\times 2\times 3\times 4^5+6\times 7\times 89.$
\item [] $9883=1+2\times 3\times 4^5+6\times 7\times 89.$
\item [] $9884=(1+2\times 3)\times (4\times 5\times 67+8\times 9).$
\item [] $9885=(12+3)\times (4+5\times (6\times 7+89)).$
\item [] $\mathit{9886=1-2\times 3+(4^5+67+8)\times 9.}$
\item [] $\mathit{9887=(1\times 234\times 5+67)\times 8-9.}$
\item [] $9888=(1\times 23\times 4+5+6)\times (7+89).$
\item [] $9889=(1+2\times 3\times 4)\times 56\times 7+89.$
\item [] $\mathit{9890=1\times 2-3+(4^5+67+8)\times 9.}$
\item [] $9891=(1^{23}\times 4^5+67+8)\times 9.$
\item [] $9892=1^{23}+(4^5+67+8)\times 9.$
\item [] $9893=12\times (3\times 45\times 6+7)+89.$
\item [] $9894=1\times (23\times 4+5)\times (6+7+89).$
\item [] $9895=1^2+3+(4^5+67+8)\times 9.$
\item [] $9896=1\times 2+3+(4^5+67+8)\times 9.$
\item [] $9897=123+(4^5+6+7\times 8)\times 9.$
\item [] $9898=1+2\times 3+(4^5+67+8)\times 9.$
\item [] $9899=1\times 2^3+(4^5+67+8)\times 9.$
\item [] $9900=(12+3)\times (4+567+89).$
\item [] $9901=1+(2+3)\times 4\times (5+6\times 7+8)\times 9.$
\item [] $\mathit{9902=1\times 2-(3\times (4-56)\times 7-8) \times  9.}$
\item [] $9903=(1+2^3)\times 4^5+678+9.$
\item [] $\mathit{9904=(-1+(2+3)^4-5)\times (6-7+8+9).}$
\item [] $9905=(1\times 234\times 5+67)\times 8+9.$
\item [] $9906=12+3+(4^5+67+8)\times 9.$
\item [] $9907=1+2\times (3^4+56\times (78+9)).$
\item [] $9908=1234\times 5+6\times 7\times 89.$
\item [] $9909=12\times (3\times 45\times 6+7+8)+9.$
\item [] $\mathit{9910=-1+2^3\times 4\times 5\times (6+7\times 8)-9.}$
\item [] $9911=((1^2+3)\times 4+567)\times (8+9).$
\item [] $9912=12\times (3+4+5\times 6+789).$
\item [] $9913=1234+(5+6)\times 789.$
\item [] $9914=(1^{23}+4)^5+6789.$
\item [] $9915=1+(2+3)^4\times 5+6789.$
\item [] $\mathit{9916=12^3+4\times (5\times 6-7)\times 89.}$
\item [] $9917=(1+2\times 3^4)\times 56+789.$
\item [] $9918=1\times 2\times (3\times 45+67\times 8\times 9).$
\item [] $9919=(1+(2+3)^4)\times 5+6789.$
\item [] $9920=1\times 2+(3+4^5+67+8)\times 9.$
\item [] $9921=(12+3\times 45)\times 67+8\times 9.$
\item [] $9922=(1^2+3^4)\times (56+7\times 8+9).$
\item [] $\mathit{9923=1-23+45\times (6+7)\times (8+9).}$
\item [] $9924=(1+2^3)\times 4^5+6+78\times 9.$
\item [] $\mathit{9925=-1+2\times (3+4)\times (-5+6\times 7\times (8+9)).}$
\item [] $\mathit{9926=1\times 2\times (3+4)\times (-5+6\times 7\times (8+9)).}$
\item [] $9927=12\times 3+(4^5+67+8)\times 9.$
\item [] $9928=1\times 2^{(3+4)}\times (5+6)\times 7+8\times 9.$
\item [] $9929=1+2\times 34\times (5+6+(7+8)\times 9).$
\item [] $9930=1+2^3\times 4\times 5\times (6+7\times 8)+9.$
\item[]$\mbox{Decreasing order}$
\item [] $\mathit{9861=9876-5-4-3-2-1.}$
\item [] $9862=9+8\times 76+5\times 43^2\times 1.$
\item [] $9863=9+8\times 76+5\times 43^2+1.$
\item [] $9864=9\times 8\times (76+54+3\times 2+1).$
\item [] $9865=9\times 8\times 7\times (6+5)+4321.$
\item [] $9866=9+(8+7)\times (654+3)+2\times 1.$
\item [] $9867=9+8\times 76+5\times (43^2+1).$
\item [] $\mathit{9868=9876-5-4+3-2\times 1.}$
\item [] $9869=(9+8)\times 7+6\times 5\times (4+321).$
\item [] $9870=9+87+6\times 543\times (2+1).$
\item [] $9871=(98+7)\times ((6+5)\times 4+3)\times 2+1.$
\item [] $9872=(9+8\times 76)\times (5+4+3\times 2+1).$
\item [] $9873=9+(87+6\times 54)\times (3+21).$
\item [] $9874=9+8+7\times (6+5)\times 4\times 32+1.$
\item [] $9875=(98+7)\times 6+5\times 43^2\times 1.$
\item [] $9876=(98+7)\times 6+5\times 43^2+1.$
\item [] $9877=(9+8)\times ((7+6)\times 5\times 4+321).$
\item [] $\mathit{9878=9876+5\times 4+3-21.}$
\item [] $9879=98+7+6\times 543\times (2+1).$
\item [] $9880=9+8+7\times (65+4^3\times 21).$
\item [] $9881=9+8\times (7\times 6\times 5+4^(3+2)\times 1).$
\item [] $9882=9+(8+7)\times 654+3\times 21.$
\item [] $9883=9\times (87\times 6+(5+4)\times 3)\times 2+1.$
\item [] $\mathit{9884=9876+5-4+3\times 2+1.}$
\item [] $9885=9+(8+7)\times (654+3)+21.$
\item [] $9886=(9\times 87+65\times 4^3)\times 2\times 1.$
\item [] $9887=(9\times 87+65\times 4^3)\times 2+1.$
\item [] $9888=(9+8\times (7+6\times 5)+4)\times 32\times 1.$
\item [] $9889=(9+8\times (7+6\times 5)+4)\times 32+1.$
\item [] $9890=9876+5+4+3+2\times 1.$
\item [] $9891=9876+5+4+3+2+1.$
\item [] $9892=9876+5+4+3\times 2+1.$
\item [] $9893=(9+8)\times 7+6\times 543\times (2+1).$
\item [] $9894=9876+5+4+3^2\times 1.$
\item [] $9895=9876+5+4+3^2+1.$
\item [] $9896=9876+5+4\times 3+2+1.$
\item [] $9897=9\times 87\times (6+5)+4\times 321.$
\item [] $9898=98\times (7+6\times 5+43+21).$
\item [] $9899=98\times ((7+6)\times 5+4+32)+1.$
\item [] $9900=(9+8+7+6)\times (5+4+321).$
\item [] $9901=9876+5\times 4+3+2\times 1.$
\item [] $9902=9876+5\times 4+3+2+1.$
\item [] $9903=9876+5\times 4+3\times 2+1.$
\item [] $\mathit{9904=9876-5+4\times 3+21.}$
\item [] $9905=9876+5+4\times 3\times 2\times 1.$
\item [] $9906=9876+5+4\times 3\times 2+1.$
\item [] $\mathit{9907=9876-5+4+32\times 1.}$
\item [] $\mathit{9908=9876+5-4+32-1.}$
\item [] $9909=9876+5+4+3+21.$
\item [] $9910=9876+(5+4\times 3)\times 2\times 1.$
\item [] $9911=9876+(5+4\times 3)\times 2+1.$
\item [] $9912=9876+(5+4+3)\times (2+1).$
\item [] $9913=9876+5\times (4+3)+2\times 1.$
\item [] $9914=9876+5+4\times 3+21.$
\item [] $\mathit{9915=9876-5+43+2-1.}$
\item [] $\mathit{9916=9876+5+4+32-1.}$
\item [] $9917=9876+5+4+32\times 1.$
\item [] $9918=9876+5+4+32+1.$
\item [] $9919=98\times (7+6\times 5+4^3)+21.$
\item [] $9920=9876+5\times 4+3+21.$
\item [] $9921=9876+5\times (4+3+2)\times 1.$
\item [] $9922=9876+5\times (4+3+2)+1.$
\item [] $9923=9876+(5\times 4+3)\times 2+1.$
\item [] $9924=9876+(5+4)\times 3+21.$
\item [] $9925=(9+8+76\times 5)\times (4\times 3\times 2+1).$
\item [] $9926=9876+5+43+2\times 1.$
\item [] $9927=9876+5+43+2+1.$
\item [] $9928=9876+5\times 4+32\times 1.$
\item [] $9929=9876+5\times 4+32+1.$
\item [] $9930=9876+5+(4+3)^2\times 1.$
\item[]$\mbox{Increasing order}$
\item [] $\mathit{9931=12^3\times 4-5+6\times 7\times 8\times 9.}$
\item [] $\mathit{9932=-12\times 3+(45+67)\times 89.}$
\item [] $\mathit{9933=-12+3\times (45+6)\times (7\times 8+9).}$
\item [] $\mathit{9934=-1\times 2+3\times (45-6+7)\times 8\times 9.}$
\item [] $9935=(1+2^{(3+4)}\times (5+6))\times 7+8\times 9.$
\item [] $9936=12\times (3\times 4+5+67+8)\times 9.$
\item [] $9937=1+2^3\times (4+56+78)\times 9.$
\item [] $9938=(12+3\times 45)\times 67+89.$
\item [] $\mathit{9939=-1\times 2\times 3+45\times (6+7)\times (8+9).}$
\item [] $9940=1+2\times 3\times 4\times 5\times (6+7\times 8+9).$
\item [] $9941=12^3\times 4+5+6\times 7\times 8\times 9.$
\item [] $9942=1\times 2\times 3\times (4\times 56\times 7+89).$
\item [] $9943=1+2\times 3\times (4\times 56\times 7+89).$
\item [] $9944=1\times (2+3)\times 4^5+67\times 8\times 9.$
\item [] $9945=1+(2+3)\times 4^5+67\times 8\times 9.$
\item [] $9946=1^{23}+45\times (6+7)\times (8+9).$
\item [] $9947=1\times 2+3\times (45+6)\times (7\times 8+9).$
\item [] $9948=1^2\times 3+45\times (6+7)\times (8+9).$
\item [] $9949=1^2+3+45\times (6+7)\times (8+9).$
\item [] $9950=1\times 2+3+45\times (6+7)\times (8+9).$
\item [] $9951=1+2+3+45\times (6+7)\times (8+9).$
\item [] $9952=1+2\times 3+45\times (6+7)\times (8+9).$
\item [] $9953=1\times 2^3+45\times (6+7)\times (8+9).$
\item [] $9954=(12\times 3^4+56+78)\times 9.$
\item [] $9955=1+2\times (3\times 4+5+67\times 8)\times 9.$
\item [] $\mathit{9956=1-2-3-4\times 5\times (6-7\times 8\times 9).}$
\item [] $9957=12+3\times (45+6)\times (7\times 8+9).$
\item [] $\mathit{9958=-1+(2+34\times 5+6)\times 7\times 8-9.}$
\item [] $9959=1\times 23\times (4+5\times (6+78)+9).$
\item [] $9960=12+3+45\times (6+7)\times (8+9).$
\item [] $9961=1+2^3\times (456+789).$
\item [] $9962=(123+456+7)\times (8+9).$
\item [] $9963=(1\times 2^3+4^5+67+8)\times 9.$
\item [] $9964=1+(2^3+(4^5+67+8))\times 9.$
\item [] $\mathit{9965=-1+(234+5)\times 6\times 7-8\times 9.}$
\item [] $\mathit{9966=(12\times 3\times 4+5)\times 67-8-9.}$
\item [] $\mathit{9967=-1+23+45\times (6+7)\times (8+9).}$
\item [] $9968=1^{23}\times (45+67)\times 89.$
\item [] $9969=1+23+45\times (6+7)\times (8+9).$
\item [] $\mathit{9970=1-2+3+(45+67)\times 89.}$
\item [] $9971=1^2\times 3+(45+67)\times 89.$
\item [] $9972=12\times (3\times 4+5\times 6+789).$
\item [] $9973=1\times 2+3+(45+67)\times 89.$
\item [] $9974=1+2+3+(45+67)\times 89.$
\item [] $9975=(1\times 2\times 3)^4+(5+6)\times 789.$
\item [] $9976=1\times 2^3+(45+67)\times 89.$
\item [] $9977=1+2^3+(45+67)\times 89.$
\item [] $9978=1+(2+34\times 5+6)\times 7\times 8+9.$
\item [] $9979=1+2\times 3+(4^5+6+78)\times 9.$
\item [] $9980=123\times (4+(5+6)\times 7)+8+9.$
\item [] $9981=12^3\times 4+(5+6\times 7\times 8)\times 9.$
\item [] $9982=1\times 2\times (3+4)\times (5+6+78\times 9).$
\item [] $9983=12+3+(45+67)\times 89.$
\item [] $9984=(12+3^4+5+6)\times (7+89).$
\item [] $9985=1+(2+3\times (4+5\times 6))\times (7+89).$
\item [] $\mathit{9986=(12+3)^4/5-67-8\times 9.}$
\item [] $9987=12+3+(4^5+6+78)\times 9.$
\item [] $9988=1\times 2\times (34\times 5+67\times 8\times 9).$
\item [] $9989=12\times 3\times (45\times 6+7)+8+9.$
\item [] $9990=1\times (23+45+6)\times (7+8)\times 9.$
\item [] $9991=1\times 23+(45+67)\times 89.$
\item [] $9992=1+23+(45+67)\times 89.$
\item [] $9993=1\times 2\times (34+5\times 6)\times 78+9.$
\item [] $9994=1+2\times (34+5\times 6)\times 78+9.$
\item [] $9995=1\times 23+(4^5+6+78)\times 9.$
\item [] $9996=1+23+(4^5+6+78)\times 9.$
\item [] $9997=1+2\times (3\times 4+5\times 6)\times 7\times (8+9).$
\item [] $9998=1\times 2+(3+45\times (6+7))\times (8+9).$
\item [] $9999=(1+23\times 45+67+8)\times 9.$
\item [] $10000=(12\times 3\times 4+5)\times 67+8+9.$
\item[]$\mbox{Decreasing order}$
\item [] $9931=9876+5+(4+3)^2+1.$
\item [] $9932=9876+5\times (4+3)+21.$
\item [] $9933=(9\times (8+7\times 6)+5\times 4+3)\times 21.$
\item [] $\mathit{9934=9876+54+3+2-1.}$
\item [] $9935=9876+54+3+2\times 1.$
\item [] $9936=9876+54+3+2+1.$
\item [] $9937=9876+54+3\times 2+1.$
\item [] $9938=9876+5\times 4\times 3+2\times 1.$
\item [] $9939=9876+5\times 4\times 3+2+1.$
\item [] $9940=9876+54+3^2+1.$
\item [] $9941=9876+5\times (4+3^2\times 1).$
\item [] $9942=9876+5\times (4+3^2)+1.$
\item [] $9943=(9\times 8\times ((7+6)\times 5+4)+3)\times 2+1.$
\item [] $9944=((987+6)\times 5+4+3)\times 2\times 1.$
\item [] $9945=9876+5+43+21.$
\item [] $9946=9876+5\times (4+3)\times 2\times 1.$
\item [] $9947=9876+5+4^3+2\times 1.$
\item [] $9948=9876+5+4+3\times 21.$
\item [] $\mathit{9949=(9+8-7+6)\times (5^4-3)-2-1.}$
\item [] $9950=9+8+7\times (6+5)\times 43\times (2+1).$
\item [] $9951=9876+5\times (4\times 3+2+1).$
\item [] $\mathit{9952=(9+8\times 76+5)\times 4\times (3+2-1).}$
\item [] $\mathit{9953=98+7\times (6+5)\times 4\times 32-1.}$
\item [] $9954=9876+54+3+21.$
\item [] $9955=98+7\times (6+5)\times 4\times 32+1.$
\item [] $\mathit{9956=9876-5+43\times 2-1.}$
\item [] $9957=9876+5\times 4\times 3+21.$
\item [] $9958=9876+(5+4)\times 3^2+1.$
\item [] $9959=9876+5\times 4+3\times 21.$
\item [] $9960=(9+8)\times 7\times 6+5\times 43^2+1.$
\item [] $9961=98+7\times (65+4^3\times 21).$
\item [] $9962=9876+54+32\times 1.$
\item [] $9963=9876+54+32+1.$
\item [] $9964=(987+6\times 5\times 4)\times 3^2+1.$
\item [] $\mathit{9965=(9\times 8\times 7-6)\times 5\times 4+3+2\times 1.}$
\item [] $9966=9876+5+4^3+21.$
\item [] $9967=9876+5+43\times 2\times 1.$
\item [] $9968=9876+5+43\times 2+1.$
\item [] $9969=(9+8\times 7+6\times 543)\times (2+1).$
\item [] $\mathit{9970=(9\times 8\times 7-6)\times 5\times 4+3^2+1.}$
\item [] $\mathit{9971=9876+(5+43)\times 2-1.}$
\item [] $9972=9876+(5+43)\times 2\times 1.$
\item [] $9973=9876+(5+43)\times 2+1.$
\item [] $\mathit{9974=9\times (87-6)+5\times 43^2\times 1.}$
\item [] $9975=((9+8\times (7+6)+5)\times 4+3)\times 21.$
\item [] $9976=9876+5\times 4\times (3+2\times 1).$
\item [] $9977=9876+5+4\times (3+21).$
\item [] $\mathit{9978=9876+(54-3)\times 2\times 1.}$
\item [] $9979=(9+8)\times (7\times 6+543+2\times 1).$
\item [] $9980=(987+6+5)\times (4+3+2+1).$
\item [] $9981=9876+5\times (4+3)\times (2+1).$
\item [] $9982=((9+8)\times (7+6)\times 5+4)\times 3^2+1.$
\item [] $9983=(9\times 8+7\times (6+5))\times (4+3\times 21).$
\item [] $9984=(98+7\times 6\times 5+4)\times 32\times 1.$
\item [] $9985=9+87\times 65+4321.$
\item [] $\mathit{9986=9876+5\times (43-21).}$
\item [] $9987=9\times 8\times 76+5\times 43\times 21.$
\item [] $\mathit{9988=(-9+8\times (7-6)\times 5^4+3)\times 2\times 1.}$
\item [] $9989=9876+5+4\times 3^{(2+1)}.$
\item [] $9990=9876+(54+3)\times 2\times 1.$
\item [] $9991=9876+(54+3)\times 2+1.$
\item [] $\mathit{9992=(9+8-7+6)\times 5^4-3^2+1.}$
\item [] $9993=9876+54+3\times 21.$
\item [] $9994=9+8\times (7+6)\times (5+43)\times 2+1.$
\item [] $\mathit{9995=9876+5\times 4\times 3\times 2-1.}$
\item [] $9996=9876+5\times 4\times 3\times 2\times 1.$
\item [] $9997=9876+5\times 4\times 3\times 2+1.$
\item [] $9998=9+8\times 7+(6+5)\times 43\times 21.$
\item [] $9999=98\times (7\times 6+5\times 4\times 3)+2+1.$
\item [] $10000=(9\times 8+7\times (6+5)\times 4^3)\times 2\times 1.$
\item[]$\mbox{Increasing order}$
\item [] $10001=(1\times 23\times (4+5)\times 6+7)\times 8+9.$
\item [] $10002=1\times 2\times ((34+5\times 6)\times 78+9).$
\item [] $10003=(1+2\times 3)\times (4\times 5\times 67+89).$
\item [] $10004=12\times 3+(45+67)\times 89.$
\item [] $10005=(12+3\times 4+5)\times (6\times 7\times 8+9).$
\item [] $10006=1+23\times (45+6\times (7\times 8+9)).$
\item [] $\mathit{10007=(1234+5+6+7)\times 8-9.}$
\item [] $10008=12\times (34+5+6+789).$
\item [] $10009=1+2^3\times (4+5)\times (67+8\times 9).$
\item [] $10010=(1\times 2+3\times 4)\times (5+6)\times (7\times 8+9).$
\item [] $10011=(1+2)\times 3\times 4^5+6+789.$
\item [] $\mathit{10012=1+(23\times (4+5)+6)\times (7\times 8-9).}$
\item [] $10013=(12+3+4)\times (5+6\times (78+9)).$
\item [] $10014=123+(4^5+67+8)\times 9.$
\item [] $10015=1+2\times (3+(4\times 5+67\times 8)\times 9).$
\item [] $\mathit{10016=-1\times 2+(3+4)^5-6789.}$
\item [] $10017=1\times (2^3+45)\times (6+7+8)\times 9.$
\item [] $10018=1+(2^3+45)\times (6+7+8)\times 9.$
\item [] $\mathit{10019=-1+2+(3+4)^5-6789.}$
\item [] $10020=12\times (34+(5+6+78)\times 9).$
\item [] $\mathit{10021=1+2+(3+4)^5-6789.}$
\item [] $10022=(1+2^{(3+4)})\times (5+6)\times 7+89.$
\item [] $10023=(1+2)\times (3\times 4^5+6)+789.$
\item [] $\mathit{10024=1^2-3+(4^5+6\times (7+8))\times 9.}$
\item [] $10025=(1234+5+6+7)\times 8+9.$
\item [] $10026=(12+3+4^5+67+8)\times 9.$
\item [] $10027=1^{23}+(4^5+6\times (7+8))\times 9.$
\item [] $10028=1\times 23\times 4\times (5\times 6+7+8\times 9).$
\item [] $10029=1+23\times 4\times (5\times 6+7+8\times 9).$
\item [] $10030=1^2\times 34\times 5\times (6\times 7+8+9).$
\item [] $10031=1^2+34\times 5\times (6\times 7+8+9).$
\item [] $10032=1\times 2+34\times 5\times (6\times 7+8+9).$
\item [] $10033=(1+2+34\times 5+6)\times 7\times 8+9.$
\item [] $10034=1\times 2^3+(4^5+6\times (7+8))\times 9.$
\item [] $10035=(1+2\times 3+4^5+6+78)\times 9.$
\item [] $\mathit{10036=-1+23\times 4\times (5+(6+7)\times 8)+9.}$
\item [] $10037=(12+34)\times (5\times 6\times 7+8)+9.$
\item [] $10038=1+23\times 4\times (5+(6+7)\times 8)+9.$
\item [] $\mathit{10039=12^3\times 4+56\times 7\times 8-9.}$
\item [] $\mathit{10040=1+(234+5)\times 6\times 7-8+9.}$
\item [] $10041=12+3+(4^5+6\times (7+8))\times 9.$
\item [] $10042=12+34\times 5\times (6\times 7+8+9).$
\item [] $10043=(1\times 2+3^4)\times (56+7\times 8+9).$
\item [] $10044=12\times 3\times (45\times 6+7)+8\times 9.$
\item [] $10045=1+(2^3+4^5+6+78)\times 9.$
\item [] $\mathit{10046=1-2+3\times (4+5\times (678-9)).}$
\item [] $10047=12\times 34+567\times (8+9).$
\item [] $10048=1+(2\times 3\times 4+567)\times (8+9).$
\item [] $10049=1\times 23+(4^5+6\times (7+8))\times 9.$
\item [] $10050=1+23+(4^5+6\times (7+8))\times 9.$
\item [] $10051=(12+3+4)\times (5\times (6+7)\times 8+9).$
\item [] $10052=123\times (4+(5+6)\times 7)+89.$
\item [] $10053=(12+3^4)\times (5\times 6+78)+9.$
\item [] $10054=1^2+(3+4^5+6\times (7+8))\times 9.$
\item [] $10055=1\times (234+5)\times 6\times 7+8+9.$
\item [] $10056=1+(234+5)\times 6\times 7+8+9.$
\item [] $10057=12^3\times 4+56\times 7\times 8+9.$
\item [] $10058=1+(2\times (3+4\times 5)+67)\times 89.$
\item [] $10059=(1^2+3)\times 4^5+67\times 89.$
\item [] $10060=1+2^{(3+4+5)}+67\times 89.$
\item [] $10061=12\times 3\times (45\times 6+7)+89.$
\item [] $10062=1\times 2\times (3+4\times 5+67\times 8)\times 9.$
\item [] $10063=1+2\times (3+4\times 5+67\times 8)\times 9.$
\item [] $10064=(1\times 2^3)^4+5+67\times 89.$
\item [] $10065=((1+2)^3\times 45+6\times 7)\times 8+9.$
\item [] $10066=1\times 2\times (3+4)\times (5+6\times 7\times (8+9)).$
\item [] $10067=1+2\times (3+4)\times (5+6\times 7\times (8+9)).$
\item [] $10068=12\times (3^4+56+78\times 9).$
\item [] $\mathit{10069=12\times 3^4\times (5+6)-7\times 89.}$
\item [] $\mathit{10070=1-2+3\times 4\times 56\times (7+8)-9.}$
\item[]$\mbox{Decreasing order}$
\item [] $10001=9876+5\times (4\times 3\times 2+1).$
\item [] $10002=9\times (8+76)+5\times 43^2+1.$
\item [] $10003=(9+8\times (7+6)\times (5+43))\times 2+1.$
\item [] $\mathit{10004=(9+8-7+6)\times 5^4+3+2-1.}$
\item [] $10005=(9+8\times 7\times 6)\times (5+4\times 3\times 2\times 1).$
\item [] $10006=9\times (8+76)+5\times (43^2+1).$
\item [] $\mathit{10007=(9+8-7+6)\times 5^4+3\times 2+1.}$
\item [] $10008=9\times 8\times (7+65+4+3\times 21).$
\item [] $10009=9876+5+4\times 32\times 1.$
\item [] $10010=9876+5+4\times 32+1.$
\item [] $10011=(9\times 8+7+6\times 543)\times (2+1).$
\item [] $10012=9\times 8+7+(6+5)\times 43\times 21.$
\item [] $10013=9876+5+4\times (32+1).$
\item [] $10014=9+87\times (6\times 5+4^3+21).$
\item [] $10015=9876+(5+4^3)\times 2+1.$
\item [] $10016=9876+5\times (4+3+21).$
\item [] $10017=9\times 8+765\times (4+3^2\times 1).$
\item [] $10018=9\times 8+765\times (4+3^2)+1.$
\item [] $\mathit{10019=9876+(5+4+3)^2-1.}$
\item [] $10020=9876+(5+43)\times (2+1).$
\item [] $10021=9876+(5+4+3)^2+1.$
\item [] $\mathit{10022=987\times 6+5+4^(3\times 2)-1.}$
\item [] $10023=987\times 6+5+4^{(3+2+1)}.$
\item [] $10024=(98\times 7+6\times 5)\times (4+3)\times 2\times 1.$
\item [] $10025=9876+5+(4\times 3)^2\times 1.$
\item [] $10026=9876+5+(4\times 3)^2+1.$
\item [] $10027=9\times (87\times 6+5\times (4+3))\times 2+1.$
\item [] $10028=9876+5+(4+3)\times 21.$
\item [] $10029=9+87+(6+5)\times 43\times 21.$
\item [] $10030=(9+8\times (7+6)+5)\times (4^3+21).$
\item [] $10031=98+7\times (6+5)\times 43\times (2+1).$
\item [] $10032=(9+8\times 7+6+5)\times 4\times (32+1).$
\item [] $10033=(9\times 8+7)\times (6+5\times 4\times 3\times 2+1).$
\item [] $10034=9\times 87+6+5\times 43^2\times 1.$
\item [] $10035=9\times 87+6+5\times 43^2+1.$
\item [] $10036=(9\times 87\times 6+5\times 4^3)\times 2\times 1.$
\item [] $10037=(9\times 87\times 6+5\times 4^3)\times 2+1.$
\item [] $10038=98+7+(6+5)\times 43\times 21.$
\item [] $10039=9\times 87+6+5\times (43^2+1).$
\item [] $10040=9876+54\times 3+2\times 1.$
\item [] $10041=9876+54\times 3+2+1.$
\item [] $10042=9+(87+65)\times (4^3+2)+1.$
\item [] $10043=98+765\times (4+3^2)\times 1.$
\item [] $10044=9\times 876+5\times 432\times 1.$
\item [] $10045=9\times 876+5\times 432+1.$
\item [] $10046=98\times 7+65\times (4\times 3)^2\times 1.$
\item [] $10047=987+6^5+4\times 321.$
\item [] $\mathit{10048=(-9-8+7+6\times 54)\times 32\times 1.}$
\item [] $10049=9\times 876+5\times (432+1).$
\item [] $10050=(9+8+7+6)\times 5\times (4+3\times 21).$
\item [] $10051=((9+8+76)\times 54+3)\times 2+1.$
\item [] $10052=(9+8)\times 7+(6+5)\times 43\times 21.$
\item [] $10053=(98+7\times (6+5)\times 4^3)\times 2+1.$
\item [] $\mathit{10054=9\times (8+7)\times 6+5\times 43^2-1.}$
\item [] $10055=9\times (8+7)\times 6+5\times 43^2\times 1.$
\item [] $10056=9876+5\times 4\times 3^2\times 1.$
\item [] $10057=9876+5\times 4\times 3^2+1.$
\item [] $10058=9+8765+4\times 321.$
\item [] $10059=9876+54\times 3+21.$
\item [] $10060=9\times (8+7)\times 6+5\times (43^2+1).$
\item [] $10061=9876+5\times (4+32+1).$
\item [] $10062=(9+87+6\times 543)\times (2+1).$
\item [] $10063=9\times (8\times 7+6\times 5)\times (4+3^2)+1.$
\item [] $10064=9\times (8+7+6+5)\times 43+2\times 1.$
\item [] $10065=9\times (8+7+6+5)\times 43+2+1.$
\item [] $10066=((9+8)\times 7\times 6+5)\times (4+3)\times 2\times 1.$
\item [] $10067=((9+8)\times 7\times 6+5)\times (4\times 3+2)+1.$
\item [] $10068=9\times (8+7)+(6+5)\times 43\times 21.$
\item [] $10069=(9+(87+6)\times 54+3)\times 2+1.$
\item [] $\mathit{10070=9876+5\times 43-21.}$
\item[]$\mbox{Increasing order}$
\item [] $10071=(1\times 23\times 45+6+78)\times 9.$
\item [] $10072=(12\times 3\times 4+5)\times 67+89.$
\item [] $10073=1^2\times 34\times (5\times 6+7)\times 8+9.$
\item [] $10074=1^2+34\times (5\times 6+7)\times 8+9.$
\item [] $10075=1\times 2+34\times (5\times 6+7)\times 8+9.$
\item [] $10076=1+2+34\times (5\times 6+7)\times 8+9.$
\item [] $10077=12\times (3\times (45\times 6+7)+8)+9.$
\item [] $\mathit{10078=-1+(23\times 4+5)\times (6+7)\times 8-9.}$
\item [] $\mathit{10079=(1+2\times (34+56)\times 7)\times 8-9.}$
\item [] $10080=(1\times 2\times 34+5+67)\times 8\times 9.$
\item [] $10081=1+2\times 3\times 4\times 5\times (67+8+9).$
\item [] $10082=1+(2^3+45\times (6+7))\times (8+9).$
\item [] $\mathit{10083=12+3\times (4+56)\times 7\times 8-9.}$
\item [] $\mathit{10084=(12+3)^4/5+6-7\times 8+9.}$
\item [] $10085=12+34\times (5\times 6+7)\times 8+9.$
\item [] $10086=(1+2)\times (3\times (4^5+67)+89).$
\item [] $\mathit{10087=-1\times 2+3\times 4\times 56\times (7+8)+9.}$
\item [] $\mathit{10088=1-2+3\times 4\times 56\times (7+8)+9.}$
\item [] $10089=1\times 2\times (34+56)\times 7\times 8+9.$
\item [] $10090=1+2\times (34+56)\times 7\times 8+9.$
\item [] $10091=123+(45+67)\times 89.$
\item [] $10092=1+2+3\times 4\times 56\times (7+8)+9.$
\item [] $\mathit{10093=12\times 3+(4\times 5\times 6-7)\times 89.}$
\item [] $\mathit{10094=-1+2\times (3^4-5)\times 67-89.}$
\item [] $10095=123+(4^5+6+78)\times 9.$
\item [] $\mathit{10096=-1+(23\times 4+5)\times (6+7)\times 8+9.}$
\item [] $10097=(1+234+5)\times 6\times 7+8+9.$
\item [] $10098=(1\times 23+4+567)\times (8+9).$
\item [] $10099=1+2\times (3+(4+5)\times (6+7\times 8))\times 9.$
\item [] $\mathit{10100=1\times (2+3)\times 4\times (-5+6+7\times 8\times 9).}$
\item [] $10101=12+3\times 4\times 56\times (7+8)+9.$
\item [] $\mathit{10102=1+(2+56\times (3\times 4))\times (7+8)-9.}$
\item [] $10103=12\times 345+67\times 89.$
\item [] $10104=(1^2+3\times 4\times 56)\times (7+8)+9.$
\item [] $10105=(1\times 2+3\times (4+56)\times 7)\times 8+9.$
\item [] $10106=1+(2+3\times (4+56)\times 7)\times 8+9.$
\item [] $10107=(1+23+4^5+67+8)\times 9.$
\item [] $10108=1^2+3\times (4\times 56\times (7+8)+9).$
\item [] $10109=1\times 2+3\times (4\times 56\times (7+8)+9).$
\item [] $10110=1\times (234+5)\times 6\times 7+8\times 9.$
\item [] $10111=1+(234+5)\times 6\times 7+8\times 9.$
\item [] $10112=1\times 2\times (34+5\times 6)\times (7+8\times 9).$
\item [] $10113=1\times 2\times 3\times (4+5\times 6\times 7\times 8)+9.$
\item [] $10114=1+2\times 3\times (4+5\times 6\times 7\times 8)+9.$
\item [] $10115=(1+23+4+567)\times (8+9).$
\item [] $10116=12\times (3+45+6+789).$
\item [] $10117=(1+(234+5)\times 6)\times 7+8\times 9.$
\item [] $\mathit{10118=1\times 2+3\times 45\times (67+8)-9.}$
\item [] $10119=(1\times 2+3\times 4\times 56)\times (7+8)+9.$
\item [] $10120=1+(2+3\times 4\times 56)\times (7+8)+9.$
\item [] $\mathit{10121=(12+3)^4/5+6+7-8-9.}$
\item [] $10122=(1+2)\times (3+4)\times (5+6\times 78+9).$
\item [] $10123=(1+2\times 3)^4+(5+6)\times 78\times 9.$
\item [] $\mathit{10124=1-2+3\times 45\times (6+78-9).}$
\item [] $10125=(12+3+4+56)\times (7+8)\times 9.$
\item [] $10126=1+(2\times 3+4+5)\times (67+8)\times 9.$
\item [] $10127=1\times (234+5)\times 6\times 7+89.$
\item [] $10128=1+(234+5)\times 6\times 7+89.$
\item [] $\mathit{10129=-1\times 23+(4^5+(6+7)\times 8)\times 9.}$
\item [] $\mathit{10130=-1+(2+3)\times (4\times 5+6)\times 78-9.}$
\item [] $10131=123\times 4+567\times (8+9).$
\item [] $10132=1\times 2\times 34\times (5\times 6+7\times (8+9)).$
\item [] $10133=12\times (3+4\times 5\times 6\times 7)+8+9.$
\item [] $10134=1^2\times 3\times 45\times (67+8)+9.$
\item [] $10135=1^2+3\times 45\times (67+8)+9.$
\item [] $10136=1\times 2+3\times 45\times (67+8)+9.$
\item [] $10137=(1+23\times 4)\times (5\times 6+7+8\times 9).$
\item [] $10138=1^2+3\times (4+5\times (67+8)\times 9).$
\item [] $10139=1\times 2+3\times (4+5\times (67+8)\times 9).$
\item [] $10140=(1\times 2^3+4)\times (56+789).$
\item[]$\mbox{Decreasing order}$
\item [] $10071=9+(87+6\times 5)\times 43\times 2\times 1.$
\item [] $10072=9+(87+6\times 5)\times 43\times 2+1.$
\item [] $10073=9876+5+4^3\times (2+1).$
\item [] $10074=9+(8\times 7+6)\times 54\times 3+21.$
\item [] $10075=(9+8\times 7)\times (6+5+(4\times 3)^2\times 1).$
\item [] $10076=9876+5\times 4\times (3^2+1).$
\item [] $\mathit{10077=-987\times 6+(5\times 4)^3\times 2-1.}$
\item [] $\mathit{10078=-987\times 6+(5\times 4)^3\times 2\times 1.}$
\item [] $\mathit{10079=-987\times 6+(5\times 4)^3\times 2+1.}$
\item [] $10080=(9+87+6\times 54)\times (3+21).$
\item [] $10081=9\times (8+7+65)\times (4\times 3+2)+1.$
\item [] $10082=9\times (87+6)+5\times 43^2\times 1.$
\item [] $10083=9\times (87+6)+5\times 43^2+1.$
\item [] $\mathit{10084=-9+87\times (6\times 5+43\times 2)+1.}$
\item [] $\mathit{10085=(-9+8+7)\times (6+5\times (4+3))^2-1.}$
\item [] $\mathit{10086=9+8\times 7\times (6+54)\times 3-2-1.}$
\item [] $10087=9\times (87+6)+5\times (43^2+1).$
\item [] $\mathit{10088=9876+5\times 43-2-1.}$
\item [] $10089=9+(8+76)\times 5\times 4\times 3\times 2\times 1.$
\item [] $10090=9+(8+76)\times 5\times 4\times 3\times 2+1.$
\item [] $10091=9+8\times 7\times (6+54)\times 3+2\times 1.$
\item [] $10092=9+8\times 7\times (6+54)\times 3+2+1.$
\item [] $10093=9876+5\times 43+2\times 1.$
\item [] $10094=9876+5\times 43+2+1.$
\item [] $10095=98\times (7+6\times 5+4^3+2)+1.$
\item [] $10096=98\times (76+(5+4)\times 3)+2\times 1.$
\item [] $10097=9+8+7\times (6+54)\times (3+21).$
\item [] $10098=(9\times 8+76+5)\times (4^3+2\times 1).$
\item [] $10099=(9\times 8+76+5)\times (4^3+2)+1.$
\item [] $10100=(98\times 7+6\times 54)\times (3^2+1).$
\item [] $10101=9876+5\times (43+2\times 1).$
\item [] $10102=9876+5\times (43+2)+1.$
\item [] $\mathit{10103=(9\times (87+6)+5)\times 4\times 3-2+1.}$
\item [] $10104=(9+8+7)\times (6\times 5\times (4+3)\times 2+1).$
\item [] $10105=9+8\times (7\times (6+54)\times 3+2\times 1).$
\item [] $10106=9876+5\times (43+2+1).$
\item [] $10107=987+6^5+4^3\times 21.$
\item [] $\mathit{10108=-98+7\times 6\times (5+4)\times 3^{(2+1)}.}$
\item [] $10109=(9+8\times 7\times (6+54))\times 3+2\times 1.$
\item [] $10110=9+8\times 7\times (6+54)\times 3+21.$
\item [] $10111=98\times 7+65\times ((4\times 3)^2+1).$
\item [] $10112=9876+5\times 43+21.$
\item [] $10113=9+(8\times 7\times 6\times 5+4)\times 3\times 2\times 1.$
\item [] $10114=9+(8\times 7\times 6\times 5+4)\times 3\times 2+1.$
\item [] $10115=98\times (76+(5+4)\times 3)+21.$
\item [] $10116=(9+8\times 7\times (6+54)+3)\times (2+1).$
\item [] $\mathit{10117=-9+8-7+(6+5+4)^3\times (2+1).}$
\item [] $10118=9+8765+4^3\times 21.$
\item [] $10119=9\times (8+7)\times 65+4^3\times 21.$
\item [] $10120=(98\times 7+6\times (5+4)^3)\times 2\times 1.$
\item [] $10121=9876+5\times (4+3)^2\times 1.$
\item [] $10122=9876+5\times (4+3)^2+1.$
\item [] $10123=(9+(8\times 7\times 6\times 5+4)\times 3)\times 2+1.$
\item [] $\mathit{10124=-9876+5^4\times 32\times 1.}$
\item [] $10125=9+8+76\times (5+4^3\times 2\times 1).$
\item [] $10126=9\times (8\times 7+65+4)\times 3^2+1.$
\item [] $10127=9\times (8+7)\times (6+5+4^3)+2\times 1.$
\item [] $10128=9876+(5+4+3)\times 21.$
\item [] $10129=9\times 8\times 7+6^5+43^2\times 1.$
\item [] $10130=9\times 87\times 6+5432\times 1.$
\item [] $10131=9\times 87\times 6+5432+1.$
\item [] $10132=(9+8)\times (7\times (65+4\times (3+2))+1).$
\item [] $10133=9876+5+4\times 3\times 21.$
\item [] $10134=(9+(8+76)\times 5\times 4)\times 3\times 2\times 1.$
\item [] $10135=9+876+5\times (43^2+1).$
\item [] $\mathit{10136=-987+6\times (5+43^2)-1.}$
\item [] $\mathit{10137=(9\times 8\times 7+6)\times 5\times 4-3\times 21.}$
\item [] $\mathit{10138=-987+6\times (5+43^2)+1.}$
\item [] $\mathit{10139=(9+8\times 7)\times (6+5\times 4)\times 3\times 2-1.}$
\item [] $10140=9+(8+7)\times 654+321.$
\item[]$\mbox{Increasing order}$
\item [] $10141=1^2+3\times 4\times (56+789).$
\item [] $10142=1\times 2+3\times 4\times (56+789).$
\item [] $10143=1+2+3\times 4\times (56+789).$
\item [] $10144=1+23\times (4+5\times 6+7+8)\times 9.$
\item [] $10145=(1^2+3\times (4+56))\times 7\times 8+9.$
\item [] $10146=12+3\times 45\times (67+8)+9.$
\item [] $10147=(1\times 2\times 34+5)\times (67+8\times 9).$
\item [] $10148=1\times (2+34\times 5)\times (6\times 7+8+9).$
\item [] $10149=12+3\times (4+5\times (67+8)\times 9).$
\item [] $10150=1+(2+3)\times (4\times 5+6)\times 78+9.$
\item [] $\mathit{10151=(-1+23)\times 456+7\times (8+9).}$
\item [] $10152=(1+234+5)\times 6\times 7+8\times 9.$
\item [] $10153=1+(2+34)\times (5\times 6\times 7+8\times 9).$
\item [] $10154=1\times 2+3\times (45\times (67+8)+9).$
\item [] $10155=1+2+3\times (45\times (67+8)+9).$
\item [] $10156=1^2+3+(4^5+(6+7)\times 8)\times 9.$
\item [] $10157=1\times 2+3+(4^5+(6+7)\times 8)\times 9.$
\item [] $10158=1\times 2\times (345+6\times 789).$
\item [] $10159=1+2\times (345+6\times 789).$
\item [] $10160=(123+4)\times (56+7+8+9).$
\item [] $10161=(1\times 23+4)\times (5+6\times 7)\times 8+9.$
\item [] $10162=1+(23+4)\times (5+6\times 7)\times 8+9.$
\item [] $\mathit{10163=1-2+3\times (4+(5+6\times 7)\times 8\times 9).}$
\item [] $10164=12\times (3+4)\times (56+7\times 8+9).$
\item [] $10165=1^2+3\times (4+(5+6\times 7)\times 8\times 9).$
\item [] $10166=(1^{23}+45)\times (6+7)\times (8+9).$
\item [] $10167=12+3+(4^5+(6+7)\times 8)\times 9.$
\item [] $10168=1\times 2+34\times (5\times 6\times 7+89).$
\item [] $10169=(1+234+5)\times 6\times 7+89 .$
\item [] $10170=((1+23)\times 45+6\times 7+8)\times 9.$
\item [] $\mathit{10171=-1\times 2+3\times (4+5\times 678)-9.}$
\item [] $\mathit{10172=1\times 2\times (3^4\times (56+7)-8-9).}$
\item [] $10173=12\times (3+4+56\times (7+8))+9.$
\item [] $\mathit{10174=-1+(2+3\times 4+5)\times 67\times 8-9.}$
\item [] $10175=(1+2\times 3\times 4)\times (5\times 67+8\times 9).$
\item [] $10176=(12+3\times 4)\times (5\times 67+89).$
\item [] $10177=(1234+5\times 6+7)\times 8+9.$
\item [] $10178=12+34\times (5\times 6\times 7+89).$
\item [] $10179=(1+2+3+4+5)\times 678+9.$
\item [] $10180=1+(2\times 3+4+5)\times 678+9.$
\item [] $10181=1\times 2+(3+4^5+(6+7)\times 8)\times 9.$
\item [] $10182=1+2+(3+4^5+(6+7)\times 8)\times 9.$
\item [] $10183=(1\times 2^3\times 4+567)\times (8+9).$
\item [] $10184=1+(2^3\times 4+567)\times (8+9).$
\item [] $10185=12\times (34\times 5+678)+9.$
\item [] $10186=1+(2+3)\times ((4\times 5+6)\times 78+9).$
\item [] $\mathit{10187=1-2+3\times 4\times (56\times (7+8)+9).}$
\item [] $10188=(1+23+4^5+6+78)\times 9.$
\item [] $10189=1+(2^3+4)\times (56\times (7+8)+9).$
\item [] $10190=1\times 2+3\times 4\times (56\times (7+8)+9).$
\item [] $10191=1^2\times 3\times (4+5\times 678)+9.$
\item [] $10192=1^2+3\times (4+5\times 678)+9.$
\item [] $10193=1\times 2+3\times (4+5\times 678)+9.$
\item [] $10194=1+2+3\times (4+5\times 678)+9.$
\item [] $10195=1+2\times 3\times (4^5+(67+8)\times 9).$
\item [] $\mathit{10196=1-2-3+4\times 5\times (6+7\times 8\times 9).}$
\item [] $10197=(1\times 2+3+4^5+(6+7)\times 8)\times 9.$
\item [] $10198=1\times 2^{(3+4+5)}+678\times 9.$
\item [] $10199=1+2^{(3+4+5)}+678\times 9.$
\item [] $10200=12\times (3+4\times 56+7\times 89).$
\item [] $10201=1+2\times 34\times 5\times (6+7+8+9).$
\item [] $10202=1+(2+3\times (4+56))\times 7\times 8+9.$
\item [] $10203=12+3\times (4+5\times 678)+9.$
\item [] $10204=1^2+3+4\times 5\times (6+7\times 8\times 9).$
\item [] $10205=12\times (3+4\times 5\times 6\times 7)+89.$
\item [] $10206=1^2\times 3^4\times (5\times 6+7+89).$
\item [] $10207=1^2+3^4\times (5\times 6+7+89).$
\item [] $10208=1\times 2+3^4\times (5\times 6+7+89).$
\item [] $10209=1^2\times 3\times (4+5\times 678+9).$
\item [] $10210=1^2+3\times (4+5\times 678+9).$
\item[]$\mbox{Decreasing order}$
\item [] $10141=(9+8\times 7)\times 6\times (5\times 4+3\times 2)+1.$
\item [] $\mathit{10142=9\times 87+65\times (4\times 3)^2-1.}$
\item [] $10143=9\times 87+65\times (4\times 3)^2\times 1.$
\item [] $10144=9\times 87+65\times (4\times 3)^2+1.$
\item [] $10145=9+8\times (7\times 6+(5\times (4+3))^2\times 1).$
\item [] $10146=9876+54\times (3+2)\times 1.$
\item [] $10147=9876+54\times (3+2)+1.$
\item [] $10148=(9+8\times (7+6)+5)\times 43\times 2\times 1.$
\item [] $10149=(9+8\times (7+6)+5)\times 43\times 2+1.$
\item [] $10150=((9+8\times 7)\times 6+5^4)\times (3^2+1).$
\item [] $\mathit{10151=(9\times 8+7\times 6\times 5)\times 4\times 3^2-1.}$
\item [] $10152=9\times (8+7+6+543)\times 2\times 1.$
\item [] $10153=9\times 8\times (76+5)+4321.$
\item [] $10154=9\times 8\times (7\times (6+5)+4^3)+2\times 1.$
\item [] $10155=9\times 8\times (7\times (6+5)+4^3)+2+1.$
\item [] $\mathit{10156=-9+8+7\times (-6+(5+4)^3\times 2-1).}$
\item [] $\mathit{10157=(9+8\times 7\times 6\times 5+4)\times 3\times 2-1.}$
\item [] $10158=(9+8\times 7\times 6\times 5+4)\times 3\times 2\times 1.$
\item [] $10159=(9+8\times 7\times 6\times 5+4)\times 3\times 2+1.$
\item [] $\mathit{10160=(9-87\times 6+5)\times (4-3-21).}$
\item [] $10161=9+8\times (76\times 5+43)\times (2+1).$
\item [] $10162=9+8\times (7+6+(5^4+3)\times 2)+1.$
\item [] $10163=9+8+(7+(6+5+4)^3)\times (2+1).$
\item [] $10164=9876+(5+4)\times 32\times 1.$
\item [] $10165=9876+(5+4)\times 32+1.$
\item [] $10166=9876+(5+4\times 3)^2+1.$
\item [] $10167=(9+8)\times (7+6)\times (5\times 4+3)\times 2+1.$
\item [] $\mathit{10168=(9\times 8\times 7+6)\times 5\times 4-32\times 1.}$
\item [] $10169=9+8\times (7+6\times 5\times 4)\times (3^2+1).$
\item [] $10170=9\times (87+6+5\times 4)\times (3^2+1).$
\item [] $10171=(987+6\times 5)\times (4+3\times 2)+1.$
\item [] $\mathit{10172=-9+8\times 7+(6+5+4)^3\times (2+1).}$
\item [] $10173=9876+(5+4)\times (32+1).$
\item [] $\mathit{10174=9\times (8+7+6)\times 54-32\times 1.}$
\item [] $\mathit{10175=9876+5\times 4^3-21.}$
\item [] $10176=(9+87)\times (6+5\times 4\times (3+2)\times 1).$
\item [] $10177=9+(8\times 7+6)\times (54\times 3+2\times 1).$
\item [] $10178=(9\times 87+65)\times 4\times 3+2\times 1.$
\item [] $10179=(9\times 87+65)\times 4\times 3+2+1.$
\item [] $10180=9\times 8+76\times (5+4\times 32)\times 1.$
\item [] $10181=9\times 8+76\times (5+4\times 32)+1.$
\item [] $10182=9\times 8\times (7+6)+5\times 43^2+1.$
\item [] $10183=(9+8)\times (7\times 65+(4\times 3)^2\times 1).$
\item [] $10184=(9+8)\times (7\times 65+(4\times 3)^2)+1.$
\item [] $10185=(98+76\times 5+4+3)\times 21.$
\item [] $10186=(98+7)\times (6+5+43\times 2)+1.$
\item [] $\mathit{10187=9\times 876+(5+43)^2-1.}$
\item [] $10188=9\times 876+(5+43)^2\times 1.$
\item [] $10189=9\times 876+(5+43)^2+1.$
\item [] $10190=9+8\times 7+(6+5+4)^3\times (2+1).$
\item [] $10191=(9\times 8+7)\times (65+43+21).$
\item [] $10192=98\times (7+6+5+43\times 2\times 1).$
\item [] $10193=98\times (7+6+5+43\times 2)+1.$
\item [] $10194=(9\times 8+7)\times (65+4^3)+2+1.$
\item [] $\mathit{10195=9876+5\times 4^3-2+1.}$
\item [] $10196=9876+5\times (43+21).$
\item [] $10197=(9\times 87+65)\times 4\times 3+21.$
\item [] $10198=9876+5\times 4^3+2\times 1.$
\item [] $10199=9876+5\times 4^3+2+1.$
\item [] $10200=9876+54\times 3\times 2\times 1.$
\item [] $10201=9876+54\times 3\times 2+1.$
\item [] $10202=(9+8+7+65+4\times 3)^2+1.$
\item [] $\mathit{10203=9\times (8+7+6)\times 54-3\times (2-1).}$
\item [] $10204=9\times 8+7+(6+5+4)^3\times (2+1).$
\item [] $10205=(9\times 8\times 7+6)\times 5\times 4+3+2\times 1.$
\item [] $10206=9876+5+4+321.$
\item [] $10207=98+76\times (5+4\times 32)+1.$
\item [] $10208=9\times 87+65\times ((4\times 3)^2+1).$
\item [] $10209=(9\times 8\times 7+6)\times 5\times 4+3^2\times 1.$
\item [] $10210=9+(8+76+5+4\times 3)^2\times 1.$
\item[]$\mbox{Increasing order}$
\item [] $10211=1\times 2+3\times (4+5\times 678+9).$
\item [] $10212=1+2+3\times (4+5\times 678+9).$
\item [] $\mathit{10213=1-23\times (4+56\times (-7+8-9)).}$
\item [] $\mathit{10214=(1+(2-3+4)^5\times 6)\times 7-8+9.}$
\item [] $10215=(12\times 3+4^5+67+8)\times 9.$
\item [] $10216=(1^2+3)\times (4+5\times (6+7\times 8\times 9)).$
\item [] $10217=1^2\times (34+567)\times (8+9).$
\item [] $10218=12+3^4\times (5\times 6+7+89).$
\item [] $10219=1+(2\times 3)^4\times 5+6\times 7\times 89.$
\item [] $10220=1+2+(34+567)\times (8+9).$
\item [] $10221=12+3\times (4+5\times 678+9).$
\item [] $10222=1+2\times 3\times (4^5+678)+9.$
\item [] $10223=1\times 2\times 3^4\times (56+7)+8+9.$
\item [] $10224=12\times 3\times 45\times 6+7\times 8\times 9.$
\item [] $10225=(1+2\times 3\times 4)\times (56\times 7+8+9).$
\item [] $\mathit{10226=1-2\times (3-4^5\times (6+7-8))-9.}$
\item [] $10227=((1+2\times 3\times 4)\times 5+6)\times 78+9.$
\item [] $\mathit{10228=-1\times 23+(4+5)\times 67\times (8+9).}$
\item [] $10229=12+(34+567)\times (8+9).$
\item [] $10230=1\times 2\times (3\times (4^5+678)+9).$
\item [] $10231=1+2\times (3\times (4^5+678)+9).$
\item [] $10232=1+2\times (3^4\times (56+7)+8)+9.$
\item [] $10233=12\times 3\times 4\times (56+7+8)+9.$
\item [] $10234=(1^2+34+567)\times (8+9).$
\item [] $10235=(1^2\times 3+45+67)\times 89.$
\item [] $10236=12\times (34+5\times 6+789).$
\item [] $10237=1\times 2+(3+45+67)\times 89.$
\item [] $10238=1+2+(3+45+67)\times 89.$
\item [] $10239=(1+(2+3)^4+56)\times (7+8)+9.$
\item [] $10240=1\times 2^{(3+4)}\times (56+7+8+9).$
\item [] $10241=1+2^{(3+4)}\times (56+7+8+9).$
\item [] $10242=12\times 345+678\times 9.$
\item [] $10243=1\times 2^{(3\times 4)}+(5+678)\times 9.$
\item [] $10244=1+2^{(3\times 4)}+(5+678)\times 9.$
\item [] $10245=(12+3)\times (4+56+7\times 89).$
\item [] $\mathit{10246=1-2\times 3+(4+5)\times 67\times (8+9).}$
\item [] $10247=12+(3+45+67)\times 89.$
\item [] $10248=1\times 2\times (3+4)\times (5\times 6+78\times 9).$
\item [] $10249=1+2\times (3+4)\times (5\times 6+78\times 9).$
\item [] $10250=1+2^{(3+4)}\times (5+67+8)+9.$
\item [] $10251=(1+23\times 4\times 5+678)\times 9.$
\item [] $10252=1^{23}+(4+5)\times 67\times (8+9).$
\item [] $\mathit{10253=1-2+3+(4+5)\times 67\times (8+9).}$
\item [] $10254=(1+2+3\times 4)\times (5+678)+9.$
\item [] $10255=1^2+3+(4+5)\times 67\times (8+9).$
\item [] $10256=1\times 2+3+(4+5)\times 67\times (8+9).$
\item [] $10257=(1234+5+6\times 7)\times 8+9.$
\item [] $10258=1+2\times 3+(4+5)\times 67\times (8+9).$
\item [] $10259=1\times 2^3+(4+5)\times 67\times (8+9).$
\item [] $10260=(12+34+5\times 6)\times (7+8)\times 9.$
\item [] $10261=1^2+3\times 4\times (5+6\times (7+8))\times 9.$
\item [] $10262=1\times 2+3\times 4\times (5+6\times (7+8))\times 9.$
\item [] $10263=1+2+3\times (4+(5+6\times 7)\times 8)\times 9.$
\item [] $10264=1\times (2+3)^4+567\times (8+9).$
\item [] $10265=1+(2+3)^4+567\times (8+9).$
\item [] $10266=1\times 2\times 3\times (4^5+678+9).$
\item [] $10267=1+2\times 3\times (4^5+678+9).$
\item [] $10268=(1+2+34+567)\times (8+9).$
\item [] $10269=123\times 45+6\times 789.$
\item [] $10270=1\times (2+3)\times (4\times 5+6)\times (7+8\times 9).$
\item [] $10271=1+(2+3)\times (4\times 5+6)\times (7+8\times 9).$
\item [] $10272=(1+23)\times (4+5\times 67+89).$
\item [] $\mathit{10273=-1+23+(4+5)\times 67\times (8+9).}$
\item [] $10274=1\times 23+(4+5)\times 67\times (8+9).$
\item [] $10275=1+23+(4+5)\times 67\times (8+9).$
\item [] $\mathit{10276=1+(2+3\times 45)\times (6+78-9).}$
\item [] $\mathit{10277=-1+2\times 3^4\times (56+7)+8\times 9.}$
\item [] $10278=((12+3+4)\times 56+78)\times 9.$
\item [] $10279=1+2\times 3^4\times (56+7)+8\times 9.$
\item [] $10280=1\times 2+(3\times (4+5)\times 6\times 7+8)\times 9.$
\item[]$\mbox{Decreasing order}$
\item [] $10211=9876+5\times (4+3\times 21).$
\item [] $10212=9\times (8+7+6)\times 54+3\times 2\times 1.$
\item [] $10213=9\times (8+7+6)\times 54+3\times 2+1.$
\item [] $10214=(9\times 8+76)\times (5+4^3)+2\times 1.$
\item [] $10215=9+(8+7+6)\times 54\times 3^2\times 1.$
\item [] $10216=9+(8+7+6)\times 54\times 3^2+1.$
\item [] $10217=9876+5\times 4+321.$
\item [] $10218=9\times 8+(7+(6+5+4)^3)\times (2+1).$
\item [] $10219=9+8+(7+6\times 5+4^3)^2+1.$
\item [] $\mathit{10220=(9\times 8+7-6)\times 5\times (4+3+21).}$
\item [] $10221=9+87+(6+5+4)^3\times (2+1).$
\item [] $\mathit{10222=(-9+8\times 7\times (65-4))\times 3+2-1.}$
\item [] $10223=9+8+7\times 6\times (5+4)\times 3^{(2+1)}.$
\item [] $10224=(9\times 8\times 7+6)\times 5\times 4+3+21.$
\item [] $10225=9+8\times 765+4^{(3\times 2)}\times 1.$
\item [] $10226=9+8\times 765+4^{(3\times 2)}+1.$
\item [] $10227=(98+7\times (6+5)\times 43)\times (2+1).$
\item [] $\mathit{10228=-98\times 7+(6\times 5+4)\times 321.}$
\item [] $\mathit{10229=(-9+87\times 6)\times 5\times 4-32+1.}$
\item [] $10230=9\times (8+7+6)\times 54+3+21.$
\item [] $\mathit{10231=(9\times 8\times 7+6)\times 5\times 4+32-1.}$
\item [] $10232=(9\times 8\times 7+6)\times 5\times 4+32\times 1.$
\item [] $10233=9876+(5+4\times 3)\times 21.$
\item [] $10234=(9\times 8+7\times 6+5)\times 43\times 2\times 1.$
\item [] $10235=(9\times 8+7\times 6+5)\times 43\times 2+1.$
\item [] $10236=9\times 8+7\times (6+5)\times 4\times (32+1).$
\item [] $\mathit{10237=9\times (8+7+6)\times 54+32-1.}$
\item [] $10238=987+6+5\times 43^2\times 1.$
\item [] $10239=987+6+5\times 43^2+1.$
\item [] $\mathit{10240=(9\times 8+7+6-5)\times 4\times 32\times 1.}$
\item [] $10241=98\times 7+65\times (4+3)\times 21.$
\item [] $10242=9\times (8\times 76+(5\times 4+3)^2+1).$
\item [] $10243=987+6+5\times (43^2+1).$
\item [] $10244=(9+8)\times 7+(6+5+4)^3\times (2+1).$
\item [] $\mathit{10245=((9\times 8+7)\times 65-4\times 3)\times 2-1.}$
\item [] $\mathit{10246=(9-8\times 7)\times (-6-5\times 43+2+1).}$
\item [] $\mathit{10247=-9+8+7\times (6+(5+4)^3\times 2\times 1).}$
\item [] $10248=987\times 6+5+4321.$
\item [] $10249=9+(8+7+65)\times 4\times 32\times 1.$
\item [] $10250=9+(8+7+65)\times 4\times 32+1.$
\item [] $10251=9876+54+321.$
\item [] $10252=(9+8)\times (7+6+54)\times 3^2+1.$
\item [] $\mathit{10253=(-9+87\times 6)\times 5\times 4-3\times 2-1.}$
\item [] $10254=9876+54\times (3\times 2+1).$
\item [] $10255=(9+8\times 7\times (6+5\times 4))\times (3\times 2+1).$
\item [] $10256=9\times 8+76\times (5+4\times 32+1).$
\item [] $10257=9+8\times (7+6+5+43)\times 21.$
\item [] $10258=(9+(8+7+65)\times 4^3)\times 2\times 1.$
\item [] $10259=(9+(8+7+65)\times 4^3)\times 2+1.$
\item [] $10260=(9\times 8+7\times 6)\times (5+4^3+21).$
\item [] $10261=9\times (87\times 6+5+43)\times 2+1.$
\item [] $10262=9\times (8\times 7+6\times 54)\times 3+2\times 1.$
\item [] $10263=(9\times 8\times 7+6)\times 5\times 4+3\times 21.$
\item [] $\mathit{10264=-987+6\times 5^4\times 3+2-1.}$
\item [] $10265=9+8\times (7+6+5^4+3)\times 2\times 1.$
\item [] $10266=(98+76)\times (54+3+2\times 1).$
\item [] $10267=(98+76)\times (54+3+2)+1.$
\item [] $10268=(9+8)\times ((7+6+54)\times 3^2+1).$
\item [] $10269=9\times (8+7+6)\times 54+3\times 21.$
\item [] $10270=(9+8\times 7)\times (6\times 5+4\times 32\times 1).$
\item [] $10271=(9+8\times 7)\times (6\times 5+4\times 32)+1.$
\item [] $10272=(9+87+6+5)\times 4\times (3+21).$
\item [] $10273=9+8\times (76\times 5+43\times 21).$
\item [] $10274=9\times 8+(7+6\times 5+4^3)^2+1.$
\item [] $10275=(9+8\times (7+6+5^4+3))\times 2+1.$
\item [] $\mathit{10276=-98+7+6\times 54\times 32-1.}$
\item [] $\mathit{10277=-98+7+6\times 54\times 32\times 1.}$
\item [] $10278=9\times 8\times 7+6\times 543\times (2+1).$
\item [] $10279=(9+(8+7)\times 6\times (54+3))\times 2+1.$
\item [] $10280=9\times (8+7\times 6\times (5+4)\times 3)+2\times 1.$
\item[]$\mbox{Increasing order}$
\item [] $10281=1\times 2\times 3\times (4+5\times 6\times 7)\times 8+9.$
\item [] $10282=1+2\times 3\times (4+5\times 6\times 7)\times 8+9.$
\item [] $\mathit{10283=(12^3+4-5)\times 6-7-8\times 9.}$
\item [] $10284=12\times (34\times 5+678+9).$
\item [] $10285=(1\times 23\times (4\times 5+6)+7)\times (8+9).$
\item [] $10286=(1+2\times 34+5)\times (67+8\times 9).$
\item [] $10287=12^3\times 4+5\times (67+8)\times 9.$
\item [] $10288=1^2+3\times (45+6\times 7\times 8)\times 9.$
\item [] $10289=1\times 2+3\times (45+6\times 7\times 8)\times 9.$
\item [] $10290=1+2+3\times (45+6\times 7\times 8)\times 9.$
\item [] $10291=1+2\times (3\times (4+5\times 6\times 7)\times 8+9).$
\item [] $\mathit{10292=(-1-2+34)\times (5+6\times 7\times 8-9).}$
\item [] $10293=(1\times 2\times 34+5)\times (6+(7+8)\times 9).$
\item [] $10294=1+(2\times 34+5)\times (6+(7+8)\times 9).$
\item [] $10295=1\times 2\times 3^4\times (56+7)+89.$
\item [] $10296=(12\times 3+4^5+6+78)\times 9.$
\item [] $10297=1+2\times 3\times 4\times (5\times (6+78)+9).$
\item [] $\mathit{10298=1\times 2+(3^4-5+67)\times 8\times 9.}$
\item [] $10299=12+3\times (45+6\times 7\times 8)\times 9.$
\item [] $10300=1\times (2+3)\times 4\times (5+6+7\times 8\times 9).$
\item [] $10301=1+(2+3)\times 4\times (5+6+7\times 8\times 9).$
\item [] $10302=(1+2)\times (3\times 45+67)\times (8+9).$
\item [] $10303=1^2+(3+(4+5)\times 67)\times (8+9).$
\item [] $10304=(1+2^3)^4+5+6\times 7\times 89.$
\item [] $10305=1^2\times 3\times 4\times (5+6)\times 78+9.$
\item [] $10306=1+2\times (3\times 4\times 5+6)\times 78+9.$
\item [] $10307=1\times 2+3\times 4\times (5+6)\times 78+9.$
\item [] $10308=1+2+3\times 4\times (5+6)\times 78+9.$
\item [] $\mathit{10309=(12^3-4+5)\times 6-7\times 8-9.}$
\item [] $\mathit{10310=(1\times 2-3-4^5-6)\times (7-8-9).}$
\item [] $10311=12^3\times 4+5\times 678+9.$
\item [] $\mathit{10312=-1+23\times 4\times (56+7\times 8)+9.}$
\item [] $10313=1\times 23\times 4\times (56+7\times 8)+9.$
\item [] $10314=(12+3)\times (4+5+678)+9.$
\item [] $10315=1+2\times ((3\times 4\times 5+6)\times 78+9).$
\item [] $\mathit{10316=1-2+3\times (4+5\times (678+9)).}$
\item [] $10317=12+3\times 4\times (5+6)\times 78+9.$
\item [] $10318=1\times 2\times (3^4\times 56+7\times 89).$
\item [] $10319=1+2\times (3^4\times 56+7\times 89).$
\item [] $10320=12\times (3+4\times 5\times 6\times 7+8+9).$
\item [] $10321=(1+2\times (3^4+5+6)\times 7)\times 8+9.$
\item [] $10322=1\times 2\times ((3^4+5+6)\times 7\times 8+9).$
\item [] $10323=123+4\times 5\times (6+7\times 8\times 9).$
\item [] $10324=(1^2+3+45+67)\times 89.$
\item [] $10325=(1^2+34)\times 5\times (6\times 7+8+9).$
\item [] $10326=1+2+3\times (4\times (5+6)\times 78+9).$
\item [] $10327=1+2\times 3\times ((4+5\times 6\times 7)\times 8+9).$
\item [] $\mathit{10328=1+23\times (456+7\times (8-9)).}$
\item [] $10329=12+3\times (4+5\times (678+9)).$
\item [] $10330=1\times 2\times (3+(45+6+7)\times 89).$
\item [] $10331=1+2\times (3+(45+6+7)\times 89).$
\item [] $10332=123\times (4+56+7+8+9).$
\item [] $\mathit{10333=12\times (3+4\times 5\times 6)\times 7-8+9.}$
\item [] $\mathit{10334=-1\times 2+34\times (5\times (67-8)+ 9).}$
\item [] $10335=(1^2+3+4+5)\times (6+789).$
\item [] $10336=1\times 2\times 34\times (56+7+89).$
\item [] $10337=1+2\times 34\times (56+7+89).$
\item [] $10338=1\times 2\times (345+67\times 8\times 9).$
\item [] $10339=1+2\times (345+67\times 8\times 9).$
\item [] $\mathit{10340=(12^3+4+5)\times 6+7-89.}$
\item [] $10341=123\times (4+5+67+8)+9.$
\item [] $10342=1+(2+(3\times 4+5)\times 67+8)\times 9.$
\item [] $10343=12\times 3\times 45\times 6+7\times 89.$
\item [] $10344=12\times (3+4+(5+6\times (7+8))\times 9).$
\item [] $\mathit{10345=1+23\times (456-7)+8+9.}$
\item [] $\mathit{10346=(1+2)\times (3456-7)+8-9.}$
\item [] $10347=12^3\times 4+5\times (678+9).$
\item [] $\mathit{10348=-1\times 2+345\times (6+7+8+9).}$
\item [] $10349=12\times (3+4\times 5\times 6)\times 7+8+9.$
\item [] $10350=1^2\times 345\times (6+7+8+9).$
\item[]$\mbox{Decreasing order}$
\item [] $10281=9\times (8\times 7+6\times 54)\times 3+21.$
\item [] $10282=98+76\times (5+4\times 32+1).$
\item [] $10283=98+7\times (6+(5+4^3)\times 21).$
\item [] $10284=((9\times 8+7)\times 65+4+3)\times 2\times 1.$
\item [] $10285=((9\times 8+7)\times 65+4+3)\times 2+1.$
\item [] $10286=(9+8)\times (76+(5\times 4+3)^2)+1.$
\item [] $10287=9\times 87+6^5+(4\times 3)^{(2+1)}.$
\item [] $10288=(9+(8+7+6)\times 54)\times 3^2+1.$
\item [] $10289=9+8\times ((7\times 6\times 5+4)\times 3\times 2+1).$
\item [] $10290=(98+76\times 5+4\times 3)\times 21.$
\item [] $10291=(9+8\times (7\times 6\times 5+4)\times 3)\times 2+1.$
\item [] $\mathit{10292=98\times 7\times (6+5+4)+3-2+1.}$
\item [] $\mathit{10293=(9\times 8+7-6)\times (54\times 3-21).}$
\item [] $10294=(987+65\times 4^3)\times 2\times 1.$
\item [] $10295=(987+65\times 4^3)\times 2+1.$
\item [] $10296=(98+7+6\times 54)\times (3+21).$
\item [] $10297=987\times 6+5^4\times (3\times 2+1).$
\item [] $\mathit{10298=98\times 7\times (6+5+4)+3^2-1.}$
\item [] $10299=98\times 7\times (6+5+4)+3^2\times 1.$
\item [] $10300=98\times 7\times (6+5+4)+3^2+1.$
\item [] $10301=9876+5\times (4^3+21).$
\item [] $10302=(9+8)\times (7\times 6+543+21).$
\item [] $\mathit{10303=9876-5+432\times 1.}$
\item [] $10304=(9+87+65)\times (43+21).$
\item [] $10305=9\times (8\times 7+65+4^(3+2)\times 1).$
\item [] $10306=9876+5\times 43\times 2\times 1.$
\item [] $10307=9876+5\times 43\times 2+1.$
\item [] $10308=9+8+7\times 6\times 5\times (4+3)^2+1.$
\item [] $\mathit{10309=-9\times 8+7+6\times (54\times 32+1).}$
\item [] $10310=(987+(6+5)\times 4)\times (3^2+1).$
\item [] $10311=9876+5\times (43\times 2+1).$
\item [] $10312=98\times 7+6^5+43^2+1.$
\item [] $10313=9876+5+432\times 1.$
\item [] $10314=9876+5+432+1.$
\item [] $10315=(9\times 87+6\times (5+4)^3)\times 2+1.$
\item [] $10316=9+(8+7\times (6+(5+4)^3))\times 2+1.$
\item [] $10317=98\times 7\times (6+5+4)+3^{(2+1)}.$
\item [] $\mathit{10318=(9-8+7)\times 6\times 5\times 43-2\times 1.}$
\item [] $10319=(9+8)\times (7+6\times 5\times 4\times (3+2)\times 1).$
\item [] $10320=9\times 8+7\times (6+(5+4)^3\times 2\times 1).$
\item [] $10321=9\times 8+7\times (6+(5+4)^3\times 2)+1.$
\item [] $10322=98\times 7\times (6+5+4)+32\times 1.$
\item [] $10323=98\times 7\times (6+5+4)+32+1.$
\item [] $10324=9\times 87\times 6+5^4\times 3^2+1.$
\item [] $10325=(9+87+65)\times 4^3+21.$
\item [] $10326=(9+8\times (7\times 6\times 5+4))\times 3\times 2\times 1.$
\item [] $10327=9\times 8+7\times (6+(5+4)^3\times 2+1).$
\item [] $10328=9\times (87\times 6+5^4)+3+2\times 1.$
\item [] $10329=9\times (87\times 6+5^4)+3+2+1.$
\item [] $10330=9\times (87\times 6+5^4)+3\times 2+1.$
\item [] $\mathit{10331=9\times (87\times 6+5^4)+3^2-1.}$
\item [] $10332=9+(87\times 6+5^4)\times 3^2\times 1.$
\item [] $10333=9+(87\times 6+5^4)\times 3^2+1.$
\item [] $\mathit{10334=(9+8\times (7-654))\times (-3+2-1).}$
\item [] $10335=(9+8\times 7)\times (6\times 5+4\times 32+1).$
\item [] $10336=(9+87)\times 65+4^{(3\times 2)}\times 1.$
\item [] $10337=98\times 76+(5+4)\times 321.$
\item [] $10338=9\times 87+65\times (4+3)\times 21.$
\item [] $\mathit{10339=98+7\times (6+(5+4)^3\times 2-1).}$
\item [] $\mathit{10340=-9\times 8+76\times (5+4\times (32+1)).}$
\item [] $10341=9\times (8\times 7+6+543\times 2+1).$
\item [] $\mathit{10342=9+8\times 76\times (5+4\times 3)-2-1.}$
\item [] $\mathit{10343=-9-8-7+6\times 54\times 32-1.}$
\item [] $10344=(9+8+7)\times (6+5\times (4^3+21)).$
\item [] $\mathit{10345=-9-8-7+6\times 54\times 32+1.}$
\item [] $10346=98+7\times (6+(5+4)^3\times 2\times 1).$
\item [] $10347=987+65\times (4\times 3)^2\times 1.$
\item [] $10348=987+65\times (4\times 3)^2+1.$
\item [] $10349=(9\times 8+7)\times (65+4^3+2)\times 1.$
\item [] $10350=(9\times 8+7)\times (65+4^3+2)+1.$
\item[]$\mbox{Increasing order}$
\item [] $10351=1^2+345\times (6+7+8+9).$
\item [] $10352=1\times 2+345\times (6+7+8+9).$
\item [] $10353=1+2+345\times (6+7+8+9).$
\item [] $10354=1+(2\times 3+(4+5)\times 67)\times (8+9).$
\item [] $10355=(1\times 2+3\times 4+5)\times (67\times 8+9).$
\item [] $10356=1\times 2\times 3\times 4^5+6\times 78\times 9.$
\item [] $10357=1+2\times 3\times 4^5+6\times 78\times 9.$
\item [] $10358=(1+2\times 3^4)\times (56+7)+89.$
\item [] $10359=123\times 45+67\times 8\times 9.$
\item [] $10360=1+23\times (4+5)\times (6\times 7+8)+9.$
\item [] $\mathit{10361=(1+2)\times 3456+7\times (8-9).}$
\item [] $10362=12+345\times (6+7+8+9).$
\item [] $10363=(1\times 2\times 3^4+5)\times (6+7\times 8)+9.$
\item [] $10364=1+(2\times 3^4+5)\times (6+7\times 8)+9.$
\item [] $10365=1+2\times (3+4+(567+8)\times 9).$
\item [] $\mathit{10366=(12^3-4+5)\times 6-7+8-9.}$
\item [] $\mathit{10367=12\times 3\times 4\times (5+67)+8-9.}$
\item [] $10368=12\times (3+4+5+6+78)\times 9.$
\item [] $10369=1+2\times 3^4\times (5+6\times 7+8+9).$
\item [] $10370=1^2\times 34\times ((5\times 6+7)\times 8+9).$
\item [] $10371=1^2+34\times ((5\times 6+7)\times 8+9).$
\item [] $10372=1\times 2+34\times ((5\times 6+7)\times 8+9).$
\item [] $10373=(1+2^3)\times 4^5+(6+7)\times 89.$
\item [] $10374=123+(4+5)\times 67\times (8+9).$
\item [] $10375=1+2\times 3\times (4+5\times (6\times 7\times 8+9)).$
\item [] $\mathit{10376=(1+2)\times 3456+7-8+9.}$
\item [] $10377=1\times 2\times (3^4+567)\times 8+9.$
\item [] $10378=1+2\times (3^4+567)\times 8+9.$
\item [] $\mathit{10379=(12^3-4-5)\times 6-7+8\times 9.}$
\item [] $10380=(1+2+3\times 4)\times (5+678+9).$
\item [] $\mathit{10381=(12^3-4+5)\times 6-7\times (8-9).}$
\item [] $10382=1234\times 5+6\times 78\times 9.$
\item [] $10383=(1^2+3\times 4\times (5+6))\times 78+9.$
\item [] $10384=1\times 2\times (3^4\times (56+7)+89).$
\item [] $10385=(1234+56+7)\times 8+9.$
\item [] $10386=(12+3\times (4+5)\times 6\times 7+8)\times 9.$
\item [] $10387=(1+2^3+4)\times (5+6\times 7)\times (8+9).$
\item [] $10388=1+(2^3+(4+5)\times 67)\times (8+9).$
\item [] $10389=12\times (3+4+(5+6)\times 78)+9.$
\item [] $\mathit{10390=-1\times 23+(4+5)\times (6+7)\times 89.}$
\item [] $\mathit{10391=-1+23\times 456-7-89.}$
\item [] $10392=(1+2)\times 3456+7+8+9.$
\item [] $10393=1+2\times 3\times (4^5+6+78\times 9).$
\item [] $\mathit{10394=1-2+(3+4)\times (5+6)\times (7+8)\times 9.}$
\item [] $10395=(1+2^3)\times (4^5+6\times 7+89).$
\item [] $10396=(12+34)\times (5+(6+7)\times (8+9)).$
\item [] $10397=1\times 2+(3+4)\times (5+6)\times (7+8)\times 9.$
\item [] $10398=1+2+(3+4)\times (5+6)\times (7+8)\times 9.$
\item [] $\mathit{10399=((1+2)\times 3^4+5)\times 6\times 7-8-9.}$
\item [] $10400=(1+2^3+4)\times (5+6+789).$
\item [] $10401=12\times (3\times 45\times 6+7\times 8)+9.$
\item [] $10402=((1+2)^3+4)\times 5\times 67+8+9.$
\item [] $\mathit{10403=12^3-4+(5+6)\times 789.}$
\item [] $10404=12\times 3\times (4\times 56+7\times 8+9).$
\item [] $10405=1^2+34\times (5+6+7)\times (8+9).$
\item [] $10406=(1+2)\times (3456+7)+8+9.$
\item [] $10407=1+2+34\times (5+6+7)\times (8+9).$
\item [] $\mathit{10408=-1+23\times 456-7-8\times 9.}$
\item [] $10409=1\times (2+3)\times 4\times 5\times (6+7)\times 8+9.$
\item [] $10410=(1\times 2+345)\times (6+7+8+9).$
\item [] $10411=12^3+4+(5+6)\times 789.$
\item [] $\mathit{10412=1\times 2-3+(4+5)\times (6+7)\times 89.}$
\item [] $10413=(1\times 2+3+45+67)\times 89.$
\item [] $10414=1^{23}+(4+5)\times (6+7)\times 89.$
\item [] $\mathit{10415=(1+2)\times 3456+7\times 8-9.}$
\item [] $10416=12+34\times (5+6+7)\times (8+9).$
\item [] $10417=1^2+3+(4+5)\times (6+7)\times 89.$
\item [] $10418=1\times 2+3+(4+5)\times (6+7)\times 89.$
\item [] $10419=1\times 2\times 3+(4+5)\times (6+7)\times 89.$
\item [] $10420=1+2\times 3+(4+5)\times (6+7)\times 89.$
\item[]$\mbox{Decreasing order}$
\item [] $10351=(9\times 8\times (7+6)+5)\times (4+3\times 2+1).$
\item [] $10352=9\times (8+7\times 6)\times (5\times 4+3)+2\times 1.$
\item [] $10353=98\times 7\times (6+5+4)+3\times 21.$
\item [] $10354=(9+8)\times 7\times (6+(5+4)\times 3^2)+1.$
\item [] $10355=9+8\times 7+6\times 5\times (4+3)^{(2+1)}.$
\item [] $10356=9876+5\times 4\times (3+21).$
\item [] $10357=((9\times 8+7)\times 65+43)\times 2+1.$
\item [] $10358=(98+765)\times 4\times 3+2\times 1.$
\item [] $10359=9876+(5\times 4+3)\times 21.$
\item [] $10360=(9\times 8+76)\times 5\times (4+3^2+1).$
\item [] $10361=(9\times 8+76)\times 5\times (4+3)\times 2+1.$
\item [] $10362=9876+54\times 3^2\times 1.$
\item [] $10363=9876+54\times 3^2+1.$
\item [] $\mathit{10364=-9-8+7+6\times (54\times 32+1).}$
\item [] $10365=((98+765)\times 4+3)\times (2+1).$
\item [] $10366=9+8\times 76\times (5+4\times 3)+21.$
\item [] $\mathit{10367=-9+8+(76+5)\times 4\times 32\times 1.}$
\item [] $10368=(9+87)\times (6+5+43)\times 2\times 1.$
\item [] $10369=(9+87)\times (6+5+43)\times 2+1.$
\item [] $10370=(9+87)\times (65+43)+2\times 1.$
\item [] $10371=(9+87)\times (65+43)+2+1.$
\item [] $\mathit{10372=9+8-7+6\times (54\times 32-1).}$
\item [] $\mathit{10373=-9+8+7+6\times 54\times 32-1.}$
\item [] $10374=(98+76+5\times 4^3)\times 21.$
\item [] $10375=(9\times 8+7\times 6)\times (5+43\times 2)+1.$
\item [] $\mathit{10376=9-8+7+6\times 54\times 32\times 1.}$
\item [] $10377=(98+765)\times 4\times 3+21.$
\item [] $10378=9+8\times (7+6+5\times 4+3)^2+1.$
\item [] $\mathit{10379=9+8-7+6\times 54\times 32+1.}$
\item [] $10380=(98\times 7+6)\times (5+4+3+2+1).$
\item [] $10381=(98\times 7+6)\times (5+4+3\times 2)+1.$
\item [] $10382=(((9+8)\times 76+5)\times 4+3)\times 2\times 1.$
\item [] $10383=9+(8+7\times (65+4)+3)\times 21.$
\item [] $\mathit{10384=9+8+(76+5)\times 4\times 32-1.}$
\item [] $10385=9+8+(76+5)\times 4\times 32\times 1.$
\item [] $10386=9+8+(76+5)\times 4\times 32+1.$
\item [] $10387=9+8\times (7+6\times 5\times 43)+2\times 1.$
\item [] $10388=98\times (7\times 6+54+3^2+1).$
\item [] $10389=(9+87)\times (65+43)+21.$
\item [] $\mathit{10390=-98+76\times (5+4^3)\times 2\times 1.}$
\item [] $\mathit{10391=9+8+7+6\times 54\times 32-1.}$
\item [] $10392=9+8+7+6\times 54\times 32\times 1.$
\item [] $10393=9+8+7+6\times 54\times 32+1.$
\item [] $10394=((98\times 7+6)\times 5+4)\times 3+2\times 1.$
\item [] $10395=(9+87+65+4)\times 3\times 21.$
\item [] $10396=(98+7)\times (6+5)\times (4+3+2)+1.$
\item [] $10397=9\times (8+7)\times (65+4\times 3)+2\times 1.$
\item [] $10398=9+8+7+6\times (54\times 32+1).$
\item [] $10399=((9\times 8+7)\times 65+4^3)\times 2+1.$
\item [] $10400=(9+8+7\times (6+5)\times 4)\times 32\times 1.$
\item [] $10401=(9+8\times 7+65\times 4)\times 32+1.$
\item [] $10402=(9+8+(76+5)\times 4^3)\times 2\times 1.$
\item [] $10403=(9+8+(76+5)\times 4^3)\times 2+1.$
\item [] $10404=(9+87\times 6+5^4)\times 3^2\times 1.$
\item [] $10405=(9+8\times 7+6\times 5+4+3)^2+1.$
\item [] $10406=9876+(5\times 4+3)^2+1.$
\item [] $\mathit{10407=9\times 87+6^5+43^2-1.}$
\item [] $10408=9\times 87+6^5+43^2\times 1.$
\item [] $10409=9\times 87+6^5+43^2+1.$
\item [] $10410=9+8\times (7+6)\times 5\times 4\times (3+2)+1.$
\item [] $\mathit{10411=-9+8+76\times (5+4\times (32+1)).}$
\item [] $10412=987+65\times ((4\times 3)^2+1).$
\item [] $10413=((98\times 7+6)\times 5+4)\times 3+21.$
\item [] $10414=9\times (8+76+5)\times (4+3^2)+1.$
\item [] $10415=9+(8\times 7+65)\times 43\times 2\times 1.$
\item [] $10416=9876+54\times (3^2+1).$
\item [] $10417=9+8\times (76+(5\times (4+3))^2\times 1).$
\item [] $10418=9+8\times (76+(5\times (4+3))^2)+1.$
\item [] $10419=9+(8+7)\times (6+5^4+3\times 21).$
\item [] $\mathit{10420=9876+543+2-1.}$
\item[]$\mbox{Increasing order}$
\item [] $10421=(12+34+567)\times (8+9).$
\item [] $10422=12\times 3\times 45\times 6+78\times 9.$
\item [] $10423=(1+2\times 3)\times (4+(5+6)\times (7+8)\times 9).$
\item [] $10424=1\times 2+(3\times ((4+5)\times 6\times 7+8)\times 9).$
\item [] $10425=(123+45)\times (6+7\times 8)+9.$
\item [] $10426=1+((2+34)\times 5+6)\times 7\times 8+9.$
\item [] $\mathit{10427=1-2+3\times 4\times (5+6)\times (7+8\times 9).}$
\item [] $10428=12+3+(4+5)\times (6+7)\times 89.$
\item [] $10429=1^2+3\times 4\times (5+6)\times (7+8\times 9).$
\item [] $10430=1\times 2+3\times 4\times (5+6)\times (7+8\times 9).$
\item [] $10431=(1+2\times 345+6\times 78)\times 9.$
\item [] $10432=(1+2\times 3^4)\times (5+6\times 7+8+9).$
\item [] $10433=(1+2)\times 3456+7\times 8+9.$
\item [] $10434=(1+2+34)\times (5\times 6\times 7+8\times 9).$
\item [] $\mathit{10435=1-2\times 3+4\times 5\times 6\times (78+9).}$
\item [] $10436=1\times 23+(4+5)\times (6+7)\times 89.$
\item [] $10437=1+23+(4+5)\times (6+7)\times 89.$
\item [] $10438=(1\times 2+34\times ((5+6)+7))\times (8+9).$
\item [] $10439=1+(2+34\times ((5+6)+7))\times (8+9).$
\item [] $10440=(1+2)\times (3456+7+8+9).$
\item [] $10441=1^{23}+4\times 5\times 6\times (78+9).$
\item [] $10442=(12+34)\times (5\times 6\times 7+8+9).$
\item [] $10443=1^2\times 3+4\times 5\times 6\times (78+9).$
\item [] $10444=1^2+3+4\times 5\times 6\times (78+9).$
\item [] $10445=1\times 2+3+4\times 5\times 6\times (78+9).$
\item [] $10446=1\times 2\times 3+4\times 5\times 6\times (78+9).$
\item [] $10447=(1+2)\times 3456+7+8\times 9.$
\item [] $10448=1\times 2^3+4\times 5\times 6\times (78+9).$
\item [] $10449=1+2^3+4\times 5\times 6\times (78+9).$
\item [] $10450=(12+3+4)\times (5+67\times 8+9).$
\item [] $10451=1\times 23\times (4+5\times 6\times (7+8))+9.$
\item [] $10452=12^3\times 4+5\times (6+78\times 9).$
\item [] $10453=(1+2\times 3)\times (4^5+6\times 78)+9.$
\item [] $\mathit{10454=1\times 2+3\times (4-5\times (6-78\times 9)).}$
\item [] $10455=(1+2)\times 3456+78+9.$
\item [] $10456=(1^2+3)\times (4+5\times 6\times (78+9)).$
\item [] $10457=(1234+5+67)\times 8+9.$
\item [] $10458=1\times 2\times (3^4\times 5+67\times 8\times 9).$
\item [] $10459=1+2\times (3^4\times 5+67\times 8\times 9).$
\item [] $\mathit{10460=-1+(2+3\times 4\times (5+6))\times 78+9.}$
\item [] $10461=(1+2)\times (3456+7)+8\times 9.$
\item [] $10462=1+(2+3\times 4\times (5+6))\times 78+9.$
\item [] $10463=1\times 23+4\times 5\times 6\times (78+9).$
\item [] $10464=(1+2)\times 3456+7+89.$
\item [] $10465=(12\times 3^4+5\times 67)\times 8+9.$
\item [] $\mathit{10466=1-23\times (4-5\times 6\times (7+8)-9).}$
\item [] $10467=(1+2^3)\times (4^5+67+8\times 9).$
\item [] $10468=1+((2+3)\times (4\times 56+7)+8)\times 9.$
\item [] $\mathit{10469=(12^3+4+5)\times 6+7\times 8-9.}$
\item [] $\mathit{10470=(12^3-4+5)\times 6+7+89.}$
\item [] $\mathit{10471=-1+(2^3\times 4+56)\times 7\times (8+9).}$
\item [] $10472=(1\times 2^3\times 4+56)\times 7\times (8+9).$
\item [] $10473=1+(2^3\times 4+56)\times 7\times (8+9).$
\item [] $10474=((1+2)^3+4)\times 5\times 67+89.$
\item [] $\mathit{10475=-1+2\times (3^4\times 56+78\times 9).}$
\item [] $10476=1\times 2\times (3+4+567+8)\times 9.$
\item [] $10477=1+2\times (3+4+567+8)\times 9.$
\item [] $10478=(1+2)\times (3456+7)+89.$
\item [] $10479=1+2\times (3+4^5+6\times 78\times 9).$
\item [] $10480=(1^2+3)\times 4\times 5\times (6\times 7+89).$
\item [] $10481=(1+(2+34)\times 5+6)\times 7\times 8+9.$
\item [] $10482=12^3\times 4+5\times 6\times 7\times (8+9).$
\item [] $\mathit{10483=1+23\times 456-7-8+9.}$
\item [] $\mathit{10484=-1-(2+3-(4+5+6)\times 78)\times 9.}$
\item [] $10485=(12+3)\times (4+5\times (67+8\times 9)).$
\item [] $10486=(1+2\times 3)\times (4^5+6\times (7+8\times 9)).$
\item [] $10487=(1+2)\times 3456+7\times (8+9).$
\item [] $10488=(1+2\times 34)\times (56+7+89).$
\item [] $10489=1+23\times 4\times (5\times (6+7+8)+9).$
\item [] $\mathit{10490=1+23\times 456+(-7+8)^9.}$
\item[]$\mbox{Decreasing order}$
\item [] $10421=9876+543+2\times 1.$
\item [] $10422=9876+543+2+1.$
\item [] $10423=9\times (87\times 6+54+3)\times 2+1.$
\item [] $10424=(9+(8\times 7+65)\times 43)\times 2\times 1.$
\item [] $10425=9+8\times (7+6+5+4\times 321).$
\item [] $10426=9+8\times 7\times (6+5\times 4\times 3^2)+1.$
\item [] $\mathit{10427=-9+87\times 6\times 5\times 4-3-2+1.}$
\item [] $10428=(9\times 8+7)\times (65+4+3\times 21).$
\item [] $10429=(9\times 8+7)\times (6+5\times 4\times 3)\times 2+1.$
\item [] $10430=(9\times 87+65\times 4)\times (3^2+1).$
\item [] $10431=(9\times 8+7)\times (6+5)\times 4\times 3+2+1.$
\item [] $10432=(98+(7+6)\times 5)\times (43+21).$
\item [] $10433=9+8\times 7+6\times 54\times 32\times 1.$
\item [] $10434=9+8\times 7+6\times 54\times 32+1.$
\item [] $10435=(98+(7+6)\times 5)\times 4^3+2+1.$
\item [] $10436=98\times 7+6\times 5\times (4+321).$
\item [] $10437=987\times 6+5\times 43\times 21.$
\item [] $\mathit{10438=9876+5^4-3\times 21.}$
\item [] $10439=9+8\times 7+6\times (54\times 32+1).$
\item [] $10440=9876+543+21.$
\item [] $10441=9\times 8+(76+5)\times 4\times 32+1.$
\item [] $10442=(98+76)\times 5\times 4\times 3+2\times 1.$
\item [] $10443=9876+(5+4)\times 3\times 21.$
\item [] $\mathit{10444=9+87\times 6\times 5\times 4-3-2\times 1.}$
\item [] $10445=(9+8\times (7+6)\times 5\times 4)\times (3+2)\times 1.$
\item [] $10446=(9+8\times (7+6)\times 5\times 4)\times (3+2)+1.$
\item [] $10447=9\times 8+7+6\times 54\times 32\times 1.$
\item [] $10448=9\times 8+7+6\times 54\times 32+1.$
\item [] $10449=9+87\times 6\times (5+4\times 3+2+1).$
\item [] $10450=9+8\times 765+4321.$
\item [] $\mathit{10451=9+87\times 6\times 5\times 4+3-2+1.}$
\item [] $\mathit{10452=9+87\times 6\times 5\times 4+3\times (2-1).}$
\item [] $10453=9\times 8+7+6\times (54\times 32+1).$
\item [] $10454=9+87\times 6\times 5\times 4+3+2\times 1.$
\item [] $10455=9+87\times 6\times 5\times 4+3+2+1.$
\item [] $10456=9+87\times 6\times 5\times 4+3\times 2+1.$
\item [] $10457=9+8+(7+65)\times ((4\times 3)^2+1).$
\item [] $10458=9+87\times 6\times 5\times 4+3^2\times 1.$
\item [] $10459=9+87\times 6\times 5\times 4+3^2+1.$
\item [] $10460=98\times 7+6\times 543\times (2+1).$
\item [] $10461=(98+76)\times 5\times 4\times 3+21.$
\item [] $10462=9\times 87\times (6+5)+43^2\times 1.$
\item [] $10463=9\times 87\times (6+5)+43^2+1.$
\item [] $10464=9+87+6\times 54\times 32\times 1.$
\item [] $10465=9+87+6\times 54\times 32+1.$
\item [] $10466=98+(76+5)\times 4^3\times 2\times 1.$
\item [] $10467=98+(76+5)\times 4\times 32+1.$
\item [] $\mathit{10468=9876+5^4-32-1.}$
\item [] $10469=9+8\times 7+((6\times 5+4)\times 3)^2\times 1.$
\item [] $10470=9+87+6\times (54\times 32+1).$
\item [] $\mathit{10471=-9-8+76\times (5+4^3)\times 2\times 1.}$
\item [] $10472=(9+8)\times (76+54\times (3^2+1)).$
\item [] $10473=9+87\times 6\times 5\times 4+3+21.$
\item [] $10474=98+7+6\times 54\times 32+1.$
\item [] $10475=9+8+(7\times 65+43)\times 21.$
\item [] $10476=9\times (87\times 6+5\times 4\times 3)\times 2\times 1.$
\item [] $10477=9\times (87\times 6+5\times 4\times 3)\times 2+1.$
\item [] $10478=((9+8)\times 7\times (6+5)\times 4+3)\times 2\times 1.$
\item [] $10479=98+7+6\times (54\times 32+1).$
\item [] $\mathit{10480=9+87\times 6\times 5\times 4+32-1.}$
\item [] $10481=9+87\times 6\times 5\times 4+32\times 1.$
\item [] $10482=9+87\times 6\times 5\times 4+32+1.$
\item [] $10483=(9\times 87\times 6+543)\times 2+1.$
\item [] $10484=9\times 8+76\times (5+4\times (32+1)).$
\item [] $10485=(9+87\times 6\times 5)\times 4+3^2\times 1.$
\item [] $10486=(9+87\times 6\times 5)\times 4+3^2+1.$
\item [] $10487=(9+8)\times 7+6\times 54\times 32\times 1.$
\item [] $10488=(9+8)\times 7+6\times 54\times 32+1.$
\item [] $10489=9+8\times (7+6\times (5\times 43+2)+1).$
\item [] $10490=(9\times (8+76\times 5)+4)\times 3+2\times 1.$
\item[]$\mbox{Increasing order}$
\item [] $\mathit{10491=(12^3+4+5)\times 6+78-9.}$
\item [] $\mathit{10492=1+2+(3-4\times 5)\times (6-7\times 89).}$
\item [] $\mathit{10493=-1+23\times 456+7+8-9.}$
\item [] $10494=((1+2+34)\times 5\times 6+7\times 8)\times 9.$
\item [] $10495=1+(2\times (3\times 4+567)+8)\times 9.$
\item [] $10496=(1^2+3)\times 4\times (567+89).$
\item [] $10497=(1234+(5+6)\times 7)\times 8+9.$
\item [] $\mathit{10498=1\times 23\times 456-7+8+9.}$
\item [] $\mathit{10499=1+23\times 456-7+8+9.}$
\item [] $10500=12\times (3\times 45\times 6+7\times 8+9).$
\item [] $10501=(12^3+4+5)\times 6+7+8\times 9.$
\item [] $10502=(12+34+5+67)\times 89.$
\item [] $10503=(1+2)\times 3456+(7+8)\times 9.$
\item [] $10504=(1\times 2^3)^4+(5+67)\times 89.$
\item [] $10505=((1+2)\times 3^4+5)\times 6\times 7+89.$
\item [] $10506=(12\times 34+5\times 6\times 7)\times (8+9).$
\item [] $10507=(1+2\times 3)\times (4^5+6\times 78+9).$
\item [] $\mathit{10508=1-23+(4+5+6)\times 78\times 9.}$
\item [] $10509=12\times 3\times 45\times 6+789.$
\item [] $\mathit{10510=-1+23\times (4\times (56+7\times 8)+9).}$
\item [] $10511=1\times 23\times (4\times (56+7\times 8)+9).$
\item [] $10512=1\times 23\times 456+7+8+9.$
\item [] $10513=1+23\times 456+7+8+9.$
\item [] $10514=12+(3^4+5\times 6+7)\times 89.$
\item [] $10515=1+2+((3^4+5\times (6+7))\times 8)\times 9.$
\item [] $\mathit{10516=-1-23+4\times 5\times (67\times 8-9).}$
\item [] $\mathit{10517=-1\times 23+4\times 5\times (67\times 8-9).}$
\item [] $10518=(12^3+4+5)\times 6+7+89.$
\item [] $\mathit{10519=-12+34\times 5\times (6+7\times 8)-9.}$
\item [] $\mathit{10520=(1-23-4^5-6)\times (7-8-9).}$
\item [] $10521=(1\times (2\times 3)^4+5+6+7)\times 8+9.$
\item [] $10522=1+((2\times 3)^4+5+6+7)\times 8+9.$
\item [] $\mathit{10523=-1-2\times 3+(4+5+6)\times 78\times 9.}$
\item [] $10524=12+((3^4)+5\times (6+7))\times (8\times 9).$
\item [] $\mathit{10525=1-2\times 3+(4+5+6)\times 78\times 9.}$
\item [] $\mathit{10526=1-2-3+(4+5+6)\times 78\times 9.}$
\item [] $10527=(1^{23}+4\times 5\times 6)\times (78+9).$
\item [] $\mathit{10528=-1-2+34\times 5\times (6+7\times 8)-9.}$
\item [] $10529=(1+(2\times 3)^4+5+6+7)\times 8+9.$
\item [] $10530=1^{23}\times (4+5+6)\times 78\times 9.$
\item [] $10531=1^{23}+(4+5+6)\times 78\times 9.$
\item [] $10532=1\times 2+3^4\times (5+6+7\times (8+9)).$
\item [] $10533=1^2\times 3+(4+5+6)\times 78\times 9.$
\item [] $10534=1^2+3+(4+5+6)\times 78\times 9.$
\item [] $10535=1\times 2+3+(4+5+6)\times 78\times 9.$
\item [] $10536=12\times (34\times 5+6+78\times 9).$
\item [] $10537=(12+34\times 5+6)\times 7\times 8+9.$
\item [] $10538=1\times 2^3+(4+5+6)\times 78\times 9.$
\item [] $10539=1\times 234\times (5\times 6+7+8)+9.$
\item [] $10540=1+234\times (5\times 6+7+8)+9.$
\item [] $10541=(12^3+4+5)\times 6+7\times (8+9).$
\item [] $10542=12+3^4\times (5+6+7\times (8+9)).$
\item [] $\mathit{10543=12+34\times 5\times (6+7\times 8)-9.}$
\item [] $\mathit{10544=-1+2+3+4\times 5\times (67\times 8-9).}$
\item [] $10545=12+3+(4+5+6)\times 78\times 9.$
\item [] $\mathit{10546=-1+2\times 3^4\times 5\times (6+7)+8+9.}$
\item [] $10547=1\times 2\times 3^4\times 5\times (6+7)+8+9.$
\item [] $10548=12\times (34+56+789).$
\item [] $10549=1^2\times 34\times 5\times (6+7\times 8)+9.$
\item [] $10550=1^2+34\times 5\times (6+7\times 8)+9.$
\item [] $10551=1\times 2+34\times 5\times (6+7\times 8)+9.$
\item [] $10552=12^3\times 4+56\times (7\times 8+9).$
\item [] $10553=1\times 23\times 456+7\times 8+9.$
\item [] $10554=1+23\times 456+7\times 8+9.$
\item [] $10555=1\times 2\times (3^4\times 5\times (6+7)+8)+9.$
\item [] $10556=1+2\times (3^4\times 5\times (6+7)+8)+9.$
\item [] $10557=(1+2)^3+(4+5+6)\times 78\times 9.$
\item [] $10558=1+23\times (4+5+(6\times 7+8)\times 9).$
\item [] $10559=1\times 2+(3+(4+5+6)\times 78)\times 9.$
\item [] $10560=(1+2+3+4)\times (5+6)\times (7+89).$
\item[]$\mbox{Decreasing order}$
\item [] $10491=(9+8\times 76)\times (5+4\times 3)+2\times 1.$
\item [] $10492=(9+8\times 76)\times (5+4\times 3)+2+1.$
\item [] $10493=(9+8)\times 7+6\times (54\times 32+1).$
\item [] $10494=9\times (8\times 7+6\times 5\times (4+32+1)).$
\item [] $10495=9\times ((8+7+6)\times 54+32)+1.$
\item [] $10496=(9+8+(7+6)\times 5)\times 4\times 32\times 1.$
\item [] $10497=(9+8+(7+6)\times 5)\times 4^3\times 2+1.$
\item [] $\mathit{10498=1+23\times 456-(7-8)\times 9.}$
\item [] $\mathit{10499=9876+5^4-3+2-1.}$
\item [] $10500=(9+87\times 6\times 5)\times 4+3+21.$
\item [] $10501=9+87+(6\times (5+4\times 3))^2+1.$
\item [] $10502=9+8765+(4\times 3)^{(2+1)}.$
\item [] $10503=9\times (8+7)+6\times 54\times 32\times 1.$
\item [] $10504=9\times (8+7)+6\times 54\times 32+1.$
\item [] $10505=9+8+76\times (5+4^3)\times 2\times 1.$
\item [] $10506=9876+5^4+3+2\times 1.$
\item [] $10507=9876+5^4+3+2+1.$
\item [] $10508=9876+5^4+3\times 2+1.$
\item [] $10509=(9+87\times 6\times 5)\times 4+32+1.$
\item [] $10510=9876+5^4+3^2\times 1.$
\item [] $10511=9876+5^4+3^2+1.$
\item [] $10512=9+87\times 6\times 5\times 4+3\times 21.$
\item [] $10513=9\times 8\times (7+6+5+4\times 32)+1.$
\item [] $10514=(98\times 7+65)\times (4+3^2+1).$
\item [] $10515=(98\times 7+65)\times (4+3)\times 2+1.$
\item [] $10516=9876+5\times 4\times 32\times 1.$
\item [] $10517=9876+5\times 4\times 32+1.$
\item [] $10518=9+(87+6)\times (5+4\times 3^{(2+1)}).$
\item [] $10519=98\times (7+6)+5\times 43^2\times 1.$
\item [] $10520=98\times (7+6)+5\times 43^2+1.$
\item [] $10521=(9\times 8\times 7+6)\times 5\times 4+321.$
\item [] $10522=9+(8+(7+6)\times 5)\times (4\times 3)^2+1.$
\item [] $10523=9+876\times (5+4+3)+2\times 1.$
\item [] $10524=9+876\times (5+4+3)+2+1.$
\item [] $10525=9876+5^4+3+21.$
\item [] $10526=(9+8)\times (76+543)+2+1.$
\item [] $10527=9\times (8+7+6)\times 54+321.$
\item [] $10528=9876+5^4+3^{(2+1)}.$
\item [] $10529=(9+(8+7+65)\times 4)\times 32+1.$
\item [] $10530=9\times (8+76+543\times 2\times 1).$
\item [] $10531=9\times (8+7+6+5)\times (43+2)+1.$
\item [] $10532=9\times (8\times 7\times 6+54)\times 3+2\times 1.$
\item [] $10533=9876+5^4+32\times 1.$
\item [] $10534=9876+5^4+32+1.$
\item [] $\mathit{10535=9+87\times (6+5)(4\times 3)^2-1.}$
\item [] $10536=98\times 7\times 6+5\times 4\times 321.$
\item [] $10537=(9+8)\times 76+5\times 43^2\times 1.$
\item [] $10538=(9+8)\times 76+5\times 43^2+1.$
\item [] $10539=(9+87\times 6\times 5)\times 4+3\times 21.$
\item [] $10540=(9\times 8+7\times 65)\times 4\times (3+2)\times 1.$
\item [] $10541=(9\times 8+7\times 65)\times 4\times (3+2)+1.$
\item [] $10542=987+65\times (4+3)\times 21.$
\item [] $\mathit{10543=-98\times 7+6\times 5^4\times 3-21.}$
\item [] $10544=(9+8)\times (76+543)+21.$
\item [] $10545=9+8\times (7\times (6\times 5+4^3)\times 2+1).$
\item [] $10546=(9+8\times 7\times (6\times 5+4^3))\times 2\times 1.$
\item [] $10547=98+(76+5)\times 43\times (2+1).$
\item [] $10548=9\times (8\times 76+543+21).$
\item [] $10549=9+(87\times 6+5)\times 4\times (3+2)\times 1.$
\item [] $10550=9+(87\times 6+5)\times 4\times (3+2)+1.$
\item [] $10551=9\times (8\times 7\times 6+54)\times 3+21.$
\item [] $\mathit{10552=-9+(8+7)\times (6+5)\times 4^3+2-1.}$
\item [] $10553=9+8\times (7+6\times 5\times 43+21).$
\item [] $\mathit{10554=-9+(8+7)\times (6+5)\times 4^3+2+1.}$
\item [] $\mathit{10555=(9+8)\times 76\times 5+4^(3\times 2)-1.}$
\item [] $10556=(9+8)\times 76\times 5+4^{(3\times 2)}\times 1.$
\item [] $10557=9\times 87+6\times 543\times (2+1).$
\item [] $10558=(9+8)\times (76+543+2)+1.$
\item [] $\mathit{10559=9\times 8+76\times (5+4^3)\times 2-1.}$
\item [] $10560=(9+87)\times (65+43+2)\times 1.$
\item[]$\mbox{Increasing order}$
\item [] $10561=12+34\times 5\times (6+7\times 8)+9.$
\item [] $\mathit{10562=-1+23+4\times 5\times (67\times 8-9).}$
\item [] $10563=123+4\times 5\times 6\times (78+9).$
\item [] $10564=1\times 2\times (3^4\times 5\times (6+7)+8+9).$
\item [] $10565=1+2\times (3^4\times 5\times (6+7)+8+9).$
\item [] $10566=12\times 3+(4+5+6)\times 78\times 9.$
\item [] $10567=1\times 23\times 456+7+8\times 9.$
\item [] $10568=1+23\times 456+7+8\times 9.$
\item [] $10569=12+(3+(4+5+6)\times 78)\times 9.$
\item [] $\mathit{10570=1\times 23\times 456-7+89.}$
\item [] $\mathit{10571=1+23\times 456-7+89.}$
\item [] $10572=12\times (3^4+5+6+789).$
\item [] $\mathit{10573=12\times 3^4\times (5+6)-7\times (8+9).}$
\item [] $10574=1\times 2\times (3+4\times (5+6)\times 7)\times (8+9).$
\item [] $10575=1\times 23\times 456+78+9.$
\item [] $10576=1+23\times 456+78+9.$
\item [] $\mathit{10577=1\times 23\times (456+7)-8\times 9.}$
\item [] $10578=123\times (4+5\times (6+7)+8+9).$
\item [] $10579=(12\times 3\times 4+5)\times (6+7\times 8+9).$
\item [] $10580=1\times (2+3)\times 4\times (5\times (6+7)\times 8+9).$
\item [] $10581=12\times (345+67\times 8)+9.$
\item [] $\mathit{10582=-1+(23+4)\times 56\times 7+8-9.}$
\item [] $\mathit{10583=-1+23\times 456+7+89.}$
\item [] $10584=1\times 23\times 456+7+89.$
\item [] $10585=1+23\times 456+7+89.$
\item [] $10586=(1\times 2+3\times 4\times (5+6))\times (7+8\times 9).$
\item [] $10587=1+(2+3\times 4\times (5+6))\times (7+8\times 9).$
\item [] $\mathit{10588=-123+4\times 5\times 67\times 8-9.}$
\item [] $\mathit{10589=-1+2\times ((3\times 4+56)\times 78-9).}$
\item [] $\mathit{10590=(1+2)^{(3+4)}\times 5-6\times 7\times 8-9.}$
\item [] $10591=(1+2\times 3+45+67)\times 89.$
\item [] $\mathit{10592=1+(2\times 3^4-5)\times 67+8\times 9.}$
\item [] $10593=(1+2)\times (3+4+56)\times 7\times 8+9.$
\item [] $\mathit{10594=-1\times 2+3\times (4+(56-7)\times 8\times 9).}$
\item [] $10595=(1^2+3\times 4)\times (5+6\times (7+8)\times 9).$
\item [] $10596=12\times (34+56\times (7+8)+9).$
\item [] $\mathit{10597=-1\times 23+4\times 5\times (67-8)\times 9.}$
\item [] $\mathit{10598=-1+2\times (3\times 4+56)\times 78-9.}$
\item [] $\mathit{10599=1\times 2\times (3\times 4+56)\times 78-9.}$
\item [] $10600=(1+2\times 3\times 4)\times (5\times 67+89).$
\item [] $10601=1\times (23+4)\times 56\times 7+8+9.$
\item [] $10602=1+(23+4)\times 56\times 7+8+9.$
\item [] $10603=1+2\times 3^4\times 5\times (6+7)+8\times 9.$
\item [] $\mathit{10604=(-12+34)\times (5+6\times 78+9).}$
\item [] $10605=(12+3)\times (4\times 5+678+9).$
\item [] $\mathit{10606=-123+4\times 5\times 67\times 8+9.}$
\item [] $10607=1\times 23\times 456+7\times (8+9).$
\item [] $10608=1+23\times 456+7\times (8+9).$
\item [] $10609=1^2+(3+45)\times (6+7)\times (8+9).$
\item [] $10610=1\times 2+(3+45)\times (6+7)\times (8+9).$
\item [] $10611=12\times 3^4+567\times (8+9).$
\item [] $10612=(1\times 2+3\times 4)\times (56+78\times 9).$
\item [] $10613=1+2\times (3+4)\times (56+78\times 9).$
\item [] $10614=1\times (23\times 4+5\times 6)\times (78+9).$
\item [] $10615=1+(23\times 4+5\times 6)\times (78+9).$
\item [] $\mathit{10616=-1+2\times (3\times 4+56)\times 78+9.}$
\item [] $10617=1\times 2\times (3\times 4+56)\times 78+9.$
\item [] $10618=1+2\times (3\times 4+56)\times 78+9.$
\item [] $10619=1\times 2\times 3^4\times 5\times (6+7)+89.$
\item [] $10620=(1\times 2\times 3+4+5)\times (6+78\times 9).$
\item [] $10621=1+(2\times 3+4+5)\times (6+78\times 9).$
\item [] $\mathit{10622=-1+23\times 456+(7+8)\times 9.}$
\item [] $10623=1\times 23\times 456+(7+8)\times 9.$
\item [] $10624=1+23\times 456+(7+8)\times 9.$
\item [] $10625=(1+2\times 3\times 4)\times 5\times (6+7+8\times 9).$
\item [] $10626=(1+2\times 34)\times (5\times (6+7)+89).$
\item [] $10627=1+2\times ((3\times 4+56)\times 78+9).$
\item [] $\mathit{10628=(12^3+45)\times 6+7-8-9.}$
\item [] $10629=(1+2)\times (3456+78+9).$
\item [] $10630=1+(23\times (4+5+6\times 7)+8)\times 9.$
\item[]$\mbox{Decreasing order}$
\item [] $10561=(9+87)\times 65+4321.$
\item [] $10562=(9\times (8+7)+6\times 5)\times 4^3+2\times 1.$
\item [] $10563=(9+8+7\times (65+4)+3)\times 21.$
\item [] $10564=9876+5^4+3\times 21.$
\item [] $10565=(98+(76+5)\times 4^3)\times 2+1.$
\item [] $10566=9\times (87\times 6+5^4+3^{(2+1)}).$
\item [] $\mathit{10567=-98\times 7+6\times 5^4\times 3+2+1.}$
\item [] $\mathit{10568=-9\times 8+76\times 5\times 4\times (3\times 2+1).}$
\item [] $10569=9+(8+7+65)\times 4\times (32+1).$
\item [] $10570=(9+876+5^4)\times (3\times 2+1).$
\item [] $10571=9+(8+7)\times (6+5)\times 4^3+2\times 1.$
\item [] $10572=9+(8+7)\times (6+5)\times 4^3+2+1.$
\item [] $10573=98\times 76+5^4\times (3+2)\times 1.$
\item [] $10574=98\times 76+5^4\times (3+2)+1.$
\item [] $10575=(9\times (8+7)+6)\times 5\times (4\times 3+2+1).$
\item [] $10577=9+8\times (7+(654+3)\times 2\times 1).$
\item [] $10577=9+8\times (7+6\times 5+4\times 321).$
\item [] $10578=(9+8+(7+6)\times 5)\times 43\times (2+1).$
\item [] $10579=(9+8\times 7+6)\times (5+(4\times 3)^2\times 1).$
\item [] $10580=(9+8\times (7+6)\times 5)\times 4\times (3+2)\times 1.$
\item [] $10581=9+8+76\times ((5+4^3)\times 2+1).$
\item [] $\mathit{10582=9\times 8\times 7\times (6+5-4)\times 3-2\times 1.}$
\item [] $10583=9+(876+5)\times 4\times 3+2\times 1.$
\item [] $10584=(9\times 8+76+5\times 4)\times 3\times 21.$
\item [] $10585=9\times 8\times 7\times (6+5+4+3\times 2)+1.$
\item [] $10586=98+76\times (5+4^3)\times 2\times 1.$
\item [] $10587=98+76\times (5+4^3)\times 2+1.$
\item [] $\mathit{10588=(9+87\times 6)\times 5\times 4-32\times 1.}$
\item [] $\mathit{10589=(9+87\times 6)\times 5\times 4-32+1.}$
\item [] $10590=9+(8+7)\times (6+5)\times 4^3+21.$
\item [] $10591=(9+8)\times 7\times (65+4\times 3\times 2\times 1).$
\item [] $10592=9\times 8\times 76+5\times 4^(3+2)\times 1.$
\item [] $10593=9+8\times 7\times (6+54\times 3+21).$
\item [] $10594=9+(8+(7+6)\times 5)\times ((4\times 3)^2+1).$
\item [] $10595=(9\times (8+7)\times 6+5)\times (4+3^2)\times 1.$
\item [] $10596=9876+5\times (4\times 3)^2\times 1.$
\item [] $10597=9876+5\times (4\times 3)^2+1.$
\item [] $10598=98+7\times 6\times 5\times ((4+3)^2+1).$
\item [] $10599=9+(8+7)\times ((6+5)\times 4^3+2\times 1).$
\item [] $10600=((9+87)\times (6+5)+4)\times (3^2+1).$
\item [] $10601=9+8+7\times (65+4+3)\times 21.$
\item [] $10602=9+(876+5)\times 4\times 3+21.$
\item [] $10603=(9+8+76)\times (54+3)\times 2+1.$
\item [] $\mathit{10604=-9+8765+43^2-1.}$
\item [] $10605=9\times (8+76\times 5+4)\times 3+21.$
\item [] $10606=(98+7)\times (65+4+32)+1.$
\item [] $10607=9876+(5+4)^3+2\times 1.$
\item [] $10608=9876+(5+4)^3+2+1.$
\item [] $10609=9+8+(7+6\times 54)\times 32\times 1.$
\item [] $10610=9+8+(7+6\times 54)\times 32+1.$
\item [] $10611=9876+5\times (4+3)\times 21.$
\item [] $10612=987+6^5+43^2\times 1.$
\item [] $10613=987+6^5+43^2+1.$
\item [] $10614=(98+76)\times (54+3\times 2+1).$
\item [] $\mathit{10615=(9+87\times 6)\times 5\times 4-3-2\times 1.}$
\item [] $\mathit{10616=(9+87\times 6)\times 5\times 4-3-2+1.}$
\item [] $\mathit{10617=(9+87\times 6)\times 5\times 4-3\times (2-1).}$
\item [] $10618=9+(8\times 7+(6+5)\times 4+3)^2\times 1.$
\item [] $10619=98\times 7+(6+5)\times 43\times 21.$
\item [] $10620=9\times (87+6+543\times 2+1).$
\item [] $10621=(9+8+7\times 6)\times 5\times 4\times 3^2+1.$
\item [] $10622=(9+876)\times (5+4+3)+2\times 1.$
\item [] $10623=9+8765+43^2\times 1.$
\item [] $10624=9+8765+43^2+1.$
\item [] $10625=9\times (8+7)\times 65+43^2+1.$
\item [] $10626=9876+(5+4)^3+21.$
\item [] $10627=(9+87+65)\times (4^3+2)+1.$
\item [] $\mathit{10628=(9+87\times 6)\times 5\times 4+3^2-1.}$
\item [] $10629=(9+87\times 6)\times 5\times 4+3^2\times 1.$
\item [] $10630=(9+87\times 6)\times 5\times 4+3^2+1.$
\item[]$\mbox{Increasing order}$
\item [] $\mathit{10631=1\times 2^3\times 4\times 5\times 67-89.}$
\item [] $10632=1^2\times 3\times (4+5\times (6+78\times 9)).$
\item [] $10633=1^2+3\times (4+5\times (6+78\times 9)).$
\item [] $10634=1\times 2+3\times (4+5\times (6+78\times 9)).$
\item [] $10635=1+2+3\times (4+5\times (6+78\times 9)).$
\item [] $\mathit{10636=(-1+2^3\times 4\times 5)\times 67-8-9.}$
\item [] $\mathit{10637=-9\times 8+765\times (4\times 3+2)-1.}$
\item [] $10638=(1+2)\times (3\times 4^5+6\times (7+8\times 9)).$
\item [] $\mathit{10639=(-12+3+4\times 5\times 67)\times 8-9.}$
\item [] $10640=(12+3+4)\times (56+7\times 8\times 9).$
\item [] $10641=(1+2)\times (3+4+5\times (6+78\times 9)).$
\item [] $10642=(1\times 234+56\times 7)\times (8+9).$
\item [] $10643=1+(234+56\times 7)\times (8+9).$
\item [] $10644=12+3\times (4+5\times (6+78\times 9)).$
\item [] $\mathit{10645=12^3\times 4-5+6\times 7\times 89.}$
\item [] $\mathit{10646=(12^3+45)\times 6+7-8+9.}$
\item [] $10647=(1^2+3\times 4)\times (5\times 6+789).$
\item [] $\mathit{10648=1\times 2^3\times 4\times 5\times 67-8\times 9.}$
\item [] $10649=1\times 23\times (4+5\times 6\times (7+8)+9).$
\item [] $10650=1\times 2\times (3^4\times 56+789).$
\item [] $10651=1+2\times (3^4\times 56+789).$
\item [] $\mathit{10652=1-2-3+4\times (5\times 6+7)\times 8\times 9.}$
\item [] $10653=123+(4+5+6)\times 78\times 9.$
\item [] $10654=1\times 2\times (3+4)\times (5+(6+78)\times 9).$
\item [] $10655=12^3\times 4+5+6\times 7\times 89.$
\item [] $10656=(1+2)\times (3456+7+89).$
\item [] $10657=1+(23+4)\times 56\times 7+8\times 9.$
\item [] $10658=1\times 2+(3^4+5\times 6)\times (7+89).$
\item [] $10659=(1+234+56\times 7)\times (8+9).$
\item [] $10660=1^2+3+4\times (5\times 6+7)\times 8\times 9.$
\item [] $10661=1\times 2+3+4\times (5\times 6+7)\times 8\times 9.$
\item [] $10662=(12^3+45)\times 6+7+8+9.$
\item [] $10663=1+2\times 3+4\times (5\times 6+7)\times 8\times 9.$
\item [] $10664=1\times 2^3+4\times (5\times 6+7)\times 8\times 9.$
\item [] $10665=12\times (3\times 45\times 6+78)+9.$
\item [] $10666=1\times 23\times (456+7)+8+9.$
\item [] $10667=1+23\times (456+7)+8+9.$
\item [] $10668=12\times (3\times 45\times 6+7+8\times 9).$
\item [] $\mathit{10669=(1^2\times 3+4\times 56)\times (7\times 8-9).}$
\item [] $\mathit{10670=(-1+23)\times (4+56\times 7+89).}$
\item [] $10671=12+3+4\times (5\times 6+7)\times 8\times 9.$
\item [] $10672=(1^2+3)\times (4+(5\times 6+7)\times 8\times 9).$
\item [] $10673=1\times (23+4)\times 56\times 7+89.$
\item [] $10674=1+(23+4)\times 56\times 7+89.$
\item [] $10675=(1^2+34)\times ((5\times 6+7)\times 8+9).$
\item [] $\mathit{10676=1-2+(3^4+56)\times 78-9.}$
\item [] $\mathit{10677=(1^2\times 3^4+56)\times 78-9.}$
\item [] $\mathit{10678=-1+23+4\times (5\times 6+7)\times 8\times 9.}$
\item [] $10679=1\times 23+4\times (5\times 6+7)\times 8\times 9.$
\item [] $10680=12\times (345+67\times 8+9).$
\item [] $10681=1\times 23\times (45+6+7)\times 8+9.$
\item [] $10682=1+23\times (45+6+7)\times 8+9.$
\item [] $10683=(1+2)^3+4\times (5\times 6+7)\times 8\times 9.$
\item [] $10684=(1+2\times 3^4)\times 5\times (6+7)+89.$
\item [] $10685=1\times 2+(3+4\times (5\times 6+7)\times 8)\times 9.$
\item [] $10686=1+2+(3+4\times (5\times 6+7)\times 8)\times 9.$
\item [] $\mathit{10687=-1-23+4\times 5\times 67\times 8-9.}$
\item [] $\mathit{10688=-1\times 23+4\times 5\times 67\times 8-9.}$
\item [] $10689=(1+23\times (45+6+7))\times 8+9.$
\item [] $\mathit{10690=-1-2+(3\times 4+5)\times (6+7\times 89).}$
\item [] $\mathit{10691=-1+(2\times 3)^4\times 5+6\times 78\times 9.}$
\item [] $10692=12^3\times 4+5\times (6+78)\times 9.$
\item [] $10693=1+2\times 3\times (4\times 5\times 6+78)\times 9.$
\item [] $10694=1^2+(3\times 4+5)\times (6+7\times 89).$
\item [] $10695=(1^2\times 3^4+56)\times 78+9.$
\item [] $10696=1^2+(3^4+56)\times 78+9.$
\item [] $10697=1\times 2+(3^4+56)\times 78+9.$
\item [] $10698=1+2+(3^4+56)\times 78+9.$
\item [] $\mathit{10699=12\times 3^4\times (5+6)+7\times (-8+9).}$
\item [] $\mathit{10700=12\times 3^4\times (5+6)+7-8+9.}$
\item[]$\mbox{Decreasing order}$
\item [] $\mathit{10631=(9+8+7)\times (6+5+432)-1.}$
\item [] $10632=(9+8+7)\times (6+5+432)\times 1.$
\item [] $10633=9+8\times (7+654+3)\times 2\times 1.$
\item [] $10634=9+8\times (7+654+3)\times 2+1.$
\item [] $10635=9+(8+7\times 65+43)\times 21.$
\item [] $10636=9\times 8+76\times ((5+4^3)\times 2+1).$
\item [] $\mathit{10637=-9\times 8+765\times (4+3)\times 2-1.}$
\item [] $10638=9\times (87\times 6+5+4^3)\times 2\times 1.$
\item [] $10639=9\times (87\times 6+5+4^3)\times 2+1.$
\item [] $\mathit{10640=(9-8)\times 76\times 5\times (4+3+21).}$
\item [] $10641=(9+876)\times (5+4+3)+21.$
\item [] $10642=(9+8\times (7+654+3))\times 2\times 1.$
\item [] $10643=(9+8\times (7+654+3))\times 2+1.$
\item [] $10644=9\times (87\times 6+5^4)+321.$
\item [] $\mathit{10645=9-8\times 7+6\times 54\times (32+1).}$
\item [] $\mathit{10646=9\times 8\times (7+6\times 5)\times 4-3^2-1.}$
\item [] $10647=(9+8\times 7\times 6+54\times 3)\times 21.$
\item [] $10648=(9\times (8\times 7+6)+5^4)\times 3^2+1.$
\item [] $10649=9+8\times (7+(6+54+3)\times 21).$
\item [] $10650=(9+8\times 7+6)\times (5+(4\times 3)^2+1).$
\item [] $\mathit{10651=9\times 8\times (7+6\times 5)\times 4-3-2\times 1.}$
\item [] $10652=(9+87\times 6)\times 5\times 4+32\times 1.$
\item [] $10653=(9+87\times 6)\times 5\times 4+32+1.$
\item [] $10654=(9\times (8+76)+5)\times (4+3)\times 2\times 1.$
\item [] $10655=(9\times (8+76)+5)\times (4+3)\times 2+1.$
\item [] $10656=(9\times 8)\times (76+5+4+3\times 21).$
\item [] $10657=9+8\times (7\times 6+5+4\times 321).$
\item [] $10658=9\times (8+7\times (6+54\times 3))+2\times 1.$
\item [] $10659=(9+(8\times 7+6)\times 5+4)\times (32+1).$
\item [] $10660=(9+8\times 7)\times (6\times (5+4)\times 3+2\times 1).$
\item [] $10661=9\times 8\times (7+6\times 5)\times 4+3+2\times 1.$
\item [] $10662=9\times 8\times (7+6\times 5)\times 4+3\times 2\times 1.$
\item [] $10663=9\times 8\times (7+6\times 5)\times 4+3\times 2+1.$
\item [] $10664=9\times 8+(7+6\times 54)\times 32\times 1.$
\item [] $10665=9\times 8+(7+6+5\times 4)\times 321.$
\item [] $10666=9\times 8\times (7+6\times 5)\times 4+3^2+1.$
\item [] $10667=9+8+(7+6+5+4)^3+2\times 1.$
\item [] $10668=((9+8+76)\times 5+43)\times 21.$
\item [] $\mathit{10669=(9-8+7-6+5\times 4)^3+21.}$
\item [] $\mathit{10670=(98-76)\times (54\times 3^2-1).}$
\item [] $\mathit{10671=-9+87\times (6\times 5\times 4+3)-21.}$
\item [] $\mathit{10672=(9\times 8\times (7+6\times 5)+4)\times (3+2-1).}$
\item [] $10673=9+8\times (76+(5^4+3)\times 2+1).$
\item [] $10674=9\times (8+7\times 6+543)\times 2\times 1.$
\item [] $10675=9\times (8+7\times 6+543)\times 2+1.$
\item [] $10676=(9+8)\times (7+(65+4)\times 3^2\times 1).$
\item [] $10677=9+((87+6)\times 5+43)\times 21.$
\item [] $\mathit{10678=(9+876+5)\times 4\times 3-2\times 1.}$
\item [] $10679=(9+8+7\times 6)\times (5\times 4\times 3^2+1).$
\item [] $10680=9\times 8\times (7+6\times 5)\times 4+3+21.$
\item [] $10681=(9\times 8\times 7+6\times 5)\times 4\times (3+2)+1.$
\item [] $10682=98\times (7+65+4+32+1).$
\item [] $10683=(9+876+5)\times 4\times 3+2+1.$
\item [] $10684=987\times 6+(5+4^3)^2+1.$
\item [] $10685=(98+76\times (5+4^3))\times 2+1.$
\item [] $10686=9+8+(7+6+5+4)^3+21.$
\item [] $10687=9+(8\times (7\times 6+5^4)+3)\times 2\times 1.$
\item [] $10688=9\times 8\times (7+6\times 5)\times 4+32\times 1.$
\item [] $10689=9\times 8\times (7+6\times 5)\times 4+32+1.$
\item [] $10690=98+(7+6\times 54)\times 32\times 1.$
\item [] $10691=98\times (7+6)\times 5+4321.$
\item [] $10692=9+(8+7\times 6\times 5)\times (4+3)^2+1.$
\item [] $10693=(9+8)\times (7\times (6+5)\times 4+321).$
\item [] $\mathit{10694=-9-8+765\times (4+3)\times 2+1.}$
\item [] $10695=(9+8+76)\times ((54+3)\times 2+1).$
\item [] $10696=(9+8\times (7\times 6+5^4)+3)\times 2\times 1.$
\item [] $10697=98\times 76+(54+3)^2\times 1.$
\item [] $10698=98\times 76+(54+3)^2+1.$
\item [] $10699=((9+(8+7)\times 6)\times 54+3)\times 2+1.$
\item [] $\mathit{10700=9-8+7+6\times 54\times (32+1).}$
\item[]$\mbox{Increasing order}$
\item [] $10701=(1+23\times 4+5\times 6)\times (78+9).$
\item [] $10702=1^2+(3+4\times 5\times 6)\times (78+9).$
\item [] $10703=(12^3+45)\times 6+7\times 8+9.$
\item [] $10704=1+2+(3+4\times 5\times 6)\times (78+9).$
\item [] $10705=12+(3\times 4+5)\times (6+7\times 89).$
\item [] $\mathit{10706=1-2\times 3+4\times 5\times 67\times 8-9.}$
\item [] $10707=12+(3^4+56)\times 78+9.$
\item [] $10708=1\times 2\times (3^4\times 5\times (6+7)+89).$
\item [] $10709=(1+2)^3\times (4+56\times 7)+8+9.$
\item [] $10710=(12+345)\times (6+7+8+9).$
\item [] $10711=1^2+(34+56)\times 7\times (8+9).$
\item [] $10712=1\times 2+(34+56)\times 7\times (8+9).$
\item [] $10713=12\times (34+(5+6)\times 78)+9.$
\item [] $\mathit{10714=-12-3+4\times 5\times 67\times 8+9.}$
\item [] $\mathit{10715=-1+2+3+4\times 5\times 67\times 8-9.}$
\item [] $10716=12\times 3^4\times (5+6)+7+8+9.$
\item [] $10717=(12^3+45)\times 6+7+8\times 9.$
\item [] $10718=123\times (45+6\times 7)+8+9.$
\item [] $10719=(1\times 234\times 5+6+7+8)\times 9.$
\item [] $10720=1+((234\times 5+6+7)+8)\times 9.$
\item [] $10721=1\times 23\times (456+7)+8\times 9.$
\item [] $10722=1+23\times (456+7)+8\times 9.$
\item [] $10723=1^2+3\times (4+5\times 6\times 7\times (8+9)).$
\item [] $10724=1\times 2+3\times (4+5\times 6\times 7\times (8+9)).$
\item [] $10725=(12^3+45)\times 6+78+9.$
\item [] $\mathit{10726=12+3+4\times 5\times 67\times 8-9.}$
\item [] $10727=(1^2+(34+56)\times 7)\times (8+9).$
\item [] $10728=(1^2+3\times 45+6+7)\times 8\times 9.$
\item [] $10729=1^{23}\times 4\times 5\times 67\times 8+9.$
\item [] $10730=1^{23}+4\times 5\times 67\times 8+9.$
\item [] $10731=(1+2)\times (3+4+5\times 6\times 7\times (8+9)).$
\item [] $10732=1^2\times 3+4\times 5\times 67\times 8+9.$
\item [] $10733=1^2+3+4\times 5\times 67\times 8+9.$
\item [] $10734=1\times 2+3+4\times 5\times 67\times 8+9.$
\item [] $10735=1+2+3+4\times 5\times 67\times 8+9.$
\item [] $10736=1+2\times 3+4\times 5\times 67\times 8+9.$
\item [] $10737=1\times 2^3+4\times 5\times 67\times 8+9.$
\item [] $10738=1+2^3+4\times 5\times 67\times 8+9.$
\item [] $10739=1+23\times (456+7)+89.$
\item [] $\mathit{10740=-1\times 2^{(3\times 4)}+5^6-789.}$
\item [] $10741=1\times 23\times ((4+5)\times 6\times 7+89).$
\item [] $10742=1+23\times ((4+5)\times 6\times 7+89).$
\item [] $\mathit{10743=-12+(3+4\times 5)\times 6\times 78-9.}$
\item [] $10744=12+3+4\times 5\times 67\times 8+9.$
\item [] $10745=1+2\times (3\times 4+56)\times (7+8\times 9).$
\item [] $10746=(1+2)\times (3\times 4^5+6+7\times 8\times 9).$
\item [] $\mathit{10747=12\times 3+4\times 5\times 67\times 8-9.}$
\item [] $\mathit{10748=-1\times 2^3+4\times (5\times 67\times 8+9).}$
\item [] $\mathit{10749=1-2^3+4\times (5\times 67\times 8+9).}$
\item [] $10750=1^2\times (3^4+5)\times (6+7\times (8+9)).$
\item [] $10751=1^2+(3^4+5)\times (6+7\times (8+9)).$
\item [] $10752=1\times 23+4\times 5\times 67\times 8+9.$
\item [] $10753=1+23+4\times 5\times 67\times 8+9.$
\item [] $10754=1^2+(3+4\times 5\times 67)\times 8+9.$
\item [] $10755=(12+3)\times (4+5+6+78\times 9).$
\item [] $10756=(1+2)^3+4\times 5\times 67\times 8+9.$
\item [] $10757=12\times 3^4\times (5+6)+7\times 8+9.$
\item [] $10758=1\times 2\times (3+4\times 56\times (7+8+9)).$
\item [] $10759=1+2\times (3+4\times 56\times (7+8+9)).$
\item [] $10760=1^2+3+4\times (5\times 67\times 8+9).$
\item [] $10761=(1\times 23+4+5)\times 6\times 7\times 8+9.$
\item [] $10762=1+(23+4+5)\times 6\times 7\times 8+9.$
\item [] $10763=1+2\times 3+4\times (5\times 67\times 8+9).$
\item [] $10764=12\times (3\times 45\times 6+78+9).$
\item [] $10765=12\times 3+4\times 5\times 67\times 8+9.$
\item [] $10766=1\times 2\times (3+4+56\times (7+89)).$
\item [] $10767=1+2\times (3+4+56\times (7+89)).$
\item [] $\mathit{10768=(1+2^3)^4-5+6\times 78\times 9.}$
\item [] $10769=(1+2^3+45+67)\times 89.$
\item [] $10770=1+(2\times 3\times (4+5)+67)\times 89.$
\item[]$\mbox{Decreasing order}$
\item [] $10701=(9+876+5)\times 4\times 3+21.$
\item [] $10702=9\times (8\times (7+6\times 5)\times 4+3+2)+1.$
\item [] $10703=98\times (7+6\times (5+4\times 3))+21.$
\item [] $10704=9+(87+6)\times ((54+3)\times 2+1).$
\item [] $10705=9+(8+76\times 5\times 4)\times (3\times 2+1).$
\item [] $\mathit{10706=98+(76+(5+4)\times 3)^2-1.}$
\item [] $10707=98+(76+(5+4)\times 3)^2\times 1.$
\item [] $10708=98+(76+(5+4)\times 3)^2+1.$
\item [] $10709=9+8+(76+5)\times 4\times (32+1).$
\item [] $10710=(98+7)\times (65+4+32+1).$
\item [] $10711=(9+8)\times 7\times 6\times (5+4+3\times 2)+1.$
\item [] $10712=9\times 8+76\times 5\times (4+3+21).$
\item [] $10713=9+87\times (6\times 5\times 4+3)+2+1.$
\item [] $\mathit{10714=987\times (6+5)-(4\times 3)^2+1.}$
\item [] $\mathit{10715=-9+8+76\times (54\times 3-21).}$
\item [] $10716=9\times 87+(6+5)\times 43\times 21.$
\item [] $\mathit{10717=(9\times (8+7)+6)\times (-5+43)\times 2+1.}$
\item [] $\mathit{10718=(9-8\times 7)\times 6\times (5-43)+2\times 1.}$
\item [] $10719=9\times 8\times (7+6\times 5)\times 4+3\times 21.$
\item [] $10720=(9+87\times 6+5)\times 4\times (3+2)\times 1.$
\item [] $10721=(9+87\times 6+5)\times 4\times (3+2)+1.$
\item [] $10722=(987+6\times (5+4)^3)\times 2\times 1.$
\item [] $10723=(987+6\times (5+4)^3)\times 2+1.$
\item [] $10724=((9+8)\times 7\times 6\times 5+4)\times 3+2\times 1.$
\item [] $10725=(9+8\times 7+65\times 4)\times (32+1).$
\item [] $10726=(9+8\times (7+65\times 4))\times (3+2)+1.$
\item [] $10727=9+8+765\times (4\times 3+2)\times 1.$
\item [] $10728=9+8+765\times (4+3)\times 2+1.$
\item [] $10729=(9+8\times (7\times 6+5^4+3)\times 2\times 1).$
\item [] $10730=9+8\times (7\times 6+5^4+3)\times 2+1.$
\item [] $10731=9+87\times (6\times 5\times 4+3)+21.$
\item [] $\mathit{10732=9-8+7\times (6\times 5+43)\times 21.}$
\item [] $\mathit{10733=9+8+76\times (54\times 3-21).}$
\item [] $10734=(9\times 87+65\times 43)\times (2+1).$
\item [] $10735=9+(8+765\times (4+3))\times 2\times 1.$
\item [] $10736=9+(8+765\times (4+3))\times 2+1.$
\item [] $10737=987+6\times 5\times (4+321).$
\item [] $10738=98+76\times 5\times (4+3+21).$
\item [] $10739=(9+8\times (7\times 6+5^4+3))\times 2+1.$
\item [] $10740=9+(8+(7+6)\times 5)\times (4+3)\times 21.$
\item [] $10741=9\times 8+(7+6+5+4)^3+21.$
\item [] $\mathit{10742=-9+(8+7\times 6)\times 5\times 43+2-1.}$
\item [] $10743=((9+8)\times 7\times 6\times 5+4)\times 3+21.$
\item [] $10744=(9+8)\times (76\times 5+4\times 3\times 21).$
\item [] $10745=(9+8+765\times (4+3))\times 2+1.$
\item [] $10746=9\times (8\times 7+6\times (54+3))\times (2+1).$
\item [] $\mathit{10747=-9\times 8\times 7+6\times 5^4\times 3+2-1.}$
\item [] $10748=9+8+7\times (6\times 5+43)\times 21.$
\item [] $10749=98+(7+6+5+4)^3+2+1.$
\item [] $10750=((9+8)\times 7+6)\times (54+32\times 1).$
\item [] $10751=((9+8)\times 7+6)\times (54+32)+1.$
\item [] $10752=(98+7\times (6\times 5+4))\times 32\times 1.$
\item [] $10753=(98+7\times (6\times 5+4))\times 32+1.$
\item [] $\mathit{10754=9-8-7\times (6-54)\times 32+1.}$
\item [] $\mathit{10755=-9\times 8\times 7+(6\times 5^4+3)\times (2+1).}$
\item [] $\mathit{10756=9+(8+7\times 6)\times 5\times 43-2-1.}$
\item [] $10757=9+8\times 7+6\times 54\times (32+1).$
\item [] $10758=(98+(7+6)\times 5)\times (4^3+2\times 1).$
\item [] $10759=(9+8+76\times 5\times 4)\times (3\times 2+1).$
\item [] $\mathit{10760=9+(8+7\times 6)\times 5\times 43+2-1.}$
\item [] $10761=987+6\times 543\times (2+1).$
\item [] $10762=9+(8+7\times 6)\times 5\times 43+2+1.$
\item [] $\mathit{10763=(9-8\times 7)\times (-65\times 4+32-1).}$
\item [] $10764=9+(8+7)\times (654+3\times 21).$
\item [] $10765=9\times (8+(7+6+5)\times (4^3+2))+1.$
\item [] $10766=9\times (876+5\times 4^3)+2\times 1.$
\item [] $10767=(98+7+6)\times ((5+43)\times 2+1).$
\item [] $\mathit{10768=((9+8)\times 76+54)\times (3^2-1).}$
\item [] $\mathit{10769=9+8-7\times (6-54)\times 32\times 1.}$
\item [] $10770=9+87\times 6\times 5\times 4+321.$
\item[]$\mbox{Increasing order}$
\item [] $10771=12+3+4\times (5\times 67\times 8+9).$
\item [] $10772=(1^2+3)\times (4+5\times 67\times 8+9).$
\item [] $10773=123\times (45+6\times 7)+8\times 9.$
\item [] $10774=1^2+(3+4\times 5)\times 6\times 78+9.$
\item [] $10775=1\times 2+(3+4\times 5)\times 6\times 78+9.$
\item [] $10776=1+2+(3+4\times 5)\times 6\times 78+9.$
\item [] $10777=(1\times 2\times 3+4\times 5\times 67)\times 8+9.$
\item [] $10778=(1+2^3)^4+5+6\times 78\times 9.$
\item [] $10779=12\times 3^4\times (5+6)+78+9.$
\item [] $10780=1+23+4\times (5\times 67\times 8+9).$
\item [] $10781=(1+2)^3\times (4+56\times 7)+89.$
\item [] $10782=(1\times (2+3)\times 4\times 56+78)\times 9.$
\item [] $10783=(1+2)^3+4\times (5\times 67\times 8+9).$
\item [] $10784=1\times 2^3\times (4+56\times (7+8+9)).$
\item [] $10785=12+(3+4\times 5)\times 6\times 78+9.$
\item [] $\mathit{10786=-1+23\times (4+5\times (6+78+9)).}$
\item [] $10787=1\times 23\times (4+5\times (6+78+9)).$
\item [] $10788=(1\times 2\times 34+56)\times (78+9).$
\item [] $10789=1+(2\times 34+56)\times (78+9).$
\item [] $10790=123\times (45+6\times 7)+89.$
\item [] $10791=(1+(2+3)\times 4\times 56+78)\times 9.$
\item [] $10792=1\times 2^3\times 4\times 5\times 67+8\times 9.$
\item [] $10793=1+2^3\times 4\times 5\times 67+8\times 9.$
\item [] $10794=1+(2^3+4\times 5\times 67)\times 8+9.$
\item [] $10795=(1\times 2\times 34+567)\times (8+9).$
\item [] $10796=1+(2+3)\times (4\times 5\times 6+7)\times (8+9).$
\item [] $10797=12\times (3+(45+67)\times 8)+9.$
\item [] $\mathit{10798=-1\times 2+3\times (456-7\times 8)\times 9.}$
\item [] $\mathit{10799=(1\times 2\times 3\times 4\times 56+7)\times 8-9.}$
\item [] $10800=(1\times 23+4\times 5\times 6+7)\times 8\times 9.$
\item [] $10801=1+2\times (3\times 4+56+7)\times 8\times 9.$
\item [] $\mathit{10802=1\times 2+3\times (456-7\times 8)\times 9.}$
\item [] $\mathit{10803=1+2+3\times (456-7\times 8)\times 9.}$
\item [] $10804=(1+2^3\times 4\times 5)\times 67+8+9.$
\item [] $\mathit{10805=1\times 2+3\times ((456-7)\times 8+9).}$
\item [] $\mathit{10806=(1+2)^{(3+4)}\times 5+6-(7+8)\times 9.}$
\item [] $\mathit{10807=-12^3+4\times 56\times 7\times 8-9.}$
\item [] $\mathit{10808=-1+2^3\times 4\times 5\times 67+89.}$
\item [] $10809=1\times 2^3\times 4\times 5\times 67+89.$
\item [] $10810=1+2^3\times 4\times 5\times 67+89.$
\item [] $10811=12\times 3^4\times (5+6)+7\times (8+9).$
\item [] $10812=(1+2\times 34+567)\times (8+9).$
\item [] $10813=1^2+3\times (4+(56\times 7+8)\times 9).$
\item [] $10814=1\times 2+3\times (4+(56\times 7+8)\times 9).$
\item [] $10815=(1+2)\times (3+4)\times (5+6+7\times 8\times 9).$
\item [] $10816=1+(2+3)\times (4^5+67\times (8+9)).$
\item [] $10817=(1\times 2\times 3\times 4\times 56+7)\times 8+9.$
\item [] $10818=(1+2^3)^4+(5+6\times 78)\times 9.$
\item [] $10819=1+2\times (3^4+5\times (6+7)\times 8)\times 9.$
\item [] $10820=1\times 2\times (34+56\times (7+89)).$
\item [] $10821=1+2\times (34+56\times (7+89)).$
\item [] $10822=(1+2\times 3)\times (4^5+6\times (78+9)).$
\item [] $10823=12\times 3^4\times 5+67\times 89.$
\item [] $10824=12+3\times (4+(56\times 7+8)\times 9).$
\item [] $10825=(1+2\times 3\times 4\times 56+7)\times 8+9.$
\item [] $10826=1+2+(3^4+56)\times (7+8\times 9).$
\item [] $10827=12\times 3^4\times (5+6)+(7+8)\times 9.$
\item [] $\mathit{10828=1-2^{(3\times 4)}+5^6-78\times 9.}$
\item [] $10829=(1^2+34+56)\times 7\times (8+9).$
\item [] $10830=1+(2+3)\times 4\times (5+67\times 8)+9.$
\item [] $\mathit{10831=(12+3+4\times 5\times 67)\times 8-9.}$
\item [] $\mathit{10832=-12+(3+4)^5-67\times 89.}$
\item [] $10833=(1+2)\times (3\times 4^5+67\times 8)+9.$
\item [] $\mathit{10834=123+4\times 5\times 67\times 8-9.}$
\item [] $10835=12+(3^4+56)\times (7+8\times 9).$
\item [] $10836=12\times (345+(6+7\times 8)\times 9).$
\item [] $10837=1+2\times (3+4)\times (5\times 6+7\times 8)\times 9.$
\item [] $\mathit{10838=-1+2^{(3-4+5)}\times 678-9.}$
\item [] $\mathit{10839=(1+2\times 3+4+5)\times 678-9.}$
\item [] $\mathit{10840=(1+2)^{(3+4)}\times 5-(6+7)\times 8+9.}$
\item[]$\mbox{Decreasing order}$
\item [] $10771=9\times 8+7+6\times 54\times (32+1).$
\item [] $\mathit{10772=987\times (6+5)-4^3-21.}$
\item [] $10773=9\times 876+(5+4)\times 321.$
\item [] $\mathit{10774=9876-5+43\times 21.}$
\item [] $10775=9\times (8+7+6)\times (54+3)+2\times 1.$
\item [] $10776=9\times (8+7+6)\times (54+3)+2+1.$
\item [] $\mathit{10777=-9-8+7\times 6\times (5+4\times 3\times 21).}$
\item [] $\mathit{10778=(9+8)\times (7\times 6+5^4-32-1).}$
\item [] $\mathit{10779=-9\times (8+7)+(6\times 5+4)\times 321.}$
\item [] $10780=9+(8+7\times 6)\times 5\times 43+21.$
\item [] $10781=(9+8)\times 76\times 5+4321.$
\item [] $10782=9\times 8+765\times (4+3)\times 2\times 1.$
\item [] $10783=9\times 8+765\times (4+3)\times 2+1.$
\item [] $10784=9876+5+43\times 21.$
\item [] $10785=((9+8)\times 7\times 6+5)\times (4\times 3+2+1).$
\item [] $\mathit{10786=-9-8\times 76+543\times 21.}$
\item [] $10787=(9+87+65)\times (4+3\times 21).$
\item [] $10788=9+87+6\times 54\times (32+1).$
\item [] $10789=(98+76)\times (5\times 4\times 3+2)+1.$
\item [] $10790=98+(76+5)\times 4\times (32+1).$
\item [] $10791=9\times (876+5\times 4^3+2+1).$
\item [] $10792=(9+8\times 7+6)\times (5+(4+3)\times 21).$
\item [] $\mathit{10793=987\times (6+5)-43-21.}$
\item [] $10794=9\times (8+7+6)\times (54+3)+21.$
\item [] $10795=(9+8)\times (7+6\times 5\times 4)\times (3+2)\times 1.$
\item [] $10796=(9+8)\times (7+6\times 5\times 4)\times (3+2)+1.$
\item [] $10797=98+7+6\times 54\times (32+1).$
\item [] $\mathit{10798=987\times (6+5)+4-3\times 21.}$
\item [] $\mathit{10799=9\times 8\times 7\times 6\times 5-4321.}$
\item [] $10800=9\times (876+54\times 3\times 2\times 1).$
\item [] $10801=9\times (8+7\times 6)\times (5+4+3)\times 2+1.$
\item [] $\mathit{10802=(9\times (8-(7-6)^5)^4-3)/2-1.}$
\item [] $10803=9\times 8+7\times (6\times 5+43)\times 21.$
\item [] $\mathit{10804=9-8\times 76+543\times 21.}$
\item [] $\mathit{10805=(9\times 8\times 76-5-4^3)\times 2-1.}$
\item [] $\mathit{10806=(9\times 8\times 76-5-4^3)\times 2\times 1.}$
\item [] $\mathit{10807=98+765\times (4+3)\times 2-1.}$
\item [] $10808=98+765\times (4+3)\times 2\times 1.$
\item [] $10809=98+765\times (4+3)\times 2+1.$
\item [] $10810=9+(8+7)\times 6\times 5\times 4\times 3\times 2+1.$
\item [] $10811=9+8+7\times 6\times (5+4\times 3\times 21).$
\item [] $10812=(9+8)\times ((7+6\times 5\times 4)\times (3+2)+1).$
\item [] $\mathit{10813=987\times (6+5)-43-2+1.}$
\item [] $\mathit{10814=987\times (6+5)-4^3+21.}$
\item [] $10815=(9+8+7\times 65+43)\times 21.$
\item [] $\mathit{10816=987\times (6+5)-43+2\times 1.}$
\item [] $10817=9+8\times 7\times (65+4^3\times 2\times 1).$
\item [] $10818=9+8\times 7\times (65+4^3\times 2)+1.$
\item [] $10819=(9+(8+7)\times 6\times 5\times 4\times 3)\times 2+1.$
\item [] $\mathit{10820=987\times (6+5)-4-32-1.}$
\item [] $\mathit{10821=-98+7\times 65\times 4\times 3\times 2-1.}$
\item [] $10822=9876+5^4+321.$
\item [] $\mathit{10823=-98+7\times 65\times 4\times 3\times 2+1.}$
\item [] $10824=(9+8+(7+6)\times 5)\times 4\times (32+1).$
\item [] $10825=9+8\times (7\times (65+4\times 32)+1).$
\item [] $\mathit{10826=987\times (6+5)-4-3^{(2+1)}.}$
\item [] $10827=9\times (8+7)+6\times 54\times (32+1).$
\item [] $\mathit{10828=987\times (6+5)+4-32-1.}$
\item [] $10829=98+7\times (6\times 5+43)\times 21.$
\item [] $\mathit{10830=987\times (6+5)+4-32+1.}$
\item [] $10831=9+(8+765)\times (4+3)\times 2\times 1.$
\item [] $10832=9+(8+765)\times (4+3)\times 2+1.$
\item [] $10833=9+8\times (7+6\times 5+4)\times (32+1).$
\item [] $\mathit{10834=987\times (6+5)-4\times 3\times 2+1.}$
\item [] $10835=987\times 6+(5+4\times 3)^(2+1).$
\item [] $10836=98\times 7\times 6+5\times 4^3\times 21.$
\item [] $10837=(9+87+6\times 5)\times 43\times 2+1.$
\item [] $\mathit{10838=987\times (6+5)-4\times (3+2)+1.}$
\item [] $10839=((9+8)\times 7\times 6\times 5+43)\times (2+1).$
\item [] $10840=(9+(8+765)\times (4+3))\times 2\times 1.$
\item[]$\mbox{Increasing order}$
\item [] $\mathit{10841=1+2^3\times (4\times (5+6\times 7\times 8)-9).}$
\item [] $10842=1\times 23\times (456+7+8)+9.$
\item [] $10843=1+23\times (456+7+8)+9.$
\item [] $\mathit{10844=1^2\times (3+4)^5-67\times 89.}$
\item [] $10845=(12+3+4)\times 567+8\times 9.$
\item [] $10846=1^2\times 34\times (5\times (6+7\times 8)+9).$
\item [] $10847=1^2+34\times (5\times (6+7\times 8)+9).$
\item [] $10848=1\times (2+3^4+5\times 6)\times (7+89).$
\item [] $10849=(12+3+4\times 5\times 67)\times 8+9.$
\item [] $10850=((1+2)^3+4)\times (5+6\times 7\times 8+9).$
\item [] $10851=(1\times 2+3^4+56)\times 78+9.$
\item [] $10852=123+4\times 5\times 67\times 8+9.$
\item [] $10853=(1+2)\times (3^4+5)\times 6\times 7+8+9.$
\item [] $10854=(12+3)\times (45+678)+9.$
\item [] $10855=1+2\times (3\times 45+6\times 78)\times 9.$
\item [] $10856=1\times 2^3\times (4\times 5\times 67+8+9).$
\item [] $10857=(1+2\times 3+4+5)\times 678+9.$
\item [] $10858=1\times 2\times (3+45+6+7)\times 89.$
\item [] $10859=(1+2^3\times 4\times 5)\times 67+8\times 9.$
\item [] $10860=(12+3)\times (4+5\times 6\times (7+8+9)).$
\item [] $\mathit{10861=-1+(23-4)\times 567+89.}$
\item [] $10862=(12+3+4)\times 567+89.$
\item [] $10863=1^2\times 3^4\times (56+78)+9.$
\item [] $10864=1^2+3^4\times (56+78)+9.$
\item [] $10865=1\times 2+3^4\times (56+78)+9.$
\item [] $10866=1+2+3^4\times (56+78)+9.$
\item [] $\mathit{10867=(1+2)^{(3+4)}\times 5-67+8-9.}$
\item [] $\mathit{10868=(1+2)^{(3+4)}\times 5+67\times (8-9).}$
\item [] $10869=(1+2)\times ((3^4+5)\times 6\times 7+8)+9.$
\item [] $\mathit{10870=12 + (3 \times  45 - 6 - 7) \times  89.}$
\item [] $\mathit{10871=(1+(2\times 3)^4+56+7)\times 8-9.}$
\item [] $10872=(1+2)\times 3456+7\times 8\times 9.$
\item [] $10873=1234+567\times (8+9).$
\item [] $10874=1+2\times (3\times 4\times 56+7)\times 8+9.$
\item [] $10875=12+3^4\times (56+78)+9.$
\item [] $10876=(1+2^3\times 4\times 5)\times 67+89.$
\item [] $\mathit{10877=-1+2\times 3\times 4^5+6\times 789.}$
\item [] $10878=1\times 2\times 3\times 4^5+6\times 789.$
\item [] $10879=123+4\times (5\times 67\times 8+9).$
\item [] $10880=(1^2+3)\times 4\times (5+(67+8)\times 9).$
\item [] $10881=12\times 3\times (4\times 56+78)+9.$
\item [] $10882=1\times 2\times ((3\times 4\times 56+7)\times 8+9).$
\item [] $10883=1+2\times ((3\times 4\times 56+7)\times 8+9).$
\item [] $10884=12\times (3+4\times (5+(6+7)\times (8+9))).$
\item [] $10885=(1^2+3)\times 4^5+6789.$
\item [] $10886=1+2^{(3+4+5)}+6789.$
\item [] $10887=12^3\times 4+5\times (6+789).$
\item [] $\mathit{10888=(1+2)^{(3+4)}\times 5+6\times 7-89.}$
\item [] $10889=(1+(2\times 3)^4+56+7)\times 8+9.$
\item [] $10890=(1\times 2^3)^4+5+6789.$
\item [] $10891=1+2^{(3\times 4)}+5+6789.$
\item [] $\mathit{10892=(1+2)^{(3+4)}\times 5-6\times 7+8-9.}$
\item [] $\mathit{10893=-1-2\times 3+4\times 5\times (67\times 8+9).}$
\item [] $\mathit{10894=(1+2)^{(3+4)}\times 5-6\times 7-8+9.}$
\item [] $\mathit{10895=(1\times 23+4\times 5\times 67)\times 8-9.}$
\item [] $10896=1\times 2\times (3+4\times 56)\times (7+8+9).$
\item [] $10897=1+2\times (3+4\times 56)\times (7+8+9).$
\item [] $10898=1+((2\times 3)^4+5\times (6+7))\times 8+9.$
\item [] $10899=(12\times 3\times 4\times 5+6)\times (7+8)+9.$
\item [] $10900=1^{23}\times 4\times 5\times (67\times 8+9).$
\item [] $10901=1^{23}+4\times 5\times (67\times 8+9).$
\item [] $10902=(123+4+5+6)\times (7+8\times 9).$
\item [] $10903=1^2\times 3+4\times 5\times (67\times 8+9).$
\item [] $10904=1234\times 5+6\times 789.$
\item [] $10905=1\times 2+3+4\times 5\times (67\times 8+9).$
\item [] $10906=1\times 2\times 3+4\times 5\times (67\times 8+9).$
\item [] $10907=1+2\times 3+4\times 5\times (67\times 8+9).$
\item [] $10908=1\times 2^3+4\times 5\times (67\times 8+9).$
\item [] $10909=1+2^3+4\times 5\times (67\times 8+9).$
\item [] $10910=1\times 2+3\times (4+56\times 7+8)\times 9.$
\item[]$\mbox{Decreasing order}$
\item [] $10841=(9+(8+765)\times (4+3))\times 2+1.$
\item [] $\mathit{10842=987\times (6+5)-4\times 3-2-1.}$
\item [] $\mathit{10843=987\times (6+5)+4+3-21.}$
\item [] $\mathit{10844=987\times (6+5)-4\times 3-2+1.}$
\item [] $10845=9\times (8+765+432\times 1).$
\item [] $10846=9\times (8+765+432)+1.$
\item [] $\mathit{10847=-9\times 8+7\times 65\times 4\times 3\times 2-1.}$
\item [] $10848=(9\times 8+7+65\times 4)\times 32\times 1.$
\item [] $10849=(9\times 8+7+65\times 4)\times 32+1.$
\item [] $10850=98\times 76+54\times 3\times 21.$
\item [] $10851=(9+8)\times (7+6+5^4)+3+2\times 1.$
\item [] $10852=(9+8)\times (7+6+5^4)+3\times 2\times 1.$
\item [] $10853=(9+8)\times (7+6+5^4)+3\times 2+1.$
\item [] $10854=9\times (876+5+4+321).$
\item [] $10855=(9\times 87\times 6+(5+4)^3)\times 2+1.$
\item [] $10856=(9+8)\times (7+6+5^4)+3^2+1.$
\item [] $10857=9+8\times (7+65+4\times 321).$
\item [] $10858=9+(8+7+6\times 54)\times 32+1.$
\item [] $10859=9+(8+7\times 6)\times (5\times 43+2\times 1).$
\item [] $10860=9+(8+7\times 6)\times (5\times 43+2)+1.$
\item [] $\mathit{10861=987\times (6+5)+4+3-2-1.}$
\item [] $\mathit{10862=987\times (6+5)+4+3-2\times 1.}$
\item [] $10863=9\times (876+5\times (4^3+2)+1).$
\item [] $\mathit{10864=(9-8+7-6)\times 5432\times 1.}$
\item [] $10865=9+8\times (7+6\times 5\times (43+2)\times 1).$
\item [] $10866=9\times 8+7\times 6\times (5+4\times 3\times 21).$
\item [] $10867=987\times (6+5)+4+3+2+1.$
\item [] $10868=987\times (6+5)+4+3\times 2+1.$
\item [] $\mathit{10869=987\times (6+5)+4+3^2-1.}$
\item [] $10870=987\times (6+5)+4+3^2\times 1.$
\item [] $10871=987\times (6+5)+4+3^2+1.$
\item [] $10872=9\times 8\times 7+6\times 54\times 32\times 1.$
\item [] $10873=9\times 8\times 7+6\times 54\times 32+1.$
\item [] $10874=9\times 8\times ((7+6\times 5)\times 4+3)+2\times 1.$
\item [] $10875=((9+8)\times 7+6)\times (54+32+1).$
\item [] $\mathit{10876=987\times (6+5)+4\times (3+2)-1.}$
\item [] $10877=987\times (6+5)+4\times (3+2)\times 1.$
\item [] $10878=98\times (7\times 6+5+43+21).$
\item [] $10879=(9+8)\times (7+6+5^4)+32+1.$
\item [] $10880=(98+7+65)\times (43+21).$
\item [] $10881=987\times (6+5)+4\times 3\times 2\times 1.$
\item [] $10882=(98+7+65)\times 4^3+2\times 1.$
\item [] $10883=(98+7+65)\times 4^3+2+1.$
\item [] $10884=9876+(5+43)\times 21.$
\item [] $10885=987\times (6+5)+4+3+21.$
\item [] $\mathit{10886=987\times (6+5)-4+32+1.}$
\item [] $\mathit{10887=(-9-8+7\times 65\times 4\times 3)\times 2+1.}$
\item [] $10888=987\times (6+5)+4+3^{(2+1)}.$
\item [] $10889=9+(8+7\times (6+5))\times 4\times 32\times 1.$
\item [] $10890=987\times (6+5)+4\times 3+21.$
\item [] $10891=9\times (8\times 7+6+543)\times 2+1.$
\item [] $10892=98+7\times 6\times (5+4\times 3\times 21).$
\item [] $10893=987\times (6+5)+4+32\times 1.$
\item [] $10894=987\times (6+5)+4+32+1.$
\item [] $\mathit{10895=9876-5+4^(3+2)\times 1.}$
\item [] $\mathit{10896=9876-5+4^(3+2)+1.}$
\item [] $10897=987\times (6+5)+4\times (3^2+1).$
\item [] $10898=(9+(8+7\times (6+5))\times 4^3)\times 2\times 1.$
\item [] $10899=(98+7+6\times (5+4^3))\times 21.$
\item [] $10900=(9+8)\times (7+6+5^4+3)+2+1.$
\item [] $10901=(98+7+65)\times 4^3+21.$
\item [] $10902=987\times (6+5)+43+2\times 1.$
\item [] $10903=987\times (6+5)+43+2+1.$
\item [] $10904=9\times 8\times 76+5432\times 1.$
\item [] $10905=9\times 8\times 76+5432+1.$
\item [] $10906=9876+5+4^(3+2)+1.$
\item [] $10907=987\times (6+5)+(4+3)^2+1.$
\item [] $10908=9\times 87+(6+5+4)^3\times (2+1).$
\item [] $10909=(9+8)\times (7+6+5^4)+3\times 21.$
\item [] $\mathit{10910=(9+8-7)\times (6+543\times 2-1).}$
\item[]$\mbox{Increasing order}$
\item [] $10911=1+2+3\times (4+56\times 7+8)\times 9.$
\item [] $10912=((1+2)^3+4)\times (5\times 67+8+9).$
\item [] $10913=(1\times 23+4\times 5\times 67)\times 8+9.$
\item [] $10914=1\times 2\times 3\times (4^5+6+789).$
\item [] $10915=12+3+4\times 5\times (67\times 8+9).$
\item [] $10916=1\times 2+3\times (4+5\times 6\times 7)\times (8+9).$
\item [] $10917=1+2+3\times (4+5\times 6\times 7)\times (8+9).$
\item [] $\mathit{10918=(1+2)^{(3+4)}\times 5+(6-7)\times 8-9.}$
\item [] $\mathit{10919=(1+2)^{(3+4)}\times 5-6+7-8-9.}$
\item [] $10920=12+3\times (4+56\times 7+8)\times 9.$
\item [] $10921=(12\times 3^4+56\times 7)\times 8+9.$
\item [] $10922=1+2^3\times 4\times (5+6\times 7\times 8)+9.$
\item [] $10923=1\times 23+4\times 5\times (67\times 8+9).$
\item [] $10924=1+23+4\times 5\times (67\times 8+9).$
\item [] $10925=(1+2)\times (3^4+5)\times 6\times 7+89.$
\item [] $10926=(12^3+4+5)\times 6+7\times 8\times 9.$
\item [] $10927=(1+2)^3+4\times 5\times (67\times 8+9).$
\item [] $\mathit{10928=-1+(2+(3+4\times 5)\times 6)\times 78+9.}$
\item [] $10929=12\times 345+6789.$
\item [] $10930=(1+2\times 3^4)\times (5+6+7\times 8)+9.$
\item [] $10931=(123+4)\times (5\times 6+7\times 8)+9.$
\item [] $10932=12\times ((3+45)\times 6+7\times 89).$
\item [] $10933=1^2+3\times (4+56\times (7\times 8+9)).$
\item [] $10934=(12+3+4)\times (567+8)+9.$
\item [] $10935=1^2\times 3^4\times (56+7+8\times 9).$
\item [] $10936=12\times 3+4\times 5\times (67\times 8+9).$
\item [] $10937=1\times 2+3^4\times (56+7+8\times 9).$
\item [] $10938=1+2+3^4\times (56+7+8\times 9).$
\item [] $\mathit{10939=1-2\times 3+456\times (7+8+9).}$
\item [] $\mathit{10940=1-2-3+456\times (7+8+9).}$
\item [] $10941=(1+2)\times (3+4+56\times (7\times 8+9)).$
\item [] $\mathit{10942=-1+2-3+456\times (7+8+9).}$
\item [] $\mathit{10943=1\times 2-3+456\times (7+8+9).}$
\item [] $10944=1^{23}\times 456\times (7+8+9).$
\item [] $10945=1^{23}+456\times (7+8+9).$
\item [] $\mathit{10946=1-2+3+456\times (7+8+9).}$
\item [] $10947=12+3^4\times (56+7+8\times 9).$
\item [] $10948=1\times 2\times 34\times (5+67+89).$
\item [] $10949=1+2\times 34\times (5+67+89).$
\item [] $10950=1+2+3+456\times (7+8+9).$
\item [] $10951=1+2\times 3+456\times (7+8+9).$
\item [] $10952=1\times 2^3+456\times (7+8+9).$
\item [] $10953=1+2^3+456\times (7+8+9).$
\item [] $10954=1+((2\times 3)^4+5+67)\times 8+9.$
\item [] $\mathit{10955=-1+(23-45)\times (6-7\times 8\times 9).}$
\item [] $10956=123\times (4+(5+6)\times 7+8)+9.$
\item [] $\mathit{10957=(1+2)^{(3+4)}\times 5-67+89.}$
\item [] $10958= \mbox{still not available}.$
\item [] $10959=12+3+456\times (7+8+9).$
\item [] $10960=12+(3^4+5+6)\times 7\times (8+9).$
\item [] $10961=(1+2+34)\times (5\times 6+7)\times 8+9.$
\item [] $10962=12\times 3^4\times 5+678\times 9.$
\item [] $10963=1+2\times (34+567+8)\times 9.$
\item [] $10964=(1^2+3)\times (4\times (5+678)+9).$
\item [] $10965=(1+(2+34+56)\times 7)\times (8+9).$
\item [] $\mathit{10966=-1+23+456\times (7+8+9).}$
\item [] $10967=1\times 23+456\times (7+8+9).$
\item [] $10968=1\times 2\times 3\times 4^5+67\times 8\times 9.$
\item [] $10969=1+2\times 3\times 4^5+67\times 8\times 9.$
\item [] $\mathit{10970=1-2+(3+4\times 5)\times (6\times 78+9).}$
\item [] $10971=(1+2)^3+456\times (7+8+9).$
\item [] $10972=1^2+(3+4\times 5)\times (6\times 78+9).$
\item [] $10973=1\times 2+(3+4\times 5)\times (6\times 78+9).$
\item [] $10974=(12+3^4)\times (5+(6+7)\times 8+9).$
\item [] $\mathit{10975=(1+2\times 3)\times 4\times 56\times 7+8-9.}$
\item [] $\mathit{10976=(1+2)^{(3+4)}\times 5-6+7\times 8-9.}$
\item [] $\mathit{10977=(1+2\times 3)\times 4\times 56\times 7-8+9.}$
\item [] $\mathit{10978=(1+2)^{(3+4)}\times 5+6\times 7-8+9.}$
\item [] $10979=((1^2+3)^4+5)\times 6\times 7+8+9.$
\item [] $10980=12\times 3+456\times (7+8+9).$
\item[]$\mbox{Decreasing order}$
\item [] $10911=(98\times (7+6\times 5)+4)\times 3+21.$
\item [] $10912=(9+8+(76+5)\times 4)\times 32\times 1.$
\item [] $10913=(9+8+(76+5)\times 4)\times 32+1.$
\item [] $10914=(9+8)\times (7\times 6\times 5+432\times 1).$
\item [] $10915=(9+8+7\times 6)\times 5\times (4+32+1).$
\item [] $10916=(9+8)\times (7\times 6\times 5+4)\times 3+2\times 1.$
\item [] $10917=(9+8)\times (7\times 6\times 5+4)\times 3+2+1.$
\item [] $10918=(9+8)\times ((7+6+5^4)+3)+21.$
\item [] $\mathit{10919=-9+8+7\times 65\times 4\times 3\times 2\times 1.}$
\item [] $10920=987+(6+5)\times 43\times 21.$
\item [] $10921=987\times (6+5)+43+21.$
\item [] $10922=(9+8\times 7)\times (6+54\times 3)+2\times 1.$
\item [] $10923=(9+8+7)\times 65\times (4+3)+2+1.$
\item [] $10924=987\times (6+5)+4+3\times 21.$
\item [] $\mathit{10925=(98-7)\times 6\times 5\times 4+3+2\times 1.}$
\item [] $10926=9\times (8+(7+6+5)\times (4+3\times 21)).$
\item [] $10927=(9+(8+76\times 5)\times 4)\times (3\times 2+1).$
\item [] $\mathit{10928=-9-8+76\times (5+4+3)^2+1.}$
\item [] $10929=9+8\times (76+5+4\times 321).$
\item [] $10930=(9\times (8\times 7+65)+4)\times (3^2+1).$
\item [] $10931=(9+8)\times (7\times 6\times 5+432+1).$
\item [] $10932=(9+8)\times (7+6+5^4+3+2)+1.$
\item [] $\mathit{10933=(9+8\times 7\times 65-4)\times 3-2\times 1.}$
\item [] $\mathit{10934=(9+8\times 7\times 65-4)\times 3-2+1.}$
\item [] $10935=9\times (8+7+6\times 5\times 4)\times 3^2\times 1.$
\item [] $10936=9\times (8+7+6\times 5\times 4)\times 3^2+1.$
\item [] $10937=9+8+7\times 65\times 4\times 3\times 2\times 1.$
\item [] $10938=9+8+7\times 65\times 4\times 3\times 2+1.$
\item [] $\mathit{10939=(9\times 8\times 76+5-4-3)\times 2-1.}$
\item [] $10940=9+8+(7+6\times 54)\times (32+1).$
\item [] $10941=(98+76\times 5+43)\times 21.$
\item [] $10942=987\times (6+5)+4^3+21.$
\item [] $10943=987\times (6+5)+43\times 2\times 1.$
\item [] $10944=987\times (6+5)+43\times 2+1.$
\item [] $10945=(9+87)\times 6\times (5+4\times 3+2)+1.$
\item [] $10946=(9\times (87+6)+5)\times (4+3^2\times 1).$
\item [] $10947=(9\times (87+6)+5)\times (4+3^2)+1.$
\item [] $10948=(9+8+765)\times (4\times 3+2)\times 1.$
\item [] $10949=(9+8)\times 7\times (6+54+32)+1.$
\item [] $10950=9+(8\times 7\times 65+4+3)\times (2+1).$
\item [] $10951=(9\times (87+65)\times 4+3)\times 2+1.$
\item [] $\mathit{10952=-9-87\times (6-5\times 4)\times 3^2-1.}$
\item [] $10953=9\times (876+5\times 4+321).$
\item [] $10954=9+8\times 76\times (5+4+3^2)+1.$
\item [] $10955=(9+8+7\times 65\times 4\times 3)\times 2+1.$
\item [] $10956=(9\times 8\times 7\times 6+5^4+3)\times (2+1).$
\item [] $\mathit{10957=(9+8\times 7\times 65+4)\times 3-2\times 1.}$
\item [] $\mathit{10958=(9+8\times 7\times 65+4)\times 3-2+1.}$
\item [] $10959=9+(8\times 76\times (5+4)+3)\times 2\times 1.$
\item [] $10960=9+(8\times 76\times (5+4)+3)\times 2+1.$
\item [] $10961=(9+8\times 7\times 65+4)\times 3+2\times 1.$
\item [] $10962=9876+543\times 2\times 1.$
\item [] $10963=9876+543\times 2+1.$
\item [] $\mathit{10964=(-9-87+6\times 5^4)\times 3+2\times 1.}$
\item [] $10965=9+(8\times 7\times 65+4\times 3)\times (2+1).$
\item [] $\mathit{10966=(1\times 2+3)\times ((4+5-6)^7+8)-9.}$
\item [] $\mathit{10967=(9\times 8\times 76+5+4+3)\times 2-1.}$
\item [] $10968=(9\times 8\times 76+5+4+3)\times 2\times 1.$
\item [] $10969=(9\times 8\times 76+5+4+3)\times 2+1.$
\item [] $\mathit{10970=9+87\times (6+5\times 4\times 3\times 2)-1.}$
\item [] $10971=9+87\times (6\times 5\times 4+3\times 2\times 1).$
\item [] $10972=9+87\times (6+54+3)\times 2+1.$
\item [] $\mathit{10973=-98+7-6\times (5-43^2\times 1).}$
\item [] $10974=9+(8+7\times (6+5))\times 43\times (2+1).$
\item [] $\mathit{10975=(9+8+76)\times (-5+4^3)\times 2+1.}$
\item [] $10976=98\times (7+65+4\times (3^2+1)).$
\item [] $10977=9\times 8\times (7+6\times 5)\times 4+321.$
\item [] $10978=(9\times 8\times 76+5+4\times 3)\times 2\times 1.$
\item [] $10979=9+8\times 7+(6\times 5+4)\times 321.$
\item [] $10980=(9+8\times 7\times 65+4)\times 3+21.$
\item[]$\mbox{Increasing order}$
\item [] $10981=1+23\times (4+5+6\times 78)+9.$
\item [] $10982=1+(2+3^4+56)\times (7+8\times 9).$
\item [] $10983=(1+2\times 3)\times (4^5+67\times 8+9).$
\item [] $10984=1\times 2^3\times (4\times (5+6\times 7\times 8)+9).$
\item [] $10985=(1+2^3+4)\times (56+789).$
\item [] $10986=1+2^3\times 4\times (5\times 67+8)+9.$
\item [] $\mathit{10987=-1\times 2+(3\times 45+6)\times 78-9.}$
\item [] $\mathit{10988=(1+2)^{(3+4)}\times 5+6+7\times 8-9.}$
\item [] $10989=(1+234\times 5+6\times 7+8)\times 9.$
\item [] $10990=1+(23+4)\times (5\times 67+8\times 9).$
\item [] $10991=(1+2)\times 3456+7\times 89.$
\item [] $10992=1\times 23\times 456+7\times 8\times 9.$
\item [] $10993=1+23\times 456+7\times 8\times 9.$
\item [] $10994=1234\times 5+67\times 8\times 9.$
\item [] $\mathit{10995=-12+(3\times 45+6)\times 78+9.}$
\item [] $\mathit{10996=-1+234\times (5+6\times 7)+8-9.}$
\item [] $10997=(1^2+3^4)\times (56+78)+9.$
\item [] $10998=(123+4^5+67+8)\times 9.$
\item [] $10999=1+2\times (34+5)\times (6+(7+8)\times 9).$
\item [] $11000=1\times (2+3)\times 4\times (5+67\times 8+9).$
\item [] $11001=1+(2+3)\times 4\times (5+67\times 8+9).$
\item [] $\mathit{11002=(1+2)^{(3+4)}\times 5+67\times (-8+9).}$
\item [] $\mathit{11003=(1+2)^{(3+4)}\times 5+67-8+9.}$
\item [] $11004=1\times 2\times 3\times (4^5+6\times (7+8)\times 9).$
\item [] $11005=1+2\times 3\times (4^5+6\times (7+8)\times 9).$
\item [] $11006=(1+2)^(3+4)\times 5+6+7\times 8+9.$
\item [] $11007=1^2\times (3\times 45+6)\times 78+9.$
\item [] $11008=1^2+(3\times 45+6)\times 78+9.$
\item [] $11009=1\times 2+(3\times 45+6)\times 78+9.$
\item [] $11010=1+2+(3\times 45+6)\times 78+9.$
\item [] $11011=1\times 2+(3\times 456+7)\times 8+9.$
\item [] $11012=1+2+(3\times 456+7)\times 8+9.$
\item [] $\mathit{11013=(1+23)\times 456+78-9.}$
\item [] $\mathit{11014=-1+234\times (5+6\times 7)+8+9.}$
\item [] $11015=1\times 234\times (5+6\times 7)+8+9.$
\item [] $11016=1+234\times (5+6\times 7)+8+9.$
\item [] $11017=(1^2+3\times 456+7)\times 8+9.$
\item [] $11018=1\times 2+3^4\times (5+6\times 7+89).$
\item [] $11019=12+(3\times 45+6)\times 78+9.$
\item [] $11020=(1+2)^(3+4)\times 5+6+7+8\times 9.$
\item [] $11021=12+(3\times 456+7)\times 8+9.$
\item [] $11022=(1+2)\times (34+56\times (7\times 8+9)).$
\item [] $11023=(1+23)\times 456+7+8\times 9.$
\item [] $\mathit{11024=(1+23\times (4+5))\times (6+7\times 8-9).}$
\item [] $11025=(1\times 2+3\times 456+7)\times 8+9.$
\item [] $11026=1+(2+3\times 456+7)\times 8+9.$
\item [] $\mathit{11027=-1\times 2+(3\times 4\times 5^6-7)/(8+9).}$
\item [] $11028=12+3^4\times (5+6\times 7+89).$
\item [] $\mathit{11029=(1-2\times 3\times 45)\times (6-7\times 8+9).}$
\item [] $\mathit{11030=(1+2)^{(3+4)}\times 5+(6+7)\times 8-9.}$
\item [] $11031=(1+23)\times 456+78+9.$
\item [] $\mathit{11032=1\times 2\times ((3+4)\times (5-6+789)).}$
\item [] $11033=(1+2+3\times 456+7)\times 8+9.$
\item [] $11034=(1\times 2\times 345+67\times 8)\times 9.$
\item [] $11035=1+(2\times 345+67\times 8)\times 9.$
\item [] $11036=(1+2\times 3+(4+5)\times (6+7))\times 89.$
\item [] $11037=(1+2)^(3+4)\times 5+6+7+89.$
\item [] $\mathit{11038=(1-2+34)\times 5\times 67-8-9.}$
\item [] $11039=(12+3+4)\times (5+6\times (7+89)).$
\item [] $11040=(1+23)\times 456+7+89.$
\item [] $11041=1+23\times (456+7+8+9).$
\item [] $11042=1\times 2\times (3+(4\times 5+6\times 7)\times 89).$
\item [] $11043=(1+2\times 345+67\times 8)\times 9.$
\item [] $11044=1+(23+4)\times (56\times 7+8+9).$
\item [] $11045=(12^3+4+5)\times 6+7\times 89.$
\item [] $\mathit{11046=(1+2\times 3-4+5+6)\times 789.}$
\item [] $\mathit{11047=-1+2^3\times (4\times (5\times 67+8)+9).}$
\item [] $11048=(1+23+4)\times 56\times 7+8\times 9.$
\item [] $11049=12\times (3+45+67)\times 8+9.$
\item [] $11050=(1\times 2+3^4+567)\times (8+9).$
\item[]$\mbox{Decreasing order}$
\item [] $10981=(9\times 8+7)\times (6+5+4^3\times 2\times 1).$
\item [] $10982=(9\times 8+7)\times (6+5+4\times 32)+1.$
\item [] $10983=(9+8\times 7\times 65+4\times 3)\times (2+1).$
\item [] $\mathit{10984=987\times (6+5)+4\times 32-1.}$
\item [] $10985=987\times (6+5)+4\times 32\times 1.$
\item [] $10986=987\times (6+5)+4\times 32+1.$
\item [] $10987=(9+(8+7\times 65\times 4)\times 3)\times 2+1.$
\item [] $\mathit{10988=-9\times (8+7)+6\times (5+43^2)-1.}$
\item [] $10989=987\times (6+5)+4\times (32+1).$
\item [] $10990=(9\times 8\times 76+5\times 4+3)\times 2\times 1.$
\item [] $10991=(9\times 8\times 76+5\times 4+3)\times 2+1.$
\item [] $10992=9\times 8+7\times 65\times 4\times 3\times 2\times 1.$
\item [] $10993=9\times 8+7\times 65\times 4\times 3\times 2+1.$
\item [] $\mathit{10994=9+8+7^(-6+5+4)\times 32+1.}$
\item [] $10995=9\times 8+(7+6\times 54)\times (32+1).$
\item [] $\mathit{10996=9-(87-6\times 5^4)\times 3-2\times 1.}$
\item [] $10997=98+(7\times 65+4^3)\times 21.$
\item [] $10998=(9\times 8\times 76+(5+4)\times 3)\times 2\times 1.$
\item [] $10999=(9\times 8\times 76+(5+4)\times 3)\times 2+1.$
\item [] $11000=(9+8)\times (7+6+5^4+3^2)+1.$
\item [] $11001=987\times (6+5)+(4\times 3)^2\times 1.$
\item [] $11002=987\times (6+5)+(4\times 3)^2+1.$
\item [] $\mathit{11003=9\times 8\times 7\times 6+(5\times 4)^3-21.}$
\item [] $11004=987\times (6+5)+(4+3)\times 21.$
\item [] $11005=(9+8\times 7+6)\times 5\times (4+3^{(2+1)}).$
\item [] $\mathit{11006=-(9+8)\times 7+6\times (5+43^2)+1.}$
\item [] $11007=9\times (((87+65)\times 4+3)\times 2+1).$
\item [] $\mathit{11008=9\times 876+5^4\times (3+2)-1.}$
\item [] $11009=9\times 876+5^4\times (3+2\times 1).$
\item [] $11010=9+87+(6\times 5+4)\times 321.$
\item [] $11011=(9+8\times (76\times (5+4)+3))\times 2+1.$
\item [] $\mathit{11012=(98-7)\times (6\times 5\times 4+3-2)+1.}$
\item [] $\mathit{11013=-(9\times 8-7)\times 6+543\times 21.}$
\item [] $11014=(9\times 8\times 76+5\times (4+3))\times 2\times 1.$
\item [] $11015=(9\times 8\times 76+5\times (4+3))\times 2+1.$
\item [] $11016=9\times 8+76\times (5+4+3)^2\times 1.$
\item [] $11017=9\times 8+76\times (5+4+3)^2+1.$
\item [] $11018=98+7\times 65\times 4\times 3\times 2\times 1.$
\item [] $11019=98+7\times 65\times 4\times 3\times 2+1.$
\item [] $11020=(9+8\times 7+6+5)\times ((4\times 3)^2+1).$
\item [] $11021=98+(7+6\times 54)\times (32+1).$
\item [] $11022=(9+8+7\times 65\times 4)\times 3\times 2\times 1.$
\item [] $11023=((9+87+6)\times 54+3)\times 2+1.$
\item [] $11024=(9\times 87+65)\times (4+3^2\times 1).$
\item [] $11025=(9+8+76+5+4+3)^2\times 1.$
\item [] $11026=(9+8+76+5+4+3)^2+1.$
\item [] $11027=9+8\times (7\times 65+4)\times 3+2\times 1.$
\item [] $11028=9+8\times (7\times 65+4)\times 3+2+1.$
\item [] $\mathit{11029=-9-87+6\times (5+43^2)+1.}$
\item [] $\mathit{11030=-9+(8\times 7\times 6+5+4)\times 32-1.}$
\item [] $\mathit{11031=-9+8\times ((7\times 65+4)\times 3+2+1).}$
\item [] $\mathit{11032=(98-7\times 6)\times (5+4^3\times (2+1)).}$
\item [] $11033=(9+8)\times 7+(6\times 5+4)\times 321.$
\item [] $11034=9+(87+6+5+4+3)^2\times 1.$
\item [] $11035=9+(87+6+5+4+3)^2+1.$
\item [] $11036=9+(8+7)\times (6+(5+4)^3)+2\times 1.$
\item [] $11037=9\times 8\times (7\times 6+5+4)\times 3+21.$
\item [] $\mathit{11038=(9\times 8\times 7+6^5)\times 4/3-2\times 1.}$
\item [] $\mathit{11039=-98+7+6\times (5+43^2+1).}$
\item [] $11040=(9\times 8\times 76+5+43)\times 2\times 1.$
\item [] $11041=(9\times 8\times 76+5+43)\times 2+1.$
\item [] $11042=9+8+7\times (6+5+4^3)\times 21.$
\item [] $11043=98+76\times (5+4+3)^2+1.$
\item [] $\mathit{11044=(9+8\times (7\times 65+4))\times 3+2-1.}$
\item [] $11045=9\times 8\times 7\times 6+(5\times 4)^3+21.$
\item [] $11046=9+8\times (7\times 65+4)\times 3+21.$
\item [] $11047=(9\times 87+6)\times (5+4+3+2)+1.$
\item [] $11048=98\times 76+(5\times 4\times 3)^2\times 1.$
\item [] $11049=9\times (8+7)+(6\times 5+4)\times 321.$
\item [] $11050=(9+8)\times (7+6)\times (5+43+2\times 1).$
\item[]$\mbox{Increasing order}$
\item [] $11051=((1^2+3)^4+5)\times 6\times 7+89.$
\item [] $11052=12\times (345+6\times (7+89)).$
\item [] $11053=1+2\times (3+(4+5)\times 67+8)\times 9.$
\item [] $\mathit{11054=(1-2+34)\times 5\times 67+8-9.}$
\item [] $\mathit{11055=(1-2+34)\times 5\times 67\times (-8+9).}$
\item [] $\mathit{11056=(1-2+34)\times 5\times 67-8+9.}$
\item [] $\mathit{11057=((1\times 234-5)\times 6+7)\times 8+9.}$
\item [] $11058=(1\times 23\times 4+5)\times (6\times 7+8\times 9).$
\item [] $11059=1+(23\times 4+5)\times (6\times 7+8\times 9).$
\item [] $11060=(1+2+3^4+56)\times (7+8\times 9).$
\item [] $11061=1+(2+(3+4\times 5)\times 6)\times (7+8\times 9).$
\item [] $11062=(1+234)\times (5+6\times 7)+8+9.$
\item [] $11063=(1+23)\times 456+7\times (8+9).$
\item [] $11064=(1\times 2+3+456)\times (7+8+9).$
\item [] $11065=(1+23+4)\times 56\times 7+89.$
\item [] $11066=(1+2)^(3+4)\times 5+6\times 7+89.$
\item [] $11067=123+456\times (7+8+9).$
\item [] $\mathit{11068=1^2+3\times (4\times 56-7)\times (8+9).}$
\item [] $\mathit{11069=-1+234\times (5+6\times 7)+8\times 9.}$
\item [] $11070=(1+2)\times 3456+78\times 9.$
\item [] $11071=1+234\times (5+6\times 7)+8\times 9.$
\item [] $11072=(1+2^3\times 4)\times 5\times 67+8+9.$
\item [] $11073=12^3+(4+5+6)\times 7\times 89.$
\item [] $11074=(1+2)^(3+4)\times 5+67+8\times 9.$
\item [] $\mathit{11075=1+2\times (-3-4+(5+6)\times 7\times 8\times 9).}$
\item [] $11076=(1+2)^(3+4)\times 5+6+(7+8)\times 9.$
\item [] $\mathit{11077=-1+2\times ((3+4^5)\times 6-7\times 89).}$
\item [] $\mathit{11078=-1+2\times 3\times (4\times 56+7)\times 8-9.}$
\item [] $11079=123\times (4+5\times 6+7\times 8)+9.$
\item [] $11080=(1+2\times 3)^4+(5+6)\times 789.$
\item [] $11081=((1+2)^3\times (45+6)+7)\times 8+9.$
\item [] $\mathit{11082=-1+(2+3)\times 4^5+67\times 89.}$
\item [] $11083=1\times (2+3)\times 4^5+67\times 89.$
\item [] $11084=1+(2+3)\times 4^5+67\times 89.$
\item [] $11085=(1^2+3\times 45+6)\times 78+9.$
\item [] $\mathit{11086=-1+234\times (5+6\times 7)+89.}$
\item [] $11087=1\times 234\times (5+6\times 7)+89.$
\item [] $11088=1+234\times (5+6\times 7)+89.$
\item [] $11089=1+(2\times 3+456)\times (7+8+9).$
\item [] $11090=1+(2+3)\times (45\times 6+7)\times 8+9.$
\item [] $11091=(1+2)^(3+4)\times 5+67+89.$
\item [] $\mathit{11092=(1+234-5+6)\times (7\times 8-9).}$
\item [] $\mathit{11093=-1^2+(3^4+5)\times (-6+(7+8)\times 9).}$
\item [] $\mathit{11094=(1-23\times 4+5)\times (6-(7+8)\times 9).}$
\item [] $11095=12^3\times 4+(5+6\times 7)\times 89.$
\item [] $11096=(12+3+4)\times (567+8+9).$
\item [] $11097=12^3+4\times 5\times 6\times 78+9.$
\item [] $11098=1+2\times 3\times (4\times 56+7)\times 8+9.$
\item [] $11099=1\times 2+3^4\times (5\times (6+7)+(8\times 9)).$
\item [] $11100=1+2+3^4\times (5\times (6+7)+(8\times 9)).$
\item [] $11101=(1\times 2^{(3+4)}\times 5+6+7)\times (8+9).$
\item [] $11102=1\times 2\times (3+4+(5+6)\times 7\times 8\times 9).$
\item [] $11103=(1+2)\times (3\times 4^5+6+7\times 89).$
\item [] $\mathit{11104=-1+2+3\times (456+7)\times 8-9.}$
\item [] $11105=(12+3\times 456+7)\times 8+9.$
\item [] $11106=(1+2^3)^4+567\times 8+9.$
\item [] $11107=1+2\times (3\times (4\times 56+7)\times 8+9).$
\item [] $\mathit{11108=-1+23\times (4+5+6\times (7+8\times 9)).}$
\item [] $11109=(1+2\times 34)\times (5+67+89).$
\item [] $11110=1+23\times (4+5+6\times (7+8\times 9)).$
\item [] $11111=1\times 23\times 456+7\times 89.$
\item[]$\mbox{Decreasing order}$
\item [] $11051=(9+8\times 7)\times (6+54\times 3+2)+1.$
\item [] $11052=((9+8)\times (7+65)+4)\times 3^2\times 1.$
\item [] $11053=((9+8)\times (7+65)+4)\times 3^2+1.$
\item [] $11054=98\times 7+6\times 54\times 32\times 1.$
\item [] $11055=98\times 7+6\times 54\times 32+1.$
\item [] $\mathit{11056=-9+8-7-6\times (5-43^2\times 1).}$
\item [] $11057=9+8\times (7+6\times 5+4^3\times 21).$
\item [] $11058=(9\times 8\times 76+54+3)\times 2\times 1.$
\item [] $11059=(9\times 8\times 76+54+3)\times 2+1.$
\item [] $11060=98\times 7+6\times (54\times 32+1).$
\item [] $11061=(98+7\times 65)\times 4\times (3+2)+1.$
\item [] $\mathit{11062=(9\times 8\times 76-5+4^3)\times 2\times 1.}$
\item [] $\mathit{11063=(9\times 8+7\times 65\times 4\times 3)\times 2-1.}$
\item [] $11064=(9\times 8\times 76+5\times 4\times 3)\times 2\times 1.$
\item [] $11065=(9\times 8\times 76+5\times 4\times 3)\times 2+1.$
\item [] $\mathit{11066=9\times (8\times 76+5^4)-32+1.}$
\item [] $11067=(9\times 8\times 7+6+5+4\times 3)\times 21.$
\item [] $\mathit{11068=9+8-7-6\times (5-43^2+1).}$
\item [] $\mathit{11069=(-9+8\times 7-6)\times 54\times (3+2)-1.}$
\item [] $11070=(987+6\times 5\times 4)\times (3^2+1).$
\item [] $11071=9\times (8+7)\times (6+5\times (4+3))\times 2+1.$
\item [] $11072=(98\times 7+6)\times (5+4+3\times 2+1).$
\item [] $11073=9876+(54+3)\times 21.$
\item [] $11074=98\times (7\times 6+5+4^3+2\times 1).$
\item [] $11075=((9+8)\times (7+6)+5)\times (4+3)^2+1.$
\item [] $11076=9+(8+7\times 65+4^3)\times 21.$
\item [] $\mathit{11077=9-8\times 7+6\times (5+43^2\times 1).}$
\item [] $\mathit{11078=9-8\times 7+6\times (5+43^2)+1.}$
\item [] $11079=9+(8+7)\times (6+(5+4)^3+2+1).$
\item [] $\mathit{11080=9\times (8-7)\times (6+(5\times (4+3))^2)+1.}$
\item [] $\mathit{11081=(9\times 8\times 76+5+4^3)\times 2-1.}$
\item [] $11082=(9\times 8\times 76+5+4^3)\times 2\times 1.$
\item [] $11083=(9\times 8\times 76+5+4^3)\times 2+1.$
\item [] $11084=(9+8)\times (7+6\times 54+321).$
\item [] $\mathit{11085=(-9+8\times 7\times 65+4^3)\times (2+1).}$
\item [] $\mathit{11086=-9+8+7\times (6+5)\times (4\times 3)^2-1.}$
\item [] $\mathit{11087=-9+8+7\times (6+5)\times (4\times 3)^2\times 1.}$
\item [] $11088=9\times 8\times (76+54+3+21).$
\item [] $11089=(98+7\times 6\times 5)\times 4\times 3^2+1.$
\item [] $11090=98\times 7+((6\times 5+4)\times 3)^2\times 1.$
\item [] $11091=98\times 7+(6\times (5+4\times 3))^2+1.$
\item [] $11092=9\times 8+76\times ((5+4+3)^2+1).$
\item [] $\mathit{11093=(-9-8-7+6\times 5)\times 43^2-1.}$
\item [] $11094=(9+(87\times 6+5)\times (4+3))\times (2+1).$
\item [] $\mathit{11095=(-9-8-7+6\times 5)\times 43^2+1.}$
\item [] $\mathit{11096=9\times 8\times 76+5^4\times 3^2-1.}$
\item [] $11097=9\times 8\times 76+5^4\times 3^2\times 1.$
\item [] $11098=9\times 8\times 76+5^4\times 3^2+1.$
\item [] $11099=9\times (87+6\times 54)\times 3+2\times 1.$
\item [] $11100=9\times (87+6\times 54)\times 3+2+1.$
\item [] $11101=(9+8)\times (7+6+5\times 4\times 32\times 1).$
\item [] $11102=9\times (8\times 76+5^4)+3+2\times 1.$
\item [] $11103=9\times (8\times 76+5^4)+3+2+1.$
\item [] $11104=9\times (8\times 76+5^4)+3\times 2+1.$
\item [] $11105=9+8+7\times (6+5)\times (4\times 3)^2\times 1.$
\item [] $11106=9+(8\times 76+5^4)\times 3^2\times 1.$
\item [] $11107=9+(8\times 76+5^4)\times 3^2+1.$
\item [] $\mathit{11108=-9-8+7+6\times (5+43^2-1).}$
\item [] $11109=987\times (6+5)+4\times 3\times 21.$
\item [] $\mathit{11110=-987+6^5+4321.}$
\item [] $\mathit{11111=-9+8\times 7-6\times (5-43^2\times 1).}$
\end{itemize}
\end{multicols}

\section{Numbers with Subtraction and Division}

There are 1256 numbers where we have used subtraction and/or division operations. Out of them 611 are in the increasing case and 645 in decreasing case. The number 10958 is still missing. There are 178 numbers common in both the cases. Division sign is applied to find 8 numbers. These numbers are 9986, 10084, 10121 and 11027 in the increasing case and 9668, 9686, 10802 and 11038 in the decreasing case.

\section*{\textbf{Acknowledgement}}

The author is thankful to T.J. Eckman, Georgia, United States, email: jeek@jeek.net in finding some difficult representations using computerized script.

\begin{center}
-------------------------------------------------
\end{center}
\end{document}